\documentclass[oneside,a4paper,english, 11pt]{amsart}

\usepackage{adjustbox}
\usepackage[latin9]{inputenc}
\usepackage{hyperref}
\usepackage{lscape}
\usepackage{calligra,mathrsfs}  
\usepackage{enumerate}
\usepackage{amsmath}
\usepackage{amssymb}
\usepackage{amscd}
\usepackage{amsthm}

\makeatletter
\theoremstyle{plain}
\newtheorem{prop}{Proposition}[subsection]
\newtheorem{conj}[prop]{Conjecture}
\newtheorem{lem}[prop]{Lemma}
\newtheorem{thm}[prop]{Theorem}
\newtheorem{cor}[prop]{Corollary}
\newtheorem{cond}[prop]{Condition}
\newtheorem{hypo}[prop]{Hypothesis}
\newtheorem{cons}[prop]{Construction}
\theoremstyle{definition}
\newtheorem{defn}[prop]{Definition}
\newtheorem{ex}[prop]{Example}
\newtheorem{rem}[prop]{Remark}
\newtheorem{ques}[prop]{Question}

\theoremstyle{plain}

\theoremstyle{definition}

\theoremstyle{plain}

\theoremstyle{definition}

\usepackage{array,multirow}
\usepackage[arrow,matrix,tips,all,cmtip,poly,color]{xy}

\usepackage{stmaryrd}
\usepackage{url}
\usepackage{tikz-cd}
\usepackage{color}
\usepackage{caption}
\usepackage{float}
\usepackage{enumitem}

\pdfpageheight\paperheight
\pdfpagewidth\paperwidth
\newcommand{\Z}{\mathbb{Z}}

\newcommand{\Q}{\mathbb{Q}}

\newcommand{\R}{\mathbb{R}}

\newcommand{\bG}{\mathbb{G}}

\newcommand{\fa}{\mathfrak{a}}
\newcommand{\fb}{\mathfrak{b}}

\newcommand{\fg}{\mathfrak{g}}
\newcommand{\fh}{\mathfrak{h}}

\newcommand{\fl}{\mathfrak{l}}

\newcommand{\fp}{\mathfrak{p}}

\newcommand{\ft}{\mathfrak{t}}
\newcommand{\fu}{\mathfrak{u}}

\newcommand{\fm}{\mathfrak{m}}
\newcommand{\fn}{\mathfrak{n}}

\newcommand{\fs}{\mathfrak{s}}
\newcommand{\fz}{\mathfrak{z}}

\newcommand{\bA}{\mathbb{A}}

\newcommand{\bH}{\mathbb{H}}

\newcommand{\bP}{\mathbb{P}}

\newcommand{\bT}{\mathbb{T}}

\newcommand{\cA}{\mathcal{A}}
\newcommand{\cB}{\mathcal{B}}
\newcommand{\cC}{\mathcal{C}}
\newcommand{\cD}{\mathcal{D}}
\newcommand{\cE}{\mathcal{E}}
\newcommand{\cF}{\mathcal{F}}
\newcommand{\cG}{\mathcal{G}}
\newcommand{\cH}{\mathcal{H}}
\newcommand{\cI}{\mathcal{I}}

\newcommand{\cJ}{\mathcal{J}}

\newcommand{\cM}{\mathcal{M}}
\newcommand{\cN}{\mathcal{N}}
\newcommand{\cO}{\mathcal{O}}
\newcommand{\cP}{\mathcal{P}}

\newcommand{\cR}{\mathcal{R}}
\newcommand{\cS}{\mathcal{S}}
\newcommand{\cT}{\mathcal{T}}
\newcommand{\cU}{\mathcal{U}}
\newcommand{\cV}{\mathcal{V}}
\newcommand{\cW}{\mathcal{W}}
\newcommand{\cX}{\mathcal{X}}

\newcommand{\cHom}{\mathscr{H}\text{\kern -3pt{\calligra\large om}\,}}

\newcommand{\al}{\alpha}

\newcommand{\val}{\mathrm{val}}
\newcommand{\plog}{\mathrm{log}}

\newcommand{\Hom}{\mathrm{Hom}}

\newcommand{\GL}{\mathrm{GL}}

\newcommand{\un}[1]{\underline{#1}}

\newcommand{\tld}[1]{\widetilde{#1}}

\newcommand{\defeq}{\stackrel{\textrm{\tiny{def}}}{=}}

\newif\iffinalrun
\iffinalrun
  \newcommand{\mar}[1]{}
\else
  \newcommand{\mar}[1]{\marginpar{\raggedright\tiny #1}}

\DeclareMathOperator{\Dim}{dim}

\DeclareMathOperator{\dom}{dom}

\title{On $\mathrm{Ext}^{\bullet}$ between locally analytic generalized Steinberg with applications}
\author{Zicheng Qian}

\begin{document}

\maketitle

\begin{abstract}
Let $n\geq 2$ be an integer, $p$ be a prime number and $K$ be a finite extension of $\Q_p$. 
Motivated by Schraen's thesis and Gehrmann's definition of automorphic simple $\mathscr{L}$-invariants, we study the first non-vanishing extension groups between a pair of locally $K$-analytic generalized Steinberg representations of $\mathrm{GL}_n(K)$.
We study subspaces of these extension groups defined by using either relative conditions with respect to Lie subalgebras of $\fs\fl_{n}$ (isomorphic to $\fs\fl_{m}$ for some $2\leq m<n$) or maps between locally $K$-analytic generalized Steinberg representations of $\mathrm{GL}_n(K)$ with different highest weights.
The applications of these computations are two-fold. On one hand, we prove that a certain universal successive extension of filtered $(\varphi,N)$-modules can be realized as the space of homomorphisms from a suitable shift of the dual of locally $K$-analytic Steinberg representation into the de Rham complex of the Drinfeld upper-half space, generalizing one main result of Schraen's thesis from $\mathrm{GL}_{3}(\Q_p)$ to $\mathrm{GL}_{n}(K)$. On the other hand, we give a definition of higher $\mathscr{L}$-invariants for $\mathrm{GL}_n(K)$ (which we call Breuil-Schraen $\mathscr{L}$-invariants) and discuss its possible explicit relation to Fontaine-Mazur $\mathscr{L}$-invariants, using ideas from Breuil-Ding's higher $\mathscr{L}$-invariants for $\mathrm{GL}_{3}(\Q_p)$.
\end{abstract}
\tableofcontents

\section{Introduction}\label{sec: intro}
\subsection{Historical background}\label{subsec: intro history}
Let $f$ be a normalized new form of weight $2$ and level $\Gamma_0(pN)$ with $p\nmid N$.
Assume moreover that $f$ is defined over $\Q$ with $U_{p}$-eigenvalue $1$. In this case, one could attach an elliptic curve $\mathrm{E}_{/\Q}$ to $f$ which has split multiplicative reduction at $p$ and thus there exists $q\in\Q_p^{\times}$ such that $|q|<1$ and one has an isomorphism 
\[\bG_{m}^{\rm{rig}}/q^{\Z}=\mathrm{E}^{\rm{rig}}\] 
of rigid analytic groups over $\Q_p$ by Tate's $p$-adic uniformization theorem. One has $L_{p}(f,1)=0$ and the weight $2$ case of the famous Mazur-Tate-Teitelbaum conjecture (see \cite{MTT86}) claims that
\[\frac{d}{ds}~L_{p}(f,s)|_{s=1}=\mathscr{L}(f) L(f,1)\]
where the $\mathscr{L}$-invariant $\mathscr{L}(f)$ of $f$ is defined as
\[\mathscr{L}(f)=\mathscr{L}(\mathrm{E})\defeq \frac{\plog_{p}(q)}{\val_{p}(q)}.\]
This conjecture have been fully proven (see \cite{GS93}) and many alternative definitions of $\mathscr{L}(f)$ (even for weight $k\geq 2$) has been introduced and shown to be all equivalent (cf.~\cite{Ast}, \cite{Col04} and \cite{Col05} for detailed discussions).
Following the list in the introduction of \cite{Geh21}, we have the following incomplete list of definitions of $\mathscr{L}$-invariants along with their generalizations to higher dimension/rank cases.
\begin{itemize}
\item Fontaine and Mazur (see \cite{M91}) define $\mathscr{L}$-invariants for general semi-stable $p$-adic continuous Galois representations via their associated filtered $(\varphi,N)$-modules.
\item Teitelbaum (see \cite{T90}) defines $\mathscr{L}(f)$ using Jacquet-Langlands transfer, Coleman's $p$-adic integration theory and the Cerednik-Drinfeld uniformization theorem. Later Besser and de Shalit (see \cite{BdS16}) generalize this definition to some higher dimensional varieties which are $p$-adically uniformizable by Drinfeld's upper half space. 
\item Darmon (see \cite{D01}) define $\mathscr{L}(f)$ using harmonic cochains on the Bruhat-Tits tree which are invariant under a $p$-arithmetic subgroup of $\mathrm{GL}_2(\Q)$. This definition has been generalized by Orton (see \cite{O04}) to modular form of higher weights, by Spiess (see \cite{Sp14}) to Hilbert modular forms, and recently by Gehrmann (see \cite{Geh21}) to suitable automorphic representations of higher rank reductive groups.
\item Breuil (see \cite{Bre10}) defines $\mathscr{L}(f)$ by understanding the $\mathrm{GL}_{2}(\Q_p)$-action on the $f$-isotypic component of the completed cohomology of a tower of modular curves, and thus initiates the study of $\mathscr{L}$-invariants via the locally analytic representation theory of $p$-adic reductive groups. 
    Breuil's definition of $\mathscr{L}$-invariants has been generalized to the case $\mathrm{GL}_{2}(K)$ (see \cite{Schr10} and \cite{Ding16}), the case $\mathrm{GL}_{3}(\Q_p)$ (see \cite{Schr11} and more recently \cite{Bre19}, \cite{BD20} and \cite{BD19}) and partially for $\mathrm{GL}_{n}(K)$ in \cite{Ding19}.
\end{itemize}
The main results of this paper are three-folds.
\begin{itemize}
\item Motivated by Schraen's work for $\mathrm{GL}_{3}(\Q_p)$ and Gehrmann's work on automorphic simple $\mathscr{L}$-invariants, we define what we call \emph{Breuil-Schraen $\mathscr{L}$-invariants} which generalize Breuil's definition of $\mathscr{L}$-invariants to $\mathrm{GL}_{n}(K)$ for arbitrary $n$ and $K$.
\item We conjecture an explicit relation between Breuil-Schraen $\mathscr{L}$-invariants and Fontaine-Mazur $\mathscr{L}$-invariants based on some ideas from Breuil-Ding's higher $\mathscr{L}$-invariants for $\mathrm{GL}_{3}(\Q_p)$.
\item We prove that a certain universal successive extension of filtered $(\varphi,N)$-modules can be realized by mapping (suitable shift of) the dual of locally $K$-analytic Steinberg into the de Rham complex of the Drinfeld upper-half space, which generalizes a results of Schraen from $\mathrm{GL}_{3}(\Q_p)$ to $\mathrm{GL}_{n}(K)$ and relates our Breuil-Schraen $\mathscr{L}$-invariants to the geometry of Drinfeld upper-half space. In particular, we describe a natural candidate for the explicit Hyodo-Kato splitting of the de Rham complex of the Drinfeld upper-half space.
\end{itemize}
\subsection{Ideas from some previous work}\label{subsec: intro L loc an}
We recall various known approach to $\mathscr{L}$-invariants that uses locally analytic representation theory of $p$-adic reductive groups. We will assume $K=\Q_p$ throughout this section except for the general notation in \S~\ref{subsubsec: intro quick notation}.
\subsubsection{Some notation}\label{subsubsec: intro quick notation}
Let $E$ be a sufficiently large $p$-adic field and we fix an embedding $\iota: K\hookrightarrow E$.
We fix below an integer $n\geq 2$ with $G\defeq\mathrm{PGL}_{n}(K)$ and write $B$ (resp.~$B^{+}$, resp.~$T$) for its lower-triangular Borel subgroup (resp.~upper-triangular Borel subgroup, resp.~diagonal maximal torus). 
We write $\Phi^+$ for the set of positive roots associated with $(B^{+},T)$ with $\Delta\subseteq\Phi^+$ being the set of simple positive roots which we often conveniently identify with $[1,n-1]\defeq \{1,\dots,n-1\}$. For each subset $I\subseteq\Delta$, we can associate a parabolic subgroup $P_{I}\supseteq B$ as well as its standard Levi subgroup $L_{I}$ equipped with a surjection $P_{I}\twoheadrightarrow L_{I}$. 
Let $J\subseteq\Delta$.
We write $D(L_{J})$ for the unital associative $E$-algebra of locally $K$-analytic distributions, and write $\cM(L_{J})$ for the bounded derived category of the abelian category $\mathrm{Mod}_{D(L_{J})}$ of abstract $D(L_{J})$-modules.
We write $\mathrm{Rep}^{\rm{an}}_{\rm{adm}}(L_{J})$ (resp.~$\mathrm{Rep}^{\infty}_{\rm{adm}}(L_{J})$, resp.~$\mathrm{Rep}^{\rm{cont}}_{\rm{adm}}(L_{J})$) for the abelian category of admissible locally $K$-analytic (resp.~smooth, resp.~unitary Banach) $E$-representations of $L_{J}$ (see \cite{ST03}). (Throughout this paper, by locally $K$-analytic we mean locally $\iota$-analytic with coefficients in $E$.) For each $V\in \mathrm{Rep}^{\rm{an}}_{\rm{adm}}(L_{J})$, note that $V^{\vee}$ is a \emph{coadmissible $D(L_{J})$-module} in the sense of Schneider-Teitelbaum (see \cite[\S~6]{ST03}).
We have a fully faithful embedding $\mathrm{Rep}^{\infty}_{\rm{adm}}(L_{J})\hookrightarrow \mathrm{Rep}^{\rm{an}}_{\rm{adm}}(L_{J})$ and functors $\mathrm{Rep}^{\rm{cont}}_{\rm{adm}}(L_{J})\rightarrow \mathrm{Rep}^{\rm{an}}_{\rm{adm}}(L_{J})\rightarrow \mathrm{Rep}^{\infty}_{\rm{adm}}(L_{J})$ given by taking locally $K$-analytic vectors and smooth vectors.
For each $I\subseteq J$, we consider the following \emph{locally $K$-analytic principal series}
\[i_{I,J}^{\rm{an}}\defeq (\mathrm{Ind}_{P_{I}\cap L_{J}}^{L_{J}}1_{L_{I}})^{\rm{an}}\]
which comes with a natural embedding $i_{I',J}^{\rm{an}}\hookrightarrow i_{I,J}^{\rm{an}}$ for each $I\subseteq I'\subseteq J$. Finally, we consider the \emph{locally $K$-analytic generalized Steinberg representation}
\[V_{I,J}^{\rm{an}}\defeq i_{I,J}^{\rm{an}}/\sum_{I'\supsetneq I}i_{I',J}^{\rm{an}}\]
which will be a key player in (the generalization of) Breuil's definition of $\mathscr{L}$-invariants.
When $J=\Delta$, we omit it from the notation and obtain $i_{I}^{\rm{an}}$ and $V_{I}^{\rm{an}}$ for each $I\subseteq\Delta$, with $V_{\Delta}^{\rm{an}}=1_{G}$ and $\mathrm{St}_{G}^{\rm{an}}\defeq V_{\emptyset}^{\rm{an}}$ being the \emph{locally $K$-analytic Steinberg representation}.
Replacing $(\mathrm{Ind}_{P_{I}\cap L_{J}}^{L_{J}}-)^{\rm{an}}$ with $(\mathrm{Ind}_{P_{I}\cap L_{J}}^{L_{J}}-)^{\infty}$ and $(\mathrm{Ind}_{P_{I}\cap L_{J}}^{L_{J}}-)^{\rm{cont}}$ in the discussion above, we can similarly define a \emph{smooth Steinberg representation} $V_{I,J}^{\infty}$ and a \emph{continuous Steinberg representation} $V_{I,J}^{\rm{cont}}$ for each $I\subseteq J$, with $\mathrm{St}_{G}^{\infty}\defeq V_{\emptyset,\Delta}^{\infty}$ and $\mathrm{St}_{G}^{\rm{cont}}\defeq V_{\emptyset,\Delta}^{\rm{cont}}$ for short.
Note that both $i_{I,J}^{\ast}$ and $V_{I,J}^{\ast}$ are objects of $\mathrm{Rep}^{\ast}_{\rm{adm}}(L_{J})$ for each $\ast\in\{\rm{an},\infty,\rm{cont}\}$. We also know that $V_{I,J}^{\rm{an}}$ is the $E$-subspace of locally $K$-analytic vectors of $V_{I,J}^{\rm{cont}}$, and $V_{I,J}^{\infty}$ is the $E$-subspace of smooth vectors of $V_{I,J}^{\rm{an}}$, with similar facts hold for $i_{I,J}^{\ast}$.
\subsubsection{Breuil's $\mathscr{L}$-invariants for $\mathrm{PGL}_{2}(\Q_p)$}\label{subsubsec: intro Breuil L inv}
We assume that $n=2$ in this section.
We write $\mathrm{Gal}_{\Q_p}\defeq\mathrm{Gal}(\overline{\Q_p}/\Q_p)$ for the absolute Galois group of $\Q_p$ and similarly for $\mathrm{Gal}_{\Q}$.
Let $f$ be a normalized new form of weight $2$ as in \S \ref{subsec: intro history}.
In \cite{Bre04} and \cite{Bre10}, Breuil initiates the definition of $\mathscr{L}(f)$ via the theory of $p$-adic Banach (or locally analytic) representations of $G=\mathrm{GL}_{2}(\Q_p)$. We briefly recall his result as follows.
On one hand, writing $r_{f}:\mathrm{Gal}_{\Q}\rightarrow \mathrm{GL}_{2}(E)$ for the $p$-adic Galois representation associated with $f$ and $\rho_{f}\defeq r_{f}|_{\mathrm{Gal}_{\Q_p}}$, we know from our assumption on $f$ that $\rho_{f}$ has the form 
\[\left(
\begin{array}{cr}
\varepsilon & \ast\\
 0 & 1
\end{array}
\right)\] 
and thus the isomorphism class of $\rho_{f}$ is determined by an $E$-line in the $2$-dimensional space 
\[\cE_{1}\defeq \mathrm{Ext}_{\mathrm{Gal}_{\Q_p}}^1(1,\varepsilon)^{\vee}\] 
with $\varepsilon:\mathrm{Gal}_{\Q_p}\rightarrow \Z_{p}^{\times}$ being the $p$-adic cyclotomic character.
On the other hand, Breuil proves that there exists a unique (up to isomorphism) locally $\Q_p$-analytic representation $\Pi_{f}^{\rm{an}}$ of $G$ that fits into a non-split extension of the form
\[0\rightarrow \mathrm{St}_{G}^{\rm{an}}\rightarrow \Pi_{f}^{\rm{an}}\rightarrow 1_{G}\rightarrow 0\]
and moreover admits a $G$-equivariant embedding into the $f$-isotypic component of the completed cohomology of a tower of modular curves (with infinite level at $p$). Here the isomorphism class of $\Pi_{f}^{\rm{an}}$ is determined by an $E$-line in the $2$-dimensional space 
\[\mathbf{E}_{1}=\mathrm{Ext}_{G}^{1}(1_{G},\mathrm{St}_{G}^{\rm{an}}).\]
Finally, the isomorphism class of $\Pi_{f}^{\rm{an}}$ and that of $\rho_{f}$ are related via a sequence of canonical isomorphisms between $2$-dimensional $E$-vector spaces of the following form
\begin{multline}\label{equ: intro E1 regulator}
\mathbf{E}_{1}=\mathrm{Ext}_{G}^{1}(1_{G},\mathrm{St}_{G}^{\rm{an}})\buildrel\sim\over\longleftarrow \mathrm{Ext}_{G}^{1}(i_{\emptyset}^{\rm{an}},i_{\emptyset}^{\rm{an}}) \buildrel\sim\over\longleftarrow
\mathrm{Ext}_{T}^{1}(1_{T},1_{T}) \buildrel\sim\over\longleftarrow \Hom(\Q_p^{\times},E) \\
\buildrel\sim\over\longrightarrow \Hom(\mathrm{Gal}_{\Q_p},E)\buildrel\sim\over\longrightarrow H^{1}(\mathrm{Gal}_{\Q_p},1_{\mathrm{Gal}_{\Q_p}})\buildrel\sim\over\longrightarrow H^{1}(\mathrm{Gal}_{\Q_p},\varepsilon)^{\vee}=\cE_{1}.
\end{multline}
In other words, we can associate with each $E$-line $\mathscr{L}\subseteq \Hom(\Q_p^{\times},E)$ a unique (up to isomorphism) $G$-representation $\Pi_{\mathscr{L}}^{\rm{an}}$ that fits into a non-split extension
\[0\rightarrow \mathrm{St}_{G}^{\rm{an}}\rightarrow \Pi_{\mathscr{L}}^{\rm{an}}\rightarrow 1_{G}\rightarrow 0,\]
and a unique (up to isomorphism) $\rho_{\mathscr{L}}:\mathrm{Gal}_{\Q_p}\rightarrow \mathrm{GL}_{2}(E)$ that fits into
\[0\rightarrow \varepsilon\rightarrow \rho_{\mathscr{L}}\rightarrow 1_{\mathrm{Gal}_{\Q_p}}\rightarrow 0,\]
which recover $\Pi_{f}^{\rm{an}}$ and $\rho_{f}$ when we take $\mathscr{L}=E(\plog_{p}-\mathscr{L}(f)\val_{p})$.
Here $\val_{p}:\Q_p^{\times}\rightarrow \Z$ is the $p$-adic valuation with $\val_{p}(p)=1$, and $\plog_{p}:\Q_p^{\times}\rightarrow \Q_p$ is the branch of $p$-adic logarithm that satisfies $\plog_{p}(p)=0$.
The correspondence
\[
\xymatrix{
\Pi_{\mathscr{L}}^{\rm{an}} & \mathscr{L}\subseteq \Hom(\Q_p^{\times},E) \ar@{~>}[l]\ar@{~>}[r] & \rho_{\mathscr{L}} 
}
\]
is the first instance of the so-called $p$-adic local Langlands correspondence for $\mathrm{GL}_{2}(\Q_p)$.
\subsubsection{Questions from $p$-adic local-global compatibility conjectures}\label{subsubsec: intro lgc question}
To generalize Breuil's $\mathscr{L}$-invariants to $G=\mathrm{PGL}_{n}(\Q_p)$ for arbitrary $n$, we begin with the formulation of various natural questions/expectations on such generalization following the conjectural $p$-adic Langlands correspondence and its local-global compatibility.
Let $F$ be a number field with a place $v$ over $p$ such that $F_{v}=\Q_p$. Let $\cG$ be a reductive group over $F$ such that $\cG(F_{v})=G=\mathrm{PGL}_{n}(\Q_p)$. Given a compact open subgroup $U^{v}\subseteq \cG(\bA_{F}^{\infty,v})$, we can define certain universal Hecke algebra $\bT(U^{v})$ that acts on the following $E$-vector space of $p$-adic automorphic forms
\[\widehat{S}(U^{v})\defeq C^{\rm{cont}}(\cG(F)\backslash\cG(\bA_{F}^{\infty})/U^{v},\cO_{E})\otimes_{\cO_{E}}E.\]
Let $\cG^{\vee}$ be the Langlands dual group of $\cG$ and $r:\mathrm{Gal}_{F}\rightarrow \cG^{\vee}(E)$ be $p$-adic continuous Galois representation, to which we may associate a maximal ideal $\mathfrak{m}_{r}\subseteq \bT(U^{v})$ so that we may consider the following Hecke isotypic component
\[\Pi(U^{v},r)\defeq \widehat{S}(U^{v})[\mathfrak{m}_{r}]\]
whose $E$-subspace of locally $\Q_p$-analytic vectors is denoted by $\Pi(U^{v},r)^{\rm{an}}$.
The cases that are relevant to the generalization of Breuil's $\mathscr{L}$-invariants are those when $\mathrm{St}_{G}^{\infty}$ embeds into $\Pi(U^{v},r)^{\rm{an}}$ (upon twisting by some unitary character), for some favorable choices of $F$, $v$, $\cG$, $U^{v}$ and $r$ (cf.~\cite[\S 6]{BD20}). 
In this case, we would know (by compatibility with the classical local Langlands correspondence in favorable cases) that $\rho\defeq r|_{\mathrm{Gal}_{F_{v}}}$ is a semi-stable $p$-adic Galois representation of the form (again up to a determinant twist)
\begin{equation}\label{equ: intro n dim Galois}
\left(
\begin{array}{ccccrrrr}
\varepsilon^{n-1} & \ast & \cdots & \ast & \ast\\
0 & \varepsilon^{n-2} & \cdots & \ast & \ast\\
\vdots & \vdots & \ddots & \vdots & \vdots\\
0 & 0 & \cdots & \varepsilon & \ast\\
0 & 0 & \cdots & 0 & 1
\end{array}
\right)
\end{equation}
such that $D_{\rm{st}}(\rho)$ satisfies $N^{n-1}\neq 0$. Taking an arbitrary $0\neq e_0\in D_{\rm{st}}(\rho)^{\varphi=1}$ and writing $e_{k}\defeq N^{k}(e_0)$ with $\varphi(e_{k})=p^{-k}e_{k}$ for each $1\leq k\leq n-1$, we obtain a basis $\{e_{k}\}_{0\leq k\leq n-1}$ of $D_{\rm{st}}(\rho)$ and thus of $D_{\rm{dR}}(\rho)$. When the Hodge filtration on $D_{\rm{dR}}(\rho)$ lies in a sufficiently generic position (often called \emph{non-critical}), there exists a unique scalar $\mathscr{L}_{\al}=\mathscr{L}_{i,j}\in E$ for each $\al=(i,j)\in\Phi^+$ such that
\[\mathrm{Fil}^{-k}(D_{\rm{dR}}(\rho))=\mathrm{Fil}^{-k+1}(D_{\rm{dR}}(\rho))\oplus E(e_{k}+\sum_{k<j\leq n}\mathscr{L}_{k,j}e_{j})\]
for each $0\leq k\leq n-1$. The tuple of scalars
\begin{equation}\label{equ: intro FM L inv}
\un{\mathscr{L}}=\{\mathscr{L}_{\al}\}_{\al\in\Phi^+}
\end{equation}
is called the \emph{Fontaine-Mazur $\mathscr{L}$-invariants} associated with $\rho$, as defined by Mazur in \cite{M91} using Fontaine's theory.
Note that (\ref{equ: intro FM L inv}) is clearly independent of our choice of $e_0$, and in fact the isomorphism class of $\rho$ determines and depends only on (\ref{equ: intro FM L inv}).
We often write $\rho_{\un{\mathscr{L}}}$ for the unique (up to isomorphism) $\rho$ determined by the tuple $\un{\mathscr{L}}$ as above.
We use the terminology \emph{simple Fontaine-Mazur $\mathscr{L}$-invariants} for those $\mathscr{L}_{\al}$ with $\al=(i,i+1)\in\Delta$ for some $i$, and note that each such $\mathscr{L}_{\al}$ determines the isomorphism class of the $2$-dimensional subquotient of $\rho$ that fits into a non-split extension $0\rightarrow \varepsilon^{i}\rightarrow \ast\rightarrow \varepsilon^{i-1}\rightarrow 0$.
We call the rest of $\mathscr{L}_{\al}$ with $\al\in\Phi^+\setminus\Delta$ \emph{higher Fontaine-Mazur $\mathscr{L}$-invariants}.
The conjectural $p$-adic Langlands correspondence and their local-global compatibility suggest that $\Pi(U^{v},r)^{\rm{an}}$ should determine $\rho$ uniquely (and might even depends only on $\rho$ up to a finite copy).
In reality, the full structure of $\Pi(U^{v},r)^{\rm{an}}$ is very hard to understand and is only known when $n=2$ (and $K=\Q_p$).
Nevertheless, one might wish to construct a subrepresentation $\Pi(\rho)\subseteq \Pi(U^{v},r)^{\rm{an}}$ which determines $\rho$ and depends only on $\rho$. To summarize, we arrive at the following well-known question.
\begin{ques}\label{ques: intro construction Pi}
Let $\rho$ be a semi-stable $p$-adic Galois representation of the form (\ref{equ: intro n dim Galois}) such that $D_{\rm{st}}(\rho)$ satisfies $N^{n-1}\neq 0$.
How to construct $\Pi(\rho)\in\mathrm{Rep}^{\rm{an}}_{\rm{adm}}(G)$ that satisfies the following properties?
\begin{enumerate}[label=(\roman*)]
\item \label{it: intro construction Pi 1} The object $\Pi(\rho)$ is of finite length, with its semi-simplification $\Pi(\rho)^{\rm{ss}}$ being independent of the choice of $\rho$.
\item \label{it: intro construction Pi 2} The object $\Pi(\rho)$ determines and depends only on $\rho$ up to isomorphism.
\item \label{it: intro construction Pi 3} There exists a unique (up to scalars) embedding $\mathrm{St}_{G}^{\infty}\hookrightarrow\Pi(\rho)$ which identifies $\mathrm{St}_{G}^{\infty}$ with the $E$-subspace of smooth vectors of $\Pi(\rho)$.
\item \label{it: intro construction Pi 4} For each $F$, $v$, $\cG$, $U^{v}$ and $r$ such that $\Hom_{G}(\mathrm{St}_{G}^{\infty},\Pi(U^{v},r)^{\rm{an}})\neq 0$, the following cup product map
\[\Hom_{G}(\mathrm{St}_{G}^{\infty},\Pi(\rho))\otimes_E\Hom_{G}(\Pi(\rho),\Pi(U^{v},r)^{\rm{an}})\buildrel\cup\over\longrightarrow \Hom_{G}(\mathrm{St}_{G}^{\infty},\Pi(U^{v},r)^{\rm{an}})\]
is an isomorphism between $E$-vector spaces.
\end{enumerate}
\end{ques}
Using the fact that $\mathrm{St}_{G}^{\rm{cont}}$ is the universal completion of $\mathrm{St}_{G}^{\infty}$, it is easy to see that any embedding from $\mathrm{St}_{G}^{\infty}$ into $\Pi(U^{v},r)^{\rm{an}}$ extends to an embedding from $\mathrm{St}_{G}^{\rm{cont}}$ to $\Pi(U^{v},r)$ which further induces an embedding from $\mathrm{St}_{G}^{\rm{an}}$ to $\Pi(U^{v},r)^{\rm{an}}$ by taking locally analytic vectors. Hence, it is harmless for us to ask for $\Pi(\rho)$ that satisfies
\begin{equation}\label{equ: intro loc an St Pi}
\Dim_E\Hom_{G}(\mathrm{St}_{G}^{\rm{an}},\Pi(\rho))=1.
\end{equation}
\begin{rem}\label{rem: intro OS constituent}
We say a simple object of $\mathrm{Rep}^{\rm{an}}_{\rm{adm}}(G)$ is \emph{Orlik-Strauch} if it is a constituent of a locally analytic principal series representation induced from a locally algebraic character of $T$. (In fact, Orlik-Strauch have introduced a functor in \cite{OS15} that classifies all constituents of all such locally analytic principal series, whence our terminology.)
Based on the construction of various finite length objects of $\mathrm{Rep}^{\rm{an}}_{\rm{adm}}(G)$ in \cite{BQ24}, we expect that there exists an object $\Pi(\rho)\in\mathrm{Rep}^{\rm{an}}_{\rm{adm}}(G)$ (for each $\rho$) all of whose constituents being Orlik-Strauch such that all properties listed in Question~\ref{ques: intro construction Pi} hold.
\end{rem}
\subsubsection{Gehrmann's automorphic simple $\mathscr{L}$-invariants}\label{subsubsec: intro auto L inv}
We briefly recall Gehrmann's definition of automorphic $\mathscr{L}$-invariants (upon adapting our notation) to motivate our definition of Breuil-Schraen $\mathscr{L}$-invariants.
Let $F$, $v$, $\cG$ and $U^{v}$ be as in \S \ref{subsubsec: intro lgc question} with $\Gamma\defeq \cG(F)\cap(U^{v}\times\mathrm{Iw}_{v})$ with $\mathrm{Iw}_{v}\subseteq\cG(\cO_{F_{v}})$ being the Iwahori subgroup and $\Gamma^{v}\defeq \cG(F)\cap(U^{v}\times\cG(F_{v}))$.
We write $H^{\bullet}(\Gamma,-)$ and $H^{\bullet}(\Gamma^{v},-)$) for the (abstract) group cohomology with coefficients $E$ for the discrete groups $\Gamma$ and $\Gamma^{v}$, and both of them admits the action of a universal Hecke algebra $\bT(U^{v})$ by double coset operators.
We write $\mathrm{Ind}_{\Gamma}^{\Gamma^{v}}$ (resp.~$\mathrm{ind}_{\Gamma}^{\Gamma^{v}}$) for the induction functor (resp.~compact induction functor) and $\mathrm{ind}_{\mathrm{Iw}_{v}}^{\cG(F_{v})}$ for the (smooth) compact induction functor.
Then the natural maps
\[\mathrm{ind}_{\Gamma}^{\Gamma_{v}}1_{\Gamma}\rightarrow \mathrm{ind}_{\mathrm{Iw}_{v}}^{\cG(F_{v})}1_{\mathrm{Iw}_{v}}\rightarrow\mathrm{St}_{G}^{\infty}\]
induces (by taking $E$-dual) a map
\[(\mathrm{St}_{G}^{\infty})^{\vee}\rightarrow (\mathrm{ind}_{\Gamma}^{\Gamma_{v}}1_{\Gamma})^{\vee}=\mathrm{Ind}_{\Gamma}^{\Gamma_{v}}1_{\Gamma}\]
and thus the maps
\begin{equation}\label{equ: intro evaluation}
H^{\bullet}(\Gamma^v, (\mathrm{St}_{G}^{\infty})^{\vee})\rightarrow H^{\bullet}(\Gamma^v, \mathrm{Ind}_{\Gamma}^{\Gamma_{v}}1_{\Gamma})\buildrel\sim\over\longrightarrow H^{\bullet}(\Gamma, 1_{\Gamma})
\end{equation}
with the RHS isomorphism from Shapiro's Lemma for abstract group cohomologies (cf.~\cite[\S~6.3]{Wei94}).
Let $\fm\subseteq\bT(U^{v})$ be a maximal ideal.
We consider the following conditions on the choice of $F$, $v$, $\cG$, $U^{v}$, $\fm$ and a degree $k\geq 0$ (following the introduction of \cite{Geh21}).
\begin{cond}\label{cond: auto L inv}
Let $F$, $v$, $\cG$, $U^{v}$, $\fm$ and $k$ be as above.
\begin{enumerate}[label=(\roman*)]
\item \label{it: auto L inv setup 1} The $\fm$-isotypic component $H^{k}(\Gamma,1_{\Gamma})[\fm]$ is $1$-dimensional.
\item \label{it: auto L inv setup 2} The following map induced from (\ref{equ: intro evaluation})
\[H^{k}(\Gamma^v, (\mathrm{St}_{G}^{\infty})^{\vee})[\fm]\rightarrow H^{k}(\Gamma,1_{\Gamma})[\fm]\]
is an isomorphism between $1$-dimensional $E$-vector spaces.
\item \label{it: auto L inv setup 3} For each $I_1\subseteq I_0\subseteq \Delta$, the cup product map
\begin{equation}\label{equ: auto L sm cup}
\mathrm{Ext}_{G}^{\#I_0\setminus I_1}(V_{I_0}^{\infty},V_{I_1}^{\infty})\otimes_E H^{k+\#I_1}(\Gamma^v, (V_{I_1}^{\infty})^{\vee})[\fm]\buildrel\cup\over\longrightarrow H^{k+\#I_0}(\Gamma^v, (V_{I_0}^{\infty})^{\vee})[\fm]
\end{equation}
is an isomorphism between $1$-dimensional $E$-vector spaces.
\end{enumerate}
\end{cond}
A nice property of $H^{\bullet}(\Gamma^v,-)$ is that the embeddings
\[V_{I}^{\infty}\hookrightarrow V_{I}^{\rm{an}}\hookrightarrow V_{I}^{\rm{cont}}\]
induce the following isomorphisms (see \cite[Prop.~3.11]{Geh21})
\begin{equation}\label{equ: arithmetic coh completion}
H^{\bullet}(\Gamma^v, (V_{I}^{\rm{cont}})^{\vee})\buildrel\sim\over\longrightarrow H^{\bullet}(\Gamma^v, (V_{I}^{\rm{an}})^{\vee})\buildrel\sim\over\longrightarrow H^{\bullet}(\Gamma^v, (V_{I}^{\infty})^{\vee}).
\end{equation}
Suppose that the embedding $V_{I_0}^{\infty}\hookrightarrow V_{I_0}^{\rm{an}}$ and $V_{I_1}^{\infty}\hookrightarrow V_{I_1}^{\rm{an}}$ induce the following maps (which is true by results in this paper, see Lemma~\ref{lem: change left cup} and Lemma~\ref{lem: Ext subquotient})
\begin{equation}\label{equ: intro auto L inv Ext}
\mathbf{E}_{I_0,I_1}^{\infty}\defeq\mathrm{Ext}_{G}^{\#I_0\setminus I_1}(V_{I_0}^{\infty},V_{I_1}^{\infty})\hookrightarrow\mathrm{Ext}_{G}^{\#I_0\setminus I_1}(V_{I_0}^{\infty},V_{I_1}^{\rm{an}})\buildrel\sim\over\longleftarrow \mathbf{E}_{I_0,I_1}\defeq\mathrm{Ext}_{G}^{\#I_0\setminus I_1}(V_{I_0}^{\rm{an}},V_{I_1}^{\rm{an}}),
\end{equation}
we observe that (\ref{equ: arithmetic coh completion}) induces the following commutative diagram
\begin{equation}\label{equ: auto L inv diagram}
\xymatrix{
\mathbf{E}_{I_0,I_1} & \otimes_E & H^{k+\#I_1}(\Gamma^v, (V_{I_1}^{\rm{an}})^{\vee})[\fm] \ar^{\cup}[rr] \ar^{\wr}[d] & & H^{k+\#I_0}(\Gamma^v, (V_{I_0}^{\rm{an}})^{\vee})[\fm] \ar^{\wr}[d]\\
\mathbf{E}_{I_0,I_1}^{\infty} \ar@{^{(}->}[u] & \otimes_E & H^{k+\#I_1}(\Gamma^v, (V_{I_1}^{\infty})^{\vee})[\fm] \ar^{\cup}[rr] & & H^{k+\#I_0}(\Gamma^v, (V_{I_0}^{\infty})^{\vee})[\fm]
}
\end{equation}
with the bottom horizontal cup product map being an isomorphism between $1$-dimension $E$-vector spaces by \ref{it: auto L inv setup 3} of Condition~\ref{cond: auto L inv}.
We define $W_{I_0,I_1}^{k}[\fm]$ as the maximal subspace of $\mathbf{E}_{I_0,I_1}$ that annihilate $H^{k+\#I_1}(\Gamma^v, (V_{I_1}^{\rm{an}})^{\vee})[\fm]$ under the top horizontal cup product map in (\ref{equ: auto L inv diagram}).
The commutative diagram (\ref{equ: auto L inv diagram}) implies that $W_{I_0,I_1}^{k}[\fm]$ is a codimension $1$ subspace of $\mathbf{E}_{I_0,I_1}$ that satisfies $W_{I_0,I_1}^{k}[\fm]\cap \mathbf{E}_{I_0,I_1}^{\infty}=0$ (or equivalently $W_{I_0,I_1}^{k}[\fm]+\mathbf{E}_{I_0,I_1}^{\infty}=\mathbf{E}_{I_0,I_1}$).
This collection of annihilators
\begin{equation}\label{equ: intro auto annihilator}
\{W_{I_0,I_1}^{k}[\fm]\}_{I_1\subseteq I_0\subseteq\Delta}
\end{equation}
is what Gehrmann has called \emph{automorphic $\mathscr{L}$-invariants} in \cite[Def.~3.12]{Geh21}, except that he has restricted himself to the $\#I_0\setminus I_1=1$ cases possibly due to the lack of results on $\mathbf{E}_{I_0,I_1}$ for general $I_1\subseteq I_0$ at the time.
For each $I_1\subseteq I\subseteq I_0$, we consider the following cup product map
\[\kappa_{I_0,I_1}^{I}: \mathbf{E}_{I_0,I}\otimes_E \mathbf{E}_{I,I_1} \buildrel\cup\over\longrightarrow \mathbf{E}_{I_0,I_1}.\]
Given arbitrary $x\in\mathbf{E}_{I_0,I}$, $y\in\mathbf{E}_{I,I_1}$ and $v\in H^{k+\#I_1}(\Gamma^v, (V_{I_1}^{\rm{an}})^{\vee})[\fm]$, we have
\[(x\cup y)\cup v=x\cup (y\cup v)\]
with $x\cup y\in \mathbf{E}_{I_0,I_1}$ and $y\cup v\in H^{k+\#I}(\Gamma^v, (V_{I}^{\rm{an}})^{\vee})[\fm]$. This together with the definition of $W_{I_0,I_1}^{k}[\fm]$, $W_{I_0,I}^{k}[\fm]$ and $W_{I,I_1}^{k}[\fm]$ gives
\begin{equation}\label{equ: intro annihilator relation}
\kappa_{I_0,I_1}^{I}(W_{I_0,I}^{k}[\fm]\otimes_E\mathbf{E}_{I,I_1}+\mathbf{E}_{I_0,I}\otimes_EW_{I,I_1}^{k}[\fm])\subseteq W_{I_0,I_1}^{k}[\fm]\cap\mathrm{im}(\kappa_{I_0,I_1}^{I}).
\end{equation}
Suppose now that there exists some $r: \mathrm{Gal}_{F}\rightarrow \cG^{\vee}(E)$ with $\rho\defeq r|_{\mathrm{Gal}_{F_v}}$ such that $\fm=\fm_{r}$ is the maximal ideal associated with $r$. Combining an (optimistic) expectation from $p$-adic local Langlands correspondence with the (unique up to a sign) isomorphism
\begin{equation}\label{equ: intro normalize Ext}
\mathbf{E}_{I_0,I_1}\cong \mathbf{E}_{I_0\setminus I_1}\defeq \mathbf{E}_{\Delta,\Delta\setminus(I_0\setminus I_1)}=\mathrm{Ext}_{G}^{\#I_0\setminus I_1}(1_{G},V_{\Delta\setminus(I_0\setminus I_1)}^{\rm{an}})
\end{equation}
to be constructed in this paper (see \ref{it: intro cup comparison 1} of Theorem~\ref{thm: intro cup comparison} below), one might expect that $W_{I_0,I_1}^{k}[\fm]\subseteq \mathbf{E}_{I_0\setminus I_1}$, whose definition involves $F$, $v$, $\cG$, $U^{v}$, $\fm$ and $k$, to depend only on $\rho=r|_{\mathrm{Gal}_{F_v}}$ and $I_0\setminus I_1$, namely to have the form $W_{I_0\setminus I_1}(\rho)\subseteq \mathbf{E}_{I_0\setminus I_1}$ under (\ref{equ: intro normalize Ext}).
Using the isomorphism (\ref{equ: intro normalize Ext}) and similarly $\mathbf{E}_{I_0,I}\cong \mathbf{E}_{I_0\setminus I}$ as well as $\mathbf{E}_{I,I_1}\cong \mathbf{E}_{I\setminus I_1}$, we might expect the existence of a map
\[\kappa_{J,J'}:\mathbf{E}_{J}\otimes_E\mathbf{E}_{J'}\buildrel\cup\over\longrightarrow \mathbf{E}_{J\sqcup J'}\]
for each $J,J'\in\Delta$ with $J\cap J'=\emptyset$, which recovers $\kappa_{I_0,I_1}^{I}$ (up to a sign) upon taking $J=I_0\setminus I$ and $J'=I\setminus I_1$. Assuming the existence of such $\kappa_{J,J'}$, the relation (\ref{equ: intro annihilator relation}) becomes
\begin{equation}\label{equ: intro annihilator relation transform}
\kappa_{J,J'}(W_{J}(\rho)\otimes_E\mathbf{E}_{J'}+\mathbf{E}_{J}\otimes_EW_{J'}(\rho))\subseteq W_{J\sqcup J'}(\rho)\cap\mathrm{im}(\kappa_{J,J'})
\end{equation}
for each $J,J'\in\Delta$ with $J\cap J'=\emptyset$.
Such collection of annihilators
\begin{equation}\label{equ: intro auto L inv annihilator}
\{W_{I}(\rho)\subseteq\mathbf{E}_{I}\}_{I\subseteq \Delta}
\end{equation}
leads us directly to the notion of Breuil-Schraen $\mathscr{L}$-invariants (see Definition~\ref{def: intro BS L inv} below).
Suppose we can construct $\kappa_{J,J'}$ and prove its injectivity, then the inclusion (\ref{equ: intro annihilator relation transform}) would be necessarily an equality as BHS of (\ref{equ: intro annihilator relation transform}) are $E$-subspaces of codimension $1$ in $\mathrm{im}(\kappa_{J,J'})$. Exchanging $J,J'$, we optimistically expect the equality $\mathrm{im}(\kappa_{J,J'})=\mathrm{im}(\kappa_{J',J})$ as well as the following equalities 
\begin{multline}\label{equ: intro annihilator relation exchange}
\kappa_{J,J'}(W_{J}(\rho)\otimes_E\mathbf{E}_{J'}+\mathbf{E}_{J}\otimes_EW_{J'}(\rho))=W_{J\sqcup J'}(\rho)\cap\mathrm{im}(\kappa_{J,J'})\\
=W_{J\sqcup J'}(\rho)\cap\mathrm{im}(\kappa_{J',J})=\kappa_{J',J}(W_{J'}(\rho)\otimes_E\mathbf{E}_{J}+\mathbf{E}_{J'}\otimes_EW_{J}(\rho))
\end{multline}
to hold for each $\rho$ and $J,J'\in\Delta$ with $J\cap J'=\emptyset$. With $\rho$ varying, the expected equalities (\ref{equ: intro annihilator relation exchange}) further suggest the following equality
\[
E\kappa_{J,J'}(x\otimes y)=E\kappa_{J',J}(y\otimes x)
\]
for each $x\in\mathbf{E}_{J}$ and $y\in\mathbf{E}_{J'}$.
All these expectations on $\kappa_{J,J'}$ indeed hold by results of this paper, and we refer to \S~\ref{subsubsec: intro cup} for more detailed discussions.

To summarize, Gehrmann's definition of automorphic $\mathscr{L}$-invariants as annihilators (see (\ref{equ: intro auto annihilator})) justifies the importance of understanding $\mathbf{E}_{I_0,I_1}$ for each $I_1\subseteq I_0\subseteq \Delta$ as well as the cup product maps $\kappa_{I_0,I_1}^{I}$ for each $I_1\subseteq I\subseteq I_0$, with an important aspect being the construction and the study of a \emph{normalized cup product map} $\kappa_{J,J'}$ for each $J,J'\in\Delta$ with $J\cap J'=\emptyset$.
Nevertheless, the following natural question remains.
\begin{ques}\label{ques: intro normalize L inv}
How to characterize the annihilators (\ref{equ: intro auto L inv annihilator}) from $\rho$ or equivalently the Fontaine-Mazur $\mathscr{L}$-invariants of $\rho$?
\end{ques}
\subsubsection{Drinfeld realization and Schraen's $\mathscr{L}$-invariants for $\mathrm{GL}_{3}(\Q_p)$}\label{subsubsec: intro Drinfeld realization}
Now we discuss the relation between (higher) $\mathscr{L}$-invariants and Drinfeld upper-half space, with the goal being the (re)formulation of several questions including a refinement of Question~\ref{ques: intro normalize L inv}.
Throughout this discussion, we assume $K=\Q_p$ for simplicity.
The Drinfeld upper-half space is defined as $\bH_{n-1}\defeq \bP^{n-1}_{\rm{rig}}\setminus(\bigcup_{H} H)$ with $H$ running through $\Q_p$-rational hyperplane in $\bP^{n-1}_{\rm{rig}}$.
Needless to say, Drinfeld upper-half space (and its coverings) is one of the most fundamental (and relatively well-understood) object in the study of both classical and more recently $p$-adic local Langlands correspondence (cf.~\cite{Sch18}, \cite{DL17}, \cite{CDN20a}, \cite{V24}, etc).
We write $W_{K}$ for the Weil group of $K$.
In the classical $\ell$-adic setting with $\ell\neq p$, writing $R\Gamma_{c}(-,\overline{\Q}_{\ell})$ for the $\ell$-adic cohomology with compact support and $D^{b}_{\overline{\Q}_{\ell}}(G)$ for the bounded derived category of the abelian category $\mathrm{Rep}_{\overline{\Q}_{\ell}}(G)$ of smooth $\overline{\Q}_{\ell}$-representations of $G$, Dat has proved in \cite[Thm.~1.1]{Dat06} that we have the following isomorphism between continuous $\overline{\Q}_{\ell}$-representations of $W_{K}$
\begin{equation}\label{equ: intro Dat isom}
H^{\bullet}(R\Hom_{D^{b}_{\overline{\Q}_{\ell}}(G)}(R\Gamma_{c}(\bH_{n-1},\overline{\Q}_{\ell}),\pi))\buildrel\sim\over\longrightarrow \sigma_{n}(\pi)\otimes_E|\cdot|^{\frac{n-1}{2}}
\end{equation}
where $\pi=V_{I}^{\infty}$ for some $I\subseteq\Delta$, and $\sigma_{n}(\pi)$ is the $\overline{\Q}_{\ell}$-representation of $W_{K}$ associated with $\pi$ via the classical local Langlands correspondence for $\mathrm{GL}_{n}(K)$ (see \cite{HT01}, \cite{Hen00}, \cite{Sch13}).
It has been a natural question to ask for a $p$-adic analogue of (\ref{equ: intro Dat isom}).
Note that $\bH_{n-1}$ is quasi-Stein and the higher cohomologies of coherent sheaves on $\bH_{n-1}$ all vanish. Hence, we are interested in the global sections $R\Gamma_{\rm{dR}}(\bH_{n-1})$ of the de Rham complex of $\bH_{n-1}$ (with coefficients extended to $E$), which is known to be a complex of coadmissible $D(G)$-modules in the sense of Schneider-Teitelbaum. In fact, each coadmissible $D(G)$-module in the complex $R\Gamma_{\rm{dR}}(\bH_{n-1})$ is finite length and can be explicitly constructed from Orlik-Strauch functors (see \cite{ST02}, \cite{Po04}, \cite{Or08} and \cite{BQ24}).
The search for a $p$-adic analogue of (\ref{equ: intro Dat isom}) together with Question~\ref{ques: intro construction Pi} leads to the following question (see the introduction of \cite{Schr10} for more historical account).
\begin{ques}\label{ques: intro Pi Drinfeld}
Is it possible to construct $\{\Pi(\rho)\}_{\rho}$ with (\ref{equ: intro loc an St Pi}) such that
\begin{equation}\label{equ: intro Pi Drinfeld}
\Hom_{\cM(G)}(\Pi(\rho)^{\vee}[1-n], R\Gamma_{\rm{dR}}(\bH_{n-1}))\cong D_{\rm{dR}}(\rho)
\end{equation}
as filtered $(\varphi,N)$-modules (up to a twist) for each $\rho$?
\end{ques}
Here the Hodge filtration $\mathrm{Fil}_{H}^{\ell}(-)$ on LHS of (\ref{equ: intro Pi Drinfeld}) should be given by (the image of) the following maps
\[\Hom_{\cM(G)}(\Pi(\rho)^{\vee}[1-n], \sigma_{\geq \ell} R\Gamma_{\rm{dR}}(\bH_{n-1}))\rightarrow \Hom_{\cM(G)}(\Pi(\rho)^{\vee}[1-n], R\Gamma_{\rm{dR}}(\bH_{n-1}))\]
which we expect to be an embedding for each $\ell$.
The $(\varphi,N)$-action on $R\Gamma_{\rm{dR}}(\bH_{n-1})$ is more delicate and involves the Hyodo-Kato complex $R\Gamma_{\rm{HK}}(\bH_{n-1})$ of $\bH_{n-1}$ (see \cite{GK05} and \cite{CDN20b}) with coefficients extended to $E$. We have the following Hyoto-Kato isomorphism (see \cite[Thm.~0.1]{GK05})
\begin{equation}\label{equ: intro HK isom}
R\Gamma_{\rm{HK}}(\bH_{n-1})\cong R\Gamma_{\rm{dR}}(\bH_{n-1})
\end{equation}
as well as a canonical $\varphi$-equivariant splitting of $R\Gamma_{\rm{HK}}(\bH_{n-1})$
\[R\Gamma_{\rm{HK}}(\bH_{n-1})\cong \bigoplus_{k=0}^{n-1}H_{\rm{HK}}^{k}(\bH_{n-1})[-k]\]
with $\varphi$ acting on $H_{\rm{HK}}^{k}(\bH_{n-1})[-k]$ by $p^{k}$ and $N$ satisfying $N\varphi=p\varphi N$ with $N^{n-1}\neq 0$ (see \cite{dS05}).
Through $R\Gamma_{\rm{HK}}(\bH_{n-1})$, one obtains a splitting (which is an isomorphism in $\cM(G)$)
\begin{equation}\label{equ: intro dR splitting}
s_{\rm{HK}}: R\Gamma_{\rm{dR}}(\bH_{n-1})\cong \bigoplus_{k=0}^{n-1}H_{\rm{dR}}^{k}(\bH_{n-1})[-k]
\end{equation}
and thus a $(\varphi,N)$-action on $R\Gamma_{\rm{dR}}(\bH_{n-1})$ through $s_{\rm{HK}}$, with $\varphi$ acting on $H_{\rm{dR}}^{k}(\bH_{n-1})[-k]$ by $p^{k}$ for each $0\leq k\leq n-1$ and $N$ satisfying $N\varphi=p\varphi N$ as well as $N^{n-1}\neq 0$.
In fact, the existence of an abstract splitting
\begin{equation}\label{equ: intro dR abstract splitting}
s: R\Gamma_{\rm{dR}}(\bH_{n-1})\cong \bigoplus_{k=0}^{n-1}H_{\rm{dR}}^{k}(\bH_{n-1})[-k]
\end{equation}
in $\cM(G)$ is known (without invoking $R\Gamma_{\rm{HK}}(\bH_{n-1})$) and can be fully classified by the results of Dat (see \cite[Thm.~1.3]{Dat06}) and Orlik (see \cite{Or05a}). However, it is not evident how to explicitly characterize $s_{\rm{HK}}$ among all possible splittings of the form (\ref{equ: intro dR abstract splitting}). (Note that the isomorphism (\ref{equ: intro HK isom}) and thus $s_{\rm{HK}}$ depends on the choice of a uniformizer of $\Q_p$.)

We can also ask the following universal version of Question~\ref{ques: intro Pi Drinfeld} with the benefit being that we could (only for this question to make sense) ignore the problem of explicitly characterizing the splitting $s_{\rm{HK}}$.
\begin{ques}\label{ques: intro Pi Drinfeld universal}
We consider an arbitrary abstract splitting $s$ as in (\ref{equ: intro dR abstract splitting}) and equip $R\Gamma_{\rm{dR}}(\bH_{n-1})$ with a $(\varphi,N)$-action through $s$, with the $(\varphi,N)$-action on RHS of (\ref{equ: intro dR splitting}) characterized by $\varphi$ acting on $H_{\rm{dR}}^{k}(\bH_{n-1})[-k]$ by $p^{k}$ for each $0\leq k\leq n-1$ and $N$ satisfying $N\varphi=p\varphi N$ as well as $N^{n-1}\neq 0$.
Then what can we say about the structure of
\begin{equation}\label{equ: intro Pi Drinfeld universal}
\Hom_{\cM(G)}((\mathrm{St}_{G}^{\rm{an}})^{\vee}[1-n], R\Gamma_{\rm{dR}}(\bH_{n-1}))
\end{equation}
as a filtered $(\varphi,N)$-module?
\end{ques}
In \cite[Prop.~6.21]{Schr11}, Schraen has provided us with an candidate for the answer to Question~\ref{ques: intro Pi Drinfeld universal} and has confirmed it when $n=3$ (and $K=\Q_p$).
We define $\cE_{0}\defeq E$ and $\mathbf{V}_{0}\defeq 1_{\mathrm{Gal}_{\Q_p}}$. Then we define $\cE_{n}$ and $\mathbf{V}_{n}$ inductively for $n\geq 1$ by
\[\cE_{n}\defeq \mathrm{Ext}_{\mathrm{Gal}_{\Q_p}}^1(\mathbf{V}_{n-1},\varepsilon^{n})^\vee\]
and the following universal extension of finite dimensional continuous $E$-representations of $\mathrm{Gal}_{\Q_p}$
\[0\rightarrow \cE_{n}\otimes_E\varepsilon^{n}\rightarrow \mathbf{V}_{n}\rightarrow \mathbf{V}_{n-1} \rightarrow 0.\]
The following conjecture has not been explicitly stated in the literature but is somewhat hidden behind \cite[Prop.~6.21]{Schr11}
\begin{conj}\label{conj: intro Drinfeld}
For an arbitrary abstract splitting $s$ as in (\ref{equ: intro dR abstract splitting}) and the $(\varphi,N)$-action on $R\Gamma_{\rm{dR}}(\bH_{n-1})$ induced by $s$ (as described in Question~\ref{ques: intro Pi Drinfeld universal}), we have
\begin{equation}\label{equ: intro Drinfeld conj}
\Hom_{\cM(G)}((\mathrm{St}_{G}^{\rm{an}})^{\vee}[1-n], R\Gamma_{\rm{dR}}(\bH_{n-1}))\cong D_{\rm{dR}}(\mathbf{V}_{n-1}\otimes_E\varepsilon^{1-n})
\end{equation}
as filtered $(\varphi,N)$-modules.
\end{conj}
Then we can reformulate \cite[Prop.~6.21]{Schr11} as follows.
\begin{thm}[Schraen]\label{thm: intro Schraen}
Conjecture~\ref{conj: intro Drinfeld} holds for $n=3$.
\end{thm} 
One main result of this paper is to confirm Conjecture~\ref{conj: intro Drinfeld} for arbitrary $n$ (see Theorem~\ref{thm: intro Drinfeld} below).

Combining Question~\ref{ques: intro construction Pi} with Conjecture~\ref{conj: intro Drinfeld}, we arrive at the following more refined conjecture.
\begin{conj}\label{conj: intro Drinfeld dR}
There exists a family of admissible locally analytic representations $\{\Pi(\rho)\}_{\rho}$ with $\rho$ running through non-critical semi-stable $p$-adic continous Galois representations of the form (\ref{equ: intro n dim Galois}) satisfying $N^{n-1}\neq 0$ on $D_{\rm{st}}(\rho)$, such that the following statements hold.
\begin{enumerate}[label=(\roman*)]
\item \label{it: intro Drinfeld dR 1} For each $\rho$, $\Pi(\rho)$ has finite length with all constituents being Orlik-Strauch and its semi-simplification $\Pi(\rho)^{\rm{ss}}$ being independent of the choice of $\rho$.
\item \label{it: intro Drinfeld dR 2} There exists a unique (up to scalars) embedding $\mathrm{St}_{G}^{\rm{an}}\hookrightarrow \Pi(\rho)$, and the induced map
    \begin{equation}\label{equ: intro Pi rho kernel}
    \mathbf{E}_{I}\cong\mathbf{E}_{I,\emptyset}=\mathrm{Ext}_{G}^{\#I}(V_{I}^{\rm{an}},\mathrm{St}_{G}^{\rm{an}})\rightarrow \mathrm{Ext}_{G}^{\#I}(V_{I}^{\rm{an}},\Pi(\rho))
    \end{equation}
    has $1$-dimensional image with kernel denoted by $W_{I}(\rho)\subseteq \mathbf{E}_{I}$.
\item \label{it: intro Drinfeld dR 3} For each $\rho$, we have the following commutative diagram of maps between filtered $(\varphi,N)$-modules (with the $(\varphi,N)$-action on $R\Gamma_{\rm{dR}}(\bH_{n-1})$ induced from $s_{\rm{HK}}$)
\begin{equation}\label{equ: intro Drinfeld dR diagram}
\xymatrix{
\mathbf{E}_{n-1}\defeq\mathbf{E}_{\Delta}=\mathrm{Ext}_{G}^{n-1}(1_{G},\mathrm{St}_{G}^{\rm{an}}) \ar@{=}[d] & & &\\
\Hom_{\cM(G)}((\mathrm{St}_{G}^{\rm{an}})^{\vee}[1-n], H_{\rm{dR}}^{0}(\bH_{n-1})) \ar^{\sim}[rrr] \ar@{^{(}->}[d] & & & \cE_{n-1} \ar@{^{(}->}[d]\\
\Hom_{\cM(G)}((\mathrm{St}_{G}^{\rm{an}})^{\vee}[1-n], R\Gamma_{\rm{dR}}(\bH_{n-1})) \ar^{\sim}[rrr] \ar@{->>}[d] & & & D_{\rm{dR}}(\mathbf{V}_{n-1}\otimes_E\varepsilon^{1-n}) \ar@{->>}[d]\\
\Hom_{\cM(G)}((\Pi(\rho))^{\vee}[1-n], R\Gamma_{\rm{dR}}(\bH_{n-1})) \ar^{\sim}[rrr] & & & D_{\rm{dR}}(\rho\otimes_E\varepsilon^{1-n})
}
\end{equation}
with the kernel the composition of LHS vertical maps being $W_{\Delta}(\rho)$ from \ref{it: intro Drinfeld dR 2}.
\item \label{it: intro Drinfeld dR 4} The association $\rho\leadsto\Pi(\rho)\leadsto W_{\Delta}(\rho)$ is a one to one correspondence, which is uniquely determined by a unique (up to scalars) isomorphism
    \begin{equation}\label{equ: intro regulator map}
    r_{n-1}: \mathbf{E}_{n-1}\buildrel\sim\over\longrightarrow \cE_{n-1}.
    \end{equation}
\end{enumerate}
\end{conj}
Following Conjecture~\ref{conj: intro Drinfeld dR}, we have the following question which refines Question~\ref{ques: intro normalize L inv}.
\begin{ques}\label{ques: intro regulator}
Assume that Conjecture~\ref{conj: intro Drinfeld dR} holds. How to explicitly describe the isomorphism (\ref{equ: intro regulator map}) for each $n\geq 2$?
\end{ques}
Note that an explicit characterization of (\ref{equ: intro regulator map}) is closely related to the question of explicitly describing the Hyodo-Kato splitting $s_{\rm{HK}}$, and if possible, would in particular leads to an explicit relation between Breuil-Schraen $\mathscr{L}$-invariants (as in (\ref{equ: intro auto L inv annihilator})) and Fontaine-Mazur $\mathscr{L}$-invariants of $\rho$.
\subsubsection{Breuil-Ding's higher $\mathscr{L}$-invariants for $\mathrm{GL}_{3}(\Q_p)$}\label{subsubsec: intro Breuil Ding L inv}
We assume $n>2$ in this section.
When $n=3$, Question~\ref{ques: intro construction Pi} has been completely solved by Breuil and Ding in \cite{BD20}. So far, this is the only solved case of Question~\ref{ques: intro construction Pi} (assuming $n>2$), partially because Breuil-Ding's construction of the family $\{\Pi(\rho)\}_{\rho}$ uses the (family version of) $p$-adic local Langlands correspondence for $\mathrm{GL}_{2}(\Q_p)$ as a key ingredient.

We briefly recall Breuil-Ding's construction below.
Assume from now that $n=3$ with $G=\mathrm{PGL}_{3}(\Q_p)$.
For each $i\in\Delta=\{1,2\}$, we write $P_{i}\defeq P_{\{i\}}$ and $L_{i}\defeq L_{\{i\}}$ for the associated parabolic and Levi subgroups of $G$, and note that we naturally have $L_{i}\cong\mathrm{GL}_{2}(\Q_p)$.
For each $i\in\Delta$, we write $\rho_{i}$ for the unique $2$-dimensional subquotient of $\rho$ determined by $L_{i}$.
Using the $p$-adic local Langlands correspondence for $\mathrm{GL}_{2}(\Q_p)$ (and its family version), for each $i\in\Delta$ one can associate with $\rho$ an explicit finite length admissible locally analytic representation $\Pi_{L_{i}}(\rho)$ of $L_{i}\cong \mathrm{GL}_{2}(\Q_p)$, which comes with a natural isomorphism between $E$-vector spaces
\begin{equation}\label{equ: intro GL2 deform isom}
\mathrm{Ext}_{L_{i}}^{1}(\Pi_{L_{i}}(\rho),\Pi_{L_{i}}(\rho))\cong \mathrm{Ext}_{\mathrm{Gal}_{\Q_p}}^{1}(\rho_{i},\rho_{i}).
\end{equation}
One could check that $(\mathrm{Ind}_{P_{i}}^{G}\Pi_{L_{i}}(\rho))^{\rm{an}}$ has socle $V_{3-i}^{\infty}$. One defines $\Pi_{i}(\rho)^{-}$ as the unique quotient of $(\mathrm{Ind}_{P_{i}}^{G}\Pi_{L_{i}}(\rho))^{\rm{an}}$ with socle $\mathrm{St}_{G}^{\infty}$. 
By applying the functor $(\mathrm{Ind}_{P_{i}}^{G}-)^{\rm{an}}$ and then taking subquotient, one obtains a map
\begin{equation}\label{equ: intro induction subquotient}
\mathrm{Ext}_{L_{i}}^{1}(\Pi_{L_{i}}(\rho),\Pi_{L_{i}}(\rho))\rightarrow\mathrm{Ext}_{G}^{1}(V_{3-i}^{\infty},\Pi_{i}(\rho)^{-}).
\end{equation}
Using local Tate pairing (see \cite[\S 3.2.2]{BD20}), Breuil-Ding associate with $\rho$ a hyperplane in the RHS of (\ref{equ: intro GL2 deform isom}), and then use (\ref{equ: intro GL2 deform isom}) and (\ref{equ: intro induction subquotient}) to construct (see \cite[\S 3.3.4]{BD20}) a locally analytic representation $\Pi_{i}(\rho)$ of $G$ that fits into
\[0\rightarrow \Pi_{i}(\rho)^{-}\rightarrow \Pi_{i}(\rho)\rightarrow (V_{3-i}^{\infty})^{\oplus 2}.\]
Note that the locally $\Q_p$-analytic representation of $G$ constructed by Breuil-Ding in \emph{loc.cit.} is in fact a subrepresentation of $\Pi_{i}(\rho)$ which carries the same information, and it is harmless to ignore this minor difference.
Finally, they construct $\Pi(\rho)$ by taking suitable amalgamate sum of $\Pi_{1}(\rho)$ and $\Pi_{2}(\rho)$.
Under the isomorphism (\ref{equ: intro GL2 deform isom}), Breuil-Ding have emphasized in \cite[\S~5.5]{BD19} that a deformation on LHS has central character (resp.~has infinitesimal character) if and only if its corresponding deformation on RHS fixes determinant (resp.~is Hodge-Tate). We thus obtain strong constraints on the following conjectural isomorphism (from (\ref{equ: intro regulator map}) by taking $n=3$)
\begin{equation}\label{equ: intro GL3 regulator}
r_{2}:\mathbf{E}_{2}\buildrel\sim\over\longrightarrow \cE_{2}.
\end{equation}
Motivated by Breuil-Ding's results above for $\mathrm{GL}_{3}(\Q_p)$, we will construct (see \S \ref{subsubsec: intro Drinfeld} and \S \ref{subsubsec: intro BS} below) a pair of $E$-subspaces 
\begin{equation}\label{equ: intro Galois pair}
\cE_{n-1}^{\flat},\cE_{n-1}^{\sharp}\subseteq\cE_{n-1}
\end{equation}
which for $n=3$ are certain universal version of the orthogonal complement (under local Tate pairing) of the subspace of Hodge-Tate deformations and fixed determinant deformations (see the discussion in \cite[\S 5.5]{BD19}). We will also construct their automorphic counterpart, namely a pair of $E$-subspaces 
\[\mathbf{E}_{n-1}^{\flat},\mathbf{E}_{n-1}^{\sharp}\subseteq \mathbf{E}_{n-1}\]
that are expected to match the pair (\ref{equ: intro Galois pair}) under $r_{n-1}$.
We will propose a candidate to the answer of Question~\ref{ques: intro regulator} in Theorem~\ref{thm: intro regulator} below.
\subsection{Main results}\label{subsec: intro main results}
With some historical background and motivation from previous results given, we are ready to give a quick overview of the main results of this paper.
\subsubsection{Cohomology of locally analytic generalized Steinberg}\label{subsubsec: intro coh St}
Let $J\subseteq\Delta$. We recall the locally $K$-analytic generalized Steinberg representation $V_{\Delta\setminus J}^{\rm{an}}$ from \S~\ref{subsubsec: intro quick notation}. It is well-known (see \cite{OS13}) that the surjection $i_{\Delta\setminus J}^{\rm{an}}\twoheadrightarrow V_{\Delta\setminus J}^{\rm{an}}$ extends to a quasi-isomorphism of the following form
\begin{equation}\label{equ: intro Tits resolution}
\mathbf{C}_{J}\defeq [i_{\Delta}^{\rm{an}}\rightarrow \cdots \rightarrow \bigoplus_{I\supseteq \Delta\setminus J,\#I=\ell}i_{I}^{\rm{an}}\rightarrow \cdots \rightarrow i_{\Delta\setminus J}^{\rm{an}}]\buildrel\sim\over\longrightarrow V_{\Delta\setminus J}^{\rm{an}}[\#\Delta\setminus J]
\end{equation}
with $\bigoplus_{I\supseteq \Delta\setminus J,\#I=\ell}i_{I}^{\rm{an}}$ sitting in degree $-\ell$ for each $\#\Delta\setminus J\leq \ell\leq \#\Delta$.
The main actor of this paper is the following $\mathrm{Ext}$-group
\begin{equation}\label{equ: intro St coh}
\mathbf{E}_{J}\defeq\mathrm{Ext}_{G}^{\#J}(1_{G},V_{\Delta\setminus J}^{\rm{an}})\buildrel\sim\over\longleftarrow \mathrm{Ext}_{G}^{2\#J-\#\Delta}(1_{G},\mathbf{C}_{J}).
\end{equation}
We can compute (\ref{equ: intro St coh}) via a spectral sequence $E_{\bullet,J}^{\bullet,\bullet}$ (see \S~\ref{subsubsec: Tits double} for its definition and further discussions). In particular, $\mathbf{E}_{J}$ is equipped with a canonical decreasing filtration $\mathrm{Fil}^{\bullet}(\mathbf{E}_{J})$ with
\[\mathrm{gr}^{-\ell}(\mathbf{E}_{J})=E_{\infty,J}^{-\ell,\ell+2\#J-\#\Delta}\]
for each $\ell$.
We define $\cS_{J}$ as the set of all subsets $S\subseteq \Phi^+$ that satisfy 
\[\sum_{\al\in S}\al=\al_{J}\defeq \sum_{\beta\in J}\beta\in \Z_{\geq 0}\Phi^+.\]
In other words, $\cS_{J}$ is the set of all partitions of $\al_{J}$ into positive roots.
Our first main theorem focuses on the computation of the spectral sequence $E_{\bullet,J}^{\bullet,\bullet}$, with the key output being a combinatorial formula for the dimension of $\mathbf{E}_{J}$ as an $E$-vector space. 
\begin{thm}[Corollary~\ref{cor: bottom deg E2 basis}, (\ref{equ: atom to partition}), Proposition~\ref{prop: bottom deg degeneracy}]\label{thm: intro main dim}
Let $J\subseteq\Delta$. We have
\begin{equation}\label{equ: intro main dim graded}
\Dim_E\mathrm{gr}^{-\ell}(\mathbf{E}_{J})=\#\{(S,I)\mid S\in\cS_{J}, \#S=\#\Delta-\ell, I\subseteq S\cap\Delta\}
\end{equation}
for each $0\leq\ell\leq\#\Delta$. In particular, we have
\begin{equation}\label{equ: intro main dim}
\Dim_{E}\mathbf{E}_{J}=\#\{(S,I)\mid S\in\cS_{J}, I\subseteq S\cap\Delta\}.
\end{equation}
\end{thm}
The proof of the equality (\ref{equ: intro main dim graded}) for each $\ell$ has two main ingredients.
First, we explicitly construct (through combinatorial methods) a basis of $E_{2,J}^{-\ell,\ell+2\#J-\#\Delta}$ of the form
\begin{equation}\label{equ: intro main dim basis}
\{\overline{x}_{S,I}\}_{S\in\cS_{J}, \#S=\#\Delta-\ell, I\subseteq S\cap\Delta}.
\end{equation}
Then we prove
\begin{equation}\label{equ: intro E2 degenerate}
E_{2,J}^{-\ell,\ell+2\#J-\#\Delta}=E_{\infty,J}^{-\ell,\ell+2\#J-\#\Delta}
\end{equation}
which is a $E_2$-degeneracy result for this particular bi-degree $(-\ell,\ell+2\#J-\#\Delta)$.
We refer further details on both ingredients to \S~\ref{subsubsec: intro atom E2}.
We emphasize here that the construction of a special basis of $\mathbf{E}_{J}$ is indeed possible (cf.~Theorem~\ref{thm: moduli of inv}, Remark~\ref{rem: relation to FM L inv} and the discussion before Question~\ref{ques: BS FM match}) but much more involved than the combinatorial construction of the basis (\ref{equ: intro main dim basis}) of (\ref{equ: intro E2 degenerate}) (namely the graded piece $\mathrm{gr}^{-\ell}(\mathbf{E}_{J})$) for each $\ell$.

One interesting feature of $\mathbf{E}_{J}$ is that it is simultaneously equipped with several genuinely different structures yet they interact with each other in non-trivial ways.
Before passing to the second structure on $\mathbf{E}_{J}$, we need some simple facts on the combinatorics of the set 
\[\Gamma\defeq\{x\in W(G)\mid \ell(x)=\#\mathrm{Supp}(x)\}\] 
of \emph{partial-Coxeter elements} in the Weyl group $W(G)$ of $G$. By definition, $\Gamma$ consists of those elements of $W(G)$ whose arbitrary reduced decompositions (into product of simple reflections) are multiplicity free.
We equip $\Gamma$ with the partial-order $\unlhd$ such that $x\unlhd w$ for some $x,w\in\Gamma$ if and only if $w=yx$ for some $y\in\Gamma$ that satisfies $\ell(w)=\ell(x)+\ell(y)$.
For each $J\subseteq\Delta$, we define $\Gamma^{J}\defeq \{x\in\Gamma\mid\mathrm{Supp}(x)\subseteq J\}$ which inherits a partial-order from $\Gamma$. The following combinatorial fact is easy to prove yet crucial for our study of $\mathbf{E}_{J}$ (see Theorem~\ref{thm: intro coxeter filtration} below).
\begin{prop}[Proposition~\ref{prop: coxeter partition cardinality}, Proposition~\ref{prop: cardinality}]\label{prop: intro coxeter partition}
For each $J\subseteq\Delta$, we have
\begin{equation}\label{equ: intro coxeter partition}
\#\Gamma^{J}=\#\{(S,I)\mid S\in\cS_{J}, I\subseteq S\cap\Delta\}.
\end{equation}
\end{prop}

For each $x\in\Gamma^{J}$ and $I\subseteq J_{x}\defeq \Delta\setminus\mathrm{Supp}(x)$, we write $L^{I}(x)\defeq L^{I}(x\cdot 0)$ for the finite dimensional simple $U(\fl_{I})$-module with highest weight $x\cdot 0$, and $M^{I}(x)\defeq U(\fg)\otimes_{U(\fp_{I})}L^{I}(x)$ for the associated parabolic Verma module. We write $M(x)\defeq M^{\emptyset}(x)$ for short.
Note that the $U(\fl_{I})$-action on $L^{I}(x)$ extends to a $D(L_{I})$-action, so that $L^{I}(x)^{\vee}$ is naturally an \emph{algebraic representation} of $L_{I}$. For each $I\subseteq J_{x}$, we write
\[i_{x,I}^{\rm{an}}\defeq (\mathrm{Ind}_{P_{I}}^{G} L^{I}(x)^{\vee})^{\rm{an}}\]
and then define the following \emph{locally $K$-analytic generalized Steinberg representation with highest weight $x\cdot 0$}
\[
V_{x,I}^{\rm{an}}\defeq i_{x,I}^{\rm{an}}/\sum_{I\subsetneq I'\subseteq J_{x}}i_{x,I'}^{\rm{an}}
\]
Recall from \cite{OS13} that $V_{x,I}^{\rm{an}}$ fits into a quasi-isomorphism of the form
\[
\mathbf{C}_{x,I}\defeq [i_{x,J_{x}}^{\rm{an}}\rightarrow\cdots\rightarrow \bigoplus_{I\subseteq I'\subseteq J_{x},\#I'=\ell}i_{x,I'}^{\rm{an}}\rightarrow\cdots\rightarrow i_{x,I}^{\rm{an}}]\buildrel\sim\over\longrightarrow V_{x,I}^{\rm{an}}[\#I]
\]
with $\bigoplus_{I\subseteq I'\subseteq J_{x},\#I'=\ell}i_{x,I'}^{\rm{an}}$ sitting in degree $-\ell$.
Then we define
\[\mathbf{E}_{x,J}\defeq\mathrm{Ext}_{G}^{\#J}(1_{G},V_{x,\Delta\setminus J}^{\rm{an}})\]
which recovers $\mathbf{E}_{J}$ when $x=1$.
The following results summarize the relation between $V_{x,\Delta\setminus J}^{\rm{an}}$ and $\mathbf{E}_{x,J}$ for different choices of $x\in\Gamma^{J}$.
\begin{prop}\label{prop: intro coxeter map}
Let $J\subseteq\Delta$ and $x,w\in\Gamma^{J}$ with $x\unlhd w$. We have the following results.
\begin{enumerate}[label=(\roman*)]
\item \label{it: intro coxeter map 1} We have a unique (up to scalars) non-zero map
\begin{equation}\label{equ: intro coxeter St map}
V_{x,\Delta\setminus J}^{\rm{an}}\rightarrow V_{w,\Delta\setminus J}^{\rm{an}}.
\end{equation}
\item \label{it: intro coxeter map 2} The map (\ref{equ: intro coxeter St map}) induces a surjection
\begin{equation}\label{equ: intro coxeter St coh surjection}
\mathbf{E}_{x,J}\twoheadrightarrow \mathbf{E}_{w,J}.
\end{equation}
\end{enumerate}
\end{prop}

Let $J\subseteq \Delta$. We define
\begin{equation}\label{equ: intro coxeter filtration def}
\mathrm{Fil}_{x}(\mathbf{E}_{J})\defeq \bigcap_{y\in\Gamma^{J},y\not\unlhd x}\mathrm{ker}(\mathbf{E}_{J}\twoheadrightarrow \mathbf{E}_{y,J})\subseteq \mathbf{E}_{J}
\end{equation}
for each $x\in\Gamma^{J}$, and it is clear that $\mathrm{Fil}_{x}(\mathbf{E}_{J})\subseteq \mathrm{Fil}_{w}(\mathbf{E}_{J})$ for each pair of elements $x,w\in\Gamma^{J}$ that satisfy $x\unlhd w$. In other words, the collection of $E$-subspaces 
\begin{equation}\label{equ: intro coxeter filtration}
\{\mathrm{Fil}_{x}(\mathbf{E}_{J})\}_{x\in\Gamma^{J}}
\end{equation}
forms an increasing filtration of $\mathbf{E}_{J}$ with respect to the partial-order $\unlhd$ on $\Gamma^{J}$.
We define
\[\mathrm{gr}_{x}(\mathbf{E}_{J})\defeq \mathrm{Fil}_{x}(\mathbf{E}_{J})/\sum_{u\unlhd x,u\neq x}\mathrm{Fil}_{u}(\mathbf{E}_{J})\]
for each $x\in\Gamma^{J}$. 
The following is our main result on the increasing filtration (\ref{equ: intro coxeter filtration}) on $\mathbf{E}_{J}$.
\begin{thm}\label{thm: intro coxeter filtration}
Let $J\subseteq\Delta$. Then we have
\[
\sum_{u\in\Gamma^{J}}\mathrm{Fil}_{u}(\mathbf{E}_{J})=\mathbf{E}_{J}
\]
and
\begin{equation}\label{equ: intro coxeter filtration graded}
\Dim_E\mathrm{gr}_{x}(\mathbf{E}_{J})=1
\end{equation}
for each $x\in\Gamma^{J}$.
\end{thm}
The increasing filtration (\ref{equ: intro coxeter filtration}) is what we call \emph{Coxeter filtration} on $\mathbf{E}_{J}$.
Our study of it is naturally motivated by that of the Hodge filtration of the de Rham complex of the Drinfeld space (see (\ref{equ: intro auto flat}) and the discussion below it).
Note that our definition of Coxeter filtration on $\mathbf{E}_{J}$ along with our formulation of Proposition~\ref{prop: intro coxeter map} and Theorem~\ref{thm: intro coxeter filtration} is not exactly following the logical order of what we actually have done in the main text (see Theorem~\ref{thm: coxeter filtration} and Proposition~\ref{prop: coxeter filtration x surjection}).
We have decided to formulate Proposition~\ref{prop: intro coxeter map} and Theorem~\ref{thm: intro coxeter filtration} as they are to minimize the amount of technical discussions.
In reality, we first define the Coxeter filtration (\ref{equ: intro coxeter filtration}) on $\mathbf{E}_{J}$ via the layer structure of $V_{\Delta\setminus J}^{\rm{an}}$, using the combinatorial facts (\ref{equ: intro main dim}) and (\ref{equ: intro coxeter partition}) as crucial inputs (see Theorem~\ref{thm: coxeter filtration}). Then we verify that such Coxeter filtration can be equally described as (\ref{equ: intro coxeter filtration def}) in Proposition~\ref{prop: coxeter filtration x surjection}.
We refer the interested readers to \S~\ref{subsubsec: intro coxeter} for further details.

Finally, we would like to mention that $\mathbf{E}_{J}$ admits a third collection of $E$-subspaces defined via relative conditions with respect to $E$-Lie subalgebras $\fh\subseteq\fg$ of the form $\fh=[\fl_{I},\fl_{I}]$ with $I\subseteq \Delta$ being an interval (namely $\sum_{\al\in I}\al\in\Phi^+$). This third kind of structure on $\mathbf{E}_{J}$ plays a crucial (yet quite more technical) role in the proof of several of our main results.
We refer the interested readers to different parts of \S~\ref{subsec: intro technical} (and notably \S~\ref{subsubsec: relative support}) for an extensive discussion on this third structure on $\mathbf{E}_{J}$.
\subsubsection{Cup product maps}\label{subsubsec: intro cup}
Motivated by our discussion on Gehrmann's automorphic $\mathscr{L}$-invariants (see \S~\ref{subsubsec: intro auto L inv}), the second key player in our paper is a \emph{normalized cup product map} of the form
\begin{equation}\label{equ: intro main cup}
\kappa_{J,J'}: \mathbf{E}_{J}\otimes_E\mathbf{E}_{J'}\buildrel\cup\over\longrightarrow\mathbf{E}_{J\sqcup J'}
\end{equation}
for each pair of subsets $J,J'\subseteq\Delta$ with $J\cap J'=\emptyset$ (see \S~\ref{subsubsec: normalized cup} for its precise definition). 
We present below several of our main results on (\ref{equ: intro main cup}) which basically say that (\ref{equ: intro main cup}) have all the nice properties that we might expect.

For each $I_1\subseteq I_0\subseteq\Delta$, we recall
\[\mathbf{E}_{I_0,I_1}=\mathrm{Ext}_{G}^{\#I_0\setminus I_1}(V_{I_0}^{\rm{an}},V_{I_1}^{\rm{an}})\]
from (\ref{equ: intro auto L inv Ext}).
One original purpose of introducing a map of the form (\ref{equ: intro main cup}) is to streamline the study of the following cup product map
\begin{equation}\label{equ: intro main cup Ext}
\kappa_{I_0,I_1}^{I}: \mathbf{E}_{I_0,I}\otimes_E\mathbf{E}_{I,I_1}\buildrel\cup\over\longrightarrow \mathbf{E}_{I_0,I_1}
\end{equation}
for each $I_1\subseteq I\subseteq I_0$. This is confirmed by the following result.
\begin{thm}[Proposition~\ref{prop: Ext complex std seq}, Proposition~\ref{prop: Tits cup transfer}]\label{thm: intro cup comparison}
Let $I_1\subseteq I_0\subseteq \Delta$. We have the following results.
\begin{enumerate}[label=(\roman*)]
\item \label{it: intro cup comparison 1} We have an isomorphism
\[\mathbf{E}_{I_0,I_1}\buildrel\sim\over\longrightarrow\mathbf{E}_{I_0\setminus I_1}\]
which is canonically defined up to a sign.
\item \label{it: intro cup comparison 2} For each $I_1\subseteq I\subseteq I_0$, we have the following commutative diagram
\begin{equation}\label{equ: intro cup comparison}
\xymatrix{
\mathbf{E}_{I_0,I} \ar^{\wr}[d] & \otimes_E & \mathbf{E}_{I,I_1} \ar^{\kappa_{I_0,I_1}^{I}}[rrr] \ar^{\wr}[d] & & & \mathbf{E}_{I_0,I_1} \ar^{\wr}[d]\\
\mathbf{E}_{I_0\setminus I} & \otimes_E & \mathbf{E}_{I\setminus I_1} \ar^{\kappa_{I_0\setminus I, I\setminus I_1}}[rrr] & & & \mathbf{E}_{I_0\setminus I_1}
}
\end{equation}
with all the vertical maps being isomorphisms from \ref{it: intro cup comparison 1} (up to signs).
\end{enumerate}
\end{thm}
We note that the proof of \ref{it: intro cup comparison 2} of Theorem~\ref{thm: intro cup comparison} builds on the construction of a big commutative diagram (see (\ref{equ: 5 complex cup total 1})) which involves delicate truncation maps between certain $5$-complex.

It is a fairly formal consequence of the definition of (\ref{equ: intro main cup}) that it restricts to a map 
\begin{equation}\label{equ: intro main cup Fil}
\mathrm{Fil}^{-\ell_0}(\mathbf{E}_{J})\otimes_E\mathrm{Fil}^{-\ell_1}(\mathbf{E}_{J'})\buildrel\cup\over\longrightarrow \mathrm{Fil}^{-\ell_2}(\mathbf{E}_{J\sqcup J'})
\end{equation}
for each $\#\Delta\setminus J\leq \ell_0\leq \#\Delta$ and $\#\Delta\setminus J'\leq \ell_1\leq \#\Delta$ with $\ell_2=\ell_0+\ell_1-\#\Delta$, and thus further induces a map
\begin{equation}\label{equ: intro graded cup}
\mathrm{gr}(\mathbf{E}_{J})\otimes_E\mathrm{gr}(\mathbf{E}_{J'})\buildrel\cup\over\longrightarrow \mathrm{gr}(\mathbf{E}_{J\sqcup J'})
\end{equation}
with $\mathrm{gr}(\mathbf{E}_{J})\defeq\bigoplus_{\ell_0=\#\Delta\setminus J}^{\#\Delta\setminus J}\mathrm{gr}^{-\ell_0}(\mathbf{E}_{J})$ and similarly for others.
Similarly, exchanging the role of $J$ and $J'$, we see that $\kappa_{J',J}$ also induces a map
\begin{equation}\label{equ: intro graded cup prime}
\mathrm{gr}(\mathbf{E}_{J'})\otimes_E\mathrm{gr}(\mathbf{E}_{J})\buildrel\cup\over\longrightarrow \mathrm{gr}(\mathbf{E}_{J'\sqcup J}).
\end{equation}
Following Theorem~\ref{thm: intro main dim} and the discussion below it, it is natural to ask how the explicit basis of $\mathrm{gr}(\mathbf{E}_{J})$, $\mathrm{gr}(\mathbf{E}_{J'})$ and $\mathrm{gr}(\mathbf{E}_{J\sqcup J'})$ are related under (\ref{equ: intro graded cup}) and (\ref{equ: intro graded cup prime}). This is answered by the following result.
\begin{thm}[Lemma~\ref{lem: grade cup commute}, Theorem~\ref{thm: general cup}]\label{thm: intro cup graded}
Let $J,J'\subseteq \Delta$ with $J\cap J'=\emptyset$. We have the following results.
\begin{enumerate}[label=(\roman*)]
\item \label{it: intro cup graded 0} For each $x\in \mathrm{gr}(\mathbf{E}_{J})$ and $y\in \mathrm{gr}(\mathbf{E}_{J'})$, the image of $x\otimes y$ under (\ref{equ: intro graded cup}) equals the image of $y\otimes x$ under (\ref{equ: intro graded cup prime}).
\item \label{it: intro cup graded 1} Given $S\in\cS_{J}$, $I\subseteq S\cap\Delta$, $S'\in\cS_{J'}$ and $I'\subseteq S'\cap\Delta$, the image of $\overline{x}_{S,I}\otimes\overline{x}_{S',I'}$ under (\ref{equ: intro graded cup}) equals $\overline{x}_{S\sqcup S',I\sqcup I'}$ up to a sign. In particular, the map (\ref{equ: intro graded cup}) is injective.
\item \label{it: intro cup graded 2} The map (\ref{equ: intro main cup}) is injective.
\item \label{it: intro cup graded 3} If there exists $j_0\in\Delta$ such that $j<j_0<j'$ for each $j\in J$ and $j'\in J'$, then both $\kappa_{J,J'}$ and $\kappa_{J',J}$ are isomorphisms.
\end{enumerate}
\end{thm}
We refer interested readers to \S~\ref{subsubsec: normalized cup} for further discussions on the proof of Theorem~\ref{thm: intro cup graded}.

For each $J\subseteq\Delta$, we consider the following $E$-subspace spanned by \emph{decomposable elements}
\[\mathbf{E}_{J}^{<}\defeq\sum_{\emptyset\neq J'\subsetneq J}\mathrm{im}(\kappa_{J',J\setminus J})\subseteq\mathbf{E}_{J}\]
with $\overline{\mathbf{E}}_{J}\defeq \mathbf{E}_{J}/\mathbf{E}_{J}^{<}$ being the \emph{primitive quotient}. Theorem~\ref{thm: intro cup graded} has the following consequence.
\begin{cor}[Lemma~\ref{lem: E2 cup of generator}, Lemma~\ref{lem: primitive generator}]\label{cor: intro cup generator}
Let $J\subseteq\Delta$. We have the following results.
\begin{enumerate}[label=(\roman*)]
\item \label{it: intro cup generator 1} If $\#J=1$, then $\mathbf{E}_{J}^{<}=0$ and we have $\mathbf{E}_{J}=\Hom(Z_{\Delta\setminus J},E)\cong \Hom(K^{\times},E)$.
\item \label{it: intro cup generator 2} If $\#J\geq 2$, then $\mathbf{E}_{J}^{<}=\mathrm{Fil}^{-(n-3)}(\mathbf{E}_{J})$ with $\overline{\mathbf{E}}_{J}=\mathrm{gr}^{-(n-2)}(\mathbf{E}_{J})\neq 0$ if and only if $\sum_{\al\in J}\al\in\Phi^+\setminus\Delta$, in which case we have $\Dim_E\overline{\mathbf{E}}_{J}=1$.
\end{enumerate}
\end{cor}

Let $J,J'\subseteq\Delta$ with $J\cap J'=\emptyset$.
For each $x\in\Gamma^{J}$ and $y\in\Gamma^{J'}$, we define
\begin{equation}\label{equ: intro x y envelop}
\Gamma_{x,y}\defeq \{w\in\Gamma\mid \mathrm{Supp}(w)=\mathrm{Supp}(x)\sqcup\mathrm{Supp}(y), x\leq w, y\leq w\}\subseteq \Gamma^{J\sqcup J'}.
\end{equation}
Following Proposition~\ref{prop: intro coxeter map} and Theorem~\ref{thm: intro coxeter filtration}, the following result gives a fairly clean description of the relation between the Coxeter filtration on $\mathbf{E}_{J}$, $\mathbf{E}_{J'}$ and $\mathbf{E}_{J,J'}$ under the cup product map $\kappa_{J,J'}$ (see (\ref{equ: intro main cup})).
\begin{thm}[(\ref{equ: x y cup w bottom diagram}), Proposition~\ref{prop: coxeter cup comparison}, Proposition~\ref{prop: cup x y nonvanishing}]\label{thm: intro cup x y}
Let $J,J'\subseteq \Delta$ with $J\cap J'=\emptyset$. Let $x\in\Gamma^{J}$ and $y\in\Gamma^{J'}$. We have the following results.
\begin{enumerate}[label=(\roman*)]
\item \label{it: intro cup x y 1} For each $w\in\Gamma_{x,y}\subseteq \Gamma^{J\sqcup J'}$, we have a commutative diagram of the form 
    \begin{equation}\label{equ: intro cup x y diagram}
    \xymatrix{
    \mathbf{E}_{J} \ar@{->>}[d] & \otimes_E & \mathbf{E}_{J'} \ar@{->>}[d] \ar^{\kappa_{J,J'}}[rr] & & \mathbf{E}_{J\sqcup J'} \ar@{->>}[d] \\
    \mathbf{E}_{x,J} & \otimes_E & \mathbf{E}_{y,J'} \ar[rr] & & \mathbf{E}_{w,J\sqcup J'}
    }
    \end{equation}
    which is functorial with respect to the choice of $x$, $y$ and $w$, with all vertical maps being surjections from \ref{it: intro coxeter map 2} of Proposition~\ref{prop: intro coxeter map}.
\item \label{it: intro cup x y 2} The map $\kappa_{J,J'}$ restricts to a map
\begin{equation}\label{equ: intro cup x y}
\mathrm{Fil}_{x}(\mathbf{E}_{J})\otimes_E\mathrm{Fil}_{y}(\mathbf{E}_{J'})\hookrightarrow \sum_{u\in\Gamma_{x,y}}\mathrm{Fil}_{u}(\mathbf{E}_{J\sqcup J'})
\end{equation}
whose composition with
\[\sum_{u\in\Gamma_{x,y}}\mathrm{Fil}_{u}(\mathbf{E}_{J\sqcup J'})\twoheadrightarrow \mathrm{gr}_{w}(\mathbf{E}_{J\sqcup J'})\]
is surjective for each $w\in\Gamma_{x,y}$.
\end{enumerate}
\end{thm}
We refer interested readers to \S~\ref{subsubsec: interaction} for further discussions on the proof of Theorem~\ref{thm: intro cup x y}.
Taking $x=1$ in (\ref{equ: intro cup x y}) and using $\Gamma_{1,y}=\{y\}$, we obtain the inclusion
\[\kappa_{J,J'}(\mathrm{Fil}_{1}(\mathbf{E}_{J})\otimes_E\mathrm{Fil}_{y}(\mathbf{E}_{J'}))\subseteq \mathrm{Fil}_{y}(\mathbf{E}_{J\sqcup J'}).\]
Now that $\kappa_{J,J'}$ is injective (see \ref{it: intro cup graded 2} of Theorem~\ref{thm: intro cup graded}) and we have $\Dim_E\mathrm{Fil}_{1}(\mathbf{E}_{J})=1$ as well as
\[\Dim_E\mathrm{Fil}_{y}(\mathbf{E}_{J'}))=\#\{u\in\Gamma^{J'}\mid u\unlhd y\}=\Dim_E\mathrm{Fil}_{y}(\mathbf{E}_{J\sqcup J'})\]
by Theorem~\ref{thm: intro coxeter filtration}, we conclude that
\begin{equation}\label{equ: intro sm cup}
\kappa_{J,J'}(\mathrm{Fil}_{1}(\mathbf{E}_{J})\otimes_E\mathrm{Fil}_{y}(\mathbf{E}_{J'}))=\mathrm{Fil}_{y}(\mathbf{E}_{J\sqcup J'}).
\end{equation}

Theorem~\ref{thm: intro cup x y} also has the following interesting consequence.
\begin{thm}[Theorem~\ref{thm: cup top grade}]\label{thm: intro top graded}
Let $J\subseteq \Delta$ be a non-empty interval and $w\in\Gamma^{J}$ with $\mathrm{Supp}(w)=J$.
Then the composition of the following maps
\[
\mathrm{Fil}_{w}(\mathbf{E}_{J})\hookrightarrow \mathbf{E}_{J}\twoheadrightarrow \overline{\mathbf{E}}_{J}=\mathbf{E}_{J}/\mathbf{E}_{J}^{<}
\]
is surjective.
\end{thm}

Building on Theorem~\ref{thm: intro cup graded} and Theorem~\ref{thm: intro cup x y} as well as a delicate ingredient involving $E$-subspaces of $\mathbf{E}_{J}$, $\mathbf{E}_{J'}$ and $\mathbf{E}_{J\sqcup J'}$ defined using relative conditions with respect to $E$-Lie subalgebras of $\fg$, we are able to prove the following commutativity result which is predicted by our discussion of Gehrmann's automorphic $\mathscr{L}$-invariants (see the discussion around (\ref{equ: intro annihilator relation exchange})).
\begin{thm}[Theorem~\ref{thm: cup exchange}]\label{thm: intro cup commute}
Let $J,J'\subseteq\Delta$ with $J\cap J'=\emptyset$. Then we have
\begin{equation}\label{equ: intro cup commute}
\kappa_{J,J'}(x\otimes y)=\kappa_{J',J}(y\otimes x)\in\mathbf{E}_{J\sqcup J'}
\end{equation}
for each $x\in\mathbf{E}_{J}$ and $y\in\mathbf{E}_{J'}$.
\end{thm}
We refer interested readers to \S~\ref{subsubsec: commute} for an outline of the proof of Theorem~\ref{thm: intro cup commute}.

\begin{rem}\label{rem: Qp an case}
Each result listed in \S~\ref{subsubsec: intro coh St} and \S~\ref{subsubsec: intro cup} admits a locally $\Q_p$-analytic version, with the following modifications. We write $G_0$ for the underlying locally $\Q_p$-analytic group of $G$, and $\mathrm{Ext}_{G_0}^{\bullet}(-,-)$ for the corresponding locally $\Q_p$-analytic $\mathrm{Ext}$-groups between two admissible locally $\Q_p$-analytic representations. We write $i_{I,0}^{\rm{an}}$ (resp.~$V_{I,0}^{\rm{an}}$) for the locally $\Q_p$-analytic principal series (resp.~locally $\Q_p$-analytic generalized Steinberg) for each $I\subseteq\Delta$ with $\mathbf{E}_{J,0}\defeq\mathrm{Ext}_{G_0}^{\#J}(1_{G},V_{\Delta\setminus J,0}^{\rm{an}})$. 
Assume that $E$ contains the Galois closure of $K$ and write $\Sigma$ for the set of all embeddings $K\hookrightarrow E$.
We write
\[\cS_{J,0}\defeq \{S\subseteq \Phi^+\times\Sigma\mid \sum_{(\al,\sigma)\in S}\al=\sum_{\beta\in J}\beta\}\]
for each $J\subseteq\Delta$.
We summarize the locally $\Q_p$-analytic variant of our previous results as below.
\begin{enumerate}[label=(\roman*)]
\item \label{it: Qp an 1} Recall from Theorem~\ref{thm: intro main dim} that $\mathbf{E}_{J}$ admits a canonical filtration $\mathrm{Fil}^{\bullet}(\mathbf{E}_{J})$ with the (total) graded pieces $\mathrm{gr}(\mathbf{E}_{J})=\bigoplus_{\ell}\mathrm{gr}^{-\ell}(\mathbf{E}_{J})$ admitting a basis $\{\overline{x}_{\ast}\}$ indexed by
    \[\ast\in\{(S,I)\mid S\in\cS_{J},I\subseteq S\cap\Delta\}\buildrel\sim\over\longrightarrow\bigsqcup_{I\subseteq J}\cS_{J\setminus I}.\]
    Similarly, $\mathbf{E}_{J,0}$ admits a canonical filtration with the (total) graded pieces $\mathrm{gr}(\mathbf{E}_{J,0})$ admitting a basis $\{\overline{x}_{\ast}\}$ indexed by
    \[\ast\in\bigsqcup_{I\subseteq J}\cS_{J\setminus I,0}.\]
\item \label{it: Qp an 2} We write $W(G_0)\defeq W(G)\times\Sigma$. For each $\un{x}=(x_{\sigma})_{\sigma\in\Sigma}\in W(G_0)$, we write $\mathrm{Supp}(\un{x})=\bigcup_{\sigma\in\Sigma}\mathrm{Supp}(x_{\sigma})$ and $\ell(\un{x})=\sum_{\sigma\in\Sigma}\ell(x_{\sigma})$. Then we define the set of \emph{partial-Coxeter elements in $W(G_0)$} by
    \[\Gamma_{0}\defeq \{\un{x}\in W(G_0)\mid \#\mathrm{Supp}(\un{x})=\ell(\un{x})\}\]
    which is equipped with a similar partial-order $\unlhd$ as that of $\Gamma$. For each $J\subseteq\Delta$, we also define
    \[\Gamma_{0}^{J}\defeq \{\un{x}\in\Gamma_{0}\mid \mathrm{Supp}(\un{x})\subseteq J\}.\]
    For each $\un{x}\in\Gamma_{0}$ and $I\subseteq J_{\un{x}}\defeq\Delta\setminus\mathrm{Supp}(\un{x})$, we can construct an admissible locally $\Q_p$-analytic representation $V_{\un{x},I,0}^{\rm{an}}$ which recovers $V_{I,0}^{\rm{an}}$ when $\un{x}=1$.
    For each $J\subseteq\Delta$ and $\un{x}\in\Gamma_{0}^{J}$, we write $\mathbf{E}_{\un{x},J,0}\defeq\mathrm{Ext}_{G_0}^{\#J}(1_{G},V_{\un{x},\Delta\setminus J,0}^{\rm{an}})$ and have a locally $\Q_p$-analytic version of Proposition~\ref{prop: intro coxeter map}. Finally, similar to Proposition~\ref{prop: intro coxeter partition} and Theorem~\ref{thm: intro coxeter filtration}, for each $J\subseteq\Delta$ we have
    \[\#\Gamma_{0}^{J}=\sum_{I\subseteq J}\#\cS_{J\setminus I,0}=\Dim_{E}\mathbf{E}_{J,0}\]
    and $\mathbf{E}_{J,0}$ admits a Coxeter filtration of the form
    \[\{\mathrm{Fil}_{\un{x}}(\mathbf{E}_{J,0})\}_{\un{x}\in\Gamma_{0}^{J}}.\]
\item \label{it: Qp an 3} For each $J,J'\subseteq\Delta$ with $J\cap J'=\emptyset$, we have a locally $\Q_p$-analytic normalized cup product map
    \[\kappa_{J,J',0}: \mathbf{E}_{J,0}\otimes_{E}\mathbf{E}_{J,0}\buildrel\cup\over\longrightarrow \mathbf{E}_{J\sqcup J',0}\]
    and Theorem~\ref{thm: intro cup comparison}, Theorem~\ref{thm: intro cup graded}, Theorem~\ref{thm: intro cup x y} and Theorem~\ref{thm: intro cup commute} all admit locally $\Q_p$-analytic version with parallel statements.
\end{enumerate}
In particular, parallel to Definition~\ref{def: intro BS L inv} below, we have a well-defined notion of \emph{locally $\Q_p$-analytic Breuil-Schraen $\mathscr{L}$-invariants}.
\end{rem}
\subsubsection{Drinfeld realization of certain universal filtered $(\varphi,N)$-module}\label{subsubsec: intro Drinfeld}
Recall from \S~\ref{subsubsec: intro quick notation} that we have fixed an embedding $\iota: K\hookrightarrow E$.
We also recall $D(G)$ and $\cM(G)$ from \S~\ref{subsubsec: intro quick notation}.
Here we follow the notation of \cite{BQ24} which is slightly different from that of \cite{Schr11} and \S \ref{subsubsec: intro Drinfeld realization}.
We write $\bH=\bH_{n-1}=\bP^{n-1}_{\rm{rig}}\setminus(\bigcup_{H} H)$ for the Drinfeld upper-half space with $H$ running through $K$-rational hyperplane in $\bP^{n-1}_{\rm{rig}}$. The space $\bH$ is quasi-Stein so we study the global sections of de Rham complex of $\bH$
\begin{equation}\label{equ: intro dR complex}
\Omega^{\bullet}=[\Omega^0\rightarrow \cdots\rightarrow \Omega^{n-1}]
\end{equation}
which has been denoted by $R\Gamma_{\rm{dR}}(\bH_{n-1})$ in \S \ref{subsubsec: intro Drinfeld realization}.
It is a famous theorem of Schneider--Stuhler (see \cite{SS91}) that we have
\begin{equation}\label{equ: intro dR coh}
H^k(\Omega^{\bullet})\cong (V_{[1,n-1-k]}^{\infty})^\vee
\end{equation}
for each $k\geq 0$. Combining results from \cite{Dat06} (or \cite{Or05a}) and \cite[\S 6.1]{Schr11}, we know that $\Omega^{\bullet}$ \emph{splits} (abstractly and non-canonically) in $\cM(G)$, namely there exist isomorphisms in $\cM(G)$
\begin{equation}\label{equ: intro dR split}
\Omega^{\bullet}\buildrel s\over\longrightarrow \bigoplus_{k=0}^{n-1}H^k(\Omega^{\bullet})[-k]\cong \bigoplus_{k=0}^{n-1}(V_{[1,n-1-k]}^{\infty})^\vee[-k].
\end{equation}
It follows from \cite{Dat06} (or \cite{Or05a}) and \cite[\S 6.1]{Schr11} that
\begin{equation}\label{equ: intro End coh}
\mathrm{End}_{\cM(G)}(\bigoplus_{k=0}^{n-1}(V_{[1,n-1-k]}^{\infty})^\vee[-k])\cong \fb
\end{equation}
as unital $E$-algebras, with $\fb$ understood as the lower triangular matrix algebra. (In particular, there exist many different choices of isomorphisms (\ref{equ: intro dR split}).) 

We consider a pair of elements $\varphi,N$ in (\ref{equ: intro End coh}) satisfying $\varphi|_{H^k(\Omega^{\bullet})[-k]}=p^k$ for each $0\leq k\leq n-1$ and $N\varphi=p\varphi N$. Here $N$ restricts to the zero map on $H^0(\Omega^{\bullet})$, and restricts to a map 
\[H^k(\Omega^{\bullet})[-k]\rightarrow H^{k-1}(\Omega^{\bullet})[-k+1]\]
corresponding to a non-zero element of
\[\Hom_{\cM(G)}(H^k(\Omega^{\bullet})[-k],H^{k-1}(\Omega^{\bullet})[-k+1])=\mathrm{Ext}_{G}^{1}(V_{[1,n-k]}^{\infty},V_{[1,n-1-k]}^{\infty})\cong\Hom_{\infty}(K^{\times},E)\]
for each $1\leq k\leq n-1$. 
Consequently, for each choice of the splitting isomorphism $s$ as in (\ref{equ: intro dR split}), we obtain a $(\varphi,N)$-action on
\begin{equation}\label{equ: intro Hom in derived}
\Hom_{\cM(G)}(D,\Omega^{\bullet})
\end{equation}
for each $D\in \cM(G)$.

When $D=(\mathrm{St}_{G}^{\rm{an}})^\vee[1-n]$, the naive truncation $\sigma_{\geq \ell}\Omega^{\bullet}\rightarrow \Omega^{\bullet}$ induces an embedding
\[\Hom_{\cM(G)}((\mathrm{St}_{G}^{\rm{an}})^\vee[1-n],\sigma_{\geq \ell}\Omega^{\bullet})\hookrightarrow \Hom_{\cM(G)}((\mathrm{St}_{G}^{\rm{an}})^\vee[1-n],\Omega^{\bullet})\]
for each $\ell$ (see \ref{it: Hodge seq 0} of Lemma~\ref{lem: Hodge exact seq}).
Consequently, the following space
\begin{equation}\label{equ: intro main Hom}
\Hom_{\cM(G)}((\mathrm{St}_{G}^{\rm{an}})^\vee[1-n],\Omega^{\bullet})
\end{equation}
is equipped with the structure of a filtered $(\varphi,N)$-module, with its Hodge filtration given by
\[\mathrm{Fil}^{\ell}_{H}(\Hom_{\cM(G)}((\mathrm{St}_{G}^{\rm{an}})^\vee[1-n],\Omega^{\bullet}))\defeq \Hom_{\cM(G)}((\mathrm{St}_{G}^{\rm{an}})^\vee[1-n],\sigma_{\geq \ell}\Omega^{\bullet})\]
for each $\ell$.

We are interested in the de Rham $(\varphi,\Gamma)$-module associated with the filtered $(\varphi,N)$-module (\ref{equ: intro main Hom}) via Berger's functor \cite{Ber09} (see also the second paragraph at the beginning of \S~\ref{subsec: universal Galois}).
Recall that each $p$-adic continuous character $\delta: K^{\times}\rightarrow E^{\times}$ uniquely determines a rank one $(\varphi,\Gamma)$-module over $\cR=\cR_{K,E}$ denoted by $\cR(\delta)$ (cf.~\cite[Notation~6.2.2]{KPX14}).
We write $|\cdot|$ (resp.~$z$) for the continuous character $K^{\times}\rightarrow E^{\times}$ given by $x\mapsto |N_{K/\Q_p}(x)|_{p}$ (resp.~by $x\mapsto \iota(x)$) and $\theta\defeq |\cdot|z$.
We write $\mathrm{Ext}_{g}^1(-,-)$ for the \emph{de Rham extensions} between two de Rham $(\varphi,\Gamma)$-modules over $\cR$.
We define $\cE_{0}\defeq E$ and $\cD_{0}\defeq \cR(1_{K^{\times}})$. Then we define $\cE_{n}$ and $\cD_{n}$ inductively for $n\geq 1$ by
\[\cE_{n}\defeq \mathrm{Ext}_{g}^1(\cD_{n-1},\cR(\theta^{n}))^\vee\]
and the following universal de Rham extension between de Rham $(\varphi,\Gamma)$-modules
\[0\rightarrow \cE_{n}\otimes_E\cR(\theta^{n})\rightarrow \cD_{n}\rightarrow \cD_{n-1} \rightarrow 0.\]
\begin{thm}[Theorem~\ref{thm: main filtered phi N}]\label{thm: intro Drinfeld}
For each choice of the splitting isomorphism $s$ as in (\ref{equ: intro dR split}), the de Rham $(\varphi,\Gamma)$-module associated with the filtered $(\varphi,N)$-module (\ref{equ: intro main Hom}) is isomorphic to $\cD_{n-1}\otimes_{\cR} \cR(\theta^{1-n})$.
In particular, Conjecture~\ref{conj: intro Drinfeld} holds for arbitrary $n\geq 2$.
\end{thm}
If we abuse the notation $\mathrm{D}_{\rm{dR}}(-)$ for the functor that sends a de Rham $(\varphi,\Gamma)$-module to its associated filtered $(\varphi,N)$-module, then Theorem~\ref{thm: intro Drinfeld} gives an isomorphism between filtered $(\varphi,N)$-modules of the form
\begin{equation}\label{equ: intro Drinfeld phi N}
\Hom_{\cM(G)}((\mathrm{St}_{G}^{\rm{an}})^\vee[1-n],\Omega^{\bullet})\cong \mathrm{D}_{\rm{dR}}(\cD_{n-1}\otimes_{\cR} \cR(\theta^{1-n})).
\end{equation}

We give an outline of the proof of Theorem~\ref{thm: intro Drinfeld}.
We write
\[\mathbf{E}_{k,k'}\defeq \mathbf{E}_{[1,n-1-k],[1,n-1-k']}=\mathrm{Ext}_{G}^{k'-k}(V_{[1,n-1-k]}^{\rm{an}},V_{[1,n-1-k']}^{\rm{an}})\]
for each $0\leq k\leq k'\leq n-1$, and write
\[\mathbf{E}_{\leq\ell}\defeq\bigoplus_{k=0}^{\ell}\mathbf{E}_{k,\ell}\]
for each $0\leq \ell\leq n-1$.
Concerning our choice of $(\varphi,N)$-action on (\ref{equ: intro main Hom}) (using the splitting $s$), we may identify (\ref{equ: intro main Hom}) with $\mathbf{E}_{\leq n-1}$ from now on, with $\varphi$ acting by $p^{k}$ on $\mathbf{E}_{k,n-1}$ and $N$ restricting to a map $\mathbf{E}_{k,n-1}\rightarrow\mathbf{E}_{k-1,n-1}$ for each $0\leq k\leq n-1$ (with $N(\mathbf{E}_{0,n-1})=0$).
We have a natural cup product map
\[\kappa_{k,k''}^{k'}\defeq\kappa_{[1,n-1-k],[1,n-1-k'']}^{[1,n-1-k']}: \mathbf{E}_{k,k'}\otimes_E\mathbf{E}_{k',k''}\buildrel\cup\over\longrightarrow\mathbf{E}_{k,k''}\]
for each $0\leq k\leq k'\leq k''\leq n-1$, and write $\kappa^{\ell}$ for the composition of the following maps
\[
\mathbf{E}_{\leq\ell}\otimes_{E}\mathbf{E}_{\ell,n-1}=\bigoplus_{k=0}^{\ell}\mathbf{E}_{k,\ell}\otimes_E\mathbf{E}_{\ell,n-1}\\
\buildrel \bigoplus_{k=0}^{\ell}\kappa_{k,n-1}^{\ell}\over\longrightarrow \bigoplus_{k=0}^{\ell}\mathbf{E}_{k,n-1}\hookrightarrow \bigoplus_{k=0}^{n-1}\mathbf{E}_{k,n-1}=\mathbf{E}_{\leq n-1}.
\]
For each $1\leq k\leq n-1$ the map $N:\mathbf{E}_{k,n-1}\rightarrow\mathbf{E}_{k-1,n-1}$ is an embedding given by
\[\kappa_{k-1,n-1}^{k}(u_{k}\otimes-): \mathbf{E}_{k,n-1}\rightarrow \mathbf{E}_{k-1,n-1}\]
for some $0\neq u_{k}\in\mathrm{Fil}_{1}(\mathbf{E}_{k-1,k})\subseteq \mathbf{E}_{k-1,k}$.
We write $x_{k,k'}\defeq s_{n-k}s_{n-1-k}\cdots s_{n-1-k'}$ for each $0\leq k\leq k'\leq n-1$.
The key ingredient in the proof of Theorem~\ref{thm: intro Drinfeld} is to show that the Hodge filtration $\mathrm{Fil}_{H}^{\bullet}(\mathbf{E}_{\leq n-1})$ on $\mathbf{E}_{\leq n-1}$ satisfies the following conditions.
\begin{cond}\label{cond: intro Hodge}
Let $\mathbf{E}_{\leq n-1}$ be the $(\varphi,N)$-module as described above. There exists a $1$-dimensional $E$-subspace 
\begin{equation}\label{equ: intro Hodge line}
L_{\ell}\subseteq \bigoplus_{k=0}^{\ell}\mathrm{Fil}_{x_{k,\ell}}(\mathbf{E}_{k,\ell})\subseteq \bigoplus_{k=0}^{\ell}\mathbf{E}_{k,\ell}=\mathbf{E}_{\leq\ell}
\end{equation}
for each $0\leq \ell\leq n-1$ such that
\begin{enumerate}[label=(\roman*)]
\item \label{it: intro Hodge 1} the composition of
\[L_{\ell}\subseteq \bigoplus_{k=0}^{\ell}\mathrm{Fil}_{x_{k,\ell}}(\mathbf{E}_{k,\ell})\twoheadrightarrow \mathrm{Fil}_{x_{j,\ell}}(\mathbf{E}_{j,\ell})\twoheadrightarrow \mathrm{gr}_{x_{j,\ell}}(\mathbf{E}_{j,\ell})\]
is an isomorphism for each $0\leq j\leq \ell$;
\item \label{it: intro Hodge 2} we have
\[\mathrm{Fil}_{H}^{m}(\mathbf{E}_{\leq n-1})=\sum_{\ell=m}^{n-1}\kappa^{\ell}(L_{\ell}\otimes_E\mathbf{E}_{\ell,n-1})\]
for each $0\leq m\leq n-1$.
\end{enumerate}
\end{cond}
Then we prove in Proposition~\ref{prop: Galois dR Ext exact} that any filtered $(\varphi,N)$-module $\mathbf{E}_{\leq n-1}$ satisfying Condition~\ref{cond: intro Hodge} is necessarily isomorphic to $\mathrm{D}_{\rm{dR}}(\cD_{n-1}\otimes_{\cR} \cR(\theta^{1-n}))$.
Here we emphasize that the $1$-dimensional $E$-subspace (\ref{equ: intro Hodge line}) is actually given by the $E$-span of the composition of the following maps (see (\ref{equ: Hodge position composition}))
\begin{equation}\label{equ: intro explicit splitting}
\mathrm{ker}(d_{\Omega^{\bullet}}^{\ell})[-\ell]\rightarrow \sigma_{\geq\ell}\Omega^{\bullet}\rightarrow \Omega^{\bullet}\buildrel s\over\longrightarrow \bigoplus_{k=0}^{n-1}H^{k}(\Omega^{\bullet})[-k]
\end{equation}
which certainly depends on the choice of the splitting $s$, but the isomorphism class of $\mathbf{E}_{\leq n-1}$ (as a filtered $(\varphi,N)$-module satisfying Condition~\ref{cond: intro Hodge}) does not depend on the particular choice of $L_{\ell}$ (for each $0\leq \ell\leq n-1$) or the choice of $s$.
\begin{rem}\label{rem: intro s match}
Suppose now that we take $s$ to be the Hyodo-Kato splitting $s_{\rm{HK}}$. Then upon fixing the choice of an isomorphism (\ref{equ: intro dR coh}) for each $0\leq k\leq n-1$, the element $u_{k}$ (in the description of $N$) can be explicitly determined for each $1\leq k\leq n-1$ (see \cite{dS05} and \cite{V24}). To obtain an explicit full description of the filtered $(\varphi,N)$-module, the main difficulty is thus to make explicit the $1$-dimensional $E$-subspace (\ref{equ: intro Hodge line}) (which is determined by $s=s_{\rm{HK}}$ via (\ref{equ: intro explicit splitting})) for each $0\leq\ell\leq n-1$.
\end{rem}

We continue with the study of special filtered $(\varphi,N)$-submodules of (\ref{equ: intro main Hom}).
When $D=V^{\vee}[1-n]$ for some subrepresentation $V\subseteq \mathrm{St}_{G}^{\rm{an}}$, we can show that (\ref{equ: intro Hom in derived}) embeds into (\ref{equ: intro main Hom}) as filtered $(\varphi,N)$-modules (using Lemma~\ref{lem: Ext subquotient}).
Among different choices of $V\subseteq \mathrm{St}_{G}^{\rm{an}}$, we are particularly interested in the choice $V=(\Omega^{n-1})^{\vee}$ which gives the following filtered $(\varphi,N)$-module
\begin{equation}\label{equ: intro Hom special sub}
\Hom_{\cM(G)}(\Omega^{n-1}[1-n],\Omega^{\bullet}).
\end{equation}
Here we use the fact that
\[\Hom_{D(G)}((\mathrm{St}_{G}^{\rm{an}})^{\vee}, \Omega^{n-1})=1\] 
from Lemma~\ref{lem: Hodge Ext transfer} (with $\ell=n-1$ in \emph{loc.cit.}).
In fact, among all choices of subrepresentations $V\subseteq \mathrm{St}_{G}^{\rm{an}}$ with $D=V^{\vee}[1-n]$ such that $\mathrm{Fil}^{n-1}(-)$  of (\ref{equ: intro Hom in derived}) is non-zero, the choice $V=(\Omega^{n-1})^{\vee}$ is minimal (namely has minimal possible length).
Under the map $H^{0}(\Omega^{\bullet})\rightarrow\Omega^{\bullet}$, the embedding from (\ref{equ: intro Hom special sub}) to (\ref{equ: intro main Hom}) restricts to an embedding (again using Lemma~\ref{lem: Ext subquotient})
\begin{multline}\label{equ: intro auto flat}
\mathbf{E}_{n-1}^{\flat}\defeq\mathrm{Ext}_{G}^{n-1}(1_{G},(\Omega^{n-1})^{\vee})=\Hom_{\cM(G)}(\Omega^{n-1}[1-n],H^0(\Omega^{\bullet}))\\
\hookrightarrow \Hom_{\cM(G)}((\mathrm{St}_{G}^{\rm{an}})^{\vee}[1-n],H^0(\Omega^{\bullet}))=\mathrm{Ext}_{G}^{n-1}(1_{G},\mathrm{St}_{G}^{\rm{an}})=\mathbf{E}_{n-1}.
\end{multline}
Now that the explicit structure of the coadmissible $D(G)$-module $\Omega^{n-1}$ has been well-understood (see \cite{ST02}, \cite{Or08} and \cite{BQ24}), we can actually identify $\mathbf{E}_{n-1}^{\flat}$ with the following step of Coxeter filtration
\[\mathrm{Fil}_{s_{n-1}\cdots s_1}(\mathbf{E}_{\Delta})\subseteq\mathbf{E}_{\Delta}=\mathbf{E}_{n-1}\]
from Theorem~\ref{thm: intro coxeter filtration}.

We define $\cE_{n,0}^{\flat}\defeq E$, $\cD_{n,0}^{\flat}\defeq \cR(1_{K^\times})$.
Then we inductively define $\cE_{n,m}^{\flat}$ and $\cD_{n,m}^{\flat}$ for increasing $1\leq m\leq n$ by 
\[\cE_{n,m}^{\flat}\defeq \mathrm{Ext}_g^1(\cD_{n,m-1}^{\flat},\cR(|\cdot|^{m}z^{n}))^\vee\]
and the following universal de Rham extension
\[0\rightarrow \cE_{n,m}^{\flat}\otimes_E\cR(|\cdot|^{m}z^{n})\rightarrow \cD_{n,m}^{\flat}\rightarrow \cD_{n,m-1}^{\flat}\rightarrow 0.\]
We write $\cD_{n}^{\flat}\defeq \cD_{n,n}^{\flat}$ and $\cE_{n}^{\flat}\defeq \cE_{n,n}^{\flat}$ for short.
It follows from Lemma~\ref{lem: Galois sub flat} that there exists a unique (up to scalars) non-zero map $\cD_{n-1}^{\flat}\rightarrow \cD_{n-1}$ which is an embedding. Note that this embedding restricts to an embedding $\cE_{n-1}^{\flat}\otimes_E\cR(\theta^{n-1})\hookrightarrow \cE_{n-1}\otimes_E\cR(\theta^{n-1})$ and thus induces an embedding 
\begin{equation}\label{equ: intro Galois flat}
\cE_{n-1}^{\flat}\hookrightarrow \cE_{n-1}
\end{equation}
between $E$-vector spaces.
The following result describes the relation between the embedding $\cD_{n-1}^{\flat}\hookrightarrow\cD_{n-1}$ and the embedding
\begin{equation}\label{equ: intro universal phi N embedding}
\Hom_{\cM(G)}(\Omega^{n-1}[1-n],\Omega^{\bullet})\hookrightarrow \Hom_{\cM(G)}((\mathrm{St}_{G}^{\rm{an}})^\vee[1-n],\Omega^{\bullet})
\end{equation}
induced from the unique (up to scalars) surjection $(\mathrm{St}_{G}^{\rm{an}})^{\vee}\twoheadrightarrow \Omega^{n-1}$.
\begin{thm}[Theorem~\ref{thm: main filtered phi N}]\label{thm: intro Drinfeld special sub}
We have the following results.
\begin{enumerate}[label=(\roman*)]
\item \label{it: intro Drinfeld sub 1} The de Rham $(\varphi,\Gamma)$-module associated with the filtered $(\varphi,N)$-module (\ref{equ: intro Hom special sub}) is isomorphic to $\cD_{n-1}^{\flat}\otimes_{\cR} \cR(\theta^{1-n})$.
\item \label{it: intro Drinfeld sub 2} Up to non-zero scalars, the embedding (\ref{equ: intro universal phi N embedding}) between filtered $(\varphi,N)$-modules induces the unique embedding $\cD_{n-1}^{\flat}\otimes_{\cR} \cR(\theta^{1-n})\hookrightarrow \cD_{n-1}\otimes_{\cR} \cR(\theta^{1-n})$ obtained from the embedding $\cD_{n-1}^{\flat}\hookrightarrow\cD_{n-1}$ by a twist.
\end{enumerate}
\end{thm}
\subsubsection{Breuil--Schraen $\mathscr{L}$-invariants and Fontaine-Mazur $\mathscr{L}$-invariants}\label{subsubsec: intro BS}
We apply our results from \S \ref{subsubsec: intro coh St} and \S \ref{subsubsec: intro cup} to define our Breuil--Schraen $\mathscr{L}$-invariants and then propose a candidate for their explicit relation to Fontaine-Mazur $\mathscr{L}$-invariants.

For each $J\subseteq\Delta$, we write $\mathbf{E}_{J}^{\infty}\defeq\mathrm{Fil}_{1}(\mathbf{E}_{J})$ for short and note that $\mathbf{E}_{J}^{\infty}$ is a $1$-dimensional $E$-subspace of $\mathbf{E}_{J}$ by Theorem~\ref{thm: intro coxeter filtration}.
For each $J,J'\subseteq\Delta$ with $J\cap J'=\emptyset$, it follows from (\ref{equ: intro sm cup}) (by taking $y=1$ in \emph{loc.cit.}) that
\[\mathbf{E}_{J\sqcup J'}^{\infty}=\kappa_{J,J'}(\mathbf{E}_{J}^{\infty}\otimes_E\mathbf{E}_{J'}^{\infty})\] 
as $1$-dimensional $E$-subspaces of $\mathbf{E}_{J\sqcup J'}$.
We write $\mathbf{E}_{n-1}^{\infty}\defeq\mathbf{E}_{\Delta}^{\infty}\subseteq\mathbf{E}_{\Delta}=\mathbf{E}_{n-1}$ for short.
For each $J\subseteq \Delta$, we write
\[\tld{\mathbf{E}}_{J}\defeq \kappa_{\Delta\setminus J,J}(\mathbf{E}_{\Delta\setminus J}^{\infty}\otimes_E\mathbf{E}_{J})\subseteq \mathbf{E}_{\Delta}=\mathbf{E}_{n-1}\]
for short.
Now that $\mathbf{E}_{\Delta\setminus J}^{\infty}$ is $1$-dimensional and $\kappa_{\Delta\setminus J,J}$ is injective (see \ref{it: intro cup graded 2} of Theorem~\ref{thm: intro cup graded}), we know that $\kappa_{\Delta\setminus J,J}(\mathbf{E}_{\Delta\setminus J}^{\infty}\otimes_E-)$ induces a bijection between $E$-subspaces of $\mathbf{E}_{J}$ and $E$-subspaces of $\tld{\mathbf{E}}_{J}$, under which $\mathbf{E}_{J}^{\infty}$ corresponds to $\mathbf{E}_{\Delta}^{\infty}=\mathbf{E}_{n-1}^{\infty}$.
Consequently, given $J\subseteq\Delta$ and a codimension $1$ $E$-subspace $W\subseteq \mathbf{E}_{n-1}$ with $W\cap \mathbf{E}_{n-1}^{\infty}=0$, there exists a unique $E$-subspace $W_{J}\subseteq \mathbf{E}_{J}$ with codimension $1$ such that $W_{J}\cap\mathbf{E}_{J}^{\infty}=0$ and
\begin{equation}\label{equ: intro W sm cup}
\kappa_{\Delta\setminus J,J}(\mathbf{E}_{\Delta\setminus J}^{\infty}\otimes W_{J})=W\cap \tld{\mathbf{E}}_{J}
\end{equation}
as $E$-subspaces of $\mathbf{E}_{\Delta}=\mathbf{E}_{n-1}$.

The following definition generalizes the definition of $\mathscr{L}$-invariants by Breuil (when $n=2$ with $K=\Q_p$) and Schraen (when either $n=2$ with $K$ general, or $n=3$ with $K=\Q_p$).
\begin{defn}[Definition~\ref{def: total inv}]\label{def: intro BS L inv}
A \emph{Breuil-Schraen $\mathscr{L}$-invariant} for $G$ is a codimension $1$ $E$-subspace $W\subseteq \mathbf{E}_{n-1}$ that satisfies the following conditions.
\begin{enumerate}[label=(\roman*)]
\item \label{it: intro BS L inv 1} We have $W\cap \mathbf{E}_{n-1}^{\infty}=0$ and thus $W_{J}$ has codimension $1$ in $\mathbf{E}_{J}$ for each $J\subseteq \Delta$.
\item \label{it: intro BS L inv 2} For each $J,J'\subseteq \Delta$ with $J\cap J'=\emptyset$, the map $\kappa_{J,J'}$ induces the following isomorphism between $1$-dimensional $E$-vector spaces
    \[\mathbf{E}_{J}/W_{J}\otimes_E\mathbf{E}_{J'}/W_{J'}\buildrel\sim\over\longrightarrow \mathbf{E}_{J\sqcup J'}/W_{J\sqcup J'}.\]
\end{enumerate}
\end{defn}

\begin{ex}\label{ex: intro example}
We compare Definition~\ref{def: intro BS L inv} with previous definitions of $\mathscr{L}$-invariants for $n=2,3$ due to Breuil and Schraen.
\begin{enumerate}[label=(\roman*)]
\item \label{it: BS inv n=2} Assume that $n=2$. Then we have an isomorphism
\[\mathbf{E}_1=\mathbf{E}_{\Delta}\cong \Hom(K^\times,E)\]
by \ref{it: intro cup generator 1} of Corollary~\ref{cor: intro cup generator}, under which the subspace $\mathbf{E}_1^{\infty}=\mathbf{E}_{\Delta}^{\infty}$ of $\mathbf{E}_1=\mathbf{E}_{\Delta}$ is given by $E\val\subseteq \Hom(K^\times,E)$. Hence, a Breuil-Schraen $\mathscr{L}$-invariant for $G=\mathrm{PGL}_2(K)$ corresponds to a $E$-line in $\Hom(K^\times,E)$ which is different from $E\val$, which has the form $E(\log-\mathscr{L}\val)$ for some $\mathscr{L}\in E$.
\item \label{it: BS inv n=3} Assume that $n=3$ and write $\Delta=\{1,2\}$ for short. By \ref{it: intro cup generator 1} of Corollary~\ref{cor: intro cup generator} we know that $\mathbf{E}_{\{j\}}$ admits a basis of the form $\{\plog_j,\val_j\}$ for each $j\in\Delta$.
    The map
    \[\kappa_{\{1\},\{2\}}: \mathbf{E}_{\{1\}}\otimes_E \mathbf{E}_{\{2\}}\buildrel\cup\over\longrightarrow \mathbf{E}_2=\mathbf{E}_{\Delta}\]
    is injective with image $\mathbf{E}_{\Delta}^{<}$ spanned by
    \[\{\log_1\cup\log_2, \log_1\cup\val_2, \val_1\cup\log_2, \val_1\cup\val_2\}\]
    which has codimension $1$ in $\mathbf{E}_{\Delta}$ by \ref{it: intro cup generator 2} of Corollary~\ref{cor: intro cup generator}.
    In particular, $\mathbf{E}_{\Delta}$ is $5$-dimensional. We can check from Definition~\ref{def: intro BS L inv} that a Breuil-Schraen $\mathscr{L}$-invariant for $G=\mathrm{PGL}_3(K)$ is given by an $E$-line $E\val_j\neq W_{\{j\}}\subseteq\mathbf{E}_{\{j\}}$ for each $j\in\Delta$ and a $4$-dimensional subspace $W=W_{\Delta}\subseteq \mathbf{E}_{\Delta}$ such that
    \begin{equation}\label{equ: choice of hyperplane}
    \sum_{j\in\Delta}\mathbf{E}_{\Delta\setminus\{j\}}\cup W_{\{j\}}\subseteq W\neq \mathbf{E}_{\Delta}^{<}
    \end{equation}
    with LHS of (\ref{equ: choice of hyperplane}) being $3$-dimensional.
    Note that given the choice of $W_{\{j\}}\subseteq\mathbf{E}_{\{j\}}$ for each $j\in\Delta$, the choice of a $W$ satisfying (\ref{equ: choice of hyperplane}) is equivalent to the choice of an $E$-line in the $2$-dimensional quotient
    \begin{equation}\label{equ: intro 2 dim quotient}
    \mathbf{E}_{\Delta}/\big(\sum_{j\in\Delta}\mathbf{E}_{\Delta\setminus\{j\}}\cup W_{\{j\}}\big)
    \end{equation}
    which is different from the $E$-line
    \begin{equation}\label{equ: intro 2 dim quotient line}
    \mathbf{E}_{\Delta}^{<}/\big(\sum_{j\in\Delta}\mathbf{E}_{\Delta\setminus\{j\}}\cup W_{\{j\}}\big).
    \end{equation}
    In fact, Schraen's original definition of a $\mathscr{L}$-invariant for $G=\mathrm{PGL}_3(\Q_p)$ (see \cite[\S 5.5]{Schr11}) can be translated to the choice of an $E$-line $W_{\{j\}}\subseteq\mathbf{E}_{\{j\}}$ that satisfies $W_{\{j\}}\neq E\val_j$ for each $j\in\Delta$, and then the choice of a third $E$-line in (\ref{equ: intro 2 dim quotient}) which is different from (\ref{equ: intro 2 dim quotient line}).
\item Under Theorem~\ref{thm: intro cup comparison}, $W_{\{j\}}\subseteq\mathbf{E}_{\{j\}}\cong\mathbf{E}_{\{j\},\emptyset}$ is essentially Ding's simple $\mathscr{L}$-invariants in \cite{Ding19}, which is known to satisfy local-global compatibility for suitable global setup.
\end{enumerate}
\end{ex}

Now we discuss the relation between Breuil-Schraen $\mathscr{L}$-invariants and Fontaine-Mazur $\mathscr{L}$-invariants, based on the filtered $(\varphi,N)$-module $\mathbf{E}_{\leq n-1}=\bigoplus_{k=0}^{n-1}\mathbf{E}_{k,n-1}$ introduced before Condition~\ref{cond: intro Hodge}. We fix throughout this discussion a choice of isomorphism (\ref{equ: intro dR coh}) for each $0\leq k\leq n-1$.
Note that the filtered $(\varphi,N)$-module $\mathbf{E}_{\leq n-1}$ depends on the choice of a splitting $s$ (see (\ref{equ: intro dR split})) and the tuple $\un{u}=(u_{k})_{1\leq k\leq n-1}$ which give the $N$-action. 
Given a Breuil-Schraen $\mathscr{L}$-invariant $W\subseteq\mathbf{E}_{n-1}=\mathbf{E}_{\Delta}$ with associated $\{W_{J}\}_{J\subseteq\Delta}$ as in Definition~\ref{def: intro BS L inv}, we write $W_{k,n-1}\defeq W_{[1,n-1-k]}$ for short and identify it with a codimension $1$ $E$-subspace of $\mathbf{E}_{k,n-1}$ via the isomorphism $\mathbf{E}_{k,n-1}=\mathbf{E}_{[1,n-1-k],\emptyset}\cong \mathbf{E}_{[1,n-1-k]}$ from \ref{it: intro cup comparison 1} of Theorem~\ref{thm: intro cup comparison}. We write $W_{\leq n-1}\defeq\bigoplus_{k=0}^{n-1}W_{k,n-1}$ and consider the quotient
\begin{equation}\label{equ: intro phi N quotient}
\mathbf{E}_{\leq n-1}/W_{\leq n-1}=\bigoplus_{k=0}^{n-1}\mathbf{E}_{k,n-1}/W_{k,n-1}
\end{equation}
which inherits the structure of a filtered $(\varphi,N)$-module from $\mathbf{E}_{\leq n-1}$. 
We write $\cD_{s,\un{u},W}$ for the de Rham $(\varphi,\Gamma)$-module over $\cR$ associated with (\ref{equ: intro phi N quotient}) via Berger's functor, which is necessarily a quotient of $\cD_{n-1}$. In other words, each choice of $s$ and $\un{u}$ determines a matching
\begin{equation}\label{equ: intro BS FM match}
W \leftrightsquigarrow \mathbf{E}_{\leq n-1}/W_{\leq n-1} \leftrightsquigarrow \cD_{s,\un{u},W}
\end{equation}
which is entirely determined by an isomorphism of the form
\begin{equation}\label{equ: intro abstract match}
r_{n-1,s,\un{u}}: \mathbf{E}_{n-1}\buildrel\sim\over\longrightarrow \cE_{n-1}.
\end{equation}
When $s=s_{\rm{HK}}$ is the Hyodo-Kato splitting and $\un{u}$ comes from the $N$-action on $\bigoplus_{k=0}^{n-1}H_{\rm{HK}}^{k}(\bH)[-k]$, we recover (\ref{equ: intro regulator map}) from (\ref{equ: intro abstract match}), and call (\ref{equ: intro BS FM match}) the \emph{canonical matching between Breuil-Schraen $\mathscr{L}$-invariants and Fontaine-Mazur $\mathscr{L}$-invariants} in this case.
To give an explicit candidate for this canonical matching, it is thus sufficient to explicitly characterize a map of the form (\ref{equ: intro regulator map}). We devote the rest of this section to such an attempt.

For each pair of integers $m,n\geq 0$, we construct a canonical embedding (see (\ref{equ: auto stable cup}))
\begin{equation}\label{equ: intro reg auto product}
\mathbf{E}_{m}\otimes_E\mathbf{E}_{n}\rightarrow \mathbf{E}_{m+n}
\end{equation}
as well as a canonical embedding (see (\ref{equ: Galois tensor sub Ext embedding}))
\begin{equation}\label{equ: intro reg Galois product}
\cE_{m}\otimes_E\cE_{n}\rightarrow \cE_{m+n}.
\end{equation}
For each $n\geq 2$, we obtain a canonical embedding
\begin{equation}\label{equ: intro reg auto sym}
\mathrm{Sym}^{n}(\mathbf{E}_{1})\hookrightarrow \mathbf{E}_{n},
\end{equation}
which allows us to construct (see Corollary~\ref{cor: sym sharp decomposition}) a canonical surjection
\begin{equation}\label{equ: intro reg auto surjection}
\mathbf{E}_{n}\twoheadrightarrow \mathrm{Sym}^{n}(\mathbf{E}_{1})
\end{equation}
with kernel $\mathbf{E}_{n}^{\sharp}$ and whose composition with (\ref{equ: intro reg auto sym}) is the identity map.
Similarly, for each $n\geq 2$, we have a canonical embedding (see (\ref{equ: Galois sym embedding}))
\begin{equation}\label{equ: intro reg Galois sym}
\mathrm{Sym}^{n}(\cE_{1})\hookrightarrow \cE_{n}
\end{equation}
as well as a canonical surjection
\begin{equation}\label{equ: intro reg Galois surjection}
\cE_{n}\twoheadrightarrow \mathrm{Sym}^{n}(\cE_{1}).
\end{equation}
with kernel $\cE_{n}^{\sharp}$ and whose composition with (\ref{equ: intro reg Galois sym}) is the identity map.

Recall from \cite[Lem.~1.19]{Ding17} that we have a canonical isomorphism
\begin{equation}\label{equ: intro r1}
r_{1}: \mathbf{E}_{1}=\Hom(K^\times,E)\cong \mathrm{Ext}_{g}^1(\cR(1_{K^\times}),\cR(\theta))^\vee=\cE_{1}
\end{equation}
which induces a canonical isomorphism
\begin{equation}\label{equ: intro r1 sym}
\mathrm{Sym}^{n}(r_{1}): \mathrm{Sym}^{n}(\mathbf{E}_{1})\buildrel\sim\over\longrightarrow\mathrm{Sym}^{n}(\cE_{1})
\end{equation}
for each $n\geq 2$.
Upon replacing $n-1$ with $n$ in (\ref{equ: intro auto flat}) and (\ref{equ: intro Galois flat}), we have a subspace $\mathbf{E}_{n}^{\flat}\subseteq\mathbf{E}_{n}$ as well as a subspace $\cE_{n}^{\flat}\subseteq\cE_{n}$ for each $n\geq 2$.
Under a technical hypothesis (which we strongly believe to be true), the following result presents a list of explicit conditions that uniquely characterize a collection of isomorphisms between $E$-vector spaces
\begin{equation}\label{equ: intro total reg}
\{r_{n}: \mathbf{E}_{n}\buildrel\sim\over\longrightarrow \cE_{n}\}_{n\geq 1}.
\end{equation}
\begin{thm}[Theorem~\ref{thm: regulator candidate}]\label{thm: intro regulator}
Assume that Hypothesis~\ref{hypo: log projection} holds for each $n\geq 2$.
There exists at most one collection of isomorphisms (\ref{equ: intro total reg}) such that $r_{1}$ is given by (\ref{equ: intro r1}) and the following conditions hold.
\begin{enumerate}[label=(\roman*)]
\item \label{it: intro regulator 1} For each $n\geq 2$, we have the following commutative diagram
\[
\xymatrix{
\mathbf{E}_{n} \ar^{r_{n}}_{\sim}[rr] \ar@{->>}[d] & & \cE_{n} \ar@{->>}[d]\\
\mathrm{Sym}^{n}(\mathbf{E}_{1}) \ar^{\mathrm{Sym}^{n}(r_{1})}_{\sim}[rr] & & \mathrm{Sym}^{n}(\cE_{1})
}
\]
with the vertical maps from (\ref{equ: intro reg auto surjection}) and (\ref{equ: intro reg Galois surjection}).
\item \label{it: intro regulator 2} For each $m,n\geq 1$, we have the following commutative diagram
\[\xymatrix{
\mathbf{E}_{m}\otimes_E\mathbf{E}_{n} \ar^{r_{m}\otimes_E r_{n}}_{\sim}[rr] \ar[d] & & \cE_{m}\otimes_E\cE_{n} \ar[d]\\
\mathbf{E}_{m+n} \ar^{r_{m+n}}_{\sim}[rr] & & \cE_{m+n}
}\]
with the vertical maps from (\ref{equ: intro reg auto product}) and (\ref{equ: intro reg Galois product}).
\item \label{it: intro regulator 3} For each $n\geq 2$, we have $r_{n}(\mathbf{E}_{n}^{\flat})=\cE_{n}^{\flat}$.
\end{enumerate}
\end{thm}
\begin{rem}\label{rem: intro sym map}
Along the proof of Theorem~\ref{thm: intro regulator}, we see that the composition of the following maps
\begin{equation}\label{equ: intro Galois sym map}
\cE_{n}^{\flat}\hookrightarrow \cE_{n}\twoheadrightarrow \mathrm{Sym}^{n}(\cE_{1})
\end{equation}
is an isomorphism between $n+1$-dimensional $E$-vector spaces. Upon applying $-\otimes_{E}\cR(\theta^{n})$, (\ref{equ: intro Galois sym map}) is induced from the following maps between certain universal $(\varphi,\Gamma)$-modules (see Lemma~\ref{lem: universal to sym}).
\[\cD_{n}^{\flat}\hookrightarrow \cD_{n}\twoheadrightarrow \mathrm{Sym}^{n}(\cD_{1})\]
Similarly, under Hypothesis~\ref{hypo: log projection}, the composition of the following maps
\begin{equation}\label{equ: intro auto sym map}
\mathbf{E}_{n}^{\flat}\hookrightarrow \mathbf{E}_{n}\twoheadrightarrow \mathrm{Sym}^{n}(\mathbf{E}_{1})
\end{equation}
is also an isomorphism between $n+1$-dimensional $E$-vector spaces.
One might wonder if the maps (\ref{equ: intro auto sym map}) have any geometric interpretation.
\end{rem}

Following Theorem~\ref{thm: intro Drinfeld}, Theorem~\ref{thm: intro Drinfeld special sub} and Theorem~\ref{thm: intro regulator}, we arrive at the following natural question.
\begin{ques}\label{ques: intro main question}
Does the map (\ref{equ: intro regulator map}) from Conjecture~\ref{conj: intro Drinfeld dR} (which is constructed using the Hyodo-Kato splitting $s_{\rm{HK}}$) satisfies all the conditions listed in Theorem~\ref{thm: intro regulator}?
\end{ques}
We leave this question for future investigation.
\subsection{Main technical ingredients}\label{subsec: intro technical}
In this section, we discuss several main technical ingredients of the proof of the theorems listed in \S~\ref{subsec: intro main results}.
\subsubsection{Tits double complex}\label{subsubsec: Tits double}
In this section, we introduce several examples of \emph{Tits double complex}, with possibly \emph{relative conditions} or highest weights indexed by  partial-Coxeter elements, which serve as the most important technical intermediate objects in the computation of $\mathbf{E}_{J}$ for each $J\subseteq\Delta$.

Let $\fh\subseteq\fg$ be a reductive $E$-Lie subalgebra.
For each object to be defined below, we omit $\fh$ from the notation when $\fh=0$.
For each $I\subseteq \Delta$, we consider the following \emph{relative abstract $L_{I}$-cochain complex}
\[C^{\bullet}(L_{I},\fl_{I}\cap\fh)\defeq \Hom_{D(L_{I})}(\bigotimes_{U(\fl_{I}\cap\fh)}^{\bullet+1}D(L_{I})\otimes_{U(\fl_{I}\cap\fh)}1_{L_{I}}^{\vee},1_{L_{I}}^{\vee})\]
and write $H^{\bullet}(L_{I},\fl_{I}\cap\fh)$ for its cohomology.
For each $I'\subseteq I\subseteq\Delta$, we have a natural restriction map
\begin{equation}\label{equ: intro Levi restriction}
R_{I,I',\fh}^{\bullet}: C^{\bullet}(L_{I},\fl_{I}\cap\fh)\rightarrow C^{\bullet}(L_{I'},\fl_{I'}\cap\fh)
\end{equation}
which induces a restriction map between cohomologies
\begin{equation}\label{equ: intro Levi restriction coh}
\mathrm{Res}_{I,I',\fh}^{\bullet}: H^{\bullet}(L_{I},\fl_{I}\cap\fh)\rightarrow H^{\bullet}(L_{I'},\fl_{I'}\cap\fh).
\end{equation}
For each $I\subseteq\Delta$, the inclusion $D(L_{I})\subseteq D(G)$ and the embedding $1_{L_{I}}^{\vee}\hookrightarrow (i_{I}^{\rm{an}})^{\vee}$ induce a quasi-isomorphism (cf.~(\ref{equ: PS parabolic Levi resolution}))
\begin{equation}\label{equ: intro PS to Levi}
\Hom_{D(G)}(\bigotimes_{U(\fl_{I}\cap\fh)}^{\bullet+1}D(G)\otimes_{U(\fl_{I}\cap\fh)}(i_{I}^{\rm{an}})^{\vee},1_{G}^{\vee})\rightarrow C^{\bullet}(L_{I},\fl_{I}\cap\fh)
\end{equation}
with the cohomology of LHS being the following relative $\mathrm{Ext}$-groups (cf.~\cite[\S~8.7.1]{Wei94})
\[\mathrm{Ext}_{G,\fl_{I}\cap\fh}^{\bullet}(1_{G},i_{I}^{\rm{an}})\defeq \mathrm{Ext}_{D(G),U(\fl_{I}\cap\fh)}^{\bullet}((i_{I}^{\rm{an}})^{\vee},1_{G}^{\vee}).\]
Note that (\ref{equ: intro PS to Levi}) is functorial with respect to the choice of $I$ (namely with respect to the natural embedding $i_{I}^{\rm{an}}\hookrightarrow i_{I'}^{\rm{an}}$ for each $I'\subseteq I$) as well as the choice of $\fh$. 

Let $J\subseteq\Delta$.
We define what we call \emph{Tits double complex} $\cT_{J,\fh}^{\bullet,\bullet}$ as the double complex given by
\[\cT_{J,\fh}^{-\ell,\bullet}=\bigoplus_{I\supseteq\Delta\setminus J, \#I=\ell}C^{\bullet}(L_{I},\fl_{I}\cap\fh)\]
for each $\ell$, with the differential map $\cT_{J,\fh}^{-\ell,\bullet}\rightarrow \cT_{J,\fh}^{-\ell+1,\bullet}$ given by the direct sum of 
\[(-1)^{\#\{i\in I\mid i<j\}}R_{I,I\setminus\{j\},\fh}^{\bullet}: C^{\bullet}(L_{I},\fl_{I}\cap\fh)\rightarrow C^{\bullet}(L_{I\setminus\{j\}},\fl_{I\setminus\{j\}}\cap\fh)\]
over all $I\supseteq\Delta\setminus J$ and $j\in I$ with $\#I=\ell$ (see (\ref{equ: intro Levi restriction}) for $R_{I,I\setminus\{j\},\fh}^{\bullet}$).
Truncations of $\cT_{J,\fh}^{\bullet,\bullet}$ with respect to the first index give a spectral sequence $E_{\bullet,J,\fh}^{\bullet,\bullet}$ that converges to $H^{\bullet}(\mathrm{Tot}(\cT_{J,\fh}^{\bullet,\bullet}))$.
It is clear that the natural map $C^{\bullet}(L_{I},\fl_{I}\cap\fh)\rightarrow C^{\bullet}(L_{I})$ for each $I\supseteq \Delta\setminus J$ induces a map between double complex
\begin{equation}\label{equ: intro Tits double drop h}
\cT_{J,\fh}^{\bullet,\bullet}\rightarrow \cT_{J}^{\bullet,\bullet}.
\end{equation}
Upon replacing $C^{\bullet}(L_{I},\fl_{I}\cap\fh)$ with LHS of (\ref{equ: intro PS to Levi}) in the definition of $\cT_{J,\fh}^{\bullet,\bullet}$, we obtain another double complex $\cS_{J,\fh}^{\bullet,\bullet}$ which comes with a natural map
\[\cS_{J,\fh}^{\bullet,\bullet}\rightarrow \cT_{J,\fh}^{\bullet,\bullet}\]
whose associated map between spectral sequences is an isomorphism from the first page (using the fact that (\ref{equ: intro PS to Levi}) is a quasi-isomorphism for each $I\supseteq \Delta\setminus J$).
When $\fh=0$ (with $\fh$ omitted from the notation), the quasi-isomorphism (\ref{equ: intro Tits resolution}) together with the quasi-isomorphism (\ref{equ: intro PS to Levi}) for each $\Delta\setminus J\subseteq I\subseteq \Delta$ induces an isomorphism
\begin{equation}\label{equ: intro St total coh}
\mathbf{E}_{J}=\mathrm{Ext}_{G}^{2\#J-\#\Delta}(1_{G},\mathbf{C}_{J})=H^{2\#J-\#\Delta}(\mathrm{Tot}(\cS_{J}^{\bullet,\bullet}))\buildrel\sim\over\longrightarrow H^{2\#J-\#\Delta}(\mathrm{Tot}(\cT_{J}^{\bullet,\bullet})).
\end{equation} 

Now we generalize (\ref{equ: intro St total coh}) to a similar description of $\mathbf{E}_{x,J}$ for each $x\in\Gamma^{J}$.
For each $I\subseteq\Delta$, we write $D(\fg,P_{I})$ for the associative $E$-subalgebra of $D(G)$ generated by $U(\fg)$ and $D(P_{I})$.
We temporarily fix a choice of $x\in\Gamma^{J}$. 
For each $\Delta\setminus J\subseteq I\subseteq J_{x}=\Delta\setminus\mathrm{Supp}(x)$, the inclusion $D(\fg,P_{I})\subseteq D(G)$ and the embedding $M^{I}(x)\hookrightarrow (i_{x,I}^{\rm{an}})^{\vee}$ induce a quasi-isomorphism
\begin{multline}\label{equ: intro PS to parabolic x}
\Hom_{D(G)}(\bigotimes_{U(\fl_{I}\cap\fh)}^{\bullet+1}D(G)\otimes_{U(\fl_{I}\cap\fh)}(i_{x,I}^{\rm{an}})^{\vee},1_{G}^{\vee})\\
\rightarrow \Hom_{D(\fg,P_{I})}(\bigotimes_{U(\fl_{I}\cap\fh)}^{\bullet+1}D(\fg,P_{I})\otimes_{U(\fl_{I}\cap\fh)}M^I(x),1_{D(\fg,P_{I})})
\end{multline}
which is functorial with respect to the choice of $I$.
Hence, upon replacing $C^{\bullet}(L_{I},\fl_{I}\cap\fh)$ in the definition of $\cT_{J,\fh}^{\bullet,\bullet}$ with each term of the quasi-isomorphism (\ref{equ: intro PS to parabolic x}) for varying $\Delta\setminus J\subseteq I\subseteq J_{x}$, we can similarly define a map between double complex
\[
\cS_{x,J,\fh}^{\bullet,\bullet}\rightarrow \cS_{x,J,\fh,\flat}^{\bullet,\bullet}
\]
whose associated map between spectral sequences
\[
E_{\bullet,x,J,\fh}^{\bullet,\bullet}\rightarrow E_{\bullet,x,J,\fh,\flat}^{\bullet,\bullet}
\]
is an isomorphism from the first page. In particular, we obtain the following isomorphisms when $\fh=0$
\begin{equation}\label{equ: intro St total coh x}
\mathbf{E}_{x,J}=\mathrm{Ext}_{G}^{2\#J-\#\Delta}(1_{G},\mathbf{C}_{x,J})=H^{2\#J-\#\Delta}(\mathrm{Tot}(\cS_{x,J}^{\bullet,\bullet}))\buildrel\sim\over\longrightarrow H^{2\#J-\#\Delta}(\mathrm{Tot}(\cS_{x,J,\flat}^{\bullet,\bullet})).
\end{equation} 
When $x=1$ and $\fh=0$, the RHS of (\ref{equ: intro PS to parabolic x}) admits a natural quasi-isomorphism to $C^{\bullet}(L_{I})$, from which we obtain a map between double complex $\cS_{1,J,\flat}^{\bullet,\bullet}\rightarrow \cT_{J}^{\bullet,\bullet}$ whose associated map between spectral sequences is an isomorphism from the first page. In other words, (\ref{equ: intro St total coh}) factors through (\ref{equ: intro St total coh x}) with $x=1$.
\subsubsection{Atomic basis and $E_{2}$-degeneracy}\label{subsubsec: intro atom E2}
In this section, we explain the main technical ingredients in our study of the spectral sequence $E_{\bullet,J}^{\bullet,\bullet}$ and give an outline of the proof of Theorem~\ref{thm: intro main dim} when $J=\Delta$.

Assume in the rest of this section that $J=\Delta$ and write $r\defeq\#\Delta=n-1$ for short.
Note that the spectral sequence $E_{\bullet,\Delta}^{\bullet,\bullet}$ equips $\mathbf{E}_{\Delta}$ with a canonical decreasing filtration $\mathrm{Fil}^{\bullet}(\mathbf{E}_{\Delta})$ with graded piece
\[\mathrm{gr}^{-\ell}(\mathbf{E}_{\Delta})=E_{\infty,\Delta}^{-\ell,\ell+r}\]
for each $0\leq \ell\leq r$.
Following the discussion below Theorem~\ref{thm: intro main dim}, we need to construct an explicit basis
\begin{equation}\label{equ: intro comb basis}
\{\overline{x}_{S,I}\mid S\in\cS_{\Delta}, \#S=r-\ell, I\subseteq S\cap\Delta\}
\end{equation}
of $E_{2,\Delta}^{-\ell,\ell+r}$ and then prove
\begin{equation}\label{equ: intro full E2 degenerate}
E_{2,\Delta}^{-\ell,\ell+r}=E_{\infty,\Delta}^{-\ell,\ell+r}
\end{equation}
for each $0\leq \ell\leq r$.
Recall from the definition of $E_{\bullet,\Delta}^{\bullet,\bullet}$ that we have
\[E_{1,\Delta}^{-\ell,k}=\bigoplus_{I, \#I=\ell}H^{k}(L_{I})\]
for each $(-\ell,k)$, with $d_{1,\Delta}^{-\ell,k}$ given by the direct sum of 
\[(-1)^{\#\{i\in I\mid i<j\}}\mathrm{Res}_{I,I\setminus\{j\}}^{k}: H^{k}(L_{I})\rightarrow H^{k}(L_{I\setminus\{j\}})\]
over all $I$ and $j\in I$ with $\#I=\ell$.

We start with an explicit description of the restriction maps $\mathrm{Res}_{I,I'}^{\bullet}$ for each $I'\subseteq I$.
For each $I\subseteq\Delta$, we consider the normal open subgroup $L_{I}'\defeq\bigcap_{v}\mathrm{ker}(v)\subseteq L_{I}$ with $v$ running through an arbitrary basis of the $E$-vector space $\Hom_{\infty}(L_{I},E)$ of smooth additive characters of $L_{I}$.
Since $L_{I}'$ is open in $L_{I}$, we can clearly identify $\fl_{I}$ with the $E$-Lie algebra associated with $L_{I}'$.
Following the discussion around (\ref{equ: center isogeny decomposition}) and around (\ref{equ: Levi center dagger projection}), we can choose a discrete cocompact subgroup $T^{\dagger}\subseteq T$ as well as an isomorphism $\prod_{j\in\Delta}\Z_{j}\buildrel\sim\over\longrightarrow T^{\dagger}$ (with $\Z_{j}$ being a copy of $\Z$ for each $j\in\Delta$), such that for each $I\subseteq\Delta$ the composition of
\[\prod_{j\in\Delta}\Z_{j}\buildrel\sim\over\longrightarrow T^{\dagger}\hookrightarrow L_{I}\twoheadrightarrow L_{I}/L_{I}'\]
factors through
\begin{equation}\label{equ: intro T dagger projection}
\prod_{j\in\Delta}\Z_{j}\twoheadrightarrow \prod_{j\in\Delta\setminus I}\Z_{j}\buildrel\sim\over\longrightarrow Z_{I}^{\dagger}\rightarrow L_{I}/L_{I}'
\end{equation}
for some cocompact discrete subgroup $Z_{I}^{\dagger}$ of the center $Z_{I}$ of $L_{I}$, with the last map of (\ref{equ: intro T dagger projection}) being an isogeny (see Definition~\ref{def: gp isogeny}).
For each $I'\subseteq I$, we note that the inclusion $L_{I'}\subseteq L_{I}$ restricts to an inclusion
\begin{equation}\label{equ: intro Levi product}
L_{I'}'\times Z_{I'}^{\dagger}\subseteq L_{I}'\times Z_{I}^{\dagger}.
\end{equation}
For each $I\subseteq\Delta$, we write $H^{\bullet}(\fl_{I})$ for the Lie algebra cohomology (of the trivial $U(\fl_{I})$-module) and $\Hom(Z_{I}^{\dagger},E)$ for the $E$-vector space of (locally constant) additive characters of $Z_{I}^{\dagger}$.
The natural map $L_{I}'\times Z_{I}^{\dagger}\rightarrow L_{I}$ is an isogeny, which together with the inclusion $U(\fl_{I})\subseteq D(L_{I}')$ induces (see Proposition~\ref{prop: Levi decomposition}) an isomorphism
\begin{equation}\label{equ: intro Levi coh decomposition}
H^{\bullet}(L_{I})\cong\mathrm{Tot}(H^{\bullet}(\fl_{I})\otimes_E\wedge^{\bullet}\Hom(Z_{I}^{\dagger},E))
\end{equation}
for each $I\subseteq \Delta$. It follows from Koszul's theorem (cf.~Theorem~\ref{thm: Koszul thm} following \cite[Thm.~10.2, Thm.~10.3]{Kos50}) that $H^{\bullet}(\fl_{I})$ is a free wedge algebra generated by its $E$-subspace $P(\fl_{I})$ of \emph{primitive elements}, and $P(\fl_{I})$ has the form $\Hom(\fz_{I},E)\oplus\bigoplus_{\fa}P(\fa)$ with $\fz_{I}$ being the $E$-Lie algebra associated with $Z_{I}$ and $\fa$ running through all simple factors of $\fh_{I}\defeq[\fl_{I},\fl_{I}]$. Note that each $\fa$ is isomorphic to $\fs\fl_{m_{\fa}}$ for some $m_{\fa}\geq 2$, with
\[P(\fa)=\bigoplus_{m'=2}^{m_{\fa}}P^{2m'-1}(\fa)\]
and $P^{2m'-1}(\fa)=P(\fa)\cap H^{2m'-1}(\fa)$ being a $1$-dimensional $E$-vector space.
For each pair $I'\subseteq I$, we prove in Proposition~\ref{prop: Levi restriction} that we have the following commutative diagram
\begin{equation}\label{equ: intro Levi coh diagram}
\xymatrix{
H^{\bullet}(L_{I}) \ar^{\sim}[rrr] \ar^{\mathrm{Res}_{I,I'}^{\bullet}}[d] & & & \mathrm{Tot}(H^{\bullet}(\fl_{I})\otimes_E\wedge^{\bullet}\Hom(Z_{I}^{\dagger},E)) \ar[d]\\
H^{\bullet}(L_{I'}) \ar^{\sim}[rrr] & & & \mathrm{Tot}(H^{\bullet}(\fl_{I'})\otimes_E\wedge^{\bullet}\Hom(Z_{I'}^{\dagger},E))
}
\end{equation}
with both of the horizontal isomorphisms from (\ref{equ: intro Levi coh decomposition}), and the RHS vertical map being induced from the natural restriction maps $H^{\bullet}(\fl_{I})\rightarrow H^{\bullet}(\fl_{I'})$ (which is induced from a map $P(\fl_{I})\rightarrow P(\fl_{I'})$ by taking wedge product, thanks to Koszul's theorem as recalled in Theorem~\ref{thm: Koszul thm}) and the map $\wedge^{\bullet}\Hom(Z_{I}^{\dagger},E)\rightarrow \wedge^{\bullet}\Hom(Z_{I'}^{\dagger},E)$ (associated with the map $Z_{I'}^{\dagger}\rightarrow Z_{I}^{\dagger}$ induced from (\ref{equ: intro Levi product})). 
We emphasize that the map (\ref{equ: intro Levi product}) is certainly not the product of a map $L_{I'}'\rightarrow L_{I}'$ with another map $Z_{I'}^{\dagger}\rightarrow Z_{I}^{\dagger}$, and thus the description of $\mathrm{Res}_{I,I'}^{\bullet}$ as a product form in (\ref{equ: intro Levi coh diagram}) is not at all formal. The key idea in the proof of (\ref{equ: intro Levi coh diagram}) is to control $\mathrm{Res}_{I,I'}^{\bullet}$ using the following commutative diagram
\begin{equation}\label{equ: intro Levi coh relative diagram}
\xymatrix{
H^{\bullet}(L_{I},\fh) \ar[rr] \ar^{\mathrm{Res}_{I,I',\fh}^{\bullet}}[d] & & H^{\bullet}(L_{I}) \ar^{\mathrm{Res}_{I,I'}^{\bullet}}[d]\\
H^{\bullet}(L_{I'},\fl_{I'}\cap\fh) \ar[rr] & & H^{\bullet}(L_{I'})
}
\end{equation}
for varying reductive $E$-Lie subalgebras $\fh\subseteq\fl_{I}$. 
In fact, the usage of relative conditions with respect to varying $\fh$ allows us to cut out generators of $H^{\bullet}(L_{I})$ and to characterize the image of these generators under the restriction map $\mathrm{Res}_{I,I'}^{\bullet}$.

Such explicit description of $\mathrm{Res}_{I,I'}^{\bullet}$ for each pair $I'\subseteq I$ leads to an explicit and combinatorial description of the complex 
\[(E_{1,\Delta}^{\bullet,k}, d_{1,\Delta}^{\bullet,k})\]
for each $k\geq 0$. Now we construct the basis (\ref{equ: intro comb basis}) of
\[E_{2,\Delta}^{-\ell,\ell+r}=H^{-\ell}(E_{1,\Delta}^{\bullet,\ell+r})\]
for each $0\leq \ell\leq r$. We fix temporarily a choice of $S\in\cS_{\Delta}$ with $\ell\defeq r-\#S$ and $I\subseteq S\cap\Delta$. 
Note that the association
\[S\mapsto I_{S}\defeq \Delta\setminus\{i_{\al}\mid \al=(i_{\al},i_{\al}')\in S\}\]
actually induces a bijection between $\cS_{\Delta}$ and the power set of $[2,n-1]$, with $\#I_{S}=r-\#S$. The pair $I\subseteq S\cap\Delta$ of subsets of $\Delta$ determines a $\#S\cap\Delta$-dimensional $E$-subspace 
\[\Hom_{S,I}\subseteq \Hom(Z_{\Delta\setminus S},E)\hookrightarrow \Hom(Z_{I_{S}},E)\cong\Hom(\fz_{I_{S}},E)\oplus \Hom(Z_{I_{S}}^{\dagger},E),\] 
from which we can define the following $1$-dimensional $E$-subspace
\[P_{S,I}\defeq \wedge^{\#S\cap\Delta}\Hom_{S,I}\otimes_E\bigotimes_{\fa}P^{2m_{\fa}-1}(\fa)\subseteq H^{\ell+r}(L_{I_{S}})\subseteq E_{1,\Delta}^{-\ell,\ell+r}\]
where $\fa$ runs through all simple factors of $\fh_{I_{S}}$. 
We check (cf.~Lemma~\ref{lem: low deg}) that $d_{1,\Delta}^{-\ell,\ell+r}(P_{S,I})=0$ and thus can define $\overline{P}_{S,I}$ as the image of $P_{S,I}$ under $\mathrm{ker}(d_{1,\Delta}^{-\ell,\ell+r})\twoheadrightarrow E_{2,\Delta}^{-\ell,\ell+r}$. Then we can reformulate Corollary~\ref{cor: bottom deg E2 basis} as follows.
\begin{prop}[Corollary~\ref{cor: bottom deg E2 basis}]\label{prop: intro explicit E2 basis}
We have $\overline{P}_{S,I}\neq 0$ for each $S\in\cS_{\Delta}$ and $I\subseteq S\cap\Delta$. Moreover, we have
\begin{equation}\label{equ: intro explicit E2 basis}
E_{2,\Delta}^{-\ell,\ell+r}\cong \bigoplus_{S\in\cS_{\Delta}, \#S=r-\ell, I\subseteq S\cap\Delta}\overline{P}_{S,I}
\end{equation}
for each $0\leq \ell\leq r$.
\end{prop}
It is thus clear that we can construct the basis (\ref{equ: intro comb basis}) by choosing some non-zero element $\overline{x}_{S,I}\in \overline{P}_{S,I}$ for each pair $(S,I)$.

Now we move on to explain the idea of the proof of (\ref{equ: intro full E2 degenerate}). By similar combinatorial argument used in the proof of Proposition~\ref{prop: intro explicit E2 basis}, we prove (see Lemma~\ref{lem: existence of atom}) that $E_{2,\Delta}^{-\ell',k'}=0$ and thus $d_{s,\Delta}^{-\ell',k'}=0$ for each $s\geq 2$ whenever $k'-\ell'<r$. Hence, to prove (\ref{equ: intro full E2 degenerate}), it suffices to show that
\begin{equation}\label{equ: higher differential vanishing}
d_{s,\Delta}^{-\ell,\ell+r}(\overline{P}_{S,I})=0
\end{equation}
for each $s\geq 2$ and each $S\in\cS_{\Delta}$ as well as $I\subseteq S\cap\Delta$ with $\#S=r-\ell$.
The key idea in the proof of (\ref{equ: higher differential vanishing}) is the construction of a certain double complex $\cT_{\Delta,S,I}^{\bullet,\bullet}$ that comes with a map to $\cT_{\Delta}^{\bullet,\bullet}$, whose associated spectral sequence $E_{\bullet,\Delta,S,I}^{\bullet,\bullet}$ and its map to $E_{\bullet,\Delta}^{\bullet,\bullet}$ satisfy the following conditions.
\begin{itemize}
\item We have $E_{1,\Delta,S,I}^{-\ell',k'}=0$ whenever $k'-\ell'>r$.
\item The image of the map $E_{1,\Delta,S,I}^{-\ell,\ell+r}\rightarrow E_{1,\Delta}^{-\ell,\ell+r}$ is exactly $P_{S,I}$.
\end{itemize}
Note that the first condition ensures $d_{s,\Delta,S,I}^{-\ell,\ell+r}=0$ for each $s\geq 1$. This together with the second condition gives
\[d_{s,\Delta}^{-\ell,\ell+r}(P_{S,I})=0\]
for each $s\geq 1$, which implies (\ref{equ: higher differential vanishing}) for each $s\geq 2$.
When $S\cap\Delta=\emptyset$ (and in particular $I=\emptyset$), we write
\[I_{S}^{\dagger}\defeq I_{S}\setminus\{i\in I_{S}\mid i+1\notin I_{S}\}\]
and consider the $E$-Lie subalgebra $\fh_{I_{S}^{\dagger}}\subseteq\fh_{I_{S}}\cong\fl_{I_{S}}/\fz_{I_{S}}$.
Then the double complex $\cT_{\Delta,S}^{\bullet,\bullet}\defeq \cT_{\Delta,S,\emptyset}^{\bullet,\bullet}$ is given by
\[\cT_{\Delta,S}^{-\ell,\bullet}=\bigoplus_{I'\subseteq I_{S}, \#I'=\ell}C^{\bullet}(L_{I'}/Z_{I_{S}},(\fl_{I'}/\fz_{I_{S}})\cap\fh_{I_{S}^{\dagger}})\]
for each $\ell$, with the differential map $\cT_{\Delta,S}^{-\ell,\bullet}\rightarrow \cT_{\Delta,S}^{-\ell+1,\bullet}$ given by the direct sum (over all $I$ and $j\in I$ with $\#I=\ell$) of $(-1)^{\#\{i\in I'\mid i<j\}}$ times the following restriction map
\[C^{\bullet}(L_{I'}/Z_{I_{S}},(\fl_{I'}/\fz_{I_{S}})\cap\fh_{I_{S}^{\dagger}})\rightarrow C^{\bullet}(L_{I'\setminus\{j\}}/Z_{I_{S}},(\fl_{I'\setminus\{j\}}/\fz_{I_{S}})\cap\fh_{I_{S}^{\dagger}}).\]
In other words, Tits double complex with relative conditions is crucial for the proof of the degeneracy (\ref{equ: intro full E2 degenerate}).
\subsubsection{Coxeter filtration}\label{subsubsec: intro coxeter}
In this section, we fix a choice of $J\subseteq\Delta$ and give an outline of the proof of Proposition~\ref{prop: intro coxeter map} and Theorem~\ref{thm: intro coxeter filtration}.

Let $x,w\in\Gamma^{J}$ with $x\unlhd w$. We know from \cite[\S~4.2, \S~4.6]{Hum08} that there exists a unique up to scalars embedding $M(w)\hookrightarrow M(x)$. We prove that this embedding induces a non-zero map $M^{I}(w)\rightarrow M^{I}(x)$ for each $I\subseteq J_{w}\subseteq J_{x}$ (see the discussion around (\ref{equ: Lie parabolic Verma wt shift})), whose composition with $M^{I}(x)\twoheadrightarrow M^{I'}(x)$ is zero for any $I'\supseteq I$ that satisfies $I'\not\subseteq J_{w}$. Combining these facts with the Orlik-Strauch functors from \cite{OS15}, we obtain a non-zero map $i_{x,I}^{\rm{an}}\rightarrow i_{w,I}^{\rm{an}}$ for each $I\subseteq J_{w}\subseteq J_{x}$, whose composition with $i_{x,I'}^{\rm{an}}\rightarrow i_{x,I}^{\rm{an}}$ is zero for any $I'\supseteq I$ that satisfies $I'\not\subseteq J_{w}$. Consequently, we obtain a map between complex $\mathbf{C}_{x,I}\rightarrow \mathbf{C}_{w,I}$ and thus a map $V_{x,I}^{\rm{an}}\rightarrow V_{w,I}^{\rm{an}}$ which is non-zero by the discussion around (\ref{equ: St x w transfer J}).

Now we move on to the discussion on the proof of \ref{it: intro coxeter map 2} of Proposition~\ref{prop: intro coxeter map} and Theorem~\ref{thm: intro coxeter filtration}. We first clarify (formally) the input we need to construct the Coxeter filtration (\ref{equ: intro coxeter filtration})
on $\mathbf{E}_{J}$. 
For each finite length object $V\in\mathrm{Rep}^{\rm{an}}_{\rm{adm}}(G)$ and $k_0\geq 0$, we write
\[\mathrm{JH}_{G}^{k_0}(V)\defeq \{W\in\mathrm{JH}_{G}(V)\mid \mathrm{Ext}_{G}^{k_0}(1_{G},W)\neq 0\}\subseteq \mathrm{JH}_{G}(V).\]
We consider the following conditions on $V$ and $k_0$.
\begin{cond}\label{cond: intro abstract filtration}
Let $V\in\mathrm{Rep}^{\rm{an}}_{\rm{adm}}(G)$ be a finite length object and $k_0\geq 0$.
\begin{enumerate}[label=(\roman*)]
\item \label{it: intro abstract filtration 0} We have $\mathrm{Ext}_{G}^{k}(1_{G},W)=0$ for each $W\in\mathrm{JH}_{G}(V)$ and $k<k_0$.
\item \label{it: intro abstract filtration 1} Each $W\in \mathrm{JH}_{G}^{k_0}(V)$ appears with multiplicity one in $V$ and satisfies $\Dim_E\mathrm{Ext}_{G}^{k_0}(1_{G},W)=1$.
\item \label{it: intro abstract filtration 2} We have 
\begin{equation}\label{equ: intro Ext full dim}
\Dim_E\mathrm{Ext}_{G}^{k_0}(1_{G},V)=\#\mathrm{JH}_{G}^{k_0}(V).
\end{equation}
\end{enumerate}
\end{cond}
Assume for now that $V$ and $k_0$ satisfy Condition~\ref{cond: intro abstract filtration}. We write $\mathbf{E}\defeq \mathrm{Ext}_{G}^{k_0}(1_{G},V)$ for short.
Note that $\mathrm{JH}_{G}^{k_0}(V)$ is equipped with a natural partial-order given by $W'<W$ if and only if $W'\in\mathrm{JH}_{G}(\tld{W})$ with $\tld{W}$ being the unique subrepresentation of $V$ with cosocle $W$ (crucially using the fact that $W,W'\in \mathrm{JH}_{G}^{k_0}(V)$ are of multiplicity one in $V$).
Assuming Condition~\ref{cond: intro abstract filtration}, we can argue by a formal d\'evissage argument that each embedding $V'\hookrightarrow V$ induces an embedding
\[\mathrm{Ext}_{G}^{k_0}(1_{G},V')\hookrightarrow \mathrm{Ext}_{G}^{k_0}(1_{G},V)=\mathbf{E}\]
with $\Dim_E\mathrm{Ext}_{G}^{k_0}(1_{G},V')=\#\mathrm{JH}_{G}^{k_0}(V')$. Consequently, we obtain an increasing filtration indexed by the partially-ordered set $\mathrm{JH}_{G}^{k_0}(V)$
\[\{\mathrm{Fil}_{W}(\mathbf{E})\defeq \mathrm{Ext}_{G}^{k_0}(1_{G},\tld{W})\}_{W\in \mathrm{JH}_{G}^{k_0}(V)}\]
with the graded piece
\[\mathrm{gr}_{W}(\mathbf{E})\defeq \mathrm{Fil}_{W}(\mathbf{E})/\sum_{W'<W}\mathrm{Fil}_{W'}(\mathbf{E})\]
being $1$-dimensional for each $W\in \mathrm{JH}_{G}^{k_0}(V)$.
Now we consider $k_0\geq 0$ and a map $\varphi: V\rightarrow V'$ between finite length objects which together give a map
\[\varphi_{\ast}: \mathbf{E}\rightarrow \mathbf{E}'\defeq \mathrm{Ext}_{G}^{k_0}(1_{G},V').\]
Suppose now that both $(V,k_0)$ and $(V',k_0)$ satisfy Condition~\ref{cond: intro abstract filtration} and that $\mathrm{JH}_{G}^{k_0}(\mathrm{coker}(\varphi))=\emptyset$. Then we claim that $\varphi$ induces a strict map between filtered $E$-vector spaces $\mathbf{E}\twoheadrightarrow \mathbf{E}'$, with $\varphi_{\ast}(\mathrm{Fil}_{W}(\mathbf{E}))$ being non-zero if and only if $W\in \mathrm{JH}_{G}^{k_0}(V')\subseteq \mathrm{JH}_{G}^{k_0}(V)$, in which case it equals $\mathrm{Fil}_{W}(\mathbf{E}')$.

Now we consider $J\subseteq\Delta$. For each $x\in\Gamma^{J}$, it follows from Lemma~\ref{lem: x Ext with sm} that $V=V_{x,\Delta\setminus J}^{\rm{an}}$ and $k_0=\#J$ satisfy Condition~\ref{cond: intro abstract filtration} and we have a natural bijection
\[\{u\in\Gamma^{J}\mid x\unlhd u\}\buildrel\sim\over\longrightarrow \mathrm{JH}_{G}^{\#J}(V_{x,\Delta\setminus J}^{\rm{an}}): u\mapsto C_{u}^{\Delta\setminus J}\]
which respects partial-orders on BHS and is functorial with respect to the choice of $x$.
Hence, given $x,w\in\Gamma^{J}$ with $x\unlhd w$, we can apply the above discussion on $\varphi: V\rightarrow V'$ to the map $V_{x,\Delta\setminus J}^{\rm{an}}\rightarrow V_{w,\Delta\setminus J}^{\rm{an}}$ with $k_0=\#J$. We obtain a strict map between filtered $E$-vector space $\mathbf{E}_{x,J}\twoheadrightarrow \mathbf{E}_{w,J}$, with $\varphi_{\ast}(\mathrm{Fil}_{C_{u}^{\Delta\setminus J}}(\mathbf{E}_{x,J}))$ being non-zero if and only if $w\unlhd u$, in which case it equals $\mathrm{Fil}_{C_{u}^{\Delta\setminus J}}(\mathbf{E}_{w,J})$.
We recover (\ref{equ: intro coxeter filtration}) by taking
\[\mathrm{Fil}_{x}(\mathbf{E}_{J})\defeq \mathrm{Fil}_{C_{x}^{\Delta\setminus J}}(\mathbf{E}_{1,J})\]
for each $x\in\Gamma^{J}$, and Theorem~\ref{thm: intro coxeter filtration} follows from our discussion below Condition~\ref{cond: intro abstract filtration} by taking $V=V_{\Delta\setminus J}^{\rm{an}}$ and $k_0=\#J$.
\subsubsection{Normalized cup product map}\label{subsubsec: normalized cup}
In this section, we define the \emph{normalized cup product map} $\kappa_{J,J'}$ (see (\ref{equ: intro main cup})). Then we give an outline of the proof of Theorem~\ref{thm: intro cup graded}.
We write $r\defeq \#\Delta$ for short.
We fix the choice of a pair $J,J'\subseteq \Delta$ with $J\cap J'=\emptyset$ throughout this section.

For each $I\supseteq\Delta\setminus J$ and $I'\supseteq\Delta\setminus J'$, we construct in (\ref{equ: Levi cup term}) maps between complex
\begin{multline}\label{equ: intro Levi cochain cup}
\mathrm{Tot}(C^{\bullet}(L_{I})\otimes_EC^{\bullet}(L_{I'}))\\
\rightarrow \mathrm{Tot}(\Hom_{D(L_{I\cap I'})}(\bigotimes_{E}^{\bullet+1}D(L_{I\cap I'})\otimes_E(\bigotimes_E^{\bullet+1}D(L_{I\cap I'})\otimes_E1_{L_{I\cap I'}}^{\vee}),1_{L_{I\cap I'}}^{\vee})) \leftarrow C^{\bullet}(L_{I\cap I'})
\end{multline}
which are functorial with respect to the choice of $I$ and $I'$, with the leftward map being a quasi-isomorphism. Modulo a careful treatment of some standard sign twists (that show up while considering total complex of the same multi-complex taken in different orders of the indices, cf.~the discussion towards the end of \S~\ref{subsec: notation}), we obtain a map in the derived category of $E$-vector spaces (see (\ref{equ: Levi cup total}) and (\ref{equ: Levi cup Tits normalize}))
\begin{equation}\label{equ: intro main cup complex}
\mathrm{Tot}(\mathrm{Tot}(\cT_{J}^{\bullet,\bullet})[-r]\otimes_E\mathrm{Tot}(\cT_{J'}^{\bullet,\bullet})[-r])\dashrightarrow \mathrm{Tot}(\cT_{J\sqcup J'}^{\bullet,\bullet})[-r]
\end{equation}
and therefore a map
\begin{equation}\label{equ: intro main cup general} H^{k-r}(\mathrm{Tot}(\cT_{J}^{\bullet,\bullet}))\otimes_EH^{k'-r}(\mathrm{Tot}(\cT_{J'}^{\bullet,\bullet}))\buildrel\cup\over\longrightarrow H^{k+k'-r}(\mathrm{Tot}(\cT_{J\sqcup J'}^{\bullet,\bullet}))
\end{equation}
for each $k,k'\in\Z$. We obtain $\kappa_{J,J'}$ by taking $k=2\#J$ and $k'=2\#J'$ in (\ref{equ: intro main cup general}).
By truncation on the first index of $\cT_{J}^{\bullet,\bullet}$, $\cT_{J'}^{\bullet,\bullet}$ and of $\cT_{J\sqcup J'}^{\bullet,\bullet}$, we see that $\kappa_{J,J'}$ restricts to (\ref{equ: intro main cup Fil}) for each $r-\#J\leq\ell_0\leq r$ and $r-\#J'\leq\ell_1\leq r$ with $\ell_2\defeq\ell_0+\ell_1-r$, and thus further induces a map of the form (\ref{equ: intro graded cup}).
The proof of \ref{it: intro cup graded 0} of Theorem~\ref{thm: intro cup graded} is another careful treatment of the sign twists that appear between total complex of the same multi-complex taken with respect to different orders of the indices (see the proof of Lemma~\ref{lem: grade cup commute}).
As a starting point of the proof of \ref{it: intro cup graded 1} of Theorem~\ref{thm: intro cup graded}, we observe that (\ref{equ: intro main cup complex}) induces (again by truncation on the first index of $\cT_{J}^{\bullet,\bullet}$, $\cT_{J'}^{\bullet,\bullet}$ and of $\cT_{J\sqcup J'}^{\bullet,\bullet}$) a map of the form (with $\ell''\defeq \ell+\ell'-r$)
\begin{equation}\label{equ: intro E1 cup}
E_{1,J}^{-\ell,h}\otimes_EE_{1,J'}^{-\ell',h'}\rightarrow E_{1,J\sqcup J'}^{-\ell'',h+h'}
\end{equation}
which is the direct sum of maps of the form (the usual cup product map twisted by a sign $\varepsilon\in\{1,-1\}$)
\[H^{h}(L_{I})\otimes_EH^{h'}(L_{I'})\rightarrow H^{h}(L_{I\cap I'})\otimes_EH^{h'}(L_{I\cap I'})\buildrel \varepsilon\cup\over\longrightarrow H^{h+h'}(L_{I\cap I'})\]
over all $I\supseteq\Delta\setminus J$ with $\#I=\ell$ and all $I'\supseteq\Delta\setminus J'$ with $\#I'=\ell'$, with $\varepsilon$ being independent of the choice of $I$ and $I'$ (and even independent of the choice of $J$ and $J'$, see below (\ref{equ: Levi cup coh sign})).
Now we consider $S\in\cS_{J}$ with $I\subseteq S\cap\Delta$ and $\ell\defeq r-\#S$, as well as $S'\in\cS_{J'}$ with $I'\subseteq S'\cap\Delta$ and $\ell'\defeq r-\#S'$.
To check that the image of $\overline{x}_{S,I}\otimes\overline{x}_{S',I'}$ under (\ref{equ: intro graded cup}) equals $\pm\overline{x}_{S\sqcup S',I\sqcup I'}$ as claimed in \ref{it: intro cup graded 1} of Theorem~\ref{thm: intro cup graded}, it suffices to choose a lift $y_{S,I}\in E_{1,J}^{-\ell,\ell+r}$ (resp.~$y_{S',I'}\in E_{1,J'}^{-\ell',\ell'+r}$, resp.~$y_{S\sqcup S',I\sqcup I'}\in E_{1,J\sqcup J'}^{-\ell'',\ell''+r}$) of $\overline{x}_{S,I}$ (resp.~of $\overline{x}_{S',I'}$, resp.~of $\overline{x}_{S\sqcup S',I\sqcup I'}$) and then check that we have
\begin{equation}\label{equ: intro E1 explicit cup}
y_{S,I}\otimes_Ey_{S',I'} \mapsto \pm y_{S\sqcup S',I\sqcup I'}
\end{equation}
under (\ref{equ: intro E1 cup}) with $h=\ell+r$ and $h'=\ell'+r$.
Here a technical difficulty is that we do not always have convenient choices of such explicit lifts $y_{S,I}$, $y_{S',I'}$ and $y_{S\sqcup S',I\sqcup I'}$ so that (\ref{equ: intro E1 explicit cup}) is easy to check.
We write $J<J'$ if $j<j'$ for each $j\in J$ and $j'\in J'$, and write $J<<J'$ if there exists $j_0\in\Delta$ such that $j<j_0<j'$ for each $j\in J$ and $j'\in J'$ (see Definition~\ref{def: disconnected pair}). 
Our solution to the aforementioned difficulty is to make such choices of explicit lifts $y_{S,I}$, $y_{S',I'}$ and $y_{S\sqcup S',I\sqcup I'}$ and verify (\ref{equ: intro E1 explicit cup}) only when $J<J'$, and then reduce the general cases to such cases by an increasing induction on the number of pairs $(j,j')\in J\times J'$ such that $j>j'$ (see the proof of Theorem~\ref{thm: general cup}).
The injectivity of $\kappa_{J,J'}$ in \ref{it: intro cup graded 2} of Theorem~\ref{thm: intro cup graded} then follows easily from the injectivity of its graded version (\ref{equ: intro graded cup}) as in \ref{it: intro cup graded 1} of Theorem~\ref{thm: intro cup graded}.
Similarly, we prove \ref{it: intro cup graded 3} of Theorem~\ref{thm: intro cup graded} by showing first that both (\ref{equ: intro graded cup}) and (\ref{equ: intro graded cup prime}) are isomorphisms when $J<<J'$.

Last but not the least, we notice that the normalized cup product maps $\kappa_{J,J'}$ can be easily generalized to an ordered tuple of pair wise disjoint subsets $J_{1},\dots,J_{m}\subseteq\Delta$ for some $m\geq 2$, in which case we obtain a map
\begin{equation}\label{equ: intro multi cup}
\kappa_{J_1,\dots,J_{m}}: \mathbf{E}_{J_1}\otimes_E\cdots\otimes_E\mathbf{E}_{J_{m}}\buildrel\cup\over\longrightarrow\mathbf{E}_{\bigsqcup_{s=1}^{m}J_{s}}.
\end{equation}
The maps (\ref{equ: intro multi cup}) satisfy all associative relations with respect to taking union of two adjacent subsets.
In particular, given an ordered triple of pair wise disjoint subsets $J,J',J''\subseteq\Delta$, we have
\begin{equation}\label{equ: intro cup associative}
\kappa_{J\sqcup J',J''}\circ(\kappa_{J,J'}\otimes\mathrm{Id}_{\mathbf{E}_{J''}})=\kappa_{J,J',J''}=\kappa_{J,J'\sqcup J''}\circ(\mathrm{Id}_{\mathbf{E}_{J}}\otimes\kappa_{J',J''}).
\end{equation}
\subsubsection{Relative conditions and support filtration}\label{subsubsec: relative support}
Let $\fh\subseteq\fg$ be an $E$-Lie subalgebra such that $\fh=\fh_{J_{\fh}}=[\fl_{J_{\fh}},\fl_{J_{\fh}}]$ for some subinterval $J_{\fh}\subseteq\Delta$.
In this section, we study the following map
\begin{equation}\label{equ: intro relative subspace}
\mathbf{E}_{J,\fh}\defeq H^{2\#J-\#\Delta}(\cT_{J,\fh}^{\bullet,\bullet})\rightarrow H^{2\#J-\#\Delta}(\cT_{J}^{\bullet,\bullet})=\mathbf{E}_{J}
\end{equation}
induced from the map between double complex $\cT_{J,\fh}^{\bullet,\bullet}\rightarrow \cT_{J}^{\bullet,\bullet}$ from (\ref{equ: intro Tits double drop h}).

To study (\ref{equ: intro relative subspace}), we introduce what we call \emph{support filtration} 
\begin{equation}\label{equ: intro support filtration}
\{\tau^{J'}(\mathbf{E}_{J,\fh})\}_{J'\subseteq J}
\end{equation}
on $\mathbf{E}_{J,\fh}$, and similarly $\{\tau^{J'}(\mathbf{E}_{J})\}_{J'\subseteq J}$ on $\mathbf{E}_{J}$.
For each $I\subseteq\Delta$, it follows from (\ref{equ: intro T dagger projection}) that we have an isogeny
\[(\prod_{j\in\Delta\setminus I}\Z_{j})\times L_{I}'\buildrel\sim\over\longrightarrow Z_{I}^{\dagger}\times L_{I}'\rightarrow L_{I},\]
which induces the following quasi-isomorphisms
\begin{multline}\label{equ: intro Levi cochain natural}
C^{\bullet}(L_{I},\fl_{I}\cap\fh)\rightarrow C^{\bullet}((\prod_{j\in\Delta}\Z_{j})\times L_{I}',\fl_{I}\cap\fh)\\
\rightarrow C^{\bullet}(L_{I},\fl_{I}\cap\fh)_{\natural}\defeq \mathrm{Tot}(\Hom_{E}(\bigotimes_{j\in\Delta\setminus I}B_{\bullet}(\Z_{j})_{\Z_{j}},C^{\bullet}(L_{I}',\fl_{I}\cap\fh)))
\end{multline}
with the last map (induced from the shuffle product, see \cite[Prop.~8.6.13]{Wei94}) being a quasi-isomorphism by Eilenberg-Zilber's theorem (see \cite[\S~8.5]{Wei94}).
Here $B_{\bullet}(\Z_{j})\defeq \bigotimes_{E}^{\bullet+1}D(\Z_{j})$ is the Bar resolution of the trivial $D(\Z_{j})$-module with $B_{\bullet}(\Z_{j})_{\Z_{j}}$ being its $\Z_{j}$-coinvariant.
Upon replacing $C^{\bullet}(L_{I},\fl_{I}\cap\fh)$ with $C^{\bullet}(L_{I},\fl_{I}\cap\fh)_{\natural}$ in the definition of $\cT_{J,\fh}^{\bullet,\bullet}$ (see the beginning of \S~\ref{subsubsec: Tits double}), we obtain a double complex $\cT_{J,\fh,\natural}^{\bullet,\bullet}$.
We check in (\ref{equ: Levi separate discrete diagram}) that the maps (\ref{equ: intro Levi cochain natural}) are actually functorial with respect to the choice of $I$, from which we obtain a map between double complex $\cT_{J,\fh}^{\bullet,\bullet}\rightarrow \cT_{J,\fh,\natural}^{\bullet,\bullet}$ whose associated map between spectral sequences is an isomorphism from the first page.
Using the canonical truncation map $B_{\bullet}(\Z_{j})_{\Z_{j}}\rightarrow \sigma^{\leq -1}(B_{\bullet}(\Z_{j})_{\Z_{j}})=\tau^{\leq -1}(B_{\bullet}(\Z_{j})_{\Z_{j}})$ for each $j\in\Delta\setminus I$, we can define a suitable truncation $\tau^{J'}(C^{\bullet}(L_{I},\fl_{I}\cap\fh)_{\natural})$ of $C^{\bullet}(L_{I},\fl_{I}\cap\fh)_{\natural}$ which is non-zero only if $I\subseteq\Delta\setminus J'$, such that the truncation map
\[\tau^{J'}(C^{\bullet}(L_{I},\fl_{I}\cap\fh)_{\natural})\rightarrow C^{\bullet}(L_{I},\fl_{I}\cap\fh)_{\natural}\]
is functorial with respect to the choice of $\Delta\setminus J\subseteq I\subseteq \Delta\setminus J'$. Hence, we obtain a map between double complex
\[\tau^{J'}(\cT_{J,\fh,\natural}^{\bullet,\bullet})\rightarrow \cT_{J,\fh,\natural}^{\bullet,\bullet}\]
and therefore a map
\begin{equation}\label{equ: intro support map}
\tau^{J'}(\mathbf{E}_{J,\fh})\defeq H^{2\#J-r}(\mathrm{Tot}(\tau^{J'}(\cT_{J,\fh,\natural}^{\bullet,\bullet})))
\rightarrow H^{2\#J-r}(\mathrm{Tot}(\cT_{J,\fh,\natural}^{\bullet,\bullet}))\buildrel\sim\over\longleftarrow H^{2\#J-r}(\mathrm{Tot}(\cT_{J,\fh}^{\bullet,\bullet}))=\mathbf{E}_{J,\fh}.
\end{equation}
When $\fh=0$, we have the following result.
\begin{prop}[Proposition~\ref{prop: bottom deg graded}]\label{prop: intro support filtration dim}
Let $J\subseteq\Delta$. For each $J'\subseteq J$, the natural map $\tau^{J'}(\mathbf{E}_{J})\rightarrow \mathbf{E}_{J}$ is injective with
\begin{equation}\label{equ: intro support filtration dim}
\Dim_E\tau^{J'}(\mathbf{E}_{J})=\#\{(S,I)\mid S\in\cS_{J},J'\subseteq I\subseteq S\cap\Delta\}=\#\Gamma^{J\setminus J'}.
\end{equation}
\end{prop}
This finishes our definition of the support filtration $\{\tau^{J'}(\mathbf{E}_{J})\}_{J'\subseteq J}$ on $\mathbf{E}_{J}$.
In general, we have the following result concerning (\ref{equ: intro support map})
\begin{prop}[Proposition~\ref{prop: relative bottom embedding}]\label{prop: intro relative sub}
Let $J\subseteq\Delta$ and $\fh=\fh_{J_{\fh}}\subseteq\fg$ for some interval $J_{\fh}\subseteq\Delta$. Then for each $J'\subseteq J$, we have a natural cartesian diagram of the form
\begin{equation}\label{equ: intro support relative to absolute}
\xymatrix{
\tau^{J'}(\mathbf{E}_{J,\fh}) \ar[r] \ar[d] & \tau^{J'}(\mathbf{E}_{J}) \ar[d]\\
\mathbf{E}_{J,\fh} \ar[r] & \mathbf{E}_{J}
}
\end{equation}
with all maps being injective.
\end{prop}
In particular, we obtain a support filtration $\{\tau^{J'}(\mathbf{E}_{J,\fh})\}_{J'\subseteq J}$ on $\mathbf{E}_{J,\fh}$ such that the natural map $\mathbf{E}_{J,\fh}\rightarrow\mathbf{E}_{J}$ is injective and is strict with respect to the support filtration on its source and its target.
\subsubsection{Comparison between different structures}\label{subsubsec: interaction}
Let $J\subseteq \Delta$. 
So far, we have the following genuinely different resources of special $E$-subspaces of $\mathbf{E}_{J}$.
\begin{cons}\label{cons: intro subspace}
\begin{enumerate}[label=(\roman*)]
\item \label{it: intro subspace 1} Following the discussions in \S~\ref{subsubsec: normalized cup}, given an arbitrary ordered partition $J_{1},\dots,J_{m}$ of $J$, we can embed
    \begin{equation}\label{equ: intro partition tensor cup}
    \mathbf{E}_{J_{1}}\otimes_E\cdots\otimes_E\mathbf{E}_{J_{m}}
    \end{equation}
    into $\mathbf{E}_{J}$ via the normalized cup product map $\kappa_{J_1,\dots,J_{m}}$ (see (\ref{equ: intro multi cup})).
    We thus often abuse (\ref{equ: intro partition tensor cup}) for its image in $\mathbf{E}_{J}$. For each $\#\Delta\setminus J\leq\ell\leq\#\Delta$, we can recover $\mathrm{Fil}^{-\ell}(\mathbf{E}_{J})$ as the sum of (\ref{equ: intro partition tensor cup}) over all ordered partitions $J_{1},\dots,J_{m}$ that satisfy $m=\#\Delta-\ell$ and $J_{i}\neq \emptyset$ for each $1\leq i\leq m$.
\item \label{it: intro subspace 2} Following the discussions in \S~\ref{subsubsec: intro coxeter}, we have defined (using the layer structure of $V_{\Delta\setminus J}^{\rm{an}}$, see Theorem~\ref{thm: coxeter filtration}) a Coxeter filtration $\{\mathrm{Fil}_{x}(\mathbf{E}_{J})\}_{x\in\Gamma^{J}}$ on $\mathbf{E}_{J}$. This Coxeter filtration can be equally defined via the collection of surjections $\{\mathbf{E}_{J}\twoheadrightarrow\mathbf{E}_{x,J}\}_{x\in\Gamma^{J}}$.
\item \label{it: intro subspace 3} Let $\fh\subseteq\fg$ be a reductive $E$-Lie subalgebra of the form $\fh=\fh_{J_{\fh}}$ for some interval $J_{\fh}\subseteq\Delta$. Following the discussions in \S~\ref{subsubsec: relative support}, we have shown in Proposition~\ref{prop: intro relative sub} that the following natural map
    \[\mathbf{E}_{J,\fh}\rightarrow\mathbf{E}_{J}\]
    is an embedding and is strict with respect to the support filtration $\{\tau^{J'}(-)\}_{J'\subseteq J}$ on its source and target.
\end{enumerate}
\end{cons}
The comparisons between different pairs of structures in the above list are governed by several natural commutative diagrams.

Let $J\subseteq\Delta$. For each $x\in\Gamma^{J}$ and $\fh\subseteq\fg$ as in \ref{it: intro subspace 3} of Construction~\ref{cons: intro subspace}, we have the following commutative diagram
\begin{equation}\label{equ: intro relative coxeter}
\xymatrix{
\mathbf{E}_{J,\fh} \ar[r] \ar[d] & \mathbf{E}_{J} \ar[d]\\
\mathbf{E}_{x,J,\fh} \ar[r] & \mathbf{E}_{x,J}
}
\end{equation}
with
\[\mathbf{E}_{x,J,\fh}\defeq H^{2\#J-\#\Delta}(\cS_{x,J,\fh}^{\bullet,\bullet}).\]
In particular, the comparison between $\mathbf{E}_{J,\fh}$ with the Coxeter filtration $\{\mathrm{Fil}_{x}(\mathbf{E}_{J})\}_{x\in\Gamma^{J}}$ is often achieved via the computation of $\mathbf{E}_{x,J,\fh}$ for various choices of $x\in\Gamma^{J}$.

Let $J,J'\subseteq\Delta$ be a pair of subsets with $J\cap J'=\emptyset$.
Let $M\subseteq J$ and $M'\subseteq J'$ be subsets.
Let $\fh,\fh'\subseteq\Delta$ be a pair of $E$-Lie subalgebras of the form $\fh=\fh_{J_{\fh}}$ and $\fh'=\fh_{J_{\fh'}}$ for some intervals $J_{\fh},J_{\fh'}\subseteq\Delta$. Parallel to the definition of $\kappa_{J,J'}$ in \S~\ref{subsubsec: normalized cup} using (\ref{equ: intro Levi cochain cup}), we have a commutative diagram
\[
\xymatrix{
\tau^{M}(\mathbf{E}_{J,\fh}) \ar@{^{(}->}[d] & \otimes_E & \tau^{M'}(\mathbf{E}_{J',\fh'}) \ar@{^{(}->}[d] \ar^{\cup}[rr] & & \tau^{M\sqcup M'}(\mathbf{E}_{J\sqcup J',\fh\cap\fh'}) \ar@{^{(}->}[d]\\
\mathbf{E}_{J} & \otimes_E & \mathbf{E}_{J'} \ar^{\cup}[rr]& & \mathbf{E}_{J\sqcup J'}
}
\]
with the bottom horizontal map being $\kappa_{J,J'}$ and all vertical maps being natural embeddings from Proposition~\ref{prop: intro relative sub}. In particular, the comparison between \ref{it: intro subspace 1} and \ref{it: intro subspace 3} of Construction~\ref{cons: intro subspace} can be summarized as the following inclusion
\begin{equation}\label{equ: intro relative cup}
\kappa_{J,J'}(\tau^{M}(\mathbf{E}_{J,\fh})\otimes_E\tau^{M'}(\mathbf{E}_{J',\fh'}))\subseteq \tau^{M\sqcup M'}(\mathbf{E}_{J\sqcup J',\fh\cap\fh'}).
\end{equation}
Upon taking $\fh=\fh'=\fh''=0$ and $M=J$ in (\ref{equ: intro relative cup}), we have
\[\kappa_{J,J'}(\tau^{J}(\mathbf{E}_{J})\otimes_E\tau^{M'}(\mathbf{E}_{J'}))\subseteq \tau^{J\sqcup M'}(\mathbf{E}_{J\sqcup J'}).\]
Now that we have $\Dim_E\tau^{J}(\mathbf{E}_{J})=1$ and
\[\Dim_E\tau^{M'}(\mathbf{E}_{J'})=\#\Gamma^{J'\setminus M'}=\Dim_E\tau^{J\sqcup M'}(\mathbf{E}_{J\sqcup J'})\]
from Proposition~\ref{prop: intro support filtration dim}, we conclude that
\begin{equation}\label{equ: intro support cup}
\kappa_{J,J'}(\tau^{J}(\mathbf{E}_{J})\otimes_E\tau^{M'}(\mathbf{E}_{J'}))=\tau^{J\sqcup M'}(\mathbf{E}_{J\sqcup J'}).
\end{equation}

Now we consider the comparison between \ref{it: intro subspace 1} and \ref{it: intro subspace 2} of Construction~\ref{cons: intro subspace}, and then give an outline of the proof of Theorem~\ref{thm: intro cup x y}.
The key ingredient of such a comparison is the construction of a generalization of $\kappa_{J,J'}$ of the form
\begin{equation}\label{equ: intro cup x y w}
\mathbf{E}_{x,J}\otimes_E\mathbf{E}_{y,J'}\buildrel\cup\over\longrightarrow\mathbf{E}_{w,J\sqcup J'}
\end{equation}
for each $x\in\Gamma^{J}$, $y\in\Gamma^{J'}$ and $w\in\Gamma_{x,y}\subseteq\Gamma^{J\sqcup J'}$ (see (\ref{equ: x y envelop}) for the definition of $\Gamma_{x,y}$).
We fix the choice of a surjection (which is unique up to scalars) $M(1)\twoheadrightarrow M^{I}(1)$ for each $I\subseteq\Delta$.
We also fix the choice of a embedding (which is also unique up to scalars, see \cite[\S~4.2, \S~4.6]{Hum08}) $M(u)\hookrightarrow M(1)$ for each $u\in W(G)$.
The diagonal map $U(\fg)\rightarrow U(\fg)\otimes_EU(\fg)$ restricts to a diagonal map $U(\fb)\rightarrow U(\fb)\otimes_EU(\fb)$ and thus induces a diagonal map 
\begin{equation}\label{equ: intro full Verma diagonal}
M(1)\rightarrow M(1)\otimes_EM(1).
\end{equation}
For each $I\supseteq\Delta\setminus J$ and $I'\supseteq \Delta\setminus J'$, our fixed choice of surjections $M(1)\twoheadrightarrow M^{I}(1)$ and $M(1)\twoheadrightarrow M^{I'}(1)$ together with (\ref{equ: intro full Verma diagonal}) induce a map
\[M(1)\rightarrow M^{I}(1)\otimes_E M^{I'}(1)\]
which factors through a map (using \cite[Thm~9.4 (c)]{Hum08})
\begin{equation}\label{equ: intro Verma diagonal}
M^{I\cap I'}(1)\rightarrow M^{I}(1)\otimes_E M^{I'}(1).
\end{equation}
Let $x\in\Gamma^{J}$, $y\in\Gamma^{J'}$ and $w\in\Gamma_{x,y}\subseteq\Gamma^{J\sqcup J'}$. For each $\Delta\setminus J\subseteq I\subseteq J_{x}$ and $\Delta\setminus J'\subseteq I'\subseteq J_{y}$, we prove that (\ref{equ: intro Verma diagonal}) together with our fixed choices of embeddings from $M(x)$, $M(y)$ and $M(w)$ into $M(1)$ uniquely determines a commutative diagram of non-zero maps of the form (see (\ref{equ: tensor of parabolic Verma diagram}) and Lemma~\ref{lem: tensor of parabolic Verma})
\begin{equation}\label{equ: intro cup diagram Lie}
\xymatrix{
M^{I\cap I'}(w) \ar[r] \ar[d] & M^{I}(x)\otimes_E M^{I'}(y) \ar[d]\\
M^{I\cap I'}(1) \ar[r] & M^{I}(1)\otimes_E M^{I'}(1)
}
\end{equation}
which is functorial with respect to the choice of $x$, $y$, $w$, $I$ and $I'$.
Parallel to the construction of (\ref{equ: intro Levi cochain cup}) and (\ref{equ: intro main cup complex}), we can use the top horizontal map of (\ref{equ: intro cup diagram Lie}) to define the following maps between complex
\begin{multline}\label{equ: intro cochain cup x y}
\mathrm{Tot}(\Hom_{D(\fg,P_{I})}(\bigotimes_{E}^{\bullet+1}D(\fg,P_{I})\otimes_EM^{I}(x),1_{D(\fg,P_{I})})\otimes_E\Hom_{D(\fg,P_{I'})}(\bigotimes_{E}^{\bullet+1}D(\fg,P_{I'})\otimes_EM^{I'}(y),1_{D(\fg,P_{I'})}))\\
\rightarrow \mathrm{Tot}(\Hom_{D(\fg,P_{I\cap I'})}(\bigotimes_{E}^{\bullet+1}D(\fg,P_{I\cap I'})\otimes_E(M^{I}(x)\otimes_E(\bigotimes_{E}^{\bullet+1}D(\fg,P_{I\cap I'})\otimes_EM^{I'}(y))),1_{D(\fg,P_{I\cap I'})}))\\
\leftarrow \Hom_{D(\fg,P_{I\cap I'})}(\bigotimes_{E}^{\bullet+1}D(\fg,P_{I\cap I'})\otimes_E(M^{I}(x)\otimes_EM^{I'}(y)),1_{D(\fg,P_{I\cap I'})})\\
\rightarrow \Hom_{D(\fg,P_{I\cap I'})}(\bigotimes_{E}^{\bullet+1}D(\fg,P_{I\cap I'})\otimes_EM^{I\cap I'}(w),1_{D(\fg,P_{I\cap I'})})
\end{multline}
which are functorial with respect to the choice of $I$ and $I'$. Hence, we obtain the following map in the derived category of $E$-vector spaces
\begin{equation}\label{equ: intro main cup complex x y}
\mathrm{Tot}(\mathrm{Tot}(\cS_{x,J,\flat}^{\bullet,\bullet})[-r]\otimes_E\mathrm{Tot}(\cS_{y,J',\flat}^{\bullet,\bullet})[-r])\dashrightarrow \mathrm{Tot}(\cS_{w,J\sqcup J',\flat}^{\bullet,\bullet})[-r]
\end{equation}
and therefore a map
\begin{equation}\label{equ: intro cup x y general} H^{k-r}(\mathrm{Tot}(\cS_{x,J,\flat}^{\bullet,\bullet}))\otimes_EH^{k'-r}(\mathrm{Tot}(\cS_{y,J',\flat}^{\bullet,\bullet}))\buildrel\cup\over\longrightarrow H^{k+k'-r}(\mathrm{Tot}(\cS_{w,J\sqcup J',\flat}^{\bullet,\bullet}))
\end{equation}
for each $k,k'\in\Z$. We finish the construction of (\ref{equ: intro cup x y w}) by taking $k=2\#J$ and $k'=2\#J'$ in (\ref{equ: intro cup x y general}). Now that the maps in (\ref{equ: intro cochain cup x y}) are functorial with respect to the choice of $x$, $y$ and $w$, so are the maps (\ref{equ: intro main cup complex x y}) and (\ref{equ: intro cup x y general}), which together with (\ref{equ: intro St total coh x}) finish the construction of the commutative diagram (\ref{equ: intro cup x y diagram}) from \ref{it: intro cup x y 1} of Theorem~\ref{thm: intro cup x y}. 

Now we move onto the proof of \ref{it: intro cup x y 2} of Theorem~\ref{thm: intro cup x y}. 
We fix a choice of $x\in\Gamma^{J}$ and $y\in\Gamma^{J'}$. To prove (\ref{equ: intro cup x y}), it suffices to show that the composition of the following maps
\begin{equation}\label{equ: intro cup x y composition}
\mathrm{Fil}_{x}(\mathbf{E}_{J})\otimes_E\mathrm{Fil}_{y}(\mathbf{E}_{J'})\hookrightarrow \mathbf{E}_{J}\otimes_E\mathbf{E}_{J'}\buildrel\cup\over\longrightarrow \mathbf{E}_{J\sqcup J'}\twoheadrightarrow \mathbf{E}_{w',J\sqcup J'}
\end{equation}
is non-zero only if there exists $w\in\Gamma_{x,y}$ such that $w'\unlhd w$. 
By \cite[Thm.~2.2.2]{BB05} we know that there exists a unique pair of elements $x'\in\Gamma^{J}$ and $y'\in\Gamma^{J'}$ such that $w'\in\Gamma_{x',y'}$. Upon replacing $x$, $y$ and $w$ in (\ref{equ: intro cup x y diagram}) with $x'$, $y'$ and $w'$, we know that the composition of (\ref{equ: intro cup x y composition}) factors through
\[\mathrm{Fil}_{x}(\mathbf{E}_{J})\otimes_E\mathrm{Fil}_{y}(\mathbf{E}_{J'})\rightarrow  \mathbf{E}_{x',J}\otimes_E\mathbf{E}_{y',J'}\rightarrow \mathbf{E}_{w',J\sqcup J'}\]
and thus is non-zero only if $x'\unlhd x$ and $y'\unlhd y$, which forces the existence of $w\in\Gamma_{x,y}$ satisfying $w'\unlhd w$ by Lemma~\ref{lem: x y envelop partial order}.
In other words, (\ref{equ: intro cup x y}) is essentially a consequence of (\ref{equ: intro cup x y diagram}) (upon replacing $x$, $y$ and $w$ in \emph{loc.cit.} with varying $x'$, $y'$ and $w'$) and the result Lemma~\ref{lem: x y envelop partial order} on combinatorics of partial-Coxeter elements.
Concerning the proof of the second half of \ref{it: intro cup x y 2} of Theorem~\ref{thm: intro cup x y}, the main ingredient is a Lie theoretic cup product map
\begin{equation}\label{equ: intro Lie cup}
\mathrm{Ext}_{U(\fg)}^{\ell(x)}(M^{J_{x}}(x),1_{\fg})\otimes_E\mathrm{Ext}_{U(\fg)}^{\ell(y)}(M^{J_{y}}(y),1_{\fg})
\buildrel\cup\over\longrightarrow \mathrm{Ext}_{U(\fg)}^{\ell(w)}(M^{J_{w}}(w),1_{\fg})
\end{equation}
constructed from the top horizontal map of (\ref{equ: intro cup diagram Lie}).
The crucial fact here is that (\ref{equ: intro Lie cup}) is an isomorphism between $1$-dimensional $E$-vector spaces, which we prove by eventually reducing to the case $y=s_j$ and $w=s_jx>x$ for some $j\in J'$.

We finish this section by a further discussion on the support filtration $\{\tau^{J'}(\mathbf{E}_{J})\}$ on $\mathbf{E}_{J}$ as well as a refinement of Proposition~\ref{prop: intro relative sub} on the description of $\mathrm{gr}^{J'}(\mathbf{E}_{J,\fh})$ for each $J'\subseteq J$.
We write
\[\Gamma_{J}\defeq \{u\in\Gamma\mid \mathrm{Supp}(u)=J\}\]
for each $J\subseteq\Delta$ and note that $\Gamma^{J}=\bigsqcup_{J'\subseteq J}\Gamma_{J'}$.
We start with the following claim
\begin{equation}\label{equ: intro sm sub}
\tau^{J}(\mathbf{E}_{J})=\mathrm{Fil}_{1}(\mathbf{E}_{J}).
\end{equation}
In fact, the equality (\ref{equ: intro sm sub}) holds when $\#J=1$ by a direct check, and the general cases follow from an increasing induction on $\#J$ based on the equality (\ref{equ: intro support cup}) (with $M=J'$ in \emph{loc.cit.}) and the equality (\ref{equ: intro sm cup}) (with $y=1$ in \emph{loc.cit.}).
Then by taking $M=\emptyset$ in (\ref{equ: intro support cup}) and applying (\ref{equ: intro sm cup}) for each $y\in\Gamma^{J'}$, we have
\begin{multline*}
\tau^{J}(\mathbf{E}_{J\sqcup J'})=\kappa_{J,J'}(\tau^{J}(\mathbf{E}_{J})\otimes_E\mathbf{E}_{J'})=\kappa_{J,J'}(\tau^{J}(\mathbf{E}_{J})\otimes_E(\sum_{y\in\Gamma^{J'}}\mathrm{Fil}_{y}(\mathbf{E}_{J'})))\\
=\sum_{y\in\Gamma^{J'}}\kappa_{J,J'}(\tau^{J}(\mathbf{E}_{J})\otimes_E\mathrm{Fil}_{y}(\mathbf{E}_{J'}))=\sum_{y\in\Gamma^{J'}}\mathrm{Fil}_{y}(\mathbf{E}_{J\sqcup J'})
\end{multline*}
for each $J,J'\subseteq\Delta$ with $J\cap J'=\emptyset$. 
Consequently, we have
\begin{equation}\label{equ: intro support comparison}
\tau^{J'}(\mathbf{E}_{J})=\sum_{u\in\Gamma^{J\setminus J'}}\mathrm{Fil}_{u}(\mathbf{E}_{J})
\end{equation}
for each $J'\subseteq J$. In particular, the support filtration $\{\tau^{J'}(\mathbf{E}_{J})\}_{J'\subseteq J}$ on $\mathbf{E}_{J}$ can be entirely recovered from the Coxeter filtration, with graded piece given by
\begin{equation}\label{equ: intro support comparison graded}
\mathrm{gr}^{J'}(\mathbf{E}_{J})=\sum_{u\in\Gamma_{J\setminus J'}}\mathrm{gr}_{u}(\mathbf{E}_{J})
\end{equation}
for each $J'\subseteq J$. Now we are ready to state the following refinement of Proposition~\ref{prop: intro relative sub}.
\begin{prop}[Proposition~\ref{prop: grade M relative image}, Proposition~\ref{prop: graded coxeter tau comparison}]\label{prop: intro relative sub graded}
Let $J\subseteq\Delta$ and $\fh=\fh_{J_{\fh}}\subseteq\fg$ for some interval $J_{\fh}\subseteq\Delta$. Then for each $J'\subseteq J$, the natural map $\mathrm{gr}^{J'}(\mathbf{E}_{J,\fh})\rightarrow \mathrm{gr}^{J'}(\mathbf{E}_{J})$ (see Proposition~\ref{prop: intro relative sub}) is an embedding with image
\[\bigoplus_{u\in\Gamma_{J\setminus J'}, D_{L}(u)\cap J_{\fh}=\emptyset}\mathrm{gr}_{u}(\mathbf{E}_{J})\]
under (\ref{equ: intro support comparison graded}).
\end{prop}
\subsubsection{The commutativity of normalized cup product map}\label{subsubsec: commute}
In this section, we give an outline of the proof of Theorem~\ref{thm: intro cup commute}.

Combining \ref{it: intro cup graded 3} of Theorem~\ref{thm: intro cup graded} with the associative relation (\ref{equ: intro cup associative}) for each ordered disjoint triple $J,J',J''$, we can reduce the proof of Theorem~\ref{thm: intro cup commute} to the case when $J<J'$ (namely $j<j'$ for each $j\in J$ and $j'\in J'$), with both $J$ and $J'$ being intervals.

To finish the proof of Theorem~\ref{thm: intro cup commute}, we claim that it suffices to construct $1$-dimensional $E$-subspaces $\mathbf{E}_{J}^{\heartsuit}\subseteq \mathbf{E}_{J}$, $\mathbf{E}_{J'}^{\heartsuit}\subseteq \mathbf{E}_{J'}$ that satisfy the following conditions.
\begin{cond}\label{cond: intro cup heart}
Let $J,J'\subseteq\Delta$ be subintervals with $J<J'$.
\begin{enumerate}[label=(\roman*)]
\item \label{it: intro cup heart 1} We have $\mathbf{E}_{J}=\mathbf{E}_{J}^{\heartsuit}\oplus\mathbf{E}_{J}^{<}$ and $\mathbf{E}_{J'}=\mathbf{E}_{J'}^{\heartsuit}\oplus\mathbf{E}_{J'}^{<}$.
\item \label{it: intro cup heart 2} We have
\[\kappa_{J,J'}(\mathbf{E}_{J}^{\heartsuit}\otimes_E\mathbf{E}_{J'}^{\heartsuit})=\kappa_{J',J}(\mathbf{E}_{J'}^{\heartsuit}\otimes_E\mathbf{E}_{J}^{\heartsuit})\subseteq \mathbf{E}_{J\sqcup J'}.\]
\end{enumerate}
\end{cond}
In fact, it follows from Condition~\ref{cond: intro cup heart} and \ref{it: intro cup graded 0} of Theorem~\ref{thm: intro cup graded} that (\ref{equ: intro cup commute}) holds when we have $x\in \mathbf{E}_{J}^{\heartsuit}$ and $y\in \mathbf{E}_{J'}^{\heartsuit}$.
If either $x\in\mathbf{E}_{J}^{<}$ or $y\in \mathbf{E}_{J'}^{<}$, we could reduce (\ref{equ: intro cup commute}) (again using (\ref{equ: intro cup associative}) for each triple $J,J',J''$, as well as \ref{it: intro cup graded 3} of Theorem~\ref{thm: intro cup graded}) to cases with strictly smaller $J$ or $J'$. Hence, given the construction of $\mathbf{E}_{J}^{\heartsuit}$ and $\mathbf{E}_{J'}^{\heartsuit}$ satisfying Condition~\ref{cond: intro cup heart}, we can finish the proof of Theorem~\ref{thm: intro cup commute} by an increasing induction on $\#J\sqcup J'$.

We devote the rest of this section to explain the construction of $\mathbf{E}_{J}^{\heartsuit}$ and $\mathbf{E}_{J'}^{\heartsuit}$ that satisfy Condition~\ref{cond: intro cup heart}. Following the discussion around (\ref{equ: cup interval reduction}), it is harmless to assume that $J=[1,j_0]$ and $J'=[j_1,n-1]$ for some $j_0<j_1$. 

We first construct some $E$-subspaces of $\mathbf{E}_{J}$, $\mathbf{E}_{J'}$ and $\mathbf{E}_{J\sqcup J'}$ using explicit relative conditions as below.
For each $j\in\Delta$, we write $\widehat{j}\defeq\Delta\setminus\{j\}$ for short.
Recall that we write $H_{I}=[L_{I},L_{I}]$ with associated $E$-Lie algebra $\fh_{I}=[\fl_{I},\fl_{I}]$ for each $I\subseteq\Delta$.
Now we introduce a few notation following the discussion above (\ref{equ: special relative j0}).
We write $J_0\defeq\Delta\setminus J=[j_0+1,n-1]\subseteq \widehat{j_0}$, $J_1\defeq\Delta\setminus J=[1,j_1-1]\subseteq \widehat{j_1}$ and $J_2\defeq \Delta\setminus(J\sqcup J')=J_0\cap J_1=[j_0+1,j_1-1]\subseteq \widehat{j_0}\cap \widehat{j_1}$ for short.
We define $\fh_{0}\defeq (\fh_{[2,n-1]}\cap\fl_{\widehat{j_0}})/\fh_{J_0}$ (resp.~$\fh_{1}\defeq (\fh_{[1,n-2]}\cap\fl_{\widehat{j_1}})/\fh_{J_1}$, resp.~$\fh_{2}\defeq (\fh_{[2,n-2]}\cap\fl_{\widehat{j_0}\cap\widehat{j_1}})/\fh_{J_2}$) which is naturally an $E$-Lie subalgebra of the $E$-Lie algebras $\fl_{\widehat{j_0}}/\fh_{J_0}$ (resp.~$\fl_{\widehat{j_1}}/\fh_{J_1}$, resp.~$\fl_{\widehat{j_0}\cap\widehat{j_1}}/\fh_{J_2}$) associated with the locally $K$-analytic groups $L_{\widehat{j_0}}/H_{J_0}Z_{\widehat{j_0}}^{\dagger}$ (resp.~$L_{\widehat{j_1}}/H_{J_1}Z_{\widehat{j_1}}^{\dagger}$, resp.~$L_{\widehat{j_0}\cap\widehat{j_1}}/H_{J_2}Z_{\widehat{j_0}\cap\widehat{j_1}}^{\dagger}$).
Note that we actually have $H_{J_2}Z_{\widehat{j_0}\cap\widehat{j_1}}^{\dagger}\subseteq (H_{J_0}Z_{\widehat{j_0}}^{\dagger})\cap(H_{J_1}Z_{\widehat{j_1}}^{\dagger})$.
Upon replacing $C^{\bullet}(L_{I},\fl_{I}\cap\fh)$ for each $\Delta\setminus J\subseteq I\subseteq \Delta$ in the definition of $\cT_{J,\fh}^{\bullet,\bullet}$ (see the beginning of \S~\ref{subsubsec: Tits double}) with $C^{\bullet}(L_{I})_{0,\dagger}\defeq C^{\bullet}(L_{I}/H_{J_0}Z_{\widehat{j_0}}^{\dagger},(\fl_{I}/\fh_{J_0})\cap\fh_{0})$ for each $J_0\subseteq I\subseteq\widehat{j_0}$ (resp.~with $C^{\bullet}(L_{I})_{1,\dagger}\defeq C^{\bullet}(L_{I}/H_{J_1}Z_{\widehat{j_1}}^{\dagger},(\fl_{I}/\fh_{J_1})\cap\fh_{1})$ for each $J_1\subseteq I\subseteq\widehat{j_1}$, resp.~with $C^{\bullet}(L_{I})_{2,\dagger}\defeq C^{\bullet}(L_{I}/H_{J_2}Z_{\widehat{j_0}\cap \widehat{j_1}}^{\dagger},(\fl_{I}/\fh_{J_2})\cap\fh_{2})$ for each $J_2\subseteq I\subseteq\widehat{j_0}\cap\widehat{j_1}$), we can define a double complex $\cT_{0,\dagger}^{\bullet,\bullet}$ (resp.~$\cT_{1,\dagger}^{\bullet,\bullet}$, resp.~$\cT_{2,\dagger}^{\bullet,\bullet}$) which comes with a map $\cT_{0,\dagger}^{\bullet,\bullet}\rightarrow \cT_{J}^{\bullet,\bullet}$ (resp.~a map $\cT_{1,\dagger}^{\bullet,\bullet}\rightarrow \cT_{J'}^{\bullet,\bullet}$, resp.~a map $\cT_{2,\dagger}^{\bullet,\bullet}\rightarrow \cT_{J\sqcup J'}^{\bullet,\bullet}$). We write 
\[\mathbf{E}_{a,\dagger}\defeq H^{\#\Delta-2\#J_{a}}(\mathrm{Tot}(\cT_{a,\dagger}^{\bullet,\bullet}))\] 
for each $a=0,1,2$, which come with natural maps $\mathbf{E}_{0,\dagger}\rightarrow \mathbf{E}_{J}$, $\mathbf{E}_{1,\dagger}\rightarrow \mathbf{E}_{J'}$ and $\mathbf{E}_{2,\dagger}\rightarrow \mathbf{E}_{J\sqcup J'}$. Parallel to the definition of $\kappa_{J,J'}$, we can define a normalized cup product map
\begin{equation}\label{equ: intro cup relative}
\mathbf{E}_{0,\dagger}\otimes_E\mathbf{E}_{1,\dagger}\buildrel\cup\over\longrightarrow \mathbf{E}_{2,\dagger}
\end{equation}
using a natural map in the derived category of $E$-vector spaces
\[\mathrm{Tot}(C^{\bullet}(L_{I})_{0,\dagger}\otimes_EC^{\bullet}(L_{I'})_{1,\dagger})\dashrightarrow C^{\bullet}(L_{I\cap I'})_{2,\dagger}\]
which is functorial with respect to the choice of $J_0\subseteq I\subseteq\widehat{j_0}$ and $J_1\subseteq I'\subseteq\widehat{j_1}$.
Moreover, (\ref{equ: intro cup relative}) fits into a commutative diagram of the form
\begin{equation}\label{equ: intro cup relative diagram}
\xymatrix{
\mathbf{E}_{0,\dagger} \ar[d] & \otimes_E & \mathbf{E}_{1,\dagger} \ar^{\cup}[rr] \ar[d] & & \mathbf{E}_{2,\dagger} \ar[d]\\
\mathbf{E}_{J} & \otimes_E & \mathbf{E}_{J'} \ar^{\cup}[rr] & & \mathbf{E}_{J\sqcup J'}
}.
\end{equation}
We write $w_{1,j}\defeq s_1\cdots s_{j-1}s_j$ and $w_{n-1,j}\defeq s_{n-1}\cdots s_{j+1}s_j$ for each $j\in\Delta$.
We define
\[
\left\{\begin{array}{ccccccc}
\Gamma_{0,\dagger} &\defeq &\{x\in\Gamma^{J}\mid D_L(x)=\{1\}\} &= &\{w_{1,j}\mid j\in J\} & \subseteq & \Gamma^{J}\\
\Gamma_{1,\dagger} &\defeq &\{x\in\Gamma^{J'}\mid D_L(y)=\{n-1\}\} &= &\{w_{n-1,j}\mid j\in J'\} & \subseteq & \Gamma^{J'}\\
\Gamma_{2,\dagger} &\defeq &\{w\in\Gamma^{J\sqcup J'}\mid D_L(w)=\{1,n-1\}\} &= &\bigsqcup_{x\in\Gamma_{0,\dagger},y\in\Gamma_{1,\dagger}}\Gamma_{x,y} & \subseteq & \Gamma^{J\sqcup J'}
\end{array}\right.
\]
with the equalities between the middle two columns being easy exercises on the combinatorics of partial-Coxeter elements.
For each $i\geq 0$, we write
\[\tau^{i}(\mathbf{E}_{J})\defeq \sum_{M\subseteq J,\#M\leq i}\tau^{M}(\mathbf{E}_{J})\]
and write $\tau^{i}(-)$ for the induced filtration on each $E$-subspace of $\mathbf{E}_{J}$. We similarly define $\tau^{i}(-)$ for $\mathbf{E}_{J'}$ and $\mathbf{E}_{J\sqcup J'}$ as well as the induced filtration on their $E$-subspaces.
The following is our main result concerning the vertical maps in the commutative diagram (\ref{equ: intro cup relative diagram}).
\begin{prop}[Proposition~\ref{prop: relative coh E2 grade}, Proposition~\ref{prop: relative coh E2 grade dagger}, Proposition~\ref{prop: graded coxeter tau comparison}]\label{prop: intro relative image}
We have the following results.
\begin{enumerate}[label=(\roman*)]
\item \label{it: intro relative image 1} The natural map $\mathbf{E}_{0,\dagger}\rightarrow \mathbf{E}_{J}$ is injective, and the induced map between graded pieces $\mathrm{gr}^{i}(\mathbf{E}_{0,\dagger})\rightarrow\mathrm{gr}^{i}(\mathbf{E}_{J})$ has image \[\bigoplus_{x\in\Gamma_{0,\dagger},\ell(x)=\#J-i}\mathrm{gr}_{x}(\mathbf{E}_{J})=\mathrm{gr}_{w_{1,j_0-i}}(\mathbf{E}_{J})\] 
    for each $0\leq i\leq \#J$.
\item \label{it: intro relative image 2} The natural map $\mathbf{E}_{1,\dagger}\rightarrow \mathbf{E}_{J'}$ is injective, and the induced map between graded pieces $\mathrm{gr}^{i}(\mathbf{E}_{1,\dagger})\rightarrow\mathrm{gr}^{i}(\mathbf{E}_{J'})$ has image \[\bigoplus_{y\in\Gamma_{1,\dagger},\ell(y)=\#J'-i}\mathrm{gr}_{y}(\mathbf{E}_{J'})=\mathrm{gr}_{w_{n-1,j_1+i}}(\mathbf{E}_{J'})\]  
    for each $0\leq i\leq \#J'$.
\item \label{it: intro relative image 3} The natural map $\mathbf{E}_{2,\dagger}\rightarrow \mathbf{E}_{J\sqcup J'}$ is injective, and the induced map between graded pieces $\mathrm{gr}^{i}(\mathbf{E}_{2,\dagger})\rightarrow\mathrm{gr}^{i}(\mathbf{E}_{J\sqcup J'})$ has image contained in
    \[\bigoplus_{w\in\Gamma_{2,\dagger},\ell(w)=\#J\sqcup J'-i}\mathrm{gr}_{w}(\mathbf{E}_{J\sqcup J'})\]  
    for each $0\leq i\leq \#J\sqcup J'$.
\end{enumerate}
\end{prop}

Now we construct some $E$-subspaces of $\mathbf{E}_{J}$, $\mathbf{E}_{J'}$ and $\mathbf{E}_{J\sqcup J'}$ coming from the Coxeter filtration on them (cf.~Theorem~\ref{thm: intro coxeter filtration}).
Given $M\subseteq\Delta$, we say a subset $\Gamma'\subseteq\Gamma^{M}$ is \emph{closed} if any $x\in\Gamma^{M}$ satisfying $x\unlhd y$ for some $y\in\Gamma'$ also satisfy $x\in\Gamma'$.
For each closed subset $\Gamma'\subseteq \Gamma^{M}$, we write
\[\mathrm{Fil}_{\Gamma'}(\mathbf{E}_{M})\defeq\sum_{u\in\Gamma'}\mathrm{Fil}_{u}(\mathbf{E}_{M})\subseteq\mathbf{E}_{M}\]
which is an $E$-subspace of $\mathbf{E}_{M}$ of dimension $\#\Gamma'$ (see the discussion below Condition~\ref{cond: intro abstract filtration}).
We define $\Gamma_0\defeq \Gamma^{[2,j_0]}\sqcup\{w_{1,j_0}\}$, $\Gamma_1\defeq \Gamma^{[j_1,n-2]}\sqcup\{w_{n-1,j_1}\}$ and $\Gamma_{2}\defeq \bigsqcup_{x\in\Gamma_{0},y\in\Gamma_{1}}\Gamma_{x,y}$ and note that they are closed subsets of $\Gamma^{J}$, $\Gamma^{J'}$ and $\Gamma^{J\sqcup J'}$ respectively.
We define
\[\Gamma_0^+\defeq \Gamma_0\cup \Gamma_{0,\dagger}=\Gamma^{[2,j_0]}\sqcup\{w_{1,j}\mid j\in J\},\]
\[\Gamma_1^+\defeq \Gamma_1\cup \Gamma_{1,\dagger}=\Gamma^{[j_1,n-2]}\sqcup\{w_{n-1,j}\mid j\in J'\},\]
and
\[\Gamma_2^+\defeq \Gamma_2\cup \Gamma_{2,\dagger}=\bigsqcup_{x\in\Gamma_{0}^+,y\in\Gamma_{1}^+}\Gamma_{x,y}.\]
Note that the subset $\Gamma_{0,\dagger}\subseteq\Gamma^{J}$ (resp.~$\Gamma_{1,\dagger}\subseteq\Gamma^{J'}$, resp.~$\Gamma_{2,\dagger}\subseteq\Gamma^{J\sqcup J'}$) is not closed, and its closure is contained in $\Gamma_0^+$ (resp.~in $\Gamma_1^+$, resp.~in $\Gamma_2^+$). We also have the simple but crucial observation that
\begin{equation}\label{equ: intro coxeter dagger intersection}
\left\{\begin{array}{ccc}
\Gamma_0\cap \Gamma_{0,\dagger} & = & \{w_{1,j_0}\}\\
\Gamma_1\cap \Gamma_{1,\dagger} & = & \{w_{n-1,j_1}\}\\
\Gamma_2\cap \Gamma_{2,\dagger} & = & \Gamma_{w_{1,j_0},w_{n-1,j_1}}
\end{array}\right.
\end{equation}
Now we define
\[
\left\{\begin{array}{ccc}
\mathbf{E}_{J}^{\heartsuit} & \defeq & \mathbf{E}_{0,\dagger}\cap\mathrm{Fil}_{\Gamma_0}(\mathbf{E}_{J})\\
\mathbf{E}_{J'}^{\heartsuit} & \defeq & \mathbf{E}_{1,\dagger}\cap\mathrm{Fil}_{\Gamma_1}(\mathbf{E}_{J'})\\
\mathbf{E}_{J,J'}^{\heartsuit} & \defeq & \mathbf{E}_{2,\dagger}\cap\mathrm{Fil}_{\Gamma_2}(\mathbf{E}_{J\sqcup J})\cap\mathbf{E}_{J\sqcup J'}^{<},
\end{array}\right.
\]
and then endow them with the filtration $\tau^{i}(-)$ induced from $\mathbf{E}_{J}$, $\mathbf{E}_{J'}$ and $\mathbf{E}_{J\sqcup J'}$ respectively.
Based on the crucial combinatorial observation (\ref{equ: intro coxeter dagger intersection}), we prove the following result.
\begin{prop}[Proposition~\ref{prop: cup of line}]\label{prop: intro relative coxeter intersection}
Let $J=[1,j_0]$ and $J'=[j_1,n-1]$ be as above. We have the following results.
\begin{enumerate}[label=(\roman*)]
\item \label{it: intro relative coxeter 1} We have $\mathrm{Fil}_{\Gamma_0^+}(\mathbf{E}_{J})=\mathbf{E}_{0,\dagger}+\mathrm{Fil}_{\Gamma_0}(\mathbf{E}_{J})$ and $\Dim_E\mathbf{E}_{J}^{\heartsuit}=1$. Moreover, we have $\mathrm{gr}^{i}(\mathbf{E}_{J}^{\heartsuit})\neq 0$ if and only if $i=0$, in which case it is isomorphic to $\mathrm{gr}_{w_{1,j_0}}(\mathbf{E}_{J})$.
\item \label{it: intro relative coxeter 2} We have $\mathrm{Fil}_{\Gamma_1^+}(\mathbf{E}_{J'})=\mathbf{E}_{1,\dagger}+\mathrm{Fil}_{\Gamma_1}(\mathbf{E}_{J'})$ and $\Dim_E\mathbf{E}_{J'}^{\heartsuit}=1$. Moreover, we have $\mathrm{gr}^{i}(\mathbf{E}_{J'}^{\heartsuit})\neq 0$ if and only if $i=0$, in which case it is isomorphic to $\mathrm{gr}_{w_{n-1,j_1}}(\mathbf{E}_{J'})$.
\item \label{it: intro relative coxeter 3} We have $\mathrm{Fil}_{\Gamma_2^+}(\mathbf{E}_{J\sqcup J'})=\mathbf{E}_{2,\dagger}+\mathrm{Fil}_{\Gamma_2}(\mathbf{E}_{J\sqcup J'})$ and $\Dim_E\mathbf{E}_{J,J'}^{\heartsuit}=1$.
\end{enumerate}
\end{prop}
In the proof of \ref{it: intro relative coxeter 1} of Proposition~\ref{prop: intro relative coxeter intersection} (and similarly for \ref{it: intro relative coxeter 2} and \ref{it: intro relative coxeter 3} of Proposition~\ref{prop: intro relative coxeter intersection}), apart from the usage of \ref{it: intro relative image 1} of Proposition~\ref{prop: intro relative image}, a key step is to show that 
\begin{equation}\label{equ: intro dagger in coxeter}
\mathbf{E}_{0,\dagger}\subseteq \mathrm{Fil}_{\Gamma_{0}^+}(\mathbf{E}_{J}).
\end{equation}
For each $x\in\Gamma_{0,\dagger}$, we write $\zeta_{x}:\mathbf{E}_{J}\twoheadrightarrow\mathbf{E}_{x,J}$ for the associated surjection (see \ref{it: intro coxeter map 2} of Proposition~\ref{prop: intro coxeter map}) and identify $\mathrm{gr}_{x}(\mathbf{E}_{J})$ with a $1$-dimensional $E$-subspace of $\mathbf{E}_{x,J}$. Following the observation that an element $u\in\Gamma^{J}$ satisfies $u\in\Gamma_0^{+}$ if and only if there does not exist $x\in\Gamma_{0,\dagger}$ such that $x\unlhd u$ and $x\neq u$, we have
\[\mathrm{Fil}_{\Gamma_{0}^+}(\mathbf{E}_{J})=\bigcap_{x\in\Gamma_{0,\dagger}}\zeta_{x}^{-1}(\mathrm{gr}_{x}(\mathbf{E}_{J})).\]
Hence, to prove (\ref{equ: intro dagger in coxeter}), it suffices to check that the composition of
\begin{equation}\label{equ: intro dagger to x}
\mathbf{E}_{0,\dagger}\hookrightarrow\mathbf{E}_{J}\buildrel \zeta_{x}\over\longrightarrow\mathbf{E}_{x,J}
\end{equation} 
has image contained in $\mathrm{gr}_{x}(\mathbf{E}_{J})$, for each $x\in\Gamma_{0,\dagger}$. In fact, the composition of (\ref{equ: intro dagger to x}) factors through $\mathbf{E}_{J,\fh_{[2,n-1]}}$ and thus further factors through $\mathbf{E}_{x,J,\fh_{[2,n-1]}}$ (using (\ref{equ: intro relative coxeter}) for $\fh=\fh_{[2,n-1]}$), and we actually prove that the natural map
\[\mathbf{E}_{x,J,\fh_{[2,n-1]}}\rightarrow \mathbf{E}_{x,J}\]
is an embedding with image $\mathrm{gr}_{x}(\mathbf{E}_{J})$ (see Proposition~\ref{prop: 0 1 2 sm part}).

Now we are finally ready to verify Condition~\ref{cond: intro cup heart} for the $1$-dimensional $E$-subspaces $\mathbf{E}_{J}^{\heartsuit}\subseteq \mathbf{E}_{J}$ and $\mathbf{E}_{J'}^{\heartsuit}\subseteq \mathbf{E}_{J'}$ constructed above.
First, using Theorem~\ref{thm: intro top graded}, we can deduce $\mathbf{E}_{J}^{\heartsuit}\not\subseteq\mathbf{E}_{J}^{<}$ (resp.~$\mathbf{E}_{J'}^{\heartsuit}\not\subseteq\mathbf{E}_{J'}^{<}$) from $\mathrm{gr}^{0}(\mathbf{E}_{J}^{\heartsuit})=\mathrm{gr}_{w_{1,j_0}}(\mathbf{E}_{J})$ (resp.~from $\mathrm{gr}^{0}(\mathbf{E}_{J'}^{\heartsuit})=\mathrm{gr}_{w_{n-1,j_1}}(\mathbf{E}_{J'})$). This verifies \ref{it: intro cup heart 1} of Condition~\ref{cond: intro cup heart}.
It follows from \ref{it: intro cup x y 2} of Theorem~\ref{thm: intro cup x y} and the definition of $\Gamma_{0}$, $\Gamma_1$ as well as $\Gamma_2$ that we have
\begin{equation}\label{equ: intro cup commute coxeter bound}
\kappa_{J,J'}(\mathrm{Fil}_{\Gamma_0}(\mathbf{E}_{J})\otimes_E\mathrm{Fil}_{\Gamma_1}(\mathbf{E}_{J'}))\subseteq \mathrm{Fil}_{\Gamma_2}(\mathbf{E}_{J\sqcup J'}).
\end{equation}
We also have
\begin{equation}\label{equ: intro cup commute relative bound}
\kappa_{J,J'}(\mathbf{E}_{0,\dagger}\otimes_E\mathbf{E}_{1,\dagger})\subseteq \mathbf{E}_{2,\dagger}
\end{equation}
from the commutative diagram (\ref{equ: intro cup relative diagram}). Taking the intersection of BHS of (\ref{equ: intro cup commute coxeter bound}) and (\ref{equ: intro cup commute relative bound}) and then using $\mathrm{im}(\kappa_{J,J'})\subseteq\mathbf{E}_{J\sqcup J'}^{<}$, we conclude that
\begin{equation}\label{equ: intro cup commute mix bound 1}
\kappa_{J,J'}(\mathbf{E}_{J}^{\heartsuit}\otimes_E\mathbf{E}_{J'}^{\heartsuit})\subseteq \mathbf{E}_{J,J'}^{\heartsuit}.
\end{equation}
Upon exchanging $J,J'$ in the argument above, we also have
\begin{equation}\label{equ: intro cup commute mix bound 2}
\kappa_{J',J}(\mathbf{E}_{J'}^{\heartsuit}\otimes_E\mathbf{E}_{J}^{\heartsuit})\subseteq \mathbf{E}_{J,J'}^{\heartsuit}.
\end{equation}
Now that both $\kappa_{J,J'}$ and $\kappa_{J',J}$ are injective by \ref{it: intro cup graded 2} of Theorem~\ref{thm: intro cup graded}, and that $\mathbf{E}_{J}^{\heartsuit}$, $\mathbf{E}_{J'}^{\heartsuit}$ and $\mathbf{E}_{J,J'}^{\heartsuit}$ are all $1$-dimensional by Proposition~\ref{prop: intro relative coxeter intersection}, we see that both (\ref{equ: intro cup commute mix bound 1}) and (\ref{equ: intro cup commute mix bound 2}) are actually equalities and thus finish the verification of \ref{it: intro cup heart 2} of Condition~\ref{cond: intro cup heart}
\subsection{Acknowledgements}
The author is very grateful to Christophe Breuil for his guidance into this area as well as the joint work \cite{BQ24}.
We thank Benjamin Schraen for his beautiful thesis \cite{Schr11} in the $\mathrm{GL}_3(\Q_p)$ case, and the initial motivation of the present work is to generalize \cite{Schr11} to $\mathrm{GL}_n$.
We thank Lennart Gehrmann for his work \cite{Geh21} (and discussions) which offers some first evidence of the results in this work.
We thank Yiwen Ding for extensive discussions around locally analytic representations and explaining basic results on $(\varphi,\Gamma)$-modules. His joint works \cite{BD20} and \cite{BD19} with Breuil are obviously crucial to the development of this work.
We thank Florian Herzig and Stefano Morra for their interest.
We thank Yiqin He, Benchao Su, Arnaud Vanhaecke, Zhixiang Wu and many others for their comments and helpful discussions.
The author wishes to thank University of Toronto and Academy of Mathematics and Systems Science for their hospitality.
\subsection{Notation}\label{subsec: notation}
In the whole text we fix a prime $p$ and an embedding $\iota: K\rightarrow E$ of two finite extensions of $\Q_p$.

For each locally $K$-analytic group $G$, we write $\Hom(G,E)$ for the $E$-vector space of locally $K$-analytic functions (via our fixed $\iota: K\hookrightarrow E$) $f:G\rightarrow E$ satisfying $f(xy)=f(x)+f(y)$ for each $x,y\in G$.
Let $\val: K^\times\rightarrow \Z$ be the $p$-adic valuation function normalized by $\val(p)=1$. We fix a choice of locally $K$-analytic logarithm $\plog: K^\times\rightarrow K$ which satisfies $\plog(p)=0$. By viewing both $\val$ and $\plog$ as $E$-valued locally $K$-analytic functions on $K^\times$, it is a well-known fact that $\{\val,\iota\circ\plog\}$ forms a basis of $\Hom(K^\times,E)$ (and we often abuse $\plog$ for $\iota\circ\plog$ when the choice of $\iota$ is fixed).

\begin{defn}\label{def: gp isogeny}
Let $G$ and $H$ be locally $K$-analytic groups.
We say a homomorphism $\varphi: H\rightarrow G$ between locally $K$-analytic groups is an \emph{isogeny} if $\mathrm{ker}(\varphi)$ is finite and $\mathrm{im}(\varphi)$ is an open normal subgroup of $H$ with finite index.
A \emph{rational homomorphism} from $H$ to $G$, written $\psi: H\dashrightarrow G$, is by definition a diagram of the following form
\[
\xymatrix{
H & G \ar^{\psi_3}[d]\\
H_1 \ar^{\psi_1}[u] \ar^{\psi_2}[r] & G_1
}
\]
with $\psi_2$ being a homomorphism and $\psi_1$, $\psi_3$ being isogeny between locally $K$-analytic groups, modulo the evident equivalence.
If furthermore $\psi_2$ is an isogeny, we say the rational homomorphism $\psi: H\dashrightarrow G$ is a \emph{quasi-isogeny}.
\end{defn}
We observe that any isogeny $\varphi: H\rightarrow G$ induces an isomorphism
\[\varphi^{\ast}: \Hom(G,E)\buildrel\sim\over\longrightarrow \Hom(H,E).\]
Hence, any rational homomorphism $\psi: H\dashrightarrow G$ induces a well-defined map
\[\psi^{\ast}: \Hom(G,E)\rightarrow \Hom(H,E),\]
which is an isomorphism when $\psi$ is a quasi-isogeny.

Let $\mathbf{G}$ be a split reductive group scheme over $K$. We fix a split maximal torus $\mathbf{T}\subseteq \mathbf{G}$ as well as a pair of opposite Borel subgroups $\mathbf{B},\mathbf{B}^+\subseteq \mathbf{G}$ satisfying $\mathbf{B}\cap\mathbf{B}^+=\mathbf{T}$.
We write $\Phi^+$ for the set of positive roots w.r.t $(\mathbf{B}^+,\mathbf{T})$ and $\Delta\subseteq \Phi^+$ the subset of positive simple roots.
We have natural partial-order on $\Z\Phi^+=\Z\Delta$ given as follows: for each $\al,\beta\in\Z\Phi^+$, we write $\al\geq \beta$ if $\al-\beta\in\Z_{\geq 0}\Phi^+=\Z_{\geq 0}\Delta$.
Each $I\subseteq \Delta$ gives rise to a Levi subgroup $\mathbf{L}_{I}\subseteq \mathbf{G}$ as well as a pair of parabolic subgroups $\mathbf{P}_{I}\defeq \mathbf{L}_{I}\mathbf{B}$ and $\mathbf{P}_{I}^+\defeq \mathbf{L}_{I}\mathbf{B}^+$. We write $\mathbf{N}_{I}$ (resp.~$\mathbf{N}_{I}^+$) for the unipotent radical of $\mathbf{P}_{I}$ (resp.~$\mathbf{P}_{I}^+$).
We write $\mathbf{U}_{I}$ (resp.~$\mathbf{U}_{I}^+$) for the unipotent radical of $\mathbf{B}_{I}\defeq \mathbf{B}\cap\mathbf{L}_{I}$ (resp.~of $\mathbf{B}_{I}^+\defeq \mathbf{B}^+\cap\mathbf{L}_{I}$). We write $\mathbf{Z}_{I}$ for the center of $\mathbf{L}_{I}$, and $\mathbf{H}_{I}$ for the closed subgroup scheme of $\mathbf{L}_{I}$ generated by $\mathbf{U}_{I}$ and $\mathbf{U}_{I}^+$. We write $G$, $P_{I}$, $L_{I}$, etc for the sets of $K$-points of $\mathbf{G}$, $\mathbf{P}_{I}$, $\mathbf{L}_{I}$, etc, which are naturally locally $K$-analytic groups.
We write $\mathrm{Lie}(G)$, $\mathrm{Lie}(P_{I})$, $\mathrm{Lie}(L_{I})$, etc for the $K$-Lie algebras associated with $G$, $P_{I}$, $L_{I}$, etc, and $\fg\defeq \mathrm{Lie}(G)\otimes_{K,\iota}E$, $\fp_{I}\defeq \mathrm{Lie}(P_{I})\otimes_{K,\iota}E$, $\fl_{I}\defeq \mathrm{Lie}(L_{I})\otimes_{K,\iota}E$, etc for the $E$-Lie algebras obtained by an extension of scalars via $\iota: K\rightarrow E$.
Note that the natural map $H_{I}\times Z_{I}\rightarrow L_{I}$ is an isogeny and induces an isomorphism $\fh_{I}\times \fz_{I}\buildrel\sim\over\longrightarrow \fl_{I}$.

We let $W(G)\defeq N_{\mathbf{G}}(\mathbf{T})/\mathbf{T}$ be the Weyl group, $W(L_I)$ be the Weyl group of the Levi $L_I$ corresponding to $I$, $\ell(w)\in \Z_{\geq 0}$ the length of $w\in W(G)$ and $W^{I',I}$ the set of minimal length representatives of $W(L_{I'})\backslash W(G)/W(L_I)$. We let $\rho_I$ be half the sum of the roots in $\fb_I$ (when $I=\emptyset$, we just write $\rho$). In particular we have $\langle\rho_I, \alpha^\vee\rangle =-1$ for $\alpha \in I$ (where $\alpha^\vee$ is the corresponding coroot and $\langle~,~\rangle$ the usual pairing, see \cite[\S 0.2]{Hum08}) and $\rho_I$ is dominant with respect to $\fb_I$. We define the dot action $w\cdot \mu \defeq w(\mu + \rho_I) - \rho_I$ for $w\in W(L_I)$ and $\mu\in X(T)$. Note that, since $\rho - \rho_I$ is invariant under $W(L_I)$, we also have $w\cdot \mu = w(\mu + \rho) - \rho$.

Assume from now that the center $\mathbf{Z}$ of $\mathbf{G}$ is trivial.
The natural maps $Z_{I}\rightarrow L_{I}\rightarrow L_{I}/H_{I}$ (with their composition being an isogeny) induce the following isomorphisms between $E$-vector spaces
\begin{equation}\label{equ: Levi char transfer}
\Hom(L_{I}/H_{I},E)\buildrel\sim\over\longrightarrow \Hom(L_{I},E) \buildrel\sim\over\longrightarrow \Hom(Z_{I},E)
\end{equation}
and thus
\begin{equation}\label{equ: Levi char transfer sm}
\Hom_{\infty}(L_{I}/H_{I},E)\buildrel\sim\over\longrightarrow \Hom_{\infty}(L_{I},E) \buildrel\sim\over\longrightarrow \Hom_{\infty}(Z_{I},E).
\end{equation}
Note that both (\ref{equ: Levi char transfer}) and (\ref{equ: Levi char transfer sm}) are functorial with respect to the choice of $I$, namely with respect to the inclusion $L_{I'}\subseteq L_{I}$ and the induced maps $L_{I'}/H_{I'}\rightarrow L_{I}/H_{I}$ as well as $Z_{I'}\dashrightarrow Z_{I}$ for each pair $I'\subseteq I$.

For each $\al\in\Delta$, we fix an isomorphism
\[\mathbf{Z}_{\Delta\setminus\{\al\}}\buildrel\sim\over\longrightarrow \bG_m\]
between groups schemes over $K$, which induces an isomorphism
\begin{equation}\label{equ: center maximal Levi}
\varepsilon_{\al}: Z_{\Delta\setminus\{\al\}}\buildrel\sim\over\longrightarrow K^\times
\end{equation}
between locally $K$-analytic groups.
We set $\val_{\al}\defeq \val\circ\varepsilon_{\al}$ and $\plog_{\al}\defeq \plog\circ\varepsilon_{\al}$, and note that $\{\val_{\al},\plog_{\al}\}$ (resp.~$\{\val_{\al}\}$) forms a basis of $\Hom(Z_{\Delta\setminus\{\al\}},E)$ (resp.~of $\Hom_{\infty}(Z_{\Delta\setminus\{\al\}},E)$).

Taking product of the rational maps $Z_{I}\dashrightarrow Z_{\Delta\setminus\{\al\}}$ for each $\al\in\Delta$, we obtain a quasi-isogeny
\begin{equation}\label{equ: center isogeny decomposition}
Z_{I}\dashrightarrow \prod_{\al\in\Delta\setminus I}Z_{\Delta\setminus\{\al\}}
\end{equation}
which induces the following isomorphism
\begin{equation}\label{equ: Levi char decomposition}
\prod_{\al\in\Delta\setminus I}\Hom(Z_{\Delta\setminus\{\al\}},E)\buildrel\sim\over\longrightarrow \Hom(Z_{I},E)
\end{equation}
between $E$-vector spaces. 
In particular
\[\mathbf{Log}_{I}\defeq \{\val_{\al},\plog_{\al}\mid \al\in\Delta\setminus I\}\]
forms a basis of (\ref{equ: Levi char transfer}), and
\[\mathbf{Log}_{I}^{\infty}\defeq \{\val_{\al}\mid \al\in\Delta\setminus I\}\subseteq \mathbf{Log}_{I}\]
forms a basis of (\ref{equ: Levi char transfer sm}).
We fix a total order on $\mathbf{Log}_{\emptyset}$ which induces a total order on $\mathbf{Log}_{I}\subseteq \mathbf{Log}_{\emptyset}$ for each $I\subseteq \Delta$. We consider a subset $v=\{x_1<\cdots<x_k\}\subseteq \mathbf{Log}_{I}$ (equipped with the total order inherited from $\mathbf{Log}_{I}$), then we have $x_1\wedge\cdots\wedge x_k\in \wedge^k\Hom(Z_{I},E)$. This way, upon fixing a total order on $\mathbf{Log}_{\emptyset}$, we may identify $\wedge^k\mathbf{Log}_{I}\defeq \{v\subseteq \mathbf{Log}_{I}\mid \#v=k\}$ with a basis of $\wedge^k\Hom(Z_{I},E)$.
Note that each $v\subseteq \mathbf{Log}_{I}$ determines a unique maximal $I_v\subseteq \Delta$ containing $I$ such that $v\subseteq\mathbf{Log}_{I_v}$. Throughout the text, we will abuse the notation $v$ for both a subset of $\mathbf{Log}_{I}\subseteq \mathbf{Log}_{\emptyset}$ and the corresponding element of $\wedge^{\#v}\Hom(Z_{I},E)\subseteq \wedge^{\#v}\Hom(T,E)$.

For each $I\subseteq \Delta$ and $v\in \mathbf{Log}_{I}^{\infty}$, we write $\mathrm{ker}(v)_{I}\defeq \{g\in L_{I}\mid v(g)=0\}$ which is an open subgroup of $L_{I}$ (as $v(L_{I})\subseteq E$ is discrete).
We define
\begin{equation}\label{equ: trivial val gp}
L_{I}'\defeq \bigcap_{v\in \mathbf{Log}_{I}^{\infty}}\mathrm{ker}(v)_{I}\subseteq L_{I}
\end{equation}
which is an open subgroup of $L_{I}$ with $L_{I}/L_{I}'$ being a free abelian group of rank $\#\Delta\setminus I$.
We choose an arbitrary discrete cocompact subgroup $T''$ of $T=Z_{\emptyset}$.
For each $I\subseteq \Delta$, we write $Z_{I}'$ (resp.~$Z_{I}''$) for the unique maximal compact open subgroup of $Z_{I}$ (resp.~for a finite index subgroup of $T''\cap Z_{I}$), and thus the natural map $Z_{I}'\times Z_{I}''\rightarrow Z_{I}$ is an isogeny.
We also have natural maps $H_{I}\times Z_{I}'\rightarrow L_{I}'$ and $L_{I'}\times Z_{I}''\rightarrow L_{I}$ which are both isogeny.
The quasi-isogeny (\ref{equ: center isogeny decomposition}) induces a quasi-isogeny
\begin{equation}\label{equ: center discrete isogeny decomposition}
Z_{I}''\dashrightarrow \prod_{\al\in\Delta\setminus I}Z_{\Delta\setminus\{\al\}}''
\end{equation}
for each $I\subseteq \Delta$. The isogeny $Z_{I}''\rightarrow Z_{I}/Z_{I}'$ induces an isomorphism
\[\Hom_{\infty}(Z_{I},E)\buildrel\sim\over\longrightarrow \Hom(Z_{I}'',E)\]
for each $I\subseteq \Delta$.

When $\mathbf{G}=\mathrm{PGL}_{n/K}$, we choose $\mathbf{T}$ (resp.~$\mathbf{B}$, resp.~$\mathbf{B}^+$) to be the diagonal maximal torus (resp.~the lower-triangular Borel, resp.~the upper-triangular Borel).
We abuse $i$ for its corresponding element $(i,i+1)\in\Delta$, for each $1\leq j\leq n-1$.
A subset $I\subseteq \Delta=\{1,\dots,n-1\}$ is called an \emph{interval} if either $I=\emptyset$ or $I=[i,j]\defeq \{j'\in\Delta\mid i\leq j'\leq j\}$ for some $1\leq i\leq j\leq n-1$.
For each $I\subseteq \Delta$, we write
\begin{equation}\label{equ: sum of roots}
\al_I\defeq\sum_{i\in I}(i,i+1)\in\Z \Phi^+.
\end{equation}
Note that each $\al\in\Phi^+$ determines a unique non-empty subinterval $I_\al\subseteq \Delta$ such that $\al=\al_{I_\al}$. Conversely, we have $\al_I\in\Phi^+$ for each non-empty subinterval $I\subseteq \Delta$.

Given $I\subseteq \Delta$ with $r_I\defeq n-\#I$, there exists a unique tuple of integers $n_d$ for $1\leq d\leq r_I$ such that $\sum_{d'=1}^{r_I}n_{d'}=n$ and
\begin{equation}\label{equ: explicit complement}
\Delta\setminus I=\{\sum_{d'=1}^{d}n_{d'}\mid 1\leq d\leq r_I-1\}.
\end{equation}
For each $1\leq d\leq r_I$, we define 
\[
I^d\defeq \{i\mid \sum_{d'=1}^{d-1}n_{d'}< i<\sum_{d'=1}^{d}n_{d'}\}
\]
with $n_d=1+\#I^d$. Note from (\ref{equ: explicit complement}) that we have a partition $I=\bigsqcup_{d=1}^{r_I}I^d$, with $I^d$ being either empty or a maximal subinterval of $I$ for each $1\leq d\leq r_I$.
It is clear that $i<i'$ for each $i\in I^d$ and $i'\in I^{d'}$ satisfying $1\leq d<d'\leq r_I$.

For an abelian category $\cA$ and a finite length object $D\in\cA$, we write $\mathrm{JH}_{\cA}(D)$ for the set of (isomorphism classes of) constituents of $D$. If $D$ is multiplicity free, then we equip $\mathrm{JH}_{\cA}(D)$ with the following natural partial-order: two different constituents $C_1,C_2\in\mathrm{JH}_{\cA}(D)$ satisfies $C_1<C_2$ if $C_1$ shows up in the unique subobject of $D$ with cosocle $C_2$.

Let $\cC$ be an additive category in which finite direct sum exists. 
Let $\cT^{\bullet,\bullet'}$ be a bounded $\cC$-valued double complex, and we number the index $\bullet$ (resp.~$\bullet'$) by $1$ (resp.~by $2$) for convenience.
We thus write $d_{1}^{\bullet,\bullet'}$ and $d_{2}^{\bullet,\bullet'}$ for the first and second differential maps respectively.
Given the ordered tuple $\ast=(1,2)$ (resp.~$\ast=(2,1)$), we define its associated total complex $\mathrm{Tot}_{\ast}(\cT^{\bullet,\bullet'})$ by
\[\mathrm{Tot}_{\ast}(\cT^{\bullet,\bullet'})^{m}\defeq \bigoplus_{k+k'=m}\cT^{k,k'},\]
with the differential map
\[\mathrm{Tot}_{\ast}(\cT^{\bullet,\bullet'})^{m}\rightarrow \mathrm{Tot}_{\ast}(\cT^{\bullet,\bullet'})^{m+1}\]
given by $d_{1}^{k,k'}+(-1)^{k}d_{2}^{k,k'}$ (resp.~by $(-1)^{k'}d_{1}^{k,k'}+d_{2}^{k,k'}$) on the direct summand $\cT^{k,k'}$ of $\mathrm{Tot}_{\ast}(\cT^{\bullet,\bullet'})^{m}$, for each $k$, $k'$ and $m$ satisfying $k+k'=m$.
Note that the scalar automorphism
\[(-1)^{kk'}: \cT^{k,k'}\buildrel\sim\over\longrightarrow \cT^{k,k'}\]
for each $k,k'\in\Z$ induces an isomorphism
\begin{equation}\label{equ: double total swap}
\mathrm{Tot}_{(1,2)}(\cT^{\bullet,\bullet'})\buildrel\sim\over\longrightarrow \mathrm{Tot}_{(2,1)}(\cT^{\bullet,\bullet'})
\end{equation}
between complex. We call (\ref{equ: double total swap}) the \emph{standard isomorphism} between $\mathrm{Tot}_{\ast}(\cT^{\bullet,\bullet'})$ for different choices of $\ast$.

More generally, we consider a $\cC$-valued $N$-complex $\cT^{\bullet,\dots,\bullet}$ with its indices numbered by $1,\dots,N$ for some $N\geq 2$ and $i$-th differential map denoted by $d_{i}^{\bullet,\dots,\bullet}$ for each $1\leq i\leq N$.
We consider a permutation $\sigma\in\mathfrak{S}_{N}$ and view it as an ordered tuple of integers as $(\sigma(1),\cdots,\sigma(N))$.
We define the \emph{$\sigma$-twisted total complex} $\mathrm{Tot}_{\sigma}(\cT^{\bullet,\dots,\bullet})$ of $\cT^{\bullet,\dots,\bullet}$ as
\[\mathrm{Tot}_{\sigma}(\cT^{\bullet,\dots,\bullet})^{m}\defeq \bigoplus_{\sum_{i=1}^{N}k_i=m}\cT^{k_1,\dots,k_N}\]
with the differential map given by
\[\mathrm{Tot}_{\sigma}(\cT^{\bullet,\dots,\bullet})^{m}\rightarrow \mathrm{Tot}_{\sigma}(\cT^{\bullet,\dots,\bullet})^{m+1}\]
given by
\begin{equation}\label{equ: general total differential}
\sum_{i=1}^{N}(-1)^{\sum_{j=1}^{i-1}k_{\sigma(j)}}d_{\sigma(i)}^{k_1,\dots,k_N}
\end{equation}
on the direct summand $\cT^{k_1,\dots,k_N}$ of $\mathrm{Tot}_{\sigma}(\cT^{\bullet,\dots,\bullet})^{m}$, for each tuple of integers $k_1,\dots,k_N$ satisfying $\sum_{i=1}^{N}k_i=m$. 
We can interpret $\mathrm{Tot}_{\sigma}(\cT^{\bullet,\dots,\bullet})$ as successive partial total complex of $\cT^{\bullet,\dots,\bullet}$. For example, we have natural identifications
\[\mathrm{Tot}_{-,\sigma(3)}(\mathrm{Tot}_{\sigma(1),\sigma(2)}(\cT^{\bullet,\bullet,\bullet}))=\mathrm{Tot}_{\sigma}(\cT^{\bullet,\bullet,\bullet})=\mathrm{Tot}_{\sigma(1),-}(\mathrm{Tot}_{\sigma(2),\sigma(3)}(\cT^{\bullet,\bullet,\bullet}))\]
when $N=3$, and similar facts hold for general $N$.
Now we consider some $1\leq i_0\leq N-1$ and $\sigma'\defeq \sigma\cdot(i_0,i_0+1)$, with the ordered tuple associated with $\sigma'$ being
\[(\sigma'(1),\dots,\sigma'(N))=(\sigma(1),\dots,\sigma(i_0-1),\sigma(i_0+1),\sigma(i_0),\sigma(i_0+2),\dots,\sigma(N)).\]
It is then easy to check that the scalar automorphism
\[(-1)^{k_{i_0}k_{i_0+1}}: \cT^{k_1,\dots,k_N}\buildrel\sim\over\longrightarrow \cT^{k_1,\dots,k_N}\]
for each $k_1,\dots,k_N\in\Z$ induces an isomorphism
\[
\zeta_{\sigma,\sigma'}: \mathrm{Tot}_{\sigma}(\cT^{\bullet,\dots,\bullet})\buildrel\sim\over\longrightarrow \mathrm{Tot}_{\sigma'}(\cT^{\bullet,\dots,\bullet})
\]
between complex. Note that the isomorphism $\zeta_{\sigma,\sigma'}$ constructed above restricts to the identity map on $\cT^{0,\dots,0}$, and in fact this is the unique such isomorphism which is functorial with respect to the choice of arbitrary $\cC$ and $\cT^{\bullet,\dots,\bullet}$.

Given a general pair of elements $\sigma,\sigma'\in\mathfrak{S}_{N}$, there exists a sequence of elements 
\begin{equation}\label{equ: total sigma seq}
\sigma=\sigma_0,\dots,\sigma_{t}=\sigma' 
\end{equation}
such that $\sigma_{t'-1}\sigma_{t'}$ is a simple reflection of the form $(i_{t'},i_{t'}+1)$ for some $1\leq i_{t'}\leq N-1$ and for each $1\leq t'\leq t$. We thus obtain an isomorphism
\[\zeta_{\sigma,\sigma'}\defeq \zeta_{\sigma_{t-1},\sigma_{t}}\circ\cdots\circ \zeta_{\sigma_{1},\sigma_{0}}.\]
This $\zeta_{\sigma,\sigma'}$ for our general pair $\sigma,\sigma'$ is again the unique isomorphism between complex which restricts to the identity map on $\cT^{0,\dots,0}$ and is functorial with respect to the choice of arbitrary $\cC$ and $\cT^{\bullet,\dots,\bullet}$. In particular, $\zeta_{\sigma,\sigma'}$ is independent of the choice of the sequence (\ref{equ: total sigma seq}). 
To summarize, we have constructed a \emph{standard isomorphism} $\zeta_{\sigma,\sigma'}$ for each pair of elements $\sigma,\sigma'\in\mathfrak{S}_{N}$, such that $\zeta_{\sigma,\sigma}$ is the identity map for each $\sigma\in\mathfrak{S}_{N}$, and we have
\[\zeta_{\sigma,\sigma''}=\zeta_{\sigma',\sigma''}\circ\zeta_{\sigma,\sigma'}\]
for each triple of elements $\sigma,\sigma',\sigma''\in \mathfrak{S}_{N}$. 
When $\sigma=1$, we often use the shortened notation 
\[\mathrm{Tot}(\cT^{\bullet,\dots,\bullet})\defeq \mathrm{Tot}_{1}(\cT^{\bullet,\dots,\bullet}).\]

\section{Preliminary on extensions between locally analytic representations}\label{sec: preliminary}
We collect some useful results on $U(\fg)$-modules (\S \ref{subsec: g category O}), smooth representation theory (\S \ref{subsec: sm rep}), locally analytic representation theory (\S \ref{subsec: loc an rep}) and locally analytic group cohomologies (\S \ref{subsec: gp coh}). In \S \ref{subsec: OS} we recall from \cite{BQ24} an important spectral sequence (Proposition~\ref{prop: St key seq}) and its consequences, followed by further application to the study of locally analytic generalized Steinberg representations.
\subsection{Partial-Coxeter elements and partition of positive roots}\label{subsec: coxeter partition}
We prove some basic results on the combinatorics of partial-Coxeter elements and partition of elements in $\Phi^+$.

For each $w\in W(G)$, we write $\mathrm{Supp}(w)\subseteq \Delta$ as the set of $\al\in\Delta$ such that the simple reflection $s_{\al}$ appears in one (or equivalently all) reduced decomposition of $w$.
We consider the following set of \emph{partial-Coxeter} elements in $W(G)$
\begin{equation}\label{equ: partial coxeter}
\Gamma\defeq \{w\in W(G)\mid \ell(w)=\#\mathrm{Supp}(w)\}.
\end{equation}
For each $I\subseteq \Delta$, we define
\begin{equation}\label{equ: coxeter element support}
\Gamma_I\defeq \{w\in\Gamma\mid \mathrm{Supp}(w)=I\}
\end{equation}
and
\begin{equation}\label{equ: coxeter upper bound}
\Gamma^I\defeq \{w\in\Gamma\mid \mathrm{Supp}(w)\subseteq I\}=\bigsqcup_{I'\subseteq I}\Gamma_{I'}.
\end{equation}
For each $I'\subseteq I$, we have a projection map
\begin{equation}\label{equ: coxeter projection}
\Gamma_{I}\rightarrow\Gamma_{I'}
\end{equation}
that sends $w\in\Gamma_{I}$ to the unique $w'\in\Gamma_{I'}$ that satisfies $w'\leq w$ (cf.~\cite[Thm.~2.2.2]{BB05}).
For each $I,J\subseteq \Delta$, we also define
\begin{equation}\label{equ: coxeter left}
\Gamma_{I}(J)\defeq \{w\in\Gamma_{I}\mid D_L(w)\subseteq I\setminus J\}.
\end{equation}

Given $x,w\in W(G)$, we write $x\unlhd w$ if there exists $y\in W(G)$ such that $w=yx$ and $\ell(w)=\ell(y)+\ell(x)$. This defines a partial-order on $W(G)$ which restricts to a partial-order on $\Gamma$.
When $x,w\in\Gamma$, we see that $x\unlhd w$ if and only if $w=yx$ for some $y\in\Gamma$ satisfying $\mathrm{Supp}(w)=\mathrm{Supp}(y)\sqcup\mathrm{Supp}(x)$.

For each $I\subseteq \Delta$, we write $\cS_{I}$ for the set of all subsets $S\subseteq \Phi^+$ satisfying $\sum_{\beta\in S}\beta=\sum_{\al\in I}\al\in\Z\Phi^+$.
For each $I'\subseteq I$, we have a projection map
\begin{equation}\label{equ: partition projection}
\cS_{I}\rightarrow\cS_{I'}
\end{equation}
which sends $S\in\cS_{I}$ to the unique $S'\in\cS_{I'}$ satisfying the following condition: for each $\beta'\in S'$ there exists a unique $\beta\in S$ such that $I_{\beta'}$ is a connected component of $I_{\beta}\cap I'$.
For each $I,J\subseteq \Delta$, we write $\cS_{I}^{J}\subseteq \cS_{I}$ for the subset consisting of those $S$ that satisfies $I_{\beta}\not\subseteq J$ for each $\beta\in S$.

For each $I,I'\subseteq \Delta$ with $I\cap I'=\emptyset$, we have a product map
\begin{equation}\label{equ: partition product}
\cS_{I}\times \cS_{I'}\rightarrow \cS_{I\sqcup I'}
\end{equation}
which sends $(S,S')$ to $S\sqcup S'\in \cS_{I\sqcup I'}$ for each $S\in\cS_{I}$ and $S'\in\cS_{I'}$.
It is clear that (\ref{equ: partition product}) is an injection.

We write $w_0\in W(G)$ for the longest Weyl element. The standard outer involution $G\rightarrow G: A\mapsto \overline{A}\defeq w_0(A^{t})^{-1}w_0$ induces an involution $W(G)\rightarrow W(G): w\mapsto \overline{w}\defeq w_0ww_0$ as well as an involution $\Phi^+\rightarrow \Phi^+: \al\mapsto\overline{\al}\defeq -w_0(\al)$. The aforementioned involution of $W(G)$ restricts to bijections 
\begin{equation}\label{equ: coxeter involution}
\Gamma^{I}\buildrel\sim\over\longrightarrow\Gamma^{w_0(I)}: w\mapsto \overline{w}
\end{equation}
and $\Gamma_{I}\buildrel\sim\over\longrightarrow \Gamma_{w_0(I)}$ for each $I\subseteq\Delta$, while the aforementioned involution of $\Phi^+$ induces a bijection 
\begin{equation}\label{equ: partition involution}
\cS_{I}\buildrel\sim\over\longrightarrow \cS_{w_0(I)}: S\mapsto\overline{S}\defeq\{\overline{\al}\mid \al\in S\}.
\end{equation}

\begin{lem}\label{lem: coxeter partition interval}
Let $I\subseteq \Delta$ be a non-empty interval. Then we have
\begin{equation}\label{equ: coxeter partition interval}
\#\Gamma_I=2^{\#I-1}=\#\cS_I.
\end{equation}
\end{lem}
\begin{proof}
We write $I=[i,i']$ for some $1\leq i\leq i'\leq n-1$ and set $I^{-}\defeq I\setminus\{i\}$.
We write $\cP(I^{-})$ for the set of all subsets of $I^{-}$.
For each $\beta\in\Phi^+$ we write $\beta=(i_{\beta},i_{\beta}')$ for some $1\leq i_{\beta}<i_{\beta}'\leq n$.
On one hand, the map
\begin{equation}\label{equ: partition to index}
\cS_{I}\rightarrow \cP(I^{-}): S\mapsto\{i_{\beta}\mid \beta\in S, i_{\beta}\neq i\}
\end{equation}
is clearly a bijection.
On the other hand, the map
\begin{equation}\label{equ: coxeter to index}
\Gamma_I\rightarrow\cP(I^{-}): w\mapsto\{j\in I^{-}\mid s_{j-1}s_j\leq w\}
\end{equation}
is again a bijection by \cite[Thm.~2.2.2]{BB05}.
We conclude (\ref{equ: coxeter partition interval}) from the bijections (\ref{equ: partition to index}) and (\ref{equ: coxeter to index}) together with the fact that $\#\cP(I^{-})=2^{\#I^{-}}=2^{\#I-1}$.
\end{proof}

\begin{prop}\label{prop: coxeter partition cardinality}
Let $I,J\subseteq \Delta$. We have the following results.
\begin{enumerate}[label=(\roman*)]
\item \label{it: coxeter partition 1} We have $\#\Gamma_{I}=\#\cS_{I}$.
\item \label{it: coxeter partition 2} We have $\#\Gamma_{I}(J)=\#\cS_{I}^{J}$.
\end{enumerate}
\end{prop}
\begin{proof}
We write $I=\bigsqcup_{t'=1}^{t}I_{t'}$ for the decomposition of $I$ into its connected components.
It is clear that the projection map $\Gamma_{I}\rightarrow \Gamma_{I_{t'}}$ (see (\ref{equ: coxeter projection})) and the projection map $\cS_{I}\rightarrow \cS_{I_{t'}}$ (see (\ref{equ: partition projection})) for each $1\leq t'\leq t$ induce a bijection 
\begin{equation}\label{equ: component Gamma bijection}
\Gamma_{I}\buildrel\sim\over\longrightarrow\prod_{t'=1}^{t}\Gamma_{I_{t'}}
\end{equation} 
and a bijection 
\begin{equation}\label{equ: component S bijection}
\cS_{I}\buildrel\sim\over\longrightarrow\prod_{t'=1}^{t}\cS_{I_{t'}}, 
\end{equation}
which together with Lemma~\ref{lem: coxeter partition interval} finishes the proof of \ref{it: coxeter partition 1}.

Now we prove \ref{it: coxeter partition 2}.\\
If we write $I=\bigsqcup_{t'=1}^{t}I_{t'}$ for the decomposition of $I$ into its connected components, then (\ref{equ: component Gamma bijection}) restricts to a bijection
\[\Gamma_{I}(J)\buildrel\sim\over\longrightarrow\prod_{t'=1}^{t}\Gamma_{I_{t'}}(J),\]
and (\ref{equ: component S bijection}) restricts to a bijection
\[\cS_{I}^{J}\buildrel\sim\over\longrightarrow\prod_{t'=1}^{t}\cS_{I_{t'}}^{J}.\]
Hence, it is harmless to assume in the rest of the proof that $I$ is connected (namely an interval of the form $[j_{I},j_{I}']$).
Our fixed $I,J\subseteq \Delta$ determines a unique subset $\Sigma_{I,J}\subseteq \Phi_{I}^+$ such that $\al\in\Sigma_{I,J}$ if and only if $I_{\al}$ is a connected component of $I\cap J$.
For each $\al=(i_{\al},i_{\al}')\in\Phi_{I}^+$, we write $\tld{\al}=(i_{\tld{\al}},i_{\tld{\al}}')\in\Phi_{I}^+$ for the unique root that satisfies $i_{\tld{\al}}=\max\{i_{\al}-1,j_{I}\}$ and $i_{\tld{\al}}'=\min\{i_{\al}'+1,j_{I}'+1\}$.
We write $\Phi^+_{I,\flat}\subseteq \Phi_{I}^+$ for the subset consisting of those $\al=(i_{\al},i_{\al}')\in\Phi_{I}^+$ that satisfies $j_{I}<i_{\al}<i_{\al}'<j_{I}'+1$, and then write $\Sigma_{I,J,\flat}\defeq\Sigma_{I,J}\cap\Phi^+_{I,\flat}$ for short.
For each $w\in \Gamma_{I}(J)$ and $\al\in\Sigma_{I,J}$, we write $w_{\al}\in\Gamma_{I_{\tld{\al}}}$ for the image of $w$ under the projection map $\Gamma_{I}\rightarrow \Gamma_{I_{\tld{\al}}}$ (see (\ref{equ: coxeter projection})).
For each $\al\in\Sigma_{I,J}$, we write $\Gamma(J)_{I,\al}\subseteq \Gamma_{I_{\tld{\al}}}$ for the subset consisting of those $x$ that satisfies $D_L(x)\cap I_{\al}=\emptyset$. The definition of $\Gamma_{I}(J)$ ensures that the projection map $\Gamma_{I}\rightarrow \Gamma_{I_{\tld{\al}}}$ restricts to a well-defined map
\begin{equation}\label{equ: I coxeter projection 1}
\Gamma_{I}(J)\rightarrow \Gamma(J)_{I,\al}.
\end{equation}
Note that the projection map $\Gamma_{I}\rightarrow \Gamma_{I\setminus J}$ restricts to map
\begin{equation}\label{equ: I coxeter projection 2}
\Gamma_{I}(J)\rightarrow \Gamma_{I\setminus J}
\end{equation}
and we write $w_{I\setminus J}$ for the image of $w\in\Gamma_{I}(J)$ under (\ref{equ: I coxeter projection 2}).
Using (\ref{equ: I coxeter projection 1}) for each $\al\in\Sigma_{I,J}$ as well as (\ref{equ: I coxeter projection 2}), we obtain a map
\begin{equation}\label{equ: I coxeter glue}
\Gamma_{I}(J)\rightarrow \Gamma_{I\setminus J}\times\prod_{\al\in\Sigma_{I,J}}\Gamma(J)_{I,\al}: w\mapsto (w_{I\setminus J}, w_{\al}\mid \al\in \Sigma_{I,J})
\end{equation}
which is a bijection by \cite[Thm.~2.2.2]{BB05}.
Note that $\Gamma(J)_{I,\al}$ consists of exactly one element when either $i_{\tld{\al}}=i_{\al}=j_{I}$ or $i_{\tld{\al}}'=i_{\al}'=j_{I}'+1$. When $j_{I}<i_{\al}<i_{\al}'<j_{I}'+1$ namely $\al\in\Sigma_{I,J,\flat}$, each $w_{\al}\in\Gamma(J)_{I,\al}$ is uniquely determined by $D_{R}(w_{\al})$ which is a single element of $I_{\tld{\al}}$. In other words, we have
\begin{equation}\label{equ: I coxeter al part}
\#\Gamma(J)_{I,\al}=\#I_{\tld{\al}}=\#I_{\al}+2
\end{equation}
for each $\al\in\Sigma_{I,J,\flat}$.
For each $\al\in\Sigma_{I,J}$, we define $\cS^{J}_{I,\al}\subseteq \cS_{I_{\tld{\al}}}$ as the subset consisting of those $S$ that satisfy $I_{\beta}\not\subseteq I_{\al}$ for each $\beta\in S$.
When either $i_{\tld{\al}}=i_{\al}=j_{I}$ or $i_{\tld{\al}}'=i_{\al}'=j_{I}'+1$, we note that $\cS^{J}_{I,\al}=\{\tld{\al}\}$ consists of exactly one element. When $j_{I}<i_{\al}<i_{\al}'<j_{I}'+1$ namely $\al\in\Sigma_{I,J,\flat}$, $\cS^{J}_{I,\al}$ coincides with the subset of $\cS_{I_{\tld{\al}}}$ consisting of those $S$ satisfying $\#S\leq 2$. In particular, we have
\begin{equation}\label{equ: I partition al part}
\#\cS^{J}_{I,\al}=\#I_{\tld{\al}}=\#I_{\al}+2
\end{equation}
for each $\al\in\Sigma_{I,J,\flat}$. The definition of $\cS_{I}^{J}$ ensures that the projection map $\cS_{I}\rightarrow \cS_{I_{\tld{\al}}}$ restricts to a map
\begin{equation}\label{equ: I partition projection 1}
\cS_{I}^{J}\rightarrow \cS^{J}_{I,\al}
\end{equation}
for each $\al\in\Sigma_{I,J}$.
Note that the projection map $\cS_{I}\rightarrow \cS_{I\setminus J}$ restricts to map
\begin{equation}\label{equ: I partition projection 2}
\cS_{I}^{J}\rightarrow \cS_{I\setminus J}.
\end{equation}
Using (\ref{equ: I partition projection 1}) for each $\al\in\Sigma_{I,J}$ as well as (\ref{equ: I partition projection 2}), we obtain a map
\begin{equation}\label{equ: I partition glue}
\cS_{I}^{J}\rightarrow \cS_{I\setminus J}\times\prod_{\al\in\Sigma_{I,J}}\cS^{J}_{I,\al}
\end{equation}
which is also a bijection.
Combining the bijections (\ref{equ: I coxeter glue}) and (\ref{equ: I partition glue}) with (\ref{equ: I coxeter al part}), (\ref{equ: I partition al part}) and \ref{it: coxeter partition 1} (with $I\setminus J$ replacing $I$ in \emph{loc.cit.}), we conclude that
\[\#\cS_{I}^{J}=\#\cS_{I\setminus J}\prod_{\al\in\Sigma_{I,J,\flat}}\#I_{\tld{\al}}=\#\Gamma_{I\setminus J}\prod_{\al\in\Sigma_{I,J,\flat}}\#I_{\tld{\al}}=\#\Gamma_{I}(J)\]
and thus finish the proof of \ref{it: coxeter partition 2}.
\end{proof}

\begin{rem}\label{rem: atom partition}
Let $I\subseteq \Delta$. We write $I^{-}\defeq \{i\in I\mid i-1\in I\}$ for short. For example, we have $\Delta^{-}=[2,n-1]\subseteq \Delta=[1,n-1]$. Using the bijection (\ref{equ: partition to index}) by replacing $I$ in \emph{loc.cit.} with each connected component of our $I$ here, we obtain a bijection
\begin{equation}\label{equ: atom partition bijection}
\cS_{I}\buildrel\sim\over\longrightarrow \cP(I^{-})
\end{equation}
with $\cP(I^{-})$ being the set of all subsets of $I^{-}$.
\end{rem}

\begin{defn}\label{def: disconnected pair}
Let $I,I'\subseteq \Delta$. We have the following notions.
\begin{enumerate}[label=(\roman*)]
\item \label{it: disconnected 0} We say that $I$ is an \emph{interval} if either $I=\emptyset$ or there exists $j,j'\in\Delta$ such that $I=[j,j']\defeq\{j,j+1,\dots,j'\}$.
\item \label{it: disconnected 1} We say $J\subseteq I$ is a \emph{connected component} of $I$ if it is a non-empty maximal subinterval.
\item \label{it: disconnected 2} We say that $I$ is \emph{smaller than} $I'$, written $I<I'$, if we have
\[\max\{i\mid i\in I\}<\min\{i'\mid i'\in I'\}.\]
We write $\delta(I,I')\in\Z_{\geq 0}$ for the number of pairs $(J,J')$ satisfying $J>J'$, with $J$ (resp.~$J'$) being a connected component of $I$ (resp.~$I'$). In particular, $\delta(I,I')=0$ if and only if $I<I'$.
\item \label{it: disconnected 3} We write $I<<I'$ if there exists $i\in\Delta$ such that $I<\{i\}<I'$. We say that \emph{$I$ and $I'$ do not connect} if $|i-i'|\geq 2$ for any $i\in I$ and $i'\in I'$. In particular, $I$ and $I'$ do not connect if and only if $I\cap I'=\emptyset$ and $\al+\al'$ is not a root for any $\al\in I$ and $\al'\in I'$, if and only if for each choice of connected components $J$ of $I$ and $J'$ of $I'$ we have either $J<<J'$ or $J'<<J$.
\end{enumerate}
\end{defn}

Let $x\in\Gamma$. We say $x'\in\Gamma$ is a \emph{connected component} of $x$ if $x'\leq x$ and $\mathrm{Supp}(x')$ is a connected component of $\mathrm{Supp}(x)$ in the sense of \ref{it: disconnected 1} of Definition~\ref{def: disconnected pair}. Thanks to \cite[Thm.~2.2.2]{BB05}, we know that, given $x\in\Gamma$, the set of connected components of $x$ is in natural bijection with the set of connected components of $\mathrm{Supp}(x)$. Moreover, we know by the usual braided relations that distinct connected components of $x$ commute with each other.

Let $x,y\in\Gamma$ with $\mathrm{Supp}(x)\cap\mathrm{Supp}(y)=\emptyset$.
We define
\begin{equation}\label{equ: x y envelop}
\Gamma_{x,y}\defeq \{w\in\Gamma\mid \mathrm{Supp}(w)=\mathrm{Supp}(x)\sqcup\mathrm{Supp}(y), x\leq w, y\leq w\}.
\end{equation}
It is clear that $\Gamma_{x,y}=\Gamma_{y,x}$ depends only on the set $\{x,y\}$.
More generally, given a set $\{x_1,\dots,x_m\}$ of elements in $\Gamma$ with pair wise distinct support, we can similarly define
\begin{equation}\label{equ: tuple envelop}
\Gamma_{x_1,\dots,x_m}\defeq \{w\in\Gamma\mid \mathrm{Supp}(w)=\bigsqcup_{k=1}^{m}\mathrm{Supp}(x_k), x_k\leq w \ { for \ each \ } 1\leq k\leq m\}.
\end{equation}
\begin{lem}\label{lem: x y envelop partial order}
Let $x,y\in\Gamma$ with $\mathrm{Supp}(x)\cap\mathrm{Supp}(y)=\emptyset$. Then for each $x',y'\in\Gamma$ satisfying $x'\unlhd x$ and $y'\unlhd y$ and $w'\in\Gamma_{x',y'}$, there exists $w\in\Gamma_{x,y}$ such that $w'\unlhd w$.
\end{lem}
\begin{proof}
Since $x'\unlhd x$ and $y'\unlhd y$, there exists $x'',y''\in\Gamma$ such that $x=x''x'$ with $\mathrm{Supp}(x)=\mathrm{Supp}(x')\sqcup\mathrm{Supp}(x'')$ and $y=y''y'$ with $\mathrm{Supp}(y)=\mathrm{Supp}(y')\sqcup\mathrm{Supp}(y'')$. It is clear that $x'$, $x''$, $y'$ and $y''$ have pair wise disjoint support, and thus $w\defeq x''y''w'$ satisfies 
\begin{multline*}
\mathrm{Supp}(w)=\mathrm{Supp}(x'')\sqcup\mathrm{Supp}(y'')\sqcup\mathrm{Supp}(w')\\
=\mathrm{Supp}(x'')\sqcup\mathrm{Supp}(y'')\sqcup\mathrm{Supp}(x')\sqcup\mathrm{Supp}(y')=\mathrm{Supp}(x)\sqcup\mathrm{Supp}(y)
\end{multline*}
and thus $w'\unlhd w$. Since $x'\leq w'$ and $y'\leq w'$, we have $x=x''x'\leq x''y''w'$ and $y=y''y'\leq x''y''w'$ by \cite[Thm.~2.2.2]{BB05}. In other words, we have checked that $w=x''y''w'$ satisfies $w'\unlhd w$ and $w\in\Gamma_{x,y}$. The proof is thus finished.
\end{proof}

\begin{lem}\label{lem: interval Bruhat order}
Let $I_0,I_0'\subseteq \Delta$ be non-empty subsets that satisfy $I_0<I_0'$. Let $x\in\Gamma_{I_0}$, $y\in\Gamma_{I_0'}$ and $w\in\Gamma_{I_0\sqcup I_0'}$ with $x,y\leq w$. Then $w$ must be either $xy$ or $yx$, with $xy=yx$ if and only if $I_0<<I_0'$.
\end{lem}
\begin{proof}
Note that $I_0$ contains a unique connected component $I$ such that $(I_0\setminus I)<<I$, and that $I_0'$ contains a unique connected component $I'$ such that $I'<<(I_0'\setminus I')$. Thanks to \cite[Thm.~2.2.2]{BB05} and braided relations, there exist unique $x'\in\Gamma_{I}$, $x''\in\Gamma_{I_0\setminus I}$, $y'\in\Gamma_{I'}$ and $y''\in\Gamma_{I_0'\setminus I'}$ such that $x=x''x'=x'x''$, $y=y'y''=y''y'$, with moreover $x''$ commuting with $y'$ and $y''$, and $y''$ commuting with $x'$ and $x''$.
If $I_0<<I_0'$ or equivalently $I<<I'$, then $x'$ commutes with $y'$ and we have $xy=yx$ being the unique choice of $w\in\Gamma_{I_0\sqcup I_0'}$ with $x,y\leq w$ by \cite[Thm.~2.2.2]{BB05} and braided relations. Assume from now that $I_0<<I_0'$ fails, or equivalently $I\sqcup I'$ is an interval and thus there exists $j,j_1,j_2\in\Delta$ such that $j_1\leq j< j_2$, $I=[j_1,j]$ and $I'=[j+1,j_2]$. Given $j'\in[j_1,j]$ and $j''\in[j+1,j_2]$, it is clear that $s_{j'}s_{j''}\neq s_{j''}s_{j'}$ if and only if $j'=j$ and $j''=j+1$.
For each choice of $w\in\Gamma_{I_0\sqcup I_0'}$ with $x,y\leq w$, it is clear from \cite[Thm.~2.2.2]{BB05} and braided relations that there exists $w'\in\Gamma_{I\sqcup I'}$ with $x',y'\leq w'$ such that $w=x''w'y''$ with $x''$, $w'$ and $y''$ commuting with each other.
By \cite[Thm.~2.2.2]{BB05} we know that either $s_js_{j+1}\leq w'$ or $s_{j+1}s_j\leq w'$ and exactly one of them holds. We fix an arbitrary reduced decomposition of $w'$. If $s_js_{j+1}\leq w'$, then we could move $s_{j'}$ (resp.~$s_{j''}$) in the reduced decomposition of $w'$ to the left of $s_{j+1}$ (resp.~to the right of $s_j$) for each $j_1\leq j'<j$ (resp.~for each $j+1<j''\leq j_2$), and thus we deduce $w'=x'y'$. Similar argument shows that $w'=y'x'$ whenever $s_{j+1}s_j\leq w$. 
It is clear that $x'y'\neq y'x'$ and thus $xy\neq yx$ in this case.
The proof is thus finished.
\end{proof}

\begin{defn}\label{def: adjacent pair}
Let $I\subseteq \Delta$ be a non-empty interval. A pair $w_1,w_2\in\Gamma_{I}$ is called \emph{adjacent} if $w_1\neq w_2$ and there exist a reduced decomposition of $w_1$ and a reduced decomposition of $w_2$ that differ from each other by permuting two adjacent simple reflections.
\end{defn}

\begin{lem}\label{lem: interval adjacent pair}
Let $I\subseteq \Delta$ be a non-empty interval and $w_1,w_2\in\Gamma_{I}$ be an pair of elements. Then the pair $w_1,w_2$ is adjacent if and only if there exists a partition $I=I_0\sqcup I_0'$ as well as $x\in\Gamma_{I_0}$ and $y\in\Gamma_{I_0'}$ such that both $I_0$ and $I_0'$ are non-empty intervals and we have $w_1=xy$ and $w_2=yx$.
\end{lem}
\begin{proof}
As $I$ is a non-empty interval, there exists $j_1,j_2\in\Delta$ such that $j_1\leq j_2$ and $I=[j_1,j_2]$.
We only prove the `only if' direction, and `if' direction follows by reversing the argument.
Assume that $w_1,w_2$ are adjacent, then by Definition~\ref{def: adjacent pair} there exists $j,j'\in I=[j_1,j_2]$ such that $w_1=w's_js_{j'}w''$ and $w_2=w's_{j'}s_jw''$ with both expressions being reduced. As $w_1\neq w_2$, we must have $s_js_{j'}\neq s_{j'}s_j$ which is equivalent to $|j-j'|=1$. Upon exchanging $w_1$ and $w_2$ (and thus exchanging $x$ and $y$ below), it is harmless to assume that $j'=j+1$.
As $w_1,w_2\in\Gamma_{I}$, we know that $\mathrm{Supp}(w')\cap\mathrm{Supp}(w'')=\emptyset$ and that $\mathrm{Supp}(w')\cup\mathrm{Supp}(w'')=[j_1,j_2]\setminus \{j,j'\}$. By braided relation we may arrange the choice of $w',w''$ so that $\mathrm{Supp}(w')=[j_1,j]\setminus\{j\}$ and $\mathrm{Supp}(w'')=[j',j_2]\setminus\{j'\}$, and then set $x\defeq w's_j$, $y\defeq s_{j'}w''$ ,with $I_0\defeq [j_1,j]$ and $I_0'\defeq [j',j_2]$ being non-empty intervals.
The proof is thus finished.
\end{proof}

\subsection{$U(\fg)$-modules and category $\cO$}\label{subsec: g category O}
We recall some well-known results on $U(\fg)$-modules from \cite{Hum08} and \cite{BQ24}, and then prove all results we need on $\mathrm{Ext}_{U(\fg)}^{\bullet}(-,-)$.
We also prove several basic results on generalized Verma modules (cf.~Lemma~\ref{lem: Lie Ext Verma vanishing} and Lemma~\ref{lem: tensor of parabolic Verma}) and generalized Steinberg modules (see Proposition~\ref{prop: Lie St Ext w decomposition}) which have important applications in later sections. 

For each $k\geq 1$, the diagonal map $\fg\rightarrow \fg^{k+1}$ induces a diagonal map
\begin{equation}\label{equ: envelop alg diagonal}
D_{\fg}^{k}: U(\fg)\rightarrow \bigotimes_{E}^{k+1}U(\fg).
\end{equation}
We start with a discussion on resolutions of objects in $\mathrm{Mod}_{U(\fg)}$.
\begin{cons}\label{cons: resolution}
Let $M\in\mathrm{Mod}_{U(\fg)}$. We have the following two different versions of standard resolutions of $M$ by free $U(\fg)$-modules.
\begin{enumerate}[label=(\roman*)]
\item \label{it: free resolution 1} We define 
\begin{equation}\label{equ: free resolution 1}
B_{\bullet}(\fg,M)\defeq \bigotimes_{E}^{\bullet+1}U(\fg)\otimes_EM=U(\fg)\otimes_E\cdots\otimes_EU(\fg)\otimes_EM
\end{equation}
with $U(\fg)$ acting by left multiplication on the first tensor factor. Note that (\ref{equ: free resolution 1}) is a simplicial resolution of $M$ by free $U(\fg)$-modules, with face maps given by multiplying adjacent tensor factors and degeneration maps given by inserting $1$ as an extra tensor factor (cf.~\cite[\S~8.6]{Wei94}). 
\item \label{it: resolution 2} We define 
\begin{equation}\label{equ: free resolution 2}
\tld{B}_{\bullet}(\fg,M)\defeq \bigotimes_{E}^{\bullet+1}U(\fg)\otimes_EM=U(\fg)\otimes_E\cdots\otimes_EU(\fg)\otimes_EM
\end{equation}
with the action of $\delta\in U(\fg)$ on $\tld{B}_{m}(\fg,M)$ given by the multiplication of $D_{\fg}^{m+1}(\delta)\in \bigotimes_{E}^{m+2}U(\fg)$. Note that (\ref{equ: free resolution 1}) is a simplicial resolution of $M$ by $U(\fg)$-modules, with face maps given by omitting one of the tensor factors and degeneration maps induced from the diagonal map $D_{\fg}^1:U(\fg)\rightarrow U(\fg)\otimes_EU(\fg)$.
One can check that the following formula
\begin{equation}\label{equ: left diagonal transfer formula}
\delta_0\otimes_E\cdots\otimes_E\delta_{m}\otimes_Ea\mapsto D_{\fg}^{m+1}(\delta_0)\cdots D_{\fg}^{1}(\delta_{m})\cdot a
\end{equation}
gives a $U(\fg)$-equivariant isomorphism
\[B_{m}(\fg,M)\buildrel\sim\over\longrightarrow \tld{B}_{m}(\fg,M)\]
for each $m\geq 0$, which is moreover compatible with the face and degeneration maps on both sides.
In other words, we obtain the following isomorphism
\begin{equation}\label{equ: left diagonal transfer}
B_{\bullet}(\fg,M)\buildrel\sim\over\longrightarrow \tld{B}_{\bullet}(\fg,M)
\end{equation}
between simplicial resolutions of $M$.
In particular, we see that (\ref{equ: free resolution 2}) is also a simplicial resolution of $M$ by free $U(\fg)$-modules.
\item \label{it: resolution 3} We define 
\begin{equation}\label{equ: free resolution 3}
\tld{B}_{\bullet}(\fg,M)_{\dagger}\defeq U(\fg)\otimes_E(\bigotimes_{E}^{\bullet}U(\fg)\otimes_EM)
\end{equation}
with the action of $\delta\in U(\fg)$ given by the left multiplication of the first tensor factor. 
We equip (\ref{equ: free resolution 3}) with the unique possible face and degeneration maps such that 
\begin{equation}\label{equ: partial left diagonal transfer formula}
\delta_0\otimes_E\cdots\otimes_E\delta_{m}\otimes_Ea\mapsto \delta_0\cdot D_{\fg}^{m}(\delta_{1})\cdots D_{\fg}^{1}(\delta_{m})\cdot a
\end{equation}
gives an isomorphism
\begin{equation}\label{equ: partial left diagonal transfer}
B_{\bullet}(\fg,M)\buildrel\sim\over\longrightarrow \tld{B}_{\bullet}(\fg,M)_{\dagger}
\end{equation}
between simplicial resolutions of $M$.
It is clear that (\ref{equ: left diagonal transfer}) naturally factors through (\ref{equ: partial left diagonal transfer}).
\end{enumerate}
\end{cons}
We define
\[
C^{\bullet}(\fg,M)\defeq\Hom_{U(\fg)}(B_{\bullet}(\fg,M),1_{U(\fg)})
\]
and similarly for $\tld{C}^{\bullet}(\fg,M)$ and $\tld{C}^{\bullet}(\fg,M)_{\dagger}$, and thus obtain isomorphisms between cosimplicial complex of $E$-vector spaces (from (\ref{equ: left diagonal transfer}) and (\ref{equ: partial left diagonal transfer}))
\begin{equation}\label{equ: Lie coh left diagonal transfer}
\tld{C}^{\bullet}(\fg,M)\buildrel\sim\over\longrightarrow \tld{C}^{\bullet}(\fg,M)_{\dagger}\buildrel\sim\over\longrightarrow C^{\bullet}(\fg,M).
\end{equation}
We will frequently abuse the notation $C^{\bullet}(\fg,M)$ for its associated (unnormalized) cochain complex and write $d_{\fg,M}^{\bullet}$ for its differential map. We use similar convention for $\tld{C}^{\bullet}(\fg,M)$ and $\tld{C}^{\bullet}(\fg,M)_{\dagger}$.

For each $m\geq 0$, we have an embedding
\begin{equation}\label{equ: wedge to free}
\wedge^{m}\fg\otimes_EM\rightarrow \bigotimes_{E}^{m+1}U(\fg)\otimes_EM
\end{equation}
given by
\begin{equation}\label{equ: wedge to free formula}
(x_{1}\wedge\cdots\wedge x_{m})\otimes_Eu\mapsto \sum_{\sigma\in\mathfrak{S}_{m}}\varepsilon(\sigma)1\otimes_Ex_{\sigma(1)}\otimes_E\cdots\otimes_Ex_{\sigma(m)}\otimes_Eu
\end{equation}
for each $x_{i}\in\fg$ for $1\leq i\leq m$ and $u\in M$. Here $\mathfrak{S}_{m}$ is the permutation group of $\{1,\dots,m\}$ and $\varepsilon(\sigma)\in\{1,-1\}$ is the sign associated with each permutation $\sigma\in\mathfrak{S}_{m}$.
Recall from \cite[\S 1.1]{BW80} that we have the following Chevalley-Eilenberg complex
\begin{equation}\label{equ: CE g M}
\mathrm{CE}^{\bullet}(\fg,M)\defeq \Hom_{E}(\wedge^{\bullet}\fg\otimes_EM,E)=\Hom_{E}(\wedge^{\bullet}\fg,\Hom_{E}(M,E)).
\end{equation}
We write $d_{\mathrm{CE},\fg,M}^{\bullet}$ for the differential maps of (\ref{equ: CE g M}) which are explicitly given in \cite[\S 1.1]{BW80}.
For each $m\geq 0$, we notice that the map (\ref{equ: wedge to free}) induces the following maps 
\begin{equation}\label{equ: left free wedge transfer}
C^{m}(\fg,M)\hookrightarrow \Hom_{E}(B_{m}(\fg,M),E)\twoheadrightarrow \Hom_{E}(\wedge^{m}\fg\otimes_EM,E)=\mathrm{CE}^{m}(\fg,M)
\end{equation}
and similar maps from $\tld{C}^{m}(\fg,M)$ and $\tld{C}^{m}(\fg,M)_{\dagger}$ to $\mathrm{CE}^{m}(\fg,M)$.
The following result is well-known but we recall it here for convenience.
\begin{lem}\label{lem: Lie resolution left diagonal transfer}
We have the following commutative diagram of maps between complex
\begin{equation}\label{equ: Lie resolution left diagonal transfer}
\xymatrix{
\tld{C}^{\bullet}(\fg,M) \ar^{\sim}[r] \ar[d] & \tld{C}^{\bullet}(\fg,M)_{\dagger} \ar^{\sim}[r] \ar[d] & C^{\bullet}(\fg,M) \ar[d]\\
\mathrm{CE}^{\bullet}(\fg,M) \ar@{=}[r] & \mathrm{CE}^{\bullet}(\fg,M) \ar@{=}[r] & \mathrm{CE}^{\bullet}(\fg,M)
}
\end{equation}
with all three vertical maps induced from (\ref{equ: wedge to free}) for each $m\geq 0$, and the top horizontal maps being the isomorphisms from (\ref{equ: Lie coh left diagonal transfer}).
\end{lem}
\begin{proof}
Dually, it suffices to check the commutativity of a diagram of the form
\begin{equation}\label{equ: left diagonal free diagram}
\xymatrix{
\wedge^{m}\fg\otimes_EM \ar@{=}[r] \ar[d] & \wedge^{m}\fg\otimes_EM \ar@{=}[r] \ar[d] & \wedge^{m}\fg\otimes_EM \ar[d]\\
B_{m}(\fg,M)\otimes_{U(\fg)}E \ar^{\sim}[r] & \tld{B}_{m}(\fg,M)_{\dagger}\otimes_{U(\fg)}E \ar^{\sim}[r] & \tld{B}_{m}(\fg,M)\otimes_{U(\fg)}E
}
\end{equation}
for each $m\geq 0$ and their compatibility with differential maps.
We divide the proof into two steps.

\textbf{Step $1$}: We prove the commutativity of 
\begin{equation}\label{equ: left diagonal free resolution diagram}
\xymatrix{
\wedge^{m}\fg\otimes_EM \ar@{=}[r] \ar[d] & \wedge^{m}\fg\otimes_EM \ar@{=}[r] \ar[d] & \wedge^{m}\fg\otimes_EM \ar[d]\\
B_{m}(\fg,M) \ar^{\sim}[r] & \tld{B}_{m}(\fg,M)_{\dagger} \ar^{\sim}[r] & \tld{B}_{m}(\fg,M)
}
\end{equation}
for each $m\geq 0$ (which implies the commutativity of (\ref{equ: left diagonal free diagram}) for each $m\geq 0$).\\
Let $x_{i}\in\fg$ for each $1\leq i\leq m$ and $u\in M$. All vertical maps of (\ref{equ: left diagonal free resolution diagram}) are given by (\ref{equ: wedge to free formula}). We need to show that the isomorphism in the bottom row of (\ref{equ: left diagonal free resolution diagram}) leave the term
\begin{equation}\label{equ: image of wedge sum}
\sum_{\sigma\in\mathfrak{S}_{m}}\varepsilon(\sigma)1\otimes_Ex_{\sigma(1)}\otimes_E\cdots\otimes_Ex_{\sigma(m)}\otimes_Eu
\end{equation}
unchanged.
For each $\sigma\in\mathfrak{S}_{m}$, the image of $1\otimes_Ex_{\sigma(1)}\otimes_E\cdots\otimes_Ex_{\sigma(m)}\otimes_Eu$ under either the left bottom horizontal map of (\ref{equ: left diagonal free resolution diagram}) or under the composition of the bottom row of (\ref{equ: left diagonal free resolution diagram}) is given by (see (\ref{equ: left diagonal transfer formula}))
\begin{equation}\label{equ: successive left diagonal transfer}
\varepsilon(\sigma)D_{\fg}^{m}(x_{\sigma(1)})\cdots D_{\fg}^{1}(x_{\sigma(m)})\cdot u,
\end{equation}
where we abuse the following element (with $x$ in the $j$-th copy of $U(\fg)$ for the $j$-th term of the sum)
\[D_{\fg}^{i}(x)\in \sum_{j=0}^{i}1\otimes_E\cdots\otimes_Ex\otimes_E\cdots\otimes_E1\in\bigotimes_{E}^{i+1}U(\fg)\]
for its image in $\bigotimes_{E}^{m+2}U(\fg)$ via the embedding by tensor $m+1-i$ copies of $1\in U(\fg)$ from the left.
Note that we can naturally express (\ref{equ: successive left diagonal transfer}) as a sum of elements of the form
\begin{equation}\label{equ: left diagonal transfer term}
\varepsilon(\sigma) 1\otimes_Ey_{1}\otimes_E\cdots\otimes_Ey_{m}\otimes_E(y_{m+1}\cdot v)
\end{equation}
with each $y_{k}$ being the product of $x_{\sigma(i)}$ for some $i\leq k$ with respect to the increasing order of $i$ such that $x_{\sigma(i)}$ appears in exactly one $y_{k}$ for each $1\leq i\leq m$. It is then a purely combinatorial observation all such terms (\ref{equ: left diagonal transfer term}) cancel while we sum up (\ref{equ: successive left diagonal transfer}) over $\sigma\in\mathfrak{S}_{m}$, with the only exception being those terms (\ref{equ: left diagonal transfer term}) satisfying $y_{k}=x_{\sigma(k)}$ for each $1\leq k\leq m$ and $y_{m+1}=1$. In other words, upon summing up (\ref{equ: successive left diagonal transfer}) over $\sigma\in\mathfrak{S}_{m}$, we recover (\ref{equ: image of wedge sum}).

\textbf{Step $2$}: We prove the commutativity of
\begin{equation}\label{equ: left free wedge diagram}
\xymatrix{
\wedge^{m}\fg\otimes_EM \ar[r] \ar[d] &  \wedge^{m-1}\fg\otimes_EM \ar[d]\\
B_{m}(\fg,M)\otimes_{U(\fg)}E \ar[r] &  B_{m-1}(\fg,M)\otimes_{U(\fg)}E
}
\end{equation}
for each $m\geq 1$.\\
Let $x_{i}\in\fg$ for each $1\leq i\leq m$ and $u\in M$. The image of (\ref{equ: image of wedge sum}) under $B_{m}(\fg,M)\rightarrow B_{m-1}(\fg,M)$ equals 
\begin{equation}\label{equ: X sigma total}
\sum_{\sigma\in\mathfrak{S}_{m}}\varepsilon(\sigma)(X_{\sigma,0}+X_{\sigma,1}+X_{\sigma,2}), 
\end{equation}
where
\[
\left\{\begin{array}{ccc}
X_{\sigma,0} & \defeq & x_{\sigma(1)}\otimes_E\cdots\otimes_Ex_{\sigma(m)}\otimes_Eu\\
X_{\sigma,1} & \defeq & \sum_{i=1}^{m-1}(-1)^{i}1\otimes_Ex_{\sigma(1)}\otimes_E\cdots \otimes_Ex_{\sigma(i)}x_{\sigma(i+1)}\otimes_E\cdots\otimes_Ex_{\sigma(m)}\otimes_Eu\\
X_{\sigma,2} & \defeq & (-1)^{m}1\otimes_Ex_{\sigma(1)}\otimes_E\cdots\otimes_Ex_{\sigma(m-1)}\otimes_E(x_{\sigma(m)}\cdot u)
\end{array}\right.
\]
It is clear that $X_{\sigma,0}$ has image zero in $B_{m-1}(\fg,M)\otimes_{U(\fg)}E$. 
For each $1\leq \ell\leq m$, we write $\mathfrak{S}_{m,\ell}$ for the set of all bijections from $[1,m-1]$ to $[1,m]\setminus\{\ell\}$, which comes with a standard sign function $\varepsilon(-)$ with respect to the unique order preserving bijection from $[1,m-1]$ to $[1,m]\setminus\{\ell\}$. For each $\sigma\in\mathfrak{S}_{m}$ with $\sigma(m)=\ell$, we have $\sigma'\defeq\sigma|_{[1,m-1]}\in\mathfrak{S}_{m,\ell}$ and $\varepsilon(\sigma)\varepsilon(\sigma')=(-1)^{m-\ell}$. This together with the formula of $X_{\sigma,2}$ for each $\sigma\in\mathfrak{S}_{m}$ implies that
\begin{equation}\label{equ: X sigma 2}
\sum_{\sigma\in\mathfrak{S}_{m}}\varepsilon(\sigma)X_{\sigma,2}=\sum_{\ell=1}^{m}(-1)^{\ell}\sum_{\sigma'\in\mathfrak{S}_{m,\ell}}1\otimes_Ex_{\sigma'(1)}\otimes_E\cdots\otimes_Ex_{\sigma'(m-1)}\otimes_E(x_{\ell}\cdot u)\in B_{m-1}(\fg,M).
\end{equation}
Similarly, for each $1\leq i\leq m-1$ and $1\leq j\neq k\leq m$, we define $\mathfrak{S}_{m,i,j,k}$ as the set of bijections from $[1,m]\setminus\{i,i+1\}$ to $[1,m]\setminus\{j,k\}$ which comes with a standard sign function (and is independent of the order of the pair $j,k$). 
Then we check that
\begin{equation}\label{equ: X sigma 1}
\sum_{\sigma\in\mathfrak{S}_{m}}\varepsilon(\sigma)X_{\sigma,1}=\sum_{1\leq j\neq k\leq m}(-1)^{j+k+i}\sum_{\sigma''\in\mathfrak{S}_{m,i,j,k}}1\otimes_Ex_{\sigma''(1)}\otimes_E\cdots\otimes_Ex_{j}x_{k}\otimes_E\cdots\otimes_Ex_{\sigma''(m)}\otimes_Eu \in B_{m-1}(\fg,M).
\end{equation}
Combining (\ref{equ: X sigma 2}) and (\ref{equ: X sigma 1}) with our claim on $X_{\sigma,0}$ for each $\sigma\in\mathfrak{S}_{m}$, we conclude that the image of (\ref{equ: X sigma total}) in $B_{m-1}(\fg,M)\otimes_{U(\fg)}E$ equals that of $(x_1\wedge\cdots\wedge x_{m})\otimes_Eu$ under the composition of
\[\wedge^{m}\fg\otimes_EM\rightarrow \wedge^{m-1}\fg\otimes_EM\rightarrow B_{m-1}(\fg,M)\otimes_{U(\fg)}E.\]
The proof is thus finished.
\end{proof}

Let $\fh\subseteq\fg$ be an $E$-Lie subalgebra. We define
\begin{equation}\label{equ: relative left resolution}
B_{m}(\fg,\fh,M)\defeq \bigotimes_{U(\fh)}^{m+1}U(\fg)\otimes_{U(\fh)}M=U(\fg)\otimes_{U(\fh)}\cdots\otimes_{U(\fh)}U(\fg)\otimes_{U(\fh)}M,
\end{equation}
\begin{equation}\label{equ: relative partial left diagonal resolution}
\tld{B}_{m}(\fg,\fh,M)_{\dagger}\defeq U(\fg)\otimes_{U(\fh)}(\bigotimes_{E}^{m}(U(\fg)\otimes_{U(\fh)}E)\otimes_EM),
\end{equation}
and
\begin{equation}\label{equ: relative diagonal resolution}
\tld{B}_{m}(\fg,\fh,M)\defeq \bigotimes_{E}^{m+1}(U(\fg)\otimes_{U(\fh)}E)\otimes_{E}M
\end{equation}
for each $m\geq 0$, with $U(\fg)$ acting on (\ref{equ: relative left resolution}) and (\ref{equ: relative partial left diagonal resolution}) by left multiplication on the first tensor factor and acting on (\ref{equ: relative diagonal resolution}) diagonally via $D_{\fg}^{m+1}: U(\fg)\rightarrow\bigotimes_{E}^{m+2}U(\fg)$.
We have natural surjections between $U(\fg)$-modules $B_{m}(\fg,M)\twoheadrightarrow B_{m}(\fg,\fh,M)$ and the simplicial structure on $B_{\bullet}(\fg,M)$ descends to one on $B_{\bullet}(\fg,\fh,M)$. We similarly obtain the simplicial complex $\tld{B}_{\bullet}(\fg,\fh,M)$ and $\tld{B}_{\bullet}(\fg,\fh,M)_{\dagger}$.
We set
\[C^{\bullet}(\fg,\fh,M)\defeq \Hom_{U(\fg)}(B_{m}(\fg,\fh,M), 1_{U(\fg)})\]
and similarly for $\tld{C}^{\bullet}(\fg,\fh,M)$ and $\tld{C}^{\bullet}(\fg,\fh,M)_{\dagger}$, and arrive the following commutative diagram of maps between cosimplicial complex of $E$-vector spaces
\begin{equation}\label{equ: relative left diagonal diagram}
\xymatrix{
\tld{C}_{\bullet}(\fg,\fh,M) \ar^{\sim}[rr] \ar[d] & & \tld{C}_{\bullet}(\fg,\fh,M)_{\dagger} \ar^{\sim}[rr] \ar[d] & & C^{\bullet}(\fg,\fh,M) \ar[d] \\
\tld{C}_{\bullet}(\fg,M) \ar^{\sim}[rr] & & \tld{C}_{\bullet}(\fg,M)_{\dagger} \ar^{\sim}[rr] & & C^{\bullet}(\fg,M) 
}.
\end{equation}
with all horizontal maps being isomorphisms. 
Now that the composition of
\[\tld{C}_{\bullet}(\fg,\fh,M)_{\dagger}\rightarrow \tld{C}_{\bullet}(\fg,M)_{\dagger}\rightarrow \mathrm{CE}^{\bullet}(\fg,M)\]
clearly factors through the following Chevalley-Eilenberg
\begin{equation}\label{equ: relative CE M}
\mathrm{CE}^{\bullet}(\fg,\fh,M)\defeq \Hom_{U(\fh)}(\wedge^{\bullet}(\fg/\fh)\otimes_EM, 1_{U(\fh)})=\Hom_{U(\fg)}(U(\fg)\otimes_{U(\fh)}(\wedge^{\bullet}(\fg/\fh)\otimes_EM), 1_{U(\fg)}),
\end{equation}
we obtain the following commutative diagram of maps between complex of $E$-vector spaces
\begin{equation}\label{equ: relative std to CE}
\xymatrix{
\tld{C}_{\bullet}(\fg,\fh,M)_{\dagger} \ar^{\sim}[rr] \ar[d]& & C^{\bullet}(\fg,\fh,M) \ar[rr] & & C^{\bullet}(\fg,M) \ar[d]\\
\mathrm{CE}^{\bullet}(\fg,\fh,M) \ar[rrrr] & & & & \mathrm{CE}^{\bullet}(\fg,M)
}
\end{equation}
Now that $U(\fg)\otimes_E(\wedge^{\bullet}\fg\otimes_EM)\rightarrow B_{\bullet}(\fg,M)$ is a map between resolution of $M$ by free $U(\fg)$-modules that lifts the identity map of $M$, and that $U(\fg)\otimes_{U(\fh)}(\wedge^{\bullet}(\fg/\fh)\otimes_EM)\rightarrow \tld{B}_{\bullet}(\fg,\fh,M)_{\dagger}$ is a map between resolution of $M$ by $U(\fg)$-modules which are relative projective with respect to the pair $U(\fh)\subseteq U(\fg)$ that lifts the identity map of $M$, we know that each vertical map of (\ref{equ: relative std to CE}) is a quasi-isomorphism whose associated map between cohomology is the identity map, and both horizontal maps of (\ref{equ: relative std to CE}) induce by taking cohomology the following map
\[\mathrm{Ext}_{U(\fg),U(\fh)}^{\bullet}(M,1_{\fg})\rightarrow \mathrm{Ext}_{U(\fg)}^{\bullet}(M,1_{\fg}).\]
We also note that all vertical maps of (\ref{equ: Lie resolution left diagonal transfer}) are quasi-isomorphisms, and the commutative diagram (\ref{equ: Lie resolution left diagonal transfer}) induces by taking cohomology the commutative diagram of identity maps between $\mathrm{Ext}_{U(\fg)}^{\bullet}(M,1_{\fg})$.

Let $M_i\in\mathrm{Mod}_{U(\fg)}$ for $i=0,1,2$ together with a map $M_2\rightarrow M_1\otimes_EM_0$ between $U(\fg)$-modules, with $U(\fg)$ acting diagonally on $M_1\otimes_EM_0$ by $D_{\fg}^{1}: U(\fg)\rightarrow U(\fg)\otimes_EU(\fg)$.
We consider the composition of
\begin{multline}\label{equ: Lie cup composition}
C^{m_0}(\fg,M_0)\otimes_EC^{m_1}(\fg,M_1)=\Hom_{U(\fg)}(B_{m_0}(\fg,M_0),1_{U(\fg)})\otimes_E\Hom_{U(\fg)}(B_{m_1}(\fg,M_1),1_{U(\fg)})\\
\rightarrow \Hom_{U(\fg)}(M_1\otimes_EB_{m_0}(\fg,M_0),M_1)\otimes_E\Hom_{U(\fg)}(B_{m_1}(\fg,M_1),1_{U(\fg)})\\
\rightarrow \Hom_{U(\fg)}(B_{m_1}(\fg,M_1\otimes_EB_{m_0}(\fg,M_0)),B_{m_1}(\fg,M_1))\otimes_E\Hom_{U(\fg)}(B_{m_1}(\fg,M_1),1_{U(\fg)})\\
\rightarrow \Hom_{U(\fg)}(B_{m_1}(\fg,M_1\otimes_EB_{m_0}(\fg,M_0)),1_{U(\fg)})
\end{multline}
for each $m_0,m_1\geq 0$, which altogether give a map between double complex
\begin{equation}\label{equ: Lie cup left double}
C^{\bullet}(\fg,M_0)\otimes_EC^{\bullet'}(\fg,M_1)\rightarrow \Hom_{U(\fg)}(B_{\bullet'}(\fg,M_1\otimes_EB_{\bullet}(\fg,M_0)),1_{U(\fg)}).
\end{equation}
The map $B_{\bullet}(\fg,M_0)\rightarrow M_0$ induces the following map between resolutions of $M_1\otimes_EM_0$ by free $U(\fg)$-modules (which lifts the identity map of $M_1\otimes_EM_0$)
\[\mathrm{Tot}(B_{\bullet'}(\fg,M_1\otimes_EB_{\bullet}(\fg,M_0)))\rightarrow B_{\bullet'}(\fg,M_1\otimes_EM_0)\]
and further induces a quasi-isomorphism
\[C^{\bullet}(\fg,M_1\otimes_EM_0)\rightarrow \Hom_{U(\fg)}(B_{\bullet'}(\fg,M_1\otimes_EB_{\bullet}(\fg,M_0)),1_{U(\fg)}).\]
Combining this with (\ref{equ: Lie cup left double}) and the map $M_2\rightarrow M_1\otimes_EM_0$, we obtain maps
\begin{multline*}
\mathrm{Tot}(C^{\bullet}(\fg,M_0)\otimes_EC^{\bullet'}(\fg,M_1))\\
\rightarrow \mathrm{Tot}(\Hom_{U(\fg)}(B_{\bullet'}(\fg,M_1\otimes_EB_{\bullet}(\fg,M_0)),1_{U(\fg)}))\leftarrow C^{\bullet}(\fg,M_1\otimes_EM_0)\rightarrow C^{\bullet}(\fg,M_2)
\end{multline*}
with the middle leftward map being a quasi-isomorphism that induce the identity map of $\mathrm{Ext}_{U(\fg)}^{\bullet}(M_1\otimes_EM_0,1_{\fg})$ by taking cohomology. We obtain this way a map
\begin{equation}\label{equ: Lie cup left}
\kappa_{\sharp}: \mathrm{Tot}(C^{\bullet}(\fg,M_0)\otimes_EC^{\bullet}(\fg,M_1))\dashrightarrow C^{\bullet}(\fg,M_2)
\end{equation}
in the derived category of $E$-vector spaces.
Replacing $C^{\bullet}(\fg,M_i)$ with $\tld{C}^{\bullet}(\fg,M_i)$ for each $i=0,1,2$, a parallel construction gives a map
\begin{equation}\label{equ: Lie cup diagonal 1}
\kappa_{\dagger}: \mathrm{Tot}(\tld{C}^{\bullet}(\fg,M_0)\otimes_E\tld{C}^{\bullet}(\fg,M_1))\dashrightarrow \tld{C}^{\bullet}(\fg,M_2).
\end{equation}
For each $m\geq 0$, $\tld{B}_{m}(\fg,M_1)\otimes_E\tld{B}_{m}(\fg,M_0)$ admits a natural $U(\fg)\otimes_EU(\fg)$-action under which we can identify it with $\tld{B}_{m}(\fg\times\fg,M_1\otimes_EM_0)$. 
This way, we obtain a natural identification (with $\mathrm{Diag}(-)$ standing for the simplicial diagonal of a simplicial bi-complex)
\begin{equation}\label{equ: Lie simplicial diagonal}
\mathrm{Diag}(\tld{B}_{\bullet}(\fg,M_0)\otimes_E\tld{B}_{\bullet}(\fg,M_1))=\tld{B}_{\bullet}(\fg\times\fg,M_1\otimes_EM_0).
\end{equation}
We consider the following Alexander-Whitney map (see \cite[\S 8.5]{Wei94}) (associated with the simplicial bi-complex $\tld{B}_{\bullet}(\fg,M_0)\otimes_E\tld{B}_{\bullet}(\fg,M_1)$)
\begin{equation}\label{equ: Lie Diag Tot}
\mathrm{Diag}(\tld{B}_{\bullet}(\fg,M_0)\otimes_E\tld{B}_{\bullet}(\fg,M_1))\rightarrow \mathrm{Tot}(\tld{B}_{\bullet}(\fg,M_0)\otimes_E\tld{B}_{\bullet}(\fg,M_1))
\end{equation}
which is a quasi-isomorphism between complex of $U(\fg)\otimes_EU(\fg)$-modules by Eilenberg-Zilber theorem (cf.~\cite[\S 8.5]{Wei94}). 
Upon applying $\Hom_{U(\fg)\otimes_EU(\fg)}(-,1_{U(\fg)\otimes_EU(\fg)})$ to both (\ref{equ: Lie simplicial diagonal}) and (\ref{equ: Lie Diag Tot}), we obtain the following sequence of maps
\begin{multline}\label{equ: Lie Diag Tot co}
\mathrm{Tot}(\tld{C}^{\bullet}(\fg,M_0)\otimes_E\tld{C}^{\bullet}(\fg,M_1))\\
\rightarrow
\Hom_{U(\fg)\otimes_EU(\fg)}(\mathrm{Tot}(\tld{B}_{\bullet}(\fg,M_0)\otimes_E\tld{B}_{\bullet}(\fg,M_1)),1_{U(\fg)\otimes_EU(\fg)})\\
\rightarrow \Hom_{U(\fg)\otimes_EU(\fg)}(\mathrm{Diag}(\tld{B}_{\bullet}(\fg,M_0)\otimes_E\tld{B}_{\bullet}(\fg,M_1)),1_{U(\fg)\otimes_EU(\fg)})\\
=\tld{C}(\fg\times\fg,M_1\otimes_EM_0)\rightarrow \tld{C}(\fg,M_2)
\end{multline}
with the last map induced from the diagonal map $D_{\fg}^{1}:U(\fg)\rightarrow U(\fg)\otimes_EU(\fg)$ and the map $M_2\rightarrow M_1\otimes_EM_0$. Taking composition of these maps, we obtain another map in the derived category of $E$-vector spaces
\begin{equation}\label{equ: Lie cup diagonal 2}
\kappa_{\ddagger}: \mathrm{Tot}(\tld{C}^{\bullet}(\fg,M_0)\otimes_E\tld{C}^{\bullet}(\fg,M_1))\dashrightarrow \tld{C}^{\bullet}(\fg,M_2).
\end{equation}
Finally, the wedge product map
\[\Hom_{E}(\wedge^{m_0}\fg,E)\otimes_E\Hom_{E}(\wedge^{m_1}\fg,E)\rightarrow \Hom_{E}(\wedge^{m_0+m_1}\fg,E)\]
together with our given map $M_2\rightarrow M_1\otimes_EM_0$ induces a map (see (\ref{equ: CE g M}))
\[\mathrm{CE}^{m_0}(\fg,M_0)\otimes_E\mathrm{CE}^{m_1}(\fg,M_1)\rightarrow \mathrm{CE}^{m_0+m_1}(\fg,M_2)\]
for each $m_0,m_1\geq 0$, which altogether give a map between complex of $E$-vector spaces
\begin{equation}\label{equ: Lie cup wedge}
\kappa_{\flat}: \mathrm{Tot}(\mathrm{CE}^{\bullet}(\fg,M_0)\otimes_E\mathrm{CE}^{\bullet}(\fg,M_1))\rightarrow \mathrm{CE}^{\bullet}(\fg,M_2).
\end{equation}
For each $\ast\in\{\sharp,\dagger,\ddagger,\flat\}$, we have the following associated map on cohomology
\[\overline{\kappa}_{\ast}: \mathrm{Tot}(\mathrm{Ext}_{U(\fg)}^{\bullet}(M_0,1_{\fg})\otimes_E\mathrm{Ext}_{U(\fg)}^{\bullet}(M_1,1_{\fg}))\rightarrow \mathrm{Ext}_{U(\fg)}^{\bullet}(M_2,1_{\fg}),\]
using the discussion below (\ref{equ: relative std to CE}).
The following result is again well-known but we record it for our convenience.
\begin{lem}\label{lem: Lie cup diagram}
Let $M_i\in\mathrm{Mod}_{U(\fg)}$ for $i=0,1,2$ together with a map $M_2\rightarrow M_1\otimes_EM_0$ between $U(\fg)$-modules. 
The maps $\overline{\kappa}_{\ast}$ are equal for $\ast\in\{\sharp,\dagger,\ddagger,\flat\}$.
\end{lem}
\begin{proof}
The equality $\overline{\kappa}_{\sharp}=\overline{\kappa}_{\dagger}$ is clear from the isomorphism $C^{\bullet}(\fg,M_i)\cong\tld{C}^{\bullet}(\fg,M_i)$ as the definition of $\kappa_{\sharp}$ and $\kappa_{\dagger}$ are parallel with the only difference between the usage of $C^{\bullet}(\fg,M_i)$ or $\tld{C}^{\bullet}(\fg,M_i)$ for $i=0,1,2$.

We have the following natural identification between simplicial bi-complex
\begin{equation}\label{equ: simplicial bi complex transfer}
\tld{B}_{\bullet}(\fg,M_1\otimes_E\tld{B}_{\bullet}(\fg,M_0))=\tld{B}_{\bullet}(\fg,M_1)\otimes_E\tld{B}_{\bullet}(\fg,M_0),
\end{equation}
with $\mathrm{Tot}(-)$ of BHS being resolutions of $M_1\otimes_EM_0$ in $\mathrm{Mod}_{U(\fg)\otimes_EU(\fg)}$ by free $U(\fg)\otimes_EU(\fg)$-modules, and thus in particular resolutions of $M_1\otimes_EM_0$ in $\mathrm{Mod}_{U(\fg)}$ by free $U(\fg)$-modules via the diagonal map $D_{\fg}^{1}: U(\fg)\rightarrow U(\fg)\otimes_EU(\fg)$.
In particular, we obtain the following sequence of maps between resolutions of $M_1\otimes_EM_0$ in $\mathrm{Mod}_{U(\fg)}$ by free $U(\fg)$-modules
\begin{multline}\label{equ: Lie free resolution transfer}
\tld{B}_{\bullet}(\fg,M_1\otimes_EM_0)\leftarrow \mathrm{Tot}(\tld{B}_{\bullet}(\fg,M_1\otimes_E\tld{B}_{\bullet}(\fg,M_0)))\\
=\mathrm{Tot}(\tld{B}_{\bullet}(\fg,M_1)\otimes_E\tld{B}_{\bullet}(\fg,M_0))
\leftarrow \mathrm{Diag}(\tld{B}_{\bullet}(\fg,M_1)\otimes_E\tld{B}_{\bullet}(\fg,M_0))\\
=\tld{B}_{\bullet}(\fg\times\fg,M_1\otimes_EM_0)\leftarrow \tld{B}_{\bullet}(\fg,M_1\otimes_EM_0)
\end{multline}
which lift the identity map of $M_1\otimes_EM_0$ (with the leftward map in the middle row of (\ref{equ: Lie free resolution transfer}) being the Alexander-Whitney map). Hence, upon applying $\Hom_{U(\fg)}(-,1_{U(\fg)})$ and taking cohomology, (\ref{equ: Lie free resolution transfer}) induces identity maps of $\mathrm{Ext}_{U(\fg)}^{\bullet}(M_1\otimes_EM_0,1_{\fg})$. This implies the equality $\overline{\kappa}_{\dagger}=\overline{\kappa}_{\ddagger}$.

Finally, under (\ref{equ: Lie simplicial diagonal}), we have the following commutative diagram
\begin{equation}\label{equ: Lie diagonal wedge diagram}
\xymatrix{
\tld{C}(\fg\times\fg,M_1\otimes_EM_0) \ar[rr] \ar@/^/[d]^{\nabla^{\ast}} & & \mathrm{CE}^{\bullet}(\fg\times\fg,M_1\otimes_EM_0) \ar^{\wr}[d]\\
\Hom_{U(\fg)\otimes_EU(\fg)}(\mathrm{Tot}(\tld{B}_{\bullet}(\fg,M_0)\otimes_E\tld{B}_{\bullet}(\fg,M_1)),1_{U(\fg)\otimes_EU(\fg)}) \ar[rr] \ar@/^/[u] & & \mathrm{Tot}(\mathrm{CE}^{\bullet}(\fg,M_0)\otimes_E\mathrm{CE}^{\bullet}(\fg,M_1)) \ar@{=}[d]\\
\mathrm{Tot}(\tld{C}^{\bullet}(\fg,M_0)\otimes_E\tld{C}^{\bullet}(\fg,M_1)) \ar[rr] \ar[u] & & \mathrm{Tot}(\mathrm{CE}^{\bullet}(\fg,M_0)\otimes_E\mathrm{CE}^{\bullet}(\fg,M_1))
}
\end{equation}
with the right top isomorphism induced by taking $\Hom_{E}(-,E)$ from the isomorphism
\begin{equation}\label{equ: Lie wedge decomposition}
\mathrm{Tot}((\wedge^{\bullet}\fg\otimes_EM_0)\otimes_E(\wedge^{\bullet}\fg\otimes_EM_1))\cong \wedge^{\bullet}(\fg\times\fg)\otimes_E(M_0\otimes_EM_1), 
\end{equation}
the left top upward maps from (\ref{equ: Lie Diag Tot co}), and the left top downward map $\nabla^{\ast}$ being induced from the shuffle product map (see \cite[Prop.~8.6.13]{Wei94})
\begin{equation}\label{equ: Lie Diag Tot shuffle}
\mathrm{Tot}(\tld{B}_{\bullet}(\fg,M_0)\otimes_E\tld{B}_{\bullet}(\fg,M_1))\rightarrow \mathrm{Diag}(\tld{B}_{\bullet}(\fg,M_0)\otimes_E\tld{B}_{\bullet}(\fg,M_1)).
\end{equation}
Now that the two sided compositions of (\ref{equ: Lie Diag Tot shuffle}) and (\ref{equ: Lie Diag Tot}) induce identity maps on the cohomology (see the end of \cite[\S~8.5]{Wei94}), and that (\ref{equ: Lie Diag Tot shuffle}) restricts to the natural isomorphism (\ref{equ: Lie wedge decomposition}) by a direct check on explicit formula (cf.~\cite[Ex.~8.6.5]{Wei94}), we conclude the equality $\overline{\kappa}_{\ddagger}=\overline{\kappa}_{\flat}$ from (\ref{equ: Lie diagonal wedge diagram}).
\end{proof}
Now that $B_{\bullet}(\fg,M_0)$ is a resolution of $M_0$ by free $U(\fg)$-modules and $M_1\otimes_EB_{\bullet}(\fg,M_0)$ is actually a resolution of $M_1\otimes_EM_0$ by free $U(\fg)$-modules, we deduce from (\ref{equ: Lie cup composition}) that the associated map between cohomology of total complex (\ref{equ: Lie cup left double}) factors through the following maps
\begin{multline*}
\mathrm{Tot}(\mathrm{Ext}_{U(\fg)}^{\bullet}(M_0,1_{\fg})\otimes_E\mathrm{Ext}_{U(\fg)}^{\bullet}(M_1,1_{\fg}))\\
\rightarrow \mathrm{Tot}(\mathrm{Ext}_{U(\fg)}^{\bullet}(M_1\otimes_EM_0,M_1)\otimes_E\mathrm{Ext}_{U(\fg)}^{\bullet}(M_1,1_{\fg}))
\rightarrow \mathrm{Ext}_{U(\fg)}^{\bullet}(M_1\otimes_EM_0,1_{\fg}).
\end{multline*}
Consequently, combining with the map $M_2\rightarrow M_1\otimes_EM_0$, we see that $\overline{\kappa}_{\ast}$ (for $\ast\in\{\sharp,\dagger,\ddagger,\flat\}$ also fits into the following commutative diagram
\begin{equation}\label{equ: Lie cup Ext as composition}
\xymatrix{
\mathrm{Tot}(\mathrm{Ext}_{U(\fg)}^{\bullet}(M_0,1_{\fg})\otimes_E\mathrm{Ext}_{U(\fg)}^{\bullet}(M_1,1_{\fg})) \ar^{\overline{\kappa}_{\ast}}[rrr] \ar[d] & & & \mathrm{Ext}_{U(\fg)}^{\bullet}(M_2,1_{\fg}) \ar@{=}[d]\\
\mathrm{Tot}(\mathrm{Ext}_{U(\fg)}^{\bullet}(M_2,M_1)\otimes_E\mathrm{Ext}_{U(\fg)}^{\bullet}(M_1,1_{\fg})) \ar[rrr] & & & \mathrm{Ext}_{U(\fg)}^{\bullet}(M_2,1_{\fg})
}
\end{equation}

Let $s\geq 1$ and $M_i\in\mathrm{Mod}_{U(\fg)}$ for $0\leq i\leq s$. 
Note that $\bigotimes_{i=0}^{s}M_{i}$ is a $U(\fg)$-module via the diagonal map $D_{\fg}^{s}: U(\fg)\rightarrow \bigotimes_{E}^{s+1}U(\fg)$.
Let $M\in\mathrm{Mod}_{U(\fg)}$ and 
\begin{equation}\label{equ: g mod tensor}
M\rightarrow \bigotimes_{i=0}^{s}M_{i}
\end{equation}
be a map between $U(\fg)$-modules. Similar to $\kappa_{\rm{left}}$ from (\ref{equ: Lie cup left}), we can use the resolution
\[B_{\bullet}(\fg,M_{s}\otimes_EB_{\bullet}(\fg,M_{s-1}\otimes_E\cdots\otimes_EB_{\bullet}(\fg,M_0)))\]
of $\bigotimes_{i=0}^{s}M_{i}$ by free $U(\fg)$-modules and define a map
\begin{equation}\label{equ: Lie multi cup left}
\mathrm{Tot}(C^{\bullet}(\fg,M_{s})\otimes_{E}\cdots\otimes_EC^{\bullet}(\fg,M_{0}))\rightarrow C^{\bullet}(\fg,M).
\end{equation}
Similar to $\kappa_{\rm{CE}}$ from (\ref{equ: Lie cup wedge}), we can also define a map
\begin{equation}\label{equ: Lie multi cup CE}
\mathrm{Tot}(\mathrm{CE}^{\bullet}(\fg,M_{s})\otimes_{E}\cdots\otimes_E\mathrm{CE}^{\bullet}(\fg,M_{0}))\rightarrow \mathrm{CE}^{\bullet}(\fg,M).
\end{equation}
\begin{prop}\label{prop: Lie multi cup}
Let $s\geq 1$, $M_i\in\mathrm{Mod}_{U(\fg)}$ for $0\leq i\leq s$ and $M\in\mathrm{Mod}_{U(\fg)}$ be as above. 
Both the map (\ref{equ: Lie multi cup left}) and the map (\ref{equ: Lie multi cup CE}) induce by taking cohomology the same map
\begin{equation}\label{equ: Lie multi cup Ext}
\mathrm{Ext}_{U(\fg)}^{m_{s}}(M_{s},1_{U(\fg)})\otimes_{E}\cdots\otimes_E\mathrm{Ext}_{U(\fg)}^{m_{0}}(M_{0},1_{U(\fg)})\rightarrow \mathrm{Ext}_{U(\fg)}^{m}(M,1_{U(\fg)})
\end{equation}
for each $m_i\geq 0$ and $0\leq i\leq s$ with $m\defeq\sum_{i=0}^{s}m_{i}$.
Moreover, exchanging the order of $M_0,\dots,M_{s}$ only change the map (\ref{equ: Lie multi cup Ext}) up to a sign.
\end{prop}
\begin{proof}
The first claim follows directly from Lemma~\ref{lem: Lie cup diagram} by an increasing induction on $s\geq 1$. The second claim follows from the fact that (\ref{equ: Lie multi cup CE}) is by its very definition independent of the order of $M_{i}$ up to a sign.
\end{proof}

Given $M,M'\in\mathrm{Mod}_{U(\fg)}$ with $\Hom_{E}(M,M')$ viewed as $U(\fg)$-module via the diagonal map $D_{\fg}^1$, we define $\mathrm{Ext}_{\fg,\fh}^{\bullet}(M,M')$ as the cohomology of the following relative Chevalley-Eilenberg complex (see \cite[\S I.1.2 (1)]{BW80})
\begin{equation}\label{equ: relative CE Ext}
\Hom_{U(\fh)}(\wedge^{\bullet}(\fg/\fh),\Hom_{E}(M,M'))=\Hom_{U(\fh)}(\wedge^{\bullet}(\fg/\fh)\otimes_EM,M').
\end{equation}
Note that (\ref{equ: relative CE M}) can be recovered from (\ref{equ: relative CE Ext}) by taking $M'=1_{U(\fg)}$.
When $M=1_{U(\fg)}$, we also write 
\[H^{\bullet}(\fg,\fh,M')\defeq \mathrm{Ext}_{\fg,\fh}^{\bullet}(1_{U(\fg)},M')\]
for short.
When $\fh$ is reductive and both $M$ and $M'$ are $\fh$-finite, $\mathrm{Ext}_{\fg,\fh}^{\bullet}(M,M')$ can be alternatively interpreted as the $\mathrm{Ext}$-groups computed in the abelian category $\mathrm{Mod}_{U(\fg)}^{\fh-\rm{fin}}$ of locally $\fh$-finite $U(\fg)$-modules (cf.~\cite[\S 2.1]{BW80}). 
If both $M$ and $M'$ are furthermore locally $\fb$-finite and $Z(\fg)$-finite, then we have $\mathrm{Ext}_{\fg,\fh}^{k}(M,M')\neq 0$ for some $k\geq 0$ only if there exists an infinitesimal character $\xi: Z(\fg)\rightarrow E$ such that $M_{\xi}\neq 0\neq M'_{\xi}$.

Let $\fh_{i}\subseteq \fg_{i}$ be a pair of $E$-Lie algebra and $M_{i}'$ be a $U(\fg_{i})$-module for $i=0,1$. Then we have an natural isomorphism between complex of $U(\fh_0)\otimes_EU(\fh_1)=U(\fh_0\times\fh_1)$-modules
\[\mathrm{Tot}(\wedge^{\bullet}(\fg_0/\fh_0)\otimes_E\wedge^{\bullet}(\fg_1/\fh_0))\cong\wedge^{\bullet}(\fg_0\times\fg_1/\fh_0\times\fh_1)\]
which induces an identification between complex of $E$-vector spaces
\begin{multline*}
\mathrm{Tot}(\Hom_{U(\fh_0)}(\wedge^{\bullet}(\fg_0/\fh_0),M_0')\otimes_E\Hom_{U(\fh_1)}(\wedge^{\bullet}(\fg_1/\fh_1),M_1')\\
\cong\Hom_{U(\fh_0\times\fh_1)}(\wedge^{\bullet}(\fg_0\times\fg_1/\fh_0\times\fh_1),M_0'\otimes_EM_1')
\end{multline*}
whose upon taking cohomology recovers the following K\"unneth formula (cf.~\cite[\S 1.3]{BW80})
\begin{equation}\label{equ: relative Lie Kunneth}
\mathrm{Tot}(H^{\bullet}(\fg_0,\fh_0,M_0')\otimes_EH^{\bullet}(\fg_1,\fh_1,M_1'))\cong H^{\bullet}(\fg_0\times\fg_1,\fh_0\times\fh_1,M_0'\otimes_EM_1').
\end{equation}

We define $\cO_{\rm{alg}}^{\fb}\subseteq \mathrm{Mod}_{U(\fg)}$ as the full subcategory of $M$ which are locally $\fb$-finite, $\ft$-semi-simple with $M=\bigoplus_{\mu\in\Lambda}M_{\mu}$ and finitely generated as $U(\fg)$-modules. Note that $\cO_{\rm{alg}}^{\fb}$ is the full subcategory of the category $\cO$ of \cite{Hum08} of objects with integral (equivalently algebraic) weights.

For each $\mu\in \Lambda$, we have a Verma module $M(\mu)\defeq U(\fg)\otimes_{U(\fb)}\mu \in \cO_{\rm{alg}}^{\fb}$, which has an irreducible cosocle denoted by $L(\mu)$ (\cite[Thm.~1.2(f)]{Hum08}). Moreover recall that each simple object of $\cO_{\rm{alg}}^{\fb}$ has the form $L(\mu)$ for some $\mu\in \Lambda$ (\cite[\S 1.3]{Hum08}) and that each object of $\cO_{\rm{alg}}^{\fb}$ has finite length (\cite[\S 1.11]{Hum08}). In the sequel, we write $N(\mu)$ for the kernel of the surjection $M(\mu)\twoheadrightarrow L(\mu)$.

For $I\subseteq \Delta$, we consider the full subcategory $\cO_{\rm{alg}}^{\fp_{I}}\subseteq \cO_{\rm{alg}}^{\fb}$ of those $M$ which are locally $\fp_{I}$-finite, i.e.\ equal to the union of their finite dimensional $U(\fp_{I})$-submodules. As $\fp_{I}=\fl_{I}\oplus \fn_{I}$ and the category of finite dimensional $U(\fl_{I})$-modules is semi-simple, $M\in \cO_{\rm{alg}}^{\fb}$ is locally $\fp_{I}$-finite if and only if the (underlying) $U(\fl_I)$-module $M$ is a direct sum of (simple) finite dimensional $U(\fl_{I})$-modules.

Replacing $\fg$ with $\fl_{I}$, we can define an analogous full subcategory $\cO_{\fl_{I},\rm{alg}}^{\fb_{I}}$ of $\mathrm{Mod}_{U(\fl_{I})}$. We also define $\cO_{\fl_{I},\rm{alg}}^{\fl_{I}\cap \fp_{I'}}\hookrightarrow \cO_{\fl_{I},\rm{alg}}^{\fb_{I}}$ for each $I'\subseteq I$ (note that $\fl_{I}\cap \fp_{I'}$ is a parabolic subalgebra of $\fl_{I}$). We write $L^I(\mu)\in \cO_{\fl_{I},\rm{alg}}^{\fb_{I}}$ for the unique simple quotient of $U(\fl_{I})\otimes_{U(\fb_{I})}\mu$, and set
\begin{equation}\label{belongtoOb}
M^I(\mu)\defeq U(\fg)\otimes_{U(\fp_{I})}L^I(\mu).
\end{equation}
Be careful that we allow $M^I(\mu)$ and $L^I(\mu)$ to be defined for any $\mu\in \Lambda$, in particular $L^I(\mu)$ can be infinite dimensional. As $M^I(\mu)$ is a quotient of
\[M(\mu)=U(\fg)\otimes_{U(\fb)}\mu\cong U(\fg)\otimes_{U(\fp_{I})}(U(\fl_{I})\otimes_{U(\fb_{I})}\mu),\]
we see that $M^I(\mu)$ is in $\cO_{\rm{alg}}^{\fb}$ and has $L(\mu)$ as unique simple quotient. We write $N^I(\mu)$ for the kernel of the surjection $M^I(\mu)\twoheadrightarrow L(\mu)$. Moreover it follows from \cite[Prop.\ 9.3(e)]{Hum08} and \cite[Thm.\ 9.4]{Hum08} that $M^I(\mu)$ is in $\cO_{\rm{alg}}^{\fp_{I}}$ if and only if $L(\mu)$ is in $\cO_{\rm{alg}}^{\fp_{I}}$ if and only if $\mu\in \Lambda_I^{\rm dom}$ if and only if $L^I(\mu)$ is finite dimensional.

For $w\in W(G)$ and $I\subseteq \Delta$ we set $M(w)\defeq M(w\cdot 0)$, $L(w)\defeq L(w\cdot 0)$, $N(w)\defeq N(w\cdot 0)$, $L^I(w)\defeq L^I(w\cdot 0)$, $M^I(w)\defeq M^I(w\cdot 0)$ and $N^I(w)\defeq N^I(w\cdot 0)$. Note that all Jordan--H\"older factors of $N(w)$ have the form $L(w')$ for some $w'>w$ (\cite[\S\S 5.1,5.2]{Hum08} and \cite[\S 8.3(a)]{Hum08}).

The following result is well-known (cf.~\cite[Lem.~3.1.1]{BQ24}).
\begin{lem}\label{lem: dominance and left set}
Let $w\in W(G)$ be an element and $I\subseteq \Delta$ be a subset. Then $L(w)$ is in $\cO_{\rm{alg}}^{\fp_{I}}$ if and only if $I\cap D_L(w)=\emptyset$.
\end{lem}

The following result is also well-known (cf.~\cite[\S 1.10]{Hum08} and \cite[Lem.~3.1.3]{BQ24}).
\begin{lem}\label{lem: g block}
Let $I\subseteq \Delta$ and $M,M'\in\cO^{\fb_I}_{\fl_I,\rm{alg}}$. If $\mathrm{Ext}_{U(\fl_I)}^k(M,M')\neq 0$ for some $k\geq 0$, then there exists an infinitesimal character $\xi: Z(\fl_I)\rightarrow E$ such that $M_{\xi}\neq 0\neq M'_{\xi}$.
\end{lem}

Let $M$ in $\mathrm{Mod}_{U(\fg)}$ and $I\subseteq \Delta$.
We consider the following Chevalley--Eilenberg complex (as a special case of (\ref{equ: relative CE Ext}))
\begin{equation}\label{equ: CE complex}
M\rightarrow \Hom_{E}(\fn_{I}, M)\rightarrow\cdots\rightarrow \Hom_{E}(\wedge^k\fn_{I},M)\rightarrow\cdots
\end{equation}
with $M$ in degree zero (and $H^k(\fn_{I},M)$ as the cohomology group of this complex in degree $k\geq 0$). As the complex (\ref{equ: CE complex}) is $U(\fl_{I})$-equivariant, $H^k(\fn_{I},M)$ is naturally a $U(\fl_{I})$-module, and thus in particular a $U(\ft)$-module.
Let $M_{I}$ in $\mathrm{Mod}_{U(\fl_{I})}$ and $\fh\subseteq \fl_{I}$ be a reductive $E$-Lie subalgebra.
Assume moreover that both $M$ and $M_{I}$ are locally $\fh$-finite and in particular $H^{k}(\fn_{I},M)$ is locally $\fh$-finite for $k\geq 0$. Then we recall from \cite[Theorem~I.6.5]{BW80} the following Hochschild--Serre spectral sequence
\begin{equation}\label{equ: g spectral seq}
\mathrm{Ext}_{\fl_{I},\fh}^{\ell_1}(M_I, H^{\ell_2}(\fn_{I},M))\implies \mathrm{Ext}_{\fp_{I},\fh}^{\ell_1+\ell_2}(M_I, M)\cong \mathrm{Ext}_{\fg,\fh}^{\ell_1+\ell_2}(U(\fg)\otimes_{U(\fp_{I})}M_I, M)
\end{equation}
where the last isomorphism is Shapiro's lemma for relative Lie algebra cohomology (from exactness of $U(\fg)\otimes_{U(\fp_{I})}-$).
In particular, we have a canonical isomorphism
\begin{equation}\label{equ: g spectral seq 0}
\Hom_{U(\fl_{I})}(M_I, H^0(\fn_{I},M))\cong \Hom_{U(\fg)}(U(\fg)\otimes_{U(\fp_{I})}M_I, M).
\end{equation}
If furthermore $M_{I}\in\cO^{\fb_I}_{\fl_I,\rm{alg}}$ and $M\in \cO^{\fb}_{\rm{alg}}$ with $Z(\fl_{I})$ acting on $M_{I}$ by some infinitesimal character $\xi$, then we have $H^{\ell}(\fn_{I},M)\in \cO^{\fb_I}_{\fl_I,\rm{alg}}$ for each $\ell\geq 0$ (cf. \cite[Prop.~3.1.5]{BQ24}) and can replace $H^{\ell_2}(\fn_{I},M)$ in (\ref{equ: g spectral seq}) with $H^{\ell_2}(\fn_{I},M)_{\xi}$.

Recall from (\ref{equ: partial coxeter}) that $\Gamma\subseteq W(G)$ is the set of partial-Coxeter elements.
\begin{lem}\label{lem: g Ext1}
Let $x,w\in W(G)$. We have the following results.
\begin{enumerate}[label=(\roman*)]
\item \label{it: g Ext1 1} We have a natural isomorphism $\mathrm{Ext}_{U(\fg)}^1(L(x),L(w))\cong \mathrm{Ext}_{U(\fg)}^1(L(w),L(x))$.
\item \label{it: g Ext1 2} If $x<w$ and $\ell(w)=\ell(x)+1$, then we have $\Dim_E \mathrm{Ext}_{U(\fg)}^1(L(w),L(x))=1$.
\item \label{it: g Ext1 3} If $x<w$, $\ell(w)\geq \ell(x)+2$ and $w\in\Gamma$, then we have $\mathrm{Ext}_{U(\fg)}^1(L(w),L(x))=0$.
\end{enumerate}
\end{lem}
\begin{proof}
\ref{it: g Ext1 1} is \cite[(117)]{BQ24}. \ref{it: g Ext1 2} is \cite[Lem.~3.2.2 (iii)]{BQ24}.
\ref{it: g Ext1 3} follows from \cite[Lem.~3.2.2 (ii)]{BQ24} and the proof of \cite[Lem.~A.12 (i)]{BQ24} (which implies that $P_{x,w}=P_{w_0xw_0,w_0ww_0}=1$ and thus $x\prec w$ if and only if $\ell(w)=\ell(x)+1$).
\end{proof}

\begin{lem}\label{lem: n coh collection}
Let $I\subseteq \Delta$ and $x\in W(G)$. We have the following results.
\begin{enumerate}[label=(\roman*)]
\item \label{it: n coh 1} We have $H^0(\fn_I,L(x))\cong L^I(x)$.
\item \label{it: n coh 2} If $L(x)\in\cO^{\fp_{J}}_{\rm{alg}}$ for some $J\subseteq \Delta$, then $H^k(\fn_I,L(x))\in \cO^{\fl_I\cap\fp_J}_{\fl_I,\rm{alg}}$ for each $k\geq 0$.
\item \label{it: n coh 3} For each $k\geq 0$, $H^k(\fn_I,L(x))$ is a finite length $U(\fl_I)$-module with each constituent of the form $L^I(x')$ for some $x'\in W(G)$ satisfying $x'\geq x$.
\item \label{it: n coh 4} For each $k\geq 1$ and each constituent $L^I(x')$ of $H^k(\fn_I,L(x))$ (for some $x'\in W(G)$), we have $x'>x$ and $x'\notin W(L_I)x$.
\item \label{it: n coh 5} We have $H^k(\fn_I,L(1))\cong \bigoplus_{x'\in W^{I,\emptyset},\ell(x')=k}L^I(x')$.
\item \label{it: n coh 6} Let $w\in W(L_I)x$ such that either $x<w$ with $\ell(w)=\ell(x)+1$ or $x>w$ with $\ell(x)=\ell(w)+1$. Let $M$ (resp.~$M_I$) be the unique length $2$ $U(\fg)$-module with socle $L(x)$ and cosocle $L(w)$ (resp.~the unique length $2$ $U(\fl_I)$-module with socle $L^I(x)$ and cosocle $L^I(w)$) by Lemma~\ref{lem: g Ext1}. Then we have $H^0(\fn_I,M)\cong M_I$.
\end{enumerate}
\end{lem}
\begin{proof}
\ref{it: n coh 1} is \cite[Lem.~3.1.8 (ii)]{BQ24}. \ref{it: n coh 2} follows from \cite[Lem.~3.1.8 (i), Prop.~3.1.5]{BQ24}. \ref{it: n coh 5} is a well-known theorem of Kostant (cf.~\cite[Th\'eor\`eme~4.10]{Schr11}) 

We prove \ref{it: n coh 3}.
It follows from \cite[Lem.~3.1.4, Prop.~3.1.5]{BQ24} that for each $k\geq 0$ $H^k(\fn_I,L(x))$ is a finite length $U(\fl_I)$-module with each constituent of the form $L^I(x')$ for some $x'\in W(G)$. As $x'\cdot 0-x\cdot 0\in\Z_{\geq 0}\Phi^+$ if and only if $x'\geq x$, we deduce $x'\geq x$ from \cite[Lem.~3.1.2]{BQ24} and thus finish the proof of \ref{it: n coh 3}.

We prove \ref{it: n coh 4}.
Let $\xi: Z(\fl_I)\rightarrow E$ be the unique infinitesimal character such that $L^I(x)_{\xi}\neq 0$. It follows from \cite[Lem.~3.1.8 (iii)]{BQ24} and \ref{it: n coh 1} above that $H^k(\fn_I,L(x))_{\xi}=0$. In particular, we have $L^I(x')_{\xi}=0$ and thus $x'\notin W(L_I)x$ by \cite[Thm.~1.10]{Hum08}. This together with $x'\geq x$ from \ref{it: n coh 3} above forces $x'>x$.

We prove \ref{it: n coh 5}.
Let $\xi: Z(\fl_I)\rightarrow E$ be the unique infinitesimal character such that $L^I(x)_{\xi}\neq 0\neq L^I(w)_{\xi}$.
By \ref{it: n coh 4} (and its proof) we know that $H^k(\fn_I,L(x))_{\xi}=0=H^k(\fn_I,L(w))_{\xi}=0$ for $k\geq 1$ and in particular $H^0(\fn_I,M)$ fits into the following short exact sequence (using $H^0(\fn_I,L(x))\cong L^I(x)$ and $H^0(\fn_I,L(w))\cong L^I(w)$ from \ref{it: n coh 1})
\begin{equation}\label{equ: n coh short}
0\rightarrow L^I(x)\rightarrow H^0(\fn_I,M)\rightarrow L^I(w)\rightarrow 0.
\end{equation}
If (\ref{equ: n coh short}) is non-split, then we must have $H^0(\fn_I,M)\cong M_I$ by Lemma~\ref{lem: g Ext1}. Assume on the contrary that (\ref{equ: n coh short}) is split, then by (\ref{equ: g spectral seq 0}) we have
\[\Hom_{U(\fg)}(M^I(w),M)\cong \Hom_{U(\fl_I)}(L^I(w),H^0(\fn_I,M))\neq 0.\]
Since $[M^I(w):L(w)]=1=[M:L(w)]$ with both $M^I(w)$ and $M$ having cosocle $L(w)$, any non-zero map $M^I(w)\rightarrow M$ is a surjection, and thus $L(x)\in\mathrm{JH}_{U(\fg)}(M^I(w))$ which (by our assumption on $x,w$) forces $x>w$ and $\ell(x)=\ell(w)+1$. But this together with \cite[Lem.~3.3.10]{BQ24} contradicts our assumption that (\ref{equ: n coh short}) is split. The proof is thus finished.
\end{proof}

\begin{lem}\label{lem: Lie Ext Verma vanishing}
Let $I_0,I_1\subseteq \Delta$. If $I_1\not\subseteq I_0$, then we have $\mathrm{Ext}_{U(\fg)}^k(M^{I_1}(1),M^{I_0}(1))=0$ for each $k\geq 0$
\end{lem}
\begin{proof}
Let $\xi: Z(\fl_{I_1})\rightarrow E$ be the unique infinitesimal character such that $L^{I_1}(1)_{\xi}\neq 0$.
It follows from \ref{it: n coh 4} of Lemma~\ref{lem: n coh collection} that we have $H^k(\fn_{I_1},L(x))_{\xi}=0$ for each $L(x)\in\mathrm{JH}_{U(\fg)}(M^{I_0}(1))$ and $k\geq 1$. A d\'evissage with respect to $\mathrm{JH}_{U(\fg)}(M^{I_0}(1))$ gives $H^k(\fn_{I_1},M^{I_0}(1))_{\xi}=0$ for $k\geq 1$, which together with (\ref{equ: g spectral seq}) gives
\begin{equation}\label{equ: Lie Ext Verma seq}
\mathrm{Ext}_{U(\fg)}^k(M^{I_1}(1),M^{I_0}(1))\cong \mathrm{Ext}_{U(\fl_{I_1})}^k(L^{I_1}(1),H^0(\fn_{I_1},M^{I_0}(1))_{\xi})
\end{equation}
for each $k\geq 0$.
For each subset $\Lambda'\subseteq \Lambda$ and $M\in\mathrm{Mod}_{U(\fg)}$, we write $M_{\Lambda'}\defeq \bigoplus_{\mu\in\Lambda'}M_{\mu}$ for short.
We have the following decomposition
\begin{equation}\label{equ: Levi wt decomposition}
M^{I_0}(1)\cong M^{I_0}(1)_{\Z\Phi_{I_1}^+}\oplus M^{I_0}(1)_{\Lambda\setminus \Z\Phi_{I_1}^+}
\end{equation}
which is clearly $U(\fl_{I_1})$-equivariant.
Under the $U(\ft)$-equivariant isomorphism $M^{I_0}(1)\cong U(\fn_{I_0}^+)$, we have
\[M^{I_0}(1)_{\Z\Phi_{I_1}^+}=U(\fl_{I_1}\cap\fn_{I_0}^+)\cong U(\fl_{I_1})\otimes_{U(\fl_{I_1}\cap\fp_{I_0})}L^{I_0\cap I_1}(1).\]
In other words, the embedding $U(\fl_{I_1})\hookrightarrow U(\fg)$ induces the following isomorphism
\begin{equation}\label{equ: Levi Verma wt}
U(\fl_{I_1})\otimes_{U(\fl_{I_1}\cap\fp_{I_0})}L^{I_0\cap I_1}(1)\buildrel\sim\over\longrightarrow M^{I_0}(1)_{\Z\Phi_{I_1}^+}.
\end{equation}
Given a $U(\fl_{I_1})$-module $M_{I_1}$, we note that $M_{I_1,\xi}\neq 0$ only if $M_{I_1,\mu}\neq 0$ for some $\mu\in\Z\Phi_{I_1}^+$. Hence, we have $(M^{I_0}(1)_{\Lambda\setminus \Z\Phi_{I_1}^+})_{\xi}=0$, which together with (\ref{equ: Levi wt decomposition}), (\ref{equ: Levi Verma wt}) and the fact that
\[M^{I_0}(1)_{\Z\Phi_{I_1}^+} \subseteq H^0(\fn_{I_1},M^{I_0}(1))\]
forces
\[H^0(\fn_{I_1},M^{I_0}(1))_{\xi}\cong U(\fl_{I_1})\otimes_{U(\fl_{I_1}\cap\fp_{I_0})}L^{I_0\cap I_1}(1).\]
Concerning the isomorphism (\ref{equ: Lie Ext Verma seq}), it is harmless to assume from now that $I_1=\Delta\neq I_0$.
Since $I_0\neq \Delta$, we know that $\Hom_{E}(\wedge^{\Dim_E\fp_{I_0}}\fp_{I_0},E)\cong L^{I_0}(\mu)$ for some $0\neq \mu\in \Lambda_{I_0}^{\dom}$.
Note that $\mathrm{Ext}_{U(\fg)}^{\bullet}(-,U(\fg))$ defines an anti-involution of $D^{b}(\mathrm{Mod}_{U(\fg)})$, with
\[\mathrm{Ext}_{U(\fg)}^{\bullet}(L(1),U(\fg))\cong L(1)[-\Dim_E\fg]\]
and
\[\mathrm{Ext}_{U(\fg)}^{\bullet}(M^{I_0}(1),U(\fg))\cong M^{I_0}(\mu)[-\Dim_E\fp_{I_0}].\]
Since $0\neq \mu\in \Lambda_{I_0}^{\dom}$, we know that $L(1)$ and $M^{I_0}(\mu)$ have distinct infinitesimal character, and thus the anti-involution $\mathrm{Ext}_{U(\fg)}^{\bullet}(-,U(\fg))$ induces an isomorphism
\[\mathrm{Ext}_{U(\fg)}^{k}(L(1),M^{I_0}(1))\cong \mathrm{Ext}_{U(\fg)}^{k-\Dim_E\fg+\Dim_E\fp_{I_0}}(M^{I_0}(\mu),L(1))=0\]
for $k\in \Z$. The proof is thus finished.
\end{proof}

\begin{lem}\label{lem: Lie Verma bottom Ext}
Let $u\in W(G)$. We have the following results.
\begin{enumerate}[label=(\roman*)]
\item \label{it: Lie Verma bottom Ext 1} The $E$-vector space $\mathrm{Ext}_{U(\fg)}^{k}(L(u),L(1))$ is zero when $k<\ell(u)$ and $1$-dimensional when $k=\ell(u)$.
\item \label{it: Lie Verma bottom Ext 2} For each $I'\subseteq I\subseteq \Delta$, the surjection $M^{I'}(u)\twoheadrightarrow M^{I}(u)$ induces an isomorphism
    \begin{equation}\label{equ: Lie 2 Verma bottom Ext}
    \mathrm{Ext}_{U(\fg)}^{\ell(u)}(M^{I}(u),L(1))\buildrel\sim\over\longrightarrow \mathrm{Ext}_{U(\fg)}^{\ell(u)}(M^{I'}(u),L(1))
    \end{equation}
    between $1$-dimensional spaces.
\end{enumerate}
\end{lem}
\begin{proof}
Both \ref{it: Lie Verma bottom Ext 1} and \ref{it: Lie Verma bottom Ext 2} are well-known.\\
Recall from Bott's formula (cf.~\cite[\S 6.6]{Hum08}) that
\[H^{\ell}(\fu,L(1))\cong\bigoplus_{w\in W(G),\ell(w)=\ell}w\cdot 0,\]
which together with (\ref{equ: g spectral seq}) implies that the $E$-vector space $\mathrm{Ext}_{U(\fg)}^{k}(M(u),L(1))$ is zero when $k<\ell(u)$ and $1$-dimensional when $k=\ell(u)$. Now that each $L(u')\in\mathrm{JH}_{U(\fg)}(N(u))$ satisfies $u'>u$, a decreasing induction on $\ell(u)$ together with d\'evissage with respect to $0\rightarrow N(u)\rightarrow M(u)\rightarrow L(u)\rightarrow 0$ as well as d\'evissage with respect to $\mathrm{JH}_{U(\fg)}(N(u))$ gives \ref{it: Lie Verma bottom Ext 1}. 
Using d\'evissage with respect to $0\rightarrow N^{I}(u)\rightarrow M^{I}(u)\rightarrow L(u)\rightarrow 0$ and the fact that each $L(u')\in\mathrm{JH}_{U(\fg)}(N^{I}(u))$ satisfies $u'>u$, it is clear from \ref{it: Lie Verma bottom Ext 1} that $\mathrm{Ext}_{U(\fg)}^{\ell(u)}(M^{I}(u),L(1))$ is $1$-dimensional for each $I\subseteq\Delta$. We deduce the isomorphism (\ref{equ: Lie 2 Verma bottom Ext}) from \ref{it: Lie Verma bottom Ext 1} and the fact that each constituent $L(u')$ of $\mathrm{ker}(M^{I'}(u)\twoheadrightarrow M^{I}(u))$ satisfies $u'>u$.
\end{proof}

\begin{lem}\label{lem: g induction quotient}
Let $w\in W(L_I)x$ such that either $x<w$ with $\ell(w)=\ell(x)+1$ or $x>w$ with $\ell(x)=\ell(w)+1$. Let $M$ (resp.~$M_I$) be the unique length $2$ $U(\fg)$-module with socle $L(x)$ and cosocle $L(w)$ (resp.~the unique length $2$ $U(\fl_I)$-module with socle $L^I(x)$ and cosocle $L^I(w)$) by Lemma~\ref{lem: g Ext1}. Then we have a surjection
\begin{equation}\label{equ: g induction quotient}
q: U(\fg)\otimes_{U(\fp_I)}M_I\rightarrow M.
\end{equation}
\end{lem}
\begin{proof}
It follows from \ref{it: n coh 6} of Lemma~\ref{lem: n coh collection} and (\ref{equ: g spectral seq 0}) that
\[\Hom_{U(\fg)}(U(\fg)\otimes_{U(\fp_I)}M_I,M)\cong \Hom_{U(\fl_J)}(M_I,H^0(\fn_I,M))\neq 0\]
namely we have a non-zero map of the form (\ref{equ: g induction quotient}).
As we have
\[\Hom_{U(\fg)}(U(\fg)\otimes_{U(\fp_I)}M_I,L(x))\cong \Hom_{U(\fl_I)}(M_I,H^0(\fn_I,L(x)))=\Hom_{U(\fl_I)}(M_I,L^I(x))=0\]
from \ref{it: n coh 1} of Lemma~\ref{lem: n coh collection} and (\ref{equ: g spectral seq 0}), we see that $q$ must be a surjection.
\end{proof}

\begin{lem}\label{lem: Ext Verma w}
Let $I\subseteq \Delta$ and $x,w\in W(G)$ with $D_L(x)\cap I=\emptyset$ and $w^{I}$ being the unique minimal element in $W(L_{I})w$. If $w^{I}\not\leq x$, then
\begin{equation}\label{equ: Ext Verma w}
\mathrm{Ext}_{U(\fg)}^{k}(M^{I}(x),L(w))
\end{equation}
is zero for each $k\geq 0$.
\end{lem}
\begin{proof}
Thanks to (\ref{equ: g spectral seq}) we have the following spectral sequence
\begin{equation}\label{equ: Ext Verma w seq}
\mathrm{Ext}_{U(\fl_{I})}^{k}(L^{I}(x),H^{\ell}(\fn_{I},L(w)))\implies \mathrm{Ext}_{U(\fp_{I})}^{k+\ell}(L^{I}(x),L(w))\cong \mathrm{Ext}_{U(\fg)}^{k+\ell}(M^{I}(x),L(w)).
\end{equation}
We write $\xi: Z(\fl_{I})\rightarrow E$ for the unique infinitesimal character such that $L^{I}(x)_{\xi}\neq 0$.
Assume on the contrary that (\ref{equ: Ext Verma w}) is non-zero for some $k\geq 0$, then by Lemma~\ref{lem: g block} we know that $H^{\ell}(\fn_{I},L(w))_{\xi}\neq 0$ for some $\ell\geq 0$. In particular, there exists $L^{I}(w')\in\mathrm{JH}_{U(\fl_{I})}(H^{\ell}(\fn_{I},L(w)))$ that satisfies $L^{I}(w')_{\xi}\neq 0$ or equivalently $w'\in W(L_{I})x$ by \cite[Thm.~1.10]{Hum08}. It follows from \ref{it: n coh 3} of Lemma~\ref{lem: n coh collection} that any $L^{I}(w')\in\mathrm{JH}_{U(\fl_{I})}(H^{\ell}(\fn_{I},L(w)))$ satisfies $w'\geq w$, which together with $w'\in W(L_{I})x$ (and $x$ being the unique minimal element in $W(L_{I})x$ by $D_L(x)\cap I=\emptyset$) forces $x\geq w^{I}$, a contradiction.
\end{proof}

\begin{rem}\label{rem: Ext dominant Verma w}
Note that $w^{I}=1$ if and only if $w\in W(L_{I})$ if and only if $\mathrm{Supp}(w)\subseteq I$.
Consequently, when $x=1$, Lemma~\ref{lem: Ext Verma w} specializes to the result that
\[\mathrm{Ext}_{U(\fg)}^{k}(M^{I}(1),L(w))=0\]
for each $k\geq 0$ whenever $\mathrm{Supp}(w)\not\subseteq I$.
\end{rem}

\begin{lem}\label{lem: coxeter subquotient}
Let $I\subseteq \Delta$. We have the following results.
\begin{enumerate}[label=(\roman*)]
\item \label{it: g coxeter 1} For each $x\in\Gamma$, we have $[M^I(1):L(x)]\neq 0$ if and only if $\mathrm{Supp}(x)\subseteq \Delta\setminus I$, in which case $[M^I(1):L(x)]=1$.
\item \label{it: g coxeter 2} For each $x,w\in\Gamma$ with $x\leq w$ and $\mathrm{Supp}(x),\mathrm{Supp}(w)\subseteq \Delta\setminus I$, $M^I(1)$ admits a unique subquotient $M$ with socle $L(w)$ and cosocle $L(x)$ such that $L(x')\in\mathrm{JH}_{U(\fg)}(M)$ for some $x'\in\Gamma$ if and only if $x\leq x'\leq w$. Moreover, $M$ has length $2$ if $\ell(w)=\ell(x)+1$.
\item \label{it: g coxeter 3} The $M$ in \ref{it: g coxeter 2} is a quotient of $M^I(x)$ and depends only on the choice of $x$ and $w$, but not on the choice of $I$.
\end{enumerate}
\end{lem}
\begin{proof}
We divide the proof into the following steps.

\textbf{Step $1$}: We prove \ref{it: g coxeter 1}.\\
It follows from \cite[Lem.~A.12 (i)]{BQ24} that, given $x'\in W(G)$, we have $[M(x'):L(x)]\neq 0$ if and only if $x'\leq x$ in which case $[M(x'):L(x)]=1$. This together with \cite[Prop.~9.6]{Hum08} implies that
\[[M^I(1):L(x)]=\sum_{x'\in W(L_I)}(-1)^{\ell(x')}[M(x'):L(x)]=\sum_{x'\in W(L_I), x'\leq x}(-1)^{\ell(x')}.\]
As $x\in\Gamma$, any $x'\leq x$ is uniquely determined by its support $\mathrm{Supp}(x')$ (together with $x$) by \cite[Thm.~2.2.2]{BB05}. Hence $[M^I(1):L(x)]\neq 0$ if and only if $x'=1$ is the only $x'\in W(L_I)$ that satisfies $x'\leq x$, if and only if $\mathrm{Supp}(x)\cap I=\emptyset$.

\textbf{Step $2$}: We prove that $M(1)$ admits a unique subquotient $M$ with socle $L(w)$ and cosocle $L(x)$ such that $L(x')\in\mathrm{JH}_{U(\fg)}(M)$ for some $x'\in\Gamma$ if and only if $x\leq x'\leq w$, and moreover that $M$ has length $2$ when $\ell(w)=\ell(x)+1$.\\
By \cite[Thm.~5.1]{Hum08} there exist embeddings $M(w)\subseteq M(x)\subseteq M(1)$.
As we have $[M(x):L(w)]=1$ from \cite[Lem.~A.12 (i)]{BQ24}, we may define $M$ as the unique quotient of $M(x)$ with socle $L(w)$ (and cosocle $L(x)$), which is naturally a subquotient of $M(1)$. As $[M(1):L(x)]=[M(1):L(w)]=1$ from \emph{loc.cit.}, $M$ is also the unique subquotient of $M(1)$ with socle $L(w)$ and cosocle $L(x)$.
For each $x'\in\Gamma$, using \cite[Thm.~5.1]{Hum08} as well as $[M(1):L(x')]=1$ from \cite[Lem.~A.12 (i)]{BQ24}, we see that $M(x')$ is the unique $U(\fg)$-submodule of $M(1)$ with cosocle $L(x')$. From \cite[Thm.~5.1]{Hum08} we see that $M(w)\subseteq M(x')$ (resp.~$M(x')\subseteq M(x)$) if and only if $w\geq x'$ (resp.~$x'\geq x$). Consequently, as $M$ is the unique subquotient of $M(1)$ with socle $L(w)$ and cosocle $L(x)$, we have $L(x')\in\mathrm{JH}_{U(\fg)}(M)$ if and only if $M(w)\subseteq M(x')\subseteq M(x)\subseteq M(1)$ if and only if $x\leq x'\leq w$.

\textbf{Step $3$}: We prove \ref{it: g coxeter 2} and \ref{it: g coxeter 3}.\\
Recall from \cite[Thm.~9.4 (c)]{Hum08} that $M^I(1)$ is a quotient of $M(1)$. As $M$ is the unique subquotient of $M(1)$ with socle $L(w)$ and cosocle $L(x)$ and we have $[M^I(1):L(x)]=[M^I(1):L(w)]=1$ from \ref{it: g coxeter 1}, $M$ is also necessarily the unique subquotient of $M^I(1)$ with socle $L(w)$ and cosocle $L(x)$, which together with \textbf{Step $2$} gives \ref{it: g coxeter 2}. As our definition of $M$ in \textbf{Step $2$} is independent of the choice of $I$, we deduce \ref{it: g coxeter 3}.
\end{proof}

Let $x\in W(G)$ with $J_{x}\defeq\Delta\setminus\mathrm{Supp}(x)$.
For each $I\subseteq \Delta$ and $I''\subseteq I'\subseteq I\cap J_{x}$, there exists a surjection (cf.~\cite[Thm.~9.4(c)]{Hum08})
\begin{equation}\label{equ: Verma transfer surjection}
q^{x,I}_{I'',I'}: U(\fl_{I})\otimes_{U(\fl_{I}\cap\fp_{I''})}L^{I''}(x)\twoheadrightarrow U(\fl_{I})\otimes_{U(\fl_{I}\cap\fp_{I'})}L^{I'}(x),
\end{equation}
For each $I\subseteq \Delta$ and $I'\subseteq I\cap J_{x}$, we set
\begin{equation}\label{equ: Lie St x}
\mathfrak{v}_{x,I',I}\defeq \bigcap_{I'\subsetneq I''\subseteq I\cap J_{x}}\mathrm{ker}(q^{x,I}_{I',I''})\subseteq U(\fl_{I})\otimes_{U(\fl_{I}\cap\fp_{I'})}L^{I'}(x).
\end{equation}
When $x=1$ (resp.~$I=\Delta$), we omit it from the notation and obtain $\mathfrak{v}_{I',I}$ (resp.~$\mathfrak{v}_{x,I'}$).

Let $x\in W(G)$, $I\subseteq \Delta$ and $I_0\subseteq I_1\subseteq I\cap J_{x}$.
For each $I_0\subseteq I'\subseteq I_1$ and $j\in I'\setminus I_0$, we define
\begin{equation}\label{equ: differential sign}
m(I',j)\defeq \#\{j'\in I'\mid j'<j\}\in\Z_{\geq 0}.
\end{equation}
We consider a complex in $\mathrm{Mod}_{U(\fl_{I})}$ of the form
\begin{equation}\label{equ: Lie x Tits complex}
\mathfrak{c}^{x,I}_{I_0,I_1}\defeq [U(\fl_{I})\otimes_{U(\fl_{I}\cap\fp_{I_0})}L^{I_0}(x)\rightarrow \cdots  \rightarrow U(\fl_{I})\otimes_{U(\fl_{I}\cap\fp_{I_1})}L^{I_1}(x)]
\end{equation}
whose degree $\ell$ term is given by
\[\bigoplus_{I_0\subseteq I'\subseteq I_1, \#I'=\ell}U(\fl_{I})\otimes_{U(\fl_{I}\cap\fp_{I'})}L^{I'}(x),\]
with the differential map at degree $\ell-1$ restricts to
\[(-1)^{m(I',j)}q^{x,I}_{I'\setminus\{j\},I'}: U(\fl_{I})\otimes_{U(\fl_{I}\cap\fp_{I'\setminus\{j\}})}L^{I'\setminus\{j\}}(x)\rightarrow U(\fl_{I})\otimes_{U(\fl_{I}\cap\fp_{I'})}L^{I'}(x)\]
for each $I_1\subseteq I'\subseteq I_0$ and $j\in I'\setminus I_1$ satisfying $\#I'=\ell$.
Then we have
\begin{equation}\label{equ: Lie general Tits resolution}
\mathfrak{c}^{x,I}_{I_0,I_1}\cong U(\fl_{I})\otimes_{U(\fl_{I}\cap\fp_{I_1})}(\mathfrak{v}_{x,I_0,I_1})[-\#I_0]
\end{equation}
in $D^{b}(\mathrm{Mod}_{U(\fl_{I})})$. 
When $x=1$ (resp.~$I=\Delta$), we omit it from the notation and obtain $\mathfrak{c}^{I}_{I_0,I_1}$ (resp.~$\mathfrak{c}^{x}_{I_0,I_1}$).
We define $\mathfrak{c}_{I_0,I_1}\defeq \mathfrak{c}^{1,\Delta}_{I_0,I_1}$ similarly.

\begin{rem}\label{rem: choice of sign}
We consider the complex $[L(1)\buildrel\mathrm{Id}\over\longrightarrow L(1)]$ with degree supported in $[0,1]$. By our definition of $\mathfrak{c}_{\emptyset,\Delta}$, we have an embedding between complex
\[\mathfrak{c}_{\emptyset,\Delta}\rightarrow \mathrm{Tot}([L(1)\buildrel\mathrm{Id}\over\longrightarrow L(1)]^{\otimes \#\Delta})\]
induced from the natural surjection $q^{1,\Delta}_{I,\Delta}: M^{I}(1)\twoheadrightarrow L(1)$ for each $I\subseteq \Delta$.
\end{rem}

\begin{lem}\label{lem: mult suppot}
Let $x\in W(G)$. Then the quantity $[M^{I}(1):L(x)]$ depends only on $I\cup J_{x}$.
\end{lem}
\begin{proof}
We prove by an increasing induction on $\ell(x)$. It is clear that $[M^{I}(1):L(1)]=1$ for each $I\subseteq \Delta=J_{1}$.
In general, we notice that
\[\Dim_E M^{I}(1)_{x\cdot0}=\sum_{y\leq x}[M^{I}(1):L(y)]\Dim_EL(y)_{x\cdot 0}\]
which together with $\Dim_EL(x)_{x\cdot 0}=1$ gives
\begin{equation}\label{equ: mult supprt formula}
[M^{I}(1):L(x)]=\Dim_E M^{I}(1)_{x\cdot0}-\sum_{y<x}[M^{I}(1):L(y)]\Dim_EL(y)_{x\cdot 0}.
\end{equation}
For each $y<x$, by our induction hypothesis we know that $[M^{I}(1):L(y)]$ depends only on $I\cup J_{y}$, and in particular depends only on $I\cup J_{x}$. 
Using PBW basis, we see that the surjection $M^{I}(1)\twoheadrightarrow M^{I\cup J_{x}}(1)$ induces an isomorphism
\[M^{I}(1)_{x\cdot0}\buildrel\sim\over\longrightarrow M^{I\cup J_{x}}(1)_{x\cdot0}\]
between $E$-vector spaces and thus $\Dim_E M^{I}(1)_{x\cdot0}$ depends only on $I\cup J_{x}$.
We thus conclude from (\ref{equ: mult supprt formula}) that $[M^{I}(1):L(x)]$ depends only on $I\cup J_{x}$.
\end{proof}

\begin{lem}\label{lem: coxeter parabolic Verma}
Let $x,w\in\Gamma$ and $I\subseteq J_{x}$ with $x\unlhd w$ and $I\not\subseteq J_{w}$.
Then we have $L(w)\notin\mathrm{JH}_{U(\fg)}(M^{I}(x))$.
\end{lem}
\begin{proof}
Since we have $I\subseteq J_{x}$ and $I\not\subseteq J_{w}$, we can choose an arbitrary $j\in I\cap(\mathrm{Supp}(w)\setminus\mathrm{Supp}(x))$. Our assumption $x\unlhd w$ together with \cite[Thm.~2.2.2]{BB05} ensures that $x<s_jx\leq w$ and $j\in D_L(s_jx)\cap I$ (with $s_jx\in\Gamma$). Now that $M^{I}(x)$ is locally $\fp_{I}$-finite but $L(s_jx)$ is not (see Lemma~\ref{lem: dominance and left set}), we see that $L(s_jx)\notin\mathrm{JH}_{U(\fg)}(M^{I}(x))$. Now that $M(s_jx)$ has cosocle $L(s_jx)$, the composition of
\[M(s_jx)\hookrightarrow M(x)\twoheadrightarrow M^{I}(x)\]
must be zero and thus $M^{I}(x)$ is a quotient of $M(x)/M(s_jx)$. Now that we have
\[[M(s_jx):L(w)]=[M(x):L(w)]=[M(1):L(w)]=1\]
by Lemma~\ref{lem: coxeter subquotient}, we deduce that $[M(x)/M(s_jx):L(w)]=0$ and thus $L(w)\notin\mathrm{JH}_{U(\fg)}(M^{I}(x))$ (as $M^{I}(x)$ is a quotient of $M(x)/M(s_jx)$).
The proof is thus finished.
\end{proof}

\begin{lem}\label{lem: Lie St JH}
Let $x,w\in W(G)$ and $I_0\subseteq J_{x}$. 
We have the following results.
\begin{enumerate}[label=(\roman*)]
\item \label{it: Lie St JH 1} If $[\mathfrak{v}_{x,I_0}:L(w)]\neq 0$, then we have $\mathrm{Supp}(w)\supseteq \Delta\setminus I_0$. If furthermore $\ell(w)=\#\Delta\setminus I_0$, then we have $w\in\Gamma_{\Delta\setminus I_0}$ with $x\unlhd w$.
\item \label{it: Lie St JH 2} If $x=1$ and $w\in\Gamma$, then we have $[\mathfrak{v}_{I_0}:L(w)]\neq 0$ if and only if $\mathrm{Supp}(w)=\Delta\setminus I_0$, in which case $[\mathfrak{v}_{I_0}:L(w)]=1$.
\end{enumerate}
\end{lem}
\begin{proof}
Let $I_0\subseteq I_1\subseteq J_{x}$ and write $\xi: Z(\fl_{I_1})\rightarrow E$ for the unique infinitesimal character that satisfies $L^{I_1}(x)_{\xi}\neq 0$.
We divide the proof into the following steps.

\textbf{Step $1$}: We prove that $[\mathfrak{v}_{x,I_0,I_1}:L^{I_1}(w)]\neq 0$ only if $x\unlhd w$ and $\mathrm{Supp}(w)\supseteq \mathrm{Supp}(x)\sqcup(I_1\setminus I_0)$.\\
Now that $\mathfrak{v}_{x,I_0,I_1}$ is a $U(\fl_{I_1})$-submodule of $U(\fl_{I_1})\otimes_{U(\fl_{I_1}\cap\fp_{I_0})}L^{I_0}(x)$, $Z(\fl_{I_1})$ acts on $\mathfrak{v}_{x,I_0,I_1}$ by $\xi$. Hence, it suffices to consider $w\in W(G)$ that satisfies $L^{I_1}(w)_{\xi}\neq 0$, which (by \cite[Thm.~1.10]{Hum08} and $I_1\subseteq J_{x}$) holds if and only if $w\in W(L_{I_1})x$ if and only if $x\unlhd w$ and $\mathrm{Supp}(w)\subseteq \mathrm{Supp}(x)\sqcup I_1$. 
Upon applying the translation functor from the block $\cO^{\fb\cap\fl_{I_1}}_{\rm{alg},\fl_{I_1},\xi}$ to the block $\cO^{\fb\cap\fl_{I_1}}_{\rm{alg},\fl_{I_1},0}$, it is harmless to assume in the rest of \textbf{Step $1$} that $x=1$ and $w\in W(L_{I_1})$. 
By the very definition of $\mathfrak{v}_{I_0,I_1}=\mathfrak{v}_{1,I_0,I_1}$, we have
\begin{equation}\label{equ: Lie x St JH mult}
[\mathfrak{v}_{I_0,I_1}: L^{I_1}(w)]=\sum_{I_0\subseteq I\subseteq I_1}(-1)^{\#I\setminus I_0}[U(\fl_{I_1})\otimes_{U(\fl_{I_1\cap\fp_{I}})}L^{I}(1):L^{I_1}(w)].
\end{equation}
If $\mathrm{Supp}(w)\not\supseteq I_1\setminus I_0$ or equivalently $J_{w}\cap(I_1\setminus I_0)\neq \emptyset$, then we have
\begin{equation}\label{equ: Lie x St JH mult cancel}
\sum_{J\subseteq I\subseteq J\sqcup (J_{w}\cap(I_1\setminus I_0))}(-1)^{\#I\setminus I_0}[U(\fl_{I_1})\otimes_{U(\fl_{I_1\cap\fp_{I}})}L^{I}(1):L^{I_1}(w)]=0
\end{equation}
for each $I_0\subseteq J\subseteq I_0\sqcup(I_1\setminus J_{w})$, and thus (\ref{equ: Lie x St JH mult}) equals zero (by summing up (\ref{equ: Lie x St JH mult cancel}) over all possible $J$).
In other words, we have shown that (\ref{equ: Lie x St JH mult}) is non-zero only if $\mathrm{Supp}(w)\supseteq I_1\setminus I_0$ (when $x=1$).

\textbf{Step $2$}: We prove that $[\mathfrak{v}_{x,I_0,I_1}:L^{I_1}(w)]\neq 0$ for some $w\in\Gamma$ if and only if $x\unlhd w$ and $\mathrm{Supp}(w)=\mathrm{Supp}(x)\sqcup(I_1\setminus I_0)$.\\
It is harmless to assume that $x=1$ and $w\in W(L_{I_1})$ as in \textbf{Step $1$}.
For each $I_0\subseteq I\subseteq I_1$, we deduce from \ref{it: g coxeter 1} of Lemma~\ref{lem: coxeter subquotient} that 
\begin{equation}\label{equ: Lie x Verma JH}
[U(\fl_{I_1})\otimes_{U(\fl_{I_1\cap\fp_{I}})}L^{I}(1):L^{I_1}(w)]
\end{equation}
is non-zero if and only if $\mathrm{Supp}(w)\subseteq I_1\setminus I$, in which case it equals one.
If $\mathrm{Supp}(w)\not\subseteq I_1\setminus I_0$, then (\ref{equ: Lie x Verma JH}) is zero for each $I_0\subseteq I\subseteq I_1$ and thus (\ref{equ: Lie x St JH mult}) is zero.
If $\mathrm{Supp}(w)\subseteq I_1\setminus I_0$, then (\ref{equ: Lie x St JH mult}) equals 
\[\sum_{I_0\sqcup\mathrm{Supp}(w)\subseteq I\subseteq I_1}(-1)^{\#I\setminus I_0}\]
which is non-zero if and only if $\mathrm{Supp}(w)=I_1\setminus I_0$, in which case (\ref{equ: Lie x St JH mult}) equals one.

Note that \ref{it: Lie St JH 2} is the special case of \textbf{Step $2$} when $x=1$ and $I_1=J_{1}=\Delta$.

\textbf{Step $3$}: We finish the proof of \ref{it: Lie St JH 1}.\\
It follows from
\[\mathfrak{v}_{x,I_0}\cong U(\fg)\otimes_{U(\fl_{J_{x}})}\mathfrak{v}_{x,I_0,J_{x}}\]
and the exactness of $U(\fg)\otimes_{U(\fl_{J_{x}})}-$ that
\[[\mathfrak{v}_{x,I_0}:L(w)]=\sum_{x'}[\mathfrak{v}_{x,I_0,J_{x}}:L^{J_{x}}(x')][M^{J_{x}}(x'):L(w)].\]
In other words, for each $w\in W(G)$ satisfying $[\mathfrak{v}_{x,I_0}:L(w)]$, there exists $x'$ such that $[\mathfrak{v}_{x,I_0,J_{x}}:L^{J_{x}}(x')]\neq 0$ and $[M^{J_{x}}(x'):L(w)]\neq 0$.
Now that $[\mathfrak{v}_{x,I_0,J_{x}}:L^{J_{x}}(x')]\neq 0$ forces $x\unlhd x'$ and $\mathrm{Supp}(x')\supseteq \Delta\setminus I_0$ by \textbf{Step $1$} (replacing $w$ and $I_1$ with $x'$ and $J_{x}$ in \emph{loc.cit.}) and $[M^{J_{x}}(x'):L(w)]\neq 0$ forces $x'\leq w$, we conclude that $\mathrm{Supp}(w)\supseteq\mathrm{Supp}(x')\supseteq \Delta\setminus I_0$.
If furthermore $\ell(w)=\#\Delta\setminus I_0$, we deduce from $\mathrm{Supp}(w)\supseteq \Delta\setminus I_0$ that $w\in\Gamma_{\Delta\setminus I_0}$, which together with $x\unlhd x'\leq w$ (and $\ell(x')\geq \#\mathrm{Supp}(x')\geq \#\Delta\setminus I_0$) forces $x\unlhd x'=w\in\Gamma_{\Delta\setminus I_0}$.
\end{proof}

\begin{lem}\label{lem: Lie Tits cosocle}
Let $x\in\Gamma$ and $I_0\subseteq J_{x}=\Delta\setminus\mathrm{Supp}(x)$. Then we have
\begin{equation}\label{equ: Lie Tits cosocle}
\mathrm{cosoc}_{U(\fg)}(\mathfrak{v}_{x,I_0})=\bigoplus_{w\in\Gamma_{\Delta\setminus I_0}, x\unlhd w}L(w).
\end{equation}
\end{lem}
\begin{proof}
We consider an arbitrary $L(w)\in\mathrm{JH}_{U(\fg)}(\mathfrak{v}_{x,I_0})$.
Recall from (\ref{equ: Lie general Tits resolution}) that we have
\[\mathfrak{c}^{x,\Delta}_{I_0,J_{x}}\cong\mathfrak{v}_{x,I_0}[-\#I_0]\]
in $D^{b}(\mathrm{Mod}_{U(\fg)})$ and thus we have a spectral sequence of the form
\begin{equation}\label{equ: Levi Lie St resolution Ext}
\bigoplus_{I_0\subseteq I\subseteq J_{x}, \#I=\ell}\mathrm{Ext}_{U(\fg)}^{k}(M^{I_0}(x),L(w))\implies \mathrm{Ext}_{U(\fg)}^{k-\ell}(\mathfrak{c}^{x,\Delta}_{I_0,J_{x}},L(w))\cong\mathrm{Ext}_{U(\fg)}^{k-\ell+\#I_0}(\mathfrak{v}_{x,I_0},L(w)).
\end{equation}
It follows from \ref{it: Lie St JH 1} of Lemma~\ref{lem: Lie St JH} that $\mathrm{Supp}(w)\cup I_0=\Delta$ and thus
\begin{equation}\label{equ: Levi Lie St index union}
(\mathrm{Supp}(w)\setminus\mathrm{Supp}(x))\cup I_0=J_{x}.
\end{equation}
Note that any $I_0\subseteq I\subseteq J_{x}$ satisfies $w^{I}\leq x$ only if $\mathrm{Supp}(w)\setminus\mathrm{Supp}(x)\subseteq I$ which together with (\ref{equ: Levi Lie St index union}) forces $I=J_{x}$. Hence, we deduce from Lemma~\ref{lem: Ext Verma w} that
\[\mathrm{Ext}_{U(\fg)}^{k}(M^{I}(x),L(w))=0\]
for each $k\geq 0$ and $I_0\subseteq I\subsetneq J_{x}$, which together with the spectral sequence (\ref{equ: Levi Lie St resolution Ext}) gives an isomorphism
\begin{equation}\label{equ: Levi Lie St cosocle}
\Hom_{U(\fg)}(\mathfrak{v}_{x,I_0},L(w))\cong \mathrm{Ext}_{U(\fg)}^{\#J_{x}\setminus I_0}(M^{J_{x}}(x),L(w)).
\end{equation}
We write $\xi_{x}: Z(\fl_{J_{x}})\rightarrow E$ for the unique infinitesimal character such that $L^{J_{x}}(x)_{\xi_{x}}\neq 0$.
Recall from (\ref{equ: g spectral seq}) that we have the following spectral sequence
\begin{equation}\label{equ: Lie St cosocle spectral seq}
\mathrm{Ext}_{U(\fl_{J_{x}})}^{k}(L^{J_{x}}(x), H^{\ell}(\fn_{J_{x}},L(w))_{\xi_{x}})\implies \mathrm{Ext}_{U(\fg)}^{k+\ell}(M^{J_{x}}(x),L(w)).
\end{equation}
Suppose that $H^{\ell}(\fn_{J_{x}},L(w))_{\xi_{x}}\neq 0$ for some $\ell\geq 0$ and consider an arbitrary \[L^{J_{x}}(y)\in\mathrm{JH}_{U(\fl_{J_{x}})}(H^{\ell}(\fn_{J_{x}},L(w))_{\xi_{x}})\] 
which necessarily satisfies $y\in W(L_{J_{x}})x$.
If $\ell\geq 1$, then by \ref{it: n coh 4} of Lemma~\ref{lem: n coh collection} we have $w<y$ and $y\notin W(L_{J_{x}})w$, which together with $y\in W(L_{J_{x}})x$ and \cite[Prop.~2.5.1]{BB05} forces $w^{I}<x$ and thus $\mathrm{Supp}(w^{J_{x}})\subsetneq \mathrm{Supp}(x)$ (as $x\in\Gamma$). But this together with $w\in W(L_{J_{x}})w^{J_{x}}$ forces $\mathrm{Supp}(w)\not\supseteq\mathrm{Supp}(x)$ which contradicts $\mathrm{Supp}(w)\supseteq \Delta\setminus I_0\supseteq\mathrm{Supp}(x)$.
Hence, we must have $\ell=0$ and $L^{J_{x}}(w)\cong H^{0}(\fn_{J_{x}},L(w))_{\xi_{x}}\neq 0$ (see \ref{it: n coh 1} of Lemma~\ref{lem: n coh collection}), or equivalently $\ell=0$ and $w\in W(L_{J_{x}})x$ (see \cite[Thm.~1.10]{Hum08}).
This together with (\ref{equ: Lie St cosocle spectral seq}) and Lemma~\ref{lem: Lie Verma bottom Ext} implies that (\ref{equ: Levi Lie St cosocle}) is non-zero (for given $x\in\Gamma$, $I_0\subseteq J_{x}$ and $L(w)\in\mathrm{JH}_{U(\fg)}(\mathfrak{v}_{x,I_0})$ which necessarily satisfies $\mathrm{Supp}(w)\supseteq\Delta\setminus I_0$) if and only if $w\in W(L_{J_{x}})x$ and $\ell(w)=\ell(x)+\#J_{x}\setminus I_0=\#\Delta\setminus I_0$, if and only if $w\in\Gamma_{\Delta\setminus I_0}$ and $x\unlhd w$, in which case (\ref{equ: Levi Lie St cosocle}) is $1$-dimensional.
\end{proof}

Let $x,w\in \Gamma$ with $x\leq w$ and $I\subseteq J_{w}$.
It follows from \cite[\S 4.2, \S 4.6, Thm~9.4(c)]{Hum08} that we have the following commutative diagram
\begin{equation}\label{equ: Lie Verma wt shift}
\xymatrix{
M(w) \ar[r] \ar[d] & M(x) \ar[r] \ar[d] & M(1) \ar[d]\\
M^{I}(w) \ar[r] & M^{I}(x) \ar[r] & M^{I}(1)
}
\end{equation}
with all vertical maps being surjections and horizontal maps in the top row being embeddings.
Thanks to \ref{it: g coxeter 1} of Lemma~\ref{lem: coxeter subquotient} we know that $[M(1):L(w)]=1=[M^{I}(1):L(w)]$, which forces the composition $M(w)\rightarrow M^{I}(1)$ of (\ref{equ: Lie Verma wt shift}) to be non-zero. In particular, we obtain a non-zero map
\begin{equation}\label{equ: Lie parabolic Verma wt shift}
M^{I}(w) \rightarrow M^{I}(x).
\end{equation}
For each $L(x')\in\mathrm{JH}_{U(\fg)}(M^{I}(x))$, we note from \ref{it: n coh 1} of Lemma~\ref{lem: n coh collection} that
\begin{equation}\label{equ: Lie Verma wt shift Hom}
\Hom_{U(\fl_{I})}(L^{I}(w),H^0(\fn_{I},L(x')))\neq 0
\end{equation}
if and only if $x'=w$, in which case (\ref{equ: Lie Verma wt shift Hom}) is $1$-dimensional. A d\'evissage with respect to $L(x')\in\mathrm{JH}_{U(\fg)}(M^{I}(x))$ together with $[M^{I}(x):L(w)]=1$ gives
\[\Dim_E \Hom_{U(\fl_{I})}(L^{I}(w),H^0(\fn_{I},M^{I}(x)))\leq 1,\]
which together with (\ref{equ: g spectral seq 0}) gives
\[\Dim_E \Hom_{U(\fg)}(M^{I}(w),M^{I}(x))\leq 1.\]
Hence, the map (\ref{equ: Lie parabolic Verma wt shift}) constructed from (\ref{equ: Lie Verma wt shift}) is actually the unique  (up to scalars) such map between $U(\fg)$-modules.
\begin{lem}\label{lem: Lie St wt shift}
Let $x,w\in \Gamma$ with $x\unlhd w$ and $I_0\subseteq J_{w}$.
Then the map $M^{I_0}(w)\rightarrow M^{I_0}(x)$ restricts to a non-zero map
\begin{equation}\label{equ: Lie St wt shift}
\mathfrak{v}_{w,I_0}\rightarrow \mathfrak{v}_{x,I_0}
\end{equation}
which further induces an embedding
\begin{equation}\label{equ: Lie St wt shift cosocle}
\mathrm{cosoc}_{U(\fg)}(\mathfrak{v}_{w,I_0})\rightarrow \mathrm{cosoc}_{U(\fg)}(\mathfrak{v}_{x,I_0}).
\end{equation}
Under (\ref{equ: Lie Tits cosocle}), the embedding (\ref{equ: Lie St wt shift cosocle}) corresponds to the following evident embedding
\begin{equation}\label{equ: Lie St wt shift embedding}
\bigoplus_{u\in\Gamma_{\Delta\setminus I_0}, w\unlhd u}L(u)\hookrightarrow \bigoplus_{u\in\Gamma_{\Delta\setminus I_0}, x\unlhd u}L(u).
\end{equation}
\end{lem}
\begin{proof}
Upon fixing a non-zero map $M(w)\rightarrow M(x)$, there exists a unique choice of non-zero map $M^{I}(w)\rightarrow M^{I}(x)$ that fits into (\ref{equ: Lie Verma wt shift}) for each $I_0\subseteq I\subseteq J_{w}$, from which we obtain a map between complex
\begin{equation}\label{equ: Lie Tits wt shift 1}
\mathfrak{c}^{w,\Delta}_{I_0,J_{w}}\rightarrow \mathfrak{c}^{x,\Delta}_{I_0,J_{w}}.
\end{equation}
For each $I\subseteq J_{x}$ which satisfies $I\not\subseteq J_{w}$, the composition of
\[M^{I\cap J_{w}}(w)\rightarrow M^{I\cap J_{w}}(x)\twoheadrightarrow M^{I}(x)\]
is zero as $M^{I\cap J_{w}}(w)$ has cosocle $L(w)$ and $[M^{I}(x):L(w)]=0$ by Lemma~\ref{lem: coxeter parabolic Verma}. Consequently, the map (\ref{equ: Lie Tits wt shift 1}) factors through a map
\[\mathfrak{c}^{w,\Delta}_{I_0,J_{w}}\rightarrow \mathfrak{c}^{x,\Delta}_{I_0,J_{x}}\]
and thus induces a map of the form (\ref{equ: Lie St wt shift}).
Recall from Lemma~\ref{lem: Lie Tits cosocle} that we have
\begin{equation}\label{equ: Lie Tits cosocle w}
\mathrm{cosoc}_{U(\fg)}(\mathfrak{v}_{w,I_0})=\bigoplus_{u\in\Gamma_{\Delta\setminus I_0}, w\unlhd u}L(u)
\end{equation}
and
\begin{equation}\label{equ: Lie Tits cosocle x}
\mathrm{cosoc}_{U(\fg)}(\mathfrak{v}_{x,I_0})=\bigoplus_{u\in\Gamma_{\Delta\setminus I_0}, x\unlhd u}L(u).
\end{equation}
Now we consider an arbitrary $u\in\Gamma_{\Delta\setminus I_0}$ that satisfies $w\unlhd u$.
We have the following commutative diagram
\begin{equation}\label{equ: Lie St wt shift diagram}
\xymatrix{
M(u) \ar[r] \ar[d] & M(w) \ar[r] \ar[d] & M(x) \ar[d]\\
M^{I_0}(u) \ar[r] & M^{I_0}(w) \ar[r] & M^{I_0}(x)\\
 & \mathfrak{v}_{w,I_0} \ar[r] \ar[u] & \mathfrak{v}_{x,I_0} \ar[u]
}
\end{equation}
with the vertical maps from the third row to the second row being embeddings.
Now that $M(u)$ and $M^{I_0}(u)$ have cosocle $L(u)$, the horizontal maps in the top row of (\ref{equ: Lie St wt shift diagram}) are embeddings, and each term in (\ref{equ: Lie St wt shift diagram}) admits $L(u)$ as a constituent with multiplicity one (using the fact that each term in (\ref{equ: Lie St wt shift diagram}) is a subquotient of $M(1)$ and $[M(1):L(u)]=1$ from \ref{it: g coxeter 1} of Lemma~\ref{lem: coxeter subquotient}), we see that (\ref{equ: Lie St wt shift diagram}) induces maps
\begin{equation}\label{equ: u Verma to St}
M^{I_0}(u)\rightarrow \mathfrak{v}_{w,I_0}\rightarrow \mathfrak{v}_{x,I_0}
\end{equation}
which maps $\mathrm{cosoc}_{U(\fg)}(M^{I_0}(u))=L(u)$ isomorphically to the direct summand $L(u)$ of $\mathrm{cosoc}_{U(\fg)}(\mathfrak{v}_{w,I_0})$ and $\mathrm{cosoc}_{U(\fg)}(\mathfrak{v}_{x,I_0})$ under (\ref{equ: Lie Tits cosocle w}) and (\ref{equ: Lie Tits cosocle x}) respectively.
In particular, the map (\ref{equ: Lie St wt shift cosocle}) corresponds to the evident embedding (\ref{equ: Lie St wt shift embedding}) under (\ref{equ: Lie Tits cosocle w}) and (\ref{equ: Lie Tits cosocle x}).
The proof is thus finished.
\end{proof}

Let $I_0\subseteq \Delta$, $x\in\Gamma^{\Delta\setminus I_0}$ and $u\in\Gamma_{\Delta\setminus I_0}$ with $x\unlhd u$.
Thanks to (\ref{equ: Lie Tits cosocle x}) and (\ref{equ: u Verma to St}), we have non-zero maps 
\begin{equation}\label{equ: u Verma to St to simple}
M^{I_0}(u)\rightarrow \mathfrak{v}_{x,I_0}\twoheadrightarrow L(u)
\end{equation}
whose composition is the unique (up to scalars) surjection 
\begin{equation}\label{equ: u Verma to simple}
M^{I_0}(u)\twoheadrightarrow L(u).
\end{equation}
It follows from \ref{it: Lie Verma bottom Ext 2} of Lemma~\ref{lem: Lie Verma bottom Ext} that (\ref{equ: u Verma to simple}) induces the following isomorphism between $1$-dimensional $E$-vector spaces
\begin{equation}\label{equ: Lie Ext u Verma}
\mathrm{Ext}_{U(\fg)}^{\#\Delta\setminus I_0}(L(u),L(1))\buildrel\sim\over\longrightarrow \mathfrak{e}_{I_0,u}\defeq \mathrm{Ext}_{U(\fg)}^{\#\Delta\setminus I_0}(M^{I_0}(u),L(1)).
\end{equation}
\begin{prop}\label{prop: Lie St Ext w decomposition}
Let $I_0\subseteq \Delta$ and $x,w\in\Gamma^{\Delta\setminus I_0}$ with $x\unlhd w$. We have the following results.
\begin{enumerate}[label=(\roman*)]
\item \label{it: Lie St Ext w 2} The maps (\ref{equ: u Verma to St to simple}) for each $u\in\Gamma_{\Delta\setminus I_0}$ satisfying $x\unlhd u$ induce isomorphisms
\begin{equation}\label{equ: Lie St Ext w 2}
\bigoplus_{u\in\Gamma_{\Delta\setminus I_0}, x\unlhd u}\mathrm{Ext}_{U(\fg)}^{\#\Delta\setminus I_0}(L(u),L(1))\buildrel\sim\over\longrightarrow \mathrm{Ext}_{U(\fg)}^{\#\Delta\setminus I_0}(\mathfrak{v}_{x,I_0},L(1)) \buildrel\sim\over\longrightarrow\bigoplus_{u\in\Gamma_{\Delta\setminus I_0},x\unlhd u}\mathfrak{e}_{I_0,u}
\end{equation}
between $\#\Gamma_{J_{x}\setminus I_0}$-dimensional $E$-vector spaces, whose composition is given by the direct sum of (\ref{equ: Lie Ext u Verma}) over all $u\in\Gamma_{\Delta\setminus I_0}$ satisfying $x\unlhd u$.
\item \label{it: Lie St Ext w 3} The non-zero map (\ref{equ: Lie St wt shift}) induces the following commutative diagram
\begin{equation}\label{equ: Lie St Ext w 3}
\xymatrix{
\mathrm{Ext}_{U(\fg)}^{\#\Delta\setminus I_0}(\mathfrak{v}_{x,I_0},L(1)) \ar^{\sim}[rr] \ar[d] & & \bigoplus_{u\in\Gamma_{\Delta\setminus I_0},x\unlhd u}\mathfrak{e}_{I_0,u} \ar[d]\\
\mathrm{Ext}_{U(\fg)}^{\#\Delta\setminus I_0}(\mathfrak{v}_{w,I_0},L(1)) \ar^{\sim}[rr] & & \bigoplus_{u\in\Gamma_{\Delta\setminus I_0},w\unlhd u}\mathfrak{e}_{I_0,u}
}
\end{equation}
with the horizontal maps being isomorphisms from \ref{equ: Lie St Ext w 2}, and the right vertical map being the evident projection.
\end{enumerate}
\end{prop}
\begin{proof}
We write
\[\Gamma_{x}\defeq \{u\in \Gamma_{\Delta\setminus I_0}\mid x\unlhd u\}\]
for short.
Note that $y\mapsto yx$ gives a bijection
\[\Gamma_{J_{x}\setminus I_0}\buildrel\sim\over\longrightarrow \Gamma_{x}\]
and in particular $\#\Gamma_{x}=\#\Gamma_{J_{x}\setminus I_0}$.
By \ref{it: Lie St JH 1} of Lemma~\ref{lem: Lie St JH} and Lemma~\ref{lem: Lie Tits cosocle} we know that each $L(u')\in\mathrm{JH}_{U(\fg)}(\mathrm{rad}_{U(\fg)}(\mathfrak{v}_{x,I_0}))$ satisfies $\ell(u')>\#\Delta\setminus I_0$ and thus we know from Lemma~\ref{lem: Lie Verma bottom Ext} that $\mathrm{Ext}_{U(\fg)}^{k}(L(u'),L(1))$ is zero for each $k\leq \#\Delta\setminus I_0$. A d\'evissage with respect to $\mathrm{JH}_{U(\fg)}(\mathrm{rad}_{U(\fg)}(\mathfrak{v}_{x,I_0}))$ gives
\[\mathrm{Ext}_{U(\fg)}^{k}(\mathrm{rad}_{U(\fg)}(\mathfrak{v}_{x,I_0}),L(1))=0\]
for each $k\leq \#\Delta\setminus I_0$, which together with a further d\'evissage with respect to
\[0\rightarrow \mathrm{rad}_{U(\fg)}(\mathfrak{v}_{x,I_0})\rightarrow \mathfrak{v}_{x,I_0}\rightarrow \bigoplus_{u\in\Gamma_{x}}L(u)\rightarrow 0\]
gives the LHS isomorphism of (\ref{equ: Lie St Ext w 2}). This together with (\ref{equ: Lie Ext u Verma}) gives the RHS isomorphism of (\ref{equ: Lie St Ext w 2}) as well as the rest of \ref{it: Lie St Ext w 2}.
Finally, \ref{it: Lie St Ext w 3} follows directly from \ref{it: Lie St Ext w 2} (as well as its variant with $w$ replacing $x$) and Lemma~\ref{lem: Lie Tits cosocle}.
\end{proof}

Let $I_0\subseteq\Delta$ and $x\in\Gamma^{\Delta\setminus I_0}$.
Let $J\subseteq\Delta$ with $I_0\cup J=\Delta$ and thus $\mathrm{Supp}(x)\subseteq J$.
For each $I_0\subseteq I\subseteq J_{x}$, we have a natural surjection $M^{I\cap J}(x)\twoheadrightarrow M^{I}(x)$ (cf.~\cite[Thm.~9.4(c)]{Hum08}) which is functorial with respect to the choice of $I$. Up to a sign twist on these surjections, we obtain a map between complex $\mathfrak{c}^{x}_{I_0\cap J,J_{x}\cap J}[-\#J]\rightarrow \mathfrak{c}^{x}_{I_0,J_{x}}$ which together with (\ref{equ: Lie general Tits resolution}) induces a map
\begin{equation}\label{equ: Lie Tits J shift}
U(\fg)\otimes_{U(\fp_{J})}\mathfrak{v}_{x,I_0\cap J,J}\rightarrow \mathfrak{v}_{x,I_0}.
\end{equation}
Similar to the discussion following (\ref{equ: Lie St wt shift diagram}), we have a commutative diagram of non-zero maps
\begin{equation}\label{equ: Lie Tits J shift u diagram}
\xymatrix{
M^{I_0\cap J}(u) \ar[rr] \ar@{->>}[d] & &  U(\fg)\otimes_{U(\fp_{J})}\mathfrak{v}_{x,I_0\cap J,J} \ar[d]\\
M^{I_0}(u) \ar[rr] & & \mathfrak{v}_{x,I_0}
}
\end{equation}
for each $u\in\Gamma^{\Delta\setminus I_0}$ with $x\unlhd u$, with the top horizontal map of (\ref{equ: Lie Tits J shift u diagram}) equally obtained from the unique (up to scalars) non-zero map
\begin{equation}\label{equ: Lie Tits J u}
U(\fl_{J})\otimes_{U(\fl_{J}\cap\fp_{I_0\cap J})}L^{I_0\cap J}(u)\rightarrow \mathfrak{v}_{x,I_0\cap J,J}
\end{equation}
by applying $U(\fg)\otimes_{U(\fp_{J})}-$.
Taking direct sum of $u\in\Gamma^{\Delta\setminus I_0}$ satisfying $x\unlhd u$, we obtain the following commutative diagram
\begin{equation}\label{equ: Lie Tits J shift u diagram sum}
\xymatrix{
\bigoplus_{u\in\Gamma^{\Delta\setminus I_0}, x\unlhd u} M^{I_0\cap J}(u) \ar@{->>}[rr] \ar@{->>}[d] & &  U(\fg)\otimes_{U(\fp_{J})}\mathfrak{v}_{x,I_0\cap J,J} \ar@{->>}[d]\\
\bigoplus_{u\in\Gamma^{\Delta\setminus I_0}, x\unlhd u} M^{I_0}(u) \ar@{->>}[rr] & & \mathfrak{v}_{x,I_0}
}
\end{equation}
Here the bottom horizontal map of (\ref{equ: Lie Tits J shift u diagram sum}) is surjective by Lemma~\ref{lem: Lie Tits cosocle} and the discussion below (\ref{equ: Lie St wt shift diagram}). Similarly, upon replacing $\mathfrak{v}_{x,I_0}$ with $\mathfrak{v}_{x,I_0\cap J,J}$ in Lemma~\ref{lem: Lie Tits cosocle} and the discussion below (\ref{equ: Lie St wt shift diagram}), we see that the direct sum of (\ref{equ: Lie Tits J u}) over $u\in\Gamma^{\Delta\setminus I_0}$ with $x\unlhd u$ gives a surjection, which by applying $U(\fg)\otimes_{U(\fp_{J})}-$ implies that the top horizontal map of (\ref{equ: Lie Tits J shift u diagram sum}) is also a surjection.
Following the proof of \ref{it: Lie St Ext w 2} of Proposition~\ref{prop: Lie St Ext w decomposition}, we see that the commutative diagram (\ref{equ: Lie Tits J shift u diagram sum}) induces the following commutative diagram
\begin{equation}\label{equ: Lie Tits J shift u diagram Ext}
\xymatrix{
\mathrm{Ext}_{U(\fg)}^{\#\Delta\setminus I_0}(\mathfrak{v}_{x,I_0},1_{\fg}) \ar^{\sim}[rrr] \ar^{\wr}[d] & & & \bigoplus_{u\in\Gamma^{\Delta\setminus I_0}, x\unlhd u}\mathfrak{e}_{I_0,u} \ar^{\wr}[d]\\
\mathrm{Ext}_{U(\fg)}^{\#\Delta\setminus I_0}(U(\fg)\otimes_{U(\fp_{J})}\mathfrak{v}_{x,I_0\cap J,J},1_{\fg}) \ar^{\sim}[rrr] & & & \bigoplus_{u\in\Gamma^{\Delta\setminus I_0}, x\unlhd u}\mathfrak{e}_{I_0\cap J,u}
}
\end{equation}
with the right vertical map being the direct sum over $u$ of the isomorphism $\mathfrak{e}_{I_0,u}\buildrel\sim\over\longrightarrow \mathfrak{e}_{I_0\cap J,u}$ between $1$-dimensional $E$-vector spaces from \ref{it: Lie Verma bottom Ext 2} of Lemma~\ref{lem: Lie Verma bottom Ext}.

Recall from \cite[\S 4.2,4.6]{Hum08} that we have a unique (up to scalars) embedding $M(u)\hookrightarrow M(u')$ for each $u,u'\in W(G)$ with $u'\leq u$. We thus view $M(u)$ as a $U(\fg)$-submodule of $M(1)$ for each $u\in W(G)$, with $M(u)\subseteq M(u')$ for each $u,u'\in W(G)$ with $u'\leq u$.
\begin{lem}\label{lem: Verma intersection}
Let $x,y\in W(G)$. Assume that there exists $j_0\in\Delta$ such that $j<j_0$ for each $j\in\mathrm{Supp}(x)$ and that $j>j_0$ for each $j\in\mathrm{Supp}(y)$. Then we have $xy=yx$ and $M(xy)=M(x)\cap M(y)$.
\end{lem}
\begin{proof}
Under our assumption, we evidently have $xy=yx$. We write $I\defeq \mathrm{Supp}(x)$ and $I'\defeq \mathrm{Supp}(y)$ for short and note that $\fh_{I\sqcup I'}\cong \fh_{I}\times\fh_{I'}$. 
We write $M_{I}(u)\defeq U(\fh_{I})\otimes_{U(\fh_{I}\cap\fb)}u\cdot0\in\mathrm{Mod}_{U(\fh_{I})}$ for each $u\in W(L_{I})$. We define $M_{I'}(v)$ for each $v\in W(L_{I'})$ and $M_{I\sqcup I'}(w)$ for each $w\in W(L_{I\sqcup I'})$ similarly.
We have an inclusion of $U(\fh_{I})$-modules $M_{I}(x)\subseteq M_{I}(1)$ as well as an inclusion of $U(\fh_{I'})$-modules $M_{I'}(y)\subseteq M_{I'}(1)$, whose tensor over $E$ gives
\begin{multline}\label{equ: Verma tensor Levi factor}
M_{I\sqcup I'}(xy)=M_{I}(x)\otimes_E M_{I'}(y)=(M_{I}(x)\otimes_E M_{I'}(1))\cap (M_{I}(1)\otimes_E M_{I'}(y))\\
=M_{I\sqcup I'}(x)\cap M_{I\sqcup I'}(y)\subseteq M_{I}(1)\otimes_E M_{I'}(1)=M_{I\sqcup I'}(1).
\end{multline}
Applying $U(\fg)\otimes_{U(\fp_{I\sqcup I'})}-$ to (\ref{equ: Verma tensor Levi factor}) gives
\[M(xy)=M(x)\cap M(y)\subseteq M(1).\]
\end{proof}

The diagonal map $U(\fg)\rightarrow U(\fg)\otimes_EU(\fg)$ restricts to a map $U(\fp_{I\cap I'})\rightarrow U(\fp_{I})\otimes_EU(\fp_{I'})$ and thus induces a map
\begin{equation}\label{equ: Verma diagonal}
M^{I\cap I'}(1)\rightarrow M^{I}(1)\otimes_E M^{I'}(1).
\end{equation}
for each $I,I'\subseteq \Delta$.
\begin{lem}\label{lem: tensor of Verma}
Let $x,y,w\in\Gamma$ with $\mathrm{Supp}(x)\cap\mathrm{Supp}(y)=\emptyset$ and $w\in\Gamma_{x,y}$ (see (\ref{equ: x y envelop})). Then we have a unique (up to scalars) non-zero map that fits into the following commutative diagram
\begin{equation}\label{equ: tensor of Verma diagram}
\xymatrix{
M(w) \ar[r] \ar[d] & M(x)\otimes_EM(y) \ar[d]\\
M(1) \ar[r] & M(1)\otimes_EM(1)
}
\end{equation}
with the bottom horizontal map being (\ref{equ: Verma diagonal}) with $I=I'=\emptyset$.
Moreover, (\ref{equ: tensor of Verma diagram}) is functorial with respect to the choice of $x$, $y$ and $w$.
\end{lem}
\begin{proof}
Under the unique (up to scalars) embedding $M(x)\subseteq M(1)$, $M(y)\subseteq M(1)$ and $M(w)\subseteq M(1)$ as well as the diagonal embedding $M(1)\subseteq M(1)\otimes_EM(1)$, the commutative diagram (\ref{equ: tensor of Verma diagram}) is equivalent to the following inclusion
\begin{equation}\label{equ: tensor of Verma inclusion}
M(w)\subseteq (M(x)\otimes_EM(y))\cap M(1)\subseteq M(1)\otimes_E M(1).
\end{equation}
We divide the proof of (\ref{equ: tensor of Verma inclusion}) into the following steps.

\textbf{Step $1$}: We prove (\ref{equ: tensor of Verma inclusion}) when there exists $j_0\in\Delta$ such that $j<j_0$ for each $j\in\mathrm{Supp}(x)$ and that $j>j_0$ for each $j\in\mathrm{Supp}(y)$.\\
We have $w=xy=yx$ in this case.
We temporarily borrow notation from the proof of Lemma~\ref{lem: Verma intersection}. 
On one hand, the $U(\fp_{I})$-equivariant embedding $M_{I}(x)\subseteq M(x)$ and $U(\fp_{I'})$-equivariant embedding $M_{I'}(y)\subseteq M(y)$ together induce a $U(\fp_{I\sqcup I'})$-equivariant embedding
\[M_{I\sqcup I'}(w)=M_{I}(x)\otimes_E M_{I'}(y)\subseteq M(x)\otimes_EM(y).\]
On the other hand, the $U(\fp_{I})$-equivariant embedding $M_{I}(x)\subseteq M_{I}(1)$ and $U(\fp_{I'})$-equivariant embedding $M_{I'}(y)\subseteq M_{I'}(1)$ together induce a $U(\fp_{I\sqcup I'})$-equivariant embedding
\[M_{I\sqcup I'}(w)=M_{I}(x)\otimes_E M_{I'}(y)\subseteq M_{I}(1)\otimes_E M_{I'}(1)=M_{I\sqcup I'}(1)\subseteq M(1).\]
We thus obtain a $U(\fp_{I\sqcup I'})$-equivariant embedding
\[M_{I\sqcup I'}(w)=M_{I}(x)\otimes_E M_{I'}(y)\subseteq (M(x)\otimes_EM(y))\cap M(1).\]
This together with
\[\Hom_{U(\fp_{I\sqcup I'})}(M_{I\sqcup I'}(w),-)\cong \Hom_{U(\fg)}(M(w),-)\]
from (\ref{equ: g spectral seq 0}) (and that the unique (up to scalars) non-zero map $M(w)\rightarrow M(1)$ is an embedding) gives the inclusion (\ref{equ: tensor of Verma inclusion}).

\textbf{Step $2$}: We prove (\ref{equ: tensor of Verma inclusion}) when both $\mathrm{Supp}(x)$ and $\mathrm{Supp}(y)$ are intervals.\\
If either $x=1$ or $y=1$, we conclude (\ref{equ: tensor of Verma inclusion}) from \textbf{Step $1$}.
Upon exchanging $x,y$, it is harmless to assume that $\mathrm{Supp}(x)=[j_x,j_x']\neq\emptyset$ and $\mathrm{Supp}(y)=[j_y,j_y']\neq\emptyset$ with $j_x'<j_y$.
We write $x'\leq x$ (resp.~$y'\leq y$) for the unique element that satisfies $\mathrm{Supp}(x')=[j_x,j_x'-1]$ (resp.~$\mathrm{Supp}(x')=[j_y+1,j_y']$). Applying \textbf{Step $1$} to the pair $x,y'$ and the pair $x',y$, we deduce that
$M(xy')\subseteq (M(x)\otimes_EM(y'))\cap M(1)$ and $M(x'y)\subseteq (M(x')\otimes_EM(y))\cap M(1)$, which together with Lemma~\ref{lem: Verma intersection} gives
\begin{multline*}
M(w)\subseteq M(xy')\cap M(x'y)\subseteq (M(x)\otimes_EM(y'))\cap (M(x')\otimes_EM(y))\cap M(1)\\
=((M(x)\cap M(x'))\otimes_E(M(y)\cap M(y')))\cap M(1)=(M(x)\otimes_E M(y))\cap M(1).
\end{multline*}

\textbf{Step $3$}: We prove (\ref{equ: tensor of Verma inclusion}) for general $x$, $y$ and $w$.\\
We define $W_{x}$ as the set of connected components of $x$, namely the set of $x'\in \Gamma$ with $x'\leq x$ and $\mathrm{Supp}(x')$ being an non-empty subinterval of $\mathrm{Supp}(x)$. We define $W_{y}$ similarly. 
For each $x'\in W_{x}$ and $y'\in W_{y}$, there exists a unique $w_{x',y'}\leq w$ such that $\mathrm{Supp}(w_{x',y'})=\mathrm{Supp}(x')\sqcup\mathrm{Supp}(y')$ (cf.~\cite[Thm.~2.2.2]{BB05}). Note that our assumption on $x$, $y$ and $w$ ensures that $x',y'\leq w_{x',y'}$ for each $x'\in W_{x}$ and $y'\in W_{y}$. Now that both $\mathrm{Supp}(x')$ and $\mathrm{Supp}(y')$ are non-empty intervals, we have
\[M(w_{x',y'})\subseteq (M(x')\otimes_E M(y'))\cap M(1)\]
by \textbf{Step $2$}.
Taking intersection over all $x'\in W_{x}$ and $y'\in W_{y}$, we obtain the inclusions
\begin{multline*}
M(w)\subseteq \bigcap_{x'\in W_{x},y'\in W_{y}}M(w_{x',y'})\subseteq \bigcap_{x'\in W_{x},y'\in W_{y}}(M(x')\otimes_E M(y'))\cap M(1)\\
=\big((\bigcap_{x'\in W_{x}}M(x'))\otimes_E(\bigcap_{y'\in W_{y}}M(y'))\big)\cap M(1)=(M(x)\otimes_E M(y))\cap M(1)
\end{multline*}
with the last equality using Lemma~\ref{lem: Verma intersection}.
The proof is thus finished.
\end{proof}

Let $x,y,w\in\Gamma$ with $\mathrm{Supp}(x)\cap\mathrm{Supp}(y)=\emptyset$ and $w\in\Gamma_{x,y}$ (see (\ref{equ: x y envelop})). Let $I\subseteq \Delta\setminus\mathrm{Supp}(x)$ and $I'\subseteq \Delta\setminus\mathrm{Supp}(y)$.
Now that $M^{I}(x)$ is $\fp_{I}$-finite and $M^{I'}(y)$ is $\fp_{I'}$-finite, we know that $M^{I}(x)\otimes_EM^{I'}(y)$ is $\fp_{I\cap I'}$-finite and thus the composition of (with LHS map from Lemma~\ref{lem: tensor of Verma})
\[M(w)\hookrightarrow M(x)\otimes_E M(y)\twoheadrightarrow M^{I}(x)\otimes_EM^{I'}(y)\]
factors through $M^{I\cap I'}(w)$ by \cite[Thm.~9.4(c)]{Hum08}. We thus obtain a well-defined map
\begin{equation}\label{equ: parabolic Verma x y w}
M^{I\cap I'}(w)\rightarrow M^{I}(x)\otimes_EM^{I'}(y)
\end{equation}
which together with (\ref{equ: Verma diagonal}) forms the following commutative diagram
\begin{equation}\label{equ: tensor of parabolic Verma diagram}
\xymatrix{
M^{I\cap I'}(w) \ar[r] \ar[d] & M^{I}(x)\otimes_EM^{I'}(y) \ar[d]\\
M^{I\cap I'}(1) \ar[r] & M^{I}(1)\otimes_EM^{I'}(1)
}.
\end{equation}
Note that (\ref{equ: tensor of parabolic Verma diagram}) is clearly functorial with respect to the choice of $x$, $y$, $w$, $I$ and $I'$.

\begin{lem}\label{lem: tensor of Verma basic}
Take $y=s_j$ and $w=s_jx$ in (\ref{equ: parabolic Verma x y w}) for some $j\in J_{x}$ with $I\subseteq J_{x}$ and $I'\subseteq \Delta\setminus\{j\}$.
We have the following results.
\begin{enumerate}[label=(\roman*)]
\item \label{it: tensor of Verma basic 1} The resulting map
\begin{equation}\label{equ: tensor of Verma basic}
M^{I\cap I'}(s_jx)\rightarrow M^{I}(x)\otimes_EM^{I'}(s_j)
\end{equation}
is non-zero.
\item \label{it: tensor of Verma basic 2} The map (\ref{equ: tensor of Verma basic}) induces an isomorphism
\begin{equation}\label{equ: tensor of Verma Ext1}
\mathrm{Ext}_{U(\fg)}^1(M^{I'}(s_j),L(1))\buildrel\sim\over\longrightarrow\mathrm{Ext}_{U(\fg)}^1(M^{I\cap I'}(s_jx),M^{I}(x))
\end{equation}
between $1$ dimensional $E$-vector spaces.
\end{enumerate}
\end{lem}
\begin{proof}
We prove \ref{it: tensor of Verma basic 1}.\\
Now that the map (\ref{equ: tensor of Verma basic}) is functorial with respect to the choice of $I$ and $I'$ by its construction, we have the following commutative diagram
\begin{equation}\label{equ: tensor of Verma basic diagram}
\xymatrix{
M(s_jx) \ar[r] \ar[d] &  M(x)\otimes_EM(s_j) \ar[d]\\
M^{I\cap I'}(s_jx) \ar[r] & M^{I}(x)\otimes_EM^{I'}(s_j)
}
\end{equation}
with both vertical maps being surjections and the top horizontal map being an embedding.
We choose arbitrary generators $0\neq x_j\in\fu_{\al}$ and $0\neq y_j\in\fu_{\al}^+$ with $\al=(j,j+1)$.
We choose arbitrary highest weight vectors $0\neq v_{s_j}\in M(s_j)_{s_j\cdot 0}$ and $0\neq v_{x}\in M(x)_{x\cdot 0}$, and write $\overline{v}_{s_j}$ (resp.~$\overline{v}_{x}$) for its image in $M^{I'}(s_j)$ (resp.~$M^{I}(x)$) which is a highest weight vector inside.
Now that $(s_jx)\cdot0-x\cdot0-s_j\cdot0=N\al$ for some $N\geq 0$, there exists $c_k\in E$ for $0\leq k\leq N$ such that the image of $M(s_jx)_{s_jx\cdot0}$ under the top horizontal map of (\ref{equ: tensor of Verma basic diagram}) is spanned by
\begin{equation}\label{equ: parabolic Verma tensor vector}
v_{s_jx}\defeq \sum_{k=0}^{N}c_k (x_j^{N-k}\cdot v_{x})\otimes_E(x_j^{k}\cdot v_{s_j})\neq 0.
\end{equation}
We write $L$ for the kernel of the (fixed) surjection $M(x)\twoheadrightarrow M^{I}(x)$. It is clear that $L\subseteq N(x)$ and thus $L_{\mu}\neq 0$ for some $\mu\in\Lambda$ only if $\mu-y\cdot 0\in\Z_{\geq 0}\Phi^+$ for some $y>x$. Hence, either $y\neq s_jx$ and $\mu-x\cdot0-s_{j'}\cdot 0\in\Z_{\geq 0}\Phi^+$ for some $j'\in\mathrm{Supp}(y)\setminus\{j\}$, or $y=s_jx$ and $\mu-s_jx\cdot 0\in \Z_{\geq 0}\Phi^+$. In particular, we have $\mu-x\cdot0\neq k\al$ for each $0\leq k\leq N$, and thus $\{x_j^{N-k}\cdot\overline{v}_{x}\}_{0\leq k\leq N}$ is linearly independent in $M^{I}(x)=M(x)/L$. Now that $\{x_j^{k}\cdot\overline{v}_{s_j}\}_{0\leq k\leq N}$ is clearly linearly independent in $M^{I'}(s_j)$ (using $j\notin I'$), we know that
\[\{(x_j^{N-k}\cdot\overline{v}_{x})\otimes_E(x_j^{k}\cdot\overline{v}_{s_j})\}_{0\leq k\leq N}\]
is linearly independent in $M^{I}(x)\otimes_EM^{I'}(s_j)$ and in particular the image $\overline{v}_{s_jx}$ of $v_{s_jx}$ in $M^{I}(x)\otimes_EM^{I'}(s_j)$ is non-zero and necessarily a highest weight vector. In other words, we have shown that the composition of
\[M(s_jx)\twoheadrightarrow M^{I\cap I'}(s_jx)\rightarrow M^{I}(x)\otimes_EM^{I'}(s_j)\]
is non-zero, and in particular (\ref{equ: tensor of Verma basic}) is non-zero.

We prove \ref{it: tensor of Verma basic 2}.\\
The functor $M^{I}(x)\otimes_E-$ clearly induces a map
\[\mathrm{Ext}_{U(\fg)}^1(M^{I'}(s_j),L(1))\rightarrow \mathrm{Ext}_{U(\fg)}^1(M^{I}(x)\otimes_E M^{I'}(s_j),M^{I}(x))\]
which together with the map (\ref{equ: tensor of Verma basic}) induces the desired map
\begin{equation}\label{equ: tensor of Verma Ext1 prime}
\mathrm{Ext}_{U(\fg)}^1(M^{I'}(s_j),L(1))\rightarrow\mathrm{Ext}_{U(\fg)}^1(M^{I\cap I'}(s_jx),M^{I}(x)).
\end{equation}
It follows from Lemma~\ref{lem: g Ext1} that both the source and the target of (\ref{equ: tensor of Verma Ext1 prime}) are $1$ dimensional.
Hence, it suffices to show that the map (\ref{equ: tensor of Verma Ext1 prime}) is non-zero. Let $M$ be the unique (up to isomorphism) $U(\fg)$-module that fits into a non-split short exact sequence 
\begin{equation}\label{equ: tensor of Verma Ext1 seq}
0\rightarrow L(1)\rightarrow M\rightarrow M^{I'}(s_j)\rightarrow 0, 
\end{equation}
which clearly split over $U(\ft)$ with $M\twoheadrightarrow M^{I'}(s_j)$ inducing an isomorphism $L(1)_{0}\buildrel\sim\over\longrightarrow M_{0}$ and an isomorphism $M_{\mu}\buildrel\sim\over\longrightarrow M^{I'}(s_j)_{\mu}$ for each $\mu\neq 0$. We continue to use the notation in the proof of \ref{it: tensor of Verma basic 1}, and write $\overline{v}_{s_j}'\in M_{s_j\cdot 0}$ for the unique element that lifts $0\neq\overline{v}_{s_j}\in M^{I'}(s_j)_{s_j\cdot0}$. 
For each $\beta\in\Phi^+\setminus\{\al\}$, we have $M_{s_j\cdot 0-\beta}=M^{I'}(s_j)_{s_j\cdot0-\beta}=0$ and thus $u\cdot \overline{v}_{s_j}'=0$ for each $u\in\fu_{\beta}$. If $y_j\cdot\overline{v}_{s_j}'=0$, then we have $u\cdot \overline{v}_{s_j}'=0$ for each $u\in\fu$ which gives a non-zero map $M(s_j)\rightarrow M$ that necessarily factors through $M^{I'}(s_j)$ as $M$ is locally $\fp_{I'}$-finite. This forces (\ref{equ: tensor of Verma Ext1 seq}) to split and thus a contradiction.
Hence, we must have $y_j\cdot\overline{v}_{s_j}'\neq 0$ and spans $L(1)\subseteq M$.
We set
\[\overline{v}_{s_jx}'\defeq \sum_{k=0}^{N}c_k(x_j^{N-k}\cdot\overline{v}_{x})\otimes_E(x_j^{k}\cdot\overline{v}_{s_j}').\]
Now that $v_{s_jx}$ is the highest weight vector of $M(s_jx)\subseteq M(x)\otimes_EM(s_j)$, we know that $u\cdot \overline{v}_{s_jx}=0$ for each $u\in\fu$, and therefore 
\[y_j\cdot \overline{v}_{s_jx}'=c_0 (x_j^{N}\cdot\overline{v}_{x})\otimes_E(y_j\cdot\overline{v}_{s_j}')\in M^{I}(x)\otimes_EL(1)\] 
with $u\cdot\overline{v}_{s_jx}'=0$ for each $\beta\in\Phi^+\setminus\{\al\}$ and $u\in\fu_{\beta}$.
Consequently, we deduce that 
\[y_j^{N+1}\cdot\overline{v}_{s_jx}'=y_j^{N}\cdot(y_j\cdot \overline{v}_{s_jx}')=c_0 \overline{v}_{x}\otimes_E(y_j\cdot\overline{v}_{s_j}')\]
is a highest weight vector of $M^{I}(x)\otimes_EL(1)=M^{I}(x)$. This forces the following short exact sequence (induced from (\ref{equ: tensor of Verma Ext1 seq}) by applying $M^{I}(x)\otimes_E-$)
\[0\rightarrow M^{I}(x)\rightarrow M^{I}(x)\otimes_EM\rightarrow M^{I}(x)\otimes_EM^{I'}(s_j)\rightarrow 0\]
to remain non-split after pulling back along the map $M^{I\cap I'}(s_jx)\rightarrow M^{I}(x)\otimes_EM^{I'}(s_j)$.
We conclude that the map (\ref{equ: tensor of Verma Ext1 prime}) is non-zero and finish the proof.
\end{proof}

\begin{lem}\label{lem: tensor of parabolic Verma}
Let $x,y,w\in\Gamma$ with $w\in\Gamma_{x,y}$ (see (\ref{equ: x y envelop})). Let $I\subseteq J_{x}$ and $I'\subseteq J_{y}$. Then the map (\ref{equ: parabolic Verma x y w}) is non-zero.
\end{lem}
\begin{proof}
Now that (\ref{equ: parabolic Verma x y w}) is functorial with respect to the choice of $I$ and $I'$, it is harmless to prove this for $I=J_{x}$ and $I'=J_{y}$ (with $I\cap I'=J_{w}$). Now that the map
\begin{equation}\label{equ: parabolic Verma full tensor}
M^{J_{w}}(w)\rightarrow \bigotimes_{j\in\mathrm{Supp}(w)}M^{\Delta\setminus\{j\}}(s_j)
\end{equation} 
factors through (\ref{equ: parabolic Verma x y w}), it suffices to show that (\ref{equ: parabolic Verma full tensor}) is non-zero. 
Now that (\ref{equ: parabolic Verma full tensor}) factors through
\[M^{J_{w}}(w)\rightarrow M^{J_{s_{j_0}w}}(s_{j_0}w)\otimes_E M^{\Delta\setminus\{j_0\}}(s_{j_0}) \rightarrow\bigotimes_{j\in\mathrm{Supp}(w)}M^{\Delta\setminus\{j\}}(s_j)\]
for each $j_0\in D_L(w)$, we conclude the non-vanishing of (\ref{equ: parabolic Verma full tensor}) from \ref{it: tensor of Verma basic 1} of Lemma~\ref{lem: tensor of Verma basic} and an increasing induction on $w$ with respect to the partial-order $\unlhd$.
\end{proof}

\begin{lem}\label{lem: weak dominance}
Let $x,w\in\Gamma$ with $x<w$. Let $M$ be the unique subquotient of $M(1)$ with socle $L(w)$ and cosocle $L(x)$ (see \ref{it: g coxeter 2} of Lemma~\ref{lem: coxeter subquotient}).
Let $j\in\mathrm{Supp}(w)\setminus\mathrm{Supp}(x)$ such that there exists $j_1\in\mathrm{Supp}(x)$ that satisfies $|j-j_1|=1$ and $s_{j_1}s_j\leq w$.
Then we have $M\in\cO^{\fp_{\{j\}}}_{\rm{alg}}$.
\end{lem}
\begin{proof}
Let $w_1$ be the unique element satisfying $x<w_1<w$ and $\mathrm{Supp}(w_1)=\mathrm{Supp}(w)\setminus\{j\}$ (see \cite[Thm.~2.2.2]{BB05}). We write $J\defeq \{j\}$ for short and let $\xi_J: Z(\fl_J)\rightarrow E$ be the unique infinitesimal character such that $L^J(x)_{\xi_J}\neq 0$. It follows from \cite[Thm.~1.10]{Hum08} that $L^J(\mu)_{\xi_J}\neq 0$ for some $\mu\in \Lambda$ if and only if $\mu\in\{x\cdot 0, s_jx\cdot 0\}$.
It follows from \cite[Thm.~2.2.2]{BB05} and $w\in\Gamma$ that $s_{j_1}s_j$ is the unique element $w'\leq w$ that satisfies $\mathrm{Supp}(w')=\{j,j_1\}$. In particular, we have $s_js_{j_1}\not\leq w$, which together with $j_1\in\mathrm{Supp}(x)$ (and thus $s_{j_1}\leq x$) forces $s_jx\not\leq w$.
Note that $j\in\mathrm{Supp}(w)\setminus\mathrm{Supp}(x)$ and $x\in\Gamma$ implies that $x<s_jx\in\Gamma$.
Hence, it follows from Lemma~\ref{lem: coxeter subquotient} (by taking $I=\emptyset$ in \emph{loc.cit.}) that $[M:L(x)]=[M:L(w)]=1$, $M$ is a quotient of $M(x)$ and thus $L(s_jx)\notin\mathrm{JH}_{U(\fg)}(M)$ (using $s_jx\in\Gamma$ but $s_jx\not\leq w$). Hence, each $L(x')\in\mathrm{JH}_{U(\fg)}(\mathrm{rad}_{U(\fg)}(M))$ satisfies $x'>x$ and $x'\not>s_jx$, which together with \ref{it: n coh 3} of Lemma~\ref{lem: n coh collection} implies that $H^k(\fn_J,L(x'))_{\xi_J}=0$ for each $k\geq 0$. By \ref{it: n coh 1} and \ref{it: n coh 4} of Lemma~\ref{lem: n coh collection} we also have $H^0(\fn_J,L(x))_{\xi_J}\cong H^0(\fn_J,L(x))\cong L^J(x)$ and that $H^k(\fn_J,L(x))_{\xi_J}=0$ for $k\geq 1$. Hence, by a d\'evissage with respect to $\mathrm{JH}_{U(\fg)}(M)$, we deduce that $H^0(\fn_J,M)_{\xi_J}\cong H^0(\fn_J,L(x))_{\xi_J}\cong L^J(x)$. This together with (\ref{equ: g spectral seq}) gives
\[\Hom_{U(\fg)}(M^J(x),M)\cong \Hom_{U(\fl_J)}(L^J(x),H^0(\fn_J,M))\cong \Hom_{U(\fl_J)}(L^J(x),H^0(\fn_J,M)_{\xi_J})\neq 0.\]
In particular, we obtain a non-zero map $M^J(x)\rightarrow M$ which must be a surjection as $L(x)$ appears with multiplicity $1$ in both $M^J(x)$ and $M$, and only in the cosocle of both $M^J(x)$ and $M$. This together with $M^J(x)\in\cO^{\fp_J}_{\rm{alg}}$ forces $M\in\cO^{\fp_J}_{\rm{alg}}$.
\end{proof}

\begin{lem}\label{lem: Ext with dominant relation}
Let $x,w\in W(G)$ with $x=s_jw$ for some $j\in D_L(w)$. Let $M$ be the unique length $2$ $U(\fg)$-module with socle $L(x)$ and cosocle $L(w)$ (see Lemma~\ref{lem: g Ext1}). Then we have
\begin{equation}\label{equ: Ext relation vanishing}
\mathrm{Ext}_{U(\fg)}^k(M,L(1))=0
\end{equation}
for $k\leq \ell(x)$, and the cup product map
\begin{equation}\label{equ: Ext relation cup}
\mathrm{Ext}_{U(\fg)}^1(L(w),L(x))\otimes_E\mathrm{Ext}_{U(\fg)}^{\ell(x)}(L(x),L(1))\buildrel\cup\over\longrightarrow \mathrm{Ext}_{U(\fg)}^{\ell(w)}(L(w),L(1))
\end{equation}
is an isomorphism.
\end{lem}
\begin{proof}
The short exact sequence $0\rightarrow L(x)\rightarrow M\rightarrow L(w)\rightarrow 0$ induces the following exact sequence
\begin{equation}\label{equ: Ext with domiant seq}
\mathrm{Ext}_{U(\fg)}^k(L(w),L(1))\rightarrow \mathrm{Ext}_{U(\fg)}^k(M,L(1))\rightarrow \mathrm{Ext}_{U(\fg)}^k(L(x),L(1))\rightarrow \mathrm{Ext}_{U(\fg)}^{k+1}(L(w),L(1))
\end{equation}
for each $k\geq 0$.
Recall from Lemma~\ref{lem: Lie Verma bottom Ext} that $\mathrm{Ext}_{U(\fg)}^k(L(x),L(1))=0$ for $k<\ell(x)$, $\mathrm{Ext}_{U(\fg)}^k(L(w),L(1))=0$ for $k<\ell(w)=\ell(x)+1$, and that
\[\Dim_E \mathrm{Ext}_{U(\fg)}^{\ell(x)}(L(x),L(1))=1=\Dim_E \mathrm{Ext}_{U(\fg)}^{\ell(w)}(L(w),L(1)).\]
This together with (\ref{equ: Ext with domiant seq}) implies that (\ref{equ: Ext relation vanishing}) holds for each $k<\ell(x)$, and that (\ref{equ: Ext relation vanishing}) holds for $k=\ell(x)$ if and only if the cup product map (\ref{equ: Ext relation cup}) is an isomorphism.
We write $J\defeq \{j\}$ for short and let $M_J$ be the unique $U(\fl_J)$-module with socle $L^J(x)$ and cosocle $L^J(w)$ (see Lemma~\ref{lem: g Ext1} with $\fl_J$ here replacing $\fg$ in \emph{loc.cit.}).
Let $\xi: Z(\fl_J)\rightarrow E$ be the unique infinitesimal character such that $L^J(x)_{\xi}\neq 0$.
It follows from Lemma~\ref{lem: g induction quotient} (by replacing $I$ in \emph{loc.cit.} with $J$ here) we have a surjection $q: U(\fg)\otimes_{U(\fp_J)}M_J\twoheadrightarrow M$. As $[\mathrm{ker}(q)]=[U(\fg)\otimes_{U(\fp_J)}M_J]-[M]=[M^J(x)]+[M^J(w)]-[L(x)]-[L(w)]$ and $w>x$, we know that each $L(y)\in\mathrm{JH}_{U(\fg)}(\mathrm{ker}(q))$ satisfying $\ell(y)>\ell(x)$ and thus $\mathrm{Ext}_{U(\fg)}^k(L(y),L(1))=0$ for $k\leq \ell(x)$ by \ref{it: Lie Verma bottom Ext 1} of Lemma~\ref{lem: Lie Verma bottom Ext}. By a d\'evissage with respect to $\mathrm{JH}_{U(\fg)}(\mathrm{ker}(q))$, we have $\mathrm{Ext}_{U(\fg)}^k(\mathrm{ker}(q),L(1))=0$ for $k\leq \ell(x)$, and thus $q$ induces an isomorphism
\begin{equation}\label{equ: Ext Verma with dominant}
\mathrm{Ext}_{U(\fg)}^{\ell(x)}(M,L(1))\buildrel\sim\over\longrightarrow \mathrm{Ext}_{U(\fg)}^{\ell(x)}(U(\fg)\otimes_{U(\fp_J)}M_J,L(1)).
\end{equation}
By \ref{it: n coh 6} of Lemma~\ref{lem: n coh collection} we know that $H^{\ell}(\fn_J,L(1))_{\xi}\neq 0$ if and only if $\ell=\ell(x)$, in which case we have $H^{\ell(x)}(\fn_J,L(1))_{\xi}\cong L^J(x)$. This together with (\ref{equ: g spectral seq}) forces
\[\mathrm{Ext}_{U(\fg)}^{\ell(x)}(U(\fg)\otimes_{U(\fp_J)}M_J,L(1))\cong \Hom_{U(\fl_J)}(M_J,H^{\ell(x)}(\fn_J,L(1))_{\xi})=\Hom_{U(\fl_J)}(M_J,L^J(x))=0,\]
which together with (\ref{equ: Ext Verma with dominant}) gives (\ref{equ: Ext relation vanishing}) for $k=\ell(x)$.
\end{proof}

Let $x,y,w\in\Gamma$ with $w\in\Gamma_{x,y}$ (see (\ref{equ: x y envelop})).
Recall from (\ref{equ: Lie multi cup Ext}) that the map $M^{J_{w}}(w)\rightarrow M^{J_{x}}(x)\otimes_EM^{J_{y}}(y)$ (see Lemma~\ref{lem: tensor of parabolic Verma}) induce a map of the following form (see (\ref{equ: Lie Ext u Verma}) for notation)
\begin{equation}\label{equ: Lie bottom cup isom}
\mathfrak{e}_{J_{x},x}\otimes_E\mathfrak{e}_{J_{y},y}\buildrel\cup\over\longrightarrow 
\mathfrak{e}_{J_{w},w}.
\end{equation}
\begin{prop}\label{prop: Lie bottom cup isom}
Let $x,y,w\in\Gamma$ with $w\in\Gamma_{x,y}$. Then the map (\ref{equ: Lie bottom cup isom}) is an isomorphism between $1$-dimensional $E$-vector spaces.
\end{prop}
\begin{proof}
Note that both the source and the target of (\ref{equ: Lie bottom cup isom}) are known to be $1$-dimensional by Lemma~\ref{lem: Lie Verma bottom Ext}).
Recall from Proposition~\ref{prop: Lie multi cup} that we have the following map
\begin{equation}\label{equ: Lie multi cup nonvanishing}
\bigotimes_{j\in\mathrm{Supp}(w)}\mathfrak{e}_{\Delta\setminus\{j\},s_j}\buildrel\cup\over\longrightarrow\mathfrak{e}_{J_{w},w}
\end{equation}
between $1$-dimensional $E$-vector spaces (see Lemma~\ref{lem: Lie Verma bottom Ext}), which is well-defined (upon changing the order between tensor factors) up to a sign and moreover factors through (\ref{equ: Lie bottom cup isom}). Hence, to prove that (\ref{equ: Lie bottom cup isom}) is an isomorphism, it suffices to show that the map (\ref{equ: Lie multi cup nonvanishing}) is an isomorphism. Upon using an increasing induction on $w$ with respect to the partial-order $\unlhd$, we are reduced to show that the following map
\begin{equation}\label{equ: Lie cup nonvanishing induction}
\mathfrak{e}_{J_{s_{j}},s_{j}}\otimes_E\mathfrak{e}_{J_{u},u}\buildrel\cup\over\longrightarrow 
\mathfrak{e}_{J_{w},w}
\end{equation}
is an isomorphism for $u\defeq s_jw<w$ and each $j\in D_L(w)$.
Thanks to (\ref{equ: Lie cup Ext as composition}), we have the following commutative diagram
\begin{equation}\label{equ: Lie cup x y nonvanishing diagram}
\xymatrix{
\mathfrak{e}_{J_{s_{j}},s_{j}} \ar^{\wr}[d] & \otimes_E & \mathfrak{e}_{J_{u},u} \ar^{\cup}[rr] \ar@{=}[d] & & \mathfrak{e}_{J_{w},w} \ar@{=}[d]\\
\mathrm{Ext}_{U(\fg)}^{1}(M^{J_{w}}(w),M^{J_{u}}(u)) & \otimes_E & \mathfrak{e}_{J_{u},u} \ar^{\cup}[rr] & & \mathfrak{e}_{J_{w},w}
}
\end{equation}
with the left vertical map being an isomorphism between $1$-dimensional $E$-vector spaces by \ref{it: tensor of Verma basic 2} of Lemma~\ref{lem: tensor of Verma basic} and the bottom horizontal map being an isomorphism between $1$-dimensional $E$-vector spaces by Lemma~\ref{lem: g Ext1}, Lemma~\ref{lem: Lie Verma bottom Ext} and Lemma~\ref{lem: Ext with dominant relation}.
The proof is thus finished.
\end{proof}

\begin{lem}\label{lem: Ext special Verma}
Let $x,w\in \Gamma$ with $x<w$ and $\ell(w)=\ell(x)+1$ and $M$ be the unique length $2$ $U(\fg)$-module with socle $L(x)$ and cosocle $L(w)$ (see Lemma~\ref{lem: g Ext1}). Then for each $I\subseteq I_x\cap I_w$, there exists a finite dimensional $U(\fp_I)$-module $X$ such that we have a surjection $q: U(\fg)\otimes_{U(\fp_I)}X\rightarrow M$ and each $L(\mu)\in\mathrm{ker}(q)$ satisfies $\mathrm{Ext}_{U(\fg)}^k(L(\mu), L(1))=0$ for $k\leq \ell(x)$. In particular, we have
\begin{equation}\label{equ: Ext special Verma}
\mathrm{Ext}_{U(\fg)}^k(U(\fg)\otimes_{U(\fp_I)}X, L(1))=0
\end{equation}
for each $k\leq \ell(x)$.
\end{lem}
\begin{proof}
Let $\xi: Z(\fg)\rightarrow E$ be the unique infinitesimal character such that $L(x)_{\xi}\neq 0$.
As $I\subseteq I_x\cap I_w$, we know from Lemma~\ref{lem: dominance and left set} that both $L(x)$ and $L(w)$ are locally $\fp_I$-finite, and so is $M$. As $M$ is a finitely generated $U(\fg)$-module which is locally $\fp_{I}$-finite, there exists a finite dimensional $U(\fp_I)$-submodule $X\subseteq M$ such that the induced map $q: U(\fg)\otimes_{U(\fp_I)}X\rightarrow M$ is a surjection.
Note that $\Dim_E M_{x\cdot 0}=1$ and that $M_{\mu}\neq 0$ only if $\mu-x\cdot0\in \Z_{\geq 0}\Phi^+$, we deduce from $q$ being a surjection that $\Dim_E X_{x\cdot 0}=1=\Dim_E (U(\fg)\otimes_{U(\fp_I)}X)_{x\cdot 0}$ and that $(U(\fg)\otimes_{U(\fp_I)}X)_{\mu}\neq 0$ only if $\mu-x\cdot0\in \Z_{\geq 0}\Phi^+$.
In particular, we deduce that each $L(\mu)\in\mathrm{ker}(q)$ satisfies either $L(\mu)_{\xi}=0$ or $\mu=y\cdot 0$ for some $y>x$, and thus $\mathrm{Ext}_{U(\fg)}^k(L(\mu), L(1))=0$ for $k\leq \ell(x)$ in both cases. In particular, the surjection $q$ induces an isomorphism
\[\mathrm{Ext}_{U(\fg)}^k(M, L(1))\buildrel\sim\over\longrightarrow \mathrm{Ext}_{U(\fg)}^k(U(\fg)\otimes_{U(\fp_I)}X, L(1))\]
for $k\leq \ell(x)$.
It remains to show that
\begin{equation}\label{equ: Ext special Verma length two}
\mathrm{Ext}_{U(\fg)}^k(M, L(1))=0
\end{equation}
for each $k\leq \ell(x)$. 
Let $j\in\Delta$ be the unique element of $\mathrm{Supp}(w)\setminus\mathrm{Supp}(x)$. Thanks to (\ref{equ: Lie cup Ext as composition}), we have the following commutative diagram
\begin{equation}\label{equ: Ext special Verma diagram}
\xymatrix{
\mathfrak{e}_{J_{s_{j}},s_{j}} \ar[d] & \otimes_E & \mathfrak{e}_{J_{x},x} \ar^{\cup}[rr] \ar@{=}[d] & & \mathfrak{e}_{J_{w},w} \ar@{=}[d]\\
\mathrm{Ext}_{U(\fg)}^{1}(M^{J_{w}}(w),M^{J_{x}}(x)) & \otimes_E & \mathfrak{e}_{J_{x},x} \ar^{\cup}[rr] & & \mathfrak{e}_{J_{w},w}
}
\end{equation}
with the top horizontal map being an isomorphism by Proposition~\ref{prop: Lie bottom cup isom}. This forces the bottom horizontal map of (\ref{equ: Ext special Verma diagram}) to be non-zero. Now that the surjection $M^{J_{x}}(x)\twoheadrightarrow L(x)$ induces an isomorphism between $1$-dimensional $E$-vector spaces (see \ref{it: Lie Verma bottom Ext 2} of Lemma~\ref{lem: Lie Verma bottom Ext})
\[\mathrm{Ext}_{U(\fg)}^{\ell(x)}(L(x),L(1))\buildrel\sim\over\longrightarrow \mathfrak{e}_{J_{x},x},\]
we see that bottom horizontal map of (\ref{equ: Ext special Verma diagram}) factors through
\begin{equation}\label{equ: Ext special Verma cup 1}
\mathrm{Ext}_{U(\fg)}^{1}(M^{J_{w}}(w),L(x))\otimes_E\mathrm{Ext}_{U(\fg)}^{\ell(x)}(L(x),L(1))\buildrel\cup\over\longrightarrow \mathrm{Ext}_{U(\fg)}^{\ell(w)}(M^{J_{w}}(w),L(1))=\mathfrak{e}_{J_{w},w}
\end{equation}
which forces (\ref{equ: Ext special Verma cup 1}) to be non-zero.
It follows from (\ref{equ: g spectral seq}) and \cite[Lem.~3.2.7]{BQ24} that $\mathrm{Ext}_{U(\fg)}^{1}(M^{J_{w}}(w),L(x))$ is $1$-dimensional, which forces the surjection $M^{J_{w}}(w)\twoheadrightarrow L(w)$ to induce an isomorphism between $1$-dimensional $E$-vector spaces
\[\mathrm{Ext}_{U(\fg)}^{1}(L(w),L(x))\buildrel\sim\over\longrightarrow \mathrm{Ext}_{U(\fg)}^{1}(M^{J_{w}}(w),L(x)).\]
This together with (\ref{equ: Ext special Verma cup 1}) and \ref{it: Lie Verma bottom Ext 2} of Lemma~\ref{lem: Lie Verma bottom Ext} gives the following isomorphism between $1$-dimensional $E$-vector spaces
\begin{equation}\label{equ: Ext special Verma cup 2}
\mathrm{Ext}_{U(\fg)}^{1}(L(w),L(x))\otimes_E\mathrm{Ext}_{U(\fg)}^{\ell(x)}(L(x),L(1))\buildrel\cup\over\longrightarrow \mathrm{Ext}_{U(\fg)}^{\ell(w)}(L(w),L(1)).
\end{equation}
A further d\'evissage with respect to the short exact sequence $0\rightarrow L(x)\rightarrow M\rightarrow L(w)\rightarrow 0$ together with \ref{it: Lie Verma bottom Ext 1} of Lemma~\ref{lem: Lie Verma bottom Ext} gives the vanishing (\ref{equ: Ext special Verma length two}) for each $k\leq \ell(x)$. The proof is thus finished.
\end{proof}

\subsection{Smooth representations}\label{subsec: sm rep}
We recall from \cite{BQ24} some well-known results on smooth $E$-representations of $G$ with focus on $\mathrm{Ext}_{G}^{\bullet}(-,-)^{\infty}$ within the Bernstein block that contains $1_{G}$.

Let $G$ be a locally profinite group. We write $\mathrm{Rep}^{\infty}(G)$ for the abelian category of all smooth representations of $G$ over $E$ and $\mathrm{Rep}^{\infty}_{\rm{adm}}(G)$ for the full abelian subcategory of admissible ones. For each finite length object $V\in \mathrm{Rep}^{\infty}_{\rm{adm}}(G)$, we write $\mathrm{JH}_{G}(V)\defeq \mathrm{JH}_{\mathrm{Rep}^{\infty}_{\rm{adm}}(G)}(V)$ for its set of constituents.
We write $\mathrm{Ext}_{G}^\bullet(-,-)^{\infty}$ for the $\mathrm{Ext}$-groups in the abelian category $\mathrm{Rep}^{\infty}(G)$, and write $H^\bullet(G,-)^{\infty}\defeq \mathrm{Ext}_{G}^\bullet(1_{G},-)^{\infty}$ for short, with $1_{G}$ being the trivial representation of $G$.

Let $H\subseteq G$ be a closed subgroup and $\pi^{\infty}\in \mathrm{Rep}^{\infty}(H)$. We define $(\mathrm{Ind}_{H}^{G}\pi^{\infty})^{\infty}$ to be the $E$-vector space of locally constant functions $f: G\rightarrow \pi^{\infty}$ such that $f(xh)=h^{-1}\cdot f(x)$ for $x\in G$ and $h\in H$, which is naturally a (left) smooth $G$-representation via $(g(f))(x)\defeq f(g^{-1}x)$ ($g\in G$, $x\in G$, $f\in (\mathrm{Ind}_{H}^{G}\pi^{\infty})^{\infty}$).
We also consider $(\mathrm{ind}_{H}^{G}\pi^{\infty})^{\infty}$ as the subspace of $(\mathrm{Ind}_{H}^{G}\pi^{\infty})^{\infty}$ consisting of those $f$ for which there exists a compact open subset $C_f$ of $G$ such that $f(x)=0$ for $x\notin C_fH$.
They give the so-called (unnormalized) induction and compact induction functors
\[(\mathrm{Ind}_H^G\cdot)^{\infty}, (\mathrm{ind}_H^G\cdot)^{\infty}: \mathrm{Rep}^{\infty}(H) \rightarrow \mathrm{Rep}^{\infty}(G)\]
which are both exact. Note that they do not send $\mathrm{Rep}^{\infty}_{\rm{adm}}(H)$ to $\mathrm{Rep}^{\infty}_{\rm{adm}}(G)$ in general.

Let $G$ be the $K$-points of a $p$-adic reductive algebraic group over $K$. For (the $K$-points of) a parabolic subgroup $P=L_PN_P\subseteq G$, we have the \emph{unnormalized parabolic induction functor} $(\mathrm{Ind}_{P}^{G}\cdot)^{\infty}: \mathrm{Rep}^{\infty}(L_P)\rightarrow \mathrm{Rep}^{\infty}(G)$ which in that case restricts to a functor $(\mathrm{Ind}_{P}^{G}\cdot)^{\infty}: \mathrm{Rep}^{\infty}_{\rm{adm}}(L_P)\rightarrow \mathrm{Rep}^{\infty}_{\rm{adm}}(G)$ (see \cite[\S III.2.3]{Ren}).
The functor $\mathrm{Ind}_{P}^{G}$ admits a left adjoint functor $J_{N_P}: \mathrm{Rep}^{\infty}(G)\rightarrow \mathrm{Rep}^{\infty}(L_P)$ which is exact and restricts to a functor $J_{N_P}: \mathrm{Rep}^{\infty}_{\rm{adm}}(G)\rightarrow \mathrm{Rep}^{\infty}_{\rm{adm}}(L_P)$ (cf. \cite[\S\S VI.1.1, VI.6.1]{Ren} but note that our $J_{N_P}$ is the unnormalized Jacquet functor).

Now we take $G=\mathrm{PGL}_n(K)$ as in \S \ref{subsec: notation}.
For $I'\subseteq I\subseteq \Delta$, we use the shortened notation
\[i_{I',I}^{\infty}(\cdot)\defeq (\mathrm{Ind}_{P_{I'} \cap L_{I}}^{L_{I}}\cdot)^{\infty},\ J_{I,I'}(\cdot)\defeq J_{N_{I'}\cap L_{I}}(\cdot).\]
Then for each $I''\subseteq I'\subseteq I$, we clearly have
\[i_{I'',I}^{\infty}(\cdot)\cong i_{I',I}^{\infty}(i_{I'',I'}^{\infty}(\cdot)),\ J_{I,I''}(\cdot)\cong J_{I',I''}(J_{I,I'}(\cdot)).\]
We set $i_{I',I}^{\infty}\defeq i_{I',I}^{\infty}(1_{L_{I'}})$ and note that there exists an injection
\begin{equation}\label{equ: sm PS injection}
\theta_{I,I'}^{\infty}: 1_{L_{I}}\rightarrow i_{I',I}^{\infty}=i_{I',I}^{\infty}(1_{L_{I'}}).
\end{equation}
We choose $\theta_{I,I'}^{\infty}$ so that its sends constant functions on $L_{I}$ to themselves as elements in $i_{I',I}^{\infty}(1_{L_{I'}})$. For each $I''\subseteq I'\subseteq I\subseteq \Delta$, we have an injection
\[i_{I',I}^{\infty}(\theta_{I',I''}^{\infty}): i_{I',I}^{\infty}=i_{I',I}^{\infty}(1_{L_{I'}})\rightarrow i_{I',I}^{\infty}(i_{I'',I'}^{\infty}(1_{L_{I''}}))= i_{I'',I}^{\infty}(1_{L_{I''}})=i_{I'',I}^{\infty}\]
which satisfies
\begin{equation}\label{equ: sm consistent embedding}
\theta_{I,I''}^{\infty}=i_{I',I}^{\infty}(\theta_{I',I''}^{\infty})\circ\theta_{I,I'}^{\infty}.
\end{equation}
For each $I'\subseteq I\subseteq \Delta$, we set
\begin{equation}\label{equ: sm St}
V_{I',I}^{\infty}\defeq i_{I',I}^{\infty}/\sum_{I'\subsetneq I''\subseteq I}i_{I'',I}^{\infty}.
\end{equation}

We write $\widehat{T}^{\infty}$ for the set of smooth $E$-valued characters of $T$, which is naturally an abelian group under multiplication. We let $W(G)$ act on the left on $\widehat{T}^{\infty}$ via the \emph{smooth dot action} as in \cite[(35)]{BQ24}. For each finite length $\pi^{\infty}\in\mathrm{Rep}^{\infty}_{\rm{adm}}(L_I)$, it is well-known that $J_{I,\emptyset}(\pi^{\infty})\in\mathrm{Rep}^{\infty}_{\rm{adm}}(T)$ is of finite length and we write $\cJ(\pi^{\infty})\defeq\mathrm{JH}_{T}(J_{I,\emptyset}(\pi^{\infty}))$ for short.
For each $I\subseteq \Delta$, we write $\cB^I\defeq \cB^I_{W(L_I)\cdot 1_T}$ for the full subcategory of $\mathrm{Rep}^{\infty}_{\rm{adm}}(L_I)$ consisting of finite length objects $\pi^{\infty}$ such that each $\sigma^{\infty}\in\mathrm{JH}_{L_I}(\pi^{\infty})$ satisfies $\emptyset\neq \cJ(\sigma^{\infty})\subseteq W(L_I)\cdot 1_T$. Note from \cite[Lem.~2.1.11, Lem.~2.1.15]{BQ24} that $i_{I',I}^{\infty}, V_{I',I}^{\infty}\in \cB^I$ for each $I'\subseteq I\subseteq \Delta$.

Given $I_0,I_1\subseteq \Delta$, we set $d(I_0,I_1)\defeq \#I_0\setminus I_1+\#I_1\setminus I_0$ and note that $d(I_0,I_1)=d(I_1,I_0)$.
We set $[I_0,I_1]\defeq \{I_2\subseteq \Delta\mid d(I_0,I_1)=d(I_0,I_2)+d(I_2,I_1)\}$. We may equip $[I_0,I_1]$ with the following partial-order: $I_2,I_2'\in [I_0,I_1]$ satisfies $I_2\leq I_2'$ if and only if $d(I_0,I_1)=d(I_0,I_2)+d(I_2,I_2')+d(I_2',I_1)$, or equivalently $I_2\in[I_0,I_2']$ and $I_2'\in[I_2,I_1]$.
\begin{lem}\label{lem: sm cube}
We have the following results.
\begin{enumerate}[label=(\roman*)]
\item \label{it: sm cube 1} For each $I\subseteq \Delta$ and $I_0,I_1\subseteq I$, there exists a unique multiplicity free finite length $Q_{I}(I_0,I_1)\in\mathrm{Rep}^{\infty}_{\rm{adm}}(L_I)$ which has simple socle $V_{I_0,I}^{\infty}$ and simple cosocle $V_{I_1,I}^{\infty}$. Moreover, we have
    \begin{equation}\label{equ: sm cube 1}
    \mathrm{JH}_{L_I}(Q_{I}(I_0,I_1))=\{V_{I',I}^{\infty}\mid I'\in[I_0,I_1]\}
    \end{equation}
    with $V_{I_0',I}^{\infty}\leq V_{I_1',I}^{\infty}$ inside (\ref{equ: sm cube 1}) for some $I_0',I_1'\in[I_0,I_1]$ if and only if $I_0'\leq I_1'$ inside $[I_0,I_1]$.
\item \label{it: sm cube 2} For each $I\subseteq \Delta$ and $I_0,I_1,I_0',I_1'\subseteq I$, we have
\begin{equation}\label{equ: sm cube 2}
\Hom_{L_{I}}(Q_{I}(I_0,I_1),Q_{I}(I_0',I_1'))\neq 0
\end{equation}
if and only if $I_0'\in [I_0,I_1]$ and $I_1\in [I_0',I_1']$, in which case (\ref{equ: sm cube 2}) is one dimensional and any such non-zero map has image $Q_{I}(I_0',I_1)$.
\item \label{it: sm cube 3} For each $I\subseteq I'\subseteq \Delta$ and $I_0,I_1\subseteq I$, we have $i_{I,I'}^{\infty}(Q_{I}(I_0,I_1))\cong Q_{I'}(I_0\sqcup(I'\setminus I),I_1)$.
\end{enumerate}
\end{lem}
\begin{proof}
We deduce \ref{it: sm cube 1} and \ref{it: sm cube 2} from \cite[Lem.~2.2.1]{BQ24} and \ref{it: sm cube 3} from \cite[Lem.~2.2.6]{BQ24}.
\end{proof}

Let $I\subseteq \Delta$ and $I_0,I_1\subseteq I$. Given $\pi_i^{\infty}\in\mathrm{Rep}^{\infty}_{\rm{adm}}(L_{I_i})$ for $i=0,1$, we define
\begin{equation}\label{equ: sm distance}
d_{I}(\pi_0^{\infty},\pi_1^{\infty})\defeq \min\{k\mid \mathrm{Ext}_{L_I}^k(i_{I_0,I}^{\infty}(\pi_0^{\infty}), i_{I_1,I}^{\infty}(\pi_1^{\infty}))\neq 0\}.
\end{equation}
In particular, we have $d_{I}(\pi_0^{\infty},\pi_1^{\infty})=\infty$ if $\mathrm{Ext}_{L_I}^k(i_{I_0,I}^{\infty}(\pi_0^{\infty}), i_{I_1,I}^{\infty}(\pi_1^{\infty}))=0$ for $k\geq 0$.
For each $I_0,I_1,I_2\subseteq \Delta$, we set $d(I_2,[I_0,I_1])\defeq \min\{d(I_2,I)\mid I\in[I_0,I_1]\}$.
\begin{lem}\label{lem: general sm distance}
We have the following results.
\begin{enumerate}[label=(\roman*)]
\item \label{it: general sm distance 1} For each $I\subseteq \Delta$ and $I_0,I_1,I_0',I_1'\subseteq I$, we have $d_I(Q_{I}(I_0,I_1),Q_{I}(I_0',I_1'))=d_I(Q_{I}(I_1',I_0'),Q_{I}(I_1,I_0))$
\item \label{it: general sm distance 2} For each $I\subseteq \Delta$ and $I_0,I_1\subseteq I$, we have $d_I(V_{I_0,I}^{\infty},V_{I_1,I}^{\infty})=d(I_0,I_1)$ and
    \[\mathrm{Ext}_{L_I}^{d(I_0,I_1)}(V_{I_0,I}^{\infty},V_{I_1,I}^{\infty})^{\infty}=1.\]
    If moreover $I=\Delta$, then we have $\mathrm{Ext}_{G}^{k}(V_{I_0,\Delta}^{\infty},V_{I_1,\Delta}^{\infty})^{\infty}=0$ for $k\neq d(I_0,I_1)$.
\item \label{it: general sm distance 3} For each $I\subseteq \Delta$ and $I_0,I_1,I_2\subseteq I$, we have $d_I(V_{I_2,I}^{\infty},Q_{I}(I_1,I_0))=d_I(Q_{I}(I_0,I_1),V_{I_2,I}^{\infty})<\infty$ if and only if $d(I_2,[I_0,I_1])=d(I_2,I_1)$, in which case $d_I(V_{I_2,I}^{\infty},Q_{I}(I_1,I_0))=d_I(Q_{I}(I_0,I_1),V_{I_2,I}^{\infty})=d(I_2,I_1)$.
\item \label{it: general sm distance 4} For each $I\subseteq \Delta$ and $I_0,I_1,I_0',I_1'\subseteq I$, if $d_I(Q_{I}(I_1',I_0'),Q_{I}(I_1,I_0))<\infty$, then we have
    \[\Dim_E \mathrm{Ext}_{L_I}^{d_I(Q_{I}(I_1',I_0'),Q_{I}(I_1,I_0))}(Q_{I}(I_1',I_0'),Q_{I}(I_1,I_0))=1.\]
\end{enumerate}
\end{lem}
\begin{proof}
We deduce \ref{it: general sm distance 1} from the discussion before \cite[Lem.~2.2.3]{BQ24}, \ref{it: general sm distance 2} from \cite[Lem.~2.2.3]{BQ24} (which follows from \cite[Cor.~2]{Or05a}, see also \cite{Dat06}), \ref{it: general sm distance 3} from \cite[Lem.~2.2.9]{BQ24} (combined with \ref{it: general sm distance 1} above), and \ref{it: general sm distance 4} from \cite[Lem.~2.2.4]{BQ24}.
\end{proof}

\begin{lem}\label{lem: special sm distance}
Let $I_1\subseteq I_0\subseteq \Delta$ and $J\subseteq \Delta$.
Then we have $d_{\Delta}(Q_{\Delta}(I_0,I_1),Q_{\Delta}(\Delta,J))<\infty$ if and only if $\Delta\setminus I_0\subseteq J\subseteq (\Delta\setminus I_0)\sqcup I_1$, in which case we have $d_{\Delta}(Q_{\Delta}(I_0,I_1),Q_{\Delta}(\Delta,J))=\#\Delta\setminus I_0$.
\end{lem}
\begin{proof}
For each $I\subseteq\Delta$, we observe that $d(I,[J,\Delta])=d(I,\Delta)=\#\Delta\setminus I$ if and only if $\Delta\setminus J\subseteq I$ if and only if $\Delta\setminus I\subseteq J$. Hence, for each $I\in[I_0,I_1]$ that satisfies $\Delta\setminus J\not\subseteq I$, we have $d(I,[J,\Delta])\neq d(I,\Delta)$ and thus $d_{\Delta}(V_{I}^{\infty},Q_{\Delta}(\Delta,J))=d_{\Delta}(Q_{\Delta}(J,\Delta),V_{I}^{\infty})=\infty$ by \ref{it: general sm distance 3} of Lemma~\ref{lem: general sm distance}. If $\Delta\setminus J\not\subseteq I_0$, then each $I\in[I_0,I_1]$ satisfies $\Delta\setminus J\not\subseteq I$ and thus $d_{\Delta}(V_{I}^{\infty},Q_{\Delta}(\Delta,J))=d_{\Delta}(Q_{\Delta}(J,\Delta),V_{I}^{\infty})=\infty$, which together with a d\'evissage on $\mathrm{JH}_{G}(Q_{\Delta}(I_0,I_1))$ gives $d(Q_{\Delta}(I_0,I_1),Q_{\Delta}(\Delta,J))=\infty$.
Assume from now that $\Delta\setminus J\subseteq I_0$, then $Q_{\Delta}(I_0,I_1)$ contains $Q_{\Delta}(I_0,I_1\cup(\Delta\setminus J))$ as a subrepresentation (see \ref{it: sm cube 2} of Lemma~\ref{lem: sm cube}), and each $V_{I}^{\infty}\in\mathrm{JH}_{G}(Q_{\Delta}(I_0,I_1)/Q_{\Delta}(I_0,I_1\cup(\Delta\setminus J)))$ satisfies $\Delta\setminus J\not\subseteq I$ and thus $d_{\Delta}(V_{I}^{\infty},Q_{\Delta}(\Delta,J))=\infty$ by previous discussion. Hence, by d\'evissage on $\mathrm{JH}_{G}(Q_{\Delta}(I_0,I_1)/Q_{\Delta}(I_0,I_1\cup(\Delta\setminus J)))$, we deduce $d_{\Delta}(Q_{\Delta}(I_0,I_1)/Q_{\Delta}(I_0,I_1\cup(\Delta\setminus J)),Q_{\Delta}(\Delta,J))=\infty$, which together with further d\'evissage with respect to $0\rightarrow Q_{\Delta}(I_0,I_1\cup(\Delta\setminus J))\rightarrow Q_{\Delta}(I_0,I_1)\rightarrow Q_{\Delta}(I_0,I_1)/Q_{\Delta}(I_0,I_1\cup(\Delta\setminus J))\rightarrow 0$ gives
\begin{equation}\label{equ: sm cube equal distance}
d_{\Delta}(Q_{\Delta}(I_0,I_1),Q_{\Delta}(\Delta,J))=d_{\Delta}(Q_{\Delta}(I_0,I_1\cup(\Delta\setminus J)),Q_{\Delta}(\Delta,J))
\end{equation}
Under our assumption $\Delta\setminus J\subseteq I_0$ above, we observe that $I_1\cup(\Delta\setminus J)=I_0$ if and only if $I_0\setminus I_1\subseteq \Delta\setminus J$ if and only if $J\subseteq \Delta\setminus (I_0\setminus I_1)=(\Delta\setminus I_0)\sqcup I_1$. As $I_0\setminus (I_1\cup(\Delta\setminus J))\subseteq \Delta\setminus (\Delta\setminus J)=J$, we see that
\begin{equation}\label{equ: sm cube distance end}
d(J',[I_0,I_1\cup(\Delta\setminus J)])=d(J',I_0)
\end{equation}
for each $J'\supseteq J$ (or equivalently $J'\in[\Delta,J]$).
If $I_1\cup(\Delta\setminus J)=I_0$, then our assumption $\Delta\setminus J\subseteq I_0$ above ensures $d(I_0,[J,\Delta])=d(I_0,\Delta)$, which together with \ref{it: general sm distance 3} of Lemma~\ref{lem: general sm distance} gives $d_{\Delta}(V_{I_0}^{\infty},Q_{\Delta}(\Delta,J))=d(I_0,\Delta)=\#\Delta\setminus I_0$. If $I_1\cup(\Delta\setminus J)\neq I_0$, then by (\ref{equ: sm cube distance end}) and \ref{it: general sm distance 3} of Lemma~\ref{lem: general sm distance} we have $d_{\Delta}(Q_{\Delta}(I_0,I_1\cup(\Delta\setminus J)),V_{J'}^{\infty})=\infty$ for each $J'\in[\Delta,J]$, which together with a d\'evissage with respect to $\mathrm{JH}_{G}(Q_{\Delta}(\Delta,J))$ gives $d_{\Delta}(Q_{\Delta}(I_0,I_1\cup(\Delta\setminus J)),Q_{\Delta}(\Delta,J))=\infty$.
Up to this stage, we have shown that $d_{\Delta}(Q_{\Delta}(I_0,I_1\cup(\Delta\setminus J)),Q_{\Delta}(\Delta,J))<\infty$ if and only if $I_1\cup(\Delta\setminus J)=I_0$, in which case we have $d_{\Delta}(Q_{\Delta}(I_0,I_1\cup(\Delta\setminus J)),Q_{\Delta}(\Delta,J))=\#\Delta\setminus I_0$. This together with (\ref{equ: sm cube equal distance}) finishes the proof.
\end{proof}

\subsection{Locally analytic representations}\label{subsec: loc an rep}
We collect some basic notions in the theory of locally $K$(or $\Q_p$)-analytic representations of $G$ following \cite{ST01}, \cite{ST03} and \cite{ST05}.

For a paracompact locally $K$-analytic manifold $M$ of finite dimension (\cite[\S 8]{Sch2}), we write $C^{\rm{an}}(M)$ for the space of $E$-valued locally $K$-analytic functions on $M$ equipped with its usual locally convex topology (\cite[\S 12]{Sch2}). We write $D(M)\defeq C^{\rm{an}}(M)_b^\vee$ for its continuous $E$-dual equipped with the strong topology (\cite[\S 9]{Sch1}, \cite[\S 2]{ST01}).
Let $M_0$ be the locally $\Q_p$-analytic manifold underlying $M$. We write $C^{\rm{an}}(M_0)$ for the space of $E$-valued locally $\Q_p$-analytic functions on $M_0$ equipped with its usual locally convex topology, and $D(M_0)\defeq C^{\rm{an}}(M_0)_b^\vee$ for its continuous $E$-dual equipped with the strong topology. Then we have a closed embedding $C^{\rm{an}}(M)\hookrightarrow C^{\rm{an}}(M_0)$ which induces a strict surjection $D(M_0)\twoheadrightarrow D(M)$.

Similarly, we write $C^{\infty}(M)$ for the space of $E$-valued locally constant functions on $M$ equipped with its usual locally convex topology and $D^{\infty}(M)\defeq C^{\infty}(M)^\vee_b$ for its strong continuous $E$-dual. We have a closed embedding $C^{\infty}(M)\hookrightarrow C^{\rm{an}}(M)$ (use that $C^{\infty}(M)$ is the kernel of the continuous map $C^{\rm{an}}(M)\rightarrow C^{\rm{an}}(TM)$ where $TM$ is the tangent bundle, see \cite[\S 9]{Sch2}, \cite[Def.~12.4.i]{Sch2} and \cite[Rk.~6.2]{Sch2}) which induces a strict continuous surjection $D(M)\twoheadrightarrow D^{\infty}(M)$ (\cite[\S 2]{ST}).

If $M=G$ is a locally $K$-analytic group (automatically paracompact by \cite[Cor.~18.8]{Sch2}), then $D(G)$ is a unital associative $E$-algebra (with multiplication being the convolution, see \cite[Prop.~2.3]{ST01}).
We write $G_0$ for the locally $\Q_p$-analytic group underlying $G$, then the strict surjection $D(G_0)\twoheadrightarrow D(G)$ is a homomorphism between unital $E$-algebras. Let $\fg\defeq \mathrm{Lie}(G)\otimes_{K,\iota}E$ and $\fg_0\defeq \mathrm{Lie}(G_0)\otimes_{\Q_p}E$ be the $E$-Lie algebras associated with $G$ and $G_0$ respectively, then we have $\fg_0\cong \fg\times\fg'$ with $\fg'\defeq \prod_{\iota'\neq \iota}\mathrm{Lie}(G)\otimes_{K,\iota'}E$. By \cite[Lem.~5.1]{Schm08} we know that $\fg'D(G_0)$ is a closed two-sided ideal inside $D(G_0)$ with $D(G_0)\twoheadrightarrow D(G)$ inducing a topological isomorphism $D(G_0)/\fg'D(G_0)\buildrel\sim\over\longrightarrow D(G)$ between unital $E$-algebras.

For $A$ a Fr\'echet-Stein algebra we let $\cC_A$ be the abelian category of coadmissible left $A$-modules (see \cite[\S 3]{ST03}), which is a full subcategory of $\mathrm{Mod}_{A}$ (\cite[Cor.~3.5]{ST03}). By the discussion before and after \cite[Lemma~3.6]{ST03}, each $M\in \cC_A$ carries a canonical topology as a $E$-Fr\'echet space and any $A$-linear map in $\cC_A$ is continuous and strict.

If $G$ is a compact locally $K$-analytic group, by \cite[Thm.~5.1]{ST03} the algebra $D(G)$ is Fr\'echet-Stein, and so is its quotient $D^{\infty}(G)$ by \cite[Prop.~3.7]{ST03}. Moreover the continuous surjection $D(G)\twoheadrightarrow D^{\infty}(G)$ induces a fully faithful embedding $\cC_{D^{\infty}(G)}\hookrightarrow \cC_{D(G)}$. When $G$ is not necessarily compact, one defines the category $\cC_{D(G)}$ of coadmissible $D(G)$-modules over $E$ as the full subcategory of $\mathrm{Mod}_{D(G)}$ of $D(G)$-modules which are coadmissible as $D(H)$-module for any - equivalently one - compact open subgroup $H$ of $G$ (see \cite[\S 6]{ST03} and the references there). Replacing $\mathrm{Mod}_{D(G)}$ by $\mathrm{Mod}_{D(G_0)}$ and $\mathrm{Mod}_{D^{\infty}(G)}$, one defines in a similar way $\cC_{D(G_0)}$ and $\cC_{D^{\infty}(G)}$ which fits into embeddings $\cC_{D^{\infty}(G)}\rightarrow\cC_{D(G)}\rightarrow \cC_{D(G_0)}$.

Let $G$ be an arbitrary locally $K$-analytic group. We write $\mathrm{Rep}^{\rm{an}}(G)$ for the category of continuous $G$-representations on locally convex topological spaces of compact type. Given $V\in \mathrm{Rep}^{\rm{an}}(G)$, its strong dual $V^\vee\defeq V_b^\vee$ is a Fr\'echet nuclear space equipped with a separately continuous $D(G)$-action. We say that $V$ is \emph{admissible} if $V_b^\vee\in \cC_{G}$.
We write $\mathrm{Rep}^{\rm{an}}_{\rm{adm}}(G)$ for the category of admissible locally analytic representations of $G$ over $E$. The strong dual gives an anti-equivalence $\mathrm{Rep}^{\rm{an}}_{\rm{adm}}(G)\xrightarrow{\sim}\cC_{D(G)}$ (\cite[Thm.~6.3]{ST03}). Via this equivalence, the ``inverse image'' of the full subcategory $\cC_{D^{\infty}(G)}$ recovers the abelian category $\mathrm{Rep}^{\infty}_{\rm{adm}}(G)$ (see \S\ref{subsec: sm rep}) of admissible smooth representations of $G$ over $E$ (\cite[Thm.~6.6]{ST03}). Similarly, we can define $\mathrm{Rep}^{\rm{an}}_{\rm{adm}}(G_0)$ which is anti-equivalent to $\cC_{D(G_0)}$ and contains $\mathrm{Rep}^{\rm{an}}_{\rm{adm}}(G)$ as a full subcategory.

Let $R$ be a unital associated $E$-algebra with $S\subseteq R$ being a unital associative $E$-subalgebra.
We introduce relative Bar resolution, relative extension groups and their cup product associated with this pair $S\subseteq R$ (cf.~\cite[\S~8.7.1]{Wei94}). The explicit construction of cup product will be useful in later sections such as \S \ref{subsec: cup Tits double} and \S \ref{subsec: cup tits}.
Given $D\in\mathrm{Mod}_{R}$, we write
\begin{equation}\label{equ: std relative resolution general}
B_{\bullet}(R,S,D)\defeq \bigotimes_{S}^{\bullet+1}R\otimes_{S} D=R\otimes_{S}R\otimes_{S}\cdots\otimes_{S}R\otimes_{S}D
\end{equation}
for the standard $S$-split relative projective resolution (relative Bar resolution with respect to $S\subseteq R$) of $D$, with the cohomology of 
\[\Hom_{R}(B_{\bullet}(R,S,D),D')\]
recovering the relative extension groups $\mathrm{Ext}_{R,S}^{\bullet}(D,D')$.
When $S=0$, $\mathrm{Ext}_{R,S}^{\bullet}(D,D')=\mathrm{Ext}_{R}^{\bullet}(D,D')$ is the usual $\mathrm{Ext}$-groups computed in the abelian category $\mathrm{Mod}_{R}$.
Note that (\ref{equ: std relative resolution general}) is evidently functorial with respect to the choice of $D$.
Let $D,D',D''\in \mathrm{Mod}_{R}$. The construction $B_{\bullet}(R,S,-)$ induces a map
\[\Hom_{R}(B_{\bullet}(R,S,D),D')\rightarrow \Hom_{R}(\mathrm{Tot}(B_{\bullet}(R,S,B_{\bullet}(R,S,D))),B_{\bullet}(R,S,D')),\]
which together with the composition map
\begin{multline*}
\mathrm{Tot}(\Hom_{R}(\mathrm{Tot}(B_{\bullet}(R,S,B_{\bullet}(R,S,D))),B_{\bullet}(R,S,D'))\otimes_E\Hom_{R}(B_{\bullet}(R,S,D'),D''))\\
\rightarrow \Hom_{R}(\mathrm{Tot}(B_{\bullet}(R,S,B_{\bullet}(R,S,D))),D'')
\end{multline*}
defines a map
\begin{equation}\label{equ: abstract cup product composition}
\mathrm{Tot}(\Hom_{R}(B_{\bullet}(R,S,D),D')\otimes_E\Hom_{R}(B_{\bullet}(R,S,D),D'))\rightarrow \Hom_{R}(\mathrm{Tot}(B_{\bullet}(R,S,B_{\bullet}(R,S,D))),D'').
\end{equation}
The map $B_{\bullet}(R,S,D)\rightarrow D$ induces a map
\[\mathrm{Tot}(B_{\bullet}(R,S,B_{\bullet}(R,S,D)))\rightarrow B_{\bullet}(R,S,D)\]
between $S$-split relative projective (with respect to $S\subseteq R$) resolution of $D$ and thus a quasi-isomorphism
\[\Hom_{R}(B_{\bullet}(R,S,D)),D'')\rightarrow \Hom_{R}(\mathrm{Tot}(B_{\bullet}(R,S,B_{\bullet}(R,S,D))),D'')\]
between complex of $E$-vector spaces, which together with (\ref{equ: abstract cup product composition}) gives a map
\begin{equation}\label{equ: abstract cup product}
\mathrm{Tot}(\Hom_{R}(B_{\bullet}(R,S,D),D')\otimes_E\Hom_{R}(B_{\bullet}(R,S,D),D'))\dashrightarrow \Hom_{R}(B_{\bullet}(R,S,D),D'').
\end{equation}
in the derived category of $E$-vector spaces. By taking cohomology of (\ref{equ: abstract cup product}), one recovers
\begin{equation}\label{equ: abstract cup product Ext}
\mathrm{Tot}(\mathrm{Ext}_{R,S}^{\bullet}(D,D')\otimes_E\mathrm{Ext}_{R,S}^{\bullet}(D',D''))\rightarrow \mathrm{Ext}_{R,S}^{\bullet}(D,D'').
\end{equation}

Let $R\twoheadrightarrow \overline{R}$ be a surjection between associative $E$-algebras. Let $D\in\mathrm{Mod}_{R}$ and $\overline{D}\in\mathrm{Mod}_{\overline{R}}$. Let $\cI^{\bullet}$ be an injective resolution of $D$ in $\mathrm{Mod}_{R}$, then the cohomology of
\[\Hom_{R}(\overline{D},\cI^{\bullet})=\Hom_{\overline{R}}(\overline{D},\Hom_{R}(\overline{R},\cI^{\bullet}))\]
computes $\mathrm{Ext}_{R}^{\bullet}(\overline{D},D)$. The following double complex
\[\Hom_{\overline{R}}(B_{\bullet}(\overline{R},\overline{D}),\Hom_{R}(\overline{R},\cI^{\bullet}))\]
induces a spectral sequence of the form
\begin{equation}\label{equ: abstract ST}
\mathrm{Ext}_{\overline{R}}^{k}(\overline{D},\mathrm{Ext}_{R}^{\ell}(\overline{R},D))\implies \mathrm{Ext}_{R}^{k+\ell}(\overline{D},D).
\end{equation}
Let $\overline{D}'\in\mathrm{Mod}_{\overline{R}}$ be another object.
Similar to the definition of (\ref{equ: abstract cup product composition}), we can define a map
\begin{multline}\label{equ: abstract ST cup}
\mathrm{Tot}(\Hom_{\overline{R}}(B_{\bullet}(\overline{R},\overline{D}'),\overline{D})\otimes_E \Hom_{\overline{R}}(B_{\bullet}(\overline{R},\overline{D}),\Hom_{R}(\overline{R},\cI^{\bullet})))\\
\rightarrow \mathrm{Tot}(\Hom_{\overline{R}}(B_{\bullet}(\overline{R},B_{\bullet}(\overline{R},\overline{D})),\Hom_{R}(\overline{R},\cI^{\bullet}))),
\end{multline}
which leads to the following commutative diagram
\begin{equation}\label{equ: abstract ST cup diagram}
\xymatrix{
\mathrm{Ext}_{\overline{R}}^{k_0}(\overline{D}',\overline{D}) \ar@{=}[r] & \mathrm{Ext}_{\overline{R}}^{k_0}(\overline{D}',\overline{D})\\
\otimes_E & \otimes_E \\
\mathrm{Ext}_{\overline{R}}^{k_1}(\overline{D},\mathrm{Ext}_{R}^{\ell}(\overline{R},D)) \ar@{=>}[r] \ar^{\cup}[d] & \mathrm{Ext}_{R}^{k_1+\ell}(\overline{D},D) \ar^{\cup}[d]\\
\mathrm{Ext}_{\overline{R}}^{k_0+k_1}(\overline{D},\mathrm{Ext}_{R}^{\ell}(\overline{R},D)) \ar@{=>}[r] & \mathrm{Ext}_{R}^{k_0+k_1+\ell}(\overline{D},D)
}
\end{equation}

Let $G$ be a locally $K$-analytic group $G$ with $\fg=\mathrm{Lie}(G)\otimes_{K,\iota}E$ its associated $E$-Lie algebra and $D(G)=D(G,E)$ be its associated locally $K$-analytic distribution algebra with coefficients $E$.
For each $E$-Lie subalgebra $\fh\subseteq\fg$ and $D,D'\in\mathrm{Mod}_{D(G)}$, we can associate the relative $\mathrm{Ext}$-groups $\mathrm{Ext}_{D(G),U(\fh)}^\bullet(D,D')$. If $\fh=0$, this recovers the usual $\mathrm{Ext}$-groups $\mathrm{Ext}_{D(G)}^\bullet(D,D')$ defined for the abelian category $\mathrm{Mod}_{D(G)}$.

Let $G$ and $\fh\subseteq \fg$ be as above.
For $V_0,V_1$ in $\mathrm{Rep}^{\rm{an}}_{\rm{adm}}(G)$, we use the notation
\begin{equation}\label{extdef}
\mathrm{Ext}_{G,\fh}^\bullet(V_0,V_1)\defeq \mathrm{Ext}_{D(G),U(\fh)}^\bullet(V_1^\vee,V_0^\vee).
\end{equation}
Taking $V_0=1_{G}$ in (\ref{extdef}), we write
\[
H^\bullet(G,\fh,V_1)\defeq \mathrm{Ext}^\bullet_{G,\fh}(1_{G}, V_1)=\mathrm{Ext}^\bullet_{D(G),U(\fh)}(V_1^\vee, 1_{G}^\vee).
\]
We again omit $\fh$ from the above notation when $\fh=0$.

Compared with \S \ref{subsec: sm rep}, for $V_0,V_1$ in $\mathrm{Rep}^{\infty}_{\rm{adm}}(G)$ we have canonical isomorphisms (cf.~\cite[Lem.~4.2.3]{BQ24})
\begin{equation}\label{sm Ext def}
\mathrm{Ext}_G^\bullet(V_0,V_1)^{\infty}\cong\mathrm{Ext}_{D^{\infty}(G)}^\bullet(V_1^\vee,V_0^\vee).
\end{equation}
Taking $V_0=1_{G}$ in (\ref{sm Ext def}), we see that
\[H^\bullet(G,V_1)^{\infty}\cong \mathrm{Ext}_{D^{\infty}(G)}^\bullet(V_1^\vee, 1_{G}^\vee)\]

\begin{lem}\label{lem: ST spectral seq}
Let $\fh\subseteq \fg$ be an $E$-Lie subalgebra.
For each $D^{\infty}\in\mathrm{Mod}_{D^{\infty}(G)}$ and $D\in\mathrm{Mod}_{D(G)}$, we have the following spectral sequence
\begin{equation}\label{equ: general ST seq 1}
\mathrm{Ext}_{D^{\infty}(G)}^k(D^{\infty}, H^{\ell}(\fg, \fh, D))\implies \mathrm{Ext}_{D(G),U(\fh)}^{k+\ell}(D^{\infty},D).
\end{equation}
\end{lem}
\begin{proof}
We follow the argument in \cite[p.306]{ST05} with minor modification. On one hand, $\bigotimes^{\bullet+1}_{U(\fh)}U(\fg)$ is a relative projective resolution of $1_{\fg}$ w.r.t the pair $U(\fh)\subseteq U(\fg)$, making 
\[D(G)\otimes_{U(\fh)}(\bigotimes_{U(\fh)}^{\bullet}U(\fg))=D(G)\otimes_{U(\fg)}(\bigotimes^{\bullet+1}_{U(\fh)}U(\fg))\]
a relative projective resolution of $D^{\infty}(G)=D(G)\otimes_{U(\fg)}E$ with respect to the pair $U(\fh)\subseteq D(G)$. This together with
\[\Hom_{D(G)}(D(G)\otimes_{U(\fh)}(\bigotimes_{U(\fh)}^{\bullet}U(\fg)),D_2)=\Hom_{U(\fg)}(\bigotimes^{\bullet+1}_{U(\fh)}U(\fg),D_2)\]
gives an isomorphism
\[\mathrm{Ext}_{D(G),U(\fh)}^{\bullet}(D^{\infty}(G),D_2)\cong H^{\bullet}(\fg,\fh,D_2).\]
To conclude following the argument of \cite[p.306]{ST05}, it suffices to notice that the functor
\[\Hom_{D(G)}(D^{\infty}(G),-): \mathrm{Mod}_{D(G)}\rightarrow \mathrm{Mod}_{D^{\infty}(G)}\]
sends $D(G)$-modules that are relative injective with respect to $U(\fh)\subseteq D(G)$ to injective $D^{\infty}(G)$-modules.
\end{proof}

For each bounded cochain complex
\[\mathbf{C}=[\cdots\rightarrow C^{\ell}\rightarrow C^{\ell+1}\rightarrow \cdots]\in \mathrm{Ch}^{b}(\mathrm{Rep}^{\rm{an}}_{\rm{adm}}(G)),\]
we associate the following bounded cochain complex
\[\mathbf{C}^\vee\defeq [\cdots\rightarrow (C^{\ell+1})^\vee \rightarrow (C^{\ell})^\vee \rightarrow \cdots]\in \mathrm{Ch}^{b}(\cC_{G})\subseteq \mathrm{Ch}^{b}(\mathrm{Mod}_{D(G)})\]
with $(C^{\ell})^\vee$ sitting in degree $-\ell$.
This defines a contravariant embedding $\mathrm{Ch}^{b}(\mathrm{Rep}^{\rm{an}}_{\rm{adm}}(G))\rightarrow \mathrm{Ch}^{b}(\mathrm{Mod}_{D(G)}): \mathbf{C}\mapsto \mathbf{C}^\vee$ which is fully faithful.
We write $\cM(G)$ for the bounded derived category of $\mathrm{Mod}_{D(G)}$.
Given two cochain complex $\mathbf{C}_0, \mathbf{C}_1$ in $\mathrm{Ch}^{b}(\mathrm{Rep}^{\rm{an}}_{\rm{adm}}(G))$, we say they are \emph{isomorphic in the derived sense}, written $\mathbf{C}_0\cong \mathbf{C}_1$,
if $\mathbf{C}_1^\vee\cong \mathbf{C}_0^\vee$ in $\cM(G)$.
For $\mathbf{C}_0, \mathbf{C}_1$ in $\mathrm{Ch}^{b}(\mathrm{Rep}^{\rm{an}}_{\rm{adm}}(G))$, we set
\begin{equation}\label{equ: general Ext complex}
\mathrm{Ext}_{G}^{\bullet}(\mathbf{C}_0,\mathbf{C}_1)\defeq \mathrm{Ext}_{\cM(G)}^{\bullet}(\mathbf{C}_1^\vee,\mathbf{C}_0^\vee).
\end{equation}

We abuse each $V\in\mathrm{Rep}^{\rm{an}}_{\rm{adm}}(G)$ as a bounded cochain complex sitting at degree $0$.
For $\mathbf{C}_0, \mathbf{C}_1$ in $\mathrm{Ch}^{b}(\mathrm{Rep}^{\rm{an}}_{\rm{adm}}(G))$, general properties in $\cM(G)$ together with our definition (\ref{equ: general Ext complex}) gives the following two spectral sequences
\begin{equation}\label{equ: first complex spectral seq}
\mathrm{Ext}^k_G(\mathbf{C}_0,\mathbf{C}_1^{\ell})\implies \mathrm{Ext}_{G}^{k+\ell}(\mathbf{C}_0,\mathbf{C}_1)
\end{equation}
and
\begin{equation}\label{equ: second complex spectral seq}
\mathrm{Ext}^k_G(\mathbf{C}_0^{\ell},\mathbf{C}_1)\implies \mathrm{Ext}_{G}^{k-\ell}(\mathbf{C}_0,\mathbf{C}_1)
\end{equation}

Let $H\subseteq G$ be a pair of locally $K$-analytic groups with $G/H$ being compact. For each $V\in\mathrm{Rep}^{\rm{an}}(H)$, we define
\[(\mathrm{Ind}_{H}^{G}V)^{\rm{an}}\defeq \{f\in C^{\rm{an}}(G,V)\mid f(gh)=h^{-1}\cdot f(g) \ {\rm for each } \ g\in G \ {\rm and } \ h\in H\}\]
which is an object of $\mathrm{Rep}^{\rm{an}}(G)$ as $G/H$ is compact (cf.~\cite[Prop.~2.1.1]{Eme}). It follows from \cite[Prop.~2.1.2]{Eme} that $(\mathrm{Ind}_{H}^{G}-)^{\rm{an}}$ restricts to an exact functor
\[(\mathrm{Ind}_{H}^{G}-)^{\rm{an}}: \mathrm{Rep}^{\rm{an}}_{\rm{adm}}(H)\rightarrow \mathrm{Rep}^{\rm{an}}_{\rm{adm}}(G).\]

When $G=\mathrm{PGL}_n(K)$, we use the following shortened notation
\begin{equation}\label{equ: parabolic induction functor}
i_{I',I}^{\rm{an}}(-)\defeq (\mathrm{Ind}_{L_I\cap P_{I'}}^{L_I}-)^{\rm{an}}
\end{equation}
for each $I'\subseteq I\subseteq \Delta$.
We write $i_{I',I}^{\rm{an}}\defeq i_{I',I}^{\rm{an}}(1_{L_{I'}})$.
\subsection{Results on group cohomologies}\label{subsec: gp coh}
We prove all results we need on locally $K$-analyic group cohomologies, with the main result being Proposition~\ref{prop: Levi decomposition}.

Recall from \S \ref{subsec: notation} that we fix an embedding $\iota: K\rightarrow E$.
Let $G$ be a locally $K$-analytic group with $\fg=\mathrm{Lie}(G)\otimes_{K,\iota}E$ being its associated $E$-Lie algebra.
Let $\fh\subseteq \fg$ be a $E$-Lie subalgebra.
For each $V\in\mathrm{Rep}^{\rm{an}}_{\rm{adm}}(G)$, we define (see (\ref{equ: std relative resolution general}))
\begin{equation}\label{equ: abstract relative cochain}
C^{\bullet}(G,\fh,V)\defeq \Hom_{D(G)}(B_{\bullet}(D(G),U(\fh),V^{\vee}),1_{G}^{\vee})
\end{equation}
for short. When $V=1_{G}$, we omit $V$ from the notation and write $H^{\bullet}(G,\fh,1_{G})$ for the cohomology of (\ref{equ: abstract relative cochain}).
\begin{lem}\label{lem: relative ST seq}
Let $G$ and $\fh\subseteq\fg$ be as above. Then we have a spectral sequence of the form
\begin{equation}\label{equ: relative ST seq}
H^{k}(G,1_{G})^{\infty}\otimes_EH^{\ell}(\fg,\fh,1_{\fg})\implies H^{k+\ell}(G,\fh,1_{G}).
\end{equation}
\end{lem}
\begin{proof}
This is a special case of Lemma~\ref{lem: ST spectral seq} when $D=D^{\infty}=1_{G}^{\vee}$.
\end{proof}

\begin{defn}\label{def: degenerate group}
Let $G$ be a locally $K$-analytic group, $\fg$ be its associated $E$-Lie algebra and $\fh\subseteq \fg$ be a $E$-Lie subalgebra.
We say that the pair $(G,\fh)$ is \emph{degenerate} if the spectral sequence (\ref{equ: relative ST seq}) degenerates at the second page. When $\fh=0$, we simply say that $G$ is \emph{degenerate}.
\end{defn}
Note by definition above that the pair $(G,\fh)$ is degenerate if and only if we have
\begin{equation}\label{equ: degenerate dim}
\Dim_E H^k(G,\fh,1_{G})=\sum_{k_0+k_1=k}\Dim_E H^{k_0}(G,1_{G})^{\infty} \Dim_E H^{k_1}(\fg,\fh,1_{\fg})
\end{equation}
for each $k\geq 0$.

\begin{lem}\label{lem: group coh finite difference}
Let $\varphi: G'\dashrightarrow G$ be a quasi-isogeny between locally $K$-analytic group (see Definition~\ref{def: gp isogeny}) and $\fh\subseteq \fg$ be a $E$-Lie subalgebra. Then $\varphi$ induces isomorphisms $H^{\bullet}(G,1_{G})^{\infty}\buildrel\sim\over\longrightarrow H^{\bullet}(G',1_{G'})^{\infty}$, $H^{\bullet}(\fg,\fh,1_{\fg})\buildrel\sim\over\longrightarrow H^{\bullet}(\fg',\varphi^{-1}(\fh),1_{\fg'})$ as well as $H^{\bullet}(G,\fh,1_{G})\buildrel\sim\over\longrightarrow H^{\bullet}(G',\varphi^{-1}(\fh),1_{G'})$ between graded $E$-algebras. In particular, $(G,\fh)$ is degenerate if and only if $(G',\varphi^{-1}(\fh))$ is degenerate.
\end{lem}
\begin{proof}
According to Definition~\ref{def: gp isogeny}, it is harmless to assume that $\varphi$ is an isogeny from $G'$ to $G$.
Now that $\varphi$ induces an isomorphism $\fg'\buildrel\sim\over\longrightarrow \fg$ between $E$-Lie algebras, it also induces an isomorphism $\varphi^{-1}(\fh)\buildrel\sim\over\longrightarrow \fh$ and we do not distinguish between $\varphi^{-1}(\fh)$ and $\fh$ in the following. In particular, $\varphi$ induces an isomorphism
\begin{equation}\label{equ: group coh finite difference g}
H^{\bullet}(\fg,\fh,1_{\fg})\buildrel\sim\over\longrightarrow H^{\bullet}(\fg',\fh,1_{\fg'})
\end{equation}
between graded $E$-algebras.
For each finite group $M$, the abelian category $\mathrm{Rep}^{\infty}(M)$ is semi-simple and thus we have $H^k(M,1_{M})^{\infty}=0$ for $k>0$. By Hochschild--Serre spectral sequences for smooth $E$-representations of the pair $\mathrm{ker}(\varphi)\subseteq G'$ as well as the pair $\mathrm{im}(\varphi)\subseteq G$, we see that $\varphi$ induces an isomorphism
\begin{equation}\label{equ: group coh finite difference sm}
H^{\bullet}(G,1_{G})^{\infty}\buildrel\sim\over\longrightarrow H^{\bullet}(G',1_{G'})^{\infty}
\end{equation}
between graded $E$-algebras.
Note that $\varphi$ induces a map between spectral sequences from (\ref{equ: relative ST seq}) to its variant with $(G',\fh)$ replacing $(G,\fh)$. Thanks to (\ref{equ: group coh finite difference g}) and (\ref{equ: group coh finite difference sm}), we know that this map between spectral sequences is an isomorphism at the second page, and thus induces an isomorphism between limits, which gives the isomorphism
\[H^{\bullet}(G,\fh,1_{G})\buildrel\sim\over\longrightarrow H^{\bullet}(G',\fh,1_{G'})\]
between graded $E$-algebras.
Moreover, we notice that (\ref{equ: relative ST seq}) degenerates at the second page if and only if its variant for $(G',\fh)$ degenerates at the second page. In other words, $(G,\fh)$ is degenerate if and only if $(G',\fh)$ is degenerate.
\end{proof}

Given a discrete group $M$ and $D\in\mathrm{Mod}_{D(M)}$, we write
\begin{equation}\label{equ: discrete group cochain}
C^{\bullet}(M,D)\defeq \Hom_{D(M)}(B_{\bullet}(D(M)),D)
\end{equation}
for short. When the $D(M)$-action on $D$ is trivial, we have
\begin{equation}\label{equ: discrete group cochain trivial}
C^{\bullet}(M,D)=\Hom_{E}((B_{\bullet}(D(M)))_{M},D)
\end{equation}
where $(-)_{M}$ denotes the $M$-coinvariant.
\begin{lem}\label{lem: Kunneth K discrete}
Let $G$ be locally $K$-analytic group, $\fh\subseteq \fg$ an $E$-Lie subalgebra and $M$ be a discrete group.
Assume that $H^k(M,1_{M})$ is finite dimensional for each $k\geq 0$. We have the following results.
\begin{enumerate}[label=(\roman*)]
\item \label{it: Kunneth K 1} We have a canonical isomorphism
\begin{equation}\label{equ: Kunneth K discrete}
H^{\bullet}(G\times M, \fh, 1_{G\times M})\cong H^{\bullet}(G,\fh,1_{G})\otimes_E H^{\bullet}(M,1_{M})
\end{equation}
between $E$-vector spaces.
\item \label{it: Kunneth K 2} If $(G,\fh)$ is degenerate (see Definition~\ref{def: degenerate group}), then $(G\times M,\fh)$ is also degenerate.
\end{enumerate}
\end{lem}
\begin{proof}
Note that $G$ and $M$ are naturally subgroups of $G\times M$ with $G\subseteq G\times M$ being open and thus
\begin{equation}\label{equ: product discrete gp}
D(G\times M)\cong D(M)\otimes_E D(G)
\end{equation}
as unital associative $E$-algebras.
We recall $C^{\bullet}(G,\fh)$ from (\ref{equ: abstract relative cochain}) and define $C^{\bullet}(G\times M,\fh)$ similarly.
We have a natural identification as complex of $D(M)\otimes_ED(G)$-modules
\begin{equation}\label{equ: simplicial diagonal relative}
B_{\bullet}(D(M)\otimes_ED(G),U(\fh))=\mathrm{Diag}(B_{\bullet}(D(M))\otimes_EB_{\bullet}(D(G),U(\fh)))
\end{equation}
with $\mathrm{Diag}(-)$ standing for simplicial diagonal.
We write $(-)_{M}$ for the $M$-coinvariant of a representation of $M$.
The isomorphism (\ref{equ: simplicial diagonal relative}) between $D(M)\otimes_ED(G)$-modules together with the identification (\ref{equ: product discrete gp}) between unital associative $E$-algebras induces an isomorphism between $E$-vector spaces (see (\ref{equ: discrete group cochain trivial}))
\begin{multline}\label{equ: product double Hom}
C^{k}(G\times M,\fh)=\Hom_{D(G\times M)}(B_{k}(D(G\times M),U(\fh)),1_{D(G\times M)})\\
=\Hom_{D(M)\otimes_ED(G)}(B_{k}(D(M))\otimes_EB_{k}(D(G),U(\fh)),1_{D(M)\otimes_ED(G)})
\\\cong C^{k}(M, C^k(G,\fh))=\Hom_{E}((B_{k}(D(M)))_{M},C^k(G,\fh))
\end{multline}
for each $k\geq 0$, where $C^k(G,\fh)$ is equipped with the trivial $D(M)$-action.
We consider the following two cosimplicial bicomplex
\[B^{\bullet,\bullet}\defeq C^{\bullet}(M)\otimes_EC^{\bullet}(G,\fh)\]
and
\[\widehat{B}^{\bullet,\bullet}\defeq C^{\bullet}(M,C^{\bullet}(G,\fh)),\]
which come with a natural embedding $B^{\bullet,\bullet}\rightarrow \widehat{B}^{\bullet,\bullet}$ and therefore a map
\begin{equation}\label{equ: total complex completion}
\mathrm{Tot}(B^{\bullet,\bullet})\rightarrow \mathrm{Tot}(\widehat{B}^{\bullet,\bullet})
\end{equation}
between their associated total (cochain) complex.
Since both the functor $-\otimes_EC^{\bullet}(G,\fh)$ and the functor $\Hom_{E}(-,C^{\bullet}(G,\fh))$ are exact on $E$-vector spaces, the embedding $B^{\bullet,\bullet}\rightarrow \widehat{B}^{\bullet,\bullet}$ induces a map between spectral sequence which on the first page is given by
\[H^{k'}(M,1_{M})\otimes_EC^{k}(G,\fh) \rightarrow \Hom_{E}(H^{k'}((B_{\bullet}(D(M)))_{M}), C^{k}(G,\fh))\]
for each $k,k'\geq 0$, which is an isomorphism as 
\[H^{k'}(M,1_{M})=\Hom_{E}(H^{k'}((B_{\bullet}(D(M)))_{M}),E)\] 
is a finite dimensional $E$-vector space for each $k'\geq 0$ by our assumption.
In particular, the map (\ref{equ: total complex completion}) is a quasi-isomorphism.
We write $\mathrm{Diag}(\widehat{B}^{\bullet,\bullet})$ for the cosimplicial diagonal of $\widehat{B}^{\bullet,\bullet}$ with
\[\mathrm{Diag}(\widehat{B}^{\bullet,\bullet})^k\defeq C^{k}(M, C^{k}(G,\fh))\cong C^{k}(G\times M,\fh)\]
for $k\geq 0$ by (\ref{equ: product double Hom}). In other words, we can identify (the cochain complex associated with the cosimplicial complex $\mathrm{Diag}(\widehat{B}^{\bullet,\bullet})$) with $C^{\bullet}(G\times M,\fh)$, whose cohomology recovers $H^{\bullet}(G\times M,\fh,1_{G\times M})$.
By Eilenberg-Zilber's theorem, the Alexander-Whitney map (see \cite[\S 8.5]{Wei94})
\begin{equation}\label{equ: discrete AW}
\mathrm{Tot}(\widehat{B}^{\bullet,\bullet})\rightarrow \mathrm{Diag}(\widehat{B}^{\bullet,\bullet})
\end{equation}
is a quasi-isomorphism.
Now that the cohomology of $\mathrm{Tot}(B^{\bullet,\bullet})$ is clearly $\mathrm{Tot}(H^{\bullet}(G,\fh,1_{G})\otimes_E H^{\bullet}(M,1_{M}))$, we conclude (\ref{equ: Kunneth K discrete}) from the quasi-isomorphisms (\ref{equ: total complex completion}) and (\ref{equ: discrete AW}).
Upon replacing all $C^{\bullet}(G,\fh)$ in the argument above with $C^{\bullet}(G)^{\infty}$, a parallel argument gives
\begin{equation}\label{equ: sm Kunneth K discrete}
H^{\bullet}(G\times M, 1_{G\times M})^{\infty}\cong H^{\bullet}(G,1_{G})^{\infty}\otimes_E H^{\bullet}(M,1_{M})^{\infty}
\end{equation}
with $H^{\bullet}(M,1_{M})^{\infty}=H^{\bullet}(M,1_{M})$ as $M$ is discrete.

If $(G,\fh)$ is degenerate, then we have
\[\Dim_E H^{k}(G,\fh,1_{G})=\sum_{k_0+k_1=k}\Dim_E H^{k_0}(G,1_{G})^{\infty} \Dim_E H^{k_1}(\fg,\fh,1_{\fg})\]
for each $k\geq 0$, which together with (\ref{equ: Kunneth K discrete}) and (\ref{equ: sm Kunneth K discrete}) implies
\begin{multline*}
\Dim_E H^{m}(G\times M,\fh,1_{G\times H})=\sum_{k+\ell=m}\Dim_E H^{k}(G,\fh,1_{G}) \Dim_E H^{\ell}(M,1_{M})\\
=\sum_{k_0+k_1+\ell=m}\Dim_E H^{k_0}(G,1_{G})^{\infty} \Dim_E H^{k_1}(\fg,\fh,1_{\fg}) \Dim_E H^{\ell}(M,1_{M})\\
=\sum_{m_0+k_1=m}\big(\sum_{k_0+\ell=m_0}\Dim_E H^{k_0}(G,1_{G})^{\infty}\Dim_E H^{\ell}(M,1_{M})^{\infty}\big)\Dim_E H^{k_1}(\fg,\fh,1_{\fg})\\
=\sum_{m_0+k_1=m}\Dim_E H^{m_0}(G\times M,1_{G\times M})^{\infty}\Dim_E H^{k_1}(\fg,\fh,1_{\fg})
\end{multline*}
and thus $(G\times M,\fh)$ is degenerate. Here we use the fact that $H^{\ell}(M,1_{M})=H^{\ell}(M,1_{M})^{\infty}$ as graded $E$-algebras as $M$ is discrete.
\end{proof}

\begin{lem}\label{lem: abelian gp coh}
Let $M$ be a finitely generated discrete abelian group. Then we have the following canonical isomorphisms between graded $E$-algebras
\begin{equation}\label{equ: abelian gp coh}
H^{\bullet}(M,1_{M})\cong \wedge^{\bullet}H^1(M,1_{M})\cong \wedge^{\bullet}\Hom(M,E).
\end{equation}
\end{lem}
\begin{proof}
Thanks to Lemma~\ref{lem: group coh finite difference} it is harmless to assume that $M\cong \Z^{m}$ for some $m\geq 0$.
By inductive applying Lemma~\ref{lem: Kunneth K discrete} with $M=\Z$ in \emph{loc.cit.}, it suffices to show that the graded $E$-algebra $H^{\bullet}(\Z,1_{\Z})$ is supported in degree $[0,1]$ with $H^1(\Z,1_{\Z})\cong \Hom(\Z,E)$.
In fact, $1_{Z}^\vee$ admits the following resolution by free $D(\Z)$-modules
\[0\rightarrow D(\Z) \buildrel p_1\over \rightarrow D(\Z) \buildrel p_2\over \rightarrow 1_{\Z}^\vee\rightarrow 0\]
with $p_1(\delta)(x)=\delta(x)-\delta(x-1)$ and $p_2(f)(x)=\sum_{y\in\Z}\delta(y)$ for each $\delta\in D(\Z)$ and $x\in \Z$. We conclude by the observation that $\Hom_{D(\Z)}(D(\Z),1_{\Z}^\vee)$ is $1$-dimensional and $p_1$ induces the zero endomorphism of $\Hom_{D(\Z)}(D(\Z),1_{\Z}^\vee)$.
\end{proof}

Let $\mathbf{G}$ be a reductive group scheme over $K$ and assume from now that $G=\mathbf{G}(K)$.
Let $\fh\subseteq\fg$ be an $E$-Lie algebra.
We recall other related notation from \S \ref{subsec: notation}, including the closed subgroups $H_{I}$, $Z_{I}$, $L_{I}'$, $Z_{I}'$ and $Z_{I}''$ of $L_{I}$.
\begin{lem}\label{lem: center compact}
Let $I\subseteq \Delta$, we have $H^k(L_{I}', 1_{L_{I}'})^{\infty}=0$ for each $k>0$. In particular, $(L_{I}',\fl_{I}\cap\fh)$ is degenerate and we have a canonical isomorphism
\begin{equation}\label{equ: center compact}
H^{\bullet}(L_{I}',\fl_{I}\cap\fh,1_{L_{I}'})\buildrel\sim\over\longrightarrow H^{\bullet}(\fl_{I},\fl_{I}\cap\fh,1_{\fl_{I}})
\end{equation}
between graded $E$-algebras.
\end{lem}
\begin{proof}
As $H_{I}$ is semi-simple, it follows from \cite[Thm.~1]{Or05a} (see also \cite[\S 2.2.13]{Dat06}) that $H^k(H_{I}, 1_{H_{I}})^{\infty}=0$ for $k>0$. Let $P_{\bullet}$ be an arbitrary projective resolution of $1_{H_{I}}$ in $\mathrm{Rep}^{\infty}(H_{I})$, then the cohomology of $\Hom_{H_{I}}(P_{\bullet},1_{H_{I}})$ is concentrated in degree zero. Now that $Z_{I}'$ is compact, the complex $P_{\bullet}$ equipped with trivial $Z_{I}'$-action is also a projective resolution of $1_{H_{I}\times Z_{I}'}$ in $\mathrm{Rep}^{\infty}(H_{I}\times Z_{I}')$, with the cohomology of
\[\Hom_{H_{I}\times Z_{I}'}(P_{\bullet},1_{H_{I}\times Z_{I}'})\buildrel\sim\over\longrightarrow\Hom_{H_{I}}(P_{\bullet},1_{H_{I}})\]
concentrated in degree zero. In other words, we have shown that
\begin{equation}\label{equ: sm coh product}
H^k(H_{I}\times Z_{I}', 1_{H_{I}\times Z_{I}'})^{\infty}=0
\end{equation}
for $k>0$. Now that the natural map $H_{I}\times Z_{I}'\rightarrow L_{I}'$ is an isogeny, we finish the proof using (\ref{equ: sm coh product}) and Lemma~\ref{lem: group coh finite difference}.
\end{proof}

Let $M$ be a finitely generated discrete free abelian group.
We write $B_{\bullet}(M)_{M})\defeq B_{\bullet}(D(M))_{M})$ for short below.
We consider an isogeny $\Z^{m}\rightarrow M$ and $D$ be an arbitrary $E$-vector space. By an increasing induction on $m\geq 1$, we deduce from Eilenberg-Zilber's theorem (cf.~\cite[\S 8.5]{Wei94}) that the following shuffle product map (see \cite[Prop.~8.6.13]{Wei94})
\[\mathrm{Tot}(\bigotimes_{E}^{m}B_{\bullet}(D(\Z))_{\Z})\rightarrow B_{\bullet}(D(\Z^{m}))_{\Z^{m}}\]
is a quasi-isomorphism.
This induces the following quasi-isomorphisms
\begin{equation}\label{equ: M cochain split}
C^{\bullet}(M,D)=\Hom_{E}((B_{\bullet}(D(M)))_{M},D)\rightarrow \Hom_{E}((B_{\bullet}(D(\Z^{m})))_{\Z^{m}},D)
\rightarrow \Hom_{E}(\mathrm{Tot}(\bigotimes_{E}^{m}B_{\bullet}(\Z)_{\Z}),D)
\end{equation}
which are functorial with respect to the choice of the $E$-vector space $D$. 
Let $G$ be a locally $K$-analytic group and $\fh\subseteq\fg$ be an $E$-Lie subalgebra.
We have a shuffle product map
\[\mathrm{Tot}(B_{\bullet}(D(G),U(\fh))\otimes_EB_{\bullet}(D(M)))\rightarrow B_{\bullet}(D(G)\otimes_ED(M),U(\fh))=B_{\bullet}(D(G\times M),U(\fh))\]
which upon taking $\Hom_{D(G)\otimes_ED(M)}(-,1_{D(G)\otimes_ED(M)})$ and then combined with (\ref{equ: M cochain split}) (with $D=C^{k}(G,\fh)$ for $k\geq 0$ in \emph{loc.cit.}) gives the following quasi-isomorphisms
\begin{equation}\label{equ: G M natural}
C^{\bullet}(G\times M,\fh)\rightarrow \mathrm{Tot}(C^{\bullet}(M,C^{\bullet}(G,\fh)))\rightarrow C^{\bullet}(G\times M,\fh)_{\natural}\defeq \mathrm{Tot}(\Hom_{E}(\bigotimes_{E}^{m}B_{\bullet}(\Z)_{\Z},C^{\bullet}(G,\fh))).
\end{equation}
It is clear that the maps in (\ref{equ: G M natural}) are functorial with respect to the choice of $\fh\subseteq \fg$.
Recall from (\ref{equ: abstract cup product composition}) and (\ref{equ: abstract cup product}) (by taking $R=D(G\times M)$ and $S=U(\fh)$ in \emph{loc.cit.}) that we defined the following map in the derived category of $E$-vector spaces
\begin{equation}\label{equ: G M cup}
\mathrm{Tot}(C^{\bullet}(G\times M,\fh)\otimes_EC^{\bullet}(G\times M,\fh))\dashrightarrow C^{\bullet}(G\times M,\fh).
\end{equation}
One can check from the explicit formula for shuffle product map (see \cite[Prop.~8.6.13, Ex.~8.6.5]{Wei94}) that the map (\ref{equ: G M cup}) fits into a commutative diagram of the form
\begin{equation}\label{equ: G M cup diagram}
\xymatrix{
\mathrm{Tot}(C^{\bullet}(G\times M,\fh)\otimes_EC^{\bullet}(G\times M,\fh)) \ar@{-->}[r] \ar[d] & C^{\bullet}(G\times M,\fh) \ar[d]\\
\mathrm{Tot}(C^{\bullet}(G\times M,\fh)_{\natural}\otimes_EC^{\bullet}(G\times M,\fh)_{\natural}) \ar@{-->}[r] & C^{\bullet}(G\times M,\fh)_{\natural}
}
\end{equation}
with vertical maps being quasi-isomorphisms from (\ref{equ: G M natural}).
Taking cohomology of (\ref{equ: G M cup diagram}) and then combining with Lemma~\ref{lem: abelian gp coh}, we obtain isomorphisms between graded $E$-algebras
\begin{equation}\label{equ: Kunneth graded algebra}
H^{\bullet}(G\times M,\fh,1_{G\times M})\cong H^{\bullet}(G,\fh,1_{G})\otimes_EH^{\bullet}(M,1_{M})\cong H^{\bullet}(G,\fh,1_{G})\otimes_E\wedge^{\bullet}\Hom(M,E)
\end{equation}
refining Lemma~\ref{lem: Kunneth K discrete}.

Let $G$ and $G'$ be two locally $K$-analytic groups. Let $M$ and $M'$ be two finitely generated discrete abelian groups.
We consider a homomorphism $G'\times M'\rightarrow G\times M$ between locally $K$-analytic groups and an isogeny $\Z^{m'}\rightarrow M$ such that the composition of $G'\times \Z^{m'}\rightarrow G'\times M'\rightarrow G\times M$ factors through
\begin{equation}\label{equ: loc an separate discrete}
G'\times \Z^{m'}\buildrel\sim\over\longrightarrow G'\times \Z^{m'-m}\times \Z^{m}\rightarrow G\times M
\end{equation}
for some $m\leq m'$, with $\Z^{m'}\buildrel\sim\over\longrightarrow \Z^{m'-m}\times \Z^{m}$ given by partition of the $m'$-coordinates into two subsets (with cardinality $m'-m$ and $m$ respectively), and the second map of (\ref{equ: loc an separate discrete}) being the product of a map $G'\times \Z^{m'-m}$ and an isogeny $\Z^{m}\rightarrow M$. 
Let $\fh\subseteq\fg$ and $\fh'\subseteq\fg'$ be $E$-Lie algebras such that the map $\fg'\rightarrow \fg$ (associated with $G'\rightarrow G$) induces a map $\fh'\rightarrow \fh$.
By considering the variant of (\ref{equ: G M natural}) for each term of (\ref{equ: loc an separate discrete}), we obtain the following commutative diagram
\[
\xymatrix{
C^{\bullet}(G\times M,\fh) \ar[d] & \\
C^{\bullet}(G\times \Z^{m},\fh) \ar[r] \ar[d] & \mathrm{Tot}(\Hom_{E}(\bigotimes_{E}^{m}B_{\bullet}(\Z)_{\Z},C^{\bullet}(G,\fh))) \ar[d]\\
C^{\bullet}(G'\times \Z^{m'-m}\times \Z^{m},\fh') \ar[r] \ar@{=}[d] & \mathrm{Tot}(\Hom_{E}(\bigotimes_{E}^{m}B_{\bullet}(\Z)_{\Z},C^{\bullet}(G'\times \Z^{m'-m},\fh'))) \ar[d]\\
C^{\bullet}(G'\times \Z^{m'},\fh') \ar[r] & \mathrm{Tot}(\Hom_{E}(\bigotimes_{E}^{m'}B_{\bullet}(\Z)_{\Z},C^{\bullet}(G,\fh')))
}
\]
which factors through
\begin{equation}\label{equ: G M splitting diagram}
\xymatrix{
C^{\bullet}(G\times M,\fh) \ar[r] \ar[d] & C^{\bullet}(G\times M,\fh)_{\natural} \ar[d]\\
C^{\bullet}(G'\times M',\fh') \ar[r] & C^{\bullet}(G'\times M',\fh')_{\natural}
}
\end{equation}
with both horizontal maps being quasi-isomorphisms.

Now we consider $\Z^{\#\Delta}=\prod_{i\in\Delta}\Z_{i}$ with $\Z_{i}$ standing for the copy of $\Z$ indexed by some $i\in\Delta$.
We fix the choice of an isogeny $\varsigma_{\emptyset}: \prod_{i\in\Delta}\Z_{i}\rightarrow T''$ whose composition with $T''\subseteq T\subseteq L_{I}$ factors through
\begin{equation}\label{equ: Levi center dagger projection}
\prod_{i\in\Delta}\Z_{i}\twoheadrightarrow \prod_{i\in\Delta\setminus I}\Z_{i}\buildrel\varsigma_{I}\over\longrightarrow Z_{I}''\subseteq Z_{I}\subseteq L_{I}
\end{equation}
with the first map being the evident projection and the middle map $\varsigma_{I}$ being an isogeny. We write $Z_{I}^{\dagger}$ for the image of $\varsigma_{I}$ which is a subgroup of $Z_{I}''$ of finite index. 

Let $I\subseteq J\subseteq\Delta$. 
The isogeny $L_{I}'\times(Z_{I}^{\dagger}\cap H_{J})\times Z_{J}^{\dagger}\rightarrow L_{I}$ together with the isomorphism $(\fl_{I}\cap\fh_{J})\times\fz_{J}\cong\fl_{I}$ induce the following inclusion
\begin{equation}\label{equ: Levi Lie discrete inclusion}
U(\fl_{I}\cap\fh_{J})\otimes_{E}U(\fz_{J})\otimes_{E}D(Z_{I}^{\dagger}\cap H_{J})\otimes_{E}D(Z_{J}^{\dagger})\subseteq D(L_{I})
\end{equation}
and thus a shuffle product map (see \cite[Prop.~8.6.13]{Wei94})
\[\mathrm{Tot}(B_{\bullet}(U(\fl_{I}\cap\fh_{J}))\otimes_EB_{\bullet}(U(\fz_{J}))\otimes_EB_{\bullet}(D(Z_{I}^{\dagger}\cap H_{J}))\otimes_EB_{\bullet}(D(Z_{J}^{\dagger})))\rightarrow B_{\bullet}(D(L_{I})).\]
This together with Lemma~\ref{lem: center compact} (as well as Eilenberg-Zilber's theorem, see \cite[\S 8.5]{Wei94}) induces the following quasi-isomorphism
\begin{multline}\label{equ: Levi Lie discrete cochain}
C^{\bullet}(L_{I})\rightarrow \mathrm{Tot}(\Hom_{E}(B_{\bullet}(Z_{J}^{\dagger})_{Z_{J}^{\dagger}},\Hom_{E}(B_{\bullet}(Z_{I}^{\dagger}\cap L_{J}')_{Z_{J}^{\dagger}\cap H_{J}},C^{\bullet}(\fz_{I},C^{\bullet}(\fl_{I}\cap\fh_{J})))))\\
=\mathrm{Tot}(\Hom_{E}(B_{\bullet}(Z_{J}^{\dagger})_{Z_{J}^{\dagger}},C^{\bullet}(\fz_{I},\Hom_{E}(B_{\bullet}(Z_{I}^{\dagger}\cap L_{J}')_{Z_{J}^{\dagger}\cap H_{J}},C^{\bullet}(\fl_{I}\cap\fh_{J})))))
\end{multline}
Now that the inclusion (\ref{equ: Levi Lie discrete inclusion}) factors through $D(L_{I}\cap H_{J})\otimes_{E}U(\fz_{J})\otimes_{E}D(Z_{J}^{\dagger})$, we see that the quasi-isomorphism (\ref{equ: Levi Lie discrete cochain}) factors through
\begin{equation}\label{equ: Levi Lie discrete factor}
\mathrm{Tot}(\Hom_{E}(B_{\bullet}(Z_{J}^{\dagger})_{Z_{J}^{\dagger}},C^{\bullet}(\fz_{I},C^{\bullet}(L_{I}\cap H_{J})))).
\end{equation}
Now that the natural map (induced from the shuffle product map associated with $U(\fl_{I}\cap\fh_{J})\otimes_{E}D(Z_{I}^{\dagger}\cap H_{J})\subseteq D(L_{I}\cap H_{J})$, see \cite[Prop.~8.6.13]{Wei94})
\[C^{\bullet}(L_{I}\cap H_{J})\rightarrow \Hom_{E}(B_{\bullet}(Z_{I}^{\dagger}\cap L_{J}')_{Z_{J}^{\dagger}\cap H_{J}},C^{\bullet}(\fl_{I}\cap\fh_{J}))\]
is a quasi-isomorphism again by Eilenberg-Zilber's theorem and Lemma~\ref{lem: center compact}, we deduce that the map from (\ref{equ: Levi Lie discrete factor}) to RHS of (\ref{equ: Levi Lie discrete cochain}) is a quasi-isomorphism, and thus the map from LHS of (\ref{equ: Levi Lie discrete cochain}) to (\ref{equ: Levi Lie discrete factor}) is also a quasi-isomorphism. Combining with the following quasi-isomorphism (induced from the shuffle product map associated with the isomorphism $\varsigma_{J}: \prod_{j\notin J}\Z_{j}\rightarrow Z_{J}^{\dagger}$)
\[\bigotimes_{j\notin J}B_{\bullet}(\Z_{j})_{\Z_{j}}\rightarrow B_{\bullet}(Z_{J}^{\dagger})_{Z_{J}^{\dagger}},\]
we obtain a quasi-isomorphism
\begin{equation}\label{equ: Levi Lie discrete center factor}
C^{\bullet}(L_{I})\rightarrow \mathrm{Tot}(\Hom_{E}(\bigotimes_{j\notin J}B_{\bullet}(\Z_{j})_{\Z_{j}},C^{\bullet}(L_{I}\cap H_{J})))
\end{equation}
which is easily checked to be functorial with respect to the choice of $I\subseteq J$.

Let $I'\subseteq I\subseteq\Delta$. The definition of $Z_{I}^{\dagger}$ and $Z_{I'}^{\dagger}$ as above ensures that the inclusion $L_{I'}\subseteq L_{I}$ induces a map $L_{I'}'\times Z_{I'}^{\dagger}\rightarrow L_{I}'\times Z_{I}^{\dagger}$.
Let $\fh\subseteq \fg$ be an $E$-Lie subalgebra.
Upon replacing $G$, $G'$, $M$, $M'$, $\fh$ and $\fh'$ in (\ref{equ: loc an separate discrete}) with $L_{I}'$, $L_{I'}'$, $Z_{I}^{\dagger}$, $Z_{I'}^{\dagger}$, $\fl_{I}\cap\fh$ and $\fl_{I'}\cap\fh$ and then replacing the isogeny $\Z^{m}\rightarrow M$ and $\Z^{m'}\rightarrow M'$ in \emph{loc.cit.} with $\varsigma_{I}$ and $\varsigma_{I'}$, we deduce from (\ref{equ: G M splitting diagram}) (and the fact that the isogeny $L_{I}'\times Z_{I}^{\dagger}\rightarrow L_{I}$ restricts to the isogeny $L_{I'}'\times Z_{I'}^{\dagger}\rightarrow L_{I'}$) the following commutative diagram
\begin{equation}\label{equ: Levi separate discrete diagram}
\xymatrix{
C^{\bullet}(L_{I},\fl_{I}\cap\fh) \ar[r] \ar[d] & C^{\bullet}(L_{I}'\times Z_{I}^{\dagger},\fl_{I}\cap\fh) \ar[d] \\
C^{\bullet}(L_{I'},\fl_{I'}\cap\fh) \ar[r] & C^{\bullet}(L_{I'}'\times Z_{I'}^{\dagger},\fl_{I'}\cap\fh)
}
\end{equation}
with both horizontal maps being quasi-isomorphisms.
\begin{prop}\label{prop: Levi decomposition}
Let $\fh\subseteq\fg$ be an $E$-Lie subalgebra. Let $M\subseteq Z_{I}^{\dagger}$ be a subgroup.
Let $I'\subseteq I\subseteq \Delta$. We have the following results.
\begin{enumerate}[label=(\roman*)]
\item \label{it: Levi decomposition 0} The pair $(L_{I}/M,\fl_{I}\cap\fh)$ is degenerate.
\item \label{it: Levi decomposition 1}
The isogeny $Z_{I}^{\dagger}/M\times L_{I}'\rightarrow L_{I}/M$ induces the following isomorphisms between graded $E$-algebras
\begin{multline}\label{equ: Kunneth decomposition K}
H^{\bullet}(L_{I}/M,\fl_{I}\cap\fh,1_{L_{I}/M})\cong H^{\bullet}(Z_{I}^{\dagger}/M,1_{Z_{I}^{\dagger}/M})\otimes_E H^{\bullet}(L_{I}',\fl_{I}\cap\fh,1_{L_{I}'})\\
\cong \wedge^{\bullet}\Hom(Z_{I}^{\dagger}/M,E)\otimes_E H^{\bullet}(\fl_{I},\fl_{I}\cap\fh,1_{\fl_{I}}).
\end{multline}
\item \label{it: Levi decomposition 2} For each $E$-Lie subalgebra $\fh'\subseteq\fh$ and subgroup $M'\subseteq M$, the restriction map
\[H^{\bullet}(L_{I}/M,\fl_{I}\cap\fh,1_{L_{I}})\rightarrow H^{\bullet}(L_{I}/M',\fl_{I}\cap\fh',1_{L_{I}})\]
corresponds under (\ref{equ: Kunneth decomposition K}) to the tensor between the restriction map
\[H^{\bullet}(\fl_{I},\fl_{I}\cap\fh,1_{L_{I}})\rightarrow H^{\bullet}(\fl_{I},\fl_{I}\cap\fh',1_{L_{I}})\]
and the natural map $\wedge^{\bullet}\Hom(Z_{I}^{\dagger}/M,E)\rightarrow \wedge^{\bullet}\Hom(Z_{I}^{\dagger}/M',E)$ induced from $Z_{I}^{\dagger}/M'\twoheadrightarrow Z_{I}^{\dagger}/M$.
\end{enumerate}
\end{prop}
\begin{proof}
We prove \ref{it: Levi decomposition 0}.\\
Now that $(L_{I}',\fl_{I}\cap\fh)$ is degenerate by Lemma~\ref{lem: center compact}, the pair $(L_{I}'\times (Z_{I}^{\dagger}/M),\fl_{I}\cap\fh)$ is degenerate by Lemma~\ref{lem: Kunneth K discrete} (upon taking $G=L_{I}'$ and $H=Z_{I}^{\dagger}/M$ in \emph{loc.cit.}). As the map $L_{I}'\times (Z_{I}^{\dagger}/M)\rightarrow L_{I}$ is an isogeny, we conclude from Lemma~\ref{lem: group coh finite difference} that $(L_{I}/M,\fl_{I}\cap\fh)$ is degenerate.

We prove \ref{it: Levi decomposition 1}.\\
As the natural map $L_{I}'\times (Z_{I}^{\dagger}/M)\rightarrow L_{I}$ is an isogeny, we deduce (\ref{equ: Kunneth decomposition K}) from Lemma~\ref{lem: group coh finite difference}, Lemma~\ref{lem: center compact} and (\ref{equ: Kunneth graded algebra}).

Finally \ref{it: Levi decomposition 2} follows directly from the fact that the maps
\[C^{\bullet}(L_{I}/M,\fl_{I}\cap\fh) \rightarrow C^{\bullet}(L_{I}'\times (Z_{I}^{\dagger}/M),\fl_{I}\cap\fh)\rightarrow \mathrm{Tot}(\Hom_{E}(B_{\bullet}(Z_{I}^{\dagger}/M)_{Z_{I}^{\dagger}/M},C^{\bullet}(L_{I}',\fl_{I}\cap\fh)))\]
are functorial with respect to the choice of $\fh$ and $M$.
\end{proof}

\begin{lem}\label{lem: Levi parabolic coh}
Let $I\subseteq \Delta$ and $\fh\subseteq\fl_{I}$ be a reductive $E$-Lie subalgebra. We have the following results.
\begin{enumerate}[label=(\roman*)]
\item \label{it: Levi parabolic coh 1} The homomorphisms $L_{I}\hookrightarrow P_{I}\twoheadrightarrow L_{I}$ induce the following quasi-isomorphisms
\begin{equation}\label{equ: parabolic degenerate sm}
C^{\bullet}(L_{I})^{\infty} \buildrel\sim\over\longrightarrow C^{\bullet}(P_{I})^{\infty}\buildrel\sim\over\longrightarrow C^{\bullet}(L_{I})^{\infty}
\end{equation}
as well as the following quasi-isomorphisms
\begin{equation}\label{equ: parabolic degenerate}
C^{\bullet}(L_{I},\fh) \buildrel\sim\over\longrightarrow C^{\bullet}(P_{I},\fh)\buildrel\sim\over\longrightarrow C^{\bullet}(L_{I},\fh),
\end{equation}
with the associated maps between cohomology for (\ref{equ: parabolic degenerate sm}) and (\ref{equ: parabolic degenerate}) being isomorphisms between graded $E$-algebras.
\item \label{it: Levi parabolic coh 3} The pair $(P_{I},\fh)$ is degenerate. 
\end{enumerate}
\end{lem}
\begin{proof}
Since the composition of $L_{I}\hookrightarrow P_{I}\twoheadrightarrow L_{I}$ is the identity map of $L_{I}$, it suffices to study the map $C^{\bullet}(L_{I})^{\ast}\rightarrow C^{\bullet}(P_{I})^{\ast}$ induced from $P_{I}\twoheadrightarrow L_{I}$ for $\ast\in\{,\infty\}$.
As the $J_{N_I}: \mathrm{Rep}^{\infty}(P_I)\rightarrow \mathrm{Rep}^{\infty}(L_I)$ is exact and $J_{N_I}(1_{P_I})\cong 1_{L_I}$, we obtain canonical isomorphisms
\begin{equation}\label{equ: Levi parabolic sm Ext}
H^{\bullet}(L_{I},1_{L_{I}})^{\infty}\cong \mathrm{Ext}_{L_I}^k(J_{N_I}(1_{P_I}),1_{L_I})^{\infty}\buildrel\sim\over\longrightarrow H^{\bullet}(P_{I},1_{P_{I}})^{\infty}
\end{equation}
for $k\geq 0$.
We write $\xi: Z(\fl_{I})\rightarrow E$ for the unique infinitesimal character such that $L^{I}(1)_{\xi}\neq 0$.
Recall from \ref{it: n coh 1} and \ref{it: n coh 4} of Lemma~\ref{lem: n coh collection} (together with \cite[Thm.~1.10]{Hum08}) that $H^{0}(\fn_{I},L(1))=L^{I}(1)=1_{\fl_{I}}$ and that $H^{k}(\fn_{I},L(1))_{\xi}=0$ for each $k\geq 1$, which together with (\ref{equ: g spectral seq}) implies that $\fp_{I}\twoheadrightarrow \fl_{I}$ induces the following isomorphism
\begin{equation}\label{equ: g Levi parabolic coh}
H^k(\fl_{I},\fh,1_{\fl_I})\buildrel\sim\over\longrightarrow H^k(\fp_{I},\fh,1_{\fp_{I}})
\end{equation}
for $k\geq 0$.
Note that the functor $\mathrm{Mod}_{D(L_{I})}\rightarrow\mathrm{Mod}_{D(P_{I})}$ (induced from the homomorphism $P_I\twoheadrightarrow L_I$) leads to the following map between the spectral sequences
\[
\xymatrix{
\mathrm{Ext}_{D^{\infty}(L_I)}^k(1_{L_{I}}^{\vee},\mathrm{Ext}_{D(L_{I}),U(\fh)}^{\ell}(D^{\infty}(L_{I}),1_{L_{I}}^{\vee}) \ar@{=>}[r] \ar[d] & \mathrm{Ext}_{D(L_{I}),U(\fh)}^{k+\ell}(1_{L_{I}}^{\vee},1_{L_{I}}^{\vee}) \ar[d]\\
\mathrm{Ext}_{D^{\infty}(P_I)}^k(1_{P_{I}}^{\vee},\mathrm{Ext}_{D(P_{I}),U(\fh)}^{\ell}(D^{\infty}(P_{I}),1_{P_{I}}^{\vee}) \ar@{=>}[r] & \mathrm{Ext}_{D(P_{I}),U(\fh)}^{k+\ell}(1_{P_{I}}^{\vee},1_{P_{I}}^{\vee})
}
\]
which under Lemma~\ref{lem: ST spectral seq} (upon taking $\fh=0$, $D=1_{D(G)}$ and then replacing $G$ with $L_{I}$ or $P_{I}$ in \emph{loc.cit.} respectively) becomes
\begin{equation}\label{equ: Levi parabolic seq}
\xymatrix{
H^{k+\ell}(L_I,1_{L_I})^{\infty}\otimes_E H^{\ell}(\fl_{I},\fh,1_{\fl_I}) \ar@{=>}[r] \ar[d] & H^{k+\ell}(L_I,\fh,1_{L_I}) \ar[d]\\
H^{k+\ell}(P_I,1_{P_I})^{\infty}\otimes_E H^{\ell}(\fp_{I},\fh,1_{\fp_I}) \ar@{=>}[r] & H^{k+\ell}(P_I,\fh,1_{P_I})
}.
\end{equation}
Now that the map (\ref{equ: Levi parabolic seq}) between spectral sequences is an isomorphism on the first page by (\ref{equ: g Levi parabolic coh}) and (\ref{equ: Levi parabolic sm Ext}), we finish the proof of \ref{it: Levi parabolic coh 1}.
As the top row of (\ref{equ: Levi parabolic seq}) degenerates at the first page by \ref{it: Levi decomposition 0} of Proposition~\ref{prop: Levi decomposition}, we also conclude that the bottom row of (\ref{equ: Levi parabolic seq}) must also degenerate at the first page, namely the pair $(P_I,\fh)$ is degenerate, finishing the proof of \ref{it: Levi parabolic coh 3}.
\end{proof}

Let $x\in W(G)$, $I\subseteq I_{x}=\Delta\setminus D_L(x)$ and $\fh\subseteq\fp_{I}$ be an $E$-Lie subalgebra.
Recall that $\mathrm{Ext}_{G,\fh}^{\bullet}(1_{G},i_{x,I}^{\rm{an}})$ can be computed as the cohomology of the complex
\begin{equation}\label{equ: PS standard complex}
\Hom_{D(G)}(B_{\bullet}(D(G),U(\fh),(i_{x,I}^{\rm{an}})^\vee), 1_{D(G)}).
\end{equation}
The inclusion $P_{I}\subseteq G$ induces the following map between complex of free $D(G)$-modules
\begin{multline}\label{equ: PS parabolic standard complex}
D(G)\otimes_{D(P_{I})}(B_{\bullet}(D(P_{I}),U(\fh),L^{I}(x)))=D(G)\otimes_{U(\fh)}(\bigotimes_{U(\fh)}^{\bullet}D(P_{I})\otimes_{U(\fh)}L^{I}(x))\\
=D(G)\otimes_{U(\fh)}(\bigotimes_{U(\fh)}^{\bullet+1}D(P_{I})\otimes_{D(P_{I})}L^{I}(x))\rightarrow D(G)\otimes_{U(\fh)}(\bigotimes_{U(\fh)}^{\bullet+1}D(G)\otimes_{D(P_{I})}L^{I}(x))\\
=D(G)\otimes_{U(\fh)}(\bigotimes_{U(\fh)}^{\bullet}D(G)\otimes_{U(\fh)}(D(G)\otimes_{D(P_{I})}L^{I}(x)))=B_{\bullet}(D(G),U(\fh),(i_{x,I}^{\rm{an}})^\vee).
\end{multline}
Now that we have $L^{I}(x)$ is a $D(P_{I})$-module of finite $E$-dimension and thus
\begin{equation}\label{equ: Tor vanishing}
\mathrm{Tor}_k^{D(P_{I})}(D(G),L^{I}(x))=0
\end{equation}
for each $k>0$ by the argument in the first paragraph of \cite[Lem.~4.3.3]{BQ24} (modeled on \cite[Lem.~6.3]{ST05}), we observe that $D(G)\otimes_{D(P_{I})}(B_{\bullet}(D(P_{I}),U(\fh),L^{I}(x)))$ is exact and a resolution of $D(G)\otimes_{D(P_{I})}L^{I}(x)=(i_{x,I}^{\rm{an}})^\vee$ by $D(G)$-modules which are relative projective with respect to the pair $U(\fh)\subseteq D(G)$.
Hence, the composition of (\ref{equ: PS parabolic standard complex}) is a map between resolutions of $(i_{x,I}^{\rm{an}})^\vee$ by free $D(G)$-modules, which induces the following quasi-isomorphism
\begin{multline}\label{equ: PS parabolic resolution x}
\Hom_{D(G)}(B_{\bullet}(D(G),U(\fh),(i_{x,I}^{\rm{an}})^\vee), 1_{D(G)})\\
\rightarrow \Hom_{D(G)}(D(G)\otimes_{D(P_{I})}(B_{\bullet}(D(P_{I}),U(\fh),L^{I}(x))),1_{D(G)})=\Hom_{D(P_{I})}(B_{\bullet}(D(P_{I}),U(\fh),L^{I}(x)),1_{D(P_{I})}).
\end{multline}
Similarly, we also see that the quasi-isomorphism (\ref{equ: PS parabolic resolution x}) further factors through the following quasi-isomorphism
\[
\Hom_{D(\fg,P_{I})}(B_{\bullet}(D(\fg,P_{I}),U(\fh),M^{I}(x)), 1_{D(\fg,P_{I})})\rightarrow \Hom_{D(P_{I})}(B_{\bullet}(D(P_{I}),U(\fh),L^{I}(x)),1_{D(P_{I})})\]
with $D(\fg,P_{I})\defeq U(\fg)\otimes_{U(\fp_{I})}D(P_{I})\subseteq D(G)$.
To summarize, we have the following results.
\begin{lem}\label{lem: g P Ext transfer}
Let $x\in W(G)$ and $I\subseteq I_{x}$. The inclusions $D(P_{I})\subseteq D(\fg,P_{I})\subseteq D(G)$ between unital associative $E$-algebras and the embeddings $L^{I}(x)\hookrightarrow M^{I}(x)\hookrightarrow (i_{x,I}^{\rm{an}})^{\vee}$ induce the following quasi-isomorphisms
\begin{multline}\label{equ: g P Ext transfer}
\Hom_{D(G)}(B_{\bullet}(D(G),U(\fh),(i_{x,I}^{\rm{an}})^\vee), 1_{D(G)})\\
\rightarrow \Hom_{D(\fg,P_{I})}(B_{\bullet}(D(\fg,P_{I}),U(\fh),M^{I}(x)), 1_{D(\fg,P_{I})})\\
\rightarrow \Hom_{D(P_{I})}(B_{\bullet}(D(P_{I}),U(\fh),L^{I}(x)),1_{D(P_{I})})
\end{multline}
which are functorial with respect to the choice of $I$.
\end{lem}
The $x=1$ case of Lemma~\ref{lem: g P Ext transfer} together with Lemma~\ref{lem: Levi parabolic coh} implies that the inclusions $L_{I}\subseteq P_{I}\subseteq G$ together with $1_{L_{I}}^{\vee}\hookrightarrow (i_{I}^{\rm{an}})^{\vee}$ induce the following quasi-isomorphisms
\begin{equation}\label{equ: PS parabolic Levi resolution}
\Hom_{D(G)}(B_{\bullet}(D(G),U(\fh),(i_{I}^{\rm{an}})^{\vee}),1_{D(G)})\rightarrow C^{\bullet}(P_{I},\fh)\rightarrow C^{\bullet}(L_{I},\fh)
\end{equation}
which are functorial with respect to the choice of $I$.
Taking cohomology of (\ref{equ: PS parabolic Levi resolution}) gives
\[
\mathrm{Ext}_{G,\fh}^{\bullet}(1_{G},i_{I}^{\rm{an}})
\buildrel\sim\over\longrightarrow H^{\bullet}(P_{I},\fh,1_{P_{I}}) \buildrel\sim\over\longrightarrow H^{\bullet}(L_{I},\fh,1_{L_{I}})
\]
Similarly, the inclusions $\fl_{I}\subseteq \fp_{I}\subseteq \fg$ together with $1_{\fl_{I}}\hookrightarrow M^{I}(1)$ induce the following quasi-isomorphisms
\begin{equation}\label{equ: Verma parabolic Levi resolution}
\xymatrix{
\Hom_{U(\fg)}(B_{\bullet}(U(\fg),U(\fh),M^{I}(1)),1_{\fg}) \ar[r] \ar[d] & C^{\bullet}(\fp_{I},\fh) \ar[r] \ar[d] & C^{\bullet}(\fl_{I},\fh) \ar[d]\\
\mathrm{CE}^{\bullet}(\fg,\fh,M^{I}(1)) \ar[r] & \mathrm{CE}^{\bullet}(\fp_{I},\fh) \ar[r] & \mathrm{CE}^{\bullet}(\fl_{I},\fh)
}
\end{equation}
which are functorial with respect to the choice of $I$.
Taking cohomology of both row of (\ref{equ: Verma parabolic Levi resolution}) gives
\[
\mathrm{Ext}_{U(\fg),U(\fh)}^{\bullet}(M^{I}(1),1_{\fg})
\buildrel\sim\over\longrightarrow H^{\bullet}(\fp_{I},\fh,1_{\fp_{I}}) \buildrel\sim\over\longrightarrow H^{\bullet}(\fl_{I},\fh,1_{\fl_{I}}).
\]
\begin{lem}\label{lem: parabolic restriction x}
Let $x\in W(G)$, $I\subseteq J_{x}=\Delta\setminus\mathrm{Supp}(x)$ and $\fh\subseteq\fl_{I}$ be a reductive $E$-Lie subalgebra. 
Then we have a canonical isomorphism
\begin{equation}\label{equ: parabolic restriction x}
\mathrm{Ext}_{D(P_{I}),U(\fh)}^{k}(L^{I}(x),1_{P_{I}}^{\vee})\buildrel\sim\over\longrightarrow H^{k}(L_{I},\fh,1_{L_{I}})\otimes_E\mathrm{Ext}_{U(\fp_{I})}^{\ell(x)}(L^{I}(x),1_{\fp_{I}})
\end{equation}
for each $k\geq 0$ which is functorial with respect to the choice of $I$.
\end{lem}
\begin{proof}
Since $I\subseteq J_{x}$ and $\fh\subseteq\fl_{I}$, $\fn_{J_{x}}D(P_{I})\subseteq D(P_{I})$ is a two-sided ideal with the unital associative $E$-algebra $A\defeq D(P_{I})/\fn_{J_{x}}D(P_{I})$ naturally containing $D(L_{I})$ and thus $U(\fh)$ as associative $E$-subalgebra. Similar argument as in Lemma~\ref{lem: ST spectral seq} gives the following spectral sequence
\begin{multline}\label{equ: A Ext spectral seq}
\mathrm{Ext}_{A,U(\fh)}^{k}(L^{I}(x),H^{\ell}(\fn_{J_{x}},1_{D(P_{I})}))\\ =\mathrm{Ext}_{A,U(\fh)}^{k}(L^{I}(x),\mathrm{Ext}_{D(P_{I})}^{\ell}(A,1_{D(P_{I})}))\implies\mathrm{Ext}_{D(P_{I}),U(\fh)}^{k+\ell}(L^{I}(x),1_{D(P_{I})})
\end{multline}
Note that we also have $A/(\fp_{I}/\fn_{J_{x}})A=D^{\infty}(P_{I})$, and the inclusion $D(L_{I})\subseteq A$ induces the following map between spectral sequences
\begin{equation}\label{equ: A Levi spectral seq map}
\xymatrix{
\mathrm{Ext}_{D^{\infty}(P_{I})}^{i}(1_{P_{I}}^{\vee},H^{j}(\fp_{I}/\fn_{J_{x}},\fh,D)) \ar@{=>}[r] \ar[d] & \mathrm{Ext}_{A,U(\fh)}^{i+j}(1_{A},D) \ar[d]\\
\mathrm{Ext}_{D^{\infty}(L_{I})}^{i}(1_{L_{I}}^{\vee},H^{j}(\fl_{I},\fh,D)) \ar@{=>}[r] & \mathrm{Ext}_{D(L_{I}),U(\fh)}^{i+j}(1_{L_{I}}^{\vee},D)
}
\end{equation}
for each $D\in\mathrm{Mod}_{D(L_{I})}$ viewed as an object of $\mathrm{Mod}_{A}$ via the natural surjection $A\twoheadrightarrow D(L_{I})$. Now that $I\subseteq J_{x}$ and thus $L^{I}(x)$ is a $A$-module of finite $E$-dimension, $\Hom_{E}(L^{I}(x),-)$ is a well-defined endo-functor of $\mathrm{Mod}_{A}$ that sends relative injective objects to relative injective objects (with respect to the pair $U(\fh)\subseteq A$), we have natural isomorphisms
\begin{equation}\label{equ: A Ext x swap}
\mathrm{Ext}_{A,U(\fh)}^{k}(L^{I}(x),H^{\ell}(\fn_{J_{x}},1_{D(P_{I})}))\cong \mathrm{Ext}_{A,U(\fh)}^{k}(1_{A},\Hom_{E}(L^{I}(x),H^{\ell}(\fn_{J_{x}},1_{D(P_{I})})))
\end{equation}
for $k\geq 0$. We assume in the rest of the proof that
\[D\defeq \Hom_{E}(L^{I}(x),H^{\ell}(\fn_{J_{x}},1_{D(P_{I})}))\]
for some $\ell\geq 0$. 
We write $\xi:Z(\fl_{I})\rightarrow E$ for the unique infinitesimal character such that $L^{I}(1)_{\xi}\neq 0$ (with $L^{I}(1)=1_{\fl_{I}}$).
Recall from (an evident variant of) (\ref{equ: g spectral seq}) that we have the following spectral sequence
\begin{equation}\label{equ: A Ext n coh spectral seq}
H^{p}(\fl_{I},\fh,H^{q}(\fn_{I}/\fn_{J_{x}},D)_{\xi})\implies H^{p+q}(\fp_{I}/\fn_{J_{x}},\fh,D).
\end{equation}
It follows from Lemma~\ref{lem: n coh collection} that $H^{q}(\fn_{I}/\fn_{J_{x}},D)_{\xi}\neq 0$ if and only if $q=0$ and $\ell=\ell(x)$ in which case it is canonically isomorphic to the following $E$-vector space
\[D_{\xi}=\Hom_{U(\fl_{I})}(L^{I}(x),H^{\ell}(\fn_{J_{x}},1_{D(P_{I})}))\cong\mathrm{Ext}_{U(\fp_{I})}^{\ell(x)}(L^{I}(x),1_{\fp_{I}})\]
equipped with the trivial $U(\fl_{I})$-action.
In other words, (\ref{equ: A Ext n coh spectral seq}) degenerates to an isomorphism
\[H^{p}(\fl_{I},\fh,H^{0}(\fn_{I}/\fn_{J_{x}},D)_{\xi})\buildrel\sim\over\longrightarrow H^{p}(\fp_{I}/\fn_{J_{x}},\fh,D)\]
with BHS being possibly non-zero only if $\ell=\ell(x)$. Combining this with (\ref{equ: Levi parabolic sm Ext}), (\ref{equ: A Ext spectral seq}), (\ref{equ: A Levi spectral seq map}) and (\ref{equ: A Levi spectral seq map}), we see that (\ref{equ: A Ext spectral seq}) degenerates to the isomorphisms
\begin{multline}\label{equ: A Ext isom}
\mathrm{Ext}_{A,U(\fh)}^{k}(1_{A},1_{A})\cong\mathrm{Ext}_{A,U(\fh)}^{k}(1_{A},\Hom_{E}(L^{I}(x),H^{\ell(x)}(\fn_{J_{x}},1_{D(P_{I})})))\\
\cong\mathrm{Ext}_{A,U(\fh)}^{k}(L^{I}(x),H^{\ell(x)}(\fn_{J_{x}},1_{D(P_{I})}))\buildrel\sim\over\longrightarrow\mathrm{Ext}_{D(P_{I}),U(\fh)}^{k+\ell(x)}(L^{I}(x),1_{D(P_{I})})
\end{multline}
for each $k\geq 0$, and the map between spectral sequence (\ref{equ: A Levi spectral seq map}) is an isomorphism from the first page which in particular gives an isomorphism
\begin{equation}\label{equ: A Levi isom}
\mathrm{Ext}_{A,U(\fh)}^{\bullet}(1_{A},D)\buildrel\sim\over\longrightarrow H^{\bullet}(L_{I},\fh,D)=H^{\bullet}(L_{I},\fh,1_{L_{I}})\otimes_ED.
\end{equation}
We obtain (\ref{equ: parabolic restriction x}) as the composition of (\ref{equ: A Levi isom}) with the inverse of (\ref{equ: A Ext isom}), and finish the proof by noticing that both (\ref{equ: A Ext isom}) and (\ref{equ: A Levi isom}) are functorial with respect to the choice of $I$ (as the inclusion $P_{I'}\subseteq P_{I}$ induces an inclusion $D(P_{I'})/\fn_{J_{x}}D(P_{I'})\subseteq D(P_{I})/\fn_{J_{x}}D(P_{I})$ which further restricts to the inclusion $D(L_{I'})\subseteq D(L_{I})$, for each $I'\subseteq I$).
\end{proof}
Along the proof of Lemma~\ref{lem: g P Ext transfer} and Lemma~\ref{lem: parabolic restriction x}, we also have canonical isomorphisms
\begin{equation}\label{equ: Lie parabolic restriction x}
\mathrm{Ext}_{U(\fg),U(\fh)}^{k}(M^{I}(x),1_{\fg})\buildrel\sim\over\longrightarrow\mathrm{Ext}_{U(\fp_{I}),U(\fh)}^{k}(L^{I}(x),1_{\fp_{I}})\buildrel\sim\over\longrightarrow H^{k-\ell(x)}(\fl_{I},\fh,1_{\fl_{I}})
\end{equation}
for each $k\geq 0$ which is functorial with respect to the choice of $I$. Moreover, the inclusion $U(\fg)\subseteq D(\fg,P_{I})$ (which restricts to the inclusion $U(\fl_{I})\subseteq D(L_{I})$) further induces a commutative diagram of the form
\begin{equation}\label{equ: group to Lie parabolic restriction x}
\xymatrix{
\mathrm{Ext}_{D(\fg,P_{I}),U(\fh)}^{k}(M^{I}(x),1_{D(\fg,P_{I})}) \ar^{\sim}[rr] \ar[d] & & H^{k-\ell(x)}(L_{I},\fh,1_{L_{I}}) \ar[d] \\
\mathrm{Ext}_{U(\fg),U(\fh)}^{k}(M^{I}(x),1_{\fg}) \ar^{\sim}[rr] & & H^{k-\ell(x)}(\fl_{I},\fh,1_{\fl_{I}})
}.
\end{equation}
Note that the vertical maps of (\ref{equ: group to Lie parabolic restriction x}) are isomorphisms between $1$-dimensional $E$-vector spaces when $k=\ell(x)$.
Combining (the proof of) Lemma~\ref{lem: parabolic restriction x} with (\ref{equ: abstract ST cup diagram}) (by taking $R=D(P_{I})$ and $\overline{R}=D(P_{I})\fn_{J_{x}}D(P_{I})$ in \emph{loc.cit.}), we obtain the following commutative diagram (with $L^{I}(x)^{\ast}\defeq \Hom_{E}(L^{I}(x),E)$ for short)
\begin{equation}\label{equ: x sm cup}
\xymatrix{
H^{k_0}(L_{I},1_{L_{I}})^{\infty} \ar@{=}[d] & \otimes_E & \mathrm{Ext}_{D(P_{I})}^{k_1}(1_{D(P_{I})},L^{I}(x)^{\ast}) \ar^{\wr}[d] \ar^{\cup}[r] & \mathrm{Ext}_{D(P_{I})}^{k_0+k_1}(1_{D(P_{I})},L^{I}(x)^{\ast}) \ar^{\wr}[d]\\
H^{k_0}(L_{I},1_{L_{I}})^{\infty} & \otimes_E & H^{k_1-\ell(x)}(L_{I},1_{L_{I}})\otimes_EX \ar^{\cup}[r] & H^{k_0+k_1-\ell(x)}(L_{I},1_{L_{I}})\otimes_EX
}
\end{equation}
for each $k_0,k_1\geq 0$, with $X\defeq\mathrm{Ext}_{U(\fp_{I})}^{\ell(x)}(L^{I}(x),1_{\fp_{I}})$ for short.
Here the top horizontal cup product map is defined via
\[H^{k_0}(L_{I},1_{L_{I}})^{\infty}\rightarrow \mathrm{Ext}_{D(P_{I})}^{k_0}(1_{D(P_{I})},1_{D(P_{I})})\]
via $P_{I}\twoheadrightarrow L_{I}$, and both the middle and right vertical maps are isomorphisms from Lemma~\ref{lem: parabolic restriction x}.
For technical convenience, we consider the associative $E$-subalgebra $S\defeq U(\fp_{I})\otimes_ED(Z_{I}^{\dagger})\subseteq D(P_{I})$ and note that it admits a natural diagonal map $D_{S}: S\rightarrow S\otimes_ES$, so that we could define two version of standard resolutions $B_{\bullet}(S,-)$ and $\tld{B}_{\bullet}(S,-)$ (similar to Construction~\ref{cons: resolution}) and relate them using $D_{S}$. In fact, we could define cup product map of $\mathrm{Ext}_{S}^{\bullet}(-,-)$ using either $B_{\bullet}(S,-)$ and $\tld{B}_{\bullet}(S,-)$. We thus obtain the following commutative diagram
\begin{equation}\label{equ: x sm cup Lie transfer}
\xymatrix{
\mathrm{Ext}_{D(P_{I})}^{k_0}(1_{D(P_{I})},1_{D(P_{I})})\ar^{\wr}[d] & \otimes_E &\mathrm{Ext}_{D(P_{I})}^{k_1}(1_{D(P_{I})},L^{I}(x)^{\ast}) \ar^{\cup}[r] \ar^{\wr}[d]& \mathrm{Ext}_{D(P_{I})}^{k_0+k_1}(1_{D(P_{I})},L^{I}(x)^{\ast}) \ar^{\wr}[d]\\
\mathrm{Ext}_{S}^{k_0}(1_{S},1_{S})\ar^{\wr}[d] & \otimes_E &\mathrm{Ext}_{S}^{k_1}(1_{S},L^{I}(x)^{\ast}) \ar^{\cup}[r] \ar^{\wr}[d]& \mathrm{Ext}_{S}^{k_0+k_1}(1_{S},L^{I}(x)^{\ast}) \ar^{\wr}[d]\\
\mathrm{Ext}_{S}^{k_0}(1_{S},1_{S}) & \otimes_E &\mathrm{Ext}_{S}^{k_1}(L^{I}(x),1_{S}) \ar^{\cup}[r] & \mathrm{Ext}_{S}^{k_0+k_1}(L^{I}(x),1_{S}) \\
\mathrm{Ext}_{D(P_{I})}^{k_0}(1_{D(P_{I})},1_{D(P_{I})})\ar_{\wr}[u] & \otimes_E &\mathrm{Ext}_{D(P_{I})}^{k_1}(L^{I}(x),1_{D(P_{I})}) \ar^{\cup}[r] \ar_{\wr}[u]& \mathrm{Ext}_{D(P_{I})}^{k_0+k_1}(L^{I}(x),1_{D(P_{I})}) \ar_{\wr}[u]
}
\end{equation}
with all vertical maps being natural isomorphisms, with both the isomorphisms from the first row to the second row and the isomorphisms from the fourth row to the third induced from $S\subseteq P_{I}$ (and thus $B_{\bullet}(S,-)\rightarrow B_{\bullet}(D(P_{I}),-)$), and the isomorphisms between the second row and the third row using $B_{\bullet}(S,-)\cong\tld{B}_{\bullet}(S,-)$.
The definition of the cup product map in the fourth row of (\ref{equ: x sm cup Lie transfer}) uses the maps
\begin{multline*}
\Hom_{D(P_{I})}(B_{\bullet}(D(P_{I})),1_{D(P_{I})})\rightarrow \Hom_{D(P_{I})}(L^{I}(x)\otimes_EB_{\bullet}(D(P_{I})),L^{I}(x))\\
\rightarrow \mathrm{Tot}(\Hom_{D(P_{I})}(B_{\bullet}(D(P_{I}),L^{I}(x)\otimes_EB_{\bullet}(D(P_{I}))),B_{\bullet}(D(P_{I}),L^{I}(x))))
\end{multline*}
and similar construction will be discussed into details in \S~\ref{subsec: cup Tits double} (cf.~(\ref{equ: x y cup term}))
\subsection{Orlik--Strauch functor}\label{subsec: OS}
We recall from \cite{BQ24} an important spectral sequence (see Proposition~\ref{prop: St key seq}) to compute $\mathrm{Ext}_{G}^{\bullet}(-,-)$ between admissible locally $K$-analytic representations that can be constructed using Orlik-Strauch functors (see \cite{OS15}). This spectral sequence has many useful consequences (cf.~Lemma~\ref{lem: OS socle cosocle}, Lemma~\ref{lem: Ext1 OS cube} and Lemma~\ref{lem: Ext with sm}) which in particular help us study the layer structure of locally $K$-analytic generalized Steinberg representations (see Lemma~\ref{lem: basic subquotient of St} and Lemma~\ref{lem: C coxeter factor}.)
Note that the set of constituents of a locally $K$-analytic generalized Steinberg representation has been studied in \cite{OS13}.

We assume that $G=\mathrm{PGL}_n(K)$.
Given $I\subseteq \Delta$, $M$ in $\cO_{\rm{alg}}^{\fp_I}$ and $\pi^{\infty}$ in $\mathrm{Rep}^{\infty}_{\rm{adm}}(L_I)$, Orlik and Strauch define $\cF_{P_I}^G(M,\pi^{\infty})$ in $\mathrm{Rep}^{\rm{an}}_{\rm{adm}}(G)$, see \cite{OS15}.
We collect below some main properties of the functor $\cF_{P_I}^G(-,-)$ from \cite{OS15}.
\begin{prop}\label{prop: OS property}
We have the following results.
\begin{enumerate}[label=(\roman*)]
\item \label{it: OS property 1} The functor $\cF_{P_I}^G(-,-)$ is contravariant (resp.~covariant) in the first (resp.~the second) argument, and is exact for both arguments.
\item \label{it: OS property 2} For each $I_1\subseteq I\subseteq \Delta$, $M\in\cO^{\fp_I}_{\rm{alg}}$ and $\pi_1^{\infty}\in\mathrm{Rep}^{\infty}_{\rm{adm}}(L_{I_1})$, we have a canonical isomorphism
    $$\cF_{P_{I_1}}^G(M,\pi_1^{\infty})\cong \cF_{P_I}^G(M,i_{I_1,I}^{\infty}(\pi_1^{\infty})).$$
\item \label{it: OS property 3} Let $I\subseteq \Delta$, $M\in\cO^{\fp_I}_{\rm{alg}}$ and $\pi^{\infty}\in\mathrm{Rep}^{\infty}_{\rm{adm}}(L_I)$ with $I$ maximal for $M$. If both $M$ and $\pi^{\infty}$ are simple, then so is $\cF_{P_I}^G(M,\pi^{\infty})$.
\item \label{it: OS property 4} For each $I\subseteq \Delta$, $\mu\in \Lambda_I^+$ and $\pi^{\infty}\in\mathrm{Rep}^{\infty}_{\rm{adm}}(L_I)$, we have
    $$\cF_{P_I}^G(M^I(\mu),\pi^{\infty})\cong i_I^{\rm{an}}(L^I(\mu)^\vee\otimes_E \pi^{\infty}).$$
\end{enumerate}
\end{prop}

Given $I\subseteq\Delta$, we recall the Bernstein block $\cB^{I}$ from the paragraph below (\ref{equ: sm St}).
\begin{lem}\label{lem: JH OS}
Let $I\subseteq \Delta$, $M\in\cO^{\fp_I}_{\rm{alg}}$ and $\pi^{\infty}\in\cB^I$. Let $x\in W(G)$ and $\sigma^{\infty}\in\cB^{I_x}$ be simple. Then we have \[[\cF_{P_I}^{G}(M,\pi^{\infty}):\cF_{P_{I_x}}^{G}(L(x),\sigma^{\infty})]\neq 0\]
only if $[M:L(x)]\neq 0$ (which forces $I\subseteq I_x$ by Lemma~\ref{lem: dominance and left set}), in which case we have
\begin{equation}\label{equ: JH OS}
[\cF_{P_I}^{G}(M,\pi^{\infty}):\cF_{P_{I_x}}^{G}(L(x),\sigma^{\infty})]=[M:L(x)][i_{I,I_x}^{\infty}(\pi^{\infty}):\sigma^{\infty}].
\end{equation}
\end{lem}
\begin{proof}
This follows directly from \ref{it: OS property 1} and \ref{it: OS property 2} of Proposition~\ref{prop: OS property} as well as \cite[Lem.~5.1.1]{BQ24}.
\end{proof}

\begin{lem}\label{lem: OS vers induction}
Let $I'\subseteq I\subseteq \Delta$, $M_I\in\cO^{\fl_I\cap\fp_{I'}}_{\fl_I,\rm{alg}}$ and $\pi^{\infty}\in\mathrm{Rep}^{\infty}_{\rm{adm}}(L_{I'})$. Then we have the following isomorphism in $\mathrm{Rep}^{\rm{an}}_{\rm{adm}}(G)$
\begin{equation}\label{equ: OS vers induction}
i_{I}^{\rm{an}}(\cF_{L_I\cap P_{I'}}^{L_I}(M_I,\pi^{\infty}))\cong \cF_{P_{I'}}^{G}(U(\fg)\otimes_{U(\fp_I)}M_I, \pi^{\infty}).
\end{equation}
\end{lem}
\begin{proof}
We choose a finite dimensional $U(\fl_I\cap\fp_{I'})$-module which is $\ft$-semi-simple such that we have a surjection
\[q_I: U(\fl_I)\otimes_{U(\fl_I\cap\fp_{I'})}X\twoheadrightarrow M_I\]
in $\cO^{\fl_I\cap\fp_{I'}}_{\fl_I,\rm{alg}}$, which induces a surjection
\[q: U(\fg)\otimes_{U(\fp_{I'})}X\cong U(\fg)\otimes_{U(\fp_I)}(U(\fp_I)\otimes_{U(\fp_{I'})}X)\cong U(\fg)\otimes_{U(\fp_I)}(U(\fl_I)\otimes_{U(\fl_I\cap\fp_{I'})}X)\twoheadrightarrow U(\fg)\otimes_{U(\fp_I)}M_I\]
in $\cO^{\fp_{I'}}_{\rm{alg}}$.
In the rest of this proof we freely borrow notation from \cite[\S 4.2, 4.3]{BQ24} with the only exception being that we write $G^0$, $P_I^0=G^0\cap P_I$ and $L_I^0=G^0\cap L_I$ for $G_0$, $P_{I,0}$ and $L_{I,0}$ in \emph{loc.cit.} (and this is because in this paper we use $G_0$, $P_{I,0}$ and $L_{I,0}$ for the locally $\Q_p$-analytic group underlying $G$, $P_I$ and $L_I$).
Recall that (cf.~the discussion after \cite[Lem.~4.2.2]{BQ24}) $D(G^0)\cong \varprojlim_{r\in\cI}D(G^0)_r$ is a Fr\'echet--Stein algebra, and thus one can define the abelian category of coadmissible $D(G^0)$-modules. Similar facts hold for $D(G^0)_{P_I^0}$, $D(G^0)_{P_{I'}^0}$ and $D(G^0)_1=D(G)_1$ etc.
Recall from \cite[Prop.~4.3.1, Lem.~4.3.3]{BQ24} that we have the following short exact sequence of coadmissible $D(L_I)_1$-modules
\begin{multline}\label{equ: Levi completion sequence}
0\rightarrow D(L_I)_{L_I\cap P_{I'}}\cdot\mathrm{ker}(q_I)\rightarrow D(L_I)_{L_I\cap P_{I'}}\otimes_{D(L_I\cap P_{I'})}X\\
\rightarrow \cM_I\defeq D(L_I)_{L_I\cap P_{I'}}\otimes_{D(L_I\cap P_{I'})}X/D(L_I)_{L_I\cap P_{I'}}\cdot\mathrm{ker}(q_I)\rightarrow 0
\end{multline}
with $\cM_I$, $D(L_I)_{L_I\cap P_{I'}}\otimes_{D(L_I\cap P_{I'})}X$ and $D(L_I)_{L_I\cap P_{I'}}\cdot\mathrm{ker}(q_I)$ being the canonical Fr\'echet completion of $M_I$, $U(\fl_I)\otimes_{U(\fl_I\cap\fp_{I'})}X$ and $\mathrm{ker}(q_I)$ respectively as coadmissible $D(L_I)_{L_I\cap P_{I'}}$-modules. We could write (\ref{equ: Levi completion sequence}) as the inverse limit of a short exact sequence of finitely generated $D(L_I^0)_{L_I^0\cap P_{I'}^0,r}$-modules for $r\in\cI$, and then apply $D(G^0)_{P_{I'}^0,r}\otimes_{D(P_I^0)_{P_{I'}^0,r}}-$ and finally take $\varprojlim_{r\in\cI}$. Using \cite[Lem.~4.2.4]{BQ24} we obtain a strict short exact sequence of Fr\'echet $D(G)_{P_{I'}}$-modules
\begin{multline}\label{equ: induction completion sequence}
0\rightarrow D(G)_{P_{I'}}\widehat{\otimes}_{D(P_I)_{P_{I'}},\iota}(D(L_I)_{L_I\cap P_{I'}}\cdot\mathrm{ker}(q_I))\\
\rightarrow D(G)_{P_{I'}}\widehat{\otimes}_{D(P_I)_{P_{I'}},\iota}(D(L_I)_{L_I\cap P_{I'}}\otimes_{D(L_I\cap P_{I'})}X)\rightarrow D(G)_{P_{I'}}\widehat{\otimes}_{D(P_I)_{P_{I'}},\iota}\cM_I\rightarrow 0.
\end{multline}
Note that
\begin{multline*}
D(G)_{P_{I'}}\widehat{\otimes}_{D(P_I)_{P_{I'}},\iota}(D(L_I)_{L_I\cap P_{I'}}\otimes_{D(L_I\cap P_{I'})}X)\\
\cong D(G)_{P_{I'}}\widehat{\otimes}_{D(P_I)_{P_{I'}},\iota}(D(P_I)_{P_{I'}}\otimes_{D(P_{I'})}X)\cong D(G)_{P_{I'}}\otimes_{D(P_{I'})}X
\end{multline*}
are coadmissible $D(G)_{P_{I'}}$-modules, we deduce that \[D(G)_{P_{I'}}\widehat{\otimes}_{D(P_I)_{P_{I'}},\iota}(D(L_I)_1\cdot\mathrm{ker}(q_I))\] 
and 
\[D(G)_{P_{I'}}\widehat{\otimes}_{D(P_I)_{P_{I'}},\iota}\cM_I\] 
are also coadmissible $D(G)_{P_{I'}}$-modules. Now that $\mathrm{ker}(q)=U(\fg)\otimes_{U(\fp_I)}\mathrm{ker}(q_I)$ is dense in \[D(G)_1\widehat{\otimes}_{D(P_I)_1}(D(L_I)_1\cdot\mathrm{ker}(q_I))\cong D(G)_{P_{I'}}\widehat{\otimes}_{D(P_I)_{P_{I'}},\iota}(D(L_I)_{L_I\cap P_{I'}}\cdot\mathrm{ker}(q_I))\]
as an $E$-subspace of
\[D(G)_1\widehat{\otimes}_{D(P_I)_1}(D(L_I)_1\otimes_{D(L_I\cap P_{I'})_1}X)\cong D(G)_1\otimes_{D(P_{I'})_1}X\cong D(G)_{P_{I'}}\otimes_{D(P_{I'})}X,\]
we deduce that
\[D(G)_{P_{I'}}\widehat{\otimes}_{D(P_I)_{P_{I'}},\iota}(D(L_I)_{L_I\cap P_{I'}}\cdot\mathrm{ker}(q_I))=D(G)_{P_{I'}}\cdot\mathrm{ker}(q)\]
as closed $D(G)_{P_{I'}}$-submodules of $D(G)_{P_{I'}}\otimes_{D(P_{I'})}X$, and in particular
\begin{equation}\label{equ: canonical completion induction}
D(G)_{P_{I'}}\widehat{\otimes}_{D(P_I)_{P_{I'}},\iota}\cM_I\cong D(G)_{P_{I'}}\otimes_{D(P_{I'})}X/D(G)_{P_{I'}}\cdot\mathrm{ker}(q).
\end{equation}
In other words, $D(G)_{P_{I'}}\widehat{\otimes}_{D(P_I)_{P_{I'}}}\cM_I$ is the canonical Fr\'echet completion of $U(\fg)\otimes_{U(\fp_I)}M_I$ as a coadmissible $D(G)_{P_{I'}}$-module (see \cite[Prop.~4.3.1, Lem.~4.3.3]{BQ24}).
Recall from \cite[Prop.~4.3.6]{BQ24} that we have the following isomorphism between coadmissible $D(L_I)$-modules
\[\cF_{L_I\cap P_{I'}}^{L_I}(M_I,\pi^{\infty})^\vee\cong D(L_I)\widehat{\otimes}_{D(L_I)_{L_I\cap P_{I'}},\iota}(\cM_I\widehat{\otimes}_E(\pi^{\infty})^\vee),\]
which together with \cite[Prop.~2.1.2]{Eme} (upon rephrasing the proof of \emph{loc.cit.} using distribution algebras) gives the following isomorphisms between coadmissible $D(G)$-modules
\begin{multline}\label{equ: OS induction as tensor}
i_{I}^{\rm{an}}(\cF_{L_I\cap P_{I'}}^{L_I}(M_I,\pi^{\infty}))^\vee\cong D(G)\widehat{\otimes}_{D(P_I),\iota}(\cF_{L_I\cap P_{I'}}^{L_I}(M_I,\pi^{\infty})^\vee)\\
\cong D(G)\widehat{\otimes}_{D(P_I),\iota}(D(L_I)\widehat{\otimes}_{D(L_I)_{L_I\cap P_{I'}},\iota}(\cM_I\widehat{\otimes}_E(\pi^{\infty})^\vee))\cong D(G)\widehat{\otimes}_{D(P_I),\iota}(D(P_I)\widehat{\otimes}_{D(P_I)_{P_{I'}},\iota}(\cM_I\widehat{\otimes}_E(\pi^{\infty})^\vee))\\
\cong D(G)\widehat{\otimes}_{D(P_I)_{P_{I'}},\iota}(\cM_I\widehat{\otimes}_E(\pi^{\infty})^\vee)\cong D(G)\widehat{\otimes}_{D(G)_{P_{I'}},\iota}(D(G)_{P_{I'}}\widehat{\otimes}_{D(P_I)_{P_{I'}},\iota}(\cM_I\widehat{\otimes}_E(\pi^{\infty})^\vee))
\end{multline}
Using similar argument as in \cite[Lem.~4.3.5 (i)]{BQ24}, we have the following isomorphism between coadmissible $D(G)_{P_{I'}}$-modules
\begin{equation}\label{equ: OS move tensor}
D(G)_{P_{I'}}\widehat{\otimes}_{D(P_I)_{P_{I'}},\iota}(\cM_I\widehat{\otimes}_E(\pi^{\infty})^\vee)\cong (D(G)_{P_{I'}}\widehat{\otimes}_{D(P_I)_{P_{I'}},\iota}\cM_I)\widehat{\otimes}_E(\pi^{\infty})^\vee.
\end{equation}
Now that $D(G)_{P_{I'}}\widehat{\otimes}_{D(P_I)_{P_{I'}},\iota}\cM_I$ is the canonical Fr\'echet completion of $U(\fg)\otimes_{U(\fp_{I'})}M_I\in \cO^{\fp_{I'}}_{\rm{alg}}$ as a coadmissible $D(G)_{P_{I'}}$-module, we conclude (\ref{equ: OS vers induction}) from (\ref{equ: OS induction as tensor}), (\ref{equ: OS move tensor}) and \cite[Prop.~4.3.6]{BQ24}.
\end{proof}

Let $x\in W(G)$. For $I'\subseteq I\subseteq I_x=\Delta\setminus D_L(x)$, we define (using \ref{it: OS property 4} of Proposition~\ref{prop: OS property}, by replacing $G$, $P_I$ and $L^I(\mu)$ with $L_{I}$, $L_{I}\cap P_{I'}$ and $L^{I'}(x)$ here)
\[i_{x,I',I}^{\rm{an}}\defeq i_{I',I}^{\rm{an}}(L^{I'}(x)^\vee)\cong \cF_{L_{I}\cap P_{I'}}^{L_{I}}(U(\fl_{I})\otimes_{U(\fl_{I}\cap\fp_{I'})}L^{I'}(x),1_{L_{I'}}).\]
For each $I''\subseteq I'\subseteq I\subseteq I_x$, the surjection (\ref{equ: Verma transfer surjection}) together with the injection $\theta_{I',I''}^{\infty}: 1_{L_{I'}}\hookrightarrow i_{I'',I'}^{\infty}(1_{L_{I''}})$ (see (\ref{equ: sm PS injection})) as well as \ref{it: OS property 1} and \ref{it: OS property 2} of Proposition~\ref{prop: OS property} induces the following embeddings
\begin{multline}\label{equ: x loc an PS embedding}
i_{x,I',I}^{\rm{an}}=\cF_{L_{I}\cap P_{I'}}^{L_{I}}(U(\fl_{I})\otimes_{U(\fl_{I}\cap\fp_{I'})}L^{I'}(x),1_{L_{I'}})\hookrightarrow \cF_{L_{I}\cap P_{I'}}^{L_{I}}(U(\fl_{I})\otimes_{U(\fl_{I}\cap\fp_{I'})}L^{I'}(x),i_{I'',I'}^{\infty}(1_{L_{I''}}))\\
\cong\cF_{L_{I}\cap P_{I''}}^{L_{I}}(U(\fl_{I})\otimes_{U(\fl_{I}\cap\fp_{I'})}L^{I'}(x),1_{L_{I''}})
\hookrightarrow \cF_{L_{I}\cap P_{I''}}^{L_{I}}(U(\fl_{I})\otimes_{U(\fl_{I}\cap\fp_{I''})}L^{I''}(x),1_{L_{I''}})=i_{x,I'',I}^{\rm{an}}.
\end{multline}
We write $\theta^{x,I}_{I',I''}$ for the composition of (\ref{equ: x loc an PS embedding}). 
Note that we have $q^{x,I}_{I''',I'}=q^{x,I}_{I'',I'}\circ q^{x,I}_{I''',I''}$ for each $I'''\subseteq I''$ and that $q^{x,J}_{I'',I'}$ restricts to $q^{x,I}_{I'',I'}$ for each $I\subseteq J\subseteq I_x$.
This together with (\ref{equ: sm consistent embedding}) (by replacing $I,I',I''$ in \emph{loc.cit.} with $I',I'',I'''$ here) implies that
\begin{equation}\label{equ: consistent PS embedding}
\theta^{x,I}_{I',I'''}=\theta^{x,I}_{I'',I'''}\circ \theta^{x,I}_{I',I''}
\end{equation}
for each $I'''\subseteq I''\subseteq I'\subseteq I\subseteq I_{x}$.
For each $J\supseteq I$, it follows from Lemma~\ref{lem: OS vers induction} (by replacing $G$, $M_I$ and $\pi^{\infty}$ in \emph{loc.cit.} with $L_J$, $U(\fl_I)\otimes_{U(\fl_I\cap\fp_{I'})}L^{I'}(x)$ and $1_{L_{I'}}$) that
\begin{equation}\label{equ: induction x PS}
i_{I,J}^{\rm{an}}(i_{x,I',I}^{\rm{an}})\cong i_{x,I',J}^{\rm{an}}
\end{equation}
and we have the following equality between maps
\begin{equation}\label{equ: induction x PS embedding}
i_{I,J}^{\rm{an}}(\theta^{x,I}_{I',I''})=\theta^{x,J}_{I',I''}: i_{x,I',J}^{\rm{an}}\hookrightarrow i_{x,I'',J}^{\rm{an}}.
\end{equation}
For each $I\subseteq \Delta$ and $I'\subseteq I\cap J_{x}=I\setminus\mathrm{Supp}(x)$, we set
\begin{equation}\label{equ: loc an St x}
V_{x,I',I}^{\rm{an}}\defeq i_{x,I',I}^{\rm{an}}/\sum_{I'\subsetneq I''\subseteq J_{x}}i_{x,I'',I}^{\rm{an}}.
\end{equation}
When $x=1$ (resp.~$I=\Delta$), we omit it from the notation.
In particular, we obtain $i_{I',I}^{\rm{an}}$, $V_{I',I}^{\rm{an}}$, $\theta^{I}_{I',I''}$, $i_{x,I'}^{\rm{an}}$ and $V_{x,I'}^{\rm{an}}$ etc.

Let $x\in W(G)$, $I\subseteq \Delta$ and $I_0\subseteq I_1\subseteq I\cap I_{x}=I\setminus D_L(x)$.
We consider a complex in $\mathrm{Rep}^{\rm{an}}_{\rm{adm}}(L_{I})$ of the form
\begin{equation}\label{equ: x Tits complex}
\mathbf{C}^{x,I}_{I_0,I_1}\defeq [i_{x,I_1,I}^{\rm{an}}\rightarrow \cdots \rightarrow i_{x,I_0,I}^{\rm{an}}],
\end{equation}
whose degree $-\ell$ term is given by
\[\bigoplus_{I_0\subseteq I'\subseteq I_1, \#I'=\ell}i_{x,I',I}^{\rm{an}},\]
with the differential map at degree $-\ell$ restricts to
\[(-1)^{m(I',j)}\theta^{x,I}_{I',I'\setminus\{j\}}: i_{x,I',I}^{\rm{an}}\rightarrow i_{x,I'\setminus\{j\},I}^{\rm{an}}\]
for each $I_0\subseteq I'\subseteq I_1$ and $j\in I'\setminus I_0$ satisfying $\#I'=\ell$ (see (\ref{equ: differential sign}) for $m(I',j)$).
Note that the $U(\fl_{I})$-equivariant embedding $U(\fl_{I})\otimes_{U(\fl_{I}\cap\fp_{I'})}L^{I'}(x)\hookrightarrow (i_{x,I',I}^{\rm{an}})^\vee$ for each $I_0\subseteq I'\subseteq I_1$ induces an embedding (see (\ref{equ: Lie x Tits complex}))
\begin{equation}\label{equ: Lie to loc an Tits}
\mathfrak{c}^{x,I}_{I_0,I_1}\hookrightarrow (\mathbf{C}^{x,I}_{I_0,I_1})^\vee
\end{equation}
between two complex of $U(\fl_{I})$-modules supported in degree $[\#I_0,\#I_1]$.
For each $x\in W(G)$, $I\subseteq \Delta$, $I_0\subseteq I_1\subseteq I\cap J_{x}$ and $I_0\subseteq I_0',I_1'\subseteq I_1$, we have obvious truncation maps
\begin{equation}\label{equ: general Tits truncation}
\mathbf{C}^{x,I}_{I_0,I_1'}\rightarrow \mathbf{C}^{x,I}_{I_0,I_1}\rightarrow \mathbf{C}^{x,I}_{I_0',I_1}.
\end{equation}
Given $J\supseteq I$, it follows from (\ref{equ: induction x PS}) and (\ref{equ: induction x PS embedding}) that we have the following isomorphism of complex
\begin{equation}\label{equ: induction Tits complex}
i_{I,J}^{\rm{an}}(\mathbf{C}^{x,I}_{I_0,I_1})\cong \mathbf{C}^{x,J}_{I_0,I_1}.
\end{equation}
It follows from the last paragraph of the proof of \cite[Thm.~4.2]{OS13} that $\mathbf{C}^{x,I_1}_{I_0,I_1}\cong V_{x,I_0,I_1}^{\rm{an}}[\#I_0]$, which together with (\ref{equ: induction Tits complex}) gives the following isomorphism
\begin{equation}\label{equ: general Tits resolution}
\mathbf{C}^{x,I}_{I_0,I_1}\cong i_{I_1,I}^{\rm{an}}(V_{x,I_0,I_1}^{\rm{an}})[\#I_0]
\end{equation}
in the derived sense.
If $I=\Delta$, we write $\mathbf{C}^{x}_{I_0,I_1}\defeq \mathbf{C}^{x,\Delta}_{I_0,I_1}$ for short. If furthermore $x=1$, we write $\mathbf{C}_{I_0,I_1}\defeq \mathbf{C}^{1}_{I_0,I_1}$ for short.

Let $x\in W(G)$ and $I_0\subseteq J_{x}$.
For each $I_0\subseteq I',J\subseteq J_{x}$, the representation $i_{x,I'}^{\rm{an}}$ contains a subrepresentation
\[\tau^{J}(i_{x,I'}^{\rm{an}})\defeq \cF_{P_{I'}}^{G}(M^{J\cup I'}(x),1_{L_{I'}})\cong \cF_{P_{J\cup I'}}^{G}(M^{J\cup I'}(x),i_{I',J\cup I'}^{\infty}),\]
so that the complex $\mathbf{C}^{x}_{I_0,J_{x}}$ contains a subcomplex $\tau^{J}(\mathbf{C}^{x}_{I_0,J_{x}})$ whose degree $-\ell$ term is given by
\[\bigoplus_{I_0\subseteq I'\subseteq J_{x}, \#I'=\ell}\tau^{J}(i_{x,I'}^{\rm{an}}).\]
Now that $M^{I'\cup J}(x)$ is naturally isomorphic to the coproduct of $M^{I_0\cup J}(x)$ and $M^{I'}(x)$ over $M^{I_0}(x)$ by \cite[Thm~9.4(c)]{Hum08}, we deduce from Proposition~\ref{prop: OS property} that
\[\tau^{J}(i_{x,I'}^{\rm{an}})=\tau^{J}(i_{x,I_0}^{\rm{an}})\cap i_{x,I'}^{\rm{an}}\subseteq i_{x,I_0}^{\rm{an}}\]
for each $I_0\subseteq I'\subseteq J_{x}$, and thus
\[\tau^{J}(\mathbf{C}^{x}_{I_0,J_{x}})\cong \tau^{J}(V_{x,I_0}^{\rm{an}})[\#I_0]\]
in the derived sense, where $\tau^{J}(V_{x,I_0}^{\rm{an}})$ is the image of the composition of \[\tau^{J}(i_{x,I_0}^{\rm{an}})\hookrightarrow i_{x,I_0}^{\rm{an}}\twoheadrightarrow V_{x,I_0}^{\rm{an}}.\]
We write
\[\mathrm{gr}^{J}(i_{x,I'}^{\rm{an}})\defeq \tau^{J}(i_{x,I'}^{\rm{an}})/\sum_{J'\supsetneq }\tau^{J'}(i_{x,I'}^{\rm{an}})\]
for short, and note that it is zero if $I'\not\subseteq J$, and is isomorphic to
\[\cF_{P_{I'}}^{G}(\mathfrak{v}_{x,J},1_{L_{I'}})\cong \cF_{P_{J}}^{G}(\mathfrak{v}_{x,J},i_{I',J}^{\infty})\]
if $I'\subseteq J$ (see (\ref{equ: Lie St x}) for $\mathfrak{v}_{x,J}$ here).
In particular,
\[\mathrm{gr}^{J}(\mathbf{C}^{x}_{I_0,J_{x}})\defeq \tau^{J}(\mathbf{C}^{x}_{I_0,J_{x}})/\sum_{J'\supsetneq J}\tau^{J'}(\mathbf{C}^{x}_{I_0,J_{x}})\]
is a complex whose degree $-\ell$ term is given by
\[\bigoplus_{I_0\subseteq I'\subseteq J, \#I'=\ell}\cF_{P_{J}}^{G}(\mathfrak{v}_{x,J},i_{I',J}^{\infty}),\]
which moreover satisfies
\[\mathrm{gr}^{J}(\mathbf{C}^{x}_{I_0,J_{x}})\cong \cF_{P_{J}}^{G}(\mathfrak{v}_{x,J},V_{I_0,J}^{\infty})[\#I_0]\]
in the derived sense.
To summarize, we obtain a decreasing filtration
\begin{equation}\label{equ: loc an St J filtration}
\{\tau^{J}(V_{x,I_0}^{\rm{an}})\}_{I_0\subseteq J\subseteq J_{x}}
\end{equation}
which satisfies
\begin{equation}\label{equ: loc an St J grade}
\mathrm{gr}^{J}(V_{x,I_0}^{\rm{an}})\cong \cF_{P_{J}}^{G}(\mathfrak{v}_{x,J},V_{I_0,J}^{\infty})
\end{equation}
for each $I_0\subseteq J\subseteq J_{x}$.

Let $x,w\in W(G)$ with $x\leq w$ and $I_0\subseteq I_1\subseteq J_{w}$. The embedding $M(w)\hookrightarrow M(x)$ induces by \cite[Thm.~9.4(c)]{Hum08} a (possibly zero) map $M^{I'}(w)\rightarrow M^{I'}(x)$ and thus a (possibly zero) map $i_{x,I'}^{\rm{an}}\rightarrow i_{w,I'}^{\rm{an}}$ which is functorial with respect to the choice of $I_0\subseteq I'\subseteq I_1$. Hence, we obtain a map between complex
\begin{equation}\label{equ: naive Tits x w transfer}
\mathbf{C}^{x}_{I_0,I_1}\rightarrow \mathbf{C}^{w}_{I_0,I_1}.
\end{equation}

Let $I_0\subseteq \Delta$ and $x,w\in\Gamma^{\Delta\setminus I_0}$ with $x\unlhd w$.
For each $I_0\subseteq I'\subseteq J_{x}$ with $I'\not\subseteq J_{w}$, we recall from Lemma~\ref{lem: coxeter parabolic Verma} that $[M^{I'}(x):L(w)]=0$, which together with $\mathrm{cosoc}_{U(\fg)}(M^{I'\cap J_{w}}(w))=L(w)$ forces the composition of
\[M^{I'\cap J_{w}}(w)\rightarrow M^{I'\cap J_{w}}(x)\twoheadrightarrow M^{I'}(x)\]
to be zero. Now that the composition of
\[i_{x,I'}^{\rm{an}}\hookrightarrow i_{x,I'\cap J_{w}}^{\rm{an}}\rightarrow i_{w,I'\cap J_{w}}^{\rm{an}}\]
factors through
\begin{multline*}
i_{x,I'}^{\rm{an}}=\cF_{P_{I'}}^{G}(M^{I'}(x),1_{L_{I'}})\hookrightarrow \cF_{P_{I'}}^{G}(M^{I'}(x),i_{I'\cap J_{w},I'}^{\infty})=\cF_{P_{I'\cap J_{w}}}^{G}(M^{I'}(x),1_{L_{I'\cap J_{w}}})\\
\hookrightarrow \cF_{P_{I'\cap J_{w}}}^{G}(M^{I'\cap J_{w}}(x),1_{L_{I'\cap J_{w}}})=i_{x,I'\cap J_{w}}^{\rm{an}}\rightarrow \cF_{P_{I'\cap J_{w}}}^{G}(M^{I'\cap J_{w}}(w),1_{L_{I'\cap J_{w}}})=i_{w,I'\cap J_{w}}^{\rm{an}},
\end{multline*}
it is also zero. Consequently, the map between complex $\mathbf{C}^{x}_{I_0,J_{w}}\rightarrow \mathbf{C}^{w}_{I_0,J_{w}}$ (see (\ref{equ: naive Tits x w transfer}) with $I_1=J_{w}$ in \emph{loc.cit.}) factors through $\mathbf{C}^{x}_{I_0,J_{x}}$ and gives a map
\begin{equation}\label{equ: C Tits x w transfer}
\mathbf{C}^{x}_{I_0,J_{x}}\rightarrow \mathbf{C}^{w}_{I_0,J_{w}}
\end{equation}
which further induces a map
\begin{equation}\label{equ: St x w transfer}
V_{x,I_0}^{\rm{an}}\rightarrow V_{w,I_0}^{\rm{an}}
\end{equation}
under (\ref{equ: general Tits resolution}). Note that (\ref{equ: St x w transfer}) restricts to a map
\[
\tau^{J}(V_{x,I_0}^{\rm{an}})\rightarrow \tau^{J}(V_{w,I_0}^{\rm{an}})
\]
and thus induces a map
\begin{equation}\label{equ: St x w transfer J}
\mathrm{gr}^{J}(V_{x,I_0}^{\rm{an}})\rightarrow \mathrm{gr}^{J}(V_{w,I_0}^{\rm{an}})
\end{equation}
for each $J\subseteq J_{w}\subseteq J_{x}$. When $J=\emptyset$, the map (\ref{equ: St x w transfer J}) is induced from the non-zero map $\mathfrak{v}_{w,I_0}\rightarrow \mathfrak{v}_{x,I_0}$ (see Lemma~\ref{lem: Lie St wt shift} and its proof) by applying $\cF_{P_{I_0}}^{G}(-,1_{L_{I_0}})$. In particular, we know that the map (\ref{equ: St x w transfer}) is non-zero.

\begin{defn}\label{def: isotypic factor}
Let $J\subseteq \Delta$ and $x\in W(G)$. A simple object $V\in\mathrm{Rep}^{\rm{an}}_{\rm{adm}}(L_J)$ is called \emph{$x$-typic} if there exists a simple object $\pi^{\infty}\in\mathrm{Rep}^{\infty}_{\rm{adm}}(L_{I_x\cap J})$ such that $V\cong \cF_{P_{I_x}\cap L_J}^{L_J}(L^J(x),\pi^{\infty})$.
\end{defn}

For each $x\in W(G)$, $J\subseteq \Delta$ and $I\subseteq I_x\cap J$, we set
\begin{equation}\label{equ: loc an factor St}
C_{x,I,J}\defeq \cF_{P_{I_x}\cap L_{J}}^{L_{J}}(L^{J}(x),V_{I,I_x\cap J}^{\infty}) \in \mathrm{Rep}^{\rm{an}}_{\rm{adm}}(L_{J})
\end{equation}
which is irreducible by \ref{it: OS property 3} of Proposition~\ref{prop: OS property}.
Note that $C_{x,I,J}$ is $x$-typic in the sense of Definition~\ref{def: isotypic factor}
When $J=\Delta$, we write $C_{x,I}\defeq C_{x,I,\Delta}$ for short.

The following result is well-known (see \cite[Thm.~4.6]{OS13}).
\begin{lem}\label{lem: JH mult PS}
Let $I_1\subseteq \Delta$, $x\in\Gamma$ and $I\subseteq I_x$. We have the following results
\begin{enumerate}[label=(\roman*)]
\item \label{it: JH mult PS 1} We have $[i_{I_1}^{\rm{an}}: C_{x,I}]$ if and only if $\mathrm{Supp}(x)\subseteq \Delta\setminus I_1$ and $I_1\subseteq I$, in which case $[i_{I_1}^{\rm{an}}: C_{x,I}]=1$.
\item \label{it: JH mult PS 2} We have $[V_{I_1}^{\rm{an}}: C_{x,I}]$ if and only if $\mathrm{Supp}(x)\subseteq \Delta\setminus I_1$ and $I_1\subseteq I\subseteq I_1\sqcup(\mathrm{Supp}(x)\setminus D_L(x))$, in which case $[V_{I_1}^{\rm{an}}: C_{x,I}]=1$.
\end{enumerate}
\end{lem}
\begin{proof}
Recall from \ref{it: OS property 4} of Proposition~\ref{prop: OS property} that we have $i_{I_1}^{\rm{an}}=i_{I_1}^{\rm{an}}(1_{L_{I_1}})\cong \cF_{P_{I_1}}^{G}(M^{I_1}(1),1_{L_{I_1}})$.
It follows from Lemma~\ref{lem: JH OS} that $[i_{I_1}^{\rm{an}}: C_{x,I}]\neq 0$ if and only if $[M^{I_1}(1):L(x)]\neq 0$ (with $I_1\subseteq I_x$) and $[i_{I_1,I_x}^{\infty}(1_{L_{I_1}}): V_{I,I_x}^{\infty}]\neq 0$, in which case we have
\begin{equation}\label{equ: JH mult PS Verma}
[i_{I_1}^{\rm{an}}: C_{x,I}]=[M^{I_1}(1):L(x)][i_{I_1,I_x}^{\infty}(1_{L_{I_1}}): V_{I,I_x}^{\infty}].
\end{equation}
Recall from Lemma~\ref{it: g coxeter 1} of Lemma~\ref{lem: coxeter subquotient} that $[M^{I_1}(1):L(x)]\neq 0$ if and only if $\mathrm{Supp}(x)\subseteq \Delta\setminus I_1$, in which case $[M^{I_1}(1):L(x)]=1$.
We also recall from Lemma~\ref{lem: sm cube} that $i_{I_1,I_x}^{\infty}(1_{L_{I_1}})\cong Q_{I_x}(I_x,I_1)$ and that $[Q_{I_x}(I_x,I_1): V_{I,I_x}^{\infty}]\neq 0$ if and only if $I_1\subseteq I\subseteq I_x$, in which case $[Q_{I_x}(I_x,I_1): V_{I,I_x}^{\infty}]=1$. These together with (\ref{equ: JH mult PS Verma}) imply that $[i_{I_1}^{\rm{an}}: C_{x,I}]\neq 0$ if and only if $\mathrm{Supp}(x)\subseteq \Delta\setminus I_1$ and $I_1\subseteq I\subseteq I_x$, or equivalently $I_1\subseteq I\setminus\mathrm{Supp}(x)$, in which case we have $[i_{I_1}^{\rm{an}}: C_{x,I}]=1$.
It follows from (\ref{equ: general Tits resolution}) (upon replacing $x$, $I$, $I_0$ and $I_1$ with $1$, $\Delta$, $I_1$ and $\Delta$ here) that
\[[V_{I_1}^{\rm{an}}: C_{x,I}]=\sum_{I_1\subseteq I_1'}(-1)^{\#I_1'\setminus I_1}[i_{I_1'}^{\rm{an}}: C_{x,I}]=\sum_{I_1\subseteq I_1'\subseteq I\setminus\mathrm{Supp}(x)}(-1)^{\#I_1'\setminus I_1}\]
is non-zero if and only if $I_1=I\setminus\mathrm{Supp}(x)$, or equivalently $\mathrm{Supp}(x)\subseteq \Delta\setminus I_1$ and $I_1\subseteq I\subseteq I_1\sqcup(\mathrm{Supp}(x)\setminus D_L(x))$. The proof is thus finished.
\end{proof}

The following result is a special case of \cite[Prop.~4.5.13 (iv)]{BQ24}.
\begin{prop}\label{prop: St key seq}
We consider $V_i=\cF_{P_{I_i}}^G(M_i,\pi_i^{\infty})$ such that $I_i\subseteq \Delta$, $M_i\in\cO^{\fp_{I_i}}_{\rm{alg}}$ and $\pi_i^{\infty}\in\cB^{I_i}$ for $i=0,1$. Assume moreover that $M_1\cong U(\fg)\otimes_{U(\fp_{I_1})}X$ for some $\ft$-semi-simple finite dimensional $U(\fp_{I_1})$-module $X$.
Then we have a spectral sequence
\begin{equation}\label{equ: St key seq}
\mathrm{Ext}_{U(\fg)}^{\ell}(M_1, M_0)\otimes_E \mathrm{Ext}_{L_{I_1}}^k(i_{I_0\cap I_1,I_1}^{\infty}(J_{I_0,I_0\cap I_1}(\pi_0^{\infty})), \pi_1^{\infty})^{\infty}
\implies \mathrm{Ext}_{G}^{k+\ell}(V_0,V_1).
\end{equation}
In particular, we have an isomorphism between $E$-vector spaces
\begin{equation}\label{equ: St key seq Hom}
\Hom_{U(\fg)}(M_1, M_0)\otimes_E \Hom_{L_{I_1}}(i_{I_0\cap I_1,I_1}^{\infty}(J_{I_0,I_0\cap I_1}(\pi_0^{\infty})), \pi_1^{\infty})\cong\Hom_{G}(V_0,V_1).
\end{equation}
\end{prop}
\begin{proof}
Our definition of $\cB^{I_i}$ ensures that $\Sigma_i\defeq W(L_{I_i})\cdot\cJ(\pi_i^{\infty})= W(L_{I_i})\cdot 1_T$. Since we also have $1_T\in \Sigma_0\cap \Sigma_1\neq \emptyset$ and $\pi_i^{\infty}\in \cB^{I_i}=\cB^{I_i}_{\Sigma_i}$ (with $\Sigma_i=W(L_{I_i})\cdot 1_T$ being \emph{$G$-regular}, see the paragraph below \cite[(37)]{BQ24}) for $i=0,1$, we can apply \cite[Prop.~4.5.13 (iv)]{BQ24} with $w=1$ and deduce (\ref{equ: St key seq}). Here we use the identifications of functors $i_{I_0,I_1,1}^{\infty}(-)=i_{I_0\cap I_1,I_1}^{\infty}(-)$ and $J_{I_0,I_1,1}(-)=J_{I_0,I_0\cap I_1}(-)$ from \cite[(43),(45)]{BQ24}. Strictly speaking, our assumption $M_1\cong U(\fg)\otimes_{U(\fp_{I_1})}X$ here is weaker than that of \cite[Prop.~4.5.13 (iv)]{BQ24} (which says $M_1\cong U(\fg)\otimes_{U(\fp_{I_1})}L^{I_1}(\mu)$ for some $\mu\in\Lambda_{I_1}^+$, see the beginning of \cite[\S 4.5]{BQ24}), but the reader could check that all arguments in \cite[\S 4.5]{BQ24} work for such more general $M_1$.
\end{proof}

\begin{rem}\label{rem: trivial block adjunction}
We continue to use notation in the proof of Proposition~\ref{prop: St key seq}. Since $1_{T}\in\Sigma_0\cap\Sigma_1$, we know that $\Sigma_1\cap w^{-1}\cdot\Sigma_0\neq 0$ for some $w\in W^{I_0,I_1}$ if and only if $w=1$, and thus \[J_{I_1}(i_{I_0}^{\infty}(\pi_0^{\infty}))_{\cB^{I_1}_{\Sigma_1}}\cong i_{I_0,I_1,1}^{\infty}(J_{I_0,I_1,1}(\pi_0^{\infty}))_{\cB^{I_1}_{\Sigma_1}}\cong i_{I_0\cap I_1,I_1}^{\infty}(J_{I_0,I_0\cap I_1}(\pi_0^{\infty}))_{\cB^{I_1}_{\Sigma_1}}\]
by \cite[Lem~2.1.18]{BQ24}, which together with adjunction (cf.~\cite[(31)]{BQ24}) (and property of the Bernstein block $\cB^{I_1}_{\Sigma_1}$) gives isomorphisms
\begin{equation}\label{equ: trivial block adjunction}
\mathrm{Ext}_{G}^k(i_{I_0}^{\infty}(\pi_0^{\infty}), i_{I_1}^{\infty}(\pi_1^{\infty}))^{\infty}\cong \mathrm{Ext}_{L_{I_1}}^k(J_{I_1}(i_{I_0}^{\infty}(\pi_0^{\infty}))_{\cB^{I_1}_{\Sigma_1}},\pi_1^{\infty})^{\infty}\cong \mathrm{Ext}_{L_{I_1}}^k(i_{I_0\cap I_1,I_1}^{\infty}(J_{I_0,I_0\cap I_1}(\pi_0^{\infty})), \pi_1^{\infty})^{\infty}
\end{equation}
\end{rem}

\begin{lem}\label{lem: sm distance Ext vanishing}
We consider $V_i=\cF_{P_{I_i}}^G(M_i,\pi_i^{\infty})$ such that $I_i\subseteq \Delta$, $M_i\in\cO^{\fp_{I_i}}_{0}$ and $\pi_i^{\infty}\in\cB^{I_i}$ for $i=0,1$. Assume that
\[d_{\Delta}(\pi_0^{\infty},\pi_1^{\infty})=d_{\Delta}(i_{I_0}^{\infty}(\pi_0^{\infty}), i_{I_1}^{\infty}(\pi_1^{\infty}))=\infty.\]
Then we have $\mathrm{Ext}_{G}^k(V_0,V_1)=0$ for each $k\geq 0$.
\end{lem}
\begin{proof}
This follows directly from \cite[Lem.~4.5.17 (i)]{BQ24}.
\end{proof}

\begin{lem}\label{lem: Ext between sm}
Let $I_0,I_1\subseteq \Delta$. We have a natural isomorphism
\begin{equation}\label{equ: Ext between sm}
\mathrm{Ext}_{G}^k(V_{I_0}^{\infty},V_{I_1}^{\infty})^{\infty}\buildrel\sim\over\longrightarrow \mathrm{Ext}_{G}^k(V_{I_0}^{\infty},V_{I_1}^{\infty})
\end{equation}
for each $k\leq d(I_0,I_1)+2$. In particular, we have $\mathrm{Ext}_{G}^k(V_{I_0}^{\infty},V_{I_1}^{\infty})=0$ when either $k<d(I_0,I_1)$ or $k-d(I_0,I_1)\in\{1,2\}$, and $\Dim_E \mathrm{Ext}_{G}^{d(I_0,I_1)}(V_{I_0}^{\infty},V_{I_1}^{\infty})=1$.
\end{lem}
\begin{proof}
The isomorphisms (\ref{equ: Ext between sm}) for $k\leq d(I_0,I_1)+2$ follow from \cite[(369)]{BQ24} together with the fact that $\mathrm{Ext}_{U(\fg)}^k(L(1),L(1))=0$ for $k=1,2$ as $\fg$ is semi-simple. The rest of the claims then follow from (\ref{equ: Ext between sm}) combined with \ref{it: general sm distance 2} and \ref{it: general sm distance 4} of Lemma~\ref{lem: general sm distance}.
\end{proof}

\begin{lem}\label{lem: OS layer}
Let $x\in W(G)$ and $\pi^{\infty}\in\cB^{I_x}$ be multiplicity free. Then $\cF_{P_{I_x}}^{G}(L(x),-)$ induces an order-preserving bijection
\begin{equation}\label{equ: OS layer}
\mathrm{JH}_{L_{I_x}}(\pi^{\infty})\buildrel\sim\over\longrightarrow \mathrm{JH}_{G}(\cF_{P_{I_x}}^{G}(L(x),\pi^{\infty}))
\end{equation}
between partially-ordered sets.
\end{lem}
\begin{proof}
This follows from \cite[Lem.~5.1.2]{BQ24} by taking $I=I_x$ and $\Sigma=W(L_{I_x})\cdot 1_T$ in \emph{loc.cit.}.
\end{proof}

\begin{lem}\label{lem: OS socle cosocle}
Let $I\subseteq \Delta$, $M\in\cO^{\fp_I}_{\rm{alg}}$ and $\pi^{\infty}\in\cB^I$. Then we have
\begin{equation}\label{equ: OS socle}
\mathrm{soc}_{G}(\cF_{P_I}^{G}(M,\pi^{\infty})=\mathrm{soc}_{G}(\cF_{P_I}^{G}(\mathrm{cosoc}_{U(\fg)}(M),\pi^{\infty})
\end{equation}
and
\begin{equation}\label{equ: OS cosocle}
\mathrm{cosoc}_{G}(\cF_{P_I}^{G}(M,\pi^{\infty})=\mathrm{cosoc}_{G}(\cF_{P_I}^{G}(\mathrm{soc}_{U(\fg)}(M),\pi^{\infty}).
\end{equation}
\end{lem}
\begin{proof}
The equality (\ref{equ: OS socle}) can be checked by applying \cite[Lem.~5.1.7 (i)]{BQ24} with $I_1=I$, $M_1=M$ and an arbitrary irreducible $V_0=\cF_{P_{I_0}}^{G}(L(x_0),\pi_0^{\infty})$.
Similarly, The equality (\ref{equ: OS cosocle}) can be checked by applying \cite[Lem.~5.1.7 (ii)]{BQ24} with $I_1=I$, $M_1=M$ and an arbitrary irreducible $V_0=\cF_{P_{I_0}}^{G}(L(x_0),\pi_0^{\infty})$.
\end{proof}

\begin{lem}\label{lem: Ext1 OS cube}
Let $x_i\in W(G)$ and $I_i,I_i'\subseteq I_{x_i}$ with $V_i=\cF_{P_{I_{x_i}}}^{G}(L(x_i),Q_{I_{x_i}}(I_i,I_i'))$ for $i=0,1$ with $x_0\neq x_1$.
Then we have a canonical isomorphism
\begin{equation}\label{equ: Ext1 OS cube}
\mathrm{Ext}_{G}^1(V_0, V_1)\cong \mathrm{Ext}_{U(\fg)}^1(L(x_1),L(x_0))\otimes_E \Hom_{G}(Q_{\Delta}(I_0\sqcup D_L(x_0),I_0'), Q_{\Delta}(I_1\sqcup D_L(x_1),I_1')).
\end{equation}
\end{lem}
\begin{proof}
We use the shortened notation $J_i\defeq I_{x_i}$ and note that $\Delta\setminus J_i=D_L(x_i)$ for $i=0,1$. It follows from \ref{it: sm cube 3} of Lemma~\ref{lem: sm cube} (by taking $I=J_i$, $I'=\Delta$, $I_0=I_i$ and $I_1=I_i'$) that
\begin{equation}\label{equ: Ext1 OS induction cube}
i_{J_i,\Delta}^{\infty}(Q_{J_i}(I_i,I_i'))\cong Q_{\Delta}(I_i\sqcup(\Delta\setminus J_i),I_i')=Q_{\Delta}(I_i\sqcup D_L(x_i),I_i')
\end{equation}
for $i=0,1$.
By taking $k=0$ and replacing $I_i$ and $\pi_i^{\infty}$ in (\ref{equ: trivial block adjunction}) with $J_i$ and $Q_{J_i}(I_i,I_i')$ and then applying (\ref{equ: Ext1 OS induction cube}), we obtain a canonical isomorphism
\begin{equation}\label{equ: Ext1 OS sm Hom}
\Hom_{G}(Q_{\Delta}(I_0\sqcup D_L(x_0),I_0'), Q_{\Delta}(I_1\sqcup D_L(x_1),I_1'))\cong \Hom_{L_{J_1}}(i_{J_0\cap J_1,J_1}^{\infty}(J_{J_0,J_0\cap J_1}(Q_{J_0}(I_0,I_0'))),Q_{J_1}(I_1,I_1')).
\end{equation}
It follows from \cite[Prop.~5.1.14]{BQ24} (by replacing $I_i$ and $\pi_i^{\infty}$ in \emph{loc.cit.} with $J_i$ and $Q_{J_i}(I_i,I_i')$ here) that we have a canonical isomorphism
\[\mathrm{Ext}_{G}^1(V_0, V_1)\cong \mathrm{Ext}_{U(\fg)}^1(L(x_1),L(x_0))\otimes_E \Hom_{L_{J_1}}(i_{J_0\cap J_1,J_1}^{\infty}(J_{J_0,J_0\cap J_1}(Q_{J_0}(I_0,I_0'))),Q_{J_1}(I_1,I_1')),\]
which together with (\ref{equ: Ext1 OS sm Hom}) gives (\ref{equ: Ext1 OS cube}).
\end{proof}

\begin{rem}\label{rem: length 2 g OS}
Let $x_i\in W(G)$ for $i=0,1$ that satisfy either $x_0<x_1$ with $\ell(x_1)=\ell(x_0)+1$, or $x_0>x_1$ with $\ell(x_0)=\ell(x_1)+1$. Then there exists a unique length $2$ $U(\fg)$-module $M$ with socle $L(x_0)$ and cosocle $L(x_1)$. Let $I'\subseteq I\subseteq I_{x_0}\cap I_{x_1}$. By \ref{it: OS property 1} of Proposition~\ref{prop: OS property} we know that $\cF_{P_I}^{G}(M,V_{I',I}^{\infty})$ fits into the short exact sequence
\begin{equation}\label{equ: length 2 g OS}
0\rightarrow \cF_{P_I}^{G}(L(x_1),V_{I',I}^{\infty})\rightarrow \cF_{P_I}^{G}(M,V_{I',I}^{\infty})\rightarrow \cF_{P_I}^{G}(L(x_0),V_{I',I}^{\infty})\rightarrow 0.
\end{equation}
It follows from Lemma~\ref{lem: OS socle cosocle} that $\mathrm{soc}_{G}(\cF_{P_I}^{G}(M,V_{I',I}^{\infty}))=\mathrm{soc}_{G}(\cF_{P_I}^{G}(L(x_1),V_{I',I}^{\infty}))$ and thus (\ref{equ: length 2 g OS}) is non-split. Therefore by Lemma~\ref{lem: Ext1 OS cube} we know that $\cF_{P_I}^{G}(M,V_{I',I}^{\infty})\in\mathrm{Rep}^{\rm{an}}_{\rm{adm}}(G)$ is the unique object that fits into a non-split extension of the form
\[0\rightarrow \cF_{P_I}^{G}(L(x_1),V_{I',I}^{\infty})\rightarrow \ast \rightarrow \cF_{P_I}^{G}(L(x_0),V_{I',I}^{\infty})\rightarrow 0.\]
\end{rem}

\begin{lem}\label{lem: Ext1 OS simple}
Let $x_i\in W(G)$ and $I_i\subseteq I_{x_i}$ for $i=0,1$ with $x_0\neq x_1$. Then we have a canonical isomorphism
\begin{equation}\label{equ: Ext1 OS simple}
\mathrm{Ext}_{G}^1(C_{x_0,I_0}, C_{x_1,I_1})\cong \mathrm{Ext}_{U(\fg)}^1(L(x_1),L(x_0))\otimes_E \Hom_{G}(Q_{\Delta}(I_0\sqcup D_L(x_0),I_0), Q_{\Delta}(I_1\sqcup D_L(x_1),I_1)).
\end{equation}
In particular, (\ref{equ: Ext1 OS simple}) is non-zero if and only if $\mathrm{Ext}_{U(\fg)}^1(L(x_1),L(x_0))\neq 0$ and $I_0=I_1\sqcup(D_L(x_1)\setminus D_L(x_0))$, in which case it has dimension $\Dim_E \mathrm{Ext}_{U(\fg)}^1(L(x_1),L(x_0))$.
\end{lem}
\begin{proof}
The isomorphism (\ref{equ: Ext1 OS simple}) follows directly from (\ref{equ: Ext1 OS cube}) by taking $I_0'=I_0$ and $I_1'=I_1$. It suffices to show that
\begin{equation}\label{equ: Ext1 OS simple sm}
\Hom_{G}(Q_{\Delta}(I_0\sqcup D_L(x_0),I_0), Q_{\Delta}(I_1\sqcup D_L(x_1),I_1))\neq 0
\end{equation}
if and only if $I_0=I_1\sqcup(D_L(x_1)\setminus D_L(x_0))$, in which case (\ref{equ: Ext1 OS simple sm}) has dimension one. It follows from \ref{it: sm cube 2} of Lemma~\ref{lem: sm cube} that (\ref{equ: Ext1 OS simple sm}) holds if and only if $I_1\sqcup D_L(x_1)\in [I_0\sqcup D_L(x_0),I_0]$ and $I_0\in [I_1\sqcup D_L(x_1),I_1]$, if and only if $I_1\subseteq I_0\subseteq I_1\sqcup D_L(x_1)\subseteq I_0\sqcup D_L(x_0)$ if and only if $I_0=I_1\sqcup(D_L(x_1)\setminus D_L(x_0))$, in which case (\ref{equ: Ext1 OS simple sm}) is one dimensional. The proof is thus finished.
\end{proof}

\begin{lem}\label{lem: Ext with sm general}
Let $I\subseteq \Delta$, $\pi_0^{\infty}\in\cB^{\Delta}$ and $\pi_1^{\infty}\in \cB^I$.
\begin{enumerate}[label=(\roman*)]
\item \label{it: Ext with sm general 1} For each $w\in W(G)$ with $L(w)\in\cO^{\fp_I}_{\rm{alg}}$, the $E$-vector space $\mathrm{Ext}_{G}^k(\pi_0^{\infty}, \cF_{P_I}^{G}(L(w),\pi_1^{\infty}))$ vanishes for $k<d_{\Delta}(\pi_0^{\infty}, \pi_1^{\infty})+\ell(w)$ and is isomorphic to $\mathrm{Ext}_{G}^{d_{\Delta}(\pi_0^{\infty},\pi_1^{\infty})}(\pi_0^{\infty},i_{I,\Delta}^{\infty}(\pi_1^{\infty}))^{\infty}$ for $k=d_{\Delta}(\pi_0^{\infty}, \pi_1^{\infty})+\ell(w)$.
\item \label{it: Ext with sm general 2} Let $x,w\in W(G)$ with $x<w$, $\ell(w)=\ell(x)+1$ such that the unique length $2$ $U(\fg)$-module $M$ with socle $L(x)$ and cosocle $L(w)$ (see Lemma~\ref{lem: g Ext1}) satisfies $M\in\cO^{\fp_I}_{\rm{alg}}$. Then $\mathrm{Ext}_{G}^k(\pi_0^{\infty}, \cF_{P_I}^{G}(M,\pi_1^{\infty}))$ vanishes for $k\leq d_{\Delta}(\pi_0^{\infty}, \pi_1^{\infty})+\ell(x)$.
\end{enumerate}
\end{lem}
\begin{proof}
By the definition of $d_{\Delta}(\pi_0^{\infty}, \pi_1^{\infty})$ and \ref{it: general sm distance 4} of Lemma~\ref{lem: general sm distance} we know that $\mathrm{Ext}_{G}^k(\pi_0^{\infty},i_{I_w,\Delta}^{\infty}(\pi_1^{\infty}))^{\infty}$ vanishes for $k<d_{\Delta}(\pi_0^{\infty}, \pi_1^{\infty})$ and is $1$-dimensional for $k=d_{\Delta}(\pi_0^{\infty}, \pi_1^{\infty})$. For each finite dimensional $U(\fp_I)$-module $X$, we have the following spectral sequence
\begin{equation}\label{equ: Ext with sm general seq}
\mathrm{Ext}_{U(\fg)}^k(U(\fg)\otimes_{U(\fp_I)}X,L(1))\otimes_E \mathrm{Ext}_{G}^{\ell}(\pi_0^{\infty}, i_{I,\Delta}^{\infty}(\pi_1^{\infty}))^{\infty}\Rightarrow \mathrm{Ext}_{G}^{k+\ell}(\pi_0^{\infty}, \cF_{P_I}^{G}(U(\fg)\otimes_{U(\fp_I)}X,\pi_1^{\infty})).
\end{equation}
by Proposition~\ref{prop: St key seq} and Remark~\ref{rem: trivial block adjunction}.
If we take $X=L^I(w)$ for some $w\in W^{I,\emptyset}$ with $U(\fg)\otimes_{U(\fp_I)}X=M^I(w)$, then by \ref{it: n coh 5} of Lemma~\ref{lem: n coh collection} together with (\ref{equ: g spectral seq}) we know that $\mathrm{Ext}_{U(\fg)}^k(M^I(w),L(1))$ vanishes for $k<\ell(w)$ and is $1$-dimensional for $k=\ell(w)$, which together with (\ref{equ: Ext with sm general seq}) and the previous discussion on $\mathrm{Ext}_{G}^{\bullet}(\pi_0^{\infty},i_{I_w,\Delta}^{\infty}(\pi_1^{\infty}))^{\infty}$ implies that $\mathrm{Ext}_{G}^k(\pi_0^{\infty}, \cF_{P_I}^{G}(M^I(w),\pi_1^{\infty}))$ vanishes for $k<d_{\Delta}(\pi_0^{\infty}, \pi_1^{\infty})+\ell(w)$ and is $1$-dimensional for $k=d_{\Delta}(\pi_0^{\infty}, \pi_1^{\infty})+\ell(w)$.
Similarly, if we take $X$ as in Lemma~\ref{lem: Ext special Verma}, we know that $\mathrm{Ext}_{U(\fg)}^k(U(\fg)\otimes_{U(\fp_I)}X,L(1))$ vanishes for $k\leq \ell(x)$, which together with (\ref{equ: Ext with sm general seq}) and the previous discussion on $\mathrm{Ext}_{G}^{\bullet}(\pi_0^{\infty},i_{I_w,\Delta}^{\infty}(\pi_1^{\infty}))^{\infty}$ implies that $\mathrm{Ext}_{G}^k(\pi_0^{\infty}, \cF_{P_I}^{G}(U(\fg)\otimes_{U(\fp_I)}X,\pi_1^{\infty}))$ vanishes for $k\leq d_{\Delta}(\pi_0^{\infty}, \pi_1^{\infty})+\ell(x)$.

We prove \ref{it: Ext with sm general 1} by a decreasing induction on $\ell(w)$.\\
If $w$ is the maximal length element in $W^{I,\emptyset}$, then we have $L(w)=M^I(w)=U(\fg)\otimes_{U(\fp_I)}L^I(w)$ and conclude \ref{it: Ext with sm general 1} from the discussion on $\mathrm{Ext}_{G}^{\bullet}(\pi_0^{\infty}, \cF_{P_I}^{G}(M^I(w),\pi_1^{\infty}))$ in previous paragraph. Otherwise, we have a surjection $q_1: M^I(w)\twoheadrightarrow L(w)$ with each $L(y)\in\mathrm{JH}_{U(\fg)}(\mathrm{ker}(q_1))$ satisfying $\ell(y)>\ell(w)$ and thus $\mathrm{Ext}_{G}^k(\pi_0^{\infty}, \cF_{P_I}^{G}(L(y),\pi_1^{\infty}))=0$ for each $k\leq d_{\Delta}(\pi_0^{\infty}, \pi_1^{\infty})+\ell(w)<d_{\Delta}(\pi_0^{\infty}, \pi_1^{\infty})+\ell(y)$ by our induction hypothesis, which together with a d\'evissage with respect to $\mathrm{JH}_{U(\fg)}(\mathrm{ker}(q_1))$ gives $\mathrm{Ext}_{G}^k(\pi_0^{\infty}, \cF_{P_I}^{G}(\mathrm{ker}(q_1),\pi_1^{\infty}))=0$ for $k\leq d_{\Delta}(\pi_0^{\infty}, \pi_1^{\infty})+\ell(w)$. A further d\'evissage with respect to
\[0\rightarrow \cF_{P_I}^{G}(L(w),\pi_1^{\infty})\rightarrow \cF_{P_I}^{G}(M^I(w),\pi_1^{\infty})\rightarrow \cF_{P_I}^{G}(\mathrm{ker}(q_1),\pi_1^{\infty})\rightarrow 0\]
gives an isomorphism
\[\mathrm{Ext}_{G}^k(\pi_0^{\infty}, \cF_{P_I}^{G}(L(w),\pi_1^{\infty}))\buildrel\sim\over\longrightarrow\mathrm{Ext}_{G}^k(\pi_0^{\infty}, \cF_{P_I}^{G}(M^I(w),\pi_1^{\infty}))\]
for each $k\leq d_{\Delta}(\pi_0^{\infty}, \pi_1^{\infty})+\ell(w)$, which together with the discussion on $\mathrm{Ext}_{G}^{\bullet}(\pi_0^{\infty}, \cF_{P_I}^{G}(M^I(w),\pi_1^{\infty}))$ in previous paragraph gives \ref{it: Ext with sm general 1}.

We prove \ref{it: Ext with sm general 2}.\\
Let $\xi: Z(\fg)\rightarrow E$ be the unique infinitesimal character such that $L(x)_{\xi}\neq 0$.
Let $X$ as in Lemma~\ref{lem: Ext special Verma} equipped with a surjection $q: U(\fg)\otimes_{U(\fp_I)}X\twoheadrightarrow M$ with each $L(\mu)\in\mathrm{ker}(q)$ satisfying either $L(\mu)_{\xi}=0$ or $\mu=y\cdot 0$ for some $y>x$. In the second case, we deduce
\begin{equation}\label{equ: Ext with sm general 2 vanishing}
\mathrm{Ext}_{G}^k(\pi_0^{\infty}, \cF_{P_I}^{G}(L(\mu),\pi_1^{\infty}))=0
\end{equation}
for $k\leq d_{\Delta}(\pi_0^{\infty}, \pi_1^{\infty})+\ell(x)<d_{\Delta}(\pi_0^{\infty}, \pi_1^{\infty})+\ell(y)$ by \ref{it: Ext with sm general 1}. In the first case, the vanishing (\ref{equ: Ext with sm general 2 vanishing}) holds for $k\geq 0$ by \cite[Lem.~4.5.18]{BQ24}.
By a d\'evissage with respect to $\mathrm{JH}_{U(\fg)}(\mathrm{ker}(q))$ we have
\[\mathrm{Ext}_{G}^k(\pi_0^{\infty}, \cF_{P_I}^{G}(\mathrm{ker}(q),\pi_1^{\infty}))=0\]
for $k\leq d_{\Delta}(\pi_0^{\infty}, \pi_1^{\infty})+\ell(x)$, which together with a further d\'evissage with respect to
\[0\rightarrow \cF_{P_I}^{G}(M,\pi_1^{\infty})\rightarrow \cF_{P_I}^{G}(U(\fg)\otimes_{U(\fp_I)}M_I,\pi_1^{\infty})\rightarrow \cF_{P_I}^{G}(\mathrm{ker}(q),\pi_1^{\infty})\rightarrow 0\]
as well as the vanishing of $\mathrm{Ext}_{G}^k(\pi_0^{\infty}, \cF_{P_I}^{G}(U(\fg)\otimes_{U(\fp_I)}X,\pi_1^{\infty}))$ for $k\leq d_{\Delta}(\pi_0^{\infty}, \pi_1^{\infty})+\ell(x)$ in the first paragraph finishes the proof of \ref{it: Ext with sm general 2}.
\end{proof}

Recall the distance $d(-,-)$ between two subsets of $\Delta$ from the discussion before Lemma~\ref{lem: sm cube}.
\begin{lem}\label{lem: Ext with sm}
Let $I_0\subseteq \Delta$, $w\in W(G)$ and $I\subseteq I_w$.
\begin{enumerate}[label=(\roman*)]
\item \label{it: Ext with sm 1} If $I_0\not\supseteq D_L(w)$, then we have $\mathrm{Ext}_{G}^k(V_{I_0}^{\infty}, C_{w,I})=0$ for $k\geq 0$.
\item \label{it: Ext with sm 2} If $I_0\supseteq D_L(w)$, then $\mathrm{Ext}_{G}^k(V_{I_0}^{\infty}, C_{w,I})$ vanishes for $k<d(I_0\setminus D_L(w),I)+\ell(w)$, and is $1$-dimensional for $k=d(I_0\setminus D_L(w),I)+\ell(w)$.
\item \label{it: Ext with sm 3} Let $x\in W(G)$ with $x<w$, $\ell(w)=\ell(x)+1$ and $D_L(x)=D_L(w)$ (with $I_x=I_w$). Let $M$ be the unique length $2$ $U(\fg)$-module with socle $L(x)$ and cosocle $L(w)$. Then we have
    \begin{equation}\label{equ: Ext with sm OS vanishing}
    \mathrm{Ext}_{G}^k(V_{I_0}^{\infty}, \cF_{P_{I_w}}^{G}(M,V_{I,I_w}^{\infty}))=0
    \end{equation}
    for $k\leq d(I_0\setminus D_L(x),I)+\ell(x)$.
\end{enumerate}
\end{lem}
\begin{proof}
We fix a choice of $w\in W(G)$ through the proof.
Recall from \ref{it: sm cube 3} of Lemma~\ref{lem: sm cube} that $i_{I_w,\Delta}^{\infty}(V_{I_w,I}^{\infty})\cong Q_{\Delta}(I\sqcup D_L(w),I)$. By \ref{it: general sm distance 3} of Lemma~\ref{lem: general sm distance}, we observe that
\begin{equation}\label{equ: Ext with sm distance}
d_{\Delta}(V_{I_0}^{\infty},Q_{\Delta}(I\sqcup D_L(w),I))=d_{\Delta}(Q_{\Delta}(I,I\sqcup D_L(w)),V_{I_0}^{\infty})<\infty
\end{equation}
if and only if $d(I_0,[I,I\sqcup D_L(w)])=d(I_0,I\sqcup D_L(w))$ if and only if $D_L(w)\subseteq I_0$, in which case (\ref{equ: Ext with sm distance}) equals $d(I_0\setminus D_L(w),I)$ and 
\[\Dim_E \mathrm{Ext}_{G}^{d(I_0\setminus D_L(w),I)}(V_{I_0}^{\infty},Q_{\Delta}(I\sqcup D_L(w),I))^{\infty}=1\] 
(see \ref{it: general sm distance 4} of Lemma~\ref{lem: general sm distance}).
These together with \ref{it: Ext with sm general 1} of Lemma~\ref{lem: Ext with sm general} (resp.~\ref{it: Ext with sm general 2} of Lemma~\ref{lem: Ext with sm general}) give \ref{it: Ext with sm 1} and \ref{it: Ext with sm 2} (resp.~gives \ref{it: Ext with sm 3}).
The proof is thus finished.
\end{proof}

\begin{lem}\label{lem: x Ext with sm}
Let $I_1\subseteq I_0\subseteq\Delta$ and $x\in\Gamma^{I_0\setminus I_1}$. Let $w\in W(G)$ and $I\subseteq I_w$ such that $[V_{x,I_1}^{\rm{an}}:C_{w,I}]\neq 0$.
Then we have $\mathrm{Ext}_{G}^k(V_{I_0}^{\infty}, C_{w,I})=0$ for $k<\#I_0\setminus I_1$, and
\begin{equation}\label{equ: x Ext with sm minimal deg}
\mathrm{Ext}_{G}^{\#I_0\setminus I_1}(V_{I_0}^{\infty}, C_{w,I})\neq 0
\end{equation}
if and only if $x\unlhd w$, $w\in\Gamma^{I_0\setminus I_1}$ and $I=I_1\sqcup(\mathrm{Supp}(w)\setminus D_L(w))$, in which case (\ref{equ: x Ext with sm minimal deg}) is $1$-dimensional.
\end{lem}
\begin{proof}
If $I_0\not\supseteq D_L(w)$, then 
\begin{equation}\label{equ: x Ext with sm}
\mathrm{Ext}_{G}^k(V_{I_0}^{\infty}, C_{w,I})
\end{equation}
is zero for $k\geq 0$ by \ref{it: Ext with sm 1} of Lemma~\ref{lem: Ext with sm}. It is thus harmless to assume in the rest of the proof that $D_L(w)\subseteq I_0$.
We write $k_0\defeq d(I_0\setminus D_L(w),I)+\ell(w)$ for short.
It follows from \ref{it: Ext with sm 2} of Lemma~\ref{lem: Ext with sm} that (\ref{equ: x Ext with sm}) is zero when $k<k_0$ and $1$ dimensional when $k=k_0$. We devote the rest of the proof to show that $k_0\geq \#I_0\setminus I_1$, with the equality holds if and only if $x\unlhd w$, $w\in\Gamma^{I_0\setminus I_1}$ and $I=I_1\sqcup(\mathrm{Supp}(w)\setminus D_L(w))$.

Recall from (\ref{equ: loc an St J filtration}) and (\ref{equ: loc an St J grade}) (with $I_1$ replacing $I_0$ in \emph{loc.cit.}) that $V_{x,I_1}^{\rm{an}}$ admits a decreasing filtration 
\[\{\tau^{J}(V_{x,I_1}^{\rm{an}})\}_{I_1\subseteq J\subseteq J_{x}}\]
with graded pieces
\[\mathrm{gr}^{J}(V_{x,I_1}^{\rm{an}})\cong \cF_{P_{J}}^{G}(\mathfrak{v}_{x,J},V_{I_1,J}^{\infty})\]
for each $I_1\subseteq J\subseteq J_{x}$. Now that $[V_{x,I_1}^{\rm{an}}:C_{w,I}]\neq 0$, there exists $I_1\subseteq J\subseteq J_{x}$ such that
\[[\cF_{P_{J}}^{G}(\mathfrak{v}_{x,J},V_{I_1,J}^{\infty}):C_{w,I}]\neq 0,\]
or equivalently 
\begin{equation}\label{equ: x J St w}
[\mathfrak{v}_{x,J}:L(w)]\neq 0
\end{equation}
and 
\begin{equation}\label{equ: x J St sm}
[i_{J,I_w}^{\infty}(V_{I_1,J}^{\infty}):V_{I,I_w}^{\infty}]\neq 0.
\end{equation}
We deduce from (\ref{equ: x J St w}) and \ref{it: Lie St JH 1} of Lemma~\ref{lem: Lie St JH} that $\mathrm{Supp}(w)\supseteq\Delta\setminus J$ and thus $\ell(w)\geq \#\Delta\setminus J$, with the equality holds if and only if $w\in\Gamma_{\Delta\setminus J}$ and $x\unlhd w$.
We deduce from (\ref{equ: x J St sm}) and \ref{it: sm cube 3} of Lemma~\ref{lem: sm cube} that $I_1\subseteq I\subseteq I_1\sqcup(I_w\setminus J)$, and thus
\begin{multline*}
d(I_0\setminus D_L(w),I)=d(I_0,I\sqcup D_L(w))\geq \#I_0\setminus I_1-d(I_1,I\sqcup D_L(w))=\#I_0\setminus I_1-\#(I\sqcup D_L(w))+\#I_1\\
\geq \#I_0\setminus I_1-\#(I_1\sqcup(I_w\setminus J)\sqcup D_L(w))+\#I_1=\#I_0\setminus I_1-\#\Delta\setminus J,
\end{multline*}
with equality holds if and only if $I\sqcup D_L(w)\subseteq I_0$ and $I=I_1\sqcup(I_w\setminus J)$, if and only if $I\sqcup D_L(w)=I_1\sqcup(\Delta\setminus J)\subseteq I_0$.
Consequently, we see that
\[k_0=d(I_0\setminus D_L(w),I)+\ell(w)\geq \#I_0\setminus I_1-\#\Delta\setminus J+\#\Delta\setminus J=\#I_0\setminus I_1,\]
with equality holds if and only if $w\in\Gamma_{\Delta\setminus J}$, $x\unlhd w$ and $I\sqcup D_L(w)=I_1\sqcup(\Delta\setminus J)\subseteq I_0$, if and only if $x\unlhd w$, $w\in\Gamma^{I_0\setminus I_1}$ and $I=I_1\sqcup(\mathrm{Supp}(w)\setminus D_L(w))$ (with $J$ determined by $w$ via $\mathrm{Supp}(w)=\Delta\setminus J$).
\end{proof}

\begin{lem}\label{lem: special Ext with sm}
Let $j\in\Delta$ and $I\subseteq\widehat{j}\defeq \Delta\setminus\{j\}$. Then there exists a unique $V\in\mathrm{Rep}^{\rm{an}}_{\rm{adm}}(G)$ that fits into a non-split extension of the form $0\rightarrow \cF_{P_{\widehat{j}}}^{G}(L(s_j),V_{I,\widehat{j}}^{\infty})\rightarrow V\rightarrow V_{I\sqcup\{j\}}^{\infty}\rightarrow 0$. Moreover, we have
\begin{equation}\label{equ: special Ext with sm}
\mathrm{Ext}_{G}^k(V_{J}^{\infty},V)=0
\end{equation}
for each $J\supseteq I\sqcup\{j\}$ and $k\leq \#J\setminus I-1$.
\end{lem}
\begin{proof}
The existence and uniqueness of $V$ is clear from Lemma~\ref{lem: Ext1 OS simple} by taking $x_0=1$, $x_1=s_j$, $I_0=I\sqcup\{j\}$ and $I_1=I$. We prove the vanishing (\ref{equ: special Ext with sm}) for $k\leq \#J\setminus I-1=\#J\setminus(I\sqcup\{j\})$ by an increasing induction on $\#J$. If $J=I\sqcup\{j\}$, we clearly have $\Hom_{G}(V_{I\sqcup\{j\}}^{\infty},V)=0$ by the definition of $V$.
Assume from now that $J\supsetneq I\sqcup\{j\}$ and choose an arbitrary $j'\in J\setminus(I\sqcup\{j\})$ with $J'\defeq J\setminus\{j'\}$. It follows from \ref{it: general sm distance 2} of Lemma~\ref{lem: general sm distance} and Lemma~\ref{lem: Ext between sm} that $\Dim_E \mathrm{Ext}_{G}^1(V_{J}^{\infty},V_{J'}^{\infty})=1$, which together with \ref{it: sm cube 1} of Lemma~\ref{lem: sm cube} implies that $Q_{\Delta}(J',J)\in\mathrm{Rep}^{\rm{an}}_{\rm{adm}}(G)$ is the unique length $2$ object with socle $V_{J'}^{\infty}$ and cosocle $V_{J}^{\infty}$.
Note from \ref{it: sm cube 3} of Lemma~\ref{lem: sm cube} that $i_{\widehat{j},\Delta}^{\infty}(V_{I,\widehat{j}}^{\infty})\cong Q_{\Delta}(I\sqcup\{j\},I)$. It is thus clear from \ref{it: general sm distance 3} of Lemma~\ref{lem: general sm distance} that $d_{\Delta}(Q_{\Delta}(J',J),V_{I\sqcup\{j\}}^{\infty})=\infty=d_{\Delta}(Q_{\Delta}(J',J),V_{I}^{\infty})$ which by a d\'evissage with respect to $0\rightarrow V_{I\sqcup\{j\}}^{\infty}\rightarrow Q_{\Delta}(I\sqcup\{j\},I)\rightarrow V_{I}^{\infty}\rightarrow 0$ gives 
\[d_{\Delta}(Q_{\Delta}(J',J),V_{I,\widehat{j}}^{\infty})=d_{\Delta}(Q_{\Delta}(J',J),Q_{\Delta}(I\sqcup\{j\},I))=\infty.\] 
These together with Lemma~\ref{lem: sm distance Ext vanishing} implies that
\[\mathrm{Ext}_{G}^k(Q_{\Delta}(J',J),V_{I\sqcup\{j\}}^{\infty})=0=\mathrm{Ext}_{G}^k(Q_{\Delta}(J',J),\cF_{P_{\widehat{j}}}^{G}(L(s_j),V_{I,\widehat{j}}^{\infty}))\]
for $k\geq 0$, which by d\'evissage with respect to $0\rightarrow \cF_{P_{\widehat{j}}}^{G}(L(s_j),V_{I,\widehat{j}}^{\infty})\rightarrow V\rightarrow V_{I\sqcup\{j\}}^{\infty}\rightarrow 0$ gives
\begin{equation}\label{equ: special Ext with sm vanishing}
\mathrm{Ext}_{G}^k(Q_{\Delta}(J',J),V)=0
\end{equation}
for $k\geq 0$. Hence, the short exact sequence $0\rightarrow V_{J'}^{\infty}\rightarrow Q_{\Delta}(J',J)\rightarrow V_{J}^{\infty}\rightarrow 0$ induces an isomorphism
\begin{equation}\label{equ: special Ext with sm induction}
\mathrm{Ext}_{G}^k(V_{J'}^{\infty},V)\buildrel\sim\over\longrightarrow \mathrm{Ext}_{G}^{k+1}(V_{J}^{\infty},V)
\end{equation}
for each $k\geq 0$. By our inductive hypothesis we know that LHS of (\ref{equ: special Ext with sm induction}) is zero for $k\leq \#J'\setminus I-1$, which forces the RHS of \emph{loc.cit.} to be zero for $k\leq \#J\setminus I-1$. The proof is thus finished.
\end{proof}

Let $I\subseteq \Delta$ and $x\in\Gamma$ with $\mathrm{Supp}(x)\subseteq \Delta\setminus I$.
We set
\begin{equation}\label{equ: x factor}
C_{x}^{I}\defeq \cF_{P_{I_x}}^{G}(L(x),V_{I\sqcup(\mathrm{Supp}(x)\setminus D_L(x)),I_x}^{\infty})
\end{equation}
and
\begin{equation}\label{equ: x cube}
Q_{x}^{I}\defeq \cF_{P_{I_x}}^{G}(L(x),Q_{I_x}(I\sqcup(\mathrm{Supp}(x)\setminus D_L(x)),I)).
\end{equation}
As we have $\mathrm{soc}_{L_{I_x}}(Q_{I_x}(I\sqcup(\mathrm{Supp}(x)\setminus D_L(x)),I))\cong V_{I\sqcup(\mathrm{Supp}(x)\setminus D_L(x)),I_x}^{\infty}$ from \ref{it: sm cube 1} of Lemma~\ref{lem: sm cube}, we deduce from Lemma~\ref{lem: OS layer} (as well as (\ref{equ: x factor}) and (\ref{equ: x cube})) that
\begin{equation}\label{equ: x cube socle}
\mathrm{soc}_{G}(Q_{x}^{I})\cong C_{x}^{I}.
\end{equation}
We temporarily fix our choice of $I$ and $x$.
For each $I\subseteq I'\subseteq I\sqcup(\mathrm{Supp}(x)\setminus D_L(x))$, as $[V_{I}^{\rm{an}}: C_{x,I'}]=1$ by Lemma~\ref{lem: JH mult PS}, we write $\tld{C}_{x,I'}$ for the unique subrepresentation of $V_{I}^{\rm{an}}$ with cosocle $C_{x,I'}$.
In particular, we write $\tld{C}_{x}^{I}\defeq \tld{C}_{x,I\sqcup(\mathrm{Supp}(x)\setminus D_L(x))}$ (with cosocle $C_{x}^{I}$) and $\tld{Q}_{x}^{I}\defeq \tld{C}_{x,I}$ (with cosocle $C_{x,I}$) for short.
For each $I\subseteq I'\subseteq I''\subseteq I\sqcup(\mathrm{Supp}(x)\setminus D_L(x))$, it follows from \ref{it: sm cube 1} of Lemma~\ref{lem: sm cube} (applied to $Q_{I_x}(I\sqcup(\mathrm{Supp}(x)\setminus D_L(x)),I)$) that $V_{I'',I_x}^{\infty}\leq V_{I',I_x}^{\infty}$ inside $\mathrm{JH}_{L_{I_x}}(Q_{I_x}(I\sqcup(\mathrm{Supp}(x)\setminus D_L(x)),I))$ and thus $C_{x,I''}\leq C_{x,I'}$ inside $\mathrm{JH}_{G}(Q_{x}^{I})$ by (\ref{equ: x cube}) and Lemma~\ref{lem: OS layer}, which further implies $\tld{C}_{x,I''}\subseteq \tld{C}_{x,I'}$.
By taking $I''=I\sqcup(\mathrm{Supp}(x)\setminus D_L(x))$ and $I'=I$, we have
\[\tld{C}_{x}^{I}\subseteq \tld{Q}_{x}^{I}.\]
\begin{lem}\label{lem: basic subquotient of St}
Let $I\subseteq \Delta$ and $x\in\Gamma$ with $\mathrm{Supp}(x)\subseteq \Delta\setminus I$. We have the following results.
\begin{enumerate}[label=(\roman*)]
\item \label{it: basic quotient 1} The subrepresentation $\tld{Q}_{x}^{I}\subseteq V_{I}^{\rm{an}}$ admits a unique quotient that is isomorphic to $Q_{x}^{I}$, and this subquotient of $V_{I}^{\rm{an}}$ exhausts all $x$-typic constituents of $V_{I}^{\rm{an}}$.
\item \label{it: basic quotient 2} For each $x'\leq x$, we have $\tld{Q}_{x'}^{I}\subseteq \tld{Q}_{x}^{I}$.
\item \label{it: basic quotient 3} Given $w\in\Gamma$ and $I'\subseteq I_w$, we have $C_{w,I'}\in\mathrm{JH}_{G}(\tld{Q}_{x}^{I})$ if and only if $w\leq x$ and $C_{w,I'}\in\mathrm{JH}_{G}(Q_{w}^{I})$.
\end{enumerate}
\end{lem}
\begin{proof}
We fix our choice of $I$ and $x$ throughout the proof.
We prove \ref{it: basic quotient 1}. \\
It follows from \ref{it: sm cube 1} of Lemma~\ref{lem: sm cube} that $Q_{I_x}(I\sqcup(\mathrm{Supp}(x)\setminus D_L(x)),I)$ is multiplicity free with
\[\mathrm{JH}_{L_{I_x}}(Q_{I_x}(I\sqcup(\mathrm{Supp}(x)\setminus D_L(x)),I))=\{V_{I',I_x}^{\infty}\mid I\subseteq I'\subseteq I\sqcup(\mathrm{Supp}(x)\setminus D_L(x))\},\]
which together with Lemma~\ref{lem: JH OS} and our definition of $Q_{x}^{I}$ in (\ref{equ: x cube}) implies that $Q_{x}^{I}$ is multiplicity free with
\begin{equation}\label{equ: JH x cube}
\mathrm{JH}_{G}(Q_{x}^{I})=\{C_{x,I'}\mid I\subseteq I'\subseteq I\sqcup(\mathrm{Supp}(x)\setminus D_L(x))\}.
\end{equation}
Recall from  Lemma~\ref{lem: JH mult PS} that the set of $x$-typic constituents of $V_{I}^{\rm{an}}$ equals (\ref{equ: JH x cube}) and each $x$-typic constituent of $V_{I}^{\rm{an}}$ appears with multiplicity one. Note that $i_{\emptyset}^{\rm{an}}=\cF_{B}^{G}(M(1),1_{T})$ admits $\cF_{B}^{G}(L(x),1_{T})\cong \cF_{P_{I_x}}^{G}(L(x),Q_{I_x}(I_x,\emptyset))$ as a subquotient, which further admits $Q_{x}^{I}$ as a subquotient (as $Q_{I_x}(I\sqcup(\mathrm{Supp}(x)\setminus D_L(x)),I)$ is a subquotient of $Q_{I_x}(I_x,\emptyset)$ by \ref{it: sm cube 2} of Lemma~\ref{lem: sm cube}). Since $V_{I}^{\rm{an}}$ is a subquotient of $i_{\emptyset}^{\rm{an}}=\cF_{B}^{G}(M(1),1_{T})$ and each element of (\ref{equ: JH x cube}) appears with multiplicity one in both $V_{I}^{\rm{an}}$ and $i_{\emptyset}^{\rm{an}}$, the subquotient $Q_{x}^{I}$ is necessarily a subquotient of $V_{I}^{\rm{an}}$. Since $Q_{I_x}(I\sqcup(\mathrm{Supp}(x)\setminus D_L(x)),I)$ has cosocle $V_{I,I_x}^{\infty}$ (see \ref{it: sm cube 1} of Lemma~\ref{lem: sm cube}), we deduce from Lemma~\ref{lem: JH OS} that $Q_{x}^{I}$ has cosocle $C_{x,I}$, which together with the definition of $\tld{Q}_{x}^{I}$ forces $Q_{x}^{I}$ to be a quotient $\tld{Q}_{x}^{I}$.

We prove \ref{it: basic quotient 2}.\\
We fix a choice of $x'\leq x$. By the definition of $\tld{Q}_{x'}^{I}$ and $\tld{Q}_{x}^{I}$, it suffices to show that $C_{x',I}\in\mathrm{JH}_{G}(\tld{Q}_{x}^{I})$. Let $M$ be the unique quotient of $M^{I}(1)$ with socle $L(x)$ (by replacing $x$ and $w$ in \ref{it: g coxeter 2} of Lemma~\ref{lem: coxeter subquotient} with $1$ and $x$ here). By \emph{loc.cit.} we know that $L(x')\in\mathrm{JH}_{U(\fg)}(M)$ and thus
\[\cF_{P_{I}}^{G}(L(x'),1_{L_{I}})\cong \cF_{P_{I_{x'}}}^{G}(L(x'),i_{I,I_{x'}}^{\infty}(1_{L_{I}}))\cong \cF_{P_{I_{x'}}}^{G}(L(x'),Q_{I_{x'}}(I_{x'},I))\]
is a subquotient of $\cF_{P_{I}}^{G}(M,1_{L_{I}})$. As $M$ has socle $L(x)$, by Lemma~\ref{lem: OS socle cosocle} and Lemma~\ref{lem: JH OS} we know that
\begin{multline}\label{equ: cosocle M sub}
\mathrm{cosoc}_{G}(\cF_{P_{I}}^{G}(M,1_{L_{I}}))=\mathrm{cosoc}_{G}(\cF_{P_{I}}^{G}(L(x),1_{L_{I}}))\\
=\mathrm{cosoc}_{G}(\cF_{P_{I_x}}^{G}(L(x),Q_{I_x}(I_x,I)))=\cF_{P_{I_x}}^{G}(L(x),\mathrm{cosoc}_{L_{I_x}}(Q_{I_x}(I_x,I)))= C_{x,I}.
\end{multline}
The surjection $M^{I}(1)\twoheadrightarrow M$ and the embedding $1_{L_{I}}\hookrightarrow i_{\emptyset,I}^{\infty}(1_{T})$ naturally induce an embedding
\[\cF_{P_{I}}^{G}(M,1_{L_{I}})\hookrightarrow \cF_{P_{I}}^{G}(M^{I}(1),1_{L_{I}})=i_{I}^{\rm{an}},\]
which together with (\ref{equ: cosocle M sub}) and $C_{x,I}$ having multiplicity one in $i_{I}^{\rm{an}}$ (see  Lemma~\ref{lem: JH mult PS}) forces $\cF_{P_{I}}^{G}(M,1_{L_{I}})$ to be the unique subrepresentation of $i_{I}^{\rm{an}}$ with cosocle $C_{x,I}$. It follows from Lemma~\ref{lem: JH OS}, $[M:L(x')]=1$ from \ref{it: g coxeter 1} of Lemma~\ref{lem: coxeter subquotient}, and $i_{I,I_{x'}}^{\infty}(1_{L_{I}})\cong Q_{I_{x'}}(I_{x'},I)$ with $[Q_{I_{x'}}(I_{x'},I):V_{I,I_{x'}}^{\infty}]=1$ from Lemma~\ref{lem: sm cube} that
\[[\cF_{P_{I}}^{G}(M,1_{L_{I}}):C_{x',I}]=[M:L(x')][i_{I,I_{x'}}^{\infty}(1_{L_{I}}),V_{I,I_{x'}}^{\infty}]=1,\]
which together with $[i_{I}^{\rm{an}},C_{x',I}]=[V_{I}^{\rm{an}},C_{x',I}]=1$ from  Lemma~\ref{lem: JH mult PS} forces $C_{x',I}$ to appear with multiplicity one in the image of the composition $\cF_{P_{I}}^{G}(M,1_{L_{I}})\hookrightarrow i_{I}^{\rm{an}}\twoheadrightarrow V_{I}^{\rm{an}}$ which is nothing but $\tld{Q}_{x}^{I}$ by its own definition and (\ref{equ: cosocle M sub}).

We prove \ref{it: basic quotient 3}. Let $w\in\Gamma$ and $I'\subseteq I_w$ such that $C_{w,I'}\in\mathrm{JH}_{G}(\tld{Q}_{x}^{I})$. By the proof of \ref{it: basic quotient 2} we see that $C_{w,I'}\in\mathrm{JH}_{G}(\cF_{P_{I}}^{G}(M,1_{L_{I}}))$ and thus $[M:L(w)]\neq 0$ (see Lemma~\ref{lem: JH OS}). This together with $w\in\Gamma$ and \ref{it: g coxeter 1} of Lemma~\ref{lem: coxeter subquotient} forces $w\leq x$. Note that $Q_{w}^{I}$ is a quotient of $\tld{Q}_{w}^{I}$ by \ref{it: basic quotient 1} and thus a subquotient of $\tld{Q}_{x}^{I}$ by \ref{it: basic quotient 2}. It follows from (\ref{equ: x cube}) (with $w$ replacing $x$) and Lemma~\ref{lem: JH mult PS} that $\mathrm{JH}_{G}(Q_{w}^{I})=\{C_{w,I'}\mid I\subseteq I\subseteq I\sqcup(\mathrm{Supp}(w)\setminus D_L(w))\}$ exhausts all $w$-typic constituents of $V_{I}^{\rm{an}}$. This finishes the proof of \ref{it: basic quotient 3}.
\end{proof}

\begin{lem}\label{lem: length 2 subquotient}
Let $I\subseteq \Delta$, $x_i\in\Gamma$ and $I_i\subseteq I_{x_i}$ for $i=0,1$ with $x_0\neq x_1$.
Then $V_{I}^{\rm{an}}$ admits a length $2$ subquotient with socle $C_{x_1,I_1}$ and cosocle $C_{x_0,I_0}$ if and only if $I\subseteq I_i\subseteq I\sqcup(\mathrm{Supp}(x_i)\setminus D_L(x_i))$ (with $\mathrm{Supp}(x_i)\subseteq \Delta\setminus I$) for $i=0,1$, $x_0>x_1$ with $\ell(x_0)=\ell(x_1)+1$, and $I_0=I_1\sqcup(D_L(x_1)\setminus D_L(x_0))$. In particular, if $x_0=s_jx_1$ for some $j\in\mathrm{Supp}(x_0)\setminus\mathrm{Supp}(x_1)$, then $V_{I}^{\rm{an}}$ admits a unique length $2$ subquotient with socle $C_{x_1}^{I}$ and cosocle $C_{x_0}^{I}$.
\end{lem}
\begin{proof}
It follows from Lemma~\ref{lem: JH mult PS} that, for $i=0,1$, we have $C_{x_i,I_i}\in\mathrm{JH}_{G}(V_{I}^{\rm{an}})$ if and only if $I\subseteq I_i\subseteq I\sqcup(\mathrm{Supp}(x_i)\setminus D_L(x_i))$ with $\mathrm{Supp}(x_i)\subseteq \Delta\setminus I$, in which case we have
\begin{equation}\label{equ: length 2 subquotient mult 1}
[i_{I}^{\rm{an}}:C_{x_i,I_i}]=[V_{I}^{\rm{an}}:C_{x_i,I_i}]=[\cF_{P_I}^{G}(L(x_i),1_{L_I}):C_{x_i,I_i}]=1.
\end{equation}
Hence, we assume $I\subseteq I_i\subseteq I\sqcup(\mathrm{Supp}(x_i)\setminus D_L(x_i))$ throughout the rest of the proof.

We prove the `only if' direction.\\
We write $V$ for the subquotient of $V_{I}^{\rm{an}}$ with socle $C_{x_1,I_1}$ and cosocle $C_{x_0,I_0}$.
Let $M_0$ be the unique quotient of $M^I(1)$ with socle $L(x_0)$ (see \ref{it: g coxeter 2} of Lemma~\ref{lem: coxeter subquotient}) and note that $[M_0:L(x_1)]\neq 0$ if and only if $x_1<x_0$ (using our assumption $x_0\neq x_1$), in which case $[M_0:L(x_1)]=1$. By \ref{it: OS property 1} of Proposition~\ref{prop: OS property} we know that $\cF_{P_I}^{G}(M_0,1_{L_I})$ is a subrepresentation of $i_{I}^{\rm{an}}=\cF_{P_I}^{G}(M^I(1),1_{L_I})$. It follows from $[M_0:L(x_0)]=1$ and (\ref{equ: length 2 subquotient mult 1}) that $V$ must be a subquotient of $\cF_{P_I}^{G}(M_0,1_{L_I})$, which together with Lemma~\ref{lem: JH OS} forces $[M_0:L(x_1)]\neq 0$ and thus $x_1<x_0$ by our previous discussion.
The existence of $V$ forces $\mathrm{Ext}_{G}^1(C_{x_0,I_0},C_{x_1,I_1})\neq 0$, which together with Lemma~\ref{lem: Ext1 OS simple} gives $\mathrm{Ext}_{U(\fg)}^1(L(x_0),L(x_1))\neq 0$ and $I_0=I_1\sqcup(D_L(x_1)\setminus D_L(x_0))$. Note that $\mathrm{Ext}_{U(\fg)}^1(L(x_0),L(x_1))\neq 0$ together with $x_1<x_0$ and $x_0\in\Gamma$ forces $\ell(x_0)=\ell(x_1)+1$ by Lemma~\ref{lem: g Ext1}.

We prove the `if' direction.\\
As $x_1<x_0$ with $\ell(x_0)=\ell(x_1)+1$, there exists a unique length $2$ subquotient $M$ of $M^I(1)$ with socle $L(x_0)$ and cosocle $L(x_1)$ (see \ref{it: g coxeter 2} of Lemma~\ref{lem: coxeter subquotient}). It is clear from \ref{it: OS property 1} of Proposition~\ref{prop: OS property} that $\cF_{P_I}^{G}(M,1_{L_I})$ is a subquotient of $i_{I}^{\rm{an}}=\cF_{P_I}^{G}(M^I(1),1_{L_I})$. Our assumption $I_0=I_1\sqcup(D_L(x_1)\setminus D_L(x_0))$ together with the proof of Lemma~\ref{lem: Ext1 OS simple} (see (\ref{equ: Ext1 OS simple sm})) implies that $d_{\Delta}(V_{I_0,I_{x_0}}^{\infty},V_{I_1,I_{x_1}}^{\infty})=0$. Consequently, we can apply \cite[Lem.~5.1.19]{BQ24} with $M$ and $I$ as above and
\[\sigma_{x_i}^{\infty}\defeq V_{I_i,I_{x_i}}^{\infty}\in\mathrm{JH}_{L_{I_{x_i}}}(i_{I,I_{x_i}}^{\infty}(1_{L_I}))\] for $i=0,1$ and deduce that $\cF_{P_I}^{G}(M,1_{L_I})$ admits a unique length $2$ subquotient $V$ with socle $C_{x_1,I_1}$ and cosocle $C_{x_0,I_0}$. Since both $\cF_{P_I}^{G}(M,1_{L_I})$ and $V_{I}^{\rm{an}}$ are subquotients of $i_{I}^{\rm{an}}$ and we have $[i_{I}^{\rm{an}}:C_{x_i,I_i}]=1$ for $i=0,1$ from (\ref{equ: length 2 subquotient mult 1}), we see that the length $2$ subquotient $V$ of $\cF_{P_I}^{G}(M,1_{L_I})$ must be a subquotient of $V_{I}^{\rm{an}}$.

We prove the last claim.\\
We take $I_i\defeq I\sqcup(\mathrm{Supp}(x_i)\setminus D_L(x_i))$ for $i=0,1$.
It suffices to check $I_0=I_1\sqcup(D_L(x_1)\setminus D_L(x_0))$. Since $\mathrm{Supp}(x_0)=\mathrm{Supp}(x_1)\sqcup\{j\}$ and $j\in D_L(x_0)$, we deduce that $\mathrm{Supp}(x_0)\setminus D_L(x_0)=\mathrm{Supp}(x_1)\setminus D_L(x_0)=(\mathrm{Supp}(x_1)\setminus D_L(x_1))\sqcup (D_L(x_1)\setminus D_L(x_0))$, which clearly gives $I_0=I_1\sqcup(D_L(x_1)\setminus D_L(x_0))$.
\end{proof}

Recall the partial-order $\unlhd$ on $\Gamma$ from the discussion below (\ref{equ: coxeter left}).
\begin{lem}\label{lem: inclusion order}
Let $I\subseteq \Delta$ and $x,w\in\Gamma$ with $\mathrm{Supp}(x),\mathrm{Supp}(w)\subseteq \Delta\setminus I$. Then we have $\tld{C}_{x}^{I}\subseteq \tld{C}_{w}^{I}$ if and only if $x\unlhd w$.
\end{lem}
\begin{proof}
According to the definition of $\tld{C}_{x}^{I}$ and $\tld{C}_{w}^{I}$, it is clear that $\tld{C}_{x}^{I}\subseteq \tld{C}_{w}^{I}$ if and only if $C_{x}^{I}\in\mathrm{JH}_{G}(\tld{C}_{w}^{I})$ if and only if $C_{w}^{I}$ appears in the unique quotient of $V_{I}^{\rm{an}}$ with socle $C_{x}^{I}$. We divide the proof into the following steps.

\textbf{Step $1$}: We prove that $\tld{C}_{x}^{I}\subseteq \tld{C}_{w}^{I}$ when $x\unlhd w$.\\
We can choose a sequence $x=y_0\lhd y_1\lhd \cdots\lhd y_t=w$ with $t=\ell(w)-\ell(x)$ and prove $\tld{C}_{y_{t'-1}}^{I}\subseteq \tld{C}_{y_{t'}}^{I}$ for each $1\le t'\leq t$. Hence, it is harmless to assume that $x\lhd w$ with $\ell(w)=\ell(x)+1$, namely $w=s_jx$ with $\mathrm{Supp}(w)=\mathrm{Supp}(x)\setminus\{j\}$. In this case, it follows from Lemma~\ref{lem: length 2 subquotient} that $V_{I}^{\rm{an}}$ admits a (unique) length $2$ subquotient with socle $C_{x}^{I}$ and cosocle $C_{w}^{I}$, which clearly implies $\tld{C}_{x}^{I}\subseteq \tld{C}_{w}^{I}$.

\textbf{Step $2$}: We prove that $\tld{C}_{x}^{I}\not\subseteq \tld{C}_{w}^{I}$ when $x\not \leq w$.\\
Let $M$ be the unique quotient of $M^I(1)$ with socle $L(w)$ and we have $L(x)\notin\mathrm{JH}_{U(\fg)}(M)$ by \ref{it: g coxeter 2} of Lemma~\ref{lem: coxeter subquotient}. By \ref{it: OS property 1} of Proposition~\ref{prop: OS property} we know that $\cF_{P_I}^{G}(M,1_{L_I})$ is a subrepresentation of $i_{I}^{\infty}=\cF_{P_I}^{G}(M^I(1),1_{L_I})$. As $[M:L(w)]=1$ and $[M:L(x)]=0$ (see \ref{it: g coxeter 2} of Lemma~\ref{lem: coxeter subquotient}), we know that $\cF_{P_I}^{G}(M,1_{L_I})$ does not admit any $x$-typic constituent yet exhausts $w$-typic constituents of $i_{I}^{\rm{an}}$. As both $C_{x}^{I}$ and $C_{w}^{I}$ appears with multiplicity one in $i_{I}^{\rm{an}}$ and $V_{I}^{\rm{an}}$ by  Lemma~\ref{lem: JH mult PS}, we deduce that the image of the composition
\[\cF_{P_I}^{G}(M,1_{L_I})\rightarrow i_{I}^{\rm{an}}\rightarrow V_{I}^{\rm{an}}\]
admits $C_{w}^{I}$ but not $C_{x}^{I}$ as constituents. By definition of $\tld{C}_{x}^{I}$ and $\tld{C}_{w}^{I}$ we have $\tld{C}_{x}^{I}\not\subseteq \tld{C}_{w}^{I}$.

\textbf{Step $3$}: We prove that $\tld{C}_{x}^{I}\not\subseteq \tld{C}_{w}^{I}$ when $x\leq w$ and $x\not\unlhd w$.\\
We assume on the contrary that $\tld{C}_{x}^{I}\subseteq \tld{C}_{w}^{I}$.
As $x\leq w$ and $x\not\unlhd w$, there exists $j\in\mathrm{Supp}(x)$ and $j'\in\mathrm{Supp}(w)\setminus\mathrm{Supp}(x)$ such that $|j-j'|=1$ and $s_js_{j'}\leq w$ by \cite[Thm.~2.2.2]{BB05}.
We write $I_{+}\defeq I\sqcup\{j'\}$ and $I'\defeq I\sqcup(\mathrm{Supp}(w)\setminus (D_L(w)\sqcup\{j'\}))$ for short.
Let $M''$ be the unique subquotient of $M^I(1)$ with socle $L(w)$ and cosocle $L(x)$ (see \ref{it: g coxeter 2} of Lemma~\ref{lem: coxeter subquotient}), and note that we have $M''\in\cO^{\fp_{I_{+}}}_{\rm{alg}}$ from Lemma~\ref{lem: weak dominance}.
By \ref{it: OS property 1} and \ref{it: OS property 2} of Proposition~\ref{prop: OS property} (as well as $V_{I,I_{+}}^{\infty}\in\mathrm{JH}_{L_{I_{+}}}(Q_{I_{+}}(I_{+},I))$ from \ref{it: sm cube 1} of Lemma~\ref{lem: sm cube}) we know that $V''\defeq \cF_{P_{I_{+}}}^{G}(M'',V_{I,I_{+}}^{\infty})$ is a subquotient of $\cF_{P_I}^{G}(M'',1_{L_I})\cong\cF_{P_{I_{+}}}^{G}(M'',Q_{I_{+}}(I_{+},I))$ and thus a subquotient of $i_{I}^{\rm{an}}=\cF_{P_I}^{G}(M^I(1),1_{L_I})$.
Note that $i_{I_{+},I_x}^{\infty}(V_{I,I_{+}}^{\infty})\cong Q_{I_x}(I\sqcup(I_x\setminus I_{+}),I)$ by \ref{it: sm cube 3} of Lemma~\ref{lem: sm cube}, which together with $V_{I\sqcup(\mathrm{Supp}(x)\setminus D_L(x)),I_x}^{\infty}\in \mathrm{JH}_{L_{I_x}}(Q_{I_x}(I\sqcup(I_x\setminus I_{+}),I))$ (from $I\subseteq I\sqcup(\mathrm{Supp}(x)\setminus D_L(x))\subseteq I\sqcup(I_x\setminus I_{+})$ and \ref{it: sm cube 1} of Lemma~\ref{lem: sm cube}) gives
\[V_{I\sqcup(\mathrm{Supp}(x)\setminus D_L(x)),I_x}^{\infty}\in \mathrm{JH}_{L_{I_x}}(i_{I_{+},I_x}^{\infty}(V_{I,I_{+}}^{\infty}))\]
and thus $C_{x}^{I}=\cF_{P_{I_x}}^{G}(L(x),V_{I\sqcup(\mathrm{Supp}(x)\setminus D_L(x)),I_x}^{\infty})\in \mathrm{JH}_{G}(V'')$ by Lemma~\ref{lem: JH OS}.
Similarly, it follows from $I\subseteq I'\subseteq I\sqcup(I_w\setminus I_{+})$ (resp.~$j'\in I\sqcup(\mathrm{Supp}(w)\setminus D_L(w))$ but $j'\notin I\sqcup(I_w\setminus I_{+})$) that $V_{I',I_w}^{\infty}\in \mathrm{JH}_{L_{I_w}}(Q_{I_w}(I\sqcup(I_w\setminus I_{+}),I))$ and $C_{w,I'}\in\mathrm{JH}_{G}(V'')$ (resp.~$V_{I\sqcup(\mathrm{Supp}(w)\setminus D_L(w)),I_w}^{\infty}\notin \mathrm{JH}_{L_{I_w}}(Q_{I_w}(I\sqcup(I_w\setminus I_{+}),I))$ and $C_{w}^{I}\notin\mathrm{JH}_{G}(V'')$).
Recall from  Lemma~\ref{lem: JH mult PS} that $C_{x}^{I}$, $C_{w}^{I}$ and $C_{w,I'}$ all have multiplicity one in $i_{I}^{\rm{an}}$ and $V_{I}^{\rm{an}}$, which forces $C_{x}^{I}$ and $C_{w,I'}$ to have multiplicity one in $V''$ (with $C_{w}^{I}\notin\mathrm{JH}_{G}(V'')$).
Recall from \ref{it: basic quotient 1} of Lemma~\ref{lem: basic subquotient of St} that $V_{I}^{\rm{an}}$ admits $Q_{w}^{I}=\cF_{P_{I_w}}^{G}(L(w),Q_{I_w}(I\sqcup(\mathrm{Supp}(w)\setminus D_L(w)),I))$ as a subquotient. It follows from \ref{it: sm cube 1} of Lemma~\ref{lem: sm cube} that $Q_{I_w}(I\sqcup(\mathrm{Supp}(w)\setminus D_L(w)),I)$ admits a unique length $2$ subrepresentation with socle $V_{I\sqcup(\mathrm{Supp}(w)\setminus D_L(w)),I_w}^{\infty}$ and cosocle $V_{I',I_w}^{\infty}$, which together with Lemma~\ref{lem: OS layer} implies that $Q_{w}^{I}$ admits a unique length $2$ subrepresentation $C_{w,I'}^+$ with socle $C_{w}^{I}$ and cosocle $C_{w,I'}$.
We write $\tld{C}_{w,I'}$ for the unique subrepresentation of $V_{I}^{\rm{an}}$ with cosocle $C_{w,I'}$.
As $C_{w,I'}^+$ is a subrepresentation of $Q_{w}^{I}$, it is a subquotient of $V_{I}^{\rm{an}}$, which together with the definition of $\tld{C}_{w}^{I}$ and $\tld{C}_{w,I'}$ forces $\tld{C}_{w}^{I}\subseteq \tld{C}_{w,I'}$.
This together with our assumption $\tld{C}_{x}^{I}\subseteq \tld{C}_{w}^{I}$ implies that $V_{I}^{\rm{an}}$ admits a subquotient $W$ with socle $C_{x}^{I}$ and cosocle $C_{w,I'}$ which moreover satisfies $C_{w}^{I}\in\mathrm{JH}_{G}(W)$. As $C_{x}^{I}$ and $C_{w,I'}$ have multiplicity one in $V_{I}^{\rm{an}}$, $i_{I}^{\rm{an}}$ and $V''$ (with both $V''$ and $V_{I}^{\rm{an}}$ being subquotients of $i_{I}^{\rm{an}}$), the subquotient $W$ above must also be a subquotient of $V''$. We find a contradiction between $C_{w}^{I}\in\mathrm{JH}_{G}(W)$ and $C_{w}^{I}\notin\mathrm{JH}_{G}(V'')$.

We finish the proof by combining \textbf{Step $1$}, \textbf{Step $2$} and \textbf{Step $3$}.
\end{proof}

\begin{lem}\label{lem: C coxeter factor}
Let $I\subseteq \Delta$ and $x,w\in\Gamma$ with $\mathrm{Supp}(x),\mathrm{Supp}(w)\subseteq \Delta\setminus I$. Let $I'\subseteq I_x$. Then we have $C_{x,I'}\in\mathrm{JH}_{G}(\tld{C}_{w}^{I})$ if and only if $x\unlhd w$ and $C_{x,I'}=C_{x}^{I}$.
\end{lem}
\begin{proof}
Recall that we write $\tld{C}_{x,I'}$ for the unique subrepresentation of $V_{I}^{\rm{an}}$ with cosocle $C_{x,I'}$, and note that $C_{x,I'}\in\mathrm{JH}_{G}(\tld{C}_{w}^{I})$ if and only if $\tld{C}_{x,I'}\subseteq\tld{C}_{w}^{I}$.
If $x\unlhd w$ and $C_{x,I'}=C_{x}^{I}$, then we clearly have $\tld{C}_{x,I'}=\tld{C}_{x}^{I}\subseteq \tld{C}_{w}^{I}$ by Lemma~\ref{lem: inclusion order}.
We fix from now on a choice of $I$, $x$ and $w$ as above as well as $I'\subseteq I_x$ such that $\tld{C}_{x,I'}\subseteq \tld{C}_{w}^{I}$. Our assumption $C_{x,I'}\in\mathrm{JH}_{G}(V_{I}^{\rm{an}})$ together with Lemma~\ref{lem: JH mult PS} forces $I\subseteq I'\subseteq I\sqcup(\mathrm{Supp}(x)\setminus D_L(x))$.
Recall from the discussion right before Lemma~\ref{lem: basic subquotient of St} that we have $C_{x,I''}\leq C_{x,I'}$ inside $\mathrm{JH}_{G}(Q_{x}^{I})$ and thus $\tld{C}_{x,I''}\subseteq \tld{C}_{x,I'}$ for each $I''$ satisfying $I'\subseteq I''\subseteq I\sqcup(\mathrm{Supp}(x)\setminus D_L(x))$.
In particular, we have $\tld{C}_{x}^{I}\subseteq \tld{C}_{x,I'}$ and thus $\tld{C}_{x}^{I}\subseteq \tld{C}_{w}^{I}$, which together with Lemma~\ref{lem: inclusion order} forces $x\unlhd w$.
Assume on the contrary that $C_{x,I'}\neq C_{x}^{I}$ namely $I'\neq I\sqcup(\mathrm{Supp}(x)\setminus D_L(x))$, then we may choose an arbitrary $j\in(\mathrm{Supp}(x)\setminus D_L(x))\setminus I'$. By \cite[Thm.~2.2.2]{BB05} there exists a unique $x'<x$ such that $\mathrm{Supp}(x')=\mathrm{Supp}(x)\setminus\{j\}$. As $j\in\mathrm{Supp}(x)\setminus D_L(x)$ and $\mathrm{Supp}(x)=\mathrm{Supp}(x')\sqcup\{j\}$, there exists $j_0\in \mathrm{Supp}(x')$ such that $|j_0-j|=1$ and $s_{j_0}s_j\leq x\leq w$, which forces $x'\not\unlhd x$ and $x'\not\unlhd w$.
It follows from Lemma~\ref{lem: length 2 subquotient} that $V_{I}^{\rm{an}}$ admits a unique length $2$ subquotient with socle $C_{x'}^{I}$ and cosocle $C_{x,I''}$ with $I''\defeq I\sqcup(\mathrm{Supp}(x)\setminus (D_L(x)\cup\{j\}))$, which forces $\tld{C}_{x'}^{I}\subseteq \tld{C}_{x,I''}$. As $I\subseteq I'\subseteq I''\subseteq I\sqcup(\mathrm{Supp}(x)\setminus D_L(x))$, we deduce that $\tld{C}_{x,I''}\subseteq \tld{C}_{x,I'}$ (see the discussion right before Lemma~\ref{lem: basic subquotient of St}) and thus $\tld{C}_{x'}^{I}\subseteq \tld{C}_{x,I''}\subseteq \tld{C}_{x,I'}\subseteq \tld{C}_{w}^{I}$, which together with Lemma~\ref{lem: inclusion order} forces $x'\unlhd w$, a contradiction to $x'\not\unlhd w$ from previous discussion.
Consequently, we must have $C_{x,I'}=C_{x}^{I}$ and the proof is thus finished.
\end{proof}

\begin{rem}\label{rem: coxeter layer structure}
Let $I\subseteq \Delta$. We consider the following subset of $\mathrm{JH}_{G}(V_{I}^{\rm{an}})$ (see \ref{it: JH mult PS 2} of Lemma~\ref{lem: JH mult PS})
\begin{equation}\label{equ: coxeter bottom constituent}
\{C_{x}^{I}\mid x\in\Gamma^{\Delta\setminus I}\}
\end{equation}
with each $C_{x}^{I}$ appearing in $V_{I}^{\rm{an}}$ with multiplicity one. The layer structure of $V_{I}^{\rm{an}}$ equips (\ref{equ: coxeter bottom constituent}) with a natural partial-order, namely $C_{x}^{I}\leq C_{w}^{I}$ if and only if $\tld{C}_{x}^{I}\subseteq\tld{C}_{w}^{I}$. Then Lemma~\ref{lem: inclusion order} says that this natural partial-order on (\ref{equ: coxeter bottom constituent}) equals the one induced from the partial-order $\unlhd$ on $\Gamma^{\Delta\setminus I}$. We also know from Lemma~\ref{lem: C coxeter factor} that, given $x\in\Gamma^{\Delta\setminus I}$, the set of all constituents of $\tld{C}_{x}^{I}$ which are $y$-typic for some $y\in\Gamma$ is precisely given by
\[\{C_{x'}^{I}\mid x'\unlhd x\}.\]
Thanks to Lemma~\ref{lem: x Ext with sm}, the partially-ordered set (\ref{equ: coxeter bottom constituent}) will play a central role in our discussion on \emph{Coxeter filtration} on $\mathrm{Ext}$-groups such as $\mathrm{Ext}_{G}^{\#\Delta\setminus I}(1_{G},V_{I}^{\rm{an}})$, see \S~\ref{subsec: Coxeter Fil}.
\end{rem}
\subsection{Some vanishing results}\label{subsec: Ext vanishing}
We prove several vanishing results on $\mathrm{Ext}_{G}^{\bullet}(-,-)$ between a Orlik--Strauch representation and a (parabolically induced) principal series, and discuss applications of such results.

For each $x\in W(G)$ and $I\subseteq I_x$, we recall from \S \ref{subsec: OS} that $i_{x,I}^{\rm{an}}=i_{I}^{\rm{an}}(L^I(x)^\vee)\cong \cF_{P_I}^{G}(M^I(x),1_{L_I})$.
\begin{lem}\label{lem: vanishing by sm}
Let $I_0,I_1\subseteq \Delta$, $M_0\in\cO^{\fp_{I_0}}_{\rm{alg}}$, $I\subseteq I_0$ and $x_1\in W(G)$ with $I_1\subseteq I_{x_1}$.
If $I_1\not\subseteq I_0$ or $I_0\setminus I\not\subseteq I_1$, then we have
\begin{equation}\label{equ: vanishing by sm}
\mathrm{Ext}_{G}^k(\cF_{P_{I_0}}^{G}(M_0,V_{I,I_0}^{\infty}),i_{x_1,I_1}^{\rm{an}})=0
\end{equation}
for each $k\geq 0$.
\end{lem}
\begin{proof}
It follows from \ref{it: sm cube 3} of Lemma~\ref{lem: sm cube} that $i_{I_0,\Delta}^{\infty}(V_{I,I_0}^{\infty})\cong Q_{\Delta}(I\sqcup(\Delta\setminus I_0),I)$ and $i_{I_1,\Delta}^{\infty}(1_{L_{I_1}})\cong Q_{\Delta}(\Delta,I_1)$.
It follows from Lemma~\ref{lem: special sm distance} (by replacing $I_0$, $I_1$ and $J$ in \emph{loc.cit.} with $I\sqcup(\Delta\setminus I_0)$, $I$ and $I_1$ respectively) that $d_{\Delta}(Q_{\Delta}(I\sqcup(\Delta\setminus I_0),I),Q_{\Delta}(\Delta,I_1))<\infty$ if and only if $\Delta\setminus (I\sqcup(\Delta\setminus I_0))\subseteq I_1\subseteq \Delta\setminus (I\sqcup(\Delta\setminus I_0))\sqcup I$, if and only if $I_0\setminus I\subseteq I_1\subseteq I_0\cup I=I_0$ (using $I\subseteq I_0$). Since we assume either $I_1\not\subseteq I_0$ or $I_0\setminus I\not\subseteq I_1$, we deduce that
\[d_{\Delta}(V_{I,I_0}^{\infty},1_{L_{I_1}})=d_{\Delta}(i_{I_0,\Delta}^{\infty}(V_{I,I_0}^{\infty}),i_{I_1,\Delta}^{\infty}(1_{L_{I_1}}))=\infty,\]
which together with Lemma~\ref{lem: sm distance Ext vanishing} finishes the proof.
\end{proof}

\begin{lem}\label{lem: two PS vanishing}
Let $x_i\in W(G)$ and $I_i\subseteq I_{x_i}$ for $i=0,1$. If $I_1\not\subseteq I_0$, then we have
\begin{equation}\label{equ: two PS vanishing}
\mathrm{Ext}_{G}^k(i_{x_0,I_0}^{\rm{an}}, i_{x_1,I_1}^{\rm{an}})=0
\end{equation}
for each $k\geq 0$.
\end{lem}
\begin{proof}
This follows directly from Lemma~\ref{lem: vanishing by sm} by taking $M_0$ and $I$ in \emph{loc.cit.} to be $M^{I_0}(x_0)$ and $I_0$ respectively.
\end{proof}

\begin{lem}\label{lem: Hom x w St}
Let $I_0\subseteq\Delta$ and $x,w\in\Gamma^{\Delta\setminus I_0}$ with $x\unlhd w$. Then we have
\begin{equation}\label{equ: Hom x w St}
\Dim_E\Hom_{G}(V_{x,I_0}^{\rm{an}},V_{w,I_0}^{\rm{an}})=1.
\end{equation}
\end{lem}
\begin{proof}
It follows from Proposition~\ref{prop: St key seq}, (\ref{equ: Lie parabolic Verma wt shift}) and Lemma~\ref{it: sm cube 2} of Lemma~\ref{lem: sm cube} that
\[\Hom_{G}(i_{x,I_0}^{\rm{an}},i_{w,I_0}^{\rm{an}})\cong \Hom_{G}(i_{I_0}^{\infty},i_{I_0}^{\infty})\otimes_E\Hom_{U(\fg)}(M^{I_0}(w),M^{I_0}(x))\]
is $1$ dimensional. Since we have
\[\mathrm{Ext}_{G}^k(i_{x,I_0}^{\rm{an}},i_{w,I}^{\rm{an}})=0\]
for each $k\geq 0$ and each $I_0\subsetneq I\subseteq J_{w}=\Delta\setminus\mathrm{Supp}(w)$ by Lemma~\ref{lem: two PS vanishing}, we deduce from (see (\ref{equ: general Tits resolution}))
\[\mathbf{C}^{w}_{I_0,J_{w}}\cong V_{w,I_0}^{\rm{an}}[\#I_0]\]
that the surjection $i_{w,I_0}^{\rm{an}}\twoheadrightarrow V_{w,I_0}^{\rm{an}}$ induces an isomorphism
\[\Hom_{G}(i_{x,I_0}^{\rm{an}},i_{w,I_0}^{\rm{an}})\buildrel\sim\over\longrightarrow \Hom_{G}(i_{x,I_0}^{\rm{an}},V_{w,I_0}^{\rm{an}})\]
with the RHS being necessarily $1$-dimensional. Now that the surjection $i_{x,I_0}^{\rm{an}}\twoheadrightarrow V_{x,I_0}^{\rm{an}}$ evidently induces an embedding
\[\Hom_{G}(V_{x,I_0}^{\rm{an}},V_{w,I_0}^{\rm{an}})\hookrightarrow \Hom_{G}(i_{x,I_0}^{\rm{an}},V_{w,I_0}^{\rm{an}})\]
and we know the existence of a non-zero map $V_{x,I_0}^{\rm{an}}\rightarrow V_{w,I_0}^{\rm{an}}$ from (\ref{equ: St x w transfer}) and the discussion below it, we conclude (\ref{equ: Hom x w St}).
\end{proof}

For each $x\in W(G)$ and $I\subseteq \Delta$, we write $x^I\in W^{I,\emptyset}$ for the unique minimal element in $W(L_I)x$.
\begin{lem}\label{lem: vanishing by g}
Let $x_i\in W(G)$, $I_i\subseteq I_{x_i}$ and $\pi_i^{\infty}\in\cB^{I_i}$ for $i=0,1$. If $x_0^{I_1}\not\leq x_1$, then we have
\begin{equation}\label{equ: vanishing by g}
\mathrm{Ext}_{G}^k(\cF_{P_{I_0}}^{G}(L(x_0),\pi_0^{\infty}), \cF_{P_{I_1}}^{G}(M^{I_1}(x_1),\pi_1^{\infty}))=0
\end{equation}
for each $k\geq 0$.
\end{lem}
\begin{proof}
By applying Proposition~\ref{prop: St key seq} and Remark~\ref{rem: trivial block adjunction} with $M_0=L(x_0)$, $M_1=M^{I_1}(x_1)$ and then $I_0$, $I_1$, $\pi_0^{\infty}$, $\pi_1^{\infty}$ as in \emph{loc.cit.}, we obtain the spectral sequence
\begin{equation}\label{equ: vanishing by g seq}
\mathrm{Ext}_{U(\fg)}^{\ell}(M^{I_1}(x_1), L(x_0))\otimes_E \mathrm{Ext}_{G}^k(i_{I_0}^{\infty}(\pi_0^{\infty}), i_{I_1}^{\infty}(\pi_1^{\infty}))^{\infty}
\implies \mathrm{Ext}_{G}^{k+\ell}(\cF_{P_{I_0}}^{G}(L(x_0),\pi_0^{\infty}),\cF_{P_{I_1}}^{G}(M^{I_1}(x_1),\pi_1^{\infty})).
\end{equation}
Hence, it suffices to show that
\begin{equation}\label{equ: g vanishing by g}
\mathrm{Ext}_{U(\fg)}^{\ell}(M^{I_1}(x_1), L(x_0))=0
\end{equation}
for each $\ell\geq 0$, which follows from Lemma~\ref{lem: Ext Verma w} upon replacing $x$, $w$ and $I$ in \emph{loc.cit.} with $x_1$, $x_0$ and $I_1$ here.
\end{proof}

\begin{lem}\label{lem: factor PS vanishing}
Let $x_0,x_1\in W(G)$, $I_1\subseteq I_{x_1}$ and $I\subseteq I_0\subseteq I_{x_0}$. If $x_0\not\leq x_1$, then we have
\begin{equation}\label{equ: factor PS vanishing}
\mathrm{Ext}_{G}^k(\cF_{P_{I_0}}^{G}(L(x_0),V_{I,I_0}^{\infty}),i_{x_1,I_1}^{\rm{an}})=0
\end{equation}
for $k\geq 0$.
\end{lem}
\begin{proof}
If $I_1\not\subseteq I_0$, then (\ref{equ: factor PS vanishing}) holds for $k\geq 0$ by applying Lemma~\ref{lem: vanishing by sm} with $M_0=L(x_0)$. If $I_1\subseteq I_0\subseteq I_{x_0}$, then $x_0^{I_1}=x_0\not\leq x_1$ and thus (\ref{equ: factor PS vanishing}) holds for $k\geq 0$ by applying Lemma~\ref{lem: vanishing by g} with $\pi_0^{\infty}=V_{I,I_0}^{\infty}$ and $\pi_1^{\infty}=1_{L_{I_1}}$.
\end{proof}

Fir each $x\in W(G)$ and $I_0\subseteq I_1\subseteq I_x$, we recall the complex $\mathbf{C}^{x}_{I_0,I_1}$ from (\ref{equ: x Tits complex}) with $I=\Delta$ in \emph{loc.cit.}
\begin{lem}\label{lem: x loc an St vanishing}
Let $x,w\in W(G)$ with $w\not\leq x$. Let $I_0\subseteq I_1\subseteq I_x$ and $I\subseteq I_w$.
\begin{enumerate}[label=(\roman*)]
\item \label{it: x loc an St vanishing 1} We have
\begin{equation}\label{equ: x complex vanishing}
\mathrm{Ext}_{G}^k(C_{w,I},\mathbf{C}^{x}_{I_0,I_1})=0
\end{equation}
for each $k\in\Z$.
\item \label{it: x loc an St vanishing 2} The representation $i_{I_1}^{\rm{an}}(V_{x,I_0,I_1}^{\rm{an}})$ has simple socle $C_{x,I_0\sqcup(I_x\setminus I_1)}$.
\end{enumerate}
\end{lem}
\begin{proof}
We prove \ref{it: x loc an St vanishing 1}. We apply Lemma~\ref{lem: factor PS vanishing} (by replacing $x_0$, $x_1$, $I_0$, $I_1$ with $w$, $x$, $I_w$ and $I'$) and deduce that $\mathrm{Ext}_{G}^k(C_{w,I},i_{x,I'}^{\rm{an}})=0$ for each $k\in\Z$ and $I'\subseteq I_x$, which together with the definition of $\mathbf{C}^{x}_{I_0,I_1}$ gives (\ref{equ: x complex vanishing}) for each $k\in\Z$.

We prove \ref{it: x loc an St vanishing 2}. For each $I_0\subseteq I'\subseteq I_1$, by \ref{it: OS property 1} of Proposition~\ref{prop: OS property} we know that $i_{x,I',I_1}^{\rm{an}}=\cF_{L_{I_1}\cap P_{I'}}^{L_{I_1}}(U(\fl_{I_1})\otimes_{U(\fl_{I_1}\cap\fp_{I'})}L^{I'}(x),1_{L_{I'}})$ contains the subrepresentation
$\cF_{L_{I_1}\cap P_{I'}}^{L_{I_1}}(L^{I_1}(x),1_{L_{I'}})$ that exhausts its $x$-typic constituents.
By \ref{it: sm cube 3} of Lemma~\ref{lem: sm cube} we have $i_{I',I_x\cap I_1}^{\infty}(1_{L_{I'}})\cong Q_{I_x\cap I_1}(I_x\cap I_1,I')$ and thus
\begin{equation}\label{equ: Levi JH mult x}
[\cF_{L_{I_1}\cap P_{I'}}^{L_{I_1}}(L^{I_1}(x),1_{L_{I'}})]=\sum_{I'\subseteq I''\subseteq I_x\cap I_1}[C_{x,I'',I_1}]
\end{equation}
It follows from the isomorphism $V_{x,I_0,I_1}^{\rm{an}}\cong \mathbf{C}^{x,I_1}_{I_0,I_1}[-\#I_0]$ (see (\ref{equ: general Tits resolution}) with $I=I_1$) that
\[[V_{x,I_0,I_1}^{\rm{an}}]=\sum_{I_0\subseteq I'\subseteq I_1}(-1)^{\#I'\setminus I_1}[i_{x,I',I_1}^{\rm{an}}],\]
which together with (\ref{equ: Levi JH mult x}) implies that $C_{x,I_0,I_1}$ is the only $x$-typic constituent of $V_{x,I_0,I_1}^{\rm{an}}$ and it appears with multiplicity $1$. As $[U(\fl_{I_1})\otimes_{U(\fl_{I_1}\cap\fp_{I'})}L^{I'}(x):L^{I_1}(x')]\neq 0$ for some $x'\neq x$ only if $x'>x$, by Lemma~\ref{lem: JH OS} we know that any constituent $W$ of $i_{x,I_0,I_1}^{\rm{an}}$ which is not $x$-typic is $x'$-typic for some $x'>x$. Hence, by applying Lemma~\ref{lem: OS vers induction} to $i_{I_1}^{\rm{an}}(W)$ for each $W\in\mathrm{JH}_{L_{I_1}}(V_{x,I_0,I_1}^{\rm{an}})$, we see that each constituent of $i_{I_1}^{\rm{an}}(V_{x,I_0,I_1}^{\rm{an}}/C_{x,I_0,I_1})$ are $x'$-typic for some $x'>x$.
As we have
\[i_{I_1}^{\rm{an}}(C_{x,I_0,I_1})=i_{I_1}^{\rm{an}}(\cF_{P_{I_x}\cap L_{I_1}}^{L_{I_1}}(L^{I_1}(x),V_{I_0,I_1\cap I_x}^{\infty}))=\cF_{P_{I_1\cap I_x}}^{G}(M^{I_1}(x),V_{I_0,I_1\cap I_x}^{\infty})\]
by Lemma~\ref{lem: OS vers induction}, we see that each constituent of $i_{I_1}^{\rm{an}}(C_{x,I_0,I_1})/V_0$ is $x'$-typic for some $x'>x$, where
\begin{equation}\label{equ: x typic sub}
V_0\defeq \cF_{P_{I_1\cap I_x}}^{G}(L(x),V_{I_0,I_1\cap I_x}^{\infty})\cong \cF_{P_{I_x}}^{G}(L(x),i_{I_1\cap I_x,I_x}^{\infty}(V_{I_0,I_1\cap I_x}^{\infty})).
\end{equation}
As $i_{I_1\cap I_x,I_x}^{\infty}(V_{I_0,I_1\cap I_x}^{\infty})\cong Q_{I_x}(I_0\sqcup(I_x\setminus I_1),I_0)$ has socle $V_{I_0\sqcup(I_x\setminus I_1),I_x}^{\infty}$ by Lemma~\ref{lem: sm cube}, we deduce from Lemma~\ref{lem: OS layer} that (\ref{equ: x typic sub}) has socle $C_{x,I_0\sqcup(I_x\setminus I_1)}$. As any constituent $V$ of $i_{I_1}^{\rm{an}}(V_{x,I_0,I_1}^{\rm{an}})/V_0$ is necessarily $x'$-typic for some $x'>x$ by previous discussion, we deduce $\Hom_{G}(V, i_{I_1}^{\rm{an}}(V_{x,I_0,I_1}^{\rm{an}}))=0$ from \ref{it: x loc an St vanishing 1}, and in particular $V$ does not appear in the socle of $i_{I_1}^{\rm{an}}(V_{x,I_0,I_1}^{\rm{an}})$. Hence, we must have
\[\mathrm{soc}_{G}(i_{I_1}^{\rm{an}}(V_{x,I_0,I_1}^{\rm{an}}))=\mathrm{soc}_{G}(V_0)=C_{x,I_0\sqcup(I_x\setminus I_1)}.\]
The proof is thus finished.
\end{proof}

\begin{lem}\label{lem: change left cup}
Let $I_0,I_1\subseteq \Delta$. We have the following results.
\begin{enumerate}[label=(\roman*)]
\item \label{it: change left cup 1} We have $\mathrm{Ext}_{G}^k(W,V_{I_1}^{\rm{an}})=0$ for each $W\in\mathrm{JH}_{G}(V_{I_0}^{\rm{an}}/V_{I_0}^{\infty})$ and $k\geq 0$.
\item \label{it: change left cup 2} Let $V'\subseteq V\subseteq V_{I_0}^{\rm{an}}$ be subrepresentations. Then we have $\Dim_E \Hom_{G}(V',V)=1$, and the following cup product map
\begin{equation}\label{equ: change left cup}
\Hom_{G}(V',V)\otimes_E\mathrm{Ext}_{G}^k(V,V_{I_1}^{\rm{an}})\buildrel\cup\over\longrightarrow \mathrm{Ext}_{G}^k(V',V_{I_1}^{\rm{an}})
\end{equation}
is an isomorphism for each $k\geq 0$.
\end{enumerate}
\end{lem}
\begin{proof}
We have $[V_{I_0}^{\rm{an}}:V_{I_0}^{\infty}]=1$ from Lemma~\ref{lem: JH mult PS}. It also follows from \ref{it: x loc an St vanishing 2} of Lemma~\ref{lem: x loc an St vanishing} (by replacing $x$, $I_0$, $I_1$ in \emph{loc.cit.} with $1$, $\Delta$ and $I_0$) that $\mathrm{soc}_{G}(V_{I_0}^{\rm{an}})\cong V_{I_0}^{\infty}$. These together force $[V:V_{I_0}^{\infty}]=1=[V':V_{I_0}^{\infty}]$ and $\mathrm{soc}_{G}(V)\cong V_{I_0}^{\infty}\cong \mathrm{soc}_{G}(V')$. In particular, we have $[V'/V_{I_0}^{\infty}:V_{I_0}^{\infty}]=0$ and $\Hom_{G}(V'/V_{I_0}^{\infty},V)=0$, which together with a d\'evissage with respect to $0\rightarrow V_{I_0}^{\infty}\rightarrow V'\rightarrow V'/V_{I_0}^{\infty}\rightarrow 0$ gives an embedding
\begin{equation}\label{equ: change left cup embedding}
\Hom_{G}(V',V)\hookrightarrow \Hom_{G}(V_{I_0}^{\infty},V).
\end{equation}
As we clearly have $\Hom_{G}(V',V)\neq 0$ and $\Dim_E \Hom_{G}(V_{I_0}^{\infty},V)=1$, the embedding (\ref{equ: change left cup embedding}) is an isomorphism between $1$ dimensional spaces.
In particular, all non-zero maps $V'\rightarrow V$ are embeddings and differ from each other by a scalar. To prove that the cup product map (\ref{equ: change left cup}) is an isomorphism, using a d\'evissage with respect to $0\rightarrow V'\rightarrow V\rightarrow V/V'\rightarrow 0$, it suffices to show that
\begin{equation}\label{equ: change left cup vanishing}
\mathrm{Ext}_{G}^k(V/V',V_{I_1}^{\rm{an}})=0
\end{equation}
for $k\geq 0$. But since $V_{I_0}^{\infty}$ is the only smooth (or equivalently $1$-typic) constituent of $V_{I_0}^{\rm{an}}$ (see Lemma~\ref{lem: JH mult PS}) and $\mathrm{soc}_{G}(V)\cong V_{I_0}^{\infty}\cong \mathrm{soc}_{G}(V')$, we deduce each constituent of $V/V'$ has the form $C_{w,I}$ for some $w>1$ and $I\subseteq I_w$. We thus deduce from \ref{it: x loc an St vanishing 1} of Lemma~\ref{lem: x loc an St vanishing} (by taking $x$, $I_0$, $I_1$ in \emph{loc.cit.} to be $1$, $I_1$ and $I_1$) that
\[\mathrm{Ext}_{G}^k(C_{w,I},V_{I_1}^{\rm{an}})=0\]
for $k\geq 0$, which proves \ref{it: change left cup 1}. This together with a d\'evissage with respect to $\mathrm{JH}_{G}(V/V')$ gives (\ref{equ: change left cup vanishing}) which further implies \ref{it: change left cup 2}.
\end{proof}

\section{Preliminary on extensions between Tits complex}\label{sec: Ext PS}
We introduce what we call \emph{Tits double complex} to compute the cohomology of a Tits complex. Then we discuss how to reduce the computation of extensions between two Tits complex to the cohomology of a single Tits complex (see Proposition~\ref{prop: Ext complex std seq}).
\subsection{Tits construction}\label{subsec: Tits construction}
We write $\cP(\Delta)$ for the power set of $\Delta$.
We view $\cP(\Delta)$ as a category with morphisms given by inclusion between subsets $I'\subseteq I$ of $\Delta$.
Let $\cC$ be an additive category in which finite direct sum exists and $\cT: \cP(\Delta)^{\rm{op}}\rightarrow \cC$ be a functor. Concretely, a such functor $\cT$ is given by the following data
\begin{itemize}
\item an object $M_{I}\defeq\cT(I)\in \cC$ for each $I\subseteq\Delta$;
\item a morphism $\iota_{I,I'}\defeq \cT(I'\subseteq I): M_{I}\rightarrow M_{I'}$ in $\cC$ for each pair $I'\subseteq I$, such that $\iota_{I,I''}=\iota_{I',I''}\circ\iota_{I,I'}$ for each triple $I''\subseteq I'\subseteq I$.
\end{itemize}
We identify $\cT$ with a bounded $\cC$-valued multi-complex indexed by $\cP(\Delta)=\{0,1\}^{\#\Delta}$. More precisely, for each $j\in\Delta$, the $j$-degree of $M_{I}$ is $0$ (resp.~is $1$) if and only if $j\in I$ (resp.~$j\notin I$), with differential maps in $j$-th direction given by $\iota_{I,I\setminus\{j\}}: M_{I}\rightarrow M_{I\setminus\{j\}}$ if $j\in I$ and given by zero otherwise.

Given a tuple of signs $\un{c}=(c_{I})_{I\subseteq\Delta}$ in $\{-1,1\}$, we define the twist of $\cT$ by $\un{c}$, which is another functor written $\cT[\un{c}]: \cP(\Delta)^{\rm{op}}\rightarrow \cC$, by exactly the same collection of objects $M_{I}$ in $\cC$ for all $I\subseteq\Delta$, but with the morphisms $\iota_{I,I'}$ replaced with $c_{I}^{-1}c_{I'}\iota_{I,I'}$ for each pair $I'\subseteq I$ of subsets of $\Delta$. In particular, the scalar isomorphism $c_{I}: M_{I}\rightarrow M_{I}$ for each $I\subseteq\Delta$ induces a natural transform 
\begin{equation}\label{equ: twist transform}
\un{c}: \cT\rightarrow \cT[\un{c}].
\end{equation}
It is evident that the complex $\cT[\un{c}]$ is unchanged if we multiply each $c_{I}$ by a common sign in $\{-1,1\}$.

Recall from the discussion around (\ref{equ: general total differential}) that we can associate with $\cT$ its total complex $\mathrm{Tot}(\cT)$ with the differential map $\mathrm{Tot}(\cT)^{m}\rightarrow \mathrm{Tot}(\cT)^{m+1}$ given by
\[
\bigoplus_{I\subseteq\Delta, \#\Delta\setminus I=m, j\in I}(-1)^{m'(I,j)}\iota_{I,I\setminus\{j\}}
\]
where $m'(I,j)\defeq \#\{i\in\Delta\setminus I\mid i<j\}$ for each $I\subseteq\Delta$ and $j\in I$.
We write $\cT[1]\defeq\cT[1,\dots,1]$ for short. We can also associate with $\cT[1]$ its total complex $\mathrm{Tot}(\cT[1])$ with the differential map $\mathrm{Tot}(\cT[1])^{m}\rightarrow \mathrm{Tot}(\cT[1])^{m+1}$ given by
\[
\bigoplus_{I\subseteq\Delta, \#I=-m, j\in I}(-1)^{m(I,j)}\iota_{I,I\setminus\{j\}}
\]
where $m(I,j)\defeq \#\{i\in I\mid i<j\}$ for each $I\subseteq\Delta$ and $j\in I$ (see also (\ref{equ: differential sign})). It is clear that we have $m'(I,j)+m(I,j)=\#\{i\in\Delta\mid i<j\}=j-1$ for each $I\subseteq\Delta$ and $j\in I$.
We consider the tuple of signs $\un{c}=(c_{I})_{I\subseteq\Delta}$ defined by $c_{I}\defeq (-1)^{\sum_{j\in I}(j-1)}$ for each $I\subseteq\Delta$. It is clear that 
\[c_{I}^{-1}c_{I\setminus\{j\}}=(-1)^{j-1}=(-1)^{m(I,j)+m'(I,j)}\]
for each $I\subseteq\Delta$ and $j\in I$, and thus we have a natural identification $\mathrm{Tot}(\cT)=\mathrm{Tot}(\cT[1][\un{c}])[-\#\Delta]$ between complex, with the scalar automorphism $c_{I}: M_{I}\rightarrow M_{I}$ inducing the following isomorphism between complex
\begin{equation}\label{equ: Tits total convention transfer}
\mathrm{Tot}(\cT[1])[-\#\Delta]\buildrel\sim\over\longrightarrow \mathrm{Tot}(\cT[1][\un{c}])[-\#\Delta]=\mathrm{Tot}(\cT)
\end{equation}
and thus the following isomorphism between cohomolgy
\begin{equation}\label{equ: Tits total convention coh transfer}
H^{k-\#\Delta}(\mathrm{Tot}(\cT[1]))\buildrel\sim\over\longrightarrow H^{k}(\mathrm{Tot}(\cT[1][\un{c}])[-\#\Delta])=H^{k}(\mathrm{Tot}(\cT))
\end{equation}
for each $k\in \Z$.

As the first interesting example of $\cC$ and $\cT$, we can recover our construction of the Tits complex.
Assume for the moment that $\cC=\mathrm{Rep}^{\rm{an}}_{\rm{adm}}(G)$. Given $I_0\subseteq I_1\subseteq\Delta$, we consider $\cT$ given by $M_{I}=i_{I}^{\rm{an}}$ if $I_0\subseteq I\subseteq I_1$ and zero otherwise, with $\iota_{I,I'}$ being the natural inclusion $i_{I}^{\rm{an}}\hookrightarrow i_{I'}^{\rm{an}}$ of $I_0\subseteq I'\subseteq I\subseteq I_1$ and being zero otherwise. Then we have
\[\mathbf{C}_{I_0,I_1}=\mathrm{Tot}(\cT[1])\in\mathrm{Ch}^{b}(\mathrm{Rep}^{\rm{an}}_{\rm{adm}}(G)).\]

Now we move on to some more interesting examples of $\cT$ which are closely related to the computation of $\mathrm{Ext}_{G}^{\bullet}(1_{G},\mathbf{C}_{I_0,I_1})$.
Assume from now that $\cC$ is the additive category of complex of $E$-vector spaces.
The choice of a functor $\cT$ amounts to the choice of a collection $\Sigma=\{M^{\bullet}_{I}\}_{I\subseteq\Delta}$ of complex $E$-vector spaces as well as maps between complex $\iota_{I,I'}: M^{\bullet}_{I}\rightarrow M^{\bullet}_{I'}$ satisfying $\iota_{I,I''}=\iota_{I',I''}\circ\iota_{I,I'}$ for each $I''\subseteq I'\subseteq I$. For each $k\in\Z$, we can associate with $\cT$ the functor $\cT^{k}$ from $\cP(\Delta)^{\rm{op}}$ to the category of $E$-vector spaces by considering the collection $\{M^{k}_{I}\}_{I\subseteq\Delta}$ and the morphisms $\iota^{k}_{I,I'}: M^{k}_{I}\rightarrow M^{k}_{I'}$ for each $I'\subseteq I$.
We define a double complex $\cT_{\Sigma}^{\bullet,\bullet}$, called \emph{abstract Tits double complex}, by
\begin{equation}\label{equ: general double complex}
\cT_{\Sigma}^{\bullet,k}\defeq \mathrm{Tot}(\cT^{k}[1])
\end{equation}
for each $k\in\Z$, with the differential map $\cT_{\Sigma}^{\bullet,k}\rightarrow \cT_{\Sigma}^{\bullet,k+1}$ being the evident one induced from the differential map $M^{k}_{I}\rightarrow M^{k+1}_{I}$ for each $I\subseteq\Delta$.
We can define another double complex $\cN\cT_{\Sigma}^{\bullet,\bullet}$, called \emph{normalized abstract Tits double complex}, \begin{equation}\label{equ: general double complex normalized}
\cN\cT_{\Sigma}^{\bullet,k}\defeq \mathrm{Tot}(\cT^{k})
\end{equation}
for each $k\in\Z$, with the differential map $\cN\cT_{\Sigma}^{\bullet,k}\rightarrow \cN\cT_{\Sigma}^{\bullet,k+1}$ being the evident one induced from the differential map $M^{k}_{I}\rightarrow M^{k+1}_{I}$ for each $I\subseteq\Delta$.
It follows from (\ref{equ: Tits total convention transfer}) and (\ref{equ: Tits total convention coh transfer}) that we have a standard isomorphism between double complex
\begin{equation}\label{equ: Tits double complex normalize}
\cT_{\Sigma}^{\bullet,\bullet}[-\#\Delta,0]\buildrel\sim\over\longrightarrow \cN\cT_{\Sigma}^{\bullet,\bullet}
\end{equation}
which induces the following isomorphism
\begin{equation}\label{equ: Tits double complex normalize coh}
H^{k-\#\Delta}(\mathrm{Tot}(\cT_{\Sigma}^{\bullet,\bullet}))\buildrel\sim\over\longrightarrow H^{k}(\mathrm{Tot}(\cN\cT_{\Sigma}^{\bullet,\bullet}))
\end{equation}
for each $k\in\Z$.

The following notions will be useful later.
\begin{defn}\label{def: permutation preserve order}
Let $\sigma\in\mathfrak{S}_{\Delta}$ be a permutation. We have the following notions.
\begin{enumerate}[label=(\roman*)]
\item \label{it: permutation 1} Given a subset $I\subseteq\Delta$, we say that $\sigma$ is \emph{compatible} with $I$ if $\sigma^{-1}(i)<\sigma^{-1}(j)$ for each $i,j\in I$ satisfying $i<j$.
\item \label{it: permutation 2} Given a ordered pair of subsets $I,I'\subseteq\Delta$ with $I\cap I'=\emptyset$, we say that $\sigma$ is \emph{well-posed} for the pair $I,I'$ if $\sigma^{-1}(j)<\sigma^{-1}(j')$ for each $j\in I$ and $j'\in I'$.
\end{enumerate}
\end{defn}

Let $I_0\subseteq I_1\subseteq\Delta$ and assume from now on that $M^{\bullet}_{I}\neq 0$ only if $I_0\subseteq I\subseteq I_1$.
The truncation on the first index of $\cT_{\Sigma}^{\bullet,\bullet}$ induces a spectral sequence $E_{\bullet,\Sigma}^{\bullet,\bullet}$ that converges to $H^{\bullet}(\mathrm{Tot}(\cT_{\Sigma}^{\bullet,\bullet}))$ with
\begin{equation}\label{equ: abstract E1 sum}
E_{1,\Sigma}^{-\ell,k}=\bigoplus_{I_0\subseteq I\subseteq I_1, \#I=\ell}H^k(M^{\bullet}_{I})
\end{equation}
with $d_{1,\Sigma}^{-\ell,k}$ given by the direct sum of 
\[(-1)^{m(I,j)}H^{k}(\iota^{\bullet}_{I,I'}): H^k(M^{\bullet}_{I})\rightarrow H^k(M^{\bullet}_{I\setminus\{j\}})\]
for each $I_0\subseteq I\subseteq I_1$ with $\#I=\ell$ and $j\in I$.
We say that $E_{\bullet,\Sigma}^{\bullet,\bullet}$ \emph{degenerates at the second page for $(-\ell,k)$} if $E_{2,\Sigma}^{-\ell,k}=E_{\infty,\Sigma}^{-\ell,k}$.

We consider two collections of complex $\Sigma=\{M^{\bullet}_{I}\}_{I_0\subseteq I\subseteq I_1}$ and $\Sigma'=\{N^{\bullet}_{I}\}_{I_0\subseteq I\subseteq I_1}$.
Given maps between complex $M^{\bullet}_{I}\rightarrow N^{\bullet}_{I}$ which are functorial with respect to the choice of $I_0\subseteq I\subseteq I_1$, we can naturally associate with them a map between double complex
\begin{equation}\label{equ: abstract map between double}
\cT_{\Sigma}^{\bullet,\bullet}\rightarrow \cT_{\Sigma'}^{\bullet,\bullet}
\end{equation}
which further induces a map between spectral sequences
\begin{equation}\label{equ: abstract map between double spectral seq}
E_{\bullet,\Sigma}^{\bullet,\bullet}\rightarrow E_{\bullet,\Sigma'}^{\bullet,\bullet}.
\end{equation}
We say that (\ref{equ: abstract map between double}) induces an \emph{isomorphism on the first page} if (\ref{equ: abstract map between double spectral seq}) is an isomorphism on the first page, in which case the map between total complex associated with (\ref{equ: abstract map between double}) is necessarily a quasi-isomorphism.

More generally, we write 
\begin{equation}\label{equ: abstract double quasi map}
\cT_{\Sigma}^{\bullet,\bullet}\dashrightarrow\cT_{\Sigma'}^{\bullet,\bullet},
\end{equation}
and call it a \emph{quasi map} between double complex, if there exists $t\geq 1$ and collections $\Sigma=\Sigma_0,\dots,\Sigma_{2t}=\Sigma'$ as well as maps between double complex
\[\cT_{\Sigma}^{\bullet,\bullet}\leftarrow \cT_{\Sigma_1}^{\bullet,\bullet}\rightarrow\cdots \rightarrow\cT_{\Sigma_{2t'-2}}^{\bullet,\bullet}\leftarrow \cT_{\Sigma_{2t'-1}}^{\bullet,\bullet}\rightarrow \cT_{\Sigma_{2t'}}^{\bullet,\bullet}\leftarrow \cdots \leftarrow \cT_{\Sigma_{2t-1}}^{\bullet,\bullet}\rightarrow \cT_{\Sigma'}^{\bullet,\bullet}\]
such that $\cT_{\Sigma_{2t'}}^{\bullet,\bullet}\rightarrow \cT_{\Sigma_{2t'-1}}^{\bullet,\bullet}$ is an isomorphism on the first page for each $1\leq t'\leq t$.
In this case, we obtain maps
\[E_{r,\Sigma}^{\bullet,\bullet}\buildrel\sim\over\longleftarrow E_{r,\Sigma_1}^{\bullet,\bullet}\rightarrow\cdots \rightarrow E_{1,\Sigma_{2t'-2}}^{\bullet,\bullet}\buildrel\sim\over\longleftarrow E_{r,\Sigma_{2t'-1}}^{\bullet,\bullet}\rightarrow E_{r,\Sigma_{2t'}}^{\bullet,\bullet}\buildrel\sim\over\longleftarrow \cdots \buildrel\sim\over\longleftarrow E_{r,\Sigma_{2t-1}}^{\bullet,\bullet}\rightarrow E_{r,\Sigma'}^{\bullet,\bullet}\]
and thus a map
\[E_{r,\Sigma}^{\bullet,\bullet}\rightarrow E_{r,\Sigma'}^{\bullet,\bullet}\]
for each $r\geq 1$.
Note that the map $H^{h}(\mathrm{Tot}(\cT_{\Sigma_{2t'-1}}^{\bullet,\bullet}))\rightarrow H^{h}(\mathrm{Tot}(\cT_{\Sigma_{2t'-2}}^{\bullet,\bullet}))$ is an isomorphism for each $1\leq t'\leq t$ and $h\in\Z$. Hence, we also have a map
\[H^{h}(\mathrm{Tot}(\cT_{\Sigma}^{\bullet,\bullet}))\rightarrow H^{h}(\mathrm{Tot}(\cT_{\Sigma'}^{\bullet,\bullet}))\]
which is compatible with the canonical filtration on both sides, whose induced map between $\mathrm{gr}^{-\ell}(\cdot)$ recovers
\[E_{\infty,\Sigma}^{-\ell,\ell+h}\rightarrow E_{\infty,\Sigma'}^{-\ell,\ell+h}.\]

We give below some abstract criterion towards the partial-degeneracy of the spectral sequence at the second page (see Lemma~\ref{lem: abstract E2 degenerate}, Lemma~\ref{lem: abstract E2 injection} and Lemma~\ref{lem: abstract E2 seq}).
\begin{lem}\label{lem: abstract E2 degenerate}
Let $\Sigma=\{M^{\bullet}_{I}\}_{I_0\subseteq I\subseteq I_1}$ be a collection of complex as above to which we associate the double complex $\cT_{\Sigma}^{\bullet,\bullet}$ and the spectral sequence $E_{\bullet,\Sigma}^{\bullet,\bullet}$. Let $\#I_0\leq \ell\leq \#I_1$, $k\geq 0$ and $\Psi$ be a basis of $E_{2,\Sigma}^{-\ell,k}$. Assume that
\begin{enumerate}[label=(\roman*)]
\item \label{it: abstract E2 degenerate 1} we have $d_{r,\Sigma}^{-\ell-r,k+(r-1)}=0$ for each $r\geq 2$;
\item \label{it: abstract E2 degenerate 2} for each $x\in\Psi$, there exists a collection of complex $\Sigma_{x}$ as well as a quasi map between double complex $\cT_{\Sigma_{x}}^{\bullet,\bullet}\dashrightarrow \cT_{\Sigma}^{\bullet,\bullet}$ such that $E_{1,\Sigma_{x}}^{-\ell',k'}=0$ whenever $k'-\ell'>k-\ell$ and that the map $E_{2,\Sigma_{x}}^{-\ell,k}\rightarrow E_{2,\Sigma}^{-\ell,k}$ has image $Ex$.
\end{enumerate}
Then we have
\begin{equation}\label{equ: abstract E2 degenerate}
E_{2,\Sigma}^{-\ell,k}=E_{\infty,\Sigma}^{-\ell,k}.
\end{equation}
\end{lem}
\begin{proof}
Thanks to \ref{it: abstract E2 degenerate 1}, we have inclusions
\[E_{\infty,\Sigma}^{-\ell,k}\subseteq \cdots\subseteq E_{r,\Sigma}^{-\ell,k}\subseteq\cdots\subseteq E_{2,\Sigma}^{-\ell,k}\]
for $r\geq 2$. Now that $\Psi$ is a basis of $E_{2,\Sigma}^{-\ell,k}$, it suffices to show that
\begin{equation}\label{equ: abstract r differential vanishing}
d_{r,\Sigma}^{-\ell,k}(x)=0
\end{equation}
for each $x\in\Psi$ and $r\geq 2$. Fix from now a choice of $x\in\Psi$ to which we associate $\Sigma_{x}$ as in \ref{it: abstract E2 degenerate 2}. Now that $E_{1,\Sigma_{x}}^{-\ell',k'}=0$ whenever $k'-\ell'>k-\ell$, we have $d_{r,\Sigma_{x}}^{-\ell,k}=0$ for each $r\geq 1$.
We choose an arbitrary $y\in E_{2,\Sigma_{x}}^{-\ell,k}$ whose image in $E_{2,\Sigma}^{-\ell,k}$ is $x$.
For each $r\geq 2$, we abuse $y$ for its image in $E_{r,\Sigma_{x}}^{-\ell,k}$, and then deduce (\ref{equ: abstract r differential vanishing}) from the vanishing $d_{r,\Sigma_{x}}^{-\ell,k}(y)=0$.
The proof is thus finished.
\end{proof}

\begin{lem}\label{lem: abstract E2 injection}
Let $\Sigma$ and $\Sigma'$ be collections of complex with a quasi map between double complex $\cT_{\Sigma'}^{\bullet,\bullet}\dashrightarrow \cT_{\Sigma}^{\bullet,\bullet}$. Let $h\in\Z$, $\#I_0\leq \ell\leq \#I_1$ and $k\geq 0$. Assume that
\begin{enumerate}[label=(\roman*)]
\item \label{it: abstract E2 injection 1} $d_{r,\Sigma}^{-\ell,k}=0$ for each $(-\ell,k)$ with $k-\ell=h$;
\item \label{it: abstract E2 injection 2} the map $E_{2,\Sigma'}^{-\ell,k}\rightarrow E_{2,\Sigma}^{-\ell,k}$ is injective for each $(-\ell,k)$ with $k-\ell\leq h+1$;
\item \label{it: abstract E2 injection 3} $E_{2,\Sigma}^{-\ell,k}=0$ for each $(-\ell,k)$ with $k-\ell<h$.
\end{enumerate}
Then we have $E_{2,\Sigma'}^{-\ell,k}=0$ for each $(-\ell,k)$ with $k-\ell<h$, and $E_{2,\Sigma'}^{-\ell,k}=E_{\infty,\Sigma'}^{-\ell,k}$ for each $(-\ell,k)$ with $k-\ell=h$.
In particular, the map
\begin{equation}\label{equ: abstract E2 injection total}
H^{h}(\mathrm{Tot}(\cT_{\Sigma'}^{\bullet,\bullet}))\rightarrow H^{h}(\mathrm{Tot}(\cT_{\Sigma}^{\bullet,\bullet}))
\end{equation}
is injective and strict with respect to the canonical filtration on both terms, whose $\mathrm{gr}^{-\ell}(\cdot)$ recovers the injection $E_{\infty,\Sigma'}^{-\ell,\ell+h}\hookrightarrow E_{\infty,\Sigma}^{-\ell,\ell+h}$ for each $\#I_0\leq \ell\leq \#I_1$.
\end{lem}
\begin{proof}
We fix throughout the proof an arbitrary choice of $(-\ell_0,k_0)$ with $k_0-\ell_0=h$.
It is clear from \ref{it: abstract E2 injection 2} and \ref{it: abstract E2 injection 3} that $E_{2,\Sigma'}^{-\ell,k}=0$ for each $(-\ell,k)$ with $k-\ell<h$. In particular, we have $E_{2,\Sigma'}^{-\ell_0-r,k_0+(r-1)}=0$ and thus
\begin{equation}\label{equ: abstract E2 injection 1}
d_{r,\Sigma'}^{-\ell_0-r,k_0+(r-1)}=0=d_{r,\Sigma}^{-\ell_0-r,k_0+(r-1)}
\end{equation}
for each $r\geq 2$.
We prove by an increasing induction on $r\geq 2$ that
\begin{equation}\label{equ: abstract E2 injection 2}
d_{r,\Sigma'}^{-\ell_0,k_0}=0.
\end{equation}
We consider the following commutative diagram
\begin{equation}\label{equ: abstract E2 injection 3}
\xymatrix{
E_{r,\Sigma'}^{-\ell_0,k_0} \ar^{d_{r,\Sigma'}^{-\ell_0,k_0}}[r] \ar[d] & E_{r,\Sigma'}^{-\ell_0+r,k_0-(r-1)} \ar[d]\\
E_{r,\Sigma}^{-\ell_0,k_0} \ar^{d_{r,\Sigma}^{-\ell_0,k_0}}[r] & E_{r,\Sigma}^{-\ell_0+r,k_0-(r-1)}
}
\end{equation}
with the bottom horizontal map $d_{r,\Sigma}^{-\ell_0,k_0}$ known to be zero by \ref{it: abstract E2 injection 1}.
When $r=2$, both vertical maps of (\ref{equ: abstract E2 injection 3}) are injective by \ref{it: abstract E2 injection 2}, which forces $d_{2,\Sigma'}^{-\ell_0,k_0}=0$ and gives (\ref{equ: abstract E2 injection 2}) for $r=2$.
Assume from now that $r>2$.
Since $d_{r',\Sigma}^{-\ell_0,k_0}=0$ for each $r'\geq 2$ by \ref{it: abstract E2 injection 1}, we have inclusions
\[E_{\infty,\Sigma}^{-\ell_0+r,k_0-(r-1)}\subseteq \cdots\subseteq E_{r,\Sigma}^{-\ell_0+r,k_0-(r-1)}\subseteq \cdots\subseteq E_{2,\Sigma}^{-\ell_0+r,k_0+(r-1)}.\]
By our inductive hypothesis we also have $d_{r',\Sigma'}^{-\ell_0,k_0}=0$ for each $2\leq r'\leq r-1$ which gives inclusions
\[E_{r,\Sigma'}^{-\ell_0+r,k_0-(r-1)}\subseteq \cdots\subseteq E_{2,\Sigma'}^{-\ell_0+r,k_0-(r-1)}.\]
In particular, the injection $E_{2,\Sigma'}^{-\ell_0+r,k_0-(r-1)}\hookrightarrow E_{2,\Sigma}^{-\ell_0+r,k_0-(r-1)}$ from \ref{it: abstract E2 injection 2} restricts to an injection $E_{r,\Sigma'}^{-\ell_0+r,k_0-(r-1)}\hookrightarrow E_{r,\Sigma}^{-\ell_0+r,k_0-(r-1)}$, which together with (\ref{equ: abstract E2 injection 3}) gives (\ref{equ: abstract E2 injection 2}).
Finally, we conclude $E_{2,\Sigma'}^{-\ell_0,k_0}=E_{\infty,\Sigma'}^{-\ell_0,k_0}$ from (\ref{equ: abstract E2 injection 1}) and (\ref{equ: abstract E2 injection 2}).
The last claim on (\ref{equ: abstract E2 injection total}) is clear.
\end{proof}

\begin{lem}\label{lem: abstract E2 seq}
Let $\Sigma$, $\Sigma'$ and $\Sigma''$ be collections of complex with quasi maps between double complex $\cT_{\Sigma'}^{\bullet,\bullet}\dashrightarrow \cT_{\Sigma}^{\bullet,\bullet}\dashrightarrow \cT_{\Sigma''}^{\bullet,\bullet}$. Let $h\in\Z$, $\#I_0\leq \ell\leq \#I_1$ and $k\geq 0$.
Assume that
\begin{enumerate}[label=(\roman*)]
\item \label{it: abstract E2 seq 1} $d_{r,\Sigma}^{-\ell,k}=0$ for each $(-\ell,k)$ with $k-\ell=h$;
\item \label{it: abstract E2 seq 2} $E_{2,\Sigma'}^{-\ell,k}\rightarrow E_{2,\Sigma}^{-\ell,k}\rightarrow E_{2,\Sigma''}^{-\ell,k}$ is a short exact sequence for each $(-\ell,k)$ with $k-\ell\leq h+1$;
\item \label{it: abstract E2 seq 3} $E_{2,\Sigma}^{-\ell,k}=0$ for each $(-\ell,k)$ with $k-\ell<h$.
\end{enumerate}
Then we have $E_{2,\Sigma'}^{-\ell,k}=0=E_{2,\Sigma''}^{-\ell,k}$ for each $(-\ell,k)$ with $k-\ell<h$, and $E_{2,\Sigma'}^{-\ell,k}=E_{\infty,\Sigma'}^{-\ell,k}$ as well as $E_{2,\Sigma''}^{-\ell,k}=E_{\infty,\Sigma''}^{-\ell,k}$ for each $(-\ell,k)$ with $k-\ell=h$.
In particular, we have a short exact sequence
\begin{equation}\label{equ: abstract E2 seq total}
0\rightarrow H^{h}(\mathrm{Tot}(\cT_{\Sigma'}^{\bullet,\bullet}))\rightarrow H^{h}(\mathrm{Tot}(\cT_{\Sigma}^{\bullet,\bullet}))\rightarrow H^{h}(\mathrm{Tot}(\cT_{\Sigma''}^{\bullet,\bullet}))\rightarrow 0
\end{equation}
which is strict with respect to the canonical filtration on each term, whose $\mathrm{gr}^{-\ell}(\cdot)$ recovers the short exact sequence
\[0\rightarrow E_{\infty,\Sigma'}^{-\ell,\ell+h}\rightarrow E_{\infty,\Sigma}^{-\ell,\ell+h}\rightarrow E_{\infty,\Sigma''}^{-\ell,\ell+h}\rightarrow 0\]
for each $\#I_0\leq \ell\leq \#I_1$.
\end{lem}
\begin{proof}
We fix for now a choice of $(-\ell_0,k_0)$ with $k_0-\ell_0=h$. The surjection $E_{2,\Sigma}^{-\ell_0,k_0}\twoheadrightarrow E_{2,\Sigma''}^{-\ell_0,k_0}$ from \ref{it: abstract E2 seq 2} together with \ref{it: abstract E2 seq 1} clearly gives $d_{r,\Sigma''}^{-\ell_0,k_0}=0$ for each $r\geq 2$. For each $r\geq 2$, we have $E_{2,\Sigma}^{-\ell_0-r,k+(r-1)}=0$ from \ref{it: abstract E2 seq 2} and \ref{it: abstract E2 seq 3} and thus $d_{r,\Sigma''}^{-\ell_0-r,k+(r-1)}=0$. These together give $E_{2,\Sigma''}^{-\ell_0,k_0}=E_{\infty,\Sigma''}^{-\ell_0,k_0}$. Thanks to Lemma~\ref{lem: abstract E2 injection}, we also have $E_{2,\Sigma'}^{-\ell_0,k_0}=E_{\infty,\Sigma'}^{-\ell_0,k_0}$, under which the short exact sequence
\[0\rightarrow E_{2,\Sigma'}^{-\ell_0,k_0}\rightarrow E_{2,\Sigma}^{-\ell_0,k_0}\rightarrow E_{2,\Sigma''}^{-\ell_0,k_0}\rightarrow 0\]
from \ref{it: abstract E2 seq 2} induces a short exact sequence
\begin{equation}\label{equ: abstract infty page seq}
0\rightarrow E_{\infty,\Sigma'}^{-\ell_0,k_0}\rightarrow E_{\infty,\Sigma}^{-\ell_0,k_0}\rightarrow E_{\infty,\Sigma''}^{-\ell_0,k_0}\rightarrow 0.
\end{equation}
The LHS injection (resp.~RHS surjection) of (\ref{equ: abstract infty page seq}) for each $(-\ell_0,k_0)$ with $k_0-\ell_0=h$ ensures the injectivity (resp.~surjectivity) of the LHS map (resp.~RHS map) of (\ref{equ: abstract E2 seq total}). The exactness of (\ref{equ: abstract E2 seq total}) in the middle follows again from (\ref{equ: abstract infty page seq}) (for each $(-\ell_0,k_0)$ with $k_0-\ell_0=h$) by counting dimension.
\end{proof}

We devote the next section to a list of concrete examples of $\cT_{\Sigma}^{\bullet,\bullet}$.
\subsection{Tits double complex}\label{subsec: std seq}
We introduce the Tits double complex $\cT_{I_0,I_1}^{\bullet,\bullet}$ associated with a pair $I_0\subseteq I_1\subseteq\Delta$ whose total cohomology computes the cohomology of the Tits complex $\mathbf{C}_{I_0,I_1}$ (see (\ref{equ: Tits double complex total coh})). 
We also introduce a list of variants of $\cT_{I_0,I_1}^{\bullet,\bullet}$ (by allowing general highest weights or by adding relative conditions) which will be important in later applications.

Let $x\in W(G)$ and $I_0,I_1\subseteq J_{x}=\Delta\setminus\mathrm{Supp}(x)$.
Recall from the paragraph before Definition~\ref{def: isotypic factor} that we have the following complex
\begin{equation}\label{equ: general Tits complex}
\mathbf{C}^{x}_{I_0,I_1}:~i_{x,I_1}^{\rm{an}}\rightarrow\cdots\rightarrow\bigoplus_{I_0\subseteq I\subseteq I_1, \#I=\ell}i_{x,I}^{\rm{an}}\rightarrow\cdots\rightarrow i_{x,I_0}^{\rm{an}}
\end{equation}
with $i_{x,I_0}^{\rm{an}}$ placed at degree $-\#I_0$. (We use the convention $\mathbf{C}^{x}_{I_0,I_1}\defeq 0$ if $I_0\not\subseteq I_1$.)
By taking $I=\Delta$ in (\ref{equ: general Tits resolution}), we deduce that
\begin{equation}\label{equ: complex as induction}
\mathbf{C}^{x}_{I_0,I_1}\cong i_{I_1,\Delta}^{\rm{an}}(V_{x,I_0,I_1}^{\rm{an}})[\#I_0].
\end{equation}
We define two double complex $\cS_{x,I_0,I_1}^{\bullet,\bullet}$ and $\cS_{x,I_0,I_1,\flat}^{\bullet,\bullet}$ as $\cT_{\Sigma}^{\bullet,\bullet}$ by taking $M^{\bullet}_{I}$ to be
\[\Hom_{D(G)}(B_{\bullet}(D(G),(i_{x,I}^{\rm{an}})^\vee),1_{D(G)})\]
and
\[\Hom_{D(\fg,P_{I})}(B_{\bullet}(D(\fg,P_{I}),M^{I}(x),1_{D(\fg,P_{I})})\]
respectively.
For each $I_0\subseteq I\subseteq I_1$, recall from Lemma~\ref{lem: g P Ext transfer} that the inclusion $D(\fg,P_{I})\subseteq D(G)$ and $M^{I}(x)\hookrightarrow (i_{x,I}^{\rm{an}})^{\vee}$ induces a quasi-isomorphism
\[\Hom_{D(G)}(B_{\bullet}(D(G),(i_{x,I}^{\rm{an}})^\vee),1_{D(G)})\rightarrow \Hom_{D(\fg,P_{I})}(B_{\bullet}(D(\fg,P_{I}),M^{I}(x),1_{D(\fg,P_{I})})\]
which is functorial with respect to the choice of $I$. Hence, we obtain a map between double complex
\begin{equation}\label{equ: x Tits double flat}
\cS_{x,I_0,I_1}^{\bullet,\bullet}\rightarrow \cS_{x,I_0,I_1,\flat}^{\bullet,\bullet}
\end{equation}
which induces isomorphisms on the first page with the cohomology of their total complex gives
\[\mathrm{Ext}_{G}^{\bullet}(1_{G},\mathbf{C}^{x}_{I_0,I_1})=\mathrm{Ext}_{G}^{\bullet}(1_{G},i_{I_1,\Delta}^{\rm{an}}(V_{x,I_0,I_1}^{\rm{an}})[\#I_0]).\]
More generally, given a reductive $E$-Lie subalgebra $\fh\subseteq\fg$, we define a double complex $\cS_{x,I_0,I_1,\fh}^{\bullet,\bullet}$ as $\cT_{\Sigma}^{\bullet,\bullet}$ by taking $M^{\bullet}_{I}$ to be
\[\Hom_{D(G)}(B_{\bullet}(D(G),U(\fl_{I}\cap\fh),(i_{x,I}^{\rm{an}})^\vee),1_{G}^\vee)\]
and similarly for $\cS_{x,I_0,I_1,\fh,\flat}^{\bullet,\bullet}$, and these double complex are functorial with respect to the choice of $\fh$.
When $x=1$ (resp.~$\fh=0$, resp.~$I_1=\Delta\setminus\mathrm{Supp}(x)$) we omit $x$ (resp.~$\fh$, resp.~$I_1$)from the notation.

Let $I_0\subseteq I_1\subseteq \Delta$.
We define a double complex $\cT_{I_0,I_1}^{\bullet,\bullet}$ by taking $M^{\bullet}_{I}$ to be $C^{\bullet}(L_{I})$ and $\iota_{I,I'}$ to be the natural restriction map $C^{\bullet}(L_{I})\rightarrow C^{\bullet}(L_{I'})$.
Let $v\subseteq \mathbf{Log}_{\emptyset}$ (see the discussion below (\ref{equ: Levi char decomposition})) which comes with subsets $I_{v}\subseteq \Delta$. We abuse the notation $v$ for its corresponding element of $\wedge^{\#v}\Hom(Z_{I_{v}},E)$ as in \emph{loc.cit.}.
We use the shortened notation $\widehat{j}\defeq\Delta\setminus\{j\}$ for each $j\in\Delta$.
For each $I_0\subseteq I\subseteq I_1\cap I_{v}$, recall from (\ref{equ: Levi Lie discrete center factor}) (upon replacing $J$ in \emph{loc.cit.} with $I_{v}$ here) that the inclusion $U(\fz_{I_{v}})\otimes_ED(Z_{I_{v}}^{\dagger})\otimes_ED(L_{I}\cap H_{I_{v}})\subseteq D(L_{I})$ together with the isomorphism between (discrete) groups $\prod_{j\notin I_{v}}\mathbf{Z}_{j}\cong Z_{I_{v}}^{\dagger}$ induces a quasi-isomorphism
\[C^{\bullet}(L_{I})\rightarrow \mathrm{Tot}(\mathrm{CE}^{\bullet}(\fz_{I_{v}})\otimes_E\Hom_{E}(\bigotimes_{j\notin I_{v}}B_{\bullet}(\Z_{j})_{\Z_{j}},C^{\bullet}(L_{I}\cap H_{I_{v}})))\]
which is functorial with respect to the choice of $I$. Using the quasi-isomorphism
\[\wedge^{\bullet}\Hom(\Z_{j},E)\rightarrow \Hom_{E}(B_{\bullet}(\Z_{j})_{\Z_{j}},E)\]
for each $j\notin I_{v}$, as well as the isomorphisms
\[\wedge^{\bullet}\Hom(Z_{I_{v}},E)\cong\mathrm{Tot}(\wedge^{\bullet}\Hom(\fz_{I_{v}},E)\otimes_E\wedge^{\bullet}\Hom(Z_{I_{v}}^{\dagger},E))\cong
\mathrm{Tot}(\wedge^{\bullet}\Hom(\fz_{I_{v}},E)\otimes_E\bigotimes_{j\notin I_{v}}\wedge^{\bullet}\Hom(\Z_{j},E)),\]
we obtain an isomorphism in the derived category of $E$-vector spaces
\[\mathrm{Tot}(\wedge^{\bullet}\Hom(Z_{I_{v}},E)\otimes_EC^{\bullet}(L_{I}\cap H_{I_{v}}))\dashrightarrow C^{\bullet}(L_{I})\]
and in particular a morphism in the derived category of $E$-vector spaces
\[v\otimes_EC^{\bullet-\#v}(L_{I}\cap H_{I_{v}}))\dashrightarrow C^{\bullet}(L_{I})\]
which is functorial with respect to the choice of $I$.
We thus obtain the following (quasi) maps
\begin{equation}\label{equ: double fix v}
\cT_{I_0,I_1,v}^{\bullet,\bullet}\dashrightarrow \cT_{I_0,I_1\cap I_{v}}^{\bullet,\bullet}\rightarrow \cT_{I_0,I_1}^{\bullet,\bullet}
\end{equation}
with $\cT_{I_0,I_1,v}^{\bullet,\bullet}$ being the double complex $\cT_{\Sigma}^{\bullet,\bullet}$ from (\ref{equ: general double complex}) with $M^{\bullet}_{I}$ taken to be $C^{\bullet-\#v}(L_{I}/Z_{I_{v}})\otimes_E v$ when $I_0\subseteq I\subseteq I_1\cap I_{v}$, and taken to be zero otherwise.
Further given $J\subseteq I_1\cap I_{v}$ and a (closed) normal subgroup $M\subseteq L_{J}$ containing $Z_{I_{v}}$ with associated $E$-Lie algebra $\fm$, and a reductive $E$-Lie subalgebra $\fh\subseteq \fg$ satisfying $\fh\cap\fm=0$, we define a double complex $\cT_{I_0,I_1,v,\fh,J,M}^{\bullet,\bullet}$ as $\cT_{\Sigma}^{\bullet,\bullet}$ from (\ref{equ: general double complex}) by taking $M^{\bullet}_{I}$ to be $C^{\bullet-\#v}(L_{I}/M,(\fl_{I}/\fm)\cap\fh)\otimes_E v$ (where we identify $\fh$ with a $E$-Lie subalgebra of $\fl_{I}/\fm$ using $\fh\cap\fm=0$) when $I_0\subseteq I\subseteq I_1\cap I_{v}\cap J$, and to be zero otherwise.
We naturally obtain a map between double complex
\begin{equation}\label{equ: double v relative general}
\cT_{I_0,I_1,v,\fh,J,M}^{\bullet,\bullet}\rightarrow \cT_{I_0,I_1,v}^{\bullet,\bullet}.
\end{equation}
When $v=\emptyset$ (resp.~$\fh=\fg$, resp.~$M=Z_{I_{v}}$, resp.~$J\supseteq I_1\cap I_{v}$), we omit this subscript $v$ (resp.~$\fh$, resp.~$M$, resp.~$J$) from the notation $\cT_{I_0,I_1,v,\fh,J,M}^{\bullet,\bullet}$.
Note that (\ref{equ: double fix v}) and (\ref{equ: double v relative general}) induce the following maps between spectral sequences
\begin{equation}\label{equ: double v relative spectral seq}
E_{\bullet,I_0,I_1,v,\fh,J,M}^{\bullet,\bullet}\rightarrow E_{\bullet,I_0,I_1,v}^{\bullet,\bullet}\rightarrow E_{\bullet,I_0,I_1\cap I_{v}}^{\bullet,\bullet}\rightarrow E_{\bullet,I_0,I_1}^{\bullet,\bullet}
\end{equation}
Note that we have
\begin{equation}\label{equ: general double E1}
E_{1,I_0,I_1,v,\fh,J,M}^{-\ell,k}=\bigoplus_{I_0\subseteq I\subseteq I_1\cap J,\#I=\ell}H^{k-\#v}(L_{I}/M,(\fl_{I}/\fm)\cap\fh,1_{L_{I}/M})\otimes_Ev
\end{equation}
for each $(-\ell,k)$.

Recall from Lemma~\ref{lem: g P Ext transfer} and (\ref{equ: PS parabolic Levi resolution}) that we have the following quasi-isomorphisms
\begin{equation}\label{equ: PS to Levi resolution I}
\Hom_{D(G)}(B_{\bullet}(D(G),(i_{I}^{\rm{an}})^{\vee}),1_{D(G)})\rightarrow \Hom_{D(\fg,P_{I})}(B_{\bullet}(D(\fg,P_{I}),M^{I}(x)),1_{D(\fg,P_{I})})\rightarrow C^{\bullet}(L_{I})
\end{equation}
which are functorial with respect to the choice of $I$. Hence, we obtain maps between double complex
\begin{equation}\label{equ: coh Tits Levi}
\cS_{I_0,I_1}^{\bullet,\bullet}\rightarrow \cS_{I_0,I_1,\flat}^{\bullet,\bullet}\rightarrow \cT_{I_0,I_1}^{\bullet,\bullet}
\end{equation}
which induce isomorphisms on the first page.
In particular, we have a spectral sequence (see also (\ref{equ: first complex spectral seq}) with $\mathbf{C}_0=1_{G}$ and $\mathbf{C}_1=\mathbf{C}_{I_0,I_1}$)
\begin{equation}\label{equ: K Tits seq}
E_{\bullet,I_0,I_1}^{\bullet,\bullet}
\end{equation}
converging to
\begin{multline}\label{equ: Tits double complex total coh}
\mathrm{Ext}_{G}^{\bullet}(1_{G},i_{I_1}^{\rm{an}}(V_{I_0,I_1}^{\rm{an}})[\#I_0])\cong\mathrm{Ext}_{G}^{\bullet}(1_{G},\mathbf{C}_{I_0,I_1})\\
=H^{\bullet}(\mathrm{Tot}(\cS_{I_0,I_1}^{\bullet,\bullet}))\buildrel\sim\over\longrightarrow H^{\bullet}(\mathrm{Tot}(\cS_{I_0,I_1,\flat}^{\bullet,\bullet}))\buildrel\sim\over\longrightarrow H^{\bullet}(\mathrm{Tot}(\cT_{I_0,I_1}^{\bullet,\bullet}))
\end{multline}
whose first page has the form
\begin{multline}\label{equ: first page as Levi coh K}
E_{1,I_0,I_1}^{-\ell,k}=\bigoplus_{I_0\subseteq I\subseteq I_1, \#I=\ell}\mathrm{Ext}_{G}^k(1_{G},i_{I}^{\rm{an}})\\
\buildrel\sim\over\longrightarrow \bigoplus_{I_0\subseteq I\subseteq I_1, \#I=\ell}\mathrm{Ext}_{D(\fg,P_{I})}^k(M^{I}(1),1_{D(\fg,P_{I})}) \buildrel\sim\over\longrightarrow \bigoplus_{I_0\subseteq I\subseteq I_1, \#I=\ell} H^k(L_{I},1_{L_{I}}).
\end{multline}

Let $I\subseteq \Delta$.
We recall the isogeny $\varsigma_{I}:\prod_{j\in\Delta\setminus I}\Z_{j}\rightarrow Z_{I}''$ with image $Z_{I}^{\dagger}$ in the paragraph before (\ref{equ: Levi separate discrete diagram}). We use the shortened notation $B_{\bullet,j}\defeq B_{\bullet}(\Z_{j})_{\Z_{j}}$ for each $j\in\Delta\setminus I$ and write $B_{\bullet}^{I}\defeq \bigotimes_{j\in\Delta\setminus I}B_{\bullet,j}$ for short.
For each $j\in\Delta\setminus I$, we have natural identifications $\sigma_{\geq 0}(B_{\bullet,j})=\tau_{\geq 0}(B_{\bullet,j})$ and $\sigma_{\leq -1}(B_{\bullet,j})=\tau_{\leq -1}(B_{\bullet,j})$ and in particular a natural splitting
\[B_{\bullet,j}=\sigma_{\geq 0}(B_{\bullet,j})\oplus \sigma_{\leq -1}(B_{\bullet,j})=\tau_{\geq 0}(B_{\bullet,j})\oplus \tau_{\leq -1}(B_{\bullet,j}).\]
For each $J\subseteq \Delta\setminus I$, we define
\[\tau^{J}(B_{\bullet}^{I})\defeq
(\bigotimes_{j\in\Delta\setminus (I\sqcup J)}B_{\bullet,j})\otimes_E(\bigotimes_{j\in J}\sigma_{\leq -1}(B_{\bullet,j}))
\]
with $\tau^{\emptyset}(B_{\bullet}^{I})=B_{\bullet}^{I}$ and
\[\mathrm{gr}^{J}(B_{\bullet}^{I})\defeq (\bigotimes_{j\in\Delta\setminus (I\sqcup J)}\sigma_{\geq 0}(B_{\bullet,j}))\otimes_E(\bigotimes_{j\in J}\sigma_{\leq -1}(B_{\bullet,j}))\cong \bigotimes_{j\in J}\sigma_{\leq -1}(B_{\bullet,j}).
\]
Note that we have natural truncation map between complex $\tau^{J}(B_{\bullet}^{I})\rightarrow \tau^{J'}(B_{\bullet}^{I})$ for each $J\subseteq J'$.
For each $i\geq 0$, we define
\[\tau^{i}(B_{\bullet}^{I})\defeq \sum_{J\subseteq\Delta\setminus I, \#J\geq i}\tau^{J}(B_{\bullet}^{I})=\sum_{J\subseteq\Delta\setminus I, \#J=i}\tau^{J}(B_{\bullet}^{I})\]
and
\[\mathrm{gr}^{i}(B_{\bullet}^{I})\defeq \sum_{J\subseteq\Delta\setminus I, \#J=i}\mathrm{gr}^{J}(B_{\bullet}^{I}).\]
Recall from Lemma~\ref{lem: abelian gp coh} that we have
\[H^{\bullet}(\Hom_{E}(B_{\bullet}^{I},E))=\mathrm{Tot}(\bigotimes_{j\in\Delta\setminus I}H^{\bullet}(\Z_{j},1_{\Z_{j}}))
=\mathrm{Tot}(\bigotimes_{j\in\Delta\setminus I}\wedge^{\bullet}\Hom(\Z_{j},E))=\wedge^{\bullet}\Hom(\prod_{j\in\Delta\setminus I}\Z_{j},E)\]
which can be identified with
\[H^{\bullet}(Z_{I}^{\dagger},1_{Z_{I}^{\dagger}})=\wedge^{\bullet}\Hom(Z_{I}^{\dagger},E)=\wedge^{\bullet}\Hom_{\infty}(Z_{I},E)\]
under the isogeny $\varsigma_{I}$. Under these identifications, we write
\begin{equation}\label{equ: sm Hom J filtration}
\tau^{J}(\wedge^{\bullet}\Hom_{\infty}(Z_{I},E))=\tau^{J}(\wedge^{\bullet}\Hom(Z_{I}^{\dagger},E))\defeq H^{\bullet}(\Hom_{E}(\tau^{J}(B_{\bullet}^{I}),E))
\end{equation}
and similarly for the others. Note that we have
\begin{equation}\label{equ: sm Hom J graded}
\mathrm{gr}^{J}(\wedge^{\bullet}\Hom_{\infty}(Z_{I},E))=\bigotimes_{j\in J}\Hom(\Z_{j},E)\cong\bigotimes_{j\in J}\Hom(Z_{\Delta\setminus\{j\}}^{\dagger},E)\cong\wedge^{\#J}\Hom(Z_{\Delta\setminus J}^{\dagger},E)
\end{equation}
under the isogeny $\varsigma_{\Delta\setminus\{j\}}$ for each $j\in J$ and the isogeny $\varsigma_{\Delta\setminus J}$.

Let $I_0\subseteq I_1\subseteq\Delta$ and $\fh\subseteq \fg$ be an $E$-Lie subalgebra.
For each $I\subseteq\Delta$, recall from (\ref{equ: Levi separate discrete diagram}) that we have a quasi-isomorphism
\begin{equation}\label{equ: Levi to natural}
C^{\bullet}(L_{I},\fl_{I}\cap\fh)
\rightarrow C^{\bullet}(L_{I},\fl_{I}\cap\fh)_{\natural}=\mathrm{Tot}(\Hom_{E}(B_{\bullet}^{I},C^{\bullet}(L_{I}',\fl_{I}\cap\fh)))
\end{equation}
which is functorial with respect to the choice of $I$ and the map $L_{I'}'\times Z_{I'}^{\dagger}\rightarrow L_{I}'\times Z_{I}^{\dagger}$ induced from $L_{I'}\rightarrow L_{I}$.
We define a double complex $\cT_{I_0,I_1,\fh,\natural}^{\bullet,\bullet}$ by taking $M^{\bullet}_{I}$ to be $C^{\bullet}(L_{I},\fl_{I}\cap\fh)_{\natural}$, which comes with a map
\begin{equation}\label{equ: double Levi to natural}
\cT_{I_0,I_1,\fh}^{\bullet,\bullet}\rightarrow \cT_{I_0,I_1,\fh,\natural}^{\bullet,\bullet}
\end{equation}
between double complex which induces an isomorphism on the first page.
For each $J\subseteq\Delta$, we define double complex $\tau^{J}(\cT_{I_0,I_1,\fh}^{\bullet,\bullet})$ and $\mathrm{gr}^{J}(\cT_{I_0,I_1,\fh}^{\bullet,\bullet})$ by taking $M^{\bullet}_{I}$ to be zero if $J\not\subseteq\Delta\setminus I$, and to be
\begin{equation}\label{equ: Levi cochain natural J}
\tau^{J}(C^{\bullet}(L_{I},\fl_{I}\cap\fh)_{\natural})\defeq \mathrm{Tot}(\Hom_{E}(\tau^{J}(B_{\bullet}^{I}),C^{\bullet}(L_{I}',\fl_{I}\cap\fh)))
\end{equation}
and
\[\mathrm{gr}^{J}(C^{\bullet}(L_{I},\fl_{I}\cap\fh)_{\natural})\defeq \mathrm{Tot}(\Hom_{E}(\mathrm{gr}^{J}(B_{\bullet}^{I}),C^{\bullet}(L_{I}',\fl_{I}\cap\fh)))\]
respectively, for each $I_0\subseteq I\subseteq I_1$.
We can similarly define $\tau^{i}(\cT_{I_0,I_1,\fh}^{\bullet,\bullet})$ and $\mathrm{gr}^{i}(\cT_{I_0,I_1,\fh}^{\bullet,\bullet})$, with an identification
\begin{equation}\label{equ: grade J direct sum}
\mathrm{gr}^{i}(\cT_{I_0,I_1,\fh}^{\bullet,\bullet})=\bigoplus_{J\subseteq \Delta, \#J=i}\mathrm{gr}^{J}(\cT_{I_0,I_1,\fh}^{\bullet,\bullet}).
\end{equation}
For each $J\subseteq J'\subseteq \Delta$, we have an embedding between double complex
\begin{equation}\label{equ: double J embedding}
\tau^{J'}(\cT_{I_0,I_1,\fh}^{\bullet,\bullet})\hookrightarrow \tau^{J}(\cT_{I_0,I_1,\fh}^{\bullet,\bullet}).
\end{equation}
For each $i\geq 0$, we have the following short exact sequence between double complex
\begin{equation}\label{equ: double i seq}
0\rightarrow\tau^{i+1}(\cT_{I_0,I_1,\fh}^{\bullet,\bullet})\hookrightarrow \tau^{i}(\cT_{I_0,I_1,\fh}^{\bullet,\bullet})\twoheadrightarrow \mathrm{gr}^{i}(\cT_{I_0,I_1,\fh}^{\bullet,\bullet})\rightarrow 0.
\end{equation}
We write $E_{\bullet,I_0,I_1,\fh,\natural}^{\bullet,\bullet}$ for the spectral sequence associated with $\cT_{I_0,I_1,\fh,\natural}^{\bullet,\bullet}$. Note that (\ref{equ: double Levi to natural}) induces an isomorphism
\[E_{\bullet,I_0,I_1,\fh}^{\bullet,\bullet}\buildrel\sim\over\longrightarrow E_{\bullet,I_0,I_1,\fh,\natural}^{\bullet,\bullet}\]
between spectral sequences from the first page, so we do not distinguish between these two spectral sequences from now.
For each $J\subseteq \Delta$ and $i\geq 0$, we define $\tau^{J}(E_{\bullet,I_0,I_1,\fh}^{\bullet,\bullet})$ as the spectral sequence associated with $\tau^{J}(\cT_{I_0,I_1,\fh}^{\bullet,\bullet})$, and similarly for the others.
For each $J\subseteq \Delta$, we have (using Lemma~\ref{lem: center compact} and (\ref{equ: sm Hom J filtration}))
\[\tau^{J}(E_{1,I_0,I_1,\fh}^{-\ell,k})=\bigoplus_{I_0\subseteq I\subseteq I_1\setminus J, \#I=\ell}\mathrm{Tot}(H^{\bullet}(\fl_{I},\fl_{I}\cap\fh,1_{\fl_{I}})\otimes_E\tau^{J}(\wedge^{\bullet}\Hom(Z_{I}^{\dagger},E)))^{k}\]
for each $(-\ell,k)$. Similar formula hold for $\mathrm{gr}^{J}(E_{1,I_0,I_1,\fh}^{-\ell,k})$, $\tau^{i}(E_{1,I_0,I_1,\fh}^{-\ell,k})$ and $\mathrm{gr}^{i}(E_{1,I_0,I_1,\fh}^{-\ell,k})$ for each $i\geq 0$ and $(-\ell,k)$.
For each $J\subseteq J'$, (\ref{equ: double J embedding}) induces an embedding
\[\tau^{J'}(E_{1,I_0,I_1,\fh}^{-\ell,k})\hookrightarrow \tau^{J}(E_{1,I_0,I_1,\fh}^{-\ell,k})\]
given by the direct sum over all $I_0\subseteq I\subseteq I_1\setminus J$ with $\#I=\ell$ of the following map
\[\mathrm{Tot}(H^{\bullet}(\fl_{I},\fl_{I}\cap\fh,1_{\fl_{I}})\otimes_E\tau^{J'}(\wedge^{\bullet}\Hom(Z_{I}^{\dagger},E)))^{k}\hookrightarrow \mathrm{Tot}(H^{\bullet}(\fl_{I},\fl_{I}\cap\fh,1_{\fl_{I}})\otimes_E\tau^{J}(\wedge^{\bullet}\Hom(Z_{I}^{\dagger},E)))^{k}.\]

Let $I_0\subseteq I_1\subseteq \Delta$, $J\subseteq\Delta$ and $\fh\subseteq\fg$ be a reductive $E$-Lie subalgebra.
We define double complex $\mathrm{CE}_{I_0,I_1,\fh}^{\bullet,\bullet}$ by taking $M^{\bullet}_{I}$ to be $\mathrm{CE}^{\bullet}(\fl_{I},\fl_{I}\cap\fh)$, with $\iota_{I,I'}$ being zero when $I\cap J\neq\emptyset$ and being to be the restriction map 
\[\mathrm{CE}^{\bullet}(\fl_{I},\fl_{I}\cap\fh)\rightarrow \mathrm{CE}^{\bullet}(\fl_{I'},\fl_{I'}\cap\fh).\]
We omit $\fh$ from the notation when $\fh=0$.
We also similarly define a double complex $\mathrm{Inv}_{I_0,I_1}^{\bullet,\bullet}$ by taking $M^{\bullet}_{I}$ to be $\mathrm{Inv}^{\bullet}(\fl_{I})=\mathrm{CE}^{\bullet}(\fl_{I})^{\fl_{I}}$.
It follows from \cite[\S 9, \S 10]{Kos50}) that $\mathrm{Inv}^{\bullet}(\fl_{I})$ has zero differential maps and the embedding (which is clearly functorial with respect to the choice of $I$)
\[\mathrm{Inv}^{\bullet}(\fl_{I})\hookrightarrow \mathrm{CE}^{\bullet}(\fl_{I})\]
is a quasi-isomorphism. We thus obtain an embedding between double complex
\begin{equation}\label{equ: double Inv to CE}
\mathrm{Inv}_{I_0,I_1}^{\bullet,\bullet}\hookrightarrow \mathrm{CE}_{I_0,I_1}^{\bullet,\bullet}
\end{equation}
which induces an isomorphism on the first page.
Since $\mathrm{Inv}^{\bullet}(\fl_{I})$ has zero differential maps for each $I_0\subseteq I\subseteq I_1$, we have the following isomorphism between complex
\[\mathrm{Tot}(\mathrm{Inv}_{I_0,I_1}^{\bullet,\bullet})\cong\bigoplus_{k\in\Z}\mathrm{Inv}_{I_0,I_1}^{\bullet,k}[-k]\]
and therefore
\begin{equation}\label{equ: Lie Tits coh decomposition} H^{m}(\mathrm{Tot}(\mathrm{CE}_{I_0,I_1}^{\bullet,\bullet}))\buildrel\sim\over\longleftarrow H^{m}(\mathrm{Tot}(\mathrm{Inv}_{I_0,I_1}^{\bullet,\bullet}))
\cong \bigoplus_{k\geq 0}H^{m-k}(\mathrm{Inv}_{I_0,I_1}^{\bullet,k})\cong \bigoplus_{\ell=\#I_0}^{\#I_1}H^{-\ell}(\mathrm{Inv}_{I_0,I_1}^{\bullet,\ell+m})
\end{equation}
for each $m\in\Z$.

For each $J\subseteq\Delta$ and $I_0\subseteq I\subseteq I_1\setminus J$, we deduce from Lemma~\ref{lem: center compact} that the inclusion $U(\fl_{I})\subseteq D(L_{I}')$ together with Lemma~\ref{lem: Lie resolution left diagonal transfer} induces the following quasi-isomorphisms
\begin{multline*}
\mathrm{gr}^{J}(C^{\bullet}(L_{I},\fl_{I}\cap\fh)_{\natural})=\mathrm{Tot}(\Hom_{E}(\mathrm{gr}^{J}(B_{\bullet}^{I}),C^{\bullet}(L_{I}',\fl_{I}\cap\fh)))\\
\rightarrow \mathrm{Tot}(\Hom_{E}(\mathrm{gr}^{J}(B_{\bullet}^{I}),\mathrm{CE}^{\bullet}(\fl_{I},\fl_{I}\cap\fh))) =\mathrm{Tot}(\mathrm{CE}^{\bullet}(\fl_{I},\fl_{I}\cap\fh)\otimes_E\Hom_{E}(\mathrm{gr}^{J}(B_{\bullet}^{I}),E))\\
\leftarrow\mathrm{CE}^{\bullet-\#J}(\fl_{I},\fl_{I}\cap\fh)\otimes_E \wedge^{\#J}\Hom(Z_{\Delta\setminus J}^{\dagger},E)
\end{multline*}
which is functorial with respect to the choice of $I$.
Hence, we obtain a quasi-map between double complex
\begin{equation}\label{equ: grade J Z to CE}
\mathrm{gr}^{J}(\cT_{I_0,I_1,\fh}^{\bullet,\bullet})\dashrightarrow \mathrm{CE}_{I_0,I_1\setminus J,\fh}^{\bullet,\bullet-\#J}\otimes_E\wedge^{\#J}\Hom(Z_{\Delta\setminus J}^{\dagger},E)
\end{equation}
which induces an isomorphism on the first page.
In particular, taking $J=\emptyset$, we obtain the following maps between double complex
\begin{equation}\label{equ: double Levi to CE}
\cT_{I_0,I_1,\fh}^{\bullet,\bullet}\rightarrow \cT_{I_0,I_1,\fh,\natural}^{\bullet,\bullet}\rightarrow \mathrm{gr}^{\emptyset}(\cT_{I_0,I_1,\fh}^{\bullet,\bullet}) \rightarrow \mathrm{CE}_{I_0,I_1,\fh}^{\bullet,\bullet}.
\end{equation}
with both the first map and the third map inducing isomorphisms on the first page.
We will frequently identify the spectral sequence attached to $\mathrm{CE}_{I_0,I_1,\fh}^{\bullet,\bullet}$ with $\mathrm{gr}^{\emptyset}(E_{\bullet,I_0,I_1,\fh}^{\bullet,\bullet})$. In particular, when $\fh=0$, (\ref{equ: Lie Tits coh decomposition}) induces an isomorphism
\begin{equation}\label{equ: Lie Tits coh splitting}
H^{m}(\mathrm{Tot}(\mathrm{CE}_{I_0,I_1}^{\bullet,\bullet}))
\cong\bigoplus_{\ell=\#I_0}^{\#I_1}\mathrm{gr}^{\emptyset}(E_{2,I_0,I_1}^{-\ell,\ell+m})
\end{equation}

Let $x\in \Gamma$, $I_0\subseteq I_1\subseteq J_{x}=\Delta\setminus\mathrm{Supp}(x)$ and $\fh\subseteq\fg$ be a reductive $E$-Lie algebra.
We define a double complex $\tld{\mathrm{CE}}_{x,I_0,I_1,\fh}^{\bullet,\bullet}$ by taking $M^{\bullet}_{I}$ to be $\mathrm{CE}^{\bullet}(\fg,\fl_{I}\cap\fh,M^{I}(x))$ for each $I_0\subseteq I\subseteq I_1$. We omit the notation $x$ when $x=1$.
If $x=1$, the inclusion $\fl_{I}\subseteq\fg$ (fixed) surjection $M^{I}(1)\twoheadrightarrow L(1)$ induces a quasi-isomorphism (see (\ref{equ: Verma parabolic Levi resolution}))
\[\mathrm{CE}^{\bullet}(\fg,\fl_{I}\cap\fh,M^{I}(1))\rightarrow \mathrm{CE}^{\bullet}(\fl_{I},\fl_{I}\cap\fh)\]
which is functorial with respect to the choice of $I$. We thus obtain a map between double complex
\begin{equation}\label{equ: Tits Verma to Levi}
\tld{\mathrm{CE}}_{I_0,I_1,\fh}^{\bullet,\bullet}\rightarrow \mathrm{CE}_{I_0,I_1,\fh}^{\bullet,\bullet}
\end{equation}
which induces isomorphisms on the first page. We clearly have a commutative diagram that maps (\ref{equ: coh Tits Levi}) (upon omitting its middle term) to (\ref{equ: Tits Verma to Levi}).

Let $x,w\in \Gamma$ with $x\unlhd w$ and $I_0\subseteq J_{w}=\Delta\setminus\mathrm{Supp}(w)$. Recall from (\ref{equ: C Tits x w transfer}) that we have a map between complex $\mathbf{C}^{x}_{I_0}\rightarrow \mathbf{C}^{w}_{I_0}$ whose dual restricts to a map between complex of $U(\fg)$-modules $\mathfrak{c}^{w}_{I_0}\rightarrow \mathfrak{c}^{x}_{I_0}$. We thus obtain the following commutative diagram of maps between double complex
\begin{equation}\label{equ: double x w diagram}
\xymatrix{
\cS_{x,I_0,\fh}^{\bullet,\bullet} \ar[r] \ar[d] & \cS_{x,I_0,\fh,\flat}^{\bullet,\bullet} \ar[r] \ar[d] & \tld{\mathrm{CE}}_{x,I_0,\fh}^{\bullet,\bullet} \ar[d]\\
\cS_{w,I_0,\fh}^{\bullet,\bullet} \ar[r] & \cS_{w,I_0,\fh,\flat}^{\bullet,\bullet} \ar[r] & \tld{\mathrm{CE}}_{w,I_0,\fh}^{\bullet,\bullet}
}
\end{equation}
for each reductive $E$-Lie subalgebra $\fh\subseteq \fg$.
\subsection{Reduce extension to cohomology}\label{subsec: Ext Tits preliminary}
We show how to reduce the computation of extensions between two Tits complex to that of the cohomology of a single Tits complex, with the main result being Proposition~\ref{prop: Ext complex std seq}.

\begin{lem}\label{lem: Ext PS change}
Let $x\in W(G)$, $I_0\subseteq\Delta$ and $I_1\subseteq I_{x}$. We have the following results.
\begin{enumerate}[label=(\roman*)]
\item \label{it: Ext PS change 1} If $I_0\subseteq I_{x}$ and $I_1\not\subseteq I_0$, then we have 
    \begin{equation}\label{equ: Ext PS change vanishing}
    \mathrm{Ext}^k_{G}(i_{x,I_0}^{\rm{an}},i_{x,I_1}^{\rm{an}})=0
    \end{equation}
    for each $k\geq 0$.
\item \label{it: Ext PS change 1 sm} If $I_1\not\subseteq I_0$, then we have 
    \begin{equation}\label{equ: Ext PS change vanishing sm}
    \mathrm{Ext}^k_{G}(i_{I_0}^{\infty},i_{x,I_1}^{\rm{an}})=0
    \end{equation}
    for each $k\geq 0$.
\item \label{it: Ext PS change 2} If $I_1\subseteq I_0\subseteq I_{x}$, then the embedding $\cF_{P_{I_0}}^{G}(L(x),1_{L_{I_0}})\hookrightarrow i_{I_0}^{\rm{an}}$ induces an isomorphism
    \begin{equation}\label{equ: Ext PS sm change}
    \mathrm{Ext}^k_{G}(i_{I_0}^{\rm{an}},i_{I_1}^{\rm{an}})\buildrel\sim\over\longrightarrow \mathrm{Ext}^k_{G}(\cF_{P_{I_0}}^{G}(L(x),1_{L_{I_0}}),i_{I_1}^{\rm{an}})
    \end{equation}
    for each $k\geq 0$.
\item \label{it: Ext PS change 3} If $I_1\subseteq I_0\subseteq I_{x}$, then for each $I_0\subseteq I_0'$ the embedding $i_{x,I_0'}^{\rm{an}}\hookrightarrow i_{x,I_0}^{\rm{an}}$ induces an isomorphism
    \begin{equation}\label{equ: Ext PS change}
    \mathrm{Ext}^k_{G}(i_{x,I_0}^{\rm{an}},i_{x,I_1}^{\rm{an}})\buildrel\sim\over\longrightarrow \mathrm{Ext}^k_{G}(i_{x,I_0'}^{\rm{an}},i_{x,I_1}^{\rm{an}})
    \end{equation}
    for each $k\geq 0$.
\item \label{it: Ext PS change 3 sm} If $I_1\subseteq I_0$, then for each $I_0\subseteq I_0'$ the embedding $i_{I_0'}^{\infty}\hookrightarrow i_{I_0}^{\infty}$ induces an isomorphism
    \begin{equation}\label{equ: Ext PS change sm}
    \mathrm{Ext}^k_{G}(i_{I_0}^{\infty},i_{x,I_1}^{\rm{an}})\buildrel\sim\over\longrightarrow \mathrm{Ext}^k_{G}(i_{I_0'}^{\infty},i_{x,I_1}^{\rm{an}})
    \end{equation}
    for each $k\geq 0$
\end{enumerate}
\end{lem}
\begin{proof}
\ref{it: Ext PS change 1} follows from Lemma~\ref{lem: two PS vanishing} by taking $x_0=x_1=x$ in \emph{loc.cit.}, and \ref{it: Ext PS change 1 sm} follows from Lemma~\ref{lem: vanishing by sm} by taking $M_0=L(1)$, $I=I_0$ and $x_1=x$ in \emph{loc.cit.}.
Assume from now that $I_1\subseteq I_0$.

We prove \ref{it: Ext PS change 2}.\\
We write $V\defeq \cF_{P_{I_0}}^{G}(L(x),1_{L_{I_0}})$ for short.
As $i_{x,I_0}^{\rm{an}}/V\cong \cF_{P_{I_0}}^{G}(N^{I_0}(x),1_{L_{I_0}})$ from Proposition~\ref{prop: OS property}, each $W\in\mathrm{JH}_{G}(i_{x,I_0}^{\rm{an}}/V)$ has the form $W=C_{w,I}$ for some $w>x$ and $I\subseteq I_w$, which satisfies $\mathrm{Ext}_G^k(W,i_{x,I_1}^{\rm{an}})=0$ for $k\geq 0$ by Lemma~\ref{lem: factor PS vanishing}. Upon a d\'evissage with respect to $\mathrm{JH}_{G}(i_{x,I_0}^{\rm{an}}/V)$ gives $\mathrm{Ext}_G^k(i_{x,I_0}^{\rm{an}}/V,i_{x,I_1}^{\rm{an}})=0$ for $k\geq 0$. A further d\'evissage with respect to $0\rightarrow V\rightarrow i_{I_0}^{\rm{an}}\rightarrow i_{I_0}^{\rm{an}}/V\rightarrow 0$ gives \ref{it: Ext PS change 2}.

We prove \ref{it: Ext PS change 3}.\\
Let $W\in\mathrm{JH}_{G}(i_{x,I_0}^{\rm{an}}/i_{x,I_0'}^{\rm{an}})$ and write $W=C_{w,I}$ for some $I\subseteq I_{w}$. If $w>x$, we have $\mathrm{Ext}_G^k(W,i_{x,I_1}^{\rm{an}})=0$ for $k\geq 0$ by Lemma~\ref{lem: factor PS vanishing}. If $w=x$, then we have $I_1\subseteq I_0\subseteq I\subsetneq I_{x}$ which gives $\mathrm{Ext}_G^k(W,i_{x,I_1}^{\rm{an}})=0$ for $k\geq 0$ using Lemma~\ref{lem: vanishing by sm}.
A d\'evissage with respect to $\mathrm{JH}_{G}(i_{x,I_0}^{\rm{an}}/i_{x,I_0'}^{\rm{an}})$ thus gives $\mathrm{Ext}_G^k(i_{x,I_0}^{\rm{an}}/i_{x,I_0'}^{\rm{an}},i_{x,I_1}^{\rm{an}})=0$ for each $k\geq 0$, which together with a further d\'evissage with respect to $0\rightarrow i_{x,I_0'}^{\rm{an}}\rightarrow i_{x,I_0}^{\rm{an}}\rightarrow i_{x,I_0}^{\rm{an}}/i_{x,I_0'}^{\rm{an}}\rightarrow 0$ gives (\ref{equ: Ext PS change}).

We prove \ref{it: Ext PS change 3 sm}.\\ 
Note that each $W\in\mathrm{JH}_{G}(i_{I_0}^{\infty}/i_{I_0'}^{\infty})$ has the form $V_{I}^{\infty}$ for some $I_1\subseteq I_0\subseteq I\subsetneq \Delta$ and thus
$d_{\Delta}(V_{I}^{\infty},Q_{\Delta}(\Delta,I_1))=\infty$
by Lemma~\ref{lem: special sm distance}, which together with Lemma~\ref{lem: sm distance Ext vanishing} gives $\mathrm{Ext}_G^k(W,i_{x,I_1}^{\rm{an}})=0$ for $k\geq 0$. A d\'evissage with respect to $W\in\mathrm{JH}_{G}(i_{I_0}^{\infty}/i_{I_0'}^{\infty})$ gives $\mathrm{Ext}_G^k(i_{I_0}^{\infty}/i_{I_0'}^{\infty},i_{x,I_1}^{\rm{an}})=0$ for $k\geq 0$, which together with a further d\'evissage with respect to 
$0\rightarrow i_{I_0'}^{\infty}\rightarrow i_{I_0}^{\infty}\rightarrow i_{I_0}^{\infty}/i_{I_0'}^{\infty}\rightarrow 0$ gives (\ref{equ: Ext PS change sm}).
\end{proof}

Given $I_0\subseteq I_1\subseteq \Delta$, we define a complex $\mathbf{C}_{I_0,I_1}^{\infty}$ by replacing $i_{I}^{\rm{an}}$ in the definition of $\mathbf{C}_{I_0,I_1}$ (for each $I_0\subseteq I\subseteq I_1$) with $i_{I}^{\infty}$. Similar to (\ref{equ: general Tits resolution}), we have
\[\mathbf{C}_{I_0,I_1}^{\infty}\buildrel\sim\over\longrightarrow i_{I_0}^{\infty}(V_{I_0,I_1}^{\infty})[\#I_0].\]
\begin{lem}\label{lem: acyclic complex}
Let $x\in W(G)$ with $J_{x}=\Delta\setminus\mathrm{Supp}(x)$. Let $I_0\subseteq I_1\subseteq\Delta$ and $I\subseteq J_{x}$. We have the following results.
\begin{enumerate}[label=(\roman*)]
\item \label{it: acyclic complex 1} If $I_1\subseteq J_{x}$ and $I_0\cup I\neq I_1$, then we have 
\[
    \mathrm{Ext}_{G}^{k}(\mathbf{C}^{x}_{I_0,I_1},i_{x,I}^{\rm{an}})=0
\]
for $k\in\Z$.
\item \label{it: acyclic complex 1 sm} If $I_0\cup I\neq I_1$, then we have 
\[
   \mathrm{Ext}_{G}^{k}(\mathbf{C}_{I_0,I_1}^{\infty},i_{x,I}^{\rm{an}})=0
\]
for $k\in\Z$.
\item \label{it: acyclic complex 2} If $I_1\subseteq J_{x}$ and $I_0\cup I=I_1$, then the truncation map $\mathbf{C}^{x}_{I_0,I_1}\rightarrow i_{x,I_1}^{\rm{an}}[\#I_1]$ induces an isomorphism
\begin{equation}\label{equ: acyclic complex reduction}
\mathrm{Ext}_{G}^{k-\#I_1}(i_{x,I_1}^{\rm{an}},i_{x,I}^{\rm{an}})\buildrel\sim\over\longrightarrow \mathrm{Ext}_{G}^{k}(\mathbf{C}^{x}_{I_0,I_1},i_{x,I}^{\rm{an}}) 
\end{equation}
for $k\in \Z$. 
\item \label{it: acyclic complex 2 sm} If $I_0\cup I=I_1$, then the truncation map $\mathbf{C}_{I_0,I_1}^{\infty}\rightarrow i_{I_1}^{\infty}[\#I_1]$ induces an isomorphism
\begin{equation}\label{equ: acyclic complex reduction sm}
\mathrm{Ext}_{G}^{k-\#I_1}(i_{I_1}^{\infty},i_{x,I}^{\rm{an}})\buildrel\sim\over\longrightarrow \mathrm{Ext}_{G}^{k}(\mathbf{C}_{I_0,I_1}^{\infty},i_{x,I}^{\rm{an}}) 
\end{equation}
for $k\in \Z$.
\end{enumerate}
\end{lem}
\begin{proof}
We prove \ref{it: acyclic complex 1} and \ref{it: acyclic complex 2}.\\
Recall from \ref{it: Ext PS change 1} of Lemma~\ref{lem: Ext PS change} that we have $\mathrm{Ext}_{G}^{k}(i_{x,I'}^{\rm{an}},i_{x,I}^{\rm{an}})=0$ for $k\in \Z$ whenever $I\not\subseteq I'$. A d\'evissage thus implies that the truncation map $\mathbf{C}^{x}_{I_0,I_1}\rightarrow \mathbf{C}^{x}_{I_0\cup I,I_1}$ induces an isomorphism
\begin{equation}\label{equ: acyclic complex truncation}
\mathrm{Ext}_{G}^{k}(\mathbf{C}^{x}_{I_0\cup I,I_1},i_{x,I}^{\rm{an}})\buildrel\sim\over\longrightarrow \mathrm{Ext}_{G}^{k}(\mathbf{C}^{x}_{I_0,I_1},i_{x,I}^{\rm{an}})
\end{equation}
for each $k\in \Z$. We write $E_{I}^{\bullet}\defeq \mathrm{Ext}_{G}^{\bullet}(i_{x,I}^{\rm{an}},i_{x,I}^{\rm{an}})$ for short and recall from \ref{it: Ext PS change 3} of Lemma~\ref{lem: Ext PS change} that the inclusion $i_{x,I'}^{\rm{an}}\hookrightarrow i_{x,I}^{\rm{an}}$ induces an isomorphism
\[E_{I}^{k}=\mathrm{Ext}_{G}^{k}(i_{x,I}^{\rm{an}},i_{x,I}^{\rm{an}})\buildrel\sim\over\longrightarrow \mathrm{Ext}_{G}^{k}(i_{x,I'}^{\rm{an}},i_{x,I}^{\rm{an}})\]
for each $k\in\Z$.
It thus follows from (\ref{equ: second complex spectral seq}) that we can compute LHS of (\ref{equ: acyclic complex truncation}) by a spectral sequence $E_{\bullet}^{\bullet,\bullet}$ whose first page has the form
\begin{equation}\label{equ: acyclic complex seq}
E_{1}^{\ell,k}=\bigoplus_{I_0\cup I\subseteq I'\subseteq I_1, \#I'=\ell}E_{I}^{k}
\end{equation}
with one copy of $E_{I}^{k}$ for each choice of $I'$, and with the differential map $d_{1}^{\ell-1,k}: E_{1}^{\ell-1,k}\rightarrow E_{1}^{\ell,k}$ given by $(-1)^{m(I',i)}$ times the identity map from the copy of $E_{I}^{k}$ indexed by $I'\setminus\{i\}$ to the copy indexed by $I'$, for each $I_0\cup I\subseteq I'\subseteq I_1$ with $\#I'=\ell$ and each $i\in I'\setminus(I_0\cup I)$.
Consequently, the complex $E_{1}^{\bullet,k}$ is acyclic when $I_0\cup I\neq I_1$, and has single term $E_{1}^{\#I_1,k}\cong E_{I}^{k}$ when $I_0\cup I=I_1$. This gives both \ref{it: acyclic complex 1} and \ref{it: acyclic complex 2}.
Finally, the proofs of \ref{it: acyclic complex 1 sm} and \ref{it: acyclic complex 2 sm} are parallel upon replacing the usage of \ref{it: Ext PS change 1} and \ref{it: Ext PS change 3} of Lemma~\ref{lem: Ext PS change} in the proofs of \ref{it: acyclic complex 1} and \ref{it: acyclic complex 2} with \ref{it: Ext PS change 1 sm} and \ref{it: Ext PS change 3 sm} of Lemma~\ref{lem: Ext PS change}.
\end{proof}

\begin{prop}\label{prop: Ext complex std seq}
Let $x\in W(G)$ with $J_{x}=\Delta\setminus\mathrm{Supp}(x)$. Let $I_0\subseteq I_1\subseteq \Delta$ and $I_2\subseteq I_3\subseteq J_{x}$ be two pairs. We have the following results.
\begin{enumerate}[label=(\roman*)]
\item \label{it: Ext complex std seq} If $I_1\subseteq J_{x}$, then the truncation maps $\mathbf{C}^{x}_{I_2,I_3}\rightarrow \mathbf{C}^{x}_{I_2\cup(I_1\setminus I_0),I_3}\leftarrow \mathbf{C}^{x}_{I_2\cup(I_1\setminus I_0),I_1\cap I_3}$ and $\mathbf{C}^{x}_{I_0,I_1}\rightarrow i_{x,I_1}^{\rm{an}}[\#I_1]$ induce the following isomorphisms
\begin{multline}\label{equ: Ext complex std truncation}
\mathrm{Ext}_{G}^k(\mathbf{C}^{x}_{I_0,I_1}, \mathbf{C}^{x}_{I_2,I_3})\buildrel\sim\over\longrightarrow \mathrm{Ext}_{G}^k(\mathbf{C}^{x}_{I_0,I_1}, \mathbf{C}^{x}_{I_2\cup(I_1\setminus I_0),I_3})
\buildrel\sim\over\longleftarrow \mathrm{Ext}_{G}^k(\mathbf{C}^{x}_{I_0,I_1}, \mathbf{C}^{x}_{I_2\cup(I_1\setminus I_0),I_1\cap I_3})\\
\buildrel\sim\over\longleftarrow \mathrm{Ext}_{G}^k(i_{x,I_1}^{\rm{an}}[\#I_1], \mathbf{C}^{x}_{I_2\cup(I_1\setminus I_0),I_1\cap I_3})=\mathrm{Ext}_{G}^{k-\#I_1}(i_{x,I_1}^{\rm{an}}, \mathbf{C}^{x}_{I_2\cup(I_1\setminus I_0),I_1\cap I_3})
\end{multline}
for $k\in \Z$. 
\item \label{it: Ext complex std seq sm}
The truncation maps $\mathbf{C}^{x}_{I_2,I_3}\rightarrow \mathbf{C}^{x}_{I_2\cup(I_1\setminus I_0),I_3}\leftarrow \mathbf{C}^{x}_{I_2\cup(I_1\setminus I_0),I_1\cap I_3}$ and $\mathbf{C}_{I_0,I_1}^{\infty}\rightarrow i_{I_1}^{\infty}[\#I_1]$ induce the following isomorphisms
\begin{multline}\label{equ: Ext complex std truncation sm}
\mathrm{Ext}_{G}^k(\mathbf{C}_{I_0,I_1}^{\infty}, \mathbf{C}^{x}_{I_2,I_3})\buildrel\sim\over\longrightarrow \mathrm{Ext}_{G}^k(\mathbf{C}_{I_0,I_1}^{\infty}, \mathbf{C}^{x}_{I_2\cup(I_1\setminus I_0),I_3})
\buildrel\sim\over\longleftarrow \mathrm{Ext}_{G}^k(\mathbf{C}_{I_0,I_1}^{\infty}, \mathbf{C}^{x}_{I_2\cup(I_1\setminus I_0),I_1\cap I_3})\\
\buildrel\sim\over\longleftarrow \mathrm{Ext}_{G}^k(i_{I_1}^{\infty}[\#I_1], \mathbf{C}^{x}_{I_2\cup(I_1\setminus I_0),I_1\cap I_3})=\mathrm{Ext}_{G}^{k-\#I_1}(i_{I_1}^{\infty}, \mathbf{C}^{x}_{I_2\cup(I_1\setminus I_0),I_1\cap I_3})
\end{multline}
for $k\in \Z$. 
\end{enumerate}
\end{prop}
\begin{proof}
We prove \ref{it: Ext complex std seq}.\\
Note that $I\subseteq J_{x}$ satisfies both $I_2\subseteq I\subseteq I_3$ and $I\cup I_0=I_1$ if and only if $I_2\cup(I_1\setminus I_0)\subseteq I\subseteq I_1\cap I_3$, which together with \ref{it: acyclic complex 1} and \ref{it: acyclic complex 2} of Lemma~\ref{lem: acyclic complex} gives the isomorphisms on the first row of (\ref{equ: Ext complex std truncation}). As $I'\subseteq J_{x}$ satisfies $I_0\subseteq I'\subseteq I_1$ and $I'\supseteq I_2\cup(I_1\setminus I_0)$ only if $I'=I_1$, the isomorphism
\[\mathrm{Ext}_{G}^k(i_{x,I_1}^{\rm{an}}[\#I_1], \mathbf{C}^{x}_{I_2\cup(I_1\setminus I_0),I_1\cap I_3})\buildrel\sim\over\longrightarrow\mathrm{Ext}_{G}^k(\mathbf{C}^{x}_{I_0,I_1}, \mathbf{C}^{x}_{I_2\cup(I_1\setminus I_0),I_1\cap I_3}) \]
follows from \ref{it: Ext PS change 1} of Lemma~\ref{lem: Ext PS change}.

The proof of \ref{it: Ext complex std seq sm} is parallel using \ref{it: acyclic complex 1 sm} and \ref{it: acyclic complex 2 sm} of Lemma~\ref{lem: acyclic complex} as well as \ref{it: Ext PS change 1 sm} of Lemma~\ref{lem: Ext PS change}.
\end{proof}

\begin{rem}\label{rem: Hom Tits}
Under the further assumption $I_0\subseteq I_2\subseteq I_1\subseteq I_3$ and $k=0$, then \ref{it: Ext complex std seq} of Proposition~\ref{prop: Ext complex std seq} says that the truncation maps $\mathbf{C}^{x}_{I_2,I_3}\rightarrow \mathbf{C}^{x}_{I_1,I_3}\leftarrow i_{x,I_1}^{\rm{an}}[\#I_1]$ and $\mathbf{C}^{x}_{I_0,I_1}\rightarrow i_{x,I_1}^{\rm{an}}[\#I_1]$ induce the following isomorphisms
\begin{multline*}
\Hom_{G}(\mathbf{C}^{x}_{I_0,I_1}, \mathbf{C}^{x}_{I_2,I_3})\buildrel\sim\over\longrightarrow \Hom_{G}(\mathbf{C}^{x}_{I_0,I_1}, \mathbf{C}^{x}_{I_1,I_3})
\buildrel\sim\over\longleftarrow \Hom_{G}(\mathbf{C}^{x}_{I_0,I_1}, i_{x,I_1}^{\rm{an}}[\#I_1])\\
\buildrel\sim\over\longleftarrow \Hom_{G}(i_{x,I_1}^{\rm{an}}[\#I_1], i_{x,I_1}^{\rm{an}}[\#I_1])=\Hom_{G}(i_{x,I_1}^{\rm{an}},i_{x,I_1}^{\rm{an}})
\end{multline*}
with the last term being $1$-dimensional by for example Proposition~\ref{prop: St key seq} (see also \cite[Lem.~5.1.1]{BQ24}).
In particular, $\Hom_{G}(\mathbf{C}^{x}_{I_0,I_1}, \mathbf{C}^{x}_{I_2,I_3})$ is $1$-dimensional and is spanned by the composition of the truncation maps
\[\mathbf{C}^{x}_{I_0,I_1}\rightarrow \mathbf{C}^{x}_{I_2,I_1} \rightarrow \mathbf{C}^{x}_{I_2,I_3}\]
\end{rem}

Let $I_2\subseteq I_0\subseteq \Delta$ with $x\in\Gamma^{I_0\setminus I_2}$ (namely $J_{x}\supseteq J_2\defeq I_2\sqcup(\Delta\setminus I_0)$). We write $r\defeq \#\Delta=n-1$ for short. Recall from (\ref{equ: naive Tits x w transfer}) that $i_{I_2}^{\rm{an}}\rightarrow i_{x,I_2}^{\rm{an}}$ induces a map $\mathbf{C}_{I_2,\Delta}\rightarrow\mathbf{C}^{x}_{I_2,J_{x}}$. This together with standard truncation maps between various (variant of) Tits complex give the following commutative diagram
\begin{equation}\label{equ: Ext Tits transfer x}
\xymatrix{
\mathrm{Ext}_{G}^{\bullet}(\mathbf{C}_{I_0,\Delta},\mathbf{C}_{I_2,\Delta}) \ar[r] \ar^{\wr}[d] & \mathrm{Ext}_{G}^{\bullet}(\mathbf{C}_{I_0,\Delta}^{\infty},\mathbf{C}_{I_2,\Delta}) \ar[r] \ar^{\wr}[d] & \mathrm{Ext}_{G}^{\bullet}(\mathbf{C}_{I_0,\Delta}^{\infty},\mathbf{C}^{x}_{I_2,J_{x}}) \ar^{\wr}[d]\\
\mathrm{Ext}_{G}^{\bullet}(\mathbf{C}_{I_0,\Delta},\mathbf{C}_{J_2,\Delta}) \ar[r] & \mathrm{Ext}_{G}^{\bullet}(\mathbf{C}_{I_0,\Delta}^{\infty},\mathbf{C}_{J_2,\Delta}) \ar[r] & \mathrm{Ext}_{G}^{\bullet}(\mathbf{C}_{I_0,\Delta}^{\infty},\mathbf{C}^{x}_{J_2,J_{x}})\\
\mathrm{Ext}_{G}^{\bullet-r}(1_{G},\mathbf{C}_{J_2,\Delta}) \ar@{=}[r] \ar_{\wr}[u] & \mathrm{Ext}_{G}^{\bullet-r}(1_{G},\mathbf{C}_{J_2,\Delta}) \ar[r] \ar_{\wr}[u]& \mathrm{Ext}_{G}^{\bullet-r}(1_{G},\mathbf{C}^{x}_{J_2,J_{x}}) \ar_{\wr}[u]
}
\end{equation}
with all vertical maps in the first column (resp.~in the second and third column) being isomorphisms by \ref{it: Ext complex std seq} of Proposition~\ref{prop: Ext complex std seq} with $x=1$ in \emph{loc.cit.} (resp.~by \ref{it: Ext complex std seq sm} of Proposition~\ref{prop: Ext complex std seq}).
Using $\mathbf{C}_{I,\Delta}\cong V_{I}^{\rm{an}}[\#I]$ and $\mathbf{C}_{I,\Delta}^{\infty}\cong V_{I}^{\infty}[\#I]$ for each $I\subseteq\Delta$ and $\mathbf{C}^{x}_{I',J_{x}}\cong V_{x,I'}^{\rm{an}}[\#I']$ for each $I'\subseteq J_{x}$, we can rewrite (\ref{equ: Ext Tits transfer x}) as follows
\begin{equation}\label{equ: Ext Tits transfer x prime}
\xymatrix{
\mathrm{Ext}_{G}^{\bullet}(V_{I_0}^{\rm{an}},V_{I_2}^{\rm{an}}) \ar^{\sim}[r] \ar^{\wr}[d] & \mathrm{Ext}_{G}^{\bullet}(V_{I_0}^{\infty},V_{I_2}^{\rm{an}}) \ar[r] \ar^{\wr}[d] & \mathrm{Ext}_{G}^{\bullet}(V_{I_0}^{\infty},V_{x,I_2}^{\rm{an}}) \ar^{\wr}[d]\\
\mathrm{Ext}_{G}^{\bullet+\#\Delta\setminus I_0}(V_{I_0}^{\rm{an}},V_{J_2}^{\rm{an}}) \ar^{\sim}[r] & \mathrm{Ext}_{G}^{\bullet+\#\Delta\setminus I_0}(V_{I_0}^{\infty},V_{J_2}^{\rm{an}}) \ar[r] & \mathrm{Ext}_{G}^{\bullet+\#\Delta\setminus I_0}(V_{I_0}^{\infty},V_{x,J_2}^{\rm{an}})\\
\mathrm{Ext}_{G}^{\bullet}(1_{G},V_{J_2}^{\rm{an}}) \ar@{=}[r] \ar_{\wr}[u] & \mathrm{Ext}_{G}^{\bullet}(1_{G},V_{J_2}^{\rm{an}}) \ar[r] \ar_{\wr}[u]& \mathrm{Ext}_{G}^{\bullet}(1_{G},V_{x,J_2}^{\rm{an}}) \ar_{\wr}[u]
}
\end{equation}
Note that the horizontal maps from the first column to the second column of (\ref{equ: Ext Tits transfer x prime}) are isomorphisms by Lemma~\ref{lem: change left cup}.
We will see in Proposition~\ref{prop: coxeter filtration x surjection} that the horizontal maps from the second column to the third column of (\ref{equ: Ext Tits transfer x prime}) are actually surjections.

\section{Combinatorics of Tits complex}\label{sec: comb Tits}
Let $I_1\subseteq I_0\subseteq \Delta$ be a pair.
We study the spectral sequence $E_{\bullet,I_0,I_1}^{\bullet,\bullet}$ associated with the double complex $\cT_{I_0,I_1}^{\bullet,\bullet}$ introduced in \S \ref{subsec: std seq}.
We develop a combinatorial technique to compute the second page of $E_{\bullet,I_0,I_1}^{\bullet,\bullet}$ (see Proposition~\ref{prop: atom basis}) and then prove that it degenerates at the second page for the bottom non-vanishing degree (see Proposition~\ref{prop: bottom deg degeneracy}).
We also study several variants of $\cT_{I_0,I_1}^{\bullet,\bullet}$ all of which have been introduced in \S \ref{subsec: std seq} (see \S \ref{subsec: truncation atom} and \S \ref{subsec: relative Tits}).
\subsection{Basis of group cohomology}\label{subsec: std bases}
We construct a standard basis for $H^{\bullet}(L_{I},1_{L_{I}})$ for each $I\subseteq \Delta$, using primitive classes in $H^{\bullet}(\fs\fl_{m,K},1_{\fs\fl_{m,K}})$ for $m\geq 2$ and the basis $\{\log,\val\}$ of $\Hom(K^\times,E)$.
Given $I'\subseteq I\subseteq\Delta$, we also describe the restriction map 
\[H^{\bullet}(L_{I},1_{L_{I}})\rightarrow H^{\bullet}(L_{I'},1_{L_{I'}})\]
under such standard bases of the source and the target (see Proposition~\ref{prop: Levi restriction}).

Let $\fg$ be a (split) reductive $E$-Lie algebra. It follows from \cite[Thm.~9.3]{Kos50} that we have an isomorphism of graded $E$-algebra
\[H_{\bullet}(\fg,1_{\fg})\cong (\wedge^{\bullet}\fg)^\fg.\]
Let $D(\fg)$ be the subspace of $(\wedge^{\bullet}\fg)^\fg$ spanned by elements of the form $u\wedge v$ (namely decomposable). Note that we have a pairing $H_{\bullet}(\fg,1_{\fg})\times H^{\bullet}(\fg,1_{\fg})\rightarrow E$ and let $P(\fg)$ be the subspace of $H^{\bullet}(\fg,1_{\fg})$ given by the orthogonal complement of $D(\fg)$. We set $P^k(\fg)\defeq P(\fg)\cap H^k(\fg,1_{\fg})$ for each $k\geq 0$.
We recall the following classical result from \cite[Thm.~10.2, Thm.~10.3]{Kos50}.
\begin{thm}\label{thm: Koszul thm}
We have the following isomorphism between graded $E$-algebras
\begin{equation}\label{equ: Koszul isom}
\wedge(P(\fg)^\prime)\cong  H^{\bullet}(\fg,1_{\fg}).
\end{equation}
Moreover, (\ref{equ: Koszul isom}) is functorial in the following sense: for each morphism between (split) reductive $E$-Lie algebra $\varphi:~\fh\rightarrow\fg$, we have an induced morphism $\varphi^\ast:~P(\fg)\rightarrow P(\fh)$ which determines the morphism $\varphi^\ast:~ H^{\bullet}(\fg,1_{\fg})\rightarrow  H^{\bullet}(\fh,1_{\fh})$ completely, via (\ref{equ: Koszul isom}).
\end{thm}

For each $m\geq 1$, we write $\mathbf{G}_{m}\defeq \mathrm{PGL}_{m/K}$, $G_{m}\defeq \mathbf{G}_{m}(K)$ and $\fg_{m}\defeq \mathrm{Lie}(G_{m})\otimes_KE=\fs\fl_{m,K}\otimes_KE$.
For each $2\leq m'\leq m$, we have a standard embedding 
\begin{equation}\label{equ: std Levi block embedding}
\iota_{m',m}: \fg_{m'}\hookrightarrow \fg_{m}
\end{equation}
by mapping $\fg_{m'}$ into the left-upper block of $\fg_m$.
We fix $n\geq 1$ in the following discussion and write $G\defeq G_n$ as well as $\fg\defeq \fg_{n}$.
The following results are classical (cf.~\cite{CE48} and \cite{Kos50}).
\begin{thm}\label{thm: primitive class}
Let $2\leq m'\leq m$. We have the following results.
\begin{enumerate}[label=(\roman*)]
\item \label{it: primitive class 1} We have $\dim_EP^k(\fg_m)=1$ if $k=2a-1$ for some $2\leq a\leq m$ and $P^k(\fg_m)=0$ otherwise.
\item \label{it: primitive class 2} The standard embedding $\iota_{m',m}$ from (\ref{equ: std Levi block embedding}) induces an isomorphism $P^k(\fg_{m})\buildrel\sim\over\longrightarrow P^k(\fg_{m'})$ for each $k\leq 2m'-1$ and zero map otherwise.
\item \label{it: primitive class 3} Upon identifying $\fg_{m'}$ with an $E$-Lie subalgebra of $\fg_{m}$ via $\iota_{m',m}$, there exists a subspace $P(\fg_{m},\fg_{m'})\subseteq H^{\bullet}(\fg_{m},\fg_{m'},1_{\fg_{m}})$ such that
    \[H^{\bullet}(\fg_{m},\fg_{m'},1_{\fg_{m}})\cong\wedge P(\fg_{m},\fg_{m'})\]
    and that the natural map
    \begin{equation}\label{equ: relative Lie change}
    H^{\bullet}(\fg_{m},\fg_{m'},1_{\fg_{m}})\rightarrow H^{\bullet}(\fg_{m},1_{\fg_{m}})
    \end{equation}
    is the embedding induced (via wedge product) from an embedding
    \begin{equation}\label{equ: relative Lie change primitive}
    P(\fg_{m},\fg_{m'})\hookrightarrow P(\fg_{m})
    \end{equation}
    whose image is $\bigoplus_{a=m'+1}^{m}P^{2a-1}(\fg_{m})$. In particular, if $m'=m-1$, then we have \[H^{\bullet}(\fg_{m},\fg_{m'},1_{\fg_{m}})\buildrel\sim\over\longrightarrow\wedge P^{2m-1}(\fg_{m})\] 
    via the embedding (\ref{equ: relative Lie change}).
\end{enumerate}
\end{thm}
\begin{proof}
All \ref{it: primitive class 1}, \ref{it: primitive class 2} and \ref{it: primitive class 3} are well-known. Thanks to \cite[Thm.~10.2, Thm.~10.3]{Kos50}, we reduce ourselves to prove parallel results for $\fs\fl_{m,\Q}$ and $\fs\fl_{m',\Q}$ (replacing $\fg_m$ and $\fg_{m'}$) and further to the case for $\fs\fl_{m,\R}$ and $\fs\fl_{m',\R}$. Then by \cite[Thm.~15.1, Thm.~15.2]{CE48} it suffices to prove corresponding results for the following commutative diagram of maps
\begin{equation}\label{equ: primitive class top diagram}
\xymatrix{
H^{\bullet}_{\rm{dR}}(SU(m)/SU(m'),\R) \ar[r] \ar^{\wr}[d] & H^{\bullet}_{\rm{dR}}(SU(m),\R) \ar[r] \ar^{\wr}[d] & H^{\bullet}_{\rm{dR}}(SU(m'),\R) \ar^{\wr}[d]\\
H^{\bullet}_{\rm{sing}}(SU(m)/SU(m'),\R) \ar[r] & H^{\bullet}_{\rm{sing}}(SU(m),\R) \ar[r] & H^{\bullet}_{\rm{sing}}(SU(m'),\R)
}
\end{equation}
with the vertical maps being isomorphisms from classical comparison theorems.
Note that $H^{\bullet}_{\rm{sing}}(SU(m)/SU(m'),\R)$ can be computed using the Leray spectral sequence associated with the fiber bundle $SU(m)/SU(m')\rightarrow SU(m)/SU(m-1)\cong S^{2m-1}$ (with fiber $SU(m-1)/SU(m')$) for increasing $m$ (with $m>m'$ and $SU(m')$ being the trivial group when $m'=1$), the desired claims on the bottom row of (\ref{equ: primitive class top diagram}) follow from a comparison between these spectral sequences.
\end{proof}

\begin{rem}\label{rem: primitive class G inv}
Now that $H^{\bullet}(\fg_{m},1_{\fg_{m}})$ is $\mathrm{Ad}(g)$-invariant for each $g\in G_{m}$, \ref{it: primitive class 2} and \ref{it: primitive class 3} of Theorem~\ref{thm: primitive class} remain true if we replace $\iota_{m',m}$ with 
\begin{equation}\label{equ: conj std embedding}
\mathrm{Ad}(g)^{\ast}(\iota_{m',m}): \mathrm{Ad}(g)(\fg_{m'})\rightarrow \fg_{m}
\end{equation}
for some $g\in G_{m}$. Note that maps of the form (\ref{equ: conj std embedding}) including all embeddings of $\fg_{m'}$ into $\fg_{m}$ as an arbitrary standard Levi block (not necessarily the left-upper block in our definition of $\iota_{m',m}$).
\end{rem}

For each $I\subseteq\Delta$ and each $1\leq d\leq r_I$, we have an embedding
\begin{equation}\label{equ: std Levi factor}
\fg_{n_d}\hookrightarrow \fl_{I}\hookrightarrow \fg=\fg_{n}
\end{equation}
which induces the following map
\begin{equation}\label{equ: restriction to Levi}
\mathrm{Res}_{n,I}^{d}:~P(\fg)=P(\fg_{n})\rightarrow P(\fg_{n_d}).
\end{equation}
Similar argument as in Remark~\ref{rem: primitive class G inv} shows that the map (\ref{equ: restriction to Levi}) equals $\iota_{n_{d},n}^{\ast}$.

Let $m\geq 2$ and $\tau_{m}$ be the involution of $G_{m}$ given by $A\mapsto (A^{t})^{-1}$. Note that $\tau_{m}$ induces an involution of $\fg_{m}$ and thus an involution of $P(\fg_{m})$ as well as $H^{k}(\fg_{m},1_{\fg_{m}})$ for each $k\geq 0$ by Theorem~\ref{thm: Koszul thm}. We abuse the notation $\tau_{m}$ for these induced actions.
\begin{lem}\label{lem: primitive class outer}
For each $m\geq 2$, $\tau_{m}$ acts on $P^{2m-1}(\fg_{m})$ by $(-1)^{m}$.
\end{lem}
\begin{proof}
We fix a choice of $m\geq 2$ throughout the proof. We write $B_{m}^+$ (resp.~$T_{m}$) for the upper-triangular Borel subgroup (resp.~diagonal maximal torus) of $G_{m}$ with associated $E$-Lie subalgebra $\fb_{m}^{+}$ (resp.~$\ft_{m}$) and the set of positive roots $\Phi_{m}^+$. We similarly have $\fb_{m-1}$, $\ft_{m-1}$ and $\Phi_{m-1}^+$.
We embed $G_{m-1}$ in $G_{m}$ via $\iota_{m-1,m}$. Note that the involution $\tau_{m}$ of $G_{m}$ stablizes $G_{m-1}$ and restricts to the involution $\tau_{m-1}$ on $G_{m-1}$. Hence, $\tau_{m}$ induces an involution on $P(\fg_{m},\fg_{m-1})$ and $H^{k}(\fg_{m},\fg_{m-1},1_{\fg_{m}})$ for each $k\geq 0$, and by \ref{it: primitive class 3} of Theorem~\ref{thm: primitive class} we have an isomorphism
\begin{equation}\label{equ: primitive class outer}
P^{2m-1}(\fg_{m},\fg_{m-1})\buildrel\sim\over\longrightarrow P^{2m-1}(\fg_{m})
\end{equation}
between $1$-dimensional $E$-vector spaces which is compatible with the involution $\tau_{m}$ on BHS. 
Note that we have 
\[\fg_{m}/\fg_{m-1}=\ft_{m}/\ft_{m-1}\oplus\bigoplus_{\al\in\Phi_{m}^+\setminus\Phi_{m-1}^{+}}\fu_{\al}\oplus\fu_{-\al}\] 
with the transpose exchanging $\fu_{\al}$ and $\fu_{-\al}$ for each $\al\in\Phi_{m}^+\setminus\Phi_{m-1}^{+}$.
Now that $\tau_{m}$ acts on $\fg_{m}/\fg_{m-1}$ by the composition of transpose with $-\mathrm{Id}_{\fg_{m}/\fg_{m-1}}$, we conclude by the observation that
\[P^{2m-1}(\fg_{m},\fg_{m-1})=\wedge^{2m-1}(\fg_{m}/\fg_{m-1})=(\ft_{m}/\ft_{m-1})\wedge\bigwedge_{\al\in\Phi_{m}^+\setminus\Phi_{m-1}^{+}}(\fu_{\al}\wedge \fu_{-\al}).\]
\end{proof}

For each odd integer $3\leq k\leq 2m-1$, we choose an element $0\neq v_{m}^{k}\in P^k(\fg_{m})$ (and thus $P^k(\fg_{m})=E v_{m}^{k}$, see \ref{it: primitive class 1} of Theorem~\ref{thm: primitive class}). We also write $v_{m}^{0}\defeq 1\in  H^{0}(\fg_{m},1_{\fg_{m}})$ for convenience. We define $\Sigma_{m,k}$ as the set of subsets
\[\Lambda=\{m_1,\dots, m_r\}\subseteq \{3,5,\dots,2m-3,2m-1\}\]
for some $r\geq 1$ such that $\sum_{s=1}^rm_s=k$. 
Throughout the context, we may and do choose $v_{m}^{k}$ for varying $m$ such that
\[\iota_{m',m}^{\ast}(v_{m}^{k})=v_{m'}^{k}\]
for each $2\leq m'\leq m$ and each odd integer $3\leq k\leq 2m'-1$.
For each $\Lambda\in\Sigma_{m,k}$, we define $P^{\Lambda}$ as the one dimensional subspace of $\wedge P(\fg_{m})$ given by
\[P^{m_1}(\fg_{m})\wedge\cdots\wedge P^{m_r}(\fg_{m}),\]
which clearly does not depend on the choice of order on $\Lambda$. We write $\tld{\Lambda}$ for an ordered set $\{m_1,\cdots,m_r\}$ whose underlying set is $\Lambda$, and then set $v_{m}^{\tld{\Lambda}}\defeq v^{m_1}\wedge\cdots\wedge v^{m_r}$. So different $v_{m}^{\tld{\Lambda}}$ (for different choices of ordered sets $\tld{\Lambda}$ with the same underlying set $\Lambda$) differ by an explicit sign. Our convention includes $\Sigma_{m,0}\defeq \{\emptyset\}$ and $v_{m}^{\emptyset}\defeq v_{m}^{0}=1\in H^{0}(\fg_{m},1_{\fg_{m}})$.
Without further comment, we fix the choice of the order on $\tld{\Lambda}$ to be the one with decreasing order on integers, and write $v_{m}^{\Lambda}$ instead of $v_{m}^{\tld{\Lambda}}$.
\begin{cor}\label{cor: decomposition of coh}
We have natural isomorphisms
\begin{equation}\label{equ: canonical decomposition}
H^k(\fg_{m},1_{\fg_{m}})\cong \bigoplus_{\Lambda\in\Sigma_{m,k}}P^{\Lambda}
\end{equation}
for each $k\geq 0$. In particular, $H^k(\fg_{m},1_{\fg_{m}})=0$ if $k>(n^2-1)$.
\end{cor}
\begin{proof}
The isomorphism (\ref{equ: canonical decomposition}) follows from Theorem~\ref{thm: Koszul thm} and Theorem~\ref{thm: primitive class}. The vanishing is clear from $\dim \fg_{m}=(m^2-1)$, and can also be seen from the fact that
\[m^2-1=\sum_{s=2}^m2s-1\]
which implies that $\Sigma_{m,k}=\emptyset$ when $k>(m^2-1)$.
\end{proof}

Let $I\subseteq \Delta$. We recall from the discussion around (\ref{equ: explicit complement}) that we have $r_{I}=n-\#I$, $I=\bigsqcup_{d=1}^{r_{I}}I^{d}$ and $n_{d}=\#I^{d}+1$ with $\fh_{I^{d}}\cong\fg_{n_{d}}$ for each $1\leq d\leq r_{I}$, so that we have the following isomorphisms between $E$-Lie algebras
\begin{equation}\label{equ: adjoint isom g}
\fl_{I}\cong\fz_{I}\times\fh_{I}\cong \fz_{I}\times\prod_{d=1}^{r_{I}}\fh_{I^{d}}\cong\fz_{I}\times\prod_{d=1}^{r_{I}}\fg_{n_{d}}.
\end{equation}
For each subinterval $J\subseteq I$ with $\fh=\fh_{J}$, the isomorphism (\ref{equ: adjoint isom g}) restricts to an isomorphism
\[\fl_{I}\cap\fh\cong\fh_{I\cap J}\cong \prod_{d=1}^{r_{I}}\fh_{I^{d}\cap J}\cong \prod_{d=1}^{r_{I}}\fg_{1+\#I^{d}\cap J}.\]
By applying the K\"unneth formula for relative Lie algebra cohomology (cf.~(\ref{equ: relative Lie Kunneth})) to (\ref{equ: adjoint isom g}), we obtain the following isomorphism between graded $E$-algebras
\begin{multline}\label{equ: ss group coh g}
H^{\bullet}(\fl_{I},\fh,1_{\fl_{I}})
\cong \wedge^{\bullet}\Hom(\fz_{I},E)\otimes_E H^{\bullet}(\fh_{I},\fh,1_{\fh_{I}})\\
\cong \wedge^{\bullet}\Hom(\fz_{I},E)\otimes_E \bigotimes_{d=1}^{r_{I}}H^{\bullet}(\fh_{I^{d}},\fh_{I^{d}\cap J},1_{\fh_{I^{d}}})
\cong \wedge^{\bullet}\Hom(\fz_{I},E)\otimes_E \bigotimes_{d=1}^{r_{I}}H^{\bullet}(\fg_{n_{d}},\fg_{1+\#I^{d}\cap J},1_{\fg_{n_{d}}}).
\end{multline}
which together with \ref{it: Levi decomposition 1} of Proposition~\ref{prop: Levi decomposition} gives the following isomorphisms between graded algebras
\begin{multline}\label{equ: Levi Kunneth}
H^{\bullet}(L_{I},\fh,1_{L_{I}})\cong \wedge^{\bullet}\Hom(Z_{I}^{\dagger},E)\otimes_E H^{\bullet}(\fl_{I},\fh,1_{\fl_{I}})\\
\cong\wedge^{\bullet}\Hom(Z_{I}^{\dagger},E)\otimes_E\wedge^{\bullet}\Hom(\fz_{I},E)\otimes_E H^{\bullet}(\fh_{I},\fh,1_{\fh_{I}})\\
\cong \wedge^{\bullet}\Hom(Z_{I},E)\otimes_E\bigotimes_{d=1}^{r_I} H^{\bullet}(\fg_{n_d},\fg_{1+\#I^{d}\cap J},1_{\fg_{n_d}}).
\end{multline}

Given $I\subseteq \Delta$ and $\un{k}=\{k_d\}_{0\leq d\leq r_I}$ a tuple of non-negative integers satisfying $|\un{k}|\defeq \sum_{d=0}^{r_I}k_d=k$. We choose an element $\Lambda_d\in\Sigma_{n_d,k_d}$ for each $1\leq d\leq r_I$ and then write $\un{\Lambda}=\{\Lambda_d\}_{1\leq d\leq r_I}$ for the tuple. We set
\begin{equation}\label{equ: Lambda index}
\Sigma_{I,\un{k}}\defeq \prod_{d=1}^{r_I}\Sigma_{n_d,k_d}
\end{equation}
and
\begin{equation}\label{equ: std basis tensor}
v_{I,\un{k}}^{\un{\Lambda}}\defeq \bigotimes_{d=1}^{r_I}v_{n_d}^{\Lambda_d}\in\bigotimes_{d=1}^{r_I} H^{k_d}(\fg_{n_d},1_{\fg_{n_d}})
\end{equation}
for each $\un{\Lambda}\in\Sigma_{I,\un{k}}$, with the tensor taken over increasing $d$.
It is clear that $\{v_{I,\un{k}}^{\un{\Lambda}}\}_{\un{\Lambda}\in\Sigma_{I,\un{k}}}$ is a basis of $\bigotimes_{d=1}^{r_I} H^{k_d}(\fg_{n_d},1_{\fg_{n_d}})$. By (\ref{equ: Levi Kunneth}) and Corollary~\ref{cor: decomposition of coh} we know that
\begin{equation}\label{equ: std basis}
\{v\otimes_Ev_{I,\un{k}}^{\un{\Lambda}}\mid v\subseteq \mathbf{Log}_{I},~ \#v=k_0,~ \un{\Lambda}\in\Sigma_{I,\un{k}},~|\un{k}|=k\}
\end{equation}
forms a basis of $H^k(L_{I},1_{L_{I}})$ for each $I\subseteq \Delta$ and $k\geq 0$.

Now that a basis of $H^{\bullet}(L_{I},1_{L_{I}})$ has been constructed for each $I$, we need to describe the restriction maps using these bases.
\begin{prop}\label{prop: Levi restriction}
For each $I'\subseteq I$, the restriction map 
\[\mathrm{Res}^{\bullet}_{I,I'}: H^{\bullet}(L_{I},1_{L_{I}})\rightarrow H^{\bullet}(L_{I'},1_{L_{I'}})\] 
corresponds under (\ref{equ: Levi Kunneth}) to the map
\begin{equation}\label{equ: Kunneth decomposition K restriction}
\mathrm{res}^{\bullet}_{I,I'}: \mathrm{Tot}(\wedge^{\bullet}\Hom(Z_{I}^{\dagger},E)\otimes_E H^{\bullet}(\fl_{I},1_{\fl_{I}})) \rightarrow \mathrm{Tot}(\wedge^{\bullet}\Hom(Z_{I'}^{\dagger},E)\otimes_E H^{\bullet}(\fl_{I'},1_{\fl_{I'}}))
\end{equation}
induced from the inclusion $\fl_{I'}\subseteq \fl_{I}$ and the map $Z_{I'}^{\dagger}\rightarrow Z_{I}^{\dagger}$.
\end{prop}
\begin{proof}
Using the equalities
\[\mathrm{Res}^{\bullet}_{I,I''}=\mathrm{Res}^{\bullet}_{I,I'}\circ \mathrm{Res}^{\bullet}_{I',I''}\]
and
\[\mathrm{res}^{\bullet}_{I,I''}=\mathrm{res}^{\bullet}_{I,I'}\circ \mathrm{res}^{\bullet}_{I',I''}\]
for arbitrary $I''\subseteq I'\subseteq I$, it is harmless to assume throughout this proof that $\#I\setminus I'=1$.
We write $H'\defeq L_{I'}'\cap H_{I}$ and $H^{\dagger}\defeq Z_{I'}^{\dagger}\cap H_{I}$ for short with $H'\times H^{\dagger}\rightarrow L_{I'}\cap H_{I}$ being an isogeny and $H^{\dagger}$ being a discrete abelian group which is free of rank $1$ (by our assumption that $\#I\setminus I'=1$).
Using the quasi-isomorphism (\ref{equ: Levi Lie discrete center factor}) which is functorial with respect to the choice of $I$ in \emph{loc.cit.}, we only need to establish a commutative diagram of the form
\begin{equation}\label{equ: Levi restriction flat diagram}
\xymatrix{
H^{\bullet}(H_{I},1_{H_{I}}) \ar^{\mathrm{Res}^{\bullet}_{I,I',\flat}}[rrr] \ar^{\wr}[d] & & & H^{\bullet}(H'\times H^{\dagger},1_{H'\times H^{\dagger}}) \ar^{\wr}[d]\\
H^{\bullet}(\fh_{I},1_{\fh_{I}}) \ar^{\mathrm{res}^{\bullet}_{I,I',\flat}}[rrr] & & & \wedge^{\bullet}\Hom(H^{\dagger},E)\otimes_E H^{\bullet}(\fl_{I'}\cap\fh_{I},1_{\fl_{I'}\cap\fh_{I}})
}
\end{equation}
with both vertical isomorphisms from evident variant of Proposition~\ref{prop: Levi decomposition}, and the bottom horizontal map being the tensor of the restriction map $H^{\bullet}(\fh_{I},1_{\fh_{I}})\rightarrow H^{\bullet}(\fl_{I'}\cap\fh_{I},1_{\fl_{I'}\cap\fh_{I}})$ (induced from the inclusion $\fl_{I'}\cap\fh_{I}\subseteq \fh_{I}$) and the evident map $E\hookrightarrow \wedge^{\bullet}\Hom(H^{\dagger},E)$.
Similar to (\ref{equ: ss group coh g}) and (\ref{equ: Levi Kunneth}), we have
\[H^{\bullet}(H_{I},1_{H_{I}})\buildrel\sim\over\longrightarrow H^{\bullet}(\fh_{I},1_{\fh_{I}})=\wedge P(\fh_{I})\]
with $P(\fh_{I})=\bigoplus_{d=1}^{r_{I}}P(\fh_{I^{d}})$. It is thus sufficient to show that
\begin{equation}\label{equ: Levi restriction flat primitive}
\mathrm{Res}^{\bullet}_{I,I',\flat}|_{P^{2m-1}(\fh_{I^{d}})}=\mathrm{res}^{\bullet}_{I,I',\flat}|_{P^{2m-1}(\fh_{I^{d}})}.
\end{equation}
for each $1\leq d\leq r_{I}$ with $I^{d}\neq \emptyset$ (namely $n_{d}\geq 2$) and $2\leq m\leq n_{d}$.

We fix from now such a choice of $d$ and $m$ and prove (\ref{equ: Levi restriction flat primitive}).\\
For each subinterval $J\subseteq I$ with $\fh_{I'}\cap\fh_{J}=\fh_{I'\cap J}$, we define
\begin{multline*}
H^{\bullet}_{J}\defeq H^{\bullet}(H'\times H^{\dagger},\fh_{I'\cap J},1_{H'\times H^{\dagger}})\\
\buildrel\sim\over\longrightarrow \mathrm{Tot}(\wedge^{\bullet}\Hom(H^{\dagger},E)\otimes_E \wedge^{\bullet}\Hom(\fz_{I'}\cap\fh_{I},E)\otimes_E \bigotimes_{d'=1}^{r_{I'}}H^{\bullet}(\fh_{(I')^{d'}},\fh_{(I')^{d'}\cap J},1_{(I')^{d'}})),
\end{multline*}
and omit $J$ from the notation when $J=\emptyset$.
For each $1\leq d'\leq r_{I'}$, we have either $(I')^{d'}\subseteq I^{d}$ or $(I')^{d'}\cap I^{d}=\emptyset$, and write $H^{\bullet}_{d'}\defeq H^{\bullet}_{(I')^{d'}}$ for short.
For each subinterval $J\subseteq I$, we clearly have a commutative diagram of the form
\[
\xymatrix{
H^{\bullet}(H_{I},\fh_{J},1_{H_{I}}) \ar^{\varphi_{J}}[rr] \ar@{^{(}->}[d] & & H^{\bullet}_{J} \ar@{^{(}->}[d]\\
H^{\bullet}(H_{I},1_{H_{I}}) \ar^{\mathrm{Res}^{\bullet}_{I,I',\flat}}[rr] & & H^{\bullet}
}
\]
with both vertical maps being embeddings by variants of (\ref{equ: Levi Kunneth}) as well as \ref{it: primitive class 3} of Theorem~\ref{thm: primitive class}.
If $m>\#I^{d}\cap J+1$, then by \ref{it: primitive class 3} of Theorem~\ref{thm: primitive class} we have
\[P^{2m-1}(\fh_{I^{d}},\fh_{I^{d}}\cap\fh_{J})\buildrel\sim\over\longrightarrow P^{2m-1}(\fh_{I^{d}}),\]
and thus
\[\mathrm{Res}^{2m-1}_{I,I',\flat}(P^{2m-1}(\fh_{I^{d}}))=\varphi_{J}(P^{2m-1}(\fh_{I^{d}},\fh_{I^{d}\cap J}))\subseteq H^{2m-1}_{J}\subseteq H^{2m-1}.\]
We have the following two possibilities.
\begin{itemize}
\item If either $(I')^{d'}\cap I^{d}=\emptyset$ or $\#(I')^{d'}+1<m$, then we can take $J$ to be $(I')^{d'}$ which satisfies $m>\#I^{d}\cap J+1$ and thus
\begin{equation}\label{equ: Kunneth H bound 1}
\mathrm{Res}^{2m-1}_{I,I',\flat}(P^{2m-1}(\fh_{I^{d}}))\subseteq H^{2m-1}_{J}=H^{2m-1}_{d'}\subseteq H^{2m-1}.
\end{equation}
\item If $(I')^{d'}\subseteq I^{d}$ and $\#(I')^{d'}+1\geq m$, then we can choose $J\subsetneq (I')^{d'}\subseteq I^{d}$ so that $\#J=m-2$ (and thus $m>\#J+1=\#I^{d}\cap J+1$). Thanks to \ref{it: primitive class 3} of Theorem~\ref{thm: primitive class}, we deduce from $J\subsetneq (I')^{d'}$ and $\#J=m-2$ that
    \[H^{2m-1}(\fh_{{I'}^{d'}},\fh_{J},1_{\fh_{{I'}^{d'}}})=P^{2m-1}(\fh_{{I'}^{d'}})\]
    and $H^{k}(\fh_{{I'}^{d'}},\fh_{J},1_{\fh_{{I'}^{d'}}})=0$ for each $0<k<2m-1$. Consequently, we have
\begin{multline}\label{equ: Kunneth H bound 2}
\mathrm{Res}^{2m-1}_{I,I',\flat}(P^{2m-1}(\fh_{I^{d}}))\subseteq H^{2m-1}_{J}\\
=\bigoplus_{k=0}^{2m-1}H^{k}(\fh_{{I'}^{d'}},\fh_{J},1_{\fh_{{I'}^{d'}}})\otimes_EH^{2m-1-k}_{d'}\subseteq P^{2m-1}(\fh_{{I'}^{d'}})\oplus H^{2m-1}_{d'}\subseteq H^{2m-1}.
\end{multline}
\end{itemize}
We use the convention $P^{2m-1}(\fh_{{I'}^{d'}})=0$ for each $1\leq d'\leq r_{I'}$ that satisfies $\#(I')^{d'}+1<m$.
We set 
\[H^{2m-1}_{\natural}\defeq \bigcap_{d',(I')^{d'}\cap I^{d}=\emptyset}H^{2m-1}_{d'}\] 
and note that $P^{2m-1}(\fh_{{I'}^{d'}})\subseteq H^{2m-1}_{\natural}$ for each $d'$ that satisfies $(I')^{d'}\subseteq I^{d}$. It thus follows from (\ref{equ: Kunneth H bound 1}) and (\ref{equ: Kunneth H bound 2}) that
\begin{multline}\label{equ: Kunneth H bound intersection}
\mathrm{Res}^{2m-1}_{I,I',\flat}(P^{2m-1}(\fh_{I^{d}}))\subseteq H^{2m-1}_{\natural}\cap\bigcap_{d',(I')^{d'}\subseteq I^{d}}(P^{2m-1}(\fh_{{I'}^{d'}})\oplus H^{2m-1}_{d'})\\
=\bigcap_{d',(I')^{d'}\subseteq I^{d}}(P^{2m-1}(\fh_{{I'}^{d'}})\oplus (H^{2m-1}_{d'}\cap H^{2m-1}_{\natural})).
\end{multline}
If $I\setminus I'\subseteq I^{d}$, then we know that $(I')^{d'}\subseteq I^{d}$ if and only if $d'\in\{d,d+1\}$, with $P^{2m-1}(\fh_{{I'}^{d}})\subseteq H^{2m-1}_{d+1}\cap H^{2m-1}_{\natural}$, $P^{2m-1}(\fh_{{I'}^{d+1}})\subseteq H^{2m-1}_{d}\cap H^{2m-1}_{\natural}$, $P^{2m-1}(\fh_{{I'}^{d}})\cap P^{2m-1}(\fh_{{I'}^{d+1}})=0$, and
\begin{multline*}
(H^{2m-1}_{d}\cap H^{2m-1}_{\natural})\cap(H^{2m-1}_{d+1}\cap H^{2m-1}_{\natural})=\bigcap_{d'=1}^{r_{I'}}H^{2m-1}_{d'}\\
=\mathrm{Tot}(\wedge^{\bullet}\Hom(H^{\dagger},E)\otimes_E\wedge^{\bullet}\Hom(\fz_{I'}\cap\fh_{I},E))^{2m-1}=0
\end{multline*}
since $2m-1\geq 3$ (with $m\geq 2$) yet both $\Hom(H^{\dagger},E)$ and $\Hom(\fz_{I'}\cap\fh_{I},E)$ are $1$-dimensional (by our assumption $\#I\setminus I'=1$). This together with (\ref{equ: Kunneth H bound intersection}) gives
\begin{multline}\label{equ: Kunneth H bound sum}
\mathrm{Res}^{2m-1}_{I,I',\flat}(P^{2m-1}(\fh_{I^{d}}))\subseteq H^{2m-1}_{\natural}\cap\bigcap_{d',(I')^{d'}\subseteq I^{d}}(P^{2m-1}(\fh_{{I'}^{d'}})\oplus H^{2m-1}_{d'})\\
=\bigoplus_{d',(I')^{d'}\subseteq I^{d}}P^{2m-1}(\fh_{{I'}^{d'}})\subseteq P^{2m-1}(\fl_{I'}\cap\fh_{I^{d}}).
\end{multline}
If $I^{d}\subseteq I'$, then there exists a unique $1\leq d'\leq r_{I'}$ such that $(I')^{d'}=I^{d}$ and a similar argument still gives (\ref{equ: Kunneth H bound sum}).
In other words, we have shown that $\mathrm{Res}^{2m-1}_{I,I',\flat}$ always restricts to a map
\[
P^{2m-1}(\fh_{I^{d}})\rightarrow P^{2m-1}(\fl_{I'}\cap\fh_{I^{d}}).
\]
This together with the following commutative diagram (induced from $U(\fh_{I^{d}})\subseteq U(\fh_{I})\subseteq D(H_{I})$)
\[
\xymatrix{
H^{\bullet}(H_{I},1_{H_{I}}) \ar[r] \ar[d] & H^{\bullet}(H'\times H^{\dagger},1_{H'\times H^{\dagger}}) \ar[d]\\
H^{\bullet}(\fh_{I^{d}},1_{\fh_{I^{d}}}) \ar[r] & H^{\bullet}(\fl_{I'}\cap\fh_{I^{d}},1_{\fl_{I'}\cap\fh_{I^{d}}})
}
\]
gives (\ref{equ: Levi restriction flat primitive}).
The proof is thus finished.
\end{proof}

Let $I\subseteq\Delta$. The isomorphisms between $E$-Lie algebras induce the following identifications
\[P(\fl_{I})=\Hom(\fz_{I},E)\oplus P(\fh_{I})=\Hom(\fz_{I},E)\oplus\bigoplus_{d=1}^{r_{I}}P(\fh_{I^{d}})=\Hom(\fz_{I},E)\oplus\bigoplus_{d=1}^{r_{I}}P(\fg_{n_{d}}).\]
Let $J\subseteq\Delta$ be an interval and $\fh\defeq\fh_{J}$.
We set
\begin{equation}\label{equ: explicit relative primitive}
P(\fl_{I},\fl_{I}\cap\fh)\defeq \Hom(\fz_{I\cup J},E)\oplus\bigoplus_{d=1}^{r_{I}}P(\fh_{I^{d}},\fh_{I^{d}\cap J})=\Hom(\fz_{I\cup J},E)\oplus\bigoplus_{d=1}^{r_{I}}\bigoplus_{a=\#I^{d}\cap J+2}^{n_{d}}P^{2a-1}(\fg_{n_{d}})\subseteq P(\fl_{I}).
\end{equation}
\begin{lem}\label{lem: relative interval embedding Lie}
Let $I\subseteq\Delta$ and $\fh=\fh_{J}$ for some interval $J\subseteq\Delta$.
Then there exists a natural isomorphism between graded $E$-algebras
\begin{equation}\label{equ: relative interval embedding Lie isom}
H^{\bullet}(\fl_{I},\fl_{I}\cap\fh,1_{\fl_{I}})\buildrel\sim\over\longrightarrow \wedge P(\fl_{I},\fl_{I}\cap\fh)
\end{equation}
which is functorial with respect to the choice of $I$ and $\fh$, such that the natural map
\begin{equation}\label{equ: relative interval embedding Lie}
H^{\bullet}(\fl_{I},\fl_{I}\cap\fh,1_{\fl_{I}}) \rightarrow H^{\bullet}(\fl_{I},1_{\fl_{I}})=\wedge P(\fl_{I})
\end{equation}
is the embedding between graded $E$-algebras given by the wedge product of the inclusion (\ref{equ: explicit relative primitive}).
\end{lem}
\begin{proof}
If $J\subseteq I$, then there exists a unique $1\leq d\leq r_{I}$ such that $J\subseteq I^{d}$, and the desired results follows directly from (\ref{equ: relative Lie Kunneth}). Hence, we assume in the rest of the proof that $J\not\subseteq I$.
We write $J_{+}$ for the unique subinterval of $I\cup J$ that contains $J$. The isomorphism between $E$-Lie algebras $\fh_{I\setminus J_{+}}\times\fz_{I\cup J}\times(\fl_{I}\cap\fh_{J_{+}})$ together with $\fl_{I}\cap \fh\subseteq \fl_{I}\cap\fh_{J_{+}}$ and (\ref{equ: relative Lie Kunneth}) gives the following commutative diagram of maps between graded $E$-algebras
\begin{equation}\label{equ: relative interval embedding Lie reduction}
\xymatrix{
H^{\bullet}(\fl_{I},\fl_{I}\cap\fh,1_{\fl_{I}}) \ar^{\sim}[rrr] \ar[d] & & & \mathrm{Tot}(H_{\natural}^{\bullet}\otimes_EH^{\bullet}(\fl_{I}\cap\fh_{J_{+}},\fl_{I}\cap \fh,1_{\fl_{I}\cap\fh_{J_{+}}})) \ar[d]\\
H^{\bullet}(\fl_{I},1_{\fl_{I}}) \ar^{\sim}[rrr] & & & \mathrm{Tot}(H_{\natural}^{\bullet}\otimes_EH^{\bullet}(\fl_{I}\cap\fh_{J_{+}},1_{\fl_{I}\cap\fh_{J_{+}}})) 
}
\end{equation}
with
\[H_{\natural}^{\bullet}\defeq H^{\bullet}(\fh_{I\setminus J_{+}}\times\fz_{I\cup J},1_{\fh_{I\setminus J_{+}}\times\fz_{I\cup J}})\cong \mathrm{Tot}(H^{\bullet}(\fh_{I\setminus J_{+}},1_{\fh_{I\setminus J_{+}}})\otimes_EH^{\bullet}(\fz_{I\cup J},1_{\fz_{I\cup J}}))\]
and the horizontal maps of (\ref{equ: relative interval embedding Lie reduction}) being isomorphisms from (\ref{equ: relative Lie Kunneth}). Using (\ref{equ: relative interval embedding Lie reduction}) and the equality 
\[P(\fl_{I},\fl_{I}\cap\fh)=P(\fh_{I\setminus J_{+}})\oplus\Hom(\fz_{I\cup J},E)\oplus P(\fl_{I}\cap\fh_{J_{+}},\fl_{I}\cap\fh),\]
we see that (\ref{equ: relative interval embedding Lie}) is an embedding with image $\wedge P(\fl_{I},\fl_{I}\cap\fh)$ if and only if the following map
\[H^{\bullet}(\fl_{I}\cap\fh_{J_{+}},\fl_{I}\cap \fh,1_{\fl_{I}\cap\fh_{J_{+}}})\rightarrow H^{\bullet}(\fl_{I}\cap\fh_{J_{+}},1_{\fl_{I}\cap\fh_{J_{+}}})=\wedge P(\fl_{I}\cap\fh_{J_{+}})\]
is an embedding with image $\wedge P(\fl_{I}\cap\fh_{J_{+}},\fl_{I}\cap\fh)$. Hence, upon replacing $\fl_{I}\subseteq\fg$ with $\fl_{I}\cap\fh_{J_{+}}\subseteq \fh_{J_{+}}$, it is harmless to assume further that $J_{+}=\Delta$ and thus $I\cup J=\Delta$ in the rest of the proof.
Recall that we have assumed that $J\not\subseteq I$ and thus $I\subsetneq \Delta$ with $r_{I}\geq 2$.
Now that $I\cup J=\Delta$ and $J$ is an interval, we see that $I^{d}\subseteq J$ for each $1<d<r_{I}$ and thus
\begin{equation}\label{equ: relative interval embedding Lie primitive}
P(\fl_{I},\fl_{I}\cap\fh)=P(\fh_{I^{-}},\fh_{I^{-}\cap J})\oplus P(\fh_{I^{+}},\fh_{I^{+}\cap J})
\end{equation}
with $I^{-}\defeq I^{1}$, $I^{+}\defeq I^{r_{I}}$ and $I^{\dagger}\defeq I^{-}\sqcup I^{+}$ for short.
Now that we have $\fh_{I^{\dagger}}=\fh_{I^{-}}\times\fh_{I^{+}}$ with $\fh\cap\fh_{I^{\dagger}}=\fh_{I^{\dagger}\cap J}=\fh_{I^{-}\cap J}\times\fh_{I^{+}\cap J}$, the inclusion $\fh_{I^{\dagger}}\subseteq\fl_{I}$ together with (\ref{equ: relative Lie Kunneth}) induces a map between graded $E$-algebras
\begin{multline}\label{equ: relative interval embedding Lie key}
H^{\bullet}(\fl_{I},\fl_{I}\cap\fh,1_{\fl_{I}})\rightarrow H^{\bullet}(\fh_{I^{\dagger}},\fh_{I^{\dagger}\cap J},1_{\fh_{I^{\dagger}}})\\
=\mathrm{Tot}(H^{\bullet}(\fh_{I^{-}},\fh_{I^{-}\cap J},1_{\fh_{I^{-}}})\otimes_EH^{\bullet}(\fh_{I^{+}},\fh_{I^{+}\cap J},1_{\fh_{I^{+}}}))
=\wedge(P(\fh_{I^{-}},\fh_{I^{-}\cap J})\oplus P(\fh_{I^{+}},\fh_{I^{+}\cap J})).
\end{multline}
We divide the rest of the proof into the following steps.

\textbf{Step $1$}: We prove that the map (\ref{equ: relative interval embedding Lie key}) is an isomorphism.\\
We consider the following natural map between (relative) Chevalley-Eilenberg complex (induced from the inclusion $\fh_{I^{\dagger}}\subseteq\fl_{I}$)
\begin{equation}\label{equ: relative interval embedding Lie CE}
\Hom_{U(\fl_{I}\cap\fh)}(\wedge^{\bullet}(\fl_{I}/\fl_{I}\cap\fh),E)\rightarrow \Hom_{U(\fh_{I^{\dagger}\cap J})}(\wedge^{\bullet}(\fh_{I^{\dagger}}/\fh_{I^{\dagger}\cap J}),E)
\end{equation}
whose associated map between cohomology recovers the map (\ref{equ: relative interval embedding Lie key}). It is easy to see that the inclusion $\fh_{I^{\dagger}}\subseteq\fl_{I}$ induces an isomorphism $\fh_{I^{\dagger}}/\fh_{I^{\dagger}\cap J}\buildrel\sim\over\longrightarrow \fl_{I}/\fl_{I}\cap\fh$ between $E$-vector spaces, and thus an isomorphism between $E$-vector spaces (for each $k\geq 0$)
\begin{equation}\label{equ: relative interval embedding Lie CE comparison}
\Hom_{E}(\wedge^{k}(\fl_{I}/\fl_{I}\cap\fh),E)\buildrel\sim\over\longrightarrow \Hom_{E}(\wedge^{k}(\fh_{I^{\dagger}}/\fh_{I^{\dagger}\cap J}),E),
\end{equation}
with an adjoint action of $\fl_{I}\cap\fh$ on LHS and an adjoint action of $\fh_{I^{\dagger}\cap J}$ on RHS that are compatible with each other. On one hand, we note that $\fh_{I^{d}}\subseteq \fh$ acts trivially on LHS of (\ref{equ: relative interval embedding Lie CE comparison}) for each $1<d<r_{I}$. On the other hand, we note that $\fz_{I}\subseteq\fl_{I}$ admits a natural adjoint action on LHS of (\ref{equ: relative interval embedding Lie CE comparison}) which is again trivial.
Now that there exists an $E$-Lie subalgebra $\ft'\subseteq \ft\cap \fh_{I^{\dagger}}$ such that
\[\ft'\times\fh_{I^{\dagger}\cap J}\times\fz_{I}\times\prod_{1<d<r_{I}}\fh_{I^{d}}=(\fl_{I}\cap\fh)\fz_{I}\subseteq \fl_{I},\]
we conclude that the LHS of (\ref{equ: relative interval embedding Lie CE}) is the subcomplex of the RHS of (\ref{equ: relative interval embedding Lie CE}) given by invariants under the adjoint $\ft'$-action.
Consequently, (\ref{equ: relative interval embedding Lie CE}) is a quasi-isomorphism and thus (\ref{equ: relative interval embedding Lie key}) is an isomorphism. (In terms of invariant differential forms on a homogenous space $H_1/H_2$ where $H_2\subseteq H_1$ is a pair of compact Lie groups with a maximal torus $S\subseteq H_1$ normalizing $H_2$, we consider the average map $\omega\mapsto\int_{S}\mathrm{Ad}(s)^{\ast}\omega$.)

Thanks to \textbf{Step $1$}, for each $\ast\in\{+,-\}$, we write $P_{\ast,\fh}=\bigoplus_{a=\#I^{\ast}\cap J+2}^{\#I^{\ast}+1}P_{\ast,\fh}^{2a-1}$ for the $E$-subspace of $H^{\bullet}(\fl_{I},\fl_{I}\cap\fh,1_{\fl_{I}})$ that corresponds to $P(\fh_{I^{\ast}},\fh_{I^{\ast}\cap J})$ under the isomorphism (\ref{equ: relative interval embedding Lie key}).

\textbf{Step $2$}: For each $\ast\in\{+,-\}$ and $\#I^{\ast}\cap J+2\leq a\leq \#I^{\ast}+1$, we prove that the map (\ref{equ: relative interval embedding Lie}) sends $P_{\ast,\fh}^{2a-1}$ to $P^{2a-1}(\fh_{I^{\ast}})\subseteq P(\fl_{I})\subseteq H^{\bullet}(\fl_{I},1_{\fl_{I}})$.\\
We assume from now $\ast=-$ with the $\ast=+$ case being similar. We write $\fh'\defeq \fh_{J'}$ with $J'\defeq [n_{1}-a+2,n-1]$ and note that $J'$ is the maximal subinterval of $\Delta$ such that $P_{-,\fh'}^{2a-1}$ is defined.
Now that the natural map
\[H^{\bullet}(\fl_{I},\fl_{I}\cap\fh',1_{\fl_{I}})\rightarrow H^{\bullet}(\fl_{I},\fl_{I}\cap\fh,1_{\fl_{I}})\]
restricts to an isomorphism $P_{-,\fh'}^{2a-1}\buildrel\sim\over\longrightarrow P_{-,\fh}^{2a-1}$ between $1$-dimensional $E$-vector spaces by \textbf{Step $1$} (with (\ref{equ: relative interval embedding Lie key}) being functorial with respect to the choice of $\fh$), it suffices to show that the natural map
\[H^{\bullet}(\fl_{I},\fl_{I}\cap\fh',1_{\fl_{I}})\rightarrow H^{\bullet}(\fl_{I},1_{\fl_{I}})\]
sends $P_{-,\fh'}^{2a-1}$ to $P^{2a-1}(\fh_{I^{-}})\subseteq P(\fl_{I})\subseteq H^{\bullet}(\fl_{I},1_{\fl_{I}})$. Upon replacing $\fh$ with $\fh'$, we assume in the rest of this step that $J=[n_{1}-a+2,n-1]$. We write $\fz'\defeq \fz_{\Delta\setminus\{n_1\}}$ and $\fl_{I}'\defeq \fh_{I^{-}}\times\fz'$ for short, and note that the isomorphism $\fl_{I}\cong \fl_{I}'\times(\fl_{I}\cap\fh_{[n_1+1,n-1]})$ restricts to an isomorphism $\fl_{I}\cap\fh\cong (\fl_{I}'\cap\fh)\times(\fl_{I}\cap\fh_{[n_1+1,n-1]})$. Hence, the surjection $\fl_{I}\twoheadrightarrow\fl_{I}'$ induces the following commutative diagram
\begin{equation}\label{equ: relative interval embedding diagram}
\xymatrix{
H^{\bullet}(\fl_{I}',\fl_{I}'\cap\fh,1_{\fl_{I}}) \ar[r] \ar^{\wr}[d] & H^{\bullet}(\fl_{I}',1_{\fl_{I}'}) \ar@{^{(}->}[d] \\
H^{\bullet}(\fl_{I},\fl_{I}\cap\fh,1_{\fl_{I}}) \ar[r] & H^{\bullet}(\fl_{I},1_{\fl_{I}})
}
\end{equation}
Note that the top horizontal map of (\ref{equ: relative interval embedding diagram}) factors through
\[H^{\bullet}(\fl_{I}',\fh_{I^{-}\cap J},1_{\fl_{I}'})\cong H^{\bullet}(\fh_{I^{-}},\fh_{I^{-}\cap J},1_{\fh_{I^{-}}})\otimes_EH^{\bullet}(\fz',1_{\fz'})=\wedge(P(\fh_{I^{-}},\fh_{I^{-}\cap J})\oplus\Hom(\fz',E))\]
with the first isomorphism from (\ref{equ: relative Lie Kunneth}). In particular, as $\fz'$ is $1$-dimensional and $2a-1\geq 3$, we have
\[H^{2a-1}(\fl_{I}',\fh_{I^{-}\cap J},1_{\fl_{I}'})=P^{2a-1}(\fh_{I^{-}},\fh_{I^{-}\cap J}),\]
which together with (\ref{equ: relative interval embedding diagram}) forces the image of $P_{-,\fh}^{2a-1}$ under (\ref{equ: relative interval embedding Lie}) to be contained in $P^{2a-1}(\fh_{I^{-}})\subseteq P(\fl_{I})\subseteq H^{\bullet}(\fl_{I},1_{\fl_{I}})$.

Finally, we conclude our desired results on (\ref{equ: relative interval embedding Lie}) from \textbf{Step $1$}, \textbf{Step $2$} as well as the following commutative diagram
\[
\xymatrix{
H^{\bullet}(\fl_{I},\fl_{I}\cap\fh,1_{\fl_{I}}) \ar^{\sim}[rr] \ar[d] & & H^{\bullet}(\fh_{I^{\dagger}},\fh_{I^{\dagger}\cap J},1_{\fh_{I^{\dagger}}}) \ar@{=}[r] \ar@{^{(}->}[d] & \wedge(P(\fh_{I^{-}},\fh_{I^{-}\cap J})\oplus P(\fh_{I^{+}},\fh_{I^{+}\cap J})) \ar@{^{(}->}[d]\\
H^{\bullet}(\fl_{I},1_{\fl_{I}}) \ar@{->>}[rr] & & H^{\bullet}(\fh_{I^{\dagger}},1_{\fh_{I^{\dagger}}}) \ar@{=}[r] & \wedge(P(\fh_{I^{-}})\oplus P(\fh_{I^{+}}))
}
\]
induced from the inclusion $\fh_{I^{-}}\times\fh_{I^{+}}=\fh_{I^{\dagger}}\subseteq\fl_{I}$.
\end{proof}

\begin{prop}\label{prop: relative interval embedding}
Let $I\subseteq\Delta$, $M\subseteq Z_{I}^{\dagger}$ be a subgroup and $\fh=\fh_{J}$ for some interval $J\subseteq\Delta$.
Then there exists a natural isomorphism between graded $E$-algebras
\begin{equation}\label{equ: relative interval embedding isom}
H^{\bullet}(L_{I}/M,\fl_{I}\cap\fh,1_{\fl_{I}})\buildrel\sim\over\longrightarrow \wedge (\Hom(Z_{I}^{\dagger}/M,E)\oplus P(\fl_{I},\fl_{I}\cap\fh))
\end{equation}
which is functorial with respect to the choice of $I$, $\fh$ and $M$, such that the natural map
\begin{equation}\label{equ: relative interval embedding}
H^{\bullet}(L_{I}/M,\fl_{I}\cap\fh,1_{L_{I}/M}) \rightarrow H^{\bullet}(L_{I},1_{\fl_{I}})
\end{equation}
is the embedding between graded $E$-algebras given by the wedge product of the inclusion
\[\Hom(Z_{I}^{\dagger}/M,E)\oplus P(\fl_{I},\fl_{I}\cap\fh)\subseteq \Hom(Z_{I}^{\dagger},E)\oplus P(\fl_{I})\]
under the identification (\ref{equ: Levi Kunneth}).
\end{prop}
\begin{proof}
This follows directly from Lemma~\ref{lem: relative interval embedding Lie} and \ref{it: Levi decomposition 2} of Proposition~\ref{prop: Levi decomposition} (by taking $M'=0$ and $\fh'=0$ in \emph{loc.cit.}).
\end{proof}

\begin{rem}\label{rem: h tuple}
When $M=0$, Proposition~\ref{prop: relative interval embedding} implies that the image of (\ref{equ: relative interval embedding}) admits a basis of the form $\{x_{\Theta}\}$ with $\Theta=(v,I,\un{k},\un{\Lambda})$ running through all tuples with $v\subseteq\mathbf{Log}_{I}$ satisfying $v\setminus\mathbf{Log}_{I}^{\infty}\subseteq\mathbf{Log}_{I\cup J}$ and 
\[\Lambda_{d}\subseteq \{2a-1\mid \#I^{d}\cap J+2\leq a\leq n_{d}\}\]
for each $1\leq d\leq r_{I}$. We use the term \emph{$\fh$-tuple} for tuples that satisfy the conditions above.
\end{rem}

\begin{cor}\label{cor: relative Levi restriction}
Let $I'\subseteq I\subseteq\Delta$, $M\subseteq Z_{I}^{\dagger}$ be a subgroup, and $\fh=\fh_{J}$ for some interval $J\subseteq \Delta$. Then the following restriction map 
\[\mathrm{Res}^{\bullet}_{I,I',\fh,M}: H^{\bullet}(L_{I}/M,\fl_{I}\cap\fh,1_{L_{I}/M})\rightarrow H^{\bullet}(L_{I'}/M,\fl_{I'}\cap\fh,1_{L_{I'}/M})\] 
corresponds under (\ref{equ: relative interval embedding isom}) and its variant with $I'$ replacing $I$ to the map
\begin{multline}\label{equ: relative Kunneth decomposition K restriction}
\mathrm{res}^{\bullet}_{I,I',\fh,M}: \mathrm{Tot}(\wedge^{\bullet}\Hom(Z_{I}^{\dagger}/M,E)\otimes_E H^{\bullet}(\fl_{I},\fl_{I}\cap\fh,1_{\fl_{I}})) \\
\rightarrow \mathrm{Tot}(\wedge^{\bullet}\Hom(Z_{I'}^{\dagger}/M,E)\otimes_E H^{\bullet}(\fl_{I'},\fl_{I'}\cap\fh,1_{\fl_{I'}}))
\end{multline}
induced from the map between pairs $(\fl_{I'},\fl_{I'}\cap\fh)\rightarrow(\fl_{I},\fl_{I}\cap\fh)$ and the map $Z_{I'}^{\dagger}/M\rightarrow Z_{I}^{\dagger}/M$.
\end{cor}
\begin{proof}
Thanks to Proposition~\ref{prop: Levi restriction}, Lemma~\ref{lem: relative interval embedding Lie} and Proposition~\ref{prop: relative interval embedding}, we conclude
\[\mathrm{Res}^{\bullet}_{I,I',\fh,M}=\mathrm{res}^{\bullet}_{I,I',\fh,M}\]
as the restriction of
\[\mathrm{Res}^{\bullet}_{I,I'}=\mathrm{res}^{\bullet}_{I,I'}\]
to
\[H^{\bullet}(L_{I}/M,\fl_{I}\cap\fh,1_{L_{I}/M})\cong \wedge^{\bullet}\Hom(Z_{I}^{\dagger}/M,E)\otimes_E H^{\bullet}(\fl_{I},\fl_{I}\cap\fh,1_{\fl_{I}})\]
via (\ref{equ: relative interval embedding}) and (\ref{equ: relative interval embedding Lie}).
\end{proof}
\subsection{Atomic subcomplex}\label{subsec: comb subcomplex}
We index the standard basis of $E_{1,I_0,I_1}^{\bullet,\bullet}$ by certain tuples $\Theta=(v,I,\un{k},\un{\Lambda})$ (see (\ref{equ: E1 std basis})). We define an equivalence relation on the set of all tuples (see Definition~\ref{def: partial order}) and isolate out the so-called $(I_0,I_1)$-atomic equivalence classes (see Definition~\ref{def: atom element}). We associate an element $x_{\Omega}\in E_{1,I_0,I_1}^{-\ell,k}$ with each $(I_0,I_1)$-atomic equivalence classes $\Omega$ with bidegree $(-\ell,k)$, and then define $E_{1,I_0,I_1,\diamond}^{-\ell,k}\subseteq E_{1,I_0,I_1}^{-\ell,k}$ as the $E$-subspace spanned by $x_{\Omega}$ for all $(I_0,I_1)$-atomic $\Omega$ with bidegree $(-\ell,k)$. We prove that $E_{1,I_0,I_1,\diamond}^{\bullet,k}\subseteq E_{1,I_0,I_1}^{\bullet,k}$ is an embedding between complex (see Lemma~\ref{lem: atomic subcomplex}) and moreover a quasi-isomorphism (with cohomologies $E_{2,I_0,I_1}^{\bullet,k}$, see Proposition~\ref{prop: quasi isom}) for each $k\geq 0$.
We will frequently use notation introduced around (\ref{equ: explicit complement}).

Let $I_0\subseteq I\subseteq I_1$ be a subset satisfying $\#I=\ell$ (with $r_{I}=n-\ell$) and $\un{k}=\{k_d\}_{1\leq d\leq r_I}$ be a tuple satisfying $|\un{k}|=k$. Let $v\subseteq \mathbf{Log}_{I}$ with $\#v=k_0$ and $\un{\Lambda}\in\Sigma_{I,\un{k}}$ (see (\ref{equ: Lambda index})), we use in the following the shortened notation $\Theta=(v,I,\un{k},\un{\Lambda})$ and write $x_\Theta\defeq v\otimes v_{I,\un{k}}^{\un{\Lambda}}$. Note that we have $I_0\subseteq I\subseteq I_v\cap I_1$ for each $v\subseteq \mathbf{Log}_{I}$.
It follows from (\ref{equ: std basis}) and (\ref{equ: first page as Levi coh K}) that the space $E_{1,I_0,I_1}^{-\ell,k}$ admits a basis of the form
\begin{equation}\label{equ: E1 std basis}
\{x_\Theta\mid \Theta=(v,I,\un{k},\un{\Lambda}),~v\subseteq \mathbf{Log}_{\emptyset},~\un{\Lambda}\in\Sigma_{I,\un{k}},~I_0\subseteq I\subseteq I_v\cap I_1,~|\un{k}|=k,~\#I=\ell,~|v|=k_0\}.
\end{equation}
Hence, we have a decomposition (with $E_{1,I_0,I_1,v}^{-\ell,k}$ being the $v$-isotypic direct summand spanned by $x_{\Theta}$ for a fixed $v$)
\begin{equation}\label{equ: E1 v decomposition}
E_{1,I_0,I_1}^{-\ell,k}=\bigoplus_{v\subseteq \mathbf{Log}_{\emptyset}} E_{1,I_0,I_1,v}^{-\ell,k}
\end{equation}
which is compatible with the differential $d_{1,I_0,I_1}^{-\ell,k}$.
For each $x\in E_{1,I_0,I_1,v}^{-\ell,k}$ and each $\Theta=(v,I,\un{k},\un{\Lambda})$, we write $c_{\Theta}(x)$ for the coefficient of $x$ attached to $x_\Theta$.

Note that the sets $I\subseteq I_v\cap I_1$ determines a sequence
\begin{equation}\label{equ: v tuple integers}
0=r_{v,I_1,I}^0<r_{v,I_1,I}^1<\cdots<r_{v,I_1,I}^{r_{I_v\cap I_1}}=r_I
\end{equation}
characterized by
$\Delta\setminus (I_v\cap I_1)=\{\sum_{d=1}^{r_{v,I_1,I}^s}n_d\mid 1\leq s\leq r_{I_v\cap I_1}-1\}$.
Given $1\leq d\leq r_I$, note that $\sum_{d'=1}^{d-1}n_{d'}\in \{0\}\sqcup\Delta\setminus (I_v\cap I_1)$ if and only if there exists $1\leq s\leq r_{I_v\cap I_1}$ such that $d=r_{v,I_1,I}^{s-1}+1$, and similarly $\sum_{d'=1}^{d}n_{d'}\in \{n\}\sqcup\Delta\setminus (I_v\cap I_1)$ if and only if there exists $1\leq s\leq r_{I_v\cap I_1}$ such that $d=r_{v,I_1,I}^{s}$.

For each tuple $\Theta=(v,I,\un{k},\un{\Lambda})$ and $1\leq d\leq r_I$, we write $i_{\Theta,d}\defeq \sum_{d'=1}^dn_{d'}\in\Delta\setminus I$.
\begin{cons}\label{cons: tuple operation}
We have the following two fundamental constructions which (if defined) send a given tuple $\Theta=(v,I,\un{k},\un{\Lambda})$ to a new tuple.
\begin{enumerate}[label=(\roman*)]
\item \label{it: operation 1} For each $i\in\Delta\setminus I$, there exists a unique $1\leq d\leq r_I$ such that $i=i_{\Theta,d}$, and we define $p_i^+(\Theta)\defeq (v,I\sqcup\{i\},\un{k}',\un{\Lambda}')$ by the condition that
 $\Lambda'_{d'}=\Lambda_{d'}$ for each $1\leq d'\leq d-1$, $\Lambda'_{d'}=\Lambda_{d'+1}$ for each $d+1\leq d'\leq r_I-1$ and $\Lambda'_{d}=\Lambda_d\sqcup \Lambda_{d+1}$. So $p_i^+(\Theta)$ is well defined if and only if $i\in\Delta\setminus I$ and $\Lambda_d\cap \Lambda_{d+1}=\emptyset$.
\item \label{it: operation 2} For each $i\in I$, there exists a unique $1\leq d\leq r_I$ such that $i\in I^d$, and we define $p_i^-(\Theta)\defeq (v,I\setminus\{i\},\un{k}',\un{\Lambda}')$ by the condition that
 $\Lambda'_{d'}=\Lambda_{d'}$ for each $1\leq d'\leq d-1$, $\Lambda'_{d'}=\Lambda_{d'-1}$ for each $d+2\leq d'\leq r_I+1$, $\Lambda'_{d}=\Lambda_d$ and $\Lambda'_{d+1}=\emptyset$. So $p_i^-(\Theta)$ is well defined if and only if $i\in I$ and $\max\{m\mid m\in \Lambda_d\}\leq 2(i-i_{\Theta,d-1})-1$.
\end{enumerate}
\end{cons}

We have the following equivalence relation on the set of all tuples.
\begin{defn}\label{def: partial order}
Let $\Theta=(v,I,\un{k},\un{\Lambda})$, $\Theta'=(v',I',\un{k}',\un{\Lambda}')$ be two tuples.
\begin{enumerate}[label=(\roman*)]
\item \label{it: tuple partial order 1} We say that \emph{$\Theta'$ is an improvement of $\Theta$ with level $d$} for some $1\leq d\leq r_I$, if $\Theta'=p_{i_{\Theta,d}}^+(p_{i_{\Theta',d}}^-(\Theta))$ (with $i_{\Theta,d}\in I'\subseteq I_v\cap I_1$). Note that this forces $v=v'$, $\un{k}=\un{k}'$, $\un{\Lambda}=\un{\Lambda}'$ and $r_{v,I_1,I}^s=r_{v,I_1,I'}^s$ for each $0\leq s\leq r_{I_v\cap I_1}$.
\item \label{it: tuple partial order 2} We say that \emph{$\Theta$ is smaller than $\Theta'$ with level $\geq d$} for some $1\leq d\leq r_I$, written $\Theta<_{d}\Theta'$, if there exists a sequence of tuples $\Theta=\Theta^0,\dots,\Theta^m=\Theta'$ such that $\Theta^{m'}$ is an improvement of $\Theta^{m'-1}$ with level $\geq d$ for each $1\leq m'\leq m$. We say that $\Theta$ is \emph{maximal} if it is not smaller than any other tuple $\Theta'$ with any level.
\item \label{it: tuple partial order 3} We say that $\Theta$ and $\Theta'$ are \emph{equivalent} if there exists $\Theta''$ such that $\Theta<_{1}\Theta''$ and $\Theta'<_{1}\Theta''$.
\end{enumerate}
\end{defn}

Let $\Omega$ be an equivalence class of tuples and let $\Theta=(v,I,\un{k},\un{\Lambda})\in\Omega$.
The \emph{bidegree} of $\Omega$ is defined to be $(-\ell_\Omega,k_\Omega)\defeq (-\#I,|\un{k}|)$ which does not depend on the choice of $\Theta\in\Omega$.
Similarly, the integer $r_{v,I_1,I}^s$ does not depend on the choice of $\Theta\in\Omega$, so we write
\begin{equation}\label{equ: v class integer}
r_{\Omega}^s\defeq r_{v,I_1,I}^s
\end{equation}
for each $0\leq s\leq r_{I_v\cap I_1}$.
Note that $\sum_{d=1}^{r_{\Omega}^s}n_d\in\Delta\setminus I_v$ depends only on the choice of $v$ and $s$ but not on the choice of $I$.
We also note that $v$ and $\un{\Lambda}=\{\Lambda_d\}_{1\leq d\leq n-\ell_{\Omega}}$ does not depend on the choice of $\Theta\in\Omega$.
For each $\Theta\in\Omega$ define its \emph{sign} as
\begin{equation}\label{equ: sign of tuple}
\varepsilon(\Theta)\defeq (-1)^{\sum_{i\in I}i}.
\end{equation}
We consider the following element
\begin{equation}\label{equ: equivalence class element}
x_\Omega\defeq \sum_{\Theta=(v,I,\un{k},\un{\Lambda})\in\Omega}\varepsilon(\Theta)x_\Theta \in E_{1,I_0,I_1}^{-\ell,k}.
\end{equation}

\begin{lem}\label{lem: unique maximal tuple}
Let $\Omega$ be an equivalence class of tuples and $\Theta=(v,I,\un{k},\un{\Lambda})\in\Omega$. We have the following results.
\begin{enumerate}[label=(\roman*)]
\item \label{it: unique maximal tuple 1} The equivalence class $\Omega$ contains a unique maximal element.
\item \label{it: unique maximal tuple 2} The equivalence class $\Omega$ depends only on the choice of $v$, $\ell_{\Omega}$, $\un{\Lambda}$ and $r_{\Omega}^s$ for each $0\leq s\leq r_{I_v\cap I_1}$.
\end{enumerate}
\end{lem}
\begin{proof}
Note that $i_{\Theta,r_{\Omega}^s}\in\Delta\setminus I_v$ does not depend on the choice of $\Theta\in\Omega$ for each $0\leq s\leq r_{I_v\cap I_1}$.
We construct an explicit tuple $\Theta^{\max}=(v,I^{\max},\un{k},\un{\Lambda})\in\Omega$ (or equivalently an explicit $I_0\subseteq I^{\max}\subseteq I_v\cap I_1$ satisfying $r_{v,I_1,I^{\max}}^s=r_{\Omega}^s$ for each $0\leq s\leq r_{I_v\cap I_1}$). We write $\{n^{\max}_d\}_{1\leq d\leq r_{I^{\max}}}$ for the tuple of integers corresponding to $I^{\max}$, and note that the choice of $I^{\max}$ is equivalent to the choice of $i_{\Theta^{\max},d}=sum_{d'=1}^{d}n^{\max}_{d'}$ for each $1\leq d\leq r_{I^{\max}}=n-\ell_{\Omega}$.
For each $1\leq s\leq r_{I_v\cap I_1}$, we necessarily have $i_{\Theta^{\max},r_{\Omega}^{s-1}}=i_{\Theta,r_{\Omega}^{s-1}}\in \{0\}\sqcup\Delta\setminus I_v$, and we inductively define $i_{\Theta^{\max},d}$ as the smallest integer in $\Delta\setminus I_0$ such that $i_{\Theta^{\max},d}-i_{\Theta^{\max},d-1}=n^{\max}_d\geq \frac{\max\{m\mid m\in\Lambda_d\}+1}{2}$, for increasing $r_{\Omega}^{s-1}+1\leq d\leq r_{\Omega}^{s}-1$. Thanks to the given tuple of integers $\{i_{\Theta,d}\}_{1\leq d\leq n-\ell_{\Omega}}$, we see that $i_{\Theta^{\max},d}$ necessarily exists and satisfies $i_{\Theta^{\max},d}\leq i_{\Theta,d}$ for each $1\leq s\leq r_{I_v\cap I_1}$ and inductively for increasing $r_{\Omega}^{s-1}+1\leq d\leq r_{\Omega}^{s}-1$.
Moreover, we may pass from $\Theta$ to $\Theta^{\max}$ by successively applying $p_{i_{\Theta,d}}^+(p_{i_{\Theta^{\max},d}}^-(\cdot))$ for each increasing $1\leq d\leq r_I$ that satisfies $i_{\Theta^{\max},d}<i_{\Theta,d}$, and in particular we have $\Theta<_{1}\Theta^{\max}$.
To summarize, we have constructed an explicit tuple $\Theta^{\max}\in\Omega$ which must be the unique maximal tuple in $\Omega$. Since the construction of $\Theta^{\max}$ depends only on the choice of $v$, $\ell_{\Omega}$, $\un{\Lambda}$ and $r_{\Omega}^s$ for each $0\leq s\leq r_{I_v\cap I_1}$, so does $\Omega$. This finishes the proof of both \ref{it: unique maximal tuple 1} and \ref{it: unique maximal tuple 2}.
\end{proof}

\begin{lem}\label{lem: equivalence glue}
Let $\Omega$ be an equivalence class of tuples and $\Theta,\Theta'$ be two elements inside. Then for each $1\leq s\leq r_{I_v\cap I_1}$ and $r_{\Omega}^{s-1}+1\leq d\leq r_{\Omega}^s-1$, the tuples $p_{i_{\Theta,d}}^+(\Theta)$ and $p_{i_{\Theta',d}}^+(\Theta')$ are equivalent.
\end{lem}
\begin{proof}
According to Definition~\ref{def: partial order}, it is harmless to assume that $\Theta'$ is an improvement of $\Theta$ with level $d_1$, namely $\Theta'=p_{i_{\Theta,d_1}}^+(p_{i_{\Theta',d_1}}^-(\Theta))$, for some $1\leq s\leq r_{I_v\cap I_1}$ and $r_{\Omega}^{s-1}+1\leq d_1\leq r_{\Omega}^s$. If $d_1=d$, then we clearly have $p_{i_{\Theta,d}}^+(\Theta)=p_{i_{\Theta',d}}^+(\Theta')$. If $d_1\neq d$, then we have $i\defeq i_{\Theta',d}=i_{\Theta,d}$ and $p_i^+(\Theta')=p_{i_{\Theta,d_1}}^+(p_{i_{\Theta',d_1}}^-(p_i^+(\Theta)))$, namely $p_i^+(\Theta')$ is an improvement of $p_i^+(\Theta)$ with level $d_1$.
The proof is thus finished.
\end{proof}

\begin{lem}\label{lem: tuple sum of two}
Let $\Theta=(v,I,\un{k},\un{\Lambda})$ be a tuple with bidegree $(-\ell,k)$, and $1\leq d\leq r_I$ be an integer that satisfies $\Lambda_d=\emptyset$ and $i_{\Theta,d-1},i_{\Theta,d}\in (I_v\cap I_1)\setminus I$. Then we have
\begin{equation}\label{equ: tuple sum of two}
\varepsilon(p_{i_{\Theta,d-1}}^+(\Theta))c_{\Theta}(d_{1,I_0,I_1}^{-\ell,k}(x_{p_{i_{\Theta,d-1}}^+(\Theta)}))+\varepsilon(p_{i_{\Theta,d}}^+(\Theta))c_{\Theta}(d_{1,I_0,I_1}^{-\ell,k}(x_{p_{i_{\Theta,d}}^+(\Theta)}))=0.
\end{equation}
\end{lem}
\begin{proof}
It follows from (\ref{equ: sign of tuple}) that $\varepsilon(p_{i_{\Theta,d-1}}^+(\Theta))=(-1)^{i_{\Theta,d-1}}\varepsilon(\Theta)$ and $\varepsilon(p_{i_{\Theta,d}}^+(\Theta))=(-1)^{i_{\Theta,d}}\varepsilon(\Theta)$, and thus
\begin{equation}\label{equ: tuple sum of two sign 1}
\varepsilon(p_{i_{\Theta,d}}^+(\Theta))=(-1)^{i_{\Theta,d}-i_{\Theta,d-1}}\varepsilon(p_{i_{\Theta,d-1}}^+(\Theta))=(-1)^{n_d}\varepsilon(p_{i_{\Theta,d-1}}^+(\Theta)).
\end{equation}
Recall from the definition of the differential map $d_{1,I_0,I_1}^{-\ell,k}$ (see the discussion around (\ref{equ: abstract E1 sum}) and (\ref{equ: differential sign})) that \[c_{\Theta}(d_{1,I_0,I_1}^{-\ell,k}(x_{p_{i_{\Theta,d-1}}^+(\Theta)}))=(-1)^{m(I\sqcup\{i_{\Theta,d-1}\},i_{\Theta,d-1})}\] 
and 
\[c_{\Theta}(d_{1,I_0,I_1}^{-\ell,k}(x_{p_{i_{\Theta,d}}^+(\Theta)}))=(-1)^{m(I\sqcup\{i_{\Theta,d}\},i_{\Theta,d})},\] 
and thus
\begin{equation}\label{equ: tuple sum of two sign 2}
c_{\Theta}(d_{1,I_0,I_1}^{-\ell,k}(x_{p_{i_{\Theta,d}}^+(\Theta)}))=(-1)^{\#I^d}c_{\Theta}(d_{1,I_0,I_1}^{-\ell,k}(x_{p_{i_{\Theta,d}}^+(\Theta)}))=(-1)^{n_d-1}c_{\Theta}(d_{1,I_0,I_1}^{-\ell,k}(x_{p_{i_{\Theta,d}}^+(\Theta)})).
\end{equation}
The equalities (\ref{equ: tuple sum of two sign 1}) and (\ref{equ: tuple sum of two sign 2}) together give (\ref{equ: tuple sum of two}).
\end{proof}

\begin{defn}\label{def: std element}
For each tuple $\Theta=(v,I,\un{k},\un{\Lambda})$, we define $d_{I_0,I_1,\Theta}$ as the maximal integer, if exists, such that $\Lambda^{d_{I_0,I_1,\Theta}}=\emptyset$ and $d_{I_0,I_1,\Theta}\neq r_{v,I_1,I}^{s-1}+1$ for each $1\leq s\leq r_{I_v\cap I_1}$.
We say that $\Theta$ is \emph{$(I_0,I_1)$-standard} if there does not exist $i\in I$ such that $i>i_{\Theta,d_{I_0,I_1,\Theta}-1}$ and $p_{i}^{-}(\Theta)$ is defined (see Construction~\ref{cons: tuple operation}).
\end{defn}

\begin{defn}\label{def: atom element}
Let $\Theta=(v,I,\un{k},\un{\Lambda})$ be a tuple.
We say that $\Theta$ is \emph{maximally $(I_0,I_1)$-atomic} if
\begin{enumerate}[label=(\roman*)]
\item \label{it: atom element 1} $\Theta$ is maximal (see \ref{it: tuple partial order 2} of Definition~\ref{def: partial order});
\item \label{it: atom element 2} if $\Lambda_d=\emptyset$ for some $1\leq d\leq r_I$, then there exists $1\leq s\leq r_{I_v\cap I_1}$ such that $d=r_{v,I_1,I}^{s-1}+1$; and
\item \label{it: atom element 3} for each $1\leq s\leq r_{I_v\cap I_1}$, there does not exist $i\in I^{r_{v,I_1,I}^{s}}$ such that $p_i^-(\Theta)$ is defined.
\end{enumerate}
We say that an equivalence class $\Omega$ is \emph{$(I_0,I_1)$-atomic} if it contains a maximally $(I_0,I_1)$-atomic element.
We say that a tuple $\Theta$ is \emph{$(I_0,I_1)$-atomic} if it lies in an $(I_0,I_1)$-atomic equivalence class.
\end{defn}

\begin{rem}\label{rem: atom element bis}
Let $\Theta=(v,I,\un{k},\un{\Lambda})$ be a tuple. We have the following observations.
\begin{enumerate}[label=(\roman*)]
\item \label{it: atom element bis 1} The tuple $\Theta$ satisfies \ref{it: atom element 1} and \ref{it: atom element 3} of Definition~\ref{def: atom element} if and only if $p_{i}^{-}(\Theta)$ is not defined for any $i\in I$.
\item \label{it: atom element bis 2} If $\Theta$ satisfies \ref{it: atom element 1} and \ref{it: atom element 3} of Definition~\ref{def: atom element}, then we have $I^{d}\subseteq I_0$ for each $1\leq d\leq r_{I}$ that satisfies $\Lambda_{d}=\emptyset$.
\end{enumerate}
\end{rem}

By Lemma~\ref{lem: unique maximal tuple} we know that each $(I_0,I_1)$-atomic equivalence class contains a unique maximally $(I_0,I_1)$-atomic element which is its unique maximal element.

\begin{lem}\label{lem: atom property}
Let $\Theta=(v,I,\un{k},\un{\Lambda})$ be an $(I_0,I_1)$-atomic tuple. Then $\Theta$ satisfies both \ref{it: atom element 2} and \ref{it: atom element 3} of Definition~\ref{def: atom element}.
\end{lem}
\begin{proof}
Let $\Omega$ be the equivalence class that contains $\Theta$, and $\Theta'=(v,I',\un{k},\un{\Lambda})$ be the unique maximally $(I_0,I_1)$-atomic tuple in $\Omega$ (knowing that $\un{\Lambda}$ and $r_{v,I_1,I}^{s}$ do not depend on the choice of $\Theta\in\Omega$). Hence, $\Theta$ satisfies \ref{it: atom element 2} of Definition~\ref{def: atom element} since $\Theta'$ does.
Recall from the proof of Lemma~\ref{lem: unique maximal tuple} we can explicitly construct the maximal element $\Theta'$ in $\Omega$ which satisfies $i_{\Theta',d}\leq i_{\Theta,d}$ for each $1\leq d\leq r_I=n-\ell_{\Omega}$ and $i_{\Theta,r_{\Omega}^s}=i_{\Theta',r_{\Omega}^s}\notin I_v\cap I_1$ for each $0\leq s\leq r_{I_v\cap I_1}$. In particular, for each $1\leq s\leq r_{I_v\cap I_1}$, we have $I^{r_{\Omega}^s}\subseteq (I')^{r_{\Omega}^s}$. Assume on the contrary that $p_i^-(\Theta)$ is defined for some $i\in I^{r_{v,I_1,I}^{s}}$, then $p_i^-(\Theta')$ is necessarily defined as $I^{r_{\Omega}^s}\subseteq (I')^{r_{\Omega}^s}$, a contradiction to the fact that $\Theta'$ satisfies \ref{it: atom element 3} of Definition~\ref{def: atom element}. The proof is thus finished.
\end{proof}

\begin{lem}\label{lem: atom glue}
Let $\Theta=(v,I,\un{k},\un{\Lambda})$ be a tuple, $1\leq s_0\leq r_{I_v\cap I_1}$ and $r_{v,I_1,I}^{s_0-1}+1\leq d_0\leq r_{v,I_1,I}^{s_0}-1$ (with $i_{\Theta,d_0}\in I_v\cap I_1$). Assume that $\Lambda_{d_0+1}\neq \emptyset$ if $d_0=r_{v,I_1,I}^{s_0-1}+1$, and that $\Lambda_{d_0}\neq \emptyset\neq \Lambda_{d_0+1}$ if $d_0>r_{v,I_1,I}^{s_0-1}+1$. If $p_{i_{\Theta,d_0}}^+(\Theta)$ is $(I_0,I_1)$-atomic, then so is $\Theta$.
\end{lem}
\begin{proof}
Let $\Theta'$ be the unique maximal tuple in the equivalence class of $\Theta$ (see Lemma~\ref{lem: unique maximal tuple} and Definition~\ref{def: atom element}). By Lemma~\ref{lem: equivalence glue} we know that $p_{i_{\Theta,d_0}}^+(\Theta)$ and $p_{i_{\Theta',d_0}}^+(\Theta')$ are equivalent. As $p_{i_{\Theta,d_0}}^+(\Theta)$ is $(I_0,I_1)$-atomic by our assumption, we deduce that $p_{i_{\Theta',d_0}}^+(\Theta')$ is $(I_0,I_1)$-atomic.
We prove that $\Theta'$ is maximally $(I_0,I_1)$-atomic by checking the conditions Definition~\ref{def: atom element}.
\ref{it: atom element 1} of Definition~\ref{def: atom element} holds by our choice of $\Theta'$.
We use the shortened notation $\Theta''\defeq p_{i_{\Theta',d_0}}^+(\Theta')=(v,I'',\un{k}'',\un{\Lambda}'')$ with $I''\defeq I'\sqcup\{i_{\Theta',d_0}\}$.

We check \ref{it: atom element 2} of Definition~\ref{def: atom element} for $\Theta'$.\\
Assume on the contrary that there exists $1\leq d_1\leq r_{I'}$ such that $\Lambda_{d_1}=\emptyset$ and $d_1\neq r_{v,I_1,I'}^{s-1}+1$ for each $1\leq s\leq r_{I_v\cap I_1}$.
(Note that $d_1\neq r_{v,I_1,I'}^{s-1}+1$ for each $1\leq s\leq r_{I_v\cap I_1}$ if and only if $i_{\Theta',d_1}\in I_v\cap I_1$.)
Recall from our assumption that $\Lambda_{d_0+1}\neq \emptyset$ if $d_0=r_{v,I_1,I}^{s_0-1}+1$, and that $\Lambda_{d_0}\neq \emptyset\neq \Lambda_{d_0+1}$ if $d_0>r_{v,I_1,I}^{s_0-1}+1$. Since $\Lambda_{d_1}=\emptyset$ and $d_1\neq r_{v,I_1,I'}^{s-1}+1=r_{v,I_1,I}^{s-1}+1$ for each $1\leq s\leq r_{I_v\cap I_1}$, we deduce that $d_1\notin \{d_0,d_0+1\}$.
We set $d_2\defeq d_1$ if $d_1<d_0$, and $d_2\defeq d_1-1$ if $d_1>d_0+1$. Then by the definition of $\Theta''=p_{i_{\Theta',d_0}}^+(\Theta')$ in Construction~\ref{cons: tuple operation} we have $\Lambda_{d_2}''=\Lambda_{d_1}=\emptyset$ with $i_{\Theta'',d_2}=i_{\Theta',d_1}\in I_v\cap I_1$ (namely $d_2\neq r_{v,I_1,I''}^{s-1}+1$ for each $1\leq s\leq r_{I_v\cap I_1}$), which contradicts the fact that $\Theta''$ is $(I_0,I_1)$-atomic (see Lemma~\ref{lem: atom property}).

We check \ref{it: atom element 3} of Definition~\ref{def: atom element} for $\Theta'$.\\
Assume on the contrary that there exists $1\leq s\leq r_{I_v\cap I_1}$ and $i\in I^{r_{v,I_1,I'}^{s}}$ such that $p_i^-(\Theta')$ is defined. Then we notice that $p_i^-(\Theta'')=p_{i_{\Theta',d_0}}^+(p_i^-(\Theta'))$ is defined, with $i\in I^{r_{v,I_1,I'}^{s}}\subseteq I^{r_{v,I_1,I''}^{s}}$, contradicting the fact that $\Theta''$ is $(I_0,I_1)$-atomic (see again Lemma~\ref{lem: atom property}).
\end{proof}

Let $\#I_0\leq \ell\leq \#I_1$ and $k\geq 0$. For each $v\subseteq \mathbf{Log}_{\emptyset}$, we define $E_{1,I_0,I_1,\diamond,v}^{-\ell,k}\subseteq E_{1,I_0,I_1,v}^{-\ell,k}$ as the subspace spanned by $x_{\Omega}$ with $\Omega$ running through equivalence classes of tuples with the fixed $v$ and bidegree $(-\ell_{\Omega},k_{\Omega})=(-\ell,k)$. We also set
\[
E_{1,I_0,I_1,\diamond}^{-\ell,k}\defeq \bigoplus_{v\subseteq \mathbf{Log}_{\emptyset}}E_{1,I_0,I_1,\diamond,v}^{-\ell,k}.
\]
As different equivalence classes of tuples do not intersect, it is clear that
$\{x_{\Omega}\}$ forms a basis of $E_{1,I_0,I_1,\diamond}^{-\ell,k}$ where $\Omega$ runs through all $(I_0,I_1)$-atomic equivalence classes with bidegree $(-\ell,k)$.
\begin{lem}\label{lem: atomic subcomplex}
Let $k\geq 0$ and $v\subseteq \mathbf{Log}_{\emptyset}$.
Then we have $d_{1,I_0,I_1}^{-\ell,k}(E_{1,I_0,I_1,\diamond,v}^{-\ell,k})\subseteq E_{1,I_0,I_1,\diamond,v}^{-\ell+1,k}$ for each $\ell$. Consequently, $E_{1,I_0,I_1,\diamond,v}^{\bullet,k}$ is a subcomplex of $E_{1,I_0,I_1,v}^{\bullet,k}$, and $E_{1,I_0,I_1,\diamond}^{\bullet,k}$ is a subcomplex of $E_{1,I_0,I_1}^{\bullet,k}$.
\end{lem}
\begin{proof}
It suffices to show that
\begin{equation}\label{equ: atom differential}
d_{1,I_0,I_1}^{-\ell,k}(x_\Omega)\in E_{1,I_0,I_1,\diamond,v}^{-\ell+1,k}
\end{equation}
for each equivalence class $\Omega$ of $(I_0,I_1)$-atomic tuples with fixed $v$ and bidegree $(-\ell_\Omega,k_\Omega)=(-\ell,k)$.
We consider a tuple $\Theta=(v,I,\un{k},\un{\Lambda})\in\Omega$, $i\in I$ and another tuple $\Theta'=(v,I',\un{k}',\un{\Lambda}')$ such that
\begin{equation}\label{equ: atom differential tuple}
c_{\Theta'}(\mathrm{Res}_{I,I\setminus\{i\}}^{k}(x_\Theta))\neq 0.
\end{equation}
Note that (\ref{equ: atom differential tuple}) holds if and only if $\Theta=p_i^+(\Theta')$.
Note also that there exists a unique $1\leq s\leq r_{I_v\cap I_1}$ and $r_{v,I_1,I}^{s-1}+1\leq d_0\leq r_{v,I_1,I}^s$ such that $i=i_{\Theta',d_0}\in I^{d_0}$, and $\Theta=p_i^+(\Theta')$ ensures that $\Lambda'_d=\Lambda_d$ if $d<d_0$, $\Lambda'_{d+1}=\Lambda_d$ if $d>d_0$, and that $\Lambda'_{d_0}\sqcup\Lambda'_{d_0+1}=\Lambda_{d_0}$ (with $\Lambda'_{d_0}\cap\Lambda'_{d_0+1}=\emptyset$).
We write $\Omega'$ for the equivalence class containing $\Theta'$.
We have the following cases.

\textbf{Case $1$}: Assume that $\Lambda'_{d_0+1}=\emptyset$. We have $\Theta'=p_i^-(\Theta)$ in this case (see Construction~\ref{cons: tuple operation}). Since $\Theta$ is $(I_0,I_1)$-atomic and $p_i^-(\Theta)$ is defined, we deduce from Lemma~\ref{lem: atom property} that $d_0<r_{v,I_1,I}^s$ and thus $i_{\Theta,d_0}=i_{\Theta',d_0+1}\in I_v\cap I_1$. Hence, $\Theta''\defeq p_{i_{\Theta',d_0+1}}^+(\Theta')=p_{i_{\Theta',d_0+1}}^+(p_i^-(\Theta))$ is defined and we have $\Theta<_{1}\Theta''$ (so $\Theta''\in\Omega$) by Definition~\ref{def: partial order}.
Now we consider an arbitrary $1\leq d\leq r_{I'}$ that satisfies $i_{\Theta',d}\in I_v\cap I_1$ so that $\Theta'''=(v,I''',\un{k}''',\un{\Lambda}''')\defeq p_{i_{\Theta',d}}^+(\Theta')$ is defined. 
We set $d_1\defeq d_0$ if $d<d_0$ and $d_1\defeq d_0+1$ if $d>d_0+1$. Then we have $\Lambda'''_{d_1}=\Lambda'_{d_0+1}=\emptyset$ with $r_{v,I_1,I'''}^{s-1}+1<d_1\leq r_{v,I_1,I'''}^s$ and thus $\Theta'''$ cannot be $(I_0,I_1)$-atomic (see Lemma~\ref{lem: atom property}) and thus cannot be an element of $\Omega$.
In other words, among $p_{i'}^+(\Theta')$ for all possible choices of $i'\in (I_v\cap I_1)\setminus I'$, $\Theta=p_i^+(\Theta')$ and $\Theta''=p_{i_{\Theta',d_0+1}}^+(\Theta')$ are the only tuples that belong to $\Omega$ (and share the common $\un{\Lambda}$).
We thus deduce from Lemma~\ref{lem: tuple sum of two} (with $\Theta'$ and $d_0+1$ here replacing $\Theta$ and $d$ in \emph{loc.cit.}) that $c_{\Theta'}(d_{1,I_0,I_1}^{-\ell,k}(x_\Omega))=0$.
Since $\Lambda'_{d_0+1}=\emptyset$ with $r_{v,I_1,I}^{s-1}+1<d_0+1\leq r_{v,I_1,I}^s$, it follows from Lemma~\ref{lem: atom property} that $\Theta'$ is not $(I_0,I_1)$-atomic.

\textbf{Case $2$}: Assume that $\Lambda'_{d_0}=\emptyset\neq \Lambda'_{d_0+1}$ and $d_0>r_{v,I_1,I'}^{s-1}+1$.
In this case, we know that $\Theta'$ is not $(I_0,I_1)$-atomic (see Lemma~\ref{lem: atom property}) and that $i_{\Theta',d_0-1}\in I_v\cap I_1$. Hence, $\Theta''\defeq p_{i_{\Theta',d_0-1}}^+(\Theta')$ is defined with $\Theta'=p_{i_{\Theta',d_0-1}}^-(\Theta'')$ and $\Theta''<_{1}\Theta$. Similar to \textbf{Case $1$}, among $p_{i'}^+(\Theta')$ for all possible choices of $i'\in (I_v\cap I_1)\setminus I'$, $\Theta=p_i^+(\Theta')$ and $\Theta''=p_{i_{\Theta',d_0-1}}^+(\Theta')$ are the only tuples that belong to $\Omega$ (and share the common $\un{\Lambda}$).
We thus deduce from Lemma~\ref{lem: tuple sum of two} (with $\Theta'$ and $d_0$ here replacing $\Theta$ and $d$ in \emph{loc.cit.}) that $c_{\Theta'}(d_{1,I_0,I_1}^{-\ell,k}(x_\Omega))=0$.

\textbf{Case $3$}: Assume that $\Lambda'_{d_0}=\emptyset\neq \Lambda'_{d_0+1}$ and $d_0=r_{v,I_1,I'}^{s-1}+1$.
As $\Theta=p_i^+(\Theta')$ is $(I_0,I_1)$-atomic, $\Lambda'_{d_0+1}\neq \emptyset$ and $d_0=r_{v,I_1,I'}^{s-1}+1$, we deduce from Lemma~\ref{lem: atom glue} (with $\Theta'$ and $i$ here replacing $\Theta$ and $i_{\Theta,d_0}$ in \emph{loc.cit.}) that $\Theta'$ is $(I_0,I_1)$-atomic. Note that we have $\Lambda_{d_0}=\Lambda'_{d_0}\sqcup\Lambda'_{d_0+1}\neq \emptyset$, from which we observe that among all possible choices of $i'\in (I_v\cap I_1)\setminus I'$, $p_{i'}^+(\Theta')\in\Omega$ if and only if $p_{i'}^+(\Theta')=\Theta$ with $i'=i$.

\textbf{Case $4$}: Assume that $\Lambda'_{d_0}\neq \emptyset\neq \Lambda'_{d_0+1}$.
As $\Theta=p_i^+(\Theta')$ is $(I_0,I_1)$-atomic and $\Lambda'_{d_0}\neq \emptyset\neq \Lambda'_{d_0+1}$, we deduce from Lemma~\ref{lem: atom glue} (with $\Theta'$ and $i$ here replacing $\Theta$ and $i_{\Theta,d_0}$ in \emph{loc.cit.}) that $\Theta'$ is $(I_0,I_1)$-atomic.
Now we consider $\Theta''=(v,I'',\un{k}'',\un{\Lambda}'')\defeq p_{i'}^+(\Theta')$ for some $i'\in (I_v\cap I_1)\setminus I'$, and write $i'=i_{\Theta',d_1}$ for some $1\leq d_1\leq r_{I'}$. If $d_1<d_0$ (resp.~$d_1>d_0$), then we have $\Lambda''_{d_0}=\Lambda'_{d_0+1}\subsetneq \Lambda_{d_0}$ (resp.~$\Lambda''_{d_0}=\Lambda'_{d_0}\subsetneq \Lambda_{d_0}$) and thus $\Theta''\notin \Omega$.

Now we consider an arbitrary equivalence class $\Omega'$ and a tuple $\Theta'\in\Omega'$ such that $c_{\Theta'}(d_{1,I_0,I_1}^{-\ell,k}(x_\Omega))\neq 0$. By definition of $x_\Omega$, there exists $\Theta\in\Omega$ such that $c_{\Theta'}(d_{1,I_0,I_1}^{-\ell,k}(x_{\Theta}))\neq 0$, and thus the pair $\Theta',\Theta$ falls into one of the four cases above. In \textbf{Case $1$} and \textbf{Case $2$}, we have $c_{\Theta'}(d_{1,I_0,I_1}^{-\ell,k}(x_\Omega))=0$, a contradiction to our choice of $\Theta'$. Hence, the pair $\Theta',\Theta$ must fall into either \textbf{Case $3$} or \textbf{Case $4$}.
In both \textbf{Case $3$} and \textbf{Case $4$}, the tuple $\Theta'$ is $(I_0,I_1)$-atomic and we have
$c_{\Theta'}(d_{1,I_0,I_1}^{-\ell,k}(x_\Omega))=\varepsilon(\Theta)c_{\Theta'}(d_{1,I_0,I_1}^{-\ell,k}(x_{\Theta}))$ which implies
\begin{multline}\label{equ: atom differential common}
\varepsilon(\Theta')^{-1}c_{\Theta'}(d_{1,I_0,I_1}^{-\ell,k}(x_\Omega))=\varepsilon(\Theta')^{-1}\varepsilon(\Theta)c_{\Theta'}(d_{1,I_0,I_1}^{-\ell,k}(x_{\Theta}))\\
=(-1)^{i_{\Theta',d_0}}(-1)^{m(I,i_{\Theta',d_0})}=(-1)^{i_{\Theta',d_0}-m(I,i_{\Theta',d_0})}
\end{multline}
with $m(I,i_{\Theta',d_0})=\#\{i\in I\mid i<i_{\Theta',d_0}\}$. We observe that
\begin{multline}\label{equ: absolute sign}
i_{\Theta',d_0}-m(I,i_{\Theta',d_0})=1+\#\{i\in \Delta\mid i<i_{\Theta',d_0}\}-\#\{i\in I\mid i<i_{\Theta',d_0}\}\\
=1+\#\{i\in \Delta\setminus I\mid i<i_{\Theta',d_0}\}=1+(d_0-1)=d_0,
\end{multline}
which depends only on the equivalence classes $\Omega$ and $\Omega'$. Hence, the sign (\ref{equ: atom differential common}) depends only on $\Omega$ and $\Omega'$ and we denote it by $\varepsilon(\Omega,\Omega')\in\{1,-1\}$. (In fact, $\Theta'\mapsto \Theta=p_{i_{\Theta',d_0}}^+(\Theta')$ defines a map $\Omega'\rightarrow \Omega$.) Now that we have $\varepsilon(\Theta')^{-1}c_{\Theta'}(d_{1,I_0,I_1}^{-\ell,k}(x_\Omega))=\varepsilon(\Omega,\Omega')$, we conclude that
\begin{equation}\label{equ: atom differential sum}
d_{1,I_0,I_1}^{-\ell,k}(x_\Omega)=\sum_{\Omega'}\varepsilon(\Omega,\Omega')x_{\Omega'}\in E_{1,I_0,I_1,\diamond,v}^{-\ell+1,k}
\end{equation}
where $\Omega'$ runs through all $(I_0,I_1)$-atomic equivalence classes satisfying the condition that there exists $\Theta'=(v,I',\un{k}',\un{\Lambda}')\in\Omega'$ and $i\in (I_v\cap I_1)\setminus I'$ such that $p_i^+(\Theta')\in\Omega$.
The proof is thus finished.
\end{proof}

\begin{defn}\label{def: simple partial order}
We consider all sets of the form $\mathbf{I}_{>d}\defeq (I^{d''})_{d''>d}$ for some $I_0\subseteq I\subseteq I_v\cap I_1$ with $\#I=\ell$ and $1\leq d\leq r_I=n-\ell$. Then the set $\{\mathbf{I}_{>d}\}_{I,d}$ admits a natural partial order described as follows. Given $\mathbf{I}_{>d}$ and $\mathbf{I}'_{>d'}$ for some $(I,d)$ and $(I',d')$, we say that $\mathbf{I}'_{>d'}<\mathbf{I}_{>d}$ if exactly one of the following holds
\begin{itemize}
\item $d'<d$ and $(I')^{d''}=I^{d''}$ for each $d''>d$;
\item there exists $d_\flat\geq \max\{d,d'\}$ such that $(I')^{d_\flat}\subsetneq I^{d_\flat}$ and $(I')^{d''}=I^{d''}$ for each $d''>d_\flat$.
\end{itemize}
It is not difficult to check that this is a well defined partial order on $\{\mathbf{I}_{>d}\}_{I,d}$.
\end{defn}

For each $(I_0,I_1)$-standard $\Theta$ (see Definition~\ref{def: std element}), we set $i_{I_0,I_1,\Theta}\defeq \sum_{d=1}^{d_{I_0,I_1,\Theta}-1}n_d\in (I_v\cap I_1)\setminus I$ and consider 
\[\Theta^+\defeq p_{i_{I_0,I_1,\Theta}}^+(\Theta)\] 
which is clearly well defined. For each $(I_0,I_1)$-standard $\Theta$, we associate the set $\mathbf{I}_\Theta\defeq (I^d)_{d>d_{I_0,I_1,\Theta}}$.
\begin{lem}\label{lem: independent elements}
For each $v\subseteq \mathbf{Log}_{\emptyset}$ with $\#v=k_0$, the following subset of $E_{1,I_0,I_1,v}^{-\ell,k}$
\begin{equation}\label{equ: std independent}
\{d_{1,I_0,I_1}^{-\ell-1,k}(x_{\Theta^+})\mid \Theta=(v,I,\un{k},\un{\Lambda})\text{ is }(I_0,I_1)\text{-standard }, |\un{k}|=k, \#I=\ell\}
\end{equation}
is linearly independent.
\end{lem}
\begin{proof}
We fix the choice of $k$, $\ell$, $v$, and write $S$ for the set of $(I_0,I_1)$-standard tuples $\Theta=(v,I,\un{k},\un{\Lambda})$ with bi-degree $(-\ell,k)$. We consider the subspace $E_S\subseteq E_{1,I_0,I_1,v}^{-\ell,k}$ spanned by the linearly independent set $\{x_{\Theta}\}_{\Theta\in S}$. We have a natural projection map $\mathrm{pr}_S: E_{1,I_0,I_1,v}^{-\ell,k}\twoheadrightarrow E_S$ given by $x\mapsto \sum_{\Theta\in S}c_{\Theta}(x)\Theta$.
For $\Theta,\Theta'\in S$, we write $\Theta \prec \Theta'$ if $\mathbf{I}_{\Theta}<\mathbf{I}_{\Theta'}$ under the partial-order introduced in Definition~\ref{def: simple partial order}. This defines a partial-order on $S$ and we fix the choice of a total-order on $S$ that refines this partial-order $\prec$. This total-order on $S$ together with the (unordered) set $\{x_{\Theta}\}_{\Theta\in S}$ determines a basis of $E_S$.

Now we consider two different tuples $\Theta=(v,I,\un{k},\un{\Lambda})$ and $\Theta'=(v,I',\un{k}',\un{\Lambda}')$ in $S$ that satisfy
\[c_{\Theta'}\left(d_{1,I_0,I_1}^{-\ell-1,k}(x_{\Theta^+})\right)\neq 0.\]
Then we have the following possibilities
\begin{itemize}
\item If $d_{I_0,I_1,\Theta'}\leq d_{I_0,I_1,\Theta}-1$, then we have $\mathbf{I}_{\Theta'}=\mathbf{I}_{\Theta}\sqcup\{(I')^{d}\}_{d_{I_0,I_1,\Theta'}<d\leq d_{I_0,I_1,\Theta}}$.
\item If $d_{I_0,I_1,\Theta'}\geq d_{I_0,I_1,\Theta}$, then $I^{d_{I_0,I_1,\Theta'}+1}\cap I_0\neq \emptyset$, and $\Theta^+=p_i^+(\Theta')$ with $i\in I^{d_{I_0,I_1,\Theta'}+1}$. Moreover, the tuple $\mathbf{I}_{\Theta'}=\{(I')^{d}\mid d> d_{I_0,I_1,\Theta'}\}$ satisfies $(I')^{d_{I_0,I_1,\Theta'}+1}\subsetneq I^{d_{I_0,I_1,\Theta'}+1}$ and $(I')^{d}=I^{d}$ for each $d>d_{I_0,I_1,\Theta'}+1$.
\end{itemize}
Hence, we always have $\mathbf{I}_{\Theta'}<\mathbf{I}_\Theta$ for partial order introduced in Definition~\ref{def: simple partial order}. In other words, if we consider the matrix that represents the set $\{\mathrm{pr}_S(d_{1,I_0,I_1}^{-\ell-1,k}(x_{\Theta^+}))\}_{\Theta\in S}$ using the set $\{x_{\Theta}\}_{\Theta\in S}$ (both ordered using the fixed total order on $S$), this matrix is unipotent and thus invertible. Hence, $\{\mathrm{pr}_S(d_{1,I_0,I_1}^{-\ell-1,k}(x_{\Theta^+}))\}_{\Theta\in S}$ is another basis of $E_S$, and in particular the set (\ref{equ: std independent}) is linearly independent.
\end{proof}

\begin{lem}\label{lem: vanishing of std}
For each $v\subseteq\mathbf{Log}_{\emptyset}$ and $0\neq x\in E_{1,I_0,I_1,v}^{-\ell,k}$, there exists $x'\in E_{1,I_0,I_1,v}^{-\ell-1,k}$ such that
\begin{equation}\label{equ: vanishing of std}
c_{\Theta}(x-d_{1,I_0,I_1}^{-\ell-1,k}(x'))=0
\end{equation}
for each $(I_0,I_1)$-standard $\Theta=(v,I,\un{k},\un{\Lambda})$.
\end{lem}
\begin{proof}
This follows immediately from Lemma~\ref{lem: independent elements} and its proof. In fact, we can construct $x'$ as a linear combination of various $x_{\Theta^+}$ for $(I_0,I_1)$-standard tuples $\Theta$, by an induction on the partial order introduced in Definition~\ref{def: simple partial order}.
\end{proof}

For each tuple $\Theta$, we define an integer
\[e_\Theta\defeq \sum_{d=1}^{r_I}\max\{m\mid m\in\Lambda_d\}\]
which will be useful in our later induction argument. (Here we use the convention $\max\{m\mid m\in\emptyset\}=0$.) It has the following simple property.
\begin{lem}\label{lem: numerical induction}
Let $\Theta$ be a tuple. We have the following results.
\begin{enumerate}[label=(\roman*)]
\item \label{it: numerical induction 1} If $i\in (I_v\cap I_1)\setminus I$ and $p_i^+(\Theta)$ is defined, then $e_{p_i^+(\Theta)}\leq e_\Theta$. Moreover, $e_{p_i^+(\Theta)}=e_\Theta$ if and only if either $\Lambda_d=\emptyset$ or $\Lambda_{d+1}=\emptyset$, where $1\leq d\leq r_I$ satisfies $i_{\Theta,d}=i$.
\item \label{it: numerical induction 2} If $i\in I$ and $p_i^-(\Theta)$ is defined, then $e_{p_i^-(\Theta)}= e_\Theta$.
\end{enumerate}
\end{lem}
\begin{proof}
This is immediate from the definition of $p_i^+(\Theta)$ (resp.~$p_i^-(\Theta)$). We use the fact that, if $\Lambda_d\cap\Lambda_{d+1}=\emptyset$, then
\[\max\{m\mid m\in\Lambda_d\sqcup\Lambda_{d+1}\}\leq\max\{m\mid m\in\Lambda_d\}+\max\{m\mid m\in\Lambda_{d+1}\},\]
and the equality holds if and only if either $\Lambda_d=\emptyset$ or $\Lambda_{d+1}=\emptyset$.
\end{proof}

\begin{defn}\label{def: coarse partial order}
Let $(-\ell,k)$ be a bi-degree and $\Theta=(v,I,\un{k},\un{\Lambda})$, $\Theta'=(v,I',\un{k}',\un{\Lambda}')$ be two tuples with bi-degree $(-\ell,k)$.
We write $\Theta\dashv\Theta'$ if exactly one of the following conditions holds.
\begin{enumerate}[label=(\roman*)]
\item \label{it: coarse partial order 1} There exists $1\leq s\leq r_{I_1\cap I_v}$ such that $r_{v,I_1,I}^{s'}=r_{v,I_1,I'}^{s'}$ for each $s'>s$ and $r_{v,I_1,I}^{s}>r_{v,I_1,I'}^{s}$.
\item \label{it: coarse partial order 2} We have $r_{v,I_1,I}^{s}=r_{v,I_1,I'}^{s}$ for each $1\leq s\leq r_{I_1\cap I_v}$ and $e_{\Theta'}<e_{\Theta}$.
\item \label{it: coarse partial order 3} We have $\un{\Lambda}=\un{\Lambda}'$ (with $e_{\Theta}=e_{\Theta'}$) and $\Theta<_{1}\Theta'$ in the sense of Definition~\ref{def: partial order}.
\end{enumerate}
It is easy to see that this defines a partial-order on the set of all tuples with bi-degree $(-\ell,k)$.
\end{defn}

\begin{lem}\label{lem: tuple reduction}
Let $v\subseteq\mathbf{Log}_{\emptyset}$ and $x\in E_{1,I_0,I_1,v}^{-\ell,k}$ such that $d_{1,I_0,I_1}^{-\ell,k}(x)\in E_{1,I_0,I_1,\diamond,v}^{-\ell+1,k}$ and $c_{\Theta}(x)=0$ for each $\Theta=(v,I,\un{k},\un{\Lambda})$ which is either $(I_0,I_1)$-standard or maximally $(I_0,I_1)$-atomic. Then $x=0$.
\end{lem}
\begin{proof}
Let $\Theta=(v,I,\un{k},\un{\Lambda})$ be a tuple with bi-degree $(-\ell,k)$.
We prove $c_{\Theta}(x)=0$ by a decreasing induction on the partial-order $\dashv$ introduced in Definition~\ref{def: coarse partial order}.
We assume inductively that $c_{\Theta'}(x)=0$ for all $\Theta'$ satisfying $\Theta\dashv\Theta'$.
It is harmless to assume that $\Theta$ is neither $(I_0,I_1)$-standard nor maximally $(I_0,I_1)$-atomic. Hence, there exists $1\leq s_1\leq r_{I_v\cap I_1}$ and $r_{v,I_1,I}^{s_1-1}+1\leq d_1\leq r_{v,I_1,I}^{s_1}$ and $i_1\in I^{d_1}$ such that $p_{i_1}^-(\Theta)$ is defined. We choose $d_1$ and $i_1$ to be maximal possible, and thus $p_{i_1}^-(\Theta)$ is $(I_0,I_1)$-standard. If $c_\Theta(x)\neq 0$, then there exists $i_2\in (I_v\cap I_1)\setminus I$ such that $c_{p_{i_2}^+p_{i_1}^-(\Theta)}(x)\neq 0$. 
Let $1\leq s_2\leq r_{I_v\cap I_1}$ and $r_{v,I_1,I}^{s_2-1}+1\leq d_2< r_{v,I_1,I}^{s_2}$ such that $i_2=i_{\Theta,d_2}$.
We have the following possibilities.
\begin{itemize}
\item If $d_2<d_1$ (with $i_2<i_1$), then $p_{i_2}^+p_{i_1}^-(\Theta)$ is $(I_0,I_1)$-standard with bi-degree $(-\ell,k)$. 
\item If $s_2>s_1$, then we have $r_{v,I_1,I}^{s}=r_{v,I_1,I'}^{s}$ for each $s\geq s_2$ and $r_{v,I_1,I}^{s_2-1}>r_{v,I_1,I'}^{s_2-1}$ and thus $\Theta\dashv\Theta'$ by \ref{it: coarse partial order 1} of Definition~\ref{def: coarse partial order}.
\item If $s_2=s_1$ and $d_2>d_1$, then as $\Lambda_{d}\neq \emptyset$ for each $r_{v,I_1,I}^{s_1-1}+1<d\leq r_{v,I_1,I}^{s_1}$ and thus $\Lambda_{d_2}\neq\emptyset\neq\Lambda_{d_2+1}$, we deduce from Lemma~\ref{lem: numerical induction} that $e_{p_{i_2}^+p_{i_1}^-(\Theta)}<e_\Theta$ and thus $\Theta\dashv\Theta'$ by \ref{it: coarse partial order 2} of Definition~\ref{def: coarse partial order}. 
\item If $d_2=d_1$, then we clearly have $\Theta<p_{i_2}^+p_{i_1}^-(\Theta)$ and thus $\Theta\dashv\Theta'$ by \ref{it: coarse partial order 3} of Definition~\ref{def: coarse partial order}. 
\end{itemize}
Now that we arrive at a contradiction in each case above, we conclude $c_\Theta(x)=0$ and finish the proof of the induction step.
\end{proof}

\begin{rem}\label{rem: atom key}
The fundamental idea behind Lemma~\ref{lem: independent elements} and Lemma~\ref{lem: tuple reduction} is the following observation. Let $(-\ell,k)$ be a bi-degree and $\Theta=(v,I,\un{k},\un{\Lambda})$ be a tuple. Then either $\Theta$ is maximally $(I_0,I_1)$-atomic, or $\Theta$ is $(I_0,I_1)$-standard, or there exists a maximal possible $i\in I$ such that the tuple $p_{i}^{-}(\Theta)$ is defined (with bi-degree $(-\ell+1,k)$), with moreover $p_{i}^{-}(\Theta)$ being $(I_0,I_1)$-standard.
\end{rem}

\begin{prop}\label{prop: quasi isom}
The embedding between complex
\begin{equation}\label{equ: atom complex embedding}
E_{1,I_0,I_1,\diamond,v}^{\bullet,k}\rightarrow E_{1,I_0,I_1,v}^{\bullet,k}
\end{equation}
is a quasi-isomorphism for each $v\subseteq\mathbf{Log}_{\emptyset}$ and $k\geq 0$.
In particular, the embedding between complex
\[E_{1,I_0,I_1,\diamond}^{\bullet,k}\rightarrow E_{1,I_0,I_1}^{\bullet,k}\]
is a quasi-isomorphism for each $k\geq 0$.
\end{prop}
\begin{proof}
The embedding (\ref{equ: atom complex embedding}) induces a map
\begin{equation}\label{equ: atom coh map}
H^{-\ell}(E_{1,I_0,I_1,\diamond,v}^{\bullet,k})\rightarrow H^{-\ell}(E_{1,I_0,I_1,v}^{\bullet,k})
\end{equation}
for each $\#I_0\leq \ell\leq \#I_1$ and $k\geq 0$.
We divide the proof into the following steps.

\textbf{Step $1$}: We prove that the map (\ref{equ: atom coh map}) is surjective, namely for each $x\in E_{1,I_0,I_1,v}^{-\ell,k}$ satisfying $d_{1,I_0,I_1,v}^{-\ell,k}(x)=0$, there exists $x'\in E_{1,I_0,I_1,v}^{-\ell-1,k}$ and $x''\in E_{1,I_0,I_1,\diamond, v}^{-\ell,k}$ such that $x=x''+d_{1,I_0,I_1}^{-\ell-1,k}(x')$.\\
Given such $x$, we can choose $x'$ as in Lemma~\ref{lem: vanishing of std}, and then take
\[x''\defeq \sum_{\Omega}c_{\Theta}(x-d_{1,I_0,I_1}^{-\ell-1,k}(x'))\varepsilon(\Theta)x_\Omega\in E_{1,I_0,I_1,\diamond, v}^{-\ell,k}\]
where $\Omega$ runs through all equivalence classes of $(I_0,I_1)$-atomic tuples with the fixed bidegree $(-\ell,k)$ and $v$ as above, and $\Theta\in\Omega$ is maximally $(I_0,I_1)$-atomic with sign $\varepsilon(\Theta)\in\{1,-1\}$ as in Definition~\ref{def: std element}. The equality $x=x''+d_{1,I_0,I_1}^{-\ell-1,k}(x')$ then follows from Lemma~\ref{lem: tuple reduction}.

\textbf{Step $2$}: We prove that the map (\ref{equ: atom coh map}) is injective, namely for each $y\in E_{1,I_0,I_1,\diamond,v}^{-\ell,k}$, and $y'\in E_{1,I_0,I_1,v}^{-\ell-1,k}$ satisfying $y=d_{1,I_0,I_1}^{-\ell-1,k}(y')$, there exists $y''\in E_{1,I_0,I_1,\diamond,v}^{-\ell-1,k}$ such that $y=d_{1,I_0,I_1}^{-\ell-1,k}(y'')$.\\
By Lemma~\ref{lem: vanishing of std} we may assume that $c_{\Theta}(y')=0$ for each $(I_0,I_1)$-standard $\Theta$. Then we take
\[y''\defeq \sum_{\Omega}c_{\Theta}(y')\varepsilon(\Theta)x_\Omega\in E_{1,I_0,I_1,\diamond, v}^{-\ell-1,k}\]
where $\Omega$ runs through all equivalence classes of $(I_0,I_1)$-atomic tuples with the fixed bidegree $(-\ell-1,k)$ and $v$ as above and $\Theta\in\Omega$ is maximally $(I_0,I_1)$-atomic with sign $\varepsilon(\Theta)\in\{1,-1\}$ as in Definition~\ref{def: std element}. Then we have
\[d_{1,I_0,I_1}^{-\ell-1,k}(y'-y'')=y-d_{1,I_0,I_1}^{-\ell-1,k}(y'')\in E_{1,I_0,I_1,\diamond, v}^{-\ell,k},\]
which together with Lemma~\ref{lem: tuple reduction} (with bidegree $(-\ell,k)$ there replaced with $(-\ell-1,k)$) implies that $y'-y''=0$ and thus $y=d_{1,I_0,I_1}^{-\ell-1,k}(y'')$. The proof is thus finished.
\end{proof}
\subsection{Computation of $E_{2}$ in bottom degree}\label{subsec: bottom deg}
We prove that $(I_0,I_1)$-atomic equivalence class with bidegree $(-\ell,k)$ exists only when $k\geq \ell+\#I_1-2\#I_0$, and then explicitly classify all $(I_0,I_1)$-atomic equivalence classes with bidegree $(-\ell,k)$ that satisfies $k-\ell\in\{\#I_1-2\#I_0,\#I_1-2\#I_0+1\}$.
As an application, we obtain an explicit basis of $E_{2,I_0,I_1}^{-\ell,k}$ for such $(-\ell,k)$ in Proposition~\ref{prop: atom basis}.

For each $1\leq d\leq r_I$ satisfying $I^d\cap I_0\neq\emptyset$, there exists a unique maximal possible $1\leq d'\leq r_{I_0}$ such that $\emptyset\neq I_0^{d'}\subseteq I^d$ and $I^d\setminus I_0^{d'}=I^{d,-}\sqcup I^{d,+}$ with $I^{d,-}\defeq \{i\in I^d\mid i<i' \text{ for each }i'\in I_0^{d'}\}$ and $I^{d,+}\defeq \{i\in I^d\mid i>i' \text{ for each }i'\in I_0^{d'}\}$. So we always have $I^{d,+}\cap I_0=\emptyset$ but $I^{d,-}\cap I_0$ might be non-empty. 

For each $i\in\Delta$, the inclusion $L_{\Delta\setminus\{i\}}\subseteq G$ induces a map
\begin{equation}\label{equ: restriction to maximal}
\mathrm{Res}_{n,i}^{\bullet}:~ H^{\bullet}(G,1_{G})\rightarrow H^{\bullet}(L_{\Delta\setminus\{i\}},1_{L_{\Delta\setminus\{i\}}}).
\end{equation}

For each $\Lambda\in\Sigma_{n,k}$ and $i\in\Delta$, we write $\mathbf{D}_{n,\Lambda,i}$ for those subsets $\Lambda'\subseteq \Lambda$ satisfying $\max\{m\mid m\in\Lambda'\}\leq 2i-1$ and $\max\{m\mid m\in\Lambda\setminus\Lambda'\}\leq 2(n-i)-1$. Note that both $\Lambda'$ and $\Lambda\setminus\Lambda'$ are allowed to be empty.
\begin{lem}\label{lem: restriction to maximal Levi}
Let $\Lambda\in\Sigma_{n,k}$ and $i\in\Delta$. Then we have $\mathrm{Res}_{n,i}^{k}(P^{\Lambda}_n)\neq 0$ if and only if $\mathbf{D}_{n,\Lambda,i}\neq \emptyset$, in which case we have
\[\mathrm{Res}_{n,i}^{k}(v_n^{\Lambda})=\sum_{\Lambda'\in\mathbf{D}_{n,\Lambda,i}}\varepsilon(\Lambda,\Lambda')v_{i}^{\Lambda'}\otimes_Ev_{n-i}^{\Lambda\setminus\Lambda'}\]
with $\varepsilon(\Lambda,\Lambda')\in\{1,-1\}$ being a sign determined by the pair $\Lambda'\subseteq \Lambda$.
\end{lem}
\begin{proof}
This follows from Theorem~\ref{thm: Koszul thm} and Theorem~\ref{thm: primitive class}.
\end{proof}

\begin{lem}\label{lem: vanishing of res}
Let $\Theta=(v,I,\un{k},\un{\Lambda})$ be a tuple and $1\leq d\leq r_I$ be an integer.
We have the following results.
\begin{enumerate}[label=(\roman*)]
\item \label{it: vanishing res 1} If $I^d\cap I_0=\emptyset$, then the following statements are equivalent
\begin{itemize}
\item there exists $i\in I^d$ such that $\mathrm{Res}_{I,I\setminus\{i\}}^{k}(x_\Theta)\neq 0$;
\item $n_d\geq 2$ and $\mathrm{Res}_{I,I\setminus\{i\}}^{k}(x_\Theta)\neq 0$ for $i\in \{i_{\Theta,d}-1, i_{\Theta,d-1}+1\}$;
\item $n_d\geq 2$ and $2n_d-1\notin\Lambda_d$.
\end{itemize}
\item \label{it: vanishing res 2} If $I^d\cap I_0\neq \emptyset$, then the following statements are equivalent
\begin{itemize}
\item there exists $i\in I^{d,+}$ such that $\mathrm{Res}_{I,I\setminus\{i\}}^{k}(x_\Theta)\neq 0$;
\item $I^{d,+}\neq \emptyset$ and $\mathrm{Res}_{I,I\setminus\{i\}}^{k}(x_\Theta)\neq 0$ for $i=i_{\Theta,d}-1$;
\item $I^{d,+}\neq \emptyset$ and $2n_d-1\notin\Lambda_d$.
\end{itemize}
\item \label{it: vanishing res 3} If $I^d\cap I_0\neq \emptyset$ and $I^{d,+}=\emptyset$, then the following statements are equivalent
\begin{itemize}
\item there exists $i\in I^{d,-}$ such that $c_{p_i^-(\Theta)}(\mathrm{Res}_{I,I\setminus\{i\}}^{k}(x_\Theta))\neq 0$;
\item $I^{d,-}\neq \emptyset$ and $c_{p_i^-(\Theta)}(\mathrm{Res}_{I,I\setminus\{i\}}^{k}(x_\Theta))\neq 0$ for $i=i_{\Theta,d-1}+\#I^{d,-}$;
\item $I^{d,-}\neq \emptyset$ and $\max\{m\mid m\in\Lambda_d\}<2\#I^{d,-}+1$.
\end{itemize}
\end{enumerate}
\end{lem}
\begin{proof}
We only prove \ref{it: vanishing res 1} and the other two cases are similar.
Let $i\in I^d$ and note that we have $1\leq i-i_{\Theta,d-1}\leq n_d-1$ and $n_d-(i-i_{\Theta,d-1})=i_{\Theta,d}-i$.
By applying Lemma~\ref{lem: restriction to maximal Levi} (with $G_{n_d}$ and $i-i_{\Theta,d-1}$ replacing $G$ and $i$ in \emph{loc.cit.}) we observe that $\mathrm{Res}_{I,I\setminus\{i\}}^{k}(x_\Theta)\neq 0$ if and only if $\mathbf{D}_{n_d,\Lambda_d,i-i_{\Theta,d-1}}\neq \emptyset$. Recall that $\mathbf{D}_{n_d,\Lambda_d,i-i_{\Theta,d-1}}\neq \emptyset$ means that there exists $\Lambda'\subseteq \Lambda$ such that $\max\{m\mid m\in\Lambda'\}\leq 2(i-i_{\Theta,d-1})-1\leq 2n_d-3$ and $\max\{m\mid m\in\Lambda_d\setminus\Lambda'\}\leq 2(n_d-(i-i_{\Theta,d-1}))-1\leq 2n_d-3$. Hence, $\mathbf{D}_{n_d,\Lambda_d,i-i_{\Theta,d-1}}\neq \emptyset$ only if $\max\{m\mid m\in\Lambda_d\}\leq 2n_d-3$, which is further equivalent to $2n_d-1\notin\Lambda_d$. As $\max\{m\mid m\in\Lambda_d\}\leq 2n_d-3$ clearly implies $\mathrm{Res}_{I,I\setminus\{i\}}^{k}(x_\Theta)\neq 0$ for $i\in \{i_{\Theta,d}-1, i_{\Theta,d-1}+1\}$, the proof is thus finished.
\end{proof}

\begin{lem}\label{lem: maximally atomic}
Let $\Theta=(v,I,\un{k},\un{\Lambda})$ be a maximally $(I_0,I_1)$-atomic tuple and  $1\leq d\leq r_I$. Recall the integers $r_{v,I_1,I}^{s}$ from (\ref{equ: v tuple integers}).
Then we have
\begin{enumerate}[label=(\roman*)]
\item \label{it: maximal atom 1} If $I^d\cap I_0=\emptyset$ and $2n_d-1\notin\Lambda_d$, then $n_d=1$, $\Lambda_d=\emptyset$ and there exists $1\leq s\leq r_{I_v\cap I_1}$ such that $d=r_{v,I_1,I}^{s-1}+1$.
\item \label{it: maximal atom 2} If $I^d\cap I_0\neq \emptyset$ and $I^{d,+}\neq \emptyset$, then we have $2n_d-1\in\Lambda_d$;
\item \label{it: maximal atom 3} If $I^d\cap I_0\neq \emptyset$, $I^{d,+}=\emptyset$ and $\max\{m\mid m\in\Lambda_d\}<2\#I^{d,-}+1$, then $I^{d,-}=\emptyset$, $\Lambda_d=\emptyset$ and there exists $1\leq s\leq r_{I_v\cap I_1}$ such that $d=r_{v,I_1,I}^{s-1}+1$.
\end{enumerate}
\end{lem}
\begin{proof}
We may assume throughout that $d\neq r_{v,I_1,I}^s$ for each $1\leq s\leq r_{I_v\cap I_1}$ (see Definition~\ref{def: std element}).

We prove \ref{it: maximal atom 1}.
Assume that $I^d\cap I_0=\emptyset$ and $2n_d-1\notin\Lambda_d$. If $n_d\geq 2$, then by \ref{it: vanishing res 1} of Lemma~\ref{lem: vanishing of res} we have $\Theta<p_{i_{\Theta,d}}^+p_{i_{\Theta,d}-1}^-(\Theta)$ which contradicts the maximality of $\Theta$. Hence, we have $n_d=1$, $k_d=0$ and $\Lambda_d=\emptyset$, which gives a $1\leq s\leq r_{I_v\cap I_1}$ such that $d=r_{v,I_1,I}^{s-1}+1$ by \ref{it: atom element 2} of Definition~\ref{def: atom element}.

We prove \ref{it: maximal atom 2}.
Assume that $I^d\cap I_0\neq \emptyset$ and $I^{d,+}\neq \emptyset$. If $2n_d-1\notin\Lambda_d$, then by \ref{it: vanishing res 2} of Lemma~\ref{lem: vanishing of res} we have $\Theta<p_{i_{\Theta,d}}^+p_{i_{\Theta,d}-1}^-(\Theta)$, which contradicts the maximality of $\Theta$.

We prove \ref{it: maximal atom 3}.
Assume that $I^d\cap I_0\neq \emptyset$, $I^{d,+}=\emptyset$ and $\max\{m\mid m\in\Lambda_d\}< 2\#I^{d,-}+1$. If $I^{d,-}\neq \emptyset$, then by \ref{it: vanishing res 3} of Lemma~\ref{lem: vanishing of res} we have \[\Theta<p_{i_{\Theta,d}}^+p_{i_{\Theta,d-1}+\#I^{d,-}}^-(\Theta),\]
which contradicts the maximality of $\Theta$. Hence, we must have $I^{d,-}=\emptyset$, which together with $\max\{m\mid m\in\Lambda_d\}< 2\#I^{d,-}+1$ forces $k_d=0$ (namely $\Lambda_d=\emptyset$), and thus there exists $1\leq s\leq r_{I_v\cap I_1}$ such that $d=r_{v,I_1,I}^{s-1}+1$ again by \ref{it: atom element 2} of Definition~\ref{def: atom element}.
\end{proof}

Let $v\subseteq \mathbf{Log}_{\emptyset}$ be a subset.
\begin{cond}\label{cond: bottom deg v}
We consider the following list of conditions.
\begin{enumerate}[label=(\roman*)]
\item \label{it: bottom atom v} We have $\#v\cap \mathbf{Log}_{\Delta\setminus\{i\}}=1$ for each $i\in\Delta\setminus I_v$.
\item \label{it: bottom atom v bis} There exists a unique $i_0\in\Delta\setminus I_v$ such that $\#v\cap \mathbf{Log}_{\Delta\setminus\{i_0\}}=2$ and $\#v\cap \mathbf{Log}_{\Delta\setminus\{i\}}=1$ for each $i_0\neq i\in\Delta\setminus I_v$.
\end{enumerate}
\end{cond}

Let $\Theta=(v,I,\un{k},\un{\Lambda})$ be a tuple, $1\leq s\leq r_{I_v\cap I_1}$ and $r_{v,I_1,I}^{s-1}+1\leq d\leq r_{v,I_1,I}^{s}$.
\begin{cond}\label{cond: bottom deg s}
We consider the following conditions.
\begin{enumerate}[label=(\roman*)]
\item \label{it: bottom atom 1} If $I^d\cap I_0=\emptyset$, then $\Lambda_d\neq \{2n_d-1\}$ if and only if $n_d=1$, $\Lambda_d=\emptyset$ and $d=r_{v,I_1,I}^{s-1}+1$.
\item \label{it: bottom atom 2} If $I^d\cap I_0\neq \emptyset$, then $I^d\cap I_0$ is an interval with $I^{d,+}=\emptyset$, and we have $\Lambda_d\neq \{2\#I^{d,-}+1\}$ if and only if $I^{d,-}=\emptyset$, $\Lambda_d=\emptyset$ and $d=r_{v,I_1,I}^{s-1}+1$.
\end{enumerate}
\end{cond}

Let $\Theta$, $s$ and $d$ be as above.
\begin{cond}\label{cond: bottom deg s bis}
We consider the following conditions.
\begin{enumerate}[label=(\roman*)]
\item \label{it: bottom atom 1 bis} If $I^d\cap I_0=\emptyset$, then we have $\Lambda_d=\{2n_d-1\}$.
\item \label{it: bottom atom 2 bis} If $I^d\cap I_0\neq \emptyset$, then $I^d\cap I_0$ is an interval with $I^{d,+}=\emptyset\neq I^{d,-}$ and $\Lambda_d=\{2\#I^{d,-}+1\}$.
\end{enumerate}
\end{cond}

Now we classify $(I_0,I_1)$-atomic tuples when $k-\ell\in\{\#I_1-2\#I_0,\#I_1-2\#I_0+1\}$.
\begin{lem}\label{lem: existence of atom}
Let $\Theta=(v,I,\un{k},\un{\Lambda})$ be a $(I_0,I_1)$-atomic tuple with bidegree $(-\ell,k)$ and $1\leq d\leq r_I$ be an integer. We have the following results.
\begin{enumerate}[label=(\roman*)]
\item \label{it: existence atom 1} We have
\begin{equation}\label{equ: atom deg bound}
\#I_1-2\#I_0\leq k-\ell\leq n^2-n.
\end{equation}
\item \label{it: existence atom 2} If $\Theta$ is maximal with $k-\ell=\#I_1-2\#I_0$, then we have $I_v\cup I_1=\Delta$, and $\Theta$ satisfies \ref{it: bottom atom v} of Condition~\ref{cond: bottom deg v} as well as Condition~\ref{cond: bottom deg s} for each $1\leq s\leq r_{I_v\cap I_1}$.
\item \label{it: existence atom 3} If $\Theta$ is maximal with $k-\ell=\#I_1-2\#I_0+1$, then either we have $\#I_v\cup I_1=n-2$ with $\Theta$ satisfying \ref{it: bottom atom v} of Condition~\ref{cond: bottom deg v} as well as Condition~\ref{cond: bottom deg s} for each $1\leq s\leq r_{I_v\cap I_1}$, or we have $I_v\cup I_1=\Delta$ with $\Theta$ satisfying \ref{it: bottom atom v bis} of Condition~\ref{cond: bottom deg v} as well as Condition~\ref{cond: bottom deg s} for each $1\leq s\leq r_{I_v\cap I_1}$, or we have $I_v\cup I_1=\Delta$ with $\Theta$ satisfying \ref{it: bottom atom v} of Condition~\ref{cond: bottom deg v}, Condition~\ref{cond: bottom deg s bis} for a unique choice of $1\leq s_0\leq r_{I_v\cap I_1}$ and Condition~\ref{cond: bottom deg s} for each $1\leq s\neq s_0\leq r_{I_v\cap I_1}$
\end{enumerate}
\end{lem}
\begin{proof}
It is harmless to assume that our $\Theta$ is maximally $(I_0,I_1)$-atomic (see Definition~\ref{def: atom element}).
Note that $v\cap \mathbf{Log}_{\Delta\setminus\{i\}}\neq \emptyset$ if and only if $i\in\Delta\setminus I_v$, which together with $\#\mathbf{Log}_{\Delta\setminus\{i\}}=2$ and $\#\Delta\setminus I_v=r_{I_v}-1$ (as well as $I\subseteq I_v$ which forces $r_{I_v}\leq r_I$) implies that
\begin{equation}\label{equ: v deg bound}
r_{I_v}-1 \leq k_0\leq 2(r_{I_v}-1)\leq 2(r_I-1).
\end{equation}
As $\sum_{d=1}^{r_I}n_d=n$, we notice that
\[\sum_{d=1}^{r_I}n_d^2=(\sum_{d=1}^{r_I}n_d)^2-2\sum_{1\leq d'<d''\leq r_I}n_{d'}n_{d''}\leq n^2-2(r_I-1),\]
which together with (\ref{equ: v deg bound}) implies that
\begin{multline*}
\sum_{d=0}^{r_I}k_d\leq k_0+\sum_{d=1}^{r_I}(n_d^2-1)= k_0-r_I+\sum_{d=1}^{r_I}n_d^2\\
\leq 2(r_I-1)-r_I+n^2-2(r_I-1)=n^2-r_I=n^2-n+\ell.
\end{multline*}
We have the following observations from Lemma~\ref{lem: maximally atomic}.
\begin{itemize}
\item If $I^d\cap I_0=\emptyset$ and there does not exist $1\leq s\leq r_{I_v\cap I_1}$ such that $d=r_{v,I_1,I}^{s-1}+1$, then $k_d\geq 2n_d-1$ and the equality holds if and only if $\Lambda_d=\{2n_d-1\}$.
\item If $I^d\cap I_0=\emptyset$ and there exists $1\leq s\leq r_{I_v\cap I_1}$ such that $d=r_{v,I_1,I}^{s-1}+1$, then either $n_d=1$ with $k_d=0$, or $k_d\geq 2n_d-1$. So we always have $k_d\geq 2n_d-2$ and the equality holds if and only if $n_d=1$.
\item If $I^d\cap I_0\neq \emptyset$ and there does not exist $1\leq s\leq r_{I_v\cap I_1}$ such that $d=r_{v,I_1,I}^{s-1}+1$, then $k_d\geq 2\#I^{d,-}+1\geq 2n_d-2\#I^d\cap I_0-1$ and the equality holds if and only if $I^d\cap I_0$ is an interval, $I^{d,+}=\emptyset$ and $\Lambda_d=\{2n_d-2\#I^d\cap I_0-1\}$. Note that $I^{d,+}\neq \emptyset$ would force $k_d\geq 2n_d-1>2\#I^{d,-}+1$.
\item If $I^d\cap I_0\neq \emptyset$ and there exists $1\leq s\leq r_{I_v\cap I_1}$ such that $d=r_{v,I_1,I}^{s-1}+1$, then either $I^{d,+}=I^{d,-}=\emptyset$ with $k_d\geq 0$, or $I^{d,+}=\emptyset\neq I^{d,-}$ with $k_d\geq 2\#I^{d,-}+1\geq 2n_d-2\#I^d\cap I_0-1$, or $I^{d,+}\neq \emptyset$ with $k_d\geq 2n_d-1$. So we always have $k_d\geq 2n_d-2\#I^d\cap I_0-2$ and the equality holds if and only if $I^{d,+}=I^{d,-}=\emptyset$.
\end{itemize}
Consequently, for each $1\leq s\leq r_{I_v\cap I_1}$, we sum up $k_d$ over $r_{v,I_1,I}^{s-1}+1\leq d\leq r_{v,I_1,I}^{s}$ and obtain
\begin{equation}\label{equ: atom s part}
\sum_{d=r_{v,I_1,I}^{s-1}+1}^{r_{v,I_1,I}^s}k_d\geq -1+\sum_{d=r_{v,I_1,I}^{s-1}+1}^{r_{v,I_1,I}^s}(2n_d-2\#I^d\cap I_0-1).
\end{equation}
Furthermore, by summing up (\ref{equ: atom s part}) over $1\leq s\leq r_{I_v\cap I_1}$ and using $\sum_{d=1}^{r_I}n_d=n$ as well as $I_0=\bigsqcup_{d=1}^{r_I}I^d\cap I_0$, we deduce that
\[\sum_{d=1}^{r_I}k_d\geq -r_{I_v\cap I_1}+2n-r_I-2\#I_0,\]
which together with $k_0\geq r_{I_v}-1$ implies that
\begin{multline}\label{equ: atom lower bound}
k=\sum_{d=0}^{r_I}k_d\geq (r_{I_v}-1)-r_{I_v\cap I_1}+2n-r_I-2\#I_0=(n-\#I_v)-1-(n-\#I_v\cap I_1)+2n-(n-\ell)-2\#I_0\\
=\ell+n-1+\#I_v\cap I_1-\#I_v-2\#I_0=\ell+n-1+\#I_1-\#I_v\cup I_1-2\#I_0\geq \ell+\#I_1-2\#I_0.
\end{multline}
We have the following observations from the proof of (\ref{equ: atom lower bound}).
\begin{itemize}
\item We have $k=\ell+\#I_1-2\#I_0$ if and only if $I_v\cup I_1=\Delta$, $k_0=r_{I_v}-1$ and (\ref{equ: atom s part}) is an equality for each $1\leq s\leq r_{I_v\cap I_1}$.
\item We have $k=\ell+\#I_1-2\#I_0+1$ if and only if either $\#I_v\cup I_1=n-2$, $k_0=r_{I_v}-1$ and (\ref{equ: atom s part}) is an equality for each $1\leq s\leq r_{I_v\cap I_1}$, or $I_v\cup I_1=\Delta$, $k_0=r_{I_v}$ and (\ref{equ: atom s part}) is an equality for each $1\leq s\leq r_{I_v\cap I_1}$, or $I_v\cup I_1=\Delta$, $k_0=r_{I_v}-1$ and there exists a unique $1\leq s_0\leq r_{I_v\cap I_1}$ such that (\ref{equ: atom s part}) is an equality for each $1\leq s\neq s_0\leq r_{I_v\cap I_1}$ and moreover
    \[\sum_{d=r_{v,I_1,I}^{s_0-1}+1}^{r_{v,I_1,I}^{s_0}}k_d=\sum_{d=r_{v,I_1,I}^{s_0-1}+1}^{r_{v,I_1,I}^{s_0}}(2n_d-2\#I^d\cap I_0-1).\]
\end{itemize}
According to the proof of (\ref{equ: v deg bound}) and (\ref{equ: atom s part}), we have the following observations.
\begin{itemize}
\item The equality $k_0=r_{I_v}-1$ holds if and only if \ref{it: bottom atom v} of Condition~\ref{cond: bottom deg v} holds. We have $k_0=r_{I_v}$ if and only if \ref{it: bottom atom v bis} of Condition~\ref{cond: bottom deg v} holds.
\item Let $1\leq s\leq r_{I_v\cap I_1}$. The inequality (\ref{equ: atom s part}) is an equality if and only if Condition~\ref{cond: bottom deg s} holds. We have
    \[\sum_{d=r_{v,I_1,I}^{s-1}+1}^{r_{v,I_1,I}^s}k_d=\sum_{d=r_{v,I_1,I}^{s-1}+1}^{r_{v,I_1,I}^s}(2n_d-2\#I^d\cap I_0-1)\]
    if and only if Condition~\ref{cond: bottom deg s bis} holds.
\end{itemize}
The proofs of \ref{it: existence atom 2} and \ref{it: existence atom 3} are thus finished.
\end{proof}

\begin{rem}\label{rem: bottom deg atom}
Let $\Theta=(v,I,\un{k},\un{\Lambda})$ be an $(I_0,I_1)$-atomic tuple with bidegree $(-\ell,k)$. By applying Lemma~\ref{lem: existence of atom} to the unique maximal element (see Lemma~\ref{lem: unique maximal tuple}) in the equivalence class of $\Theta$ (which share the same tuple of sets $\un{\Lambda}$ with $\Theta$), we have the following results.
\begin{enumerate}[label=(\roman*)]
\item \label{it: bottom deg atom tuple 1} If $k-\ell=\#I_1-2\#I_0$, then for each $1\leq s\leq r_{I_v\cap I_1}$ we have $\Lambda_{r_{v,I_1,I}^{s-1}+1}=\emptyset$ and $\#\Lambda_d=1$ for $r_{v,I_1,I}^{s-1}+1<d\leq r_{v,I_1,I}^{s}$.
\item \label{it: bottom deg atom tuple 2} If $k-\ell=\#I_1-2\#I_0+1$, then we have $\#\Lambda_d=1$ for $r_{v,I_1,I}^{s-1}+1<d\leq r_{v,I_1,I}^{s}$ and $1\leq s\leq r_{I_v\cap I_1}$. Moreover, there exists at most one $1\leq s_0\leq r_{I_v\cap I_1}$ such that $\Lambda_{r_{v,I_1,I}^{s_0-1}+1}\neq \emptyset$, in which case $\#\Lambda_{r_{v,I_1,I}^{s_0-1}+1}=1$.
\end{enumerate}
\end{rem}

\begin{lem}\label{lem: low deg}
Let $\Omega$ be an equivalence class of $(I_0,I_1)$-atomic tuples with $(-\ell_\Omega,k_\Omega)=(-\ell,k)$. Recall the integers $r_{\Omega}^{s}$ from (\ref{equ: v class integer}).
We have the following results.
\begin{enumerate}[label=(\roman*)]
\item \label{it: low deg 1} If $k-\ell=\#I_1-2\#I_0$, then we have $d_{1,I_0,I_1}^{-\ell,k}(x_{\Omega})=0$.
\item \label{it: low deg 2} If $k-\ell=\#I_1-2\#I_0+1$, then $d_{1,I_0,I_1}^{-\ell,k}(x_{\Omega})\neq 0$ if and only if $I_v\cup I_1=\Delta$ and there exists a unique $1\leq s_0\leq r_{I_v\cap I_1}$ such that $\Lambda_{r_{\Omega}^{s_0-1}+1}\neq \emptyset$ and $I^d\cap I_0\neq \emptyset$ for some $r_{\Omega}^{s_0-1}+1\leq d\leq r_{\Omega}^{s_0}$.
    Moreover, if $d_{1,I_0,I_1}^{-\ell,k}(x_{\Omega})\neq 0$, then there exists a unique equivalence class $\Omega'$ of $(I_0,I_1)$-atomic tuples such that $d_{1,I_0,I_1}^{-\ell,k}(x_{\Omega})=\pm x_{\Omega'}$, and $\Omega'$ also determines $\Omega$ uniquely.
\end{enumerate}
\end{lem}
\begin{proof}
Assume for the moment that $d_{1,I_0,I_1}^{-\ell,k}(x_\Omega)\neq 0$, by (\ref{equ: atom differential sum}) we know that there exists a $(I_0,I_1)$-atomic equivalence class $\Omega'$ (with bidegree $(-\ell+1,k)$), $\Theta'=(v,I',\un{k}',\un{\Lambda}')\in\Omega'$ and $i'\in(I_v\cap I_1)\setminus I'$ such that $\Theta=(v,I,\un{k},\un{\Lambda})\defeq p_{i'}^+(\Theta')\in\Omega$.
Thanks to Lemma~\ref{lem: equivalence glue} it is harmless to assume that $\Theta'$ is maximally $(I_0,I_1)$-atomic.
By the discussion before (\ref{equ: atom differential sum}), we know that the pair $\Theta',\Theta$ falls into either \textbf{Case $3$} or \textbf{Case $4$} in the proof of Lemma~\ref{lem: atomic subcomplex} (as $\Theta'$ is $(I_0,I_1)$-atomic).
We borrow the notation $d_0$ and $s$ from \emph{loc.cit.}.
If $k=\ell+\#I_1-2\#I_0$ and the pair $\Theta',\Theta$ falls into \textbf{Case $3$} (resp.~\textbf{Case $4$}), then by discussion in \emph{loc.cit.} we have $\Lambda_{r_{\Omega}^{s-1}+1}\neq \emptyset$ (resp.~$\#\Lambda_{d_0}\geq 2$), which contradicts \ref{it: bottom deg atom tuple 1} of Remark~\ref{rem: bottom deg atom}.
If $k=\ell+\#I_1-2\#I_0+1$ and the pair $\Theta',\Theta$ falls into \textbf{Case $4$}, then by discussion in \emph{loc.cit.} we have $\#\Lambda_{d_0}\geq 2$, which contradicts \ref{it: bottom deg atom tuple 2} of Remark~\ref{rem: bottom deg atom}.
We assume from now on that $k=\ell+\#I_1-2\#I_0+1$ and that the pair $\Theta',\Theta$ falls into \textbf{Case $3$}, in which case we have $d_0=r_{\Omega}^{s-1}+1$, $\Lambda'_{r_{\Omega}^{s-1}+1}=\emptyset$ and $\Lambda_{r_{\Omega}^{s-1}+1}=\Lambda'_{r_{\Omega}^{s-1}+2}\neq \emptyset$. By \ref{it: bottom deg atom tuple 2} of Remark~\ref{rem: bottom deg atom} we know that such $s$ (with $\Lambda_{r_{v,I_1,I}^{s-1}+1}\neq \emptyset$) is uniquely determined by $\Omega$, and thus $\Omega'$ is uniquely determined by $\Omega$.
If $I^d\cap I_0=\emptyset$ for all $r_{\Omega}^{s-1}+1\leq d\leq r_{\Omega}^s$, then we have $\Lambda_d=\{2n_{d}-1\}$ for each $r_{\Omega}^{s-1}+1\leq d\leq r_{\Omega}^s$ by Lemma~\ref{lem: unique maximal tuple} and \ref{it: existence atom 3} of Lemma~\ref{lem: existence of atom}, and in particular there does not exist $\Theta'=(v,I',\un{k}',\un{\Lambda}')$ such that $\Theta=p_{i_{\Theta',d_0}}^+(\Theta')$.
If $I^{d_1}\cap I_0\neq\emptyset$ for some $r_{\Omega}^{s-1}+1\leq d_1\leq r_{\Omega}^s$, then we define $\Omega'$ as the unique equivalence class that satisfies $p_{i_{\Theta',r_{\Omega}^{s-1}+1}}^{+}(\Theta')\in\Omega$ for each $\Theta'\in\Omega'$ (see \ref{it: unique maximal tuple 2} of Lemma~\ref{lem: unique maximal tuple}).
In fact, writing $\Theta_0$ for the unique maximal element in $\Omega$ and choosing $r_{\Omega}^{s-1}+1\leq d_1\leq r_{\Omega}^s$ to be the minimal integer satisfying $I^{d_1}\cap I_0\neq\emptyset$, we observe that $\Omega'$ is the equivalence class that contains
\begin{equation}\label{equ: low deg tuple transfer}
p_{i_{\Theta_0,r_{\Omega}^{s-1}}+1}^-\circ\cdots \circ p_{i_{\Theta_0,d}}^+\circ p_{i_{\Theta_0,d}+1}^- \circ\cdots \circ p_{i_{\Theta_0,d_1}}^+\circ p_{i_{\Theta_0,d_1}+1}^-(\Theta_0)
\end{equation}
as its unique maximal element.
For each $r_{\Omega}^{s-1}+1\leq d\leq r_{\Omega}^s$, recall from \ref{it: existence atom 3} of Lemma~\ref{lem: existence of atom} that we have either $I^d\cap I_0=\emptyset$ with $\Lambda_{d}=\{2n_{d}-1\}$ or $I^d\cap I_0\neq\emptyset$ with $I^{d,-}\neq \emptyset$ and $\Lambda_d=\{2\#I^{d,-}+1\}$. In particular, we have $i_{\Theta^{\max},d-1}+1\in I^{\max,d}$ for each $r_{\Omega}^{s-1}+1\leq d\leq r_{\Omega}^s$ and can check that (\ref{equ: low deg tuple transfer}) is well-defined.
We also note that such $\Omega'$ as defined above can be characterized by the following properties (see \ref{it: unique maximal tuple 2} of Lemma~\ref{lem: unique maximal tuple})
\begin{itemize}
\item $r_{\Omega'}^{s'}=r_{\Omega}^{s'}$ for each $1\leq s'\leq s-1$ and $r_{\Omega'}^{s'}=r_{\Omega}^{s'}+1$ for each $s\leq s'\leq r_{I_v\cap I_1}$;
\item $\Lambda'_{d'}=\Lambda_{d'}$ for each $1\leq d'\leq r_{\Omega'}^{s-1}$, $\Lambda'_{r_{\Omega'}^{s-1}+1}=\emptyset$ and $\Lambda'_{d'}=\Lambda_{d'-1}$ for each $r_{\Omega'}^{s-1}+2\leq d'\leq r_{I'}$.
\end{itemize}
Now that $s$ is the unique integer for $\Omega'$ to fail Condition~\ref{cond: bottom deg s}, we see that $\Omega'$ also determines $s$ and thus $\Omega$ uniquely. The proof is thus finished.
\end{proof}

Let $\Omega$ be an equivalence class of $(I_0,I_1)$-atomic tuples, $\Theta=(v,I,\un{k},\un{\Lambda})\in\Omega$ be an element and $1\leq s\leq r_{I_v\cap I_1}$ be an integer.
We consider the following condition on $\Theta$ and $s$.
\begin{cond}\label{cond: bottom deg s bis bis}
If $\Lambda_{r_{\Omega}^{s-1}+1}\neq \emptyset$, then we have $I^d\cap I_0=\emptyset$ and $\Lambda_d=\{2n_d-1\}$ for each $r_{\Omega}^{s-1}+1\leq d\leq r_{\Omega}^{s}$.
\end{cond}
Note that Condition~\ref{cond: bottom deg s bis bis} is actually independent of the choice of $\Theta\in\Omega$, and is purely a condition on $\Omega$ and $s$.

Recall that $E_{2,I_0,I_1}^{-\ell,k}=\mathrm{ker}(d_{1,I_0,I_1}^{-\ell,k})/\mathrm{im}(d_{1,I_0,I_1}^{-\ell-1,k})$ for each $\#I_0\leq \ell\leq \#I_1$ and $k\geq 0$.
\begin{prop}\label{prop: atom basis}
Let $\#I_0\leq \ell\leq \#I_1$ and $k\geq 0$. We have the following results.
\begin{enumerate}[label=(\roman*)]
\item \label{it: atom basis 1} If $k=\ell+\#I_1-2\#I_0$, then we have $d_{1,I_0,I_1}^{-\ell,k}(E_{1,I_0,I_1,\diamond}^{-\ell,k})=0$, and the composition of
\begin{equation}\label{equ: atom basis E2 1}
E_{1,I_0,I_1,\diamond}^{-\ell,k}\hookrightarrow \mathrm{ker}(d_{1,I_0,I_1}^{-\ell,k})\twoheadrightarrow E_{2,I_0,I_1}^{-\ell,k}
\end{equation}
is an isomorphism.
\item \label{it: atom basis 2} If $k=\ell+\#I_1-2\#I_0+1$, then the composition of
\begin{equation}\label{equ: atom basis E2 2}
E_{1,I_0,I_1,\diamond}^{-\ell,k}\cap \mathrm{ker}(d_{1,I_0,I_1}^{-\ell,k})\hookrightarrow \mathrm{ker}(d_{1,I_0,I_1}^{-\ell,k})\twoheadrightarrow E_{2,I_0,I_1}^{-\ell,k}
\end{equation}
is an isomorphism. Moreover, the intersection $E_{1,I_0,I_1,\diamond}^{-\ell,k}\cap \mathrm{ker}(d_{1,I_0,I_1}^{-\ell,k})$ is spanned by $x_{\Omega}$ with $\Omega$ running through equivalence classes of $(I_0,I_1)$-atomic tuples from \ref{it: existence atom 3} of Lemma~\ref{lem: existence of atom} which moreover satisfy Condition~\ref{cond: bottom deg s bis bis}.
\end{enumerate}
\end{prop}
\begin{proof}
We prove \ref{it: atom basis 1}.\\
Assume that $k=\ell+\#I_1-2\#I_0$.
We have $d_{1,I_0,I_1}^{-\ell,k}(E_{1,I_0,I_1,\diamond}^{-\ell,k})=0$ by \ref{it: low deg 1} of Lemma~\ref{lem: low deg}.
It follows from Lemma~\ref{lem: existence of atom} (and Lemma~\ref{lem: unique maximal tuple}) that there does not exist $(I_0,I_1)$-atomic equivalence class of bidegree $(-\ell-1,\ell+\#I_1-2\#I_2)$ and thus $E_{1,I_0,I_1,\diamond}^{-\ell-1,k}=0$ by definition. This together with Proposition~\ref{prop: quasi isom} finishes the proof.

We prove \ref{it: atom basis 2}.\\
Assume that $k=\ell+\#I_1-2\#I_0+1$.
The claim on $E_{1,I_0,I_1,\diamond}^{-\ell,k}\cap \mathrm{ker}(d_{1,I_0,I_1}^{-\ell,k})$ follows from \ref{it: low deg 2} of Lemma~\ref{lem: low deg}. This together with $d_{1,I_0,I_1}^{-\ell-1,k}(E_{1,I_0,I_1,\diamond}^{-\ell-1,k})=0$ from \ref{it: low deg 1} of Lemma~\ref{lem: low deg} (or \ref{it: atom basis 1} above) and Proposition~\ref{prop: quasi isom} finishes the proof.
\end{proof}

For $\#I_0\leq \ell\leq \#I_1$, we write $\Psi_{I_0,I_1,\ell}$ (resp.~$\Psi_{I_0,I_1,\ell}'$) for the set of equivalence classes $\Omega$ of $(I_0,I_1)$-atomic tuples satisfying $(-\ell_\Omega,k_\Omega)=(-\ell,\ell+\#I_1-2\#I_0)$ (resp.~satisfying $(-\ell_\Omega,k_\Omega)=(-\ell,\ell+\#I_1-2\#I_0+1)$ and $d_{1,I_0,I_1}^{-\ell,\ell+\#I_1-2\#I_0+1}(x_{\Omega})=0$, see \ref{it: atom basis 2} of Proposition~\ref{prop: atom basis}).
For each $\Omega\in \Psi_{I_0,I_1,\ell}$ (resp.~$\Omega\in \Psi_{I_0,I_1,\ell}'$), we write $\overline{x}_{\Omega}$ for the image of $x_{\Omega}$ under $\mathrm{ker}(d_{1,I_0,I_1}^{-\ell,k_{\Omega}})\twoheadrightarrow E_{2,I_0,I_1}^{-\ell,k_{\Omega}}$.
Proposition~\ref{prop: atom basis} has the following immediate consequence.
\begin{cor}\label{cor: bottom deg E2 basis}
Let $\#I_0\leq \ell\leq \#I_1$ and $k\geq 0$. We have the following results.
\begin{enumerate}[label=(\roman*)]
\item \label{it: bottom deg E2 basis 1} If $k-\ell=\#I_1-2\#I_0$, then $E_{2,I_0,I_1}^{-\ell,k}$ admits a basis of the form
\[\{\overline{x}_{\Omega}\}_{\Omega\in \Psi_{I_0,I_1,\ell}}.\]
\item \label{it: bottom deg E2 basis 2} If $k-\ell=\#I_1-2\#I_0+1$, then $E_{2,I_0,I_1}^{-\ell,k}$ admits a basis of the form
\[\{\overline{x}_{\Omega}\}_{\Omega\in \Psi_{I_0,I_1,\ell}'}.\]
\end{enumerate}
\end{cor}

In fact, we can index the basis of $E_{2,I_0,I_1}^{-\ell,\ell+\#I_1-2\#I_0}$ for varying $\#I_0\leq \ell\leq \#I_1$ (see \ref{it: bottom deg E2 basis 1} of Corollary~\ref{cor: bottom deg E2 basis}) using the set $\cS_{I_1\setminus I_0}$ introduced in \S \ref{subsec: coxeter partition}, which is the set of all subsets $S\subseteq \Phi^+$ that satisfy $\sum_{\al\in S}\al=\sum_{\al\in I_1\setminus I_0}\al$.
We consider $\#I_0\leq \ell\leq \#I_1$ and $\Omega\in\Psi_{I_0,I_1,\ell}$ with associated $v\subseteq\mathbf{Log}_{\emptyset}$ and $\Theta$ being the unique maximal element of $\Omega$. We define $I\defeq\{\al\in\Delta\mid \val_{\al}\in v\}$. We define $S\subseteq \Phi^+$ by requiring that $\al\in S\cap\Delta$ if and only if $v\cap\mathbf{Log}_{\Delta\setminus\{\al\}}\neq\emptyset$, and that $\al\in S\setminus \Delta$ if and only if $I_{\al}=\{i_{\Theta,d-1}\}\sqcup(J^{d}\setminus I_0)$ for some $1\leq d\leq n-\ell$ with $J^{d}\not\subseteq I_0$ (or equivalently $\Lambda_{d}\neq\emptyset$).
It is clear that $\Omega$ (or equivalently $\Theta=(v,J,\un{k},\un{\Lambda})$ with $\ell=\#J$) is uniquely determined by $I_0\subseteq I_1\subseteq \Delta$, the set $I\subseteq \Delta$, and the subset $S\subseteq \Phi^+$, with $\#S=n-1-\ell$.
In other words, we have constructed a bijection
\begin{equation}\label{equ: atom to partition}
\bigsqcup_{\ell=\#I_0}^{\#I_1}\Psi_{I_0,I_1,\ell}\buildrel\sim\over\longrightarrow\{(S,I)\mid S\in \cS_{I_0\setminus I_1}, I\subseteq S\cap\Delta\}
\end{equation}
which sends $\Omega$ to $(S,I)$ as above.
For each $\#I_0\leq \ell\leq \#I_1$, we define $\mathrm{gr}^{J}(\Psi_{I_0,I_1,\ell})$ (resp.~$\mathrm{gr}^{i}(\Psi_{I_0,I_1,\ell})$) as the subset of $\Psi_{I_0,I_1,\ell}$ consisting of those $\Omega$ whose associated $(S,I)$ satisfies $I=J$ (resp.~satisfies $\#I=i$).
It follows from (\ref{equ: E1 v decomposition}) that 
\[\mathrm{gr}^{J}(E_{2,I_0,I_1}^{-\ell,\ell+\#I_1-2\#I_0})=\bigoplus_{v\cap\mathbf{Log}_{\emptyset}^{\infty}=\mathbf{Log}_{\Delta\setminus J}^{\infty}}E_{2,I_0,I_1,v}^{-\ell,\ell+\#I_1-2\#I_0}\] 
admits a basis of the form
\[\{\overline{x}_{\Omega}\}_{\Omega\in\mathrm{gr}^{J}(\Psi_{I_0,I_1,\ell})}.\]
Similar fact holds if we replace $\mathrm{gr}^{J}(-)$ with $\mathrm{gr}^{i}(-)$ above.

Let $I_0'\subseteq I_0\subseteq I_1\subseteq \Delta$ and $J\supseteq I_1\setminus I_0'$. 
For each $I_0\subseteq I\subseteq I_1$, we write $L_{I,\flat}\defeq L_{I\setminus I_0'}$ and $L_{I,\flat\flat}\defeq L_{I\setminus I_0'}\cap H_{J}$ for short, with $\fl_{I,\flat}$ and $\fl_{I,\flat\flat}$ being their associated $E$-Lie algebra.
We also write $Z_{I,\flat}^{\dagger}\defeq Z_{I\setminus I_0'}^{\dagger}$ and $Z_{I,\flat\flat}^{\dagger}\defeq Z_{I\setminus I_0'}^{\dagger}\cap H_{J}$ for short.
We define $\cT_{I_0,I_1,\flat}^{\bullet,\bullet}$ and $\cT_{I_0,I_1,\flat\flat}^{\bullet,\bullet}$ as the double complex $\cT_{\Sigma}^{\bullet,\bullet}$ from (\ref{equ: general double complex}) by taking $M^{\bullet}_{I}$ to be $C^{\bullet}(L_{I,\flat})$ and $C^{\bullet}(L_{I,\flat\flat})$.
The natural restriction maps
\begin{equation}\label{equ: Levi restriction J}
C^{\bullet}(L_{I})\rightarrow C^{\bullet}(L_{I,\flat})\rightarrow C^{\bullet}(L_{I,\flat\flat})
\end{equation}
for each $I_0\subseteq I\subseteq I_1$ induces the following maps between double complex
\[\cT_{I_0,I_1}^{\bullet,\bullet}\rightarrow \cT_{I_0,I_1,\flat}^{\bullet,\bullet}\rightarrow \cT_{I_0,I_1,\flat\flat}^{\bullet,\bullet}\] 
and thus maps between the second page of their associated spectral sequences
\begin{equation}\label{equ: Tits shift E2}
E_{2,I_0,I_1}^{-\ell,k}\rightarrow E_{2,I_0,I_1,\flat}^{-\ell,k}\rightarrow E_{2,I_0,I_1,\flat\flat}^{-\ell,k}.
\end{equation}
We omit $I_0'$ from the notation when $I_0'=I_0$.
\begin{prop}\label{prop: Tits bottom shift}
Let $I_0'\subseteq I_0\subseteq I_1\subseteq\Delta$ with $\#I_0\leq \ell\leq \#I_1$ and $k\geq 0$. We have the following results.
\begin{enumerate}[label=(\roman*)]
\item \label{it: Tits bottom shift 1} If $k-\ell=\#I_1-2\#I_0$, then both maps in (\ref{equ: Tits shift E2}) are isomorphisms.
\item \label{it: Tits bottom shift 2} If $k-\ell=\#I_1-2\#I_0+1$, then the LHS map of (\ref{equ: Tits shift E2}) is an injection.
\end{enumerate}
\end{prop}
\begin{proof}
The inclusion $L_{I,\flat\flat}\subseteq L_{I,\flat}$ induces the following commutative diagram of maps between graded $E$-algebras (similar to the discussions around (\ref{equ: Levi Lie discrete center factor}))
\[
\xymatrix{
H^{\bullet}(L_{I,\flat},1_{L_{I,\flat}}) \ar^{\sim}[rrr] \ar[d] & & & \mathrm{Tot}(\wedge^{\bullet}\Hom(Z_{I,\flat}^{\dagger},E)\otimes_EH^{\bullet}(\fl_{I,\flat},1_{\fl_{I,\flat}})) \ar[d]\\
H^{\bullet}(L_{I,\flat\flat},1_{L_{I,\flat\flat}}) \ar^{\sim}[rrr] & & & \mathrm{Tot}(\wedge^{\bullet}\Hom(Z_{I,\flat\flat}^{\dagger},E)\otimes_EH^{\bullet}(\fl_{I,\flat\flat},1_{\fl_{I,\flat\flat}}))
}
\]
with the right vertical map being surjective. In particular, $H^{\bullet}(L_{I,\flat\flat},1_{L_{I,\flat\flat}})$ admits a basis indexed by $(I_0\setminus I_0',I_1\setminus I_0')$-tuples $\Theta_{\flat}=(v^{\flat},I,\un{k}^{\flat},\un{\Lambda}^{\flat})$ that satisfy $v^{\flat}\cap\mathbf{Log}_{J}=\emptyset$.

Recall from Proposition~\ref{prop: atom basis} that $E_{2,I_0,I_1}^{-\ell,k}$ admits a basis of the form
\begin{equation}\label{equ: Tits Levi restriction E2 basis}
\{\overline{x}_{\Omega}\}_{\Omega\in \star}
\end{equation}
with $\star=\Psi_{I_0,I_1,\ell}$ when $k-\ell=\#I_1-2\#I_0$, and $\star=\Psi_{I_0,I_1,\ell}'$ when $k-\ell=\#I_1-2\#I_0+1$.
Similarly, when either $\ast=\flat$ or $\ast=\flat\flat$, we can define (using variants of Lemma~\ref{lem: existence of atom} and Lemma~\ref{lem: low deg}) $\Psi_{I_0,I_1,\ast,\ell}$ and $\Psi_{I_0,I_1,\ast,\ell}'$ so that $\{\overline{x}_{\Omega_{\ast}}\}_{\Omega_{\ast}\in \Psi_{I_0,I_1,\ast,\ell}}$ forms a basis of $E_{2,I_0,I_1,\flat}^{-\ell,\ell+\#I_1-2\#I_0}$ and $\{\overline{x}_{\Omega_{\ast}}\}_{\Omega_{\ast}\in \Psi_{I_0,I_1,\ast,\ell}'}$ forms a basis of $E_{2,I_0,I_1,\flat}^{-\ell,\ell+\#I_1-2\#I_0+1}$.
Following the discussion around (\ref{equ: twist transform}), there exists a tuple of signs $\un{c}$ such that $\cT_{I_0,I_1,\flat}^{\bullet,\bullet}=\cT_{I_0\setminus I_0',I_1\setminus I_0'}^{\bullet,\bullet}[\#I_0',0][\un{c}]$ and thus $E_{2,I_0,I_1,\flat}^{\bullet,\bullet}=E_{2,I_0\setminus I_0',I_1\setminus I_0'}^{\bullet+\#I_0',\bullet}[\un{c}]$ under which we have natural bijections $\Psi_{I_0,I_1,\flat,\ell}\buildrel\sim\over\longrightarrow \Psi_{I_0\setminus I_0',I_1\setminus I_0',\ell-\#I_0'}$ and $\Psi_{I_0,I_1,\flat,\ell}'\buildrel\sim\over\longrightarrow \Psi_{I_0\setminus I_0',I_1\setminus I_0',\ell-\#I_0'}'$. We identify these pairs of sets via these natural bijections from now on.
Following the variants of Lemma~\ref{lem: existence of atom} and Lemma~\ref{lem: low deg} upon replacing $(G,I_0,I_1)$ in \emph{loc.cit.} with $(G,I_0\setminus I_0',I_1\setminus I_0')$ and $(H_{J},I_0\setminus I_0',I_1\setminus I_0')$, we see that there is a natural bijection $\Psi_{I_0,I_1,\flat,\ell}=\Psi_{I_0\setminus I_0',I_1\setminus I_0',\ell-\#I_0'}\buildrel\sim\over\longrightarrow \Psi_{I_0,I_1,\flat\flat,\ell}$ and thus an isomorphism 
\[E_{2,I_0,I_1,\flat}^{-\ell,\ell+\#I_1-2\#I_0}\buildrel\sim\over\longrightarrow E_{2,I_0,I_1,\flat\flat}^{-\ell,\ell+\#I_1-2\#I_0}\]
for each $\ell$.
It thus suffices to show that the LHS map of (\ref{equ: Tits shift E2}) is an isomorphism when $k-\ell=\#I_1-2\#I_0$ and is an injection when $k-\ell=\#I_1-2\#I_0+1$.

We fix a choice of $\Omega\in\star$ as in (\ref{equ: Tits Levi restriction E2 basis}), and write $\Theta^{\max}=(v,J,\un{k},\un{\Lambda})$ for the unique (see Lemma~\ref{lem: unique maximal tuple}) maximal $(I_0,I_1)$-atomic tuple in $\Omega$.
We also write $\Theta=(v,I,\un{k},\un{\Lambda})$ for an arbitrary tuple in $\Omega$. We write $\{n^{J}_{d}\}_{1\leq d\leq r_J=n-\ell}$ (resp.~$\{n_{d}\}_{1\leq d\leq r_I=n-\ell}$) for the tuple of integers associated with $J$ (resp.~with $I$).
Now that $\Theta\in\Omega$ and $\Theta^{\max}$ is the unique maximal element in $\Omega$, we have $i_{\Theta,d}\geq i_{\Theta^{\max},d}$ for each $1\leq d\leq n-\ell$ with the inequality being equality when $d=r_{\Omega}^{s}$ for each $1\leq s\leq r_{I_v\cap I_1}$.
We divide the rest of the proof into several steps.

\textbf{Step $1$}: We prove that the image of $x_{\Theta}$ under
\begin{equation}\label{equ: atom Levi shift restriction}
H^{k}(L_{I},1_{L_{I}})\rightarrow H^{k}(L_{I\setminus I_0},1_{L_{I\setminus I_0}})
\end{equation}
is non-zero only if $\Theta=\Theta^{\max}$.\\
Let $\Theta$ be a tuple in $\Omega$ such that the image of $x_{\Theta}$ under (\ref{equ: atom Levi shift restriction}) is non-zero.
We fix an arbitrary $1\leq s\leq r_{I_v\cap I_1}$ and prove by decreasing induction on $r_{\Omega}^{s-1}+1\leq d\leq r_{\Omega}^{s}$ that $i_{\Theta,d}=i_{\Theta^{\max},d}$.
It is clear that $i_{\Theta,r_{\Omega}^{s}}=i_{\Theta^{\max},r_{\Omega}^{s}}$.
We consider some $r_{\Omega}^{s-1}+1\leq d< r_{\Omega}^{s}$ and note that we have $i_{\Theta,d+1}=i_{\Theta^{\max},d+1}$ by induction hypothesis. Now that $d+1\geq r_{\Omega}^{s-1}+2$, we have $\Lambda_{d+1}\neq \emptyset$ by Condition~\ref{cond: bottom deg s bis}.
Since the image of $x_{\Theta}$ under (\ref{equ: atom Levi shift restriction}) is non-zero, the element of $H^{k_{d+1}}(G_{n_{d+1}},1_{G_{n_{d+1}}})$ corresponding to $\Lambda_{d+1}$ has non-zero image under the following restriction map
\begin{equation}\label{equ: Levi block shift restriction}
H^{k_{d+1}}(G_{n_{d+1}},1_{G_{n_{d+1}}})\rightarrow H^{k_{d+1}}(G_{n_{d+1}}\cap L_{I\setminus I_0},1_{G_{n_{d+1}}\cap L_{I\setminus I_0}}).
\end{equation}
More precisely, we have the following cases.
\begin{itemize}
\item If $J^{d+1}\cap I_0=\emptyset$ (or equivalently $J^{d+1}\setminus I_0=J^{d+1}$), then we have $\Lambda_{d+1}=\{2n^{J}_{d+1}-1\}$ by \ref{it: bottom atom 1} of Condition~\ref{cond: bottom deg s}, which together with \ref{it: primitive class 2} of Theorem~\ref{thm: primitive class} forces $n_{d+1}\geq n^{J}_{d+1}$. This together with $i_{\Theta,d+1}=i_{\Theta^{\max},d+1}$ and $i_{\Theta,d}\geq i_{\Theta^{\max},d}$ forces $i_{\Theta,d}= i_{\Theta^{\max},d}$ (and $n_{d+1}= n^{J}_{d+1}$).
\item If $J^{d+1}\cap I_0\neq\emptyset$, then $J^{d+1}\cap I_0$ is an interval with $J^{d+1,+}=\emptyset\neq J^{d+1,-}$ and $\Lambda_{d+1}=\{2\#J^{d+1,-}+1\}$ by \ref{it: bottom atom 2} of Condition~\ref{cond: bottom deg s}. Hence, $G_{n_{d+1}}\cap L_{I_1\setminus I_0}$ can be identified with the subgroup $\mathrm{PGL}_{\#I^{d+1,-}+1}\times(K^\times)^{\#I_0\cap I^{d+1}}$ inside $G_{n_{d+1}}\cong \mathrm{PGL}_{n_{d+1}}(K)$ with $I^{d+1,-}\subseteq J^{d+1,-}$ (from $i_{\Theta,d}\geq i_{\Theta^{\max},d}$). It thus follows again from \ref{it: primitive class 2} of Theorem~\ref{thm: primitive class} that the image of $P^{2\#J^{d+1,-}+1}(\fg_{n_{d+1}})$ under (\ref{equ: Levi block shift restriction}) is non-zero if and only if $\#I^{d+1,-}=\#J^{d+1,-}$, if and only if $I^{d+1,-}=J^{d+1,-}$, if and only if $i_{\Theta,d}=i_{\Theta^{\max},d}$.
\end{itemize}
This finishes the proof of \textbf{Step $1$}.

\textbf{Step $2$}: We prove that the map
\begin{equation}\label{equ: Tits shift E2 minimal}
E_{2,I_0,I_1}^{-\ell,k}\rightarrow E_{2,I_0,I_1,\flat}^{-\ell,k}
\end{equation}
is an injection when $k-\ell\leq \#I_1-2\#I_0+1$.\\
Now that the map (\ref{equ: Tits shift E2 minimal}) for $I_0'=I_0$ factors through the map (\ref{equ: Tits shift E2 minimal}) for general $I_0'\subseteq I_0$, it suffices to prove the injectivity of (\ref{equ: Tits shift E2 minimal}) when $I_0'=I_0$ and $k-\ell\leq \#I_1-2\#I_0+1$. We assume in the rest of this step that $I_0'=I_0$.
When $\Theta=\Theta^{\max}$, the image of $x_{\Theta}$ under (\ref{equ: atom Levi shift restriction}) is $\pm x_{\Theta_{\flat}}$ for a maximal $(\emptyset,I_1\setminus I_0)$-atomic tuple $\Theta_{\flat}=(v,J\setminus I_0,\un{k}^{\flat},\un{\Lambda}^{\flat})$ uniquely determined by $\Theta$ via the following condition: if $G_{n_{d}^{J}}\cap L_{I\setminus I_0}$ is the $e$-th Levi block of $L_{I\setminus I_0}$, then we have $\Lambda^{\flat}_{e}\defeq \Lambda_{d}$.
Moreover, the $(\emptyset,I_1\setminus I_0)$-equivalence class $\Omega_{\flat}$ containing $\Theta_{\flat}$ actually has $\Theta_{\flat}$ as its unique element. The association $\Theta^{\max}\mapsto\Theta_{\flat}$ determines an association $\Omega\mapsto\Omega_{\flat}$ which gives an injection
\[\Psi_{I_0,I_1,\ell}\hookrightarrow\Psi_{\emptyset,I_1\setminus I_0,\ell-\#I_0}\]
when $k=\ell+\#I_1-2\#I_0$, and an injection
\[\Psi_{I_0,I_1,\ell}'\hookrightarrow\Psi_{\emptyset,I_1\setminus I_0,\ell-\#I_0}'\]
when $k=\ell+\#I_1-2\#I_0+1$, with the image of $x_{\Omega}$ under
\[E_{1,I_0,I_1}^{-\ell,k}\rightarrow E_{1,\emptyset,I_1\setminus I_0}^{-\ell+\#I_0,k}\]
given by $\pm x_{\Omega^{\flat}}$. Now that the map (\ref{equ: Tits shift E2 minimal}) (when $I_0'=I_0$) sends the basis (\ref{equ: Tits Levi restriction E2 basis}) of $E_{2,I_0,I_1}^{-\ell,k}$ to the linearly independent subset
\[\{\overline{x}_{\Omega_{\flat}}\}_{\Omega\in \star}\]
of $E_{2,\emptyset,I_1\setminus I_0}^{-\ell+\#I_0,k}$, it must be an injection.

\textbf{Step $3$}: We prove that the map (\ref{equ: Tits shift E2 minimal}) is an isomorphism when $k-\ell=\#I_1-2\#I_0$.\\
Now that the map (\ref{equ: Tits shift E2 minimal}) factors through $E_{2,I_0,I_1,I_0',\flat}^{-\ell,k}$, we deduce from \textbf{Step $2$} that (\ref{equ: Tits shift E2 minimal}) is injective when $k-\ell\leq \#I_1-2\#I_0+1$. 
Upon replacing the pair $(I_0,I_1)$ with $(I_0\setminus I_0',I_1\setminus I_0')$ in (\ref{equ: atom to partition}), we have
\[
\sum_{\ell=\#I_0}^{\#I_1}\Dim_EE_{2,I_0,I_1}^{-\ell,\ell+\#I_1-2\#I_0}=\#\{(S,I)\mid S\in\cS_{I_1\setminus I_0},I\subseteq S\cap\Delta\}
=\sum_{\ell=\#I_0}^{\#I_1}\Dim_EE_{2,I_0,I_1,I_0',\flat}^{-\ell,\ell+\#I_1-2\#I_0}.
\]
This forces the map (\ref{equ: Tits shift E2 minimal}) (which is an injection by \textbf{Step $2$}) to be an isomorphism for each $\ell$. The proof is thus finished.
\end{proof}

\subsection{Degeneracy at $E_{2}$ in bottom degree}\label{subsec: bottom deg degeneracy}
We prove that $E_{\bullet,I_0,I_1}^{\bullet,\bullet}$ degenerates at the second page for bottom non-vanishing degree (see Proposition~\ref{prop: bottom deg degeneracy}).

Let $I_0\subseteq I_1\subseteq \Delta$ and $\Omega$ be an equivalence class of $(I_0,I_1)$-atomic tuples with $\Theta=(v,J,\un{k},\un{\Lambda})$ being its unique maximal element (see \ref{it: unique maximal tuple 1} of Lemma~\ref{lem: unique maximal tuple}). We consider the following condition.
\begin{cond}\label{cond: unique atom}
For each $1\leq d\leq r_{J}=n-\ell_{\Omega}$, we have $J^{d}\cap I_0\neq \emptyset$ only if $\Lambda_{d}=\emptyset$.
\end{cond}
\begin{lem}\label{lem: unique atom}
Let $\Omega$ and $\Theta$ be as above. If $\Omega$ satisfies Condition~\ref{cond: unique atom}, then we have $\Omega=\{\Theta\}$.
\end{lem}
\begin{proof}
We write $\{n_{d}\}_{1\leq d\leq r_{J}}$ for the tuple of integers associated with $J$.
Assume on the contrary that there exists $\Theta'=(v,J',\un{k},\un{\Lambda})\in\Omega$ such that $\Theta'<\Theta$. Then there exists $1\leq d\leq r_{J'}=r_{J}=n-\ell_{\Omega}$ such that $i_{\Theta',d}>i_{\Theta,d}$. Now that we have $i_{\Theta',r_{\Omega}^{s}}=i_{\Theta,r_{\Omega}^{s}}$ for each $1\leq s\leq r_{I_v\cap I_1}$, We can choose such $d$ to be maximal possible and there exists a unique $1\leq s\leq r_{I_v\cap I_1}$ such that $r_{\Omega}^{s-1}< d < r_{\Omega}^{s}$. This together with $\Theta$ being $(I_0,I_1)$-atomic (see \ref{it: atom element 2} of Definition~\ref{def: atom element}) forces $\Lambda_{d+1}\neq \emptyset$ with $d+1\leq r_{\Omega}^{s}$. Now that $\Omega$ satisfies Condition~\ref{cond: unique atom} by our assumption, we know that $J^{d+1}\cap I_0=\emptyset$ and $2n_{d+1}-1\in\Lambda_{d+1}$. The maximality of $d$ ensures that $i_{\Theta',d+1}=i_{\Theta,d+1}$ and thus $i_{\Theta',d+1}-i_{\Theta',d}<i_{\Theta,d+1}-i_{\Theta,d}=n_{d}$, contradicting the fact that $2n_{d+1}-1\in\Lambda_{d+1}$. The proof is thus finished.
\end{proof}

We fix temporarily an equivalence class $\Omega$ of $(I_0,I_1)$-atomic tuples with bi-degree $(-\ell_{\Omega},k_{\Omega})$. We assume that $\Omega$ satisfies Condition~\ref{cond: unique atom} and moreover satisfies either $\Omega\in\Psi_{I_0,I_1,\ell_{\Omega}}$ or $\Omega\in\Psi_{I_0,I_1,\ell_{\Omega}}'$, with $\Theta=(v,J,\un{k},\un{\Lambda})$ being the unique (maximal) element of $\Omega$ (see Lemma~\ref{lem: unique atom}).
We write $\{n_{d}\}_{1\leq d\leq r_{J}}$ for the tuple of integers associated with $J$ and $k_0=\#v$.
Our assumption on $\Omega$ ensures the following facts (see Lemma~\ref{lem: existence of atom} and Proposition~\ref{prop: atom basis}).
\begin{itemize}
\item If $\Lambda_{d}=\emptyset$, then there exists $1\leq s\leq r_{I_v\cap I_1}$ such that $d=r_{\Omega}^{s-1}+1$ and $J^{d}\subseteq I_0$.
\item If $\Lambda_{d}\neq\emptyset$, then we have $J^{d}\cap I_0=\emptyset$ and $\Lambda_{d}=\{2n_{d}-1\}$ with $n_{d}=1+\#J^{d}\geq 2$.
\end{itemize}

Let $1\leq s\leq r_{I_v\cap I_1}$ and $r_{\Omega}^{s-1}+1\leq d\leq r_{\Omega}^{s}$.
We write $H_{\Omega,d}\defeq H_{J^{d}}$ and $\fh_{\Omega,d}$ for its associated $E$-Lie algebra.
Now that each connected component of $I_0$ has the form $J^{d}$ for some $d$ satisfying $\Lambda_{d}=\emptyset$ (and thus $d=r_{\Omega}^{s-1}+1$ for some $1\leq s\leq r_{I_v\cap I_1}$), we deduce that
\begin{equation}\label{equ: Omega empty block}
H_{I_0}\cong \prod_{d,\Lambda_{d}=\emptyset}H_{\Omega,d},
\end{equation}
which is a normal subgroup of $L_{J}$.

For each $s$ and $d$ satisfying $\Lambda_{d}=\{2n_{d}-1\}$ (with $J^{d}\cap I_0=\emptyset$ and $n_{d}=1+\#J^{d}\geq 2$), we write $H_{\Omega,d,\dagger}\defeq H_{J^{d}\setminus\{i_{\Theta,d}-1\}}$ for short with $\fh_{\Omega,d,\dagger}$ being its associated $E$-Lie algebra. We define 
\[\fh_{\Omega,\dagger}\defeq \prod_{d,\Lambda_{d}\neq\emptyset}\fh_{\Omega,d,\dagger}\subseteq \fh_{\Omega}\defeq \prod_{d,\Lambda_{d}\neq\emptyset}\fh_{\Omega,d}\cong \fl_{J}/\fz_{J}\fh_{I_0}.\]

Let $I\subseteq J$ be a subset and $\{m_{e}\}_{1\leq e\leq r_{I}}$ be the associated tuple of integers.
Let $1\leq d\leq r_{J}$ with $\Lambda_{d}\neq \emptyset$. There exists a unique $1\leq e\leq r_{I}$ such that $i_{\Theta,d}=\sum_{e'=1}^{e}m_{e'}$.
If $m_{e}=1$, then $\ft\subseteq \fl_{I}$ induces the following isomorphism
\begin{equation}\label{equ: relative Lie torus}
\ft\cap \fh_{\Omega,d}/\ft\cap \fh_{\Omega,d,\dagger}\buildrel\sim\over\longrightarrow \fl_{I}\cap \fh_{\Omega,d}/\fl_{I}\cap \fh_{\Omega,d,\dagger}
\end{equation}
between $1$-dimensional $E$-vector spaces.
If $m_{e}=1+\#I^{e}\geq 2$ (with $I^{e}\neq \emptyset$), then $\fg_{m_{e}}\cong \fh_{I^{e}}\subseteq \fl_{I}$ induce the following isomorphism
\begin{equation}\label{equ: relative Lie ss}
\fg_{m_{e}}/\fg_{m_{e}-1}\buildrel\sim\over\longrightarrow \fl_{I}\cap \fh_{\Omega,d}/\fl_{I}\cap \fh_{\Omega,d,\dagger}
\end{equation}
between $E$-vector spaces.

For each $I_0\subseteq I\subseteq J$, we view $\fl_{I}\cap \fh_{\Omega,\dagger}\subseteq \fl_{I}\cap \fh_{\Omega}$ as $E$-Lie subalgebras of $\fl_{I}/\fh_{I_0}\fz_{J}$.
Now we define $\cT_{\Omega}^{\bullet,\bullet}$ as $\cT_{I_0,I_1,v,\fh_{\Omega,\dagger},J,Z_{J}H_{I_0}}^{\bullet,\bullet}$ (see (\ref{equ: double v relative general})).
In other words, we define $\cT_{\Omega}^{\bullet,\bullet}$ as $\cT_{\Sigma}^{\bullet,\bullet}$ by taking $M^{\bullet}_{I}$ to be zero when $I\not\subseteq J$, and to be
\[C^{\bullet}(L_{I}/H_{I_0}Z_{J},\fl_{I}\cap \fh_{\Omega,\dagger},1_{L_{I}/H_{I_0}Z_{J}})\otimes_E v\]
when $I_0\subseteq I\subseteq J$.
We deduce from (\ref{equ: double fix v}), (\ref{equ: double v relative general}) and (\ref{equ: double v relative spectral seq}) the following (quasi) maps between double complex
\begin{equation}\label{equ: Omega double transfer}
\cT_{\Omega}^{\bullet,\bullet}\rightarrow \cT_{I_0,I_1,v}^{\bullet,\bullet}\dashrightarrow \cT_{I_0,I_v\cap I_1}^{\bullet,\bullet}\rightarrow \cT_{I_0,I_1}^{\bullet,\bullet}
\end{equation}
which induce maps between spectral sequences
\begin{equation}\label{equ: Omega spectral seq transfer}
E_{\bullet,\Omega}^{\bullet,\bullet}\rightarrow E_{\bullet,I_0,I_1,v}^{\bullet,\bullet}\rightarrow E_{\bullet,I_0,I_v\cap I_1}^{\bullet,\bullet}\rightarrow E_{\bullet,I_0,I_1}^{\bullet,\bullet}.
\end{equation}
\begin{lem}\label{lem: Omega seq lift}
Let $\Omega$ be as above. We have the following results.
\begin{enumerate}[label=(\roman*)]
\item \label{it: Omega seq lift 1} We have $E_{1,\Omega}^{-\ell,k}=0$ for each $(-\ell,k)$ satisfying $k-\ell>k_{\Omega}-\ell_{\Omega}$. Moreover, we have
\begin{equation}\label{equ: Omega infinite page}
E_{1,\Omega}^{-\ell_{\Omega},k_{\Omega}}=E_{\infty,\Omega}^{-\ell_{\Omega},k_{\Omega}}
\end{equation}
\item \label{it: Omega seq lift 2} The natural map (from the composition of (\ref{equ: Omega spectral seq transfer}))
\begin{equation}\label{equ: Omega seq lift E1}
E_{1,\Omega}^{-\ell_{\Omega},k_{\Omega}}\rightarrow E_{1,I_0,I_1}^{-\ell_{\Omega},k_{\Omega}}
\end{equation}
is an embedding with image $Ex_{\Omega}$.
\end{enumerate}
\end{lem}
\begin{proof}
We prove \ref{it: Omega seq lift 1}.\\
Note that we have
\begin{equation}\label{equ: Omega E1}
E_{1,\Omega}^{-\ell,k}=v\otimes_E \bigoplus_{I\subseteq J, \#I=\ell}H^{k}(L_{I}/H_{I_0}Z_{J},\fl_{I}\cap\fh_{\Omega,\dagger},1_{L_{I}/H_{I_0}Z_{J}})
\end{equation}
for each $0\leq \ell\leq \#J$ and $k\geq 0$.
Let $I\subseteq J$ with $\#I=\ell$ and $\{m_{e}\}_{1\leq e\leq r_{I}}$ be the tuple of integers associated with $I$.
Recall from Proposition~\ref{prop: Levi decomposition} that we have the following isomorphism between graded $E$-algebras
\begin{multline}\label{equ: Omega Lie sm}
H^{\bullet}(L_{I}/H_{I_0}Z_{J},\fl_{I}\cap\fh_{\Omega},1_{L_{I}/H_{I_0}Z_{J}})\\
\cong \wedge^{\bullet}\Hom(Z_{I}^{\dagger}/Z_{J}^{\dagger},E)\otimes_E H^{\bullet}(\fl_{I}/\fh_{I_0}\fz_{J},\fl_{I}\cap\fh_{\Omega,\dagger},1_{\fl_{I}/\fh_{I_0}\fz_{J}}).
\end{multline}
Note that we have
\begin{equation}\label{equ: Omega sm dim}
\Dim_E \Hom(Z_{I}^{\dagger}/Z_{J}^{\dagger},E)=\#J\setminus I=\#J-\#I=\ell_{\Omega}-\ell
\end{equation}
and in particular $\wedge^{k}\Hom_{\infty}(Z_{I}/Z_{J},E)=0$ when $k>\ell_{\Omega}-\ell$.
It follows from the definition of $\fh_{\Omega}$ and $\fh_{\Omega,\dagger}$ (as well as the K\"unneth formula for relative Lie algebra cohomology, see (\ref{equ: relative Lie Kunneth})) that we have the following isomorphisms between graded $E$-algebras
\begin{equation}\label{equ: Omega Lie decomposition}
H^{\bullet}(\fl_{I}/\fh_{I_0}\fz_{J},\fl_{I}\cap\fh_{\Omega,\dagger},1_{\fl_{I}/\fh_{I_0}\fz_{J}})
\cong \bigotimes_{d,\Lambda_{d}\neq\emptyset}H^{\bullet}(\fl_{I}\cap\fh_{\Omega,d},\fl_{I}\cap\fh_{\Omega,d,\dagger},1_{\fl_{I}\cap\fh_{\Omega,d}})
\end{equation}
For each $1\leq d\leq r_{J}$ with $\Lambda_{d}\neq\emptyset$, there exists a unique $1\leq e\leq r_{I}$ such that $i_{\Theta,d}=\sum_{e'=1}^{e}m_{e'}$, and we write $(I\cap J^{d})^{+}\defeq I^{e}$ with $\#I^{e}=m_{e}-1$.
It follows from the definition of $\fh_{\Omega,d}$ and $\fh_{\Omega,d,\dagger}$ as well as \ref{it: primitive class 3} of Theorem~\ref{thm: primitive class} that
\[H^{k}(\fl_{I}\cap\fh_{\Omega,d},\fl_{I}\cap\fh_{\Omega,d,\dagger},1_{\fl_{I}\cap\fh_{\Omega,d}})\cong H^{k}(\fg_{m_{e}},\fg_{m_{e}-1},1_{\fg_{m_{e}}})\]
is $1$-dimensional when $k=0,2m_{e}-1$ and is zero otherwise.
Note that we have
\begin{equation}\label{equ: index relation}
k_{\Omega}'\defeq k_{\Omega}-\#v=\sum_{d,\Lambda_{d}\neq \emptyset}(2\#J^{d}+1)=2\#J\setminus I_0+\#\{d\mid \Lambda_{d}\neq \emptyset\}.
\end{equation}
Consequently, we deduce that (\ref{equ: Omega Lie decomposition}) is zero for degree $k'$ when
\[k'>\sum_{d,\Lambda_{d}\neq \emptyset}(2\#(I\cap J^{d})^{+}+1),\]
which together with (\ref{equ: Omega sm dim}) and (\ref{equ: index relation}) implies that (\ref{equ: Omega Lie sm}) is zero for degree $k'$ when
\begin{multline*}
k'-\ell>\delta(I)\defeq -\#I+\#J\setminus I+\sum_{d,\Lambda_{d}\neq \emptyset}(2\#(I\cap J^{d})^{+}+1)\\
=-\#I+\#J\setminus I+\#\{d\mid \Lambda_{d}\neq \emptyset\}+2\sum_{d,\Lambda_{d}\neq \emptyset}\#(I\cap J^{d})^{+}\\
=k_{\Omega}'-2\#J\setminus I_0-\#I+\#J\setminus I+2\sum_{d,\Lambda_{d}\neq \emptyset}\#(I\cap J^{d})^{+}\\
=k_{\Omega}'-\#J-2\#I\setminus I_0+2\sum_{d,\Lambda_{d}\neq \emptyset}\#(I\cap J^{d})^{+}\\
=k_{\Omega}'-\ell_{\Omega}-2\sum_{d,\Lambda_{d}\neq \emptyset}\big(\#(I\cap J^{d})-\#(I\cap J^{d})^{+}\big),
\end{multline*}
and is $1$-dimensional when $k'-\ell=\delta(I)$.
Note that we have
\[\delta(I)\leq k_{\Omega}'-\ell_{\Omega}\]
with this equality holds if and only if $(I\cap J^{d})^{+}=I\cap J^{d}$ for each $1\leq d\leq r_{J}$ with $\Lambda_{d}\neq\emptyset$.
Hence, we deduce from (\ref{equ: Omega E1}) that $E_{1,\Omega}^{-\ell,k}=0$ whenever $k-\ell>k_{\Omega}-\ell_{\Omega}$ (with $k'=k-\#v$). In particular, we deduce that
\[d_{r,\Omega}^{-\ell,k}=0\]
whenever $k-\ell\geq k_{\Omega}-\ell_{\Omega}$. Now that we also clearly have
\[d_{r,\Omega}^{-\ell_{\Omega}-r,k_{\Omega}+(r-1)}=0\]
for each $r\geq 1$ (with $\ell_{\Omega}+r>\ell_{\Omega}=\#J$), we conclude (\ref{equ: Omega infinite page}).

We prove \ref{it: Omega seq lift 2}.\\
It is clear from (\ref{equ: Omega E1}) that we have
\[E_{1,\Omega}^{-\ell_{\Omega},k_{\Omega}}=v\otimes_E H^{k_{\Omega}'}(L_{J}/H_{I_0}Z_{J},\fl_{J}\cap\fh_{\Omega,\dagger},1_{L_{J}/H_{I_0}Z_{J}}),\]
and thus the map (\ref{equ: Omega seq lift E1}) factors through
\begin{equation}\label{equ: Omega relative to absolute}
v\otimes_E H^{k_{\Omega}'}(L_{J}/H_{I_0}Z_{J},\fl_{J}\cap\fh_{\Omega,\dagger},1_{L_{J}/H_{I_0}Z_{J}})
\rightarrow v\otimes_E H^{k_{\Omega}'}(L_{J}/Z_{J},1_{L_{J}/Z_{J}}).
\end{equation}
This together with the following isomorphism between graded $E$-algebras
\begin{equation}\label{equ: Omega Levi coh}
H^{\bullet}(L_{J}/Z_{J},1_{L_{J}/Z_{J}})\cong H^{\bullet}(\fl_{J}/\fz_{J},1_{\fl_{J}/\fz_{J}})
\end{equation}
as well as (\ref{equ: Omega Lie sm}), (\ref{equ: Omega Lie decomposition}) and the discussion below it implies that (\ref{equ: Omega seq lift E1}) factors through
\begin{multline*}
E_{1,\Omega}^{-\ell_{\Omega},k_{\Omega}}
\cong v\otimes_E \bigotimes_{d,\Lambda_{d}\neq \emptyset}H^{2\#J^{d}+1}(\fh_{\Omega,d},\fh_{\Omega,d,\dagger},1_{\fh_{\Omega,d}})\\
\cong v\otimes_E \bigotimes_{d,\Lambda_{d}\neq \emptyset}P^{2\#J^{d}+1}(\fh_{\Omega,d})
\subseteq v\otimes_E \bigotimes_{d,\Lambda_{d}\neq \emptyset}H^{2\#J^{d}+1}(\fh_{\Omega,d},1_{\fh_{\Omega,d}})\\
\subseteq v\otimes_E H^{k_{\Omega}'}(\fl_{J}/\fz_{J},1_{\fl_{J}/\fz_{J}})
\subseteq v\otimes_E H^{k_{\Omega}'}(L_{J}/Z_{J},1_{L_{J}/Z_{J}})
\end{multline*}
with the image being clearly $x_{\Omega}$. This finishes the proof of \ref{it: Omega seq lift 2}.
\end{proof}

\begin{lem}\label{lem: bottom deg differential vanishing}
Let $0\leq \ell\leq \#I_1$ and $k\geq 0$.
Assume that $k-\ell\leq\#I_1+1$. Then we have
\begin{equation}\label{equ: bottom deg infty page}
E_{2,\emptyset,I_1}^{-\ell,k}=E_{\infty,\emptyset,I_1}^{-\ell,k}.
\end{equation}
\end{lem}
\begin{proof}
Recall from \ref{it: existence atom 1} of Lemma~\ref{lem: existence of atom} that $E_{2,\emptyset,I_1}^{-\ell,k}=0$ when $k<\ell+\#I_1$, so we only need to prove (\ref{equ: bottom deg infty page}) for $k\in\{\ell+\#I_1,\ell+\#I_1+1\}$.
If $k=\ell+\#I_1$, we have $E_{2,\emptyset,I_1}^{-\ell-r,k+(r-1)}=0$ by the discussion above and thus $d_{2,\emptyset,I_1}^{-\ell-r,k+(r-1)}=0$ for each $r\geq 2$. We thus conclude (\ref{equ: bottom deg infty page}) from \ref{it: Omega seq lift 1} and \ref{it: Omega seq lift 2} of Lemma~\ref{lem: Omega seq lift} by applying Lemma~\ref{lem: abstract E2 degenerate} with $\Psi=\Psi_{I_0,I_1,\ell}$ and the quasi map $\cT_{\Sigma_{x}}^{\bullet,\bullet}\dashrightarrow \cT_{\Sigma}^{\bullet,\bullet}$ being $\cT_{\Omega}^{\bullet,\bullet}\rightarrow \cT_{I_0,I_1}^{\bullet,\bullet}$.
Now we prove (\ref{equ: bottom deg infty page}) when $k=\ell+\#I_1+1$. For each $r\geq 2$, the following fact (with $k+(r-1)=\ell+r+\#I_1$)
\[E_{2,\emptyset,I_1}^{-\ell-r,k+(r-1)}=E_{\infty,\emptyset,I_1}^{-\ell-r,k+(r-1)}\]
gives $d_{r,\emptyset,I_1}^{-\ell-r,k+(r-1)}=0$. Hence, we can again conclude (\ref{equ: bottom deg infty page}) in this case from \ref{it: Omega seq lift 1} and \ref{it: Omega seq lift 2} of Lemma~\ref{lem: Omega seq lift} by applying Lemma~\ref{lem: abstract E2 degenerate} with $\Psi=\Psi_{I_0,I_1,\ell}'$ and the quasi map $\cT_{\Sigma_{x}}^{\bullet,\bullet}\dashrightarrow \cT_{\Sigma}^{\bullet,\bullet}$ being $\cT_{\Omega}^{\bullet,\bullet}\rightarrow \cT_{I_0,I_1}^{\bullet,\bullet}$.
\end{proof}

\begin{rem}\label{rem: general degeneracy}
It is not impossible that the assumption $k-\ell\leq\#I_1+1$ in Lemma~\ref{lem: bottom deg differential vanishing} might be removed with further inputs from the normalized cup product maps to be introduced in \S~\ref{subsec: cup Tits double}.
\end{rem}

\begin{prop}\label{prop: bottom deg degeneracy}
Let $I_0\subseteq I_1\subseteq \Delta$ and $h=\#I_1-2\#I_0$.
Then we have
\begin{equation}\label{equ: bottom deg degeneracy}
E_{\infty,I_0,I_1}^{-\ell,\ell+h}=E_{2,I_0,I_1}^{-\ell,\ell+h}
\end{equation}
for each $\#I_0\leq \ell\leq \#I_1$.
\end{prop}
\begin{proof}
Recall from \ref{it: existence atom 1} of Lemma~\ref{lem: existence of atom} that $E_{2,\emptyset,I_1\setminus I_0}^{-\ell'+\#I_0,k'}=0$ for each $(\ell',k')$ satisfying $k'-\ell'<h$.
Recall from Lemma~\ref{lem: bottom deg differential vanishing} that $d_{2,\emptyset,I_1\setminus I_0}^{-\ell+\#I_0,k}=0$ for each $(-\ell,k)$ with $k-\ell=h$. We also recall from Proposition~\ref{prop: Tits bottom shift} (and its proof) that the map
\[E_{2,I_0,I_1}^{-\ell,k}\rightarrow E_{2,I_0,I_1,\flat}^{-\ell,k}\cong E_{2,\emptyset,I_1\setminus I_0}^{-\ell+\#I_0,k}\]
is injective for each $(-\ell,k)$ with $k-\ell\leq h+1$.
Consequently, we conclude (\ref{equ: bottom deg degeneracy}) by applying Lemma~\ref{lem: abstract E2 injection} to the map $\cT_{I_0,I_1}^{\bullet,\bullet}\rightarrow \cT_{I_0,I_1,\flat}^{\bullet,\bullet}\cong\cT_{\emptyset,I_1\setminus I_0}^{\bullet+\#I_0,\bullet}$.
\end{proof}

We recall from the discussion below (\ref{equ: double Levi to natural}) for the notation in the following result. For each $(-\ell,k)$, comparing \emph{loc.cit.} with 
\[E_{2,I_0,I_1}^{-\ell,k}\cong \bigoplus_{v\subseteq\mathbf{Log}_{\emptyset}}E_{2,I_0,I_1,v}^{-\ell,k}\] 
from (\ref{equ: E1 v decomposition}), we have
\begin{equation}\label{equ: tau J E2}
\tau^{J}(E_{2,I_0,I_1}^{-\ell,k})\cong\bigoplus_{v\supseteq\mathbf{Log}_{\Delta\setminus J}^{\infty}}E_{2,I_0,I_1,v}^{-\ell,k}
\end{equation}
as well as
\[\mathrm{gr}^{J}(E_{2,I_0,I_1}^{-\ell,k})\cong\bigoplus_{v\cap\mathbf{Log}_{\infty}^{\infty}=\mathbf{Log}_{\Delta\setminus J}^{\infty}}E_{2,I_0,I_1,v}^{-\ell,k}\]
for each $J\subseteq\Delta$, and similarly we have
\begin{equation}\label{equ: tau i E2}
\tau^{i}(E_{2,I_0,I_1}^{-\ell,k})\cong\bigoplus_{\#v\cap\mathbf{Log}_{\emptyset}^{\infty}\geq i}E_{2,I_0,I_1,v}^{-\ell,k}
\end{equation}
as well as
\[\mathrm{gr}^{i}(E_{2,I_0,I_1}^{-\ell,k})\cong\bigoplus_{\#v\cap\mathbf{Log}_{\emptyset}^{\infty}= i}E_{2,I_0,I_1,v}^{-\ell,k}\]
for each $i\geq 0$.
\begin{prop}\label{prop: bottom deg graded}
Let $I_0\subseteq I_1\subseteq \Delta$, $h=\#I_1-2\#I_0$. We have the following results.
\begin{enumerate}[label=(\roman*)]
\item \label{it: bottom deg graded 1} We have $\tau^{J}(E_{2,I_0,I_1}^{-\ell,\ell+h})=\tau^{J}(E_{\infty,I_0,I_1}^{-\ell,\ell+h})$ and $\mathrm{gr}^{J}(E_{2,I_0,I_1}^{-\ell,\ell+h})=\mathrm{gr}^{J}(E_{\infty,I_0,I_1}^{-\ell,\ell+h})$ for each $\#I_0\leq \ell\leq \#I_1$ and $J\subseteq \Delta$. The quasi-map $\tau^{J}(\cT_{I_0,I_1}^{\bullet,\bullet})\dashrightarrow \cT_{I_0,I_1}^{\bullet,\bullet}$ induces an embedding
    \[H^{h}(\mathrm{Tot}(\tau^{J}(\cT_{I_0,I_1}^{\bullet,\bullet})))\hookrightarrow H^{h}(\mathrm{Tot}(\cT_{I_0,I_1}^{\bullet,\bullet}))\]
    for each $J\subseteq\Delta$, so that
    \begin{equation}\label{equ: bottom deg J filtration}
    \{H^{h}(\mathrm{Tot}(\tau^{J}(\cT_{I_0,I_1}^{\bullet,\bullet})))\}_{J\subseteq\Delta}
    \end{equation}
    forms a separated exhaustive decreasing filtration on $H^{h}(\mathrm{Tot}(\cT_{I_0,I_1}^{\bullet,\bullet}))$ with graded piece \[H^{h}(\mathrm{Tot}(\mathrm{gr}^{J}(\cT_{I_0,I_1}^{\bullet,\bullet})))=0\] 
    when $J\not\subseteq I_1\setminus I_0$, and
    \[H^{h}(\mathrm{Tot}(\mathrm{gr}^{J}(\cT_{I_0,I_1}^{\bullet,\bullet})))\cong H^{h-\#J}(\mathrm{Tot}(\mathrm{CE}_{I_0,I_1\setminus J}^{\bullet,\bullet}))\otimes_E\wedge^{\#J}\Hom(Z_{\Delta\setminus J}^{\dagger},E)\]
    for each $J\subseteq I_1\setminus I_0$.
\item \label{it: bottom deg graded 2} We have $\tau^{i}(E_{2,I_0,I_1}^{-\ell,\ell+h})=\tau^{i}(E_{\infty,I_0,I_1}^{-\ell,\ell+h})$ and $\mathrm{gr}^{i}(E_{2,I_0,I_1}^{-\ell,\ell+h})=\mathrm{gr}^{i}(E_{\infty,I_0,I_1}^{-\ell,\ell+h})$ for each $\#I_0\leq \ell\leq \#I_1$ and $i\geq 0$. The quasi-map $\tau^{i}(\cT_{I_0,I_1}^{\bullet,\bullet})\dashrightarrow \cT_{I_0,I_1}^{\bullet,\bullet}$ induces an embedding
    \[H^{h}(\mathrm{Tot}(\tau^{i}(\cT_{I_0,I_1}^{\bullet,\bullet})))\hookrightarrow H^{h}(\mathrm{Tot}(\cT_{I_0,I_1}^{\bullet,\bullet}))\]
    for each $i\geq 0$, so that
    \begin{equation}\label{equ: bottom deg i filtration}
    \{H^{h}(\mathrm{Tot}(\tau^{i}(\cT_{I_0,I_1}^{\bullet,\bullet})))\}_{i\geq 0}
    \end{equation}
    forms a separated exhaustive decreasing filtration on $H^{h}(\mathrm{Tot}(\cT_{I_0,I_1}^{\bullet,\bullet}))$ with graded piece
    \begin{multline*}
    H^{h}(\mathrm{Tot}(\mathrm{gr}^{i}(\cT_{I_0,I_1}^{\bullet,\bullet})))\cong \bigoplus_{J\subseteq I_1\setminus I_0, \#J=i} H^{h}(\mathrm{Tot}(\mathrm{gr}^{J}(\cT_{I_0,I_1}^{\bullet,\bullet})))\\
    \cong \bigoplus_{J\subseteq I_1\setminus I_0, \#J=i}H^{h-i}(\mathrm{Tot}(\mathrm{CE}_{I_0,I_1\setminus J}^{\bullet,\bullet}))\otimes_E\wedge^{i}\Hom(Z_{\Delta\setminus J}^{\dagger},E)
    \end{multline*}
    for each $i\geq 0$.
\end{enumerate}
\end{prop}
\begin{proof}
Recall from (\ref{equ: double Levi to natural}) that the map $\cT_{I_0,I_1}^{\bullet,\bullet}\rightarrow \cT_{I_0,I_1,\natural}^{\bullet,\bullet}$ induces an isomorphism on the first page. Now that we have $E_{2,I_0,I_1}^{-\ell,k}=0$ when $k-\ell<h$ by \ref{it: existence atom 1} of Lemma~\ref{lem: existence of atom} and $E_{2,I_0,I_1}^{-\ell,\ell+h}=E_{\infty,I_0,I_1}^{-\ell,\ell+h}$ for each $\#I_0\leq\ell\leq \#I_1$ by Proposition~\ref{prop: bottom deg degeneracy}, \ref{it: bottom deg graded 1} follows by applying Lemma~\ref{lem: abstract E2 injection} to
\[\tau^{J}(\cT_{I_0,I_1}^{\bullet,\bullet})\rightarrow \cT_{I_0,I_1,\natural}^{\bullet,\bullet},\]
and then applying Lemma~\ref{lem: abstract E2 seq} to
\[\sum_{J'\supsetneq J}\tau^{J'}(\cT_{I_0,I_1}^{\bullet,\bullet})\rightarrow \tau^{J}(\cT_{I_0,I_1}^{\bullet,\bullet})\rightarrow \mathrm{gr}^{J}(\cT_{I_0,I_1}^{\bullet,\bullet}),\]
using (\ref{equ: tau J E2}).
Similarly, \ref{it: bottom deg graded 2} follows by applying Lemma~\ref{lem: abstract E2 seq} to
\[\tau^{i}(\cT_{I_0,I_1}^{\bullet,\bullet})\rightarrow \cT_{I_0,I_1}^{\bullet,\bullet},\]
and then applying Lemma~\ref{lem: abstract E2 seq} to
\[\tau^{i+1}(\cT_{I_0,I_1}^{\bullet,\bullet})\rightarrow \tau^{i}(\cT_{I_0,I_1}^{\bullet,\bullet})\rightarrow \mathrm{gr}^{i}(\cT_{I_0,I_1}^{\bullet,\bullet}),\]
using (\ref{equ: tau i E2}).
Note that the description of $H^{h}(\mathrm{Tot}(\mathrm{gr}^{J}(\cT_{I_0,I_1}^{\bullet,\bullet})))$ and $H^{h}(\mathrm{Tot}(\mathrm{gr}^{i}(\cT_{I_0,I_1}^{\bullet,\bullet})))$ uses (\ref{equ: grade J Z to CE}) and (\ref{equ: grade J direct sum}).
\end{proof}

Let $J\subseteq\Delta$ with $m\defeq\#J$.
We write $\Delta_{m}$ for the (standard) positive simple root associated with $G_{m+1}$ and $\cT_{J_0,J_1,m}^{\bullet,\bullet}$ for the Tits double complex associated with $G_{m+1}$ and $J_0\subseteq J_1\subseteq\Delta_{m}$.
The isomorphism $H_{J}\cong G_{m+1}$ induces an identification (see the discussion around (\ref{equ: general double complex normalized}) and Proposition~\ref{prop: Tits bottom shift} with $I_1=\Delta$ and $I_0'=I_0=\Delta\setminus J$ in \emph{loc.cit.} for the notation)
\[\cN\cT_{\Delta\setminus J,\Delta,\flat\flat}^{\bullet,\bullet}=\cN\cT_{\emptyset,\Delta_{m},m}^{\bullet,\bullet}.\]
Thanks to (\ref{equ: Tits double complex normalize}) and \ref{it: Tits bottom shift 1} of Proposition~\ref{prop: Tits bottom shift}, we have the following maps between double complex
\begin{equation}\label{equ: Tits double Levi restriction}
\cT_{\Delta\setminus J,\Delta}^{\bullet,\bullet}\rightarrow \cT_{\Delta\setminus J,\Delta,\flat\flat}^{\bullet,\bullet}\buildrel\sim\over\longrightarrow \cT_{\emptyset,\Delta_{m},m}^{\bullet+\#\Delta\setminus J,\bullet}
\end{equation}
which induce the following isomorphisms on the second page
\[E_{2,\Delta\setminus J,\Delta}^{-\ell,\ell+2\#J-\#\Delta}\buildrel\sim\over\longrightarrow E_{2,\Delta\setminus J,\Delta,\flat\flat}^{-\ell,\ell+2\#J-\#\Delta}\buildrel\sim\over\longrightarrow E_{2,\emptyset,\Delta_{m},m}^{-\ell+\#\Delta\setminus J,\ell+2\#J-\#\Delta}\]
for each $\#\Delta\setminus J\leq \ell\leq \#\Delta$. Further combined with Proposition~\ref{prop: bottom deg degeneracy}, we know that the maps (\ref{equ: Tits double Levi restriction}) induce the following isomorphism
\begin{equation}\label{equ: Tits interval reduction}
\mathbf{E}_{J}\defeq H^{2\#J-\#\Delta}(\mathrm{Tot}(\cT_{\Delta\setminus J,\Delta}^{\bullet,\bullet}))\buildrel\sim\over\longrightarrow \mathbf{E}_{m}\defeq H^{m}(\mathrm{Tot}(\cT_{\emptyset,\Delta_{m},m}^{\bullet,\bullet})).
\end{equation}
\subsection{Extensions between locally analytic generalized Steinberg}\label{subsec: Ext St}
By combining results from \S \ref{subsec: std seq}, \S \ref{subsec: comb subcomplex}, \S \ref{subsec: bottom deg} and \S \ref{subsec: bottom deg degeneracy}, we construct a canonical filtration on the bottom degree $\mathrm{Ext}_{G}^{\bullet}$ between two Tits complex and give explicit bases for graded pieces of this filtration (see Theorem~\ref{thm: ext complex}).

Recall from Proposition~\ref{prop: Ext complex std seq} that we have a canonical isomorphism
\begin{equation}\label{equ: Ext between two Tits}
\mathrm{Ext}_{G}^{k}(\mathbf{C}_{I_0,I_1},\mathbf{C}_{I_2,I_3})\cong \mathrm{Ext}_{G}^{k-\#I_1}(1_{G},\mathbf{C}_{J_0,J_1})
\end{equation}
for each $k\in\Z$, with $J_0\defeq I_2\cup (I_1\setminus I_0)$ and $J_1\defeq I_1\cap I_3$.
Note that $J_0\subseteq J_1$ if and only if $I_2\subseteq I_1\subseteq I_0\cup I_3$.
The spectral sequence $E_{\bullet,J_0,J_1}^{\bullet,\bullet}$ induces a canonical filtration on (\ref{equ: Ext between two Tits}) with graded pieces given by
\begin{equation}\label{equ: ext complex graded}
\mathrm{gr}^{-\ell}(\mathrm{Ext}_{G}^{k}(\mathbf{C}_{I_0,I_1},\mathbf{C}_{I_2,I_3}))\cong E_{\infty,J_0,J_1}^{-\ell,\ell+k-\#I_1}.
\end{equation}

This together with Proposition~\ref{prop: atom basis}, Lemma~\ref{lem: bottom deg differential vanishing} and Proposition~\ref{prop: bottom deg degeneracy} gives the following result.
\begin{thm}\label{thm: ext complex}
Let $I_0\subseteq I_1\subseteq \Delta$ and $I_2\subseteq I_3\subseteq \Delta$ be two pairs and let $k\in\Z$. We have the following results.
\begin{enumerate}[label=(\roman*)]
\item \label{it: ext complex 1} If (\ref{equ: Ext between two Tits}) is non-zero, then we have $I_2\subseteq I_1\subseteq I_0\cup I_3$ and $\#J_1-2\#J_0+\#I_1\leq k\leq n^2-n+\#I_1$.
\item \label{it: ext complex 2} If $k=\#J_1-2\#J_0+\#I_1$, then we have
\[E_{\infty,J_0,J_1}^{-\ell,\ell+k-\#I_1}=E_{2,J_0,J_1}^{-\ell,\ell+k-\#I_1}\]
which admits a basis of the form $\{\overline{x}_{\Omega}\mid \Omega\in \Psi_{J_0,J_1,\ell}\}$ for each $\#J_0\leq \ell\leq \#J_1$.
\item \label{it: ext complex 3} If $J_0=\emptyset$ and $k=\#J_1+\#I_1+1$, then we have
\[E_{\infty,\emptyset,J_1}^{-\ell,\ell+k-\#I_1}=E_{2,\emptyset,J_1}^{-\ell,\ell+k-\#I_1}\]
which admits a basis of the form $\{\overline{x}_{\Omega}\mid \Omega\in \Psi_{\emptyset,J_1,\ell}'\}$ for each $\#J_0\leq \ell\leq \#J_1$.
\end{enumerate}
\end{thm}

Recall that $\mathbf{C}_{I,\Delta}\cong V_{I}^{\rm{an}}[\#I]$ for each $I\subseteq \Delta$. Taking $I_1=I_3=\Delta$, we have $J_0=I_2\cup (\Delta \setminus I_0)\subseteq \Delta =J_1$ and
\[\mathrm{Ext}_{G}^{\bullet}(\mathbf{C}_{I_0,\Delta},\mathbf{C}_{I_2,\Delta})\cong \mathrm{Ext}_{G}^{\bullet+\#I_2-\#I_0}(V_{I_0}^{\rm{an}},V_{I_2}^{\rm{an}}).\]
Note that we have
\begin{multline*}
(\#J_1-2\#J_0+\#I_1)+\#I_2-\#I_0=2\#\Delta-2\#I_2\cup (\Delta \setminus I_0)+\#I_2-\#I_0\\
=2\#\Delta+2\#I_2\cap (\Delta \setminus I_0)-2\#I_2-2\#\Delta \setminus I_0+\#I_2-\#I_0=2\#I_2\setminus I_0+\#I_0-\#I_2\\
=\#I_0+\#I_2-2\#I_0\cap I_2=\#I_0\setminus I_2+\#I_2\setminus I_0=d(I_0,I_2).
\end{multline*}
We write
\begin{equation}\label{equ: bottom deg Ext St}
\mathbf{E}_{I_0,I_2}\defeq \mathrm{Ext}_{G}^{d(I_0,I_2)}(V_{I_0}^{\rm{an}},V_{I_2}^{\rm{an}})
\end{equation}
and
\[\mathbf{E}_{n-1}'\defeq \mathrm{Ext}_{G}^{n}(1_{G},V_{\emptyset}^{\rm{an}})\]
for short.
We obtain the following results by taking $I_1=I_3=\Delta$ in Theorem~\ref{thm: ext complex}.
\begin{cor}\label{cor: Ext between St}
Let $I_0,I_2\subseteq \Delta$ and $k\geq 0$. We have the following results.
\begin{enumerate}[label=(\roman*)]
\item \label{it: Ext between St 1} If $\mathrm{Ext}_{G}^k(V_{I_0}^{\rm{an}},V_{I_2}^{\rm{an}})\neq 0$, then $d(I_0,I_2)\leq k\leq n^2-1+\#I_2-\#I_0$.
\item \label{it: Ext between St 2} The space $\mathbf{E}_{I_0,I_2}$ admits a canonical decreasing filtration $\mathrm{Fil}^{\bullet}(\mathbf{E}_{I_0,I_2})$ with $\mathrm{gr}^{-\ell}(\mathbf{E}_{I_0,I_2})\cong E_{2,J_0,\Delta}^{-\ell,\ell+\#\Delta-2\#J_0}$ admitting a basis of the form $\{\overline{x}_{\Omega}\mid \Omega\in \Psi_{J_0,\Delta,\ell}\}$ for each $\#J_0\leq \ell\leq \#\Delta$.
\item \label{it: Ext between St 3} The space $\mathbf{E}_{n-1}'$ admits a canonical decreasing filtration
    $\mathrm{Fil}^{\bullet}(\mathbf{E}_{n-1}')$ with graded piece $\mathrm{gr}^{-\ell}(\mathbf{E}_{n-1}')\cong E_{2,\emptyset,\Delta}^{-\ell,\ell+\#\Delta+1}$ admitting a basis of the form $\{\overline{x}_{\Omega}\mid \Omega\in \Psi_{\emptyset,\Delta,\ell}'\}$ for each $0\leq \ell\leq \#\Delta$.
\end{enumerate}
\end{cor}
\subsection{Twisted $(I_0,I_1)$-atomic tuples}\label{subsec: twisted tuples}
We continue to fix a pair $I_0\subseteq I_1\subseteq \Delta$ as well as a choice of $\#I_0\leq \ell\leq \#I_1$ and $k\geq 0$. We also fix a choice of $v\subseteq \mathbf{Log}_{\emptyset}$ as well as a sequence of integers $0=r^0<r^1<\cdots<r^{r_{I_v\cap I_1}}=n-\ell$.
If $k=\ell+\#I_1-2\#I_0$ (resp.~$k=\ell+\#I_1-2\#I_0+1$), we recall from \ref{it: atom basis 1} of Proposition~\ref{prop: atom basis} (resp.~\ref{it: atom basis 2} of Proposition~\ref{prop: atom basis}) that $\{\overline{x}_{\Omega}\}_{\Omega\in \Psi_{I_0,I_1,\ell}}$ (resp.~$\{\overline{x}_{\Omega}\}_{\Omega\in \Psi_{I_0,I_1,\ell}'}$) forms a basis of $E_{2,I_0,I_1}^{-\ell,k}$.
If $k=\ell+\#I_1-2\#I_0$ (resp.~$k=\ell+\#I_1-2\#I_0+1$), we consider a choice of $\Omega\in \Psi_{I_0,I_1,\ell}$ (resp.~$\Omega\in \Psi_{I_0,I_1,\ell}'$) with the fixed $v$.
Then for each choice of $1\leq s_0\leq r_{I_v\cap I_1}$ that satisfies Condition~\ref{cond: bottom deg s} (with respect to the choice of $\Omega$ above) and each $r^{s_0-1}+1\leq d_0\leq r^{s_0}$, we construct below the \emph{$(s_0,d_0)$-twist} of $\Omega$ which is another set of tuples denoted $\Omega^{s_0,d_0}$.
We will prove in Proposition~\ref{prop: atom twist common class} that $x_{\Omega^{s_0,d_0}}\defeq \sum_{\Theta\in\Omega^{s_0,d_0}}\varepsilon(\Theta)x_{\Theta}$ satisfies
\[x_{\Omega}-x_{\Omega^{s_0,d_0}}\in d_{1,I_0,I_1}^{-\ell-1,k}(E_{1,I_0,I_1}^{-\ell-1,k}).\]
In particular, we will conclude that $x_{\Omega^{s_0,d_0}}\in \mathrm{ker}(d_{1,I_0,I_1}^{-\ell,k})$ and that $x_{\Omega}$ and $x_{\Omega^{s_0,d_0}}$ share the same image under $\mathrm{ker}(d_{1,I_0,I_1}^{-\ell,k})\twoheadrightarrow E_{2,I_0,I_1}^{-\ell,k}$.
The existence of such $(s_0,d_0)$-twist $\Omega^{s_0,d_0}$ of $\Omega$ is crucial to the proof of Lemma~\ref{lem: separated cup injection} and Lemma~\ref{lem: separated cup disconnected isom} and thus of Theorem~\ref{thm: general cup}.

\begin{defn}\label{def: atom twist}
Let $\Omega$ be an equivalence class with the fixed $v$ and bidegree $(-\ell,k)$ which moreover satisfies $r_{\Omega}^s=r^s$ for each $1\leq s\leq r_{I_v\cap I_1}$. Let $\Theta=(v,I,\un{k},\un{\Lambda})\in\Omega$ be a tuple. For each $1\leq s_0\leq r_{I_v\cap I_1}$ and $r^{s_0-1}+1\leq d_0\leq r^{s_0}$, we say that $\Theta$ is \emph{$(s_0,d_0)$-atomic} if the following conditions hold
\begin{enumerate}[label=(\roman*)]
\item \label{it: atom twist 1} $I^{d_0}\subseteq I_0$ and $\Lambda_{d_0}=\emptyset$;
\item \label{it: atom twist 2} $\Lambda_d=\{2n_d-1\}$ for each $r_{v,I_1,I}^{s_0-1}+1\leq d\leq r_{v,I_1,I}^{s_0}$ satisfying $d\neq d_0$ and $I^d\cap I_0=\emptyset$;
\item \label{it: atom twist 3} $I^d\cap I_0$ is an interval for each $r_{v,I_1,I}^{s_0-1}+1\leq d\leq r_{v,I_1,I}^{s_0}$ satisfying $I^d\cap I_0\neq \emptyset$;
\item \label{it: atom twist 4} $I^{d,+}=\emptyset$, $I^{d,-}\neq \emptyset$ and $\Lambda_d=\{2\#I^{d,-}+1\}$ for each $d_0< d\leq r_{v,I_1,I}^{s_0}$ satisfying $I^d\cap I_0\neq \emptyset$;
\item \label{it: atom twist 5} $I^{d,-}=\emptyset$, $I^{d,+}\neq \emptyset$ and $\Lambda_d=\{2\#I^{d,+}+1\}$ for each $r_{v,I_1,I}^{s_0-1}+1\leq d< d_0$ satisfying $I^d\cap I_0\neq \emptyset$;
\item \label{it: atom twist 6} for each $1\leq s\leq r_{I_v\cap I_1}$ with $s\neq s_0$, $\Theta$ satisfies either Condition~\ref{cond: bottom deg s} or Condition~\ref{cond: bottom deg s bis bis}.
\end{enumerate}
\end{defn}

In the following, we fix our choice of $(I_0,I_1)$-atomic equivalence class $\Omega$, $1\leq s_0\leq r_{I_v\cap I_1}$ and $r^{s_0-1}+1\leq d_0\leq r^{s_0}$ and assume that $\Omega$ and $s_0$ satisfy
\begin{cond}\label{cond: twist atom set up}
\begin{enumerate}[label=(\roman*)]
\item \label{it: twist atom set up 1} $\Omega$ has the fixed $v$ and bidegree $(-\ell,k)$, and thus satisfies $r_{\Omega}^s=r^s$ for each $1\leq s\leq r_{I_v\cap I_1}$;
\item \label{it: twist atom set up 2} the unique maximal element $\Theta$ in $\Omega$ (see Lemma~\ref{lem: unique maximal tuple}) is $(s_0,1)$-atomic.
\end{enumerate}
\end{cond}
We use the shortened notation
\[i_{\Theta,d}'\defeq \max\{i\in I\setminus I_0\mid i<i_{\Theta,d}\}\]
for each $2\leq d\leq r_I$.
We fix a choice of $1\leq s_0\leq r_{I_v\cap I_1}$.
For each $r^{s_0-1}+2\leq d_0\leq r^{s_0}$, we set
\begin{equation}\label{equ: maximal atom twist}
\Theta^{s_0,d_0}=(v,I^{s_0,d_0},\un{k}^{s_0,d_0},\un{\Lambda}^{s_0,d_0})\defeq p_{i_{\Theta,d_0}'}^-\circ p_{i_{\Theta,d_0-1}}^+\circ\cdots \circ p_{i_{\Theta,r^{s_0-1}+2}'}^-\circ p_{i_{\Theta,r^{s_0-1}+1}}^+(\Theta).
\end{equation}
We use the convention $\Theta^{s_0,r^{s_0-1}+1}\defeq \Theta$ for convenience.

\begin{lem}\label{lem: explicit twist atom}
Let $1\leq s_0\leq r_{I_v\cap I_1}$. We have the following results.
\begin{enumerate}[label=(\roman*)]
\item \label{it: explicit twist atom 1} The tuple $\Theta^{s_0,d_0}$ is $(s_0,d_0)$-atomic.
\item \label{it: explicit twist atom 2} For each tuple $\Theta'$ which is $(s_0,d_0)$-atomic, there exists a $(s_0,1)$-atomic tuple $\Theta$ such that $\Theta'=\Theta^{s_0,d_0}$.
\end{enumerate}
\end{lem}
\begin{proof}
We first prove \ref{it: explicit twist atom 1}.\\

Since $\Theta$ satisfies either Condition~\ref{cond: bottom deg s} or Condition~\ref{cond: bottom deg s bis bis} for each $1\leq s\leq r_{I_v\cap I_1}$ satisfying $s\neq s_0$, so is $\Theta^{s_0,d_0}$ as we have $I^d=(I^{s_0,d_0})^d$ and $\Lambda_d=\Lambda^{s_0,d_0}_d$ for each $1\leq s\leq r_{I_v\cap I_1}$ and $r^{s-1}+1\leq d\leq r^{s}$ with $s\neq s_0$.

We prove that $\Theta^{s_0,d_0}$ is $(s_0,d_0)$-atomic by an increasing induction on $r^{s_0-1}+1\leq d_0\leq r^{s_0}$. The case $d_0=r^{s_0-1}+1$ is obvious. Assume now that $r^{s_0-1}+2\leq d_0\leq r^{s_0}$ and that $\Theta^{s_0,d_0-1}$ is $(s_0,d_0-1)$-atomic. Recall from (\ref{equ: maximal atom twist}) that
\begin{equation}\label{equ: explicit twist atom induction}
\Theta^{s_0,d_0}=p_{i_{\Theta,d_0}'}^-(p_{i_{\Theta,d_0-1}}^+(\Theta^{s_0,d_0-1})).
\end{equation}
We consider an arbitrary $r^{s_0-1}+1\leq d\leq r^{s_0}$. If either $d<d_0-1$ or $d>d_0$, then we have $(I^{s_0,d_0})^d=(I^{s_0,d_0-1})^d$ and $\Lambda^{s_0,d_0}_d=\Lambda^{s_0,d_0-1}_d$, and the condition to be checked on $(I^{s_0,d_0})^d$ and $\Lambda^{s_0,d_0}_d$ follows directly from that of $(I^{s_0,d_0-1})^d$ and $\Lambda^{s_0,d_0-1}_d$.
For each $d\in \{d_0-1,d_0\}$, we note that $(I^{s_0,d_0})^d\cap I_0=(I^{s_0,d_0-1})^d\cap I_0$ is an interval if non-empty by our inductive assumption, which gives \ref{it: atom twist 3} of Definition~\ref{def: atom twist} when $d\in \{d_0-1,d_0\}$.
It follows from (\ref{equ: explicit twist atom induction}) that $i_{\Theta^{s_0,d_0},d_0}=i_{\Theta^{s_0,d_0-1},d_0}=i_{\Theta,d_0}$, which together with the definition of $i_{\Theta,d_0}'$ forces $(I^{s_0,d_0})^{d_0}=[i_{\Theta,d_0}'+1,i_{\Theta,d_0}-1]\subseteq I_0$ (with the convention $[i_{\Theta,d_0}'+1,i_{\Theta,d_0}-1]=\emptyset$ if $i_{\Theta,d_0}=i_{\Theta,d_0}'+1$), and thus \ref{it: atom twist 1} of Definition~\ref{def: atom twist} follows.
If $(I^{s_0,d_0})^{d_0-1}\cap I_0=(I^{s_0,d_0-1})^{d_0-1}\cap I_0$ is non-empty and thus an interval, then we deduce from $(I^{s_0,d_0-1})^{d_0-1}\subseteq I_0$ (by our inductive assumption on $\Theta^{s_0,d_0-1}$) and $i_{\Theta^{s_0,d_0},d_0-1}=i_{\Theta^{s_0,d_0-1},d_0-1}$ (by (\ref{equ: explicit twist atom induction})) that $(I^{s_0,d_0})^{d_0-1,-}=\emptyset$.
Our inductive assumption gives $\Lambda^{s_0,d_0-1}_{d_0-1}=\emptyset$ and $\Lambda^{s_0,d_0-1}_{d_0}=\{2(i_{\Theta,d_0}'-i_{\Theta,d_0-1})+1\}$, which gives $\Lambda^{s_0,d_0}_{d_0}=\emptyset$ and $\Lambda^{s_0,d_0}_{d_0-1}=\{2(i_{\Theta,d_0}'-i_{\Theta,d_0-1})+1\}$ by (\ref{equ: explicit twist atom induction}). This together with
\begin{equation}\label{equ: atom twist block relation}
\left\{\begin{array}{ccc}
(I^{s_0,d_0})^{d_0-1} & = & (I^{s_0,d_0-1})^{d_0-1}\sqcup [i_{\Theta,d_0-1},i_{\Theta,d_0}'-1]\\ 
(I^{s_0,d_0-1})^{d_0} & = & (I^{s_0,d_0})^{d_0}\sqcup[i_{\Theta,d_0-1}+1,i_{\Theta,d_0}']
\end{array}\right.
\end{equation}
(both equalities again by (\ref{equ: explicit twist atom induction})) gives the rest of conditions in \ref{it: atom twist 2}, \ref{it: atom twist 4} and \ref{it: atom twist 5} of Definition~\ref{def: atom twist} for $\Theta^{s_0,d_0}$ to be $(s_0,d_0)$-atomic.

Now we prove \ref{it: explicit twist atom 2} by an increasing induction on $d_0$.\\
Given $\Theta'$ that satisfy our assumption, we can check that there exists a unique tuple $\Theta''$ that satisfies
\begin{equation}\label{equ: reverse atom twist}
\Theta'=p_{i_{\Theta',d_0-1}}^-(p_{i_{\Theta',d_0-2}''}^+(\Theta''))
\end{equation}
with $i_{\Theta',d_0-2}''\defeq \min\{i\notin I_0\mid i>i_{\Theta',d_0-2}\}$. In fact, the argument in the proof of \ref{it: explicit twist atom 1} (upon replacing $\Theta^{s_0,d_0-1}$ and $\Theta^{s_0,d_0}$ there with $\Theta''$ and $\Theta'$ here) shows that $\Theta''$ is $(s_0,d_0-1)$-atomic. By our inductive assumption there exists a unique $(s_0,1)$-atomic tuple $\Theta$ such that $\Theta''=\Theta^{s_0,d_0-1}$, and (\ref{equ: reverse atom twist}) becomes (\ref{equ: explicit twist atom induction}) for this unique $\Theta$, with $\Theta'=\Theta^{s_0,d_0}$.
\end{proof}

\begin{defn}\label{def: twist partial order}
Let $\Theta'=(v,I',\un{k}',\un{\Lambda}')$ and $\Theta''=(v,I'',\un{k}'',\un{\Lambda}'')$ be two tuples satisfying $\un{\Lambda}'=\un{\Lambda}''$ and $r_{v,I_1,I'}^{s}=r_{v,I_1,I''}^{s}=r^{s}$ for each $1\leq s\leq r_{I_v\cap I_1}$.
Given $1\leq s_0\leq r_{I_v\cap I_1}$ and $r^{s_0-1}+1\leq d_0\leq r^{s_0}$, we have the following notions.
\begin{enumerate}[label=(\roman*)]
\item \label{it: twist partial order 1} Let $1\leq d\leq r_{I'}=r_{I''}$. We say that $\Theta''$ is a \emph{$(s_0,d_0)$-improvement} of $\Theta'$ if exactly one of the following holds
\begin{itemize}
\item $\Theta''$ is an improvement of $\Theta'$ with level $d$ (see \ref{it: tuple partial order 1} of Definition~\ref{def: partial order}) for some $1\leq d\leq r^{s_0-1}-1$ or some $r^{s_0}+1\leq d\leq r_{I'}=r_{I''}$;
\item there exists $d\geq d_0$ such that $i_{\Theta'',d}\in (I')^d\setminus I_0$, $i_{\Theta',d}>i_{\Theta'',d}\geq i_{\Theta^{s_0,d_0},d_0}$ and $\Theta''=p_{i_{\Theta',d}}^+(p_{i_{\Theta'',d}}^-(\Theta'))$;
\item there exists $d\leq d_0-1$ such that $i_{\Theta',d}\in (I'')^d\setminus I_0$, $i_{\Theta',d}<i_{\Theta'',d}\leq i_{\Theta^{s_0,d_0},d_0-1}$ and $\Theta'=p_{i_{\Theta'',d}}^+(p_{i_{\Theta',d}}^-(\Theta''))$.
\end{itemize}
\item \label{it: twist partial order 2} We say that $\Theta'$ is \emph{$(s_0,d_0)$-smaller} than $\Theta''$ if there exist $t\geq 1$ and tuples $\Theta_{t'}'$ for $0\leq t'\leq t$ such that
$\Theta_{0}'=\Theta'$, $\Theta_{t}'=\Theta''$ and that $\Theta_{t'}'$ is a $(s_0,d_0)$-improvement of $\Theta_{t'-1}'$ for each $1\leq t'\leq t$. Note that $\Theta'$ being $(s_0,d_0)$-smaller than $\Theta''$ implies $i_{\Theta',d}\leq i_{\Theta'',d}\leq i_{\Theta^{s_0,d_0},d_0-1}$ for each $r^{s_0-1}+1\leq d\leq d_0-1$ and $i_{\Theta',d}\geq i_{\Theta'',d}\geq i_{\Theta^{s_0,d_0},d_0}$ for each $d_0\leq d\leq r^{s_0}$.
\end{enumerate}
It is clear that our notion above depends on the integers $i_{\Theta^{s_0,d_0},d_0-1}$ and $i_{\Theta^{s_0,d_0},d_0}$ which further depends on our fixed choice of $\Omega$ (or equivalently its maximal element $\Theta$).
\end{defn}

We define $\Omega^{s_0,d_0}$ as the set of all tuples that are $(s_0,d_0)$-smaller than $\Theta^{s_0,d_0}$.
We set
\begin{equation}\label{equ: twist atom basis}
x_{\Omega^{s_0,d_0}}\defeq \sum_{\Theta'\in\Omega^{s_0,d_0}}\varepsilon(\Theta')x_{\Theta'} \in E_{1,I_0,I_1}^{-\ell,k}.
\end{equation}

Let $\Theta'=(v,I',\un{k}',\un{\Lambda}')$ be a tuple.
\begin{cond}\label{cond: twist smaller}
We consider the following conditions on $\Theta'$
\begin{enumerate}[label=(\roman*)]
\item \label{it: twist smaller 1} $r_{v,I_1,I'}^s=r^s$ for each $1\leq s\leq r_{I_v\cap I_1}$ and $(I')^d=(I^{s_0,d_0})^d=I^d$ for each $1\leq d\leq r_{I'}$ satisfying either $d\leq r^{s_0-1}$ or $d\geq r^{s_0}+1$;
\item \label{it: twist smaller 2} $\un{\Lambda}'=\un{\Lambda}^{s_0,d_0}$;
\item \label{it: twist smaller 3} $i_{\Theta',d_0-1}\leq i_{\Theta^{s_0,d_0},d_0-1}$ and $i_{\Theta',d_0}\geq i_{\Theta^{s_0,d_0},d_0}$.
\end{enumerate}
\end{cond}
It follows from Definition~\ref{def: twist partial order} that each $\Theta'=(v,I',\un{k}',\un{\Lambda}')\in \Omega^{s_0,d_0}$ satisfies Condition~\ref{cond: twist smaller}.
The converse is also true, namely we have the following results.
\begin{lem}\label{lem: twist atom criterion}
Let $\Theta'=(v,I',\un{k}',\un{\Lambda}')$ be a tuple that satisfies Condition~\ref{cond: twist smaller}. Then we have $\Theta'\in\Omega^{s_0,d_0}$.
\end{lem}
\begin{proof}
By Definition~\ref{def: twist partial order} we know that any $(s_0,d_0)$-improvement of $\Theta'$ still satisfies Condition~\ref{cond: twist smaller}. Hence, it is harmless to assume moreover that $\Theta'$ is \emph{$(s_0,d_0)$-maximal} (namely $\Theta'$ does not admit any $(s_0,d_0)$-improvement) and prove that $\Theta'=\Theta^{s_0,d_0}$.
Assume on the contrary that $\Theta'\neq \Theta^{s_0,d_0}$, then by \ref{it: twist smaller 1} of Condition~\ref{cond: twist smaller} there exists $r^{s_0-1}+1\leq d\leq r^{s_0}-1$ such that $i_{\Theta',d}\neq i_{\Theta^{s_0,d_0},d}$.
We have the following cases.
\begin{itemize}
\item If $d=d_0-1$, then we have $i_{\Theta',d}<i_{\Theta^{s_0,d_0},d}$ and $i_{\Theta',d+1}\geq i_{\Theta^{s_0,d_0},d+1}$ with $\Lambda'_{d_0}=\Lambda^{s_0,d_0}_{d_0}=\emptyset$, and thus there exists a unique tuple $\Theta''$ such that $\Theta'=p_{i_{\Theta^{s_0,d_0},d}}^+(p_{i_{\Theta',d}}^-(\Theta''))$. In other words, $\Theta''$ is a $(s_0,d_0)$-improvement of $\Theta'$, contradicting our assumption on $\Theta'$.
\item If $d=d_0$, then we have $i_{\Theta',d}>i_{\Theta^{s_0,d_0},d}$ and $i_{\Theta',d-1}\leq i_{\Theta^{s_0,d_0},d-1}$ with $\Lambda'_{d_0}=\Lambda^{s_0,d_0}_{d_0}=\emptyset$, and thus $\Theta''=p_{i_{\Theta',d}}^+(p_{i_{\Theta^{s_0,d_0},d}}^-(\Theta'))$ is a $(s_0,d_0)$-improvement of $\Theta'$, contradicting our assumption on $\Theta'$.
\item If $d<d_0-1$ and $i_{\Theta',d_0-1}=i_{\Theta^{s_0,d_0},d_0-1}$, then we may choose the maximal possible such $d$ so that we have $i_{\Theta',d}\neq i_{\Theta^{s_0,d_0},d}$ and $i_{\Theta',d+1}=i_{\Theta^{s_0,d_0},d+1}$.
    If $i_{\Theta',d}< i_{\Theta^{s_0,d_0},d}$, then by $\Lambda'_{d+1}=\Lambda^{s_0,d_0}_{d+1}$ there exists a unique tuple $\Theta''$ such that $\Theta'=p_{i_{\Theta^{s_0,d_0},d}}^+(p_{i_{\Theta',d}}^-(\Theta''))$, namely $\Theta''$ is a $(s_0,d_0)$-improvement of $\Theta'$ which contradicts our assumption on $\Theta'$.
    If $i_{\Theta',d}> i_{\Theta^{s_0,d_0},d}$, then we deduce from $i_{\Theta',d+1}=i_{\Theta^{s_0,d_0},d+1}$, $i_{\Theta',d}\notin I_0$ as well as Definition~\ref{def: atom twist} (more precisely its \ref{it: atom twist 2}, \ref{it: atom twist 3} and \ref{it: atom twist 5}) that $n'_{d+1}=i_{\Theta',d+1}-i_{\Theta',d}<\frac{m_{d+1}+1}{2}$ with $\Lambda'_{d+1}=\Lambda^{s_0,d_0}_{d+1}=\{m_{d+1}\}$, another contradiction.
\item If $d>d_0$ and $i_{\Theta',d_0}=i_{\Theta^{s_0,d_0},d_0}$, then we may choose the minimal possible such $d$ so that we have $i_{\Theta',d}\neq i_{\Theta^{s_0,d_0},d}$ and $i_{\Theta',d-1}=i_{\Theta^{s_0,d_0},d-1}$.
    If $i_{\Theta',d}> i_{\Theta^{s_0,d_0},d}$, then by $\Lambda'_{d}=\Lambda^{s_0,d_0}_{d}$ we know that $\Theta''=p_{i_{\Theta',d}}^+(p_{i_{\Theta^{s_0,d_0},d}}^-(\Theta'))$ is a $(s_0,d_0)$-improvement of $\Theta'$ which contradicts our assumption on $\Theta'$.
    If $i_{\Theta',d}< i_{\Theta^{s_0,d_0},d}$, then we deduce from $i_{\Theta',d-1}=i_{\Theta^{s_0,d_0},d-1}$, $i_{\Theta',d}\notin I_0$ as well as Definition~\ref{def: atom twist} (more precisely its \ref{it: atom twist 2}, \ref{it: atom twist 3} and \ref{it: atom twist 4}) that $n'_{d}=i_{\Theta',d}-i_{\Theta',d-1}<\frac{m_{d}+1}{2}$ with $\Lambda'_{d}=\Lambda^{s_0,d_0}_{d}=\{m_{d}\}$, another contradiction.
\end{itemize}
We finally conclude that $\Theta'=\Theta^{s_0,d_0}$ and finish the proof.
\end{proof}

For each $1\leq s_0\leq r_{I_v\cap I_1}$ and $r^{s_0-1}+2\leq d_0\leq r^{s_0}$, we define \[\Theta^{s_0,d_0,+}=(v,I^{s_0,d_0,+},\un{k}^{s_0,d_0,+},\un{\Lambda}^{s_0,d_0,+})\defeq p_{i_{\Theta^{s_0,d_0},d_0-1}}^+(\Theta^{s_0,d_0})\]
and then $\Omega^{s_0,d_0,+}$ as the set of all tuples of the form $p_{i_{\Theta',d_0-1}}^+(\Theta')$ for some $\Theta'\in\Omega^{s_0,d_0}$.
\begin{lem}\label{lem: atom twist +}
Let $s_0$, $d_0$ and $\Omega^{s_0,d_0,+}$ as above.
Then a tuple $\Theta''=(v,I'',\un{k}'',\un{\Lambda}'')$ with bidegree $(-\ell-1,k)$ belongs to $\Omega^{s_0,d_0,+}$ if and only if it satisfies
\begin{enumerate}[label=(\roman*)]
\item \label{it: twist smaller + 1} $r_{v,I_1,I''}^s=r^s$ for each $1\leq s\leq s_0-1$ and $r_{v,I_1,I''}^s=r^s-1$ for each $s_0\leq s\leq r_{I_v\cap I_1}$;
\item \label{it: twist smaller + 2} $(I'')^d=(I^{s_0,d_0,+})^d$ for each $1\leq d\leq r_{I''}$ satisfying either $d\leq r_{v,I_1,I''}^{s_0-1}$ or $d\geq r_{v,I_1,I''}^{s_0}+1$;
\item \label{it: twist smaller + 3} $\un{\Lambda}''=\un{\Lambda}^{s_0,d_0,+}$;
\item \label{it: twist smaller + 4} $i_{\Theta'',d_0-2}\leq i_{\Theta^{s_0,d_0,+},d_0-2}$ and $i_{\Theta'',d_0-1}\geq i_{\Theta^{s_0,d_0,+},d_0-1}$.
\end{enumerate}
\end{lem}
\begin{proof}
We first prove the `only if' direction.
Assume that $\Theta''\in \Omega^{s_0,d_0,+}$, namely $\Theta''=p_{i_{\Theta',d_0-1}}^+(\Theta')$ for some $\Theta'=(v,I',\un{k}',\un{\Lambda}')\in\Omega^{s_0,d_0}$. Since $\Theta'$ is $(s_0,d_0)$-smaller than $\Theta^{s_0,d_0}$, we have $\un{\Lambda}'=\un{\Lambda}^{s_0,d_0}$, $(I')^d=(I^{s_0,d_0})^d$ for each $1\leq d\leq r_{I'}$ satisfying either $d\leq r^{s_0-1}$ or $d\geq r^{s_0}+1$, $i_{\Theta',d}\leq i_{\Theta^{s_0,d_0},d}$ for each $r^{s_0-1}+1\leq d\leq d_0-1$, and $i_{\Theta',d}\geq i_{\Theta^{s_0,d_0},d}$ for each $d_0\leq d\leq r^{s_0}$. This together with $\Theta''=p_{i_{\Theta',d_0-1}}^+(\Theta')$ and $\Theta^{s_0,d_0,+}=p_{i_{\Theta^{s_0,d_0},d_0-1}}^+(\Theta^{s_0,d_0})$ gives \ref{it: twist smaller + 1}, \ref{it: twist smaller + 2}, \ref{it: twist smaller + 3} and \ref{it: twist smaller + 4}.

Now we prove the `if' direction.
By \ref{it: twist smaller + 4} we know that $\Theta'\defeq p_{i_{\Theta^{s_0,d_0},d_0-1}}^-(\Theta'')$ is defined and satisfy $\Theta''=p_{i_{\Theta',d_0-1}}^+(\Theta')$. It suffices to notice that \ref{it: twist smaller + 1}, \ref{it: twist smaller + 2}, \ref{it: twist smaller + 3} and \ref{it: twist smaller + 4} together ensures that $\Theta'$ satisfies Condition~\ref{cond: twist smaller}, and thus $\Theta'\in\Omega^{s_0,d_0}$ by Lemma~\ref{lem: twist atom criterion}.
\end{proof}

For each $1\leq s_0\leq r_{I_v\cap I_1}$ and $r^{s_0-1}+2\leq d_0\leq r^{s_0}$, we set
\[x_{\Omega^{s_0,d_0,+}}\defeq \sum_{\Theta''\in \Omega^{s_0,d_0,+}}\varepsilon(\Theta'')x_{\Theta''}.\]
\begin{lem}\label{lem: atom twist transfer}
For each $1\leq s_0\leq r_{I_v\cap I_1}$ and $r^{s_0-1}+2\leq d_0\leq r^{s_0}$, we have
\begin{equation}\label{equ: atom twist transfer}
d_{1,I_0,I_1}^{-\ell-1,k}(x_{\Omega^{s_0,d_0,+}})=(-1)^{d_0-1}(x_{\Omega^{s_0,d_0}}+x_{\Omega^{s_0,d_0-1}}).
\end{equation}
\end{lem}
\begin{proof}
We consider an arbitrary tuple $\Theta'=(v,I',\un{k}',\un{\Lambda}')$ such that there exists $i'\in (I_v\cap I_1)\setminus I'$ that satisfies $\Theta''=(v,I'',\un{k}'',\un{\Lambda}'')\defeq p_{i'}^+(\Theta')\in \Omega^{s_0,d_0,+}$. It suffices to determine
\begin{equation}\label{equ: atom twist differential}
c_{\Theta'}(d_{1,I_0,I_1}^{-\ell-1,k}(x_{\Omega^{s_0,d_0,+}}))
\end{equation}
for each such tuple $\Theta'$. We write $1\leq d_1\leq r_{I''}$ for the unique integer such that $i'\in (I'')^{d_1}$. In particular, we have $i'=i_{\Theta',d_1}$, $\Lambda''_{d_1}=\Lambda'_{d_1}\sqcup\Lambda'_{d_1+1}$, $\Lambda''_d=\Lambda'_d$ for each $d<d_1$, and $\Lambda''_d=\Lambda'_{d+1}$ for each $d>d_1$.
Note that all tuples in $\Omega^{s_0,d_0,+}$ share the common $\un{\Lambda}''$.
If either $d_1\leq r^{s_0-1}$ or $d_1\geq r^{s_0}+1$, then similar argument as in \textbf{Case $1$} and \textbf{Case $2$} of the proof of Lemma~\ref{lem: atomic subcomplex} (which uses Lemma~\ref{lem: tuple sum of two}) shows that (\ref{equ: atom twist differential}) is zero.

We assume from now that $r^{s_0-1}+1\leq d_1\leq r^{s_0}$, which forces $\#\Lambda''_{d_1}=1$ by the definition of $\Omega^{s_0,d_0}$. Hence, we have either $\Lambda'_{d_1}=\emptyset$ with $\Lambda'_{d_1+1}=\Lambda''_{d_1}$, or $\Lambda'_{d_1}=\Lambda''_{d_1}$ with $\Lambda'_{d_1+1}=\emptyset$.
We also note that $\Lambda'_d=\Lambda''_d$ for each $r^{s_0-1}+1\leq d<d_1$ and $\Lambda'_d=\Lambda''_{d-1}$ for each $d_1+1<d\leq r^{s_0}$, with $\#\Lambda'_d=1$ in both cases.

If $d_1\neq d_0-1$, then among all choices of $1\leq d\leq r_{I'}$ satisfying $i_{\Theta',d}\in (I_v\cap I_1)\setminus I'$, we can check that $\Theta'''=(v,I''',\un{k}''',\un{\Lambda}''')\defeq p_{i_{\Theta',d}}^+(\Theta')$ satisfies $\un{\Lambda}'''=\un{\Lambda}''$ if and only if either $d=d_1$, or $d=d_1-1$ with $\Lambda'_{d_1}=\emptyset$, or $d=d_1+1$ with $\Lambda'_{d_1+1}=\emptyset$.
Note that $\Omega''\in \Omega^{s_0,d_0,+}$ satisfies \ref{it: twist smaller + 1}, \ref{it: twist smaller + 2}, \ref{it: twist smaller + 3} and \ref{it: twist smaller + 4} of Lemma~\ref{lem: atom twist +}.
If $\Lambda'_{d_1}=\emptyset$ (resp.~$\Lambda'_{d_1+1}=\emptyset$), then we can that check that $p_{i_{\Theta',d_1-1}}^+(\Theta')$ (resp.~$p_{i_{\Theta',d_1+1}}^+(\Theta')$) still satisfies \ref{it: twist smaller + 1}, \ref{it: twist smaller + 2}, \ref{it: twist smaller + 3} and \ref{it: twist smaller + 4} of \emph{loc.cit.}, which forces $p_{i_{\Theta',d_1-1}}^+(\Theta')\in\Omega^{s_0,d_0,+}$ (resp.~$p_{i_{\Theta',d_1+1}}^+(\Theta')\in\Omega^{s_0,d_0,+}$) by \emph{loc.cit.}.
Consequently, we deduce from Lemma~\ref{lem: tuple sum of two} (by taking $\Theta$ \emph{loc.cit.} to be $\Theta'$ and $d$ in \emph{loc.cit.} to be the unique element in $\{d_1,d_1+1\}$ that satisfies $\Lambda'_d=\emptyset$) that (\ref{equ: atom twist differential}) is zero.

We assume in the rest of the proof that $d_1=d_0-1$.
Recall from Lemma~\ref{lem: atom twist +} that $i_{\Theta'',d_0-2}\leq i_{\Theta^{s_0,d_0,+},d_0-2}$ and $i_{\Theta'',d_0-1}\geq i_{\Theta^{s_0,d_0,+},d_0-1}$.
Recall from (\ref{equ: explicit twist atom induction}) and the definition of $\Theta^{s_0,d_0,+}$ that
\[p_{i_{\Theta^{s_0,d_0-1},d_0-1}}^+(\Theta^{s_0,d_0-1})=\Theta^{s_0,d_0,+}=p_{i_{\Theta^{s_0,d_0},d_0-1}}^+(\Theta^{s_0,d_0})\]
and thus we have
\[i_{\Theta^{s_0,d_0-1},d_0-2}=i_{\Theta^{s_0,d_0,+},d_0-2}=i_{\Theta^{s_0,d_0},d_0-2}\]
and
\[i_{\Theta^{s_0,d_0-1},d_0}=i_{\Theta^{s_0,d_0,+},d_0-1}=i_{\Theta^{s_0,d_0,+},d_0}.\]
Recall from \ref{it: atom twist 1} of Definition~\ref{def: atom twist} that $(I^{s_0,d_0-1})^{d_0-1}\subseteq I_0$ and $(I^{s_0,d_0})^{d_0}\subseteq I_0$, which together with (\ref{equ: atom twist block relation}) implies that $i_{\Theta^{s_0,d_0-1},d_0-1}<i_{\Theta^{s_0,d_0},d_0-1}$ and
\[i_{\Theta^{s_0,d_0-1},d_0-1}=\min\{i\in (I_v\cap I_1)\setminus I_0\mid i>i_{\Theta^{s_0,d_0-1},d_0-2}\}\]
and
\[i_{\Theta^{s_0,d_0},d_0-1}=\max\{i\in (I_v\cap I_1)\setminus I_0\mid i<i_{\Theta^{s_0,d_0},d_0}\}.\]
Now that $\Theta''=p_{i_{\Theta',d_0-1}}^+(\Theta')$ with $i_{\Theta',d_0-1}\in (I_v\cap I_1)\setminus I'\subseteq (I_v\cap I_1)\setminus I_0$, we have either $i_{\Theta',d_0-1}\leq i_{\Theta^{s_0,d_0-1},d_0-2}$, or $i_{\Theta^{s_0,d_0-1},d_0-1}\leq i_{\Theta',d_0-1}\leq i_{\Theta^{s_0,d_0},d_0-1}$, or $i_{\Theta',d_0-1}\geq i_{\Theta^{s_0,d_0},d_0}$.
Since $\Lambda'_{d_0-1}\sqcup\Lambda'_{d_0}=\Lambda''_{d_0-1}=\Lambda^{s_0,d_0,+}_{d_0-1}$ consist of exactly one element, we have either $\Lambda'_{d_0-1}=\emptyset$ with $\Lambda'_{d_0}=\Lambda^{s_0,d_0,+}_{d_0-1}$, or $\Lambda'_{d_0}=\emptyset$ with $\Lambda'_{d_0-1}=\Lambda^{s_0,d_0,+}_{d_0-1}$.
We write $d_2\in\{d_0-1,d_0\}$ for the unique element such that $\Lambda'_{d_2}=\emptyset$.
Then among $1\leq d\leq r_{I'}$ that satisfies $i_{\Theta',d}\in (I_v\cap I_1)\setminus I'$, we can check that $\Theta'''=(v,I''',\un{k}''',\un{\Lambda}''')\defeq p_{i_{\Theta',d}}^+(\Theta')$ satisfies $\un{\Lambda}'''=\un{\Lambda}''=\un{\Lambda}^{s_0,d_0,+}$ and $r_{v,I_1,I'''}^s=r_{v,I_1,I''}^s$ for each $1\leq s\leq r_{I_v\cap I_1}$ if and only if $d\in \{d_2-1,d_2\}$. Note that we have
\begin{equation}\label{equ: atom twist differential prime}
c_{\Theta'}(d_{1,I_0,I_1}^{-\ell-1,k}(x_{\Omega^{s_0,d_0,+}}))=\sum_{\Theta'''}\varepsilon(\Theta''')c_{\Theta'}(d_{1,I_0,I_1}^{-\ell-1,k}(x_{\Theta'''}))
\end{equation}
with $\Theta'''$ running through all tuples of the form $p_{i_{\Theta',d}}^+(\Theta')$ (for $d\in \{d_2-1,d_2\}$) that belongs to $\Omega^{s_0,d_0,+}$.
We have the following possibilities.
\begin{itemize}
\item If $d_2=d_0-1$ and $i_{\Theta',d_0-1}\leq i_{\Theta^{s_0,d_0-1},d_0-2}=i_{\Theta^{s_0,d_0,+},d_0-2}$, then we have \[\left\{\begin{array}{c}
    i_{p_{i_{\Theta',d_2-1}}^+(\Theta'),d_0-2}=i_{\Theta',d_0-1}\leq i_{\Theta^{s_0,d_0,+},d_0-2}\\
    i_{p_{i_{\Theta',d_2-1}}^+(\Theta'),d_0-1}=i_{\Theta',d_0}=i_{\Theta'',d_0-1}\geq i_{\Theta^{s_0,d_0,+},d_0-2}
    \end{array}\right.\] 
    and thus $p_{i_{\Theta',d_2-1}}^+(\Theta')\in\Omega^{s_0,d_0,+}$ by Lemma~\ref{lem: atom twist +}. This together with $p_{i_{\Theta',d_2}}^+(\Theta')=\Theta''\in\Omega^{s_0,d_0,+}$, (\ref{equ: atom twist differential prime}) and Lemma~\ref{lem: tuple sum of two} (with $\Theta$ and $d$ in \emph{loc.cit.} taken to be $\Theta'$ and $d_2$ here) gives \[c_{\Theta'}(d_{1,I_0,I_1}^{-\ell-1,k}(x_{\Omega^{s_0,d_0,+}}))=0.\]
\item If $d_2=d_0-1$ and $i_{\Theta',d_0-1}\geq i_{\Theta^{s_0,d_0-1},d_0-1}$, then we have \[i_{p_{i_{\Theta',d_2-1}}^+(\Theta'),d_0-2}=i_{\Theta',d_0-1}\geq i_{\Theta^{s_0,d_0-1},d_0-1}>i_{\Theta^{s_0,d_0,+},d_0-2}\] and thus $p_{i_{\Theta',d_2-1}}^+(\Theta')\notin\Omega^{s_0,d_0,+}$ by Lemma~\ref{lem: atom twist +}. This together with $p_{i_{\Theta',d_2}}^+(\Theta')=\Theta''\in\Omega^{s_0,d_0,+}$, (\ref{equ: atom twist differential prime}) and (\ref{equ: absolute sign}) (with $d_0-1$ here replacing $d_0$ in \emph{loc.cit.}) gives
    \begin{multline}\label{equ: atom twist differential transfer 1}
    \varepsilon(\Theta')^{-1}c_{\Theta'}(d_{1,I_0,I_1}^{-\ell-1,k}(x_{\Omega^{s_0,d_0,+}}))=\varepsilon(\Theta')^{-1}\varepsilon(\Theta'')(-1)^{m(I'',i_{\Theta',d_0-1})}\\
    =(-1)^{i_{\Theta',d_0-1}}(-1)^{m(I'',i_{\Theta',d_0-1})}=(-1)^{d_0-1}.
    \end{multline}
\item If $d_2=d_0$ and $i_{\Theta',d_0-1}\geq i_{\Theta^{s_0,d_0},d_0}=i_{\Theta^{s_0,d_0,+},d_0-1}$, then we have \[\left\{\begin{array}{c}
    i_{p_{i_{\Theta',d_2}}^+(\Theta'),d_0-1}=i_{\Theta',d_0-1}\geq i_{\Theta^{s_0,d_0,+},d_0-1}\\
    i_{p_{i_{\Theta',d_2}}^+(\Theta'),d_0-2}=i_{\Theta',d_0-2}=i_{\Theta'',d_0-2}\leq i_{\Theta^{s_0,d_0,+},d_0-2}, 
    \end{array}\right.\]
    and thus $p_{i_{\Theta',d_2}}^+(\Theta')\in\Omega^{s_0,d_0,+}$ by Lemma~\ref{lem: atom twist +}. This together with $p_{i_{\Theta',d_2-1}}^+(\Theta')=\Theta''\in\Omega^{s_0,d_0,+}$, (\ref{equ: atom twist differential prime}) and Lemma~\ref{lem: tuple sum of two} (with $\Theta$ and $d$ in \emph{loc.cit.} taken to be $\Theta'$ and $d_2$ here) gives \[c_{\Theta'}(d_{1,I_0,I_1}^{-\ell-1,k}(x_{\Omega^{s_0,d_0,+}}))=0.\]
\item If $d_2=d_0$ and $i_{\Theta',d_0-1}\leq i_{\Theta^{s_0,d_0},d_0-1}$, then we have \[i_{p_{i_{\Theta',d_2}}^+(\Theta'),d_0-1}=i_{\Theta',d_0-1}\leq i_{\Theta^{s_0,d_0},d_0-1}<i_{\Theta^{s_0,d_0,+},d_0-1}\] and thus $p_{i_{\Theta',d_2}}^+(\Theta')\notin\Omega^{s_0,d_0,+}$ by Lemma~\ref{lem: atom twist +}. This together with $p_{i_{\Theta',d_2-1}}^+(\Theta')=\Theta''\in\Omega^{s_0,d_0,+}$, (\ref{equ: atom twist differential prime}) and (\ref{equ: absolute sign}) (with $d_0-1$ here replacing $d_0$ in \emph{loc.cit.}) gives
    \begin{multline}\label{equ: atom twist differential transfer 2}
    \varepsilon(\Theta')^{-1}c_{\Theta'}(d_{1,I_0,I_1}^{-\ell-1,k}(x_{\Omega^{s_0,d_0,+}}))=\varepsilon(\Theta')^{-1}\varepsilon(\Theta'')(-1)^{m(I'',i_{\Theta',d_0-1})}\\
    =(-1)^{i_{\Theta',d_0-1}}(-1)^{m(I'',i_{\Theta',d_0-1})}=(-1)^{d_0-1}.
    \end{multline}
\end{itemize}
We conclude (\ref{equ: atom twist transfer}) by combining the four cases above.
\end{proof}

Recall that we have fixed a choice of $(I_0,I_1)$-atomic equivalence class $\Omega$ (with maximal element $\Theta=(v,I,\un{k},\un{\Lambda})$) with bidegree $(-\ell,k)$ satisfying $k=\ell+\#I_1-2\#I_0$.
\begin{prop}\label{prop: atom twist common class}
Let $\Omega$ be a $(I_0,I_1)$-atomic equivalence class and $1\leq s_0\leq r_{I_v\cap I_1}$ be an integer that satisfy Condition~\ref{cond: twist atom set up}.
Then for each $r^{s_0-1}+1\leq d_0\leq r^{s_0}$, we have
\begin{equation}\label{equ: atom twist common class}
x_{\Omega^{s_0,d_0}}-(-1)^{d_0-r^{s_0-1}}x_\Omega\in d_{1,I_0,I_1}^{-\ell-1,k}(E_{1,I_0,I_1}^{-\ell-1,k})
\end{equation}
\end{prop}
\begin{proof}
We fix a choice of $1\leq s_0\leq r_{I_v\cap I_1}$.
It follows from Lemma~\ref{lem: atom twist transfer} that
\[x_{\Omega^{s_0,d_0}}+x_{\Omega^{s_0,d_0-1}}\in d_{1,I_0,I_1}^{-\ell-1,k}(E_{1,I_0,I_1}^{-\ell-1,k})\]
for each $r^{s_0-1}+2\leq d_0\leq r^{s_0}$. This together with $x_{\Omega^{s_0,1}}=x_\Omega$ and an increasing induction on $d_0$ gives (\ref{equ: atom twist common class}).
\end{proof}

We continue to fix $I_0\subseteq I_1$, $v$, the integers $r^{s}$ for $0\leq s\leq r_{I_v\cap I_1}$ as well as a $(I_0,I_1)$-atomic equivalence class $\Omega$ with bidegree $(-\ell,k)$ and its unique maximal element of the form $\Theta=(v,I,\un{k},\un{\Lambda})$.
We fix below a choice of $1\leq s_0\leq r_{I_1\cap I_v}-1$ that satisfies $i_0\defeq i_{\Theta,r_{\Omega}^{s_0}}\in\Delta\setminus I_1$, and then write $I_1'\defeq I_1\sqcup\{i_0\}$.
For each $(I_0,I_1)$-tuple $\Theta'=(v,I',\un{k}',\un{\Lambda}')$ with the fixed $v$ and bidegree $(-\ell,k)$, we write $r_{i_0}(\Theta')$ for the associated $(I_0,I_1')$-tuple with the fixed $v$ and bidegree $(-\ell,k)$
For each set $\Omega'$ of $(I_0,I_1)$-tuples, we also write
\[r_{i_0}(\Omega')\defeq \{r_{i_0}(\Theta')\mid \Theta'\in\Omega'\}\]
for the associated set of $(I_0,I_1')$-tuples, and
\[\tld{r}_{i_0}(\Omega')\defeq \{p_{i_0}^{+}(r_{i_0}(\Theta'))\mid \Theta'\in\Omega'\}.\]
\begin{lem}\label{lem: atom twist vanishing}
Let $\Omega$, $\Theta$, $s_0$ and $i_0$ be as above.
Assume that $\Theta$ satisfies Condition~\ref{cond: bottom deg s} for both $s_0$ and $s_0+1$, and that $\Theta$ satisfies either Condition~\ref{cond: bottom deg s} or Condition~\ref{cond: bottom deg s bis bis} for $s\neq s_0,s_0+1$. Then we have
\begin{equation}\label{equ: atom twist vanishing}
x_{r_{i_0}(\Omega)}\in d_{1,I_0,I_1'}^{-\ell-1,k}(E_{1,I_0,I_1'}^{-\ell-1,k}).
\end{equation}
\end{lem}
\begin{proof}
We write $d_0\defeq r_{\Omega}^{s_0}$ for short.
We consider the $(s_0,d_0)$-twist $\Omega^{s_0,d_0}$ of $\Omega$ (see above (\ref{equ: twist atom basis})) and recall from Proposition~\ref{prop: atom twist common class} that
\[x_{\Omega}-(-1)^{r_{\Omega}^{s_0}-r_{\Omega}^{s_0-1}}x_{\Omega^{s_0,r_{\Omega}^{s_0}}}\in d_{1,I_0,I_1}^{-\ell-1,k}(E_{1,I_0,I_1}^{-\ell-1,k})\]
and therefore
\[x_{r_{i_0}(\Omega)}-(-1)^{r_{\Omega}^{s_0}-r_{\Omega}^{s_0-1}}x_{r_{i_0}(\Omega^{s_0,r_{\Omega}^{s_0}})}\in d_{1,I_0,I_1'}^{-\ell-1,k}(E_{1,I_0,I_1'}^{-\ell-1,k}).\]
Hence, in order to prove (\ref{equ: atom twist vanishing}), it suffices to show that
\begin{equation}\label{equ: atom twist vanishing shift}
x_{r_{i_0}(\Omega^{s_0,r^{s_0}})}=(-1)^{r^{s_0}-1}d_{1,I_0,I_1'}^{-\ell-1,k}(x_{\tld{r}_{i_0}(\Omega^{s_0,r^{s_0}})}).
\end{equation}
We consider an arbitrary tuple $\Theta'=(v,I',\un{k}',\un{\Lambda}')$ such that there exists $i'\in (I_v\cap I_1)\setminus I'$ that satisfies $\Theta''=(v,I'',\un{k}'',\un{\Lambda}'')\defeq p_{i'}^+(\Theta')\in \tld{r}_{i_0}(\Omega^{s_0,r^{s_0}})$. It suffices to determine
\begin{equation}\label{equ: atom twist differential change v}
\varepsilon(\Theta')c_{\Theta'}(d_{1,I_0,I_1'}^{-\ell-1,k}(x_{\tld{r}_{i_0}(\Omega^{s_0,r^{s_0}})}))
\end{equation}
for each such tuple $\Theta'$.
If $i'\neq i_0$, then similar argument as in \textbf{Case $1$} and \textbf{Case $2$} of the proof of Lemma~\ref{lem: atomic subcomplex} (which uses Lemma~\ref{lem: tuple sum of two}) shows that (\ref{equ: atom twist differential change v}) is zero.
Hence, (\ref{equ: atom twist differential change v}) is non-zero if and only if $\Theta'\in r_{i_0}(\Omega^{s_0,r^{s_0}})$ in which case (\ref{equ: atom twist differential change v}) equals (see (\ref{equ: atom twist differential transfer 2}) for similar computation)
\[\varepsilon(\Theta')\varepsilon(p_{i_0}^+(\Theta'))c_{\Theta'}(d_{1,I_0,I_1}^{-\ell-1,k}(x_{p_{i_0}^+(\Theta')}))=(-1)^{r^{s_0}-1}.\]
\end{proof}
\subsection{Outer involution on bottom degree cohomology}\label{subsec: outer involution}
In this section, we study the involution on $\mathbf{E}_{\Delta}=H^{n-1}(\mathrm{Tot}(\cT_{\emptyset,\Delta}^{\bullet,\bullet}))$ induced from the outer involution of $G$ given by $A\mapsto\overline{A}\defeq w_0(A^{t})^{-1}w_0$.
We write $r\defeq \#\Delta=n-1$ for short throughout this section.

For each $I\subseteq\Delta$, we note that the involution $A\mapsto\overline{A}$ of $G$ restricts to an isomorphism $P_{I}\buildrel\sim\over\longrightarrow P_{w_0(I)}$ as well as an isomorphism $L_{I}\buildrel\sim\over\longrightarrow L_{w_0(I)}$ between locally $K$-analytic groups. This involution on $G$ induces an involution on $\mathrm{Rep}^{\rm{an}}_{\rm{adm}}(G)$ that sends $i_{I}^{\rm{an}}$ (resp.~$V_{I}^{\rm{an}}$ to $i_{w_0(I)}^{\rm{an}}$ (resp.~$V_{w_0(I)}^{\rm{an}}$, and further induces the following commutative diagram
\begin{equation}\label{equ: cochain outer involution diagram}
\xymatrix{
\Hom_{D(G)}(B(D(G),(i_{w_0(I)}^{\rm{an}})^{\vee}),1_{G}^{\vee}) \ar[r] \ar^{\wr}[d] & C^{\bullet}(P_{w_0(I)}) \ar[r] \ar^{\wr}[d]& C^{\bullet}(L_{w_0(I)}) \ar^{\wr}[d]\\ 
\Hom_{D(G)}(B(D(G),(i_{I}^{\rm{an}})^{\vee}),1_{G}^{\vee}) \ar[r] & C^{\bullet}(P_{I}) \ar[r] & C^{\bullet}(L_{I}) \\
}
\end{equation}
between complex of $E$-vector spaces, with horizontal maps being quasi-isomorphisms from (\ref{equ: PS parabolic Levi resolution}). It is clear that (\ref{equ: cochain outer involution diagram}) is functorial w.r.t the choice of $I$.
We define double complex $\overline{\cS}_{\emptyset,\Delta}^{\bullet,\bullet}$ and $\overline{\cT}_{\emptyset,\Delta}^{\bullet,\bullet}$ as $\cT_{\Sigma}^{\bullet,\bullet}$ from (\ref{equ: general double complex}) by taking $M^{\bullet}_{I}$ to be $\Hom_{D(G)}(B(D(G),(i_{w_0(I)}^{\rm{an}})^{\vee}),1_{G}^{\vee})$ and $C^{\bullet}(L_{w_0(I)})$ respectively for each $I\subseteq\Delta$. 
We omit $\emptyset,\Delta$ from the notation and write $\cS^{\bullet,\bullet}$, $\cT^{\bullet,\bullet}$, $\overline{\cS}^{\bullet,\bullet}$ and $\overline{\cT}^{\bullet,\bullet}$ instead.

Given $I\subseteq\Delta$ and $j\in I$ with $I'\defeq I\setminus\{j\}$ and $\ell\defeq \#I$, we recall $m(I,j)\in\Z$ from (\ref{equ: differential sign}) and consider the natural restriction map
\[\mathrm{Res}_{I,I'}^{\bullet}: C^{\bullet}(L_{I})\rightarrow C^{\bullet}(L_{I'})\]
induced from the inclusion $L_{I'}\subseteq L_{I}$.
By the definition of $\cT^{\bullet,\bullet}$ and $\overline{\cT}^{\bullet,\bullet}$, for each $k\geq 0$, we observe that $(-1)^{m(I,j)}\mathrm{Res}_{I,I'}^{k}$ appears as a direct summand of the differential map $\cT^{-\ell,k}\rightarrow \cT^{-\ell+1,k}$, while $(-1)^{m(w_0(I),w_0(j))}\mathrm{Res}_{I,I'}^{k}$ appears as a direct summand of the differential map $\overline{\cT}^{-\ell,k}\rightarrow \overline{\cT}^{-\ell+1,k}$.
Now that we have $m(I,j)+m(w_0(I),w_0(j))=\#I'=\ell-1$ by (\ref{equ: differential sign}), we see that $(-1)^{m(I,j)}\mathrm{Res}_{I,I'}^{k}$ and $(-1)^{m(w_0(I),w_0(j))}\mathrm{Res}_{I,I'}^{k}$ differ by a sign $(-1)^{\ell-1}$. 
Writing $\varepsilon(i)\defeq \frac{i(i+1)}{2}$ for each $i\in\Z$, we see that the scalar automorphism
\[\varepsilon(\ell-1): \bigoplus_{I\subseteq\Delta,\#I=\ell}C^{\bullet}(L_{I})\buildrel\sim\over\longrightarrow \bigoplus_{I\subseteq\Delta,\#I=\ell}C^{\bullet}(L_{I})\]
induces an isomorphism between double complex $\cT^{\bullet,\bullet}\buildrel\sim\over\longrightarrow\overline{\cT}^{\bullet,\bullet}$.
We can also construct an isomorphism  between double complex $\cS^{\bullet,\bullet}\buildrel\sim\over\longrightarrow\overline{\cS}^{\bullet,\bullet}$ by a parallel argument. Combining these with the commutative diagram (\ref{equ: cochain outer involution diagram}), we obtain the following commutative diagram of maps between double complex
\begin{equation}\label{equ: Tits outer involution diagram}
\xymatrix{
\cS^{\bullet,\bullet} \ar^{\sim}[r] \ar[d] & \overline{\cS}^{\bullet,\bullet} \ar^{\sim}[r] \ar[d] & \cS^{\bullet,\bullet} \ar[d]\\
\cT^{\bullet,\bullet} \ar^{\sim}[r] & \overline{\cT}^{\bullet,\bullet} \ar^{\sim}[r] & \cT^{\bullet,\bullet}
}
\end{equation}
with all horizontal maps being isomorphisms between double complex, and all vertical maps inducing isomorphisms at the first page.
We write $\theta$ for the involution of the double complex $\cT^{\bullet,\bullet}$ given by the composition of the bottom horizontal maps of (\ref{equ: Tits outer involution diagram}), and abuse the same notation $\theta$ for the induced involution on $H^{r}(\mathrm{Tot}(\cT^{\bullet,\bullet}))$ as well as the involution on the spectral sequence $E_{\bullet,\emptyset,\Delta}^{\bullet,\bullet}$ associated with $\cT^{\bullet,\bullet}=\cT_{\emptyset,\Delta}^{\bullet,\bullet}$.
Note that (\ref{equ: Tits outer involution diagram}) together with $\mathbf{C}_{\emptyset,\Delta}\cong\mathrm{St}_{G}^{\rm{an}}$ implies that the isomorphism (see (\ref{equ: Tits double complex total coh}) with $I_0=\emptyset$ and $I_1=\Delta$ in \emph{loc.cit.})
\[\mathrm{Ext}_{G}^{r}(1_{G},\mathrm{St}_{G}^{\rm{an}})=\mathbf{E}_{\Delta}\cong H^{r}(\mathrm{Tot}(\cT^{\bullet,\bullet}))\]
is compatible with the involution $\theta$ on BHS, with the involution on LHS induced from that of $\mathrm{Rep}^{\rm{an}}_{\rm{adm}}(G)$ as discussed before (\ref{equ: cochain outer involution diagram}).

We define $\mathbf{E}_{\Delta}^{\theta=\pm1}\subseteq \mathbf{E}_{\Delta}$ as the eigen-subspace of $\theta$ with respect to the eigenvalue $\pm 1$, and define $(E_{\bullet,\emptyset,\Delta}^{\bullet,\bullet})^{\theta=\pm1}$ similarly.
Parallel the definition of the involution $\theta$ of $\cT^{\bullet,\bullet}$, we can similarly define involutions on $\cT_{\emptyset,\Delta,\natural}^{\bullet,\bullet}$, $\tau^{i}(\cT^{\bullet,\bullet})\defeq \tau^{i}(\cT_{\emptyset,\Delta}^{\bullet,\bullet})$ and $\mathrm{gr}^{i}(\cT^{\bullet,\bullet})\defeq \mathrm{gr}^{i}(\cT_{\emptyset,\Delta}^{\bullet,\bullet})$ for each $i\geq 0$, which we again denote by $\theta$. We also define various $(-)^{\theta=\pm1}$ similarly.
The following result is clear.
\begin{lem}\label{lem: outer involution Fil}
Let $\theta$ be the involution on $\mathbf{E}_{\Delta}$ as constructed above. We have the following results.
\begin{enumerate}[label=(\roman*)]
\item \label{it: outer involution Fil 0} We have $\mathbf{E}_{\Delta}=\mathbf{E}_{\Delta}^{\theta=1}\oplus\mathbf{E}_{\Delta}^{\theta=-1}$.
\item \label{it: outer involution Fil 1} The involution $\theta$ is strict with respect to the canonical decreasing filtration $\{\mathrm{Fil}^{-\ell}(\mathbf{E}_{\Delta})\}_{\ell}$ on $\mathbf{E}_{\Delta}$. In particular, $\mathbf{E}_{\Delta}^{\theta=\pm1}$ inherits a canonical decreasing filtration
    \[\{\mathbf{E}_{\Delta}^{\theta=\pm1}\cap \mathrm{Fil}^{-\ell}(\mathbf{E}_{\Delta})=\mathrm{Fil}^{-\ell}(\mathbf{E}_{\Delta})^{\theta=\pm1}\}_{\ell}\]
    with graded pieces
    \[\mathrm{gr}^{-\ell}(\mathbf{E}_{\Delta})^{\theta=\pm1}=(E_{\infty,\emptyset,\Delta}^{-\ell,\ell+r})^{\theta=\pm1}=(E_{2,\emptyset,\Delta}^{-\ell,\ell+r})^{\theta=\pm1}.\]
\item \label{it: outer involution Fil 2} The involution $\theta$ is strict with respect to the decreasing filtration \[\{\tau^{i}(\mathbf{E}_{\Delta})=H^{r}(\mathrm{Tot}(\tau^{i}(\cT^{\bullet,\bullet})))\}_{i\geq 0}\] 
    on $\mathbf{E}_{\Delta}$. In particular, $\mathbf{E}_{\Delta}^{\theta=\pm1}$ inherits a decreasing filtration
    \[\{\mathbf{E}_{\Delta}^{\theta=\pm1}\cap \tau^{i}(\mathbf{E}_{\Delta})=\tau^{i}(\mathbf{E}_{\Delta})^{\theta=\pm1}\}_{i\geq 0}\]
    with graded pieces
    \[\mathrm{gr}^{i}(\mathbf{E}_{\Delta})^{\theta=\pm1}=H^{r}(\mathrm{Tot}(\tau^{i}(\cT^{\bullet,\bullet})))^{\theta=\pm1}.\]
\end{enumerate}
\end{lem}

Recall from (\ref{equ: atom to partition}) that we have a bijection
\begin{equation}\label{equ: atom to partition full}
\bigsqcup_{\ell=0}^{r}\Psi_{\emptyset,\Delta,\ell}\buildrel\sim\over\longrightarrow \{(S,I)\mid S\in\cS_{\Delta}, I\subseteq S\cap\Delta\}
\end{equation}
which restricts to a bijection
\[\Psi_{\emptyset,\Delta,\ell}\buildrel\sim\over\longrightarrow \{(S,I)\mid S\in\cS_{\Delta}, I\subseteq S\cap\Delta, \#S=r-\ell\}\]
for each $0\leq \ell\leq r$. 
We write $\overline{x}_{S,I}\defeq\overline{x}_{\Omega}$ for the element of $E_{2,\emptyset,\Delta}^{-\ell,\ell+r}=E_{\infty,\emptyset,\Delta}^{-\ell,\ell+r}$ for each $S\in\cS_{\Delta}$ and $I\subseteq S\cap\Delta$ with $\#S=r-\ell$, where $\Omega\in \Psi_{\emptyset,\Delta,\ell}$ is the equivalence class corresponding to $(S,I)$ under (\ref{equ: atom to partition full}).
Recall from (\ref{equ: partition involution}) that $\cS_{\Delta}$ admits an involution $S\mapsto\overline{S}$.
\begin{lem}\label{lem: outer involution E2 basis}
Let $\theta$ be the involution on $\mathbf{E}_{\Delta}$ as in Lemma~\ref{lem: outer involution Fil}.
We have the following results.
\begin{enumerate}[label=(\roman*)]
\item \label{it: outer involution E2 0} For each $S\in\cS_{\Delta}$ and $I\subseteq S\cap\Delta$ with $\#S=r-\ell$, we have
\begin{equation}\label{equ: outer involution E2 0}
\theta(\overline{x}_{S,I})=\varepsilon(S,I)\overline{x}_{\overline{S},w_0(I)}\in E_{2,\emptyset,\Delta}^{-\ell,\ell+r}
\end{equation}
for some sign $\varepsilon(S,I)\in\{1,-1\}$.
\item \label{it: outer involution E2 1} When $\#S=1$ with $S=\overline{S}$, we have
\begin{equation}\label{equ: outer involution E2 1}
\theta(\overline{x}_{S,\emptyset})=-\varepsilon(r)\overline{x}_{S,\emptyset}\in E_{2,\emptyset,\Delta}^{-r+1,2r-1}.
\end{equation}
\item \label{it: outer involution E2 2} For each $I\subseteq\Delta$, we have
\begin{equation}\label{equ: outer involution E2 2}
\theta(\overline{x}_{\Delta,I})=(-1)^{r}\varepsilon(r)\overline{x}_{\Delta,w_0(I)}\in E_{2,\emptyset,\Delta}^{0,r}.
\end{equation}
\end{enumerate}
\end{lem}
\begin{proof}
Let $S\in\cS_{\Delta}$ and $I\subseteq S\cap\Delta$ with $\#S=r-\ell$. We write $\Omega\in \Psi_{\emptyset,\Delta,\ell}$ for the equivalence class of $(\emptyset,\Delta)$-tuples corresponding to $(S,I)$ under (\ref{equ: atom to partition full}), and $\Theta=(v,I,\un{k},\un{\Lambda})\in\Omega$ for the unique (maximal) element inside (see Lemma~\ref{lem: unique atom}). 
Similarly, we write $\Omega'\in \Psi_{\emptyset,\Delta,\ell}$ for the equivalence class of $(\emptyset,\Delta)$-tuples corresponding to $(\overline{S},w_0(I))$ under (\ref{equ: atom to partition full}), and $\Theta'=(v',I',\un{k}',\un{\Lambda}')\in \Omega'$ be the unique (maximal) element inside. 
We write $\tau$ for the automorphism of $G$ given by $A\mapsto (A^{t})^{-1}$ which clearly restricts to an automorphism of $L_{w_0(I)}$.
It follows from the definition of $\theta$ that $\theta(x_{\Theta})$ equals the image of $x_{\Theta}$ under the following isomorphisms
\begin{equation}\label{equ: outer involution Levi coh}
H^{\ell+r}(L_{I},1_{L_{I}})\buildrel\sim\over\longrightarrow H^{\ell+r}(L_{I},1_{L_{I}})\buildrel\sim\over\longrightarrow H^{\ell+r}(L_{w_0(I)},1_{L_{w_0(I)}})\buildrel\sim\over\longrightarrow H^{\ell+r}(L_{w_0(I)},1_{L_{w_0(I)}}),
\end{equation}
with first isomorphism (resp.~second isomorphism, resp.~third isomorphism) given by the scalar automorphism $\varepsilon(\ell-1)$ (resp.~induced from $A\mapsto w_0Aw_0$, resp.~induced from $\tau$). 
In particular, using Lemma~\ref{lem: primitive class outer}, we have 
\[\theta(x_{\Theta})=\varepsilon(S,I)'x_{\Theta''}\] 
for some sign $\varepsilon(S,I)'\in\{1,-1\}$ and a tuple $\Theta''=(v'',I'',\un{k}'',\un{\Lambda}'')$ characterized by the following conditions.
\begin{itemize}
\item We have $\plog_{i}\in v''$ (resp.~$\val_{i}\in v''$) if and only if $\plog_{n-i}\in v$ (resp.~$\val_{n-i}\in v$).
\item We have $I''=w_0(I)$ with $r_{I''}=r_{I}$. We have $\Lambda''_{d}=\Lambda_{r_{I}+1-d}$ and $k''_{d}=k_{r_{I}+1-d}$ for each $1\leq d\leq r_{I}$.
\end{itemize}
Then we have the crucial observation that $v'=v''$ and
\[\Theta''=(\cdots(\Theta')^{1,r_{\Omega'}^{1}}\cdots)^{r_{I_{v'}},r_{\Omega'}^{r_{I_{v'}}}},\]
which together with Proposition~\ref{prop: atom twist common class} implies that
\[\varepsilon(\Theta'')x_{\Theta''}-(-1)^{r_{I_{v'}}}\varepsilon(\Theta')x_{\Theta'}\in d_{1,\emptyset,\Delta}^{-\ell-1,\ell+r}(E_{1,\emptyset,\Delta}^{-\ell-1,\ell+r}).\]
Now that $\theta(\varepsilon(\Theta)x_{\Theta})=\varepsilon(\Theta)\varepsilon(S,I)'x_{\Theta''}$ (resp.~$\varepsilon(\Theta')x_{\Theta'}$) has image $\theta(\overline{x}_{S,I})$ (resp.~image $\overline{x}_{\overline{S},w_0(I)}$) in $E_{2,\emptyset,\Delta}^{-\ell,\ell+r}$, we conclude (\ref{equ: outer involution E2 0}) with 
\begin{equation}\label{equ: outer involution total sign}
\varepsilon(S,I)\defeq \varepsilon(\Theta)\varepsilon(\Theta'')^{-1}\varepsilon(S,I)'(-1)^{r_{I_{v'}}}.
\end{equation}
This finishes the proof of \ref{equ: outer involution E2 0}.

Now we check \ref{it: outer involution E2 1} and \ref{it: outer involution E2 2} separately by making explicit the signs involved in the proof of \ref{it: outer involution E2 0}.\\
We write $\{n''_{d}\}_{1\leq d\leq r_{I''}=r_{I}}$ for the integers associated with $I''$, which satisfies $n''_{d}=n_{r_{I}+1-d}$ for each $1\leq d\leq r_{I}$.
Thanks to Lemma~\ref{lem: primitive class outer} and (\ref{equ: Levi Kunneth}), we know that $\tau$ acts by 
\[(-1)^{\#v}\prod_{1\leq d\leq r_{I''},n''_{d}\geq 2}(-1)^{n''_{d}}=(-1)^{r}\]
which is a sign independent of the choice of $\ell$.
\begin{itemize}
\item Assume that $\#S=1$ (with $S=\overline{S}$ and $I=\emptyset$). Then the first map (resp.~the second map) of (\ref{equ: outer involution Levi coh}) is the scalar automorphism $\varepsilon(r-2)$ (resp.~sends $x_{\Theta}$ to $x_{\Theta''}$) and thus $\varepsilon(S,I)'=\varepsilon(r-2)(-1)^{r}$. Now that we have $I_{v'}=\Delta$ with $r_{I_{v'}}=1$ and $\varepsilon(\Theta)\varepsilon(\Theta'')^{-1}=(-1)^{r-1}$, we conclude from (\ref{equ: outer involution total sign}) that
    \[\varepsilon(S,I)=(-1)^{r-1}\varepsilon(r-2)(-1)^{r}(-1)=-\varepsilon(r).\]
\item Assume that $S=\Delta$ (with $S=\overline{S}$). Then we have $\Theta'=\Theta''$ and thus $\varepsilon(S,I)=\varepsilon(S,I)'$ by (\ref{equ: outer involution total sign}). For each $i\in\Delta$, we define $x_{i}$ (resp.~$x_{i}'$) as the unique element of $v\cap\mathbf{Log}_{\Delta\setminus\{i\}}$ (resp.~$v'\cap\mathbf{Log}_{\Delta\setminus\{i\}}$). We note that first map of (\ref{equ: outer involution Levi coh}) is the identity map, and that the second map of (\ref{equ: outer involution Levi coh}) sends $x_{\Theta}=x_1\cup\cdots\cup x_{r}$ to
    \[x_{r}'\cup\cdots\cup x_{1}'=\varepsilon(r)x_{1}'\cup\cdots\cup x_{r}'.\]
    We thus conclude that $\varepsilon(S,I)=\varepsilon(S,I)'=(-1)^{r}\varepsilon(r)$ in this case.
\end{itemize}
\end{proof}

We set $\theta_{n-1}\defeq \varepsilon(r)\theta$ which is an involution on $\mathbf{E}_{n-1}=\mathbf{E}_{\Delta}$. It thus follows from \ref{it: outer involution E2 1} of Lemma~\ref{lem: outer involution E2 basis} (as well as Lemma~\ref{lem: outer involution Fil}) that the canonical filtration $\{\mathrm{Fil}^{-\ell}(\mathbf{E}_{n-1}^{\theta_{n-1}=-1})\}_{\ell}$ satisfies
\begin{equation}\label{equ: outer top grade negative 1}
\mathrm{gr}^{1-r}(\mathbf{E}_{n-1}^{\theta_{n-1}=-1})=\mathrm{gr}^{1-r}(\mathbf{E}_{n-1})=E\overline{x}_{\{\al_{\Delta}\},\emptyset}
\end{equation}
with $\al_{\Delta}=\sum_{\al\in\Delta}\al\in\Phi^+$ (see (\ref{equ: sum of roots})).
\subsection{Truncation maps and $(I_0,I_1)$-atomic tuples}\label{subsec: truncation atom}
We fix a choice of $I_0\subseteq I_1\subseteq I_1'\Delta$ and $J\subseteq \Delta$ with $I_1=I_1'\cap J$ throughout this section.
Upon taking $M_{I}^{\bullet}$ to be $C^{\bullet}(L_{I}/Z_{J})$ for each $I_0\subseteq I\subseteq I_1$ and zero otherwise in the definition of $\cT_{\Sigma}^{\bullet,\bullet}$ (see (\ref{equ: general double complex})), we obtain a double complex $\cT_{I_0,I_1,J}^{\bullet,\bullet}$ which comes with maps
\[\cT_{I_0,I_1,J}^{\bullet,\bullet}\rightarrow \cT_{I_0,I_1}^{\bullet,\bullet}\rightarrow \cT_{I_0,I_1'}^{\bullet,\bullet}\]
which induce maps between spectral sequences
\begin{equation}\label{equ: center char map seq}
E_{r,I_0,I_1,J}^{-\ell,k}\rightarrow E_{r,I_0,I_1}^{-\ell,k}\rightarrow E_{r,I_0,I_1'}^{-\ell,k}.
\end{equation}
We write $\psi_{r,\ell}$ for the composition of (\ref{equ: center char map seq}) with $k=\#I_1'-2\#I_0$.
In this section, We study the map $\psi_{r,\ell}$ for $r=1,2$ and each $\ell$, with the main result being Proposition~\ref{prop: J truncation map E2} and Proposition~\ref{prop: J truncation total basis}.

For each $i\in J$, we fix an isomorphism $\mathbf{Z}_{J\setminus\{i\}}/\mathbf{Z}_{J}\buildrel\sim\over\longrightarrow \bG_m$ between groups schemes over $K$, which induces an isomorphism
\[\varepsilon_{i}^{J}: Z_{J\setminus\{i\}}/Z_{J}\buildrel\sim\over\longrightarrow K^\times\]
between locally $K$-analytic groups.
We set $\val_{i}^{J}\defeq \val\circ\varepsilon_{i}^{J}$ and $\plog_{i}^{J}\defeq \plog\circ\varepsilon_{i}^{J}$ for each $i\in J$, and note that $\mathbf{Log}_{J\setminus\{i\}}^{J}\defeq \{\val_{i}^{J},\plog_{i}^{J}\}$ (resp.~$\{\val_{i}^{J}\}$) forms a basis of $\Hom(Z_{J\setminus\{i\}}/Z_{J},E)$ (resp.~of $\Hom_{\infty}(Z_{J\setminus\{i\}}/Z_{J},E)$).
Let $I\subseteq J$.
Similar to (\ref{equ: center isogeny decomposition}), the rational maps $Z_{I}\dashrightarrow Z_{J\setminus\{i\}}$ for $i\in J\setminus I$ induces the following quasi-isogeny
\[Z_{I}/Z_{J} \dashrightarrow \prod_{i\in J\setminus I}Z_{J\setminus\{i\}}/Z_{J}\]
which further induces the following isomorphism
\begin{equation}\label{equ: Levi char decomposition J}
\prod_{i\in J\setminus I}\Hom(Z_{J\setminus\{i\}}/Z_{J},E) \buildrel\sim\over\longrightarrow \Hom(Z_{I}/Z_{J},E)
\end{equation}
between $E$-vector spaces. 
In particular
\[\mathbf{Log}_{I}^{J}\defeq \{\val_{i}^{J},\plog_{i}^{J}\mid i\in J\setminus I\}\]
forms a basis of (\ref{equ: Levi char decomposition J}), and
\[\mathbf{Log}_{I}^{J,\infty}\defeq \{\val_{i}^{J}\mid i\in J\setminus I\}\subseteq \mathbf{Log}_{I}^{J}\]
forms a basis of $\Hom_{\infty}(Z_{I}/Z_{J},E)$.
Note that we have a natural injection
\begin{equation}\label{equ: center index transfer}
\mathbf{Log}_{\emptyset}^{J}\hookrightarrow \mathbf{Log}_{\emptyset}
\end{equation}
which sends $\val_{i}^{J}$ (resp.~$\plog_{i}^{J}$) to $\val_{i}$ (resp.~$\plog_{i}$) for each $i\in J$. Thus the fixed order on $\mathbf{Log}_{\emptyset}$ induces via (\ref{equ: center index transfer}) an order on $\mathbf{Log}_{\emptyset}^{J}$, through which we identify $\wedge^k\mathbf{Log}_{I}^{J}\defeq \{v\subseteq \mathbf{Log}_{I}^{J}\mid \#v=k\}$ with a basis of $\wedge^k\Hom(Z_{I}/Z_{J},E)\cong \wedge^k\Hom(Z_{I}\cap H_{J},E)$.
For each $v\subseteq \mathbf{Log}_{\emptyset}^{J}$, we define $I_{v}^{J}\subseteq J$ as the maximal subset such that $v\subseteq \mathbf{Log}_{I_{v}^{J}}^{J}$. If we abuse the notation $v$ for its image under (\ref{equ: center index transfer}), then we have $I_{v}^{J}=I_{v}\cap J$.

Upon replacing $L_{I}$ with $L_{I}/Z_{J}$ in (\ref{equ: E1 std basis}), we may define $x_{\Theta}$ for each \emph{$J$-tuple} $\Theta=(v,I,\un{k},\un{\Lambda})$ with $v\subseteq \mathbf{Log}_{\emptyset}^{J}$ and $I\subseteq I_{v}^{J}$, so that $E_{1,I_0,I_1,J}^{-\ell,k}$ admits a basis of the form $\{x_{\Theta}\}$ with $\Theta=(v,I,\un{k},\un{\Lambda})$ running through all $J$-tuples that satisfy $I_0\subseteq I\subseteq I_{v}^{J}\cap I_1$, $\#I=\ell$ and $|\un{k}|=k$.
Similar to Definition~\ref{def: partial order}, there is a natural equivalence relation on the set of all $J$-tuples, and we can define $x_{\Omega}\in E_{1,I_0,I_1,J}^{-\ell,k}$ for each such equivalence class $\Omega$ with bidegree $(-\ell,k)$.
Similar to Definition~\ref{def: atom element}, we have the notion of \emph{$(I_0,I_1)$-atomic $J$-tuples}. For each bidegree $(-\ell,k)$, we define $E_{1,I_0,I_1,J,\diamond}^{-\ell,k}$ as the $E$-span of $x_{\Omega}$ with $\Omega$ running through all equivalence classes of $(I_0,I_1)$-atomic $J$-tuples with bidegree $(-\ell,k)$.
Parallel argument as in Lemma~\ref{lem: atomic subcomplex} and Proposition~\ref{prop: quasi isom} shows that $E_{1,I_0,I_1,J,\diamond}^{\bullet,k}$ defines a subcomplex of $E_{1,I_0,I_1,J}^{\bullet,k}$ and the embedding
\[E_{1,I_0,I_1,J,\diamond}^{\bullet,k}\rightarrow E_{1,I_0,I_1,J}^{\bullet,k}\]
is a quasi-isomorphism for each $k\in\Z$.

We assume in the rest of this section that $J=\Delta\setminus\{i_0\}$ for some $i_0\in\Delta\setminus I_1$ and set $I_1'\defeq I_1\sqcup\{i_0\}$.
For each $i\in\Delta$, we consider the following commutative diagram
\begin{equation}\label{equ: J truncation center relation}
\xymatrix{
Z_{J\setminus\{j\}}/Z_{J} \ar^{\varepsilon_{j}^{J}}[d] & Z_{J\setminus\{i\}} \ar@{->>}[l] \ar@{-->}[r] & Z_{\Delta\setminus\{j\}}\times Z_{J} \ar^{\varepsilon_{j}\times \varepsilon_{i_0}}[d]\\
K^{\times} & & K^{\times}\times K^{\times}
}
\end{equation}
Since the rational map $Z_{J\setminus\{j\}}\dashrightarrow Z_{\Delta\setminus\{j\}}\times Z_{J}$ is a quasi-isogeny, there exists $m,a,b\in \Z$ such that
\[(\varepsilon_{j}^{J})^{m}=\varepsilon_{j}^{a}\varepsilon_{i_0}^{b}\]
as rational maps $Z_{J\setminus\{j\}}\rightarrow K^\times$ (under (\ref{equ: J truncation center relation})), which implies
\begin{equation}\label{equ: J truncation log}
\plog_{j}^{J}=a_{j}^{J}\plog_{j}+b_{j}^{J}\plog_{i_0}
\end{equation}
and
\begin{equation}\label{equ: J truncation val}
\val_{j}^{J}=a_{j}^{J}\val_{j}+b_{j}^{J}\val_{i_0}
\end{equation}
as elements in $\Hom(Z_{J\setminus\{j\}},E)$, with $a_{j}^{J}\defeq \frac{a}{m}$ and $b_{j}^{J}\defeq \frac{b}{m}$ being rational numbers.

\begin{lem}\label{lem: central char transfer}
For each $j\in J$, we have 
\begin{equation}\label{equ: central char transfer}
a_{j}^{J}\neq 0\neq b_{j}^{J}.
\end{equation}
\end{lem}
\begin{proof}
We fix our choice of $j\in J$ and it is harmless to assume that $j>i_0$.
It suffices to prove (\ref{equ: central char transfer}) for an explicit choice of $\varepsilon_{j}^{J}$ and $\varepsilon_{j'}$ for each $j'\in\Delta=[1,n-1]$.
Under the standard isomorphism
\[\Hom(T,E)\cong \{(x_1,\dots,x_n)\in E^{n}\mid \sum_{k=1}^{n}x_k=0\},\]
we may choose $\varepsilon_{j}$ so that
\begin{equation}\label{equ: explicit epsilon j}
\varepsilon_{j'}=(n-j',\dots,n-j',-j',\dots,-j')
\end{equation}
for each $j'\in\Delta$, with $n-j'$ appearing $j'$ times and $-j'$ appearing $n-j'$ times.
Since $j>i_0$, we may choose $\varepsilon_{j}^{J}$ so that
\begin{equation}\label{equ: explicit epsilon J i}
\varepsilon_{j}^{J}=(0,\dots,0,n-j,\dots,n-j,i_0-j,\dots,i_0-j)
\end{equation}
with $0$ appearing $i_0$ times, $n-j$ appearing $j-i_0$ times, and $i_0-j$ appearing $n-j$ times. It is thus clear from (\ref{equ: explicit epsilon j}) that
\[\varepsilon_{j}^{n-i_0}\varepsilon_{i_0}^{-(n-j)}=(0,\dots,0,n(n-j),\dots,n(n-j),n(i_0-j),\dots,n(i_0-j))=(\varepsilon_{j}^{J})^{n},\]
from which we deduce that $a_{j}^{J}=\frac{n-i_0}{n}\neq 0$ and $b_{j}^{J}=-\frac{n-j}{n}\neq 0$.
The proof is thus finished.
\end{proof}

Recall from the beginning of \S~\ref{subsec: Ext St} that we write $\Psi_{I_0,I_1',\ell}$ for the set of equivalent classes $\Omega$ of $(I_0,I_1')$-atomic tuples with bidegree $(-\ell,\ell+\#I_1'-2\#I_0)$.
We write $\Psi_{I_0,I_1,J,\ell}'$ for the set of equivalent classes $\Omega_J$ of $(I_0,I_1)$-atomic $J$-tuples with bidegree $(-\ell,\ell+\#I_1'-2\#I_0)=(-\ell,\ell+\#I_1-2\#I_0+1)$ which moreover satisfies $d_{1,I_0,I_1,J}^{-\ell,\ell+\#I_1-2\#I_0+1}(x_{\Omega_J})=0$.
It follows from Corollary~\ref{cor: bottom deg E2 basis} (and its variant for $J$-tuples) that $\{\overline{x}_{\Omega}\}_{\Omega\in \Psi_{I_0,I_1',\ell}}$ forms a basis of $E_{2,I_0,I_1'}^{-\ell,\ell+\#I_1'-2\#I_0}$, and $\{\overline{x}_{\Omega_J}\}_{\Omega_J\in \Psi_{J,I_0,I_1,\ell}'}$ forms a basis of $E_{2,I_0,I_1,J}^{-\ell,\ell+\#I_1'-2\#I_0}$.
Our goal in the rest of this section is to determine the image of $\psi_{2,\ell}(\overline{x}_{\Omega_J})$ for each $\Omega_J\in \Psi_{I_0,I_1,J,\ell}'$ and $\#I_0\leq \ell\leq \#I_1$.

Note that each $v\subseteq \mathbf{Log}_{\emptyset}^{J}$ uniquely determines a sequence $0=r_{v}^{J,0}<r_{v}^{J,1}<\cdots<r_{v}^{J,r_{I_{v}^{J}\cap I_1}}=n$ such that
\[\Delta\setminus(I_{v}^{J}\cap I_1)=\{r_{v}^{J,s}\mid 1\leq s\leq r_{I_{v}^{J}\cap I_1}-1\}.\]
Note that each $\Omega_J\in \Psi_{J,I_0,I_1,\ell}'$ uniquely determines a subset $v\subseteq \mathbf{Log}_{\emptyset}^{J}$ and thus uniquely determines an integer $1\leq s_0\leq r_{I_{v}^{J}\cap I_1}-1$ such that $i_0=r_{v}^{J,s_0}$.

For each $v\subseteq \mathbf{Log}_{\emptyset}^{J}$, we often abuse the same notation for its image in $\mathbf{Log}_{\emptyset}$ under the injection (\ref{equ: center index transfer}), and thus each $J$-tuple $\Theta_J=(v,I,\un{k},\un{\Lambda})$ with respect to $(I_0,I_1)$ uniquely determines a tuple $\Theta=(v,I,\un{k},\un{\Lambda})$ w.r.t $(I_0,I_1)$, to which we may further associate the tuple $r_{i_0}(\Theta)=(v,I,\un{k},\un{\Lambda})$ with respect to $(I_0,I_1')$.
We say that $\Theta$ is the \emph{underlying $(I_0,I_1)$-tuple} of $\Theta_{J}$, and $r_{i_0}(\Theta)$ is the \emph{underlying $(I_0,I_1')$-tuple} of $\Theta_{J}$.
Similarly, for each set $\Omega_{J}$ of $J$-tuples with the fixed $v$, we write $\Omega$ for the underlying $(I_0,I_1)$-tuples of the $J$-tuples in $\Omega_{J}$, and $r_{i_0}(\Omega)$ for the underlying $(I_0,I_1')$-tuples of the $J$-tuples in $\Omega_{J}$.

Let $v\subseteq \mathbf{Log}_{\emptyset}$ with $i_0\in I_{v}$ (namely $v\cap \mathbf{Log}_{\Delta\setminus\{i_0\}}=\emptyset$). For each $(I_0,I_1)$-tuple $\Theta=(v,I,\un{k},\un{\Lambda})$, we define a new $(I_0,I_1)$-tuple $u_{j}^{\plog}(\Theta)\defeq ((v\setminus\{\plog_{j}\})\sqcup\{\plog_{i_0}\},I,\un{k},\un{\Lambda})$ if $\plog_{j}\in v$, and define a new $(I_0,I_1)$-tuple $u_{j}^{\val}(\Theta)\defeq ((v\setminus\{\val_{j}\})\sqcup\{\val_{i_0}\},I,\un{k},\un{\Lambda})$ if $\val_{j}\in v$.
If $\Omega$ is a set of $(I_0,I_1)$-tuples with the fixed $v$ above, we similarly define $u_{j}^{\plog}(\Omega)$ if $\plog_{j}\in v$, and define $u_{j}^{\val}(\Omega)$ if $\val_{j}\in v$.
If $v\cap\mathbf{Log}_{\Delta\setminus\{j\}}=1$, then we define $u_{j}(\Theta)$ as $u_{j}^{\plog}(\Theta)$ (resp.~as $u_{j}^{\val}(\Theta)$) when $\plog_{j}\in v$ (resp.~when $\val_{j}\in v$), and we similarly define $u_{j}(\Omega)$ for each $\Omega$ which is a set of $(I_0,I_1)$-tuples with this fixed $v$.

\begin{cond}\label{cond: J atom classification}
Let $\#I_0\leq \ell\leq \#I_1$.
Upon replacing $G$ with $L_{J}/Z_{J}$ in \ref{it: existence atom 3} of Lemma~\ref{lem: existence of atom} and \ref{it: atom basis 2} of Proposition~\ref{prop: atom basis}, we see that $\Omega_J\in \Psi_{I_0,I_1,J,\ell}'$ if and only if the unique maximal $J$-tuple $\Theta_{J}^{\max}=(v,I,\un{k},\un{\Lambda})$ in $\Omega_{J}$ satisfies exactly one of the following conditions.
\begin{enumerate}[label=(\roman*)]
\item \label{it: J atom case 1} We have $\#I_{v}^{J}\cup I_1=\#J-1$, $\#v\cap\mathbf{Log}_{J\setminus\{j\}}^{J}\leq 1$ for each $j\in J$. Moreover, $\Theta_{J}^{\max}$ satisfies the $L_{J}/Z_{J}$-variant of Condition~\ref{cond: bottom deg s} for each $1\leq s\leq r_{I_{v}^{J}\cap I_1}$.
\item \label{it: J atom case 2} We have $I_{v}^{J}\cup I_1=J$ and there exists a unique $i_1\in J$ such that $\{\val_{i_1}^{J},\plog_{i_1}^{J}\}\subseteq v$ with $\#v\cap\mathbf{Log}_{J\setminus\{i\}}^{J}\leq 1$ for each $i_1\neq j\in J$. Moreover, $\Theta_{J}^{\max}$ satisfies the $L_{J}/Z_{J}$-variant of Condition~\ref{cond: bottom deg s} for each $1\leq s\leq r_{I_{v}^{J}\cap I_1}$.
\item \label{it: J atom case 3} We have $I_{v}^{J}\cup I_1=J$ and $\#v\cap\mathbf{Log}_{J\setminus\{j\}}^{J}\leq 1$ for each $j\in J$. Moreover, there exists a unique $1\leq s_1\leq r_{I_{v}^{J}\cap I_1}$ such that $\Theta_{J}^{\max}$ satisfies the $L_{J}/Z_{J}$-variant of Condition~\ref{cond: bottom deg s bis bis} for $s_1$, as well as the $L_{J}/Z_{J}$-variant of Condition~\ref{cond: bottom deg s} for each $1\leq s\leq r_{I_{v}^{J}\cap I_1}$ with $s\neq s_1$.
\end{enumerate}
\end{cond}
We define $\Psi_{I_0,I_1,J,\ell}^{\natural}$ (resp.~$\Psi_{I_0,I_1,J,\ell}^{\dagger}$, resp.~$\Psi_{I_0,I_1,J,\ell}^{\ddagger}$) as the subset of $\Psi_{I_0,I_1,J,\ell}'$ consisting of those $\Omega_{J}$ that satisfying \ref{it: J atom case 1} (resp.~\ref{it: J atom case 2}, resp.~\ref{it: J atom case 3}) of Condition~\ref{cond: J atom classification}.

\begin{prop}\label{prop: J truncation map E2}
Let $i_0\in\Delta$, $I_0\subseteq I_1\subseteq J=\Delta\setminus\{i_0\}$ and $I_1'=I_1\sqcup\{i_0\}$ as above. Let $\#I_0\leq \ell\leq \#I_1$ and $\Omega_J\in \Psi_{I_0,I_1,J,\ell}'$, with $\Theta_{J}^{\max}$ being the unique maximal $J$-tuple in $\Omega_J$, and $\Omega$ be the set of underlying $(I_0,I_1)$-tuples of the $J$-tuples in $\Omega_{J}$. Let $v\subseteq\mathbf{Log}_{\emptyset}^{J}$ be the set determined by $\Omega_{J}$ with $i\defeq \#v\cap\mathbf{Log}_{\emptyset}^{J,\infty}$. Let $1\leq s_0\leq r_{I_{v}^{J}\cap I_1}-1$ be the unique integer such that $i_0=i_{\Theta_{J}^{\max},r_{v}^{J,s_0}}$. We have the following results.
\begin{enumerate}[label=(\roman*)]
\item \label{it: J truncation map E2 1} If $\Omega_J\in \Psi_{I_0,I_1,J,\ell}^{\natural}$, then we have
\begin{equation}\label{equ: J truncation map E2 1}
\psi_{2,\ell}(\overline{x}_{\Omega_J})=0.
\end{equation}
\item \label{it: J truncation map E2 2} If $\Omega_J\in \Psi_{I_0,I_1,J,\ell}^{\dagger}$, then we have
\begin{equation}\label{equ: J truncation map E2 2}
\psi_{2,\ell}(\overline{x}_{\Omega_J})=\pm a_{i_1}^{J}b_{i_1}^{J}(\overline{x}_{r_{i_0}(u_{i_1}^{\plog}(\Omega))}-\overline{x}_{r_{i_0}(u_{i_1}^{\val}(\Omega))})
\end{equation}
with $a_{i_1}^{J}b_{i_1}^{J}\neq 0$.
\item \label{it: J truncation map E2 3} If $\Theta_{J}^{\max}$ satisfies the $L_{J}/Z_{J}$-variant of Condition~\ref{cond: bottom deg s bis bis} for $s_0$, then there exists a unique $\tld{\Omega}\in \Psi_{I_0,I_1',\ell}$ such that $r_{i_0}(\Omega)=\tld{\Omega}^{s_0,r_{v}^{J,s_0}+1}$ and
\begin{equation}\label{equ: J truncation map E2 3}
\psi_{2,\ell}(\overline{x}_{\Omega_J})=\overline{x}_{\tld{\Omega}}.
\end{equation}
\item \label{it: J truncation map E2 4} If $\Theta_{J}^{\max}$ satisfies the $L_{J}/Z_{J}$-variant of Condition~\ref{cond: bottom deg s bis bis} for $s_0+1$, then we have $r_{i_0}(\Omega)\in \Psi_{I_0,I_1',\ell}$ and
\begin{equation}\label{equ: J truncation map E2 4}
\psi_{2,\ell}(\overline{x}_{\Omega_J})=\overline{x}_{r_{i_0}(\Omega)}.
\end{equation}
\item \label{it: J truncation map E2 5} If $\Theta_{J}^{\max}$ satisfies the $L_{J}/Z_{J}$-variant of Condition~\ref{cond: bottom deg s bis bis} for some $1\leq s_1\leq r_{I_v\cap I_1}$ with $s_1\notin\{s_0,s_0+1\}$, $i_1\defeq i_{\Theta_{J}^{\max},r_{v}^{J,s_1}}$ and $i_1'\defeq i_{\Theta_{J}^{\max},r_{v}^{J,s_1-1}}$, then there exists a unique $\tld{\Omega}\in \Psi_{I_0,I_1',\ell}$ such that $r_{i_0}(u_{i_1}(\Omega))=\tld{\Omega}^{s_1,r_{v}^{J,s_1}+1}$ (when $s_1<r_{I_v\cap I_1}$) and we have
\begin{equation}\label{equ: J truncation map E2 5}
\psi_{2,\ell}(\overline{x}_{\Omega_J})=\pm b_{i_1'}^{J}\overline{x}_{r_{i_0}(u_{i_1'}(\Omega))}\pm b_{i_1}^{J}\overline{x}_{\tld{\Omega}}.
\end{equation}
Here we omit $\pm b_{i_1'}^{J}\overline{x}_{r_{i_0}(u_{i_1'}(\Omega))}$ (resp.~$\pm b_{i_1}^{J}\overline{x}_{\tld{\Omega}}$) from (\ref{equ: J truncation map E2 5}) when $s_1=1$ (resp.~when $s_1=r_{I_v\cap I_1}$).
\end{enumerate}
\end{prop}
\begin{proof}
We prove \ref{it: J truncation map E2 1}.\\
For each $j\in J\setminus I_v$, we have $\#v\cap\mathbf{Log}_{\Delta\setminus\{j\}}=1$ by assumption and thus $u_{j}(\Omega)$ (resp.~$r_{i_0}(u_{j}(\Omega))$) is a well-defined $(I_0,I_1)$-tuple (resp.~$(I_0,I_1')$-tuple) with a fixed $v_{j}$.
Now that $\Theta_{J}^{\max}$ satisfies the $L_{J}/Z_{J}$-variant of Condition~\ref{cond: bottom deg s} for each $1\leq s\leq r_{I_{v}^{J}\cap I_1}$, we see that its underlying $(I_0,I_1)$-tuple $\Theta^{\max}$ is $(I_0,I_1)$-atomic and satisfies Condition~\ref{cond: bottom deg s} for each $1\leq s\leq r_{I_{v}\cap I_1}$, with $\Omega$ being the equivalence class of $(I_0,I_1)$-tuples that contains $\Theta^{\max}$. It thus follows from Lemma~\ref{lem: atom twist vanishing} that
\begin{equation}\label{equ: J truncation map E2 1 vanishing 1}
\overline{x}_{r_{i_0}(\Omega)}=0.
\end{equation}
Similarly, for each $j\in J\setminus I_v$, we check that there exists a unique equivalence class $\Omega_j$ of $(I_0,I_1'\setminus\{j\})$-tuples such that $r_{i_0}(u_{j}(\Omega))=r_{j}(\Omega_j)$ with the unique maximal element $\Theta_j^{\max}$ in $\Omega_j$ satisfying Condition~\ref{cond: bottom deg s} for each $1\leq s\leq r_{I_{v_j}\cap I_1}$, which together with Lemma~\ref{lem: atom twist vanishing} (upon replacing $\Omega$ and $i_0$ in \emph{loc.cit.} with $\Omega_j$ and $j$) gives
\begin{equation}\label{equ: J truncation map E2 1 vanishing 2}
\overline{x}_{r_{i_0}(u_{j}(\Omega))}=\overline{x}_{r_{j}(\Omega_j)}=0.
\end{equation}
Now we observe from (\ref{equ: J truncation log}) and (\ref{equ: J truncation val}) that $\psi_{1,\ell}(x_{\Omega_{J}})$ is a linear combination of $\overline{x}_{r_{i_0}(\Omega)}$ and $\overline{x}_{r_{i_0}(u_{j}(\Omega))}$ for various $j\in J\setminus I_v$, which together with (\ref{equ: J truncation map E2 1 vanishing 1}) and (\ref{equ: J truncation map E2 1 vanishing 2}) gives (\ref{equ: J truncation map E2 1}).

We prove \ref{it: J truncation map E2 2}.\\
According to (\ref{equ: J truncation log}) and (\ref{equ: J truncation val}), we have
\begin{equation}\label{equ: J truncation map E2 2 E1}
\psi_{1,\ell}(x_{\Omega_{J}})=\pm a_{i_1}^{J}b_{i_1}^{J}(x_{r_{i_0}(u_{i_1}^{\plog}(\Omega))}-x_{r_{i_0}(u_{i_1}^{\val}(\Omega))}))+\ast
\end{equation}
with $\ast$ being a linear combination of $x_{r_{i_0}(\Omega)}$, $x_{r_{i_0}(u_{i_1}^{\val}(u_{i_1}^{\plog}(\Omega)))}$ and $x_{r_{i_0}(u_{j}(\Omega))}$ for $j\in J\setminus (I_v\sqcup\{i_1\})$.
Parallel argument as that of (\ref{equ: J truncation map E2 1 vanishing 1}) and (\ref{equ: J truncation map E2 1 vanishing 2}) gives $\overline{x}_{r_{i_0}(\Omega)}=0$ and $\overline{x}_{r_{i_0}(u_{j}(\Omega))}=0$ for $j\in J\setminus (I_v\sqcup\{i_1\})$.
Now we observe that there exists a unique equivalence class $\Omega_{i_1}$ of $(I_0,I_1'\setminus\{i_1\})$-tuples (with a fixed $v_{i_1}$) such that
\[r_{i_0}(u_{i_1}^{\val}(u_{i_1}^{\plog}(\Omega)))=r_{i_1}(\Omega_{i_1})\]
with the unique maximal element $\Theta_{i_1}^{\max}$ in $\Omega_i$ satisfying Condition~\ref{cond: bottom deg s} for each $1\leq s\leq r_{I_{v_{i_1}}\cap I_1}$, which together with Lemma~\ref{lem: atom twist vanishing} (upon replacing $\Omega$ and $i_0$ in \emph{loc.cit.} with $\Omega_{i_1}$ and $i_1$) gives
\[\overline{x}_{r_{i_0}(u_{i_1}^{\val}(u_{i_1}^{\plog}(\Omega)))}=0.\]
This together with $\overline{x}_{r_{i_0}(\Omega)}=0$, $\overline{x}_{r_{i_0}(u_{j}(\Omega))}=0$ for $j\in J\setminus (I_v\sqcup\{i_1\})$ as well as (\ref{equ: J truncation map E2 2 E1}) gives (\ref{equ: J truncation map E2 2}).

We prove \ref{it: J truncation map E2 5} and the proof of \ref{it: J truncation map E2 3} and \ref{it: J truncation map E2 4} are simpler.\\
Recall that we write $\Theta^{\max}$ for the underlying $(I_0,I_1)$-tuple of $\Theta_{J}^{\max}$.
Our assumption on $\Theta_{J}^{\max}$ ensures that $\Theta^{\max}$ satisfies Condition~\ref{cond: bottom deg s bis bis} for $s_1$, and satisfies Condition~\ref{cond: bottom deg s} for each $1\leq s\leq r_{I_v\cap I_1}$ with $s\neq s_1$.
When $s_1<r_{I_v\cap I_1}$, we observe that $r_{i_0}(u_{i_1}(\Theta^{\max}))$ satisfies Condition~\ref{cond: bottom deg s} for each $1\leq s\leq r_{I_v\cap I_1'}$, and $r_{i_0}(u_{i_1}(\Omega))$ is exactly the unique equivalence class of $(I_0,I_1')$-tuples that contains $r_{i_0}(u_{i_1}(\Theta^{\max}))$.
When $s_1>1$, we observe that $r_{i_0}(u_{i_1'}(\Theta^{\max}))$ is $(s_1,r_{v}^{J,s_1}+1)$-atomic (see Definition~\ref{def: atom twist}) and thus there exists a unique equivalence class $\tld{\Omega}$ of $(I_0,I_1')$-atomic tuples with maximal element $\tld{\Theta}^{\max}$ such that $r_{i_0}(u_{i_1}(\Theta^{\max}))=\tld{\Theta}^{s_1,r_{v}^{J,s_1}+1}$. Now that the fixed $\tld{v}$ for $\tld{\Omega}$ satisfies $I_{\tld{v}}=(I_v\setminus\{i_0\})\sqcup\{i_1'\}$ (and thus $I_{\tld{v}}\cup I_1'=\Delta$) and $\#\tld{v}\cap \mathrm{Log}_{\Delta\setminus\{j\}}=1$ for each $j\in\Delta\setminus I_{\tld{v}}$, we easily check that $\tld{\Omega}\in \Psi_{I_0,I_1',\ell}$.
Since $\Theta^{\max}$ satisfies Condition~\ref{cond: bottom deg s bis bis} for $s_1$, we see that $\Theta$ is smaller than $\Theta^{\max}$ as a $(I_0,I_1)$-tuple if and only if $r_{i_0}(u_{i_1}(\Theta))$ is $(s_1,r_{v}^{J,s_1}+1)$-smaller than $r_{i_0}(u_{i_1}(\Theta^{\max}))$ as $(I_0,I_1')$-tuples. This together with Lemma~\ref{lem: twist atom criterion} implies that
\begin{equation}\label{equ: J truncation atom twist}
r_{i_0}(u_{i_1}(\Omega))=\tld{\Omega}^{s_1,r_{v}^{J,s_1}+1}.
\end{equation}
Now it follows from (\ref{equ: J truncation log}) and (\ref{equ: J truncation val}) that
\begin{equation}\label{equ: J truncation map E2 2 E1 prime}
\psi_{1,\ell}(x_{\Omega_{J}})=
\pm b_{i_1'}^{J}x_{r_{i_0}(u_{i_1'}(\Omega))}\pm b_{i_1}^{J}x_{r_{i_0}(u_{i_1}(\Omega))}+\ast
\end{equation}
with $\ast$ being a linear combination of $x_{r_{i_0}(\Omega)}$ and $x_{r_{i_0}(u_{j}(\Omega))}$ for $j\in J\setminus (I_v\sqcup\{i_1,i_1'\})$.
Parallel argument as that of (\ref{equ: J truncation map E2 1 vanishing 1}) and (\ref{equ: J truncation map E2 1 vanishing 2}) gives $\overline{x}_{r_{i_0}(\Omega)}=0$ and $\overline{x}_{r_{i_0}(u_{j}(\Omega))}=0$ for $j\in J\setminus (I_v\sqcup\{i_1,i_1'\})$, which together with (\ref{equ: J truncation map E2 2 E1 prime}), (\ref{equ: J truncation atom twist}) and $\overline{x}_{\tld{\Omega}}=\overline{x}_{\tld{\Omega}^{s_1,r_{v}^{J,s_1}+1}}$ from Proposition~\ref{prop: atom twist common class} imply (\ref{equ: J truncation map E2 5}).
\end{proof}

We fix the choice of a $i_0\in\Delta$ (with $\al_0\defeq (i_0,i_0+1)$) and assume $J=I_1=\Delta\setminus\{i_0\}$ in the rest of this section.
For each $\al\in\Phi^+$, we write $\al=(i_{\al},j_{\al})$ with $1\leq i_{\al}<j_{\al}\leq n$.
In particular, we have
\[L_{J}/Z_{J}=\mathrm{PGL}_{i_0}(K)\times \mathrm{PGL}_{n-i_0}(K).\]
Note that the maps between double complex
\[\cT_{\emptyset,J,J}^{\bullet,\bullet}\rightarrow \cT_{\emptyset,J}^{\bullet,\bullet}\rightarrow \cT_{\emptyset,\Delta}^{\bullet,\bullet}\]
induce a canonical map
\begin{equation}\label{equ: special J truncation total}
\vartheta_{i_0}\defeq \psi_{J}^{J}: \mathrm{Ext}_{L_{J}/Z_{J}}^{n-1}(1_{L_{J}/Z_{J}},\mathrm{St}_{L_{J}/Z_{J}}^{\rm{an}})=H^{n-1}(\mathrm{Tot}(\cT_{\emptyset,J,J}^{\bullet,\bullet}))
\rightarrow H^{n-1}(\mathrm{Tot}(\cT_{\emptyset,\Delta}^{\bullet,\bullet}))=\mathbf{E}_{n-1}.
\end{equation}
which is compatible with the canonical filtration on both the source and the target (induced from the map $E_{\bullet,\emptyset,J,J}^{\bullet,\bullet}\rightarrow E_{\bullet,\emptyset,\Delta}^{\bullet,\bullet}$ between spectral sequences) and thus induces a map
\begin{equation}\label{equ: special J truncation total graded}
\vartheta_{i_0,\ell}: \mathrm{gr}^{-\ell}(\mathrm{Ext}_{L_{J}/Z_{J}}^{n-1}(1_{L_{J}/Z_{J}},\mathrm{St}_{L_{J}/Z_{J}}^{\rm{an}}))=E_{2,\emptyset,J,J}^{-\ell,\ell+n-1}\rightarrow \mathrm{gr}^{-\ell}(\mathbf{E}_{n-1})=E_{2,\emptyset,\Delta}^{-\ell,\ell+n-1}
\end{equation}
for each $0\leq \ell\leq n-2$.

Recall from Corollary~\ref{cor: Ext between St} that
\[\mathrm{gr}^{\ell}(\mathbf{E}_{n-1})=E_{2,\emptyset,\Delta}^{-\ell,\ell+n-1}\] 
admitting a basis of the form
\begin{equation}\label{equ: full grade ell basis}
\{\overline{x}_{S,I}\mid S\in \cS_{\Delta}, \#S=n-1-\ell, I\subseteq S\cap\Delta\}
\end{equation}
for each $0\leq \ell\leq n-2$.
Now we introduce a new equivalence relation on the set of all pairs $(S,I)$ with $S\in \cS_{\Delta}$ and $I\subseteq S\cap\Delta$.
\begin{defn}\label{def: shift log basis}
We have the following notions.
\begin{enumerate}[label=(\roman*)]
\item \label{it: shift log basis 1} Let $S,S'\in\cS_{\Delta}$ be two elements that satisfy $\al_0\in S\cap S'$. We say that $S$ and $S'$ are \emph{adjacent} if $\#S\setminus S'=\#S'\setminus S=2$, $(S\setminus S')\cap\Delta\neq \emptyset\neq (S'\setminus S)\cap\Delta$ and
\begin{equation}\label{equ: shift root}
\sum_{\al\in S\setminus S'}\al=\sum_{\al\in S\setminus S'}\al\in\Phi^+.
\end{equation}
Note that any adjacent pair $S,S'$ uniquely determines a root $(i,j)\in\Phi^+$ (see (\ref{equ: shift root})) such that $j-i\geq 3$ with $S\setminus S'=\{(i,j-1),(j-1,j)\}$ and $S'\setminus S=\{(i,i+1),(i+1,j)\}$, upon exchanging $S,S'$. Note also that we have either $j\leq i_0$ or $i\geq i_0+1$.
\item \label{it: shift log basis 2} Let $(S,I)$ and $(S',I')$ be two pairs that satisfy $\al_0\in S\cap S'$. We say that $(S,I)$ and $(S',I')$ are \emph{adjacent} if $S$ and $S'$ are adjacent with $I\setminus I'=(S\setminus S')\cap I$ and $I'\setminus I=(S'\setminus S)\cap I'$. We say that $(S,I)$ and $(S',I')$ are \emph{equivalent} if there exist pairs $(S_{t'},I_{t'})$ for $0\leq t'\leq t$ with $(S,I)=(S_{0},I_{0})$ and $(S',I')=(S_{t},I_{t})$ such that $(S_{t'-1},I_{t'-1})$ and $(S_{t'},I_{t'})$ are adjacent for $1\leq t'\leq t$.
\end{enumerate}
\end{defn}

Note that, for each $0\leq \ell\leq n-2$, $\overline{x}_{S,I}$ belongs to (\ref{equ: full grade ell basis}) if and only if $\overline{x}_{S',I'}$ belongs to (\ref{equ: full grade ell basis}) for each pair $(S',I')$ that is equivalent to $(S,I)$.

We combine Proposition~\ref{prop: J truncation map E2} with Lemma~\ref{lem: E2 cup of generator} and obtain the following result describes the image of (\ref{equ: special J truncation total}).
\begin{lem}\label{lem: J truncation image generator}
Let $i_0\in\Delta$ with $\al_0=(i_0,i_0+1)$. Let $0\leq \ell\leq n-2$ and $(S,I)$ be a pair with $S\in\cS_{\Delta}$, $I\subseteq S\cap \Delta$ and $\#S=n-1-\ell$.
\begin{enumerate}[label=(\roman*)]
\item \label{it: J truncation image 1} If there exist $\beta_0\in S\setminus \Delta$ such that either $i_{\beta_0}=i_0$ or $j_{\beta_0}=i_0+1$, then we have $\overline{x}_{S,I}\in \mathrm{im}(\vartheta_{i_0,\ell})$.
\item \label{it: J truncation image 2} If $\al_0\in S$ and there exist $\beta_0\in S\setminus\Delta$ such that either $i_{\beta_0}=1$ or $j_{\beta_0}=n$, then we have $\overline{x}_{S,I}\in \mathrm{im}(\vartheta_{i_0,\ell})$.
\item \label{it: J truncation image 3} If $\al_0\in S$ and $I,I'\subseteq S\cap\Delta$ satisfy $\#I\setminus I'=\#I'\setminus I=1$, then we have
\begin{equation}\label{equ: J truncation image 3}
\overline{x}_{S,I}-\overline{x}_{S,I'}\in\mathrm{im}(\vartheta_{i_0,\ell}).
\end{equation}
\item \label{it: J truncation image 4} Let $(S',I')$ be another pair such that $\al_0\in S\cap S'$ with $(S,I)$ and $(S',I')$ being adjacent, then we have
    \begin{equation}\label{equ: J truncation image 4}
    a \overline{x}_{S',I'}+b \overline{x}_{S,I}\in  \mathrm{im}(\vartheta_{i_0,\ell})
    \end{equation}
    for some $a,b\in E^{\times}$.
\end{enumerate}
\end{lem}
\begin{proof}
Recall from Lemma~\ref{lem: E2 cup of generator} and the proof of Proposition~\ref{prop: total cup basis} that there exists a unique $\Omega_0\in \Psi_{\emptyset,\Delta,\ell_0}$ such that $\overline{x}_{S,I}=\pm\overline{x}_{\Omega_0}$ with $\ell_0\defeq n-1-\#S$. Note that $\Omega_0$ determines a $v_0\subseteq \mathbf{Log}_{\emptyset}$ and a sequence of integers $r_{\Omega_0}^{s}$ for $0\leq s\leq r_{I_{v_0}}$ with $S\cap\Delta=\{(i,i+1)\mid i\in v_0\}$. We write $\Theta_0^{\max}=(v_0,I_0,\un{k},\un{\Lambda})$ for the unique maximal tuple (with respect to $(\emptyset,\Delta)$) in $\Omega_0$, and $\{n_{d}^{\max}\}_{1\leq d\leq r_{I_0}=n-\#I_0}$ for the tuple of integers associated with $I_0$.
In each case below (according to our assumption on $(S,I)$), we will construct a maximally $(\emptyset,J)$-atomic $J$-tuple $\Theta_{J}^{\max}$ whose associated tuple with respect to $(\emptyset,J)$ is denoted by $\Theta^{\max}$.
We write $\Omega_{J}=(v,I^{J},\un{k}^{J},\un{\Lambda}^{J})$ for the equivalence class of $J$-tuples (with respect to $(\emptyset,J)$) containing $\Theta_{J}^{\max}$, whose associated equivalence class $\Omega$ of tuples (with respect to $(\emptyset,J)$) contains $\Theta^{\max}$.

We prove \ref{it: J truncation image 1}.\\
Let $\beta_0\in S\setminus\Delta$ be as in our assumption. We have the following two possibilities.
\begin{itemize}
\item If $i_{\beta_0}=i_0$, then there exists $1\leq s_0\leq r_{I_{v_0}}$ and $r_{\Omega_0}^{s_0-1}+1\leq d_0\leq r_{\Omega_0}^{s_0}-1$ such that $i_0=i_{\Theta_0^{\max},d_0}$, in which case there exists a unique $J$-tuple $\Theta_{J}^{\max}$ with respect to $(\emptyset,J)$ with its associated tuple $\Theta^{\max}$ with respect to $(\emptyset,J)$ such that $\Theta_0^{\max}=r_{i_0}(\Theta^{\max})$. Now that $\Lambda_{d}=\{2n_d-1\}$ for each $d_0+1\leq d\leq r_{\Omega_0}^{s_0}$, we see that $\Theta_{J}^{\max}$ satisfies the $L_{J}/Z_{J}$-variant of Condition~\ref{cond: bottom deg s bis bis} for $s_0+1$, and thus we have 
    \[\overline{x}_{S,I}=\pm\overline{x}_{\Omega_0}=\pm\vartheta_{i_0,\ell}(\overline{x}_{\Omega_{J}})\] 
    by \ref{it: J truncation map E2 4} of Proposition~\ref{prop: J truncation map E2}.
\item If $j_{\beta_0}=i_0+1$, then there exists $1\leq s_0\leq r_{I_{v_0}}$ and $r_{\Omega_0}^{s_0-1}+2\leq d_0\leq r_{\Omega_0}^{s_0}$ such that $i_0+1=i_{\Theta_0^{\max},d_0}$, in which case there exists a unique $J$-tuple $\Theta_{J}^{\max}$ with respect to $(\emptyset,J)$ with its associated tuple $\Theta^{\max}$ with respect to $(\emptyset,J)$ such that $(\Theta_0^{\max})^{s_0,d_0}=r_{i_0}(\Theta^{\max})$. We have \[\overline{x}_{S,I}=\pm\overline{x}_{\Omega_0}=\pm\vartheta_{i_0,\ell}(\overline{x}_{\Omega_{J}})\] 
    by \ref{it: J truncation map E2 3} of Proposition~\ref{prop: J truncation map E2}.
\end{itemize}

We prove \ref{it: J truncation image 2}.\\
Let $\beta_0\in S\setminus\Delta$ be as in our assumption. Note that $\al_0=(i_0,i_0+1)\in S$ forces $i_0\in \Delta\setminus I_{v_0}$. We have the following two possibilities.
\begin{itemize}
\item If $i_{\beta_0}=1$, then $\Lambda_{1}=\emptyset$ (with $i_{\Theta_0^{\max},1}=1$) and $\Lambda_{2}=\{2n_{2}^{\max}-1\}$ with $i_{\Theta_0^{\max},2}-1=n_{2}^{\max}\geq 2$. We set $i_1\defeq n_{2}^{\max}=i_{\Theta_0^{\max},2}-1$ and $I^{J}\defeq (I_0\sqcup\{1\})\setminus\{i_1\}$. We set $v\defeq (v_0\setminus\{\plog_{i_0}\})\sqcup\{\plog_{i_1}\}$ if $\plog_{i_0}\in v_0$, and $v\defeq (v_0\setminus\{\val_{i_0}\})\sqcup\{\val_{i_1}\}$ if $\val_{i_0}\in v_0$. We also set $\Lambda_{1}^{J}\defeq \{2i_1-1\}$, $\Lambda_{2}^{J}\defeq \emptyset$ and $\Lambda_{d}^{J}\defeq \Lambda_{d}$ for each $3\leq d\leq r_{I}=r_{I_0}$. We thus obtain a $J$-tuple $\Theta_{J}^{\max}=(v,I^{J},\un{k}^{J},\un{\Lambda}^{J})$ whose associated $\Omega_{J}$ and $\Omega$ satisfies $u_{i_1}(\Omega)=\Omega_0^{1,1}$. This together with \ref{it: J truncation map E2 5} of Proposition~\ref{prop: J truncation map E2} (with $s_1=1$ in \emph{loc.cit.}) gives 
    \[\overline{x}_{S,I}\in \mathrm{im}(\vartheta_{i_0,\ell}).\]
\item If $j_{\beta_0}=n$, then we have $\Lambda_{r_{I_0}}=\{2n_{r_{I_0}}-1\}$ (with $n_{r_{I_0}}\geq 2$) and $i_1'\defeq i_{\Theta_0^{\max},r_{I_0}-1}\in I_{v_0}$. We set $v\defeq (v_0\setminus\{\plog_{i_0}\})\sqcup\{\plog_{i_1'}\}$ if $\plog_{i_0}\in v_0$, and $v\defeq (v_0\setminus\{\val_{i_0}\})\sqcup\{\val_{i_1'}\}$ if $\val_{i_0}\in v_0$.
    We set $I^{J}\defeq I_0$ and $\Lambda_{d}^{J}\defeq \Lambda_{d}$ for each $1\leq d\leq r_{I}=r_{I_0}$. We thus obtain a $J$-tuple $\Theta_{J}^{\max}=(v,I^{J},\un{k}^{J},\un{\Lambda}^{J})$ whose associated $\Omega_{J}$ and $\Omega$ satisfies $u_{i_1}(\Omega)=\Omega_0$. This together with \ref{it: J truncation map E2 5} of Proposition~\ref{prop: J truncation map E2} (with $s_1=r_{I_v\cap I_1}$ in \emph{loc.cit.}) gives 
    \[\overline{x}_{S,I}\in \mathrm{im}(\vartheta_{i_0,\ell}).\]
\end{itemize}

We prove \ref{it: J truncation image 3}.\\
Assume now that $\al_0\in S$ and let $I,I'\subseteq S\cap\Delta$ be two subsets satisfying $\#I\setminus I'=\#I'\setminus I=1$. We write $I\setminus I'=\{i_1\}$ and $I'\setminus I=\{i_1'\}$ for short.
By Lemma~\ref{lem: E2 cup of generator} and the proof of Proposition~\ref{prop: total cup basis} we know that there exists a unique $\Omega_0'\in \Psi_{\emptyset,\Delta,\ell_0}$ such that $\overline{x}_{S,I'}=\pm\overline{x}_{\Omega_0'}$, and $\Omega_0$ determines a $v_0'\subseteq \mathbf{Log}_{\emptyset}$ with $I_{v_0'}=I_{v_0}$ and a sequence of integers $r_{\Omega_0'}^{s}$ with $r_{\Omega_0'}^{s}=r_{\Omega_0}^{s}$ for $0\leq s\leq r_{I_{v_0}}$. In fact, we have $v_0\setminus v_0'=\{\val_{i_1},\plog_{i_1'}\}$ and $v_0'\setminus v_0=\{\val_{i_1'},\plog_{i_1}\}$, and the unique maximal tuple in $\Omega_0'$ has the form $(v_0',J_0,\un{k},\un{\Lambda})$.
We have the following two possibilities.
\begin{itemize}
\item If $i_0\in\{i_1,i_1'\}$, we may assume that $i_0=i_1'$ upon exchanging $I,I'$. In this case, the data $(v_0'\setminus\{\val_{i_0}\})\sqcup\{\val_{i_1}\}$, $J_0$ and $\un{\Lambda}$ uniquely determines a $J$-tuple $\Theta_{J}^{\max}$ (with respect to $(\emptyset,J)$) whose associated equivalence class $\Omega$ of tuples with respect to $(\emptyset,J)$ satisfies $\Omega_0=r_{i_0}(u_{i_1}^{\plog}(\Omega))$ and $\Omega_0'=r_{i_0}(u_{i_1}^{\val}(\Omega))$. This together with \ref{it: J truncation map E2 2} of Proposition~\ref{prop: J truncation map E2} gives (\ref{equ: J truncation image 3}).
\item If $i_0\notin\{i_1,i_1'\}$, then we have either $i_0\in I\cap I'$ or $i_0\in S\setminus(I\cup I')$. We set $I''\defeq (I\cup I')\setminus\{i_0\}$ if $i_0\in I\cap I'$, and $I''\defeq (I\cap I')\sqcup\{i_0\}$ if $i_0\in S\setminus (I\cup I')$. It is clear that both pairs $I,I''$ and $I'',I'$ satisfy the condition of the previous paragraph, and thus we deduce that
    \[\overline{x}_{S,I}-\overline{x}_{S,I'}=(\overline{x}_{S,I}-\overline{x}_{S,I''})+(\overline{x}_{S,I''}-\overline{x}_{S,I'})\in\mathrm{im}(\vartheta_{i_0,\ell}).\]
\end{itemize}

We prove \ref{it: J truncation image 4}.\\
Parallel to the discussion in the first paragraph of the proof, we may associate a unique $\Omega_0'\in \Psi_{\emptyset,\Delta,\ell_0}$ such that $\overline{x}_{S',I'}=\pm\overline{x}_{\Omega_0'}$.
Now that $(S,I)$ and $(S',I')$ are adjacent, by \ref{it: shift log basis 1} of Definition~\ref{def: shift log basis} we may associate a pair of integers $i_1,i_1'\in\Delta$ such that
\[\sum_{\al\in S\setminus S'}\al=(i_1',i_1+1)=\sum_{\al\in S'\setminus S}\al\]
with $i_1-i_1'\geq 2$.
Upon exchanging $(S,I)$ and $(S',I')$, it is harmless to assume that $\al_1\defeq (i_1,i_1+1)\in S$ and $\al_1'\defeq (i_1',i_1'+1)\in S'$.
Note from \ref{it: shift log basis 1} of Definition~\ref{def: shift log basis} that we have $\al_0\in I$ (resp.~ $\al_1\in I$) if and only if $\al_0\in I'$ (resp.~if and only if $\al_1'\in I'$). We have the following possibilities.
\begin{itemize}
\item Assume that $\al_1\in I$, $\al_1'\in I'$ and $\al_0\in I\cap I'$. Then there exists a unique $\Omega_J$ and $1\leq s_1\leq r_{I_v\cap J}$ (with maximal element $\Theta_J^{\max}=(v,I^J,\un{k}^J,\un{\Lambda}^J)$, $i_1=i_{\Theta_J^{\max},s_1}$ and $i_1'=i_{\Theta_J^{\max},s_1-1}$) such that the associated $\Omega$ satisfies $\Omega_0=r_{i_0}(u_{i_1'}(\Omega))$ and $(\Omega_0')^{s_1,2}=r_{i_0}(u_{i_1}(\Omega))$. This together with \ref{it: J truncation map E2 5} of Proposition~\ref{prop: J truncation map E2} gives (\ref{equ: J truncation image 4}). The case when $\al_1\notin I$, $\al_1'\notin I'$ and $\al_0\notin I\cup I'$ is similar.
\item Assume that $\al_1\in I$, $\al_1'\in I'$ and $\al_0\notin I\cup I'$. We set $I''\defeq (I'\setminus\{\al_1'\})\cup\{\al_0\}$ and note that there exists a unique $\Omega_0''\in \Psi_{\emptyset,\Delta,\ell_0}$ such that $\overline{x}_{S',I''}=\pm\overline{x}_{\Omega_0''}$. Then again there exists a unique $\Omega_J$ and $1\leq s_1\leq r_{I_v\cap J}$ (with maximal element $\Theta_J^{\max}=(v,I^J,\un{k}^J,\un{\Lambda}^J)$, $i_1=i_{\Theta_J^{\max},s_1}$ and $i_1'=i_{\Theta_J^{\max},s_1-1}$) such that the associated $\Omega$ satisfies $\Omega_0=r_{i_0}(u_{i_1'}(\Omega))$ and $(\Omega_0'')^{s_1,2}=r_{i_0}(u_{i_1}(\Omega))$, which together with \ref{it: J truncation map E2 5} of Proposition~\ref{prop: J truncation map E2} gives
    \[a \overline{x}_{S',I''}+b \overline{x}_{S,I}\in  \mathrm{im}(\vartheta_{i_0,\ell})\]
    for some $a,b\in E^{\times}$.
    Now that we have $\overline{x}_{S',I''}-\overline{x}_{S',I'}\in \mathrm{im}(\vartheta_{i_0,\ell})$ from \ref{it: J truncation image 3}, we again deduce (\ref{equ: J truncation image 4}). The case when $\al_1\notin I$, $\al_1'\notin I'$ and $\al_0\in I\cap I'$ is similar.
\end{itemize}
The proof is thus finished.
\end{proof}

Assume in the rest of this section that $i_0=n-1$ and omit $i_0$ from the notation.
\begin{prop}\label{prop: J truncation total basis}
Let $0\leq \ell\leq n-2$. We have the following results.
\begin{enumerate}[label=(\roman*)]
\item \label{it: J truncation total basis 1} If $\ell>0$, then $\vartheta_{\ell}$ is surjective.
\item \label{it: J truncation total basis 2} If $\ell=0$, then the space $\mathrm{im}(\vartheta_{0})$ admits a basis of the form
\begin{equation}\label{equ: J truncation total basis 2}
\{\overline{x}_{\Delta,I}-\overline{x}_{\Delta,[n-\#I,n-1]}\mid [n-\#I,n-1]\neq I\subseteq \Delta\}.
\end{equation}
\end{enumerate}
\end{prop}
\begin{proof}
Recall that $E_{2,\emptyset,\Delta}^{-\ell,\ell+n-1}$ admits a basis of the form (\ref{equ: full grade ell basis}). When $\ell=0$, it is clear that (\ref{equ: full grade ell basis}) differs from 
\begin{equation}\label{equ: J truncation modified basis}
\{\overline{x}_{\Delta,I}-\overline{x}_{\Delta,[n-\#I,n-1]}\mid [n-\#I,n-1]\neq I\subseteq \Delta\}\sqcup\{\overline{x}_{\Delta,[n-i,n-1]}\mid 0\leq i\leq n-1\}
\end{equation}
by an invertible matrix. In particular, (\ref{equ: J truncation modified basis}) forms a basis of $E_{2,\emptyset,\Delta}^{0,n-1}$, and (\ref{equ: J truncation total basis 2}) is linearly independent.
We divide the rest of the proof into the following steps.

\textbf{Step $1$}: We prove that
\[\overline{x}_{\Delta,I}-\overline{x}_{\Delta,[n-\#I,n-1]}\in \mathrm{im}(\vartheta_{0}).\]
We write $I'\defeq [n-\#I,n-1]$ and $t\defeq \#I\setminus I'=\#I'\setminus I$ for short. We choose an arbitrary sequence of sets $J_{t'}$ for $1\leq t'\leq t$ such that $J_{0}=I$, $J_{t}=I'$ and $\#J_{t'}\setminus J_{t'-1}=\#J_{t'-1}\setminus J_{t'}=1$ for each $1\leq t'\leq t$. It follows from \ref{it: J truncation image 3} of Lemma~\ref{lem: J truncation image generator} that $\overline{x}_{\Delta,J_{t'}}-\overline{x}_{\Delta,J_{t'-1}}\in \mathrm{im}(\vartheta_{0})$, and therefore
\[\overline{x}_{\Delta,I}-\overline{x}_{\Delta,I'}=\sum_{t'=1}^{t}\overline{x}_{\Delta,J_{t'}}-\overline{x}_{\Delta,J_{t'-1}}\in \mathrm{im}(\vartheta_{0}).\]

\textbf{Step $2$}: We prove that $\mathrm{im}(\vartheta_{0})$ is spanned by (\ref{equ: J truncation total basis 2}).\\
Following the notation of Proposition~\ref{prop: J truncation map E2}, we know that $\mathrm{im}(\vartheta_{0})$ is spanned by $\vartheta_{0}(\overline{x}_{\Omega_{J}})$ for varying $\Omega_{J}\in \Psi_{\emptyset,J,J,0}'$. Now that $I_1=J=\Delta\setminus\{n-1\}$ (under the notation of Proposition~\ref{prop: J truncation map E2}), we see that $\Psi_{\emptyset,J,J,0}'=\Psi_{\emptyset,J,J,0}^{\dagger}$ (namely $\Psi_{\emptyset,J,J,0}^{\natural}=\emptyset=\Psi_{\emptyset,J,J,0}^{\ddagger}$) and from \ref{it: J truncation map E2 2} of Proposition~\ref{prop: J truncation map E2} that
\begin{equation}\label{equ: J truncation image zero piece}
\vartheta_{0}(\overline{x}_{\Omega_{J}})=c(\overline{x}_{\Delta,I'}-\overline{x}_{\Delta,I''})
\end{equation}
for some $c\in E^{\times}$ and $I',I''\subseteq \Delta$ associated with $\Omega_{J}$ with $\#I'=\#I''$. Now that (\ref{equ: J truncation image zero piece}) clearly lies in the span of (\ref{equ: J truncation total basis 2}), the proof is finished.

\textbf{Step $3$}: We prove that $\overline{x}_{S,I}\in \mathrm{im}(\vartheta_{\ell})$ for each $S\in\cS_{\Delta}$ satisfying $I\subseteq S\cap\Delta$ and $\#S=n-1-\ell<n-1$ (with $S\neq\Delta$).\\
We have the following possibilities.
\begin{itemize}
\item If there exists $\beta_0\in S\setminus\Delta$ such that $j_{\beta_0}=n=i_0+1$, then we have $\overline{x}_{S,I}\in \mathrm{im}(\vartheta_{\ell})$ by \ref{it: J truncation image 1} of Lemma~\ref{lem: J truncation image generator}.
\item If there exists $\beta_0\in S\setminus\Delta$ such that $i_{\beta_0}=1<i_0=n-1$, then we have $\overline{x}_{S,I}\in \mathrm{im}(\vartheta_{\ell})$ by \ref{it: J truncation image 2} of Lemma~\ref{lem: J truncation image generator}.
\item Assume now that $\al_0=(n-1,n)\in S$ and $(1,2)\in S$. As $S\neq \Delta$, there exists $S'\in \cS_{\Delta}$ and $I'\subseteq S'\cap\Delta$ such that $\al_0\in S'$, $(S,I)$ and $(S',I')$ are equivalent (see Definition~\ref{def: shift log basis}) and that there exists $\beta_0\in S'\setminus\Delta$ such that $i_{\beta_0}=1<i_0=n-1$. We again have $\overline{x}_{S',I'}\in \mathrm{im}(\vartheta_{\ell})$ by \ref{it: J truncation image 2} of Lemma~\ref{lem: J truncation image generator}. Now that $(S,I)$ and $(S',I')$ are equivalent, we have by definition a sequence of pairs $(S_{r'},I_{r'})$ for $1\leq r'\leq r$ with $(S_0,I_0)=(S,I)$, $(S_{r},I_{r})=(S',I')$ such that $(S_{r'},I_{r'})$ and $(S_{r'-1},I_{r'-1})$ are adjacent for $1\leq r'\leq r$.
    Hence, it follows from \ref{it: J truncation image 4} of Lemma~\ref{lem: J truncation image generator} that
    \[a_{r'} \overline{x}_{S_{r'},I_{r'}}+ b_{r'} \overline{x}_{S_{r'-1},I_{r'-1}}\in \mathrm{im}(\vartheta_{\ell})\]
    for some $a_{r'},b_{r'}\in E^{\times}$ and each $1\leq r'\leq r$, which together with $\overline{x}_{S_{r},I_{r}}=\overline{x}_{S',I'}\in \mathrm{im}(\vartheta_{\ell})$ gives $\overline{x}_{S,I}=\overline{x}_{S_{0},I_{0}}\in \mathrm{im}(\vartheta_{\ell})$.
\end{itemize}

It follows from \textbf{Step $3$} that 
\[\mathrm{im}(\vartheta_{\ell})=E_{2,\emptyset,\Delta}^{-\ell,\ell+n-1}\]
when $\ell>0$, which gives \ref{it: J truncation total basis 1}.
Combining \textbf{Step $1$} with \textbf{Step $2$}, we know that the linearly independent set (\ref{equ: J truncation total basis 2}) forms a basis of $\mathrm{im}(\vartheta_{0})$, and finish the proof of \ref{it: J truncation total basis 2}.
\end{proof}

\subsection{Relative conditions}\label{subsec: relative Tits}
Let $I_0\subseteq I_1\subseteq \Delta$ and $\fh\subseteq \fg$ be an $E$-Lie subalgebra such that $\fh=\fh_{J_{\fh}}$ for some interval $J_{\fh}\subseteq \Delta$.
We consider the maps between double complex $\cT_{I_0,I_1,\fh}^{\bullet,\bullet}\rightarrow \cT_{I_0,I_1}^{\bullet,\bullet}$ and $\mathrm{CE}_{I_0,I_1,\fh}^{\bullet,\bullet}\rightarrow \mathrm{CE}_{I_0,I_1}^{\bullet,\bullet}$ (constructed in \S \ref{subsec: std seq}) and study the associated maps between spectral sequences for the bottom non-vanishing degree with the main result being Proposition~\ref{prop: relative bottom embedding} and Proposition~\ref{prop: grade M relative image}.
Then we focus on the study of $\cT_{I_0,I_1,\fh,J,M}^{\bullet,\bullet}\rightarrow \cT_{I_0,I_1}^{\bullet,\bullet}$ (and its associated map between spectral sequences) for specific choice of $\fh$, $J$ and $M$ with main results being Proposition~\ref{prop: relative coh E2 grade} and Proposition~\ref{prop: relative coh E2 grade dagger}. These results are key ingredients of the proof of the commutativity of cup product in \S \ref{subsec: commutativity}.

Let $\fh\subseteq \fg$ be as above.
For each $\#I_0\leq \ell\leq \#I_1$, we define $\Psi_{I_0,I_1,\fh,\ell}\subseteq \Psi_{I_0,I_1,\ell}$ as the subset consisting of those $\Omega$ whose associated $(S,I)$ under (\ref{equ: atom to partition}) satisfies the following condition: each $\beta\in S\setminus I$ satisfies $I_{\beta}\not\subseteq J_{\fh}$ (see the discussion below (\ref{equ: sum of roots}) for $I_{\beta}$). 
For each $J\subseteq \Delta$ (resp.~$i\geq 0$), we also define $\mathrm{gr}^{J}(\Psi_{I_0,I_1,\fh,\ell})$ (resp.~$\mathrm{gr}^{i}(\Psi_{I_0,I_1,\fh,\ell})$) as the subset of $\Psi_{I_0,I_1,\fh,\ell}$ consisting of those $\Omega$ whose associated $(S,I)$ satisfies $I=J$ (resp.~$\#I=i$).

\begin{lem}\label{lem: relative E2 atom}
Let $I_0\subseteq I_1\subseteq \Delta$ and $\fh=\fh_{J_{\fh}}$ for some interval $J_{\fh}\subseteq \Delta$. Let $\#I_0\leq \ell\leq \#I_1$ and $\Omega\in\Psi_{I_0,I_1,\fh,\ell}$. Then there exists an element $\overline{x}_{\Omega,\fh}\in E_{2,I_0,I_1,\fh}^{-\ell,\ell+\#I_1-2\#I_0}$ whose image under
\begin{equation}\label{equ: relative E2 atom}
E_{2,I_0,I_1,\fh}^{-\ell,\ell+\#I_1-2\#I_0}\rightarrow E_{2,I_0,I_1}^{-\ell,\ell+\#I_1-2\#I_0}
\end{equation}
equals $\overline{x}_{\Omega}$.
If furthermore $\Omega\in\mathrm{gr}^{M}(\Psi_{I_0,I_1,\fh,\ell})$ for some $M\subseteq \Delta$ (uniquely determined by $\Omega$), then $\overline{x}_{\Omega,\fh}$ naturally determined an element of $\mathrm{gr}^{M}(E_{2,I_0,I_1,\fh}^{-\ell,\ell+\#I_1-2\#I_0})$ whose image under
\begin{equation}\label{equ: relative E2 atom J}
\mathrm{gr}^{M}(E_{2,I_0,I_1,\fh}^{-\ell,\ell+\#I_1-2\#I_0})\rightarrow \mathrm{gr}^{M}(E_{2,I_0,I_1}^{-\ell,\ell+\#I_1-2\#I_0})
\end{equation}
equals $\overline{x}_{\Omega}$.
\end{lem}
\begin{proof}
We write $h\defeq \#I_1-2\#I_0$ for short.
We fix a choice of $\#I_0\leq \ell\leq \#I_1$ as well as $\Omega\in\Psi_{I_0,I_1,\fh,\ell}$ throughout the proof. 
We write $J_{\fh}=[j_{\fh},j_{\fh}']$ and given an explicit construction of $\overline{x}_{\Omega,\fh}\in E_{2,I_0,I_1,\fh}^{-\ell,\ell+\#I_1-2\#I_0}$ below.

We write $\Theta=(v,J,\un{k},\un{\Lambda})$ for the unique maximal element in $\Omega$ (which is a $(I_0,I_1)$-atomic tuple), and $\{n_{d}\}_{1\leq d\leq r_{J}}$ for the tuple of integers associated with $J$.
We write $(S,I)$ for the pair associated with $\Omega$ under (\ref{equ: atom to partition}). 
Recall that $\beta=(j_{\beta},j_{\beta}')\in S\cap\Delta$ (with $j_{\beta}'=j_{\beta}+1$) if and only if $j_{\beta}\in\Delta\setminus I_{v}$ (or equivalently there exists a unique $1\leq s\leq r_{I_{v}\cap I_1}-1$ such that $j_{\beta}=i_{\Theta,r_{\Omega}^{s}}$ and $j_{\beta}'=i_{\Theta,r_{\Omega}^{s}+1}$), and that $\beta=(j_{\beta},j_{\beta}')\in S\setminus\Delta$ (with $j_{\beta}'>j_{\beta}+1$) if and only if there exists $1\leq s\leq r_{I_{v}\cap I_1}$ and $r_{\Omega}^{s-1}+2\leq d\leq r_{\Omega}^{s}$ such that $j_{\beta}=i_{\Theta,d-1}$ and $j_{\beta}'=\min\{j\in[i_{\Theta,d-1}+1,i_{\Theta,d}]\mid j\notin I_1\setminus I_0\}$ with $[j_{\beta}',i_{\Theta,d}-1]=J^{d}\cap I_0$.
In particular, note that $\beta\in S$ is in bijection with $2\leq d\leq r_{J}$.
By the very definition of $\Psi_{I_0,I_1,\fh,\ell}$ we have $I_{\beta}\not\subseteq J_{\fh}$ for each $\beta\in S\setminus I$.
If $\beta\in\Delta$, then we have $I_{\beta}=\{\beta\}$ and thus $\beta\in J_{\fh}$ only if $\beta\in I$. In other words, any $2\leq d\leq r_{J}$ that satisfies $i_{\Theta,d}=i_{\Theta,d-1}+1$ and $j_{\fh}\leq i_{\Theta,d-1}\leq j_{\fh}'$ must satisfy $\val_{i_{\Theta,d-1}}\in v$.
We consider $1\leq s_0\leq r_{I_{v}\cap I_1}$ and $r_{\Omega}^{s_0-1}+2\leq d_0\leq r_{\Omega}^{s_0}$ such that $i_{\Theta,d_0-1}+1\leq j_{\fh}$ and that $d_0$ is maximal possible. If such $d_0$ does not exist, we write $s_0\defeq 1$ and $d_0\defeq 1$ for convention.
For each $d<d_0$, we have $i_{\Theta,d}\leq i_{\Theta,d_0-1}<j_{\fh}$ and thus $J^{d}\cap J_{\fh}=\emptyset$.
For each $d>d_0$ with $\Lambda_{d}\neq\emptyset$, we have $i_{\Theta,d-1}\geq j_{\fh}$ by the maximality of $d_0$, which forces $\max\{j\in J^{d}\setminus I_0\}>j_{\fh}'$ as the root $\beta\in S\setminus\Delta$ associated with such $d$ satisfies $I_{\beta}\not\subseteq J_{\fh}$.
Now we write 
\[\Omega'\defeq (\cdots(\Omega^{s_0,d_0})^{s_0-1,r_{\Omega}^{s_0-1}}\cdots)^{1,r_{\Omega}^1}\] 
for short and consider an arbitrary $\Theta'=(v,J',\un{k}',\un{\Lambda}')\in\Omega'$.
Let $1\leq s\leq r_{I_{v}\cap I_1}$ and $r_{\Omega'}^{s-1}+1\leq d\leq r_{\Omega'}^{s}$ with $\Lambda'_{d}\neq\emptyset$. 
If $s<s_0$, then we obviously have $(J')^{d}\cap J_{\fh}=\emptyset$.
If $d>d_0$, then we have $(J')^{d}=J^{d}$ with $J^{d}\setminus I_0\not\subseteq J_{\fh}$ by previous discussion.
If $s=s_0$ and $d<d_0-1$, then $i_{\Theta',d}\leq i_{\Theta^{s_0,d_0},d}\leq i_{\Theta^{s_0,d_0},d_0-2}<i_{\Theta,d_0-1}<j_{\fh}$ and thus $(J')^{d}\cap J_{\fh}=\emptyset$.
If $s=s_0$ and $d=d_0-1$, then we have $i_{\Theta,d_0-1}\in (J')^{d_0-1}\setminus(J_{\fh}\sqcup I_0)$ and thus $(J')^{d_0-1}\not\subseteq J_{\fh}\sqcup I_0$.
To summarize, we have $(J')^{d}\not\subseteq J_{\fh}\sqcup I_0$ for each $1\leq d\leq r_{J'}$ with $\Lambda'_{d}\neq\emptyset$.
In particular, for each $\Theta'$ and $1\leq d\leq r_{J'}$ with $\Lambda'_{d}\neq\emptyset$, we obtain from Proposition~\ref{prop: relative interval embedding} and Remark~\ref{rem: h tuple} an element $x_{\Theta',\fh}\in H^{\ell+h}(L_{J'},\fl_{J'}\cap\fh,1_{L_{J'}})$ whose image under the embedding
\[H^{\ell+h}(L_{J'},\fl_{J'}\cap\fh,1_{L_{J'}})\hookrightarrow H^{\ell+h}(L_{J'},1_{L_{J'}})\]
is $x_{\Theta'}$.
We set
\begin{equation}\label{equ: relative explicit atom E1}
x_{\Omega',\fh}\defeq \sum_{\Theta'\in\Omega'}\varepsilon(\Theta')x_{\Theta',\fh}\in E_{1,I_0,I_1,\fh}^{-\ell,\ell+h}
\end{equation}
and note that its image in $E_{1,I_0,I_1,\fh}^{-\ell,\ell+h}$ is nothing but $x_{\Omega'}$.
Thanks to Proposition~\ref{prop: relative interval embedding}, we know that the following natural map
\[H^{\bullet}(L_{J'},\fl_{J'}\cap\fh,1_{L_{J'}})\rightarrow H^{\bullet}(L_{J'},1_{L_{J'}})\]
is an embedding between graded $E$-algebras for each $I_0\subseteq J'\subseteq I_1$, from which we obtain the following commutative diagram
\begin{equation}\label{equ: relative E2 atom diagram}
\xymatrix{
E_{1,I_0,I_1,\fh}^{-\ell,\ell+h} \ar^{d_{1,I_0,I_1,\fh}^{-\ell,\ell+h}}[rr] \ar@{^{(}->}[d] & & E_{1,I_0,I_1,\fh}^{-\ell+1,\ell+h} \ar@{^{(}->}[d] \\
E_{1,I_0,I_1}^{-\ell,\ell+h} \ar^{d_{1,I_0,I_1}^{-\ell,\ell+h}}[rr] & & E_{1,I_0,I_1}^{-\ell+1,\ell+h}
}
\end{equation}
with both vertical maps being injective. Here the explicit description of $d_{1,I_0,I_1,\fh}^{-\ell,\ell+h}$ using $\fh$-tuples relies on Corollary~\ref{cor: relative Levi restriction}.
Recall from \ref{it: atom basis 1} of Proposition~\ref{prop: atom basis} and Proposition~\ref{prop: atom twist common class} that $d_{1,I_0,I_1}^{-\ell,\ell+h}(x_{\Omega'})=0$ and that the image of $x_{\Omega'}$ under
\[\mathrm{ker}(d_{1,I_0,I_1}^{-\ell,\ell+h})\twoheadrightarrow E_{2,I_0,I_1}^{-\ell,\ell+h}\]
equals $\overline{x}_{\Omega}$, we deduce from (\ref{equ: relative E2 atom diagram}) (especially the injectivity of the right vertical map) that $d_{1,I_0,I_1,\fh}^{-\ell,\ell+h}(x_{\Omega',\fh})=0$ and that the image of $x_{\Omega',\fh}$ under the composition of
\begin{equation}\label{equ: relative E2 atom image}
\mathrm{ker}(d_{1,I_0,I_1,\fh}^{-\ell,\ell+h})\twoheadrightarrow E_{2,I_0,I_1,\fh}^{-\ell,\ell+h}\rightarrow E_{2,I_0,I_1}^{-\ell,\ell+h}
\end{equation}
equals $\overline{x}_{\Omega}$. We define $\overline{x}_{\Omega,\fh}$ as the image of $x_{\Omega',\fh}$ under the LHS map of (\ref{equ: relative E2 atom image}) and finish the proof.
\end{proof}

\begin{lem}\label{lem: relative E2 atom degeneracy}
Let $I_0\subseteq I_1\subseteq \Delta$, $\fh$, $\#I_0\leq \ell\leq \#I_1$, $\Omega\in\Psi_{I_0,I_1,\fh,\ell}$ and $M\subseteq\Delta$ as in Lemma~\ref{lem: relative E2 atom}. Assume that $I_0\subseteq J_{\fh}$. Then we have
\begin{equation}\label{equ: relative E2 atom vanishing}
d_{r,I_0,I_1,\fh}^{-\ell,\ell+\#I_1-2\#I_0}(\overline{x}_{\Omega,\fh})=0
\end{equation}
and
\begin{equation}\label{equ: relative E2 atom vanishing J}
\mathrm{gr}^{M}(d_{r,I_0,I_1,\fh}^{-\ell,\ell+\#I_1-2\#I_0})(\overline{x}_{\Omega,\fh})=0
\end{equation}
for each $r\geq 2$.
\end{lem}
\begin{proof}
We continue to use notation in the proof of Lemma~\ref{lem: relative E2 atom}.
In particular, we have $v\subseteq\mathbf{Log}_{\emptyset}$, $1\leq s_0\leq r_{I_{v}\cap I_1}$, $r_{\Omega}^{s_0-1}+1\leq d_0\leq r_{\Omega}^{s_0}$ and $\Omega'=\Omega^{s_0,d_0}$, with an associated element $x_{\Omega',\fh}\in E_{1,I_0,I_1,\fh}^{-\ell,\ell+h}$ as in (\ref{equ: relative explicit atom E1}).
We continue to write $\Theta=(v,I,\un{k},\un{\Lambda})$ for the unique maximal element in $\Omega$, and $\Theta'=(v,J',\un{k}',\un{\Lambda}')$ for an arbitrary element of $\Omega'$.
Let $I_0^{e}$ be a non-empty subinterval of $I_0$ (for some $1\leq e\leq r_{I_0}$), then there exists a unique $1\leq d\leq n-\ell$ such that $I_0^{e}\subseteq J^{d}$. Our assumption $I_0\subseteq J_{\fh}$ forces $\emptyset\neq I_0^{e}\subseteq J^{d}\cap J_{\fh}$.
By our discussion on $s_0$ and $d_0$ in the proof of Lemma~\ref{lem: relative E2 atom}, we have the following possibilities.
\begin{itemize}
\item If $\Lambda_{d}=\emptyset$, then we have $\Lambda'_{d}=\emptyset$ and $(J')^{d}=J^{d}=I_0^{e}\subseteq J_{\fh}$.
\item If $\Lambda_{d}\neq\emptyset$ and $d>d_0$, then we have $\max\{j\in J^{d}\setminus I_0\}>j_{\fh}'$ which contradicts $\emptyset\neq I_0^{e}\subseteq J^{d}\cap J_{\fh}$.
\item If $\Lambda_{d}\neq\emptyset$ and $d=d_0$, then we have $\Lambda'_{d_0}=\emptyset$ and $(J')^{d_0}=I_0^{e}\subseteq J_{\fh}$.
\end{itemize}
Consequently, we always have $\Lambda'_{d}=\emptyset$ and $(J')^{d}=I_0^{e}\subseteq J_{\fh}$. An argument parallel to that of Lemma~\ref{lem: unique atom} thus shows that $\Theta'=\Theta^{s_0,d_0}$ is the unique element of $\Omega'$.

In the rest of the proof, following the ideas in the proof of Lemma~\ref{lem: bottom deg differential vanishing} (which uses Lemma~\ref{lem: abstract E2 degenerate}), we construct a double complex $\cT_{\Omega',\fh}^{\bullet,\bullet}$ which comes with a quasi-map between double complex $\cT_{\Omega',\fh}^{\bullet,\bullet}\dashrightarrow\cT_{I_0,I_1,\fh}^{\bullet,\bullet}$ whose associated map between spectral sequences $E_{\bullet,\Omega',\fh}^{\bullet,\bullet}\rightarrow E_{\bullet,I_0,I_1,\fh}^{\bullet,\bullet}$ satisfy the following conditions.
\begin{itemize}
\item We have $E_{1,\Omega',\fh}^{-\ell',k'}=0$ for each $(-\ell',k')$ with $k'-\ell'>h$.
\item The image of $E_{1,\Omega',\fh}^{-\ell,\ell+h}\rightarrow E_{1,I_0,I_1,\fh}^{-\ell,\ell+h}$ has image $x_{\Omega',\fh}$.
\end{itemize}
We write $v'\defeq v\setminus \mathbf{Log}_{\emptyset}^{\infty}$ for short. We write $J_{\fh,+}=[j_{\fh,+},j_{\fh,+}']$ for the unique subinterval of $J_{\fh}\cup J'$ that contains $J_{\fh}$, and define a new interval $J_{\fh,\dagger}=[j_{\fh,\dagger},j_{\fh,\dagger}']$ that fits into $J_{\fh}\subseteq J_{\fh,\dagger}\subseteq J_{\fh,+}$ by $j_{\fh,\dagger}=\min\{j_{\fh,+}+1,j_{\fh}\}$ and $j_{\fh,\dagger}'=\max\{j_{\fh,+}'-1,j_{\fh}'\}$. For each $1\leq d\leq n-\ell$ with $(J')^{d}\cap J_{\fh}=\emptyset$ (or equivalently $(J')^{d}\cap J_{\fh,+}=\emptyset$) and $\Lambda'_{d}\neq\emptyset$, we define $\fh_{\Omega',d,\dagger}\defeq \fh_{(J')^{d}\setminus\{i_{\Theta',d}-1\}}\subseteq \fh_{\Omega',d}\defeq \fh_{(J')^{d}}$. We thus obtain the following $E$-Lie subalgebra
\[\fh_{\Omega',\fh,\dagger}\defeq\fh_{J_{\fh,\dagger}}\times\prod_{d,\Lambda'_{d}\neq\emptyset,(J')^{d}\cap J_{\fh}=\emptyset}\fh_{\Omega',d,\dagger}\subseteq \fh_{J'\cup J_{\fh}}=\fh_{J_{\fh,+}}\times\prod_{d,\Lambda'_{d}\neq\emptyset,(J')^{d}\cap J_{\fh}=\emptyset}\fh_{\Omega',d}.\]
Now we are ready to define $\cT_{\Omega',\fh}^{\bullet,\bullet}$ as the double complex $\cT_{\Sigma}^{\bullet,\bullet}$ (see (\ref{equ: general double complex})) with $M^{\bullet}_{I}$ given by
\[\mathrm{Hom}_{E}(\tau^{M}(B_{\bullet}^{I}),v'\otimes_EC^{\bullet}(L_{I}'/(L_{I}'\cap Z_{I_{v'}}),\fl_{I}\cap\fh_{\Omega',\dagger},1_{L_{I}'/(L_{I}'\cap Z_{I_{v'}})}))\]
for each $I_0\subseteq I\subseteq J'$, and given by zero otherwise. Here we understand $\fl_{I}\cap\fh_{\Omega',\dagger}$ as an $E$-Lie subalgebra of $\fl_{I}/\fz_{I_{v'}}$ for each $I_0\subseteq I\subseteq J'$.
The verification of desired conditions satisfied by $E_{\bullet,\Omega',\fh}^{\bullet,\bullet}\rightarrow E_{\bullet,I_0,I_1,\fh}^{\bullet,\bullet}$ uses parallel argument as that of Lemma~\ref{lem: Omega seq lift}, based on Proposition~\ref{prop: relative interval embedding}, Remark~\ref{rem: h tuple} and Corollary~\ref{cor: relative Levi restriction}.
\end{proof}

\begin{lem}\label{lem: relative E2 basis}
Let $I_0\subseteq I_1\subseteq \Delta$ and $\fh=\fh_{J_{\fh}}$ for some interval $I_0\subseteq J_{\fh}\subseteq \Delta$. 
Then we have the following results.
\begin{enumerate}[label=(\roman*)]
\item \label{it: relative E2 basis 1} For each $\#I_0\leq \ell\leq \#I_1$, the $E$-vector space $E_{2,I_0,I_1,\fh}^{-\ell,\ell+\#I_1-2\#I_0}$ admits a basis of the form
\begin{equation}\label{equ: relative E2 basis}
\{\overline{x}_{\Omega,\fh}\}_{\Omega\in \Psi_{I_0,I_1,\fh,\ell}}
\end{equation}
\item \label{it: relative E2 basis 2} For each $\#I_0\leq \ell\leq \#I_1$ and $M\subseteq I_1\setminus I_0$, the $E$-vector space $\mathrm{gr}^{M}(E_{2,I_0,I_1,\fh}^{-\ell,\ell+\#I_1-2\#I_0})$ admits a basis of the form
\begin{equation}\label{equ: relative E2 basis M}
\{\overline{x}_{\Omega,\fh}\}_{\Omega\in \mathrm{gr}^{M}(\Psi_{I_0,I_1,\fh,\ell})}.
\end{equation}
\end{enumerate}
\end{lem}
\begin{proof}
For each bi-degree $(-\ell,k)$ and $v\subseteq\mathbf{Log}_{I_0}$ with $v\setminus \mathbf{Log}_{I_0}^{\infty}\subseteq \mathbf{Log}_{I_0\cup J_{\fh}}$, we define $E_{1,I_0,I_1,\fh,v}^{-\ell,k}$ as the $E$-span of $x_{\Theta}$ with $\Theta=(v,I,\un{k},\un{\Lambda})$ running through all $\fh$-tuples (see Remark~\ref{rem: h tuple}) that satisfy $\#I=\ell$ and $|\un{k}|=k$ (namely with bi-degree $(-\ell,k)$). Then we have
\[E_{1,I_0,I_1,\fh}^{-\ell,k}=\bigoplus_{v}E_{1,I_0,I_1,\fh,v}^{-\ell,k}\]
for each $(-\ell,k)$, with the differential map $d_{1,I_0,I_1,\fh}^{-\ell,k}$ also being the direct sum of differential maps
\[d_{1,I_0,I_1,\fh,v}^{-\ell,k}: E_{1,I_0,I_1,\fh,v}^{-\ell,k}\rightarrow E_{1,I_0,I_1,\fh,v}^{-\ell+1,k},\]
and thus induces a decomposition
\[E_{2,I_0,I_1,\fh}^{-\ell,k}=\bigoplus_{v}E_{2,I_0,I_1,\fh,v}^{-\ell,k}\]
with $E_{2,I_0,I_1,\fh,v}^{-\ell,k}\defeq \mathrm{ker}(d_{1,I_0,I_1,\fh}^{-\ell,k})/\mathrm{im}(d_{1,I_0,I_1,\fh}^{-\ell-1,k})$ for each $(-\ell,k)$. To prove \ref{it: relative E2 basis 1} and \ref{it: relative E2 basis 2}, it suffices to fix a choice of $v$ and show that $E_{2,I_0,I_1,\fh,v}^{-\ell,\ell+h}$ admits basis of the form
\begin{equation}\label{equ: relative E2 basis v}
\{\overline{x}_{\Omega,\fh}\}_{\Omega\in \Psi_{I_0,I_1,\fh,v,\ell}}
\end{equation}
where $\Psi_{I_0,I_1,\fh,v,\ell}\subseteq \Psi_{I_0,I_1,\fh,\ell}$ is the subset consisting of those equivalence classes with our fixed $v$.
We write $J_{\fh}=[j_{\fh},j_{\fh}']$ as usual. 
For each $\fh$-tuple $\Theta$, we use the symbol $\Omega$ for its associated equivalence class.
Following the proof of Lemma~\ref{lem: relative E2 atom}, we see that our fixed choice of $v$ actually uniquely determines an integer $1\leq s_0\leq r_{I_{v}\cap I_1}$ (which is independent of the choice of $\Theta$ or its associated $\Omega$ with this fixed $v$).
Given a tuple $\Theta=(v,I,\un{k},\un{\Lambda})$, $1\leq d\leq r_{I}$ and $i\in I^{d}$, we define $q_{i}^{-}(\Theta)$ as the unique tuple (if exists) of the form $(v,I\setminus\{i\},\un{k}',\un{\Lambda}')$ with $\Lambda'_{d'}=\Lambda_{d'}$ if $d'<d$, $\Lambda'_{d}=\emptyset$, and $\Lambda'_{d'}=\Lambda_{d-1}$ if $d'>d$.
We write $\Sigma^{-\ell,k}$ for the set of all $\fh$-tuples with bi-degree $(-\ell,k)$ and our fixed $v$.
For each $\Theta\in \Sigma^{-\ell,k}$, we have a unique $i_1\in I$ (if exists) such that exactly one of the following holds
\begin{itemize}
\item $i_1<j_{\fh}$ and $i_1$ is the minimal integer such that $q_{i_1}^{-}(\Theta)$ is defined;
\item $i_1\geq j_{\fh}$, $q_{i}^{-}(\Theta)$ is not defined for each $i<j_{\fh}$, and $i_1$ is the maximal integer such that $p_{i_1}^{-}(\Theta)$ is defined.
\end{itemize}
For each $\Theta\in \Sigma^{-\ell,k}$, we have a unique $1\leq s_1\leq r_{I_{v}\cap I_1}$ and a unique $r_{\Omega}^{s_1-1}+1\leq d_1\leq r_{\Omega}^{s_1}$ (if exists) such that $\Lambda_{d_1}=\emptyset$ and exactly one of the following holds
\begin{itemize}
\item $s_1\leq s_0$, $i_{\Theta,d_1}<j_{\fh}$ and $d_1$ is the minimal such integer that $d_1\neq r_{\Omega}^{s_1}$;
\item $s_1\geq s_0$, $i_{\Theta,d_1-1}\geq j_{\fh}$ and $d_1$ is the maximal such integer that $d_1\neq r_{\Omega}^{s_1-1}+1$, with $\Lambda_{d}=\emptyset$ for some $1\leq d\leq r_{I}$ satisfying $i_{\Theta,d}<j_{\fh}$ only if $d=r_{\Omega}^{s}$ for some $1\leq s<s_0$.
\end{itemize}

Similar to Definition~\ref{def: std element}, we say that a $\fh$-tuple $\Theta=(v,I,\un{k},\un{\Lambda})$ is \emph{$(I_0,I_1,\fh)$-standard} if $d_1$ exists and either $i_1$ does not exist, or $i_{\Theta,d_1}<j_{\fh}$ with $i_{\Theta,d_1}<i_1$, or $i_{\Theta,d_1-1}\geq j_{\fh}$ with $j_{\fh}\leq i_1<i_{\Theta,d_1-1}$.
We define $\tld{\Theta}\defeq p_{i_{\Theta,d_1}}^{+}(\Theta)$ if $i_{\Theta,d_1}<j_{\fh}$, and $\tld{\Theta}\defeq p_{i_{\Theta,d_1-1}}^{+}(\Theta)$ if $i_{\Theta,d_1-1}\geq j_{\fh}$.
We write $\Sigma_{\rm{s}}^{-\ell,k}$ for the set of $(I_0,I_1,\fh)$-standard $\fh$-tuples with bi-degree $(-\ell,k)$, and set
\[\tld{\Sigma}_{\rm{s}}^{-\ell,k}\defeq \{\tld{\Theta}\mid \Theta\in \Sigma_{\rm{s}}^{-\ell+1,k}\}.\]
Similar to Definition~\ref{def: atom element}, we say a $\fh$-tuple $\Theta=(v,I,\un{k},\un{\Lambda})$ is \emph{maximally $(I_0,I_1,\fh)$-atomic} if neither $i_1$ nor $d_1$ exists.
We write $\Sigma_{\rm{a}}^{-\ell,k}$ for the set of maximally $(I_0,I_1,\fh)$-atomic $\fh$-tuples with bi-degree $(-\ell,k)$.
We divide the rest of the proof into the following steps.

\textbf{Step $1$}: We prove that
\begin{equation}\label{equ: h tuple classification}
\Sigma^{-\ell,k}=\Sigma_{\rm{s}}^{-\ell,k}\sqcup\tld{\Sigma}_{\rm{s}}^{-\ell,k}\sqcup \Sigma_{\rm{a}}^{-\ell,k}
\end{equation}
for each bi-degree $(-\ell,k)$.\\
Let $\Theta\in \Sigma^{-\ell,k}$ be a $\fh$-tuple. We assume that $\Theta\notin \Sigma_{\rm{a}}^{-\ell,k}$ (namely at least one of $i_1$ and $d_1$ exists) and show that either $\Theta\in \Sigma_{\rm{s}}^{-\ell,k}$ or $\Theta\in \tld{\Sigma}_{\rm{s}}^{-\ell,k}$ with exactly one of them holds.
We have the following disjoint possibilities.
\begin{itemize}
\item If $i_1<j_{\fh}$ and either $d_1$ does not exist or $i_1<i_{\Theta,d_1}$, then $\Theta^{\dagger}\defeq q_{i_1}^{-}(\Theta)$ is $(I_0,I_1,\fh)$-standard. Namely we have $\Theta^{\dagger}\in\Sigma_{\rm{s}}^{-\ell+1,k}$ and thus $\Theta\in\tld{\Sigma}_{\rm{s}}^{-\ell,k}$.
\item If $i_1\geq j_{\fh}$ and either $d_1$ does not exist or $i_1>i_{\Theta,d_1-1}\geq j_{\fh}$, then $\Theta^{\dagger}\defeq p_{i_1}^{-}(\Theta)$ is $(I_0,I_1,\fh)$-standard. Namely we have $\Theta^{\dagger}\in\Sigma_{\rm{s}}^{-\ell+1,k}$ and thus $\Theta\in\tld{\Sigma}_{\rm{s}}^{-\ell,k}$.
\item If $i_{\Theta,d_1}<j_{\fh}$ and either $i_1$ does not exist or $i_1>i_{\Theta,d_1}$, then $\Theta$ is $(I_0,I_1,\fh)$-standard, namely we have $\Theta\in\Sigma_{\rm{s}}^{-\ell,k}$.
\item If $i_{\Theta,d_1-1}\geq j_{\fh}$ and either $i_1$ does not exist or $j_{\fh}\leq i_1<i_{\Theta,d_1-1}$, then $\Theta$ is $(I_0,I_1,\fh)$-standard, namely we have $\Theta\in\Sigma_{\rm{s}}^{-\ell,k}$.
\end{itemize}

\textbf{Step $2$}: Let $(-\ell,k)$ be a bi-degree. We prove that the following subset
\begin{equation}\label{equ: h image of differential}
\{d_{1,I_0,I_1,\fh,v}^{-\ell,k}(x_{\tld{\Theta}})\}_{\Theta\in \Sigma_{\rm{s}}^{-\ell+1,k}}\subseteq E_{1,I_0,I_1,\fh,v}^{-\ell+1,k}
\end{equation}
is linearly independent.\\
Recall that we can uniquely associate an integer $d_1$ with each $\Theta=(v,I,\un{k},\un{\Lambda})\in \Sigma_{\rm{s}}^{-\ell+1,k}$.
We consider an arbitrary $\Theta'=(v,I',\un{k},\un{\Lambda})\in \Sigma_{\rm{s}}^{-\ell+1,k}$ with associated integer $d_1'$ such that $\Theta'\neq \Theta$ and
\begin{equation}\label{equ: h std coefficient}
c_{\Theta'}(d_{1,I_0,I_1,\fh,v}^{-\ell,k}(x_{\tld{\Theta}}))\neq 0,
\end{equation}
namely there exists $i\in (I_{v}\cap I_1)\setminus I'$ such that $\tld{\Theta}=p_{i}^{+}(\Theta')$.
We have the following possibilities.
\begin{itemize}
\item If $i<i_{\Theta,d_1}<j_{\fh}$ (with $i_{\Theta,d_1-1}+1=i_{\Theta,d_1}$), then $i\in I^{d}$ for some $d<d_1$, and $\tld{\Theta}=p_{i}^{+}(\Theta')$ forces $q_{i_{\Theta,d-1}+1}^{-}(\Theta)$ to be defined and thus contradicts $\Theta$ being $(I_0,I_1,\fh)$-standard.
\item If $i_{\Theta,d_1}<j_{\fh}$ and $i>i_{\Theta,d_1}$, then we must have $d_1'>d_1$ and $(I')^{d}=(I)^{d}$ for each $d<d_1$.
\item If $i>i_{\Theta,d_1-1}\geq j_{\fh}$, then we must have $d_1'\geq d_1$, with $(I')^{d}=I^{d}$ for each $d>d_1'+1$ and $(I')^{d_1'+1}\subsetneq I^{d_1'+1}$.
\item If $i_{\Theta,d_1-1}\geq j_{\fh}$ and $i<i_{\Theta,d_1-1}$, then we have $d_1'<d_1$ and $(I')^{d}=(I)^{d}$ for each $d>d_1$.
\end{itemize}
Similar to Definition~\ref{def: simple partial order} and the proof of Lemma~\ref{lem: independent elements}, we can define a partial-order $\prec_{\fh}$ on $\Sigma_{\rm{s}}^{-\ell+1,k}$ so that (\ref{equ: h std coefficient}) holds only if $\Theta\prec_{\fh}\Theta'$. A parallel argument as in the last paragraph of the proof of Lemma~\ref{lem: independent elements} shows that (\ref{equ: h image of differential}) is linearly independent.

\textbf{Step $3$}: We prove that $E_{2,I_0,I_1,\fh,v}^{-\ell,\ell+h}$ admits a basis of the form (\ref{equ: relative E2 basis v}).\\
Parallel to the proof of Lemma~\ref{lem: existence of atom}, we can check that $\Sigma_{\rm{a}}^{-\ell,\ell+h}$ consists of exactly those $\fh$-tuples $\Theta'$ as constructed in the proof of Lemma~\ref{lem: relative E2 atom} (for our fixed $v$) with each $\Theta'$ being the unique $\fh$-tuple in its associated equivalence class $\Omega'$ and $\varepsilon(\Theta')x_{\Theta'}=x_{\Omega',\fh}$ by definition. Consequently, we know from Lemma~\ref{lem: relative E2 atom} that the subset
\[\{x_{\Theta'}\}_{\Theta'\in \Sigma_{\rm{a}}^{-\ell,\ell+h}}\subseteq E_{1,I_0,I_1,\fh,v}^{-\ell,\ell+h}\]
is killed by $d_{1,I_0,I_1,\fh,v}^{-\ell,\ell+h}$ and induces a linearly independent subset of $E_{2,I_0,I_1,\fh,v}^{-\ell,\ell+h}$ with
\begin{equation}\label{equ: relative E2 basis lower bound}
\Dim_E E_{2,I_0,I_1,\fh,v}^{-\ell,\ell+h}\geq \#\Sigma_{\rm{a}}^{-\ell,\ell+h}.
\end{equation}
In another direction, we know from \textbf{Step $2$} that
\[\Dim_E\mathrm{im}(d_{1,I_0,I_1,\fh,v}^{-\ell'-1,\ell+h})\geq \#\Sigma_{\rm{s}}^{-\ell',\ell+h}\]
for each $\ell'\in\{\ell,\ell-1\}$, which together with \textbf{Step $1$} implies that
\begin{multline*}
\Dim_EE_{2,I_0,I_1,\fh,v}^{-\ell,\ell+h}=\Dim_EE_{1,I_0,I_1,\fh,v}^{-\ell,\ell+h}-\Dim_E\mathrm{im}(d_{1,I_0,I_1,\fh,v}^{-\ell-1,\ell+h})-\Dim_E\mathrm{im}(d_{1,I_0,I_1,\fh,v}^{-\ell,\ell+h})\\
\leq \#\Sigma^{-\ell,\ell+h}-\#\Sigma_{\rm{s}}^{-\ell,\ell+h}-\#\Sigma_{\rm{s}}^{-\ell+1,\ell+h}
=\#\Sigma^{-\ell,\ell+h}-\#\Sigma_{\rm{s}}^{-\ell,\ell+h}-\#\tld{\Sigma}_{\rm{s}}^{-\ell,\ell+h}=\#\Sigma_{\rm{a}}^{-\ell,\ell+h}.
\end{multline*}
This together with (\ref{equ: relative E2 basis lower bound}) and the discussion above it forces (\ref{equ: relative E2 basis v}) to be a basis of $E_{2,I_0,I_1,\fh,v}^{-\ell,\ell+h}$.
\end{proof}

Let $I_0\subseteq I_1\subseteq \Delta$ and $\fh=\fh_{J_{\fh}}$ for some interval $I_0\subseteq J_{\fh}\subseteq \Delta$.
We continue to write $h=\#I_1-2\#I_0$.
The map between double complex $\cT_{I_0,I_1,\fh}^{\bullet,\bullet}\rightarrow \cT_{I_0,I_1}^{\bullet,\bullet}$ induces the following map
\begin{equation}\label{equ: relative bottom embedding}
H^{h}(\mathrm{Tot}(\cT_{I_0,I_1,\fh}^{\bullet,\bullet}))\rightarrow H^{h}(\mathrm{Tot}(\cT_{I_0,I_1}^{\bullet,\bullet}))
\end{equation}
\begin{prop}\label{prop: relative bottom embedding}
Let $I_0\subseteq I_1\subseteq \Delta$ and $\fh$ be as above. Then we have the following results.
\begin{enumerate}[label=(\roman*)]
\item \label{it: relative bottom embedding 1} The map (\ref{equ: relative bottom embedding}) is an embedding which is strict with respect to the canonical filtration on its source and its target, with the induced map between $\mathrm{gr}^{-\ell}$ being the embedding
    \begin{equation}\label{equ: relative bottom embedding 1}
    E_{2,I_0,I_1,\fh}^{-\ell,\ell+h}\hookrightarrow E_{2,I_0,I_1}^{-\ell,\ell+h}
    \end{equation}
    for each $\#I_0\leq\ell\leq \#I_1$.
\item \label{it: relative bottom embedding 2} For each $M\subseteq I_1\setminus I_0$, we have the following commutative diagram
\[
\xymatrix{
H^{h}(\mathrm{Tot}(\tau^{M}(\cT_{I_0,I_1,\fh}^{\bullet,\bullet}))) \ar@{^{(}->}[r] \ar@{^{(}->}[d]& H^{h}(\mathrm{Tot}(\tau^{M}(\cT_{I_0,I_1}^{\bullet,\bullet}))) \ar@{^{(}->}[d]\\
H^{h}(\mathrm{Tot}(\cT_{I_0,I_1,\fh}^{\bullet,\bullet})) \ar@{^{(}->}[r] & H^{h}(\mathrm{Tot}(\cT_{I_0,I_1}^{\bullet,\bullet}))
}
\]
with all maps being embeddings, from which we obtain compatible filtration $\{\tau^{M}(-)\}_{M\subseteq I_1\setminus I_0}$ on BHS of (\ref{equ: relative bottom embedding}). The map (\ref{equ: relative bottom embedding}) is strict with respect to these $\tau^{M}(-)$ filtration, with the induced map between $\mathrm{gr}^{M}(-)$ being the embedding
    \begin{equation}\label{equ: relative bottom embedding grade M}
    H^{h}(\mathrm{Tot}(\mathrm{gr}^{M}(\cT_{I_0,I_1,\fh}^{\bullet,\bullet})))\hookrightarrow H^{h}(\mathrm{Tot}(\mathrm{gr}^{M}(\cT_{I_0,I_1}^{\bullet,\bullet})))
    \end{equation}
\end{enumerate}
\end{prop}
\begin{proof}
We prove \ref{it: relative bottom embedding 1}.\\
Recall from Lemma~\ref{lem: relative E2 basis} that (\ref{equ: relative E2 basis}) forms a basis of the source of (\ref{equ: relative bottom embedding 1}). Now that (\ref{equ: relative bottom embedding 1}) sends (\ref{equ: relative E2 basis}) to the linearly independent subset
\[\{\overline{x}_{\Omega}\}_{\Omega\in\Psi_{I_0,I_1,\fh,\ell}}\subseteq \{\overline{x}_{\Omega}\}_{\Omega\in\Psi_{I_0,I_1,\ell}}\]
by Lemma~\ref{lem: relative E2 atom}, we conclude that (\ref{equ: relative bottom embedding 1}) is an embedding. Now that (\ref{equ: relative E2 basis}) forms a basis of the source of (\ref{equ: relative bottom embedding 1}), we deduce from Lemma~\ref{lem: relative E2 atom degeneracy} that
\[E_{2,I_0,I_1,\fh}^{-\ell,\ell+h}=E_{\infty,I_0,I_1,\fh}^{-\ell,\ell+h},\]
which together with $E_{2,I_0,I_1}^{-\ell,\ell+h}=E_{\infty,I_0,I_1}^{-\ell,\ell+h}$ from Proposition~\ref{prop: bottom deg degeneracy} identifies the map (\ref{equ: relative bottom embedding 1}) with the map
\[
E_{\infty,I_0,I_1,\fh}^{-\ell,\ell+h}\rightarrow E_{\infty,I_0,I_1}^{-\ell,\ell+h}
\]
which is now known to be an embedding. Now that the map (\ref{equ: relative bottom embedding}) between filtered $E$-vector spaces induces an embedding between $\mathrm{gr}^{-\ell}$ for each $\ell$, we conclude that (\ref{equ: relative bottom embedding}) is an embedding which is strict with respect to the canonical filtration on its source and its target.

We prove \ref{it: relative bottom embedding 2}.\\
A parallel argument as in \ref{it: relative bottom embedding 1} using the embedding
\[\mathrm{gr}^{M}(E_{2,I_0,I_1,\fh}^{-\ell,\ell+h})\hookrightarrow \mathrm{gr}^{M}(E_{2,I_0,I_1}^{-\ell,\ell+h})\]
for each $\ell$ shows that (\ref{equ: relative bottom embedding}) is an embedding which is strict with respect to the canonical filtration $\mathrm{Fil}^{-\ell}(-)$ on its source and its target. The prove of the rest of claims of \ref{it: relative bottom embedding 2} follows similar argument as in Proposition~\ref{prop: bottom deg graded}, using Lemma~\ref{lem: abstract E2 injection} and Lemma~\ref{lem: abstract E2 seq}.
\end{proof}

\begin{lem}\label{lem: CE relative j vanishing}
Let $I_0\subseteq \Delta$, $j\in\Delta\setminus I_0$ and $x\in\Gamma$ with $\mathrm{Supp}(x)\cap(I_0\sqcup\{j\})=\emptyset$. Then we have
\begin{equation}\label{equ: CE relative j vanishing}
H^{\ell(x)+1-\#I_0}(\mathrm{Tot}(\tld{\mathrm{CE}}_{x,I_0,I_0\sqcup\{j\},\fh_{\{j\}}}^{\bullet,\bullet}))=0
\end{equation}
\end{lem}
\begin{proof}
We write $I_0'\defeq I_0\sqcup\{j\}$ and $\fh'\defeq \fh_{\{j\}}$ for short, and $\tld{E}_{\bullet,x,I_0,I_0',\fh'}^{\bullet,\bullet}$ for the spectral sequence associated with the double complex $\tld{\mathrm{CE}}_{x,I_0,I_0',\fh'}^{\bullet,\bullet}$. In particular, we have the following short exact sequence
\begin{equation}\label{equ: CE relative j seq}
0\rightarrow \tld{E}_{\infty,x,I_0,I_0',\fh'}^{-\#I_0,1+\ell(x)} \rightarrow H^{1-\#I_0}(\mathrm{Tot}(\tld{\mathrm{CE}}_{I_0,I_0',\fh'}^{\bullet,\bullet}))\rightarrow \tld{E}_{\infty,x,I_0,I_0',\fh'}^{-\#I_0',2+\ell(x)} \rightarrow 0.
\end{equation}
Thanks to (\ref{equ: Lie parabolic restriction x}), the differential map
\[\tld{E}_{1,x,I_0,I_0',\fh'}^{-\#I_0',k}\rightarrow \tld{E}_{1,x,I_0,I_0',\fh'}^{-\#I_0,k}\]
is given (up to a sign) by the restriction map
\begin{equation}\label{equ: CE relative j restriction}
H^{k-\ell(x)}(\fl_{I_0'},\fh',1_{\fl_{I_0'}}) \rightarrow H^{k-\ell(x)}(\fl_{I_0},\fh'\cap\ft,1_{\fl_{I_0}})
\end{equation}
for each $k\geq 0$. Note that the inclusions between pairs $(\fz_{I_0'},0)\subseteq(\fl_{I_0},\fh'\cap\ft)\subseteq (\fl_{I_0'},\fh')$ induce a commutative diagram that maps (\ref{equ: CE relative j restriction}) to the identity map of $H^{k-\ell(x)}(\fz_{I_0'},1_{\fz_{I_0'}})$. Now that $H^{\ell}(\fm,1_{\fm})=0$ for $\ell=1,2$ for any semi-simple $E$-Lie algebra $\fm$, we know that the natural restriction maps
\[H^{k-\ell(x)}(\fl_{I_0'},\fh',1_{\fl_{I_0'}})\rightarrow H^{k-\ell(x)}(\fz_{I_0'},1_{\fz_{I_0'}})\]
and
\[H^{k-\ell(x)}(\fl_{I_0},\fh'\cap\ft,1_{\fl_{I_0}})\rightarrow H^{k-\ell(x)}(\fz_{I_0'},1_{\fz_{I_0'}})\]
are isomorphisms for $k\in\{1+\ell(x),2+\ell(x)\}$. In other words, the map (\ref{equ: CE relative j restriction}) is an isomorphism for $k\in\{1+\ell(x),2+\ell(x)\}$ which forces
\[\tld{E}_{2,x,I_0,I_0',\fh'}^{-\#I_0,1+\ell(x)}=0=\tld{E}_{2,x,I_0,I_0',\fh'}^{-\#I_0',2+\ell(x)}\]
and of course
\[\tld{E}_{\infty,x,I_0,I_0',\fh'}^{-\#I_0,1+\ell(x)}=0=\tld{E}_{\infty,x,I_0,I_0',\fh'}^{-\#I_0',2+\ell(x)}.\]
This together with (\ref{equ: CE relative j seq}) gives (\ref{equ: CE relative j vanishing}).
\end{proof}

Let $j\in I_1\setminus I_0$ with $\fh_{\{j\}}\subseteq \fh$ (or equivalently $j\in J_{\fh}$).
Let $x\in\Gamma$ with $\mathrm{Supp}(x)=(I_1\setminus I_0)\setminus\{j\}$ and thus $\ell(x)+1=\#I_1\setminus I_0$.
Recall from (\ref{equ: Tits Verma to Levi}) and (\ref{equ: double x w diagram}) that we have the following commutative diagram of maps between double complex
\begin{equation}\label{equ: relative CE double diagram}
\xymatrix{
\tld{\mathrm{CE}}_{x,I_0,I_1,\fh_{\{j\}}}^{\bullet,\bullet}  \ar[d]  & \tld{\mathrm{CE}}_{I_0,I_1,\fh}^{\bullet,\bullet}  \ar[l] \ar[r] \ar[d] & \mathrm{CE}_{I_0,I_1,\fh}^{\bullet,\bullet} \ar[d] \\
\tld{\mathrm{CE}}_{x,I_0,I_1}^{\bullet,\bullet} & \tld{\mathrm{CE}}_{I_0,I_1}^{\bullet,\bullet}  \ar[l] \ar[r]  & \mathrm{CE}_{I_0,I_1}^{\bullet,\bullet}
}
\end{equation}
with the horizontal maps from the middle column to right column inducing isomorphism on the first page by the description of (\ref{equ: Tits Verma to Levi}). 
The commutative diagram (\ref{equ: relative CE double diagram}) together with (\ref{equ: Lie general Tits resolution}) and Proposition~\ref{prop: Lie St Ext w decomposition} induces the following commutative diagram (with $h\defeq \#I_1-2\#I_0$ for short)
\begin{equation}\label{equ: relative CE double diagram total}
\xymatrix{
H^{h}(\mathrm{Tot}(\tld{\mathrm{CE}}_{x,I_0,I_1,\fh_{\{j\}}}^{\bullet,\bullet}))  \ar[d]  & H^{h}(\mathrm{Tot}(\tld{\mathrm{CE}}_{I_0,I_1,\fh}^{\bullet,\bullet}))  \ar[l] \ar^{\sim}[r] \ar[d] & H^{h}(\mathrm{Tot}(\mathrm{CE}_{I_0,I_1,\fh}^{\bullet,\bullet})) \ar[d] \\
H^{h}(\mathrm{Tot}(\tld{\mathrm{CE}}_{x,I_0,I_1}^{\bullet,\bullet})) \ar@{=}[d] & H^{h}(\mathrm{Tot}(\tld{\mathrm{CE}}_{I_0,I_1}^{\bullet,\bullet}))  \ar[l] \ar^{\sim}[r] \ar@{=}[d] & H^{h}(\mathrm{Tot}(\mathrm{CE}_{I_0,I_1}^{\bullet,\bullet}))\\
\mathrm{Ext}_{U(\fg)}^{h}(\mathfrak{c}^{x}_{I_0,I_1},L(1)) \ar^{\wr}[d] & \mathrm{Ext}_{U(\fg)}^{h}(\mathfrak{c}_{I_0,I_1},L(1)) \ar[l] \ar^{\wr}[d]&\\
\mathrm{Ext}_{U(\fg)}^{\#I_1\setminus I_0}(U(\fg)\otimes_{U(\fp_{I_1})}\mathfrak{v}_{x,I_0,I_1},L(1)) \ar^{\wr}[d] & \mathrm{Ext}_{U(\fg)}^{\#I_1\setminus I_0}(U(\fg)\otimes_{U(\fp_{I_1})}\mathfrak{v}_{I_0,I_1},L(1)) \ar[l] \ar^{\wr}[d] &\\
\mathrm{Ext}_{U(\fg)}^{\#I_1\setminus I_0}(M^{I_0}(s_jx),L(1)) \ar@{=}[d] & \bigoplus_{w\in\Gamma_{I_1\setminus I_0}}\mathrm{Ext}_{U(\fg)}^{\#I_1\setminus I_0}(M^{I_0}(w),L(1)) \ar@{->>}[l] \ar@{=}[d]\\
\mathfrak{e}_{I_0,s_jx} & \bigoplus_{w\in\Gamma_{I_1\setminus I_0}}\mathfrak{e}_{I_0,w} \ar@{->>}[l]
}.
\end{equation}
Here all the vertical maps below the second row are isomorphisms (see Proposition~\ref{prop: Lie St Ext w decomposition}), and the horizontal map on the bottom row is given by the evident projection.
In particular, we recall from (\ref{equ: relative CE double diagram total}) the following map
\begin{equation}\label{equ: CE relative h0 transfer}
H^{h}(\mathrm{Tot}(\mathrm{CE}_{I_0,I_1,\fh}^{\bullet,\bullet}))\rightarrow H^{h}(\mathrm{Tot}(\mathrm{CE}_{I_0,I_1}^{\bullet,\bullet})).
\end{equation}
as well as the isomorphism
\begin{equation}\label{equ: CE w decomposition h0}
H^{h}(\mathrm{Tot}(\mathrm{CE}_{I_0,I_1}^{\bullet,\bullet}))\buildrel\sim\over\longrightarrow \bigoplus_{w\in\Gamma_{I_1\setminus I_0}}\mathfrak{e}_{I_0,w}
\end{equation}
Recall from (\ref{equ: coxeter left}) that
\[\Gamma_{I_1\setminus I_0}(J_{\fh})=\{w\in\Gamma_{I_1\setminus I_0}, D_L(w)\cap J_{\fh}=\emptyset\}\subseteq \Gamma_{I_1\setminus I_0}\]
\begin{lem}\label{lem: CE relative upper bound}
Let $I_0\subseteq I_1\subseteq \Delta$ and $\fh$ be as above with $h=\#I_1-2\#I_0$. Then the image of (\ref{equ: CE relative h0 transfer}) is contained in
\begin{equation}\label{equ: CE relative upper bound}
\bigoplus_{w\in\Gamma_{I_1\setminus I_0}(J_{\fh})}\mathfrak{e}_{I_0,w}
\end{equation}
under the isomorphism (\ref{equ: CE w decomposition h0}).
\end{lem}
\begin{proof}
Let $j\in J_{\fh}$ and $y\in\Gamma_{I_1\setminus I_0}$ with $j\in D_L(w)$. Then $x\defeq s_jy<y$ satisfies $\mathrm{Supp}(x)=(I_1\setminus I_0)\setminus\{j\}$. Thanks to (\ref{equ: relative CE double diagram total}), we know that the composition of (\ref{equ: CE relative h0 transfer}) and (\ref{equ: CE w decomposition h0}) with the projection
\[\bigoplus_{w\in\Gamma_{I_1\setminus I_0}}\mathfrak{e}_{I_0,w}\twoheadrightarrow \mathfrak{e}_{I_0,y}\]
factors through $H^{h}(\mathrm{Tot}(\tld{\mathrm{CE}}_{x,I_0,I_1,\fh_{\{j\}}}^{\bullet,\bullet}))$ which is zero by Lemma~\ref{lem: CE relative j vanishing}. This forces the image of (\ref{equ: CE relative h0 transfer}) to be contained in (\ref{equ: CE relative upper bound}), under the isomorphism (\ref{equ: CE w decomposition h0}).
\end{proof}

\begin{prop}\label{prop: grade M relative image}
Let $I_0\subseteq I_1\subseteq \Delta$ and $\fh$ be as above with $h=\#I_1-2\#I_0$. 
Assume that $I_0\subseteq J_{\fh}$. Then for each $M\subseteq I_1\setminus I_0$, the composition of the following maps
\begin{equation}\label{equ: grade M relative image maps}
H^{h}(\mathrm{Tot}(\mathrm{gr}^{M}(\cT_{I_0,I_1,\fh}^{\bullet,\bullet})))\hookrightarrow H^{h}(\mathrm{Tot}(\mathrm{gr}^{M}(\cT_{I_0,I_1}^{\bullet,\bullet})))
\buildrel\sim\over\longrightarrow \wedge^{\#M}\Hom(Z_{\Delta\setminus M}^{\dagger},E)\otimes_E\big(\bigoplus_{w\in\Gamma_{(I_1\setminus I_0)\setminus M}}\mathfrak{e}_{I_0,w}\big)
\end{equation}
equals
\begin{equation}\label{equ: grade M relative image}
\wedge^{\#M}\Hom(Z_{\Delta\setminus M}^{\dagger},E)\otimes_E\big(\bigoplus_{w\in\Gamma_{(I_1\setminus I_0)\setminus J}(J_{\fh})}\mathfrak{e}_{I_0,w}\big).
\end{equation}
\end{prop}
\begin{proof}
We write $\mathbf{E}_{M}^{\infty}\defeq \wedge^{\#M}\Hom(Z_{\Delta\setminus M}^{\dagger},E)$ for short.
We have the following commutative diagram
\begin{equation}\label{equ: CE relative image diagram}
\xymatrix{
H^{h}(\mathrm{Tot}(\mathrm{gr}^{M}(\cT_{I_0,I_1,\fh}^{\bullet,\bullet}))) \ar^{\sim}[r] \ar[d] & H^{h-\#M}(\mathrm{Tot}(\mathrm{CE}_{I_0,I_1\setminus M,\fh}^{\bullet,\bullet}))\otimes_E\mathbf{E}_{M}^{\infty} \ar[d]\\
H^{h}(\mathrm{Tot}(\mathrm{gr}^{M}(\cT_{I_0,I_1}^{\bullet,\bullet}))) \ar^{\sim}[r] & H^{h-\#M}(\mathrm{Tot}(\mathrm{CE}_{I_0,I_1\setminus M}^{\bullet,\bullet}))\otimes_E\mathbf{E}_{M}^{\infty}
}
\end{equation}
with the horizontal isomorphisms from (\ref{equ: grade J Z to CE}), the LHS vertical map being (\ref{equ: relative bottom embedding grade M}), and the RHS vertical map being the tensor of the natural map
\begin{equation}\label{equ: CE relative image 1}
H^{h-\#M}(\mathrm{Tot}(\mathrm{CE}_{I_0,I_1\setminus M,\fh}^{\bullet,\bullet}))\rightarrow H^{h-\#M}(\mathrm{Tot}(\mathrm{CE}_{I_0,I_1\setminus M}^{\bullet,\bullet}))
\end{equation}
with the identity map of $\mathbf{E}_{M}^{\infty}$.
We write $X\subseteq H^{h-\#M}(\mathrm{Tot}(\mathrm{CE}_{I_0,I_1\setminus M}^{\bullet,\bullet}))$ for the image of (\ref{equ: CE relative image 1}), so that $X\otimes_E\mathbf{E}_{M}^{\infty}$ is the image of (\ref{equ: relative bottom embedding grade M}) by the commutative diagram (\ref{equ: CE relative image diagram}).
On one hand, we have
\begin{multline}\label{equ: relative grade M dim}
\Dim_EX\otimes_E\mathbf{E}_{M}^{\infty}=\Dim_EH^{h}(\mathrm{Tot}(\mathrm{gr}^{M}(\cT_{I_0,I_1,\fh}^{\bullet,\bullet})))\\
=\sum_{\ell=\#I_0}^{\#I_1}\Dim_E\mathrm{gr}^{M}(E_{2,I_0,I_1,\fh}^{-\ell,\ell+h})=\sum_{\ell=\#I_0}^{\#I_1}\#\Psi_{I_0,I_1,\fh,\ell}=
\#\cS_{(I_0\setminus I_1)\setminus M}^{J_{\fh}}.
\end{multline}
using the definition of $\Psi_{I_0,I_1,\fh,\ell}$ which is based on the bijection (\ref{equ: atom to partition}).
On the other hand, we have
\begin{equation}\label{equ: CE relative upper bound M dim}
\Dim_E X\leq \#\Gamma_{(I_1\setminus I_0)\setminus M}(J_{\fh})
\end{equation}
from Lemma~\ref{lem: CE relative upper bound} (upon replacing $I_1$ in \emph{loc.cit.} with $I_1\setminus M$) with the equality holds if and only if the composition of (\ref{equ: CE relative image 1}) with
\begin{equation}\label{equ: CE M decomposition}
H^{h-\#M}(\mathrm{Tot}(\mathrm{CE}_{I_0,I_1\setminus M}^{\bullet,\bullet}))\buildrel\sim\over\longrightarrow \bigoplus_{w\in\Gamma_{(I_1\setminus I_0)\setminus M}}\mathfrak{e}_{I_0,w}
\end{equation}
exactly equals
\begin{equation}\label{equ: CE relative upper bound M}
\bigoplus_{w\in\Gamma_{(I_1\setminus I_0)\setminus M}(J_{\fh})}\mathfrak{e}_{I_0,w}.
\end{equation}
Now that we have 
\[\#\Gamma_{(I_1\setminus I_0)\setminus M}(J_{\fh})=\#\cS_{(I_0\setminus I_1)\setminus M}^{J_{\fh}}\]
from \ref{it: coxeter partition 2} of Proposition~\ref{prop: coxeter partition cardinality}, we conclude that (\ref{equ: CE relative upper bound M dim}) must be an equality and thus the composition of (\ref{equ: CE relative image 1}) with (\ref{equ: CE M decomposition}) has image exactly (\ref{equ: CE relative upper bound M}).
This together with the commutative diagram (\ref{equ: CE relative image diagram}) finishes the proof.
\end{proof}

We fix a choice of a pair $j_0,j_1\in\Delta$ with $j_0<j_1$ in the rest of this section.
We write $\widehat{j}_0\defeq \Delta\setminus\{j_0\}$, $\widehat{j}_1\defeq \Delta\setminus\{j_1\}$, $J_0\defeq [j_0+1,n-1]$, $J_1\defeq [1,j_1-1]$ and $J_2\defeq J_0\cap J_1=[j_0+1,j_1-1]$ for short.
Note that we have $J_0\cup J_1=\Delta$ and thus $\#J_2=\#J_0+\#J_1-\#\Delta$.
We write $\fh_{0}\defeq (\fh_{[2,n-1]}\cap\fl_{\widehat{j}_0})/\fh_{J_0}\subseteq \fl_{\widehat{j}_0}/\fh_{J_0}$, $\fh_{1}\defeq (\fh_{[1,n-2]}\cap\fl_{\widehat{j}_1})/\fh_{J_1}\subseteq \fl_{\widehat{j}_1}/\fh_{J_1}$ and $\fh_{2}\defeq (\fh_{[2,n-2]}\cap\fl_{\widehat{j}_0\cap\widehat{j}_1})/\fh_{J_2}\subseteq \fl_{\widehat{j}_0\cap\widehat{j}_1}/\fh_{J_2}$ for short. 
Let $Z_{\widehat{j}_0}^{\dagger}\subseteq Z_{\widehat{j}_0}$ and $Z_{\widehat{j}_1}^{\dagger}\subseteq Z_{\widehat{j}_1}$ be discrete subgroups of rank $1$ (see the discussion below (\ref{equ: Levi center dagger projection})).
Following the notation of (\ref{equ: double v relative general}), we consider the double complex 
\[\cT_{0}^{\bullet,\bullet}\defeq \cT_{J_0,\Delta,\emptyset,\fh_{0},\widehat{j}_0,H_{J_0}}^{\bullet,\bullet}\] 
and 
\[\cT_{0,\dagger}^{\bullet,\bullet}\defeq \cT_{J_0,\Delta,\emptyset,\fh_{0},\widehat{j}_0,H_{J_0}Z_{\widehat{j}_0}^{\dagger}}^{\bullet,\bullet}\] 
which fit into maps 
\begin{equation}\label{equ: special relative j0}
\cT_{0,\dagger}^{\bullet,\bullet}\rightarrow \cT_{0}^{\bullet,\bullet}\rightarrow \cT_{J_0,\Delta,\fh_{[2,n-1]}}^{\bullet,\bullet}\rightarrow \cT_{J_0,\Delta}^{\bullet,\bullet}.
\end{equation} 
Similarly, we consider the double complex 
\[\cT_{1}^{\bullet,\bullet}\defeq \cT_{J_1,\Delta,\emptyset,\fh_{1},\widehat{j}_1,H_{J_1}}^{\bullet,\bullet}\] 
and 
\[\cT_{1,\dagger}^{\bullet,\bullet}\defeq \cT_{J_1,\Delta,\emptyset,\fh_{1},\widehat{j}_1,H_{J_1}Z_{\widehat{j}_1}^{\dagger}}^{\bullet,\bullet}\] which fit into maps 
\begin{equation}\label{equ: special relative j1}
\cT_{1,\dagger}^{\bullet,\bullet}\rightarrow \cT_{1}^{\bullet,\bullet}\rightarrow \cT_{J_1,\Delta,\fh_{[1,n-2]}}^{\bullet,\bullet}\rightarrow \cT_{J_1,\Delta}^{\bullet,\bullet}.
\end{equation} 
Finally, we consider the double complex 
\[\cT_{2}^{\bullet,\bullet}\defeq \cT_{J_2,\Delta,\emptyset,\fh_{2},\widehat{j}_0\cap \widehat{j}_1,H_{J_2}}^{\bullet,\bullet}\] 
and
\begin{equation}\label{equ: special relative double j0 j1}
\cT_{2,\dagger}^{\bullet,\bullet}\defeq \cT_{J_2,\Delta,\emptyset,\fh_{2},\widehat{j}_0\cap \widehat{j}_1,H_{J_2}Z_{\widehat{j}_0}^{\dagger}Z_{\widehat{j}_1}^{\dagger}}^{\bullet,\bullet}
\end{equation}
which fit into maps 
\begin{equation}\label{equ: special relative j0 j1}
\cT_{2,\dagger}^{\bullet,\bullet}\rightarrow \cT_{2}^{\bullet,\bullet} \rightarrow \cT_{J_2,\Delta,\fh_{[2,n-2]}}^{\bullet,\bullet}\rightarrow \cT_{J_2,\Delta}^{\bullet,\bullet}.
\end{equation} 
We write $E_{\bullet,0}^{\bullet,\bullet}$ (resp.~$E_{\bullet,0,\dagger}^{\bullet,\bullet}$) for the spectral sequence associated with $\cT_{0}^{\bullet,\bullet}$ (resp.~with $\cT_{0,\dagger}^{\bullet,\bullet}$) and similarly for the others.

We will consider below $\#J_0\leq \ell_0\leq \#\Delta$ and $\#J_1\leq \ell_1\leq \#\Delta$ with $\ell_2\defeq \ell_0+\ell_1-\#\Delta$ satisfying $\#J_2\leq \ell_2\leq \#\Delta$.
We also write $\ell_0'\defeq \ell_0-\#J_0$, $\ell_1'\defeq \ell_1-\#J_1$ and $\ell_2'\defeq \ell_2-\#J_2$, which satisfy $\ell_2'=\ell_0'+\ell_1'$.

By the definition of $\mathrm{gr}^{J}(\Psi_{J_0,\Delta,\fh_{[2,n-1]},\ell_0})$ for each $J\subseteq \Delta$ and each $\#J_0\leq \ell_0\leq \#\Delta$, we see that it consists of at most one element $\Omega_{0,J}$ with associated pair $(S_{J},J)$, which exists if and only if either $\#J_0<\ell_0<\#\Delta$ with $S_{J}\setminus J=\{(1,\ell_0'+2)\}$ and $J=[\ell_0'+2,j_0]$, or $\ell_0=\#J_0$ with $S_{J}=[1,j_0]$ and $J=[2,j_0]$, or $\ell_0=\#J_0$ with $S_{J}=[1,j_0]$ and $J=[1,j_0]$.

By the definition of $\mathrm{gr}^{J}(\Psi_{J_1,\Delta,\fh_{[1,n-2]},\ell_1})$ for each $J\subseteq \Delta$ and each $\#J_1\leq \ell_1\leq \#\Delta$, we see that it consists of at most one element $\Omega_{1,J}$ with associated pair $(S_{J},J)$, which exists if and only if either $\#J_1<\ell_1<\#\Delta$ with $S_{J}\setminus J=\{(n-1-\ell_1',n)\}$ and $J=[j_1,n-2-\ell_1']$, or $\ell_1=\#J_1$ with $S_{J}=[j_1,n-1]$ and $J=[j_1,n-2]$, or $\ell_1=\#J_1$ with $S_{J}=[j_1,n-1]$ and $J=[j_1,n-1]$.

We write $\cP(I)$ for the power set of $I$ (namely the set of all subsets of $I$).
Using the bijection (\ref{equ: atom to partition}), we see that the injection 
\[\cS_{\Delta\setminus J_0}\times \cS_{\Delta\setminus J_1}\hookrightarrow \cS_{\Delta\setminus J_2}: (S,S')\mapsto S\sqcup S'\]
and the bijection
\[\cP_{\Delta\setminus J_0}\times \cP_{\Delta\setminus J_1}\buildrel\sim\over\longrightarrow \cP_{\Delta\setminus J_2}: (I,I')\mapsto I\sqcup I'\]
induce an injection
\begin{equation}\label{equ: atom partition cup injection 0 1}
\Psi_{J_0,\Delta,\ell_0}\times \Psi_{J_1,\Delta,\ell_1}\hookrightarrow \Psi_{J_2,\Delta,\ell_2}
\end{equation}
for each $\#J_0\leq \ell_0\leq \#\Delta$ and $\#J_1\leq \ell_1\leq \#\Delta$, which altogether gives an injection
\begin{equation}\label{equ: atom partition cup injection}
\big(\bigsqcup_{\ell_0=\#J_0}^{\#\Delta}\Psi_{J_0,\Delta,\ell_0}\big)\times \big(\bigsqcup_{\ell_1=\#J_1}^{\#\Delta}\Psi_{J_1,\Delta,\ell_1}\big)\hookrightarrow \bigsqcup_{\ell_2=\#J_2}^{\#\Delta}\Psi_{J_2,\Delta,\ell_2}.
\end{equation}
Now that $\fh_{[2,n-2]}=\fh_{[2,n-1]}\cap\fh_{[1,n-2]}$, the injection (\ref{equ: atom partition cup injection 0 1}) restricts to an injection
\begin{equation}\label{equ: atom partition cup injection 0 1 relative}
\Psi_{J_0,\Delta,\fh_{[2,n-1]},\ell_0}\times \Psi_{J_1,\Delta,\fh_{[1,n-2]},\ell_1}\hookrightarrow \Psi_{J_2,\Delta,\fh_{[2,n-2]},\ell_2}
\end{equation}
by the definition of each term of (\ref{equ: atom partition cup injection 0 1 relative}) based on the bijection (\ref{equ: atom to partition}).
The injection (\ref{equ: atom partition cup injection 0 1 relative}) further restricts to an injection 
\begin{equation}\label{equ: atom partition cup injection 0 1 J}
\mathrm{gr}^{J}(\Psi_{J_0,\Delta,\fh_{[2,n-1]},\ell_0})\times \mathrm{gr}^{J'}(\Psi_{J_1,\Delta,\fh_{[1,n-2]},\ell_1})\hookrightarrow \mathrm{gr}^{J\sqcup J'}(\Psi_{J_2,\Delta,\fh_{[2,n-2]},\ell_2})
\end{equation}
for each $\#J_0\leq \ell_0\leq \#\Delta$, $\#J_1\leq \ell_1\leq \#\Delta$, $J\subseteq\Delta\setminus J_0$ and $J'\subseteq \Delta\setminus J_1$ (with $J\cap J'=\emptyset$ and $J\sqcup J'\subseteq \Delta\setminus J_2$).
For each $\ell_0$, $\ell_1$, $J$ and $J'$ such that $\Omega_{0,J}\in \mathrm{gr}^{J}(\Psi_{J_0,\Delta,\fh_{[2,n-1]},\ell_0})$ and $\Omega_{1,J'}\in \mathrm{gr}^{J'}(\Psi_{J_1,\Delta,\fh_{[1,n-2]},\ell_1})$ are defined, we define $\Omega_{2,J,J'}$ as the image of $(\Omega_{0,J},\Omega_{1,J'})$ under (\ref{equ: atom partition cup injection 0 1 relative}).

We write below $h_0\defeq \#\Delta-2\#J_0$, $h_1\defeq \#\Delta-2\#J_1$ and $h_2\defeq \#\Delta-2\#J_2=h_0+h_1+\#\Delta$.

The maps (\ref{equ: special relative j0}) induce the following maps
\begin{equation}\label{equ: special relative j0 E2}
E_{2,0,\dagger}^{-\ell_0,\ell_0+h_0}\rightarrow E_{2,0}^{-\ell_0,\ell_0+h_0}\rightarrow E_{2,J_0,\Delta}^{-\ell_0,\ell_0+h_0}.
\end{equation}
We write $\mathbf{E}_{0,\dagger,\ell_0}$ (resp.~$\mathbf{E}_{0,\ell_0}$) for the image of $E_{2,0,\dagger}^{-\ell_0,\ell_0+h_0}$ (resp.~$E_{2,0}^{-\ell_0,\ell_0+h_0}$) in $E_{2,J_0,\Delta}^{-\ell_0,\ell_0+h_0}$.
Similarly, the maps (\ref{equ: special relative j1}) induce the following maps
\begin{equation}\label{equ: special relative j1 E2}
E_{2,1,\dagger}^{-\ell_1,\ell_1+h_1}\rightarrow E_{2,1}^{-\ell_1,\ell_1+h_1}\rightarrow E_{2,J_1,\Delta}^{-\ell_1,\ell_1+h_1}
\end{equation}
and the maps (\ref{equ: special relative j0 j1}) induce the following maps
\begin{equation}\label{equ: special relative j0 j1 E2}
E_{2,2,\dagger}^{-\ell_2,\ell_2+h_2}\rightarrow E_{2,2}^{-\ell_2,\ell_2+h_2}\rightarrow E_{2,J_2,\Delta}^{-\ell_2,\ell_2+h_2}.
\end{equation}
We define $\mathbf{E}_{1,\dagger,\ell_1}\subseteq \mathbf{E}_{1,\ell_1}\subseteq E_{2,J_1,\Delta}^{-\ell_1,\ell_1+h_1}$ and $\mathbf{E}_{2,\dagger,\ell_2}\subseteq \mathbf{E}_{2,\ell_2}\subseteq E_{2,J_2,\Delta}^{-\ell_2,\ell_2+h_2}$ similarly.
\begin{lem}\label{lem: relative coh E2}
Let $j_0,j_1\in\Delta$ as above with associated $J_0$, $J_1$ and $J_2=J_0\cap J_1$.
We have the following results.
\begin{enumerate}[label=(\roman*)]
\item \label{it: relative coh E2 0} Let $\#J_0\leq \ell_0\leq \#\Delta$. The natural maps (\ref{equ: special relative j0 E2}) are embeddings. We have 
    \[\mathbf{E}_{0,\dagger,\ell_0}=\mathbf{E}_{0,\ell_0}=E\overline{x}_{\Omega_{0,[\ell_0'+2,j_0]}}\]
    when $\#J_0<\ell_0<\#\Delta$, and
    \[\mathbf{E}_{0,\dagger,\#J_0}=E\overline{x}_{\Omega_{0,[2,j_0]}}\subseteq E\overline{x}_{\Omega_{0,[1,j_0]}}\oplus E\overline{x}_{\Omega_{0,[2,j_0]}}=\mathbf{E}_{0,\#J_0}.\]
\item \label{it: relative coh E2 1} Let $\#J_1\leq \ell_1\leq \#\Delta$. The natural maps (\ref{equ: special relative j1 E2}) are embeddings. We have
    \[\mathbf{E}_{1,\dagger,\ell_1}=\mathbf{E}_{1,\ell_1}=E\overline{x}_{\Omega_{1,[j_1,n-2-\ell_1']}}\]
    when $\#J_1<\ell_1<\#\Delta$, and
    \[\mathbf{E}_{1,\dagger,\#J_1}=E\overline{x}_{\Omega_{1,[j_1,n-2]}}\subseteq E\overline{x}_{\Omega_{1,[j_1,n-1]}}\oplus E\overline{x}_{\Omega_{1,[j_1,n-1]}}=\mathbf{E}_{1,\#J_1}.\]
\item \label{it: relative coh E2 2} Let $\#J_2\leq \ell_2\leq \#\Delta$. The natural maps (\ref{equ: special relative j0 j1 E2}) are embeddings. The space $\mathbf{E}_{2,\ell_2}$ is spanned by $\overline{x}_{\Omega_{2,J,J'}}$ for all $\Omega_{2,J,J'}$ with bidegree $(-\ell_2,\ell_2+h_2)$, and its subspace $\mathbf{E}_{2,\dagger,\ell_2}$ is spanned by those $\overline{x}_{\Omega_{2,J,J'}}$ which further satisfy $J\neq[1,j_0]$ and $J'\neq [j_1,n-1]$.
\end{enumerate}
\end{lem}
\begin{proof}
We only prove \ref{it: relative coh E2 2}, and the proof of \ref{it: relative coh E2 0} and \ref{it: relative coh E2 1} are similar and simpler.
Let $J_2\subseteq I\subseteq \widehat{j}_0\cap\widehat{j}_1$.
We write $I^{-}\defeq I^{1}$ and $n_{-}\defeq n_{1}$, and then $P_{-}(I)\defeq P^{2n_{-}-1}(\fh_{I^{-}})\cong P^{2n_{-}-1}(\fg_{n_{-}})$ if $n_{-}>1$, and $P_{-}(I)\defeq \Hom(\fz_{[2,n-1]},E)$ if $n_{-}=1$.
Similarly, we write $I^{+}\defeq I^{r_{I}}$ and $n_{+}\defeq n_{r_{I}}$, and then $P_{+}(I)\defeq P^{2n_{+}-1}(\fh_{I^{+}})\cong P^{2n_{+}-1}(\fg_{n_{+}})$ if $n_{+}>1$, and $P_{+}(I)\defeq \Hom(\fz_{[1,n-2]},E)$ if $n_{+}=1$.
Note that the degree of $P_{\ast}(I)$ is always $2n_{\ast}-1$ for $\ast\in\{+,-\}$.
Following Proposition~\ref{prop: relative interval embedding} and Remark~\ref{rem: h tuple}, we have
\[
E_{1,2}^{-\ell,k}=\bigoplus_{J_2\subseteq I\subseteq \widehat{j}_0\cap\widehat{j}_1,\#I=\ell}H^{k}(L_{I}/H_{J_2},\fl_{I}\cap\fh_{2},1_{L_{I}/H_{J_2}})
\]
with the following isomorphism between graded $E$-algebras
\begin{equation}\label{equ: relative coh E1 2 product}
H^{\bullet}(L_{I}/H_{J_2},\fl_{I}\cap\fh_{2},1_{L_{I}/H_{J_2}})\buildrel\sim\over\longrightarrow H^{\bullet}(L_{I},\fl_{I}\cap\fh_{2},1_{L_{I}})
\cong \wedge (\Hom(Z_{I}^{\dagger},E)\oplus P_{-}(I)\oplus P_{+}(I)).
\end{equation}
and the differential map $d_{1,2}^{-\ell,k}$ described using Corollary~\ref{cor: relative Levi restriction}.
Similar to Remark~\ref{rem: h tuple}, (\ref{equ: relative coh E1 2 product}) admits a basis indexed by $\fh_{[2,n-2]}$-tuples $\Theta$, and we write $y_{\Theta}$ for the element of (\ref{equ: relative coh E1 2 product}) whose image in $H^{\bullet}(L_{I},1_{L_{I}})$ is $x_{\Theta}$.
We write $j_0'\defeq n_{-}=n_1$ and $j_1'\defeq n-n_{+}=n-n_{r_{I}}$.
We fix from now a choice of $\#J_2\leq \ell_2\leq\#\widehat{j}_0\cap\widehat{j}_1$ and a choice of $J_2\subseteq I\subseteq \widehat{j}_0\cap\widehat{j}_1$ such that $\#I=\ell_2$ and
\begin{equation}\label{equ: relative coh E1 direct summand}
H^{\ell_2+h_2}(L_{I}/H_{J_2},\fl_{I}\cap\fh_{2},1_{L_{I}/H_{J_2}})\neq 0.
\end{equation}
Thanks to (\ref{equ: relative coh E1 2 product}), there exists $k_0\in\{0,2j_0'-1\}$, $k_1\in\{0,2(n-j_1')-1\}$ and $0\leq k_2\leq \#\Delta\setminus I=n-1-\ell_2$ such that $k_0+k_1+k_2-\ell_2=h_2$.
For each $y_{\Omega}$ in (\ref{equ: relative coh E1 direct summand}) that satisfies $d_{1,2}^{-\ell_2,\ell_2+h_2}(y_{\Omega})=0$, we write $\overline{y}_{\Omega}$ for the image of $y_{\Omega}$ in $E_{2,2}^{-\ell_2,\ell_2+h_2}$.
We have the following possibilities.

\textbf{Case $1$} Assume that $k_0=k_1=0$ with $J_2\subseteq I$. Then $\#\Delta-2\ell_2\geq k_2-\ell_2=h_2=\#\Delta-2\#J_2$ together with $J_2\subseteq I$ forces $\ell_2=\#J_2$, $I=J_2$ and $k_2=\#\Delta\setminus J_2$, in which case (\ref{equ: relative coh E1 direct summand}) is spanned by $y_{\Theta}$ with $\Theta=(v,I,\un{k},\un{\Lambda})$ being the unique tuple that satisfies $I=J_2$ and $v=\mathbf{Log}_{J_2}^{\infty}$ (with $d_{1,2}^{-\#J_2,\#J_2+h_2}=0$). We have \[\overline{x}_{\Theta}=\overline{x}_{\Omega_{2,[1,j_0],[j_1,n-1]}}\] 
in this case.
        
\textbf{Case $2$} Assume that $k_0=0$ and $k_1=2n_{+}-1$ with $I^{+}\sqcup J_2\subseteq I$. If $I^{+}\sqcup J_2\subsetneq I$, then we have
\begin{multline*}
k_0+k_1+k_2-\ell_2\leq 2n_{+}-1+\#\Delta-2\#I\\
=\#\Delta-2\#I+2\#I^{+}+1\\
\leq \#\Delta-2\#J_2-2I\setminus(J_2\sqcup I^{+})+1<\#\Delta-2\#J_2=h_{2},
\end{multline*}
a contradiction. Hence, we must have $I=J_2\sqcup I^{+}$ with $I^{-}=\emptyset$ (or equivalently $j_0'=1$), in which case $k_0+k_1+k_2-\ell_2=h_2$ forces $k_2=\#\Delta\setminus I-1$. In particular, (\ref{equ: relative coh E1 direct summand}) is spanned by $y_{\Theta}$ where $\Theta=(v,I,\un{k},\un{\Lambda})$ runs through tuples that are characterized by $\#v\cap \mathbf{Log}_{I}^{\infty}=\#\Delta\setminus I-1$ and then either $\{\plog_{n-1}\}\subseteq v\subseteq \mathbf{Log}_{I}^{\infty}\sqcup\{\plog_{n-1}\}$ with $j_1'=n-1$, or $v\subseteq\mathbf{Log}_{I}^{\infty}$ and $\Lambda_{r_{I}}=\{2(n-j_1')-1\}$ with $j_1'<n-1$. It is clear that $d_{1,2}^{-\ell_2,\ell_2+h_2}(y_{\Theta})=0$ for all such $\Theta$, and thus (\ref{equ: relative coh E1 direct summand}) is annihilated by $d_{1,2}^{-\ell_2,\ell_2+h_2}$. If there exists $j\in[1,j_1'-1]\setminus J_2$ such that $\val_{j}\notin v$, then we have $v\cap\{\plog_{j},\val_{j}\}=\emptyset$ and
\[y_{\Theta}=\pm d_{1,2}^{-\ell_2-1,\ell_2+h_2}(y_{p_{j}^+(\Theta)})\in\mathrm{im}(d_{1,2}^{-\ell_2-1,\ell_2+h_2}).\]
Note that such $j$ does not exist if and only if either $v=\mathbf{Log}_{I\sqcup\{n-1\}}^{\infty}\sqcup\{\plog_{n-1}\}$ with $j_1'=n-1$, or $v=\mathbf{Log}_{I\sqcup\{j_1'\}}^{\infty}$ and $\Lambda_{r_{I}}=\{2(n-j_1')-1\}$ with $j_1'<n-1$ (with $\val_{j_1'}\notin v$ in both cases).
We thus conclude that the image of (\ref{equ: relative coh E1 direct summand}) in $E_{2,2}^{-\ell_2,\ell_2+h_2}$ is spanned by $\overline{y}_{\Theta}$, with the image $\overline{x}_{\Theta}$ of $\overline{y}_{\Theta}$ in $E_{2,J_2,\Delta}^{-\ell_2,\ell_2+h_2}$ identified with $\overline{x}_{\Omega_{2,[1,j_0],[j_1,j_1'-1]}}$.

\textbf{Case $3$} Assume that $k_0=2n_{-}-1$ and $k_1=0$ with $I^{-}\sqcup J_2\subseteq I$. A parallel discussion as in the last case shows that $I=I^{-}\sqcup J_2$ with $j_1'=n-1$ and $k_2=\#\Delta\setminus I-1$. Moreover, (\ref{equ: relative coh E1 direct summand}) is annihilated by $d_{1,2}^{-\ell_2,\ell_2+h_2}$, and its image in $E_{2,2}^{-\ell_2,\ell_2+h_2}$ spanned by various $\overline{y}_{\Theta}$, with the image $\overline{x}_{\Theta}$ of $\overline{y}_{\Theta}$ in $E_{2,J_2,\Delta}^{-\ell_2,\ell_2+h_2}$ identified with $\overline{x}_{\Omega_{2,[j_0'+1,j_0],[j_1,n-1]}}$.

\textbf{Case $4$} Assume that $k_0=2n_{-}-1$ and $k_1=2n_{+}-1$ with $I^{-}\sqcup I^{+}\sqcup J_2\subseteq I$. Similar discussion based on $k_0+k_1+k_2-\ell_2=h_2$ shows that we either have $I\setminus(I^{-}\sqcup I^{+}\sqcup J_2)=\{j\}$ for some $j$ with $k_2=\#\Delta\setminus I$, or we have $I=I^{-}\sqcup I^{+}\sqcup J_2$ with $k_2=\#\Delta\setminus I-2$. If $k_2=\#\Delta\setminus I$, then (\ref{equ: relative coh E1 direct summand}) is spanned by $y_{\Theta}$ which satisfies
\[d_{1,2}^{-\ell_2,\ell_2+h_2}(y_{\Theta})=\pm y_{p_{j}^{-}(\Theta)}.\]
If $k_2=\#\Delta\setminus I-2$, then (\ref{equ: relative coh E1 direct summand}) is spanned by $y_{\Theta}$ with $\Theta=(v,I,\un{k},\un{\Lambda})$ runs through tuples that satisfies $\#v\cap \mathbf{Log}_{I}^{\infty}=\#\Delta\setminus I-2$, $v\setminus\mathbf{Log}_{I}^{\infty}\subseteq\{\plog_{1},\plog_{n-1}\}$, with $\plog_{1}\in v$ when $j_0'=1$ and $\Lambda_{1}=\{2j_0'-1\}$ when $j_0'>1$, and with $\plog_{n-1}\in v$ when $j_1'=n-1$ and $\Lambda_{r_{I}}=\{2(n-j_1')-1\}$ when $j_1'<n-1$.
If there exists $j'\in[j_0'+1,j_1'-1]\setminus J_2$ such that $\val_{j'}\notin v$, then we have $v\cap\{\log_{j'},\val_{j'}\}=\emptyset$ and
\[y_{\Theta}=\pm d_{1,2}^{-\ell_2-1,\ell_2+h_2}(y_{p_{j'}^+(\Theta)})\in\mathrm{im}(d_{1,2}^{-\ell_2-1,\ell_2+h_2}).\]
Otherwise we have $\val_{j'}\in v$ for each $[j_0'+1,j_1'-1]\setminus J_2$, which forces 
\[v\cap \mathbf{Log}_{I}^{\infty} =\mathbf{Log}_{I\sqcup\{j_0',j_1'\}}^{\infty}.\] 
Consequently, when $k_2=\#\Delta\setminus I-2$, (\ref{equ: relative coh E1 direct summand}) is annihilated by $d_{1,2}^{-\ell_2,\ell_2+h_2}$, with its image in $E_{2,2}^{-\ell_2,\ell_2+h_2}$ is spanned by $\overline{y}_{\Theta}$, with the image of $\overline{x}_{\Theta}$ of $\overline{y}_{\Theta}$ in $E_{2,J_2,\Delta}^{-\ell_2,\ell_2+h_2}$ identified with $\overline{x}_{\Omega_{2,[j_0'+1,j_0],[j_1,j_1'-1]}}$.

To summarize, we have shown that $E_{2,2}^{-\ell_2,\ell_2+h_2}$ is spanned by a set $\{\overline{y}_{\Theta}\}$ for particular tuples $\Theta$ as described above, with its image $\{\overline{x}_{\Theta}\}$ in $E_{2,J_2,\Delta}^{-\ell_2,\ell_2+h_2}$ being linearly independent (and each $\overline{x}_{\Theta}$ identified with certain $\overline{x}_{\Omega_{2,J,J'}}$ for some $J,J'$ determined by $\Theta$). This proves that the RHS map of (\ref{equ: special relative j0 j1 E2}) is an embedding, with $\mathbf{E}_{2,\ell_2}$ admitting a basis as described in \ref{it: relative coh E2 2}.

Now we study the LHS map of (\ref{equ: special relative j0 j1 E2}).
Again by Proposition~\ref{prop: relative interval embedding} we have
\[
E_{1,2,\dagger}^{-\ell,k}=\bigoplus_{J_2\subseteq I\subseteq \widehat{j}_0\cap\widehat{j}_1,\#I=\ell}H^{k}(L_{I}/H_{J_2}Z_{\widehat{j}_0}^{\dagger}Z_{\widehat{j}_1}^{\dagger},\fl_{I}\cap\fh_{2},1_{L_{I}/H_{J_2}Z_{\widehat{j}_0}^{\dagger}Z_{\widehat{j}_1}^{\dagger}})
\]
with the following isomorphism between graded $E$-algebras
\begin{equation}\label{equ: relative coh E1 2 product dagger}
H^{\bullet}(L_{I}/H_{J_2}Z_{\widehat{j}_0}^{\dagger}Z_{\widehat{j}_1}^{\dagger},\fl_{I}\cap\fh_{2},1_{L_{I}/H_{J_2}Z_{\widehat{j}_0}^{\dagger}Z_{\widehat{j}_1}^{\dagger}})\buildrel\sim\over\longrightarrow \wedge (\Hom(Z_{I}^{\dagger}/Z_{\widehat{j}_0}^{\dagger}Z_{\widehat{j}_1}^{\dagger},E)\oplus P_{-}(I)\oplus P_{+}(I))
\end{equation}
and the differential map $d_{1,2,\dagger}^{-\ell,k}$ described using Corollary~\ref{cor: relative Levi restriction}.
It is clear that (\ref{equ: relative coh E1 2 product dagger}) is zero when
\[k>k_{I}\defeq (\#\Delta\setminus I-2)+(2n_{-}-1)+(2n_{+}-1)=\#\Delta\setminus I+2\#I^{-}+2\#I^{+}\]
and is $1$-dimensional when $k=k_{I}$. Note that
\begin{equation}\label{equ: relative coh deg upper bound}
k_{I}-\#I=\#\Delta-2\#I\setminus(I^{-}\sqcup I^{+})\leq h_2=\#\Delta-2\#J_2
\end{equation}
with equality holds if and only if $I=J_2\sqcup I^{-}\sqcup I^{+}$. 
Consequently, we have
\[
E_{1,2,\dagger}^{-\ell_2,\ell_2+h_{2}-\varepsilon}
\cong \bigoplus_{I=J_2\sqcup I^{-}\sqcup I^{+},\#I=\ell_2}
P_{-}(I)\otimes_E P_{+}(I)\otimes_E\wedge^{n-3-\ell_2-\varepsilon}\Hom(Z_{I}^{\dagger}/Z_{\widehat{j}_0}^{\dagger}Z_{\widehat{j}_1}^{\dagger},E)+\ast_{-}+\ast_{+},
\]
for $\varepsilon=0,1$, where $\ast_{-}\neq 0$ if and only if $\varepsilon=1$ and $I_{-}=\emptyset$ with 
\[\ast_{-}=P_{+}(I)\otimes_E\wedge^{n-3-\ell_2}\Hom(Z_{I}^{\dagger}/Z_{\widehat{j}_0}^{\dagger}Z_{\widehat{j}_1}^{\dagger},E),\]
and $\ast_{+}\neq 0$ if and only if $\varepsilon=1$ and $I_{+}=\emptyset$ with 
\[\ast_{+}=P_{-}(I)\otimes_E\wedge^{n-3-\ell_2}\Hom(Z_{I}^{\dagger}/Z_{\widehat{j}_0}^{\dagger}Z_{\widehat{j}_1}^{\dagger},E).\]
In particular, we deduce that
\[d_{1,2,\dagger}^{-\ell_2,\ell_2+h_2-\varepsilon}=0\]
for each $\#J_2\leq \ell_2\leq \#\widehat{j}_0\cap\widehat{j}_1$ and $\varepsilon=0,1$.
Now that $d_{1,2,\dagger}^{-\ell_2,\ell_2+h_2}=0=d_{1,2,\dagger}^{-\ell_2-1,\ell_2+h_2}=0$, we have $E_{1,2,\dagger}^{-\ell_2,\ell_2+h_{2}}=E_{2,2,\dagger}^{-\ell_2,\ell_2+h_{2}}$ and it suffices to understand the image of the composition of the following maps
\begin{multline}\label{equ: explicit E1 drop dagger}
P_{-}(I)\otimes_E P_{+}(I)\otimes_E\wedge^{n-3-\ell_2}\Hom(Z_{I}^{\dagger}/Z_{\widehat{j}_0}^{\dagger}Z_{\widehat{j}_1}^{\dagger},E)\\
\rightarrow 
P_{-}(I)\otimes_E P_{+}(I)\otimes_E\wedge^{n-3-\ell_2}\Hom(Z_{I}^{\dagger},E)
\rightarrow E_{2,2}^{-\ell_2,\ell_2+h_2}
\end{multline}
for each $j_0'\in[1,j_0]$ and $j_1'\in[j_1,n-1]$ with $I=J_2\sqcup I^{-}\sqcup I^{+}=J_2\sqcup [1,j_0'-1]\sqcup[j_1'+1,n-1]$.
Using (\ref{equ: J truncation val}) and Lemma~\ref{lem: central char transfer}, we know that the image of
\[\wedge^{n-3-\ell_2}\Hom(Z_{I}^{\dagger}/Z_{\widehat{j}_0}^{\dagger}Z_{\widehat{j}_1}^{\dagger},E)\rightarrow \wedge^{n-3-\ell_2}\Hom(Z_{I}^{\dagger},E)\]
is spanned by
\begin{equation}\label{equ: explicit val wedge 2 dagger}
v_{\dagger}\defeq \wedge_{j\in[j_0',j_0-1]}(\val_{j}-c_{j}\val_{j_0})\otimes_E \wedge_{j\in[j_1+1,j_1']}(\val_{j}-c_{j}\val_{j_1})
\end{equation}
for some $c_{j}\in E^{\times}$ for each $j\in (\widehat{j}_0\cap\widehat{j}_1)\setminus I=[j_0',j_0-1]\sqcup [j_1+1,j_1']$. Recall from our previous discussion on $E_{2,2}^{-\ell_2,\ell_2+h_2}$ that the image of
\[P_{-}(I)\otimes_E P_{+}(I)\otimes_E\otimes_Ev\]
in $E_{2,J_2,\Delta}^{-\ell_2,\ell_2+h_2}$ is non-zero for some $v\subseteq\mathbf{Log}_{I}^{\infty}$ with $\#v=n-3-\ell_2$ if and only if $v=\mathbf{Log}_{I\sqcup\{j_0',j_1'\}}^{\infty}$ (with $j_0'>0$ and $j_1'<n-1$). Hence, the image of the composition of (\ref{equ: explicit E1 drop dagger}) equals the image of
\[P_{-}(I)\otimes_E P_{+}(I)\otimes_E\otimes_Ev_{\dagger}\]
in $E_{2,J_2,\Delta}^{-\ell_2,\ell_2+h_2}$, which equals $Ec_{j_0'}c_{j_1'}\overline{x}_{\Theta}=E\overline{x}_{\Omega_{2,[j_0'+1,j_0],[j_1,j_1'-1]}}$ with $\Theta=(v,I,\un{k},\un{\Lambda})$ being a $(J_2,\Delta)$-tuple satisfying $I=[1,j_0'-1]\sqcup[j_1'+1,n-1]\sqcup J_2$ and $v=\mathbf{Log}_{I\sqcup\{j_0',j_1'\}}^{\infty}$.
This suffices to conclude that the LHS map of (\ref{equ: special relative j0 j1 E2}) is an embedding, with $\mathbf{E}_{2,\dagger,\ell_2}$ spanned by $\overline{x}_{\Omega_{2,[j_0'+1,j_0],[j_1,j_1'-1]}}$ for all $j_0'\in[1,j_0]$ and $j_1'\in[j_1,n-1]$ such that $I=[1,j_0'-1]\sqcup[j_1'+1,n-1]\sqcup J_2$ satisfies $\#I=\ell_2$.
This finishes the proof of \ref{it: relative coh E2 2}.
\end{proof}

\begin{lem}\label{lem: relative coh degeneracy}
Let $j_0,j_1\in\Delta$ as above with associated $J_0$, $J_1$ and $J_2=J_0\cap J_1$.
Then for each $0\leq a\leq 2$, $\#J_a\leq \ell_a\leq \#\Delta$ and $\ast\in\{,\dagger\}$, we have
\begin{equation}\label{equ: relative coh degeneracy}
E_{2,a,\ast}^{-\ell_a,\ell_a+h_a}=E_{\infty,a,\ast}^{-\ell_a,\ell_a+h_a}.
\end{equation}
\end{lem}
\begin{proof}
We only prove (\ref{equ: relative coh degeneracy}) for $a=2$, and the cases $a=0,1$ are similar and simpler.
Now that the maps (\ref{equ: special relative j0 j1 E2}) are embeddings by \ref{it: relative coh E2 2} of Lemma~\ref{lem: relative coh E2} and we have 
\[E_{2,J_2,\Delta}^{-\ell_2,\ell_2+h_2}=E_{\infty,J_2,\Delta}^{-\ell_2,\ell_2+h_2}\]
by Proposition~\ref{prop: bottom deg degeneracy}, it suffices to show that 
\begin{equation}\label{equ: relative coh degeneracy vanishing 1}
d_{r,2,\ast}^{-\ell_2,\ell_2+h_2}=0
\end{equation}
for each $r\geq 2$ and $\ast\in \{,\dagger\}$.
When $\ast=\dagger$, (\ref{equ: relative coh degeneracy vanishing 1}) follows from the fact that $E_{1,2,\dagger}^{-\ell,k}=0$ whenever $k-\ell>h_2$ (see the discussion around (\ref{equ: relative coh deg upper bound})).
It thus remains to show that
\begin{equation}\label{equ: relative coh degeneracy vanishing 2}
d_{r,2}^{-\ell_2,\ell_2+h_2}(\overline{y}_{\Theta})=0
\end{equation}
for each $r\geq 2$ and the unique tuple $\Theta=(v,I,\un{k},\un{\Lambda})$ described in the proof of \ref{it: relative coh E2 2} of Lemma~\ref{lem: relative coh E2} with $I=[1,j_0'-1]\sqcup[j_1'+1,n-1]\sqcup J_2$ satisfying $\#I=\ell_2$, $v=\mathbf{Log}_{(I\sqcup\{j_0',j_1'\})\cap\Delta}^{\infty}$ and the image of $\overline{y}_{\Theta}$ in $E_{2,J_2,\Delta}^{-\ell_2,\ell_2+h_2}$ being $\overline{x}_{\Omega_{2,[j_0'+1,j_0],[j_1,j_1'-1]}}$.
If $j_0'>0$ and $j_1'<n$, then $\overline{y}_{\Omega}$ lies in the image of the embedding $E_{2,2,\dagger}^{-\ell_2,\ell_2+h_2}\hookrightarrow E_{2,2}^{-\ell_2,\ell_2+h_2}$ and we know (\ref{equ: relative coh degeneracy vanishing 2}) from previous discussion. If $j_0'=0$ and $j_1'=n$, then (\ref{equ: relative coh degeneracy vanishing 2}) is also evident.
It remains to treat the cases when either $j_0'=0$ with $j_1'<n$, or $j_0'>0$ with $j_1'=n$.
we consider $\fh_{2,0}\defeq \fh_{J_0\cap\widehat{j}_1}/\fh_{J_2}$ and $\fh_{2,1}\defeq \fh_{J_1\cap\widehat{j}_0}/\fh_{J_2}$, and then consider the following double complex
\[\cT_{2,\dagger,0}^{\bullet,\bullet}\defeq \cT_{J_2,J_0,\mathbf{Log}_{J_0}^{\infty},\fh_{2,0},J_0\cap \widehat{j}_1,H_{J_2}Z_{\widehat{j}_1}^{\dagger}}^{\bullet,\bullet}\]
and
\[\cT_{2,\dagger,1}^{\bullet,\bullet}\defeq \cT_{J_2,J_1,\mathbf{Log}_{J_1}^{\infty},\fh_{2,1},J_1\cap \widehat{j}_0,H_{J_2}Z_{\widehat{j}_0}^{\dagger}}^{\bullet,\bullet}\]
both equipped with natural maps to $\cT_{2}^{\bullet,\bullet}$.
Similar to the discussion around (\ref{equ: relative coh deg upper bound}) and Lemma~\ref{lem: relative coh E2}, we can show that $E_{1,2,\dagger,0}^{-\ell,k}=0$ whenever $k-\ell>h_2$, and that the map $\cT_{2,\dagger,0}^{\bullet,\bullet}\rightarrow \cT_{2}^{\bullet,\bullet}$ induced an embedding $E_{2,2,\dagger,0}^{-\ell_2,h_2}\rightarrow E_{2,2}^{-\ell_2,\ell_2+h_2}$ with image $E\overline{y}_{\Theta}$.
Here the tuple $\Theta$ is characterized by the fact that the image of $\overline{y}_{\Theta}$ in $E_{2,J_2,\Delta}^{-\ell_2,\ell_2+h_2}$ is $\pm\overline{x}_{\Omega_{2,[1,j_0],[j_1,j_1'-1]}}$ for the unique $j_1'\in[j_1,n-1]$ such that $I=[j_1'+1,n-1]\sqcup J_2$ satisfies $\#I=\ell_2$. We thus clearly have $d_{r,2,\dagger,0}^{-\ell_2,\ell_2+h_2}=0$ for each $r\geq 2$, which gives (\ref{equ: relative coh degeneracy vanishing 2}) for each $r\geq 2$ and the tuple $\Theta$ just characterized (which satisfies $j_0'=0$ and $j_1'<n$).
Similarly, using the map $\cT_{2,\dagger,1}^{\bullet,\bullet}\rightarrow \cT_{2}^{\bullet,\bullet}$ we can prove (\ref{equ: relative coh degeneracy vanishing 2}) for each $r\geq 2$ and those tuple $\Theta$ characterized by $j_0'>0$ and $j_1'=n$.
The proof is thus finished.
\end{proof}

It follows from Lemma~\ref{lem: relative coh E2} and Lemma~\ref{lem: relative coh degeneracy} that the maps (\ref{equ: special relative j0}) induce embeddings
\begin{equation}\label{equ: relative total a}
H^{h_a}(\mathrm{Tot}(\cT_{a,\dagger}^{\bullet,\bullet}))\hookrightarrow H^{h_a}(\mathrm{Tot}(\cT_{a}^{\bullet,\bullet}))\hookrightarrow H^{h_a}(\mathrm{Tot}(\cT_{J_{a},\Delta}^{\bullet,\bullet}))
\end{equation}
for each $a=0,1,2$.
For each $J\subseteq \Delta\setminus J_0=[1,j_0]$ and $J'\subseteq \Delta\setminus J_1=[j_1,n-1]$, it follows from Proposition~\ref{prop: relative bottom embedding} that the RHS map of (\ref{equ: relative total a}) is strict with respect to the $\tau^{\bullet}(-)$ filtration on its source and its target, and thus together with Proposition~\ref{prop: Lie St Ext w decomposition} induces (for $a=0,1,2$ respectively) the following maps
\begin{multline}\label{equ: relative total grade J 0}
\mathrm{gr}^{J}(H^{h_0}(\mathrm{Tot}(\cT_{0}^{\bullet,\bullet})))\hookrightarrow
\mathrm{gr}^{J}(H^{h_0}(\mathrm{Tot}(\cT_{J_{0},\Delta}^{\bullet,\bullet})))\\
=H^{h_0-\#J}(\mathrm{CE}_{J_0,\Delta\setminus J}^{\bullet,\bullet})\otimes_E\wedge^{\#J}\Hom(Z_{\Delta\setminus J}^{\dagger},E)
\buildrel\sim\over\longrightarrow \bigoplus_{w\in\Gamma_{\Delta\setminus(J_0\sqcup J)}}\mathfrak{e}_{J_{0},w}\otimes_E\wedge^{\#J}\Hom(Z_{\Delta\setminus J}^{\dagger},E),
\end{multline}
\begin{multline}\label{equ: relative total grade J 1}
\mathrm{gr}^{J'}(H^{h_1}(\mathrm{Tot}(\cT_{1}^{\bullet,\bullet})))\hookrightarrow
\mathrm{gr}^{J'}(H^{h_1}(\mathrm{Tot}(\cT_{J_{1},\Delta}^{\bullet,\bullet})))\\
=H^{h_1-\#J'}(\mathrm{CE}_{J_1,\Delta\setminus J'}^{\bullet,\bullet})\otimes_E\wedge^{\#J'}\Hom(Z_{\Delta\setminus J'}^{\dagger},E)
\buildrel\sim\over\longrightarrow \bigoplus_{w\in\Gamma_{\Delta\setminus(J_{1}\sqcup J')}}\mathfrak{e}_{J_{1},w}\otimes_E\wedge^{\#J'}\Hom(Z_{\Delta\setminus J'}^{\dagger},E)
\end{multline}
and
\begin{multline}\label{equ: relative total grade J 2}
\mathrm{gr}^{J\sqcup J'}(H^{h_2}(\mathrm{Tot}(\cT_{2}^{\bullet,\bullet})))\hookrightarrow
\mathrm{gr}^{J\sqcup J'}(H^{h_2}(\mathrm{Tot}(\cT_{J_{2},\Delta}^{\bullet,\bullet})))\\
=H^{h_2-\#J-\#J'}(\mathrm{CE}_{J_0,\Delta\setminus J}^{\bullet,\bullet})\otimes_E\wedge^{\#J\sqcup J'}\Hom(Z_{\Delta\setminus(J\sqcup J')}^{\dagger},E)\\
\buildrel\sim\over\longrightarrow \bigoplus_{w\in\Gamma_{\Delta\setminus(J_2\sqcup J\sqcup J')}}\mathfrak{e}_{J_{2},w}\otimes_E\wedge^{\#J\sqcup J'}\Hom(Z_{\Delta\setminus(J\sqcup J')}^{\dagger},E).
\end{multline}

\begin{prop}\label{prop: relative coh E2 grade}
Let $j_0,j_1\in\Delta$ as above with associated $J_0$, $J_1$ and $J_2=J_0\cap J_1$.
Let $h_0\defeq \#\Delta-2\#J_0$, $h_1\defeq \#\Delta-2\#J_1$ and $h_2\defeq \#\Delta-2\#J_2=h_0+h_1+\#\Delta$.
Let $J\subseteq \Delta\setminus J_0$ and $J'\subseteq \Delta\setminus J_1$ with $J\sqcup J'\subseteq \Delta\setminus J_2$.
We have the following results.
\begin{enumerate}[label=(\roman*)]
\item \label{it: relative coh E2 grade 0} We have $\mathrm{gr}^{J}(H^{h_0}(\mathrm{Tot}(\cT_{0}^{\bullet,\bullet})))\neq 0$ if and only if $J=[j_{J}+1,j_0]$ for some $0\leq j_{J}\leq j_0$, in which case the composition of (\ref{equ: relative total grade J 0}) has image 
    \[\mathfrak{e}_{J_0,w_{1,j_{J}}}\otimes_E\wedge^{\#J}\Hom(Z_{\Delta\setminus J}^{\dagger},E).\]
\item \label{it: relative coh E2 grade 1} We have $\mathrm{gr}^{J'}(H^{h_1}(\mathrm{Tot}(\cT_{1}^{\bullet,\bullet})))\neq 0$ if and only if $J'=[j_1,j_{J'}-1]$ for some $j_1\leq j_{J}\leq n$, in which case the composition of (\ref{equ: relative total grade J 1}) has image 
    \[\mathfrak{e}_{J_{1},w_{n-1,j_{J'}}}\otimes_E\wedge^{\#J'}\Hom(Z_{\Delta\setminus J'}^{\dagger},E).\]
\item \label{it: relative coh E2 grade 2} We have $\mathrm{gr}^{J\sqcup J'}(H^{h_2}(\mathrm{Tot}(\cT_{2}^{\bullet,\bullet})))\neq 0$ if and only if $J=[j_{J},j_0]$ for some $1\leq j_{J}\leq j_0+1$ and $J'=[j_1,j_{J'}]$ for some $j_1-1\leq j_{J}\leq n-1$, in which case the image of the composition of (\ref{equ: relative total grade J 2}) is a $1$-dimensional $E$-subspace of
    \begin{equation}\label{equ: relative coh E2 grade 2 image}
    \sum_{w\in \Gamma_{w_{1,j_{J}},w_{n-1,j_{J'}}}}\mathfrak{e}_{J_2,w}\otimes_E\wedge^{\#J\sqcup J'}\Hom(Z_{\Delta\setminus(J\sqcup J')}^{\dagger},E).
    \end{equation}
    In particular, if $J\sqcup J'\sqcup J_2\neq \emptyset$ (namely $J$, $J'$ and $J_2$ not all empty), then (\ref{equ: relative coh E2 grade 2 image}) equals 
    \[
    \mathfrak{e}_{J_2,w_{1,j_{J}}w_{n-1,j_{J'}}}\otimes_E\wedge^{\#J\sqcup J'}\Hom(Z_{\Delta\setminus(J\sqcup J')}^{\dagger},E).
    \]
\end{enumerate}
\end{prop}
\begin{proof}
We only treat \ref{it: relative coh E2 grade 2} as the proof of \ref{it: relative coh E2 grade 0} and \ref{it: relative coh E2 grade 1} are easier.
It follows from \ref{it: relative coh E2 2} of Lemma~\ref{lem: relative coh E2} that $\mathrm{gr}^{J\sqcup J'}(E_{2,2}^{-\ell_2,\ell_2+h_2})\neq 0$ for some $J\subseteq \Delta\setminus J_0$ and $J'\subseteq \Delta\setminus J_1$ if and only if $J=[j_{J},j_0]$ for some $1\leq j_{J}\leq j_0+1$ with $\ell_0=\#J_0+\max\{0,j_{J}-2\}$ and $J'=[j_1,j_{J'}]$ for some $j_1-1\leq j_{J}\leq n-1$ with $\ell_1=\#J_1+\max\{0,n-2-j_{J'}\}$, in which case the map (with $\ell_2=\ell_0+\ell_1-\#\Delta$)
\[\mathrm{gr}^{J\sqcup J'}(E_{2,2}^{-\ell_2,\ell_2+h_2})\rightarrow \mathrm{gr}^{J\sqcup J'}(E_{2,J_{2},\Delta}^{-\ell_2,\ell_2+h_2})\]
is an embedding with image $E\overline{x}_{\Omega_{2,J,J'}}$. 
Similar to Lemma~\ref{lem: relative coh degeneracy}, we can show that
\[\mathrm{gr}^{J\sqcup J'}(E_{2,2}^{-\ell_2,\ell_2+h_2})=\mathrm{gr}^{J\sqcup J'}(E_{\infty,2}^{-\ell_2,\ell_2+h_2})\]
and thus the map
\begin{equation}\label{equ: relative total grade J 2 prime}
H^{h_2}(\mathrm{gr}^{J\sqcup J'}(\mathrm{Tot}(\cT_{2}^{\bullet,\bullet})))\rightarrow H^{h_2}(\mathrm{Tot}(\mathrm{gr}^{J\sqcup J'}(\cT_{J_{2},\Delta}^{\bullet,\bullet})))
\end{equation}
is an embedding, with the image being non-zero if and only if $J$ and $J'$ satisfy the aforementioned condition, in which case the image is given by $E\overline{x}_{\Omega_{2,J,J'}}$.
Now that the map between double complex $\mathrm{gr}^{J\sqcup J'}(\cT_{2}^{\bullet,\bullet})\rightarrow \mathrm{gr}^{J\sqcup J'}(\cT_{J_{2},\Delta}^{\bullet,\bullet})$ factors through $\mathrm{gr}^{J\sqcup J'}(\cT_{J_{2},\Delta,\fh_{[2,n-2]}}^{\bullet,\bullet})$, the map (\ref{equ: relative total grade J 2 prime}) factors through $H^{h_2}(\mathrm{gr}^{J\sqcup J'}(\mathrm{Tot}(\cT_{J_{2},\Delta,\fh_{[2,n-2]}}^{\bullet,\bullet})))$.
We thus conclude the rest of \ref{it: relative coh E2 grade 2} from Proposition~\ref{prop: grade M relative image} by taking $I_0$, $I_1$, $\fh$ and $M$ to be $J_{2}$, $\Delta$, $\fh_{[2,n-2]}$ and $J\sqcup J'$ in \emph{loc.cit.}.
\end{proof}

\begin{prop}\label{prop: relative coh E2 grade dagger}
Let $j_0,j_1\in\Delta$ as above with associated $J_0$, $J_1$ and $J_2=J_0\cap J_1$.
Let $h_0\defeq \#\Delta-2\#J_0$, $h_1\defeq \#\Delta-2\#J_1$ and $h_2\defeq \#\Delta-2\#J_2=h_0+h_1+\#\Delta$.
We have the following results.
\begin{enumerate}[label=(\roman*)]
\item \label{it: relative coh E2 grade dagger 0} The composition of the following maps
\[
H^{h_0}(\mathrm{Tot}(\cT_{0,\dagger}^{\bullet,\bullet}))\hookrightarrow H^{h_0}(\mathrm{Tot}(\cT_{0}^{\bullet,\bullet}))
\twoheadrightarrow H^{h_0}(\mathrm{Tot}(\cT_{0}^{\bullet,\bullet}))/\tau^{\Delta\setminus J_0}(H^{h_0}(\mathrm{Tot}(\cT_{0}^{\bullet,\bullet})))
\]
is an isomorphism.
\item \label{it: relative coh E2 grade dagger 1} The composition of the following maps
\[
H^{h_1}(\mathrm{Tot}(\cT_{1,\dagger}^{\bullet,\bullet}))\hookrightarrow H^{h_1}(\mathrm{Tot}(\cT_{1}^{\bullet,\bullet}))
\twoheadrightarrow H^{h_1}(\mathrm{Tot}(\cT_{1}^{\bullet,\bullet}))/\tau^{\Delta\setminus J_1}(H^{h_1}(\mathrm{Tot}(\cT_{1}^{\bullet,\bullet})))
\]
is an isomorphism.
\item \label{it: relative coh E2 grade dagger 2} The composition of the following maps
\begin{multline}\label{equ: relative coh E2 dagger composition}
H^{h_2}(\mathrm{Tot}(\cT_{2,\dagger}^{\bullet,\bullet}))\hookrightarrow H^{h_2}(\mathrm{Tot}(\cT_{2}^{\bullet,\bullet}))\\
\twoheadrightarrow H^{h_2}(\mathrm{Tot}(\cT_{2}^{\bullet,\bullet}))/\big(\tau^{\Delta\setminus J_0}(H^{h_2}(\mathrm{Tot}(\cT_{2}^{\bullet,\bullet})))+\tau^{\Delta\setminus J_1}(H^{h_2}(\mathrm{Tot}(\cT_{2}^{\bullet,\bullet})))\big)
\end{multline}
is an isomorphism.
\end{enumerate}
\end{prop}
\begin{proof}
We only treat \ref{it: relative coh E2 grade dagger 2} with \ref{it: relative coh E2 grade dagger 0} and \ref{it: relative coh E2 grade dagger 1} being similar and simpler.
We equip the last term of (\ref{equ: relative coh E2 dagger composition}) with the canonical decreasing filtration $\{\mathrm{Fil}^{-\ell}(-)\}_{\ell}$ induced from that of $H^{h_2}(\mathrm{Tot}(\cT_{2}^{\bullet,\bullet}))$, and write its graded piece $\mathrm{gr}^{-\ell}(-)$ as $\mathbf{E}_{2,\ell}'$ for short. Then the maps in (\ref{equ: relative coh E2 dagger composition}) are strict with respect to the canonical filtration on each term, and thus induce the maps between associated graded pieces
\begin{equation}\label{equ: relative coh E2 dagger graded composition}
\mathbf{E}_{2,\dagger,\ell}\hookrightarrow \mathbf{E}_{2,\ell}
\twoheadrightarrow \mathbf{E}_{2,\ell}'
\end{equation}
for each $\#J_2\leq\ell\leq \#\Delta$ (using the degeneracy result Lemma~\ref{lem: relative coh degeneracy}). It then follows from \ref{it: relative coh E2 2} of Lemma~\ref{lem: relative coh E2} and \ref{it: relative coh E2 grade 2} of Proposition~\ref{prop: relative coh E2 grade} that the composition of (\ref{equ: relative coh E2 dagger graded composition}) is an isomorphism for each $\ell$, which forces the composition of (\ref{equ: relative coh E2 dagger composition}) to be an isomorphism.
This finishes the proof of \ref{it: relative coh E2 grade dagger 2}.
\end{proof}

Let $u\in\Gamma$ with $J_{0}\subseteq J_{u}$ being an interval. For each $I\subseteq J_{u}$, it follows from Proposition~\ref{prop: Levi decomposition}, Lemma~\ref{lem: g P Ext transfer}, Lemma~\ref{lem: parabolic restriction x} and Corollary~\ref{cor: relative Levi restriction} that we have the following isomorphisms
\begin{multline}\label{equ: g P to Levi u}
\mathrm{Ext}_{D(\fg,P_{I}),U(\fl_{I}\cap\fh_{J_{0}})}^{\bullet+\ell(u)}(M^{I}(u),1_{D(\fg,P_{I})})\\
\cong H^{\bullet}(L_{I},\fl_{I}\cap\fh_{J_{0}},1_{L_{I}})\cong\mathrm{Tot}(\wedge^{\bullet}\Hom(Z_{I}^{\dagger},E)\otimes_EH^{\bullet}(\fl_{I},\fl_{I}\cap\fh_{J_{0}},1_{\fl_{I}}))
\end{multline}
for $k\geq 0$ which are functorial with respect to the choice of $I$ and $J_{0}$. Taking $J_0=\emptyset$ and $I=J$ for some $J\subseteq J_{u}$ in (\ref{equ: g P to Levi u}), we obtain the following maps (with $J'\defeq \Delta\setminus(J_{u}\setminus J)=J\sqcup\mathrm{Supp}(u)$ for short)
\begin{multline}\label{equ: g P u bottom}
\wedge^{\#J_{u}\setminus J}\Hom(Z_{J'}^{\dagger},E)=H^{\#J_{u}\setminus J}(L_{J'},1_{L_{J'}})^{\infty}\hookrightarrow H^{\#J_{u}\setminus J}(L_{J},1_{L_{J}})^{\infty}\hookrightarrow H^{\#J_{u}\setminus J}(L_{J},1_{L_{J}})\\
\cong \mathrm{Ext}_{D(\fg,P_{I})}^{\#\Delta\setminus J}(M^{I}(u),1_{D(\fg,P_{I})})
=H^{\#\Delta-2\#J}(\mathrm{Tot}(\cS_{u,J,J,\flat}^{\bullet,\bullet}))\rightarrow H^{\#\Delta-2\#J}(\mathrm{Tot}(\cS_{u,J,J_{u},\flat}^{\bullet,\bullet})).
\end{multline}
We write $H^{\#\Delta-2\#J}(\mathrm{Tot}(\cS_{u,J,J_{u},\flat}^{\bullet,\bullet}))^{\infty}$ for the image of the composition of (\ref{equ: g P u bottom}).
\begin{lem}\label{lem: S u sm sub}
Let $u\in\Gamma$ with $J_{u}$ being an interval and $J\subseteq J_{u}$.
Then the image of 
\begin{equation}\label{equ: S u sm sub}
H^{\#\Delta-2\#J}(\mathrm{Tot}(\cS_{u,J,J_{u},\fh_{J_{u}},\flat}^{\bullet,\bullet}))\rightarrow H^{\#\Delta-2\#J}(\mathrm{Tot}(\cS_{u,J,J_{u},\flat}^{\bullet,\bullet}))
\end{equation}
equals $H^{\#\Delta-2\#J}(\mathrm{Tot}(\cS_{u,J,J_{u},\flat}^{\bullet,\bullet}))^{\infty}$ which is $1$-dimensional.
\end{lem}
\begin{proof}
We write $E_{\bullet,u,\vartriangle}^{\bullet,\bullet}$ (resp.~$E_{\bullet,u}^{\bullet,\bullet}$) for the spectral sequence associated with the double complex $\cS_{u,J,J_{u},\fh_{J_{u}},\flat}^{\bullet,\bullet}$ (resp.~the double complex $\cS_{u,J,J_{u},\flat}^{\bullet,\bullet}$). The map between double complex $\cS_{u,J,J_{u},\fh_{J_{u}},\flat}^{\bullet,\bullet}\rightarrow \cS_{u,J,J_{u},\flat}^{\bullet,\bullet}$ induces a map between spectral sequences
\[E_{\bullet,u,\vartriangle}^{\bullet,\bullet+\ell(u)}\rightarrow E_{\bullet,u}^{\bullet,\bullet}\]
which on the first page is given by the map from
\[E_{1,u,\vartriangle}^{-\ell,k+\ell(u)}=\bigoplus_{J\subseteq I\subseteq J_{u}, \#I=\ell}\mathrm{Tot}(\wedge^{\bullet}\Hom(Z_{I}^{\dagger},E)\otimes_EH^{\bullet}(\fl_{I},\fl_{I}\cap\fh_{J_{u}},1_{\fl_{I}}))^{k}\]
to
\[E_{1,u}^{-\ell,k+\ell(u)}=\bigoplus_{J\subseteq I\subseteq J_{u}, \#I=\ell}\mathrm{Tot}(\wedge^{\bullet}\Hom(Z_{I}^{\dagger},E)\otimes_EH^{\bullet}(\fl_{I},1_{\fl_{I}}))^{k}\]
for each $\#J\leq\ell\leq \#J_{u}$ and $k\geq 0$, induced from the natural map $H^{\bullet}(\fl_{I},\fl_{I}\cap\fh_{J_{u}},1_{\fl_{I}})\rightarrow H^{\bullet}(\fl_{I},1_{\fl_{I}})$ for each $J\subseteq I\subseteq J_{u}$. Since $I\subseteq J_{u}$, the $E$-Lie subalgebra $\fl_{I}\cap\fh_{J_{u}}\subseteq\fl_{I}$ is an ideal with $\fl_{I}/\fl_{I}\cap\fh_{J_{u}}\cong \fz_{J_{u}}$ and we have the following canonical isomorphism
\[H^{\bullet}(\fl_{I},\fl_{I}\cap\fh_{J_{u}},1_{\fl_{I}}))\cong \wedge^{\bullet}\Hom(\fz_{J_{u}},E)\]
between graded $E$-algebra. This way, we obtain the following commutative diagram
\begin{equation}\label{equ: E1 u reduction}
\xymatrix{
E_{1,u,\vartriangle}^{-\ell,k+\ell(u)} \ar^{\sim}[rr] \ar@{^{(}->}[d] & & E_{1,J,J_{u},\vartriangle}^{-\ell,k} \ar@{^{(}->}[d]\\
E_{1,u}^{-\ell,k+\ell(u)} \ar^{\sim}[rr] & & E_{1,J,J_{u}}^{-\ell,k}
}
\end{equation}
with vertical maps being embeddings (see the discussion around (\ref{equ: g P to Levi u}) and then Lemma~\ref{lem: relative interval embedding Lie}) and
\begin{equation}\label{equ: E1 u wedge}
E_{1,J,J_{u},\vartriangle}^{-\ell,k}\defeq \bigoplus_{J\subseteq I\subseteq J_{u}, \#I=\ell}\mathrm{Tot}(\wedge^{\bullet}\Hom(Z_{I}^{\dagger},E)\otimes_E\wedge^{\bullet}\Hom(\fz_{J_{u}},E))^{k}.
\end{equation}
For each $k\geq 0$, we define a complex $E_{1}^{\bullet,k}$ by
\begin{equation}\label{equ: E1 u sm wedge}
E_{1}^{-\ell,k}=\bigoplus_{J\subseteq I\subseteq J_{u},\#I=\ell}\wedge^{k}\Hom(Z_{I}^{\dagger},E)
\end{equation}
with the differential map $d_{1}^{-\ell,k}: E_{1}^{-\ell,k}\rightarrow E_{1}^{-\ell+1,k}$ given by the direct sum of $(-1)^{m(I,j)}$ (see (\ref{equ: differential sign})) times the map 
\[\wedge^{k}\Hom(Z_{I}^{\dagger},E)\rightarrow \wedge^{k}\Hom(Z_{I\setminus\{j\}}^{\dagger},E)\]
induced from the map $Z_{I\setminus\{j\}}^{\dagger}\rightarrow Z_{I}^{\dagger}$, for each $J\subseteq I\subseteq J_{u}$ with $\#I=\ell$ and $j\in I\setminus J$. 
For each $v\subseteq \mathbf{Log}_{J}^{\infty}$, we define $E_{v}^{\bullet}$ as the complex given by
\[E_{v}^{-\ell}=\bigoplus_{J\subseteq I\subseteq I_{v}\cap J_{u},\#I=\ell}Ev\]
with a copy of $Ev$ for each $I$ and the differential map given by the direct sum of $(-1)^{m(I,j)}\mathrm{Id}_{Ev}$ from the copy indexed by $I$ to the copy indexed by $I\setminus\{j\}$ for each $j\in I\setminus J$.
It is clear that $E_{v}^{\bullet}$ is not acyclic if and only if $I_{v}\cap J_{u}=J$ if and only if $I_{v}\cap (J_{u}\setminus J)=\emptyset$, in which case $E_{v}^{\bullet}$ is concentrated in degree $-\#J$.
Then we observe that
\[E_{1}^{\bullet,k}=\bigoplus_{v\subseteq \mathbf{Log}_{J}^{\infty}, \#v=k}E_{v}^{\bullet}\]
for each $k\geq 0$, and thus $E_{2}^{-\ell,k}\defeq H^{-\ell}(E_{1}^{\bullet,k})$ is non-zero for some $k,\ell$ with $k-\ell\leq \#J_{u}\setminus J-\#J$ if and only if $\ell=\#J$ and $k=\#J_{u}\setminus J$, in which case the inclusion
\[\wedge^{\#J_{u}\setminus J}\Hom(Z_{J'}^{\dagger},E)=E_{\mathbf{Log}_{J'}^{\infty}}^{\bullet}\hookrightarrow E_{1}^{\bullet,\#J_{u}\setminus J}\]
inducing an isomorphism between $1$-dimensional $E$-vector spaces
\begin{equation}\label{equ: E2 u sm}
\wedge^{\#J_{u}\setminus J}\Hom(Z_{J'}^{\dagger},E)\buildrel\sim\over\longrightarrow E_{2}^{-\#J,\#J_{u}\setminus J}.
\end{equation}
It follows from (\ref{equ: E1 u wedge}) and (\ref{equ: E1 u sm wedge}) that we have
\[E_{2,u,\vartriangle}^{-\ell,k+\ell(u)}=E_{2,J,J_{u},\vartriangle}^{-\ell,k}\defeq \bigoplus_{k'+k''=k}E_{2}^{-\ell,k'}\otimes_E\wedge^{k''}\Hom(\fz_{J_{u}},E),\]
which is non-zero for some $(-\ell,k)$ satisfying $k-\ell=\#J_{u}\setminus J-\#J$ if and only if $\ell=\#J$, $k'=\#J_{u}\setminus J$ and $k''=0$, in which case we have
\[E_{2,u,\vartriangle}^{-\#J,\#J_{u}\setminus J+\ell(u)}=E_{2}^{-\#J,\#J_{u}\setminus J}\subseteq E_{2,J,J_{u}}^{-\#J,\#J_{u}\setminus J}\cong E_{2,u}^{-\#J,\#J_{u}\setminus J+\ell(u)}.\]
In other words, we have shown that the $E$-vector space $E_{2,u,\vartriangle}^{-\ell,k+\ell(u)}$ is non-zero for some $(-\ell,k)$ satisfying $k-\ell=\#J_{u}\setminus J-\#J$ if and only if $\ell=\#J$ and $k=\#J_{u}\setminus J$, in which case it embeds into $E_{2,u}^{-\ell,k+\ell(u)}$ with image (\ref{equ: E2 u sm}). The proof is thus finished.
\end{proof}

Now we move on to the application of Lemma~\ref{lem: S u sm sub}.
Let $x\in\Gamma^{\Delta\setminus J_0}$ be an element satisfying $D_L(x)=\{1\}$, and $y\in\Gamma^{\Delta\setminus J_1}$ be an element satisfying $D_L(y)=\{n-1\}$.
Note that there exists a unique $j\in[1,j_0]$ such that $x=w_{1,j}\defeq s_1\cdots s_{j}$, and a unique $j'\in[j_1,n-1]$ such that $y=w_{n-1,j'}\defeq s_{n-1}\cdots s_{j'}$.
We also consider an element $w\in\Gamma_{x,y}$ (see (\ref{equ: x y envelop})).
Given $x$, $y$ and $w$ as above, we have the following commutative diagrams (using the discussion around (\ref{equ: x Tits double flat}), (\ref{equ: coh Tits Levi}) and (\ref{equ: double x w diagram}))
\begin{equation}\label{equ: 0 x diagram}
\xymatrix{
\cT_{0,\dagger}^{\bullet,\bullet} \ar[r] &\cT_{J_0,\Delta,\fh_{[2,n-1]}}^{\bullet,\bullet} \ar[d] & \cS_{J_0,\Delta,\fh_{[2,n-1]}}^{\bullet,\bullet} \ar[l] \ar[r] \ar[d]& \cS_{x,J_0,J_{x},\fh_{[2,n-1]}}^{\bullet,\bullet} \ar[r] \ar[d] & \cS_{x,J_0,J_{x},\fh_{[2,n-1]},\flat}^{\bullet,\bullet} \ar[d]\\
&\cT_{J_0,\Delta}^{\bullet,\bullet} & \cS_{J_0,\Delta}^{\bullet,\bullet} \ar[l] \ar[r] & \cS_{x,J_0,J_{x}}^{\bullet,\bullet} \ar[r] & \cS_{x,J_0,J_{x},\flat}^{\bullet,\bullet}
},
\end{equation}
\begin{equation}\label{equ: 1 y diagram}
\xymatrix{
\cT_{1,\dagger}^{\bullet,\bullet} \ar[r] & \cT_{J_1,\Delta,\fh_{[1,n-2]}}^{\bullet,\bullet} \ar[d] & \cS_{J_1,\Delta,\fh_{[1,n-2]}}^{\bullet,\bullet} \ar[l] \ar[r] \ar[d]& \cS_{y,J_1,J_{y},\fh_{[1,n-2]}}^{\bullet,\bullet} \ar[r] \ar[d]& \cS_{y,J_1,J_{y},\fh_{[1,n-2]},\flat}^{\bullet,\bullet} \ar[d]\\
& \cT_{J_1,\Delta}^{\bullet,\bullet} & \cS_{J_1,\Delta}^{\bullet,\bullet} \ar[l] \ar[r] & \cS_{y,J_1,J_{y}}^{\bullet,\bullet} \ar[r]& \cS_{y,J_1,J_{y},\flat}^{\bullet,\bullet}
},
\end{equation}
and
\begin{equation}\label{equ: 2 w diagram}
\xymatrix{
\cT_{2,\dagger}^{\bullet,\bullet} \ar[r] & \cT_{J_2,\Delta,\fh_{[2,n-2]}}^{\bullet,\bullet} \ar[d] & \cS_{J_2,\Delta,\fh_{[2,n-2]}}^{\bullet,\bullet} \ar[l] \ar[r] \ar[d] & \cS_{w,J_2,J_{w},\fh_{[2,n-2]}}^{\bullet,\bullet} \ar[r] \ar[d] & \cS_{w,J_2,J_{w},\fh_{[2,n-2]},\flat}^{\bullet,\bullet} \ar[d]\\
& \cT_{J_2,\Delta}^{\bullet,\bullet} & \cS_{J_2,\Delta}^{\bullet,\bullet} \ar[l] \ar[r] & \cS_{w,J_2,J_{w}}^{\bullet,\bullet} \ar[r] & \cS_{w,J_2,J_{w},\flat}^{\bullet,\bullet}
}.
\end{equation}
Recall from the discussion around (\ref{equ: x Tits double flat}) and (\ref{equ: coh Tits Levi}) that the maps between double complex of the form $\cS_{\ast}^{\bullet,\bullet}\rightarrow \cT_{\ast}^{\bullet,\bullet}$ and of the form $\cS_{\ast}^{\bullet,\bullet}\rightarrow \cS_{\ast,\flat}^{\bullet,\bullet}$ from (\ref{equ: 0 x diagram}), (\ref{equ: 1 y diagram}) and (\ref{equ: 2 w diagram}) all induce isomorphisms on the first page.
Consequently, we obtain the following maps
\begin{equation}\label{equ: 0 1 2 sm part 0}
\mathbf{E}_{0,\dagger}=H^{\#\Delta-2\#J_0}(\mathrm{Tot}(\cT_{0,\dagger}^{\bullet,\bullet}))\rightarrow H^{\#\Delta-2\#J_0}(\mathrm{Tot}(\cS_{x,J_0,J_{x},\flat}^{\bullet,\bullet})),
\end{equation}
\begin{equation}\label{equ: 0 1 2 sm part 1}
\mathbf{E}_{1,\dagger}=H^{\#\Delta-2\#J_1}(\mathrm{Tot}(\cT_{1,\dagger}^{\bullet,\bullet}))\rightarrow H^{\#\Delta-2\#J_1}(\mathrm{Tot}(\cS_{y,J_1,J_{y},\flat}^{\bullet,\bullet}))
\end{equation}
and
\begin{equation}\label{equ: 0 1 2 sm part 2}
\mathbf{E}_{2,\dagger}=H^{\#\Delta-2\#J_2}(\mathrm{Tot}(\cT_{2,\dagger}^{\bullet,\bullet}))\rightarrow H^{\#\Delta-2\#J_2}(\mathrm{Tot}(\cS_{w,J_2,J_{w},\flat}^{\bullet,\bullet})).
\end{equation}
The following result follows directly from Lemma~\ref{lem: S u sm sub} by taking the pair $(u,J)$ to be $(x,J_0)$ (resp.~$(y,J_1)$, resp.~$(w,J_2)$) and then using the commutative diagram (\ref{equ: 0 x diagram}) (resp.~(\ref{equ: 1 y diagram}), resp.~(\ref{equ: 2 w diagram})).
\begin{prop}\label{prop: 0 1 2 sm part}
Let $x\in\Gamma^{\Delta\setminus J_0}$ withy $D_L(x)=\{1\}$, $y\in\Gamma^{\Delta\setminus J_1}$ with $D_L(y)=\{n-1\}$ and $w\in\Gamma_{x,y}$. Then we have the following results.
\begin{enumerate}[label=(\roman*)]
\item \label{it: 0 1 2 sm part 0} The image of the map (\ref{equ: 0 1 2 sm part 0}) is contained in \[H^{\#\Delta-2\#J_{0}}(\mathrm{Tot}(\cS_{x,J_{0},J_{x},\flat}^{\bullet,\bullet}))^{\infty}.\]
\item \label{it: 0 1 2 sm part 1} The image of the map (\ref{equ: 0 1 2 sm part 1}) is contained in
\[H^{\#\Delta-2\#J_{1}}(\mathrm{Tot}(\cS_{y,J_{1},J_{y},\flat}^{\bullet,\bullet}))^{\infty}.\]
\item \label{it: 0 1 2 sm part 2} The image of the map (\ref{equ: 0 1 2 sm part 2}) is contained in
\[H^{\#\Delta-2\#J_{2}}(\mathrm{Tot}(\cS_{w,J_{2},J_{w},\flat}^{\bullet,\bullet}))^{\infty}.\]
\end{enumerate}
\end{prop}

\section{Coxeter filtration and cup product}\label{sec: coxeter}
For each $I\subseteq\Delta$, we introduce our \emph{Coxeter filtration} on
\[\mathbf{E}_{I}=H^{2\#I-\#\Delta}(\mathrm{Tot}(\cT_{\Delta\setminus I,\Delta}^{\bullet,\bullet}))=\mathrm{Ext}_{G}^{\#I}(1_{G},V_{\Delta\setminus I}^{\rm{an}})\]
indexed by $\Gamma^{I}$.
Then we introduce our \emph{normalized cup product map}
\[\kappa_{I,I'}:\mathbf{E}_{I}\otimes_E\mathbf{E}_{I'}\buildrel\cup\over\longrightarrow\mathbf{E}_{I\sqcup I'}\]
for each $I,I'\subseteq\Delta$ with $I\cap I'=\emptyset$. We study the relation between Coxeter filtration on $\mathbf{E}_{I}$, $\mathbf{E}_{I'}$ and $\mathbf{E}_{I\sqcup I'}$ under $\kappa_{I,I'}$ (see Proposition~\ref{prop: coxeter cup comparison} and Proposition~\ref{prop: cup x y nonvanishing}), as well as the relation between $\kappa_{I,I'}$ and $\kappa_{I',I}$ (see Theorem~\ref{thm: cup exchange}).
Finally, for each $I_1\subseteq I\subseteq I_0$, we reduce the study of the usual cup product map
\[\kappa_{I_0,I_1}^{I}: \mathrm{Ext}_{G}^{\#I_0\setminus I}(V_{I_0}^{\rm{an}},V_{I}^{\rm{an}})\otimes_E\mathrm{Ext}_{G}^{\#I\setminus I_1}(V_{I}^{\rm{an}},V_{I_1}^{\rm{an}})\buildrel\cup\over\longrightarrow \mathrm{Ext}_{G}^{\#I_0\setminus I_1}(V_{I_0}^{\rm{an}},V_{I_1}^{\rm{an}})\]
to that of $\kappa_{I_0\setminus I, I\setminus I_1}$ (see Corollary~\ref{cor: bottom Tits cup transfer}).
\subsection{Coxeter filtration}\label{subsec: Coxeter Fil}
For each $I_1\subseteq I_0\subseteq \Delta$, we define an increasing filtration $\{\mathrm{Fil}_{x}(\mathbf{E}_{I_0,I_1})\}_{x\in\Gamma^{I_0\setminus I_1}}$ on $\mathbf{E}_{I_0,I_1}$ (see Theorem~\ref{thm: coxeter filtration}) crucially using Lemma~\ref{lem: x Ext with sm} and Corollary~\ref{cor: Ext between St}.

For each $I\subseteq \Delta$, we recall the subsets $\Gamma_I\subseteq \Gamma^I\subseteq \Gamma\subseteq W(G)$ and the set $\cS_{I}$ of partitions of $\sum_{\al\in I}\al$ from \S \ref{subsec: coxeter partition}.
\begin{prop}\label{prop: cardinality}
For each $I_1\subseteq I_0\subseteq \Delta$, we have
\begin{equation}\label{equ: cardinality}
\Dim_E\mathbf{E}_{I_0,I_1}=\#\Gamma^{I_0\setminus I_1}.
\end{equation}
\end{prop}
\begin{proof}
We write $J_0\defeq \Delta\setminus(I_0\setminus I_1)$ for short.
Recall from \ref{it: Ext between St 2} of Corollary~\ref{cor: Ext between St} that $\mathbf{E}_{I_0,I_1}$ admits a canonical decreasing filtration $\mathrm{Fil}^{\bullet}(\mathbf{E}_{I_0,I_1})$ with grade piece
\[\mathrm{gr}^{-\ell}(\mathbf{E}_{I_0,I_1})\cong E_{2,J_0,\Delta}^{-\ell,\ell+\#\Delta-2\#J_0}\]
admitting a basis of the form $\{\overline{x}_{\Omega}\mid \Omega\in \Psi_{J_0,\Delta,\ell}\}$ for each $\#J_0\leq \ell\leq \#\Delta$.
This together with the bijection (see the discussion around (\ref{equ: atom to partition}))
\[\bigsqcup_{\ell=\#J_0}^{\#\Delta}\Psi_{J_0,\Delta,\ell}\buildrel\sim\over\longrightarrow \{(S,I)\mid S\in \cS_{I_0\setminus I_1}, I\subseteq S\cap\Delta\}\]
implies that
\begin{multline}\label{equ: total dim count}
\Dim_E\mathbf{E}_{I_0,I_1}=\sum_{\ell=\#J_0}^{\#\Delta}\Dim_E E_{2,J_0,\Delta}^{-\ell,\ell+\#\Delta-2\#J_0}=\sum_{\ell=\#J_0}^{\#\Delta}\#\Psi_{J_0,\Delta,\ell}\\
=\#\{(S,I)\mid S\in \cS_{I_0\setminus I_1}, I\subseteq S\cap\Delta\}
=\#\{(S,I)\mid I\subseteq I_0\setminus I_1, I\subseteq S\in \cS_{I_0\setminus I_1}\}\\
=\#\{(S',I)\mid I\subseteq I_0\setminus I_1, S'\in \cS_{(I_0\setminus I_1)\setminus I}\}
=\#\{(S',I')\mid I'\subseteq I_0\setminus I_1, S'\in \cS_{I'}\}=\sum_{I'\subseteq I_0\setminus I_1}\#\cS_{I'}.
\end{multline}
Recall that we have $\Gamma^{I_0\setminus I_1}=\bigsqcup_{I\subseteq I_0\setminus I_1}\Gamma_I$ by definition and thus
\begin{equation}\label{equ: card of Gamma}
\#\Gamma^{I_0\setminus I_1}=\sum_{I\subseteq I_0\setminus I_1}\#\Gamma_I.
\end{equation}
We thus conclude (\ref{equ: cardinality}) by combining (\ref{equ: total dim count}) and (\ref{equ: card of Gamma}) with \ref{it: coxeter partition 1} of Proposition~\ref{prop: coxeter partition cardinality}.
\end{proof}

Let $I_1\subseteq I_0\subseteq \Delta$.
Thanks to Lemma~\ref{lem: change left cup}, the unique (up to scalars) embedding $V_{I_0}^{\infty}\hookrightarrow V_{I_0}^{\rm{an}}$ induces an isomorphism
\begin{equation}\label{equ: Ext St to Ext with sm}
\mathbf{E}_{I_0,I_1}\buildrel\sim\over\longrightarrow \mathrm{Ext}_{G}^{\#I_0\setminus I_1}(V_{I_0}^{\infty},V_{I_1}^{\rm{an}}).
\end{equation}
We consider a pair of subrepresentations $V\subseteq V'\subseteq V_{I_1}^{\rm{an}}$.
For each $k<\#I_0\setminus I_1$, it follows from Lemma~\ref{lem: x Ext with sm} (with $x=1$ in \emph{loc.cit.}) that $\mathrm{Ext}_{G}^k(V_{I_0}^{\infty}, W)=0$ for each $W\in\mathrm{JH}_{G}(V_{I_1}^{\rm{an}})$, and thus $\mathrm{Ext}_{G}^k(V_{I_0}^{\infty}, V'/V)=0=\mathrm{Ext}_{G}^k(V_{I_0}^{\infty}, V_{I_1}^{\rm{an}}/V')$ by a simple d\'evissage. In particular, the inclusions $V\subseteq V'\subseteq V_{I_1}^{\rm{an}}$ induce embeddings
\begin{equation}\label{equ: subspace of Ext}
\mathrm{Ext}_{G}^{\#I_0\setminus I_1}(V_{I_0}^{\infty},V)\hookrightarrow \mathrm{Ext}_{G}^{\#I_0\setminus I_1}(V_{I_0}^{\infty},V')\hookrightarrow\mathrm{Ext}_{G}^{\#I_0\setminus I_1}(V_{I_0}^{\infty},V_{I_1}^{\rm{an}})\buildrel\sim\over\longleftarrow\mathbf{E}_{I_0,I_1}
\end{equation}
with the last leftward isomorphism from (\ref{equ: Ext St to Ext with sm}).
For each subrepresentation $V\subseteq V_{I_1}^{\rm{an}}$, we write
\[\mathrm{Fil}_{V}(\mathbf{E}_{I_0,I_1})\defeq \mathrm{Ext}_{G}^{\#I_0\setminus I_1}(V_{I_0}^{\infty},V)\hookrightarrow \mathbf{E}_{I_0,I_1}\]
and
\begin{equation}\label{equ: Coxeter filtration index}
\Gamma_{V}\defeq \{x\in\Gamma^{I_0\setminus I_1}\mid C_{x}^{I_1}\in\mathrm{JH}_{G}(V)\}\subseteq \Gamma^{I_0\setminus I_1}
\end{equation}
for short.
\begin{lem}\label{lem: Ext subquotient}
Let $I_1\subseteq I_0\subseteq \Delta$. We have the following results.
\begin{enumerate}[label=(\roman*)]
\item \label{it: Ext subquotient 1} For each subquotient $V$ of $V_{I_1}^{\rm{an}}$, we have
\begin{equation}\label{equ: Ext subquotient dim}
\Dim_E \mathrm{Ext}_{G}^{\#I_0\setminus I_1}(V_{I_0}^{\infty},V)=\#\Gamma_{V}.
\end{equation}
\item \label{it: Ext subquotient 2} For each triple of subrepresentations $V_1\subseteq V_2\subseteq V_3\subseteq V_{I_1}^{\rm{an}}$, the short exact sequence $0\rightarrow V_2/V_1\rightarrow V_3/V_1 \rightarrow V_3/V_2\rightarrow 0$ induces a short exact sequence
\begin{equation}\label{equ: Ext subquotient seq}
0\rightarrow \mathrm{Ext}_{G}^{\#I_0\setminus I_1}(V_{I_0}^{\infty},V_2/V_1)\rightarrow \mathrm{Ext}_{G}^{\#I_0\setminus I_1}(V_{I_0}^{\infty},V_3/V_1)\rightarrow \mathrm{Ext}_{G}^{\#I_0\setminus I_1}(V_{I_0}^{\infty},V_3/V_2)\rightarrow 0.
\end{equation}
\end{enumerate}
\end{lem}
\begin{proof}
By Lemma~\ref{lem: x Ext with sm} (with $x=1$ in \emph{loc.cit.}) and a simple d\'evissage on $\mathrm{JH}_{G}(V)$ we deduce that
\begin{equation}\label{equ: Ext subquotient dim bound}
\Dim_E \mathrm{Ext}_{G}^{\#I_0\setminus I_1}(V_{I_0}^{\infty},V)\leq \#\Gamma_{V}
\end{equation}
for each subquotient $V$ of $V_{I_1}^{\rm{an}}$.
Now we consider a triple $V_1\subseteq V_2\subseteq V_3$ as in \ref{it: Ext subquotient 2}. As (\ref{equ: Ext subquotient dim bound}) holds for $V\in \{V_1,V_2/V_1,V_3/V_1,V_3/V_2,V_{I_1}^{\rm{an}}/V_3\}$, by simple d\'evissage we have
\begin{multline}\label{equ: Ext subquotient devissage}
\Dim_E \mathrm{Ext}_{G}^{\#I_0\setminus I_1}(V_{I_0}^{\infty},V_{I_1}^{\rm{an}})\\
\leq \Dim_E \mathrm{Ext}_{G}^{\#I_0\setminus I_1}(V_{I_0}^{\infty},V_1)+\Dim_E \mathrm{Ext}_{G}^{\#I_0\setminus I_1}(V_{I_0}^{\infty},V_3/V_1)+\Dim_E \mathrm{Ext}_{G}^{\#I_0\setminus I_1}(V_{I_0}^{\infty},V_{I_1}^{\rm{an}}/V_3)\\
\leq \Dim_E \mathrm{Ext}_{G}^{\#I_0\setminus I_1}(V_{I_0}^{\infty},V_1)+\Dim_E \mathrm{Ext}_{G}^{\#I_0\setminus I_1}(V_{I_0}^{\infty},V_2/V_1)\\
+\mathrm{Ext}_{G}^{\#I_0\setminus I_1}(V_{I_0}^{\infty},V_3/V_2)+\Dim_E \mathrm{Ext}_{G}^{\#I_0\setminus I_1}(V_{I_0}^{\infty},V_{I_1}^{\rm{an}}/V_3)\\
\leq \#\Gamma_{V_1}+\#\Gamma_{V_2/V_1}+\#\Gamma_{V_3/V_2}+\#\Gamma_{V_{I_1}^{\rm{an}}/V_3}=\#\Gamma_{V_{I_1}^{\rm{an}}}=\#\Gamma^{I_0\setminus I_1}.
\end{multline}
Since we have
\[\Dim_E \mathrm{Ext}_{G}^{\#I_0\setminus I_1}(V_{I_0}^{\infty},V_{I_1}^{\rm{an}})=\Dim_E \mathbf{E}_{I_0,I_1}=\#\Gamma^{I_0\setminus I_1}\]
by Proposition~\ref{prop: cardinality}, the inequalities in (\ref{equ: Ext subquotient devissage}) must all be equalities, which clearly gives the short exact sequence (\ref{equ: Ext subquotient seq}). As any subquotient $V$ of $V_{I_1}^{\rm{an}}$ has the form $V_2/V_1$ for some $V_1\subseteq V_2\subseteq V_{I_1}^{\rm{an}}$, the argument above also shows that (\ref{equ: Ext subquotient dim bound}) is always an equality.
\end{proof}

Let $I_1\subseteq I_0\subseteq \Delta$. For each $x\in \Gamma^{I_0\setminus I_1}$, we recall from  Lemma~\ref{lem: JH mult PS} that $[V_{I_1}^{\rm{an}}:C_{x}^{I_1}]=1$ and from discussion before Lemma~\ref{lem: basic subquotient of St} that $\tld{C}_{x}^{I_1}$ is the unique subrepresentation of $V_{I_1}^{\rm{an}}$ with cosocle $C_{x}^{I_1}$.
We set
\begin{equation}\label{equ: x filtration}
\mathrm{Fil}_{x}(\mathbf{E}_{I_0,I_1})\defeq \mathrm{Fil}_{\tld{C}_{x}^{I_1}}(\mathbf{E}_{I_0,I_1})\subseteq \mathbf{E}_{I_0,I_1}.
\end{equation}
More generally, for each subset $\Sigma\subseteq\Gamma^{I_0\setminus I_1}$, we write
\begin{equation}\label{equ: Coxeter subset filtration}
\mathrm{Fil}_{\Sigma}(\mathbf{E}_{I_0,I_1})\defeq \sum_{x\in\Sigma}\mathrm{Fil}_{x}(\mathbf{E}_{I_0,I_1})
\end{equation}
for short.
Recall from Lemma~\ref{lem: inclusion order} that $\tld{C}_{x'}^{I_1}\subseteq \tld{C}_{x}^{I_1}$ if and only if $x'\unlhd x$. This together with the discussion before Lemma~\ref{lem: Ext subquotient} implies that the following collection of subspaces of $\mathbf{E}_{I_0,I_1}$
\begin{equation}\label{equ: second filtration}
\{0, \mathbf{E}_{I_0,I_1}\}\sqcup\{\mathrm{Fil}_{x}(\mathbf{E}_{I_0,I_1})\}_{x\in \Gamma^{I_0\setminus I_1}}
\end{equation}
gives an increasing $\Gamma^{I_0\setminus I_1}$-filtration on $\mathbf{E}_{I_0,I_1}$ with respect to the partial-order $\unlhd$ on $\Gamma^{I_0\setminus I_1}$ (with $0$ being the zero subspace of $\mathbf{E}_{I_0,I_1}$).
We write $\mathrm{gr}_{x}(\mathbf{E}_{I_0,I_1})\defeq \mathrm{Fil}_{x}(\mathbf{E}_{I_0,I_1})/\sum_{x'\lhd x}\mathrm{Fil}_{x'}(\mathbf{E}_{I_0,I_1})$ for each $x\in \Gamma^{I_0\setminus I_1}$.
\begin{thm}\label{thm: coxeter filtration}
Let $I_1\subseteq I_0\subseteq \Delta$.
The increasing filtration (\ref{equ: second filtration}) is separated and exhaustive, and we have a natural isomorphism
\begin{equation}\label{equ: second filtration graded}
\mathrm{gr}_{x}(\mathbf{E}_{I_0,I_1})\buildrel\sim\over\longrightarrow \mathrm{Ext}_{G}^{\#I_0\setminus I_1}(V_{I_0}^{\infty}, C_{x}^{I_1})
\end{equation}
between $1$ dimensional spaces for each $x\in \Gamma^{I_0\setminus I_1}$.
\end{thm}
\begin{proof}
The filtration (\ref{equ: second filtration}) is obviously separated and exhaustive by definition.
We fix a choice of $x\in \Gamma^{I_0\setminus I_1}$ and define $V_1\defeq 0$, $V_2\defeq \sum_{x'\lhd x}\tld{C}_{x'}^{I_1}$ and $V_3\defeq \tld{C}_{x}^{I_1}$ as subrepresentations of $V_{I_1}^{\rm{an}}$.
It is clear from the discussion above Lemma~\ref{lem: Ext subquotient} that
\[\mathrm{Ext}_{G}^{\#I_0\setminus I_1}(V_{I_0}^{\infty},V_2)=\sum_{x'\lhd x}\mathrm{Ext}_{G}^{\#I_0\setminus I_1}(V_{I_0}^{\infty},\tld{C}_{x}^{I_1})=\sum_{x'\lhd x}\mathrm{Fil}_{x'}(\mathbf{E}_{I_0,I_1})\]
and that $\mathrm{Ext}_{G}^{\#I_0\setminus I_1}(V_{I_0}^{\infty},V_3)=\mathrm{Fil}_{x}(\mathbf{E}_{I_0,I_1})$, which together with (\ref{equ: Ext subquotient seq}) gives a natural isomorphism
\begin{equation}\label{equ: second filtration graded prime}
\mathrm{gr}_{x}(\mathbf{E}_{I_0,I_1})\buildrel\sim\over\longrightarrow \mathrm{Ext}_{G}^{\#I_0\setminus I_1}(V_{I_0}^{\infty}, V_3/V_2).
\end{equation}
Note that $V_3=\tld{C}_{x}^{I_1}$ has cosocle $C_{x}^{I_1}$, and we have $\Gamma_{V_3/V_2}=\{x\}$ by Lemma~\ref{lem: inclusion order}. In particular, the surjection $V_3/V_2\twoheadrightarrow C_{x}^{I_1}$ induces an isomorphism between $1$ dimensional spaces
\[\mathrm{Ext}_{G}^{\#I_0\setminus I_1}(V_{I_0}^{\infty}, V_3/V_2)\buildrel\sim\over\longrightarrow \mathrm{Ext}_{G}^{\#I_0\setminus I_1}(V_{I_0}^{\infty}, C_{x}^{I_1})\]
by Lemma~\ref{lem: Ext subquotient} (and Lemma~\ref{lem: x Ext with sm}, with $x=1$ in \emph{loc.cit.}).
The proof is thus finished.
\end{proof}

\begin{rem}\label{rem: outer coxeter}
The standard outer involution $G\rightarrow G: A\mapsto \overline{A}=w_0(A^{t})^{-1}w_0$ induces an involution on $\mathrm{Rep}^{\rm{an}}_{\rm{adm}}(G)$ that sends $V_{I}^{\infty}\subseteq V_{I}^{\rm{an}}$ to $V_{w_0(I)}^{\infty}\subseteq V_{w_0(I)}^{\rm{an}}$ for each $I\subseteq\Delta$. 
Given $I_1\subseteq I_0\subseteq\Delta$, this outer involution of $G$ induces an isomorphism
\begin{equation}\label{equ: outer Ext isom}
\mathbf{E}_{I_0,I_1}\buildrel\sim\over\longrightarrow\mathbf{E}_{w_0(I_0),w_0(I_1)}
\end{equation}
as well as a bijection (see (\ref{equ: coxeter involution}))
\[
\Gamma^{I_0\setminus I_1}\buildrel\sim\over\longrightarrow \Gamma^{w_0(I_0)\setminus w_0(I_1)}, x\mapsto\overline{x}.
\]
Now that the outer involution of $G$ sends $\tld{C}_{x}^{I_1}\subseteq V_{I_1}^{\rm{an}}$ to $\tld{C}_{\overline{x}}^{w_0(I_1)}\subseteq V_{w_0(I_1)}^{\rm{an}}$ for each $x\in \Gamma^{I_0\setminus I_1}$, we conclude that the isomorphism (\ref{equ: outer Ext isom}) restricts to an isomorphism
\[\mathrm{Fil}_{x}(\mathbf{E}_{I_0,I_1})\buildrel\sim\over\longrightarrow \mathrm{Fil}_{\overline{x}}(\mathbf{E}_{w_0(I_0),w_0(I_1)})\]
for each $x\in \Gamma^{I_0\setminus I_1}$.
\end{rem}

\begin{cor}\label{cor: tau J Coxeter filtration}
Let $I_1\subseteq I_0\subseteq\Delta$ and $J\subseteq\Delta$ with $\tau^{J}(V_{I_1}^{\rm{an}})$ from (\ref{equ: loc an St J filtration}). Then the following map
\begin{equation}\label{equ: tau J Coxeter filtration}
\mathrm{Ext}_{G}^{\#I_0\setminus I_1}(V_{I_0}^{\infty},\tau^{J}(V_{I_1}^{\rm{an}}))\hookrightarrow 
\mathrm{Ext}_{G}^{\#I_0\setminus I_1}(V_{I_0}^{\infty},V_{I_1}^{\rm{an}})=\mathbf{E}_{I_0,I_1}
\end{equation}
induced from $\tau^{J}(V_{I_1}^{\rm{an}})\hookrightarrow V_{I_1}^{\rm{an}}$ is an embedding.
Moreover, the image of (\ref{equ: tau J Coxeter filtration}) is non-zero only if $J\supseteq \Delta\setminus(I_0\setminus I_1)$, in which case it equals
\begin{equation}\label{equ: tau J Coxeter filtration image}
\sum_{x\in\Gamma^{(I_0\setminus I_1)\cap J}}\mathrm{Fil}_{x}(\mathbf{E}_{I_0,I_1}).
\end{equation}
\end{cor}
\begin{proof}
It follows from Theorem~\ref{thm: coxeter filtration} (together with the notation (\ref{equ: Coxeter filtration index}) and (\ref{equ: Coxeter subset filtration})) that the map (\ref{equ: tau J Coxeter filtration}) is an embedding with image
\[\sum_{x\in\Gamma_{\tau^{J}(V_{I_1}^{\rm{an}})}}\mathrm{Fil}_{x}(\mathbf{E}_{I_0,I_1}).\]
Now that $\tau^{J}(V_{I_1}^{\rm{an}})$ admits a filtration with graded pieces being $\mathrm{gr}^{J'}(V_{I_1}^{\rm{an}})$ for each $J'\supseteq J$, it suffices to show that
\begin{equation}\label{equ: Coxeter index grade J}
\Gamma_{\mathrm{gr}^{J'}(V_{I_1}^{\rm{an}})}
\end{equation}
is non-empty only if $J'\supseteq \Delta\setminus(I_0\setminus I_1)$, in which case it equals $\Gamma_{(I_0\setminus I_1)\cap J}$.
Given $x\in\Gamma^{I_0\setminus I_1}$, we deduce from (\ref{equ: loc an St J grade}) that
\[
[\mathrm{gr}^{J'}(V_{I_1}^{\rm{an}}):C_{x}^{I_1}]=[\cF_{P_{J'}}^{G}(\mathfrak{v}_{J'},V_{I_1,J'}^{\infty}):C_{x}^{I_1}]
=[\mathfrak{v}_{J'}:L(x)][i_{J',I_{x}}^{\infty}(V_{I_1,J'}^{\infty}):V_{I_1\sqcup(\mathrm{Supp}(x)\cap I_{x}),I_{x}}^{\infty}]
\]
which is non-zero (and equals $1$) if and only if $\mathrm{Supp}(x)=\Delta\setminus J'$ (using Lemma~\ref{lem: JH OS}, \ref{it: Lie St JH 2} of Lemma~\ref{lem: Lie St JH}, as well as \ref{it: sm cube 1} and \ref{it: sm cube 3} of Lemma~\ref{lem: sm cube}).
Hence, $x\in\Gamma^{I_0\setminus I_1}$ belongs to (\ref{equ: Coxeter index grade J}) if and only if $\mathrm{Supp}(x)=\Delta\setminus J'$, if and only if $J'\supseteq\Delta\setminus(I_0\setminus I_1)$ and $\mathrm{Supp}(x)=(I_0\setminus I_1)\cap J'$. The proof is thus finished. 
\end{proof}

Let $I_1\subseteq I_0\subseteq \Delta$ and $x\in\Gamma^{I_0\setminus I_1}$. Recall from Lemma~\ref{lem: Hom x w St} (upon replacing $x$ and $w$ in \emph{loc.it.} with $1$ and $x$ respectively) that there exists a unique (up to scalars) non-zero map 
\begin{equation}\label{equ: x St surjection}
V_{I_1}^{\rm{an}}\rightarrow V_{x,I_1}^{\rm{an}}.
\end{equation}
\begin{prop}\label{prop: coxeter filtration x surjection}
Then map (\ref{equ: x St surjection}) induces a surjection
\begin{equation}\label{equ: coxeter filtration x surjection}
\zeta_{I_0,I_1}^{x}: \mathbf{E}_{I_0,I_1}\cong\mathrm{Ext}_{G}^{\#I_0\setminus I_1}(V_{I_0}^{\infty},V_{I_1}^{\rm{an}})\twoheadrightarrow \mathrm{Ext}_{G}^{\#I_0\setminus I_1}(V_{I_0}^{\infty},V_{x,I_1}^{\rm{an}}).
\end{equation}
which satisfies
\begin{equation}\label{equ: coxeter filtration x surjection kernel}
\mathrm{ker}(\zeta_{I_0,I_1}^{x})=\sum_{y\in\Gamma^{I_0\setminus I_1}, x\not\unlhd y}\mathrm{Fil}_{y}(\mathbf{E}_{I_0,I_1}).
\end{equation}
\end{prop}
\begin{proof}
We write $V$ for the image of (\ref{equ: x St surjection}). It follows from \ref{it: x loc an St vanishing 2} of Lemma~\ref{lem: x loc an St vanishing} that we have
\[\mathrm{soc}_{G}(V)=\mathrm{soc}_{G}(V_{x,I_1}^{\rm{an}})=C_{x}^{I_1}.\]
This together with $[V_{I_1}^{\rm{an}}:C_{x}^{I_1}]=1$ (see \ref{it: JH mult PS 2} of Lemma~\ref{lem: JH mult PS}) ensures that $V$ is the unique quotient of $V_{I_1}^{\rm{an}}$ with socle $C_{x}^{I_1}$, which satisfies 
\[\Gamma_{V}=\{y\in\Gamma^{I_0\setminus I_1}\mid x\unlhd y\}\]
by Lemma~\ref{lem: inclusion order}. We thus deduce from Lemma~\ref{lem: Ext subquotient} that the surjection $V_{I_1}^{\rm{an}}\twoheadrightarrow V$ induces a surjection
\begin{equation}\label{equ: coxeter filtration x surjection prime}
\mathrm{Ext}_{G}^{\#I_0\setminus I_1}(V_{I_0}^{\infty},V_{I_1}^{\rm{an}})\twoheadrightarrow \mathrm{Ext}_{G}^{\#I_0\setminus I_1}(V_{I_0}^{\infty},V)
\end{equation}
whose kernel is (\ref{equ: coxeter filtration x surjection kernel}). Thanks to Lemma~\ref{lem: x Ext with sm}, we know that
\[\mathrm{Ext}_{G}^{k}(V_{I_0}^{\infty},W)=0\]
for each $W\in\mathrm{JH}_{G}(V_{x,I_1}^{\rm{an}}/V)$ and $k\leq \#I_0\setminus I_1$, and thus
\[\mathrm{Ext}_{G}^{k}(V_{I_0}^{\infty},V_{x,I_1}^{\rm{an}}/V)=0\]
for each $k\leq \#I_0\setminus I_1$ by a d\'evissage with respect to $\mathrm{JH}_{G}(V_{x,I_1}^{\rm{an}}/V)$.
A further d\'evissage with respect to $0\rightarrow V\rightarrow V_{x,I_1}^{\rm{an}}\rightarrow V_{x,I_1}^{\rm{an}}/V\rightarrow 0$ ensures that the embedding $V\hookrightarrow V_{x,I_1}^{\rm{an}}$ induces an isomorphism
\[\mathrm{Ext}_{G}^{\#I_0\setminus I_1}(V_{I_0}^{\infty},V)\buildrel\sim\over\longrightarrow \mathrm{Ext}_{G}^{\#I_0\setminus I_1}(V_{I_0}^{\infty},V_{x,I_1}^{\rm{an}}).\]
This together with our description of (\ref{equ: coxeter filtration x surjection prime}) finishes the proof.
\end{proof}

Proposition~\ref{prop: coxeter filtration x surjection} has the following immediate consequences.
\begin{cor}\label{cor: coxeter filtration as kernel}
Let $I_1\subseteq I_0\subseteq \Delta$ and $x\in\Gamma^{I_0\setminus I_1}$. Then we have
\begin{equation}\label{equ: coxeter filtration as kernel}
\mathrm{Fil}_{x}(\mathbf{E}_{I_0,I_1})=\bigcap_{y\not\unlhd x}\mathrm{ker}(\zeta_{I_0,I_1}^{y}).
\end{equation}
\end{cor}

\begin{cor}\label{cor: coxeter filtration x w surjection}
Let $I_1\subseteq I_0\subseteq \Delta$ and $x,w\in\Gamma^{I_0\setminus I_1}$ with $x\unlhd w$. Then the non-zero map $V_{x,I_1}^{\rm{an}}\rightarrow V_{w,I_1}^{\rm{an}}$ (see (\ref{equ: St x w transfer})) induces a surjection
\begin{equation}\label{equ: coxeter filtration x w surjection}
\mathrm{Ext}_{G}^{\#I_0\setminus I_1}(V_{I_0}^{\infty},V_{x,I_1}^{\rm{an}})\twoheadrightarrow \mathrm{Ext}_{G}^{\#I_0\setminus I_1}(V_{I_0}^{\infty},V_{w,I_1}^{\rm{an}}).
\end{equation}
\end{cor}
\begin{proof}
Now that the non-zero map $V_{I_1}^{\rm{an}}\rightarrow V_{w,I_1}^{\rm{an}}$ factors through $V_{x,I_1}^{\rm{an}}$ (by the discussion above (\ref{equ: St x w transfer})), we know that $\zeta_{I_0,I_1}^{w}$ factors as the composition of $\zeta_{I_0,I_1}^{x}$ with (\ref{equ: coxeter filtration x w surjection}). The surjectivity of (\ref{equ: coxeter filtration x w surjection}) thus follows from that of $\zeta_{I_0,I_1}^{w}$ by Proposition~\ref{prop: coxeter filtration x surjection} (upon replacing $x$ in \emph{lod.cit.} with $w$).
\end{proof}

Let $I_1\subseteq I_0\subseteq \Delta$ and $x\in\Gamma^{I_0\setminus I_1}$.
Recall from (\ref{equ: Ext Tits transfer x prime}) (upon replacing $I_2$ and $J_2$ in \emph{loc.cit.} with $I_1$ and $I_1\sqcup(\Delta\setminus I_0)=\Delta\setminus(I_0\setminus I_1)$), we see that the unique (up to scalars) non-zero maps $V_{I_1}^{\rm{an}}\rightarrow $ (cf.~Lemma~\ref{lem: Hom x w St}) induce the following commutative diagram
\begin{equation}\label{equ: Ext Tits transfer x bottom}
\xymatrix{
\mathbf{E}_{I_0,I_1} \ar[r] \ar[d] & \mathrm{Ext}_{G}^{\#I_0\setminus I_1}(V_{I_0}^{\infty},V_{I_1}^{\rm{an}}) \ar[r] \ar[d] & \mathrm{Ext}_{G}^{\#I_0\setminus I_1}(V_{I_0}^{\infty},V_{x,I_1}^{\rm{an}}) \ar[d]\\
\mathrm{Ext}_{G}^{\#\Delta\setminus I_1}(V_{I_0}^{\rm{an}},V_{\Delta\setminus(I_0\setminus I_1)}^{\rm{an}}) \ar[r] & \mathrm{Ext}_{G}^{\#\Delta\setminus I_1}(V_{I_0}^{\infty},V_{\Delta\setminus(I_0\setminus I_1)}^{\rm{an}}) \ar[r] & \mathrm{Ext}_{G}^{\#\Delta\setminus I_1}(V_{I_0}^{\infty},V_{x,\Delta\setminus(I_0\setminus I_1)}^{\rm{an}})\\
\mathbf{E}_{I_0\setminus I_1} \ar[r] \ar[u] & \mathrm{Ext}_{G}^{\#I_0\setminus I_1}(1_{G},V_{\Delta\setminus(I_0\setminus I_1)}^{\rm{an}}) \ar[r] \ar[u]& \mathbf{E}_{x,I_0\setminus I_1} \ar[u]
}
\end{equation}
with all vertical maps being isomorphisms, the horizontal maps from the first column to the second column being isomorphisms, and the horizontal maps from the second column to the third column being surjections by Proposition~\ref{prop: coxeter filtration x surjection}. Using (\ref{equ: coxeter filtration as kernel}) as well as
\[\mathrm{Fil}_{x}(\mathbf{E}_{I_0\setminus I_1})=\mathrm{Fil}_{x}(\mathbf{E}_{\Delta,\Delta\setminus(I_0\setminus I_1)})=\bigcap_{y\not\unlhd x}\mathrm{ker}(\zeta_{\Delta,\Delta\setminus(I_0\setminus I_1)}^{y}),\]
we deduce from (\ref{equ: Ext Tits transfer x bottom}) the following result.
\begin{cor}\label{cor: coxeter x transfer}
Let $I_1\subseteq I_0\subseteq \Delta$. The isomorphism $\mathbf{E}_{I_0,I_1}\cong\mathbf{E}_{I_0\setminus I_1}$ from \ref{it: Ext complex std seq} of Proposition~\ref{prop: Ext complex std seq} (see also the first column of (\ref{equ: Ext Tits transfer x bottom})) restricts to an isomorphism
\begin{equation}\label{equ: coxeter x transfer}
\mathrm{Fil}_{x}(\mathbf{E}_{I_0,I_1})\cong\mathrm{Fil}_{x}(\mathbf{E}_{I_0\setminus I_1})
\end{equation}
for each $x\in\Gamma^{I_0\setminus I_1}$.
\end{cor}

Let $I_1\subseteq\Delta$ and $x\in\Gamma_{\Delta\setminus I_1}$. Now that $J_{x}=I_1$, we have $V_{x,I_1}^{\rm{an}}=i_{x,I_1}^{\rm{an}}$ and thus obtain a unique (up to scalars) non-zero map $V_{I_1}^{\rm{an}}\rightarrow i_{x,I_1}^{\rm{an}}$ (see (\ref{equ: x St surjection})) which fits into the following commutative diagram
\begin{equation}\label{equ: x St surjection Lie diagram}
\xymatrix{
M^{I_1}(x) \ar[r] \ar@{^{(}->}[d] & \mathfrak{v}_{I_1} \ar@{^{(}->}[d]\\
(i_{x,I_1}^{\rm{an}})^{\vee} \ar[r] & (V_{I_1}^{\rm{an}})^{\vee}
}
\end{equation}
with the vertical maps being the natural $U(\fg)$-equivariant embeddings (cf.~(\ref{equ: Lie to loc an Tits})).
\begin{lem}\label{lem: grade 0 Lie transfer}
Let $I_1\subseteq\Delta$ with $h\defeq\#\Delta\setminus I_1$. Then the commutative diagram (\ref{equ: x St surjection Lie diagram}) for each $x\in\Gamma_{\Delta\setminus I_1}$ together with $U(\fg)\subseteq D(G)$ induces a commutative diagram of the following form
\begin{equation}\label{equ: grade 0 Lie transfer}
\xymatrix{
\mathbf{E}_{\Delta\setminus I_1} \ar@{->>}[r] \ar@{->>}[d] & \bigoplus_{x\in \Gamma_{\Delta\setminus I_1}}\mathrm{Ext}_{G}^{h}(1_{G},i_{x,I_1}^{\rm{an}}) \ar^{\wr}[d]\\
\mathrm{Ext}_{U(\fg)}^{h}(\mathfrak{v}_{I_1},1_{\fg}) \ar^{\sim}[r] & \bigoplus_{x\in \Gamma_{\Delta\setminus I_1}}\mathfrak{e}_{I_1,x}
}.
\end{equation}
\end{lem}
\begin{proof}
It is clear that the commutative diagram (\ref{equ: x St surjection Lie diagram}), upon taking direct sum over $x\in \Gamma_{\Delta\setminus I_1}$, induces a commutative diagram of the form (\ref{equ: grade 0 Lie transfer}) without knowing the maps in \emph{loc.cit.} being surjections or isomorphisms. The bottom horizontal map of (\ref{equ: grade 0 Lie transfer}) is an isomorphism by \ref{it: Lie St Ext w 2} of Proposition~\ref{prop: Lie St Ext w decomposition}, upon replacing $I_0$, $x$ and $u$ in \emph{loc.cit.} with $I_1$, $1$ and $x$ here. The top horizontal map of (\ref{equ: grade 0 Lie transfer}) is an surjection that factors through
\[\mathbf{E}_{\Delta\setminus I_1}\twoheadrightarrow \bigoplus_{x\in \Gamma_{\Delta\setminus I_1}}\mathrm{gr}_{x}(\mathbf{E}_{\Delta\setminus I_1})\buildrel\sim\over\longrightarrow \bigoplus_{x\in \Gamma_{\Delta\setminus I_1}}\mathrm{Ext}_{G}^{h}(1_{G},i_{x,I_1}^{\rm{an}})\]
by Proposition~\ref{prop: coxeter filtration x surjection} (upon taking $I_0=\Delta$ in \emph{loc.cit.}).
The right vertical map of (\ref{equ: grade 0 Lie transfer}) is the direct sum over $x\in \Gamma_{\Delta\setminus I_1}$ of the composition of isomorphisms
\[\mathrm{Ext}_{G}^{h}(1_{G},i_{x,I_1}^{\rm{an}})\buildrel\sim\over\longrightarrow \mathrm{Ext}_{D(\fg,P_{I_1})}^{h}(M^{I_1}(x),1_{D(\fg,P_{I_1})})\buildrel\sim\over\longrightarrow\mathfrak{e}_{I_1,x}\]
thanks to Lemma~\ref{lem: g P Ext transfer} and the spectral sequence
\[H^{k}(P_{I_1},1_{P_{I_1}})^{\infty}\otimes_E\mathrm{Ext}_{U(\fg)}^{\ell}(M^{I_1}(x),1_{\fg})\implies \mathrm{Ext}_{D(\fg,P_{I_1})}^{k+\ell}(M^{I_1}(x),1_{D(\fg,P_{I_1})}).\]
The proof is thus finished.
\end{proof}
\subsection{Cup product of Tits double complex}\label{subsec: cup Tits double}
Let $I_0,I_0'\subseteq \Delta$ with $I_0\cap I_0'=\emptyset$. Let $x\in\Gamma^{I_0}$, $y\in\Gamma^{I_0'}$ and $w\in\Gamma_{x,y}\subseteq\Gamma^{I_0\sqcup I_0'}$. We define our \emph{normalized cup product map}
\[
\xymatrix{
\mathbf{E}_{x,I_0} \ar@{=}[d] & \otimes_E & \mathbf{E}_{y,I_0'} \ar^{\cup}[rr] \ar@{=}[d] & &\mathbf{E}_{w,I_0\sqcup I_0'} \ar@{=}[d]\\
\mathrm{Ext}_{G}^{\#I_0}(1_{G},V_{x,\Delta\setminus I_0}^{\rm{an}}) & \otimes_E & \mathrm{Ext}_{G}^{\#I_0'}(1_{G},V_{y,\Delta\setminus I_0'}^{\rm{an}}) \ar^{\cup}[rr] & &\mathrm{Ext}_{G}^{\#I_0\sqcup I_0'}(1_{G},V_{w,\Delta\setminus(I_0\sqcup I_0')}^{\rm{an}})
}
\]
which is functorial with respect to the choice of $x$, $y$ and $w$ above under the $\unlhd$ partial-order (see (\ref{equ: x y cup w bottom diagram})). As an immediate consequence, we obtain some constraint between the cup product map $\kappa_{I_0,I_0'}$ and the Coxeter filtration on $\mathbf{E}_{I_0}$, $\mathbf{E}_{I_0'}$ and $\mathbf{E}_{I_0\sqcup I_0'}$ (see Proposition~\ref{prop: coxeter cup comparison}).

For each $I\subseteq \Delta$, we write $R_{I}\defeq D(\fg,P_{I})$ for short and use the shortened notation
\[C_{I}^{\bullet}(D)\defeq \Hom_{R_{I}}(B_{\bullet}(R_{I},D),1_{R_{I}})\]
for each $D\in\mathrm{Mod}_{R_{I}}$.
Let $x\in W(G)$ and $I_0\subseteq I_1\subseteq J_{x}=\Delta\setminus\mathrm{Supp}(x)$. 
Recall from (\ref{equ: x Tits double flat}) that we have a natural map between double complex
\begin{equation}\label{equ: double g P Ext transfer}
\cS_{x,I_0,I_1}^{\bullet,\bullet} \rightarrow \cS_{x,I_0,I_1,\flat}^{\bullet,\bullet}
\end{equation}
which induces an isomorphism on the first page. In particular, (\ref{equ: double g P Ext transfer}) induces an isomorphism
\begin{equation}\label{equ: double g P Ext transfer total}
H^{k}(\mathrm{Tot}(\cS_{x,I_0,I_1}^{\bullet,\bullet})) \buildrel\sim\over\longrightarrow H^{k}(\mathrm{Tot}(\cS_{x,I_0,I_1,\flat}^{\bullet,\bullet}))
\end{equation}
for each $k\in\Z$.

In the rest of this section, we fix a choice of $I_0,I_0'\subseteq\Delta$ with $I_0\cap I_0'=\Delta$ and write $J_0\defeq \Delta\setminus I_0$ and $J_0'\defeq\Delta\setminus I_0'$ for short. We write $r\defeq\#\Delta=n-1$ for short.
For each $u\in W(G)$ and $I\subseteq I_{u}$, we fix the choice of a surjection (which is unique up to scalars) 
\begin{equation}\label{equ: fix u Verma surjection}
M(u)\twoheadrightarrow M^{I}(u)
\end{equation}
for convenience.

Let $x\in\Gamma^{I_0}$, $y\in\Gamma^{I_0'}$ and $w\in\Gamma_{x,y}$ (see (\ref{equ: x y envelop})) with $\mathrm{Supp}(w)=\mathrm{Supp}(x)\sqcup\mathrm{Supp}(y)$. 
We write $S_{x,J_0}\defeq \{I\mid J_0\subseteq I\subseteq J_{x}\}$, $S_{y,J_0'}\defeq \{I'\mid J_0'\subseteq I'\subseteq J_{y}\}$ and $S_{w,J_0\cap J_0'}\defeq \{I''\mid J_0\cap J_0'\subseteq I''\subseteq J_{w}\}$ for short. We omit $x$ (resp.~$y$, resp.~$w$) from the notation when $x=1$ (resp.~when $y=1$, resp.~when $w=1$).
Note that the following map
\begin{equation}\label{equ: x y w index bijection}
S_{x,J_0}\times S_{y,J_0'}\rightarrow S_{w,J_0\cap J_0'}: (I,I')\mapsto I\cap I'
\end{equation}
is a well-defined bijection.
We fix the choice of an embedding from $M(x)$, $M(y)$ and $M(w)$ into $M(1)$. This together with our fixed surjections (\ref{equ: fix u Verma surjection}) determine a unique (possibly zero) map $M^{I}(x)\rightarrow M^{I}(1)$ (resp.~$M^{I'}(y)\rightarrow M^{I'}(1)$, resp.~$M^{I''}(w)\rightarrow M^{I''}(1)$) for each $I\in S_{x,J_0}$ (resp.~for each $I'\in S_{y,J_0'}$, resp.~for each $I''\in S_{w,J_0\cap J_0'}$).
Following Lemma~\ref{lem: tensor of parabolic Verma}, we obtain a commutative diagram of the form
\begin{equation}\label{equ: coxeter Verma tensor diagram}
\xymatrix{
M^{I\cap I'}(w) \ar[r] \ar[d] & M^{I}(x)\otimes_EM^{I'}(y) \ar[d]\\
M^{I\cap I'}(1) \ar[r] & M^{I}(1)\otimes_EM^{I'}(1)
}
\end{equation}
with all maps being actually non-zero (and unique based on our fixed choices of (\ref{equ: fix u Verma surjection}) and embedding from $M(x)$, $M(y)$ and $M(w)$ into $M(1)$).
Moreover, (\ref{equ: coxeter Verma tensor diagram}) is functorial with respect to the choice of $x$, $y$, $z$ as well as $I\in S_{x,J_0}$ and $I'\in S_{y,J_0'}$.

Let $I\in S_{x,J_0}$ and $I'\in S_{y,J_0'}$. 
Recall from (\ref{equ: std relative resolution general}) that
$B_{\bullet}(R_{I\cap I'},M^{I'}(y))$ is a resolution of $M^{I'}(y)$ by free $R_{I\cap I'}$-modules. This induces by $M^{I}(x)\otimes_E-$ a resolution $M^{I}(x)\otimes_EB_{\bullet}(R_{I\cap I'},M^{I'}(y))$ of $M^{I}(x)\otimes_EM^{I'}(y)$ as $R_{I\cap I'}$-modules, which further induces the following resolution
\[B_{\bullet}(R_{I\cap I'},M^{I}(x)\otimes_EB_{\bullet}(R_{I\cap I'},M^{I'}(y)))\]
of $M^{I}(x)\otimes_EM^{I'}(y)$ by free $R_{I\cap I'}$-modules. The map $B_{\bullet}(R_{I\cap I'},M^{I'}(y))\rightarrow M^{I'}(y)$ induces the following map between two resolutions of $M^{I}(x)\otimes_EM^{I'}(y)$ by free $R_{I\cap I'}$-modules
\begin{equation}\label{equ: x y free resolution transfer}
\mathrm{Tot}(B_{\bullet}(R_{I\cap I'},M^{I}(x)\otimes_EB_{\bullet}(R_{I\cap I'},M^{I'}(y))))\rightarrow B_{\bullet}(R_{I\cap I'},M^{I}(x)\otimes_EM^{I'}(y)).
\end{equation}
We write
\[
C_{I\cap I'}^{\bullet,\bullet}(M^{I}(x),M^{I'}(y))\defeq
\Hom_{R_{I\cap I'}}(B_{\bullet}(R_{I\cap I'},M^{I}(x)\otimes_EB_{\bullet}(R_{I\cap I'},M^{I'}(y))),1_{R_{I\cap I'}})
\]
for short.
Consequently, we obtain the following maps between complex
\begin{multline*}
C_{I'}^{\bullet}(M^{I'}(y))
\rightarrow C_{I\cap I'}^{\bullet}(M^{I'}(y))
\rightarrow \Hom_{R_{I\cap I'}}(M^{I}(x)\otimes_EB_{\bullet}(R_{I\cap I'},M^{I'}(y)),M^{I}(x))\\
\rightarrow \Hom_{R_{I\cap I'}}(\mathrm{Tot}(B_{\bullet}(R_{I\cap I'},M^{I}(x)\otimes_EB_{\bullet}(R_{I\cap I'},M^{I'}(y)))),B_{\bullet}(R_{I\cap I'},M^{I}(x))).
\end{multline*}
This together with (\ref{equ: x y free resolution transfer}), the restriction map
\[C_{I}^{\bullet}(M^{I}(x)) \rightarrow C_{I\cap I'}^{\bullet}(M^{I}(x))\]
and the composition with
\[\Hom_{R_{I\cap I'}}(B_{\bullet}(R_{I\cap I'},M^{I}(x)),1_{R_{I\cap I'}})\]
towards
\[\mathrm{Tot}(C_{I\cap I'}^{\bullet,\bullet}(M^{I}(x),M^{I'}(y)))\]
defines maps of the following form
\begin{equation}\label{equ: x y cup term}
\mathrm{Tot}(C_{I'}^{\bullet}(M^{I'}(y))\otimes_E C_{I}^{\bullet}(M^{I}(x)))
\rightarrow \mathrm{Tot}(C_{I\cap I'}^{\bullet,\bullet}(M^{I}(x),M^{I'}(y)))
\leftarrow C_{I\cap I'}^{\bullet}(M^{I}(x)\otimes_EM^{I'}(y))
\end{equation}
with the last map being a quasi-isomorphism and all maps in (\ref{equ: x y cup term}) being functorial with respect to the choice of $x$, $y$, $I$ and $I'$.
We define $\cN\cS_{x,y,J_0,J_0',\flat}^{\bullet,\bullet,\bullet,\bullet}$ as the $4$-complex whose $(\ell,\ell',k,k')$-term has the form
\[
\cN\cS_{x,y,J_0,J_0',\flat}^{\ell,\ell',k,k'}\defeq
\bigoplus_{I\in S_{x,J_0},I'\in S_{y,J_0'},\#\Delta\setminus I=\ell,\#\Delta\setminus I'=\ell'}C_{I\cap I'}^{k,k'}(M^{I}(x),M^{I'}(y)).
\]
We define $\cN\cS_{x,y,J_0,J_0'}^{\bullet,\bullet,\bullet}$ as the $3$-complex whose $(\ell,\ell',k)$-term has the form
\[\cN\cS_{x,y,J_0,J_0',\flat}^{\ell,\ell',k}\defeq\bigoplus_{I\in S_{x,J_0},I'\in S_{y,J_0'},\#\Delta\setminus I=\ell,\#\Delta\setminus I'=\ell'}C_{I\cap I'}^{k}(M^{I}(x)\otimes_EM^{I'}(y)).\]
Then (\ref{equ: x y cup term}) induces the following maps between complex
\begin{equation}\label{equ: x y cup total}
\mathrm{Tot}(\mathrm{Tot}(\cN\cS_{x,J_0,\flat}^{\bullet,\bullet})\otimes_E\mathrm{Tot}(\cN\cS_{y,J_0',\flat}^{\bullet,\bullet}))
\rightarrow \mathrm{Tot}(\cN\cS_{x,y,J_0,J_0',\flat}^{\bullet,\bullet,\bullet,\bullet})\leftarrow \mathrm{Tot}(\cN\cS_{x,y,J_0,J_0',\flat}^{\bullet,\bullet,\bullet})
\end{equation}
with the last map being a quasi-isomorphism, and both maps being functorial with respect to the choice of $x$ and $y$.
Following the discussion at the beginning of \S~\ref{subsec: Tits construction}, we define two functors from $\cP(\Delta)^{\rm{op}}$ to the additive category of complex of $E$-vector spaces $\cT_{x,y,J_0,J_0'}$ and $\cT_{w,J_0\cap J_0'}$ as $\cT$ in \emph{loc.cit.} by taking $\Sigma=(M_{I''}^{\bullet})_{I''\subseteq\Delta}$ so that $M_{I''}^{\bullet}=0$ when $I''\notin S_{J_0\cap J_0'}$ and $M_{I''}^{\bullet}$ equals $C_{I''}^{\bullet}(M^{I''\cup J_0}(x)\otimes_EM^{I''\cup J_0'}(y))$ and $C_{I''}^{\bullet}(M^{I''}(w))$ respectively when $I''\in S_{J_0\cap J_0'}$.
It is clear that the map $M^{I\cap I'}(w)\rightarrow M^{I}(x)\otimes_EM^{I'}(y)$ from (\ref{equ: coxeter Verma tensor diagram}) (which is functorial with respect to the choice of $x$, $y$, $w$, $I$ and $I'$) induces a natural transform between functors
\begin{equation}\label{equ: x y w functor transfer}
\cT_{x,y,J_0,J_0'}\rightarrow \cT_{w,J_0\cap J_0'}
\end{equation}
which is functorial with respect to the choice of $x$, $y$ and $w$.
Let $\sigma\in\mathfrak{S}_{\Delta}$ be an element that is compatible with both $I_0$ and $I_0'$ is moreover well-posed for the ordered pair $I_0,I_0'$ (see Definition~\ref{def: permutation preserve order}), then we have the following maps between complex (using (\ref{equ: x y w functor transfer}))
\begin{multline}\label{equ: x y w functor transfer complex}
\mathrm{Tot}(\cN\cS_{x,y,J_0,J_0',\flat}^{\bullet,\bullet,k})=\mathrm{Tot}_{\sigma}(\cT_{x,y,J_0,J_0'}^{k})\buildrel\sim\over\longrightarrow 
\mathrm{Tot}_{1}(\cT_{x,y,J_0,J_0'}^{k})\rightarrow \mathrm{Tot}_{1}(\cT_{w,J_0\cap J_0'}^{k})=\cN\cS_{w,J_0\cap J_0',\flat}^{\bullet,k}
\end{multline}
for each $k\geq 0$ (with the middle isomorphism given by $\zeta_{\sigma,1}:\mathrm{Tot}_{\sigma}(-)\rightarrow\mathrm{Tot}_{1}(-)$ from the discussion around (\ref{equ: general total differential})), which induces the following map between complex
\begin{equation}\label{equ: x y cup to w}
\mathrm{Tot}(\cN\cS_{x,y,J_0,J_0',\flat}^{\bullet,\bullet,\bullet})\rightarrow \mathrm{Tot}(\cN\cS_{w,J_0\cap J_0',\flat}^{\bullet,\bullet}).
\end{equation}
Combining (\ref{equ: x y cup total}) with (\ref{equ: x y cup to w}) and (\ref{equ: Tits double complex normalize coh}) (upon replacing $\cT_{\Sigma}^{\bullet,\bullet}$ in \emph{loc.cit.} with $\cS_{x,J_0,\flat}^{\bullet,\bullet}$, $\cS_{y,J_0',\flat}^{\bullet,\bullet}$ and $\cS_{w,J_0\cap J_0',\flat}^{\bullet,\bullet}$ respectively), we obtain the following maps for each $k_0,k_1\in\Z$
\begin{multline*}
H^{k_0-r}(\mathrm{Tot}(\cS_{x,J_0,\flat}^{\bullet,\bullet}))\otimes_E H^{k_1-r}(\mathrm{Tot}(\cS_{y,J_0',\flat}^{\bullet,\bullet}))
\cong H^{k_0}(\mathrm{Tot}(\cN\cS_{x,J_0,\flat}^{\bullet,\bullet}))\otimes_E H^{k_1}(\mathrm{Tot}(\cN\cS_{y,J_0',\flat}^{\bullet,\bullet}))\\
\rightarrow H^{k_0+k_1}(\mathrm{Tot}(\cN\cS_{x,y,J_0,J_0',\flat}^{\bullet,\bullet,\bullet,\bullet}))\buildrel\sim\over\longleftarrow 
H^{k_0+k_1}(\mathrm{Tot}(\cN\cS_{x,y,J_0,J_0',\flat}^{\bullet,\bullet,\bullet}))\\
\rightarrow H^{k_0+k_1}(\mathrm{Tot}(\cS_{w,J_0\cap J_0',\flat}^{\bullet,\bullet}))
\cong H^{k_0+k_1-r}(\mathrm{Tot}(\cS_{w,J_0\cap J_0',\flat}^{\bullet,\bullet})).
\end{multline*}
whose composition gives the following map
\begin{equation}\label{equ: x y cup w Tits}
H^{k_0-r}(\mathrm{Tot}(\cS_{x,J_0,\flat}^{\bullet,\bullet}))\otimes_E H^{k_1-r}(\mathrm{Tot}(\cS_{y,J_0',\flat}^{\bullet,\bullet}))\rightarrow H^{k_0+k_1-r}(\mathrm{Tot}(\cS_{w,J_0\cap J_0',\flat}^{\bullet,\bullet}))
\end{equation}
which is functorial with respect to the choice of $x$, $y$ and $w$.
In particular, we obtain the following commutative diagram
\begin{equation}\label{equ: x y cup w Tits diagram}
\xymatrix{
H^{k_0-r}(\mathrm{Tot}(\cS_{J_0,\flat}^{\bullet,\bullet})) \ar[d] & \otimes_E & H^{k_1-r}(\mathrm{Tot}(\cS_{J_0',\flat}^{\bullet,\bullet})) \ar^{\cup}[r] \ar[d] & H^{k_0+k_1-r}(\mathrm{Tot}(\cS_{J_0\cap J_0',\flat}^{\bullet,\bullet})) \ar[d] \\
H^{k_0-r}(\mathrm{Tot}(\cS_{x,J_0,\flat}^{\bullet,\bullet})) & \otimes_E & H^{k_1-r}(\mathrm{Tot}(\cS_{y,J_0',\flat}^{\bullet,\bullet})) \ar^{\cup}[r] & H^{k_0+k_1-r}(\mathrm{Tot}(\cS_{w,J_0\cap J_0',\flat}^{\bullet,\bullet}))
}
\end{equation}
Given three reductive $E$-Lie subalgebras $\fh_0$, $\fh_0'$ and $\fh_0''$ of $\fg$ satisfying $\fh_0''\subseteq\fh_0\cap\fh_0'$, we can define a cup product map
\begin{equation}\label{equ: x y cup w Tits relative}
H^{k_0-r}(\mathrm{Tot}(\cS_{x,J_0,\fh_0,\flat}^{\bullet,\bullet}))\otimes_E H^{k_1-r}(\mathrm{Tot}(\cS_{y,J_0',\fh_0',\flat}^{\bullet,\bullet}))\rightarrow H^{k_0+k_1-r}(\mathrm{Tot}(\cS_{w,J_0\cap J_0',\fh_0'',\flat}^{\bullet,\bullet}))
\end{equation}
by a parallel argument, and moreover we have a commutative diagram mapping (\ref{equ: x y cup w Tits relative}) to (\ref{equ: x y cup w Tits}).
Thanks to Lemma~\ref{lem: g P Ext transfer} and (\ref{equ: double g P Ext transfer total}), we obtain maps having parallel form as (\ref{equ: x y cup w Tits}), (\ref{equ: x y cup w Tits diagram}) and (\ref{equ: x y cup w Tits relative}) upon removing the subscript $\flat$ everywhere.

By taking $M_0=M^{I'}(y)$, $M_1=M^{I}(x)$ and $M_2=M^{I\cap I'}(w)$ in (\ref{equ: Lie cup left}), we have the following map in the derived category of $E$-vector spaces
\[
\mathrm{Tot}(C^{\bullet}(\fg,M^{I'}(y))\otimes_EC^{\bullet}(\fg,M^{I}(x))) \dashrightarrow C^{\bullet}(\fg,M^{I\cap I'}(w))
\]
which by its very definition fits into the following commutative diagram (induced from the inclusions of $U(\fg)$ into $R_{I}=D(\fg,P_{I})$, $R_{I'}=D(\fg,P_{I'})$ and $R_{I\cap I'}=D(\fg,P_{I\cap I'})$)
\begin{equation}\label{equ: x y pass to Lie diagram}
\xymatrix{
\mathrm{Tot}(C_{I'}^{\bullet}(M^{I'}(y))\otimes_EC_{I}^{\bullet}(M^{I}(x))) \ar@{-->}[r] \ar[d] & C_{I\cap I'}^{\bullet}(M^{I\cap I'}(w)) \ar[d]\\
\mathrm{Tot}(C^{\bullet}(\fg,M^{I'}(y))\otimes_EC^{\bullet}(\fg,M^{I}(x))) \ar@{-->}[r] & C^{\bullet}(\fg,M^{I\cap I'}(w))
}
\end{equation}

Similar to (\ref{equ: x y cup term}) (see also (\ref{equ: abstract cup product composition}) and (\ref{equ: abstract cup product})), we also have the following maps
\begin{multline}\label{equ: Levi cup term}
\mathrm{Tot}(C^{\bullet}(L_{I'})\otimes_E C^{\bullet}(L_{I}))\\
=\mathrm{Tot}(B_{\bullet}(D(L_{I'})),1_{D(L_{I'})})\otimes_E \Hom_{D(L_{I})}(B_{\bullet}(D(L_{I})),1_{D(L_{I})}))\\
\rightarrow \Hom_{D(L_{I\cap I'})}(\mathrm{Tot}(B_{\bullet}(D(L_{I\cap I'}),B_{\bullet}(D(L_{I\cap I'}))),1_{D(L_{I\cap I'})})\\
\leftarrow \Hom_{D(L_{I\cap I'})}(B_{\bullet}(D(L_{I\cap I'})),1_{D(L_{I\cap I'})})=C^{\bullet}(L_{I\cap I'})
\end{multline}
with the last leftward map being a quasi-isomorphism and all maps in (\ref{equ: Levi cup term}) being functorial with respect to the choice of $I$ and $I'$.
For each $N\subseteq\Delta$, we use the shortened notation $\cT_{N}^{\bullet,\bullet}\defeq \cT_{\Delta\setminus N,\Delta}^{\bullet,\bullet}$, $E_{\bullet,N}^{\bullet,\bullet}\defeq E_{\bullet,\Delta\setminus N,\Delta}^{\bullet,\bullet}$ (with differential maps $d_{\bullet,N}^{\bullet,\bullet}\defeq d_{\bullet,\Delta\setminus N,\Delta}^{\bullet,\bullet}$) and finally $\cN\cT_{N}^{\bullet,\bullet}\defeq \cN\cT_{\Delta\setminus N,\Delta}^{\bullet,\bullet}$, whose associated spectral sequence is denoted by $\tld{E}_{\bullet,N}^{\bullet,\bullet}$ with differential maps $\tld{d}_{\bullet,N}^{\bullet,\bullet}$.
We define $\cN\cT_{I_0,I_0'}^{\bullet,\bullet,\bullet,\bullet}$ as the $4$-complex whose $(\ell,\ell',k,k')$-term has the form
\[\cN\cT_{I_0,I_0'}^{\ell,\ell',k,k'}=\bigoplus_{I\in S_{J_0},I'\in S_{J_0'},\#\Delta\setminus I=\ell,\#\Delta\setminus I'=\ell'}\Hom_{D(L_{I\cap I'})}(\mathrm{Tot}(B_{k}(D(L_{I\cap I'}),B_{k'}(D(L_{I\cap I'}))),1_{D(L_{I\cap I'})}).\]
We define $\cN\cT_{I_0,I_0'}^{\bullet,\bullet,\bullet}$ as the $3$-complex whose $(\ell,\ell',k)$-term has the form
\[\cN\cT_{I_0,I_0'}^{\ell,\ell',k}=\bigoplus_{I\in S_{J_0},I'\in S_{J_0'},\#\Delta\setminus I=\ell,\#\Delta\setminus I'=\ell'}C^{k}(L_{I\cap I'}).\]
Then (\ref{equ: Levi cup term}) induces the following maps between complex
\begin{equation}\label{equ: Levi cup total}
\mathrm{Tot}(\mathrm{Tot}(\cN\cT_{I_0}^{\bullet,\bullet})\otimes_E\mathrm{Tot}(\cN\cT_{I_0'}^{\bullet,\bullet}))
\rightarrow \mathrm{Tot}(\cN\cT_{I_0,I_0'}^{\bullet,\bullet,\bullet,\bullet})\leftarrow \mathrm{Tot}(\cN\cT_{I_0,I_0'}^{\bullet,\bullet,\bullet})\cong \mathrm{Tot}(\cN\cT_{I_0\sqcup I_0'}^{\bullet,\bullet})
\end{equation}
with the leftward map being a quasi-isomorphism, and the last isomorphism defined via similar argument as in the construction (\ref{equ: x y cup to w}), using the same $\sigma\in\mathfrak{S}_{\Delta}$ as in \emph{loc.cit.}.
Upon replacing $\cT_{\Sigma}^{\bullet,\bullet}$ in (\ref{equ: Tits double complex normalize}) with $\cT_{I_0}^{\bullet,\bullet}$, $\cT_{I_0'}^{\bullet,\bullet}$ and $\cT_{I_0\sqcup I_0'}^{\bullet,\bullet}$, we obtain an isomorphism between double complex
\begin{equation}\label{equ: Levi cup Tits normalize}
\cT_{\ast}^{\bullet,\bullet}[-r,0]\cong \cN\cT_{\ast}^{\bullet,\bullet}
\end{equation}
with $\ast$ being either $I_0$ or $I_0'$ or $I_0\sqcup I_0'$. Note that (\ref{equ: Levi cup Tits normalize}) induces an isomorphism between the first pages of the associated spectral sequences
\begin{equation}\label{equ: Levi cup Tits normalize E1}
E_{1,\ast}^{\bullet,k}[-r]\cong \tld{E}_{1,\ast}^{\bullet,k}
\end{equation}
which is the direct sum of the scalar automorphism $(-1)^{\sum_{j\in I}(j-1)}$ on $H^{k}(L_{I},1_{L_{I}})$ for each $k\geq 0$ and various $I$.
Combining (\ref{equ: Levi cup total}) with (\ref{equ: Levi cup Tits normalize}), we obtain the following maps
\begin{multline}\label{equ: Levi cup Tits}
H^{k_0-r}(\mathrm{Tot}(\cT_{I_0}^{\bullet,\bullet}))\otimes_E H^{k_1-r}(\mathrm{Tot}(\cT_{I_0'}^{\bullet,\bullet}))\\
\cong H^{k_0}(\mathrm{Tot}(\cN\cT_{I_0}^{\bullet,\bullet}))\otimes_E H^{k_1}(\mathrm{Tot}(\cN\cT_{I_0'}^{\bullet,\bullet}))\\
\rightarrow H^{k_0+k_1}(\mathrm{Tot}(\cN\cT_{I_0,I_0'}^{\bullet,\bullet,\bullet,\bullet}))\buildrel\sim\over\longleftarrow H^{k_0+k_1}(\mathrm{Tot}(\cN\cT_{I_0,I_0'}^{\bullet,\bullet,\bullet}))
\buildrel\sim\over\longrightarrow H^{k_0+k_1}(\mathrm{Tot}(\cN\cT_{I_0\sqcup I_0'}^{\bullet,\bullet}))\\
\cong H^{k_0+k_1-r}(\mathrm{Tot}(\cT_{I_0\sqcup I_0'}^{\bullet,\bullet}))
\end{multline}
for each $k_0,k_1\in\Z$.
For each $\#J_0\leq \ell_0\leq r$ and $\#J_0'\leq\ell_1\leq r$ with $\ell_2\defeq \ell_0+\ell_1-r$, $\ell\defeq r-\ell_0$ and $\ell'\defeq r-\ell_1$, the map (\ref{equ: Levi cup total}) induces the following commutative diagram
\begin{equation}\label{equ: Levi cup E1}
\xymatrix{
E_{1,I_0}^{-\ell_0,\ell_0+k_0-r} \ar[d] & \otimes_E & E_{1,I_0'}^{-\ell_1,\ell_1+k_1-r} \ar[d] \ar^{\cup}[rr]
& & E_{1,I_0\sqcup I_0'}^{-\ell_2,\ell_2+k_0+k_1-r}\\
\tld{E}_{1,I_0}^{\ell,k_0-\ell} & \otimes_E & \tld{E}_{1,I_0'}^{\ell',k_1-\ell'}
\ar^{\cup}[rr] & & \tld{E}_{1,I_0\sqcup I_0'}^{\ell+\ell',k_0+k_1-(\ell+\ell')} \ar[u]
}
\end{equation}
with all vertical maps from (\ref{equ: Levi cup Tits normalize E1}) and the bottom horizontal cup product map nothing but the direct sum over all $I\in S_{J_0}$ with $\#I=\ell_0$ and all $I'\in S_{J_0'}$ with $\#I'=\ell_1$ of the composition of the following maps
\begin{multline*}
H^{k_0-\ell}(L_{I},1_{L_{I}})\otimes_EH^{k_1-\ell'}(L_{I'},1_{L_{I'}})\\
\rightarrow H^{k_0-\ell}(L_{I\cap I'},1_{L_{I\cap I'}})\otimes_EH^{k_1-\ell'}(L_{I\cap I'},1_{L_{I\cap I'}})\buildrel\cup\over\longrightarrow H^{k_0+k_1-(\ell+\ell')}(L_{I\cap I'},1_{L_{I\cap I'}}).
\end{multline*}
For each $I\in S_{J_0}$ and $I'\in S_{J_0'}$, we note that $I\cup I'\supseteq J_0\cup J_0'=\Delta$ and thus
\[(-1)^{\sum_{j\in I}(j-1)}(-1)^{\sum_{j\in I'}(j-1)}=\varepsilon_{\Delta}(-1)^{\sum_{j\in I\cap I'}(j-1)}\]
with $\varepsilon_{\Delta}\defeq (-1)^{\sum_{j\in \Delta}(j-1)}=(-1)^{\frac{(n-1)(n-2)}{2}}$.
Hence, the top horizontal cup product map is the direct sum over all $I\in S_{J_0}$ with $\#I=\ell_0$ and all $I'\in S_{J_0'}$ with $\#I'=\ell_1$ of the composition of the following maps
\begin{multline}\label{equ: Levi cup coh sign}
H^{k_0-\ell}(L_{I},1_{L_{I}})\otimes_EH^{k_1-\ell'}(L_{I'},1_{L_{I'}})\\
\rightarrow H^{k_0-\ell}(L_{I\cap I'},1_{L_{I\cap I'}})\otimes_EH^{k_1-\ell'}(L_{I\cap I'},1_{L_{I\cap I'}})\rightarrow H^{k_0+k_1-(\ell+\ell')}(L_{I\cap I'},1_{L_{I\cap I'}})
\end{multline}
with the last map differ from the usual cup product by the sign $\varepsilon_{\Delta}$ (which is independent of the choice of $I\in S_{J_0}$ and $I'\in S_{J_0'}$).
Note that (\ref{equ: Levi cup E1}) restricts to a commutative diagram of maps between various $\mathrm{ker}(d_{1,\ast}^{-,-})$ and $\mathrm{ker}(\tld{d}_{1,\ast}^{-,-})$, which further induces a commutative of the form
\begin{equation}\label{equ: Levi cup E2}
\xymatrix{
E_{2,I_0}^{-\ell_0,\ell_0+k_0-r} \ar[d] & \otimes_E & E_{2,I_0'}^{-\ell_1,\ell_1+k_1-r} \ar[d] \ar^{\cup}[rr]
& & E_{2,I_0\sqcup I_0'}^{-\ell_2,\ell_2+k_0+k_1-r}\\
\tld{E}_{2,I_0}^{\ell,k_0-\ell} & \otimes_E & \tld{E}_{2,I_0'}^{\ell',k_1-\ell'}
\ar^{\cup}[rr] & & \tld{E}_{2,I_0\sqcup I_0'}^{\ell+\ell',k_0+k_1-(\ell+\ell')} \ar[u]
}
\end{equation}
\begin{lem}\label{lem: grade cup commute}
Let $\#J_0\leq\ell_0\leq r$, $\#J_0'\leq\ell_1\leq r$ with $\ell_2\defeq\ell_0+\ell_1-r$ and $k_0,k_1\in\Z$. Then we have
\begin{equation}\label{equ: grade cup commute E1}
u\cup v=(-1)^{k_0k_1}v\cup u\in E_{1,I_0\sqcup I_0'}^{\ell_2,\ell_2+k_0+k_1-r}
\end{equation}
for each $u\in E_{1,I_0}^{-\ell_0,\ell_0+k_0-r}$ and $v\in E_{1,I_0'}^{-\ell_1,\ell_1+k_1-r}$.
In particular, we have
\begin{equation}\label{equ: grade cup commute E2}
x\cup y=(-1)^{k_0k_1}y\cup x\in E_{2,I_0\sqcup I_0'}^{\ell_2,\ell_2+k_0+k_1-r}
\end{equation}
for each $x\in E_{2,I_0}^{-\ell_0,\ell_0+k_0-r}$ and $y\in E_{2,I_0'}^{-\ell_1,\ell_1+k_1-r}$.
\end{lem}
\begin{proof}
Now that we have the commutative diagram (\ref{equ: Levi cup E2}) and its variant with $J_0$ and $J_0'$ exchanged (with a common sign $\varepsilon_{\Delta}$ appearing the top horizontal cup product map of both commutative diagrams, c.f.~(\ref{equ: Levi cup coh sign})), it suffices to show that we have 
\[
u\cup v=(-1)^{k_0k_1}v\cup u\in \tld{E}_{1,I_0\sqcup I_0'}^{\ell+\ell',k_0+k_1-(\ell+\ell')}
\]
for each $u\in \tld{E}_{1,I_0}^{\ell,k_0-\ell}$ and $v\in \tld{E}_{1,I_0'}^{\ell',k_1-\ell'}$, with $\ell\defeq r-\ell_0$ and $\ell'\defeq r-\ell_1$.
We view $\tld{E}_{1,I_0}^{\bullet,\bullet}$, $\tld{E}_{1,I_0'}^{\bullet,\bullet}$ and $\tld{E}_{1,I_0\sqcup I_0'}^{\bullet,\bullet}$ as double complex with differential maps on the $k$ direction being zero.
We write $H^{\bullet}_{I}\defeq H^{\bullet}(L_{I},1_{L_{I}})$ for short for each $I\subseteq \Delta$.
For each $I\in S_{J_0}$ and $I'\in S_{J_0'}$, we have the following maps
\begin{equation}\label{equ: Levi coh pair cup}
\mathrm{Tot}(H^{\bullet}_{I}\otimes_E H^{\bullet}_{I'})\rightarrow \mathrm{Tot}(H^{\bullet}_{I\cap I'}\otimes_E H^{\bullet}_{I\cap I'})\rightarrow H^{\bullet}_{I\cap I'}
\end{equation}
that sends $(x,x')$ to 
\[x\cup x'\defeq \mathrm{Res}^{k}_{I,I\cap I'}(x)\cup \mathrm{Res}^{k'}_{I',I\cap I'}(x')\in H^{k+k'}_{I\cap I'}\] 
for each $x\in H^{k}_{I}$, $x'\in H^{k'}_{I'}$ and $k,k'\geq 0$. Thanks to Proposition~\ref{prop: Levi decomposition} and Proposition~\ref{prop: Levi restriction}, we know that
\begin{equation}\label{equ: Levi coh cup exchange}
x\cup x'=(-1)^{kk'}x'\cup x
\end{equation}
for each $x\in H^{k}_{I}$, $x'\in H^{k'}_{I'}$ and $k,k'\geq 0$.
We use the shortened notation $\bigoplus_{I}$ for $\bigoplus_{I\in S_{J_0}}$, and $\bigoplus_{I'}$ for $\bigoplus_{I'\in S_{J_0'}}$.
We define a $4$-complex $H_{I_0,I_0'}^{\bullet,\bullet,\bullet,\bullet}$ by
\begin{equation}\label{equ: double Levi coh}
H_{I_0,I_0'}^{\ell,\ell',k,k'}=\bigoplus_{I, I', \#\Delta\setminus I=\ell,\#\Delta\setminus I'=\ell'}H^{k}_{I}\otimes_E H^{k'}_{I'}
\end{equation}
for each $\ell$, $\ell'$, $k$ and $k'$.
We use symbol $1$, $2$, $3$ and $4$ to label the index $\ell$, $\ell'$, $k$ and $k'$ as the formula (\ref{equ: double Levi coh}) indicates so that we could conveniently talk about various partial total complex (see the end of \S \ref{subsec: notation}) without ambiguity.
It is clear that we have identification between $4$-complex
\[H_{I_0,I_0'}^{\bullet_1,\bullet_2,\bullet_3,\bullet_4}=\tld{E}_{1,I_0}^{\bullet_1,\bullet_3}\otimes_E\tld{E}_{1,I_0'}^{\bullet_2,\bullet_4}\]
which leads to the following identifications between successive total complex
\[
\mathrm{Tot}(\mathrm{Tot}(\tld{E}_{1,I_0'}^{\bullet_2,\bullet_4})\otimes_E\mathrm{Tot}(\tld{E}_{1,I_0}^{\bullet_1,\bullet_3}))=\mathrm{Tot}_{2,4,1,3}(H_{I_0,I_0'}^{\bullet_1,\bullet_2,\bullet_3,\bullet_4})
\]
and
\[
\mathrm{Tot}(\mathrm{Tot}(\tld{E}_{1,I_0}^{\bullet_1,\bullet_3})\otimes_E\mathrm{Tot}(\tld{E}_{1,I_0'}^{\bullet_2,\bullet_4}))=\mathrm{Tot}_{1,3,2,4}(H_{I_0,I_0'}^{\bullet_1,\bullet_2,\bullet_3,\bullet_4}).
\]
We define a $3$-complex $H_{I_0,I_0'}^{\bullet,\bullet,\bullet}$ by
\[H_{I_0,I_0'}^{\ell,\ell',k}\defeq \bigoplus_{I, I', \#\Delta\setminus I=\ell, \#\Delta\setminus I'=\ell'}H^{k}_{I\cap I'}\]
for each $\ell$, $\ell'$ and $k$.
We have a natural map between $3$-complex
\[\mathrm{Tot}_{3,4}(H_{I_0,I_0'}^{\bullet,\bullet,\bullet,\bullet})\rightarrow H_{I_0,I_0'}^{\bullet,\bullet,\bullet}\]
by (\ref{equ: Levi coh pair cup}), and similarly another natural map between $3$-complex
\[\mathrm{Tot}_{4,3}(H_{I_0,I_0'}^{\bullet,\bullet,\bullet,\bullet})\rightarrow H_{I_0,I_0'}^{\bullet,\bullet,\bullet}\]
upon exchanging $I$ and $I'$ in (\ref{equ: Levi coh pair cup}).
We have the following commutative diagrams
\begin{equation}\label{equ: Levi coh exchange cup diagram 1}
\xymatrix{
\mathrm{Tot}(\mathrm{Tot}(\tld{E}_{1,I_0}^{\bullet,\bullet})\otimes_E\mathrm{Tot}(\tld{E}_{1,I_0'}^{\bullet,\bullet})) \ar@{=}[d] &
\mathrm{Tot}(\mathrm{Tot}(\tld{E}_{1,I_0}^{\bullet,\bullet})\otimes_E\mathrm{Tot}(\tld{E}_{1,I_0'}^{\bullet,\bullet})) \ar@{=}[d]\\
\mathrm{Tot}(\mathrm{Tot}_{1,3,2,4}(H_{I_0,I_0'}^{\bullet,\bullet,\bullet,\bullet})) \ar@{}[r]|*=0[@]{\cong} \ar@{}[d]|*=0[@]{\cong} &
\mathrm{Tot}(\mathrm{Tot}_{2,4,1,3}(H_{I_0,I_0'}^{\bullet,\bullet,\bullet,\bullet})) \ar@{}[d]|*=0[@]{\cong}\\
\mathrm{Tot}(\mathrm{Tot}_{1,2,3,4}(H_{I_0,I_0'}^{\bullet,\bullet,\bullet,\bullet})) \ar@{}[r]|*=0[@]{\cong} \ar[d] &
\mathrm{Tot}(\mathrm{Tot}_{2,1,4,3}(H_{I_0,I_0'}^{\bullet,\bullet,\bullet,\bullet})) \ar[d]\\
\mathrm{Tot}(\mathrm{Tot}_{1,2}(H_{I_0,I_0'}^{\bullet,\bullet,\bullet})) \ar@{}[r]|*=0[@]{\cong} &
\mathrm{Tot}(\mathrm{Tot}_{2,1}(H_{I_0,I_0'}^{\bullet,\bullet,\bullet}))
}
\end{equation}
with the isomorphisms in the middle commutative square of (\ref{equ: Levi coh exchange cup diagram 1}) being natural isomorphisms between differential successive partial total complex (see $\zeta_{\sigma,\sigma'}$ for a pair of elements $\sigma,\sigma'\in\mathfrak{S}_{N}$ discussed at the end of \S \ref{subsec: notation}). 
Here the commutativity of the bottom square of (\ref{equ: Levi coh exchange cup diagram 1}) crucially uses (\ref{equ: Levi coh cup exchange}) and the fact that the isomorphism
\[\mathrm{Tot}_{3,4}(H_{I_0,I_0'}^{\bullet,\bullet,\bullet,\bullet})\cong \mathrm{Tot}_{4,3}(H_{I_0,I_0'}^{\bullet,\bullet,\bullet,\bullet})\]
is given by the direct sum of the automorphism
\[(-1)^{kk'}: H_{I}^{k}\otimes_EH_{I'}^{k'}\buildrel\sim\over\longrightarrow H_{I}^{k}\otimes_EH_{I'}^{k'}\]
over all $I$, $I'$, $k$ and $k'$.
Now we consider two elements $\sigma,\sigma'\in\mathfrak{S}_{\Delta}$ which are both compatible with $I_0$ and $I_0'$, with $\sigma$ (resp.~$\sigma'$) being well-posed with respect to the ordered pair $I_0,I_0'$ (resp.~the ordered pair $I_0',I_0$), following the terminology of Definition~\ref{def: permutation preserve order}.
Taking $\cT$ in (\ref{equ: general double complex normalized}) to be the functor with $\Sigma=(M_{I}^{\bullet})_{I\subseteq\Delta}$ in \emph{loc.cit.} characterized by $M_{I}^{\bullet}=0$ if $I\notin S_{J_0\cap J_0'}$ and $M_{I}^{\bullet}=H_{I}^{\bullet}$ if $I\in S_{J_0\cap J_0'}$, then we have
\[\mathrm{Tot}_{1,2}(H_{I_0,I_0'}^{\bullet,\bullet,k})=\mathrm{Tot}_{\sigma}(\cT^{k}),\]
\[\mathrm{Tot}_{2,1}(H_{I_0,I_0'}^{\bullet,\bullet,k})=\mathrm{Tot}_{\sigma'}(\cT^{k}),\]
and
\[\tld{E}_{1,I_0\sqcup I_0'}^{\bullet,k}=\mathrm{Tot}_{1}(\cT^{k})\]
for each $k\geq 0$, which means we have the following commutative diagram
\begin{equation}\label{equ: double complex comparison}
\xymatrix{
\mathrm{Tot}_{1,2}(H_{I_0,I_0'}^{\bullet,\bullet,k}) \ar^{\zeta_{\sigma,1}}[rr] \ar^{\zeta_{\sigma,\sigma'}}[d] & & \tld{E}_{1,I_0\sqcup I_0'}^{\bullet,k} \ar@{=}[d]\\
\mathrm{Tot}_{2,1}(H_{I_0,I_0'}^{\bullet,\bullet,k}) \ar^{\zeta_{\sigma',1}}[rr] & & \tld{E}_{1,I_0\sqcup I_0'}^{\bullet,k}
}
\end{equation}
for $k\geq 0$, with the isomorphism $\zeta_{\sigma,\sigma'}$ given by the direct sum of the automorphism
\[(-1)^{\ell\ell'}: H_{I''}^{k}\buildrel\sim\over\longrightarrow H_{I''}^{k}\]
over all $I''\in S_{J_0\cap J_0'}$ satisfying $\#\Delta\setminus(I''\cup J_0)=\ell$ and $\#\Delta\setminus(I''\cup J_0')=\ell'$.
To summarize, we obtain from (\ref{equ: Levi coh exchange cup diagram 1}) and (\ref{equ: double complex comparison}) a map between complex
\[\mathrm{Tot}(\mathrm{Tot}(\tld{E}_{1,I_0}^{\bullet,\bullet})\otimes_E\mathrm{Tot}(\tld{E}_{1,I_0'}^{\bullet,\bullet})) \buildrel\cup\over\longrightarrow \mathrm{Tot}(\tld{E}_{1,I_0\sqcup I_0'}^{\bullet,\bullet})\]
as well as a map between complex
\[\mathrm{Tot}(\mathrm{Tot}(\tld{E}_{1,I_0'}^{\bullet,\bullet})\otimes_E\mathrm{Tot}(\tld{E}_{1,I_0}^{\bullet,\bullet})) \buildrel\cup\over\longrightarrow \mathrm{Tot}(\tld{E}_{1,I_0\sqcup I_0'}^{\bullet,\bullet})\]
which satisfy
\[u\cup v=(-1)^{k_0k_1}v\cup u \in (\mathrm{Tot}(\tld{E}_{1,I_0\sqcup I_0'}^{\bullet,\bullet}))^{k_0+k_1}\]
for each 
\[u\in (\mathrm{Tot}(\tld{E}_{1,I_0}^{\bullet,\bullet}))^{k_0}=\bigoplus_{\ell=0}^{\#I_0}\tld{E}_{1,I_0}^{\ell,k_0-\ell}\] 
and 
\[v\in (\mathrm{Tot}(\tld{E}_{1,I_0'}^{\bullet,\bullet}))^{k_1}=\bigoplus_{\ell'=0}^{\#I_0'}\tld{E}_{1,I_0'}^{\ell',k_1-\ell'}.\]
The proof is thus finished.
\end{proof}

When $k_0=2\#I_0$ and $k_1=2\#I_0'$, we recall from Proposition~\ref{prop: bottom deg degeneracy}) that $E_{2,I_0}^{-\ell_0,\ell_0+k_0-r}=E_{\infty,I_0}^{-\ell_0,\ell_0+k_0-r}$, $E_{2,I_0'}^{-\ell_1,\ell_1+k_1-r}=E_{\infty,I_0'}^{-\ell_1,\ell_1+k_1-r}$ and $E_{2,I_0\sqcup I_0'}^{-\ell_2,\ell_2+k_0+k_1-r}=E_{\infty,I_0\sqcup I_0'}^{-\ell_2,\ell_2+k_0+k_1-r}$ for each $\#J_0\leq \ell_0\leq r$ and $\#J_0'\leq \ell_1\leq r$ with $\ell_2=\ell_0+\ell_1-r$, and thus the composition of (\ref{equ: Levi cup Tits}) actually restricts to a map
\[\mathrm{Fil}^{-\ell_0}(H^{k_0-r}(\mathrm{Tot}(\cT_{I_0}^{\bullet,\bullet})))\otimes_E \mathrm{Fil}^{-\ell_1}(H^{k_1-r}(\mathrm{Tot}(\cT_{I_0'}^{\bullet,\bullet})))\rightarrow \mathrm{Fil}^{-\ell_2}(H^{k_0+k_1-r}(\mathrm{Tot}(\cT_{I_0\sqcup I_0'}^{\bullet,\bullet})))\]
for each $\#J_0\leq \ell_0\leq r$ and $\#J_0'\leq \ell_1\leq r$ with $\ell_2\defeq \ell_0+\ell_1-r$, whose induced maps between graded pieces are given by the top horizontal map of (\ref{equ: Levi cup E2}).

We use the shortened notation $A_{I}\defeq D(L_{I}')$ for each $I\subseteq\Delta$ below.
Similar to (\ref{equ: Levi cup term}), we also have the following map in the derived category of $E$-vector spaces
\begin{multline}\label{equ: Levi cup term natural}
\mathrm{Tot}(C^{\bullet}(L_{I'})_{\natural}\otimes_EC^{\bullet}(L_{I})_{\natural})\rightarrow \mathrm{Tot}(C^{\bullet}(L_{I\cap I'})_{\natural}\otimes_EC^{\bullet}(L_{I\cap I'})_{\natural})\\
\rightarrow \mathrm{Tot}(\Hom_{E}(\bigotimes_{j\in I\sqcup I'}B_{\bullet}(\Z_{j},B_{\bullet}(\Z_{j}))_{\Z_{j}},\Hom_{A_{I\cap I'}}(B_{\bullet}(A_{I\cap I'},B_{\bullet}(A_{I\cap I'})),1_{A_{I\cap I'}})))\\
\leftarrow C^{\bullet}(L_{I\cap I'})_{\natural}
\end{multline}
with all maps in (\ref{equ: Levi cup term natural}) being functorial with respect to the choice of $I$ and $I'$, and the last map in (\ref{equ: Levi cup term natural}) being a quasi-isomorphism.
Hence, we obtain a map
\begin{equation}\label{equ: Levi cup Tits natural}
H^{k_0-r}(\mathrm{Tot}(\cT_{J_0,\Delta,\natural}^{\bullet,\bullet}))\otimes_E H^{k_1-r}(\mathrm{Tot}(\cT_{J_0',\Delta,\natural}^{\bullet,\bullet}))\rightarrow H^{k_0+k_1-r}(\mathrm{Tot}(\cT_{J_0\cap J_0',\Delta,\natural}^{\bullet,\bullet}))
\end{equation}
for each $k_0,k_1\in\Z$.
Now that the map (see (\ref{equ: abstract relative cochain}) and (\ref{equ: Levi to natural}) with $\fh=0$ in \emph{loc.cit.}) $C^{\bullet}(L_{I})\rightarrow C^{\bullet}(L_{I})_{\natural}$ is the evaluation at $1_{L_{I}}$ of the following natural transform between functors
\[C^{\bullet}(L_{I},-)\rightarrow \mathrm{Tot}(\Hom_{E}(B_{\bullet}^{I},C^{\bullet}(L_{I}',-)))\]
and similar facts hold with $I'$ and $I\cap I'$ replacing $I$, we obtain a natural map from (\ref{equ: Levi cup term natural}) to (\ref{equ: Levi cup term}) which altogether forms a commutative diagram that is functorial with respect to the choice of $I$ and $I'$. In particular, we obtain a commutative diagram that maps the composition of (\ref{equ: Levi cup Tits}) to (\ref{equ: Levi cup Tits natural}).
Since the maps $\cT_{J_0,\Delta}^{\bullet,\bullet}\rightarrow \cT_{J_0,\Delta,\natural}^{\bullet,\bullet}$ induces an isomorphism on the first page (see the discussion around (\ref{equ: double Levi to natural}) with $\fh=0$ in \emph{loc.cit.}) and similar facts hold with $J_0'$ and $J_0\cap J_0'$ replacing $J_0$, the aforementioned map from (\ref{equ: Levi cup Tits}) to (\ref{equ: Levi cup Tits natural}) is actually an isomorphism.
For each $J\subseteq I_0$ and $J'\subseteq I_0'$, we recall from (\ref{equ: Levi cochain natural J}) that $\tau^{J}(C^{\bullet}(L_{I})_{\natural}))$ is defined as certain natural truncation of $C^{\bullet}(L_{I})_{\natural})$ and similar facts hold if we replace $I$ and $J$ with $I'$ and $J'$, or with $I\cap I'$ and $J\sqcup J'$.
Consequently, the maps (\ref{equ: Levi cup term natural}) restrict to maps
\begin{multline}\label{equ: Levi cup term natural J}
\mathrm{Tot}(\tau^{J'}(C^{\bullet}(L_{I'})_{\natural})\otimes_E\tau^{J}(C^{\bullet}(L_{I})_{\natural}))\rightarrow \mathrm{Tot}(\tau^{J'}(C^{\bullet}(L_{I\cap I'})_{\natural})\otimes_E\tau^{J}(C^{\bullet}(L_{I\cap I'})_{\natural}))\\
\rightarrow 
\mathrm{Tot}(\Hom_{E}(\bigotimes_{j\in I\sqcup I'}\sigma^{\leq \ast_j}\mathrm{Tot}(B_{\bullet}(\Z_{j},B_{\bullet}(\Z_{j}))_{\Z_{j}}),\Hom_{A_{I\cap I'}}(B_{\bullet}(A_{I\cap I'},B_{\bullet}(A_{I\cap I'})),1_{A_{I\cap I'}})))\\
\leftarrow \tau^{J\sqcup J'}(C^{\bullet}(L_{I\cap I'})_{\natural})
\end{multline}
with $\ast_j= -1$ if $j\in J\sqcup J'$ and $\ast_{j}= 0$ otherwise.
Note that all maps in (\ref{equ: Levi cup term natural J}) are functorial with respect to the choice of $I$ and $I'$, and the last map in (\ref{equ: Levi cup term natural J}) is a quasi-isomorphism.
Given $J\subseteq I_0$ and $J'\subseteq I_0'$, we use the shortened notation $\tau^{J}(\cT_{I_0}^{\bullet,\bullet})$ (resp.~$\tau^{J'}(\cT_{I_0'}^{\bullet,\bullet})$, resp.~$\tau^{J\sqcup J'}(\cT_{I_0\sqcup I_0'}^{\bullet,\bullet})$) for $\tau^{J}(\cT_{J_0,\Delta,\natural}^{\bullet,\bullet})$ (resp.~$\tau^{J'}(\cT_{J_0',\Delta,\natural}^{\bullet,\bullet})$, resp.~$\tau^{J\sqcup J'}(\cT_{J_0\cap J_0',\Delta,\natural}^{\bullet,\bullet})$) and similarly with $\mathrm{gr}^{\ast}(-)$ replacing $\tau^{\ast}(-)$.
Consequently, we deduce from (\ref{equ: Levi cup term natural J}) the following map
\begin{equation}\label{equ: Levi cup Tits natural J}
H^{k_0-r}(\mathrm{Tot}(\tau^{J}(\cT_{I_0}^{\bullet,\bullet})))\otimes_E H^{k_1-r}(\mathrm{Tot}(\tau^{J'}(\cT_{I_0'}^{\bullet,\bullet})))\rightarrow H^{k_0+k_1-r}(\mathrm{Tot}(\tau^{J\sqcup J'}(\cT_{I_0\sqcup I_0'}^{\bullet,\bullet})))
\end{equation}
for each $k_0,k_1\in\Z$. 
It is clear that we have a commutative diagram which maps each term of (\ref{equ: Levi cup Tits natural J}) to that of (\ref{equ: Levi cup Tits natural}). 
We also consider the maps
\begin{multline}\label{equ: Levi cup term natural graded J}
\mathrm{Tot}(\mathrm{gr}^{J'}(C^{\bullet}(L_{I'})_{\natural})\otimes_E\mathrm{gr}^{J}(C^{\bullet}(L_{I})_{\natural}))\rightarrow \mathrm{Tot}(\mathrm{gr}^{J'}(C^{\bullet}(L_{I\cap I'})_{\natural})\otimes_E\mathrm{gr}^{J}(C^{\bullet}(L_{I\cap I'})_{\natural}))\\
\rightarrow 
\mathrm{Tot}(\Hom_{E}(\bigotimes_{j\in J\sqcup J'}\sigma^{\leq -1}\mathrm{Tot}(B_{\bullet}(\Z_{j},B_{\bullet}(\Z_{j}))_{\Z_{j}}),\Hom_{A_{I\cap I'}}(B_{\bullet}(A_{I\cap I'},B_{\bullet}(A_{I\cap I'})),1_{A_{I\cap I'}})))\\
\leftarrow \mathrm{gr}^{J\sqcup J'}(C^{\bullet}(L_{I\cap I'})_{\natural})
\end{multline}
which induce a map of the form
\begin{equation}\label{equ: Levi cup Tits natural graded J}
H^{k_0-r}(\mathrm{Tot}(\mathrm{gr}^{J}(\cT_{I_0}^{\bullet,\bullet})))\otimes_E H^{k_1-r}(\mathrm{Tot}(\mathrm{gr}^{J'}(\cT_{I_0'}^{\bullet,\bullet})))\rightarrow H^{k_0+k_1-r}(\mathrm{Tot}(\mathrm{gr}^{J\sqcup J'}(\cT_{I_0\sqcup I_0'}^{\bullet,\bullet})))
\end{equation}
for each $J\subseteq I_0$, $J'\subseteq I_0'$ and each $k_0,k_1\in\Z$. 
Here we omit the tensor factor $\sigma^{\geq 0}\mathrm{Tot}(B_{\bullet}(\Z_{j},B_{\bullet}(\Z_{j}))_{\Z_{j}})\cong E$ for each $j\in (I\sqcup I')\setminus(J\sqcup J')$ from the expression of the middle term of (\ref{equ: Levi cup term natural graded J}).
Note that there is a commutative diagram that maps each term of (\ref{equ: Levi cup Tits natural J}) to (\ref{equ: Levi cup Tits natural graded J}).
When $k_0=2\#I_0$ and $k_1=2\#I_0'$, we recall from \ref{it: bottom deg graded 1} of Proposition~\ref{prop: bottom deg graded} that
\[\{H^{k_0-r}(\mathrm{Tot}(\tau^{J}(\cT_{I_0}^{\bullet,\bullet})))\}_{J\subseteq I_0}\]
forms a separated exhaustive decreasing filtration on $H^{k_0-r}(\mathrm{Tot}(\cT_{I_0}^{\bullet,\bullet}))$ with graded pieces given by $H^{k_0-r}(\mathrm{Tot}(\mathrm{gr}^{J}(\cT_{I_0}^{\bullet,\bullet})))$ for each $J\subseteq I_0$, and similar facts hold with $(I_0',k_1)$ and $(I_0\sqcup I_0',k_0+k_1)$ replacing $(I_0,k_0)$.
Consequently, assuming $k_0=2\#I_0$ and $k_1=2\#I_0'$, we know that (\ref{equ: Levi cup Tits natural}) restricts to (\ref{equ: Levi cup Tits natural J}) for each $J,J'$, which further induces the map (\ref{equ: Levi cup Tits natural graded J}) for each $J,J'$.

Assume temporarily that $J'=I_0'$ and $J=\emptyset$ in (\ref{equ: Levi cup term natural graded J}).
For each $\Delta\setminus I_0\subseteq I\subseteq\Delta$, the embedding between complex
\[E\cong \tau^{\leq 0}(C^{\bullet}(L_{I\setminus I_0'}')\hookrightarrow C^{\bullet}(L_{I\setminus I_0'}')\]
induces the embedding
\begin{equation}\label{equ: sm embeds Levi}
\Hom_{E}(\mathrm{gr}^{I_0'}(B_{\bullet}^{I}),E)\cong \Hom_{E}(\mathrm{gr}^{I_0'}(B_{\bullet}^{I\setminus I_0'}),E)\hookrightarrow \mathrm{gr}^{I_0'}(C^{\bullet}(L_{I\setminus I_0'})_{\natural}).
\end{equation}
Under the embedding (\ref{equ: sm embeds Levi}), the maps (\ref{equ: Levi cup term natural graded J}) (with $I'=\Delta\setminus I_0'$) pull back to
\begin{multline}\label{equ: Levi cup term natural graded J sm}
\mathrm{Tot}(\Hom_{E}(\mathrm{gr}^{I_0'}(B_{\bullet}^{I}),E)\otimes_EC^{\bullet}(L_{I}'))\\
\rightarrow \mathrm{Tot}(\Hom_{E}(\mathrm{gr}^{I_0'}(B_{\bullet}^{I\setminus I_0'}),E)\otimes_EC^{\bullet}(L_{I\setminus I_0'}'))\\
\rightarrow 
\mathrm{Tot}(\Hom_{E}(\bigotimes_{j\in I_0'}\sigma^{\leq -1}\mathrm{Tot}(B_{\bullet}(\Z_{j},B_{\bullet}(\Z_{j}))_{\Z_{j}}),C^{\bullet}(L_{I\setminus I_0'}')))\\
\leftarrow \mathrm{gr}^{I_0'}(C^{\bullet}(L_{I\setminus I_0'})_{\natural})=\Hom_{E}(\mathrm{gr}^{I_0'}(B_{\bullet}^{I\setminus I_0'}),C^{\bullet}(L_{I\setminus I_0'}')),
\end{multline}
using $\mathrm{gr}^{\emptyset}(C^{\bullet}(L_{I})_{\natural})\cong C^{\bullet}(L_{I}')$ and $\mathrm{gr}^{\emptyset}(C^{\bullet}(L_{I\setminus I_0'})_{\natural})\cong C^{\bullet}(L_{I\setminus I_0'}')$. 
We define $\mathrm{gr}^{\emptyset}(\cT_{I_0,\flat}^{\bullet,\bullet})$ (resp.~$\mathrm{CE}_{I_0,\flat}^{\bullet,\bullet}$) by replacing $C^{\bullet}(L_{I}')$ (resp.~$\mathrm{CE}^{\bullet}(\fl_{I})$) in the definition of $\mathrm{gr}^{\emptyset}(\cT_{I_0}^{\bullet,\bullet})$ (resp.~of $\mathrm{CE}_{I_0}^{\bullet,\bullet}\defeq\mathrm{CE}_{\Delta\setminus I_0,\Delta}^{\bullet,\bullet}$) with $C^{\bullet}(L_{I\setminus I_0'}')$ (resp.~$\mathrm{CE}^{\bullet}(\fl_{I\setminus I_0'})$).
We deduce from (\ref{equ: Levi cup term natural graded J}) the following maps between complex
\begin{multline}\label{equ: Levi cup Tits natural graded J sm}
\mathrm{Tot}(\Hom_{E}(\mathrm{gr}^{I_0'}(B_{\bullet}^{I}),E)\otimes_E\mathrm{gr}^{\emptyset}(\cT_{I_0}^{\bullet,\bullet}))
\rightarrow \mathrm{Tot}(\Hom_{E}(\mathrm{gr}^{I_0'}(B_{\bullet}^{I\setminus I_0'}),E)\otimes_E\mathrm{gr}^{\emptyset}(\cT_{I_0,\flat}^{\bullet,\bullet}))\\
\rightarrow \mathrm{Tot}(\Hom_{E}(\bigotimes_{j\in I_0'}\sigma^{\leq -1}\mathrm{Tot}(B_{\bullet}(\Z_{j},B_{\bullet}(\Z_{j}))_{\Z_{j}}),\mathrm{gr}^{\emptyset}(\cT_{I_0,\flat}^{\bullet,\bullet})))
\leftarrow \mathrm{gr}^{I_0'}(\cT_{I_0\sqcup I_0'}^{\bullet,\bullet}).
\end{multline}
Under the quasi-isomorphisms $C^{\bullet}(L_{I}')\rightarrow\mathrm{CE}^{\bullet}(\fl_{I})$ and $C^{\bullet}(L_{I\setminus I_0'}')\rightarrow\mathrm{CE}^{\bullet}(\fl_{I\setminus I_0'})$ (see Lemma~\ref{lem: center compact}) for each $\Delta\setminus I_0\subseteq I\subseteq\Delta$, the maps (\ref{equ: Levi cup Tits natural graded J sm}) admit quasi-isomorphisms to the maps
\begin{multline}\label{equ: Levi cup Tits natural graded J sm Lie}
\mathrm{Tot}(\Hom_{E}(\mathrm{gr}^{I_0'}(B_{\bullet}^{I}),E)\otimes_E\mathrm{CE}_{I_0}^{\bullet,\bullet})
\rightarrow \mathrm{Tot}(\Hom_{E}(\mathrm{gr}^{I_0'}(B_{\bullet}^{I\setminus I_0'}),E)\otimes_E\mathrm{CE}_{I_0,\flat}^{\bullet,\bullet})\\
\rightarrow \mathrm{Tot}(\Hom_{E}(\bigotimes_{j\in I_0'}\sigma^{\leq -1}\mathrm{Tot}(B_{\bullet}(\Z_{j},B_{\bullet}(\Z_{j}))_{\Z_{j}}),\mathrm{CE}_{I_0,\flat}^{\bullet,\bullet}))
\leftarrow \mathrm{Tot}(\Hom_{E}(\mathrm{gr}^{I_0'}(B_{\bullet}^{I\setminus I_0'}),\mathrm{CE}_{I_0,\flat}^{\bullet,\bullet})).
\end{multline}
To summarize, we have the following result.
\begin{lem}\label{lem: graded tensor as cup}
Let $I_0,I_0'\subseteq\Delta$ with $I_0\cap I_0'=\emptyset$, $h\defeq 2\#I_0-r$, $h'\defeq 2\#I_0\sqcup I_0'-r$ and $\mathbf{E}_{I_0'}^{\infty}\defeq \wedge^{\#I_0'}\Hom(Z_{\Delta\setminus I_0'}^{\dagger},E)$. We have the following commutative diagram of maps
\begin{equation}\label{equ: graded tensor as cup}
\xymatrix{
\mathrm{gr}^{\emptyset}(\mathbf{E}_{I_0}) \ar^{\sim}[rrr] & & & H^{h}(\mathrm{Tot}(\mathrm{CE}_{\Delta\setminus I_0,\Delta}^{\bullet,\bullet})) \\
\otimes_E & & & \otimes_E \\
\mathrm{gr}^{I_0'}(\mathbf{E}_{I_0'}) \ar[d] & & & \mathbf{E}_{I_0'}^{\infty} \ar_{\sim}[lll] \ar[d]\\
\mathrm{gr}^{I_0'}(\mathbf{E}_{I_0\sqcup I_0'}) \ar^{\sim}[rrr] & & & H^{h'-\#I_0'}(\mathrm{Tot}(\mathrm{CE}_{\Delta\setminus(I_0\sqcup I_0'),\Delta\setminus I_0'}^{\bullet,\bullet}))\otimes_E\mathbf{E}_{I_0'}^{\infty}
}
\end{equation}
with all the horizontal isomorphisms from (\ref{equ: grade J Z to CE}), the left vertical map from (\ref{equ: Levi cup Tits natural graded J}), and the right vertical map induced from the following maps between double complex
\begin{equation}\label{equ: graded tensor as cup double}
\mathrm{CE}_{\Delta\setminus I_0,\Delta}^{\bullet,\bullet}=\mathrm{CE}_{I_0}^{\bullet,\bullet}\rightarrow \mathrm{CE}_{I_0,\flat}^{\bullet,\bullet}\cong \mathrm{CE}_{\Delta\setminus(I_0\sqcup I_0'),\Delta\setminus I_0'}^{\bullet+\#I_0',\bullet}.
\end{equation}
\end{lem}
\begin{proof}
Now that (\ref{equ: Levi cup term natural graded J sm}) is constructed as the pull back of (\ref{equ: Levi cup term natural graded J}) (with $J=\emptyset$, $J'=I_0'$ and $I'=\Delta\setminus I_0'$ in \emph{loc.cit.}) via the embedding (\ref{equ: sm embeds Levi}), we know that the maps (\ref{equ: Levi cup Tits natural graded J sm}) (or rather their associated maps between cohomology) are compatible with (\ref{equ: Levi cup Tits natural graded J}). Since the map from (\ref{equ: Levi cup Tits natural graded J sm}) to (\ref{equ: Levi cup Tits natural graded J sm Lie}) is term by term quasi-isomorphism, we conclude the commutative diagram (\ref{equ: graded tensor as cup}). 
\end{proof}

When $x=y=w=1$, the inclusion $D(L_{I})\subseteq D(\fg,P_{I})$ and its variants with $I'$ and $I\cap I'$ replacing $I$ gives a commutative diagram that maps (\ref{equ: x y cup term}) to (\ref{equ: Levi cup term}) and this commutative diagram is functorial with respect to the choice of $I$ and $I'$.
Now that the natural maps $\cS_{J_0,\flat}^{\bullet,\bullet}\rightarrow \cT_{I_0}^{\bullet,\bullet}$, $\cS_{J_0'\flat}^{\bullet,\bullet}\rightarrow \cT_{I_0'}^{\bullet,\bullet}$ and $\cS_{J_0\cap J_0'\flat}^{\bullet,\bullet}\rightarrow \cT_{I_0\sqcup I_0'}^{\bullet,\bullet}$ all induce isomorphisms on the first page (using Lemma~\ref{lem: g P Ext transfer}), we obtain the following commutative diagram
\begin{equation}\label{equ: S T cup diagram}
\xymatrix{
H^{k_0-r}(\mathrm{Tot}(\cS_{J_0,\flat}^{\bullet,\bullet})) \ar^{\wr}[d] & \otimes_E & H^{k_1-r}(\mathrm{Tot}(\cS_{J_0',\flat}^{\bullet,\bullet})) \ar^{\cup}[r] \ar^{\wr}[d] & H^{k_0+k_1-r}(\mathrm{Tot}(\cS_{J_0\cap J_0',\flat}^{\bullet,\bullet})) \ar^{\wr}[d]\\
H^{k_0-r}(\mathrm{Tot}(\cT_{I_0}^{\bullet,\bullet})) & \otimes_E & H^{k_1-r}(\mathrm{Tot}(\cT_{I_0'}^{\bullet,\bullet})) \ar^{\cup}[r] & H^{k_0+k_1-r}(\mathrm{Tot}(\cT_{I_0\sqcup I_0'}^{\bullet,\bullet})) 
}
\end{equation}
for each $k_0,k_1\in\Z$, with all vertical maps being isomorphisms.

By taking $k_0=2\#I_0$, $k_1=2\#I_0'$ (with $k_0+k_1=2\#I_0\sqcup I_0'$) in (\ref{equ: x y cup w Tits}), we obtain the following map
\begin{equation}\label{equ: x y cup w bottom}
\mathbf{E}_{x,I_0}\otimes_E\mathbf{E}_{y,I_0'}\buildrel\cup\over\longrightarrow \mathbf{E}_{w,I_0\sqcup I_0'}
\end{equation}
which is functorial with respect to the choice of $x$, $y$ and $w$.
In particular, we obtain a map
\begin{equation}\label{equ: main bottom cup}
\kappa_{I_0,I_0'}: \mathbf{E}_{I_0}\otimes_E\mathbf{E}_{I_0'}\buildrel\cup\over\longrightarrow \mathbf{E}_{I_0\sqcup I_0'}
\end{equation}
which fits into the following commutative diagram (see also (\ref{equ: x y cup w Tits diagram}) and (\ref{equ: S T cup diagram}))
\begin{equation}\label{equ: x y cup w bottom diagram}
\xymatrix{
\mathbf{E}_{I_0} \ar@{->>}[d] & \otimes_E & \mathbf{E}_{I_0'} \ar^{\cup}[r] \ar@{->>}[d]& \mathbf{E}_{I_0\sqcup I_0'} \ar@{->>}[d]\\
\mathbf{E}_{x,I_0} & \otimes_E & \mathbf{E}_{y,I_0'} \ar^{\cup}[r] & \mathbf{E}_{w,I_0\sqcup I_0'}
}
\end{equation}
with all vertical maps being surjective by Proposition~\ref{prop: coxeter filtration x surjection}.
Taking $I=J_0=\Delta\setminus I_0$ and $I'=J_0'=\Delta\setminus I_0'$ with $I\cap I'=J_0\cap J_0'=\Delta\setminus(I_0\sqcup I_0')$ in (\ref{equ: x y pass to Lie diagram}) and assuming $\mathrm{Supp}(x)=I_0$ and $\mathrm{Supp}(y)=I_0'$ with $\mathrm{Supp}(w)=I_0\sqcup I_0'$, we have $\mathbf{E}_{x,I_0}\cong H^{\#I_0}(C_{J_0}^{\bullet}(M^{J_0}(x)))$, $\mathbf{E}_{y,I_0'}\cong H^{\#I_0'}(C_{J_0'}^{\bullet}(,M^{J_0'}(x)))$ and $H^{\#I_0\sqcup I_0'}(C_{J_0\cap J_0'}^{\bullet}(M^{J_0\cap J_0'}(w)))$, and obtain a commutative diagram of the form
\begin{equation}\label{equ: x y bottom cup to Lie}
\xymatrix{
\mathbf{E}_{x,I_0} \ar^{\wr}[d] & \otimes_E & \mathbf{E}_{y,I_0'} \ar^{\cup}[r] \ar^{\wr}[d] & \mathbf{E}_{w,I_0\sqcup I_0'} \ar^{\wr}[d]\\
\mathfrak{e}_{J_0,x} & \otimes_E & \mathfrak{e}_{J_0',y} \ar^{\cup}[r] & \mathfrak{e}_{J_0\cap J_0',w}
}
\end{equation}
with all vertical maps being isomorphisms between $1$-dimensional $E$-vector spaces by (\ref{equ: group to Lie parabolic restriction x}) (upon taking $\fh=0$ in \emph{loc.cit.} and taking $(x,I)$ in \emph{loc.cit.} to be $(x,J_0)$, $(y,J_0')$ or $(w,J_0\cap J_0')$).
Thanks to Lemma~\ref{lem: Lie cup diagram}, we know that the bottom horizontal map of (\ref{equ: x y bottom cup to Lie}) admits alternative definition through Chevalley-Eilenberg complex.

Let $\Delta_{-}\subseteq\Delta$ be an interval containing $I_0\sqcup I_0'$ with $m\defeq \#\Delta_{-}$.
Under the isomorphism $H_{\Delta_{-}}\cong G_{m+1}$ we have a bijection 
\begin{equation}\label{equ: Delta - bijection}
\Delta_{-}\buildrel\sim\over\longrightarrow \Delta_{m}=[1,m].
\end{equation}
For each $I\subseteq\Delta$, $I\cap\Delta_{-}$ corresponds to a subset of $\Delta_{m}$ under (\ref{equ: Delta - bijection}) and write $L_{m,I}\subseteq G_{m+1}$ for its associated Levi subgroup of $G_{m+1}$. In particular, the maps $G_{m+1}\cong H_{\Delta_{-}}\subseteq G$ induces the maps $L_{m,I}\cong L_{I}\cap H_{\Delta_{-}}\subseteq L_{I}$ for each $I\subseteq\Delta$.
Using the following commutative diagram
\[
\xymatrix{
\mathrm{Tot}(C^{\bullet}(L_{I'})\otimes_E C^{\bullet}(L_{I})) \ar@{-->}[r] \ar[d] & C^{\bullet}(L_{I\cap I'}) \ar[d]\\
\mathrm{Tot}(C^{\bullet}(L_{m,I'})\otimes_E C^{\bullet}(L_{m,I})) \ar@{-->}[r] & C^{\bullet}(L_{m,I\cap I'})
}
\]
which is functorial with respect to the choice of $I\in S_{J_0}$ and $I'\in S_{J_0'}$, we can construct a commutative diagram of the form  
\begin{equation}\label{equ: Tits interval reduction cup}
\xymatrix{
\mathbf{E}_{I_0} \ar[d] & \otimes_E & \mathbf{E}_{I_0'} \ar^{\cup}[r] \ar[d] & \mathbf{E}_{I_0\sqcup I_0'} \ar[d]\\
\mathbf{E}_{m,I_0} & \otimes_E & \mathbf{E}_{m,I_0'} \ar^{\cup}[r] & \mathbf{E}_{m,I_0\sqcup I_0'}
}
\end{equation}
with all vertical maps being isomorphisms from (\ref{equ: Tits interval reduction})

Given a tuple of disjoint subsets $J_1,\dots,J_m$ of $\Delta$, an element $x_k\in\Gamma^{J_{k}}$ for each $1\leq k\leq m$ and an element $w\in\Gamma_{x_1,\dots,x_m}$ (see (\ref{equ: tuple envelop})), we can similarly define a map
\begin{equation}\label{equ: tuple cup}
\mathbf{E}_{x_1,J_1}\otimes_E\cdots\otimes_E\mathbf{E}_{x_m,J_m}\buildrel\cup\over\longrightarrow \mathbf{E}_{w,\bigsqcup_{k=1}^{m}J_k}
\end{equation}
which is functorial with respect to the choice of $x_1,\dots,x_m$ and $w$.
In particular, when $x_1=\cdots=x_m=1$, we have a map
\begin{equation}\label{equ: main tuple cup}
\kappa_{J_1,\dots,J_m}: \mathbf{E}_{J_1}\otimes_E\cdots\otimes_E\mathbf{E}_{J_m}\buildrel\cup\over\longrightarrow \mathbf{E}_{\bigsqcup_{k=1}^{m}J_k}.
\end{equation}
In particular, the maps (\ref{equ: x y cup w bottom}) and (\ref{equ: main bottom cup}) satisfy natural associativity such as the following equality between maps from $\mathbf{E}_{I_0}\otimes_E\mathbf{E}_{I_0'}\otimes_E\mathbf{E}_{I_0''}$ to $\mathbf{E}_{I_0\sqcup I_0'\sqcup I_0''}$
\[\kappa_{I_0\sqcup I_0',I_0''}\circ(\kappa_{I_0,I_0'}\otimes_E\mathrm{Id}_{\mathbf{E}_{I_0''}})=\kappa_{I_0,I_0',I_0''}=\kappa_{I_0,I_0'\sqcup I_0''}\circ(\mathrm{Id}_{\mathbf{E}_{I_0}}\otimes_E\kappa_{I_0',I_0''})\]
for each triple $I_0$, $I_0'$ and $I_0''$ of pair wise disjoint subsets of $\Delta$.

For each $x,y\in\Gamma$ with $\mathrm{Supp}(x)\cap\mathrm{Supp}(y)=\emptyset$, we recall the subset $\Gamma_{x,y}\subseteq \Gamma$ from (\ref{equ: x y envelop}).
\begin{prop}\label{prop: coxeter cup comparison}
Let $I_0,I_0'\subseteq \Delta$ with $I_0\cap I_0'=\emptyset$. Let $x\in\Gamma^{I_0}$ and $y\in\Gamma^{I_0'}$.
Then the composition of the following map
\begin{equation}\label{equ: coxeter cup comparison}
\mathrm{Fil}_{x}(\mathbf{E}_{I_0})\otimes_E\mathrm{Fil}_{y}(\mathbf{E}_{I_0'})\hookrightarrow \mathbf{E}_{I_0}\otimes_E\mathbf{E}_{I_0'}\buildrel\kappa_{I_0,I_0'}\over\longrightarrow \mathbf{E}_{I_0\sqcup I_0'}
\end{equation}
is contained in
\begin{equation}\label{equ: coxeter cup comparison image}
\sum_{w\in\Gamma_{x,y}}\mathrm{Fil}_{w}(\mathbf{E}_{I_0\sqcup I_0'}).
\end{equation}
\end{prop}
\begin{proof}
We consider an arbitrary choice of $w'\in\Gamma^{I_0\sqcup I_0'}$, which uniquely determines $x'\in\Gamma^{I_0}$ and $y'\in\Gamma^{I_0'}$ such that $w'\in\Gamma_{x',y'}$ by \cite[Thm.~2.2.2]{BB05}. 
Upon replacing $x$, $y$ and $w$ in (\ref{equ: x y cup w bottom diagram}) with $x'$, $y'$ and $w'$, we obtain the following commutative diagram
\begin{equation}\label{equ: x y cup w diagram prime}
\xymatrix{
\mathbf{E}_{I_0} \ar@{->>}[d] & \otimes_E & \mathbf{E}_{I_0'} \ar[r] \ar@{->>}[d]& \mathbf{E}_{I_0\sqcup I_0'}\ar@{->>}[d]\\
\mathbf{E}_{x',I_0} & \otimes_E & \mathbf{E}_{y',I_0'} \ar[r] & \mathbf{E}_{w',I_0\sqcup I_0'}
}.
\end{equation}
If $w'\not\unlhd w$ for each $w\in\Gamma_{x,y}$, then by Lemma~\ref{lem: x y envelop partial order} we have either $x'\not\unlhd x$ or $y'\not\unlhd y$.
If $x'\not\unlhd x$, then the composition of $\mathrm{Fil}_{x}(\mathbf{E}_{I_0})\hookrightarrow \mathbf{E}_{I_0}\twoheadrightarrow \mathbf{E}_{x',I_0}$ is zero by Theorem~\ref{thm: coxeter filtration} and Proposition~\ref{prop: coxeter filtration x surjection}, and thus the composition of
\begin{equation}\label{equ: x y Fil w composition}
\mathrm{Fil}_{x}(\mathbf{E}_{I_0})\otimes_E\mathrm{Fil}_{y}(\mathbf{E}_{I_0'})\rightarrow \mathbf{E}_{I_0\sqcup I_0'}\twoheadrightarrow \mathbf{E}_{w',I_0\sqcup I_0'}
\end{equation}
is zero. Similarly, if $y'\not\unlhd y$, then the composition of $\mathrm{Fil}_{y}(\mathbf{E}_{I_0'})\hookrightarrow \mathbf{E}_{I_0'}\twoheadrightarrow \mathbf{E}_{y',I_0'}$ is zero again by Theorem~\ref{thm: coxeter filtration} and Proposition~\ref{prop: coxeter filtration x surjection}, and thus the composition of (\ref{equ: x y Fil w composition}) is again zero.
To summarize, we have shown that the composition of (\ref{equ: x y Fil w composition}) is zero for each $w'\in\Gamma^{I_0\sqcup I_0'}$ that satisfies $w'\not\unlhd w$ for each $w\in\Gamma_{x,y}$. 
This together with Theorem~\ref{thm: coxeter filtration} and Proposition~\ref{prop: coxeter filtration x surjection} implies that the image of the composition of (\ref{equ: coxeter cup comparison}) is contained in (\ref{equ: coxeter cup comparison image}).
\end{proof}
\subsection{Bottom degree cup product on graded pieces}\label{subsec: bottom cup E2}
We continue to write $r\defeq \#\Delta=n-1$ for short.
In this section, for each $I\subseteq \Delta$, we continue to use the shortened notation $\cT_{I}^{\bullet,\bullet}$ for $\cT_{\Delta\setminus I,\Delta}^{\bullet,\bullet}$ with
\[\mathbf{E}_{I}\defeq H^{2\#I-r}(\mathrm{Tot}(\cT_{I}^{\bullet,\bullet})),\]
and then define
\[\mathrm{gr}(\mathbf{E}_{I})\defeq \bigoplus_{\ell=r-\#I}^{r}\mathrm{gr}^{-\ell}(H^{2\#I-r}(\mathrm{Tot}(\cT_{I}^{\bullet,\bullet})))=\bigoplus_{\ell=r-\#I}^{r}E_{2,I}^{-\ell,\ell+2\#I-r}.\]
We fix a pair of subsets $I_0,I_0'\subseteq \Delta$ with $I_0\cap I_0'=\emptyset$, and write $J_0\defeq\Delta\setminus I_0$ as well as $J_0'\defeq\Delta\setminus I_0'$ for short.
We write $k_0\defeq 2\#I_0$ and $k_1=2\#I_0'$ for short.
For each pair $\#J_0\leq\ell_0\leq r$ and $\#J_0'\leq\ell_1\leq r$, we often write $\ell_2\defeq \ell_0+\ell_1-r$ below.
Recall from (\ref{equ: Levi cup Tits}) (and the discussion below it) that we have a map 
\begin{equation}\label{equ: bottom cup Tits}
\mathbf{E}_{I_0}\otimes_E \mathbf{E}_{I_0'}\rightarrow \mathbf{E}_{I_0\sqcup I_0'}.
\end{equation}
Note that (\ref{equ: bottom cup Tits}) restricts to a map
\begin{equation}\label{equ: bottom cup Tits filtration}
\mathrm{Fil}^{-\ell_0}(\mathbf{E}_{I_0})\otimes_E \mathrm{Fil}^{-\ell_1}(\mathbf{E}_{I_0'})\rightarrow \mathrm{Fil}^{-\ell_2}(\mathbf{E}_{I_0\sqcup I_0'})
\end{equation}
for each $\#J_0\leq\ell_0\leq r$ and $\#J_0'\leq\ell_1\leq r$, which induces the following map between graded pieces 
\begin{equation}\label{equ: bottom cup Tits graded total}
\mathrm{gr}(\mathbf{E}_{I_0})\otimes_E\mathrm{gr}(\mathbf{E}_{I_0'})\rightarrow \mathrm{gr}(\mathbf{E}_{I_0\sqcup I_0'})
\end{equation}
which is the direct sum over all $\#J_0\leq\ell_0\leq r$ and $\#J_0'\leq\ell_1\leq r$ of the following map (see the top row of (\ref{equ: Levi cup E2}))
\begin{equation}\label{equ: bottom cup Tits graded E2}
E_{2,I_0}^{-\ell_0,\ell_0+k_0-r}\otimes_E E_{2,I_0'}^{-\ell_1,\ell_1+k_1-r}\buildrel\cup\over\longrightarrow E_{2,I_0\sqcup I_0'}^{-\ell_2,\ell_2+k_0+k_1-r}.
\end{equation}
In this section, we study the map (\ref{equ: bottom cup Tits graded E2}) crucially using results from \S \ref{subsec: twisted tuples} (see Lemma~\ref{lem: separated cup injection} and Lemma~\ref{lem: separated cup disconnected isom}). As an application, we prove the injectivity of (\ref{equ: bottom cup Tits}) and give criterion for (\ref{equ: bottom cup Tits}) to be an isomorphism (see Theorem~\ref{thm: general cup}).

We start with the following simple result.
\begin{lem}\label{lem: filtered criterion}
Let
\[\varphi: (V,\mathrm{Fil}^{\bullet}(V))\rightarrow (W,\mathrm{Fil}^{\bullet}(W))\]
be a $E$-linear map between filtered $E$-vector spaces, which induces a map between the associated (total) graded pieces
\[\mathrm{gr}(\varphi): \mathrm{gr}(V)\rightarrow \mathrm{gr}(W).\]
If $\mathrm{gr}(\varphi)$ is an injection (resp.~an isomorphism), then so is $\varphi$.
\end{lem}
\begin{proof}
We choose an arbitrary basis $\Sigma$ of $V$ which is compatible with $\mathrm{Fil}^{\bullet}(V)$ and thus induces a basis of $\mathrm{gr}(V)$. If $\mathrm{gr}(\varphi)$ is an injection (resp.~an isomorphism), then $\{\varphi(v)\mid v\in\Sigma\}\subseteq W$ is linearly independent in $W$ (resp.~forms a basis of $W$), which forces $\varphi$ to be an injection (resp.~an isomorphism).
\end{proof}

Following Corollary~\ref{cor: bottom deg E2 basis}, we write $\Psi_{I_0,\ell_0}\defeq \Psi_{J_0,\Delta,\ell_0}$ for the set of $(J_0,\Delta)$-atomic equivalence classes with bidegree $(-\ell_0,\ell_0+k_0-r)$ and similarly for $\Psi_{I_0',\ell_1}\defeq \Psi_{J_0',\Delta,\ell_1}$ and $\Psi_{I_0\sqcup I_0',\ell_2}\defeq \Psi_{J_0\cap J_0',\Delta,\ell_2}$.
Given $J\subseteq \Delta$, we write $\cP(J)$ for the set of all subsets of $J$. We have a natural bijection
\begin{equation}\label{equ: power set product}
\cP(I_0)\times \cP(I_0')\buildrel\sim\over\longrightarrow \cP(I_0\sqcup I_0')
\end{equation}
that sends $(I,I')$ to $I\sqcup I'$ for each $I\in \cP(I_0)$ and $I'\in \cP(I_0')$.
Recall from (\ref{equ: partition product}) that we have an injective product map
\begin{equation}\label{equ: partition cup product}
\cS_{I_0}\times \cS_{I_0'}\hookrightarrow \cS_{I_0\sqcup I_0'}
\end{equation}
which together with the bijection (\ref{equ: power set product}) and the bijection (\ref{equ: atom to partition}) induces an injection
\begin{equation}\label{equ: atom cup product total}
\big(\bigsqcup_{\ell_0=\#J_0}^{r}\Psi_{I_0,\ell_0}\big)\times\big(\bigsqcup_{\ell_1=\#J_0'}^{r}\Psi_{I_0',\ell_1}\big)
\hookrightarrow \bigsqcup_{\ell_2=\#J_0\cap J_0'}^{r}\Psi_{I_0\sqcup I_0',\ell_2}
\end{equation}
which restricts to an injection
\begin{equation}\label{equ: atom cup product}
\Psi_{I_0,\ell_0}\times\Psi_{I_0',\ell_1}\hookrightarrow \Psi_{I_0\sqcup I_0',\ell_2}
\end{equation}
for each $\#J_0\leq\ell_0\leq r$ and $\#J_0'\leq\ell_1\leq r$ with $\ell_2=\ell_0+\ell_1-r$.
\begin{lem}\label{lem: separated cup injection}
Assume that $I_0<I_0'$ (see \ref{it: disconnected 2} of Definition~\ref{def: disconnected pair}). We have the following results.
\begin{enumerate}[label=(\roman*)]
\item \label{it: separated cup injection 1} Let $\#J_0\leq\ell_0\leq r$ and $\#J_0'\leq\ell_1\leq r$. Then for each $\Omega_0\in \Psi_{I_0,\ell_0}$ and $\Omega_1\in \Psi_{I_0',\ell_1}$, the image of $\overline{x}_{\Omega_0}\otimes_E\overline{x}_{\Omega_1}$ under (\ref{equ: bottom cup Tits graded E2}) equals $\pm\overline{x}_{\Omega_2}$, with $\Omega_2\in \Psi_{I_0\sqcup I_0',\ell_2}$ being the image of $(\Omega_0,\Omega_1)$ under (\ref{equ: atom cup product}).
\item \label{it: separated cup injection 2} The map (\ref{equ: bottom cup Tits graded E2}) is injective for each choice of $\#J_0\leq\ell_0\leq r$ and $\#J_0'\leq\ell_1\leq r$ with $\ell_2=\ell_0+\ell_1-r$. Moreover, the map (\ref{equ: bottom cup Tits}) is injective.
\end{enumerate}
\end{lem}
\begin{proof}
We first prove \ref{it: separated cup injection 1}.\\

The cases with either $I_0=\emptyset$ or $I_0'=\emptyset$ are easy. We assume in the rest of the proof that $J_0,J_0'\subsetneq \Delta$.
We fix now an arbitrary choice of $\Omega_0\in \Psi_{I_0,\ell_0}$ and $\Omega_1\in \Psi_{I_0',\ell_1}$, and write $\Theta_a^{\max}$ for the unique maximal element in $\Omega_a$ (see Lemma~\ref{lem: unique maximal tuple}) for $a=0,1$.
For $a=0,1$, our choice of $\Omega_a$ determines a $v_a\subseteq \mathbf{Log}_{\emptyset}$ which satisfies $\#v_a\cap\mathbf{Log}_{\Delta\setminus\{i\}}=1$ for each $i\in\Delta\setminus I_{v_a}$ (by applying \ref{it: bottom atom 1} of Lemma~\ref{lem: existence of atom} to the maximal element in $\Omega_a$). Since we have $J_0\subseteq I_{v_0}$ and $J_0'\subseteq I_{v_1}$, we deduce that
$\Delta\setminus I_{v_0}\subseteq I_0=\Delta\setminus J_0$, $\Delta\setminus I_{v_1}\subseteq I_0'=\Delta\setminus J_0'$ and thus $(\Delta\setminus I_{v_0})\cap (\Delta\setminus I_{v_1})\subseteq I_0\cap I_0'=\emptyset$ which forces $v_0\cap v_1=\emptyset$. We set $v_2\defeq v_0\sqcup v_1$ with $I_{v_2}=I_{v_0}\cap I_{v_1}$ and note that we have $\#v_2\cap\mathbf{Log}_{\Delta\setminus\{i\}}=1$ for each $i\in\Delta\setminus I_{v_2}=(\Delta\setminus I_{v_0})\sqcup(\Delta\setminus I_{v_1})$.
For $a=0,1$, our choice of $\Omega_a$ determines a choice of $\un{\Lambda}^a$ (and thus a choice of $\un{k}^a$) so that each $\Theta_a\in\Omega_a$ has the form $\Theta_a=(v_a, K_a, \un{k}^a,\un{\Lambda}^a)$, with $J_0\subseteq K_0$ if $a=0$, and $J_0'\subseteq K_1$ if $a=1$.

If $\Lambda^0_{d}=\emptyset$ for each $1\leq d\leq n-\ell_0$, then we set $s_0\defeq 1$ and $d_0\defeq 1$.
Otherwise, there exists a maximal possible $1\leq s_0\leq r_{I_{v_0}}$ with $d_0\defeq r_{\Omega_0}^{s_0}$ such that $\Lambda^0_{d_0}\neq \emptyset$.
Since $\Theta_0^{\max}$ is $(s,1)$-atomic for each $1\leq s\leq r_{I_{v_0}}$ by Lemma~\ref{lem: existence of atom} and Definition~\ref{def: atom twist}, we can define $(\Theta_0^{\max})^{s_0,d_0}$ as in (\ref{equ: maximal atom twist}) (with $\Theta_0^{\max}$ replacing $\Theta$ in \emph{loc.cit.}) and then define $\Omega_0'\defeq \Omega_0^{s_0,d_0}$ as the set of all tuples that are $(s_0,d_0)$-smaller than $(\Theta_0^{\max})^{s_0,d_0}$ (see Definition~\ref{def: twist partial order}). Recall from Lemma~\ref{lem: low deg} that
\[d_{1,I_0}^{-\ell_0,\ell_0+k_0-r}(x_{\Omega_0})=0,\]
which together with Proposition~\ref{prop: atom twist common class} (by replacing $\Omega$ and $\Omega^{s_0,d_0}$ in \emph{loc.cit.} with $\Omega_0$ and $\Omega_0'$) implies that
\begin{equation}\label{equ: atom differential zero 2}
d_{1,J_0,I_1}^{-\ell_0,\ell_0+k_0-r}(x_{\Omega_0'})=0
\end{equation}
and that the image of $x_{\Omega_0'}$ under
\begin{equation}\label{equ: separated cup to E2}
\mathrm{ker}(d_{1,I_0}^{-\ell_0,\ell_0+k_0-r})\twoheadrightarrow  E_{2,I_0}^{-\ell_0,\ell_0+k_0-r}
\end{equation}
equals $\overline{x}_{\Omega_0}$.
By definition of $\Omega_0'$, we know that $r_{\Omega_0'}^s=r_{\Omega_0}^s$ for each $1\leq s\leq r_{I_{v_0}}$, and each $\Theta_0'\in\Omega_0'$ has the form $(v_0,K_0',\un{k}^{0,\prime},\un{\Lambda}^{0,\prime})$ with $\Lambda^{0,\prime}_d=\Lambda^{0}_d$ for each $1\leq d\leq r_{\Omega_0}^{s_0-1}$, $\Lambda^{0,\prime}_d=\Lambda^{0}_{d+1}$ for each $r_{\Omega_0}^{s_0-1}+1\leq d\leq d_0-1$, and $\Lambda^{0,\prime}_{d}=\emptyset$ for each $d_0\leq d\leq \ell+1$.
Moreover, we notice that each tuple of the form $\Theta_0'=(v_0,K_0',\un{k}^{0,\prime},\un{\Lambda}^{0,\prime})$ (with $v_0$ and $\un{\Lambda}^{0,\prime}$ as given above, and varying $K_0'$) necessarily satisfies
\begin{equation}\label{equ: atom twist left bound}
i_{\Theta_0',r_{\Omega_0}^{s_0}-1}\leq i_{(\Theta_0^{\max})^{s_0,d_0},r_{\Omega_0}^{s_0}-1}
\end{equation}
and thus is an element of $\Omega_0'$ by Lemma~\ref{lem: twist atom criterion}.
For each choice of $\Theta_0'=(v_0,K_0',\un{k}^{0,\prime},\un{\Lambda}^{0,\prime})\in\Omega_0'$ and $\Theta_1=(v_1,K_1,\un{k}^{1},\un{\Lambda}^{1})\in\Omega_1$, our assumption $I_0<I_0'$ together with $\Delta\setminus I_{v_0}\subseteq I_0\subseteq J_0'\subseteq K_1$ and $\Delta\setminus I_{v_1}\subseteq I_0'\subseteq J_0\subseteq K_0'$ implies that there exists a unique tuple $\Theta_2'$ that satisfies
\begin{equation}\label{equ: cup of atom twist}
x_{\Theta_2'}=\mathrm{Res}_{K_0',K_0'\cap K_1}^{k_0-\ell}(x_{\Theta_0'})\cup \mathrm{Res}_{K_1,K_0'\cap K_1}^{k_1-\ell'}(x_{\Theta_1})\in E_{1,I_0\sqcup I_0'}^{-\ell_2,\ell_2+k_0+k_1-r},
\end{equation}
where we have $\Theta_2'=(v_2,K_2',\un{k}^{2,\prime},\un{\Lambda}^{2,\prime})$ with $K_2'\defeq K_0'\cap K_1$.
The data $\un{\Lambda}^{2,\prime}$ is uniquely characterized as follows.
\begin{itemize}
\item If $(K_2')^{d}\subseteq J_0\cap J_0'$, then we have $\Lambda^{2,\prime}_{d}=\emptyset$.
\item If $(K_2')^{d}\not\subseteq J_0\cap J_0'$, then there exists uniquely determined $1\leq d'\leq n-\ell_0$ and $1\leq d''\leq n-\ell_1$ such that $(K_2')^{d}\subseteq (K_0')^{d'}$ and $(K_2')^{d}\subseteq (K_1)^{d''}$, in which case we have $\Lambda^{2,\prime}_{d}\defeq \Lambda^{0,\prime}_{d'}\sqcup\Lambda^1_{d''}$. More precisely, we have either $(K_0')^{d'}\not\subseteq J_0$ and $(K_1)^{d''}\subseteq J_0'$ with $\Lambda^{0,\prime}_{d'}\neq \emptyset=\Lambda^1_{d''}$, or $(K_0')^{d'}\subseteq J_0$ and $(K_1)^{d''}\not\subseteq J_0'$ with $\Lambda^{0,\prime}_{d'}=\emptyset\neq \Lambda^1_{d''}$.
\end{itemize}
We use the shortened notation $\Theta_2'=\Theta_0'\cup \Theta_1$ to indicate that $\Theta_2'$ is the unique tuple that satisfies (\ref{equ: cup of atom twist}).

Recall that $(\Theta_0^{\max})^{s_0,d_0}$ is $(s_0,d_0)$-atomic and $(s,1)$-atomic for each $1\leq s\leq r_{I_{v_0}}$ satisfying $s\neq s_0$, and that $\Theta_1^{\max}$ is $(s,1)$-atomic for $1\leq s\leq r_{I_{v_1}}$, which together with the explicit characterization of $(\Theta_0^{\max})^{s_0,d_0}\cup \Theta_1$ above implies that $(\Theta_0^{\max})^{s_0,d_0}\cup \Theta_1$ is $(s_0,d_0)$-atomic and is $(s,1)$-atomic for each $1\leq s\leq r_{I_{v_2}}$ satisfying $s\neq s_0$.
Hence, by \ref{it: explicit twist atom 2} of Lemma~\ref{lem: explicit twist atom} we know that $(\Theta_0^{\max})^{s_0,d_0}\cup \Theta_1$ has the form $(\Theta_2^{\max})^{s_0,d_0}$ with $\Theta_2^{\max}$ being a tuple which is $(s,1)$-atomic for each $1\leq s\leq r_{I_{v_2}}$ (or equivalently maximally $(J_0\cap J_0',\Delta)$-atomic). We write $\Omega_2$ for the equivalence class containing $\Theta_2^{\max}$.
For each choice of $\Theta_0'=(v_0,K_0',\un{k}^{0,\prime},\un{\Lambda}^{0,\prime})\in\Omega_0'$ and $\Theta_1=(v_1,K_1,\un{k}^{1},\un{\Lambda}^{1})\in\Omega_1$ with $\Theta_2'=\Theta_0'\cup\Theta_1$, we have \[i_{\Theta_2',d_0-1}=i_{\Theta_0',d_0-1}\leq i_{(\Theta_0^{\max})^{s_0,d_0},d_0-1}=i_{(\Theta_2^{\max})^{s_0,d_0},d_0-1}\]
and
\[i_{\Theta_2',d_0}=i_{\Theta_1,1}\geq i_{\Theta_1^{\max},1}=i_{(\Theta_2^{\max})^{s_0,d_0},d_0},\]
which together with Lemma~\ref{lem: twist atom criterion} implies that $\Theta_2'\in \Omega_2'\defeq \Omega_2^{s_0,d_0}$.
Conversely, for each tuple $\Theta_2'=(v_2,K_2',\un{k}^2,\un{\Lambda}^{2,\prime})\in \Omega_2'$, we can recover $\Theta_0'$ and $\Theta_1$ from $\Theta_2'$ via
\begin{equation}\label{equ: separated cup recover}
K_0'\defeq K_2'\cup I_0' \ \ { and } \ \ K_1\defeq K_2'\cup I_0
\end{equation}
as well as the relation between $\un{\Lambda}^{2,\prime}$, $\un{\Lambda}^{0,\prime}$ and $\un{\Lambda}^1$ described below (\ref{equ: cup of atom twist}).
Note from (\ref{equ: separated cup recover}) that $\ell_0=\#K_0'$ and $\ell_1=\#K_1$ are uniquely determined by the choice of $\Omega_2'$.
To summarize, we have defined a bijection
\begin{equation}\label{equ: cup atom twist bijection}
\Omega_0'\times \Omega_1\rightarrow \Omega_2': (\Theta_0',\Theta_1)\mapsto \Theta_2'=\Theta_0'\cup\Theta_1.
\end{equation}
It follows from $J_0\cap J_0'\subseteq K_2'$ and (\ref{equ: separated cup recover}) that $\Delta\setminus K_2'=(\Delta\setminus K_0')\sqcup (\Delta\setminus K_1)$, which forces $\varepsilon(\Omega_0,\Omega_1)\defeq \varepsilon(\Theta_0')\varepsilon(\Theta_1)\varepsilon(\Theta_2')^{-1}$ to depend only on the choice of $\Omega_0$ and $\Omega_1$. This together with the bijection (\ref{equ: cup atom twist bijection}), the relation (\ref{equ: cup of atom twist}) ensures an equality of the form
\begin{equation}\label{equ: atom twist explicit cup}
x_{\Omega_2'}=\pm x_{\Omega_0'}\cup x_{\Omega_1}\in E_{1,I_0\sqcup I_0'}^{-\ell_2,\ell_2+k_0+k_1-r}
\end{equation}
Recall from Lemma~\ref{lem: low deg} (transferred to $\tld{E}_{\bullet,I_0\sqcup I_0'}^{\bullet,\bullet}$) that
\[d_{1,I_0\sqcup I_0'}^{-\ell_2,\ell_2+k_0+k_1-r}(x_{\Omega_2})=0,\]
which together with Proposition~\ref{prop: atom twist common class} (by replacing $\Omega$ and $\Omega^{s_0,d_0}$ in \emph{loc.cit.} with $\Omega_2$ and $\Omega_2'$) implies that
\begin{equation}\label{equ: atom differential zero 4}
d_{1,I_0\sqcup I_0'}^{-\ell_2,\ell_2+k_0+k_1-r}(x_{\Omega_2'})=0
\end{equation}
and that the image of $x_{\Omega_2'}$ under
\begin{equation}\label{equ: separated cup to E2 prime}
\mathrm{ker}(d_{1,I_0\sqcup I_0'}^{-\ell_2,\ell_2+k_0+k_1-r})\twoheadrightarrow E_{2,I_0\sqcup I_0'}^{-\ell_2,\ell_2+k_0+k_1-r}.
\end{equation}
equals $\overline{x}_{\Omega_2}$.
By Lemma~\ref{lem: low deg} we also have
\begin{equation}\label{equ: atom differential zero 5}
d_{1,I_0'}^{-\ell_1,\ell_1+k_1-r}(x_{\Omega_1})=0.
\end{equation}
and can write $\overline{x}_{\Omega_1}$ for the image of $x_{\Omega_1}$ under the surjection
\[\mathrm{ker}(d_{1,I_0'}^{-\ell_1,\ell_1+k_1-r})\twoheadrightarrow E_{1,I_0'}^{-\ell_1,\ell_1+k_1-r}.\]
Now that $\overline{x}_{\Omega_0}$ is the image of $x_{\Omega_0'}$ under (\ref{equ: separated cup to E2}), $\overline{x}_{\Omega_2}$ is the image of $x_{\Omega_2'}$ under (\ref{equ: separated cup to E2 prime}), and we have (\ref{equ: atom differential zero 2}), (\ref{equ: atom twist explicit cup}), (\ref{equ: atom differential zero 4}) as well as (\ref{equ: atom differential zero 5}), we deduce (from the definition of (\ref{equ: bottom cup Tits graded E2})) that
\[\overline{x}_{\Omega_0}\otimes_E\overline{x}_{\Omega_1} \in E_{2,I_0}^{-\ell_0,\ell_0+k_0-r}\otimes_E E_{2,I_0'}^{-\ell_1,\ell_1+k_1-r}\]
has image
\[\pm\overline{x}_{\Omega_2}\in E_{2,I_0\sqcup I_0'}^{-\ell_2,\ell_2+k_0+k_1-r}\]
under (\ref{equ: bottom cup Tits graded E2}). We write $\overline{x}_{\Omega_0}\cup\overline{x}_{\Omega_1}=\pm\overline{x}_{\Omega_2}$ in short.

Now we check that $\Omega_2$ is the image of $(\Omega_0,\Omega_1)$ under the injection (\ref{equ: atom cup product}).\\
We write $(S_{a},M_{a})$ (with $M_{a}\subseteq S_{a}\cap\Delta$) for the pair corresponding to $\Omega_{a}$ under the bijection (\ref{equ: atom to partition}) for each $a=0,1,2$. Now that $j\in M_{a}$ if and only if $\val_{j}\in v_{a}$ for each $j\in\Delta$ and $a=0,1,2$, we deduce from $v_{2}=v_{0}\sqcup v_{1}$ that $M_{2}=M_{0}\sqcup M_{1}$.
It is also clear from $v_{2}=v_{0}\sqcup v_{1}$ that 
\[S_{2}\cap\Delta=(S_{0}\cap\Delta)\sqcup (S_{1}\cap\Delta).\]
Hence, it remains to show that
\begin{equation}\label{equ: explicit partition pair}
S_{2}\setminus\Delta=(S_{0}\setminus\Delta)\sqcup (S_{1}\setminus\Delta).
\end{equation}
Following the definition of the bijection (\ref{equ: atom to partition}), we have the following characterization of elements in $S_{a}\setminus\Delta$ for each $a=0,1,2$.
\begin{itemize}
\item We have $\beta\in S_{0}\setminus\Delta$ if and only if there exists $1\leq s_1\leq r_{I_{v_0}}$ and $r_{\Omega_0}^{s_1-1}+2\leq d_1\leq r_{\Omega_0}^{s_1}$ (with $\Lambda^{0,\prime}_{d_1}\neq\emptyset$) such that $I_{\beta}=\{i_{\Theta_0^{\max},d_1-1}\}\sqcup (K_0^{d_1}\setminus J_0)$.
If $s_1<s_0$, we have $I_{\beta}=\{i_{(\Theta_0^{\max})^{s_0,d_0},d_1-1}\}\sqcup ((K_0')^{d_1}\setminus J_0)$. 
If $s_1=s_0$ (and $d_1\leq d_0=r_{\Omega_0}^{s_0}$), then we have $I_{\beta}=\{i_{(\Theta_0^{\max})^{s_0,d_0},d_1-1}\}\sqcup ((K_0')^{d_1-1}\setminus J_0)$.
\item We have $\beta\in S_{1}\setminus\Delta$ if and only if there exists $1\leq s_1\leq r_{I_{v_1}}$ and $r_{\Omega_1}^{s_1-1}+2\leq d_1\leq r_{\Omega_1}^{s_1}$ (with $\Lambda^{1}_{d_1}\neq\emptyset$) such that $I_{\beta}=\{i_{\Theta_1^{\max},d_1-1}\}\sqcup (K_1^{d_1}\setminus J_0')$.
\item We have $\beta\in S_{2}\setminus\Delta$ if and only if there exists $1\leq s_1\leq r_{I_{v_2}}$ and $r_{\Omega_2}^{s_1-1}+2\leq d_1\leq r_{\Omega_2}^{s_1}$ (with $\Lambda^{2,\prime}_{d_1}\neq\emptyset$) such that $I_{\beta}=\{i_{\Theta_2^{\max},d_1-1}\}\sqcup (K_2^{d_1}\setminus (J_0\cap J_0'))$.
If either $s_1\neq s_0$ or $d_1>d_0$, we have $I_{\beta}=\{i_{(\Theta_2^{\max})^{s_0,d_0},d_1-1}\}\sqcup ((K_2')^{d_1}\setminus (J_0\cap J_0'))$. 
If $s_1=s_0$ and $d_1\leq d_0$, then we have $I_{\beta}=\{i_{(\Theta_2^{\max})^{s_0,d_0},d_1-1}\}\sqcup ((K_2')^{d_1-1}\setminus (J_0\cap J_0'))$. 
\end{itemize}
Now we observe that, for each $1\leq d\leq n-\ell_2$ that satisfies $(K_2')^{d}\not\subseteq J_0\cap J_0'$ (or equivalently $\Lambda^{2,\prime}_d\neq\emptyset$), either we have $d<d_0$ with $(K_2')^{d}\subseteq (K_0')^{d}$ and $(K_2')^{d}\setminus(J_0\cap J_0')=(K_0')^{d}\setminus J_0$, or we have $d>d_0$ with $(K_2')^{d}\subseteq (K_1)^{d-d_0+1}$ and $(K_2')^{d}\setminus(J_0\cap J_0')=(K_0')^{d-d_0+1}\setminus J_0'$. Under our previous characterization of $S_{2}\setminus\Delta$ from $K_2'$ (and similarly characterization of $S_{0}\setminus\Delta$ and $S_{1}\setminus\Delta$ from $K_0'$ and $K_1$ respectively), we see that $\beta\in S_{2}\setminus\Delta$ if and only if either $\beta\in S_{0}\setminus\Delta$ or $\beta\in S_{1}\setminus\Delta$, which is nothing but (\ref{equ: explicit partition pair}).

We prove \ref{it: separated cup injection 2}.\\
Recall from Proposition~\ref{prop: atom basis} that $\{\overline{x}_{\Omega_0}\}_{\Omega_0\in \Psi_{I_0,\ell_0}}$ (resp.~$\{\overline{x}_{\Omega_1}\}_{\Omega_1\in \Psi_{I_0',\ell_1}}$, resp.~$\{\overline{x}_{\Omega_2}\}_{\Omega_2\in \Psi_{I_0\sqcup I_0',\ell_2}}$) forms a basis of $E_{2,I_0}^{-\ell_0,\ell_0+k_0-r}$ (resp.~ of $E_{2,I_0'}^{-\ell_1,\ell_1+k_1-r}$, resp.~ of $E_{2,I_0\sqcup I_0'}^{-\ell_2,\ell_2+k_0+k_1-r}$), which together with the injectivity of (\ref{equ: atom cup product}) gives the injectivity of (\ref{equ: bottom cup Tits graded E2}).
Now that (\ref{equ: bottom cup Tits graded total}) is the direct sum of (\ref{equ: bottom cup Tits graded E2}) over all $\#J_0\leq \ell_0\leq r$ and $\#J_0'\leq\ell_1\leq r$, we conclude that (\ref{equ: bottom cup Tits graded total}) is injective and thus (\ref{equ: bottom cup Tits}) is injective by Lemma~\ref{lem: filtered criterion}. The proof is thus finished.
\end{proof}

\begin{lem}\label{lem: separated cup disconnected isom}
Assume that $I_0<<I_0'$ (see \ref{it: disconnected 3} of Definition~\ref{def: disconnected pair}). We have the following results.
\begin{enumerate}[label=(\roman*)]
\item \label{it: disconnected isom 1} For each $\Omega_2\in \Psi_{I_0\sqcup I_0,\ell_2}$, there exists $\ell_0,\ell_1$ uniquely determined by $\Omega_2$, as well as $\Omega_0\in \Psi_{I_0,\ell_0}$ and $\Omega_1\in \Psi_{I_0',\ell_1}$ such that $\Omega_2$ is the image of $(\Omega_0,\Omega_1)$ under the map (\ref{equ: atom cup product}).
\item \label{it: disconnected isom 2} The map (\ref{equ: bottom cup Tits}) is an isomorphism.
\end{enumerate}
\end{lem}
\begin{proof}
We first prove \ref{it: disconnected isom 1}.\\

As $I_0<<I_0'$, we fix a choice of a $i\in\Delta$ that satisfies $I_0 < \{i\} < I_0'$, which necessarily satisfies $i\in J_0\cap J_0'$.
It suffices to show that the injection (\ref{equ: atom cup product}) from Lemma~\ref{lem: separated cup injection} is surjective.
We fix an arbitrary choice of $\Omega_2\in \Psi_{I_0\sqcup I_0,\ell_2}$ and give explicit construction of $\ell_0$, $\ell_1$, $\Omega_0\in \Psi_{I_0,\ell_0}$ and $\Omega_1\in \Psi_{I_0',\ell_1}$ such that $\Omega_2$ is the image of $(\Omega_0,\Omega_1)$ under (\ref{equ: atom cup product}).
Following the definition of the map (\ref{equ: atom cup product}) in the proof of Lemma~\ref{lem: separated cup injection}, we need to choose explicit $\Omega_0$, $\Omega_1$, $s_0$ and $d_0$ such that $\Omega_0'\defeq \Omega_0^{s_0,d_0}$ and $\Omega_2'\defeq \Omega_2^{s_0,d_0}$ are defined and satisfy $\Omega_2'=\Omega_0'\cup\Omega_1$ (see also (\ref{equ: atom twist explicit cup})).
Recall from \ref{it: unique maximal tuple 2} of Lemma~\ref{lem: unique maximal tuple} that the choice of $\Omega_a$ is characterized by $v_a\subseteq \mathbf{Log}_{\emptyset}$, $\un{\Lambda}^a$ and $r_{\Omega_a}^s$ for each $1\leq s\leq r_{I_{v_a}}-1$ (see the proof of Lemma~\ref{lem: separated cup injection} for notation).

Now that $i\in J_0\cap J_0'$, then there exists a unique $1\leq s_0\leq r_{I_{v_2}}$ and $r_{\Omega_2}^{s_0-1}+1\leq d_0\leq r_{\Omega_2}^{s_0}$ such that $i\in (K_2^{\max})^{d_0}$ for the unique maximal element $\Theta_2^{\max}=(v_2,K_2^{\max},\un{k}^2,\un{\Lambda}^2)\in\Omega_2$ (see Lemma~\ref{lem: unique maximal tuple}).
Then we choose
\begin{itemize}
\item $v_0\defeq v_2\cap \mathbf{Log}_{J_0}$, $v_1\defeq v_2\cap \mathbf{Log}_{J_0'}$;
\item $r_{\Omega_0}^s\defeq r_{\Omega_2}^s$ for each $1\leq s\leq r_{I_{v_0}}-1=s_0-1$, and $r_{\Omega_1}^s\defeq r_{\Omega_2}^{s+s_0-1}$ for each $1\leq s\leq r_{I_{v_1}}-1$;
\item $\Lambda^0_d\defeq \Lambda^2_d$ for each $1\leq d\leq r_{\Omega_0}^{s_0-1}$, $\Lambda^0_{r_{\Omega_0}^{s_0-1}+1}\defeq \emptyset$, and $\Lambda^0_d\defeq \Lambda^2_{d-1}$ for each $r_{\Omega_0}^{s_0-1}+2\leq d\leq r_{\Omega_0}^{s_0}=n-\ell_0$;
\item $\Lambda^1_1\defeq \emptyset$, and $\Lambda^1_d\defeq \Lambda^2_{d+d_0-1}$ for each $2\leq d\leq n-\ell_1$.
\end{itemize}
We leave it to the readers to check (following the proof of Lemma~\ref{lem: separated cup injection}) that these explicit choices work.

Now we prove \ref{it: disconnected isom 2}.\\
It follows from \ref{it: disconnected isom 1} that (\ref{equ: bottom cup Tits graded total}) is an isomorphism, which together with Lemma~\ref{lem: filtered criterion} forces (\ref{equ: bottom cup Tits}) to be an isomorphism.
\end{proof}

\begin{rem}\label{rem: cup graded not isom}
Given $\ell_2$, different choices of $\Omega_2\in \Psi_{I_0\sqcup I_0',\ell_2}$ often determine different choices of $\ell_0$ and $\ell_1$ as in \ref{it: disconnected isom 1} of Lemma~\ref{lem: separated cup disconnected isom}. Hence, even under the assumption $I_0 << I_0'$, (\ref{equ: bottom cup Tits graded E2}) is not an isomorphism in general.
\end{rem}

\begin{rem}\label{rem: opposite cup injection}
We continue to use notation from Lemma~\ref{lem: separated cup injection} and Lemma~\ref{lem: separated cup disconnected isom}.
Recall from Lemma~\ref{lem: grade cup commute} that we have
\[x\cup y=y\cup x\in E_{2,I_0\sqcup I_0'}^{-\ell_2,\ell_2+k_0+k_1-r}\]
for each $x\in E_{2,I_0}^{-\ell_0,\ell_0+k_0-r}$ and $y\in E_{2,I_0'}^{-\ell_1,\ell_1+k_1-r}$. This together with results from Lemma~\ref{lem: separated cup injection} and Lemma~\ref{lem: separated cup disconnected isom} also implies that the following cup product map
\[\mathbf{E}_{I_0'}\otimes_E\mathbf{E}_{I_0}\buildrel\cup\over\longrightarrow \mathbf{E}_{I_0\sqcup I_0'}\]
is an injection, and moreover is an isomorphism when $I_0 << I_0'$ (see \ref{it: disconnected 3} of Definition~\ref{def: disconnected pair}).
\end{rem}

\begin{thm}\label{thm: general cup}
Let $I_0,I_0'\subseteq \Delta$ with $I_0\cap I_0'=\emptyset$.
We have the following results.
\begin{enumerate}[label=(\roman*)]
\item \label{it: general cup 1} The map (\ref{equ: bottom cup Tits}) is injective.
\item \label{it: general cup 2} If $I_0$ and $I_0'$ do not connect (see \ref{it: disconnected 3} of Definition~\ref{def: disconnected pair}), then the map (\ref{equ: bottom cup Tits}) is an isomorphism.
\end{enumerate}
\end{thm}
\begin{proof}
To prove \ref{it: general cup 1} (resp.~\ref{it: general cup 2}), it suffices to show that (\ref{equ: bottom cup Tits graded total}) is injective (resp.~is an isomorphism when $I_0$ and $I_0'$ do not connect).

We prove that (\ref{equ: bottom cup Tits graded total}) is injective by an increasing induction on $\delta(I_0,I_0')$ (see \ref{it: disconnected 2} of Definition~\ref{def: disconnected pair}). 
If $\delta(I_0,I_0')=0$, then we have $I_0 < I_0'$ and thus (\ref{equ: bottom cup Tits graded total}) is injective by \ref{it: separated cup injection 2} of Lemma~\ref{lem: separated cup injection}.
Assume from now that $\delta(I_0,I_0')>0$.
Let $J\subseteq I_0$ (resp.~$J'\subseteq I_0'$) be the unique connected component such that $I\defeq I_0\setminus J << J$ and $J' << I'\defeq I_0'\setminus J'$. By \ref{it: disconnected isom 2} of Lemma~\ref{lem: separated cup disconnected isom} (upon replacing $I_0$ and $I_0'$ in \emph{loc.cit.} with $I$ and $J$) we obtain an isomorphism
\begin{equation}\label{equ: cup disconnected 1}
\mathrm{gr}(\mathbf{E}_{I})\otimes_E\mathrm{gr}(\mathbf{E}_{J})\buildrel\sim\over\longrightarrow\mathrm{gr}(\mathbf{E}_{I_0})
.
\end{equation}
Similarly, by Lemma~\ref{lem: separated cup disconnected isom} (upon replacing $I_0$ and $I_0'$ in \emph{loc.cit.} with $J'$ and $I'$) we obtain an isomorphism
\begin{equation}\label{equ: cup disconnected 2}
\mathrm{gr}(\mathbf{E}_{J'})\otimes_E\mathrm{gr}(\mathbf{E}_{I'})\buildrel\sim\over\longrightarrow \mathrm{gr}(\mathbf{E}_{I_0'}).
\end{equation}
Since $J$ and $J'$ are non-empty intervals that satisfy $J\cap J'\subseteq I_0\cap I_0'=\emptyset$, we have either $J < J'$ or $J' < J$. As $J < J'$ would force $I_0 < I_0'$ and $\delta(I_0,I_0')=0$ which is a contradiction, we must have $J' < J$.
In particular, we have
\begin{equation}\label{equ: cup induction}
\delta(I,J'), \delta(J,I'), \delta(I\sqcup J', J\sqcup I') < \delta(I\sqcup J, J'\sqcup I')=\delta(I_0,I_0'),
\end{equation}
which together with our inductive hypothesis give the following injections
\begin{equation}\label{equ: cup disconnected induction}
\left\{\begin{array}{cccccr}
\mathrm{gr}(\mathbf{E}_{I}) & \otimes_E & \mathrm{gr}(\mathbf{E}_{J'})
 & \hookrightarrow & \mathrm{gr}(\mathbf{E}_{I\sqcup J'}) &;\\
\mathrm{gr}(\mathbf{E}_{J}) & \otimes_E & \mathrm{gr}(\mathbf{E}_{I'}) & \hookrightarrow & \mathrm{gr}(\mathbf{E}_{J\sqcup I'})& ;\\
\mathrm{gr}(\mathbf{E}_{I\sqcup J'}) & \otimes_E & \mathrm{gr}(\mathbf{E}_{J\sqcup I'}) & \hookrightarrow & \mathrm{gr}(\mathbf{E}_{I_0\sqcup I_0'}) & .
\end{array}\right.
\end{equation}
By combining (\ref{equ: cup disconnected induction}) with (\ref{equ: cup disconnected 1}), (\ref{equ: cup disconnected 2}) and Lemma~\ref{lem: grade cup commute} we observe that (\ref{equ: Tits cup coh 1}) factors through (with $I_0=I\sqcup J$ and $I_0'=J'\sqcup I'$)
\begin{multline}\label{equ: cup seq injection}
\mathrm{gr}(\mathbf{E}_{I_0}) \otimes_E \mathrm{gr}(\mathbf{E}_{I_0'})
\cong \mathrm{gr}(\mathbf{E}_{I})\otimes_E \mathrm{gr}(\mathbf{E}_{J}) \otimes_E \mathrm{gr}(\mathbf{E}_{J'})\otimes_E \mathrm{gr}(\mathbf{E}_{I'})\\
\cong \mathrm{gr}(\mathbf{E}_{I})\otimes_E \mathrm{gr}(\mathbf{E}_{J'})\otimes_E \mathrm{gr}(\mathbf{E}_{J})\otimes_E \mathrm{gr}(\mathbf{E}_{I'})
\hookrightarrow \mathrm{gr}(\mathbf{E}_{I\sqcup J'})\otimes_E \mathrm{gr}(\mathbf{E}_{J\sqcup I'}) \hookrightarrow \mathrm{gr}(\mathbf{E}_{I_0\sqcup I_0'})
\end{multline}
and thus is an injection.

When $I_0$ and $I_0'$ do not connect, we prove that (\ref{equ: bottom cup Tits graded total}) is an isomorphism by an increasing induction on $\delta(I_0,I_0')$ for different choices of $I_0,I_0'$, following the argument of \ref{it: general cup 1}.\\
If $\delta(I_0,I_0')=0$, then since $I_0$ and $I_0'$ do not connect we must have $I_0 << I_0'$ and thus (\ref{equ: bottom cup Tits graded total}) is an isomorphism by \ref{it: disconnected isom 2} of Lemma~\ref{lem: separated cup disconnected isom}. If $\delta(I_0,I_0')>0$, then we define $J$, $I$, $J'$ and $I'$ as in the proof of \ref{it: general cup 1} which satisfy (\ref{equ: cup induction}). As we have $I << J$ and $J' << I'$, and $I_0=I\sqcup J$ and $I_0'=J'\sqcup I'$ do not connect, we deduce that $I$ and $J'$ (resp.~$J$ and $I'$, resp.~$I\sqcup J'$ and $J\sqcup I'$) do not connect, and thus all the injections in (\ref{equ: cup disconnected induction}) are isomorphisms by (\ref{equ: cup induction}) and our inductive hypothesis. Since (\ref{equ: bottom cup Tits graded total}) factors through (\ref{equ: cup seq injection}) which is now a sequence of isomorphisms, we conclude that (\ref{equ: bottom cup Tits graded total}) is an isomorphism.
\end{proof}

\subsection{Cup product basis}\label{subsec: cup relation}
In this section, we construct some convenient (but still non-canonical) basis of $\mathbf{E}_{I}$ for each $I\subseteq \Delta$ which is compatible with the cup product maps $\kappa_{I,I'}$ for each $I,I'\subseteq\Delta$ with $I\cap I'=\emptyset$ (see Proposition~\ref{prop: total cup basis}).

We write $r\defeq \#\Delta=n-1$ for short.
Let $I_0,I_0'\subseteq \Delta$ with $I_0\cap I_0'=\emptyset$, we recall the map $\kappa_{I_0,I_0'}$ from (\ref{equ: main bottom cup}).
For each $\#\Delta\setminus I_0\leq \ell_0\leq r$ and $\#\Delta\setminus I_0'\leq \ell_1\leq r$ with $\ell_2\defeq \ell_0+\ell_1-r$, the map $\kappa_{I_0,I_0'}$ restricts to a map
\begin{equation}\label{equ: main bottom cup Fil}
\mathrm{Fil}^{-\ell_0}(\mathbf{E}_{I_0})\otimes_E \mathrm{Fil}^{-\ell_1}(\mathbf{E}_{I_0'})\rightarrow \mathrm{Fil}^{-\ell_2}(\mathbf{E}_{I_0\sqcup I_0'}).
\end{equation}
In this section, we continue our study of the map $\kappa_{I_0,I_0'}$ (for each $I_0,I_0'\subseteq \Delta$ with $I_0\cap I_0'=\emptyset$) based on the results from Section~\ref{subsec: bottom cup E2}.

\begin{lem}\label{lem: E2 root generator}
Let $I\subseteq \Delta$. We have the following results.
\begin{enumerate}[label=(\roman*)]
\item \label{it: E2 cup generator 1} If $I=I_{\al}=\{\al\}$ for some $\al\in\Delta$, then we have canonical isomorphisms
\begin{equation}\label{equ: E2 cup generator 1}
\mathbf{E}_{I}\cong E_{2,\Delta\setminus I,\Delta}^{-(n-2),1}\cong \Hom(Z_{\Delta\setminus I},E).
\end{equation}
In particular, $\{\plog_{\al},\val_{\al}\}$ forms a basis of $\mathbf{E}_{I}$.
\item \label{it: E2 cup generator 2} If $I=I_{\al}$ for some $\al\in\Phi^+\setminus\Delta$, then we have $\Psi_{\Delta\setminus I,\Delta,n-1}=\emptyset$ and $\Psi_{\Delta\setminus I,\Delta,n-2}=\{\Omega\}$ with $\Omega$ being the unique equivalence class of $(\Delta\setminus I,\Delta)$-atomic tuples with bidegree $(-(n-2),2\#I-1)$ that satisfies $\Lambda_1=\emptyset$ and $\Lambda_2=\{2\#I-1\}$. In particular, $\{\overline{x}_{\Omega}\}$ forms a basis of $E_{2,\Delta\setminus I,\Delta}^{-(n-2),2\#I-1}$.
\end{enumerate}
\end{lem}
\begin{proof}
Both \ref{it: E2 cup generator 1} and \ref{it: E2 cup generator 2} easily follows from Corollary~\ref{cor: bottom deg E2 basis} and Proposition~\ref{prop: bottom deg degeneracy} (as well as the definition of $\mathbf{E}_{I}$).
\end{proof}

Let $J_{1},\dots,J_{s}\subseteq \Delta$ be a collection of disjoint subsets.
We consider one element $x_{i}\in\mathbf{E}_{J_{i}}$ for each $1\leq i\leq s$. There exists a unique permutation $J_{1}',\dots,J_{s}'$ of $J_{1},\dots,J_{s}$ such that $J_{1}'<\cdots<J_{s}'$. For each $1\leq i\leq s$, there exists a unique $1\leq j\leq s$ such that $J_{j}'=J_{i}$ and we set $y_{j}\defeq x_{i}$. We write
\begin{equation}\label{equ: normalized order cup}
\cup_{i=1}^{s}x_{i}\defeq y_{1}\cup\cdots\cup y_{s}\in\mathbf{E}_{\bigsqcup_{i=1}^{s}J_{i}}
\end{equation}
for short. In other words, we fix our convention to cup elements in a standard order (although we will see later in Theorem~\ref{thm: cup exchange} this choice of order is in fact unnecessary).

We recall the notation $I_{\al}\subseteq \Delta$ for $\al\in\Phi^+$ and $\al_I$ for $I\subseteq \Delta$ from the paragraph containing (\ref{equ: sum of roots}).
For each $\al\in\Delta$ (with $\#I_\al=1$), we write $x_\al^{\infty}$ (resp.~$x_{\al}$) for the elements in $\mathbf{E}_{I_\al}\cong \Hom(K^\times,E)$ corresponding to $\val_{\al}$ (resp.~$\plog_{\al}$) under \ref{it: E2 cup generator 1} of Lemma~\ref{lem: E2 root generator}.
For each $\al\in\Phi^+\setminus\Delta$ (with $\#I_\al\geq 2$), we choose an arbitrary
\begin{equation}\label{equ: cup generator choice}
x_{\al}\in\mathbf{E}_{I_\al}
\end{equation}
whose image under
\[\mathbf{E}_{I_{\al}}\twoheadrightarrow E_{2,\Delta\setminus I_{\al},\Delta}^{-\ell,\ell+\#\Delta-2\#\Delta\setminus I_{\al}}=E_{2,\Delta\setminus I_{\al},\Delta}^{-(n-2),2\#I_\al-1}\]
is given by $\overline{x}_{\Omega_{\al}}$ with $\Omega_{\al}$ being the unique element of $\Psi_{\Delta\setminus I_{\al},\Delta,n-2}$ (see \ref{it: E2 cup generator 2} of Lemma~\ref{lem: E2 root generator}).

Let $I_0\subseteq \Delta$.
For each pair $(S,I)$ with $S\in\cS_{I_0}$ and $I\subseteq S\cap\Delta$, we define
\begin{equation}\label{equ: explicit cup generator}
x_{S,I}\defeq \cup_{\al\in S}x_{\al}^{\ast}\in\mathbf{E}_{I_0}
\end{equation}
following the convention introduced in (\ref{equ: normalized order cup}), where $x_{\al}^{\ast}=x_{\al}^{\infty}$ if $\al\in I$ and $x_{\al}^{\ast}=x_{\al}$ if $\al\in S\setminus I$.
Now that we have $x_{\al}^{\ast}\in\mathrm{Fil}^{1-r}(\mathbf{E}_{I_{\al}})$ for each $\al\in\Phi^+$ (with $\ast\in\{,\infty\}$ and $\ast=\infty$ only if $\al\in\Delta$) and $x_{\al}^{\infty}\in\mathrm{\tau}^{I_{\al}}(\mathbf{E}_{I_{\al}})$ for each $\al\in\Delta$, we have (see (\ref{equ: Levi cup term natural J}) and (\ref{equ: main bottom cup Fil}))
\begin{equation}\label{equ: generator filtration}
x_{S,I}\in\mathrm{Fil}^{\#S-r}(\mathrm{\tau}^{I}(\mathbf{E}_{I_0})).
\end{equation}
Recall from (\ref{equ: atom to partition}) that we have a bijection
\begin{equation}\label{equ: cup basis atom to partition}
\bigsqcup_{\ell=\#\Delta\setminus I_0}^{\#\Delta}\Psi_{\Delta\setminus I_0,\Delta,\ell}\buildrel\sim\over\longrightarrow \{(S,I)\mid S\in\cS_{I_0}, I\subseteq S\cap\Delta\}
\end{equation}
\begin{lem}\label{lem: E2 cup of generator}
Let $I_0\subseteq \Delta$. Let $S\in\cS_{I_0}$, $I\subseteq S\cap\Delta$ and $\Omega\in\Psi_{\Delta\setminus I_0,\Delta,\ell}$ being the equivalence class associated with $(S,I)$ under (\ref{equ: cup basis atom to partition}) with $\ell=r-\#S$.
Then the image of $x_{S,I}$ under
\begin{equation}\label{equ: E2 cup of generator graded ell}
\mathrm{Fil}^{\#S-r}(\mathbf{E}_{I_0})\twoheadrightarrow \mathrm{gr}^{\#S-r}(\mathbf{E}_{I_0})\cong E_{2,\Delta\setminus I_0,\Delta}^{\#S-r,2\#I_0-\#S}
\end{equation}
equals $\pm\overline{x}_{\Omega}$.
\end{lem}
\begin{proof}
We prove this by an increasing induction on $\#S$.
If $S=\{\al\}$ with $I_0=I_{\al}$, then we have either $S=I$ with $x_{S,I}=x_{\al}^{\infty}$, or $I=\emptyset$ with $x_{S,I}=x_{\al}$ whose image under 
\begin{equation}\label{equ: E2 cup of generator quotient}
\mathrm{Fil}^{1-r}(\mathbf{E}_{I_{\al}})\twoheadrightarrow \mathrm{gr}^{1-r}(\mathbf{E}_{I_{\al}})\cong E_{2,\Delta\setminus I_{\al},\Delta}^{1-r,2\#I_{\al}-1}
\end{equation}
is $\overline{x}_{\Omega_{\al}}$ by our very choice of $x_{\al}$.
Assume from now that $\#S\geq 2$ and let $\al\in S$ be the unique element that satisfies $I_{\beta}<I_{\al}$ for each $\beta\in S'\defeq S\setminus\{\al\}$ with $I'\defeq I\setminus\{\al\}$, $I_0'\defeq I_0\setminus I_{\al}$ and $\ell'\defeq r-\#S'=\ell+1$. 
As usual, we write $x_{\al}^{\ast}\defeq x_{\al}^{\infty}$ if $\al\in I$, and $x_{\al}^{\ast}\defeq x_{\al}$ otherwise.
We write $\Omega_{\al}\in \Psi_{\Delta\setminus I_{\al},\Delta,r-1}$ for the unique equivalence class such that the image of $x_{\al}^{\ast}$ under (\ref{equ: E2 cup of generator quotient}) equals $\overline{x}_{\Omega_{\al}}$.
Now that $\#S'=\#S-1$, our induction hypothesis implies that the image of $x_{S',I'}$ under
\[
\mathrm{Fil}^{\#S'-r}(\mathbf{E}_{I_0'})\twoheadrightarrow \mathrm{gr}^{\#S'-r}(\mathbf{E}_{I_0'})\cong E_{2,\Delta\setminus I_0',\Delta}^{\#S'-r,2\#I_0'-\#S'}
\]
equals $\pm\overline{x}_{\Omega'}$ with $\Omega'\in\Psi_{\Delta\setminus I_0',\Delta,\ell'}$ being the equivalence class associated with $(S',I')$. It follows from \ref{it: separated cup injection 1} of Lemma~\ref{lem: separated cup injection} that the image of $\overline{x}_{\Omega'}\otimes_E\overline{x}_{\Omega_{\al}}$ under the map
\[E_{2,\Delta\setminus I_0',\Delta}^{\#S'-r,2\#I_0'-\#S'}\otimes_E E_{2,\Delta\setminus I_{\al},\Delta}^{1-r,2\#I_{\al}-1}\buildrel\cup\over\longrightarrow E_{2,\Delta\setminus I_0,\Delta}^{\#S-r,2\#I_0-\#S}\]
equals $\pm\overline{x}_{\Omega}$, which together with $x_{S,I}=x_{S',I'}\cup x_{\al}^{\ast}$ (by definition of $x_{S,I}$ and $x_{S',I'}$) and our previous discussion on $x_{S',I'}$ (from induction hypothesis) implies that the image of $x_{S,I}$ under (\ref{equ: E2 cup of generator graded ell}) equals $\pm\overline{x}_{\Omega}$.
\end{proof}

\begin{prop}\label{prop: total cup basis}
Let $I_0\subseteq \Delta$. Then for each $\#\Delta\setminus I_0\leq \ell\leq \#\Delta$, the $E$-vector space $\mathrm{Fil}^{-\ell}(\mathbf{E}_{I_0})$ admits a basis of the form
\begin{equation}\label{equ: cup basis Fil}
\{x_{S,I}\}_{S\in \cS_{I_0}, I\subseteq S\cap\Delta, \#S\geq r-\ell}.
\end{equation}
\end{prop}
\begin{proof}
For each $\#\Delta\setminus I_0\leq \ell\leq \#\Delta$, it follows from Lemma~\ref{lem: E2 cup of generator} that the image of the following set
\[
\{x_{S,I}\}_{S\in \cS_{I_0}, I\subseteq S\cap\Delta, \#S=r-\ell}
\]
under $\mathrm{Fil}^{-\ell}(\mathbf{E}_{I_0})\twoheadrightarrow \mathrm{gr}^{-\ell}(\mathbf{E}_{I_0})$ forms a basis of $\mathrm{gr}^{-\ell}(\mathbf{E}_{I_0})$. This finishes the proof by an increasing induction on $\ell$.
\end{proof}

For each $I_0\subseteq \Delta$, we set
\begin{equation}\label{equ: decomposable subspace}
\mathbf{E}_{I_0}^{<}\defeq \sum_{\emptyset\neq I\subsetneq I_0}\mathrm{im}(\kappa_{I_0\setminus J,J})\subseteq \mathbf{E}_{I_0}.
\end{equation}
\begin{lem}\label{lem: primitive generator}
Let $I_0\subseteq \Delta$ be subsets with $\#I_0\geq 2$. We have $\Dim_E\mathbf{E}_{I_0}/\mathbf{E}_{I_0}^{<}\leq 1$, and it is non-zero if and only if $I_0$ is an interval.
\end{lem}
\begin{proof}
Let $S\in \cS_{I_0}$ and $I\subseteq S\cap\Delta$.
If $\#S\geq 2$, there exists a unique $\al\in S$ such that $I_{\beta}<I_{\al}$ for each $\beta\in S'\defeq S\setminus\{\al\}$ with $I'\defeq I\setminus\{\al\}$, in which case we have
\[x_{S,I}=x_{S',I'}\cup x_{\{\al\},I_{\al}\cap I}\in \mathbf{E}_{I_0}^{<}\]
with $\emptyset\neq I_{\al}\subsetneq I_0$ (as $S'=S\setminus\{\al\}\neq \emptyset$).
By applying Proposition~\ref{prop: total cup basis} to $I_0$ as well as to $J$ and $I_0\setminus J$ for each $\emptyset\neq J\subsetneq I_0$, we see that $\mathbf{E}_{I_0}^{<}$ in fact admits a basis of the form
\begin{equation}\label{equ: cup basis decomposable}
\{x_{S,I}\}_{S\in \cS_{I_0}, I\subseteq S\cap\Delta, \#S\geq 2}.
\end{equation}
Consequently, $\mathbf{E}_{I_0}/\mathbf{E}_{I_0}^{<}\neq 0$ if and only if $\{\al\}\in \cS_{I_0}$ for some $\al\in\Delta$, if and only if $I_0$ is an interval, in which case $\{x_{\al}\}\subseteq \mathbf{E}_{I_0}$ induces a basis of $\mathbf{E}_{I_0}/\mathbf{E}_{I_0}^{<}$ and thus $\mathbf{E}_{I_0}/\mathbf{E}_{I_0}^{<}$ is $1$ dimensional.
\end{proof}

Let $I=[i,j]\subseteq \Delta$ be a non-empty interval.
We consider the decreasing filtration
\begin{equation}\label{equ: auto total filtration}
\mathrm{Fil}_{W}^{k}(\mathbf{E}_{I})\defeq \sum_{k+1\leq \ell\leq j}\mathrm{im}(\kappa_{[i,\ell-1],[\ell,j]})
\end{equation}
with $\mathrm{gr}_{W}^{k}(\mathbf{E}_{I})=0$ for either $k<i-1$ or $k\geq j$,
\[\mathrm{gr}_{W}^{j-1}(\mathbf{E}_{I})=\mathbf{E}_{[i,j-1]}\otimes_E\mathbf{E}_{\{j\}},\]
and being
\[\mathrm{gr}_{W}^{k}(\mathbf{E}_{I})=\mathbf{E}_{[i,k]}\otimes_E(\mathbf{E}_{[k+1,j]}/\mathbf{E}_{[k+1,j]}^{<})\]
for each $i-1\leq k<j-1$.
In particular, we have
\begin{multline}\label{equ: total dim induction}
\Dim_E \mathbf{E}_{I}=\sum_{k=i-1}^{j-1}\Dim_E \mathrm{gr}_{W}^{k}(\mathbf{E}_{I})\\
=\Dim_E \mathbf{E}_{[i,j-1]}\otimes_E\mathbf{E}_{\{j\}}+\sum_{k=i-1}^{j-2}\Dim_E \mathbf{E}_{[i,k]}\otimes_E(\mathbf{E}_{[k+1,j]}/\mathbf{E}_{[k+1,j]}^{<})\\
=2\Dim_E \mathbf{E}_{[i,j-1]}+\sum_{k=i-1}^{j-2}\Dim_E \mathbf{E}_{[i,k]}.
\end{multline}
\subsection{Cup product of Coxeter filtration}\label{subsec: cup coxeter filtration}
In this section, we refine the comparison between cup product maps and Coxeter filtration in Proposition~\ref{prop: coxeter cup comparison} by proving Proposition~\ref{prop: cup x y nonvanishing}. The main application of such refinements in the transversality result in Theorem~\ref{thm: cup top grade}.

\begin{lem}\label{lem: cup with sm}
Let $I_0,I_0'\subseteq \Delta$ with $I_0\cap I_0'=\emptyset$ and $x\in\Gamma^{I_0}$.
Then we have
\begin{equation}\label{equ: cup with sm}
\kappa_{I_0,I_0'}(\mathrm{Fil}_{x}(\mathbf{E}_{I_0})\otimes_E\mathrm{Fil}_{1}(\mathbf{E}_{I_0'}))=\mathrm{Fil}_{x}(\mathbf{E}_{I_0\sqcup I_0'})=\kappa_{I_0',I_0}(\mathrm{Fil}_{1}(\mathbf{E}_{I_0'})\otimes_E\mathrm{Fil}_{x}(\mathbf{E}_{I_0})).
\end{equation}
\end{lem}
\begin{proof}
We only prove the LHS equality of (\ref{equ: cup with sm}) and the argument for the RHS equality is symmetric.
Recall from Proposition~\ref{prop: coxeter cup comparison} that we have
\begin{equation}\label{equ: cup with sm inclusion}
\kappa_{I_0,I_0'}(\mathrm{Fil}_{x}(\mathbf{E}_{I_0})\otimes_E\mathrm{Fil}_{1}(\mathbf{E}_{I_0'}))\subseteq \mathrm{Fil}_{x}(\mathbf{E}_{I_0\sqcup I_0'}).
\end{equation}
Now that $\kappa_{I_0,I_0'}$ is injective by \ref{it: general cup 1} of Theorem~\ref{thm: general cup}, we deduce from Theorem~\ref{thm: coxeter filtration} that
\begin{multline*}
\Dim_E\mathrm{Fil}_{x}(\mathbf{E}_{I_0\sqcup I_0'})=\#\{u\in\Gamma^{I_0}\mid u\unlhd x\}=\Dim_E\mathrm{Fil}_{x}(\mathbf{E}_{I_0})\\
=\Dim_E\mathrm{Fil}_{x}(\mathbf{E}_{I_0})\Dim_E\mathrm{Fil}_{1}(\mathbf{E}_{I_0'})=\Dim_E\kappa_{I_0,I_0'}(\mathrm{Fil}_{x}(\mathbf{E}_{I_0})\otimes_E\mathrm{Fil}_{1}(\mathbf{E}_{I_0'})),
\end{multline*}
which forces the inclusion (\ref{equ: cup with sm inclusion}) to be an equality.
\end{proof}

\begin{lem}\label{lem: tau and coxeter}
Let $I_0\subseteq\Delta$. We have the following results.
\begin{enumerate}[label=(\roman*)]
\item \label{it: tau and coxeter 1} For each $J\subseteq I_0$, we have
\begin{multline}\label{equ: tau and coxeter 1}
\tau^{J}(\mathbf{E}_{I_0})=\kappa_{I_0\setminus J,J}(\mathbf{E}_{I_0\setminus J}\otimes_E\tau^{J}(\mathbf{E}_{J}))=\kappa_{I_0\setminus J,J}(\mathbf{E}_{I_0\setminus J}\otimes_E\mathrm{Fil}_{1}(\mathbf{E}_{J}))\\
=\sum_{x\in\Gamma^{I_0\setminus J}}\mathrm{Fil}_{x}(\mathbf{E}_{I_0})=\sum_{x\in\Gamma_{I_0\setminus J}}\mathrm{Fil}_{x}(\mathbf{E}_{I_0})
\end{multline}
\item \label{it: tau and coxeter 2} For each $i\geq 0$, we have
\[
\tau^{i}(\mathbf{E}_{I_0})=\sum_{x\in\Gamma^{I_0},\ell(x)\leq \#I_0-i}\mathrm{Fil}_{x}(\mathbf{E}_{I_0})=\sum_{x\in\Gamma_{I_0},\ell(x)=\#I_0-i}\mathrm{Fil}_{x}(\mathbf{E}_{I_0})
\]
\end{enumerate}
\end{lem}
\begin{proof}
Now that we have
\[\tau^{i}(\mathbf{E}_{I_0})=\sum_{J\subseteq I_0, \#J=i}\tau^{J}(\mathbf{E}_{I_0}),\]
it suffices to prove (\ref{equ: tau and coxeter 1}). Recall from (\ref{equ: Levi cup Tits natural J}) that $\kappa_{I_0\setminus J,J}$ restricts to the following cup product map
\begin{equation}\label{equ: tau and coxeter cup}
\kappa_{I_0\setminus J,J}: \mathbf{E}_{I_0\setminus J}\otimes_E\tau^{J}(\mathbf{E}_{J})\buildrel\cup\over\longrightarrow \tau^{J}(\mathbf{E}_{I_0}).
\end{equation}
On one hand, we know that (\ref{equ: tau and coxeter cup}) is injective by \ref{it: general cup 1} of Theorem~\ref{thm: general cup}. On the other hand, we know that
\[\Dim_E\tau^{J}(\mathbf{E}_{J})=\sum_{\ell=\#\Delta\setminus J}^{\Delta}\#\mathrm{gr}^{J}(\Psi_{\Delta\setminus J,\Delta,\ell})=\#\mathrm{gr}^{J}(\Psi_{\Delta\setminus J,\Delta,\#\Delta\setminus J})=1\]
and
\begin{multline}\label{equ: tau and coxeter dim}
\Dim_E\tau^{J}(\mathbf{E}_{I_0})=\sum_{\ell=\#\Delta\setminus I_0}^{\Delta}\#\mathrm{gr}^{J}(\Psi_{\Delta\setminus I_0,\Delta,\ell})=\#\{(S,I)\mid S\in\cS_{I_0}, J\subseteq I\subseteq S\cap\Delta\}\\
=\#\{(S',I')\mid S'\in\cS_{I_0\setminus J}, I'\subseteq S'\cap\Delta\}=\sum_{\ell=\#\Delta\setminus(I_0\setminus J)}^{\Delta}\#\Psi_{\Delta\setminus (I_0\setminus J),\Delta,\ell}=\Dim_E\mathbf{E}_{I_0\setminus J}
\end{multline}
from Proposition~\ref{prop: bottom deg graded} and (\ref{equ: atom to partition}) (as well as the discussion below it), with the middle equality of (\ref{equ: tau and coxeter dim}) induced from the bijection that sends $(S,I)$ to $(S',I')\defeq (S\setminus J, I\setminus J)$. 
Now that the cup product map (\ref{equ: tau and coxeter cup}) is an injection between two $E$-vector spaces of the same dimension, it must be an isomorphism. For similar reason, we see that $\tau^{J}(\mathbf{E}_{J})$ is the cup product of $\tau^{\{j\}}(\mathbf{E}_{\{j\}})=\mathrm{Fil}_{1}(\mathbf{E}_{\{j\}})$ for all $j\in J$, which together with Lemma~\ref{lem: cup with sm} gives
\[\tau^{J}(\mathbf{E}_{J})=\mathrm{Fil}_{1}(\mathbf{E}_{J}).\]
Further combining this fact with the isomorphism (\ref{equ: tau and coxeter cup}) and Lemma~\ref{lem: cup with sm}, we conclude that
\begin{multline*}
\tau^{J}(\mathbf{E}_{I_0})=\kappa_{I_0\setminus J,J}(\mathbf{E}_{I_0\setminus J}\otimes_E\tau^{J}(\mathbf{E}_{J}))
=\kappa_{I_0\setminus J,J}((\sum_{x\in\Gamma^{I_0\setminus J}}\mathrm{Fil}_{x}(\mathbf{E}_{I_0\setminus J}))\otimes_E\mathrm{Fil}_{1}(\mathbf{E}_{J}))\\
=\sum_{x\in\Gamma^{I_0\setminus J}}\kappa_{I_0\setminus J,J}(\mathrm{Fil}_{x}(\mathbf{E}_{I_0\setminus J})\otimes_E\mathrm{Fil}_{1}(\mathbf{E}_{J}))=\sum_{x\in\Gamma^{I_0\setminus J}}\mathrm{Fil}_{x}(\mathbf{E}_{I_0})
\end{multline*}
which finishes the proof of (\ref{equ: tau and coxeter 1}).
\end{proof}

\begin{rem}\label{rem: tau cup basis}
Let $J\subseteq I_0\subseteq\Delta$. It follows from (\ref{equ: tau and coxeter 1}) that the $E$-vector space $\tau^{J}(\mathbf{E}_{I_0})$ admits a basis of the form
\[
\{x_{S,I}\}_{S\in \cS_{I_0}, J\subseteq I\subseteq S\cap\Delta}.
\]
\end{rem}

\begin{prop}\label{prop: graded coxeter tau comparison}
Let $J\subseteq I_0\subseteq\Delta$. Then the following isomorphism (see (\ref{equ: grade J Z to CE}))
\begin{equation}\label{equ: graded coxeter tau comparison 1}
\mathrm{gr}^{J}(\mathbf{E}_{I_0})\buildrel\sim\over\longrightarrow \bigoplus_{x\in\Gamma_{I_0\setminus J}}\mathfrak{e}_{\Delta\setminus I_0,x}\otimes_E\wedge^{\#J}\Hom(Z_{\Delta\setminus J}^{\dagger},E)
\end{equation}
is the direct sum over $x\in\Gamma_{I_0\setminus J}$ of the following isomorphism between $1$-dimensional $E$-vector spaces
\begin{equation}\label{equ: graded coxeter tau comparison 2}
\mathrm{gr}_{x}(\mathbf{E}_{I_0})\buildrel\sim\over\longrightarrow \mathfrak{e}_{\Delta\setminus I_0,x}\otimes_E\wedge^{\#J}\Hom(Z_{\Delta\setminus J}^{\dagger},E),
\end{equation}
under the identification
\[\mathrm{gr}^{J}(\mathbf{E}_{I_0})=\bigoplus_{x\in\Gamma_{I_0\setminus J}}\mathrm{gr}_{x}(\mathbf{E}_{I_0})\]
from \ref{it: tau and coxeter 1} of Lemma~\ref{lem: tau and coxeter}.
\end{prop}
\begin{proof}
Recall from Lemma~\ref{lem: grade 0 Lie transfer} (upon taking $I_1$ in \emph{loc.cit.} to be $\Delta\setminus(I_0\setminus J)$) that we have a commutative diagram of the form
\begin{equation}\label{equ: graded coxeter tau diagram 1}
\xymatrix{
\mathrm{gr}^{\emptyset}(\mathbf{E}_{I_0\setminus J}) \ar@{=}[rrr] \ar^{\wr}[d] & & & \bigoplus_{x\in\Gamma_{I_0\setminus J}}\mathrm{gr}_{x}(\mathbf{E}_{I_0\setminus J}) \ar^{\wr}[d]\\
\mathrm{Ext}_{U(\fg)}^{\#I_0\setminus J}(\mathfrak{v}_{\Delta\setminus(I_0\setminus J)},1_{\fg}) \ar^{\sim}[rrr] & & & \bigoplus_{x\in\Gamma_{I_0\setminus J}}\mathfrak{e}_{\Delta\setminus(I_0\setminus J),x}
}
\end{equation}
We write $h\defeq 2\#(I_0\setminus J)-\#\Delta$ and $h'\defeq 2\#I_0-\#\Delta$ for short.
It follows from (\ref{equ: Tits Verma to Levi}), (\ref{equ: Lie general Tits resolution}) and (\ref{equ: Lie Tits J shift u diagram Ext}) (upon replacing $I_0$ and $J$ in \emph{loc.cit.} with $\Delta\setminus(I_0\setminus J)$ and $\Delta\setminus J$ here) that we have the following commutative diagram
\begin{equation}\label{equ: graded coxeter tau diagram 2}
\xymatrix{
H^{h}(\mathrm{Tot}(\mathrm{CE}_{\Delta\setminus(I_0\setminus J),\Delta}^{\bullet,\bullet})) \ar^{\sim}[rr] \ar^{\wr}[d] & & H^{h'-\#J}(\mathrm{Tot}(\mathrm{CE}_{\Delta\setminus I_0,\Delta\setminus J}^{\bullet,\bullet})) \ar^{\wr}[d]\\
H^{h}(\mathrm{Tot}(\tld{\mathrm{CE}}_{\Delta\setminus(I_0\setminus J),\Delta}^{\bullet,\bullet})) \ar^{\sim}[rr] \ar^{\wr}[d] & & H^{h'-\#J}(\mathrm{Tot}(\tld{\mathrm{CE}}_{\Delta\setminus I_0,\Delta\setminus J}^{\bullet,\bullet})) \ar^{\wr}[d]\\
\mathrm{Ext}_{U(\fg)}^{\#I_0\setminus J}(\mathfrak{v}_{\Delta\setminus(I_0\setminus J)},1_{\fg}) \ar^{\sim}[rr] \ar^{\wr}[d] & & \mathrm{Ext}_{U(\fg)}^{\#I_0\setminus J}(U(\fg)\otimes_{U(\fp_{\Delta\setminus J})}\mathfrak{v}_{\Delta\setminus I_0,\Delta\setminus J},1_{\fg}) \ar^{\wr}[d]\\
\bigoplus_{x\in\Gamma_{I_0\setminus J}}\mathfrak{e}_{\Delta\setminus(I_0\setminus J),x} \ar^{\sim}[rr] & & \bigoplus_{x\in\Gamma_{I_0\setminus J}}\mathfrak{e}_{\Delta\setminus I_0,x}
}
\end{equation}
with the bottom horizontal map of (\ref{equ: graded coxeter tau diagram 2}) being the direct sum over $x\in\Gamma_{I_0\setminus J}$ of the natural isomorphism between $1$-dimensional $E$-vector spaces $\mathfrak{e}_{\Delta\setminus(I_0\setminus J),x}\buildrel\sim\over\longrightarrow \mathfrak{e}_{\Delta\setminus I_0,x}$.
We conclude by combining (\ref{equ: graded coxeter tau diagram 1}), (\ref{equ: graded coxeter tau diagram 2}), Lemma~\ref{lem: cup with sm} and Lemma~\ref{lem: graded tensor as cup} (upon replacing $I_0$ and $I_0'$ in \emph{loc.cit.} with $I_0\setminus J$ and $J$ here).
\end{proof}

Let $I_0,I_0'\subseteq \Delta$ with $I_0\cap I_0'=\emptyset$ and $x\in\Gamma^{I_0}$. Now that we have
\[\kappa_{I_0,I_0'}(\mathrm{Fil}_{x'}(\mathbf{E}_{I_0})\otimes_E\mathrm{Fil}_{1}(\mathbf{E}_{I_0'}))=\mathrm{Fil}_{x'}(\mathbf{E}_{I_0\sqcup I_0'})\]
for each $x'\in\Gamma^{I_0}$ satisfying $x\not\unlhd x'$ (see Lemma~\ref{lem: cup with sm}), we have
\[\kappa_{I_0,I_0'}(\sum_{x\not\unlhd x'}\mathrm{Fil}_{x'}(\mathbf{E}_{I_0})\otimes_E\mathrm{Fil}_{1}(\mathbf{E}_{I_0'}))=\sum_{x\not\unlhd x'}\mathrm{Fil}_{x'}(\mathbf{E}_{I_0\sqcup I_0'}),\]
which together with Proposition~\ref{prop: coxeter filtration x surjection} (and the injectivity of $\kappa_{I_0,I_0'}$ from \ref{it: general cup 1} of Theorem~\ref{thm: general cup}) induces an injection
\begin{multline}\label{equ: coxeter quotient cup product 1}
\mathbf{E}_{x,I_0}\otimes_E\mathrm{Fil}_{1}(\mathbf{E}_{I_0'}))=(\mathbf{E}_{I_0}/\sum_{x\not\unlhd x'}\mathrm{Fil}_{x'}(\mathbf{E}_{I_0}))\otimes_E\mathrm{Fil}_{1}(\mathbf{E}_{I_0'})\\
\buildrel\sim\over\longrightarrow (\mathbf{E}_{I_0}\cup\mathrm{Fil}_{1}(\mathbf{E}_{I_0'}))/(\sum_{x\not\unlhd x'}\mathrm{Fil}_{x'}(\mathbf{E}_{I_0}))\cup\mathrm{Fil}_{1}(\mathbf{E}_{I_0'}))\\
=(\mathbf{E}_{I_0}\cup\mathrm{Fil}_{1}(\mathbf{E}_{I_0'}))/\sum_{x\not\unlhd x'}\mathrm{Fil}_{x'}(\mathbf{E}_{I_0\sqcup I_0'})\\
\hookrightarrow \mathbf{E}_{I_0\sqcup I_0'}/\sum_{x\not\unlhd x'}\mathrm{Fil}_{x'}(\mathbf{E}_{I_0\sqcup I_0'})=\mathbf{E}_{x,I_0\sqcup I_0'}
\end{multline}
which also comes as the restriction of the following cup product map (see (\ref{equ: x y cup w bottom diagram}))
\begin{equation}\label{equ: coxeter quotient cup product 2}
\mathbf{E}_{x,I_0}\otimes_E\mathbf{E}_{I_0'}\buildrel\cup\over\longrightarrow \mathbf{E}_{x,I_0\sqcup I_0'}.
\end{equation}

\begin{lem}\label{lem: u sm cup}
Let $u\in\Gamma$ with $J_{u}$ being an interval and let $J\subseteq J_{u}$.
Then the following $1$-dimensional $E$-subspace of $\mathbf{E}_{u,\Delta\setminus J}$ (see below (\ref{equ: g P u bottom}) and Lemma~\ref{lem: S u sm sub})
\begin{equation}\label{equ: u sm cup sub}
H^{\#\Delta-2\#J}(\mathrm{Tot}(\cS_{u,J,J_{u},\flat}^{\bullet,\bullet}))^{\infty}
\end{equation}
equals the image of the following injective map (see (\ref{equ: coxeter quotient cup product 1}) and (\ref{equ: coxeter quotient cup product 2}))
\begin{equation}\label{equ: u cup injection}
\mathbf{E}_{u,\Delta\setminus J_{u}}\otimes_E\mathrm{Fil}_{1}(\mathbf{E}_{J_{u}\setminus J})\buildrel\cup\over\longrightarrow \mathbf{E}_{u,\Delta\setminus J}.
\end{equation}
\end{lem}
\begin{proof}
We write $J'\defeq \Delta\setminus (J_{u}\setminus J)=J\sqcup\mathrm{Supp}(u)$ for short.
Note that (\ref{equ: u cup injection}) factors through
\begin{multline}\label{equ: u sm cup map seq}
\mathbf{E}_{u,\Delta\setminus J_{u}}\otimes_E\mathrm{Fil}_{1}(\mathbf{E}_{J_{u}\setminus J})\\
=\mathrm{Ext}_{D(\fg,P_{J_{u}})}^{\ell(u)}(M^{J_{u}}(u),1_{D(\fg,P_{J_{u}})})\otimes_EH^{\#J_{u}\setminus J}(L_{J'},1_{L_{J'}})^{\infty}\\
\rightarrow \mathrm{Ext}_{D(\fg,P_{J})}^{\ell(u)}(M^{J}(u),1_{D(\fg,P_{J})})\otimes_EH^{\#J_{u}\setminus J}(L_{J},1_{L_{J}})^{\infty}\\
\rightarrow\mathrm{Ext}_{D(\fg,P_{J})}^{\ell(u)}(M^{J}(u),1_{D(\fg,P_{J})})\otimes_E\mathrm{Ext}_{D(\fg,P_{J})}^{\#J_{u}\setminus J}(M^{J}(1),1_{D(\fg,P_{J})})\\
\rightarrow \mathrm{Ext}_{D(\fg,P_{J})}^{\#\Delta\setminus J}(M^{J}(u),1_{D(\fg,P_{J})})\rightarrow \mathbf{E}_{u,\Delta\setminus J}
\end{multline}
and thus its image is contained in the subspace (\ref{equ: u sm cup sub}) thanks to Lemma~\ref{lem: g P Ext transfer}, (\ref{equ: x sm cup}) and (\ref{equ: x sm cup Lie transfer}) (by taking $x$ and $I$ in \emph{loc.cit.} to be $u$ and $J$).
Now that (\ref{equ: u cup injection}) is known to be an injection by (\ref{equ: coxeter quotient cup product 1}), its image must be exactly the $1$-dimensional subspace (\ref{equ: u sm cup sub}).
\end{proof}

\begin{lem}\label{lem: Lie cup x y nonvanishing}
Let $I_0,I_0'\subseteq\Delta$ with $I_0\cap I_0'=\emptyset$. Let $x\in\Gamma_{I_0}$, $y\in\Gamma_{I_0'}$ and $w\in\Gamma_{x,y}\subseteq\Gamma_{I_0\sqcup I_0'}$. Then the following map (see (\ref{equ: x y cup w bottom}))
\begin{equation}\label{equ: Lie cup x y nonvanishing}
\mathbf{E}_{x,I_0}\otimes_E\mathbf{E}_{y,I_0'}\buildrel\cup\over\longrightarrow\mathbf{E}_{w,I_0\sqcup I_0'}
\end{equation}
is an isomorphism between $1$ dimensional $E$-vector spaces.
\end{lem}
\begin{proof}
We note that $J_x=\Delta\setminus I_0$, $J_{y}\defeq\Delta\setminus I_0'$ as well as $J_{w}=\Delta\setminus(I_0\sqcup I_0')$.
Thanks to the commutative diagram (\ref{equ: x y bottom cup to Lie}), we conclude using Proposition~\ref{prop: Lie bottom cup isom}.
\end{proof}

\begin{prop}\label{prop: cup x y nonvanishing}
Let $I_0,I_0'\subseteq \Delta$ with $I_0\cap I_0'=\emptyset$.
Let $x\in\Gamma^{I_0}$, $y\in\Gamma^{I_0'}$ and $w\in\Gamma_{x,y}\subseteq \Gamma^{I_0\sqcup I_0'}$ (see Proposition~\ref{prop: coxeter cup comparison}).
Then the composition of the following maps (see Proposition~\ref{prop: coxeter cup comparison})
\begin{equation}\label{equ: cup x y nonvanishing}
\mathrm{Fil}_{x}(\mathbf{E}_{I_0})\otimes_E\mathrm{Fil}_{y}(\mathbf{E}_{I_0'})\hookrightarrow \mathbf{E}_{I_0}\otimes_E\mathbf{E}_{I_0'}\buildrel\kappa_{I_0,I_0'}\over\longrightarrow \sum_{u\in\Gamma_{x,y}}\mathrm{Fil}_{u}(\mathbf{E}_{I_0\sqcup I_0'})\twoheadrightarrow \mathrm{gr}_{w}(\mathbf{E}_{I_0\sqcup I_0'})
\end{equation}
is surjective.
\end{prop}
\begin{proof}
We write $I_1\defeq\mathrm{Supp}(x)$ and $I_1'\defeq\mathrm{Supp}(y)$ for short with $\mathrm{Supp}(w)=I_1\sqcup I_1'$.  
Note that the composition of (\ref{equ: cup x y nonvanishing}) fits into the following commutative diagram
\begin{equation}\label{equ: cup x y nonvanishing diagram}
\xymatrix{
\mathrm{Fil}_{x}(\mathbf{E}_{I_0})\otimes_E\mathrm{Fil}_{y}(\mathbf{E}_{I_0'}) \ar^{\sim}[r] \ar[d]& \mathrm{Fil}_{1}(\mathbf{E}_{I_0\setminus I_1})\otimes_E\mathrm{Fil}_{x}(\mathbf{E}_{I_1})\otimes_E\mathrm{Fil}_{y}(\mathbf{E}_{I_1'})\otimes_E\mathrm{Fil}_{1}(\mathbf{E}_{I_0'\setminus I_1'}) \ar[d]\\
\sum_{u\in\Gamma_{x,y}}\mathrm{Fil}_{u}(\mathbf{E}_{I_0\sqcup I_0'}) \ar^{\sim}[r] \ar[d]& \sum_{u\in\Gamma_{x,y}}\mathrm{Fil}_{1}(\mathbf{E}_{I_0\setminus I_1})\otimes_E\mathrm{Fil}_{u}(\mathbf{E}_{I_1\sqcup I_1'})\otimes_E\mathrm{Fil}_{1}(\mathbf{E}_{I_0'\setminus I_1'}) \ar[d]\\
\mathrm{gr}_{w}(\mathbf{E}_{I_0\sqcup I_0'}) \ar^{\sim}[r] & \mathrm{Fil}_{1}(\mathbf{E}_{I_0\setminus I_1})\otimes_E\mathrm{gr}_{w}(\mathbf{E}_{I_1\sqcup I_1'})\otimes_E\mathrm{Fil}_{1}(\mathbf{E}_{I_0'\setminus I_1'})
}
\end{equation}
with all horizontal maps being isomorphisms by Lemma~\ref{lem: cup with sm}.
Hence, it is harmless to assume in the rest of the proof that $I_0=\mathrm{Supp}(x)$ and $I_0'=\mathrm{Supp}(y)$ (with $I_0\sqcup I_0'=\mathrm{Supp}(w)$), in which case the composition of $\mathrm{Fil}_{x}(\mathbf{E}_{I_0})\hookrightarrow\mathbf{E}_{I_0}\twoheadrightarrow\mathbf{E}_{x,I_0}$ can be identified with 
\[\mathrm{Fil}_{x}(\mathbf{E}_{I_0})\twoheadrightarrow\mathrm{gr}_{x}(\mathbf{E}_{I_0})=\mathbf{E}_{x,I_0},\]
and similar facts hold with the pairs $(y,I_0')$ and $(w,I_0\sqcup I_0')$ replacing the pair $(x,I_0)$.
In particular, under the assumption that $I_0=\mathrm{Supp}(x)$ and $I_0'=\mathrm{Supp}(y)$, the composition of (\ref{equ: cup x y nonvanishing}) factors through
\begin{equation}\label{equ: cup x y nonvanishing factor}
\mathrm{Fil}_{x}(\mathbf{E}_{I_0})\otimes_E\mathrm{Fil}_{y}(\mathbf{E}_{I_0'})\twoheadrightarrow \mathrm{gr}_{x}(\mathbf{E}_{I_0})\otimes_E\mathrm{gr}_{y}(\mathbf{E}_{I_0'})=\mathbf{E}_{x,I_0}\otimes_E\mathbf{E}_{y,I_0'}\buildrel q\over\rightarrow\mathbf{E}_{w,I_0\sqcup I_0'}=\mathrm{gr}_{w}(\mathbf{E}_{I_0\sqcup I_0'}).
\end{equation}
Now that the map $q$ is an isomorphism between $1$ dimensional $E$-vector spaces by Lemma~\ref{lem: Lie cup x y nonvanishing}, we conclude that the composition of (\ref{equ: cup x y nonvanishing factor}) is a surjection and thus finish the proof.
\end{proof}

Given an interval $I\subseteq \Delta$, we recall the $E$-subspace $\mathbf{E}_{I}^{<}\subseteq \mathbf{E}_{I}$ from (\ref{equ: decomposable subspace}) and Lemma~\ref{lem: primitive generator}.
\begin{thm}\label{thm: cup top grade}
Let $I\subseteq \Delta$ be a non-empty interval.
Then the following composition
\begin{equation}\label{equ: cup top grade}
\mathrm{Fil}_{w}(\mathbf{E}_{I})\hookrightarrow \mathbf{E}_{I}\twoheadrightarrow \mathbf{E}_{I}/\mathbf{E}_{I}^{<}
\end{equation}
is surjective for each $w\in\Gamma_{I}$.
\end{thm}
\begin{proof}
If $\#I=1$, then we have $\mathrm{Fil}_{w}(\mathbf{E}_{I})=\mathbf{E}_{I}$ and $\mathbf{E}_{I}^{<}$ with the claim of this theorem being evident.
Assume from now on that $\#I\geq 2$.
We consider the following map
\[q_{I}: \mathbf{E}_{I}\twoheadrightarrow \mathbf{E}_{I}/\sum_{x\in\Gamma^{I}\setminus \Gamma_{I}}\mathrm{Fil}_{x}(\mathbf{E}_{I})=\bigoplus_{w\in\Gamma_{I}}\mathrm{gr}_{w}(\mathbf{E}_{I}).\]
By Lemma~\ref{lem: cup with sm} we know that
\[
\mathrm{Fil}_{x}(\mathbf{E}_{I})=\kappa_{\mathrm{Supp}(x),I\setminus\mathrm{Supp}(x)}(\mathrm{Fil}_{x}(\mathbf{E}_{\mathrm{Supp}(x)})\otimes_E\mathrm{Fil}_{1}(\mathbf{E}_{I\setminus\mathrm{Supp}(x)}))\subseteq \mathbf{E}_{I}^{<}
\]
for each $x\in\Gamma^{I}\setminus \Gamma_{I}$, and thus the surjection $\mathbf{E}_{I}\twoheadrightarrow \mathbf{E}_{I}/\mathbf{E}_{I}^{<}$ factors through the following surjections
\[
\mathbf{E}_{I}\buildrel q_{I} \over\longrightarrow \bigoplus_{w\in\Gamma_{I}}\mathrm{gr}_{w}(\mathbf{E}_{I})\buildrel q_{I}' \over\longrightarrow \mathbf{E}_{I}/\mathbf{E}_{I}^{<},
\]
with $\mathrm{ker}(q_{I}')=q_{I}(\mathbf{E}_{I}^{<})$.
For each $w\in \Gamma_{I}$, we choose an arbitrary $0\neq e_w\in \mathrm{gr}_{w}(\mathbf{E}_{I})$ so that $\{e_w\mid w\in\Gamma_{I}\}$ forms a basis of $\bigoplus_{w\in\Gamma_{I}}\mathrm{gr}_{w}(\mathbf{E}_{I})$.
We consider an adjacent pair of elements $w_1,w_2\in\Gamma_{I}$. As $I$ is a non-empty interval, by Lemma~\ref{lem: interval adjacent pair} we know that there exists a partition into non-empty subintervals $I=I_0\sqcup I_0'$ as well as $x\in\Gamma_{I_0}$ and $y\in\Gamma_{I_0'}$ such that $w_1=xy$ and $w_2=yx$. Hence, it follows from Proposition~\ref{prop: coxeter cup comparison}, Proposition~\ref{prop: cup x y nonvanishing} and Lemma~\ref{lem: interval Bruhat order} that the following composition
\[\mathrm{Fil}_{x}(\mathbf{E}_{I_0})\otimes_E\mathrm{Fil}_{y}(\mathbf{E}_{I_0'})\hookrightarrow \mathbf{E}_{I_0}\otimes_E\mathbf{E}_{I_0'}\buildrel \kappa_{I_0,I_0'}\over\longrightarrow \mathbf{E}_{I}\twoheadrightarrow \mathrm{gr}_{w}(\mathbf{E}_{I})\]
is non-zero if and only if $w\in\{w_1,w_2\}=\{xy,yx\}$. Consequently, we observe that
\begin{equation}\label{equ: cup image inclusion}
q_{I}(\mathbf{E}_{I}^{<})\supseteq q_{I}(\kappa_{I_0,I_0'}(\mathrm{Fil}_{x}(\mathbf{E}_{I_0})\otimes_E\mathrm{Fil}_{y}(\mathbf{E}_{I_0'})))
\end{equation}
contains an element $e_{w_1,w_2}$ which has the form $e_{w_1,w_2}=c_{w_1,w_2}^{w_1}e_{w_1}+c_{w_1,w_2}^{w_2}e_{w_2}$ for some $c_{w_1,w_2}^{w_1},c_{w_1,w_2}^{w_2}\in E^{\times}$. This forces the two $E$-lines $Ee_{w_1}$ and $Ee_{w_2}$ to have the same image under $q_{I}'$.
Now we consider arbitrary two different elements $w,w'\in\Gamma_{I}$. It is clear that there exists $w=x_0,x_1,\dots,x_t=w'$ for some $t\geq 1$ and $x_{t'}\in\Gamma_{I}$ for $0\leq t'\leq t$ such that the pair $x_{t'-1},x_{t'}$ is adjacent for each $1\leq t'\leq t$. This implies that
\begin{equation}\label{equ: cup common line}
Eq_{I}'(e_{w})=Eq_{I}'(e_{x_0})=\cdots=Eq_{I}'(e_{x_{t'}})=\cdots=Eq_{I}'(e_{x_t})=Eq_{I}'(e_{w'}),
\end{equation}
namely $Eq_{I}'(e_{w})$ does not depend on the choice of $w\in\Gamma_{I}$.
But since $q_{I}'$ is surjective, we must have
\[\sum_{w\in\Gamma_{I}}Eq_{I}'(e_{w})=\mathbf{E}_{I}/\mathbf{E}_{I}^{<},\]
which forces
\[q_{I}'(\mathrm{gr}_{w}(\mathbf{E}_{I}))=Eq_{I}'(e_{w})=\mathbf{E}_{I}/\mathbf{E}_{I}^{<}\]
for each $w\in\Gamma_{I}$. Consequently, we conclude that (\ref{equ: cup top grade}) is surjective for each $w\in\Gamma_{I}$.
\end{proof}

\subsection{Commutativity of cup product}\label{subsec: commutativity}
We prove our main result (see Theorem~\ref{thm: cup exchange}) concerning the commutativity of our normalized cup product maps with respect to exchanging the two indices.

Let $I_0,I_0'\subseteq \Delta$ with $I_0\cap I_0'=\emptyset$.
We write $J_0\defeq \Delta\setminus I_0$, $J_0'\defeq \Delta\setminus I_0'$ and $r\defeq \#\Delta=n-1$ for short.
We continue our study of the map (see (\ref{equ: main bottom cup}))
\begin{equation}\label{equ: Tits cup coh 1}
\kappa_{I_0,I_0'}: \mathbf{E}_{I_0}\otimes_E\mathbf{E}_{I_0'}\buildrel\cup\over\longrightarrow \mathbf{E}_{I_0\sqcup I_0'}
\end{equation}
as well as the map
\begin{equation}\label{equ: Tits cup coh 2}
\kappa_{I_0',I_0}: \mathbf{E}_{I_0'}\otimes_E\mathbf{E}_{I_0}\buildrel\cup\over\longrightarrow \mathbf{E}_{I_0\sqcup I_0'}.
\end{equation}

We recall the $E$-subspaces $\mathbf{E}_{I_0}^{<}\subseteq\mathbf{E}_{I_0}$ and $\mathbf{E}_{I_0'}^{<}\subseteq\mathbf{E}_{I_0'}$ from (\ref{equ: decomposable subspace}).
\begin{prop}\label{prop: cup of line}
Assume that both $I_0$ and $I_0'$ are intervals.
Then there exists $1$-dimensional subspaces $\mathbf{E}_{I_0}^{\heartsuit}\subseteq \mathbf{E}_{I_0}$ and $\mathbf{E}_{I_0'}^{\heartsuit}\subseteq \mathbf{E}_{I_0'}$ that satisfy the following conditions
\begin{enumerate}[label=(\roman*)]
\item \label{it: cup of line 1} We have $\mathbf{E}_{I_0}^{\heartsuit}\cap\mathbf{E}_{I_0}^{<}=0$ and $\mathbf{E}_{I_0'}^{\heartsuit}\cap\mathbf{E}_{I_0'}^{<}=0$.
\item \label{it: cup of line 2} We have 
\begin{equation}\label{equ: cup of line exchange}
x\cup y=y\cup x \in \mathbf{E}_{I_0\sqcup I_0'}
\end{equation}
for each $x\in\mathbf{E}_{I_0}^{\heartsuit}$ and $y\in \mathbf{E}_{I_0'}^{\heartsuit}$.
\end{enumerate}
\end{prop}
\begin{proof}
Upon exchanging $I_0$ and $I_0'$, it is harmless to assume that $I_0<I_0'$ (see \ref{it: disconnected 2} of Definition~\ref{def: disconnected pair}).
We write $\Delta^{-}\subseteq\Delta$ for the minimal interval that contains both $I_0$ and $I_0'$ with $m\defeq \#\Delta_{-}$. 
Recall from (\ref{equ: Tits interval reduction cup}) that the maps $G_{m+1}\cong H_{\Delta_{-}}\subseteq G$ induces the following commutative diagram
\begin{equation}\label{equ: cup interval reduction}
\xymatrix{
\mathbf{E}_{I_0} \ar^{\wr}[d] & \otimes_E & \mathbf{E}_{I_0'} \ar^{\cup}[r] \ar^{\wr}[d] & \mathbf{E}_{I_0\sqcup I_0'} \ar^{\wr}[d]\\
\mathbf{E}_{m,I_0} & \otimes_E & \mathbf{E}_{m,I_0'} \ar^{\cup}[r] & \mathbf{E}_{m,I_0\sqcup I_0'}
}
\end{equation}
with all vertical maps being isomorphisms.
Hence, upon replacing $G$ with $G_{m+1}$, it is harmless to assume in the rest of the proof that $I_0=[1,j_0]$ and $I_0'=[j_1,n-1]$ for some $j_0,j_1\in\Delta$ with $j_0<j_1$. 
We write $J_0\defeq\Delta\setminus I_0=[j_0+1,n-1]$, $J_1\defeq \Delta\setminus I_0'=[1,j_1-1]$ and $J_2\defeq J_0\cap J_1$ for short.
Given such $j_0,j_1\in\Delta$, we recall from (\ref{equ: relative total a}) that the maps (\ref{equ: special relative j0}) induce embeddings
\begin{equation}\label{equ: cup line embedding 0}
\mathbf{E}_{0,\dagger}\defeq H^{2\#I_0-r}(\mathrm{Tot}(\cT_{0,\dagger}^{\bullet,\bullet}))\hookrightarrow H^{2\#I_0-r}(\mathrm{Tot}(\cT_{J_0,\Delta}^{\bullet,\bullet}))=\mathbf{E}_{I_0}.
\end{equation}
Similarly, the maps (\ref{equ: special relative j1}) induce embeddings
\begin{equation}\label{equ: cup line embedding 1}
\mathbf{E}_{1,\dagger}\defeq H^{2\#I_0'-r}(\mathrm{Tot}(\cT_{1,\dagger}^{\bullet,\bullet}))\hookrightarrow H^{2\#I_0'-r}(\mathrm{Tot}(\cT_{J_1,\Delta}^{\bullet,\bullet}))=\mathbf{E}_{I_0'},
\end{equation}
and the maps (\ref{equ: special relative j0 j1}) induce embeddings
\begin{equation}\label{equ: cup line embedding 2}
\mathbf{E}_{2,\dagger}\defeq H^{2\#I_0\sqcup I_0'-r}(\mathrm{Tot}(\cT_{2,\dagger}^{\bullet,\bullet}))\hookrightarrow H^{2\#I_0\sqcup I_0'-r}(\mathrm{Tot}(\cT_{J_2,\Delta}^{\bullet,\bullet}))=\mathbf{E}_{I_0\sqcup I_0'}.
\end{equation}
Parallel to the definition of $\kappa_{I_0,I_0'}$, we can define a map \[\mathbf{E}_{0,\dagger}\otimes_E\mathbf{E}_{1,\dagger}\buildrel\cup\over\longrightarrow \mathbf{E}_{2,\dagger}\]
which fits into the following commutative diagram
\begin{equation}\label{equ: cup line diagram}
\xymatrix{
\mathbf{E}_{0,\dagger} \ar@{^{(}->}[d] & \otimes_E & \mathbf{E}_{1,\dagger} \ar^{\cup}[r] \ar@{^{(}->}[d] & \mathbf{E}_{2,\dagger} \ar@{^{(}->}[d]\\
\mathbf{E}_{I_0} & \otimes_E & \mathbf{E}_{I_0'} \ar^{\cup}[r] & \mathbf{E}_{I_0\sqcup I_0'}
}
\end{equation}
with the vertical maps from (\ref{equ: cup line embedding 0}), (\ref{equ: cup line embedding 1}) and (\ref{equ: cup line embedding 2}) being embeddings.
Upon identifying $\mathbf{E}_{0,\dagger}$ with its image in $\mathbf{E}_{I_0}$, and similarly for $\mathbf{E}_{1,\dagger}$ and $\mathbf{E}_{2,\dagger}$, we deduce from (\ref{equ: cup line diagram}) that
\begin{equation}\label{equ: cup line dagger 1}
\kappa_{I_0,I_0'}(\mathbf{E}_{0,\dagger}\otimes_E\mathbf{E}_{1,\dagger})\subseteq \mathbf{E}_{2,\dagger}.
\end{equation}
We write $w_{1,j}\defeq s_1\cdots s_{j}$ for each $0\leq j\leq j_0$ (with convention $w_{1,0}=1$) and $w_{n-1,j}\defeq s_{n-1}\cdots s_{j}$ for each $j_1\leq j\leq n$ (with convention $w_{n-1,n}=1$). 
We define $\Gamma_0\defeq \Gamma^{[2,j_0]}\sqcup\{w_{1,j_0}\}$, $\Gamma_1\defeq \Gamma^{[j_1,n-2]}\sqcup\{w_{n-1,j_1}\}$ and $\Gamma_{2}\defeq \bigcup_{x\in\Gamma_{0},y\in\Gamma_{1}}\Gamma_{x,y}$.
We write $\mathrm{Fil}_{\Gamma_{0}}(\mathbf{E}_{I_0})\defeq \sum_{x\in\Gamma_{0}}\mathrm{Fil}_{x}(\mathbf{E}_{I_0})$ for short and similarly for $\mathrm{Fil}_{\Gamma_{1}}(\mathbf{E}_{I_0'})$ and $\mathrm{Fil}_{\Gamma_{2}}(\mathbf{E}_{I_0\sqcup I_0})$.
We define $\mathbf{E}_{I_0}^{\heartsuit}\defeq \mathbf{E}_{0,\dagger}\cap\mathrm{Fil}_{\Gamma_{0}}(\mathbf{E}_{I_0})$,
$\mathbf{E}_{I_0'}^{\heartsuit}\defeq \mathbf{E}_{1,\dagger}\cap\mathrm{Fil}_{\Gamma_{1}}(\mathbf{E}_{I_0'})$
and $\mathbf{E}_{I_0,I_0'}^{\heartsuit}\defeq \mathbf{E}_{2,\dagger}\cap\mathrm{Fil}_{\Gamma_{2}}(\mathbf{E}_{I_0\sqcup I_0'})$.
It follows from Proposition~\ref{prop: coxeter cup comparison} that
\[\kappa_{I_0,I_0'}(\mathrm{Fil}_{\Gamma_{0}}(\mathbf{E}_{I_0})\otimes_E \mathrm{Fil}_{\Gamma_{1}}(\mathbf{E}_{I_0'}))\subseteq \mathrm{Fil}_{\Gamma_{2}}(\mathbf{E}_{I_0\sqcup I_0'}),\]
which together with (\ref{equ: cup line dagger 1}) and the fact that $\mathrm{im}(\kappa_{I_0,I_0'})\subseteq\mathbf{E}_{I_0\sqcup I_0'}^{<}$ gives
\begin{equation}\label{equ: cup line bound 1}
\kappa_{I_0,I_0'}(\mathbf{E}_{I_0}^{\heartsuit}\otimes_E\mathbf{E}_{I_0'}^{\heartsuit})\subseteq \mathbf{E}_{I_0\sqcup I_0'}^{<}\cap \mathbf{E}_{I_0,I_0'}^{\heartsuit}.
\end{equation}
Symmetrically, we also have
\[\kappa_{I_0',I_0}(\mathbf{E}_{1,\dagger}\otimes_E\mathbf{E}_{0,\dagger})\subseteq \mathbf{E}_{2,\dagger},\]
which together with Proposition~\ref{prop: coxeter cup comparison} and the fact that $\mathrm{im}(\kappa_{I_0,I_0'})\subseteq\mathbf{E}_{I_0\sqcup I_0'}^{<}$ gives 
\begin{equation}\label{equ: cup line bound 2}
\kappa_{I_0',I_0}(\mathbf{E}_{I_0'}^{\heartsuit}\otimes_E\mathbf{E}_{I_0}^{\heartsuit})\subseteq \mathbf{E}_{I_0\sqcup I_0'}^{<}\cap \mathbf{E}_{I_0,I_0'}^{\heartsuit}.
\end{equation}
We divide the rest of the proof into the following steps.

\textbf{Step $1$}: We prove that $\mathbf{E}_{I_0}^{\heartsuit}$ is a $1$-dimensional $E$-subspace of $\mathbf{E}_{I_0}$ that satisfies $\mathbf{E}_{I_0}=\mathbf{E}_{I_0}^{\heartsuit}\oplus\mathbf{E}_{I_0}^{<}$, and similarly $\mathbf{E}_{I_0'}^{\heartsuit}$ is a $1$-dimensional $E$-subspace of $\mathbf{E}_{I_0'}$ that satisfies $\mathbf{E}_{I_0'}=\mathbf{E}_{I_0'}^{\heartsuit}\oplus\mathbf{E}_{I_0'}^{<}$.\\
It suffices to treat $\mathbf{E}_{I_0}^{\heartsuit}$ as the argument for $\mathbf{E}_{I_0'}^{\heartsuit}$ is parallel.
We write $\Gamma_{0}^+\defeq \Gamma^{[2,j_0]}\sqcup\{w_{1,j}\mid j\in[1,j_0]\}$ for short.
Note from \cite[Thm.~2.2.2]{BB05} that each $x\in\Gamma^{I_0}\setminus\Gamma^{[2,j_0]}$ (namely with $1\in\mathrm{Supp}(x)$) determines a unique $j\in[1,j_0]$ and $x'\in\Gamma_{[j+1,j_0]}$ such that $x=x'w_{1,j}$, namely $w_{1,j}\unlhd x$ with $x\in\Gamma_{0}^+$ if and only if $x'=1$.
In other words, an element $x\in\Gamma^{I_0}$ satisfies $x\notin\Gamma_0^+$ if and only if there exists $1\leq j<j_0$ such that $w_{1,j}\unlhd x$ and $w_{1,j}\neq x$.
Recall from \ref{it: 0 1 2 sm part 0} of Propositoin~\ref{prop: 0 1 2 sm part} (by taking $x=w_{1,j}$ in \emph{loc.cit.}) and Lemma~\ref{lem: u sm cup} that the composition of the following maps
\[\mathbf{E}_{0,\dagger}\hookrightarrow\mathbf{E}_{I_0}\buildrel \zeta_{j,I_0}\over\twoheadrightarrow\mathbf{E}_{w_{1,j},I_0}=\mathbf{E}_{I_0}/\sum_{u\in\Gamma^{I_0},w_{1,j}\not\unlhd u}\mathrm{Fil}_{u}\mathbf{E}_{I_0}\]
factors through $\mathrm{gr}_{w_{1,j}}(\mathbf{E}_{I_0})$ for each $j\in[1,j_0]$.
In other words, we have
\begin{equation}\label{equ: E 0 coxeter bound}
\mathbf{E}_{0,\dagger}\subseteq\bigcap_{j\in[1,j_0]}\zeta_{j,I_0}^{-1}(\mathrm{gr}_{w_{1,j}}(\mathbf{E}_{I_0}))
=\mathrm{Fil}_{\Gamma_{0}^+}(\mathbf{E}_{I_0})\defeq \sum_{x\in\Gamma_{0}^{+}}\mathrm{Fil}_{x}(\mathbf{E}_{I_0}).
\end{equation}
Recall from \ref{it: relative coh E2 grade dagger 0} of Proposition~\ref{prop: relative coh E2 grade dagger} that the decreasing filtration
\begin{equation}\label{equ: E 0 sm support Fil}
\{\tau^{i}(\mathbf{E}_{0,\dagger})\}_{i\geq 0}
\end{equation}
satisfies $\tau^{i}(\mathbf{E}_{0,\dagger})=\tau^{i}(\mathbf{E}_{I_0})\cap\mathbf{E}_{0,\dagger}$ for each $i\geq 0$, with $\mathrm{gr}^{i}(\mathbf{E}_{0,\dagger})\neq 0$ if and only if $0\leq i<j_0$, in which case 
\begin{equation}\label{equ: E 0 sm support grade}
\mathrm{gr}^{i}(\mathbf{E}_{0,\dagger})\buildrel\sim\over\longrightarrow\mathfrak{e}_{J_0,w_{1,j_0-i}}\otimes_E\mathbf{E}_{[j_0-i+1,j_0]}^{\infty}\subseteq \mathrm{gr}^{i}(\mathbf{E}_{I_0}).
\end{equation}
We write $\mathbf{E}_{0,\dagger}'\defeq \tau^{1}(\mathbf{E}_{0,\dagger})$ for short. It follows from (\ref{equ: E 0 sm support grade}) and Proposition~\ref{prop: graded coxeter tau comparison} that 
\begin{equation}\label{equ: E 0 2 dim intersection zero}
\mathbf{E}_{0,\dagger}'\cap \mathrm{Fil}_{\Gamma_{0}}(\mathbf{E}_{I_0})=0
\end{equation}
and thus
\[\Dim_E(\mathbf{E}_{0,\dagger}'+\mathrm{Fil}_{\Gamma_{0}}(\mathbf{E}_{I_0}))=\Dim_E\mathbf{E}_{0,\dagger}'+\Dim_E\mathrm{Fil}_{\Gamma_{0}}(\mathbf{E}_{I_0})=\Dim_E \mathrm{Fil}_{\Gamma_{0}^{+}}(\mathbf{E}_{I_0}),\]
which together with (\ref{equ: E 0 coxeter bound}) gives
\[
\mathbf{E}_{0,\dagger}+\mathrm{Fil}_{\Gamma_{0}}(\mathbf{E}_{I_0})=\mathbf{E}_{0,\dagger}'+\mathrm{Fil}_{\Gamma_{0}}(\mathbf{E}_{I_0})=\mathrm{Fil}_{\Gamma_{0}^{+}}(\mathbf{E}_{I_0}).
\]
Consequently, we have
\begin{multline}\label{equ: E 0 2 dim intersection}
\Dim_E\mathbf{E}_{I_0}^{\heartsuit}=\Dim_E\mathbf{E}_{0,\dagger}\cap\mathrm{Fil}_{\Gamma_{0}}(\mathbf{E}_{I_0})\\
=\Dim_E\mathbf{E}_{0,\dagger}+\Dim_E\mathrm{Fil}_{\Gamma_{0}}(\mathbf{E}_{I_0})-\Dim_E\mathrm{Fil}_{\Gamma_{0}^+}(\mathbf{E}_{I_0})\\
=j_0+\#\Gamma_{0}-\#\Gamma_{0}^{+}=1
\end{multline}
which together with (\ref{equ: E 0 2 dim intersection zero}) forces the composition of
\[\mathbf{E}_{I_0}^{\heartsuit}\hookrightarrow \mathbf{E}_{0,\dagger}\twoheadrightarrow \mathbf{E}_{0,\dagger}/\mathbf{E}_{0,\dagger}'=\mathrm{gr}^{\emptyset}(\mathbf{E}_{0,\dagger})\cong \mathfrak{e}_{J_0,w_{1,j_0}}\]
to be an isomorphism between $1$-dimensional $E$-vector spaces.
In other words, the composition of
\[\mathbf{E}_{I_0}^{\heartsuit}\hookrightarrow\mathbf{E}_{I_0}\twoheadrightarrow \mathrm{gr}^{\emptyset}(\mathbf{E}_{I_0})\]
is an embedding with image $\mathfrak{e}_{J_0,w_{1,j_0}}=\mathrm{gr}_{w_{1,j_0}}(\mathbf{E}_{I_0})$ (see Proposition~\ref{prop: graded coxeter tau comparison}).
Now that the composition of
\[\mathfrak{e}_{J_0,w_{1,j_0}}=\mathrm{gr}_{w_{1,j_0}}(\mathbf{E}_{I_0})\hookrightarrow\mathrm{gr}^{\emptyset}(\mathbf{E}_{I_0})\twoheadrightarrow \mathrm{gr}^{1-r}(\mathbf{E}_{I_0})=\mathbf{E}_{I_0}/\mathbf{E}_{I_0}^{<}\]
is an isomorphism by Theorem~\ref{thm: cup top grade}, we conclude that the composition of
\[\mathbf{E}_{I_0}^{\heartsuit}\hookrightarrow\mathbf{E}_{I_0}\twoheadrightarrow \mathrm{gr}^{\emptyset}(\mathbf{E}_{I_0})\twoheadrightarrow \mathrm{gr}^{1-r}(\mathbf{E}_{I_0})=\mathbf{E}_{I_0}/\mathbf{E}_{I_0}^{<}\]
is an isomorphism between $1$-dimensional $E$-vector spaces.
In other words, we have $\mathbf{E}_{I_0}=\mathbf{E}_{I_0}^{\heartsuit}\oplus\mathbf{E}_{I_0}^{<}$.

\textbf{Step $2$}: We prove that the composition of
\begin{equation}\label{equ: cup line Lie quotient}
\mathbf{E}_{I_0,I_0'}^{\heartsuit}\hookrightarrow\mathbf{E}_{I_0\sqcup I_0'}\twoheadrightarrow\mathrm{gr}^{\emptyset}(\mathbf{E}_{I_0\sqcup I_0'})\cong\bigoplus_{u\in\Gamma_{I_0\sqcup I_0'}}\mathfrak{e}_{J_{2},u}=\bigoplus_{u\in\Gamma_{I_0\sqcup I_0'}}\mathrm{gr}_{u}(\mathbf{E}_{I_0\sqcup I_0'})
\end{equation}
is an embedding with image contained in 
\begin{equation}\label{equ: cup line Lie quotient bound}
\bigoplus_{u\in\Gamma_{I_0\sqcup I_0'}\cap\Gamma_{2}}\mathfrak{e}_{J_{2},u}=\bigoplus_{u\in\Gamma_{I_0\sqcup I_0'}\cap\Gamma_{2}}\mathrm{gr}_{u}(\mathbf{E}_{I_0\sqcup I_0'}).
\end{equation}
It follows from \ref{it: relative coh E2 grade dagger 2} of Proposition~\ref{prop: relative coh E2 grade dagger} and Proposition~\ref{prop: graded coxeter tau comparison} that the decreasing filtration $\{\tau^{i}(\mathbf{E}_{I_0\sqcup I_0'})\}_{i\geq 0}$ on $\mathbf{E}_{I_0\sqcup I_0'}$ restricts to the decreasing filtration $\{\tau^{i}(\mathbf{E}_{2,\dagger})\}_{i\geq 0}$ on $\mathbf{E}_{2,\dagger}$, and the induced map
\[\mathrm{gr}^{i}(\mathbf{E}_{2,\dagger})\rightarrow \mathrm{gr}^{i}(\mathbf{E}_{I_0\sqcup I_0'})\]
is an embedding with image contained in
\begin{equation}\label{equ: cup line dagger grade bound}
\bigoplus_{u\in\Gamma^{I_0\sqcup I_0'}, D_L(u)=\{1,n-1\},\ell(u)=\#I_0\sqcup I_0'-i}\mathrm{gr}_{u}(\mathbf{E}_{I_0\sqcup I_0'}).
\end{equation}
We also recall from \ref{it: tau and coxeter 2} of Lemma~\ref{lem: tau and coxeter} that
\[\tau^{i}(\mathbf{E}_{I_0\sqcup I_0'})=\sum_{u\in\Gamma^{I_0\sqcup I_0'},\ell(u)\leq \#I_0\sqcup I_0'-i}\mathrm{Fil}_{u}(\mathbf{E}_{I_0\sqcup I_0'}),\]
and thus the decreasing filtration $\{\tau^{i}(\mathbf{E}_{I_0\sqcup I_0'})\}_{i\geq 0}$ on $\mathbf{E}_{I_0\sqcup I_0'}$ restricts to a decreasing filtration $\{\tau^{i}(\mathrm{Fil}_{\Gamma_{2}}(\mathbf{E}_{I_0\sqcup I_0'}))\}_{i\geq 0}$ on $\mathrm{Fil}_{\Gamma_{2}}(\mathbf{E}_{I_0\sqcup I_0'})$, and the induced map
\[\mathrm{gr}^{i}(\mathrm{Fil}_{\Gamma_{2}}(\mathbf{E}_{I_0\sqcup I_0'}))\rightarrow \mathrm{gr}^{i}(\mathbf{E}_{I_0\sqcup I_0'})\]
is an embedding with image given by
\begin{equation}\label{equ: cup line coxeter grade}
\bigoplus_{u\in\Gamma_{2},\ell(u)=\#I_0\sqcup I_0'-i}\mathrm{gr}_{u}(\mathbf{E}_{I_0\sqcup I_0'}).
\end{equation}
Now that $u\in\Gamma_{2}$ satisfies $D_L(u)=\{1,n-1\}$ if and only if $u\in\Gamma_{2}\cap\Gamma_{I_0\sqcup I_0'}$, we observe from (\ref{equ: cup line dagger grade bound}) and (\ref{equ: cup line coxeter grade}) that
\begin{equation}\label{equ: cup line grade intersection bound}
\mathrm{gr}^{i}(\mathbf{E}_{2,\dagger})\cap \mathrm{gr}^{i}(\mathrm{Fil}_{\Gamma_{2}}(\mathbf{E}_{I_0\sqcup I_0'}))
\end{equation}
is non-empty if and only if $i=0$, in which case it equals (\ref{equ: cup line Lie quotient bound}).
We consider an arbitrary $x\in\mathbf{E}_{I_0,I_0'}^{\heartsuit}=\mathbf{E}_{2,\dagger}\cap \mathrm{Fil}_{\Gamma_{2}}(\mathbf{E}_{I_0\sqcup I_0'})$ and let $i_0$ be the maximal integer such that $x\in\tau^{i_0}(\mathbf{E}_{I_0\sqcup I_0'})$ with image $0\neq \overline{x}\in\mathrm{gr}^{i_0}(\mathbf{E}_{I_0\sqcup I_0'})$. Our discussion above on various filtration $\{\tau^{i}(-)\}_{i\geq 0}$ implies 
\[x\in \tau^{i_0}(\mathbf{E}_{2,\dagger})\cap \tau^{i_0}(\mathrm{Fil}_{\Gamma_{2}}(\mathbf{E}_{I_0\sqcup I_0'}))\subseteq \tau^{i_0}(\mathbf{E}_{I_0\sqcup I_0'})\]
and thus
\[0\neq \overline{x}\in \mathrm{gr}^{i_0}(\mathbf{E}_{2,\dagger})\cap \mathrm{gr}^{i_0}(\mathrm{Fil}_{\Gamma_{2}}(\mathbf{E}_{I_0\sqcup I_0'}))\subseteq \mathrm{gr}^{i_0}(\mathbf{E}_{I_0\sqcup I_0'})\]
which together with our discussion around (\ref{equ: cup line grade intersection bound}) forces $i_0=0$ and that $\overline{x}$ is contained in (\ref{equ: cup line Lie quotient bound}). This finishes the proof of \textbf{Step $2$}.

\textbf{Step $3$}: We prove that $\mathbf{E}_{I_0\sqcup I_0'}^{<}\cap \mathbf{E}_{I_0,I_0'}^{\heartsuit}$ is a $1$-dimensional subspace of $\mathbf{E}_{I_0\sqcup I_0'}$.\\
Now that $\kappa_{I_0,I_0'}$ is injective by Lemma~\ref{lem: separated cup injection} (upon replacing $I_1$, $J_0$ and $J_0'$ in \emph{loc.cit.} with $\Delta$, $\Delta\setminus I_0$ and $\Delta\setminus I_0'$), we deduce from (\ref{equ: cup line bound 2}) and \textbf{Step $1$} that
\begin{equation}\label{equ: cup line image lower bound}
\Dim_E \mathbf{E}_{I_0\sqcup I_0'}^{<}\cap \mathbf{E}_{I_0,I_0'}^{\heartsuit}\geq \Dim_E \mathbf{E}_{I_0}^{\heartsuit} \Dim_E \mathbf{E}_{I_0'}^{\heartsuit}=1.
\end{equation}
It remains to show that
\begin{equation}\label{equ: cup line image upper bound}
\Dim_E \mathbf{E}_{I_0\sqcup I_0'}^{<}\cap \mathbf{E}_{I_0,I_0'}^{\heartsuit}\leq 1.
\end{equation}
If $I_0<<I_0'$ (see \ref{it: disconnected 3} of Definition~\ref{def: disconnected pair}) or equivalently $j_0<j_1-1$, then we have 
\[\Gamma_{2}\cap\Gamma_{I_0\sqcup I_0'}=\Gamma_{w_{1,j_0},w_{n-1,j_1}}=\{w_{1,j_0}w_{n-1,j_1}\}\]
and thus (\ref{equ: cup line image upper bound}) follows directly from \textbf{Step $2$}.
Otherwise, $I_0\sqcup I_0'$ is an interval and we have
\[\Gamma_{2}\cap\Gamma_{I_0\sqcup I_0'}=\Gamma_{w_{1,j_0},w_{n-1,j_1}}=\{w_{1,j_0}w_{n-1,j_1},w_{n-1,j_1}w_{1,j_0}\}.\]
Note that the surjection $\mathbf{E}_{I_0\sqcup I_0'}\twoheadrightarrow \mathbf{E}_{I_0\sqcup I_0'}/\mathbf{E}_{I_0\sqcup I_0'}^{<}$ factors through 
\[\mathrm{gr}^{\emptyset}(\mathbf{E}_{I_0\sqcup I_0'})\twoheadrightarrow \mathbf{E}_{I_0\sqcup I_0'}/\mathbf{E}_{I_0\sqcup I_0'}^{<},\]
under which the image of $\mathfrak{e}_{J_{2},u}=\mathrm{gr}_{u}(\mathbf{E}_{I_0\sqcup I_0'})$ is non-zero for each $u\in \Gamma_{2}\cap\Gamma_{I_0\sqcup I_0'}$ by Theorem~\ref{thm: cup top grade}. Consequently, if we write $\mathrm{gr}^{\emptyset}(\mathbf{E}_{I_0\sqcup I_0'}^{<})$ for the image of $\mathbf{E}_{I_0\sqcup I_0'}^{<}$ under $\mathbf{E}_{I_0\sqcup I_0'}\twoheadrightarrow\mathrm{gr}^{\emptyset}(\mathbf{E}_{I_0\sqcup I_0'})$, then the following subspace
\begin{equation}\label{equ: cup line image Lie quotient}
\mathrm{gr}^{\emptyset}(\mathbf{E}_{I_0\sqcup I_0'}^{<})\cap (\bigoplus_{u\in\Gamma_{I_0\sqcup I_0'}\cap\Gamma_{2}}\mathrm{gr}_{u}(\mathbf{E}_{I_0\sqcup I_0'}))\subseteq \mathrm{gr}^{\emptyset}(\mathbf{E}_{I_0\sqcup I_0'})
\end{equation}
is $1$-dimensional. It then follows from \textbf{Step $2$} that the composition of
\[
\mathbf{E}_{I_0\sqcup I_0'}^{<}\cap\mathbf{E}_{I_0,I_0'}^{\heartsuit}\hookrightarrow\mathbf{E}_{I_0\sqcup I_0'}\twoheadrightarrow\mathrm{gr}^{\emptyset}(\mathbf{E}_{I_0\sqcup I_0'})\cong\bigoplus_{u\in\Gamma_{I_0\sqcup I_0'}}\mathrm{gr}_{u}(\mathbf{E}_{I_0\sqcup I_0'})
\]
is an embedding with image contained in the $1$-dimensional $E$-vector space (\ref{equ: cup line image Lie quotient}), from which we conclude (\ref{equ: cup line image upper bound}) and therefore finish the proof of \textbf{Step $3$}.

\textbf{Step $4$}: We prove (\ref{equ: cup of line exchange}) for each $x\in\mathbf{E}_{I_0}^{\heartsuit}$ and $y\in \mathbf{E}_{I_0'}^{\heartsuit}$.\\
Recall from Lemma~\ref{lem: separated cup injection} and Remark~\ref{rem: opposite cup injection} that both $\kappa_{I_0,I_0'}$ and $\kappa_{I_0',I_0}$ are injective.
This together with \textbf{Step $1$} and \textbf{Step $3$} implies that both inclusions (\ref{equ: cup line bound 1}) and (\ref{equ: cup line bound 2}) must be equalities between $1$-dimensional $E$-vector spaces. In particular, we have
\begin{equation}\label{equ: exchange cup line}
\kappa_{I_0,I_0'}(\mathbf{E}_{I_0}^{\heartsuit}\otimes_E\mathbf{E}_{I_0'}^{\heartsuit})=\kappa_{I_0',I_0}(\mathbf{E}_{I_0'}^{\heartsuit}\otimes_E\mathbf{E}_{I_0}^{\heartsuit})
\end{equation}
between $1$-dimensional $E$-vector spaces.
Now we consider arbitrary $0\neq x\in\mathbf{E}_{I_0}^{\heartsuit}$ and $0\neq y\in \mathbf{E}_{I_0'}^{\heartsuit}$, with both $x\cup y$ and $y\cup x$ being non-zero elements of $\mathbf{E}_{I_0\sqcup I_0'}^{<}\cap \mathbf{E}_{I_0,I_0'}^{\heartsuit}\subseteq \mathbf{E}_{I_0\sqcup I_0'}$.
We write $\overline{x}$ (resp.~$\overline{y}$) for the image of $x$ in $\mathrm{gr}^{1-r}(\mathbf{E}_{I_0})$ (resp.~the image of $y$ in $\mathrm{gr}^{1-r}(\mathbf{E}_{I_0'})$).
Now that the composition of $\mathbf{E}_{I_0}^{\heartsuit}\hookrightarrow \mathbf{E}_{I_0}\twoheadrightarrow \mathbf{E}_{I_0}/\mathbf{E}_{I_0}^{<}=\mathrm{gr}^{1-r}(\mathbf{E}_{I_0})$ and the composition of $\mathbf{E}_{I_0'}^{\heartsuit}\hookrightarrow \mathbf{E}_{I_0'}\twoheadrightarrow \mathbf{E}_{I_0'}/\mathbf{E}_{I_0'}^{<}=\mathrm{gr}^{1-r}(\mathbf{E}_{I_0'})$ are isomorphisms by \textbf{Step $1$}, we know that $x\cup y$ and $y\cup x$ are equal as elements of $\mathbf{E}_{I_0\sqcup I_0'}$ if and only if $\overline{x}\cup\overline{y}$ and $\overline{y}\cup\overline{x}$ are equal as elements of $\mathrm{gr}^{2-r}(\mathbf{E}_{I_0\sqcup I_0'})$.
Since we do have
\[\overline{x}\cup \overline{y}=\overline{y}\cup \overline{x} \in \mathrm{gr}^{2-r}(\mathbf{E}_{I_0\sqcup I_0'})\]
by Lemma~\ref{lem: grade cup commute}, we finish the proof of (\ref{equ: cup of line exchange}) for each $x\in\mathbf{E}_{I_0}^{\heartsuit}$ and $y\in \mathbf{E}_{I_0'}^{\heartsuit}$.
\end{proof}

\begin{thm}\label{thm: cup exchange}
Let $I_0,I_0'\subseteq\Delta$ with $I_0\cap I_0'=\emptyset$. Then we have 
\begin{equation}\label{equ: cup exchange}
x\cup y=y\cup x
\end{equation}
for each $x\in \mathbf{E}_{I_0}$ and $y\in \mathbf{E}_{I_0'}$.
\end{thm}
\begin{proof}
It is harmless to assume throughout the proof that $I_0\neq \emptyset\neq I_0'$.
We will consider partition $I_0=\bigsqcup_{i=1}^{s}J_{i}$ and $I_0'=\bigsqcup_{j=1}^{t}J_{j}'$ into non-empty subsets with each $J_{i}$ and $J_{j}'$ being an interval, and moreover $J_{i}<J_{i'}$ whenever $i<i'$ and $J_{j}'< J_{j'}'$ whenever $j<j'$ (see \ref{it: disconnected 2} of Definition~\ref{def: disconnected pair}). These partitions induce the following successive cup product maps
\begin{equation}\label{equ: component cup 1}
\mathbf{E}_{J_1}\otimes_E\cdots\otimes_E \mathbf{E}_{J_s}\rightarrow \mathbf{E}_{I_0}
\end{equation}
and
\begin{equation}\label{equ: component cup 2}
\mathbf{E}_{J_1'}\otimes_E\cdots\otimes_E \mathbf{E}_{J_t'}\rightarrow \mathbf{E}_{I_0'}.
\end{equation}
We prove (\ref{equ: cup exchange}) by an increasing induction on $\#I_0\sqcup I_0'$.
If $\#I_0=1=\#I_0'$, then (\ref{equ: cup exchange}) follows from Lemma~\ref{lem: grade cup commute} (with $J_0=\Delta\setminus I_0$, $J_0'=\Delta\setminus I_0'$, $k_0=2\#I_0$ and $k_1=2\#I_0'$ in \emph{loc.cit.}).
We assume in the rest of the proof that $\#I_0\sqcup I_0'>2$.

Assume for the moment that either $I_0$ or $I_0'$ is not an interval. We can choose $I_0=\bigsqcup_{i=1}^{s}J_{i}$ (resp.~$I_0'=\bigsqcup_{j=1}^{t}J_{j}'$) so that each $J_{i}$ (resp.~$J_{j}'$) is a connected component of $I_0$ (resp.~of $I_0'$).
In particular, we have $J_{i}<<J_{i'}$ whenever $i<i'$ and $J_{j}'<< J_{j'}'$ whenever $j<j'$ (see \ref{it: disconnected 3} of Definition~\ref{def: disconnected pair}), with both (\ref{equ: component cup 1}) and (\ref{equ: component cup 2}) being isomorphisms by successively applying \ref{it: disconnected isom 2} of Lemma~\ref{lem: separated cup disconnected isom}.
This forces $\mathbf{E}_{I_0}$ to be spanned by elements of the form $x_1\cup\cdots\cup x_s$ for $x_i\in\mathbf{E}_{J_{i}}$, and $\mathbf{E}_{I_0'}$ to be spanned by elements of the form $y_1\cup\cdots\cup y_t$ for $y_j\in\mathbf{E}_{J_j'}$.
Now that either $I_0$ or $I_0'$ is not an interval, we have $\#J_i\sqcup J_j'<\#I_0\sqcup I_0'$ for each $1\leq i\leq s$ and $1\leq j\leq t$, and thus
\begin{equation}\label{equ: cup exchange induction}
x_i\cup y_j=y_j\cup x_i
\end{equation}
for each $1\leq i\leq s$ and $1\leq j\leq t$ by our induction hypothesis.
By successively exchanging $x_i$ and $y_j$ for decreasing $i$ and increasing $j$, we deduce from (\ref{equ: cup exchange induction}) that
\begin{equation}\label{equ: cup exhange partition}
(x_1\cup\cdots\cup x_s)\cup(y_1\cup\cdots\cup y_t)=(y_1\cup\cdots\cup y_t)\cup(x_1\cup\cdots\cup x_s),
\end{equation}
and thus (\ref{equ: cup exchange}) holds when either $I_0$ or $I_0'$ is not an interval.

Assume from now that both $I_0$ and $I_0'$ are intervals with $\#I_0\sqcup I_0'>2$.
Note from \ref{it: cup of line 1} of Proposition~\ref{prop: cup of line} that $\mathbf{E}_{I_0}$ is spanned by $\mathbf{E}_{I_0}^{\heartsuit}$ and $x_1\cup\cdots\cup x_s$ for $x_i\in\mathbf{E}_{J_{i}}$ with $I_0=\bigsqcup_{i=1}^{s}J_{i}$ running through all partitions that satisfy $s\geq 2$. Similarly, $\mathbf{E}_{I_0'}$ is spanned by $\mathbf{E}_{I_0'}^{\heartsuit}$ and $y_1\cup\cdots\cup y_t$ for $y_j\in\mathbf{E}_{J_{j}'}$ with $I_0'=\bigsqcup_{j=1}^{t}J_{j}'$ running through all partitions that satisfy $t\geq 2$. When $s+t>2$, we have $\#J_i\sqcup J_j'<\#I_0\sqcup I_0'$ for each $1\leq i\leq s$ and $1\leq j\leq t$, and thus (\ref{equ: cup exchange induction}) again holds for each $1\leq i\leq s$ and $1\leq j\leq t$ by our induction hypothesis. In other words, (\ref{equ: cup exchange}) is true when either $x\in\mathbf{E}_{I_0}^{<}$ or $y\in\mathbf{E}_{I_0'}^{<}$ holds. This together with \ref{it: cup of line 2} of Proposition~\ref{prop: cup of line} ensures that (\ref{equ: cup exchange}) holds for each $x\in\mathbf{E}_{I_0}$ and $y\in\mathbf{E}_{I_0'}$.
The proof is thus finished.
\end{proof}
\subsection{The subspaces $\mathbf{E}_{n-1}^{\sharp}$ and $\mathbf{E}_{n-1}^{\flat}$ of $\mathbf{E}_{n-1}$}\label{subsec: sharp subspace}
We construct a subspace $\mathbf{E}_{n-1}^{\sharp}$ of $\mathbf{E}_{n-1}$ using outer involution (studied in \S \ref{subsec: outer involution}) and our normalized cup product maps. We prove that $\mathbf{E}_{n-1}^{\sharp}$ admits a well-behaved basis (still non-canonical) in Proposition~\ref{prop: involution J total basis} and then discuss the relation between $\mathbf{E}_{n-1}^{\sharp}$ and other subspaces of $\mathbf{E}_{n-1}$ including $\mathrm{Sym}^{n-1}(\mathbf{E}_{1})$ and $\mathbf{E}_{n-1}^{\flat}$ (see Corollary~\ref{cor: sym sharp decomposition} and Proposition~\ref{prop: sharp flat intersection}).

Recall from \S~\ref{subsec: outer involution} that we have defined an involution $\theta_{n-1}$ on $\mathbf{E}_{n-1}=\mathbf{E}_{\Delta}$.
Throughout this section, we fix an arbitrary choice of
\begin{equation}\label{equ: self dual generator}
x_{m}\in\mathbf{E}_{m}^{\theta_{m}=-1}\setminus(\mathbf{E}_{m}^{\theta_{m}=-1}\cap\mathbf{E}_{m}^{<})
\end{equation}
for each $m\geq 2$.
For each interval $J\subseteq\Delta$ with $\#J=m\geq 2$, under the isomorphism $\mathbf{E}_{J}\cong\mathbf{E}_{m}$ from (\ref{equ: Tits interval reduction}), we obtain an involution $\theta_{J}$ on $\mathbf{E}_{J}$ and 
\begin{equation}\label{equ: self dual generator J}
x_{J}\in\mathbf{E}_{J}^{\theta_{J}=-1}\setminus(\mathbf{E}_{J}^{<}\cap\mathbf{E}_{J}^{\theta_{J}=-1})
\end{equation}
from the involution $\theta_{m}$ on $\mathbf{E}_{m}$ and the element (\ref{equ: self dual generator}).

\begin{lem}\label{lem: outer involution cup}
Let $J,J'\subseteq\Delta$ be subintervals with $J\cap J'=\emptyset$ and $J\sqcup J'$ being another interval. Let $w$ be the longest element of $W(L_{J\sqcup J'})$, with $w(J)$ and $w(J')$ being subintervals of $J\sqcup J'$ satisfying $w(J)\cap w(J')=\emptyset$ and $w(J)\sqcup w(J')=J\sqcup J'$. Then we have
\begin{equation}\label{equ: outer involution cup}
\theta_{J\sqcup J'}(\kappa_{J,J'}(x\otimes_Ey))=\kappa_{w(J),w(J')}(x'\otimes_Ey')
\end{equation}
for each $x\in\mathbf{E}_{J}$ and $y\in\mathbf{E}_{J'}$. Here $x'$ is the image of $\theta_{J}(x)$ under the isomorphisms $\mathbf{E}_{J}\cong\mathbf{E}_{\#J}\cong\mathbf{E}_{w(J)}$, and $y'$ is the image of $\theta_{J'}(y)$ under the isomorphisms $\mathbf{E}_{J'}\cong\mathbf{E}_{\#J'}\cong\mathbf{E}_{w(J')}$, using (\ref{equ: Tits interval reduction}).
\end{lem}
\begin{proof}
Using (\ref{equ: Tits interval reduction cup}) (upon replacing $I_0$ and $I_0'$ in \emph{loc.cit.} with $J$ and $J'$ here), it is harmless to assume that $J\sqcup J'=\Delta$ with $w=w_0$ being the longest element of $W(G)$.
Following (\ref{equ: Levi cup Tits}) (upon replacing $I_0$, $I_0'$ \emph{loc.cit.} with $J$, $J'$ or $w_0(J)$, $w_0(J')$ and taking $k_0$ and $k_1$ in \emph{loc.cit.} to be $2\#J$ and $2\#J'$ with $k_0+k_1=2r$), we have a commutative diagram of the form
\begin{equation}\label{equ: outer involution cup diagram}
\xymatrix{
\mathbf{E}_{J} \ar@{=}[d] & \otimes_E & \mathbf{E}_{J'} \ar@{=}[d] \ar^{\kappa_{J,J'}}[r] & \mathbf{E}_{\Delta} \ar@{=}[d] \\
H^{k_0}(\mathrm{Tot}(\cN\cT_{J}^{\bullet,\bullet})) \ar^{\wr}[d] & \otimes_E & H^{k_1}(\mathrm{Tot}(\cN\cT_{J'}^{\bullet,\bullet})) \ar^{\cup}[r] \ar^{\wr}[d] & H^{2r}(\mathrm{Tot}(\cN\cT_{\Delta}^{\bullet,\bullet})) \ar^{\wr}[d]\\
H^{k_0}(\mathrm{Tot}(\overline{\cN\cT}_{J}^{\bullet,\bullet})) \ar^{\wr}[d] & \otimes_E & H^{k_1}(\mathrm{Tot}(\overline{\cN\cT}_{J'}^{\bullet,\bullet})) \ar^{\cup}[r] \ar^{\wr}[d] & H^{2r}(\mathrm{Tot}(\overline{\cN\cT}_{\Delta}^{\bullet,\bullet})) \ar^{\wr}[d]\\
H^{k_0}(\mathrm{Tot}(\cN\cT_{w_0(J)}^{\bullet,\bullet})) \ar@{=}[d] & \otimes_E & H^{k_1}(\mathrm{Tot}(\cN\cT_{w_0(J')}^{\bullet,\bullet})) \ar^{\cup}[r] \ar@{=}[d] & H^{2r}(\mathrm{Tot}(\cN\cT_{\Delta}^{\bullet,\bullet})) \ar@{=}[d] \\
\mathbf{E}_{w_0(J)} & \otimes_E & \mathbf{E}_{w_0(J')} \ar^{\kappa_{w_0(J),w_0(J')}}[r] & \mathbf{E}_{\Delta} 
}
\end{equation}
with the following descriptions.
\begin{itemize}
\item The double complex $\overline{\cN\cT}_{J}^{\bullet,\bullet}$ is defined by taking $M^{\bullet}_{I}$ to be $C^{\bullet}(L_{w_0(I)})$ when $\Delta\setminus w_0(J)\subseteq I\subseteq\Delta$ and zero otherwise (see (\ref{equ: general double complex normalized})), and the double complex $\overline{\cN\cT}_{J'}^{\bullet,\bullet}$ and $\overline{\cN\cT}_{\Delta}^{\bullet,\bullet}$ are defined similarly. 
\item The composition of the vertical maps in the third column is $\theta_{\Delta}$.
\item The isomorphisms from the second row to the third row of (\ref{equ: outer involution cup diagram}) are defined using scalar automorphisms of $C^{\bullet}(L_{I})$ for various $I$ (cf.~the definition of the left bottom horizontal map of (\ref{equ: Tits outer involution diagram})).
\item The isomorphisms from the third row to the fourth row of (\ref{equ: outer involution cup diagram}) are induced from the automorphism $A\mapsto w_0Aw_0$ of $G$ which sends $L_{I}$ (resp.~$L_{I'}$) to $L_{w_0(I)}$ (resp.~$L_{w_0(I')}$) for each $\Delta\setminus w_0(J)\subseteq I\subseteq \Delta$ (resp.~for each $\Delta\setminus w_0(J')\subseteq I'\subseteq \Delta$). 
\end{itemize}
Upon replacing $C^{\bullet}(L_{I})$ with $C^{\bullet}(L_{I\cap J})$ (resp.~$C^{\bullet}(L_{I}\cap H_{J})$) for each $\Delta\setminus J\subseteq I\subseteq\Delta$ in $\cN\cT_{J}^{\bullet,\bullet}$ and $\overline{\cN\cT}_{J}^{\bullet,\bullet}$, we obtain double complex $\cN\cT_{J,\flat}^{\bullet,\bullet}$ and $\overline{\cN\cT}_{J,\flat}^{\bullet,\bullet}$ (resp.~$\cN\cT_{J,\flat\flat}^{\bullet,\bullet}$ and $\overline{\cN\cT}_{J,\flat\flat}^{\bullet,\bullet}$). We similarly define $\cN\cT_{w_0(J),\flat}^{\bullet,\bullet}$ (resp.~$\cN\cT_{w_0(J),\flat\flat}^{\bullet,\bullet}$) by replacing $C^{\bullet}(L_{w_0(I)})$ in $\cN\cT_{w_0(J)}^{\bullet,\bullet}$ with $C^{\bullet}(L_{w_0(I\cap J)})$ for each $\Delta\setminus J\subseteq I\subseteq\Delta$. We thus obtain the following commutative diagram of maps between double complex
\begin{equation}\label{equ: outer involution restriction}
\xymatrix{
\cN\cT_{J}^{\bullet,\bullet} \ar^{\sim}[r] \ar[d] & \overline{\cN\cT}_{J}^{\bullet,\bullet} \ar^{\sim}[r] \ar[d] & \cN\cT_{w_0(J)}^{\bullet,\bullet} \ar[d]\\
\cN\cT_{J,\flat}^{\bullet,\bullet} \ar^{\sim}[r] \ar[d] & \overline{\cN\cT}_{J,\flat}^{\bullet,\bullet} \ar^{\sim}[r] \ar[d] & \cN\cT_{w_0(J),\flat}^{\bullet,\bullet} \ar[d]\\
\cN\cT_{J,\flat\flat}^{\bullet,\bullet} \ar^{\sim}[r] & \overline{\cN\cT}_{J,\flat\flat}^{\bullet,\bullet} \ar^{\sim}[r] & \cN\cT_{w_0(J),\flat\flat}^{\bullet,\bullet}
}
\end{equation}
with the left top (resp.~left middle, resp.~left bottom) horizontal isomorphism defined using scalar automorphisms of $C^{\bullet}(L_{I})$ (resp.~of $C^{\bullet}(L_{I\cap J})$, resp.~of $C^{\bullet}(L_{I}\cap H_{J})$) for each $\Delta\setminus J\subseteq I\subseteq\Delta$, and the right top (resp.~right middle, resp.~right bottom) horizontal isomorphism defined using $C^{\bullet}(L_{I})\buildrel\sim\over\longrightarrow C^{\bullet}(L_{w_0(I)})$ (resp.~using $C^{\bullet}(L_{I\cap J})\buildrel\sim\over\longrightarrow C^{\bullet}(L_{w_0(I\cap J)})$, resp.~using $C^{\bullet}(L_{I}\cap H_{J})\buildrel\sim\over\longrightarrow C^{\bullet}(L_{w_0(I)}\cap H_{w_0(J)})$) induced from the automorphism $A\mapsto w_0Aw_0$ of $G$. Here the left and middle vertical maps of (\ref{equ: outer involution restriction}) are induced from the restriction maps $C^{\bullet}(L_{I})\rightarrow C^{\bullet}(L_{I\cap J})\rightarrow C^{\bullet}(L_{I}\cap H_{J})$, and the right vertical maps of (\ref{equ: outer involution restriction}) are induced from the restriction maps $C^{\bullet}(L_{w_0(I)})\rightarrow C^{\bullet}(L_{w_0(I\cap J)})\rightarrow C^{\bullet}(L_{w_0(I)}\cap H_{w_0(J)})$ (for each $\Delta\setminus w_0(J)\subseteq I\subseteq\Delta$).
Writing $w_0=w_Jw^J=(w^J)^{-1}w_J$ with $w_J$ being the longest element of $W(L_{J})$, we see that the isomorphisms 
\[L_{w_0(I\cap J)}\buildrel\sim\over\longrightarrow L_{I\cap J}: A\mapsto w_0Aw_0\]
factors as the composition of isomorphisms
\[L_{w_0(I\cap J)}\buildrel\sim\over\longrightarrow L_{I\cap J}\buildrel\sim\over\longrightarrow L_{I\cap J}: A\mapsto w^{J}A(w^J)^{-1} \mapsto w_0Aw_0\]
which is functorial with respect to the choice of $\Delta\setminus J\subseteq I\subseteq\Delta$. Consequently, the right middle and right bottom horizontal maps of (\ref{equ: outer involution restriction}) factors as
\begin{equation}\label{equ: outer involution factor}
\overline{\cN\cT}_{J,\ast}^{\bullet,\bullet}\buildrel\sim\over\longrightarrow \cN\cT_{J,\ast}^{\bullet,\bullet} \buildrel\sim\over\longrightarrow \cN\cT_{w_0(J),\ast}^{\bullet,\bullet} 
\end{equation}
for $\ast\in\{\flat,\flat\flat\}$, with the LHS map induced from $A\mapsto w_JAw_J$ and the RHS map induced from $A\mapsto w^JA(w^J)^{-1}$.
Combining with \ref{it: Tits bottom shift 1} of Proposition~\ref{prop: Tits bottom shift}, we deduce from (\ref{equ: outer involution restriction}) and (\ref{equ: outer involution factor}) the following commutative diagram
\begin{equation}\label{equ: outer involution restriction coh}
\xymatrix{
H^{k_0}(\mathrm{Tot}(\cN\cT_{J}^{\bullet,\bullet})) \ar^{\sim}[r] \ar^{\wr}[d] & H^{k_0}(\mathrm{Tot}(\overline{\cN\cT}_{J}^{\bullet,\bullet})) \ar^{\sim}[rr] \ar^{\wr}[d] & & H^{k_0}(\mathrm{Tot}(\cN\cT_{w_0(J)}^{\bullet,\bullet})) \ar^{\wr}[d]\\
H^{k_0}(\mathrm{Tot}(\cN\cT_{J,\flat}^{\bullet,\bullet})) \ar^{\sim}[r] \ar^{\wr}[d] & H^{k_0}(\mathrm{Tot}(\overline{\cN\cT}_{J,\flat}^{\bullet,\bullet})) \ar^{\sim}[r] \ar^{\wr}[d] & H^{k_0}(\mathrm{Tot}(\cN\cT_{J,\flat}^{\bullet,\bullet})) \ar^{\sim}[r] \ar^{\wr}[d] & H^{k_0}(\mathrm{Tot}(\cN\cT_{w_0(J),\flat}^{\bullet,\bullet})) \ar^{\wr}[d] \\
H^{k_0}(\mathrm{Tot}(\cN\cT_{J,\flat\flat}^{\bullet,\bullet})) \ar^{\sim}[r] \ar^{\wr}[d] & H^{k_0}(\mathrm{Tot}(\overline{\cN\cT}_{J,\flat\flat}^{\bullet,\bullet})) \ar^{\sim}[r] & H^{k_0}(\mathrm{Tot}(\cN\cT_{J,\flat\flat}^{\bullet,\bullet})) \ar^{\sim}[r] \ar^{\wr}[d] & H^{k_0}(\mathrm{Tot}(\cN\cT_{w_0(J),\flat\flat}^{\bullet,\bullet})) \ar^{\wr}[d] \\
\mathbf{E}_{\#J} \ar^{\theta_{\#J}}[rr] & & \mathbf{E}_{\#J} \ar@{=}[r] & \mathbf{E}_{\#J}
}
\end{equation}
We conclude from (\ref{equ: outer involution restriction coh}) that the composition of the left vertical maps of (\ref{equ: outer involution cup diagram}) equals the composition of
\[\mathbf{E}_{J}\cong\mathbf{E}_{\#J}\buildrel\theta_{\#J}\over\longrightarrow\mathbf{E}_{\#J}\cong \mathbf{E}_{w_0(J)}.\]
A parallel argument shows that the composition of the middle vertical maps of (\ref{equ: outer involution cup diagram}) equals the composition of
\[\mathbf{E}_{J'}\cong\mathbf{E}_{\#J'}\buildrel\theta_{\#J'}\over\longrightarrow\mathbf{E}_{\#J'}\cong \mathbf{E}_{w_0(J')}.\]
We finally finish the proof of (\ref{equ: outer involution cup}) using the commutative diagram (\ref{equ: outer involution cup diagram}) together with the above description of the composition of the vertical maps in each column of it.
\end{proof}

In the rest of this section, for each $\al\in\Phi^+\setminus\Delta$, we choose $x_{\al}$ as in (\ref{equ: cup generator choice}) to be $x_{I_{\al}}$ as in (\ref{equ: self dual generator J}) with $J=I_{\al}$ in \emph{loc.cit.}.
We again define $x_{S,I}$ for each $I_0\subseteq\Delta$, $S\in\cS_{I_0}$ and $I\subseteq S\cap\Delta$ as in (\ref{equ: explicit cup generator}), following the choices of $x_{\al}$ as above for each $\al\in\Phi^+\setminus\Delta$.
For each interval $J\subseteq\Delta$, we define
\begin{equation}\label{equ: sharp subspace}
\mathbf{E}_{J}^{\sharp}\defeq \sum_{J'}\kappa_{J',J\setminus J'}(\mathbf{E}_{J'}^{\theta_{J'}=-1}\otimes_E\mathbf{E}_{J\setminus J'}) \subseteq \mathbf{E}_{J}
\end{equation}
where $J'$ runs through subintervals of $J$ satisfying $\#J'\geq 2$.
We also write 
\begin{equation}\label{equ: sharp subspace full}
\mathbf{E}_{n-1}^{\sharp}\defeq\mathbf{E}_{\Delta}^{\sharp}\subseteq\mathbf{E}_{\Delta}=\mathbf{E}_{n-1}.
\end{equation}
\begin{lem}\label{lem: sharp subspace cup}
Let $J\subseteq\Delta$ be an interval with $\#J\geq 3$. Then we have
\begin{equation}\label{equ: sharp subspace cup}
\mathbf{E}_{J}^{<}\cap \mathbf{E}_{J}^{\theta_{J}=-1}\subseteq \sum_{J'}\kappa_{J',J\setminus J'}(\mathbf{E}_{J'}^{\theta_{J'}=-1}\otimes_E\mathbf{E}_{J\setminus J'})
\end{equation}
where $J'\subsetneq J$ runs through proper subintervals of $J$ with $\#J'\geq 2$.
\end{lem}
\begin{proof}
Recall from (\ref{equ: coxeter involution}) and (\ref{equ: partition involution}) that $\theta_{J}$ corresponds to an involution on $\cS_{J}$ and an involution on $\cP(J)$ (both still denoted by $\theta_{J}$ for simplicity).
Since we have $\theta_{I_{\al}}(x_{\al}^{\ast})=-x_{\al}^{\ast}$ (with $x_{\al}^{\ast}=x_{\al}$ when $\al\in\Phi^+\setminus\Delta$ and $x_{\al}^{\ast}\in\{\plog_{\al},\val_{\al}\}$ when $\al\in\Delta$) for each $\al\in\Phi^+$, we deduce from Lemma~\ref{lem: outer involution cup} and an increasing induction on $\#S$ that
\[\theta_{J}(x_{S,I})=(-1)^{\#S}x_{\theta_{J}(S),\theta_{J}(I)}\]
for each $S\in\cS_{J}$ and $I\subseteq S\cap\Delta$. 
In particular, the LHS of (\ref{equ: sharp subspace cup}) is spanned by elements of the form
\begin{equation}\label{equ: involution J basis}
x_{S,I}-(-1)^{\#S}x_{\theta_{J}(S),\theta_{J}(I)}
\end{equation}
for varying $S\in\cS_{J}$ with $\#S\geq 2$ and $I\subseteq S\cap\Delta$.
If $S\neq J$, then for an arbitrarily chosen $\beta\in S\setminus J$ with $S'\defeq S\setminus\{\beta\}$, $I'\defeq I\setminus\{\beta\}$, $J'\defeq I_{\beta}$ and $J''\defeq \theta_{J}(J')=I_{\theta_{J}(\beta)}$ for short, we have 
\[x_{S,I}=x_{\beta}\cup x_{S',I'}\in \kappa_{J',J\setminus J'}(\mathbf{E}_{J'}^{\theta_{J'}=-1}\otimes_E\mathbf{E}_{J\setminus J'})\] 
and
\[x_{\theta_{J}(S),\theta_J(I)}=x_{\theta_J(\beta)}\cup x_{\theta_{J}(S'),\theta_{J}(I')}\in \kappa_{J'',J\setminus J''}(\mathbf{E}_{J''}^{\theta_{J''}=-1}\otimes_E\mathbf{E}_{J\setminus J''}),\]
which altogether implies that (\ref{equ: involution J basis}) is contained in the RHS of (\ref{equ: sharp subspace cup}).
If $S=J$, then (\ref{equ: involution J basis}) has the form $x_{J,I}-x_{J,\theta_{J}(I)}$ with $\#I=\#\theta_{J}(I)$.
Given $I',I''\subseteq J$, we say that $I'$ and $I''$ are \emph{adjacent} if there exists $\beta',\beta''\in J$ such that $\beta'+\beta''\in\Phi^+$ with $I\setminus I'=\{\beta'\}$ and $I'\setminus I=\{\beta''\}$.
If $I',I''\subseteq J$ are adjacent with $\beta',\beta''\in J$ as above and $M\defeq\{\beta',\beta''\}$, then we have
\[x_{J,I'}-x_{J,I''}=(\val_{\beta'}\cup\plog_{\beta''}-\plog_{\beta'}\cup\val_{\beta''})\cup x_{J\setminus M,I\cap I'}\in\kappa_{M,J\setminus M}(\mathbf{E}_{M}^{\theta_{M}=-1}\otimes_E\mathbf{E}_{J\setminus M}).\]
In other words, $x_{J,I'}-x_{J,I''}$ is contained in the RHS of (\ref{equ: sharp subspace cup}) for each adjacent pair $I',I''\subseteq J$. Now that there exists $I_{t'}$ for $0\leq t'\leq t$ such that $I_0=I$, $I_{t}=\theta_{J}(I)$ with $I_{t'}$ and $I_{t'-1}$ being adjacent for each $1\leq t'\leq t$, we deduce that $x_{J,I}-x_{J,\theta_{J}(I)}$ is contained in the RHS of (\ref{equ: sharp subspace cup}) for each $I\subseteq J$.
This finishes the proof of (\ref{equ: sharp subspace cup}).
\end{proof}

For each interval $J=[j,j']$ and $0\leq i\leq \#J=j-j'+1$, we write $J_{i}\defeq [j'-i+1,j']\subseteq J$ for short with $\#J_i=i$.
\begin{prop}\label{prop: involution J total basis}
Let $J$ be an interval with $\#J\geq 2$. Then the subspace $\mathbf{E}_{J}^{\sharp}\subseteq\mathbf{E}_{J}$ has codimension $\#J+1$ and admits a basis of the form
\begin{equation}\label{equ: involution J total basis}
\{x_{S,I}\mid S\in\cS_{J},S\neq J,I\subseteq S\cap J\}\sqcup\bigsqcup_{i=0}^{\#J}\{x_{J,I}-x_{J,J_{i}}\mid J_{i}\neq I\subseteq J, \#I=i\}
\end{equation}
\end{prop}
\begin{proof}
We prove by an increasing induction on $\#J\geq 2$.
If $\#J=2$, namely $J=\{\beta',\beta''\}$ with $\beta'+\beta''\in\Phi^+$, then we have
\[\mathbf{E}_{J}^{\sharp}=Ex_{J}\oplus E(\val_{\beta'}\cup\plog_{\beta''}-\plog_{\beta'}\cup\val_{\beta''})\]
and the result is clear. Assume from now that $\#J\geq 3$. The very definition of $\mathbf{E}_{J}^{\sharp}$ from (\ref{equ: sharp subspace}) together with (\ref{equ: self dual generator J}) implies that 
\begin{multline*}
\mathbf{E}_{J}^{\sharp}=\mathbf{E}_{J}^{\theta_{J}=-1}+\sum_{J'\subsetneq J}\kappa_{J',J\setminus J'}(\mathbf{E}_{J'}^{\sharp}\otimes_E\mathbf{E}_{J\setminus J'})\\
=(Ex_{J}\oplus (\mathbf{E}_{J}^{<}\cap \mathbf{E}_{J}^{\theta_{J}=-1}))+\sum_{J'\subsetneq J}\kappa_{J',J\setminus J'}(\mathbf{E}_{J'}^{\sharp}\otimes_E\mathbf{E}_{J\setminus J'})
\end{multline*}
with $J'\subsetneq J$ running through proper subintervals with $\#J'\geq 2$, which together with (\ref{equ: sharp subspace cup}) gives
\begin{equation}\label{equ: involution J total basis induction}
\mathbf{E}_{J}^{\sharp}=Ex_{J}\oplus \sum_{J'\subsetneq J}\kappa_{J',J\setminus J'}(\mathbf{E}_{J'}^{\sharp}\otimes_E\mathbf{E}_{J\setminus J'}).
\end{equation}
For each $J'\subsetneq J$ (with $\#J'<\#J$), our induction hypothesis implies that $\mathbf{E}_{J'}^{\sharp}$ admits a basis of the form
\[
\{x_{S,I}\mid S\in\cS_{J'},S\neq J',I\subseteq S\cap J'\}\sqcup\bigsqcup_{i=0}^{\#J'}\{x_{J',I}-x_{J',J_{i}'}\mid J_{i}'\neq I\subseteq J', \#I=i\},
\]
and thus $\kappa_{J',J\setminus J'}(\mathbf{E}_{J'}^{\sharp}\otimes_E\mathbf{E}_{J\setminus J'})$ is contained in the span of (\ref{equ: involution J total basis}). This together with (\ref{equ: involution J total basis induction}) ensures that $\mathbf{E}_{J}^{\sharp}$ is contained in the span of (\ref{equ: involution J total basis}). 
By treating the LHS half of (\ref{equ: involution J total basis}) and the RHS half of (\ref{equ: involution J total basis}) separately following the argument of Lemma~\ref{lem: sharp subspace cup} (with $x_{\al_{J},\emptyset}=x_{J}\in\mathbf{E}_{J}^{\sharp}$ by (\ref{equ: involution J total basis induction})), we know that each element of (\ref{equ: involution J total basis}) actually belongs to $\mathbf{E}_{J}^{\sharp}$. Up to this stage, we have shown that $\mathbf{E}_{J}^{\sharp}$ equals the span of (\ref{equ: involution J total basis}). Finally, it follows from Proposition~\ref{prop: total cup basis} that (\ref{equ: involution J total basis}) is linearly independent, making it a basis of $\mathbf{E}_{J}^{\sharp}$. The proof is thus finished.
\end{proof}

Recall that $r=\#\Delta$.
Let $J\subseteq\Delta$ be an interval.
We equip $\mathbf{E}_{J}^{\sharp}$ with the canonical filtration $\{\mathrm{Fil}^{-\ell}(-)\}_{0\leq\ell\leq r}$ and the filtration $\{\tau^{i}(-)\}_{i\geq 0}$ induced from that of $\mathbf{E}_{J}$.
The following results follow directly from Proposition~\ref{prop: involution J total basis}.
\begin{cor}\label{cor: cup basis involution graded}
Let $J\subseteq\Delta$ be an interval. We have the following results.
\begin{enumerate}[label=(\roman*)]
\item \label{it: cup basis involution graded ell} The embedding 
\begin{equation}\label{equ: grade ell embedding}
\mathrm{gr}^{-\ell}(\mathbf{E}_{J}^{\sharp})\hookrightarrow \mathrm{gr}^{-\ell}(\mathbf{E}_{J})
\end{equation} is an isomorphism for each $r-\#J+1\leq \ell\leq r$. 
When $\ell=r-\#J$, the image of the embedding (\ref{equ: grade ell embedding}) has codimension $\#J+1$ and admits a basis of the form
\begin{equation}\label{equ: cup basis involution 0}
\bigsqcup_{i=0}^{\#J}\{\overline{x}_{J,I}-\overline{x}_{J,J_{i}}\}_{J_{i}\neq I\subseteq J,\#I=i}.
\end{equation}
\item \label{it: cup basis involution graded i} For each $0\leq i\leq \#J$, the image of the embedding 
\begin{equation}\label{equ: grade i embedding}
\mathrm{gr}^{i}(\mathbf{E}_{J}^{\sharp})\hookrightarrow \mathrm{gr}^{i}(\mathbf{E}_{J})
\end{equation} has codimension $1$ and admits a basis induced from
\begin{equation}\label{equ: cup basis involution graded i}
\{x_{S,I}\mid S\in\cS_{J},S\neq J,I\subseteq S\cap J,\#I=i\}\sqcup\{x_{J,I}-x_{J,J_{i}}\}_{J_{i}\neq I\subseteq J,\#I=i}\subseteq \tau^{i}(\mathbf{E}_{J}^{\sharp}).
\end{equation}
\end{enumerate}
\end{cor}

Recall that $r=\#\Delta=n-1$ and that $\mathrm{Sym}^{r}(\mathbf{E}_{1})\subseteq\mathbf{E}_{1}^{\otimes r}$ is the $E$-subspace spanned by
\begin{equation}\label{equ: sym product element}
x_{\Sigma}\defeq \sum_{\sigma\in \mathfrak{S}_{r}}x_{\sigma(1)}\otimes_E\cdots\otimes_Ex_{\sigma(r)}\in \mathbf{E}_{1}^{\otimes r}
\end{equation}
for each un-ordered tuple (multi-set) $\Sigma=\{x_1,\dots,x_r\}\subseteq \mathbf{E}_{1}$.
Recall from (\ref{equ: Tits interval reduction}) that we have an isomorphism $\mathbf{E}_{1}\cong\mathbf{E}_{\{j\}}$ for each $j\in\Delta$. 
These isomorphisms induce a natural embedding
\begin{equation}\label{equ: sym product of E1}
\mathrm{Sym}^{r}(\mathbf{E}_{1})\hookrightarrow \mathbf{E}_{r}
\end{equation}
that sends $x_{\Sigma}$ to
\begin{equation}\label{equ: sym product cup}
y_{\Sigma}\defeq \sum_{\sigma\in \mathfrak{S}_{r}}x_{\sigma(1)}\cup\cdots\cup x_{\sigma(r)}\in \mathbf{E}_{r}
\end{equation}
where we understand $x_{\sigma(j)}$ at $j$-th position as an element of $\mathbf{E}_{\{j\}}$ via the isomorphism $\mathbf{E}_{1}\cong\mathbf{E}_{\{j\}}$, for each $1\leq j\leq r$ and $\sigma\in\mathfrak{S}_{r}$.
Now that $\mathbf{E}_{1}$ admits $\{\val,\plog\}$ as a basis, $\mathrm{Sym}^{r}(\mathbf{E}_{1})$ admits a basis $x_{\Sigma_i}$ indexed by the multi-set $\Sigma_{i}$ consisting of $i$ copies of $\val$ and $r-i$ copies of $\plog$, for each $0\leq i\leq r$. Note that $x_{\Sigma_i}$ is mapped to
\[y_{\Sigma_i}=\sum_{I\subseteq\Delta,\#I=i}x_{\Delta,I}\]
under (\ref{equ: sym product of E1}) for each $0\leq i\leq r$ with $\{y_{\Sigma_i}\mid 0\leq i\leq r\}$ being linearly independent. In particular, (\ref{equ: sym product of E1}) is an embedding, and we abuse below $\mathrm{Sym}^{r}(\mathbf{E}_{1})$ as an $E$-subspace of $\mathbf{E}_{r}$ via (\ref{equ: sym product of E1}). 
The following result is a direct consequence of Proposition~\ref{prop: involution J total basis}.
\begin{cor}\label{cor: sym sharp decomposition}
We have
\begin{equation}\label{equ: sym sharp decomposition}
\mathbf{E}_{n-1}=\mathbf{E}_{n-1}^{\sharp}\oplus\mathrm{Sym}^{n-1}(\mathbf{E}_{1}).
\end{equation}
In particular, we have a unique surjection
\begin{equation}\label{equ: sym auto surjection}
\mathbf{E}_{n-1}\twoheadrightarrow \mathbf{E}_{n-1}/\mathbf{E}_{n-1}^{\sharp}\cong\mathrm{Sym}^{n-1}(\mathbf{E}_{1})
\end{equation}
which admits $\mathrm{Sym}^{n-1}(\mathbf{E}_{1})\hookrightarrow \mathbf{E}_{n-1}$ as a section.
\end{cor}

Note that the filtration
\[\tau^{1}(\mathbf{E}_{1})\subsetneq \tau^{0}(\mathbf{E}_{1})=\mathbf{E}_{1}\]
induces a $\tau^{i}(-)$ filtration on $\mathrm{Sym}^{n-1}(\mathbf{E}_{1})$ with $\mathrm{gr}^{i}(\mathrm{Sym}^{n-1}(\mathbf{E}_{1}))$ being $1$-dimensional for each $0\leq i\leq n-1$ and being zero otherwise. In particular, we have \[\tau^{n-1}(\mathrm{Sym}^{n-1}(\mathbf{E}_{1}))=\mathrm{Sym}^{n-1}(\tau^{1}(\mathbf{E}_{1}))\buildrel\sim\over\longrightarrow\tau^{n-1}(\mathbf{E}_{n-1})\subseteq\mathbf{E}_{n-1}\]
and
\begin{equation}\label{equ: sym grade 0}
\mathrm{gr}^{0}(\mathrm{Sym}^{n-1}(\mathbf{E}_{1}))=\mathrm{Sym}^{n-1}(\mathrm{gr}^{0}(\mathbf{E}_{1})).
\end{equation}
It is an easy consequence of Proposition~\ref{prop: involution J total basis} that (\ref{equ: sym sharp decomposition}) induces an isomorphism
\[\tau^{i}(\mathbf{E}_{n-1})=\tau^{i}(\mathbf{E}_{n-1}^{\sharp})\oplus\tau^{i}(\mathrm{Sym}^{n-1}(\mathbf{E}_{1}))\]
as well as an isomorphism
\begin{equation}\label{equ: sym sharp decomposition grade i}
\mathrm{gr}^{i}(\mathbf{E}_{n-1})=\mathrm{gr}^{i}(\mathbf{E}_{n-1}^{\sharp})\oplus\mathrm{gr}^{i}(\mathrm{Sym}^{n-1}(\mathbf{E}_{1})).
\end{equation}
for each $0\leq i\leq n-1$.

We write $w_{j,1}\defeq s_{j}\cdots s_1$ for each $0\leq j\leq n-1$ with convention $w_{0,1}=1$.
We set
\begin{equation}\label{equ: flat subspace}
\mathbf{E}_{n-1}^{\flat}\defeq \mathrm{Fil}_{w_{n-1,1}}(\mathbf{E}_{\Delta})\subseteq\mathbf{E}_{\Delta}=\mathbf{E}_{n-1}.
\end{equation}
Following \ref{it: tau and coxeter 2} of Lemma~\ref{lem: tau and coxeter}, we know that the decreasing filtration $\{\tau^{i}(\mathbf{E}_{n-1})_{i\geq 0}\}$ on $\mathbf{E}_{n-1}$ restricts to a decreasing filtration $\{\tau^{i}(\mathbf{E}_{n-1}^{\flat})_{i\geq 0}\}$ on $\mathbf{E}_{n-1}^{\flat}$ with 
\begin{equation}\label{equ: flat subspace grade i}
\mathrm{gr}^{i}(\mathbf{E}_{n-1}^{\flat})=\mathrm{gr}^{[n-i,n-1]}(\mathbf{E}_{n-1}^{\flat})=\mathrm{gr}_{w_{n-1-i,1}}(\mathbf{E}_{n-1})
\end{equation}
being $1$-dimensional for each $0\leq i\leq n-1$.
We write $\mathbf{E}_{n-1}^{<}\defeq\mathbf{E}_{\Delta}^{<}\subseteq \mathbf{E}_{n-1}=\mathbf{E}_{\Delta}$ (see (\ref{equ: decomposable subspace})). 
Note that $\mathbf{E}_{n-1}^{<}$ inherits a decreasing filtration $\{\tau^{i}(\mathbf{E}_{n-1}^{<})\}_{i\geq 0}$ with $\mathrm{gr}^{i}(\mathbf{E}_{n-1}^{<})$ admitting a basis induced from (cf.~the proof of Lemma~\ref{lem: primitive generator})
\[\{x_{S,I}\mid S\in\cS_{\Delta}, \#S\geq 2, I\subseteq S\cap\Delta, \#I=i\}\subseteq \tau^{i}(\mathbf{E}_{n-1}^{<})\]
for each $0\leq i\leq n-1$.
We also define 
\[\mathbf{E}_{n-1}^{\plog}\defeq Ex_{\Delta,\emptyset}\subseteq\mathrm{Sym}^{n-1}(\mathbf{E}_{1})\subseteq\mathbf{E}_{n-1},\]
which depends on our choice of $\plog$ (in the definition of $x_{\Delta,\emptyset}$, see (\ref{equ: explicit cup generator})) and satisfies $\mathbf{E}_{n-1}^{\plog}\cap \mathbf{E}_{n-1}^{\sharp}=0$ by Corollary~\ref{cor: sym sharp decomposition}. Finally we define \[\mathbf{E}_{n-1}^{\diamond}\defeq \mathbf{E}_{n-1}^{\flat}\cap (\mathbf{E}_{n-1}^{\sharp}\oplus\mathbf{E}_{n-1}^{\plog}).\]
We still do not know how to prove the following statement and therefore leave it as a hypothesis.
\begin{hypo}\label{hypo: log projection}
Let $m\geq 2$. The composition of
\begin{equation}\label{equ: log projection m}
\mathbf{E}_{m}^{\flat}\hookrightarrow\mathbf{E}_{m}\twoheadrightarrow\mathrm{Sym}^{m}(\mathrm{gr}^{0}(\mathbf{E}_{1}))
\end{equation}
is non-zero.
\end{hypo}

It follows from (\ref{equ: sym grade 0}) that the composition of (\ref{equ: log projection m}) factors through
\[\mathbf{E}_{m}^{\flat}\twoheadrightarrow\mathrm{gr}^{0}(\mathbf{E}_{m}^{\flat})\hookrightarrow\mathrm{gr}^{0}(\mathbf{E}_{m})\twoheadrightarrow \mathrm{gr}^{0}(\mathbf{E}_{m})/\mathrm{gr}^{0}(\mathbf{E}_{m}^{\sharp})\buildrel\sim\over\longrightarrow\mathrm{gr}^{0}(\mathrm{Sym}^{m}(\mathbf{E}_{1})).\]
Now that $\mathrm{gr}^{0}(\mathbf{E}_{m}^{\sharp})$ (resp.~$\mathrm{gr}^{0}(\mathbf{E}_{m}^{\flat})$) is an $E$-subspace of $\mathrm{gr}^{0}(\mathbf{E}_{m})$ of codimension $1$ by \ref{it: cup basis involution graded i} of Corollary~\ref{cor: cup basis involution graded} (resp.~of dimension $1$ by (\ref{equ: flat subspace grade i})), we see that Hypothesis~\ref{hypo: log projection} holds if and only if
\begin{equation}\label{equ: sharp flat decomposition grade 0}
\mathrm{gr}^{0}(\mathbf{E}_{m}^{\flat})\oplus\mathrm{gr}^{0}(\mathbf{E}_{m}^{\sharp})=\mathrm{gr}^{0}(\mathbf{E}_{m}).
\end{equation}

\begin{prop}\label{prop: sharp flat intersection}
Assume that Hypothesis~\ref{hypo: log projection} holds for each $2\leq m\leq n-1$. We have the following results.
\begin{enumerate}[label=(\roman*)]
\item \label{it: sharp flat intersection 1} We have
\begin{equation}\label{equ: sharp flat intersection 1}
\mathbf{E}_{n-1}^{\flat}\oplus\mathbf{E}_{n-1}^{\sharp}=\mathbf{E}_{n-1}.
\end{equation}
\item \label{it: sharp flat intersection 2} The $E$-subspace $\mathbf{E}_{n-1}^{\diamond}\subseteq \mathbf{E}_{n-1}$ is $1$-dimensional, and we have
    \begin{equation}\label{equ: sharp flat intersection 2}
    \mathbf{E}_{n-1}^{\diamond}\oplus \mathbf{E}_{n-1}^{<}=\mathbf{E}_{n-1}.
    \end{equation}
\end{enumerate}
\end{prop}
\begin{proof}
We continue to write $r=\#\Delta$.
We prove \ref{it: sharp flat intersection 1}.\\
Now that BHS of (\ref{equ: sharp flat intersection 1}) share the same dimension (using Proposition~\ref{prop: involution J total basis} and $\Dim_E\mathbf{E}_{n-1}^{\flat}=n$ by its definition), it suffices to show that
\begin{equation}\label{equ: sharp flat intersection grade i}
\mathrm{gr}^{i}(\mathbf{E}_{n-1}^{\flat})\oplus \mathrm{gr}^{i}(\mathbf{E}_{n-1}^{\sharp})=\mathrm{gr}^{i}(\mathbf{E}_{n-1})
\end{equation}
for each $0\leq i\leq n-1$. 
It follows from \ref{it: cup basis involution graded i} of Corollary~\ref{cor: cup basis involution graded} that
\begin{equation}\label{equ: sharp flat partial sm grade i}
\mathrm{gr}^{i}(\mathbf{E}_{n-1}^{\sharp})\cap\mathrm{gr}^{[n-i,n-1]}(\mathbf{E}_{n-1})
\end{equation}
admits a basis induced from
\[\{x_{S,[n-i,n-1]}\mid S\in\cS_{\Delta},S\neq \Delta,[n-i,n-1]\subseteq S\cap\Delta\}\subseteq \tau^{[n-i,n-1]}(\mathbf{E}_{n-1})\]
and thus equals
\begin{equation}\label{equ: sharp flat partial sm grade i cup 1}
\kappa_{[1,n-1-i],[n-i,n-1]}(\mathrm{gr}^{0}(\mathbf{E}_{[1,n-1-i]}^{\sharp})\otimes_E\tau^{[n-i,n-1]}(\mathbf{E}_{[n-i,n-1]})).
\end{equation}
It follows from (\ref{equ: flat subspace grade i}), Lemma~\ref{lem: cup with sm} and Lemma~\ref{lem: tau and coxeter} that
\begin{multline}\label{equ: sharp flat partial sm grade i cup 2}
\mathrm{gr}^{i}(\mathbf{E}_{n-1}^{\flat})=\mathrm{gr}^{[n-i,n-1]}(\mathbf{E}_{n-1}^{\flat})=\mathrm{gr}_{w_{n-1-i,1}}(\mathbf{E}_{n-1}^{\flat})\\
=\kappa_{[1,n-1-i],[n-i,n-1]}(\mathrm{gr}_{w_{n-1-i,1}}(\mathbf{E}_{[1,n-1-i]})\otimes_E\mathrm{Fil}_{1}(\mathbf{E}_{[n-i,n-1]}))\\
=\kappa_{[1,n-1-i],[n-i,n-1]}(\mathrm{gr}^{0}(\mathbf{E}_{[1,n-1-i]}^{\flat})\otimes_E\tau^{[n-i,n-1]}(\mathbf{E}_{[n-i,n-1]})).
\end{multline}
with $\mathbf{E}_{[1,n-1-i]}^{\flat}\defeq\mathrm{Fil}_{w_{n-1-i,1}}(\mathbf{E}_{[1,n-1-i]})$.
By our assumption and the discussion below Hypothesis~\ref{hypo: log projection}, we know that (\ref{equ: sharp flat decomposition grade 0}) holds for each $2\leq m\leq n-1$.
Taking $m=n-1-i$ and using the isomorphism $\mathbf{E}_{[1,n-1-i]}\cong\mathbf{E}_{n-1-i}$ (see (\ref{equ: Tits interval reduction})), we deduce from (\ref{equ: sharp flat decomposition grade 0}) that
\[\mathrm{gr}^{0}(\mathbf{E}_{[1,n-1-i]}^{\flat})\oplus \mathrm{gr}^{0}(\mathbf{E}_{[1,n-1-i]}^{\sharp})=\mathrm{gr}^{0}(\mathbf{E}_{[1,n-1-i]}),\]
which together with (\ref{equ: sharp flat partial sm grade i cup 2}) as well as the equality between (\ref{equ: sharp flat partial sm grade i}) and (\ref{equ: sharp flat partial sm grade i cup 1}) gives
\[\mathrm{gr}^{i}(\mathbf{E}_{n-1}^{\flat})\not\subseteq \mathrm{gr}^{i}(\mathbf{E}_{n-1}^{\sharp}).\]
Now that $\mathrm{gr}^{i}(\mathbf{E}_{n-1}^{\sharp})$ (resp.~$\mathrm{gr}^{i}(\mathbf{E}_{n-1}^{\flat})$) is an $E$-subspace of $\mathrm{gr}^{i}(\mathbf{E}_{n-1})$ of codimension $1$ by \ref{it: cup basis involution graded i} of Corollary~\ref{cor: cup basis involution graded} (resp.~of dimension $1$ by (\ref{equ: flat subspace grade i})), we finish the proof of (\ref{equ: sharp flat intersection grade i}).

We prove \ref{it: sharp flat intersection 2}.\\
It follows from \ref{it: sharp flat intersection 1} that
\begin{multline}\label{equ: E diamond dim}
\Dim_E\mathbf{E}_{n-1}^{\diamond}=\Dim_E\mathbf{E}_{n-1}^{\flat}+\Dim_E(\mathbf{E}_{n-1}^{\sharp}\oplus\mathbf{E}_{n-1}^{\plog})-\Dim_E(\mathbf{E}_{n-1}^{\flat}+(\mathbf{E}_{n-1}^{\sharp}+\mathbf{E}_{n-1}^{\plog}))\\
=\Dim_E\mathbf{E}_{n-1}^{\flat}+\Dim_E\mathbf{E}_{n-1}^{\sharp}+1-\Dim_E\mathbf{E}_{n-1}=1.
\end{multline}
It follows from \ref{it: cup basis involution graded i} of Corollary~\ref{cor: cup basis involution graded} that
\[\tau^{1}(\mathbf{E}_{n-1})\cap(\mathbf{E}_{n-1}^{\sharp}\oplus\mathbf{E}_{n-1}^{\plog})=\tau^{1}(\mathbf{E}_{n-1})\cap\mathbf{E}_{n-1}^{\sharp}=\tau^{1}(\mathbf{E}_{n-1}^{\sharp})\]
and thus
\[\tau^{1}(\mathbf{E}_{n-1}^{\flat})\cap(\mathbf{E}_{n-1}^{\sharp}\oplus\mathbf{E}_{n-1}^{\plog})\subseteq \mathbf{E}_{n-1}^{\flat}\cap\mathbf{E}_{n-1}^{\sharp}=0\]
by \ref{it: sharp flat intersection 1}. This together with (\ref{equ: E diamond dim}) forces the composition of
\begin{equation}\label{equ: diamond composition}
\mathbf{E}_{n-1}^{\diamond}\hookrightarrow \mathbf{E}_{n-1}^{\flat}\twoheadrightarrow \mathrm{gr}^{0}(\mathbf{E}_{n-1}^{\flat})
\end{equation}
to be an isomorphism between $1$-dimensional $E$-vector spaces.
Now that we clearly have $\tau^{1}(\mathbf{E}_{n-1})\subseteq \mathbf{E}_{n-1}^{<}$ (see Lemma~\ref{lem: cup with sm} and Lemma~\ref{lem: tau and coxeter}) and
\[\mathrm{gr}^{0}(\mathbf{E}_{n-1})=\mathrm{gr}^{0}(\mathbf{E}_{n-1}^{<})\oplus \mathrm{gr}^{0}(\mathbf{E}_{n-1}^{\flat})\]
by Theorem~\ref{thm: cup top grade}, we conclude \ref{it: sharp flat intersection 2} from the fact that the composition of (\ref{equ: diamond composition}) is an isomorphism.
\end{proof} 

Comparing Corollary~\ref{cor: cup basis involution graded} and Proposition~\ref{prop: J truncation total basis} (with the map $\vartheta$ from \emph{loc.cit.}), we see that
\[\mathrm{im}(\vartheta_{\ell})=\mathrm{gr}^{-\ell}(\mathbf{E}_{n-1}^{\sharp})\subseteq \mathrm{gr}^{-\ell}(\mathbf{E}_{n-1})\]
for each $0\leq \ell\leq n-1$. These equalities lead us to the following questions.
\begin{ques}\label{ques: equal subspace}
Does the following equality 
\begin{equation}\label{equ: equal subspace}
\mathrm{im}(\vartheta)=\mathbf{E}_{n-1}^{\sharp}
\end{equation}
hold?
\end{ques}
Note that the validity of (\ref{equ: equal subspace}) would (together with Proposition~\ref{prop: J truncation total basis}) in particular force the map $\vartheta$ to be strict with respect to the filtration $\{\mathrm{Fil}^{-\ell}(-)\}_{\ell}$ on both the source and the target of $\vartheta$.
One evidence for the validity of (\ref{equ: equal subspace}) will be given in Corollary~\ref{cor: det deform sym} below, which roughly says that the Galois counterpart of $\mathrm{im}(\vartheta)$ and $\mathbf{E}_{n-1}^{\sharp}$ are equal.
\subsection{Tensor product on the automorphic side}\label{subsec: auto tensor product}
We define a certain automorphic tensor product map (see (\ref{equ: auto stable cup}) and (\ref{equ: multi cup map}) below) that will be useful in \S \ref{subsec: filtered phi N} and \S \ref{subsec: map candidate}.

For each $n\geq 2$ and $I\subseteq \Delta=[1,n-1]$, we write $\mathbf{E}_{n-1,I}\defeq\mathbf{E}_{I}$ just to memorize $n$ in the notation. When $I$ is an interval, recall from (\ref{equ: Tits interval reduction}) that we have an isomorphism
\begin{equation}\label{equ: auto interval to Levi block}
\mathbf{E}_{n-1,I}\cong\mathbf{E}_{\#I}.
\end{equation}
Given $m,n\geq 0$, we define a map
\begin{equation}\label{equ: auto stable cup}
\mathbf{E}_{m}\otimes_E\mathbf{E}_{n}\rightarrow \mathbf{E}_{m+n}
\end{equation}
as the composition of the following maps
\begin{equation}\label{equ: auto stable cup seq}
\mathbf{E}_{m}\otimes_E\mathbf{E}_{n}\cong \mathbf{E}_{m+n,[1,m]}\otimes_E\mathbf{E}_{m+n,[m+1,m+n]}\buildrel\kappa_{[1,m],[m+1,m+n]}\over\longrightarrow \mathbf{E}_{m+n,[1,m+n]}\cong \mathbf{E}_{m+n}.
\end{equation}
More generally, for each (ordered) tuple of integers $m_1,\dots,m_t\geq 0$, we can similarly define a map
\begin{equation}\label{equ: multi cup map}
\mathbf{E}_{m_1}\otimes_E\cdots\otimes_E\mathbf{E}_{m_t}\rightarrow \mathbf{E}_{\sum_{\ell=1}^t m_{\ell}}
\end{equation}
by an increasing induction on $t\geq 2$.

\subsection{Cup product between $\mathrm{Ext}$ groups}\label{subsec: cup tits}
By investigating some commutative diagram (see (\ref{equ: 5 complex cup total 1})), we reduce the study of the cup product maps between $\mathrm{Ext}_{G}^{\bullet}$ between certain Tits complex to the study of the normalized cup product maps introduced in \S \ref{subsec: cup Tits double}. Our main result here is Proposition~\ref{prop: Tits cup transfer} whose applications include Corollary~\ref{cor: bottom Tits cup transfer} and Corollary~\ref{cor: bottom coxeter cup Ext}.

We continue to write $r\defeq\#\Delta=n-1$ for short.
Given a pair of subsets $N'\subseteq N\subseteq\Delta$ and $k\in\Z$, we use the shortened notation
\begin{equation}\label{equ: Tits Ext notation}
\mathbf{E}^{k}_{N,N'}\defeq \mathrm{Ext}_{G}^{k}(\mathbf{C}_{N,\Delta},\mathbf{C}_{N',\Delta})=\mathrm{Ext}_{G}^{k}(\mathbf{C}_{N,\Delta}[-r],\mathbf{C}_{N',\Delta}[-r])
\end{equation}
Let $N_2\subseteq N_1\subseteq N_0\subseteq\Delta$ be subsets and $k_0,k_1\in \Z$. We consider the following cup product map
\begin{equation}\label{equ: main Tits Ext cup}
\mathbf{E}^{k_0}_{N_0,N_1}\otimes_E \mathbf{E}^{k_1}_{N_1,N_2}\buildrel\cup\over\longrightarrow \mathbf{E}^{k_0+k_1}_{N_0,N_2}.
\end{equation}
We write $J_0\defeq N_1\sqcup(\Delta\setminus N_0)$ and $J_0'\defeq N_2\sqcup(\Delta\setminus N_1)$ for short and note that we have $J_0\cup J_0'=\Delta$ and $J_0\cap J_0'=N_2\sqcup(\Delta\setminus N_0)$.
We write $I_0\defeq\Delta\setminus J_0=N_0\setminus N_1$ and $I_0'\defeq\Delta\setminus J_0'=N_1\setminus N_2$ for short, with $r_0\defeq\#I_0$, $r_1\defeq \#I_0'$ and $r_2\defeq r_0+r_1=\#I_0\sqcup I_0'$.
\begin{lem}\label{lem: Tits cup shift}
Let $N_2\subseteq N_1\subseteq N_0\subseteq\Delta$ and $k_0,k_1\in \Z$ be as above. Then we have the following commutative diagram
\begin{equation}\label{equ: Tits cup shift}
\xymatrix{
\mathbf{E}^{k_0}_{N_0,N_1} \ar^{\wr}[d]  & \otimes_E & \mathbf{E}^{k_1}_{N_1,N_2} \ar^{\wr}[d] \ar^{\cup}[r]& \mathbf{E}^{k_0+k_1}_{N_0,N_2} \ar^{\wr}[d] \\
\mathbf{E}^{k_0}_{\Delta,J_0} & \otimes_E & \mathbf{E}^{k_1}_{J_0,J_0\cap J_0'} \ar^{\cup}[r] & \mathbf{E}^{k_0+k_1}_{\Delta,J_0\cap J_0'}\\
}
\end{equation}
with the vertical isomorphisms from Proposition~\ref{prop: Ext complex std seq}.
\end{lem}
\begin{proof}
It follows from Proposition~\ref{prop: Ext complex std seq} that we have the following commutative diagram
\begin{equation}\label{equ: Tits cup shift diagram}
\xymatrix{
\mathbf{E}^{k_0}_{N_0,N_1} \ar@{=}[d] & \otimes_E & \mathbf{E}^{k_1}_{N_1,N_2} \ar^{\wr}[d] \ar^{\cup}[r] & \mathbf{E}^{k_0+k_1}_{N_0,N_2} \ar^{\wr}[d]\\
\mathbf{E}^{k_0}_{N_0,N_1} \ar^{\wr}[d] & \otimes_E & \mathbf{E}^{k_1}_{N_1,J_0\cap J_0'} \ar^{\cup}[r] & \mathbf{E}^{k_0+k_1}_{N_0,J_0\cap J_0'} \ar@{=}[d]\\
\mathbf{E}^{k_0}_{N_0,J_0} & \otimes_E & \mathbf{E}^{k_1}_{J_0,J_0\cap J_0'} \ar_{\wr}[u] \ar@{=}[d] \ar^{\cup}[r] & \mathbf{E}^{k_0+k_1}_{N_0,J_0\cap J_0'}\\
\mathbf{E}^{k_0}_{\Delta,J_0} \ar_{\wr}[u] & \otimes_E & \mathbf{E}^{k_1}_{J_0,J_0\cap J_0'} \ar^{\cup}[r] & \mathbf{E}^{k_0+k_1}_{\Delta,J_0\cap J_0'} \ar_{\wr}[u]
}
\end{equation}
We deduce (\ref{equ: Tits cup shift}) from (\ref{equ: Tits cup shift diagram}) by inverting all the upward vertical maps and taking compositions.
\end{proof}

Thanks to Lemma~\ref{lem: Tits cup shift}, we can reduce the study of (\ref{equ: main Tits Ext cup}) to that of the following cup product map
\begin{equation}\label{equ: main Tits Ext cup reduced}
\mathbf{E}^{k_0}_{\Delta,J_0}\otimes_E \mathbf{E}^{k_1}_{J_0,J_0\cap J_0'}\buildrel\cup\over\longrightarrow \mathbf{E}^{k_0+k_1}_{\Delta,J_0\cap J_0'}
\end{equation}
Following the notation of \S \ref{subsec: cup Tits double}, we continue to write $S_{J_0}\defeq \{I\mid J_0\subseteq I\subseteq \Delta\}$ and similarly for $S_{J_0'}$ and $S_{J_0\cap J_0'}$.
We consider below several $3$-complex $\cU_{\ast}=\cU_{\ast}^{\bullet,\bullet,\bullet}$, double complex $\cV_{\ast}=\cV_{\ast}^{\bullet,\bullet}$ and $5$-complex $\cW_{\ast}=\cW_{\ast}^{\bullet,\bullet,\bullet,\bullet,\bullet}$.
The indices of $\cU_{\ast}$ (resp,~$\cV_{\ast}$, resp.~$\cW_{\ast}$) are written as integers $-\ell'$, $\ell$ and $k'$ (resp.~$\ell$ and $k$, resp.~$\ell''$, $-\ell'$, $\ell$, $k$ and $k'$).
We always write $I,I'$ for elements in $S_{J_0}$ and $I''$ for an element in $S_{J_0\cap J_0'}$, and thus will abbreviate $\bigoplus_{I\in S_{J_0}}$ as $\bigoplus_{I}$ and similarly for others.
Here we always write $\ell=\#\Delta\setminus I$, $\ell'=\#\Delta\setminus I'$ and $\ell''=\#\Delta\setminus I''$ by convention.
We write $C_{J_0}\defeq (V_{J_0}^{\rm{an}})^\vee$ and $C_{J_0\cap J_0'}\defeq (V_{J_0\cap J_0'}^{\rm{an}})^\vee$ for short. We also use the shortened notation $D_{I}\defeq (i_{I}^{\rm{an}})^{\vee}$ for each $I\subseteq \Delta$.

We define various $3$-complex $\cU_{\ast}^{\bullet,\bullet,\bullet}$ as below.
\begin{itemize}
\item We define
\[\cU^{-r_0,r_2,k'}\defeq \Hom_{D(G)}(B_{k'}(D(G),C_{J_0\cap J_0'}), C_{J_0})\]
and $\cU^{-\ell',\ell'',k'}\defeq 0$ otherwise.
\item We define
\[\cU_{+}^{-\ell',r_2,k'}\defeq \bigoplus_{I', \#\Delta\setminus I'=\ell'}\Hom_{D(G)}(B_{k'}(D(G),C_{J_0\cap J_0'}), D_{I'})\]
and $\cU_{+}^{-\ell',\ell'',k'}\defeq 0$ otherwise.
\item We define
\[\tld{\cU}_{++}^{-\ell',\ell'',k'}\defeq \bigoplus_{I', I'', \#\Delta\setminus I'=\ell', \#\Delta\setminus I''=\ell''}\Hom_{D(G)}(B_{k'}(D(G),D_{I''}), D_{I'})\]
for each $0\leq \ell'\leq r_0$, $0\leq \ell''\leq r_2$ and $k'\geq 0$.
\item We define
\[\cU_{++}^{-\ell',\ell'',k'}\defeq \bigoplus_{I', I'', \#\Delta\setminus I'=\ell', \#\Delta\setminus I''=\ell'', I''\subseteq I'}\Hom_{D(G)}(B_{k'}(D(G),D_{I''}), D_{I'})\]
for each $0\leq \ell'\leq r_0$, $0\leq \ell''\leq r_2$ and $k'\geq 0$.
\item We define
\[\cU_{\sharp}^{-\ell',\ell'',k'}\defeq \bigoplus_{I', I'', \#\Delta\setminus I'=\ell', \#\Delta\setminus I''=\ell'', I''\subseteq I'}C^{k'}(L_{I''})\]
for each $0\leq \ell'\leq r_0$, $0\leq \ell''\leq r_2$ and $k'\geq 0$.
\item We define
\[\cU_{\flat}^{0,\ell'',k'}\defeq \bigoplus_{I'',\#\Delta\setminus I''=\ell'', J_0'\subseteq I''}C^{k'}(L_{I''})\]
for each $0\leq \ell''\leq r_1$ and $k'\geq 0$, and $\cU_{\flat}^{-\ell',\ell'',k'}\defeq 0$ otherwise.
\end{itemize}

We define various double complex $\cV_{\ast}^{\bullet,\bullet}$ as below.
\begin{itemize}
\item We define
\[\cV^{r_0,k}\defeq \Hom_{D(G)}(B_{k}(D(G),C_{J_0}), 1_{G}^{\vee})\]
and $\cV^{\ell,k}\defeq 0$ otherwise.
\item We define
\[\cV_{+}^{\ell,k}\defeq \bigoplus_{I, \#\Delta\setminus I=\ell}\Hom_{D(G)}(B_{k}(D(G),D_{I}), 1_{G}^{\vee})\]
for each $0\leq \ell\leq r_0$ and $k\geq 0$.
\item We define
\[\cV_{\sharp}^{\ell,k}\defeq \bigoplus_{I, \#\Delta\setminus I=\ell}C^{k}(L_{I})\]
for each $0\leq \ell\leq r_0$ and $k\geq 0$.
\end{itemize}

We define various $5$-complex $\cW_{\ast}^{\bullet,\bullet,\bullet,\bullet,\bullet}$ as below.
\begin{itemize}
\item We define
\[\cW^{r_0,-r_0,r_2,k,k'}\defeq \Hom_{D(G)}(B_{k'}(D(G),B_{k}(D(G),C_{J_0\cap J_0'})), 1_{G}^{\vee})\]
and $\cW^{\ell,-\ell',\ell'',k,k'}\defeq 0$ otherwise.
\item We define
\[\cW_{+}^{\ell,-\ell',r_2,k,k'}\defeq \bigoplus_{I,I',\#\Delta\setminus I=\ell,\#\Delta\setminus I'=\ell',I'\subseteq I}\Hom_{D(G)}(B_{k'}(D(G),B_{k}(D(G),C_{J_0\cap J_0'})), 1_{G}^{\vee})\]
and $\cW_{+}^{\ell,-\ell',\ell'',k,k'}\defeq 0$ otherwise.
\item We define
\[\cW_{\star}^{0,0,r_2,k,k'}\defeq \Hom_{D(G)}(B_{k'}(D(G),B_{k}(D(G),C_{J_0\cap J_0'})), 1_{G}^{\vee})\]
and $\cW_{\star}^{\ell,-\ell',\ell'',k,k'}\defeq 0$ otherwise.
\item We define
\[
\cW_{++}^{\ell,-\ell',\ell'',k,k'}\defeq
\bigoplus_{I,I',I'',\#\Delta\setminus I=\ell,\#\Delta\setminus I'=\ell',\#\Delta\setminus I''=\ell'',I''\subseteq I'\subseteq I}\Hom_{D(G)}(B_{k'}(D(G),B_{k}(D(G),D_{I''})), 1_{G}^{\vee})
\]
for each $0\leq \ell,\ell'\leq r_0$, $0\leq \ell''\leq r_2$ and $k,k'\geq 0$.
\item We define
\[\cW_{\star\star}^{0,0,\ell'',k,k'}\defeq \bigoplus_{I'',\#\Delta\setminus I''=\ell''}\Hom_{D(G)}(B_{k'}(D(G),B_{k}(D(G),D_{I''})), 1_{G}^{\vee})\]
for each $0\leq \ell''\leq r_2$ and $k,k'\geq 0$, and $\cW_{\star\star}^{-\ell,\ell',-\ell'',k,k'}\defeq 0$ otherwise.
\item We define
\[
\cW_{\sharp}^{\ell,-\ell',\ell'',k,k'}\defeq
\bigoplus_{I,I',I'',\#\Delta\setminus I=\ell,\#\Delta\setminus I'=\ell',\#\Delta\setminus I''=\ell'',I''\subseteq I'\subseteq I}\Hom_{D(L_{I''})}(B_{k'}(D(L_{I''}),B_{k}(D(L_{I''}))), 1_{L_{I''}}^\vee)
\]
for each $0\leq \ell,\ell'\leq r_0$, $0\leq \ell''\leq r_2$ and $k,k'\geq 0$.
\item We define
\[
\tld{\cW}_{\flat}^{\ell,-\ell',\ell'',k,k'}\defeq
\bigoplus_{I,I',I'',\#\Delta\setminus I=\ell,\#\Delta\setminus I'=\ell',\#\Delta\setminus I''=\ell'',I''\subseteq I'}\Hom_{D(L_{I\cap I''})}(B_{k'}(D(L_{I\cap I''}),B_{k}(D(L_{I\cap I''}))), 1_{L_{I\cap I''}}^\vee)
\]
for each $0\leq \ell,\ell'\leq r_0$, $0\leq \ell''\leq r_2$ and $k,k'\geq 0$.
\item We define
\[
\cW_{\flat}^{\ell,0,\ell'',k,k'}\defeq
\bigoplus_{I,I'',\#\Delta\setminus I=\ell,\#\Delta\setminus I''=\ell'',J_0'\subseteq I''}\Hom_{D(L_{I\cap I''})}(B_{k'}(D(L_{I\cap I''}),B_{k}(D(L_{I\cap I''}))), 1_{L_{I\cap I''}}^\vee)
\]
for each $0\leq \ell\leq r_0$, $0\leq \ell''\leq r_1$ and $k,k'\geq 0$, and $\cW_{\flat}^{\ell,-\ell',\ell'',k,k'}\defeq 0$ otherwise.
\end{itemize}

We observe that the multi-complex defined above fit into the following commutative diagrams of maps between $5$-complex
\begin{equation}\label{equ: 5 complex cup}
\xymatrix{
\cU \ar[r] & \cU_{+} & \cU_{++} \ar[l] \ar[r] & \cU_{\sharp} \ar@{=}[r] & \cU_{\sharp} \ar[r] & \cU_{\flat}\\
\otimes_E & \otimes_E & \otimes_E & \otimes_E & \otimes_E & \otimes_E\\
\cV \ar[d] & \cV_{+} \ar[l] \ar@{=}[r] \ar[d] & \cV_{+} \ar[r] \ar[d] & \cV_{\sharp} \ar@{=}[r] \ar[d] & \cV_{\sharp} \ar@{=}[r] \ar[d] & \cV_{\sharp} \ar[d]\\
\cW \ar[r] & \cW_{+} & \cW_{++} \ar[l] \ar[r] & \cW_{\sharp} & \tld{\cW}_{\flat} \ar[l] \ar[r] & \cW_{\flat} \\
& \cW_{\star} \ar[u] & \cW_{\star\star} \ar[u] \ar[l]  & & &
}.
\end{equation}

Let $A^{\bullet}\rightarrow B^{\bullet}$ be a map between complex of $E$-vector spaces 
More precisely, we have the following observations.
\begin{itemize}
\item The quasi-isomorphism $C_{J_0}[r_0]\rightarrow \mathbf{C}_{J_0,\Delta}^\vee[r]$ (see (\ref{equ: complex as induction})) induces the maps $\cU\rightarrow \cU_{+}$ and $\cV_{+}\rightarrow \cV$, whose associated maps between total complex are necessarily quasi-isomorphisms.
\item The map $\cW\rightarrow \cW_{+}$ is an embedding given by a natural truncation. For each $I\in S_{J_0}$, we consider the $5$-complex $\cW_{+,I}$ given by
\[
\cW_{+,I}^{\#\Delta\setminus I,-\ell',r_2,k,k'}\defeq
\bigoplus_{I', \#\Delta\setminus I'=\ell', I'\subseteq I}\Hom_{D(G)}(B_{k'}(D(G),B_{k}(D(G),C_{J_0\cap J_0'})), 1_{G}^{\vee})
\]
and $\cW_{+,I}^{\ell,-\ell',\ell'',k,k'}\defeq 0$ otherwise. If $I\neq J_0$, then $\cW_{+,I}^{\#\Delta\setminus I,\bullet,r_2,k,k'}$ is acyclic for each $k,k'\geq 0$, and thus $\mathrm{Tot}(\cW_{+,I})$ is acyclic. Now that $\cW_{+}/\cW$ admitting a filtration with graded pieces $\cW_{+,I}$ for each $I\in S_{J_0}\setminus\{J_0\}$ with $\mathrm{Tot}(\cW_{+,I})$ being acyclic, we conclude that $\mathrm{Tot}(\cW)\rightarrow \mathrm{Tot}(\cW_{+})$ is a quasi-isomorphism.
\item The maps $\cW_{\star}\rightarrow \cW_{+}$ and $\cW_{\star\star}\rightarrow \cW_{++}$ are embeddings given by natural truncations. For each $I'\in S_{J_0}$ and $I''\in S_{J_0\cap J_0'}$ satisfying $I''\subseteq I'$, we consider the $5$-complex $\cW_{\star\star,I',I''}$ given by
\[
\cW_{\star\star,I',I''}^{\ell,-\#\Delta\setminus I',\#\Delta\setminus I'',k,k'}\defeq
\bigoplus_{I, \#\Delta\setminus I=\ell, I'\subseteq I}\Hom_{D(G)}(B_{k'}(D(G),B_{k}(D(G),D_{I''})), 1_{G}^{\vee})
\]
and $\cW_{\star\star,I',I''}^{\ell,-\ell',\ell'',k,k'}\defeq 0$ otherwise. If $I'\neq \Delta$, then $\cW_{\star\star,I',I''}^{\bullet,-\#\Delta\setminus I',\#\Delta\setminus I'',k,k'}$ is acyclic for each $k,k'\geq 0$, and thus $\mathrm{Tot}(\cW_{\star\star,I',I''})$ is acyclic. Now that $\cW_{++}/\cW_{\star\star}$ admits a filtration with graded pieces being $\cW_{\star\star,I',I''}$ for each $I'\in S_{J_0}$ and $I''\in S_{J_0\cap J_0'}$ satisfying $I'\neq \Delta$ and $I''\subseteq I'$ with $\mathrm{Tot}(\cW_{\star\star,I',I''})$ is acyclic, we conclude that $\mathrm{Tot}(\cW_{\star\star})\rightarrow \mathrm{Tot}(\cW_{++})$ is a quasi-isomorphism.
A similar argument also shows that the map $\mathrm{Tot}(\cW_{\star})\rightarrow \mathrm{Tot}(\cW_{+})$ associated with $\cW_{\star}\rightarrow \cW_{+}$ is a quasi-isomorphism.
\item The natural truncation maps
    \begin{equation}\label{equ: Tits cup Tits resolution}
    C_{J_0\cap J_0'}[r_2]\rightarrow \mathbf{C}_{J_0\cap J_0',I'}^\vee[r]\rightarrow \mathbf{C}_{J_0\cap J_0',\Delta}^\vee[r]
    \end{equation}
    for each $I'\in S_{J_0}$ induce the maps
    \begin{equation}\label{equ: Tits cup partial resolution}
    \tld{\cU}_{++}\rightarrow \cU_{++}\rightarrow \cU_{+}.
    \end{equation}
    The composition of (\ref{equ: Tits cup partial resolution}) induces a quasi-isomorphism between the associated total complex as $C_{J_0\cap J_0'}[r_2]\rightarrow \mathbf{C}_{J_0\cap J_0',\Delta}^\vee[r]$ is a quasi-isomorphism. Now that the complex
    \[\Hom_{D(G)}(B_{\bullet}(D(G),D_{I''}), D_{I'})\]
    is acyclic whenever $I''\not\subseteq I'$ (see \ref{it: Ext PS change 1} of Lemma~\ref{lem: Ext PS change}, with $I''$ and $I'$ replacing $I_1$ and $I_0$ in \emph{loc.cit.}), the map $\tld{\cU}_{++}\rightarrow \cU_{++}$ necessarily induces a quasi-isomorphism between the associated total complex. This forces $\cU_{++}\rightarrow \cU_{+}$ to also induces a quasi-isomorphism between their associated total complex. We can similarly use the maps (\ref{equ: Tits cup Tits resolution}) to define the maps $\cW_{++}\rightarrow \cW_{+}$ and $\cW_{\star\star}\rightarrow \cW_{\star}$, with the induced map $\mathrm{Tot}(\cW_{\star\star})\rightarrow \mathrm{Tot}(\cW_{\star})$ being a quasi-isomorphism.
\item For each $I'\in S_{J_0}$ and $I''\in S_{J_0\cap J_0'}$ satisfying $I''\subseteq I'$, recall from \ref{it: Ext PS change 3} of Lemma~\ref{lem: Ext PS change} (with $I_1$, $I_0$ and $I_0'$ in \emph{loc.cit.} taken to be $I''$, $I'$ and $I_1$) and (\ref{equ: PS parabolic Levi resolution}) that we have the following quasi-isomorphisms
    \[
    \Hom_{D(G)}(B_{\bullet}(D(G),D_{I''}), D_{I'})\rightarrow
    \Hom_{D(G)}(B_{\bullet}(D(G),D_{I''}), 1_{G}^{\vee})
    \rightarrow C^{\bullet}(L_{I''})
    \]
    which is functorial with respect to the choice of $I''$ (and $I'$). We thus obtain a map $\cU_{++}\rightarrow \cU_{\sharp}$ which induce quasi-isomorphisms between their associated total complex. A similar argument also gives the maps $\cV_{+}\rightarrow \cV_{\sharp}$ with the associated map $\mathrm{Tot}(\cV_{+})\rightarrow \mathrm{Tot}(\cV_{\sharp})$ being a quasi-isomorphism.
\item For each $I''\in S_{J_0\cap J_0'}$, the inclusion $D(L_{I''})\subseteq D(G)$ and the embedding $1_{L_{I''}}^{\vee}\hookrightarrow D_{I''}$ induce the following map between double complex
    \begin{equation}\label{equ: PS to Levi double}
    \Hom_{D(G)}(B_{\bullet}(D(G),B_{\bullet}(D(G),D_{I''})), 1_{G}^\vee)\rightarrow \Hom_{D(L_{I''})}(B_{\bullet}(D(L_{I''}),B_{\bullet}(D(L_{I''}))), 1_{L_{I''}}^\vee)
    \end{equation}
    which is functorial with respect to the choice of $I''$ and thus gives a map between $5$-complex $\cW_{++}\rightarrow \cW_{\sharp}$. Note that the map between total complex associated with (\ref{equ: PS to Levi double}) fits into the following commutative diagram (induced from the map $B_{\bullet}(D(G),D_{I''})\rightarrow D_{I''}$ and the map $B_{\bullet}(D(L_{I''})\rightarrow 1_{L_{I''}}^{\vee}$)
    \begin{equation}\label{equ: PS to Levi double diagram}
    \xymatrix{
    \mathrm{Tot}(\Hom_{D(G)}(B_{\bullet}(D(G),B_{\bullet}(D(G),D_{I''})), 1_{G}^{\vee})) \ar[d] & \Hom_{D(G)}(B_{\bullet}(D(G),D_{I''}), 1_{G}^{\vee})\ar[l] \ar[d]\\
    \mathrm{Tot}(\Hom_{D(L_{I''})}(B_{\bullet}(D(L_{I''}),B_{\bullet}(D(L_{I''}))), 1_{L_{I''}}^\vee)) & C^{\bullet}(L_{I''}) \ar[l]
    }.
    \end{equation}
    Note that both horizontal maps of (\ref{equ: PS to Levi double diagram}) are clearly quasi-isomorphisms, and the right vertical map of (\ref{equ: PS to Levi double diagram}) is a quasi-isomorphism by Lemma~\ref{lem: Ext PS change} and the fact that maps in (\ref{equ: PS parabolic Levi resolution}) are quasi-isomorphisms. We thus conclude that the left vertical map of (\ref{equ: PS to Levi double diagram}) is also a quasi-isomorphism. This forces the map $\mathrm{Tot}(\cW_{++})\rightarrow \mathrm{Tot}(\cW_{\sharp})$ to be a quasi-isomorphism.
\item For each $I,I'\in S_{J_0}$ with $I'\subseteq I$, recall that we have a standard surjection $D_{I'}\twoheadrightarrow D_{I}$ which induces a map
    \begin{multline*}
    \Hom_{D(G)}(B_{k'}(D(G),C_{J_0\cap J_0'}), D_{I'})\otimes_E \Hom_{D(G)}(B_{k}(D(G),D_{I}), 1_{G}^{\vee})\\
    \rightarrow \Hom_{D(G)}(B_{k}(D(G),B_{k'}(D(G),C_{J_0\cap J_0'})), 1_{G}^{\vee})
    \end{multline*}
    which is functorial with respect to the choice of $I,I'$. This defines the map $\cU_{+}\otimes_E \cV_{+}\rightarrow \cW_{+}$. Similarly, for each $I,I'\in S_{J_0}$ with $I'\subseteq I$ and $I''\in S_{J_0\cap J_0'}$ with $I''\subseteq I'$, the standard surjection $D_{I'}\twoheadrightarrow D_{I}$ induces a map
    \begin{multline}\label{equ: Tits cup PS to CE 1}
    \Hom_{D(G)}(B_{k'}(D(G),D_{I''}), D_{I'})\otimes_E\Hom_{D(G)}(B_{k}(D(G),D_{I}), 1_{G}^{\vee})\\
    \rightarrow \Hom_{D(G)}(B_{k}(D(G),B_{k'}(D(G),D_{I''})), 1_{G}^{\vee}).
    \end{multline}
    This allows us to define the map $\cU_{++}\otimes_E \cV_{++}\rightarrow \cW_{++}$ which is compatible with $\cU_{+}\otimes_E \cV_{+}\rightarrow \cW_{+}$ as in the diagram (\ref{equ: 5 complex cup}).
    Similarly, for each $I\in S_{J_0}$ and $I''\in S_{J_0\cap J_0'}$ we have a natural map
    \begin{equation}\label{equ: Tits cup PS to CE 2}
    C^{k'}(L_{I''})\otimes_E C^{k}(L_{I})
    \rightarrow \Hom_{D(L_{I\cap I''})}(B_{k}(D(L_{I\cap I''}),B_{k'}(D(L_{I\cap I''}))), 1_{L_{I\cap I''}}^\vee),
    \end{equation}
    with (\ref{equ: Tits cup PS to CE 1}) mapping to (\ref{equ: Tits cup PS to CE 2}) when $I''\subseteq I$.
    Moreover, the maps (\ref{equ: Tits cup PS to CE 1}) and (\ref{equ: Tits cup PS to CE 2}) as well as the maps between them (which form a commutative diagram) are functorial with respect to the choice of $I,I',I''$. Hence, we obtain the following commutative diagram of maps between $5$-complex
    \begin{equation}\label{equ: Tits cup PS total to CE}
    \xymatrix{
    \cU_{++}\otimes_E\cV_{+} \ar[r] \ar[d] & \cU_{\sharp}\otimes_E\cV_{\sharp} \ar@{=}[r] \ar[d] & \cU_{\sharp}\otimes_E\cV_{\sharp} \ar[r] \ar[d] & \cU_{\flat}\otimes_E\cV_{\sharp} \ar[d]\\
    \cW_{++} \ar[r] & \cW_{\sharp} & \tld{\cW}_{\flat} \ar[l] \ar[r] & \cW_{\flat}
    }
    \end{equation}
\item The map $\cU_{\sharp}\rightarrow \cU_{\flat}$ is a surjection given by a natural truncation. For each $I''\in S_{J_0\cap J_0'}$, we consider the $3$-complex $\cU_{\sharp,I''}$ given by
\[\cU_{\sharp,I''}^{-\ell',\#\Delta\setminus I'',k'}\defeq \bigoplus_{I', \#\Delta\setminus I'=\ell', I''\subseteq I'}C^{k'}(L_{I''})\]
and $\cU_{\sharp,I''}^{-\ell',\ell'',k'}\defeq 0$ otherwise. If $J_0'\not\subseteq I''$ or equivalently $I''\cup J_0\subsetneq \Delta$, then $\cU_{\sharp,I''}^{\bullet,\#\Delta\setminus I'',k'}$ is acyclic for each $k'\geq 0$, and thus $\mathrm{Tot}(\cU_{\sharp,I''})$ is acyclic. Now that the kernel of the surjection $\cU_{\sharp}\rightarrow \cU_{\flat}$ admits a filtration with graded pieces $\cU_{\sharp,I''}$ for each $I''\in S_{J_0\cap J_0'}$ satisfying $J_0'\not\subseteq I''$ with $\mathrm{Tot}(\cU_{\sharp,I''})$ being acyclic, we conclude that $\mathrm{Tot}(\cU_{\sharp})\rightarrow \mathrm{Tot}(\cU_{\flat})$ is a quasi-isomorphism.
Similar arguments show that the natural truncation map $\tld{\cW}_{\flat}\rightarrow \cW_{\flat}$ induces a quasi-isomorphism between their associated total complex.
\item The map $\tld{\cW}_{\flat}\rightarrow \cW_{\sharp}$ is a surjection given by a natural truncation. Let $I,I'\in S_{J_0}$ with $I'\not\subseteq I$ and $J\in S_{J_0\cap J_0'}$ with $J\subseteq I\cap I'$, we consider the $4$-complex $\tld{\cW}_{\sharp,I,I',J}$ given by
\[\tld{\cW}_{\sharp,I,I',J}^{\#\Delta\setminus I,-\#\Delta\setminus I',\ell'',k}\defeq \bigoplus_{I'', \#\Delta\setminus I''=\ell'', I''\subseteq I', I\cap I''=J}C^{k}(L_{J})\]
and $\tld{\cW}_{\sharp,I,I',J}^{\ell,-\ell',\ell'',k}\defeq 0$ otherwise.
For each $I''\in S_{J_0\cap J_0'}$, we note that $I''\subseteq I'$ and $I\cap I''=J$ if and only if $J\subseteq I''\subseteq J\sqcup (I'\setminus I)$ (with $I'\setminus I\neq \emptyset$). It is thus clear that $\tld{\cW}_{\sharp,I,I',J}^{\#\Delta\setminus I,-\#\Delta\setminus I',\bullet,k}$ is acyclic for each $k\geq 0$, and therefore $\mathrm{Tot}(\tld{\cW}_{\sharp,I,I',J})$ is acyclic for each $I,I'\in S_{J_0}$ with $I'\not\subseteq I$ and $J\in S_{J_0\cap J_0'}$ with $J\subseteq I\cap I'$.
Now that the kernel of the surjection $\tld{\cW}_{\flat}\rightarrow \cW_{\sharp}$ admits a filtration with each graded piece given by some $\tld{\cW}_{\sharp,I,I',J}$ as above with $\mathrm{Tot}(\tld{\cW}_{\sharp,I,I',J})$ being acyclic, we conclude that the associated map between total complex $\mathrm{Tot}(\tld{\cW}_{\flat})\rightarrow \mathrm{Tot}(\cW_{\sharp})$ is a quasi-isomorphism.
\end{itemize}

Given a multi-complex $\cC$, we use the shortened notation $H^k(\cC)\defeq H^k(\mathrm{Tot}(\cC))$ for $k\in\Z$.
Let $k_0,k_1\in \Z$ with $k_0'\defeq k_0-r$.
To summarize, the commutative diagram (\ref{equ: 5 complex cup}) induces the following commutative diagram
\begin{equation}\label{equ: 5 complex cup total 1}
\xymatrix{
H^{k_1}(\cU) \ar^{\sim}[r] & H^{k_1}(\cU_{+}) & H^{k_1}(\cU_{++}) \ar_{\sim}[l] \ar^{\sim}[r] & H^{k_1}(\cU_{\sharp}) \ar@{=}[r] & H^{k_1}(\cU_{\sharp}) \ar^{\sim}[r] & H^{k_1}(\cU_{\flat})\\
\otimes_E & \otimes_E & \otimes_E & \otimes_E & \otimes_E & \otimes_E \\
H^{k_0'}(\cV) \ar[d] & H^{k_0'}(\cV_{+}) \ar_{\sim}[l] \ar@{=}[r] \ar[d] & H^{k_0'}(\cV_{+}) \ar^{\sim}[r] \ar[d] & H^{k_0'}(\cV_{\sharp}) \ar@{=}[r] \ar[d] & H^{k_0'}(\cV_{\sharp}) \ar@{=}[r] \ar[d] & H^{k_0'}(\cV_{\sharp}) \ar[d]\\
H^{k_0'+k_1}(\cW) \ar^{\sim}[r] & H^{k_0'+k_1}(\cW_{+}) & H^{k_0'+k_1}(\cW_{++}) \ar[l] \ar^{\sim}[r]  & H^{k_0'+k_1}(\cW_{\sharp}) & H^{k_0'+k_1}(\tld{\cW}_{\flat}) \ar_{\sim}[l] \ar^{\sim}[r] & H^{k_0'+k_1}(\cW_{\flat})\\
 & H^{k_0'+k_1}(\cW_{\star}) \ar^{\wr}[u] & H^{k_0'+k_1}(\cW_{\star\star}) \ar_{\sim}[l] \ar^{\wr}[u]  & & &
}
\end{equation}
Since the natural map $S_{J_0}\times S_{J_0'}\rightarrow S_{J_0\cap J_0'}: (I,I'')\mapsto I\cap I''$ is a bijection (with $\#I\cap I''=\#I+\#I''-r$), the map
\[C^{\bullet}(L_{I\cap I''})\rightarrow \mathrm{Tot}(\Hom_{D(L_{I\cap I''})}(B_{\bullet}(D(L_{I\cap I''}),B_{\bullet}(D(L_{I\cap I''}))), 1_{L_{I\cap I''}}^\vee))\]
for each $I\in S_{J_0}$ and $I''\in S_{J_0'}$ (which is a quasi-isomorphism and is functorial with respect to the choice of $I$ and $I''$) gives a map
\[\cT_{J_0\cap J_0',\Delta}^{\bullet,\bullet}\rightarrow \mathrm{Tot}_{1,2,3}(\mathrm{Tot}_{4,5}(\cW_{\flat}))\]
which induces an isomorphism
\begin{equation}\label{equ: Tits composition transfer}
H^{k}(\mathrm{Tot}(\cT_{J_0\cap J_0',\Delta}^{\bullet,\bullet}))\buildrel\sim\over\longrightarrow H^{k}(\mathrm{Tot}(\cW_{\flat}))
\end{equation}
for each $k\in\Z$.
Now that we have $\mathrm{Tot}_{1,2}(\cU_{\flat}^{\bullet,\bullet,\bullet})=\cT_{J_0',\Delta}^{\bullet,\bullet}[-r,0]$ and $\cV_{\sharp}^{\bullet,\bullet}=\cT_{J_0,\Delta}^{\bullet,\bullet}$, we combine the right most vertical map of (\ref{equ: 5 complex cup total 1}) with (\ref{equ: Tits composition transfer}) and obtain a map
\begin{equation}\label{equ: Tits cup coh}
H^{k_1-r}(\mathrm{Tot}(\cT_{J_0',\Delta}^{\bullet,\bullet}))\otimes_E H^{k_0-r}(\mathrm{Tot}(\cT_{J_0,\Delta}^{\bullet,\bullet}))\rightarrow H^{k_0+k_1-r}(\mathrm{Tot}(\cT_{J_0\cap J_0',\Delta}^{\bullet,\bullet})).
\end{equation}

\begin{prop}\label{prop: Tits cup transfer}
Let $J_0,J_0'\subseteq \Delta$ and $k_0,k_1\in \Z$ be as above (with $J_0\cup J_0'=\Delta$). Then we have the following commutative diagram
\begin{equation}\label{equ: Tits cup transfer}
\xymatrix{
\mathbf{E}^{k_1}_{J_0,J_0\cap J_0'} \ar^{\wr}[d]  & \otimes_E & \mathbf{E}^{k_0}_{\Delta,J_0} \ar^{\wr}[d] \ar^{\cup}[r]& \mathbf{E}^{k_0+k_1}_{\Delta,J_0\cap J_0'} \ar^{\wr}[d]\\
H^{k_1-r}(\mathrm{Tot}(\cT_{J_0',\Delta}^{\bullet,\bullet})) & \otimes_E & H^{k_0-r}(\mathrm{Tot}(\cT_{J_0,\Delta}^{\bullet,\bullet})) \ar[r] & H^{k_0+k_1-r}(\mathrm{Tot}(\cT_{J_0\cap J_0',\Delta}^{\bullet,\bullet}))
}
\end{equation}
with the bottom horizontal map from (\ref{equ: Tits cup coh}), and the vertical isomorphisms equal those from Proposition~\ref{prop: Ext complex std seq} up to a sign.
\end{prop}
\begin{proof}
We deduce from (\ref{equ: 5 complex cup total 1}) a commutative diagram of the form (\ref{equ: Tits cup transfer}), with the left (resp.~middle, resp.~right) vertical map of (\ref{equ: Tits cup transfer}) being the composition of horizontal maps in the top (resp.~middle, resp.~bottom) row of (\ref{equ: 5 complex cup total 1}). It remains to check the aforementioned vertical maps in (\ref{equ: Tits cup transfer}) equal the maps from Proposition~\ref{prop: Ext complex std seq} up to a sign. The claim for the middle vertical map of (\ref{equ: Tits cup transfer}) is clear as it is constructed directly from the isomorphism $\mathbf{C}_{J_0,\Delta}\cong V_{J_0}^{\rm{an}}[\#J_0]$. We divide the rest of the proof into two steps.

\textbf{Step $1$}: We prove the claim for the left vertical map of (\ref{equ: Tits cup transfer}).\\
We define some extra $3$-complex $\cU_{\ast}^{\bullet,\bullet,\bullet}$ as below.
\begin{itemize}
\item We define
\[\tld{\cU}_{+}^{-r_0,\ell'',k'}\defeq \bigoplus_{I'', \#\Delta\setminus I''=\ell''}\Hom_{D(G)}(B_{k'}(D(G),D_{I''}), D_{J_{0}})\]
for each $0\leq \ell''\leq r_2$ and $k'\geq 0$, and $\tld{\cU}_{+}^{\ell',\ell'',k'}\defeq 0$ otherwise.
\item We define
\[\cU_{\natural}^{-r_0,\ell'',k'}\defeq \bigoplus_{I'', \#\Delta\setminus I''=\ell'', J_0'\subseteq I''}\Hom_{D(G)}(B_{k'}(D(G),D_{I''}), D_{J_{0}})\]
for each $0\leq \ell''\leq r_1$ and $k'\geq 0$, and $\tld{\cU}_{+}^{\ell',\ell'',k'}\defeq 0$ otherwise.
\item We define 
\[\tld{\cU}_{\dagger}^{-\ell',\ell'',k'}\defeq \bigoplus_{I', I'', \#\Delta\setminus I'=\ell', \#\Delta\setminus I''=\ell'', J_0'\subseteq I''}\Hom_{D(G)}(B_{k'}(D(G),D_{I''}),D_{I'})\]
and $\tld{\cU}_{\dagger}^{\ell',\ell'',k'}\defeq 0$ otherwise.
\item We define 
\[\cU_{\dagger}^{0,\ell'',k'}\defeq \bigoplus_{I'', \#\Delta\setminus I''=\ell'', J_0'\subseteq I''}\Hom_{D(G)}(B_{k'}(D(G),D_{I''}),1_{G}^{\vee})\]
and $\cU_{\dagger}^{\ell',\ell'',k'}\defeq 0$ otherwise.
\end{itemize}
We have the following commutative diagram of maps between $3$-complex
\begin{equation}\label{equ: Tits cup transfer U diagram}
\xymatrix{
\cU \ar[r] & \cU_{+} & \cU_{++} \ar[l] & \tld{\cU}_{++} \ar[l] \ar[r] \ar[ld] & \cU_{\sharp} \ar[d]\\
\tld{\cU}_{+} \ar[u] \ar[r] \ar[rrru] & \cU_{\natural} \ar[r] & \tld{\cU}_{\dagger} & \cU_{\dagger} \ar[l] \ar[r] & \cU_{\flat}
}
\end{equation}
with all associated maps between their total complex being isomorphisms. 
Note that the left vertical map of (\ref{equ: Tits cup transfer}) is nothing but the isomorphism
\begin{equation}\label{equ: Tits cup transfer U isom}
H^{k_1}(\cU)\buildrel\sim\over\longrightarrow H^{k_1}(\cU_{\flat})
\end{equation}
constructed from the first row and the right vertical map of (\ref{equ: Tits cup transfer U diagram}).
We finish the proof of \textbf{Step $1$} by observing that the other isomorphism of the form (\ref{equ: Tits cup transfer U isom}) constructed in Proposition~\ref{prop: Ext complex std seq} actually equals the composition of the various isomorphisms (or their inverses) between $H^{k_1}(\cU_{\ast})$ going through the bottom row of (\ref{equ: Tits cup transfer U diagram}).

\textbf{Step $2$}: We prove the claim for the right vertical map of (\ref{equ: Tits cup transfer}).\\
Upon replacing 
\begin{equation}\label{equ: Tits cup transfer W flat}
\Hom_{D(L_{\ast})}(B_{k'}(D(L_{\ast}),B_{k}(D(L_{\ast}))), 1_{L_{\ast}}^\vee)
\end{equation}
in the definition of $\tld{\cW}_{\flat}$ and $\cW_{\flat}$ with
\begin{equation}\label{equ: Tits cup transfer W dagger}
\Hom_{D(G)}(B_{k'}(D(G),B_{k}(D(G),D_{\ast})), 1_{G}^{\vee}),
\end{equation}
we may define another two $5$-complex $\cW_{\ddagger}$ and $\cW_{\dagger}$ which fits into the following commutative diagram of maps between $5$-complex
\begin{equation}\label{equ: Tits cup transfer W diagram}
\xymatrix{
\cW_{\star\star} \ar[r] & \cW_{++} \ar[d] & \cW_{\ddagger} \ar[l] \ar[r] \ar[d] & \cW_{\dagger} \ar[d]\\
& \cW_{\sharp} & \tld{\cW}_{\flat} \ar[l] \ar[r] & \cW_{\flat}
}
\end{equation}
with horizontal maps given by natural truncations, vertical maps induced from the map from (\ref{equ: Tits cup transfer W dagger}) to (\ref{equ: Tits cup transfer W flat}), and all maps in (\ref{equ: Tits cup transfer W diagram}) inducing quasi-isomorphisms between their associated total complex.
According to the definition of each $\cW_{\ast}$ in the top row of (\ref{equ: Tits cup transfer W diagram}), we have the following maps between $3$-complex in $\mathrm{Mod}_{D(G)}$
\begin{equation}\label{equ: Tits cup transfer X diagram}
\cX_{\star\star} \leftarrow \cX_{++} \rightarrow \cX_{\ddagger} \leftarrow \cX_{\dagger} 
\end{equation}
whose associated maps between total complex are all quasi-isomorphisms, and which upon applying 
\begin{equation}\label{equ: Tits cup transfer double resolution}
\Hom_{D(G)}(B_{\bullet}(D(G),B_{\bullet}(D(G),-)),1_{G}^{\vee}) 
\end{equation}
recover the top row of (\ref{equ: Tits cup transfer W diagram}). For each $\ast\in\{\star\star,++,\ddagger,\dagger\}$, the following canonical truncation maps
\[\mathrm{Tot}(\cX_{\ast})\rightarrow \tau_{\geq -r_2}\mathrm{Tot}(\cX_{\ast})\leftarrow \tau_{\leq -r_2}\tau_{\geq -r_2}\mathrm{Tot}(\cX_{\ast})=H^{-r_2}(\mathrm{Tot}(\cX_{\ast}))\]
induces the following isomorphism in $D^{b}(\mathrm{Mod}_{D(G)})$
\[
H^{-r_2}(\mathrm{Tot}(\cX_{\ast}))\dashrightarrow  \mathrm{Tot}(\cX_{\ast})
\]
which is functorial with respect to the choice of $\ast$. In other words, we obtain a commutative diagram in $D^{b}(\mathrm{Mod}_{D(G)})$ of the form
\begin{equation}\label{equ: Tits cup transfer X truncation}
\xymatrix{
H^{-r_2}(\mathrm{Tot}(\cX_{\star\star})) \ar[d] & H^{-r_2}(\mathrm{Tot}(\cX_{++})) \ar[l] \ar[r] \ar@{-->}[d] & H^{-r_2}(\mathrm{Tot}(\cX_{\ddagger})) \ar@{-->}[d] & H^{-r_2}(\mathrm{Tot}(\cX_{\dagger})) \ar[l] \ar[d]\\
\mathrm{Tot}(\cX_{\star\star}) & \mathrm{Tot}(\cX_{++}) \ar[l] \ar[r] & \mathrm{Tot}(\cX_{\ddagger}) & \mathrm{Tot}(\cX_{\dagger}) \ar[l]\\
}
\end{equation}
with the top row of (\ref{equ: Tits cup transfer X truncation}) being isomorphisms in $\mathrm{Mod}_{D(G)}$. Note from the definition of $\cX_{\star\star}$ and $\cX_{\dagger}$ that natural map $C_{J_0\cap J_0'}\hookrightarrow D_{J_0\cap J_0'}$ induces an isomorphism
\[C_{J_0\cap J_0'}\buildrel\sim\over\longrightarrow H^{-r_2}(\mathrm{Tot}(\cX_{\ast})\] 
for each $\ast\in\{\star\star,\dagger\}$, which together with (\ref{equ: Tits cup transfer X truncation}) gives the following commutative diagram of isomorphisms in $D^{b}(\mathrm{Mod}_{D(G)})$
\begin{equation}\label{equ: Tits cup transfer X St}
\xymatrix{
C_{J_0\cap J_0'} \ar^{\sim}[rr] \ar[d] & & H^{-r_2}(\mathrm{Tot}(\cX_{\star\star})) \ar[r] \ar@{-->}[d] & \mathrm{Tot}(\cX_{\star\star}) \ar@{-->}[d] \\
C_{J_0\cap J_0'} \ar^{\sim}[rr] & & H^{-r_2}(\mathrm{Tot}(\cX_{\dagger})) \ar[r] & \mathrm{Tot}(\cX_{\dagger})
}
\end{equation}
with the LHS vertical map being an automorphism of $C_{J_0\cap J_0'}$ which is necessarily a scalar in $E^{\times}$. 
Upon replacing $D(G)$ with the unital ring of locally constant distributions $D(G,\Z)$ on $G$ with coefficients in $\Z$, we can consider the parallel diagram of (\ref{equ: Tits cup transfer X St}) in this setting and see that the aforementioned scalar in $E^{\times}$ must be an element of $\Z^{\times}=\{1,-1\}$.
Applying (\ref{equ: Tits cup transfer double resolution}) to (\ref{equ: Tits cup transfer X St}), taking $H^{k_0'+k_1}(\mathrm{Tot}(-))$ and then combined with the left bottom corner of (\ref{equ: 5 complex cup total 1}), we obtain the following commutative diagram of isomorphisms between $E$-vector spaces
\begin{equation}\label{equ: Tits cup transfer W coh}
\xymatrix{
H^{k_0'+k_1}(\cW) \ar^{\sim}[rr] \ar[d] & & \mathrm{Ext}_{D(G)}^{k_0'+k_1}(C_{J_0\cap J_0'}[-r_2],1_{G}^{\vee}) \ar[d]\\
H^{k_0'+k_1}(\cW_{\star\star}) \ar^{\sim}[rr] \ar[d] & & \mathrm{Ext}_{D(G)}^{k_0'+k_1}(C_{J_0\cap J_0'}[-r_2],1_{G}^{\vee}) \ar[d]\\
H^{k_0'+k_1}(\cW_{\dagger}) \ar^{\sim}[rr] & & \mathrm{Ext}_{D(G)}^{k_0'+k_1}(C_{J_0\cap J_0'}[-r_2],1_{G}^{\vee})
}
\end{equation}
with RHS vertical isomorphisms of (\ref{equ: Tits cup transfer W coh}) being scalar multiplications by signs.
Now that the right vertical map of (\ref{equ: Tits cup transfer}) is given by the composition of LHS vertical isomorphisms and then bottom horizontal isomorphism of (\ref{equ: Tits cup transfer W coh}), while the other isomorphism constructed in Proposition~\ref{prop: Ext complex std seq} is simply the top horizontal isomorphism of (\ref{equ: Tits cup transfer W coh}), we conclude that they equal up to a sign and finish the proof of \textbf{Step $2$}.
\end{proof}

Upon taking $N_2\subseteq N_1\subseteq N_0$, $k_0$ and $k_1$ in (\ref{equ: main Tits Ext cup}) to be $I_1\subseteq I\subseteq I_0$, $2\#I_0\setminus I$ and $2\#I\setminus I_1$, we obtain a cup product map (see also (\ref{equ: bottom deg Ext St}) for notation)
\begin{equation}\label{equ: main Tits Ext bottom cup}
\kappa_{I_0,I_1}^{I}: \mathbf{E}_{I_0,I}\otimes_E\mathbf{E}_{I,I_1}\buildrel\cup\over\longrightarrow \mathbf{E}_{I_0,I_1}.
\end{equation}
\begin{cor}\label{cor: bottom Tits cup transfer}
Let $I_1\subseteq I\subseteq I_0\subseteq \Delta$ be as above. The map (\ref{equ: main Tits Ext bottom cup}) fits into the following commutative diagram
\begin{equation}\label{equ: bottom Tits cup transfer}
\xymatrix{
\mathbf{E}_{I_0,I} \ar[d] & \otimes_E & \mathbf{E}_{I,I_1} \ar^{\cup}[r] \ar[d] & \mathbf{E}_{I_0,I_1} \ar[d]\\
\mathbf{E}_{I_0\setminus I} & \otimes_E & \mathbf{E}_{I\setminus I_1} \ar^{\cup}[r] & \mathbf{E}_{I_0\setminus I_1}
}
\end{equation}
with the vertical maps being isomorphisms from \ref{it: Ext complex std seq} of Proposition~\ref{prop: Ext complex std seq} up to a sign.
\end{cor}
\begin{proof}
Thanks to Lemma~\ref{lem: Tits cup shift} (upon taking $N_2\subseteq N_1\subseteq N_0$, $k_0$ and $k_1$ in \emph{loc.cit.} to be $I_1\subseteq I\subseteq I_0$, $2\#I_0\setminus I$ and $2\#I\setminus I_1$), we have a commutative diagram of the form
\[
\xymatrix{
\mathbf{E}_{I_0,I} \ar[d] & \otimes_E & \mathbf{E}_{I,I_1} \ar^{\cup}[r] \ar[d] & \mathbf{E}_{I_0,I_1} \ar[d]\\
\mathbf{E}_{\Delta,I\sqcup(\Delta\setminus I_0)} & \otimes_E & \mathbf{E}_{I\sqcup(\Delta\setminus I_0),I_1\sqcup(\Delta\setminus I_0)} \ar^{\cup}[r] & \mathbf{E}_{\Delta,I\sqcup(\Delta\setminus I_0)}
}
\]
with vertical maps being isomorphisms from \ref{it: Ext complex std seq} of Proposition~\ref{prop: Ext complex std seq}.
Taking $J_0$ and $J_0'$ in Proposition~\ref{prop: Tits cup transfer}  to be $\Delta$, $I\sqcup(\Delta\setminus I_0)$ and $I_1\sqcup(\Delta\setminus I)$ (with $J_0\cap J_0'$ taken to be $I_1\sqcup(\Delta\setminus I_0)$), we deduce from (\ref{equ: Tits cup transfer}) a commutative diagram of the form
\[
\xymatrix{
\mathbf{E}_{\Delta,I\sqcup(\Delta\setminus I_0)} \ar[d] & \otimes_E & \mathbf{E}_{I\sqcup(\Delta\setminus I_0),I_1\sqcup(\Delta\setminus I_0)} \ar^{\cup}[r] \ar[d] & \mathbf{E}_{\Delta,I\sqcup(\Delta\setminus I_0)} \ar[d]\\
\mathbf{E}_{I_0\setminus I} & \otimes_E & \mathbf{E}_{I\setminus I_1} \ar^{\cup}[r] & \mathbf{E}_{I_0\setminus I_1}
}
\]
with vertical maps being isomorphisms from \ref{it: Ext complex std seq} of Proposition~\ref{prop: Ext complex std seq} up to a sign.
We conclude (\ref{equ: bottom Tits cup transfer}) by combining the two commutative diagrams above.
\end{proof}

\begin{cor}\label{cor: bottom coxeter cup Ext}
Let $I_1\subseteq I\subseteq I_0\subseteq \Delta$ be as above. Let $x\in\Gamma^{I_0\setminus I}$ and $y\in\Gamma^{I\setminus I_1}$. We have the following results on the map (\ref{equ: main Tits Ext bottom cup}).
\begin{enumerate}[label=(\roman*)]
\item \label{it: bottom coxeter cup Ext 1} We have
\begin{equation}\label{equ: bottom coxeter cup Ext bound}
\kappa_{I_0,I_1}^{I}(\mathrm{Fil}_{x}(\mathbf{E}_{I_0,I})\otimes_E\mathrm{Fil}_{y}(\mathbf{E}_{I,I_0}))\subseteq\sum_{w\in\Gamma_{x,y}}\mathrm{Fil}_{w}(\mathbf{E}_{I_0,I_1}).
\end{equation}
\item \label{it: bottom coxeter cup Ext 2} For each $w\in\Gamma_{x,y}\subseteq\Gamma^{I_0\setminus I_1}$, the composition of
\[\mathrm{Fil}_{x}(\mathbf{E}_{I_0,I})\otimes_E\mathrm{Fil}_{y}(\mathbf{E}_{I,I_0})\rightarrow \sum_{u\in\Gamma_{x,y}}\mathrm{Fil}_{u}(\mathbf{E}_{I_0,I_1})\twoheadrightarrow \mathrm{gr}_{w}(\mathbf{E}_{I_0,I_1})\]
is surjective.
\end{enumerate}
\end{cor}
\begin{proof}
This follows from Corollary~\ref{cor: bottom Tits cup transfer} by combining Proposition~\ref{prop: coxeter cup comparison} and Proposition~\ref{prop: cup x y nonvanishing}.
\end{proof}

\section{Preliminary on the Galois side}\label{sec: Galois}
In this section, we prove various results on filtered $(\varphi,N)$-modules and $(\varphi,\Gamma)$-modules that will be useful in later sections.
\subsection{Some filtered $(\varphi,N)$-modules}\label{subsec: filtered phi N}
We fix a choice of $n\geq 2$ throughout this section.
Throughout this section, we abuse the term \emph{filtered $(\varphi,N)$-module} for a finite dimensional $E$-vector space $D$ equipped with a $(\varphi,N)$-action satisfying $N\varphi=p\varphi N$ as well as decreasing filtration $\mathrm{Fil}^{\bullet}(D)$. 
We will relate this abused terminology to the usual one at the beginning of \S~\ref{subsec: universal Galois}.
We study the filtered $(\varphi,N)$-modules which are \emph{universal special with level $n-1$}, see Definition~\ref{def: universal phi N}. This notion is naturally motivated by our study of the de Rham complex of the Drinfeld upper half space (cf.~Lemma~\ref{lem: Hodge image}) and is closely related to certain successive universal de Rham extensions between rank one $(\varphi,\Gamma)$-modules (see Proposition~\ref{prop: Galois dR Ext exact}). 

For each $0\leq k\leq n-1$, we write $J^{k}\defeq [1,n-1-k]$ for short (with convention $J^{n-1}=\emptyset$).
For each $0\leq k\leq k'\leq n-1$, we write $\mathbf{E}_{k,k'}\defeq\mathbf{E}_{J^{k},J^{k'}}$, $\mathbf{E}_{k,k'}^{<}\defeq\mathbf{E}_{J^{k},J^{k'}}^{<}$, $\overline{\mathbf{E}}_{k,k'}\defeq \mathbf{E}_{k,k'}/\mathbf{E}_{k,k'}^{<}$ and $x_{k,k'}\defeq w_{n-k,n-1-k'}=s_{n-k}\cdots s_{n-1-k'}$ for short (with $\mathrm{Supp}(x_{k,k'})=[n-k,n-1-k']=J^{k'}\setminus J^{k}$ and thus $\mathrm{Fil}_{x_{k,k'}}(\mathbf{E}_{k,k'})$ is defined).
Hence, for each $0\leq k\leq k'\leq k''\leq n-1$, we have a natural cup product map
\[\kappa_{k,k''}^{k'}\defeq \kappa_{J^{k},J^{k''}}^{J^{k'}}: \mathbf{E}_{k,k'}\otimes_E\mathbf{E}_{k',k''}\buildrel\cup\over\longrightarrow\mathbf{E}_{k,k''}.\]
For each $0\leq \ell\leq n-1$, we use the following shortened notation
\begin{equation}\label{equ: phi N total cup}
\kappa^{\ell}\defeq \bigoplus_{k=0}^{\ell}\kappa_{k,n-1}^{\ell}: (\bigoplus_{k=0}^{\ell}\mathbf{E}_{k,\ell})\otimes_E\mathbf{E}_{\ell,n-1}\rightarrow \bigoplus_{k=0}^{\ell}\mathbf{E}_{k,n-1}.
\end{equation}

In the following, we equip $\bigoplus_{k=0}^{n-1}\mathbf{E}_{k,n-1}$ with the $(\varphi,N)$-action characterized by the following conditions.
\begin{itemize}
\item For each $0\leq k\leq n-1$, $\varphi$ acts on $\mathbf{E}_{k,n-1}$ by $p^k$.
\item We have $N(\mathbf{E}_{0,n-1})=0$. For each $1\leq k\leq n-1$, $N$ restricts to the map
    \[\kappa_{k-1,n-1}^{k}(c_{k}\val\otimes_E -): \mathbf{E}_{k,n-1}\rightarrow \mathbf{E}_{k-1,n-1}\] 
    for some $c_{k}\in E^{\times}$.
\end{itemize}
We have the following definition.
\begin{defn}\label{def: universal phi N}
Let $n\geq 1$ and $D$ be a filtered $(\varphi,N)$-module with $\mathrm{Fil}_{H}^0(D)=D$ and $\mathrm{Fil}_{H}^n(D)=0$. We say that $D$ is \emph{universal special with level $n-1$} if there exists an isomorphism $D\cong \bigoplus_{k=0}^{n-1}\mathbf{E}_{k,n-1}$ as $(\varphi,N)$-modules as well as an $E$-line $L_{\ell}\subseteq \bigoplus_{k=0}^{\ell}\mathrm{Fil}_{x_{k,\ell}}(\mathbf{E}_{k,\ell})$ for each $0\leq \ell\leq n-1$ such that the following conditions hold.
\begin{enumerate}[label=(\roman*)]
\item \label{it: universal phi N 2} For each $0\leq k\leq \ell$, the composition of 
\[L_{k,\ell}\hookrightarrow\mathrm{Fil}_{x_{k,\ell}}(\mathbf{E}_{k,\ell})\twoheadrightarrow \mathrm{gr}_{x_{k,\ell}}(\mathbf{E}_{k,\ell})\] 
is an isomorphism, with $L_{k,\ell}$ being the projection of $L_{\ell}$ to $\mathrm{Fil}_{x_{k,\ell}}(\mathbf{E}_{k,\ell})$.
\item \label{it: universal phi N 3} We have
\begin{equation}\label{equ: universal n step filtration}
\mathrm{Fil}_{H}^m(D)\cong\sum_{\ell=m}^{n-1}\kappa^{\ell}(L_{\ell}\otimes_E\mathbf{E}_{\ell,n-1})
\end{equation}
for each $0\leq m\leq n-1$.
\end{enumerate}
\end{defn}
\begin{lem}\label{lem: universal phi N top grade}
Let $D$ be a filtered $(\varphi,N)$-module which is universal special with level $n-1$. Then we have the following results.
Let $0\leq k\leq \ell\leq n-1$. Then we have the following results.
\begin{enumerate}[label=(\roman*)]
\item \label{it: universal phi N top 1} For each $0\leq k< \ell\leq n-1$, we have $\mathrm{gr}_{x_{k,\ell}}(\mathbf{E}_{k,\ell})=\mathrm{Fil}_{x_{k,\ell}}(\mathbf{E}_{k,\ell})/\mathrm{Fil}_{x_{k+1,\ell}}(\mathbf{E}_{k,\ell})$ with
    \begin{equation}\label{equ: universal phi N sm cup}
    \mathrm{Fil}_{x_{k+1,\ell}}(\mathbf{E}_{k,\ell})=N(\mathrm{Fil}_{x_{k+1,\ell}}(\mathbf{E}_{k+1,\ell}))\subseteq \mathbf{E}_{k,\ell}^{<}.
    \end{equation}
\item \label{it: universal phi N top 2} For each $0\leq k\leq \ell\leq n-1$, the composition of
\[L_{k,\ell}\hookrightarrow\mathrm{Fil}_{x_{k,\ell}}(\mathbf{E}_{k,\ell})\hookrightarrow \mathbf{E}_{k,\ell}\twoheadrightarrow \overline{\mathbf{E}}_{k,\ell}=\mathbf{E}_{k,\ell}/\mathbf{E}_{k,\ell}^{<}\]
is an isomorphism.
\end{enumerate}
\end{lem}
\begin{proof}
Given $0\leq k< \ell\leq n-1$, we notice that any $y\in\Gamma$ that satisfies $y\unlhd x_{k,\ell}$ and $y\neq x_{k,\ell}$ must satisfy $y\unlhd x_{k+1,\ell}$, which gives $\mathrm{gr}_{x_{k,\ell}}(\mathbf{E}_{k,\ell})=\mathrm{Fil}_{x_{k,\ell}}(\mathbf{E}_{k,\ell})/\mathrm{Fil}_{x_{k+1,\ell}}(\mathbf{E}_{k,\ell})$.
We deduce (\ref{equ: universal phi N sm cup}) from the definition of $N$ and Lemma~\ref{lem: cup with sm}, and thus finish the proof of \ref{it: universal phi N top 1}.
\ref{it: universal phi N top 2} follows from (\ref{equ: universal phi N sm cup}), \ref{it: universal phi N 2} of Definition~\ref{def: universal phi N} and Theorem~\ref{thm: cup top grade}. In fact, (\ref{equ: universal phi N sm cup}) and Theorem~\ref{thm: cup top grade} ensure that any element of
\[\mathrm{Fil}_{x_{k,\ell}}(\mathbf{E}_{k,\ell})\setminus \mathrm{Fil}_{x_{k+1,\ell}}(\mathbf{E}_{k,\ell})\]
must have non-zero image under $\mathbf{E}_{k,\ell}\twoheadrightarrow \overline{\mathbf{E}}_{k,\ell}$.
\end{proof}

For $0\leq \ell\leq n-1$, we note that the composition $\kappa^{\ell}(L_{\ell}\otimes_E\mathbf{E}_{\ell,n-1})\subseteq \bigoplus_{k=0}^{\ell}\mathbf{E}_{k,n-1}\twoheadrightarrow \mathbf{E}_{\ell,n-1}$ is an isomorphism by our assumption on $L_{\ell}$. Hence, the sum in (\ref{equ: universal n step filtration}) is a direct sum, and the inclusion $\mathrm{Fil}_{H}^m(D)\subseteq D$ induces an isomorphism between $E$-vector spaces
\begin{equation}\label{equ: complement of Hodge}
\mathrm{Fil}_{H}^m(D)\oplus\bigoplus_{\ell=0}^{m-1}\mathbf{E}_{\ell,n-1}\buildrel\sim\over\longrightarrow D
\end{equation}
for each $0\leq m\leq n-1$. Moreover, since $\kappa^{\ell}(L_{\ell}\otimes_E\mathbf{E}_{\ell,n-1})\subseteq \bigoplus_{k=0}^{\ell}\mathbf{E}_{k,n-1}\subseteq \bigoplus_{k=0}^m\mathbf{E}_{k,n-1}$ for $0\leq \ell\leq m\leq n-1$ and $\Dim_E \kappa^{\ell}(L_{\ell}\otimes_E\mathbf{E}_{\ell,n-1})=\Dim_E \mathbf{E}_{\ell,n-1}$ for $0\leq \ell\leq n-1$, we obtain an equality between subspaces of $D$
\begin{equation}\label{equ: sum of Hodge}
\bigoplus_{k=0}^m\mathbf{E}_{k,n-1}=\bigoplus_{\ell=0}^m\kappa^{\ell}(L_{\ell}\otimes_E\mathbf{E}_{\ell,n-1})
\end{equation}
for each $0\leq m\leq n-1$.

\begin{lem}\label{lem: Fil max sub phi N}
Let $D$ be filtered $(\varphi,N)$-module which is universal special with level $n-1$. Let $D^{\flat}$ be the minimal $(\varphi,N)$-submodule of $D$ that contains $\mathrm{Fil}_{H}^{n-1}(D)$. We have the following results.
\begin{enumerate}[label=(\roman*)]
\item \label{it: Fil max sub phi N 1} The isomorphism $D\cong \bigoplus_{k=0}^{n-1}\mathbf{E}_{k,n-1}$ between filtered $(\varphi,N)$-modules restricts to an isomorphism
    \begin{equation}\label{equ: Fil max sub phi N isom}
    D^{\flat}\cong \bigoplus_{k=0}^{n-1}\mathrm{Fil}_{x_{k,n-1}}(\mathbf{E}_{k,n-1})
    \end{equation}
    between $(\varphi,N)$-modules.
\item \label{it: Fil max sub phi N 2} We have
\begin{equation}\label{equ: Fil max sub phi N intersection}
\mathrm{Fil}_{H}^{\ell}(D)\cap D^{\flat}=\mathrm{Fil}_{H}^{n-1}(D)
\end{equation}
for each $1\leq \ell\leq n-1$. In particular, if we endow $D^{\flat}$ with the Hodge filtration induced from $D$, then we have $\mathrm{Fil}_{H}^{0}(D^{\flat})=D^{\flat}$, $\mathrm{Fil}_{H}^{n}(D^{\flat})=0$ and $\mathrm{Fil}_{H}^{\ell}(D^{\flat})=\mathrm{Fil}_{H}^{n-1}(D)$ for each $1\leq \ell \leq n-1$.
\end{enumerate}
\end{lem}
\begin{proof}
We prove \ref{it: Fil max sub phi N 1}.\\
Now that $\mathrm{Fil}_{H}^{n-1}(D)=L_{n-1}$ and $\varphi$ acts by $p^k$ on $L_{k,n-1}$ for each $0\leq k\leq n-1$, we see that the minimal $\varphi$-submodule of $D$ containing $\mathrm{Fil}_{H}^{n-1}(D)$ is $\bigoplus_{k=0}^{n-1}L_{k,n-1}$, and thus we have
\begin{equation}\label{equ: Fil max N decomposition}
D^{\flat}=\bigoplus_{m=0}^{n-1}\bigoplus_{k=m}^{n-1}N^{k-m}(L_{m,n-1})=\bigoplus_{k=0}^{n-1}\bigoplus_{m=0}^{k}N^{k-m}(L_{m,n-1}).
\end{equation}
It follows from \ref{it: universal phi N 2} of Definition~\ref{def: universal phi N} that we have
\begin{equation}\label{equ: Fil max N L direct sum}
\mathrm{Fil}_{x_{k,n-1}}(\mathbf{E}_{k,n-1})=L_{k,n-1}\oplus\mathrm{Fil}_{x_{k+1,n-1}}(\mathbf{E}_{k,n-1}).
\end{equation}
Now that we have (\ref{equ: universal phi N sm cup}), we deduce from (\ref{equ: Fil max N L direct sum}) and a decreasing induction on $k$ that
\[\mathrm{Fil}_{x_{k,n-1}}(\mathbf{E}_{k,n-1})=\bigoplus_{m=0}^{k}N^{k-m}(L_{m,n-1})\]
for each $0\leq k\leq n-1$, which together with (\ref{equ: Fil max N decomposition}) implies (\ref{equ: Fil max sub phi N isom}).

We prove \ref{it: Fil max sub phi N 2}.\\
Now that we have
\[\mathrm{Fil}_{H}^1(D)\cong\sum_{\ell=1}^{n-1}\kappa^{\ell}(L_{\ell}\otimes_E\mathbf{E}_{\ell,n-1})\]
from \ref{it: universal phi N 3} of Definition~\ref{def: universal phi N}, to prove (\ref{equ: Fil max sub phi N intersection}) for each $1\leq \ell\leq n-1$, it suffices to show that the following intersection
\begin{equation}\label{equ: Fil max sub phi N emptyset}
\big(\sum_{\ell=1}^{n-2}\kappa^{\ell}(L_{\ell}\otimes_E\mathbf{E}_{\ell,n-1})\big)\cap D^{\flat}
\end{equation}
is zero. 
Recall that $L_{0,\ell}$ is the projection of $L_{\ell}$ to $\mathbf{E}_{0,\ell}=\mathbf{E}_{\Delta,J^{\ell}}$ for each $1\leq \ell\leq n-1$. 
Recall from (\ref{equ: auto total filtration}) and Corollary~\ref{cor: bottom Tits cup transfer} that $\mathbf{E}_{0,n-1}$ admits a decreasing filtration with graded pieces given by $\mathbf{E}_{0,1}\otimes_E\mathbf{E}_{1,n-1}$ and $\overline{\mathbf{E}}_{0,\ell}\otimes_E\mathbf{E}_{\ell,n-1}$ for each $2\leq \ell\leq n-1$. Hence, thanks to \ref{it: universal phi N top 2} of Lemma~\ref{lem: universal phi N top grade}, we know that the following sum of $E$-subspaces
\begin{equation}\label{equ: Fil max sub direct sum}
\sum_{\ell=1}^{n-1}L_{0,\ell}\otimes_E\mathbf{E}_{\ell,n-1}\subseteq \mathbf{E}_{0,n-1}
\end{equation}
is a direct sum, and in particular the following projection map (induced from the projection $\bigoplus_{k=0}^{n-1}\mathbf{E}_{k,n-1}\twoheadrightarrow\mathbf{E}_{0,n-1}$)
\begin{equation}\label{equ: Fil max sub projection}
\sum_{\ell=1}^{n-2}\kappa^{\ell}(L_{\ell}\otimes_E\mathbf{E}_{\ell,n-1})\rightarrow \sum_{\ell=1}^{n-2}L_{0,\ell}\otimes_E\mathbf{E}_{\ell,n-1}
\end{equation}
is an injection. We consider an arbitrary element $x$ in (\ref{equ: Fil max sub phi N emptyset}). Now that the projection of $x$ to $\mathbf{E}_{0,n-1}$ is contained in
\[(\sum_{\ell=1}^{n-2}L_{0,\ell}\otimes_E\mathbf{E}_{\ell,n-1})\cap L_{0,n-1}=0,\]
we conclude that $x=0$ as the projection map (\ref{equ: Fil max sub projection}) is injective. In other words, the intersection (\ref{equ: Fil max sub phi N emptyset}) is zero and thus the proof is finished.
\end{proof}

We give here a list of simple examples of filtered $(\varphi,N)$-modules that will be technically useful in later discussions.
\begin{ex}\label{ex: filtered phi N}
\begin{enumerate}[label=(\roman*)]
\item \label{it: filtered phi N 1} For each $i\in\Z$, we write $\langle i\rangle$ (resp.~$[i]$) for the unique $1$-dimensional filtered $(\varphi,N)$-module with $\varphi=p^{-i}$, $N=0$, and $\mathrm{gr}_{H}^{0}(\langle i\rangle)\neq 0$ (resp.~$\varphi=1$, $N=0$, and $\mathrm{gr}_{H}^{-i}([i])\neq 0$). For each filtered $(\varphi,N)$-module $D$, we write $D\langle i\rangle\defeq D\otimes_E\langle i\rangle$ and $D[i]\defeq D\otimes_E[i]$ for short.
\item \label{it: filtered phi N 2} For each $n\geq m\geq 2$, we consider the $2$-dimensional filtered $(\varphi,N)$-module $\overline{S}_{n,m}=Ee_0\oplus Ee_m$ with $\varphi(e_0)=e_0$, $\varphi(e_m)=p^me_m$, $N=0$, $\mathrm{Fil}_{H}^{0}(\overline{S}_{n,m})=\overline{D}_{n,m}$, $\mathrm{Fil}_{H}^{n+1}(\overline{S}_{n,m})=0$ and $\mathrm{Fil}_{H}^{\ell}(\overline{S}_n)=L$ for $1\leq \ell\leq n$, where $L\subseteq \overline{S}_{n,m}$ is a $E$-line different from both $Ee_0$ and $Ee_m$.
    Note that $\overline{S}_{n,m}$ fits into a strict exact sequence $0\rightarrow E \rightarrow \overline{S}_{n,m}\rightarrow \langle -m\rangle[n] \rightarrow 0$ of filtered $(\varphi,N)$-modules.
\item \label{it: filtered phi N 3} Let $D$ be a filtered $(\varphi,N)$-module which is universal special with level $1$, then there exists $e_1, e_0, e_0'\in D$ such that $D=Ee_1\oplus Ee_0\oplus Ee_0'$, $\varphi(e_1)=pe_1$, $\varphi(e_0)=e_0$, $\varphi(e_0')=e_0'$, $N(e_1)=e_0$, $N(e_0)=0=N(e_0')$, $\mathrm{Fil}_{H}^{0}(D)=D$, $\mathrm{Fil}_{H}^{2}(D)=0$, and $\mathrm{Fil}_{H}^{1}(D)=E(e_1+ae_0+be_0')$ for some $a\in E$ and $b\in E^\times$. Note that the isomorphism class of $D$ does not depend on the choice of such $e_1$, $e_0$, $e_0'$ and $\mathrm{Fil}_{H}^{1}(D)$.
\item \label{it: filtered phi N 4} For each $n\geq 0$, we write $S_n=\bigoplus_{k=0}^{n}Ee_k$ for the $(\varphi,N)$-module satisfying $\varphi(e_k)=p^ke_k$ for each $0\leq k\leq n$, $N(e_0)=0$ and $N(e_k)=e_{k-1}$ for each $1\leq k\leq n$. We equip $S_n$ with the filtration such that $\mathrm{gr}_{H}^{i}(S_{n})\neq 0$ if and only if $i=0$.
\item \label{it: filtered phi N 5} For each $n\geq m\geq 1$, we write $S_{n,m}=\bigoplus_{0\leq k\leq j\leq m}Ee_{j,k}$ for the $(\varphi,N)$-module satisfying $\varphi(e_{j,k})=p^ke_{j,k}$ for each $0\leq k\leq j\leq m$, $N(e_{j,0})=0$ for each $0\leq j\leq m$, $N(e_{j,k})=e_{j,k-1}$ for each $1\leq j\leq k\leq m$, $\mathrm{Fil}_{H}^{0}(S_{n,m})=S_{n,m}$, $\mathrm{Fil}_{H}^{n+1}(S_{n,m})=0$ and $\mathrm{Fil}_{H}^{\ell}(S_{n,m})=L$ for each $1\leq \ell\leq n$ with $L$ being an $E$-line in $S_{n,m}$ whose projection to $Ee_{j,j}$ is non-zero for each $0\leq j\leq m$. We write $S_{n,0}\defeq E$ for convenience.
\end{enumerate}
\end{ex}

By comparing the proof of Lemma~\ref{lem: Fil max sub phi N} and \ref{it: filtered phi N 5} of Example~\ref{ex: filtered phi N}, we may choose $e_{j,j}$ so that $L_{j,n-1}=Ee_{j,j}$ for each $0\leq j\leq n-1$ and thus conclude that
\begin{equation}\label{equ: Fil max basis}
D^{\flat}\cong S_{n-1,n-1}
\end{equation}
as filtered $(\varphi,N)$-modules.

\begin{lem}\label{lem: Fil max phi N induction}
Let $1\leq m\leq n$. We have the following results.
\begin{enumerate}[label=(\roman*)]
\item \label{it: Fil max phi N induction 1} For each $1\leq m\leq n$, there exists strict short exact sequence
\begin{equation}\label{equ: Fil max phi N induction 1}
0\rightarrow E^{\oplus m+1}\rightarrow S_{n,m}\rightarrow S_{n,m-1}\langle -1\rangle \rightarrow 0.
\end{equation}
\item \label{it: Fil max phi N induction 2} For each $1\leq m\leq n$, there exists a strict short exact sequence
\begin{equation}\label{equ: Fil max phi N induction 2}
0\rightarrow S_{m-1}^{\oplus 2}\oplus\bigoplus_{h=1}^{n-2}S_h \rightarrow S_{n,m} \rightarrow \overline{S}_{n,m} \rightarrow 0.
\end{equation}
\end{enumerate}
\end{lem}
\begin{proof}
Under the explicit basis in \ref{it: filtered phi N 5} of Example~\ref{ex: filtered phi N}, the short exact sequence (\ref{equ: Fil max phi N induction 1}) is given by
\begin{equation}\label{equ: Fil max phi N induction basis 1}
0\rightarrow \bigoplus_{j=0}^{m}Ee_{j,0}\rightarrow \bigoplus_{0\leq k\leq j\leq m}Ee_{j,k} \rightarrow \bigoplus_{1\leq k\leq j\leq m}Ee_{j,k} \rightarrow 0
\end{equation}
and (\ref{equ: Fil max phi N induction 2}) is given by
\begin{equation}\label{equ: Fil max phi N induction basis 2}
0\rightarrow (\bigoplus_{k=0}^{m-1}Ee_{m,k})\oplus(\bigoplus_{k=0}^{m-1}Ee_{m-1,k})\oplus \bigoplus_{h=1}^{n-2}(\bigoplus_{k=0}^{h}Ee_{h,k}) \rightarrow \bigoplus_{0\leq k\leq j\leq m}Ee_{j,k} \rightarrow Ee_{m,m}\oplus Ee_{0,0} \rightarrow 0.
\end{equation}
It is straightforward to check that both (\ref{equ: Fil max phi N induction basis 1}) and (\ref{equ: Fil max phi N induction basis 2}) are actually strict short exact sequence of filtered $(\varphi,N)$-modules.
\end{proof}

Let $D\cong \bigoplus_{k=0}^{n-1}\mathbf{E}_{k,n-1}$ be a filtered $(\varphi,N)$-module which is universal special with level $n-1$ (see Definition~\ref{def: universal phi N}). Let $0\leq k\leq n-1$ and $Q\subseteq \mathbf{E}_{k,n-1}$ be an $E$-vector subspace.
We consider the following $(\varphi,N)$-submodule of $D$
\begin{equation}\label{equ: cyclic submod}
D_{Q}\defeq \bigoplus_{s=0}^{k}\kappa_{s,n-1}^{k}(\mathbf{E}_{s,k}\otimes_E Q)\subseteq \bigoplus_{s=0}^{k}\mathbf{E}_{s,n-1}
\end{equation}
and endow it with the filtration induced from $\mathrm{Fil}_{H}^{\bullet}(D)$. More precisely, we have $\mathrm{Fil}_{H}^0(D_{Q})=D_{Q}$, $\mathrm{Fil}_{H}^{k+1}(D_{Q})=0$ and
\begin{equation}\label{equ: induced Hodge filtration}
\mathrm{Fil}_{H}^{m}(D_{Q})=\bigoplus_{\ell=0}^m\kappa^{\ell}(L_{\ell}\otimes_E(\mathbf{E}_{\ell,k}\otimes_E Q))\subseteq \bigoplus_{s=0}^{m}\kappa_{s,n-1}^{k}(\mathbf{E}_{s,k}\otimes_E Q) \subseteq D_{Q}
\end{equation}
for each $1\leq m\leq k$.
Let $0\leq \ell\leq k$ and $Q'\subseteq \mathbf{E}_{\ell,n-1}$ be an $E$-vector subspace.
It follows from (\ref{equ: cyclic submod}) that $D_{Q}\cap D_{Q'}\neq 0$ if and only if $Q'\cap D_{Q}\neq 0$ if and only if $Q'\cap \kappa_{\ell,n-1}^{k}(\mathbf{E}_{k,\ell}\otimes_E Q)\neq 0$.
Assume for the moment that $\ell=k$ and $Q\cap Q'=0$.
We observe from (\ref{equ: cyclic submod}) and (\ref{equ: induced Hodge filtration}) that $D_{Q\oplus Q'}\cong D_{Q}\oplus D_{Q'}$ as  filtered $(\varphi,N)$-submodules of $D$. Moreover, any isomorphism $Q\buildrel\sim\over\longrightarrow Q'$ would induce an isomorphism $D_{Q}\buildrel\sim\over\longrightarrow D_{Q'}$ between filtered $(\varphi,N)$-modules.
When $Q$ is an $E$-line, by comparing between (\ref{equ: induced Hodge filtration}) and Definition~\ref{def: universal phi N}, we see that $D_{Q}$ is universal special with level $k$.

\begin{lem}\label{lem: universal phi N induction}
Let $D\cong \bigoplus_{k=0}^{n-1}\mathbf{E}_{k,n-1}$ be a filtered $(\varphi,N)$-module which is universal special with level $n-1$ as in Definition~\ref{def: universal phi N}.
We have the following results.
\begin{enumerate}[label=(\roman*)]
\item \label{it: universal phi N induction 1} If we endow $\mathbf{E}_{0,n-1}$ (resp.~$D/\mathbf{E}_{0,n-1}$) with the filtration induced from $D$, then we have a strict exact sequence of filtered $(\varphi,N)$-modules
    \begin{equation}\label{equ: universal phi N induction}
    0\rightarrow \mathbf{E}_{0,n-1}\rightarrow D\rightarrow D/\mathbf{E}_{0,n-1}\rightarrow 0
    \end{equation}
    with $(D/\mathbf{E}_{0,n-1})\langle 1\rangle[1]$ being universal special with level $n-2$, and $\mathrm{gr}_{H}^{i}(\mathbf{E}_{0,n-1})\neq 0$ if and only if $i=0$.
\item \label{it: universal phi N induction 2} There exists a strict exact sequence of filtered $(\varphi,N)$-modules
\begin{equation}\label{equ: universal phi N induction embedding}
0\rightarrow D_{n-2}^{\oplus 2}\oplus\bigoplus_{k=1}^{n-3}D_k \rightarrow D \rightarrow \overline{S}_{n-1,n-1} \rightarrow 0
\end{equation}
with $D_k$ being universal special with level $k$ for each $1\leq k\leq n-2$, and $\overline{S}_{n-1,n-1}$ from \ref{it: filtered phi N 2} of Example~\ref{ex: filtered phi N}.
\end{enumerate}
\end{lem}
\begin{proof}
We prove \ref{it: universal phi N induction 1}.\\
It follows from $m=1$ case of (\ref{equ: complement of Hodge}) that $\mathrm{Fil}_{H}^1(D)\subseteq D$ and $\mathbf{E}_{0,n-1}\subseteq D$ induce an isomorphism between $E$-vector spaces
\begin{equation}\label{equ: first step complement}
\mathrm{Fil}^1(D)\oplus \mathbf{E}_{0,n-1}\buildrel\sim\over\longrightarrow D.
\end{equation}
In particular, the filtration $\mathrm{Fil}_{H}^{\bullet}(D)$ on $D$ induces the (unique) filtration on $\mathbf{E}_{0,n-1}$ that satisfies $\mathrm{gr}_{H}^i(\mathbf{E}_{0,n-1})\neq 0$ if and only if $i=0$. For each $1\leq \ell\leq n-1$, we use the shortened notation
\[\kappa_{\flat}^{\ell}\defeq \bigoplus_{k=1}^{\ell}\kappa_{k,n-1}^{\ell}: (\bigoplus_{k=1}^{\ell}\mathbf{E}_{k,\ell})\otimes_E\mathbf{E}_{\ell,n-1}\rightarrow \bigoplus_{k=1}^{\ell}\mathbf{E}_{k,n-1},\]
and write $L_{\ell,\flat}\subseteq \bigoplus_{k=1}^{\ell}\mathbf{E}_{k,\ell}$ for the $E$-line which is the image of $L_{\ell}$ under the projection $\bigoplus_{k=0}^{\ell}\mathbf{E}_{k,\ell}\twoheadrightarrow \bigoplus_{k=1}^{\ell}\mathbf{E}_{k,\ell}$. Thanks to our assumption on $L_{\ell}$, we observe that the composition of 
\[L_{\ell,\flat}\hookrightarrow \mathrm{Fil}_{x_{k,\ell}}(\mathbf{E}_{k,\ell})\twoheadrightarrow \mathrm{gr}_{x_{k,\ell}}(\mathbf{E}_{k,\ell})\]
is non-zero for each $1\leq k\leq \ell$. 
We observe that the filtration $\mathrm{Fil}_{H}^{\bullet}(D)$ on $D$ induces the filtration
\[\mathrm{Fil}_{H}^{\bullet}(D/\mathbf{E}_{0,n-1})\defeq (\mathrm{Fil}_{H}^{\bullet}(D)+\mathbf{E}_{0,n-1})/\mathbf{E}_{0,n-1}\]
on $D/\mathbf{E}_{0,n-1}$, which satisfies
\begin{equation}\label{equ: induced filtration decomposition}
\mathrm{Fil}_{H}^m(D/\mathbf{E}_{0,n-1})=\bigoplus_{\ell=m}^{n-1}\kappa_{\flat}^{\ell}(L_{\ell,\flat}\otimes_E\mathbf{E}_{\ell,n-1})
\end{equation}
for each $1\leq m\leq n-1$. Comparing Definition~\ref{def: universal phi N} with (\ref{equ: induced filtration decomposition}) and the discussion on $L_{\ell,\flat}$ above, we see that $(D/\mathbf{E}_{0,n-1})\langle 1\rangle[1]$ is universal special with level $n-2$.

We prove \ref{it: universal phi N induction 2}.\\
We choose arbitrary two $E$-lines $Q_{n-2},Q_{n-2}'\subseteq\mathbf{E}_{n-2,n-1}$ such that $\mathbf{E}_{n-2,n-1}\cong Q_{n-2}\oplus Q_{n-2}'$. We also choose an $E$-line $Q_k\subseteq \mathbf{E}_{k,n-1}$ not contained in $\mathbf{E}_{k,n-1}^{<}$ for each $1\leq k\leq n-3$. By the discussion below (\ref{equ: induced Hodge filtration}) we know that the filtered $(\varphi,N)$-modules $D_{Q_{n-2}}$ and $D_{Q_{n-2}'}$ are isomorphic and universal special with level $n-2$, and $D_{Q_k}$ is universal special with level $k$ for each $1\leq k\leq n-3$.
We consider
\begin{equation}\label{equ: universal phi N sub sum}
D_{Q_{n-2}}\oplus D_{Q_{n-2}'}\oplus\bigoplus_{k=1}^{n-3}D_{Q_k}
\end{equation}
which admits a map to $D$ as filtered $(\varphi,N)$-modules. Note that the $\varphi=p^{n-1}$-eigenspace of (\ref{equ: universal phi N sub sum}) is zero, and the $\varphi=p^{\ell}$-eigenspace of (\ref{equ: universal phi N sub sum}) is
\begin{equation}\label{equ: universal phi N sub sum eigen}
\kappa_{\ell,n-1}^{n-2}(\mathbf{E}_{\ell,n-2}\otimes_E(Q_{n-2}\oplus Q_{n-2}'))\oplus \bigoplus_{k=\ell}^{n-3}\kappa_{\ell,n-1}^{k}(\mathbf{E}_{\ell,k}\otimes_E Q_k) \subseteq \mathbf{E}_{\ell,n-1}
\end{equation}
for each $0\leq \ell\leq n-2$.
Recall from (\ref{equ: auto total filtration}) (and Corollary~\ref{cor: bottom Tits cup transfer}) that for each $0\leq \ell\leq n-1$, $\mathbf{E}_{\ell,n-1}$ admits a decreasing filtration
\[\mathrm{Fil}_{W}^{k}(\mathbf{E}_{\ell,n-1})=\sum_{k'=k+1}^{n-2}\mathrm{im}(\kappa_{\ell,n-1}^{k'})\]
with graded pieces given by $\mathrm{im}(\kappa_{\ell,n-1}^{n-2}$ and $\mathbf{E}_{\ell,k}\otimes_E(\overline{\mathbf{E}}_{k,n-1})$ for $\ell\leq k<n-2$.
This implies that the map from (\ref{equ: universal phi N sub sum}) to $D$ is an isomorphism on $\varphi=p^{\ell}$-eigenspaces for each $1\leq \ell\leq n-2$, and an embedding on $\varphi=p^{n-1}$ (resp.~$\varphi=1$) eigenspaces with cokernel $\mathbf{E}_{n-1,n-1}$ (resp.~cokernel $\overline{\mathbf{E}}_{0,n-1}=\mathbf{E}_{0,n-1}/\mathbf{E}_{0,n-1}^{<}$).
In other words, the map from (\ref{equ: universal phi N sub sum}) to $D$ has cokernel $\mathbf{E}_{n-1,n-1}\oplus \overline{\mathbf{E}}_{0,n-1}$.
It is easy to check from (\ref{equ: induced Hodge filtration}) and (\ref{equ: universal n step filtration}) that the direct sum filtration on (\ref{equ: universal phi N sub sum}) equals the filtration induced from $D$.
Note that $\kappa^{\ell}(L_{\ell}\otimes_E\mathbf{E}_{\ell,n-1})$ has non-zero image under
$D\twoheadrightarrow \mathbf{E}_{n-1,n-1}\oplus \overline{\mathbf{E}}_{0,n-1}$ if and only if $\ell=n-1$ in which case this image $\overline{L}_{n-1}\subseteq \mathbf{E}_{n-1,n-1}\oplus \overline{\mathbf{E}}_{0,n-1}$ of $L_{n-1}$ has non-zero projection to both $\mathbf{E}_{n-1,n-1}$ and $\overline{\mathbf{E}}_{0,n-1}$. Consequently, the quotient $\mathbf{E}_{n-1,n-1}\oplus \overline{\mathbf{E}}_{0,n-1}$, endowed with the filtration induced from $D$, is isomorphic to $\overline{S}_{n-1,n-1}$ from \ref{it: filtered phi N 2} of Example~\ref{ex: filtered phi N}. The proof is thus finished.
\end{proof}
\subsection{Universal extensions of $(\varphi,\Gamma)$-modules}\label{subsec: universal Galois}
We study various successive universal de Rham extensions of rank one de Rham $(\varphi,\Gamma)$-modules and the maps between them.

We write $K_0$ for the maximal unramified subfield of $K$.
Recall from \cite[\S~2.2]{Bre19} that a Deligne-Fontaine module (or $(\varphi,N)$-module in short) with trivial descent data is a triple $(D,\varphi,N)$ where $D$ is a finite free $K_0\otimes_{\Q_p}E$-module, $\varphi$ is a $K_0\otimes_{\Q_p}E$-semi-linear bijective map $D\rightarrow D$, and $N$ is a $K_0\otimes_{\Q_p}E$-linear endomorphism of $D$ satisfying $N\varphi=p\varphi N$.
A \emph{filtered $(\varphi,N)$-module with trivial descent data} is a quadruplet $(D,\varphi,N,\mathrm{Fil}^{\bullet}(D_{K}))$ with $\mathrm{Fil}^{\bullet}(D_{K})$ being a decreasing filtration of $K\otimes_{\Q_p}E$-submodules of $D_{K}\defeq D\otimes_{K_0}K$.
Under the decomposition $D_{K}\cong\prod_{\sigma: K\hookrightarrow E}D_{\sigma}$ with each $D_{\sigma}$ being an $E$-vector space of dimension the rank of $D$ over $K_0\otimes_{\Q_p}E$, $\mathrm{Fil}^{\bullet}(D_{K})$ is equivalent to a decreasing filtration $\mathrm{Fil}^{\bullet}(D_{\sigma})$ of $E$-subspaces for each $\sigma: K\hookrightarrow E$.
Given our fixed $\iota: K\hookrightarrow E$, we say that $(D,\varphi,N,\mathrm{Fil}^{\bullet}(D_{K}))$ is \emph{unramified outside $\iota$} if $\mathrm{Fil}^{0}(D_{\sigma})=D_{\sigma}$ and $\mathrm{Fil}^{-1}(D_{\sigma})=0$ for each $\sigma\neq \iota$, in which case we identify it with $(D,\varphi,N,\mathrm{Fil}^{\bullet}(D_{\iota}))$.
In other words, we may identify the abused filtered $(\varphi,N)$-modules as in \S~\ref{subsec: filtered phi N} with the usual filtered $(\varphi,N)$-modules which are unramifed outside $\iota$.

Let $\cR\defeq \cR_{K,E}$ be the $E$-coefficient Robba ring for $K$.
Given a $(\varphi,\Gamma)$-module $\cD$ over $\cR$, one could attach a \emph{$B$-pair} $(W_e(\cD),W_{\rm{dR}}^+(\cD))$ with $W_e(\cD)$ (resp.~$W_{\rm{dR}}^+(\cD)$) being a semi-linear representation of $\mathrm{Gal}_K$ over $B_e=B_{\rm{cris}}^{\varphi=1}$ (resp.~over $B_{\rm{dR}}^+\otimes_{\Q_p}E$) and this functor gives an equivalence between abelian categories (see \cite[Thm.~2.2.7]{Ber09}).
When $\cD$ is \emph{de Rham}, namely when $B_{\rm{dR}}\otimes_{B_e}W_e(\cD)$ is a trivial $B_{\rm{dR}}$-representation of $\mathrm{Gal}_K$, we could associate a filtered $(\varphi,N)$-module (with descent data) $(D_{\rm{pst}}(\cD), D_{\rm{dR}}(\cD))$ (cf.~\cite[\S 2.1]{Ding24}) and we often omit $D_{\rm{pst}}(\cD)$ from the notation. Recall from \cite[Prop.~2.3.4, Thm.~2.3.5]{Ber09} (see also \cite[Thm.~1.18 (2)]{Nak09}) that $\cD\mapsto D_{\rm{dR}}(\cD)$ induces an equivalence between de Rham $(\varphi,\Gamma)$-modules over $\cR$ and filtered $(\varphi,N)$-modules (with descent data).
Given an abused filtered $(\varphi,N)$-module $D$ as in \S~\ref{subsec: filtered phi N}, we can associate with it a de Rham $(\varphi,\Gamma)$-modules $\cD$ over $\cR$ via Berger's functor, upon identifying $D$ with a usual filtered $(\varphi,N)$-module which is unramifed outside $\iota$.

Given two $(\varphi,\Gamma)$-modules $\cD$ and $\cD'$ over $\cR$, we write $\mathrm{Ext}_{\varphi,\Gamma}^1(\cD,\cD')$ for their extensions in the abelian category of all $(\varphi,\Gamma)$-modules over $\cR$. If both $\cD$ and $\cD'$ are furthermore de Rham, we write $\mathrm{Ext}_g^1(\cD,\cD')$ for the subspace consisting of de Rham extensions.

We write $|\cdot|$ (resp.~$z$) for the continuous character $K^{\times}\rightarrow E^{\times}$ given by $x\mapsto |N_{K/\Q_p}(x)|_{p}$ (resp.~by $x\mapsto \iota(x)$). We define $\theta\defeq |\cdot|z: K^{\times}\rightarrow E^{\times}$ with $\theta(x)=|N_{K/\Q_p}(x)|_{p}\iota(x)$ for each $x\in K^\times$.

For each continuous character $\delta: K^{\times}\rightarrow E^{\times}$, we can associate a rank one $(\varphi,\Gamma)$-module $\cR(\delta)$ over $\cR$ (cf.~\cite[Notation~6.2.2]{KPX14}).
Comparing with \ref{it: filtered phi N 1} of Example~\ref{ex: filtered phi N}, we see that the de Rham rank one $(\varphi,\Gamma)$-module associated with $\langle i\rangle$ (resp.~$[i]$) is $\cR(|\cdot|^i)$ (resp.~$\cR(z^i)$) for each $i\in\Z$.
\begin{lem}\label{lem: char dR Ext}
Let $\delta_i=|\cdot|^{\ell_i}z^{k_i}$ for $i=1,2$. We have the following results.
\begin{enumerate}[label=(\roman*)]
\item \label{it: char dR Ext 1} We have $\Hom_{\varphi,\Gamma}(\cR(\delta_1),\cR(\delta_2))\neq 0$ if and only if $\ell_1=\ell_2$ and $k_1\geq k_2$, in which case it is $1$-dimensional.
\item \label{it: char dR Ext 2} If $\ell_2\neq \ell_1,\ell_1+1$, then $\mathrm{Ext}_{g}^1(\cR(\delta_1),\cR(\delta_2))$ is non-zero if and only if $k_1<k_2$, in which case it is $1$-dimensional.
\item \label{it: char dR Ext 3} If $\ell_2=\ell_1+1$, then $\mathrm{Ext}_{g}^1(\cR(\delta_1),\cR(\delta_2))$ is non-zero if and only if $k_1\leq k_2$. In this case, we have $\Dim_E\mathrm{Ext}_{g}^1(\cR(\delta_1),\cR(\delta_2))=1$ when $k_1=k_2$ and $\Dim_E\mathrm{Ext}_{g}^1(\cR(\delta_1),\cR(\delta_2))=2$ when $k_1<k_2$.
\item \label{it: char dR Ext 4} If $\ell_2=\ell_1$, then $\mathrm{Ext}_{g}^1(\cR(\delta_1),\cR(\delta_2))$ is non-zero if and only if $k_1\leq k_2$. In this case, we have $\Dim_E\mathrm{Ext}_{g}^1(\cR(\delta_1),\cR(\delta_2))=1$ when $k_1=k_2$ and $\Dim_E\mathrm{Ext}_{g}^1(\cR(\delta_1),\cR(\delta_2))=2$ when $k_1<k_2$, with the push-forward map with respect to decreasing choices of $k_2$ being surjective.
\end{enumerate}
\end{lem}
\begin{proof}
We notice that \ref{it: char dR Ext 1} follows from \cite[Prop.~2.14]{Nak09}, and that \ref{it: char dR Ext 2}, \ref{it: char dR Ext 3} and \ref{it: char dR Ext 4} are special cases of \cite[Lem.~4.2, Lem.~4.3]{Nak09}.
\end{proof}

\begin{lem}\label{lem: Galois simple Hom vanishing}
Let $\delta=|\cdot|^kz^{\ell}$ for some $k,\ell\in\Z$.
We have the following results.
\begin{enumerate}[label=(\roman*)]
\item \label{it: simple vanishing 1} Let $D$ be a filtered $(\varphi,N)$-module as in \ref{it: filtered phi N 3} of Example~\ref{ex: filtered phi N} and $\cD$ be its associated de Rham $(\varphi,\Gamma)$-module. We have
\begin{equation}\label{equ: simple vanishing 1}
\Hom_{\varphi,\Gamma}(\cD,\cR(\delta))=0
\end{equation}
when either $k\neq 0,-1$ or $k=0$ and $\ell>-1$.
\item \label{it: simple vanishing 2} Let $n\geq m\geq 3$ and $\overline{S}_{n-1,m-1}$ be the filtered $(\varphi,N)$-module from \ref{it: filtered phi N 2} of Example~\ref{ex: filtered phi N} and $\overline{\cS}_{n-1,m-1}$ be its associated de Rham $(\varphi,\Gamma)$-module. We have
\begin{equation}\label{equ: simple vanishing 2}
\Hom_{\varphi,\Gamma}(\overline{\cS}_{n-1,m-1},\cR(\delta))=0
\end{equation}
when $k\neq 1-m$ and $\ell>1-n$.
\item \label{it: simple vanishing 3} Let $n\geq 2$ and $\overline{S}_{n-1}$ be the filtered $(\varphi,N)$-module from \ref{it: filtered phi N 4} of Example~\ref{ex: filtered phi N} and $\overline{\cS}_{n-1}$ be its associated de Rham $(\varphi,\Gamma)$-module. We have
\begin{equation}\label{equ: simple vanishing 3}
\Hom_{\varphi,\Gamma}(\cS_{n-1},\cR(\delta))=0
\end{equation}
when $k\neq 1-n$.
\end{enumerate}
\end{lem}
\begin{proof}
We prove \ref{it: simple vanishing 1}.\\
Note that (\ref{equ: simple vanishing 1}) follow from \ref{it: char dR Ext 1} of Lemma~\ref{lem: char dR Ext} and a d\'evissage when $k\neq 0,-1$. Suppose that (\ref{equ: simple vanishing 1}) fails for $k=0$ and some $\ell>-1$, then there exists a non-zero map $q: D\rightarrow \langle 0\rangle[\ell]$ between filtered $(\varphi,N)$-modules. Since $\ell>-1$, we have $\mathrm{Fil}^1(\langle 0\rangle[\ell])=0$ and thus $\mathrm{Fil}^1(D)\subseteq \mathrm{ker}(q)$.
Recall from \ref{it: filtered phi N 3} of Example~\ref{ex: filtered phi N} that $\mathrm{Fil}^1(D)=E(e_1+ae_0+be_0')$ for some $a\in E$ and $b\in E^{\times}$. Since $\mathrm{ker}(q)$ is stable under $(\varphi,N)$-action, we have $e_1,ae_0+be_0'\in \mathrm{ker}(q)$ and thus $e_0=N(e_1)\in \mathrm{ker}(q)$. This forces $\mathrm{ker}(q)=D$ and thus contradicts the fact that $q$ is non-zero.

The proof of \ref{it: simple vanishing 2} is similar and simpler than that of \ref{it: simple vanishing 1} and thus is omitted.

Based on the notation of \ref{it: filtered phi N 4} of Example~\ref{ex: filtered phi N}, \ref{it: simple vanishing 3} follows from the fact that any non-zero map $S_{n-1}\rightarrow \langle k\rangle[\ell]$ with $k\neq 1-n$ must contain $e_{n-1}$ in its kernel, and that $S_{n-1}$ is generated by $e_{n-1}$ as a $(\varphi,N)$-module.
\end{proof}

\begin{rem}\label{rem: extra simple vanishing}
Let $S_{n-1,1}$ be a filtered $(\varphi,N)$-module as in \ref{it: filtered phi N 5} of Example~\ref{ex: filtered phi N} and $\cS_{n-1,1}$ be its associated de Rham $(\varphi,\Gamma)$-module.
Parallel argument as in the proof of \ref{it: simple vanishing 1} of Lemma~\ref{lem: Galois simple Hom vanishing} shows that
\[\Hom_{\varphi,\Gamma}(\cS_{n-1,1},\cR(\delta))=0\]
when either $k\neq 1-n,2-n$ or $k=2-n$ with $\ell>1-n$.
\end{rem}

\begin{lem}\label{lem: Galois Hom vanishing}
Let $\delta=|\cdot|^kz^{\ell}$ for some $k,\ell\in\Z$.
Let $D$ be a filtered $(\varphi,N)$-module which is universal special with level $n-1$ (see Definition~\ref{def: universal phi N}). Let $D^{\flat}$ be the minimal $(\varphi,N)$-submodule of $D$ containing $\mathrm{Fil}^{n-1}(D)$, equipped with the filtration induced from $D$ (see Lemma~\ref{lem: Fil max sub phi N}). Let $\cD$ and $\cD^{\flat}$ be the de Rham $(\varphi,\Gamma)$-module associated with $D$ and $D^{\flat}$ respectively. We have the following results.
\begin{enumerate}[label=(\roman*)]
\item \label{it: Galois Hom vanishing 1} We have
\begin{equation}\label{equ: Galois Hom vanishing 1}
\Hom_{\varphi,\Gamma}(\cD, \cR(\delta))=0
\end{equation}
when $k\neq 1-n$ and $\ell>1-n$.
\item \label{it: Galois Hom vanishing 2} Let $n\geq m\geq 3$. We have
\begin{equation}\label{equ: Galois Hom vanishing 2}
\Hom_{\varphi,\Gamma}(\cS_{n-1,m-1},\cR(\delta))=0
\end{equation}
when $k\neq 1-m$ and $\ell>1-n$. In particular, by taking $m=n$ we have
\begin{equation}\label{equ: Galois Hom vanishing 3}
\Hom_{\varphi,\Gamma}(\cD^{\flat}, \cR(\delta))=0
\end{equation}
when $k\neq 1-n$ and $\ell>1-n$.
\end{enumerate}
\end{lem}
\begin{proof}
Now that $\cD$ and $\cD^{\flat}$ admits a filtration with graded pieces given by $\cR(|\cdot|^{a}z^{b})$ for some $1-n\leq a\leq 0$ and $b\in\Z$, (\ref{equ: Galois Hom vanishing 1}) follows directly from \ref{it: char dR Ext 1} of Lemma~\ref{lem: char dR Ext} and a d\'evissage when either $k\leq -n$ or $k\geq 1$.
We assume in the rest of the proof that $2-n\leq k\leq 0$.
For each $i\geq j\geq 0$, we write $\cS_{i}$, $\cS_{i,j}$ and $\overline{\cS}_{i,j}$ for the de Rham $(\varphi,\Gamma)$-module associated with $S_{i}$, $S_{i,j}$ and $\overline{S}_{i,j}$ respectively (see Example~\ref{ex: filtered phi N}).

We prove \ref{it: Galois Hom vanishing 1} by an increasing induction on $n\geq 2$.\\
The $n=2$ case of (\ref{equ: Galois Hom vanishing 1}) follows from \ref{it: simple vanishing 1} of Lemma~\ref{lem: Galois simple Hom vanishing}.
Assume from now that $n\geq 3$.
We deduce from \ref{it: universal phi N induction 1} of Lemma~\ref{lem: universal phi N induction} a short exact sequence
\begin{equation}\label{equ: universal phi gamma induction 1}
0\rightarrow \mathbf{E}_{0,n-1}\otimes_E \cR(1_{K^{\times}})\rightarrow \cD \rightarrow \cD'\otimes_{\cR}\cR(\theta^{-1})\rightarrow 0
\end{equation}
with $\cD'$ being the de Rham $(\varphi,\Gamma)$-module associated with a filtered $(\varphi,N)$-module which is universal special with level $n-2$.
We also deduce from \ref{it: universal phi N induction 2} of Lemma~\ref{lem: universal phi N induction} a short exact sequence
\begin{equation}\label{equ: universal phi gamma induction 2}
0\rightarrow (\cD_{n-2}')^{\oplus 2}\oplus \bigoplus_{m=1}^{n-3}\cD_{m}'\rightarrow \cD \rightarrow \overline{\cS}_{n-1,n-1}\rightarrow 0
\end{equation}
with $\cD_{m}'$ being the de Rham $(\varphi,\Gamma)$-module associated with a filtered $(\varphi,N)$-module which is universal special with level $m$, for each $1\leq m\leq n-2$.
If $k<0$, then we have $\Hom_{\varphi,\Gamma}(\cD',\cR(|\cdot|^{k+1}z^{\ell+1}))=0$ by induction hypothesis and thus
\[\Hom_{\varphi,\Gamma}(\cD'\otimes_{\cR}\cR(\theta^{-1}),\cR(\delta))=0,\]
which together with (\ref{equ: universal phi gamma induction 1}) as well as $\Hom_{\varphi,\Gamma}(\cR(1_{K^\times}),\cR(\delta))=0$ from \ref{it: char dR Ext 1} of Lemma~\ref{lem: char dR Ext} gives (\ref{equ: Galois Hom vanishing 1}).
If $k=0$, then we have $\Hom_{\varphi,\Gamma}(\cD_{m}',\cR(\delta))=0$ for each $1\leq m\leq n-2$ by induction hypothesis, which together with (\ref{equ: universal phi gamma induction 2}) as well as $\Hom_{\varphi,\Gamma}(\overline{\cS}_{n-1,n-1},\cR(\delta))=0$ from \ref{it: simple vanishing 2} of Lemma~\ref{lem: Galois simple Hom vanishing} gives (\ref{equ: Galois Hom vanishing 1}). This finishes the proof of the induction step and thus of \ref{it: Galois Hom vanishing 1}.

We fix $n$ and prove \ref{it: Galois Hom vanishing 2} by an increasing induction on $1\leq m\leq n$.\\
The $m=1$ case of (\ref{equ: Galois Hom vanishing 2}) follows from \ref{it: char dR Ext 1} of Lemma~\ref{lem: char dR Ext} and the $m=2$ case of (\ref{equ: Galois Hom vanishing 2}) follows from Remark~\ref{rem: extra simple vanishing}.
Assume from now that $m\geq 3$.
We deduce from \ref{it: Fil max phi N induction 1} of Lemma~\ref{lem: Fil max phi N induction} a short exact sequence
\begin{equation}\label{equ: Fil max phi gamma induction 1}
0\rightarrow \cR(1_{K^{\times}})^{\oplus m}\rightarrow \cS_{n-1,m-1} \rightarrow \cS_{n-1,m-2}\otimes_{\cR}\cR(|\cdot|^{-1})\rightarrow 0.
\end{equation}
We also deduce from \ref{it: Fil max phi N induction 2} of Lemma~\ref{lem: Fil max phi N induction} a short exact sequence
\begin{equation}\label{equ: Fil max phi gamma induction 2}
0\rightarrow (\cS_{n-2})^{\oplus 2}\oplus \bigoplus_{m=1}^{n-3}\cS_{m}\rightarrow \cS_{n-1,m-1} \rightarrow \overline{\cS}_{n-1,m-1}\rightarrow 0.
\end{equation}
If $k<0$, then we have $\Hom_{\varphi,\Gamma}(\cS_{n-1,m-2},\cR(|\cdot|^{k+1}z^{\ell}))=0$ by induction hypothesis and thus
\[\Hom_{\varphi,\Gamma}(\cS_{n-1,m-2}\otimes_{\cR}\cR(|\cdot|^{-1}),\cR(\delta))=0,\]
which together with (\ref{equ: Fil max phi gamma induction 1}) as well as $\Hom_{\varphi,\Gamma}(\cR(1_{K^\times}),\cR(\delta))=0$ from \ref{it: char dR Ext 1} of Lemma~\ref{lem: char dR Ext} gives (\ref{equ: Galois Hom vanishing 2}).
If $k=0$, then we have $\Hom_{\varphi,\Gamma}(\cS_{m},\cR(z^{\ell}))=0$ for each $1\leq m\leq n-2$ by \ref{it: simple vanishing 3} of Lemma~\ref{lem: Galois simple Hom vanishing}, which together with (\ref{equ: Fil max phi gamma induction 2}) as well as $\Hom_{\varphi,\Gamma}(\overline{\cS}_{n-1,m-2},\cR(z^{\ell}))=0$ from \ref{it: simple vanishing 2} of Lemma~\ref{lem: Galois simple Hom vanishing} gives (\ref{equ: Galois Hom vanishing 2}). This finishes the proof of the induction step and thus of \ref{it: Galois Hom vanishing 2}.
\end{proof}

We define $\cE_0\defeq E$ and define $\cD_0\defeq \cR(1_{K^\times})$.
For each $n\geq 1$, we inductively define $\cE_{n}\defeq \mathrm{Ext}_g^1(\cD_{n-1},\cR(\theta^{n}))^\vee$ and then $\cD_{n}$ (for increasing $n$) as the following universal de Rham extension
\begin{equation}\label{equ: universal Galois Ext}
0\rightarrow \cE_{n}\otimes_E\cR(\theta^{n})\rightarrow \cD_{n}\rightarrow \cD_{n-1}\rightarrow 0.
\end{equation}
Our inductive definition of $\cE_{n}$ and $\cD_{n}$ naturally endow $\cD_{n}$ with a decreasing filtration
\begin{equation}\label{equ: universal Galois filtration}
\cD_{n}=\mathrm{Fil}_{W}^{0}(\cD_{n})\supseteq \mathrm{Fil}_{W}^{1}(\cD_{n})\supseteq \cdots \supseteq \mathrm{Fil}_{W}^{n}(\cD_{n})\supseteq \mathrm{Fil}_{W}^{n+1}(\cD_{n})=0
\end{equation}
with $\cD_{n}/\mathrm{Fil}_{W}^{k+1}(\cD_{n})\cong \cD_{k}$ for each $0\leq k\leq n$. Note that we have
\[\mathrm{gr}_{W}^0(\cD_{n})\defeq \mathrm{Fil}_{W}^0(\cD_{n})/\mathrm{Fil}_{W}^{1}(\cD_{n})\cong \cR(1_{K^\times})\]
and
\begin{equation}\label{equ: universal Galois graded}
\mathrm{gr}_{W}^k(\cD_{n})\defeq \mathrm{Fil}_{W}^k(\cD_{n})/\mathrm{Fil}_{W}^{k+1}(\cD_{n})\cong \cE_{k}\otimes_E \cR(\theta^k)=\mathrm{Ext}_g^1(\cD_{k-1},\cR(\theta^k))^\vee\otimes_E \cR(\theta^k)
\end{equation}
for each $1\leq k\leq n$.

We define $\cE_{n,0}^{\flat}\defeq E$, $\cD_{n,0}^{\flat}\defeq \cR(1_{K^\times})$.
Then we inductively define for each $1\leq m\leq n$ that $\cE_{n,m}^{\flat}\defeq \mathrm{Ext}_g^1(\cD_{n,m-1}^{\flat},\cR(|\cdot|^{m}z^{n}))^\vee$ and then $\cD_{n,m}^{\flat}$ (for increasing $m$) as the following universal de Rham extension
\begin{equation}\label{equ: univsersal Galois Ext}
0\rightarrow \cE_{n,m}\otimes_E\cR(|\cdot|^{m}z^{n})\rightarrow \cD_{n,m}^{\flat}\rightarrow \cD_{n,m-1}^{\flat}\rightarrow 0.
\end{equation}
We write $\cE_{n}^{\flat}\defeq \cE_{n,n}^{\flat}$ and $\cD_{n}^{\flat}\defeq \cD_{n,n}^{\flat}$ for short.

\begin{prop}\label{prop: Galois dR Ext exact}
Let $n\geq 1$. We have the following results.
\begin{enumerate}[label=(\roman*)]
\item \label{it: Galois dR Ext exact 1} Let $D$ be a filtered $(\varphi,N)$-module which is universal special with level $n-1$ (see Definition~\ref{def: universal phi N}) and $\cD$ its associated de Rham $(\varphi,\Gamma)$-module. Then we have
    \begin{equation}\label{equ: universal special Ext}
    \cD\cong \cD_{n-1}\otimes_{\cR}\cR(\theta^{1-n})
    \end{equation}
    and
    \begin{equation}\label{equ: universal special Ext dim}
    \Dim_E\mathbf{E}_{n-1}=\Dim_E \cE_{n-1}.
    \end{equation}
\item \label{it: Galois dR Ext exact 2} Let $1\leq m\leq n$ and $\cS_{n-1,m-1}$ be the de Rham $(\varphi,\Gamma)$-module associated with $S_{n-1,m-1}$ (see \ref{it: filtered phi N 5} of Example~\ref{ex: filtered phi N}). Then we have
    \begin{equation}\label{equ: Fil max special Ext}
    \cS_{n-1,m-1}\cong \cD_{n-1,m-1}^{\flat}\otimes_{\cR}\cR(|\cdot|^{1-m}z^{1-n})
    \end{equation}
    and
    \begin{equation}\label{equ: Fil max special Ext dim}
    \Dim_E \cE_{n-1,m-1}^{\flat}=m.
    \end{equation}
\end{enumerate}
\end{prop}
\begin{proof}
We prove \ref{it: Galois dR Ext exact 1} by an increasing induction on $n$.\\
The $n=1$ case is obvious. We assume from now that $n\geq 2$.
It follows from \ref{it: universal phi N induction 1} of Lemma~\ref{lem: universal phi N induction} that $\cD$ fits into an exact sequence
\[0\rightarrow \mathbf{E}_{n-1}\otimes_E \cR(1_{K^\times})\rightarrow \cD \rightarrow \cD'\otimes_{\cR}\cR(\theta^{-1})\rightarrow 0\]
with $\cD'$ being the de Rham $(\varphi,\Gamma)$-module associated with a filtered $(\varphi,N)$-module which is universal special with level $n-2$. By our induction hypothesis we know that $\cD'\cong \cD_{n-2}\otimes_{\cR}\cR(\theta^{2-n})$, and thus $\cD\otimes_{\cR}\cR(\theta^{n-1})$ fits into an exact sequence
\begin{equation}\label{equ: Galois dR Ext exact 1 seq}
0\rightarrow \mathbf{E}_{n-1}\otimes_E \cR(\theta^{n-1})\rightarrow \cD\otimes_{\cR}\cR(\theta^{n-1}) \rightarrow \cD_{n-2}\rightarrow 0
\end{equation}
It follows from \ref{it: Galois Hom vanishing 1} of Lemma~\ref{lem: Galois Hom vanishing} that $\Hom_{\varphi,\Gamma}(\cD,\cR(1_{K^\times}))=0$ and thus $\Hom_{\varphi,\Gamma}(\cD\otimes_{\cR}\cR(\theta^{n-1}),\cR(\theta^{n-1}))=0$, which together with the definition of $\cD_{n-1}$ gives a surjection $\cE_{n-1}\twoheadrightarrow \mathbf{E}_{n-1}$ and thus a surjection
\begin{equation}\label{equ: Galois dR Ext surjection 1}
\cD_{n-1}\twoheadrightarrow \cD\otimes_{\cR}\cR(\theta^{n-1}).
\end{equation}
In particular, we deduce that
\begin{equation}\label{equ: Galois dR Ext dim lower bound 1}
\Dim_E \cE_{n-1}\geq \Dim_E \mathbf{E}_{n-1}.
\end{equation}
Note that (\ref{equ: Galois dR Ext dim lower bound 1}) is known to be an equality (see \ref{it: char dR Ext 3} of Lemma~\ref{lem: char dR Ext}) and thus (\ref{equ: Galois dR Ext surjection 1}) is an isomorphism when $n=2$. Assume from now that $n\geq 3$.
For each $2\leq m\leq n-1$, the short exact sequence $0\rightarrow \cE_{m-1}\otimes_E\cR(\theta)^{m-1}\rightarrow \cD_{m-1}\rightarrow \cD_{m-2}\rightarrow 0$ together with $\Hom_{\varphi,\Gamma}(\cR(\theta)^{m-1},\cR(\theta)^{n-1})=0$ from \ref{it: char dR Ext 1} of Lemma~\ref{lem: char dR Ext} gives an exact sequence
\[0\rightarrow \mathrm{Ext}_{g}^1(\cD_{m-2},\cR(\theta^{n-1}))\rightarrow \mathrm{Ext}_{g}^1(\cD_{m-1},\cR(\theta^{n-1}))\rightarrow \cE_{m-1}^{\vee}\otimes_E \mathrm{Ext}_{g}^1(\cR(\theta^{m-1}),\cR(\theta^{n-1}))\]
and thus
\begin{equation}\label{equ: Galois dR Ext exact 1 dim induction}
\Dim_E \mathrm{Ext}_{g}^1(\cD_{m-1},\cR(\theta^{n-1}))\leq \Dim_E \mathrm{Ext}_{g}^1(\cD_{m-2},\cR(\theta^{n-1}))+\Dim_E \cE_{m-1} \Dim_E \mathrm{Ext}_{g}^1(\cR(\theta^{m-1}),\cR(\theta^{n-1})).
\end{equation}
Now that we have $\Dim_E \mathrm{Ext}_{g}^1(\cR(\theta^{n-2}),\cR(\theta^{n-1}))=2$ and $\Dim_E \mathrm{Ext}_{g}^1(\cR(\theta^{m-1}),\cR(\theta^{n-1}))=1$ for each $1\leq m \leq n-2$, we deduce from (\ref{equ: Galois dR Ext exact 1 dim induction}) and an induction on $2\leq m\leq n-1$ that
\begin{equation}\label{equ: Galois dR Ext dim upper bound 1}
\Dim_E\cE_{n-1}=\Dim_E \mathrm{Ext}_{g}^1(\cD_{n-2},\cR(\theta^{n-1}))\leq 2\Dim_E \cE_{n-2}+\sum_{m=1}^{n-2}\Dim_E \cE_{m-1}.
\end{equation}
Now that we have
\[\Dim_E \mathbf{E}_{n-1}=2\Dim_E \mathbf{E}_{n-2}+\sum_{m=1}^{n-2}\Dim_E \mathbf{E}_{m-1}\]
from (\ref{equ: total dim induction}), we deduce from (\ref{equ: Galois dR Ext dim lower bound 1}) and (\ref{equ: Galois dR Ext dim upper bound 1}) that the inequality (\ref{equ: Galois dR Ext dim lower bound 1}) must be an equality, and in particular the surjection (\ref{equ: Galois dR Ext surjection 1}) is an isomorphism.

We fix $n$ and prove \ref{it: Galois dR Ext exact 2} by an increasing induction on $1\leq m\leq n$.\\
The $m=1$ case is obvious. We assume from now that $m\geq 2$.
It follows from \ref{it: Fil max phi N induction 1} of Lemma~\ref{lem: Fil max phi N induction} that $\cS_{n-1,m-1}$ fits into an exact sequence
\[0\rightarrow \cR(1_{K^\times})^{\oplus m}\rightarrow \cS_{n-1,m-1} \rightarrow \cS_{n-1,m-2}\otimes_{\cR}\cR(|\cdot|^{-1})\rightarrow 0.\]
By our induction hypothesis we know that $\cS_{n-1,m-2}\cong \cD_{n-1,m-2}^{\flat}\otimes_{\cR}\cR(|\cdot|^{2-m}z^{1-n})$, and thus $\cS_{n-1,m-1}\otimes_{\cR}\cR(|\cdot|^{m-1}z^{n-1})$ fits into an exact sequence
\begin{equation}\label{equ: Galois dR Ext exact 2 seq}
0\rightarrow \cR(|\cdot|^{m-1}z^{n-1})^{\oplus m} \rightarrow \cS_{n-1,m-1}\otimes_{\cR}\cR(|\cdot|^{m-1}z^{n-1}) \rightarrow \cD_{n-1,m-2}^{\flat}\rightarrow 0
\end{equation}
It follows from \ref{it: Galois Hom vanishing 1} of Lemma~\ref{lem: Galois Hom vanishing} that $\Hom_{\varphi,\Gamma}(\cS_{n-1,m-1},\cR(1_{K^\times}))=0$ and thus $\Hom_{\varphi,\Gamma}(\cS_{n-1,m-1}\otimes_{\cR}\cR(|\cdot|^{m-1}z^{n-1}),\cR(|\cdot|^{m-1}z^{n-1}))=0$, which together with the definition of $\cD_{n-1,m-1}^{\flat}$ gives a surjection
\begin{equation}\label{equ: Galois dR Ext surjection 2}
\cD_{n-1,m-1}^{\flat}\twoheadrightarrow \cS_{n-1,m-1}\otimes_{\cR}\cR(|\cdot|^{m-1}z^{n-1})
\end{equation}
with
\begin{equation}\label{equ: Galois dR Ext dim lower bound 2}
\Dim_E \cE_{n-1,m-1}^{\flat}\geq m.
\end{equation}
Note that (\ref{equ: Galois dR Ext dim lower bound 2}) is known to be an equality (see \ref{it: char dR Ext 3} of Lemma~\ref{lem: char dR Ext}) and thus (\ref{equ: Galois dR Ext surjection 2}) is an isomorphism when $m=2$.
Assume from now that $m\geq 2$.
For each $2\leq \ell\leq m-1$, the short exact sequence $0\rightarrow \cE_{n-1,\ell-1}^{\flat}\otimes_E\cR(|\cdot|^{\ell-1}z^{n-1})\rightarrow \cD_{n-1,\ell-1}^{\flat}\rightarrow \cD_{n-1,\ell-2}^{\flat}\rightarrow 0$ together with $\Hom_{\varphi,\Gamma}(\cR(|\cdot|^{\ell-1}z^{n-1}),\cR(|\cdot|^{m-1}z^{n-1}))=0$ from \ref{it: char dR Ext 1} of Lemma~\ref{lem: char dR Ext} gives an exact sequence
\begin{multline*}
0\rightarrow \mathrm{Ext}_{g}^1(\cD_{n-1,\ell-2}^{\flat},\cR(|\cdot|^{m-1}z^{n-1}))\rightarrow \mathrm{Ext}_{g}^1(\cD_{n-1,\ell-1}^{\flat},\cR(|\cdot|^{m-1}z^{n-1}))\\
\rightarrow \cE_{n-1,\ell-1}^{\vee}\otimes_E \mathrm{Ext}_{g}^1(\cR(|\cdot|^{\ell-1}z^{n-1}),\cR(|\cdot|^{m-1}z^{n-1})).
\end{multline*}
Now that we have $\Dim_E \mathrm{Ext}_{g}^1(\cR(|\cdot|^{m-2}z^{n-1}),\cR(|\cdot|^{m-1}z^{n-1}))=1$, $\mathrm{Ext}_{g}^1(\cR(|\cdot|^{\ell-1}z^{n-1}),\cR(|\cdot|^{m-1}z^{n-1}))=0$ for each $2\leq \ell\leq m-2$, and $\Dim_E \mathrm{Ext}_{g}^1(\cR(1_{K^\times}),\cR(|\cdot|^{m-1}z^{n-1}))=1$ from \ref{it: char dR Ext 2} of Lemma~\ref{lem: char dR Ext}, an induction on $2\leq \ell\leq m-1$ gives
\[\Dim_E\cE_{n-1,m-1}^{\flat}=\Dim_E \mathrm{Ext}_{g}^1(\cD_{n-1,m-2},\cR(|\cdot|^{m-1}z^{n-1}))\leq m.\]
Consequently, the inequality (\ref{equ: Galois dR Ext dim lower bound 2}) must be an equality, and in particular the surjection (\ref{equ: Galois dR Ext surjection 2}) is an isomorphism.
\end{proof}

Recall from the proof of \ref{it: Galois dR Ext exact 1} of Proposition~\ref{prop: Galois dR Ext exact} that the short exact sequence
\[0\rightarrow \cR(\theta^{k})\otimes_E\cE_{k}\rightarrow \cD_{k}\rightarrow \cD_{k-1}\rightarrow 0\]
induces a short exact sequence
\begin{equation}\label{equ: universal Galois Ext devissage}
0\rightarrow \mathrm{Ext}_g^1(\cR(\theta^{k}),\cR(\theta^{n}))^\vee\otimes_E\cE_{k}\rightarrow \mathrm{Ext}_g^1(\cD_{k},\cR(\theta^{n}))^\vee \rightarrow \mathrm{Ext}_g^1(\cD_{k-1},\cR(\theta^{n}))^\vee \rightarrow 0
\end{equation}
for each $1\leq k\leq n-1$.
In particular, the filtration (\ref{equ: universal Galois filtration}) induces a decreasing filtration
\[
\mathrm{Fil}_{W}^k(\cE_{n})\defeq \mathrm{Ext}_g^1(\mathrm{Fil}_{W}^k(\cD_{n-1}),\cR(\theta^{n}))^\vee\subseteq \cE_{n}
\]
with
\begin{equation}\label{equ: universal Galois Ext graded}
\mathrm{gr}_{W}^k(\cE_{n})=\cE_{k}\otimes_E \mathrm{Ext}_g^1(\cR(\theta^{k}),\cR(\theta^{n}))^\vee
\end{equation}
for each $0\leq k\leq n-1$.

\begin{defn}\label{def: special mod}
Let $\cD$ be a de Rham $(\varphi,\Gamma)$-module.
\begin{enumerate}[label=(\roman*)]
\item \label{it: special mod 1} We say that $\cD$ is \emph{weakly special with level $n$} for some $n\geq 0$ if it admits a decreasing filtration $\mathrm{Fil}_{W}^{\bullet}(\cD)$ such that $\mathrm{Fil}_{W}^{0}(\cD)=\cD$, $\mathrm{Fil}_{W}^{n+1}(\cD)=0$, $\mathrm{gr}_{W}^{0}(\cD)=\cR(1_{K^\times})$ and $\mathrm{gr}_{W}^{k}(\cD)=\cR(|\cdot|^kz^{\ell_k})^{\oplus m_k}$ for some $\ell_k\in \Z$, $m_k\geq 0$ and each $1\leq k\leq n$.
\item \label{it: special mod 2} We say that $\cD$ is \emph{special with level $n$} for some $n\geq 0$ if there exists a decreasing filtration $\mathrm{Fil}_{W}^{\bullet}(\cD)$ such that $\mathrm{Fil}_{W}^{0}(\cD)=\cD$, $\mathrm{Fil}_{W}^{n+1}(\cD)=0$ and $\mathrm{gr}_{W}^{k}(\cD)=\cR(\theta^k)$ for each $1\leq k\leq n$, and moreover $\Hom_{\varphi,\Gamma}(\mathrm{Fil}_{W}^{k-1}(\cD)/\mathrm{Fil}_{W}^{k+1}(\cD),\cR(\theta^{k}))=0$ for each $1\leq k\leq n$.
\end{enumerate}
\end{defn}

\begin{lem}\label{lem: Galois unique Hom}
Let $\cD$ and $\cD'$ be two weakly special de Rham $(\varphi,\Gamma)$-modules with level $n$ for some $n\geq 0$ (see \ref{it: special mod 1} of Definition~\ref{def: special mod}) with $\ell_{k}$ and $m_{k}$ associated with $\cD$ for each $1\leq k\leq n$. Assume that
\begin{equation}\label{equ: Galois Hom vanishing cond}
\Hom_{\varphi,\Gamma}(\cD', \cR(|\cdot|^kz^{\ell_{k}}))=0
\end{equation}
for each $1\leq k\leq n$. Then we have
\begin{equation}\label{equ: Galois Hom dim bound}
\Dim_E \Hom_{\varphi,\Gamma}(\cD',\cD)\leq 1.
\end{equation}
\end{lem}
\begin{proof}
Recall from \ref{it: char dR Ext 1} of Lemma~\ref{lem: char dR Ext} that we have $\Hom_{\varphi,\Gamma}(\cR(|\cdot|^{k}z^{\ell}),\cR(1_{K^\times}))=0$ for each $k,\ell\in \Z$ with $k\neq 0$, which together with $\Dim_E \Hom_{\varphi,\Gamma}(\cR(1_{K^\times}),\cR(1_{K^\times}))=1$ and a d\'evissage with respect to $\mathrm{Fil}_{W}^{\bullet}(\cD')$ from Definition~\ref{def: special mod} (for $\cD'$) gives
\[\Dim_E \Hom_{\varphi,\Gamma}(\cD', \cR(1_{K^\times}))=1.\]
This together with (\ref{equ: Galois Hom vanishing cond}) and a further d\'evissage with respect to $\mathrm{Fil}_{W}^{\bullet}(\cD)$ from Definition~\ref{def: special mod} (for $\cD$) gives (\ref{equ: Galois Hom dim bound}).
\end{proof}

\begin{lem}\label{lem: Galois universal Hom}
Let $\cD$ be a weakly special de Rham $(\varphi,\Gamma)$-module with level $n$ for some $n\geq 0$ and $\ell_k=k$ for each $1\leq k\leq n$ (see \ref{it: special mod 1} of Definition~\ref{def: special mod}). Assume that 
\begin{equation}\label{equ: Galois universal Hom vanishing cond}
\Hom_{\varphi,\Gamma}(\cD, \cR(\theta^{k}))=0
\end{equation}
for each $1\leq k\leq n$. Then we have 
\begin{equation}\label{equ: Galois universal Hom dim}
\Dim_E \Hom_{\varphi,\Gamma}(\cD_{n},\cD)=1
\end{equation}
and the unique (up to scalars) non-zero map $\cD_{n}\rightarrow \cD$ is a surjection.
\end{lem}
\begin{proof}
We prove our claim by an increasing induction on $n\geq 0$. When $n=0$, we have $\cD_{n}=\cR(1_{K^{\times}})=\cD$ and our claim is clear. Assume from now that $n\geq 1$. Our induction hypothesis implies that 
\[
\Dim_E \Hom_{\varphi,\Gamma}(\cD_{n-1},\cD/\mathrm{Fil}_{W}^{n}(\cD))=1
\]
and the unique (up to scalars) non-zero map $\cD_{n-1}\rightarrow \cD/\mathrm{Fil}_{W}^{n}(\cD)$ is a surjection.
The non-zero map $\cD_{n-1}\rightarrow \cD/\mathrm{Fil}_{W}^{n}(\cD)$ induces a surjection
\[
\cE_{n}=\mathrm{Ext}_{g}^{1}(\cD_{n-1},\cR(\theta^{n}))^{\vee}\twoheadrightarrow \cE\defeq \mathrm{Ext}_{g}^{1}(\cD/\mathrm{Fil}_{W}^{n}(\cD),\cR(\theta^{n}))^{\vee}
\]
and thus the following commutative diagram
\begin{equation}\label{equ: Galois universal Hom diagram}
\xymatrix{
\cE_{n}\otimes_{E}\cR(\theta^{n}) \ar@{^{(}->}[r] \ar@{->>}[d] & \cD_{n} \ar@{->>}[r] \ar[d] & \cD_{n-1} \ar@{->>}[d]\\
\cE\otimes_{E}\cR(\theta^{n}) \ar@{^{(}->}[r] \ar[d] & \cD' \ar@{->>}[r] \ar[d] & \cD/\mathrm{Fil}_{W}^{n}(\cD) \ar@{=}[d]\\
\mathrm{Fil}_{W}^{n}(\cD) \ar@{^{(}->}[r] & \cD \ar@{->>}[r] & \cD/\mathrm{Fil}_{W}^{n}(\cD)
}.
\end{equation}
with the middle row of (\ref{equ: Galois universal Hom diagram}) being the universal extension. Our assumption $\Hom_{\varphi,\Gamma}(\cD, \cR(\theta^{n}))=0$ forces the left bottom vertical map $\cE\otimes_{E}\cR(\theta^{n})\rightarrow \mathrm{Fil}_{W}^{n}(\cD)$ to be a surjection. Now that all rows of (\ref{equ: Galois universal Hom diagram}) are short exact, and all maps from the left column and the right column of (\ref{equ: Galois universal Hom diagram}) are surjective, we conclude that all maps from the middle column of (\ref{equ: Galois universal Hom diagram}) are also surjective, from which we obtain a surjection $\cD_{n}\twoheadrightarrow \cD$.
This together with Lemma~\ref{lem: Galois unique Hom} (upon replacing $\cD'$ in \emph{loc.cit.} with $\cD_{n}$) finishes the proof.
\end{proof}

\begin{lem}\label{lem: Galois sub flat}
Let $n\geq m\geq 1$. We have the following results.
\begin{enumerate}[label=(\roman*)]
\item \label{it: Galois sub flat 1} We have $\Dim_E\Hom_{\varphi,\Gamma}(\cD_{n-1,m-1}^{\flat},\cD_{m-1})=1$ and the unique (up to scalars) non-zero map $\cD_{n-1,m-1}^{\flat}\rightarrow \cD_{m-1}$ is an embedding.
\item \label{it: Galois sub flat 2} For each $\ell\geq n$, we have $\Dim_E\Hom_{\varphi,\Gamma}(\cD_{\ell-1,m-1}^{\flat},\cD_{n-1,m-1}^{\flat})=1$ and the unique (up to scalars) non-zero map $\cD_{\ell-1,m-1}^{\flat}\rightarrow \cD_{n-1,m-1}^{\flat}$ is an embedding.
\end{enumerate}
\end{lem}
\begin{proof}
We prove \ref{it: Galois sub flat 1}.\\
Let $D\cong \bigoplus_{k=0}^{n-1}\mathbf{E}_{k,n-1}$ be a filtered $(\varphi,N)$-module which is universal special with level $n-1$ (see Definition~\ref{def: universal phi N}).
Recall from Lemma~\ref{lem: Fil max sub phi N} that we have a strict embedding of filtered $(\varphi,N)$-modules $S_{n-1,n-1}\hookrightarrow D$ which induces a strict embedding of filtered $(\varphi,N)$-modules
\begin{equation}\label{equ: phi N sub embedding}
S_{n-1,m-1}\langle m-n\rangle \cong S_{n-1,n-1}/(S_{n-1,n-1}\cap \bigoplus_{k=0}^{n-m-1}\mathbf{E}_{k,n-1})\hookrightarrow D/\bigoplus_{k=0}^{n-m-1}\mathbf{E}_{k,n-1}.
\end{equation}
Now that $(D/\bigoplus_{k=0}^{n-m-1}\mathbf{E}_{k,n-1})\langle n-m\rangle[n-m]$ is universal special with level $m-1$ by inductively apply \ref{it: universal phi N induction 1} of Lemma~\ref{lem: universal phi N induction} for decreasing $m$, the (\ref{equ: phi N sub embedding}) induces an embedding between de Rham $(\varphi,\Gamma)$-modules
\[(\cD_{n-1,m-1}^{\flat}\otimes_{\cR}\cR(|\cdot|^{1-m}z^{1-n}))\otimes_{\cR}\cR(|\cdot|^{m-n}) \hookrightarrow (\cD_{m-1}\otimes_{\cR}\cR(\theta^{1-m}))\otimes_{\cR}\cR(\theta^{m-n})\]
and therefore an embedding $\cD_{n-1,m-1}^{\flat}\hookrightarrow \cD_{m-1}$. In particular, the $E$-vector space
\begin{equation}\label{equ: Galois sub flat Hom 1}
\Hom_{\varphi,\Gamma}(\cD_{n-1,m-1}^{\flat},\cD_{m-1})
\end{equation}
is non-zero. We thus deduce from Lemma~\ref{lem: Galois unique Hom} (by taking $\cD=\cD_{m-1}$ and $\cD'=\cD_{n-1,m-1}^{\flat}$ in \emph{loc.cit.}), \ref{it: Galois Hom vanishing 2} of Lemma~\ref{lem: Galois Hom vanishing} and \ref{it: Galois dR Ext exact 2} of Proposition~\ref{prop: Galois dR Ext exact} that (\ref{equ: Galois sub flat Hom 1}) is $1$ dimensional, and the aforementioned embedding $\cD_{n-1,m-1}^{\flat}\hookrightarrow \cD_{m-1}$ is the unique (up to scalars) non-zero element in (\ref{equ: Galois sub flat Hom 1}).

The proof of \ref{it: Galois sub flat 2} is very similar. We first construct an explicit embedding between filtered $(\varphi,N)$-modules
\[S_{\ell-1,m-1}\langle m-1\rangle[\ell-1]\hookrightarrow S_{n-1,m-1}\langle m-1\rangle[n-1]\]
whose associated embedding between de Rham $(\varphi,\Gamma)$-modules gives a non-zero element of
\begin{equation}\label{equ: Galois sub flat Hom 2}
\Hom_{\varphi,\Gamma}(\cD_{\ell-1,m-1}^{\flat},\cD_{n-1,m-1}^{\flat}).
\end{equation}
We thus deduce from Lemma~\ref{lem: Galois unique Hom} (by taking $\cD=\cD_{n-1,m-1}^{\flat}$ and $\cD'=\cD_{\ell-1,m-1}^{\flat}$ in \emph{loc.cit.}), \ref{it: Galois Hom vanishing 2} of Lemma~\ref{lem: Galois Hom vanishing} and \ref{it: Galois dR Ext exact 2} of Proposition~\ref{prop: Galois dR Ext exact} that (\ref{equ: Galois sub flat Hom 2}) is $1$ dimensional, and the aforementioned embedding $\cD_{\ell-1,m-1}^{\flat}\hookrightarrow \cD_{n-1,m-1}^{\flat}$ is the unique (up to scalars) non-zero element in (\ref{equ: Galois sub flat Hom 2}).
\end{proof}

Note that the unique (up to scalars) embedding $\cD_{n}^{\flat}=\cD_{n,n}^{\flat}\hookrightarrow \cD_{n}$ restricts to an embedding $\cE_{n}^{\flat}\otimes_E \cR(\theta^{n})\hookrightarrow \cE_{n}\otimes_E \cR(\theta^{n})$ and therefore induces an embedding $\cE_{n}^{\flat}\hookrightarrow \cE_{n}$ between $E$-vector spaces.
Consequently, we may identify $\cE_{n}^{\flat}$ with a canonical $E$-subspace of $\cE_{n}$ from now on.

\begin{lem}\label{lem: universal to sym}
Let $n\geq 1$. We have the following results.
\begin{enumerate}[label=(\roman*)]
\item \label{it: sym 1} We have $\Hom_{\varphi,\Gamma}(\cD_{n}^{\flat},\mathrm{Sym}^{n}(\cD_{1}))=1$ and the unique (up to scalars) non-zero map $\cD_{n}^{\flat}\rightarrow \mathrm{Sym}^{n}(\cD_{1})$ is an embedding which induces an isomorphism
    \begin{equation}\label{equ: sym isom invert t}
    \cD_{n}^{\flat}[\frac{1}{t}]\buildrel\sim\over\longrightarrow \mathrm{Sym}^{n}(\cD_{1})[\frac{1}{t}].
    \end{equation}
\item \label{it: sym 2} We have $\Hom_{\varphi,\Gamma}(\cD_{n},\mathrm{Sym}^{n}(\cD_{1}))=1$ and the unique (up to scalars) non-zero map $\cD_{n}\rightarrow \mathrm{Sym}^{n}(\cD_{1})$ is a surjection. Moreover, the composition map
    \begin{equation}\label{equ: sym composition}
    \Hom_{\varphi,\Gamma}(\cD_{n}^{\flat},\cD_{n})\otimes_E \Hom_{\varphi,\Gamma}(\cD_{n},\mathrm{Sym}^{n}(\cD_{1}))\rightarrow \Hom_{\varphi,\Gamma}(\cD_{n}^{\flat},\mathrm{Sym}^{n}(\cD_{1}))
    \end{equation}
    is an isomorphism between $1$-dimensional $E$-vector spaces.
\end{enumerate}
\end{lem}
\begin{proof}
We prove \ref{it: sym 1}.\\
Let $D=Ee_1\oplus Ee_0\oplus Ee_0'$ be a filtered $(\varphi,N)$-module as in \ref{it: filtered phi N 3} of Example~\ref{ex: filtered phi N} with $Ne_1=e_0$ and $\mathrm{Fil}_{H}^{1}(D)=E(e_1+e_0')$. Note that homogenous polynomial in variables $X,Y,Z$ with degree $n$ are spanned by $X^{n_1}Y^{n_0}Z^{n_0'}$ for all $n_1,n_0,n_0'\geq 0$ satisfying $n_1+n_0+n_0'=n$. Hence, the filtered $(\varphi,N)$-module $\mathrm{Sym}^{n}(D)$ has the form
\[\bigoplus_{n_1+n_0+n_0'=n}Ee_{n_1,n_0,n_0'},\]
with $(\varphi,N)$-actions given by $\varphi(e_{n_1,n_0,n_0'})=p^{n_1}e_{n_1,n_0,n_0'}$, $N(e_{0,n_0,n_0'})=0$ when $n_1=0$, and $N(e_{n_1,n_0,n_0'})=e_{n_1-1,n_0+1,n_0'}$ when $n_1\geq 1$. Moreover, we have $\mathrm{Fil}_{H}^{n}(\mathrm{Sym}^{n}(D))=\mathrm{Sym}^{n}(\mathrm{Fil}_{H}^{1}(D))\subseteq \mathrm{Sym}^{n}(D)$. In particular, there exists a filtered $(\varphi,N)$-module $S_{n,n}$ as in \ref{it: filtered phi N 5} of Example~\ref{ex: filtered phi N} together with an embedding $S_{n,n}\hookrightarrow \mathrm{Sym}^{n}(D)$ between filtered $(\varphi,N)$-modules which induces an isomorphism between their underlying $(\varphi,N)$-modules as well as an isomorphism $\mathrm{Fil}_{H}^{n}(S_{n,n})\buildrel\sim\over\longrightarrow \mathrm{Fil}_{H}^{n}(\mathrm{Sym}^{n}(D))$. By consider the embedding between de Rham $(\varphi,\Gamma)$-modules associated with
\[S_{n,n}\langle n\rangle[n]\hookrightarrow \mathrm{Sym}^{n}(D)\langle n\rangle[n]=\mathrm{Sym}^{n}(D\langle 1\rangle[1])\]
and using \ref{it: Galois dR Ext exact 2} of Proposition~\ref{prop: Galois dR Ext exact}, we obtain an embedding
\begin{equation}\label{equ: sym embedding}
\cD_{n}^{\flat}\hookrightarrow \mathrm{Sym}^{n}(\cD_{1})
\end{equation}
which induces an isomorphism (\ref{equ: sym isom invert t}). In particular, the $E$-vector space
\begin{equation}\label{equ: flat sym Hom}
\Hom_{\varphi,\Gamma}(\cD_{n}^{\flat},\mathrm{Sym}^{n}(\cD_{1}))
\end{equation}
is non-zero.
We thus deduce from Lemma~\ref{lem: Galois unique Hom} (by taking $\cD=\mathrm{Sym}^{n}(\cD_{1})$ and $\cD'=\cD_{n}^{\flat}$ in \emph{loc.cit.}), \ref{it: Galois Hom vanishing 2} of Lemma~\ref{lem: Galois Hom vanishing} and \ref{it: Galois dR Ext exact 2} of Proposition~\ref{prop: Galois dR Ext exact} that (\ref{equ: flat sym Hom}) is $1$ dimensional, and the embedding (\ref{equ: sym embedding}) is the unique (up to scalars) non-zero element in (\ref{equ: flat sym Hom}).
This finishes the proof of \ref{it: sym 1}.

We prove \ref{it: sym 2}.\\
Given $1\leq k\leq n$, we note that any non-zero map $\mathrm{Sym}^{n}(D)\langle n\rangle[n]\rightarrow \langle k\rangle[k]$ between filtered $(\varphi,N)$-modules must contain $\mathrm{Fil}_{H}^{0}(\mathrm{Sym}^{n}(D)\langle n\rangle[n])$ in the kernel.
Since $\mathrm{Sym}^{n}(D)\langle n\rangle[n]$ is generated by $\mathrm{Fil}_{H}^{0}(\mathrm{Sym}^{n}(D)\langle n\rangle[n])$ as a $(\varphi,N)$-module, we conclude that all maps $\mathrm{Sym}^{n}(D)\langle n\rangle[n]\rightarrow \langle k\rangle[k]$ between filtered $(\varphi,N)$-modules must be zero, or equivalently
\[\Hom_{\varphi,\Gamma}(\mathrm{Sym}^{n}(\cD_{1}),\cR(\theta^{k}))=0\]
for each $1\leq k\leq n$. This together with Lemma~\ref{lem: Galois universal Hom} (upon taking $\cD=\mathrm{Sym}^{n}(\cD_{1})$ in \emph{loc.cit.}) gives $\Hom_{\varphi,\Gamma}(\cD_{n},\mathrm{Sym}^{n}(\cD_{1}))=1$ with the unique (up to scalars) non-zero map $\cD_{n}\rightarrow \mathrm{Sym}^{n}(\cD_{1})$ being a surjection.
Finally, since the maps $\cD_{n}^{\flat}\hookrightarrow \cD_{n}\twoheadrightarrow \mathrm{Sym}^{n}(\cD_{1})$ induces isomorphisms $\cR(1_{K^\times})\buildrel\sim\over\longrightarrow \cR(1_{K^\times})\buildrel\sim\over\longrightarrow \cR(1_{K^\times})$, the composition map (\ref{equ: sym composition}) is non-zero and thus an isomorphism between $1$-dimensional $E$-vector spaces. This finishes the proof of \ref{it: sym 2}.
\end{proof}

For each $n\geq 0$, we define $\cE_{n}^{\sharp}$ as the canonical $E$-subspace of $\cE_{n}$ such that $\cE_{n}^{\sharp}\otimes_E\cR(\theta^{n})$ equals the kernel of the composition of
\[\cE_{n}\otimes_E\cR(\theta^{n})\hookrightarrow \cD_{n} \twoheadrightarrow \mathrm{Sym}^{n}(\cD_{1}).\]
In particular, we have $\cE_{n}^{\sharp}=0$ for $n\leq 1$.

\begin{prop}\label{prop: sharp flat decomposition}
We have
\begin{equation}\label{equ: sharp flat decomposition}
\cE_{n}=\cE_{n}^{\flat}\oplus \cE_{n}^{\sharp}
\end{equation}
for each $n\geq 0$.
\end{prop}
\begin{proof}
Recall from Lemma~\ref{lem: universal to sym} that we have maps
\begin{equation}\label{equ: sharp flat composition 1}
\cD_{n}^{\flat}\hookrightarrow \cD_{n}\twoheadrightarrow \mathrm{Sym}^{n}(\cD_{1})
\end{equation}
whose composition is an embedding which becomes an isomorphism after inverting $t$. By the construction of both maps in (\ref{equ: sharp flat composition 1}), we know that they restrict to the following maps
\begin{equation}\label{equ: sharp flat composition 2}
\cE_{n}^{\flat}\otimes_E\cR(\theta^{n}) \hookrightarrow \cE_{n}\otimes_E\cR(\theta^{n}) \twoheadrightarrow \mathrm{Fil}_{W}^{n}(\mathrm{Sym}^{n}(\cD_{1})).
\end{equation}
We observe that the composition of (\ref{equ: sharp flat composition 2}) must be an isomorphism as it becomes an isomorphism after inverting $t$. This together with the definition of $\cE_{n}^{\sharp}$ gives (\ref{equ: sharp flat decomposition}).
\end{proof}
\subsection{Relation with deformations under Tate duality}\label{subsec: deformation}
In this section, we relate $\cE_{n}^{\flat}$ and $\cE_{n}^{\sharp}$ to deformations of $(\varphi,\Gamma)$-modules via local Tate duality (see Corollary~\ref{cor: det deform sym}, Proposition~\ref{prop: Ext max} and Proposition~\ref{prop: special Ext max}).
We assume throughout this section that $K=\Q_p$.

\begin{lem}\label{lem: Galois Ext collection}
Let $\delta_i=|\cdot|^{\ell_i}z^{k_i}$ for $i=1,2$ and $\delta=|\cdot|^{\ell}z^{k}$. We have the following results.
\begin{enumerate}[label=(\roman*)]
\item \label{it: Galois Ext collection 0} The $E$-vector space $\mathrm{Ext}_{\varphi,\Gamma}^1(\cR(\delta_1),\cR(\delta_2))$ is $2$-dimensional if and only if either $\ell_2=\ell_1+1$ with $k_2>k_1$ or $\ell_2=\ell_1$ with $k_2\geq k_1$, and is $1$-dimensional otherwise.
\item \label{it: Galois Ext collection 1} We have $\mathrm{Ext}_{\varphi,\Gamma}^2(\cR(\delta_1),\cR(\delta_2))\neq 0$ if and only if $\ell_2=\ell_1+1$ and $k_2>k_1$, in which case it is $1$-dimensional.
\item \label{it: Galois Ext collection 2} We have $\mathrm{Ext}_{\varphi,\Gamma}^2(\cR(1_{K^{\times}}),\cD)=0$ for each de Rham $(\varphi,\Gamma)$-module $\cD$ which is special with level $1$.
\item \label{it: Galois Ext collection 3} If $\ell>\ell_1=\ell_2$ and $k_1\geq k_2$, then the embedding $\cR(\delta_1)\hookrightarrow \cR(\delta_2)$ induces an isomorphism
    \[\mathrm{Ext}_{\varphi,\Gamma}^1(\cR(\delta),\cR(\delta_1))\buildrel\sim\over\longrightarrow \mathrm{Ext}_{\varphi,\Gamma}^1(\cR(\delta),\cR(\delta_2))\]
    between $1$-dimensional $E$-vector spaces and we have $\Hom_{\varphi,\Gamma}(\cR(\delta),\cR(\delta_2)/\cR(\delta_1))=0$.
\item \label{it: Galois Ext collection 4} Let $\cD$ be a $(\varphi,\Gamma)$-module that fits into a non-split extension $0\rightarrow \cR(\delta_2)\rightarrow \cD\rightarrow \cR(\delta_1)\rightarrow 0$. If $\ell_1\neq \ell_2$ and $k_2>k_1$, then we have $\Hom_{\varphi,\Gamma}(\cR(\delta_1),\cD)=0$.
\end{enumerate}
\end{lem}
\begin{proof}
\ref{it: Galois Ext collection 0} follows from \cite[Thm.~2.9]{Col08}.
\ref{it: Galois Ext collection 1} and \ref{it: Galois Ext collection 2} follows from \cite[Prop.~2.1]{Col08}, local Tate duality and \ref{it: special mod 2} of Definition~\ref{def: special mod}.
\ref{it: Galois Ext collection 3} follows from \cite[Thm.~2.22 (i)]{Col08}.
For \ref{it: Galois Ext collection 4}, we note that any non-zero map $\cR(\delta_1)\rightarrow \cD$ is necessarily a section of the surjection $\cD\twoheadrightarrow \cR(\delta_1)$ by \cite[Prop.~2.1]{Col08}, and thus the existence of a non-zero map $\cR(\delta_1)\rightarrow \cD$ would force $\cD$ to split.
\end{proof}

\begin{lem}\label{lem: deform induction seq}
Let $n\geq 1$, $1\leq m\leq n-1$ and $\cD$ be a de Rham $(\varphi,\Gamma)$-module which is special with level $n-1$ (see \ref{it: special mod 2} of Definition~\ref{def: special mod}). Then the short exact sequence $0\rightarrow \mathrm{Fil}_{W}^{m}(\cD) \rightarrow \cD \rightarrow \cD/\mathrm{Fil}_{W}^{m}(\cD) \rightarrow 0$ induces the following short exact sequence
\begin{equation}\label{equ: deform induction seq}
0\rightarrow \mathrm{Ext}_{\varphi,\Gamma}^1(\cD/\mathrm{Fil}_{W}^{m}(\cD),\cD) \rightarrow \mathrm{Ext}_{\varphi,\Gamma}^1(\cD,\cD) \rightarrow \mathrm{Ext}_{\varphi,\Gamma}^1(\mathrm{Fil}_{W}^{m}(\cD),\cD) \rightarrow 0
\end{equation}
\end{lem}
\begin{proof}
For each $0\leq k\leq m-1$, note that $\cD$ admits a filtration with graded pieces $\mathrm{Fil}_{W}^{k}(\cD)/\mathrm{Fil}_{W}^{k+2}(\cD)$ and $\cR(\theta^{\ell})$ for each $0\leq \ell\leq n-1$ with $\ell\neq k,k+1$.
Now that we have
\[\mathrm{Ext}_{\varphi,\Gamma}^2(\cR(\theta^{k}), \mathrm{Fil}_{W}^{k}(\cD)/\mathrm{Fil}_{W}^{k+2}(\cD))\cong \mathrm{Ext}_{\varphi,\Gamma}^2(\cR(1_{K^{\times}}), (\mathrm{Fil}_{W}^{k}(\cD)/\mathrm{Fil}_{W}^{k+2}(\cD))\otimes_{\cR}\cR(\theta^{-k})=0\]
from \ref{it: Galois Ext collection 2} of Lemma~\ref{lem: Galois Ext collection}, and $\mathrm{Ext}_{\varphi,\Gamma}^2(\cR(\theta^{k}), \cR(\theta^{\ell}))=0$ for each $\ell\neq k,k+1$ from \ref{it: Galois Ext collection 1} of Lemma~\ref{lem: Galois Ext collection}, a d\'evissage gives $\mathrm{Ext}_{\varphi,\Gamma}^2(\cR(\theta^{k}), \cD)=0$
for each $0\leq k\leq m-1$, and thus
\begin{equation}\label{equ: deform induction seq Ext2}
\mathrm{Ext}_{\varphi,\Gamma}^2(\cD/\mathrm{Fil}_{W}^{m}(\cD), \cD)=0
\end{equation}
by a further d\'evissage.
For each $0\leq k\leq n-2$, we have $\Hom_{\varphi,\Gamma}(\cR(\theta^{k}),\cR(\theta^{\ell}))=0$ for each $\ell\neq k$ from \ref{it: char dR Ext 1} of Lemma~\ref{lem: char dR Ext}, and
\[\Hom_{\varphi,\Gamma}(\cR(\theta^{k}),\mathrm{Fil}_{W}^{k}(\cD)/\mathrm{Fil}_{W}^{k+2}(\cD))=0\]
from \ref{it: Galois Ext collection 4} of Lemma~\ref{lem: Galois Ext collection}, which altogether gives
\[\Hom_{\varphi,\Gamma}(\cR(\theta^{k}),\cD)=0\]
by a d\'evissage.
In particular, a further d\'evissage gives
\begin{equation}\label{equ: deform induction seq Hom}
\Hom_{\varphi,\Gamma}(\cD/\mathrm{Fil}_{W}^{k}(\cD),\cD)=0.
\end{equation}
for each $0\leq k\leq n-1$, and thus the embeddings $\cR(\theta^{n-1})\hookrightarrow \mathrm{Fil}_{W}^{m}(\cD)\hookrightarrow \cD$ induce the following embeddings
\begin{equation}\label{equ: deform induction seq Hom embedding}
\Hom_{\varphi,\Gamma}(\cD,\cD)\hookrightarrow \Hom_{\varphi,\Gamma}(\mathrm{Fil}_{W}^{m}(\cD),\cD)\hookrightarrow \Hom_{\varphi,\Gamma}(\cR(\theta^{n-1}),\cD).
\end{equation}
Now that the $\Hom_{\varphi,\Gamma}(\cD,\cD)$ is clearly non-zero and $\Hom_{\varphi,\Gamma}(\cR(\theta^{n-1}),\cD)$ is $1$-dimensional by \ref{it: char dR Ext 1} of Lemma~\ref{lem: char dR Ext} and a simple d\'evissage, we see that the embeddings in (\ref{equ: deform induction seq Hom embedding}) are actually isomorphisms between $1$-dimensional $E$-vector spaces.
The together with (\ref{equ: deform induction seq Ext2}) finishes the proof of (\ref{equ: deform induction seq}).
\end{proof}

The short exact sequence $0\rightarrow \cR(\theta^{n-1}) \rightarrow \cD \rightarrow \cD/\cR(\theta^{n-1})\rightarrow 0$ together with (\ref{equ: deform induction seq Ext2}) (with $m=n-1$ in \emph{loc.cit.}) and (\ref{equ: deform induction seq Hom}) (with $k=n-1$ in \emph{loc.cit.}) also gives the following long exact sequence
\begin{multline}\label{equ: deform induction devissage}
0\rightarrow \Hom_{\varphi,\Gamma}(\cD/\cR(\theta^{n-1}),\cD/\cR(\theta^{n-1}))\\
\rightarrow \mathrm{Ext}_{\varphi,\Gamma}^1(\cD/\cR(\theta^{n-1}),\cR(\theta^{n-1}))
\rightarrow \mathrm{Ext}_{\varphi,\Gamma}^1(\cD/\cR(\theta^{n-1}),\cD) \rightarrow \mathrm{Ext}_{\varphi,\Gamma}^1(\cD/\cR(\theta^{n-1}),\cD/\cR(\theta^{n-1}))\\
\rightarrow \mathrm{Ext}_{\varphi,\Gamma}^2(\cD/\cR(\theta^{n-1}),\cR(\theta^{n-1}))\rightarrow 0.
\end{multline}

Let $n\geq 1$ and $\cD$ be a de Rham $(\varphi,\Gamma)$-module which is special with level $n-1$.
A d\'evissage using \ref{it: Galois Ext collection 1} of Lemma~\ref{lem: Galois Ext collection} shows that the embedding $\cR(\theta^{n-1})\hookrightarrow \cD$ induces an isomorphism
\begin{equation}\label{equ: full deform Ext2}
\mathrm{Ext}_{\varphi,\Gamma}^2(\cD,\cR(\theta^{n})) \buildrel\sim\over\longrightarrow \mathrm{Ext}_{\varphi,\Gamma}^2(\cR(\theta^{n-1}),\cR(\theta^{n}))
\end{equation}
between $1$-dimensional spaces. Recall from \cite[Lem.~2.2]{BD20} that the cup product map
\begin{equation}\label{equ: Tate dual cup}
\mathrm{Ext}_{\varphi,\Gamma}^1(\cR(\theta^{n-1}),\cD) \otimes_E \mathrm{Ext}_{\varphi,\Gamma}^1(\cD,\cR(\theta^{n})) \buildrel\cup\over\longrightarrow \mathrm{Ext}_{\varphi,\Gamma}^2(\cR(\theta^{n-1}),\cR(\theta^{n}))
\end{equation}
induces a perfect pairing between $n+1$-dimensional spaces by local Tate duality.
To summarize, we have the following commutative diagram (which is the dual version of \cite[(2.2)]{BD20})
\begin{equation}\label{equ: BD pairing}
\xymatrix{
\mathrm{Ext}_{\varphi,\Gamma}^1(\cD,\cD) \ar@{->>}[d] & \otimes_E & \mathrm{Ext}_{\varphi,\Gamma}^1(\cD,\cR(\theta^{n})) \ar@{=}[d] \ar^{\cup}[r] & \mathrm{Ext}_{\varphi,\Gamma}^2(\cD,\cR(\theta^{n})) \ar^{\wr}[d]\\
\mathrm{Ext}_{\varphi,\Gamma}^1(\cR(\theta^{n-1}),\cD) & \otimes_E & \mathrm{Ext}_{\varphi,\Gamma}^1(\cD,\cR(\theta^{n})) \ar^{\cup}[r] & \mathrm{Ext}_{\varphi,\Gamma}^2(\cR(\theta^{n-1}),\cR(\theta^{n}))
}
\end{equation}
with LHS vertical map being a surjection by Lemma~\ref{lem: deform induction seq}.
We consider a de Rham $(\varphi,\Gamma)$-module $\cD^{+}$ over $\cR$ that fits into a non-split short exact sequence
\begin{equation}\label{equ: special short exact seq}
0\rightarrow \cR(\theta^{n}) \rightarrow \cD^{+} \rightarrow \cD \rightarrow 0
\end{equation}
The isomorphism class of $\cD^{+}$ is uniquely determined by a $1$-dimensional $E$-subspace
\[E[\cD^{+}]\subseteq \mathrm{Ext}_{g}^1(\cD,\cR(\theta^{n}))=\mathrm{Ext}_{\varphi,\Gamma}^1(\cD,\cR(\theta^{n})).\]
We define
\[\mathrm{Ext}_{\cD^{+}}^1(\cR(\theta^{n-1}),\cD)\subseteq \mathrm{Ext}_{\varphi,\Gamma}^1(\cR(\theta^{n-1}),\cD)\]
as the orthogonal complement of $E[\cD^{+}]$, and define
\[\mathrm{Ext}_{\cD^{+}}^1(\cD,\cD)\subseteq \mathrm{Ext}_{\varphi,\Gamma}^1(\cD,\cD)\]
as the preimage of $\mathrm{Ext}_{\cD^{+}}^1(\cR(\theta^{n-1}),\cD)$ under the LHS surjection of (\ref{equ: BD pairing}).

Let $n\geq 1$ and $\cD$ be a $(\varphi,\Gamma)$-module over $\cR$ of rank $n$. Each element of $\mathrm{Ext}_{\varphi,\Gamma}^1(\cD,\cD)$ gives rise to a $(\varphi,\Gamma)$-module $\tld{\cD}$ over $\cR$ that fits into a short exact sequence $0\rightarrow \cD\rightarrow \tld{\cD} \rightarrow \cD\rightarrow 0$. We can interpret $\tld{\cD}$ as a $(\varphi,\Gamma)$-module over $\cR_{E[\varepsilon]/\varepsilon^2}$ and consider the $(\varphi,\Gamma)$-module $\mathrm{det}(\tld{\cD})\defeq \wedge^{n}_{\cR_{E[\varepsilon]/\varepsilon^2}}\tld{\cD}$ over $\cR_{E[\varepsilon]/\varepsilon^2}$ whose underlying $(\varphi,\Gamma)$-module over $\cR$ fits into a short exact sequence
\begin{equation}\label{equ: det deform}
0\rightarrow \mathrm{det}(\cD) \rightarrow \mathrm{det}(\tld{\cD}) \rightarrow \mathrm{det}(\cD) \rightarrow 0.
\end{equation}
with $\mathrm{det}(\cD)\defeq \wedge^{n}_{\cR}\cD$.
The association $\tld{\cD}\mapsto\mathrm{det}(\tld{\cD})=\wedge^{n}_{\cR_{E[\varepsilon]/\varepsilon^2}}\tld{\cD}$ defines a map
\begin{equation}\label{equ: pass to det deform 1}
\kappa_{\mathrm{det}}: \mathrm{Ext}_{\varphi,\Gamma}^1(\cD,\cD)\rightarrow \mathrm{Ext}_{\varphi,\Gamma}^1(\mathrm{det}(\cD),\mathrm{det}(\cD))\cong \Hom(\Q_p^\times,E).
\end{equation}
We define
\[\mathrm{Ext}_{\mathrm{det}}^1(\cD,\cD)\defeq \mathrm{ker}(\kappa_{\mathrm{det}})\subseteq \mathrm{Ext}_{\varphi,\Gamma}^1(\cD,\cD)\]
which is the $E$-subspace consisting of elements such that (\ref{equ: det deform}) splits.

Given $n\geq 1$ and $\cD$ as above, we note that $\wedge^{n-1}_{\cR}\cD$ is a $(\varphi,\Gamma)$-module over $\cR$ of rank $n$, with $\wedge^{n-1}_{\cR}\cD\otimes_{\cR}\cD$ containing a canonical direct summand $\mathrm{det}(\cD)=\wedge^{n}_{\cR}\cD$. Given a deformation $\tld{\cD}$ of $\cD$ associated with an element of $\mathrm{Ext}_{\varphi,\Gamma}^1(\cD,\cD)$, the short exact sequence $0\rightarrow \cD\rightarrow \tld{\cD}\rightarrow \cD\rightarrow 0$ induces a short exact sequence
\[0\rightarrow \wedge^{n-1}_{\cR}\cD\otimes_{\cR}\cD\rightarrow \wedge^{n-1}_{\cR}\cD\otimes_{\cR}\tld{\cD}\rightarrow \wedge^{n-1}_{\cR}\cD\otimes_{\cR}\cD\rightarrow 0.\]
After a push forward along $\wedge^{n-1}_{\cR}\cD\otimes_{\cR}\cD\twoheadrightarrow \mathrm{det}(\cD)$ on LHS and a pull back along $\mathrm{det}(\cD)\hookrightarrow \wedge^{n-1}_{\cR}\cD\otimes_{\cR}\cD$ on RHS, we see that $\wedge^{n-1}_{\cR}\cD\otimes_{\cR}\tld{\cD}$ admits a natural subquotient $\mathrm{det}'(\tld{\cD})$ which is a deformation of $\mathrm{det}(\cD)$.
The association $\tld{\cD}\mapsto \mathrm{det}'(\tld{\cD})$ defines a map
\begin{equation}\label{equ: pass to det deform 2}
\kappa_{\mathrm{det}}': \mathrm{Ext}_{\varphi,\Gamma}^1(\cD,\cD)\rightarrow \mathrm{Ext}_{\varphi,\Gamma}^1(\mathrm{det}(\cD),\mathrm{det}(\cD)).
\end{equation}
We have the following result.
\begin{lem}\label{lem: det cond split}
We have
\begin{equation}\label{equ: det cond split}
\mathrm{ker}(\kappa_{\mathrm{det}}')=\mathrm{Ext}_{\mathrm{det}}^1(\cD,\cD)=\mathrm{ker}(\kappa_{\mathrm{det}}).
\end{equation}
\end{lem}
\begin{proof}
Now that $\tld{\cD}$ is a free $\cR_{E[\varepsilon]/\varepsilon^2}$-module of rank $n$, we may choose a basis
\begin{equation}\label{equ: phi gamma deform basis}
\{e_1,\dots,e_n\}
\end{equation}
of $\tld{\cD}$, namely we have $\tld{\cD}=\bigoplus_{i=1}^{n}\cR_{E[\varepsilon]/\varepsilon^2}e_i$. In particular, the underlying $\cR$-module of $\tld{\cD}$ admits
\begin{equation}\label{equ: double basis}
\{e_1,\dots,e_n\}\sqcup\{\varepsilon e_1,\dots,\varepsilon e_n\}
\end{equation}
as a basis.
There exists $A_{\varphi}\in\mathrm{GL}_{n}(\cR)$ and $B_{\varphi}\in\mathrm{M}_{n}(\cR)$ such that $A_{\varphi}+\varepsilon B_{\varphi}$ is the matrix of $\varphi_{\tld{\cD}}$ (over $\cR_{E[\varepsilon]/\varepsilon^2}$) under the basis $\{e_1,\dots,e_n\}$, with $A_{\varphi}$ being the matrix of $\varphi_{\cD}$ (over $\cR$) under the basis $\{e_1,\dots,e_n\}$. We observe (using $\varepsilon^2=0$) that
\begin{equation}\label{equ: det deform matrix 1}
\mathrm{det}(A_{\varphi}+\varepsilon B_{\varphi})=\mathrm{det}(A_{\varphi})\mathrm{det}(\mathrm{Id}_{n}+\varepsilon A_{\varphi}^{-1}B_{\varphi})=\mathrm{det}(A_{\varphi})(1+\varepsilon\mathrm{Tr}(A_{\varphi}^{-1}B_{\varphi}))
\end{equation}
is the matrix of $\varphi_{\mathrm{det}(\tld{\cD})}$ (over $\cR_{E[\varepsilon]/\varepsilon^2}$) under the basis $\{e_1\wedge\cdots\wedge e_n\}$, and thus the deformation $\mathrm{det}(\tld{\cD})$ of $\mathrm{det}(\cD)$ as $\varphi$-modules splits if and only if there exists a choice of (\ref{equ: phi gamma deform basis}) such that $\mathrm{Tr}(A_{\varphi}^{-1}B_{\varphi})=0$.
Now we notice that the matrix of $\varphi_{\tld{\cD}}$ (over $\cR$) under the basis (\ref{equ: double basis}) has the form
\[\begin{pmatrix}
  A_{\varphi} & 0 \\
  B_{\varphi} & A_{\varphi}
\end{pmatrix}.\]
Note that the basis $\{e_1,\dots,e_n\}$ induces a basis of $\wedge^{n-1}_{\cR}\cD$ so that the matrix of $\varphi_{\wedge^{n-1}_{\cR}\cD}$ (over $\cR$) under such a basis is given by $\wedge^{n-1}A_{\varphi}=\mathrm{det}(A_{\varphi})A_{\varphi}^{-1}$.
The basis (\ref{equ: double basis}) of $\tld{\cD}$ and the aforementioned basis of $\wedge^{n-1}_{\cR}\cD$ induces a basis of $\wedge^{n-1}_{\cR}\cD\otimes_{\cR}\tld{\cD}$ and thus a basis of $\tld{\mathrm{det}}'(\cD)$ (both over $\cR$), under which $\varphi_{\wedge^{n-1}_{\cR}\cD\otimes_{\cR}\tld{\cD}}$ has matrix
\[\begin{pmatrix}
  \wedge^{n-1}A_{\varphi}\otimes_E A_{\varphi} & 0 \\
  \wedge^{n-1}A_{\varphi}\otimes_E B_{\varphi} & \wedge^{n-1}A_{\varphi}\otimes_E A_{\varphi}
\end{pmatrix}\]
and $\varphi_{\tld{\mathrm{det}}'(\cD)}$ has matrix
\begin{equation}\label{equ: det deform matrix 2}
\mathrm{det}(A_{\varphi})\begin{pmatrix}
  \mathrm{Id}_{n} & 0 \\
  \mathrm{Tr}(A_{\varphi}^{-1}B_{\varphi}) & \mathrm{Id}_{n}
\end{pmatrix}.
\end{equation}
Comparing (\ref{equ: det deform matrix 1}) with (\ref{equ: det deform matrix 2}), we conclude that the deformation of $\varphi$-module $\mathrm{det}'(\tld{\cD})$ splits if and only if the deformation of $\varphi$-module $\mathrm{det}(\tld{\cD})$ splits, in which case there exists a choice of (\ref{equ: phi gamma deform basis}) such that $\mathrm{Tr}(A_{\varphi}^{-1}B_{\varphi})=0$.
Upon choosing a topological generator of $\Gamma$, we can similarly define $A_{\Gamma}\in\mathrm{GL}_{n}(\cR)$ and $B_{\Gamma}\in\mathrm{M}_{n}(\cR)$ with respect to a choice of (\ref{equ: phi gamma deform basis}), and a parallel argument shows that the deformation of $\Gamma$-module $\mathrm{det}'(\tld{\cD})$ splits if and only if the deformation of $\Gamma$-module $\mathrm{det}(\tld{\cD})$ splits, in which case we have $\mathrm{Tr}(A_{\Gamma}^{-1}B_{\Gamma})=0$ for any choice of (\ref{equ: phi gamma deform basis}).
Hence, if at least one of the deformations of $(\varphi,\Gamma)$-module $\mathrm{det}'(\tld{\cD})$ and $\mathrm{det}(\tld{\cD})$ splits, we can choose (\ref{equ: phi gamma deform basis}) such that $\mathrm{Tr}(A_{\varphi}^{-1}B_{\varphi})=0$ and $\mathrm{Tr}(A_{\Gamma}^{-1}B_{\Gamma})=0$ and thus both $\mathrm{det}'(\tld{\cD})$ and $\mathrm{det}(\tld{\cD})$ split. The proof is thus finished.
\end{proof}

In the rest of the section, thanks to Lemma~\ref{lem: det cond split}, we will often use without further explanation that
\[\mathrm{Ext}_{\mathrm{det}}^1(\cD,\cD)=\mathrm{ker}(\kappa_{\mathrm{det}}').\]

\begin{lem}\label{lem: sym det generator}
Let $n\geq 2$ and $\cB$ be a de Rham $(\varphi,\Gamma)$-module which is special with level $1$. Let $\cD\defeq\mathrm{Sym}^{n-1}(\cB)$ and $\cD^{+}\defeq\mathrm{Sym}^{n}(\cB)$.
Then we have
\begin{equation}\label{equ: sym det generator}
\zeta_{n-1,0}(\mathrm{Ext}_{\mathrm{det}}^1(\cD,\cD)\cap \mathrm{Ext}_{\cD^{+}}^1(\cD,\cD))\neq 0.
\end{equation}
\end{lem}
\begin{proof}
Note that $\cB$ has rank $2$ over $\cR$ and that $\wedge^2_{\cR}\cB\cong \cR(\theta)$. Hence, we have a canonical isomorphism
\begin{equation}\label{equ: sym decomposition}
\cB\otimes_{\cR}\cD \cong \cD^{+}\oplus \cD^{-}
\end{equation}
with $\cD^{-}\defeq (\wedge^2_{\cR}\cB)\otimes_{\cR}\mathrm{Sym}^{n-2}(\cB)=\cR(\theta)\otimes_{\cR}\mathrm{Sym}^{n-2}(\cB)$. Now we consider a deformation $\tld{\cD}$ of $\cD$ associated with an element of $\mathrm{Ext}_{\varphi,\Gamma}^1(\cD,\cD)$.
Upon pulling back along $\cR(\theta^{n-1})\hookrightarrow \cD$, we obtain a sub $\tld{\cD}'\subseteq \tld{\cD}$ that fits into
\[0\rightarrow \cD \rightarrow \tld{\cD}' \rightarrow \cR(\theta^{n-1}) \rightarrow 0.\]
We observe that $\cB\otimes_{\cR}\tld{\cD}$ is a deformation of (\ref{equ: sym decomposition}), which gives a canonical map
\[\varsigma: \mathrm{Ext}_{\varphi,\Gamma}^1(\cD,\cD)\rightarrow \mathrm{Ext}_{\varphi,\Gamma}^1(\cD^{-},\cD^{-}).\]
Using a pull-back along $\cD^{-}\hookrightarrow \cB\otimes_{\cR}\cD$ and a push-out along $\cB\otimes_{\cR}\cD\twoheadrightarrow \cD^{+}$, we obtain a subquotient $\tld{\cD}^{\dagger}$ of $\cB\otimes_{\cR}\tld{\cD}$ which fits into the short exact sequence
\[0\rightarrow \cD^{+} \rightarrow \tld{\cD}^{\dagger} \rightarrow \cD^{-} \rightarrow 0.\]
By a further pull-back along $\cR(\theta^{n-1})\hookrightarrow \cD^{-}$, we obtain a sub $\tld{\cD}^{\ddagger}\subseteq \tld{\cD}^{\dagger}$ which fits into
\[0\rightarrow \cD^{+} \rightarrow \tld{\cD}^{\ddagger} \rightarrow \cR(\theta^{n-1}) \rightarrow 0.\]
The surjection $\cB\twoheadrightarrow \cR(1_{K^{\times}})$ also induces a surjection
\begin{equation}\label{equ: sym top}
\cB\otimes_{\cR}\tld{\cD}\twoheadrightarrow \tld{\cD}.
\end{equation}
We divide the rest of the proof into the following steps.

\textbf{Step $1$}: We prove that
\begin{equation}\label{equ: sym det transverse}
\zeta_{n-1,0}(\mathrm{ker}(\varsigma))\neq 0.
\end{equation}
Let $\tld{\cD}^{-}$ be a deformation of $\cD^{-}$ associated with an element of $\mathrm{Ext}_{\varphi,\Gamma}^1(\cD^{-},\cD^{-})$. Note that we have
\[\cR(\theta^{-1})\otimes_{\cR}\cB\otimes_{\cR}\cD^{-}\cong \cD\oplus (\cR(\theta)\otimes_{\cR}\mathrm{Sym}^{n-3}(\cB))\]
with the term $\cR(\theta)\otimes_{\cR}\mathrm{Sym}^{n-3}(\cB)$ omitted when $n=2$.
Using a push-out along $\cR(\theta^{-1})\otimes_{\cR}\cB\otimes_{\cR}\cD^{-}\twoheadrightarrow \cD$ and a pull-back along $\cD\hookrightarrow \cR(\theta^{-1})\otimes_{\cR}\cB\otimes_{\cR}\cD^{-}$, we see that $\cR(\theta^{-1})\otimes_{\cR}\cB\otimes_{\cR}\tld{\cD}^{-}$ admits a natural subquotient $\tld{\cD}$ which is a deformation of $\cD$. In other words, we have constructed a map
\[\vartheta: \mathrm{Ext}_{\varphi,\Gamma}^1(\cD^{-},\cD^{-})\rightarrow \mathrm{Ext}_{\varphi,\Gamma}^1(\cD,\cD).\]
We observe that $\varsigma(\vartheta(\tld{\cD}^{-}))$ is a subquotient of $\cB\otimes_{\cR}\vartheta(\tld{\cD}^{-})$ and thus a subquotient of 
\begin{equation}\label{equ: sym det double tensor}
\cB_{\otimes_{\cR}}(\cR(\theta^{-1})\otimes_{\cR}\cB\otimes_{\cR}\tld{\cD}^{-})=(\cB_{\otimes_{\cR}}\cB\otimes_{\cR}\cR(\theta^{-1}))\otimes_{\cR}\tld{\cD}^{-}.
\end{equation}
Using
\[\cB\otimes_{\cR}\cB\otimes_{\cR}\cR(\theta^{-1})\cong (\mathrm{det}(\cB)\oplus\mathrm{Sym}^2(\cB))\otimes_{\cR}\cR(\theta^{-1})\cong  \cR(1_{K^{\times}})\oplus(\mathrm{Sym}^2(\cB)\otimes_{\cR}\cR(\theta^{-1})),\]
We see that $\varsigma(\vartheta(\tld{\cD}^{-}))$ equals the natural subquotient $\tld{\cD}^{-}$ of (\ref{equ: sym det double tensor}) obtained from a push-out along $(\cB_{\otimes_{\cR}}\cR(\theta^{-1})\otimes_{\cR}\cB)\otimes_{\cR}\tld{\cD}^{-}\twoheadrightarrow \tld{\cD}^{-}$ and a pull-back along $\tld{\cD}^{-}\hookrightarrow(\cB_{\otimes_{\cR}}\cR(\theta^{-1})\otimes_{\cR}\cB)\otimes_{\cR}\tld{\cD}^{-}$.
In other words, $\vartheta$ is a section of the map $\varsigma$ with $\varsigma$ being surjective and $\vartheta$ being injective.
Now that the composition of the injections 
\[\cR(\theta^{n-1})\hookrightarrow \cD\hookrightarrow \cR(\theta^{-1})\otimes_{\cR}\cB\otimes_{\cR}\cD^{-}\] 
factors through
\[\cR(\theta^{n-1})\buildrel\sim\over\longrightarrow \cR(\theta^{-1})\otimes_{\cR}\cR(\theta)\otimes_{\cR}\cR(\theta^{n-1})\hookrightarrow \cR(\theta^{-1})\otimes_{\cR}\cR(\theta)\otimes_{\cR}\cD^{-}\hookrightarrow \cR(\theta^{-1})\otimes_{\cR}\cB\otimes_{\cR}\cD^{-},\]
and the composition of the surjections
\[\cR(\theta^{-1})\otimes_{\cR}\cB\otimes_{\cR}\cD^{-}\twoheadrightarrow \cD\twoheadrightarrow \cR(1_{K^{\times}})\]
factors through
\[\cR(\theta^{-1})\otimes_{\cR}\cB\otimes_{\cR}\cD^{-}\twoheadrightarrow \cR(\theta^{-1})\otimes_{\cR}\cR(1_{K^{\times}})\otimes_{\cR}\cD^{-} \twoheadrightarrow \cR(\theta^{-1})\otimes_{\cR}\cR(1_{K^{\times}})\otimes_{\cR}\cR(\theta)\buildrel\sim\over\longrightarrow \cR(1_{K^{\times}}),\] 
we know that $\zeta_{n-1,0}(\varsigma(\tld{\cD}^{-}))$ is a subquotient of the following subquotient \[(\cR(\theta^{-1})\otimes_{\cR}\cR(\theta)\otimes_{\cR}\cD^{-})\oplus (\cR(\theta^{-1})\otimes_{\cR}\cR(1_{K^{\times}})\otimes_{\cR}\cD^{-})\]
of $\cR(\theta^{-1})\otimes_{\cR}\cB\otimes_{\cR}\tld{\cD}^{-}$ obtained from the push-out along $\cR(\theta^{-1})\otimes_{\cR}\cR(\theta)\otimes_{\cR}\cD^{-}\hookrightarrow \cR(\theta^{-1})\otimes_{\cR}\cB\otimes_{\cR}\cD^{-}$ and the pull-back along $\cR(\theta^{-1})\otimes_{\cR}\cB\otimes_{\cR}\cD^{-}\twoheadrightarrow \cR(\theta^{-1})\otimes_{\cR}\cR(1_{K^{\times}})\otimes_{\cR}\cD^{-}$. In other words, we have shown that 
\begin{equation}\label{equ: sym det top split}
\mathrm{im}(\vartheta)\subseteq \mathrm{ker}(\zeta_{n-1,0}).
\end{equation}
We choose an arbitrary $x\in \mathrm{Ext}_{\varphi,\Gamma}^1(\cD,\cD)\setminus \mathrm{ker}(\zeta_{n-1,0})$ and consider
\[x-\vartheta(\varsigma(x))\in \mathrm{Ext}_{\varphi,\Gamma}^1(\cD,\cD).\]
Now that we have $\varsigma\circ\vartheta=\mathrm{Id}$ and (\ref{equ: sym det top split}),
we conclude that $\varsigma(x-\vartheta(\varsigma(x)))=0$ and that $\zeta_{n-1,0}(x-\vartheta(\varsigma(x)))=\zeta_{n-1,0}(x)\neq 0$.
In other words, $x-\vartheta(\varsigma(x))$ gives a non-zero element of (\ref{equ: sym det transverse}) and the proof of \textbf{Step $1$} is finished.

\textbf{Step $2$}: We prove that
\begin{equation}\label{equ: sym det inclusion 1}
\mathrm{ker}(\varsigma)\subseteq \mathrm{Ext}_{\cD^{+}}^1(\cD,\cD).
\end{equation}

If $\tld{\cD}$ is associated with an element of $\mathrm{ker}(\varsigma)$, then $\tld{\cD}^{\dagger}$ and thus $\tld{\cD}^{\ddagger}$ are both sub of $\cB\otimes_{\cR}\tld{\cD}$. Let $\tld{\cD}''$ be the image of $\tld{\cD}^{\ddagger}$ under (\ref{equ: sym top}).
It is easy to check from \ref{it: char dR Ext 1} of Lemma~\ref{lem: char dR Ext} (together with a d\'evissage) that
\[\Dim_E \Hom_{\varphi,\Gamma}(\cR(\theta^{n-1}),\cD)=1\]
and from Lemma~\ref{lem: Galois unique Hom} that
\[\Dim_E \Hom_{\varphi,\Gamma}(\cD^{+},\cD)=1.\]
Hence, the composition of
\[\cR(\theta^{n-1})\hookrightarrow \cD^{-}\hookrightarrow \cB\otimes_{\cR}\cD \twoheadrightarrow \cD\]
is the unique (up to scalars) non-zero map $\cR(\theta^{n-1})\hookrightarrow \cD$ (which is an embedding), and the composition of
\[\cD^{+}\hookrightarrow \cB\otimes_{\cR}\cD \twoheadrightarrow \cD\]
is the unique (up to scalars) non-zero map $\cD^{+}\hookrightarrow \cD$ which is a surjection with kernel $\cR(\theta^{n})$.
Consequently, $\tld{\cD}''$ must be exactly the sub $\tld{\cD}'$ of $\tld{\cD}$ given by pull back along $\cR(\theta^{n-1})\hookrightarrow \cD$. Consequently, the very existence of $\tld{\cD}^{\ddagger}$ (which admits $\cD^{+}$ as a sub and $\tld{\cD}'$ as a quotient) ensures that the image of $\mathrm{ker}(\varsigma)$ in $\mathrm{Ext}_{\varphi,\Gamma}^1(\cR(\theta^{n-1}),\cD)$ is contained in $\mathrm{Ext}_{\cD^{+}}^1(\cR(\theta^{n-1}),\cD)$. This together with the definition of $\mathrm{Ext}_{\cD^{+}}^1(\cD,\cD)$ gives (\ref{equ: sym det inclusion 1}).

\textbf{Step $3$}: We prove that
\begin{equation}\label{equ: sym det inclusion 2}
\mathrm{ker}(\varsigma)\subseteq \mathrm{Ext}_{\mathrm{det}}^1(\cD,\cD).
\end{equation}

Using
\[\cB\otimes_{\cR}\mathrm{Sym}^{m}(\cB)=\mathrm{Sym}^{m+1}(\cB)\oplus(\mathrm{det}(\cB)\otimes_{\cR}\mathrm{Sym}^{m-1}(\cB)),\]
we see that $\cB^{\otimes m-1}\otimes_{\cR}\mathrm{det}(\cB)^{\otimes n-m}\otimes_{\cR}\mathrm{Sym}^{m-1}(\cB)$
is a canonical direct summand of $\cB^{\otimes m}\otimes_{\cR}\mathrm{det}(\cB)^{\otimes n-1-m}\otimes_{\cR}\mathrm{Sym}^{m}(\cB)$ for each $1\leq m\leq n-1$.
In particular, $\mathrm{det}(\cB)^{\otimes n-1}\cong \mathrm{det}(\cD)$ is a canonical direct summand of $\cB^{\otimes n-2}\otimes_{\cR}\cD^{-}$ which is itself a canonical direct summand of $\cB^{\otimes n-1}\otimes_{\cR}\cD$. Consequently, the map $\kappa_{\mathrm{det}}'$ from (\ref{equ: pass to det deform 2}) factors through $\mathrm{Ext}_{\varphi,\Gamma}^1(\cD^{-},\cD^{-})$. In particular, any element of $\mathrm{ker}(\varsigma)$ is mapped to zero under $\kappa_{\mathrm{det}}'$ or equivalently under $\kappa_{\mathrm{det}}$ (using Lemma~\ref{lem: det cond split}). This finishes the proof of (\ref{equ: sym det inclusion 2}).

It is clear that (\ref{equ: sym det generator}) follows from \textbf{Step $1$}, \textbf{Step $2$} and \textbf{Step $3$}.
\end{proof}

Let $n\geq 1$ and $\cD$ be a de Rham $(\varphi,\Gamma)$-module over $\cR$ which is special with level $n-1$.
Recall from \ref{it: special mod 2} of Definition~\ref{def: special mod} that $\cD$ is equipped with a decreasing filtration $\mathrm{Fil}_{W}^{\bullet}(\cD)$ with $\mathrm{gr}_{W}^{k}(\cD)\cong \cR(\theta^{k})$ for each $0\leq k\leq n-1$, and $\mathrm{gr}_{W}^{k}(\cD)=0$ otherwise.
For each $0\leq k,\ell\leq n-1$, the injection $\mathrm{Fil}_{W}^{k}(\cD)\hookrightarrow \cD$ and the surjection $\cD\twoheadrightarrow \cD/\mathrm{Fil}_{W}^{\ell+1}(\cD)$ induces a map
\[\zeta_{k,\ell}: \mathrm{Ext}_{\varphi,\Gamma}^1(\cD,\cD)\rightarrow \mathrm{Ext}_{\varphi,\Gamma}^1(\mathrm{Fil}_{W}^{k}(\cD),\cD/\mathrm{Fil}_{W}^{\ell+1}(\cD)).\]
For each $\al=(i,j)\in\Phi^+$ for some $1\leq i<j\leq n$, we define $\mathrm{Fil}_{\al}(\mathrm{Ext}_{\varphi,\Gamma}^1(\cD,\cD))$ as the intersection of $\mathrm{\ker}(\zeta_{n-h-1,n-h})$ for all $h$ satisfying either $1\leq h\leq i-1$ or $j\leq h\leq n-1$.
Let $S\subseteq \Phi^+$ be a subset such that $\beta\in S$ for each $\al\in S$ and $\beta\in\Phi^+$ satisfying $\beta\leq \al$. Then we define
\[\mathrm{Fil}_{S}(\mathrm{Ext}_{\varphi,\Gamma}^1(\cD,\cD))\defeq \sum_{\al\in S}\mathrm{Fil}_{\al}(\mathrm{Ext}_{\varphi,\Gamma}^1(\cD,\cD)).\]
When $S=\emptyset$, we have
\begin{equation}\label{equ: tri deform}
\mathrm{Fil}_{\emptyset}(\mathrm{Ext}_{\varphi,\Gamma}^1(\cD,\cD))=\mathrm{Ext}_{\mathrm{tri}}^1(\cD,\cD)
\end{equation}
with RHS being the $E$-subspace of $\mathrm{Ext}_{\varphi,\Gamma}^1(\cD,\cD)$ consisting of \emph{trianguline deformations of $\cD$}, namely those $\tld{\cD}$ admitting a decreasing filtration $\mathrm{Fil}_{W}^{\bullet}(\tld{\cD})$ with $\mathrm{gr}_{W}^{k}(\tld{\cD})$ being a deformation of $\mathrm{gr}_{W}^{k}(\cD)\cong \cR(\theta^{k})$ for each $0\leq k\leq n-1$.
For each $\al\in\Phi^+$, we write
\[\mathrm{gr}_{\al}(\mathrm{Ext}_{\varphi,\Gamma}^1(\cD,\cD))\defeq \mathrm{Fil}_{\al}(\mathrm{Ext}_{\varphi,\Gamma}^1(\cD,\cD))/\sum_{\beta<\al}\mathrm{Fil}_{\beta}(\mathrm{Ext}_{\varphi,\Gamma}^1(\cD,\cD))\]
for short.

Following the argument of \cite[Lem.~3.9]{BD20}, we see that the composition of
\[\mathrm{Ext}_{\mathrm{tri}}^1(\cD,\cD)\hookrightarrow \mathrm{Ext}_{\varphi,\Gamma}^1(\cD,\cD)\buildrel \kappa_{\mathrm{det}}\over\longrightarrow \mathrm{Ext}_{\varphi,\Gamma}^1(\mathrm{det}(\cD),\mathrm{det}(\cD))\cong\Hom(\Q_p^\times,E)\]
is surjective, and thus (\ref{equ: pass to det deform 1}) is surjective.
In particular, we have
\begin{equation}\label{equ: fix det deform dim}
\Dim_E \mathrm{Ext}_{\mathrm{det}}^1(\cD,\cD)=\Dim_E \mathrm{Ext}_{\varphi,\Gamma}^1(\cD,\cD)-2.
\end{equation}
We also define
\[\mathrm{Ext}_{\mathrm{tri},\mathrm{det}}^1(\cD,\cD)\defeq \mathrm{Ext}_{\mathrm{det}}^1(\cD,\cD)\cap \mathrm{Ext}_{\mathrm{tri}}^1(\cD,\cD)\]
and note that
\begin{equation}\label{equ: fix det tri deform dim}
\Dim_E \mathrm{Ext}_{\mathrm{tri},\mathrm{det}}^1(\cD,\cD)=\Dim_E \mathrm{Ext}_{\mathrm{tri}}^1(\cD,\cD)-2.
\end{equation}

When $\ast$ is the deformation condition given by either $\mathrm{tri}$ or $\mathrm{det}$ or $\mathrm{tri},\mathrm{det}$, we write
\begin{equation}\label{equ: ast deform subquotient}
\mathrm{Ext}_{\ast}^1(\cD/\cR(\theta^{n-1}),\cD)\defeq \mathrm{Ext}_{\varphi,\Gamma}^1(\cD/\cR(\theta^{n-1}),\cD)\cap \mathrm{Ext}_{\ast}^1(\cD,\cD)\subseteq \mathrm{Ext}_{\varphi,\Gamma}^1(\cD,\cD)
\end{equation}
and note that we have a Cartesian diagram of the form
\begin{equation}\label{equ: ast deform Cartesian square}
\xymatrix{
\mathrm{Ext}_{\ast}^1(\cD/\cR(\theta^{n-1}),\cD) \ar[r] \ar@{^{(}->}[d]& \mathrm{Ext}_{\ast}^1(\cD/\cR(\theta^{n-1}),\cD/\cR(\theta^{n-1})) \ar@{^{(}->}[d]\\
\mathrm{Ext}_{\varphi,\Gamma}^1(\cD/\cR(\theta^{n-1}),\cD) \ar[r] & \mathrm{Ext}_{\varphi,\Gamma}^1(\cD/\cR(\theta^{n-1}),\cD/\cR(\theta^{n-1}))
}.
\end{equation}
Consequently, the short exact sequence (\ref{equ: deform induction seq}) restricts to the following exact sequence
\begin{equation}\label{equ: ast deform induction seq}
0\rightarrow \mathrm{Ext}_{\ast}^1(\cD/\cR(\theta^{n-1}),\cD)\rightarrow \mathrm{Ext}_{\ast}^1(\cD,\cD)\buildrel q_{\ast} \over\longrightarrow \mathrm{Ext}_{\varphi,\Gamma}^1(\cR(\theta^{n-1}),\cR(\theta^{n-1})),
\end{equation}
and the long exact sequence (\ref{equ: deform induction devissage}) restricts to the following long exact sequence
\begin{multline}\label{equ: ast deform induction devissage}
0\rightarrow \Hom_{\varphi,\Gamma}(\cD/\cR(\theta^{n-1}),\cD/\cR(\theta^{n-1}))\\
\rightarrow \mathrm{Ext}_{\varphi,\Gamma}^1(\cD/\cR(\theta^{n-1}),\cR(\theta^{n-1}))
\rightarrow \mathrm{Ext}_{\ast}^1(\cD/\cR(\theta^{n-1}),\cD) \rightarrow \mathrm{Ext}_{\ast}^1(\cD/\cR(\theta^{n-1}),\cD/\cR(\theta^{n-1}))\\
\buildrel q_{\ast}' \over\longrightarrow \mathrm{Ext}_{\varphi,\Gamma}^2(\cD/\cR(\theta^{n-1}),\cR(\theta^{n-1})).
\end{multline}

\begin{lem}\label{lem: full deform dim}
Let $n\geq 1$ and $\cD$ be a de Rham $(\varphi,\Gamma)$-module which is special with level $n-1$. Then we have the following results.
\begin{enumerate}[label=(\roman*)]
\item \label{it: full deform dim 1} We have
\begin{equation}\label{equ: full deform dim}
\Dim_E \mathrm{Ext}_{\varphi,\Gamma}^1(\cD,\cD)=1+n^2.
\end{equation}
\item \label{it: full deform dim 2} We have
\begin{equation}\label{equ: full deform tri dim}
\Dim_E \mathrm{Ext}_{\mathrm{tri}}^1(\cD,\cD)=1+\frac{n(n+1)}{2}.
\end{equation}
\item \label{it: full deform dim 3} We have
\begin{equation}\label{equ: full deform grade}
\Dim_E \mathrm{gr}_{\al}(\mathrm{Ext}_{\varphi,\Gamma}^1(\cD,\cD))=1
\end{equation}
for each $\al\in\Phi^+$.
\end{enumerate}
\end{lem}
\begin{proof}
We prove \ref{it: full deform dim 1} by an increasing induction on $n\geq 1$.\\
The $n=1$ case is clear from $\Dim_E \mathrm{Ext}_{\varphi,\Gamma}^1(\cR(1_{K^{\times}}),\cR(1_{K^{\times}}))=2$ by \ref{it: Galois Ext collection 0} of Lemma~\ref{lem: Galois Ext collection}.
Assume from now that $n\geq 2$.
Recall from the discussion around (\ref{equ: deform induction seq Hom embedding}) (upon replacing $\cD$ in \emph{loc.cit.} with $\cD/\cR(\theta^{n-1})$) that we have
\[\Dim_E\Hom_{\varphi,\Gamma}(\cD/\cR(\theta^{n-1}),\cD/\cR(\theta^{n-1}))=1.\]
Recall from (\ref{equ: BD pairing}) and its variant with $n-1$ replacing $n$ that we have
\[
\left\{\begin{array}{cccc}
\Dim_E & \mathrm{Ext}_{\varphi,\Gamma}^2(\cD/\cR(\theta^{n-1}),\cR(\theta^{n-1})) & = & 1\\
\Dim_E & \mathrm{Ext}_{\varphi,\Gamma}^1(\cD/\cR(\theta^{n-1}),\cR(\theta^{n-1})) & = & n\\
\Dim_E & \mathrm{Ext}_{\varphi,\Gamma}^1(\cR(\theta^{n-1}),\cD) & = &n+1
\end{array}\right.
\]
We also have
\[\Dim_E \mathrm{Ext}_{\varphi,\Gamma}^1(\cD/\cR(\theta^{n-1}),\cD/\cR(\theta^{n-1}))=1+(n-1)^2\]
by our induction hypothesis.
Hence, we deduce from (\ref{equ: deform induction seq}) and (\ref{equ: deform induction devissage}) that
\begin{multline}\label{equ: deform dim induction}
\Dim_E \mathrm{Ext}_{\varphi,\Gamma}^1(\cD,\cD)=\Dim_E \mathrm{Ext}_{\varphi,\Gamma}^1(\cD/\cR(\theta^{n-1}),\cD)+\Dim_E \mathrm{Ext}_{\varphi,\Gamma}^1(\cR(\theta^{n-1}),\cD)\\
=\Dim_E \mathrm{Ext}_{\varphi,\Gamma}^1(\cD/\cR(\theta^{n-1}),\cD/\cR(\theta^{n-1}))+\Dim_E \mathrm{Ext}_{\varphi,\Gamma}^1(\cD/\cR(\theta^{n-1}),\cR(\theta^{n-1}))+\Dim_E \mathrm{Ext}_{\varphi,\Gamma}^1(\cR(\theta^{n-1}),\cD)\\
-\Dim_E \Hom_{\varphi,\Gamma}(\cD/\cR(\theta^{n-1}),\cD/\cR(\theta^{n-1}))-\Dim_E \mathrm{Ext}_{\varphi,\Gamma}^2(\cD/\cR(\theta^{n-1}),\cR(\theta^{n-1}))\\
=\Dim_E \mathrm{Ext}_{\varphi,\Gamma}^1(\cD/\cR(\theta^{n-1}),\cD/\cR(\theta^{n-1}))+2n-1=1+(n-1)^2+2n-1=1+n^2
\end{multline}
which finishes the proof of the induction step.

We write $q\defeq q_{\mathrm{tri}}$ and $q'\defeq q_{\mathrm{tri}}'$ for short (see (\ref{equ: ast deform induction seq}) and (\ref{equ: ast deform induction devissage})).
We prove \ref{it: full deform dim 2}, \ref{it: full deform dim 3} and that $q$ is surjective by an increasing induction on $n\geq 1$.\\
The $n=1$ case of \ref{it: full deform dim 2}, \ref{it: full deform dim 3} and the surjectivity of $q$ are clear. We assume from now that $n\geq 2$. By our induction hypothesis the map
\[\mathrm{Ext}_{\mathrm{tri}}^1(\cD/\cR(\theta^{n-1}),\cD/\cR(\theta^{n-1})) \rightarrow \mathrm{Ext}_{\varphi,\Gamma}^1(\cR(\theta^{n-2}),\cR(\theta^{n-2}))\]
is surjective. Now that the cup product map
\[\mathrm{Ext}_{\varphi,\Gamma}^1(\cR(\theta^{n-2}),\cR(\theta^{n-2}))\otimes_E \mathrm{Ext}_{\varphi,\Gamma}^1(\cR(\theta^{n-2}),\cR(\theta^{n-1}))\buildrel\cup\over\longrightarrow \mathrm{Ext}_{\varphi,\Gamma}^2(\cR(\theta^{n-2}),\cR(\theta^{n-1}))\]
is a perfect pairing between $2$ dimensional spaces, we deduce from (\ref{equ: BD pairing}) (upon replacing $n$ in \emph{loc.cit.} with $n-1$) that $q'$ (see (\ref{equ: ast deform induction devissage}) with $\ast=\mathrm{tri}$) is surjective. Similar to (\ref{equ: deform dim induction}), we deduce from the surjectivity of $q'$ as well as (\ref{equ: ast deform induction seq}) and (\ref{equ: ast deform induction devissage}) (with $\ast=\mathrm{tri}$ in \emph{loc.cit.}) that
\begin{equation}\label{equ: tri deform dim upper bound}
\Dim_E \mathrm{Ext}_{\mathrm{tri}}^1(\cD,\cD)\leq \Dim_E \mathrm{Ext}_{\mathrm{tri}}^1(\cD/\cR(\theta^{n-1}),\cD/\cR(\theta^{n-1}))+n=1+\frac{n(n+1)}{2}.
\end{equation}
By the definition of $\mathrm{Fil}_{\al}(\mathrm{Ext}_{\varphi,\Gamma}^1(\cD,\cD))$ and \ref{it: Galois Ext collection 0} of Lemma~\ref{lem: Galois Ext collection}, we have
\begin{equation}\label{equ: deform grade upper bound}
\Dim_E\mathrm{gr}_{\al}(\mathrm{Ext}_{\varphi,\Gamma}^1(\cD,\cD))\leq \Dim_E \mathrm{Ext}_{\varphi,\Gamma}^1(\cR(\theta^{n-i}),\cR(\theta^{n-j}))=1
\end{equation}
for each $\al=(i,j)\in\Phi^+$, and thus
\begin{multline*}
\Dim_E \mathrm{Ext}_{\varphi,\Gamma}^1(\cD,\cD)\leq \Dim_E \mathrm{Ext}_{\mathrm{tri}}^1(\cD,\cD)+\sum_{\al\in\Phi^+}\Dim_E\mathrm{gr}_{\al}(\mathrm{Ext}_{\varphi,\Gamma}^1(\cD,\cD))\\
\leq 1+\frac{n(n+1)}{2}+\frac{n(n-1)}{2}=1+n^2.
\end{multline*}
This together with \ref{it: full deform dim 1} forces both inequalities (\ref{equ: tri deform dim upper bound}) and (\ref{equ: deform grade upper bound}) to be equalities, which finishes the proof the induction step.
\end{proof}

By the definition of $\mathrm{Ext}_{\mathrm{tri}}^1(\cD,\cD)$, we notice that the RHS surjection of (\ref{equ: deform induction seq}) restricts to a surjection
\begin{equation}\label{equ: tri deform surjection}
\mathrm{Ext}_{\mathrm{tri}}^1(\cD,\cD)\twoheadrightarrow \mathrm{Ext}_{\mathrm{tri}}^1(\mathrm{Fil}_{W}^{m}(\cD),\mathrm{Fil}_{W}^{m}(\cD))
\end{equation}
for each $1\leq m\leq n-1$.
Recall that both
\begin{equation}\label{equ: fund pair 1}
\mathrm{Ext}_{\varphi,\Gamma}^1(\cR(1_{K^{\times}}),\cR(1_{K^{\times}}))\otimes_E\mathrm{Ext}_{\varphi,\Gamma}^1(\cR(1_{K^{\times}}),\cR(\theta))\buildrel\cup\over\longrightarrow \mathrm{Ext}_{\varphi,\Gamma}^2(\cR(1_{K^{\times}}),\cR(\theta))
\end{equation}
and
\begin{equation}\label{equ: fund pair 2}
\mathrm{Ext}_{\varphi,\Gamma}^1(\cR(1_{K^{\times}}),\cR(\theta))\otimes_E\mathrm{Ext}_{\varphi,\Gamma}^1(\cR(\theta),\cR(\theta))\buildrel\cup\over\longrightarrow \mathrm{Ext}_{\varphi,\Gamma}^2(\cR(1_{K^{\times}}),\cR(\theta))
\end{equation}
are perfect pairing between $2$ dimensional $E$-vector spaces, and they are related by a dual and a twist by $\cR(\theta)$.
Let $\cB$ be a de Rham $(\varphi,\Gamma)$-module which is special with level $1$. Then $\cB$ determines an $E$-line $E[\cB]\subseteq \mathrm{Ext}_{\varphi,\Gamma}^1(\cR(1_{K^{\times}}),\cR(\theta))$ and there exists an $E$-line
\[E[\cB]^{\perp}\subseteq \mathrm{Ext}_{\varphi,\Gamma}^1(\cR(1_{K^{\times}}),\cR(1_{K^{\times}}))\cong \mathrm{Ext}_{\varphi,\Gamma}^1(\cR(\theta),\cR(\theta))\]
which is both the left orthogonal complement of $E[\cB]$ under (\ref{equ: fund pair 1}) and the right orthogonal complement of $E[\cB]$ under (\ref{equ: fund pair 2}).
\begin{lem}\label{lem: tri det deform image}
Let $n\geq 2$ and $\cB$ be a de Rham $(\varphi,\Gamma)$-module which is special with level $1$ with $\cD\defeq \mathrm{Sym}^{n-1}(\cB)$. Then the following map (using (\ref{equ: tri deform surjection}) with $m=n-1$)
\begin{equation}\label{equ: tri det deform image}
\mathrm{Ext}_{\mathrm{tri},\mathrm{det}}^1(\cD,\cD)\rightarrow \mathrm{Ext}_{\varphi,\Gamma}^1(\cR(\theta^{n-1}),\cR(\theta^{n-1}))\cong \mathrm{Ext}_{\varphi,\Gamma}^1(\cR(1_{K^{\times}}),\cR(1_{K^{\times}}))
\end{equation}
has image $E[\cB]^{\perp}$.
\end{lem}
\begin{proof}
We write $q_{\flat}\defeq q_{\mathrm{tri},\mathrm{det}}$ and $q_{\flat}'\defeq q_{\mathrm{tri},\mathrm{det}}'$ for short (see (\ref{equ: ast deform induction seq}) and (\ref{equ: ast deform induction devissage})).
We prove that (\ref{equ: tri det deform image}) has image $E[\cB]^{\perp}$ by an increasing induction on $n\geq 2$.

We first treat the $n=2$ case.
Note that $\mathrm{Ext}_{\mathrm{tri},\mathrm{det}}^1(\cD/\cR(\theta^{n-1}),\cD/\cR(\theta^{n-1}))=0$ in this case, which together with (\ref{equ: ast deform induction seq}) and (\ref{equ: ast deform induction devissage}) (with $\ast=\mathrm{tri},\mathrm{det}$ in \emph{loc.cit.}) gives
\[\Dim_E \mathrm{Ext}_{\mathrm{tri},\mathrm{det}}^1(\cD/\cR(\theta^{n-1}),\cD)=1\]
and thus (using \ref{it: full deform dim 2} of Lemma~\ref{lem: full deform dim} and (\ref{equ: fix det tri deform dim}))
\begin{equation}\label{equ: tri deform image n 2}
\Dim_E \mathrm{im}(q_{\flat})=\Dim_E \mathrm{Ext}_{\mathrm{tri},\mathrm{det}}^1(\cD,\cD)-\Dim_E \mathrm{Ext}_{\mathrm{tri},\mathrm{det}}^1(\cD/\cR(\theta^{n-1}),\cD)=1.
\end{equation}
Now we consider an element of $\mathrm{Ext}_{\mathrm{tri}}^1(\cD,\cD)$ namely a trianguline deformation $\tld{\cD}$ of $\cD$. Note that $\tld{\cD}$ determines a deformation of $\cR(1_{K^{\times}})$ and $\cR(\theta)$ respectively, and thus an element
\[\psi_k\in \mathrm{Ext}_{\varphi,\Gamma}^1(\cR(\theta^{k}),\cR(\theta^{k}))\cong \Hom(\Q_p^{\times},E)\]
for $k=0,1$.
We abuse the notation $\psi_k$ for the corresponding map
\[\mathrm{Ext}_{\mathrm{tri}}^1(\cD,\cD)\rightarrow \Hom(\Q_p^{\times},E)\]
for $k=0,1$
The discussion below (\ref{equ: fund pair 1}) and (\ref{equ: fund pair 2}) shows that there exist $x_0$ (resp.~$x_1$) in $\mathrm{Ext}_{\mathrm{tri}}^1(\cD,\cD)$ that satisfies $0\neq \psi_0(x_0)\in E[\cB]^{\perp}$ and $\psi_1(x_0)=0$ (resp.~$\psi_0(x_1)=0$ and $0\neq \psi_1(x_1)\in E[\cB]^{\perp}$). We choose arbitrary $a,b\in E^{\times}$ such that $a\psi_0(x_0)+b\psi_1(x_1)=0$ and thus obtain
\[ax_0+bx_1\in \mathrm{Ext}_{\mathrm{tri},\mathrm{det}}^1(\cD,\cD)\]
that satisfies $\psi_k(ax_0+bx_1)\in E[\cB]^{\perp}$ for $k=0,1$. We conclude that $E[\cB]^{\perp}\subseteq \mathrm{im}(q_{\flat})$, which together with (\ref{equ: tri deform image n 2}) finishes the proof of the $n=2$ case.

Assume from now that $n\geq 3$.
By induction hypothesis we have
\begin{equation}\label{equ: tri det inclusion induction}
\mathrm{Ext}_{\mathrm{tri},\mathrm{det}}^1(\cD/\cR(\theta^{n-1}),\cD/\cR(\theta^{n-1}))\subseteq \mathrm{Ext}_{\cD}^1(\cD/\cR(\theta^{n-1}),\cD/\cR(\theta^{n-1})),
\end{equation}
which forces $q_{\flat}'=0$.
Consequently, we have
\[\Dim_E \mathrm{Ext}_{\mathrm{tri},\mathrm{det}}^1(\cD/\cR(\theta^{n-1}),\cD)=\Dim_E \mathrm{Ext}_{\mathrm{tri},\mathrm{det}}^1(\cD/\cR(\theta^{n-1}),\cD/\cR(\theta^{n-1}))+n-1,\]
which together with (\ref{equ: fix det deform dim}) and \ref{it: full deform dim 2} of Lemma~\ref{lem: full deform dim} forces
\[\Dim_E \mathrm{Ext}_{\mathrm{tri},\mathrm{det}}^1(\cD/\cR(\theta^{n-1}),\cD)=\Dim_E \mathrm{Ext}_{\mathrm{tri},\mathrm{det}}^1(\cD,\cD)-1.\]
Compared with (\ref{equ: ast deform induction seq}) (with $\ast=\mathrm{tri},\mathrm{det}$ in \emph{loc.cit.}), we observe that
\begin{equation}\label{equ: tri det deform image dim}
\Dim_E \mathrm{im}(q_{\flat})=1.
\end{equation}
Now we consider an arbitrary deformation $\tld{\cB}$ of $\cB$ associated with an element of $\mathrm{Ext}_{\mathrm{tri},\mathrm{det}}^1(\cB,\cB)$. Now that its image in $\mathrm{Ext}_{\mathrm{tri}}^1(\cR(\theta),\cR(\theta))$ lies in $E[\cB]^{\perp}$ by previous discussion on the $n=2$ case, there exists $\tld{\cB}^+$ that fits into $0\rightarrow \cD\rightarrow \tld{\cB}^+\rightarrow \cB\rightarrow 0$ whose push forward along $\cD\twoheadrightarrow \cB$ gives $\tld{\cB}$. We define $\tld{\cD}$ as the pull back of $\tld{\cB}^+$ along $\cD\twoheadrightarrow \cB$ and note that $\tld{\cD}$ indeed arises from an element of $\mathrm{Ext}_{\mathrm{tri},\mathrm{det}}^1(\cD,\cD)$. Now that $\cD^{\ast}\otimes_{\cR}\cR(\theta^{n-1})\cong \cD$ with $\cD^{\ast}$ being the dual of $\cD$, we see that $\tld{\cD}^{\ast}\otimes_{\cR}\cR(\theta^{n-1})$ is a deformation of $\cD$ arising from an element of $\mathrm{Ext}_{\mathrm{tri},\mathrm{det}}^1(\cD,\cD)$. Since the map (see (\ref{equ: tri deform surjection}))
\begin{equation}\label{equ: tri det deform map}
\mathrm{Ext}_{\mathrm{tri}}^1(\cD,\cD)\rightarrow \mathrm{Ext}_{\varphi,\Gamma}^1(\cR(\theta^{n-1}),\cR(\theta^{n-1}))
\end{equation}
factors through
\[\mathrm{Ext}_{\mathrm{tri}}^1(\cB\otimes_{\cR}\cR(\theta^{n-2}),\cB\otimes_{\cR}\cR(\theta^{n-2}))\cong \mathrm{Ext}_{\mathrm{tri}}^1(\cB^{\ast}\otimes_{\cR}\cR(\theta^{n-1}),\cB^{\ast}\otimes_{\cR}\cR(\theta^{n-1})),\]
we see that the image of $\tld{\cD}^{\ast}\otimes_{\cR}\cR(\theta^{n-1})$ under (\ref{equ: tri det deform map}) equals $E[\cB]^{\perp}$. In other words, we have shown that $E[\cB]^{\perp}\subseteq \mathrm{im}(q_{\flat})$, which together with (\ref{equ: tri det deform image dim}) finishes the proof of the induction step.
\end{proof}

Let $n\geq 2$ and $\cD$ be a de Rham $(\varphi,\Gamma)$-module which is special with level $n-1$ (see \ref{it: special mod 2} of Definition~\ref{def: special mod}).
We consider the following maps
\begin{equation}\label{equ: det deform image}
\mathrm{Ext}_{\mathrm{det}}^1(\cD,\cD)\hookrightarrow \mathrm{Ext}_{\varphi,\Gamma}^1(\cD,\cD)\twoheadrightarrow \mathrm{Ext}_{\varphi,\Gamma}^1(\cR(\theta^{n-1}),\cD).
\end{equation}
\begin{lem}\label{lem: sym det deform}
Let $n\geq 2$ and $\cB$ be a de Rham $(\varphi,\Gamma)$-module which is special with level $1$.
Let $\cD\defeq \mathrm{Sym}^{n-1}(\cB)$ and $\cD^{+}\defeq \mathrm{Sym}^{n}(\cB)$.
Then we have
\begin{equation}\label{equ: sym det deform}
\mathrm{Ext}_{\mathrm{det}}^1(\cD,\cD)\subseteq \mathrm{Ext}_{\cD^{+}}^1(\cD,\cD),
\end{equation}
and the composition of (\ref{equ: det deform image}) has image $\mathrm{Ext}_{\cD^{+}}^1(\cR(\theta^{n-1}),\cD)$.
\end{lem}
\begin{proof}
We write $q_{\natural}\defeq q_{\mathrm{det}}$ and $q_{\natural}'\defeq q_{\mathrm{det}}'$ for short (see (\ref{equ: ast deform induction seq}) and (\ref{equ: ast deform induction devissage})).
We prove by an increasing induction on $n\geq 2$ the inclusion (\ref{equ: sym det deform}) as well as the fact that the composition of (\ref{equ: det deform image}) has image $\mathrm{Ext}_{\cD^{+}}^1(\cR(\theta^{n-1}),\cD)$.
Now that we have
\begin{equation}\label{equ: det inclusion induction}
\mathrm{Ext}_{\mathrm{det}}^1(\cD/\cR(\theta^{n-1}),\cD/\cR(\theta^{n-1}))\subseteq \mathrm{Ext}_{\cD}^1(\cD/\cR(\theta^{n-1}),\cD/\cR(\theta^{n-1}))
\end{equation}
by our induction hypothesis, we have $q_{\natural}'=0$ and thus the following long exact sequence
\begin{multline}\label{equ: det deform induction devissage}
0\rightarrow \Hom_{\varphi,\Gamma}(\cD/\cR(\theta^{n-1}),\cD/\cR(\theta^{n-1}))\\
\rightarrow \mathrm{Ext}_{\varphi,\Gamma}^1(\cD/\cR(\theta^{n-1}),\cR(\theta^{n-1}))
\rightarrow \mathrm{Ext}_{\mathrm{det}}^1(\cD/\cR(\theta^{n-1}),\cD) \rightarrow \mathrm{Ext}_{\mathrm{det}}^1(\cD/\cR(\theta^{n-1}),\cD/\cR(\theta^{n-1}))\rightarrow 0.
\end{multline}
We deduce from (\ref{equ: det deform induction devissage}) that
\[\Dim_E \mathrm{Ext}_{\mathrm{det}}^1(\cD/\cR(\theta^{n-1}),\cD)=\Dim_E \mathrm{Ext}_{\mathrm{det}}^1(\cD/\cR(\theta^{n-1}),\cD/\cR(\theta^{n-1}))+n-1,\]
which together with (\ref{equ: fix det deform dim}) and \ref{it: full deform dim 1} of Lemma~\ref{lem: full deform dim} forces
\[\Dim_E \mathrm{Ext}_{\mathrm{det}}^1(\cD/\cR(\theta^{n-1}),\cD)=\Dim_E \mathrm{Ext}_{\mathrm{det}}^1(\cD,\cD)-n.\]
Compared with (\ref{equ: ast deform induction seq}) (with $\ast=\mathrm{det}$ in \emph{loc.cit.}), we observe that
\[
\Dim_E \mathrm{im}(q_{\natural})=n.
\]
Now we construct a basis $\{x_k\}_{0\leq k\leq n-1}$ of $\mathrm{im}(q_{\natural})$.
We define $x_{n-1}$ as an arbitrary non-zero element in $E[\cB\otimes_{\cR}\cR(\theta^{n-2})]^{\perp}$ which is contained in
\[\mathrm{im}(q_{\natural})\cap \mathrm{Ext}_{\varphi,\Gamma}^1(\cR(\theta^{n-1}),\cR(\theta^{n-1}))\subseteq \mathrm{Ext}_{\varphi,\Gamma}^1(\cR(\theta^{n-1}),\cD)\]
by Lemma~\ref{lem: tri det deform image}.
For each $0\leq k\leq n-2$, we note that $\mathrm{Fil}_{W}^{k}(\cD)\cong \mathrm{Sym}^{n-1-k}(\cB)\otimes_{\cR}\cR(\theta^{k})$ and define
\[y_k\in \mathrm{Ext}_{\varphi,\Gamma}^1(\mathrm{Fil}_{W}^{k}(\cD),\mathrm{Fil}_{W}^{k}(\cD))\cong \mathrm{Ext}_{\varphi,\Gamma}^1(\mathrm{Sym}^{n-1-k}(\cB),\mathrm{Sym}^{n-1-k}(\cB))\]
as an element in
\begin{equation}\label{equ: sym det Levi}
\mathrm{Ext}_{\mathrm{det}}^1(\mathrm{Sym}^{n-1-k}(\cB),\mathrm{Sym}^{n-1-k}(\cB))\cap \mathrm{Ext}_{\mathrm{Sym}^{n-k}(\cB)}^1(\mathrm{Sym}^{n-1-k}(\cB),\mathrm{Sym}^{n-1-k}(\cB))
\end{equation}
whose image in $\mathrm{Ext}_{\varphi,\Gamma}^1(\cR(\theta^{n-1-k}),\cR(1_{K^{\times}}))$ is non-zero, using Lemma~\ref{lem: sym det generator} upon replacing $\cD$ and $\cD^{+}$ in \emph{loc.cit.} with $\mathrm{Sym}^{n-1-k}(\cB)$ and $\mathrm{Sym}^{n-k}(\cB)$. When $k=0$, we set $\tld{y}_0\defeq y_0$.
When $k>0$, we have
\[\mathrm{Ext}_{\mathrm{det}}^1(\mathrm{Sym}^{n-1-k}(\cB),\mathrm{Sym}^{n-1-k}(\cB))\subseteq \mathrm{Ext}_{\mathrm{Sym}^{n-k}(\cB)}^1(\mathrm{Sym}^{n-1-k}(\cB),\mathrm{Sym}^{n-1-k}(\cB))\]
by our induction hypothesis. Upon taking a dual and a twist, we see that $\mathrm{Fil}_{W}^{k-1}(\cD)$ and $y_k$ are orthogonal under the top cup product map of the following commutative diagram (see \cite[(2.2)]{BD20} or the dual version of (\ref{equ: BD pairing}))
\[
\xymatrix{
\mathrm{Ext}_{\varphi,\Gamma}^1(\cR(\theta^{k-1}),\mathrm{Fil}_{W}^{k}(\cD)) \ar@{=}[d] & \otimes_E & \mathrm{Ext}_{\varphi,\Gamma}^1(\mathrm{Fil}_{W}^{k}(\cD),\mathrm{Fil}_{W}^{k}(\cD)) \ar^{\cup}[r] \ar@{->>}[d] & \mathrm{Ext}_{\varphi,\Gamma}^2(\cR(\theta^{k-1}),\mathrm{Fil}_{W}^{k}(\cD)) \ar^{\wr}[d]\\
\mathrm{Ext}_{\varphi,\Gamma}^1(\cR(\theta^{k-1}),\mathrm{Fil}_{W}^{k}(\cD)) & \otimes_E & \mathrm{Ext}_{\varphi,\Gamma}^1(\mathrm{Fil}_{W}^{k}(\cD),\cR(\theta^{k})) \ar^{\cup}[r] & \mathrm{Ext}_{\varphi,\Gamma}^2(\cR(\theta^{k-1}),\cR(\theta^{k}))
}
\]
with the middle vertical map being a surjection and the right vertical map being an isomorphism.
In particular, there exists $y_k'\in \mathrm{Ext}_{\varphi,\Gamma}^1(\mathrm{Fil}_{W}^{k-1}(\cD),\mathrm{Fil}_{W}^{k}(\cD))$ whose image under
\[\mathrm{Ext}_{\varphi,\Gamma}^1(\mathrm{Fil}_{W}^{k-1}(\cD),\mathrm{Fil}_{W}^{k}(\cD))\rightarrow \mathrm{Ext}_{\varphi,\Gamma}^1(\mathrm{Fil}_{W}^{k}(\cD),\mathrm{Fil}_{W}^{k}(\cD))\]
gives $y_k$. It follows from \ref{it: Galois Ext collection 1} of Lemma~\ref{lem: Galois Ext collection} and a d\'evissage that $\mathrm{Ext}_{\varphi,\Gamma}^2(\cD/\mathrm{Fil}_{W}^{k-1}(\cD),\mathrm{Fil}_{W}^{k}(\cD))=0$ and thus the embedding $\mathrm{Fil}_{W}^{k-1}(\cD)\hookrightarrow \cD$ induces a surjection
\begin{equation}\label{equ: sym det generator surjection}
\mathrm{Ext}_{\varphi,\Gamma}^1(\cD,\mathrm{Fil}_{W}^{k}(\cD))\twoheadrightarrow \mathrm{Ext}_{\varphi,\Gamma}^1(\mathrm{Fil}_{W}^{k-1}(\cD),\mathrm{Fil}_{W}^{k}(\cD)).
\end{equation}
Consequently, there exists
\[\tld{y}_k\in \mathrm{Ext}_{\varphi,\Gamma}^1(\cD,\mathrm{Fil}_{W}^{k}(\cD))\subseteq \mathrm{Ext}_{\varphi,\Gamma}^1(\cD,\cD)\]
whose image under (\ref{equ: sym det generator surjection}) is $y_k'$. Hence, the image of $\tld{y}_k$ in $\mathrm{Ext}_{\varphi,\Gamma}^1(\mathrm{Fil}_{W}^{k}(\cD),\mathrm{Fil}_{W}^{k}(\cD))$ gives back $y_{k}$, which forces
\begin{equation}\label{equ: tri det generator}
\tld{y}_k\in \mathrm{Ext}_{\mathrm{det}}^1(\cD,\cD).
\end{equation}
We define $x_k$ as the image of $\tld{y}_k$ under
\[\mathrm{Ext}_{\varphi,\Gamma}^1(\cD,\cD)\twoheadrightarrow \mathrm{Ext}_{\varphi,\Gamma}^1(\cR(\theta^{n-1}),\cD),\]
which by our choice of $y_k$ is an element of $\mathrm{Ext}_{\varphi,\Gamma}^1(\cR(\theta^{n-1}),\mathrm{Fil}_{W}^{k}(\cD))$ which has non-zero image in $\mathrm{Ext}_{\varphi,\Gamma}^1(\cR(\theta^{n-1}),\cR(\theta^{k}))$.
Hence, the set $\{x_k\}_{0\leq k\leq n-1}$ is linearly independent.
Recall from our definition of $x_{n-1}$ (and Lemma~\ref{lem: tri det deform image}) that $x_{n-1}$ belongs to
\begin{equation}\label{equ: tri det deform intersection}
\mathrm{Ext}_{\cD^{+}}^1(\cR(\theta^{n-1}),\cD)\cap \mathrm{im}(q_{\natural})
\end{equation}
It is clear from (\ref{equ: tri det generator}) that $x_k\in \mathrm{im}(q_{\natural})$
Since $y_k$ belongs to (\ref{equ: sym det Levi}), a push out along $\mathrm{Fil}_{W}^{k}(\cD)\hookrightarrow \cD$ (see also (\ref{equ: tri det generator})) implies that $x_k$ belongs to (\ref{equ: tri det deform intersection}) for each $0\leq k\leq n-2$.
Now that both $\mathrm{Ext}_{\cD^{+}}^1(\cR(\theta^{n-1}),\cD)$ and $\mathrm{im}(q_{\natural})$ are $n$-dimensional $E$-subspaces of $\mathrm{Ext}_{\varphi,\Gamma}^1(\cR(\theta^{n-1}),\cD)$ that contains the linearly independent set $\{x_k\}_{0\leq k\leq n-1}$, they must be equal.
This together with the definition of $\mathrm{Ext}_{\cD^{+}}^1(\cD,\cD)$ finishes the proof of the induction step.
\end{proof}

\begin{prop}\label{prop: det deform criterion}
Let $n\geq 2$ and $\cD$ be a de Rham $(\varphi,\Gamma)$-module which is special with level $n-1$ (see \ref{it: special mod 2} of Definition~\ref{def: special mod}).
Then the composition of (\ref{equ: det deform image}) is not surjective if and only if there exists a de Rham $(\varphi,\Gamma)$-module $\cB$ which is special with level $1$ such that $\cD\cong \mathrm{Sym}^{n-1}(\cB)$, in which case its image equals $\mathrm{Ext}_{\mathrm{Sym}^{n}(\cB)}^1(\cR(\theta^{n-1}),\cD)$.
\end{prop}
\begin{proof}
We prove by an increasing induction on $n\geq 2$.
When $n=2$, we have $\cD=\mathrm{Sym}^1(\cD)$ and the image of (\ref{equ: det deform image}) equals $\mathrm{Ext}_{\mathrm{Sym}^{2}(\cB)}^1(\cR(\theta),\cD)$ by Lemma~\ref{lem: sym det deform}.
Assume from now that $n\geq 3$. If $\cD\cong \mathrm{Sym}^{n-1}(\cB)$ for some $\cB$ which is special with level $1$, then the image of (\ref{equ: det deform image}) equals $\mathrm{Ext}_{\mathrm{Sym}^{n}(\cB)}^1(\cR(\theta^{n-1}),\cD)$ by Lemma~\ref{lem: sym det deform}. Assume in the rest of the proof that there does not exist a $\cB$ which is special with level $1$ such that $\cD\cong \mathrm{Sym}^{n-1}(\cB)$. We prove that (\ref{equ: det deform image}) is surjective in this case. 
We write $q_{\dagger}\defeq q_{\mathrm{det}}$ and $q_{\ddagger}\defeq q_{\mathrm{det}}'$ for short (see (\ref{equ: ast deform induction seq}) and (\ref{equ: ast deform induction devissage})).
Note that the map $q_{\ddagger}$ factors through
\begin{multline*}
\mathrm{Ext}_{\mathrm{det}}^1(\cD/\cR(\theta^{n-1}),\cD/\cR(\theta^{n-1}))\buildrel q_{\ddagger}' \over\longrightarrow \mathrm{Ext}_{\varphi,\Gamma}^1(\cR(\theta^{n-2}),\cD/\cR(\theta^{n-1})) \\
\buildrel q_{\ddagger}'' \over\longrightarrow
\mathrm{Ext}_{\varphi,\Gamma}^2(\cR(\theta^{n-2}),\cR(\theta^{n-1}))\buildrel\sim\over\longleftarrow \mathrm{Ext}_{\varphi,\Gamma}^2(\cD/\cR(\theta^{n-1}),\cR(\theta^{n-1})),
\end{multline*}
and that $q_{\ddagger}''$ is a surjection with kernel $\mathrm{Ext}_{\cD}^1(\cR(\theta^{n-2}),\cD/\cR(\theta^{n-1}))$.
If $\cD/\cR(\theta^{n-1})\cong \mathrm{Sym}^{n-2}(\cB)$ for some $\cB$ which is special with level $1$, then we have
\[\mathrm{im}(q_{\ddagger}')=\mathrm{Ext}_{\mathrm{Sym}^{n-1}(\cB)}^1(\cR(\theta^{n-2}),\cD/\cR(\theta^{n-1}))\]
by Lemma~\ref{lem: sym det deform}, which together with $\cD\not\cong \mathrm{Sym}^{n-1}(\cB)$ and the description of $q_{\ddagger}''$ above forces $q_{\ddagger}$ to be a surjection.
If $\cD/\cR(\theta^{n-1})\not\cong \mathrm{Sym}^{n-2}(\cB)$ for any $\cB$ which is special with level $1$, then $q_{\ddagger}'$ is a surjection by our induction hypothesis, which together with the description of $q_{\ddagger}''$ above also forces $q_{\ddagger}$ to be a surjection.
Since $q_{\ddagger}$ is surjective, we deduce from (\ref{equ: ast deform induction devissage}) (with $\ast=\mathrm{det}$ in \emph{loc.cit.}) as well as (\ref{equ: fix det deform dim}) and \ref{it: full deform dim 1} of Lemma~\ref{lem: full deform dim} that
\begin{multline*}
\Dim_E \mathrm{Ext}_{\mathrm{det}}^1(\cD/\cR(\theta^{n-1}),\cD)=\Dim_E \mathrm{Ext}_{\mathrm{det}}^1(\cD/\cR(\theta^{n-1}),\cD/\cR(\theta^{n-1}))+n-2\\
=\Dim_E \mathrm{Ext}_{\mathrm{det}}^1(\cD,\cD)-n-1,
\end{multline*}
which together with (\ref{equ: ast deform induction seq}) (with $\ast=\mathrm{det}$ in \emph{loc.cit.}) and
\[\Dim_E \mathrm{Ext}_{\varphi,\Gamma}^1(\cR(\theta^{n-1}),\cD)=n+1\]
forces $q_{\dagger}$ (namely the composition of (\ref{equ: det deform image})) to be surjective.
This finishes the proof of the induction step.
\end{proof}

We can slightly reformulate Proposition~\ref{prop: det deform criterion} as below. Recall the $E$-subspace $\cE_{n}^{\sharp}\subseteq \cE_{n}$ from the paragraph before Proposition~\ref{prop: sharp flat decomposition}.
\begin{cor}\label{cor: det deform sym}
Let $n\geq 2$ and $\cD^{+}$ be a de Rham $(\varphi,\Gamma)$-module which is special with level $n$ (see \ref{it: special mod 2} of Definition~\ref{def: special mod}) with $\cD\defeq \cD^{+}/\cR(\theta^{n})$ being special with level $n-1$. Then the following statements are equivalent.
\begin{enumerate}[label=(\roman*)]
\item \label{it: det deform sym 1} We have the following inclusion
\[\mathrm{Ext}^{1}_{\mathrm{det}}(\cR(\theta^{n-1}),\cD)\subseteq\mathrm{Ext}_{\cD^{+}}^{1}(\cR(\theta^{n-1}),\cD)\]
between subspaces of $\mathrm{Ext}_{\varphi,\Gamma}^{1}(\cR(\theta^{n-1}),\cD)$.
\item \label{it: det deform sym 2} There exists a de Rham $(\varphi,\Gamma)$-module $\cB$ which is special with level $1$ such that $\cD^{+}\cong\mathrm{Sym}^{n}(\cB)$ and $\cD\cong\mathrm{Sym}^{n-1}(\cB)$.
\item \label{it: det deform sym 3} The composition of
\[\cE_{n}^{\sharp}\otimes_E\cR(\theta^{n})\hookrightarrow \cD_{n}\twoheadrightarrow \cD^{+}\]
is zero.
\end{enumerate}
\end{cor}

Corollary~\ref{cor: det deform sym} clarifies the relation between the $E$-subspace $\cE_{n}^{\sharp}\subseteq \cE_{n}$ and Galois deformations under local Tate duality.
In the rest of this section, we discuss a similar relation between the $E$-subspace $\cE_{n}^{\flat}\subseteq \cE_{n}$ and Galois deformations under local Tate duality, following the ideas of \cite[\S~5.5]{BD19}.

\begin{lem}\label{lem: Tate dual property}
Let $n\geq 1$ $\cD$ be a $(\varphi,\Gamma)$-module equipped with a filtration $\mathrm{Fil}_{W}^{\bullet}(\cD)$ such that $\mathrm{Fil}_{W}^{0}(\cD)=\cD$, $\mathrm{Fil}_{W}^{n}(\cD)=0$ and that $\mathrm{gr}_{W}^{k}(\cD)$ is $\cR(|\cdot|^kz^{\ell_k})$-isotpyic for some $\ell_k\in\Z$ for each $0\leq k\leq n-1$. Then for each $1\leq k\leq n-1$, the short exact sequence
\begin{equation}\label{equ: Galois devissage seq}
0\rightarrow \mathrm{gr}_{W}^{k}(\cD)\rightarrow \cD/\mathrm{Fil}_{W}^{k+1}(\cD)\rightarrow \cD/\mathrm{Fil}_{W}^{k}(\cD)\rightarrow 0
\end{equation}
induces the following short exact sequences
\begin{equation}\label{equ: Tate dual seq 1}
0\rightarrow \mathrm{Ext}_{\varphi,\Gamma}^1(\cR(\theta^{n-1}),\mathrm{gr}_{W}^{k}(\cD))\rightarrow \mathrm{Ext}_{\varphi,\Gamma}^1(\cR(\theta^{n-1}),\cD/\mathrm{Fil}_{W}^{k+1}(\cD))\rightarrow \mathrm{Ext}_{\varphi,\Gamma}^1(\cR(\theta^{n-1}),\cD/\mathrm{Fil}_{W}^{k}(\cD)) \rightarrow 0
\end{equation}
and
\begin{equation}\label{equ: Tate dual seq 2}
0\rightarrow \mathrm{Ext}_{\varphi,\Gamma}^1(\mathrm{gr}_{W}^{k}(\cD),\cR(\theta^{n}))^{\vee}\rightarrow \mathrm{Ext}_{\varphi,\Gamma}^1(\cD/\mathrm{Fil}_{W}^{k+1}(\cD),\cR(\theta^{n}))^{\vee}\rightarrow \mathrm{Ext}_{\varphi,\Gamma}^1(\cD/\mathrm{Fil}_{W}^{k}(\cD),\cR(\theta^{n}))^{\vee} \rightarrow 0.
\end{equation}
\end{lem}
\begin{proof}
Recall from \ref{it: char dR Ext 1} of Lemma~\ref{lem: char dR Ext} that $\Hom_{\varphi,\Gamma}(\cR(\theta^{n-1}),\cR(|\cdot|^kz^{\ell_k}))=0$ for $0\leq k\leq n-2$, which together with a d\'evissage gives
\begin{equation}\label{equ: universal Hom vanishing}
\Hom_{\varphi,\Gamma}(\cR(\theta^{n-1}),\cD/\mathrm{Fil}_{W}^{k}(\cD))=0
\end{equation}
for $1\leq k\leq n-1$. It follows from \ref{it: Galois Ext collection 1} of Lemma~\ref{lem: Galois Ext collection} that
\begin{equation}\label{equ: universal Ext2 vanishing}
\mathrm{Ext}_{\varphi,\Gamma}^2(\cR(\theta^{n-1}),\cR(|\cdot|^kz^{\ell_k}))=0
\end{equation}
for $0\leq k\leq n-1$. Combining (\ref{equ: universal Hom vanishing}) and (\ref{equ: universal Ext2 vanishing}) with a d\'evissage with respect to (\ref{equ: Galois devissage seq}) for each $1\leq k\leq n-1$ gives (\ref{equ: Tate dual seq 1}) for each $1\leq k\leq n-1$.
Similarly, using \ref{it: char dR Ext 1} of Lemma~\ref{lem: char dR Ext} and \ref{it: Galois Ext collection 1} of Lemma~\ref{lem: Galois Ext collection} together with some d\'evissage, we can show that
\[\Hom_{\varphi,\Gamma}(\cR(|\cdot|^kz^{\ell_k}),\cR(\theta^{n}))=0\]
for each $0\leq k\leq n-1$ and
\[\mathrm{Ext}_{\varphi,\Gamma}^2(\cD/\mathrm{Fil}_{W}^{k}(\cD),\cR(\theta^{n}))=0\]
for each $0\leq k\leq n-2$, which together with a d\'evissage with respect to (\ref{equ: Galois devissage seq}) for each $1\leq k\leq n-1$ gives (\ref{equ: Tate dual seq 2}) for each $1\leq k\leq n-1$.
\end{proof}

\begin{lem}\label{lem: Galois Ext universal}
Let $n\geq 1$. Then we have
\begin{equation}\label{equ: full dR Ext}
\mathrm{Ext}_{g}^1(\cD,\cR(\theta^{n}))=\mathrm{Ext}_{\varphi,\Gamma}^1(\cD,\cR(\theta^{n}))
\end{equation}
if $\cD$ equals either $\cD_{n-1}$ or $\cD_{n-1}^{\flat}$.
\end{lem}
\begin{proof}
Let $\cD$ be one of $\cD_{n-1}$, $\cD_{n-1}^{\flat}$ and $\cD_{n,n-1}^{\flat}$ which is equipped with the usual decreasing filtration $\mathrm{Fil}_{W}^{\bullet}(\cD)$.
Note that such $\cD$ clearly satisfies the assumption of Lemma~\ref{lem: Tate dual property} so that (\ref{equ: Tate dual seq 2}) holds for each $1\leq k\leq n-1$.
Recall from the proof of Proposition~\ref{prop: Galois dR Ext exact} that we have a short exact sequence
\[0\rightarrow \mathrm{Ext}_{g}^1(\mathrm{gr}_{W}^{k}(\cD),\cR(\theta^{n}))^{\vee}\rightarrow \mathrm{Ext}_{g}^1(\cD/\mathrm{Fil}_{W}^{k+1}(\cD),\cR(\theta^{n}))^{\vee}\rightarrow \mathrm{Ext}_{g}^1(\cD/\mathrm{Fil}_{W}^{k}(\cD),\cR(\theta^{n}))^{\vee} \rightarrow 0\]
for each $1\leq k\leq n-1$, which together with (\ref{equ: Tate dual seq 2}) and a snake's Lemma argument gives
\[\mathrm{Ext}_{g}^1(\cD/\mathrm{Fil}_{W}^{k+1}(\cD),\cR(\theta^{n}))=\mathrm{Ext}_{\varphi,\Gamma}^1(\cD/\mathrm{Fil}_{W}^{k+1}(\cD),\cR(\theta^{n}))\]
for $0\leq k\leq n-1$ by an increasing induction on $k$. Taking $k=n-1$ gives (\ref{equ: full dR Ext}).
\end{proof}

Using Lemma~\ref{lem: Galois Ext universal}, we see that the embeddings $\cD_{n,n-1}^{\flat}\hookrightarrow \cD_{n-1}^{\flat}\hookrightarrow \cD_{n-1}$ induce a commutative diagram of the form
\begin{equation}\label{equ: Galois sub flat diagram}
\xymatrix{
\mathrm{Ext}_{\varphi,\Gamma}^1(\cR(\theta^{n-1}),\cD_{n,n-1}^{\flat}) \ar^{p_1}[r] \ar_{q_1}^{\wr}[d] & \mathrm{Ext}_{\varphi,\Gamma}^1(\cR(\theta^{n-1}),\cD_{n-1}^{\flat}) \ar@{^{(}->}[r]^{p_2} \ar_{q_2}^{\wr}[d] & \mathrm{Ext}_{\varphi,\Gamma}^1(\cR(\theta^{n-1}),\cD_{n-1}) \ar_{q_3}^{\wr}[d]\\
\mathrm{Ext}_{\varphi,\Gamma}^1(\cD_{n,n-1}^{\flat},\cR(\theta^{n}))^{\vee} \ar^{p_3}[r] \ar@{->>}[d]_{q_4} & \mathrm{Ext}_{\varphi,\Gamma}^1(\cD_{n-1}^{\flat},\cR(\theta^{n}))^{\vee} \ar@{^{(}->}[r]^{p_4} \ar@{=}[d] & \mathrm{Ext}_{\varphi,\Gamma}^1(\cD_{n-1},\cR(\theta^{n}))^{\vee} \ar@{=}[d]\\
\mathrm{Ext}_{g}^1(\cD_{n,n-1}^{\flat},\cR(\theta^{n}))^{\vee} \ar@{^{(}->}[r]^{p_5} & \mathrm{Ext}_{g}^1(\cD_{n-1}^{\flat},\cR(\theta^{n}))^{\vee} \ar@{^{(}->}[r]^{p_6} & \mathrm{Ext}_{g}^1(\cD_{n-1},\cR(\theta^{n}))^{\vee}
}
\end{equation}
with $q_1$, $q_2$ and $q_3$ being isomorphisms by local Tate duality, $q_4$ being a surjection by the definition of $\mathrm{Ext}_{g}^1$, $p_5$ and $p_6$ being embeddings by Lemma~\ref{lem: Galois sub flat}, and finally $p_2$ as well as $p_4$ being embeddings by the commutativity of (\ref{equ: Galois sub flat diagram}).
Now that $q_1$ is an isomorphism and $q_4$ is a surjection, the commutativity of (\ref{equ: Galois sub flat diagram}) gives
\begin{equation}\label{equ: flat image composition}
\cE_{n}^{\flat}=\mathrm{im}(p_6\circ p_5)=\mathrm{im}(p_4\circ p_3\circ q_1)=\mathrm{im}(q_3\circ p_2\circ p_1)
\end{equation}

The following discussion is directly motivated by the discussions on $\mathrm{Ext}_{\mathrm{HT}}^1$ in \cite[\S~5.5]{BD19}.
We consider a short exact sequence
\begin{equation}\label{equ: phi Gamma Ext}
0\rightarrow \cD'\rightarrow \cD\rightarrow \cD''\rightarrow 0
\end{equation}
of $(\varphi,\Gamma)$-modules over $\cR$ with both $\cD'$ and $\cD''$ being de Rham (and $\cD$ being almost de Rham).
By applying $D_{\rm{pdR}}(-)$ to (\ref{equ: phi Gamma Ext}) we obtain a (strict) short exact sequence of filtered $E$-vector space
\begin{equation}\label{equ: pdR Ext}
0\rightarrow D_{\rm{dR}}(\cD') \rightarrow D_{\rm{pdR}}(\cD) \rightarrow D_{\rm{dR}}(\cD'') \rightarrow 0.
\end{equation}
Moreover, $D_{\rm{pdR}}(\cD)$ is equipped with the so-called \emph{Fontaine's $N$-operator} which is a nilpotent $E$-linear endomorphism $\cN$ of $D_{\rm{pdR}}(\cD)$ compatible with the Hodge filtration $\mathrm{Fil}_{H}^{\bullet}(D_{\rm{pdR}}(\cD))$ on $D_{\rm{pdR}}(\cD)$. Now that both $\cD'$ and $\cD''$ are de Rham, $\cN$ acts by zero on both $D_{\rm{dR}}(\cD')$ and $D_{\rm{pdR}}(\cD)/D_{\rm{dR}}(\cD')\cong D_{\rm{dR}}(\cD'')$ and thus induces a $E$-linear map
\[\cN: D_{\rm{dR}}(\cD'')\rightarrow D_{\rm{dR}}(\cD').\]
We write
\[h_0\defeq \max\{h\in\Z \mid \mathrm{Fil}_{H}^{h}(D_{\rm{dR}}(\cD'))\neq 0\}\]
and $\mathrm{Fil}_{H}^{\max}(D_{\rm{dR}}(\cD'))\defeq \mathrm{Fil}_{H}^{h_0}(D_{\rm{dR}}(\cD'))$.
We define $\mathrm{Ext}_{\max}^1(\cD'',\cD')$ as the $E$-subspace of $\mathrm{Ext}_{\varphi,\Gamma}^1(\cD'',\cD')$ consisting of those $\cD$ that satisfy
\begin{equation}\label{equ: Fontaine N image}
\mathrm{im}(\cN)\subseteq \mathrm{Fil}_{H}^{\max}(D_{\rm{dR}}(\cD')).
\end{equation}
By definition we have the following exact sequence
\begin{equation}\label{equ: Fontaine N seq}
0\rightarrow \mathrm{Ext}_{g}^1(\cD'',\cD')\rightarrow \mathrm{Ext}_{\max}^1(\cD'',\cD')\rightarrow \Hom_{E}(D_{\rm{dR}}(\cD''), \mathrm{Fil}_{H}^{\max}(D_{\rm{dR}}(\cD'))).
\end{equation}
\begin{prop}\label{prop: Ext max}
Let $n\geq 1$. Recall the map $p_1$ from (\ref{equ: Galois sub flat diagram}). We have the following results.
\begin{enumerate}[label=(\roman*)]
\item \label{it: Ext max 1} The subspace
\[\mathrm{Ext}_{\max}^1(\cR(\theta^{n-1}),\cD_{n-1}^{\flat})\subseteq \mathrm{Ext}_{\varphi,\Gamma}^1(\cR(\theta^{n-1}),\cD_{n-1}^{\flat})\]
equals $\mathrm{im}(p_1)$ and fits into the following short exact sequence
\begin{equation}\label{equ: Ext max seq}
0\rightarrow \mathrm{Ext}_{g}^1(\cR(\theta^{n-1}),\cD_{n-1}^{\flat})\rightarrow \mathrm{Ext}_{\max}^1(\cR(\theta^{n-1}),\cD_{n-1}^{\flat})\rightarrow \Hom_{E}(D_{\rm{dR}}(\cR(\theta^{n-1})), \mathrm{Fil}_{H}^{0}(D_{\rm{dR}}(\cD_{n-1}^{\flat}))) \rightarrow 0.
\end{equation}
\item \label{it: Ext max 2} The commutative diagram (\ref{equ: Galois sub flat diagram}) induces a canonical isomorphism
\begin{equation}\label{equ: Ext max isom}
\mathrm{Ext}_{\max}^1(\cR(\theta^{n-1}),\cD_{n-1}^{\flat})\cong \cE_{n}^{\flat}.
\end{equation}
\end{enumerate}
\end{prop}
\begin{proof}
Recall that $\cD_{n-1}^{\flat}$ is equipped with a decreasing filtration $\mathrm{Fil}_{W}^{\bullet}(\cD_{n-1}^{\flat})$ which restricts to a decreasing filtration $\mathrm{Fil}_{W}^{\bullet}(\cD_{n,n-1}^{\flat})$ on $\cD_{n,n-1}^{\flat}$.
We write $\cD'\defeq \mathrm{Fil}_{W}^{1}(\cD_{n,n-1}^{\flat})$ and $\cD''\defeq \mathrm{Fil}_{W}^{1}(\cD_{n-1}^{\flat})$ for short.
We write $\mathrm{Fil}_{H}^{\bullet}(D_{\rm{dR}}(\cD'))$ for the Hodge filtration on $D_{\rm{dR}}(\cD')$.
For each $1\leq k\leq \ell\leq n-1$, we observe that
\[\mathrm{Fil}_{H}^{h}(D_{\rm{dR}}(\mathrm{Fil}_{W}^{k}(\cD_{n,n-1}^{\flat})/\mathrm{Fil}_{W}^{\ell+1}(\cD_{n,n-1}^{\flat})))=0\]
for each $h>-n$. Consequently, for each $(\varphi,\Gamma)$-module corresponding to an element of \[\mathrm{Ext}_{\varphi,\Gamma}^1(\cR(\theta^{n-1}),\mathrm{Fil}_{W}^{k}(\cD_{n,n-1}^{\flat})/\mathrm{Fil}_{W}^{\ell+1}(\cD_{n,n-1}^{\flat})),\] the Fontaine's $N$-operator $\cN$ must act trivially on its associated $D_{\rm{pdR}}(-)$.
In other words, we have
\begin{equation}\label{equ: Ext max dR sub 1}
\mathrm{Ext}_{g}^1(\cR(\theta^{n-1}),\mathrm{Fil}_{W}^{k}(\cD_{n,n-1}^{\flat})/\mathrm{Fil}_{W}^{\ell+1}(\cD_{n,n-1}^{\flat}))=\mathrm{Ext}_{\varphi,\Gamma}^1(\cR(\theta^{n-1}),\mathrm{Fil}_{W}^{k}(\cD_{n,n-1}^{\flat})/\mathrm{Fil}_{W}^{\ell+1}(\cD_{n,n-1}^{\flat}))
\end{equation}
for each $1\leq k\leq \ell\leq n-1$, which together with (\ref{equ: Tate dual seq 1}) (with $\cD$ in \emph{loc.cit.} taken to be $\cD'$ here) gives
\begin{multline}\label{equ: Ext max dR seq}
0\rightarrow \mathrm{Ext}_{g}^1(\cR(\theta^{n-1}),\mathrm{gr}_{W}^{k}(\cD_{n,n-1}^{\flat}))\rightarrow \mathrm{Ext}_{g}^1(\cR(\theta^{n-1}),\cD'/\mathrm{Fil}_{W}^{k+1}(\cD_{n,n-1}^{\flat}))\\
\rightarrow \mathrm{Ext}_{g}^1(\cR(\theta^{n-1}),\cD'/\mathrm{Fil}_{W}^{k}(\cD_{n,n-1}^{\flat})) \rightarrow 0
\end{multline}
for each $1\leq k\leq n-1$.
Now that $\mathrm{Ext}_{g}^1(\cR(\theta^{n-1}),\cR(1_{K^{\times}}))=0$ from \ref{it: char dR Ext 2} of Lemma~\ref{lem: char dR Ext}, we know that $\cD'\hookrightarrow \cD_{n,n-1}^{\flat}$ induces an isomorphism
\begin{equation}\label{equ: Ext max dR sub 2}
\mathrm{Ext}_{g}^1(\cR(\theta^{n-1}),\cD')\buildrel\sim\over\longrightarrow \mathrm{Ext}_{g}^1(\cR(\theta^{n-1}),\cD_{n,n-1}^{\flat}).
\end{equation}
It follows from \ref{it: char dR Ext 1} of Lemma~\ref{lem: char dR Ext} that $\Hom_{\varphi,\Gamma}(\cR(\theta^{n-1}),\cR(1_{K^{\times}}))=0$ and from \ref{it: Galois Ext collection 1} of Lemma~\ref{lem: Galois Ext collection} (with a d\'evissage) that $\mathrm{Ext}_{\varphi,\Gamma}^2(\cR(\theta^{n-1}),\cD')=0$.
Hence, the short exact sequence $0\rightarrow \cD'\rightarrow \cD_{n,n-1}^{\flat}\rightarrow \cR(1_{K^{\times}})\rightarrow 0$ induces the following short exact sequence
\begin{equation}\label{equ: Ext max seq 1}
0\rightarrow \mathrm{Ext}_{\varphi,\Gamma}^1(\cR(\theta^{n-1}),\cD')\rightarrow \mathrm{Ext}_{\varphi,\Gamma}^1(\cR(\theta^{n-1}),\cD_{n,n-1}^{\flat})\rightarrow \mathrm{Ext}_{\varphi,\Gamma}^1(\cR(\theta^{n-1}),\cR(1_{K^{\times}})) \rightarrow 0.
\end{equation}
Let $\cD^{\dagger}$ be an arbitrary $(\varphi,\Gamma)$-module represented by an element of $\mathrm{Ext}_{\varphi,\Gamma}^1(\cR(\theta^{n-1}),\cD_{n,n-1}^{\flat})$. 
Note that the filtered $E$-vector space $D_{\rm{pdR}}(\cD^{\dagger})$ admits Fontaine's $N$-operator $\cN$ as an endomorphism and moreover fits into the following strict short exact sequence
\[0\rightarrow D_{\rm{dR}}(\cD_{n,n-1}^{\flat}) \rightarrow D_{\rm{pdR}}(\cD^{\dagger})\rightarrow D_{\rm{dR}}(\cR(\theta^{n-1})) \rightarrow 0.\]
Note that $D_{\rm{dR}}(\cR(\theta^{n-1}))$ is a $1$-dimensional $E$-vector space with $\mathrm{gr}_{H}^{h}(D_{\rm{dR}}(\cR(\theta^{n-1})))\neq 0$ if and only if $h=1-n$, and the Hodge filtration on $D_{\rm{dR}}(\cD_{n,n-1}^{\flat})$ satisfies $\mathrm{gr}_{H}^{h}(D_{\rm{dR}}(\cD_{n,n-1}^{\flat}))\neq 0$ if and only if $h=0,-n$ with $\Dim_E \mathrm{gr}_{H}^{0}(D_{\rm{dR}}(\cD_{n,n-1}^{\flat}))=1$. Now that $\cN$ stablizes the Hodge filtration on $D_{\rm{pdR}}(\cD^{\dagger})$, its associated element of 
\[\Hom_{E}(D_{\rm{dR}}(\cR(\theta^{n-1})),D_{\rm{dR}}(\cD_{n,n-1}^{\flat}))\]
must be contained in
\[\Hom_{E}(D_{\rm{dR}}(\cR(\theta^{n-1})),\mathrm{Fil}_{H}^{0}(D_{\rm{dR}}(\cD_{n,n-1}^{\flat}))).\]
In other words, we have shown that
\[\mathrm{Ext}_{\max}^1(\cR(\theta^{n-1}),\cD_{n,n-1}^{\flat})=\mathrm{Ext}_{\varphi,\Gamma}^1(\cR(\theta^{n-1}),\cD_{n,n-1}^{\flat}).\]
Combining this with (\ref{equ: Ext max seq 1}), (\ref{equ: Fontaine N seq}) (by taking $\cD''=\cR(\theta^{n-1})$ and $\cD'=\cD_{n,n-1}^{\flat}$ in \emph{loc.cit.}), (\ref{equ: Ext max dR sub 1}) and (\ref{equ: Ext max dR sub 2}), we arrive at the following commutative diagram
\begin{equation}\label{equ: Ext max diagram 1}
\xymatrix{
\mathrm{Ext}_{\varphi,\Gamma}^1(\cR(\theta^{n-1}),\cD') \ar@{^{(}->}[r] \ar^{\wr}[d] & \mathrm{Ext}_{\varphi,\Gamma}^1(\cR(\theta^{n-1}),\cD_{n,n-1}^{\flat}) \ar@{->>}[r] \ar@{=}[d] & \mathrm{Ext}_{\varphi,\Gamma}^1(\cR(\theta^{n-1}),\cR(1_{K^{\times}}))  \ar^{\wr}[d]\\
\mathrm{Ext}_{g}^1(\cR(\theta^{n-1}),\cD_{n,n-1}^{\flat}) \ar@{^{(}->}[r] & \mathrm{Ext}_{\max}^1(\cR(\theta^{n-1}),\cD_{n,n-1}^{\flat}) \ar@{->>}[r] & \Hom_{E}(D_{\rm{dR}}(\cR(\theta^{n-1})), \mathrm{Fil}_{H}^{0}(D_{\rm{dR}}(\cD_{n,n-1}^{\flat})))
}
\end{equation}
with the RHS vertical isomorphism sending a $(\varphi,\Gamma)$-module corresponding to an element of the $E$-vector space $\mathrm{Ext}_{\varphi,\Gamma}^1(\cR(\theta^{n-1}),\cR(1_{K^{\times}}))$ to the Fontaine's $N$-operator $\cN$ on its associated $D_{\rm{pdR}}(-)$.
We consider the following commutative diagram induced from the embedding $\cD_{n,n-1}^{\flat}\hookrightarrow \cD_{n-1}^{\flat}$
\begin{equation}\label{equ: Ext max diagram 2}
\xymatrix{
\mathrm{Ext}_{g}^1(\cR(\theta^{n-1}),\cD_{n,n-1}^{\flat}) \ar@{^{(}->}[r] \ar[d] & \mathrm{Ext}_{\max}^1(\cR(\theta^{n-1}),\cD_{n,n-1}^{\flat}) \ar@{->>}[r] \ar[d] & \Hom_{E}(D_{\rm{dR}}(\cR(\theta^{n-1})), \mathrm{Fil}_{H}^{0}(D_{\rm{dR}}(\cD_{n,n-1}^{\flat}))) \ar[d]\\
\mathrm{Ext}_{g}^1(\cR(\theta^{n-1}),\cD_{n-1}^{\flat}) \ar@{^{(}->}[r] & \mathrm{Ext}_{\max}^1(\cR(\theta^{n-1}),\cD_{n-1}^{\flat}) \ar[r] & \Hom_{E}(D_{\rm{dR}}(\cR(\theta^{n-1})), \mathrm{Fil}_{H}^{0}(D_{\rm{dR}}(\cD_{n-1}^{\flat})))
}
\end{equation}
with both rows being exact sequences (see (\ref{equ: Fontaine N seq})).
Now that $\cD_{n,n-1}^{\flat}\hookrightarrow \cD_{n-1}^{\flat}$ induces an isomorphism between $1$-dimensional $E$-vector spaces
\[\mathrm{Fil}_{H}^{0}(D_{\rm{dR}}(\cD_{n,n-1}^{\flat}))\buildrel\sim\over\longrightarrow \mathrm{Fil}_{H}^{0}(D_{\rm{dR}}(\cD_{n-1}^{\flat})),\]
we see that the right vertical map of (\ref{equ: Ext max diagram 2}) is an isomorphism, and thus the bottom right horizontal map of (\ref{equ: Ext max diagram 2}) is surjective. The bottom horizontal row of (\ref{equ: Ext max diagram 2}) thus gives the desired (\ref{equ: Ext max seq}).
Now that the embedding $\cR(|\cdot|^{k}z^{n})\hookrightarrow \cR(|\cdot|^{k}z^{n-1})$ induces a surjection
\[\mathrm{Ext}_{g}^1(\cR(\theta^{n-1}),\cR(|\cdot|^{k}z^{n}))\twoheadrightarrow \mathrm{Ext}_{g}^1(\cR(\theta^{n-1}),\cR(|\cdot|^{k}z^{n-1}))\]
for each $1\leq k\leq n-1$, we deduce from (\ref{equ: Ext max dR seq}) and a d\'evissage that the embedding $\mathrm{Fil}_{W}^{k}(\cD_{n,n-1}^{\flat})\hookrightarrow \mathrm{Fil}_{W}^{k}(\cD_{n-1}^{\flat})$ induces a surjection
\[\mathrm{Ext}_{g}^1(\cR(\theta^{n-1}),\mathrm{Fil}_{W}^{k}(\cD_{n,n-1}^{\flat}))\twoheadrightarrow \mathrm{Ext}_{g}^1(\cR(\theta^{n-1}),\mathrm{Fil}_{W}^{k}(\cD_{n-1}^{\flat})),\]
for each $1\leq k\leq n-1$. In particular, we obtain a surjection
\begin{equation}\label{equ: Ext max surjection}
\mathrm{Ext}_{g}^1(\cR(\theta^{n-1}),\cD')\twoheadrightarrow \mathrm{Ext}_{g}^1(\cR(\theta^{n-1}),\cD''),
\end{equation}
which together with (\ref{equ: Ext max dR sub 2}) and its variant with $\cD''$ and $\cD_{n-1}^{\flat}$ replacing $\cD'$ and $\cD_{n,n-1}^{\flat}$ implies that the left vertical map of (\ref{equ: Ext max diagram 2}) is surjective. Consequently, the middle vertical map of (\ref{equ: Ext max diagram 2}) is also surjective, which together with (\ref{equ: Ext max diagram 1}) gives
\[\mathrm{im}(p_1)=\mathrm{Ext}_{\max}^1(\cR(\theta^{n-1}),\cD_{n-1}^{\flat}).\]
This together with (\ref{equ: flat image composition}) gives the isomorphism (\ref{equ: Ext max isom}).
\end{proof}

\begin{prop}\label{prop: special Ext max}
Let $n\geq 1$ and $\cD$ be a de Rham $(\varphi,\Gamma)$-module which is special with level $n-1$ (see \ref{it: special mod 2} of Definition~\ref{def: special mod}). Then we have
\begin{equation}\label{equ: Ext max Hom}
\Dim_E \Hom_{\varphi,\Gamma}(\cD_{n-1}^{\flat},\cD)=1,
\end{equation}
and the unique (up to scalars) non-zero map $\cD_{n-1}^{\flat}\rightarrow \cD$ induces a surjection
\begin{equation}\label{equ: Ext max quotient}
\mathrm{Ext}_{\max}^1(\cR(\theta^{n-1}),\cD_{n-1}^{\flat})\twoheadrightarrow \mathrm{Ext}_{\max}^1(\cR(\theta^{n-1}),\cD)
\end{equation}
with RHS being $2$-dimensional and fitting into the short exact sequence
\begin{equation}\label{equ: Ext max quotient seq}
0\rightarrow \mathrm{Ext}_{g}^1(\cR(\theta^{n-1}),\cD)\rightarrow \mathrm{Ext}_{\max}^1(\cR(\theta^{n-1}),\cD)\rightarrow \Hom_{E}(D_{\rm{dR}}(\cR(\theta^{n-1})), \mathrm{Fil}_{H}^{0}(D_{\rm{dR}}(\cD))) \rightarrow 0.
\end{equation}
\end{prop}
\begin{proof}
Thanks to Lemma~\ref{lem: Galois universal Hom} we have
\[\Dim_E \Hom_{\varphi,\Gamma}(\cD_{n-1},\cD)=1\]
with the unique (up to scalars) non-zero map $\cD_{n-1}\rightarrow \cD$ being a surjection. Now that the maps 
\begin{equation}\label{equ: special Ext max composition}
\cD_{n,n-1}^{\flat}\rightarrow \cD_{n-1}^{\flat}\rightarrow \cD_{n-1}\twoheadrightarrow \cD 
\end{equation}
induce by taking $\mathrm{gr}_{W}^{0}(-)$ isomorphisms between $\cR(1_{K^{\times}})$, we know that the composition of (\ref{equ: special Ext max composition}) is non-zero. This together with Lemma~\ref{lem: Galois unique Hom} gives (\ref{equ: Ext max Hom}) as well as
\[\Dim_E \Hom_{\varphi,\Gamma}(\cD_{n,n-1}^{\flat},\cD)=1.\]
We write $\cD^{\flat}$ (resp.~$\cD^{\flat\flat}$) for the image of the unique (up to scalars) non-zero map $\cD_{n-1}^{\flat}\rightarrow \cD$ (resp.~$\cD_{n,n-1}^{\flat}\rightarrow \cD$).
The maps $\cD^{\flat\flat}\hookrightarrow \cD^{\flat}\hookrightarrow \cD$ restrict to maps $\cR(|\cdot|^{n-1}z^{n})\hookrightarrow \cR(\theta^{n-1})=\cR(\theta^{n-1})$ and induce the following commutative diagram
\begin{equation}\label{equ: Ext max quotient diagram 1}
\xymatrix{
\mathrm{Ext}_{g}^1(\cR(\theta^{n-1}),\cR(|\cdot|^{n-1}z^{n})) \ar@{->>}[r] \ar[d] & \mathrm{Ext}_{g}^1(\cR(\theta^{n-1}),\cR(\theta^{n-1})) \ar@{=}[r] \ar[d] & \mathrm{Ext}_{g}^1(\cR(\theta^{n-1}),\cR(\theta^{n-1})) \ar^{\wr}[d]\\
\mathrm{Ext}_{g}^1(\cR(\theta^{n-1}),\cD^{\flat\flat}) \ar[r] & \mathrm{Ext}_{g}^1(\cR(\theta^{n-1}),\cD^{\flat}) \ar[r] & \mathrm{Ext}_{g}^1(\cR(\theta^{n-1}),\cD)
}
\end{equation}
Here the left top horizontal map of (\ref{equ: Ext max quotient diagram 1}) is surjective with $1$-dimensional target by \ref{it: char dR Ext 4} of Lemma~\ref{lem: char dR Ext}, and the RHS vertical map is an isomorphism by $\mathrm{Ext}_{g}^1(\cR(\theta^{n-1}),\cR(\theta^{k}))=0$ for each $0\leq k\leq n-2$ (see \ref{it: char dR Ext 2} of Lemma~\ref{lem: char dR Ext}) and a d\'evissage.
Now that $\cD^{\flat\flat}\hookrightarrow \cD^{\flat}\hookrightarrow \cD$ induce isomorphisms between $1$-dimensional $E$-vector spaces
\[\mathrm{Fil}_{H}^{0}(D_{\rm{dR}}(\cD^{\flat\flat}))\buildrel\sim\over\longrightarrow \mathrm{Fil}_{H}^{0}(D_{\rm{dR}}(\cD^{\flat})) \buildrel\sim\over\longrightarrow \mathrm{Fil}_{H}^{0}(D_{\rm{dR}}(\cD)),\]
parallel argument as below (\ref{equ: Ext max diagram 2}) gives the following commutative diagram
\begin{equation}\label{equ: Ext max quotient diagram 2}
\xymatrix{
\mathrm{Ext}_{g}^1(\cR(\theta^{n-1}),\cD^{\flat\flat}) \ar@{^{(}->}[r] \ar[d] & \mathrm{Ext}_{\max}^1(\cR(\theta^{n-1}),\cD^{\flat\flat}) \ar@{->>}[r] \ar[d] & \Hom_{E}(D_{\rm{dR}}(\cR(\theta^{n-1})), \mathrm{Fil}_{H}^{0}(D_{\rm{dR}}(\cD^{\flat\flat}))) \ar^{\wr}[d]\\
\mathrm{Ext}_{g}^1(\cR(\theta^{n-1}),\cD^{\flat}) \ar@{^{(}->}[r] \ar@{->>}[d] & \mathrm{Ext}_{\max}^1(\cR(\theta^{n-1}),\cD^{\flat}) \ar@{->>}[r] \ar@{->>}[d] & \Hom_{E}(D_{\rm{dR}}(\cR(\theta^{n-1})), \mathrm{Fil}_{H}^{0}(D_{\rm{dR}}(\cD^{\flat}))) \ar^{\wr}[d] \\
\mathrm{Ext}_{g}^1(\cR(\theta^{n-1}),\cD) \ar@{^{(}->}[r] & \mathrm{Ext}_{\max}^1(\cR(\theta^{n-1}),\cD) \ar@{->>}[r] & \Hom_{E}(D_{\rm{dR}}(\cR(\theta^{n-1})), \mathrm{Fil}_{H}^{0}(D_{\rm{dR}}(\cD)))
}.
\end{equation}
Here we recall from (\ref{equ: Ext max quotient diagram 1}) that the left bottom vertical map of (\ref{equ: Ext max quotient diagram 2}) is surjective, which together with all rows of \emph{loc.cit.} being short exact sequences and all vertical maps in the right column of \emph{loc.cit.}  being isomorphisms forces the middle bottom vertical map of \emph{loc.cit.} to be surjective as well.
Similarly, we also know that the composition of the two vertical maps in the left (resp.~middle) column of \emph{loc.cit.} is surjective.
In particular, we obtain the short exact sequence (\ref{equ: Ext max quotient seq}) as the bottom row of (\ref{equ: Ext max quotient diagram 2}), with
\begin{multline*}
\Dim_E \mathrm{Ext}_{\max}^1(\cR(\theta^{n-1}),\cD)\\
=\Dim_E \mathrm{Ext}_{g}^1(\cR(\theta^{n-1}),\cD)+\Dim_E \Hom_{E}(D_{\rm{dR}}(\cR(\theta^{n-1})), \mathrm{Fil}_{H}^{0}(D_{\rm{dR}}(\cD)))=1+1=2
\end{multline*}
Now that we have $\mathrm{Ext}_{\varphi,\Gamma}^2(\cR(\theta^{n-1}),\cR(|\cdot|^kz^n))=0$ for each $1\leq k\leq n-1$ from \ref{it: Galois Ext collection 1} of Lemma~\ref{lem: Galois Ext collection}, by d\'evissage the surjection $\cD_{n,n-1}^{\flat}\twoheadrightarrow \cD^{\flat\flat}$ induces a surjection
\[\mathrm{Ext}_{\varphi,\Gamma}^1(\cR(\theta^{n-1}),\mathrm{Fil}_{W}^{1}(\cD_{n,n-1}^{\flat}))\twoheadrightarrow \mathrm{Ext}_{\varphi,\Gamma}^1(\cR(\theta^{n-1}),\mathrm{Fil}_{W}^{1}(\cD^{\flat\flat})).\]
This together with $\mathrm{gr}_{W}^{0}(\cD_{n,n-1}^{\flat})\buildrel\sim\over\longrightarrow \mathrm{gr}_{W}^{0}(\cD^{\flat\flat})\cong\cR(1_{K^{\times}})$ and $\mathrm{Ext}_{g}^1(\cR(\theta^{n-1}),\cR(1_{K^{\times}}))=0$ (see \ref{it: char dR Ext 2} of Lemma~\ref{lem: char dR Ext}) gives a surjection
\[\mathrm{Ext}_{g}^1(\cR(\theta^{n-1}),\cD_{n,n-1}^{\flat})\twoheadrightarrow \mathrm{Ext}_{g}^1(\cR(\theta^{n-1}),\cD^{\flat\flat}),\]
and thus a surjection
\begin{equation}\label{equ: Ext max universal to flatflat}
\mathrm{Ext}_{\max}^1(\cR(\theta^{n-1}),\cD_{n,n-1}^{\flat})\twoheadrightarrow \mathrm{Ext}_{\max}^1(\cR(\theta^{n-1}),\cD^{\flat\flat})
\end{equation}
using the top row of both (\ref{equ: Ext max diagram 2}) and (\ref{equ: Ext max quotient diagram 2}).
Taking the composition of the two vertical maps in the middle column of (\ref{equ: Ext max quotient diagram 2}) with (\ref{equ: Ext max universal to flatflat}), we obtain a surjection
\[\mathrm{Ext}_{\max}^1(\cR(\theta^{n-1}),\cD_{n,n-1}^{\flat})\twoheadrightarrow \mathrm{Ext}_{\max}^1(\cR(\theta^{n-1}),\cD)\]
which necessarily factors through a surjection (\ref{equ: Ext max quotient}) (as the non-zero map $\cD_{n,n-1}^{\flat}\rightarrow \cD$ factors through $\cD_{n-1}^{\flat}$).
\end{proof}

\subsection{Tensor product on the Galois side}\label{subsec: Galois tensor product}
We define and study a certain tensor product map (see (\ref{equ: Galois tensor sub Ext embedding}) and (\ref{equ: Galois multi cup map}) below) on the Galois side that will be useful in \S \ref{subsec: map candidate}.

Let $m,n\geq 0$. When $n=0$, we have a canonical embedding $\cR(\theta^{m})\otimes_E\cE_{m}\hookrightarrow \cD_{m}$ from the construction of $\cD_m$. Suppose that we are given an embedding
\begin{equation}\label{equ: Galois tensor sub -}
\cD_{n-1}\otimes_{\cR} (\cR(\theta^{m})\otimes_E \cE_{m}) \hookrightarrow \cD_{m+n-1}
\end{equation}
for some $n\geq 1$. It follows from (\ref{equ: universal Galois Ext devissage}) (by taking $k=m+n$ and $\cD=\cD_{n-1}\otimes_{\cR}(\cR(\theta^{m})\otimes_E \cE_{m})$ in \emph{loc.cit.}) that the embedding (\ref{equ: Galois tensor sub -}) induces an embedding
\begin{multline}\label{equ: Galois tensor sub Ext embedding}
\cE_{n}\otimes_E \cE_{m}=\mathrm{Ext}_g^1(\cD_{n-1},\cR(\theta^{n}))^{\vee}\otimes_E \cE_{m}
=\mathrm{Ext}_g^1(\cD_{n-1}\otimes_{\cR} \cR(\theta^{m}),\cR(\theta^{m+n}))^{\vee}\otimes_E \cE_{m}\\
=\mathrm{Ext}_g^1(\cD_{n-1}\otimes_{\cR}(\cR(\theta^{m})\otimes_E \cE_{m}),\cR(\theta^{m+n}))^{\vee}\hookrightarrow \mathrm{Ext}_g^1(\cD_{m+n-1},\cR(\theta^{m+n}))^{\vee}=\cE_{m+n}.
\end{multline}
By consider the universal extension associated with each term of (\ref{equ: Galois tensor sub Ext embedding}), we obtain an embedding
\begin{equation}\label{equ: Galois tensor sub}
\cD_{n}\otimes_{\cR} (\cR(\theta^{m})\otimes_E \cE_{m})\hookrightarrow \cD_{m+n}.
\end{equation}
Hence, by an increasing induction on $n\geq 0$, we obtain a natural choice of (\ref{equ: Galois tensor sub Ext embedding}) and (\ref{equ: Galois tensor sub}) for each $m,n\geq 0$.
More generally, for each (ordered) tuple of integers $m_1,\dots,m_t\geq 0$, we can similarly define a map
\begin{equation}\label{equ: Galois multi cup map}
\cE_{m_1}\otimes_E\cdots\otimes_E\cE_{m_t}\rightarrow \cE_{\sum_{\ell=1}^t m_{\ell}}
\end{equation}
by an increasing induction on $t\geq 2$.
In particular, we have a natural map
\begin{equation}\label{equ: Galois sym embedding}
\mathrm{Sym}^{m}(\cE_{1})\hookrightarrow \cE_{m}
\end{equation}
for each $m\geq 1$.

Recall from \cite[Lem~1.19]{Ding17}
that $\cE_{1}\cong \Hom(K^\times,E)$ and we define $\cE_{1}^{\infty}\defeq E\val\subseteq \Hom(K^\times,E)\subseteq \cE_{1}$ and $\cE_{1}^{\plog}\defeq E\plog\subseteq \Hom(K^\times,E)\subseteq \cE_{1}$. For general $n\geq 2$, we define $\cE_{n}^{\plog}$ as the image of the composition of the following maps
\begin{equation}\label{equ: Galois log cup}
\cE_{1}^{\plog}\otimes_E\cdots\otimes_E\cE_{1}^{\plog}\hookrightarrow \cE_{1}\otimes_E\cdots\otimes_E\cE_{1}\hookrightarrow \cE_{n}
\end{equation}
with the second map from (\ref{equ: Galois multi cup map}).
\begin{lem}\label{lem: Galois Ext decomposition}
For each $n\geq 2$, we have
\begin{equation}\label{equ: Galois log intersection}
\cE_{n}^{\sharp}\cap \cE_{n}^{\plog}=0
\end{equation}
and the intersection
\begin{equation}\label{equ: Galois diamond}
\cE_{n}^{\diamond}\defeq (\cE_{n}^{\sharp}\oplus\cE_{n}^{\plog})\cap \cE_{n}^{\flat}
\end{equation}
is a $1$-dimensional $E$-subspace of $\cE_{n}$ whose image under $\cE_{n}\twoheadrightarrow \mathrm{gr}_{W}^{0}(\cE_{n})$ (see (\ref{equ: universal Galois Ext graded})) is non-zero.
\end{lem}
\begin{proof}
We define $\cB\defeq \cD_1/(\cE_1^{\infty}\otimes_E\cR(\theta))$ which is special with level $1$. The surjection $\cD_1\twoheadrightarrow \cB$ induces a surjection $\mathrm{Sym}^{n}(\cD_1)\twoheadrightarrow \mathrm{Sym}^{n}(\cB)$ and thus a surjection
\begin{equation}\label{equ: special sym}
\cD_{n}\twoheadrightarrow \mathrm{Sym}^{n}(\cB)
\end{equation}
using \ref{it: sym 2} of Lemma~\ref{lem: universal to sym}.
We prove by an increasing induction on $n\geq 1$ that the composition of (\ref{equ: special sym}) with $\cE_{n}^{\plog}\otimes_E\cR(\theta^{n})\hookrightarrow \cD_{n}$ is an embedding. The $n=1$ case is obvious from $\cE_{1}=\cE_{1}^{\infty}\oplus\cE_{1}^{\plog}$. Assume from now that $n\geq 2$.
By our induction hypothesis the composition of
\begin{equation}\label{equ: Galois composition induction}
\cE_{n-1}^{\plog}\otimes_E\cR(\theta^{n-1})\hookrightarrow \cD_{n-1}\twoheadrightarrow \mathrm{Sym}^{n-1}(\cB)
\end{equation}
is an embedding.
Recall from (\ref{equ: Galois tensor sub}) that we have an embedding
\[\cD_1\otimes_{\cR}(\cE_{n-1}\otimes_E\cR(\theta^{n-1}))\hookrightarrow \cD_{n}\]
which together with $\cD_1\twoheadrightarrow \cB$ and the composition of (\ref{equ: Galois composition induction}) induces the following commutative diagram
\begin{equation}\label{equ: special sym diagram 1}
\xymatrix{
\cD_1\otimes_{\cR}(\cE_{n-1}^{\plog}\otimes_E\cR(\theta^{n-1})) \ar@{^{(}->}[r] \ar@{^{(}->}[dd] & \cD_1\otimes_{\cR}\mathrm{Sym}^{n-1}(\cB) \ar@{->>}[d]\\
& \cB\otimes_{\cR}\mathrm{Sym}^{n-1}(\cB) \ar@{->>}[d]\\
\cD_{n} \ar@{->>}[r] & \mathrm{Sym}^{n}(\cB)
}.
\end{equation}
Now that $\cD_1\twoheadrightarrow \cB$ restricts to an embedding $\cE_1^{\plog}\otimes_E\cR(\theta)\hookrightarrow \cB$, and $\cB\otimes_{\cR}\mathrm{Sym}^{n-1}(\cB)\twoheadrightarrow \mathrm{Sym}^{n}(\cB)$ restricts to an embedding $(\cE_1^{\plog}\otimes_E\cR(\theta))\otimes_{\cR}\mathrm{Sym}^{n-1}(\cB)\hookrightarrow \mathrm{Sym}^{n}(\cB)$, the diagram restricts to
\begin{equation}\label{equ: special sym diagram 2}
\xymatrix{
\cE_{n}^{\plog}\otimes_E\cR(\theta^{n}) \ar@{^{(}->}[r] \ar@{^{(}->}[d] & (\cE_1^{\plog}\otimes_E\cR(\theta))\otimes_{\cR}\mathrm{Sym}^{n-1}(\cB) \ar@{^{(}->}[d]\\
\cD_{n} \ar@{->>}[r] & \mathrm{Sym}^{n}(\cB)
}.
\end{equation}
Upon chasing the diagram (\ref{equ: special sym diagram 2}) we see that the composition of (\ref{equ: special sym}) with $\cE_{n}^{\plog}\otimes_E\cR(\theta^{n})\hookrightarrow \cD_{n}$ is an embedding, which finishes the proof of the induction step. Now that $\cE_{n}^{\sharp}\otimes_E\cR(\theta^{n})$ is by definition the kernel of the composition of $\cE_{n}\otimes_E\cR(\theta^{n})\hookrightarrow \cD_{n}\twoheadrightarrow\mathrm{Sym}^{n}(\cD_{1})$, we conclude (\ref{equ: Galois log intersection}). Recall from Proposition~\ref{prop: sharp flat decomposition} that $\cE_{n}=\cE_{n}^{\flat}\oplus\cE_{n}^{\sharp}$, and therefore
\[\Dim_E \cE_{n}^{\diamond}=\Dim_E \cE_{n}^{\flat}+\Dim_E \cE_{n}^{\sharp}\oplus\cE_{n}^{\plog}-\Dim_E \cE_{n}=1.\]
Recall from the proof of \ref{it: sym 1} of Lemma~\ref{lem: universal to sym} that the embedding $\cD_{n}^{\flat}\hookrightarrow \mathrm{Sym}^{n}(\cD_{1})$ induces an isomorphism
\begin{equation}\label{equ: Galois sym decomposition}
\cE_{n}^{\flat}\buildrel\sim\over\longrightarrow \mathrm{Sym}^{n}(\cE_1)=\mathrm{Sym}^{n}(\cE_1^{\infty}\oplus\cE_1^{\plog})\cong (\mathrm{Sym}^{n-1}(\cE_1^{\infty}\oplus\cE_1^{\plog})\otimes_E\cE_1^{\infty})\oplus\cE_{n}^{\plog}\cong (\cE_{n-1}^{\flat}\otimes_E\cE_1^{\infty})\oplus \cE_1^{\plog}.
\end{equation}
Let $\overline{\cD}_{n}$ be the unique de Rham $(\varphi,\Gamma)$-module that fits into the non-split extension (see \ref{it: char dR Ext 2} of Lemma~\ref{lem: char dR Ext})
\[0\rightarrow \cR(\theta^{n})\rightarrow \overline{\cD}_{n}\rightarrow \cR(1_{K^{\times}}) \rightarrow 0.\]
Using \ref{it: char dR Ext 1} of Lemma~\ref{lem: char dR Ext}, we see that the surjection $\cD_{n-1}\twoheadrightarrow \cR(1_{K^{\times}})$ induces an embedding
\[\mathrm{Ext}_{g}^1(\cR(1_{K^{\times}}),\cR(\theta^{n}))\hookrightarrow \mathrm{Ext}_{g}^1(\cD_{n-1},\cR(\theta^{n}))\]
and thus a surjection $\cD_{n}\twoheadrightarrow\overline{\cD}_{n}$. We also obtain a surjection $\cD_{n}^{\flat}\twoheadrightarrow\overline{\cD}_{n}$ similarly. Note that both $\cD_{n}\twoheadrightarrow\overline{\cD}_{n}$ and $\cD_{n}^{\flat}\twoheadrightarrow\overline{\cD}_{n}$ are the unique (up to scalars) non-zero such map by Lemma~\ref{lem: Galois unique Hom}. Combined with \ref{it: Galois sub flat 1} of Lemma~\ref{lem: Galois sub flat}, we obtain the following commutative diagram
\[
\xymatrix{
\cD_{n}^{\flat} \ar@{->>}[r] \ar@{^{(}->}[d]& \overline{\cD}_{n} \ar@{=}[d]\\
\cD_{n} \ar@{->>}[r] & \overline{\cD}_{n}
}
\]
which restricts to the following diagram
\begin{equation}\label{equ: grade 0 diagram 2}
\xymatrix{
\cE_{n}^{\flat}\otimes_E\cR(\theta^{n}) \ar@{->>}[r] \ar@{^{(}->}[d]& \mathrm{gr}_{W}^0(\cE_{n})\otimes_E\cR(\theta^{n}) \ar@{=}[d]\\
\cE_{n}\otimes_E\cR(\theta^{n}) \ar@{->>}[r] & \mathrm{gr}_{W}^0(\cE_{n})\otimes_E\cR(\theta^{n})
}.
\end{equation}
Note that the kernel of the top horizontal map of (\ref{equ: grade 0 diagram 2}) contains $(\cE_{n-1}^{\flat}\otimes_E\cE_1^{\infty})\otimes_E\cR(\theta^{n})$ and thus must equal it by dimension reason. This together with (\ref{equ: Galois sym decomposition}) forces the composition of
\[\cE_{n}^{\plog}\hookrightarrow \cE_{n}^{\flat}\twoheadrightarrow \mathrm{gr}_{W}^0(\cE_{n})\]
to be an isomorphism between $1$-dimensional $E$-vector spaces. The proof is thus finished.
\end{proof}

Using \ref{it: sym 2} of Lemma~\ref{lem: universal to sym} (upon replacing $n-1$ in \emph{loc.cit.} with $n$), we see that there exists a unique surjection
\begin{equation}\label{equ: universal to sym section}
\cE_{n}\twoheadrightarrow \mathrm{Sym}^{n}(\cE_{1})
\end{equation}
which admits the composition of
\[\mathrm{Sym}^{n}(\cE_{1})\hookrightarrow \cE_{1}\otimes_E\cdots\otimes_E\cE_{1}\hookrightarrow \cE_{n}\]
as a section. Here the middle term has $n$-copies of $\cE_{1}$ and the RHS map is from (\ref{equ: Galois multi cup map}).

\section{Drinfeld realization}\label{sec: application}
We apply results from previous sections to construct a certain universal filtered $(\varphi,N)$-module from the de Rham complex of the Drinfeld upper half space.
\subsection{De Rham complex of the Drinfeld upper half space}\label{subsec: Drinfeld}
We recall some known results on the de Rham complex of the Drinfeld upper half space from \cite{SS91}, \cite{ST02}, \cite{Dat06}, \cite{Or08}, \cite{Schr11} and more recently \cite{BQ24}.

We consider the Drinfeld space $\bH_{/K}$ of dimension $n-1$ defined as $\bH\defeq \bP^{n-1}_{\rm{rig}}\setminus\bigcup_{\cH}\cH$ where $\cH$ runs through $K$-rational hyperplanes inside $\bP^{n-1}_{\rm{rig}}$. Recall from \cite[Prop.~1.4]{SS91} that $\bH$ is quasi-Stein and that for any coherent sheaf $\cF$ on $\bH$ we have $H^k(\cF)=0$ for $k>0$. We follow the convention of \cite{ST02} and \cite{Or08} and endow $\bH$ with the left action of $G=\GL_n(K)$ coming from the left action of $G$ on $\bP^{n-1}_{\rm{rig}}$ given by $(z_0,\dots,z_{n-1})\mapsto (z_0,\dots,z_{n-1})g^{-1}$ (matrix product) for $g\in G$. This action in turn induces an action of $G$ (which we will always take on the left) on the global sections $\cF(\bH)$ for any $G$-equivariant vector bundle $\cF$ on $\bH$. Moreover if $\cF$ is the restriction to $\bH$ of a $G$-equivariant vector bundle on $\bP^{n-1}_{\rm{rig}}$, then by the argument in \cite[p.593]{Or08} $\cF(\bH)$ is naturally a coadmissible $D(G,K)$-module and thus $\cF(\bH)\otimes_{K,\iota}E$ is a coadmissible $D(G)=D(G,E)$-module. For instance this applies to the sheaf of differential $k$-forms on $\bH$.\\

We obtain a complex of the following form
\begin{equation}\label{equ: Drinfeld dR}
\Omega^\bullet=[\Omega^0\longrightarrow \Omega^1\longrightarrow \cdots \longrightarrow \Omega^{n-1}]
\end{equation}
by considering the \emph{global sections} of the de Rham complex of $\bH$ and then applying $-\otimes_{K,\iota}E$.
By the above discussion this is a complex of coadmissible $D(G)$-modules. The action of $G$ induces an action on each cohomology group $H^k(\Omega^\bullet)$ for $0\leq k\leq n-1$.
We continue to use the shortened notation $J^k=[1,n-1-k]\subseteq \Delta$ for $0\leq k\leq n-1$, with $J^0=\Delta$ and $J^{n-1}=\emptyset$ (see also the beginning of \S \ref{subsec: filtered phi N}).
We recall the following seminal result of Schneider--Stuhler \cite[Thm.~3.1, Lem.~4.1]{SS91}.

\begin{thm}\label{thm: dR coh}
We have a $G$-equivariant isomorphism $H^k(\Omega^\bullet)\cong (V_{J^k}^{\infty})^\vee$ for $0\leq k\leq n-1$.
\end{thm}

For each $j,j'\in\Delta$, we define $w_{j,j'}=s_js_{j-1}\cdots s_{j'}$ if $j\geq j'$, and $w_{j,j'}=s_js_{j+1}\cdots s_{j'}$ if $j\leq j'$.
We recall the following results from \cite[(547), Thm.~5.3.12, Thm.~5.4.12]{BQ24}, which refine results in \cite{Or08} and \cite[Thm.~2.2]{Or13}. Recall from \S \ref{subsec: notation} that given a multiplicity free finite length $D(G)$-module $D$, the set $\mathrm{JH}_{D(G)}(D)$ is equipped with a natural partial-order.
\begin{thm}\label{thm: dR complex}
For each $0\leq k\leq n-1$, $\Omega^k$ is a multiplicity free finite length coadmissible $D(G)$-module with the following description of the partially-ordered set $\mathrm{JH}_{D(G)}(\Omega^k)$.
\begin{enumerate}[label=(\roman*)]
\item \label{it: dR complex 1} We have \[\mathrm{JH}_{D(G)}(\Omega^0)=\{(V_{\Delta}^{\infty})^\vee\}\sqcup\{C_{w_{j,n-1},[1,j-1]}^\vee\mid j\in\Delta\}\] with its partial-order generated by $(V_{\Delta}^{\infty})^\vee<C_{s_{n-1},[1,j-1]}^\vee$ and $C_{w_{j+1,n-1},[1,j]}^\vee<C_{w_{j,n-1},[1,j-1]}^\vee$ for each $1\leq j\leq n-2$;
\item \label{it: dR complex 2} We have \[\mathrm{JH}_{D(G)}(\Omega^{n-1})=\{(V_{\emptyset}^{\infty})^\vee\}\sqcup\{C_{w_{j,1},[1,j-1]}^\vee\mid j\in\Delta\}\] with its partial-order generated by $(V_{\emptyset}^{\infty})^\vee>C_{s_{1},\emptyset}^\vee$ and $C_{w_{j+1,1},[1,j]}^\vee<C_{w_{j,1},[1,j-1]}^\vee$ for each $1\leq j\leq n-2$;
\item \label{it: dR complex 3} For each $1\leq k\leq n-2$, we have \[\mathrm{JH}_{D(G)}(\Omega^k)=\{(V_{[1,n-1-k]}^{\infty})^\vee\}\sqcup\{C_{w_{j,n-1-k},[1,j-1]}^\vee,C_{w_{j,n-k},[1,j-1]}^\vee\mid j\in\Delta\}\]
    with its partial-order generated by $C_{s_{n-k},[1,n-1-k]}^\vee<(V_{[1,n-1-k]}^{\infty})^\vee<C_{s_{n-1-k},[n-2-k]}^\vee$, $C_{w_{j,n-k},[1,j-1]}^\vee<C_{w_{j,n-1-k},[1,j-1]}^\vee$ for each $j\in\Delta$, and $C_{w_{j+1,j'},[1,j]}^\vee<C_{w_{j,j'},[1,j-1]}^\vee$ for each $1\leq j\leq n-2$ and $j'\in\{n-1-k,n-k\}$;
\item \label{it: dR complex 4} For each $1\leq k\leq n-1$, $\Omega^{k-1}/\mathrm{ker}(d_{\Omega^{\bullet}}^{k-1})\cong\mathrm{im}(d_{\Omega^{\bullet}}^{k-1})\subseteq \Omega^k$ is uniserial of length $n-1$, with $\mathrm{JH}_{D(G)}(\mathrm{im}(d_{\Omega^{\bullet}}^{k-1}))=\{C_{w_{j,n-k},[1,j]}^\vee\mid j\in\Delta\}$ and its partial-order generated by $C_{w_{j+1,1},[n-k,j]}^\vee<C_{w_{j,1},[n-k,j-1]}^\vee$ for each $1\leq j\leq n-2$.
\end{enumerate}
\end{thm}

Recall from \S \ref{subsec: loc an rep} that $\cM(G)$ is the bounded derived category of the abelian category $\mathrm{Mod}_{D(G)}$ of abstract $D(G)$-modules.
It follows from \cite{Dat06} and \cite[\S 6.1]{Schr11} (see also Lemma~\ref{lem: Ext between sm}) that we may fix the choice of an isomorphism (splitting) in $\cM(G)$ of the following form
\begin{equation}\label{equ: abstract dR splitting}
s: \Omega^{\bullet}\buildrel\sim\over\longrightarrow \bigoplus_{k=0}^{n-1}H^k(\Omega^{\bullet})[-k],
\end{equation}
as well as the choice of isomorphisms in Theorem~\ref{thm: dR coh}, which altogether induce a map in $\cM(G)$
\begin{equation}\label{equ: abstract dR splitting ell}
s_k: \Omega^{\bullet}\rightarrow H^k(\Omega^{\bullet})[-k]\cong
(V_{J^k}^{\infty})^\vee[-k]
\end{equation}
for each $0\leq k\leq n-1$. Note that the splitting (\ref{equ: abstract dR splitting}) also induces a splitting (defined as the composition of $\bigoplus_{k=0}^{\ell}s_k$ with $\tau^{\leq \ell}\Omega^{\bullet}\rightarrow \Omega^{\bullet}$)
\begin{equation}\label{equ: abstract dR splitting truncation}
s_{\leq \ell}: \tau^{\leq \ell}\Omega^{\bullet}\buildrel\sim\over\longrightarrow \bigoplus_{k=0}^{\ell}H^k(\Omega^{\bullet})[-k]\cong
\bigoplus_{k=0}^{\ell}(V_{J^{k}}^{\infty})^\vee[-k]
\end{equation}
for each $0\leq \ell\leq n-1$.

\begin{rem}\label{rem: explicit splitting}
In \cite{BQ24}, we have constructed an explicit splitting
\[s_{n-1}: \Omega^{\bullet}\rightarrow H^{n-1}(\Omega^{\bullet})[1-n]\cong
(V_{\emptyset}^{\infty})^\vee[1-n]=(\mathrm{St}_G^{\infty})^\vee[1-n]\]
of the top degree cohomology of $\Omega^{\bullet}$.
\end{rem}

\subsection{$(\varphi, N)$-action}\label{subsec: phi N}
Based on the choice of a splitting (\ref{equ: abstract dR splitting}), we define a $(\varphi,N)$-action on $\Omega^{\bullet}$ such that $\varphi$ acts on $H^k(\Omega^{\bullet})$ by $p^k$ for $0\leq k\leq n-1$, $N^{n-1}\neq 0$ and $N\varphi=p\varphi N$.

It follows from Lemma~\ref{lem: Ext between sm} that we have
\begin{multline}\label{equ: End coh sum}
\mathrm{End}_{\cM(G)}(\bigoplus_{k=0}^{n-1}(V_{J^k}^{\infty})^\vee[-k])\\
\cong \bigoplus_{0\leq k\leq \ell\leq n-1}\Hom_{\cM(G)}((V_{J^{\ell}}^{\infty})^\vee[-\ell],(V_{J^k}^{\infty})^\vee[-k])\\
\cong\bigoplus_{0\leq k\leq \ell\leq n-1}\mathrm{Ext}^{\ell-k}(V_{J^k}^{\infty},V_{J^{\ell}}^{\infty}).
\end{multline}
We define $\varphi$ as the element of (\ref{equ: End coh sum}) that satisfies $\varphi((V_{J^k}^{\infty})^\vee[-k])=(V_{J^k}^{\infty})^\vee[-k]$ and
\begin{equation}\label{equ: varphi coh sum}
\varphi|_{(V_{J^k}^{\infty})^\vee[-k]}=p^k
\end{equation}
for each $0\leq k\leq n-1$.
Recall from Lemma~\ref{lem: Ext between sm} 
the following isomorphism
\begin{equation}\label{equ: N as Ext1}
\Hom_{\cM(G)}((V_{J^k}^{\infty})^\vee[-k],(V_{J^{k-1}}^{\infty})^\vee[-k+1])=\mathrm{Ext}_{G}^1(V_{J^{k-1}}^{\infty},V_{J^k}^{\infty})\cong \Hom_{\infty}(K^\times,E)=E\val
\end{equation}
for each $1\leq k\leq n-1$.
Hence, using the identification (\ref{equ: N as Ext1}), we may define $N$ as an arbitrary element of (\ref{equ: End coh sum}) that satisfies $N\varphi=p\varphi N$ and $N^{n-1}\neq 0$. According to (\ref{equ: varphi coh sum}), we must have $N((V_{J^0}^{\infty})^\vee)=0$ and then $N((V_{J^k}^{\infty})^\vee[-k])\subseteq (V_{J^{k-1}}^{\infty})^\vee[-k+1]$, with the restriction $N|_{(V_{J^k}^{\infty})^\vee[-k]}$ being a non-zero element of (\ref{equ: N as Ext1}),
namely $c_{k}\val$ for some $c_{k}\in E^{\times}$, for each $1\leq k\leq n-1$.
We consider the following functors
\begin{equation}\label{equ: dR functor}
\Hom_{\cM(G)}(-,\Omega^{\bullet})\buildrel s_{\ast}\over\longrightarrow \Hom_{\cM(G)}(-, \bigoplus_{k=0}^{n-1}H^k(\Omega^{\bullet})[-k])\cong \Hom_{\cM(G)}(-,\bigoplus_{k=0}^{n-1}(V_{J^k}^{\infty})^\vee[-k])
\end{equation}
from $\cM(G)$ to the abelian category of $E$-vector spaces.
The $(\varphi,N)$-action on $\bigoplus_{k=0}^{n-1}(V_{J^k}^{\infty})^\vee[-k]$ induces a $(\varphi,N)$-action on $\Omega^{\bullet}$, and therefore on $\Hom_{\cM(G)}(C,\Omega^{\bullet})$ for each object $C\in\cM(G)$. It is clear that the $(\varphi,N)$-action on $\Hom_{\cM(G)}(C,\Omega^{\bullet})$ depends on the choice of the splitting $s$ from (\ref{equ: abstract dR splitting}).

\subsection{Hodge filtration}\label{subsec: Hodge}
We begin our study of
\begin{equation}\label{equ: Hom St}
\Hom_{\cM(G)}((\mathrm{St}_G^{\rm{an}})^\vee[1-n],\Omega^{\bullet}).
\end{equation}
Our choice of $(\varphi,N)$-action on $\Omega^{\bullet}$ induces a $(\varphi,N)$-action on (\ref{equ: Hom St}).
We will prove in \ref{it: Hodge seq 0} of Lemma~\ref{lem: Hodge exact seq} that the naive truncations $\sigma_{\geq k}\Omega^{\bullet}$ induce a \emph{Hodge filtration}
\begin{equation}\label{equ: Hom St filtration}
\{\Hom_{\cM(G)}((\mathrm{St}_G^{\rm{an}})^\vee[1-n], \sigma_{\geq k}\Omega^{\bullet})\}_{0\leq k\leq n-1}
\end{equation}
on (\ref{equ: Hom St}).
The rest of this section is devoted to a detailed study of the filtered $(\varphi,N)$-module (\ref{equ: Hom St}). We will prove that it is \emph{universal special with level $n-1$} (see Lemma~\ref{lem: Hodge image} and Definition~\ref{def: universal phi N})

For $0\leq \ell\leq n-1$, we have the following truncation maps
\[\sigma_{\geq \ell}\Omega^{\bullet}\rightarrow \Omega^{\bullet}\rightarrow \tau_{\geq \ell}\Omega^{\bullet}\]
and
\[\tau^{\leq \ell}\Omega^{\bullet}\rightarrow \Omega^{\bullet}\rightarrow \sigma^{\leq \ell}\Omega^{\bullet}.\]
Since we have
\begin{multline}
[0\rightarrow \cdots \rightarrow 0\rightarrow \Omega^{\ell}/\mathrm{ker}(d_{\Omega^{\bullet}}^{\ell})\rightarrow \Omega^{\ell+1}\rightarrow \cdots\rightarrow \Omega^{n-1}]\\
\buildrel\sim\over\longrightarrow[0\rightarrow \cdots \rightarrow 0\rightarrow 0\rightarrow \Omega^{\ell+1}/\mathrm{im}(d_{\Omega^{\bullet}}^{\ell})\rightarrow \cdots\rightarrow \Omega^{n-1}]=\tau_{\geq \ell+1}\Omega^{\bullet}
\end{multline}
in $\cM(G)$, we see that the short exact sequence $0\rightarrow \mathrm{ker}(d_{\Omega^{\bullet}}^{\ell})\rightarrow \Omega^{\ell}\rightarrow \Omega^{\ell}/\mathrm{ker}(d_{\Omega^{\bullet}}^{\ell})\rightarrow 0$ in $\mathrm{Mod}_{D(G)}$ induces a distinguished triangle
\begin{equation}\label{equ: Hodge induction triangle 1}
\mathrm{ker}(d_{\Omega^{\bullet}}^{\ell})[-\ell]\rightarrow \sigma_{\geq \ell}\Omega^{\bullet} \rightarrow \tau_{\geq \ell+1}\Omega^{\bullet} \rightarrow
\end{equation}
in $\cM(G)$.
For each $1\leq \ell\leq n-1$, the short exact sequence $0\rightarrow \mathrm{im}(d_{\Omega^{\bullet}}^{\ell-1})\rightarrow \Omega^{\ell}\rightarrow \Omega^{\ell}/\mathrm{im}(d_{\Omega^{\bullet}}^{\ell-1})\rightarrow 0$ in $\mathrm{Mod}_{D(G)}$ also induces a distinguished triangle
\begin{equation}\label{equ: Hodge induction triangle 2}
\mathrm{im}(d_{\Omega^{\bullet}}^{\ell-1})[-\ell]\rightarrow \sigma_{\geq \ell}\Omega^{\bullet} \rightarrow \tau_{\geq \ell}\Omega^{\bullet} \rightarrow
\end{equation}
in $\cM(G)$.

The results in the rest of this section generalizes that of \cite[\S~6.5]{Schr11} from the case of $\mathrm{PGL}_{3}(\Q_p)$ to the general cases of $G=\mathrm{PGL}_{n}(K)$.
\begin{lem}\label{lem: Hodge exact seq}
Let $0\leq \ell\leq n-1$. We have the following results.
\begin{enumerate}[label=(\roman*)]
\item \label{it: Hodge seq 0} The maps $\mathrm{ker}(d_{\Omega^{\bullet}}^{\ell})[-\ell]\rightarrow \sigma_{\geq \ell}\Omega^{\bullet}\rightarrow \Omega^{\bullet}$ induce embeddings
\begin{multline}\label{equ: Hodge seq 0}
\Hom_{\cM(G)}((\mathrm{St}_G^{\rm{an}})^\vee[1-n],\mathrm{ker}(d_{\Omega^{\bullet}}^{\ell})[-\ell])\hookrightarrow
\Hom_{\cM(G)}((\mathrm{St}_G^{\rm{an}})^\vee[1-n],\sigma_{\geq \ell}\Omega^{\bullet})\\
\hookrightarrow \Hom_{\cM(G)}((\mathrm{St}_G^{\rm{an}})^\vee[1-n],\Omega^{\bullet})
\end{multline}
\item \label{it: Hodge seq 1} We have the following equality
\begin{equation}\label{equ: Hodge seq 1}
\Dim_E \Hom_{\cM(G)}((\mathrm{St}_G^{\rm{an}})^\vee[1-n],\tau_{\geq \ell}\Omega^{\bullet})=\sum_{\ell'=\ell}^{n-1}\Dim_E \mathbf{E}_{J^{\ell'},\emptyset}.
\end{equation}
\item \label{it: Hodge seq 2} The distinguished triangle (\ref{equ: Hodge induction triangle 2}) induces an isomorphism
\begin{equation}\label{equ: Hodge seq 2}
\Hom_{\cM(G)}((\mathrm{St}_G^{\rm{an}})^\vee[-m],\sigma_{\geq \ell}\Omega^{\ell})\buildrel\sim\over\longrightarrow \Hom_{\cM(G)}((\mathrm{St}_G^{\rm{an}})^\vee[-m],\tau_{\geq \ell}\Omega^{\ell})
\end{equation}
for each $m\in\Z$.
\item \label{it: Hodge seq 3} The surjection $\mathrm{ker}(d_{\Omega^{\bullet}}^{\ell})\twoheadrightarrow H^{\ell}(\Omega^{\bullet})\cong (V_{J^{\ell}}^{\infty})^\vee$ induces an isomorphism
\begin{equation}\label{equ: Hodge seq 3}
\Hom_{\cM(G)}((\mathrm{St}_G^{\rm{an}})^\vee[-m],\mathrm{ker}(d_{\Omega^{\bullet}}^{\ell})[-\ell])\buildrel\sim\over\longrightarrow \Hom_{\cM(G)}((\mathrm{St}_G^{\rm{an}})^\vee[-m],(V_{J^{\ell}}^{\infty})^\vee[-\ell])
\end{equation}
for each $m\in\Z$.
\item \label{it: Hodge seq 4} The distinguished triangle (\ref{equ: Hodge induction triangle 1}) induces the following short exact sequence
\begin{multline}\label{equ: Hodge seq 4}
0\rightarrow \Hom_{\cM(G)}((\mathrm{St}_G^{\rm{an}})^\vee[1-n],\mathrm{ker}(d_{\Omega^{\bullet}}^{\ell})[-\ell])\rightarrow \Hom_{\cM(G)}((\mathrm{St}_G^{\rm{an}})^\vee[1-n],\sigma_{\geq \ell}\Omega^{\bullet})\\
\rightarrow \Hom_{\cM(G)}((\mathrm{St}_G^{\rm{an}})^\vee[1-n],\tau_{\geq \ell+1}\Omega^{\bullet})\rightarrow 0.
\end{multline}
\item \label{it: Hodge seq 5} We have the following equality
\begin{multline}\label{equ: Hodge seq 5}
\Hom_{\cM(G)}((\mathrm{St}_G^{\rm{an}})^\vee[1-n],\sigma_{\geq \ell}\Omega^{\bullet})\\
= \Hom_{\cM(G)}((\mathrm{St}_G^{\rm{an}})^\vee[1-n],\mathrm{ker}(d_{\Omega^{\bullet}}^{\ell})[-\ell])\oplus \Hom_{\cM(G)}((\mathrm{St}_G^{\rm{an}})^\vee[1-n],\sigma_{\geq \ell+1}\Omega^{\bullet})
\end{multline}
as $E$-subspaces of $\Hom_{\cM(G)}((\mathrm{St}_G^{\rm{an}})^\vee[1-n],\Omega^{\bullet})$.
\end{enumerate}
\end{lem}
\begin{proof}
Let $0\leq \ell\leq n-1$ and $R$ be a constituent of $\Omega^{\ell}$. By Theorem~\ref{thm: dR complex} we know that either $R\cong (V_{J^{\ell}}^{\infty})^\vee$ or $R\cong C_{x,I}^\vee$ for some $x\neq 1$ and $I\subseteq I_x$. If $R\cong C_{x,I}^\vee$ for some $x\neq 1$ and $I\subseteq I_x$, it follows from \ref{it: change left cup 1} of Lemma~\ref{lem: change left cup} (by taking $I_0=\emptyset$, $I_1=J^{\ell}$ and $W=C_{x,I}$ in \emph{loc.cit.}) that we have
\begin{equation}\label{equ: dR non sm vanishing}
\Hom_{\cM(G)}((\mathrm{St}_G^{\rm{an}})^\vee[-m],C_{x,I}^\vee[-\ell])=0
\end{equation}
for each $m\in\Z$. If $R\cong (V_{J^{\ell}}^{\infty})^\vee$, it follows from \ref{it: change left cup 2} of Lemma~\ref{lem: change left cup} (by taking $I_1=\emptyset$, $I_0=J^{\ell}$, $V'=V_{J^{\ell}}^{\infty}$ and $V=V_{J^{\ell}}^{\rm{an}}$ in \emph{loc.cit.}) we have
\begin{equation}\label{equ: dR sm 1}
\Hom_{\cM(G)}((\mathrm{St}_G^{\rm{an}})^\vee[-m],(V_{J^{\ell}}^{\infty})^\vee[-\ell])=\mathrm{Ext}_{G}^{m-\ell}(V_{J^{\ell}}^{\infty},\mathrm{St}_G^{\rm{an}})\buildrel\sim\over\longleftarrow \mathrm{Ext}_{G}^{m-\ell}(V_{J^{\ell}}^{\rm{an}},\mathrm{St}_G^{\rm{an}}),
\end{equation}
which by \ref{it: Ext between St 1} of Corollary~\ref{cor: Ext between St} is vanishing for $m<n-1$ and gives an isomorphism
\begin{equation}\label{equ: dR sm 2}
\mathbf{E}_{J^{\ell},\emptyset}\buildrel\sim\over\longrightarrow\Hom_{\cM(G)}((\mathrm{St}_G^{\rm{an}})^\vee[1-n],(V_{J^{\ell}}^{\infty})^\vee[-\ell])
\end{equation}
when $m=n-1$.
It follows from (\ref{equ: dR non sm vanishing}) and a d\'evissage with respect to the constituents of $\mathrm{im}(d^{\ell-1}_{\Omega^{\bullet}})$ (see \ref{it: dR complex 4} of Theorem~\ref{thm: dR complex}) that we have
\[\Hom_{\cM(G)}((\mathrm{St}_G^{\rm{an}})^\vee[-m],\mathrm{im}(d^{\ell-1}_{\Omega^{\bullet}})[-\ell])=0\]
for each $m\in\Z$, which together with the short exact sequence $0\rightarrow \mathrm{im}(d^{\ell-1}_{\Omega^{\bullet}})\rightarrow \mathrm{ker}(d^{\ell}_{\Omega^{\bullet}})\rightarrow (V_{J^{\ell}}^{\infty})^\vee\rightarrow 0$ and the distinguished triangle $\mathrm{im}(d^{\ell-1}_{\Omega^{\bullet}})[-\ell]\rightarrow \sigma_{\geq}\Omega^{\bullet}\rightarrow \tau_{\geq \ell}\Omega^{\bullet}\rightarrow $ gives
\ref{it: Hodge seq 2} and \ref{it: Hodge seq 3} respectively.

We prove \ref{it: Hodge seq 0}. Note that both complex $\tau_{\geq \ell+1}\Omega^{\bullet}$ and $\sigma^{\leq \ell-1}\Omega^{\bullet}$ admit a filtration with each graded piece being $R[-\ell']$ for some constituent $R$ be a constituent of $\Omega^{\ell'}$ for some $0\leq \ell'\leq n-1$. Hence, by combining (\ref{equ: dR non sm vanishing}) and (\ref{equ: dR sm 1}) with a d\'evissage on such filtration on $\tau_{\geq \ell+1}\Omega^{\bullet}$ and $\sigma^{\leq \ell-1}\Omega^{\bullet}$, we deduce that
\begin{equation}\label{equ: Hodge seq vanishing}
\Hom_{\cM(G)}((\mathrm{St}_G^{\rm{an}})^\vee[-m],\tau_{\geq \ell+1}\Omega^{\bullet})=0=\Hom_{\cM(G)}((\mathrm{St}_G^{\rm{an}})^\vee[-m],\sigma^{\leq \ell-1}\Omega^{\bullet})
\end{equation}
for each $m<n-1$, which together with the distinguished triangles (\ref{equ: Hodge induction triangle 1}) and $\sigma_{\geq \ell}\Omega^{\bullet}\rightarrow \Omega^{\bullet}\rightarrow \sigma^{\leq \ell-1}\Omega^{\bullet}\rightarrow $ finishes the proof.

We prove \ref{it: Hodge seq 1}. The splitting $s$ from (\ref{equ: abstract dR splitting}) induces a splitting (by taking the composition of $\Omega^{\bullet}\rightarrow \tau_{\geq \ell}\Omega^{\bullet}$ with $s^{-1}$)
\[\tau_{\geq \ell}\Omega^{\bullet}\buildrel\sim\over\longleftarrow \bigoplus_{\ell'=\ell}^{n-1}H^{\ell'}(\Omega^{\bullet})[-\ell']\cong \bigoplus_{\ell'=\ell}^{n-1}(V_{J^{\ell'}}^{\infty})^\vee[-\ell'],\]
which together with (\ref{equ: dR sm 2}) (by replacing $\ell$ in \emph{loc.cit.} with each $\ell'\geq \ell$) gives an isomorphism
\[\Hom_{\cM(G)}((\mathrm{St}_G^{\rm{an}})^\vee[1-n],\tau_{\geq \ell}\Omega^{\bullet})\buildrel\sim\over\longleftarrow \bigoplus_{\ell'=\ell}^{n-1}\mathbf{E}_{J^{\ell'},\emptyset}\]
and in particular (\ref{equ: Hodge seq 1}).

We prove \ref{it: Hodge seq 4}. The distinguished triangle (\ref{equ: Hodge induction triangle 1}) together with the vanishing (\ref{equ: Hodge seq vanishing}) for $m<n-1$ induces an exact sequence
\begin{multline}\label{equ: Hodge seq 4 seq}
\Hom_{\cM(G)}((\mathrm{St}_G^{\rm{an}})^\vee[1-n],\mathrm{ker}(d_{\Omega^{\bullet}}^{\ell})[-\ell])\rightarrow \Hom_{\cM(G)}((\mathrm{St}_G^{\rm{an}})^\vee[1-n],\sigma_{\geq \ell}\Omega^{\bullet})\\
\rightarrow \Hom_{\cM(G)}((\mathrm{St}_G^{\rm{an}})^\vee[1-n],\tau_{\geq \ell+1}\Omega^{\bullet}).
\end{multline}
It follows from $m=n-1$ case of \ref{it: Hodge seq 3} and (\ref{equ: dR sm 2}) that we have
\begin{equation}\label{equ: Hodge seq 4 dim 1}
\Dim_E \Hom_{\cM(G)}((\mathrm{St}_G^{\rm{an}})^\vee[1-n],\mathrm{ker}(d_{\Omega^{\bullet}}^{\ell})[-\ell])=\Dim_E \mathbf{E}_{J^{\ell},\emptyset}.
\end{equation}
By combing \ref{it: Hodge seq 1} with \ref{it: Hodge seq 2} we have
\begin{equation}\label{equ: Hodge seq 4 dim 2}
\Dim_E \Hom_{\cM(G)}((\mathrm{St}_G^{\rm{an}})^\vee[1-n],\sigma_{\geq \ell}\Omega^{\bullet})=\Dim_E \Hom_{\cM(G)}((\mathrm{St}_G^{\rm{an}})^\vee[1-n],\tau_{\geq \ell}\Omega^{\bullet})=\sum_{\ell'=\ell}^{n-1}\Dim_E \mathbf{E}_{J^{\ell'},\emptyset}.
\end{equation}
It follows again from \ref{it: Hodge seq 1} (by replacing $\ell$ with $\ell+1$ in \emph{loc.cit.}) that
\begin{equation}\label{equ: Hodge seq 4 dim 3}
\Dim_E \Hom_{\cM(G)}((\mathrm{St}_G^{\rm{an}})^\vee[1-n],\tau_{\geq \ell+1}\Omega^{\bullet})=\sum_{\ell'=\ell+1}^{n-1}\Dim_E \mathbf{E}_{J^{\ell'},\emptyset}.
\end{equation}
By combining (\ref{equ: Hodge seq 4 dim 1}), (\ref{equ: Hodge seq 4 dim 2}) and (\ref{equ: Hodge seq 4 dim 3}), we see that the exact sequence (\ref{equ: Hodge seq 4 seq}) must be a short exact sequence. This finishes the proof of \ref{it: Hodge seq 4}.

We prove \ref{it: Hodge seq 5}. Similar to the argument of \ref{it: Hodge seq 0}, we see that the map $\sigma_{\geq \ell+1}\Omega^{\bullet}\rightarrow \sigma_{\geq \ell}\Omega^{\bullet}$ induces an embedding
\[\Hom_{\cM(G)}((\mathrm{St}_G^{\rm{an}})^\vee[1-n],\sigma_{\geq \ell+1}\Omega^{\bullet})\hookrightarrow \Hom_{\cM(G)}((\mathrm{St}_G^{\rm{an}})^\vee[1-n],\sigma_{\geq \ell}\Omega^{\bullet}).\]
It follows from \ref{it: Hodge seq 2} (by replacing $\ell$ with $\ell+1$) that the composition $\sigma_{\geq \ell+1}\Omega^{\bullet}\rightarrow \sigma_{\geq \ell}\Omega^{\bullet}\rightarrow \tau_{\geq \ell+1}\Omega^{\bullet}$ induces an isomorphism
\[\Hom_{\cM(G)}((\mathrm{St}_G^{\rm{an}})^\vee[1-n],\sigma_{\geq \ell+1}\Omega^{\bullet})\buildrel\sim\over\longrightarrow \Hom_{\cM(G)}((\mathrm{St}_G^{\rm{an}})^\vee[1-n],\tau_{\geq \ell+1}\Omega^{\bullet}),\]
which together with the short exact sequence (\ref{equ: Hodge seq 4}) gives the equality (\ref{equ: Hodge seq 5}) as $E$-subspaces of $\Hom_{\cM(G)}((\mathrm{St}_G^{\rm{an}})^\vee[1-n],\Omega^{\bullet})$.
\end{proof}

\begin{lem}\label{lem: dR Ext vanishing}
Let $0\leq \ell'\leq \ell\leq n-2$. Then we have
\begin{equation}\label{equ: dR Ext vanishing}
\mathrm{Ext}_{G}^k(V_{J^{\ell'}}^{\infty},(\Omega^{\ell})^\vee)=0
\end{equation}
for each $k\leq \ell-\ell'$.
\end{lem}
\begin{proof}
If $\ell'=\ell\leq n-2$, then by Theorem~\ref{thm: dR complex} we know that $V_{J^{\ell}}^{\infty}$ does not show up in the socle of $(\Omega^{\ell})^\vee$ and thus $\Hom_{G}(V_{J^{\ell}}^{\infty},(\Omega^{\ell})^\vee)=0$.
Assume from now that $\ell'<\ell\leq n-2$, then by \ref{it: dR complex 3} of Theorem~\ref{thm: dR complex} we know that $(\Omega^{\ell})^\vee$ admits a filtration with each graded piece $V$ being either length $2$ with socle $C_{w_{j,n-1-\ell},[1,j-1]}$ and cosocle $C_{w_{j,n-\ell},[1,j-1]}$ for some $n-\ell\leq j\leq n-1-\ell'$, or being length $2$ with socle $C_{s_{n-1-\ell},[1,n-2-\ell]}$ and cosocle $V_{J^{\ell}}^{\infty}=V_{[1,n-1-\ell]}^{\infty}$, or being $C_{w_{j,j'},[1,j-1]}$ for some $j\in[1,n-1-\ell]\sqcup[n-\ell',n-1]$ and $j'\in\{n-1-\ell,n-\ell\}$.
In the third case with $V$ being irreducible, we have
\begin{equation}\label{equ: dR Ext vanishing subquotient}
\mathrm{Ext}_{G}^k(V_{J^{\ell'}}^{\infty},V)=0
\end{equation}
for $k\leq \ell-\ell'$ by \ref{it: Ext with sm 1} and \ref{it: Ext with sm 2} of Lemma~\ref{lem: Ext with sm} (by taking $I_0=J^{\ell'}$, $w=w_{j,j'}$ and $I=[1,j-1]$ in \emph{loc.cit.}).
In the first case with $V$ being length $2$ with socle $C_{w_{j,n-1-\ell},[1,j-1]}$ and cosocle $C_{w_{j,n-\ell},[1,j-1]}$ for some $n-\ell\leq j\leq n-1-\ell'$, we have $V\cong \cF_{P_{\widehat{j}}}^{G}(M,V_{[1,j-1],\widehat{j}}^{\infty})$ (cf.~ Remark~\ref{rem: length 2 g OS}) with $\widehat{j}\defeq \Delta\setminus\{j\}$ and $M$ being the unique length $2$ $U(\fg)$-module with socle $L(w_{j,n-\ell})$ and cosocle $L(w_{j,n-1-\ell})$ (see Lemma~\ref{lem: g Ext1}), in which case (\ref{equ: dR Ext vanishing subquotient}) holds for $k\leq \ell-\ell'$ by \ref{it: Ext with sm 3} of Lemma~\ref{lem: Ext with sm} (by taking $I_0=J^{\ell'}$, $x=w_{j,n-\ell}$, $w=w_{j,n-1-\ell}$ and $I=[1,j-1]$ in \emph{loc.cit.}).
In this second case with $V$ being length $2$ with socle $C_{s_{n-1-\ell},[1,n-2-\ell]}$ and cosocle $V_{J^{\ell}}^{\infty}=V_{[1,n-1-\ell]}^{\infty}$, (\ref{equ: dR Ext vanishing subquotient}) holds for $k\leq \ell-\ell'$ due to Lemma~\ref{lem: special Ext with sm} (by taking $j=n-1-\ell$, $I=[1,n-2-\ell]$ and $J=J^{\ell'}\supseteq I\sqcup\{j\}$ in \emph{loc.cit.}).
As (\ref{equ: dR Ext vanishing subquotient}) holds for $k\leq \ell-\ell'$ for each graded piece $V$ of the aforementioned filtration on $(\Omega^{\ell})^\vee$, we deduce (\ref{equ: dR Ext vanishing}) by a d\'evissage with respect to this filtration.
\end{proof}

\begin{rem}\label{rem: geometric proof}
In a private discussion, we are informed by Benchao Su that Lemma~\ref{lem: dR Ext vanishing} admits a more geometric and conceptual proof without using the explicit structure of $\Omega^{\ell}$ as an input.
\end{rem}

We recall our notation introduced at the beginning of \S \ref{subsec: filtered phi N}, including 
\[\mathbf{E}_{k,k'}=\mathbf{E}_{J^{k},J^{k'}}\buildrel\sim\over\longrightarrow\mathrm{Ext}_{G}^{\#J^{k}\setminus J^{k'}}(V_{J^{k}}^{\infty},V_{J^{k'}}^{\rm{an}}),\] 
together with $\mathbf{E}_{k,k'}^{\infty}\defeq \mathbf{E}_{J^{k},J^{k'}}^{\infty}$, $\mathbf{E}_{k,k'}^{<}$, $\overline{\mathbf{E}}_{k,k'}$ and $x_{k,k'}$ for each $0\leq k\leq k'\leq n-1$, $\kappa_{k,k''}^{k'}$ for each $0\leq k\leq k'\leq k''\leq n-1$, and $\kappa^{\ell}$ for each $0\leq \ell\leq n-1$.
\begin{lem}\label{lem: Hodge Ext transfer}
Let $0\leq \ell'\leq \ell\leq n-1$. Then we have
\begin{equation}\label{equ: Hodge Ext transfer dim}
\Dim_E \Hom_{G}(\mathrm{ker}(d_{\Omega^{\bullet}}^{\ell})^\vee, V_{J^{\ell}}^{\rm{an}})=1
\end{equation}
and the unique up to scalar non-zero map $\mathrm{ker}(d_{\Omega^{\bullet}}^{\ell})^\vee\rightarrow V_{J^{\ell}}^{\rm{an}}$ has image $\tld{C}_{x_{0,\ell}}^{J^{\ell}}$ (cf.~ the discussions around (\ref{equ: x cube socle}) and (\ref{equ: x filtration})). Moreover, the surjection $\mathrm{ker}(d_{\Omega^{\bullet}}^{\ell})^\vee\twoheadrightarrow \tld{C}_{x_{0,\ell}}^{J^{\ell}}$ induces an isomorphism
\begin{equation}\label{equ: Hodge Ext transfer}
\mathrm{Ext}_{G}^{\ell-\ell'}(V_{J^{\ell'}}^{\infty},\mathrm{ker}(d_{\Omega^{\bullet}}^{\ell})^\vee)\buildrel\sim\over\longrightarrow
\mathrm{Ext}_{G}^{\ell-\ell'}(V_{J^{\ell'}}^{\infty},\tld{C}_{x_{0,\ell}}^{J^{\ell}}).
\end{equation}
Moreover, the RHS of (\ref{equ: Hodge Ext transfer}) equals the subspace
\[\mathrm{Fil}_{x_{\ell',\ell}}(\mathbf{E}_{\ell',\ell})\subseteq \mathbf{E}_{\ell',\ell}\]
which admits a decreasing filtration
\[\{\mathrm{Fil}_{x_{\ell'',\ell}}(\mathbf{E}_{\ell',\ell})\}_{\ell'\leq \ell''\leq \ell}\]
with graded pieces being the $1$-dimensional $E$-vector spaces $\mathrm{gr}_{x_{\ell'',\ell}}(\mathbf{E}_{\ell',\ell})$ for each $\ell'\leq \ell''\leq \ell$. In particular, (\ref{equ: Hodge Ext transfer}) is an isomorphism between $\ell-\ell'+1$-dimensional $E$-vector spaces.
\end{lem}
\begin{proof}
If $\ell=0$, then $\mathrm{ker}(d_{\Omega^{\bullet}}^{\ell})^\vee=V_{\Delta}^{\infty}=1_{G}$ and the claims are clear.
Assume from now that $1\leq \ell\leq n-1$, then it follows from the description of the multiplicity free finite length coadmissible $D(G)$-modules $\Omega^{\ell}$ and $\Omega^{\ell}/\mathrm{ker}(d_{\Omega^{\bullet}}^{\ell})\cong \mathrm{im}(d_{\Omega^{\bullet}}^{\ell})$ in Theorem~\ref{thm: dR complex} that $\mathrm{ker}(d_{\Omega^{\bullet}}^{\ell})^\vee\in\mathrm{Rep}^{\rm{an}}_{\rm{adm}}(G)$ is finite length and multiplicity free and satisfies
\[\mathrm{JH}_{G}(\mathrm{ker}(d_{\Omega^{\bullet}}^{\ell})^\vee)=\{V_{J^{\ell}}^{\infty}\}\sqcup\{C_{w_{j,n-\ell},[1,j-1]}\mid j\in\Delta\}\]
with its partial-order generated by $V_{J^{\ell}}^{\infty}<C_{s_{n-\ell},[1,n-1-\ell]}$ and $C_{w_{j,n-\ell},[1,j-1]}<C_{w_{j+1,n-\ell},[1,j]}$ for each $1\leq j\leq n-2$. We consider the unique up to scalar surjection $q: \mathrm{ker}(d_{\Omega^{\bullet}}^{\ell})^\vee\rightarrow V$ with $V$ being the unique quotient of $\mathrm{ker}(d_{\Omega^{\bullet}}^{\ell})^\vee$ with socle $V_{J^{\ell}}^{\infty}$, which satisfies
\begin{equation}\label{equ: Hodge Ext transfer JH}
\mathrm{JH}_{G}(V)=\{V_{J^{\ell}}^{\infty}\}\sqcup\{C_{w_{j,n-\ell},[1,j-1]}\mid n-\ell\leq j\leq n-1\}
\end{equation}
with its partial-order generated by $V_{J^{\ell}}^{\infty}<C_{s_{n-\ell},[1,n-1-\ell]}$ and $C_{w_{j,n-\ell},[1,j-1]}<C_{w_{j+1,n-\ell},[1,j]}$ for each $n-\ell\leq j\leq n-2$.
Recall from \ref{it: x loc an St vanishing 2} of Lemma~\ref{lem: x loc an St vanishing} (by taking $x=1$, $I_0=\Delta$ and $I_1=J^{\ell}$ in \emph{loc.cit.}) that $V_{J^{\ell}}^{\rm{an}}$ has socle $V_{J^{\ell}}^{\infty}$, and thus the surjection $q: \mathrm{ker}(d_{\Omega^{\bullet}}^{\ell})^\vee\rightarrow V$ induces an isomorphism
\begin{equation}\label{equ: Hodge Ext transfer 1}
\Hom_{G}(V,V_{J^{\ell}}^{\rm{an}})\buildrel\sim\over\longrightarrow \Hom_{G}(\mathrm{ker}(d_{\Omega^{\bullet}}^{\ell})^\vee,V_{J^{\ell}}^{\rm{an}}).
\end{equation}
For each $W\in\mathrm{JH}_{G}(\mathrm{ker}(d_{\Omega^{\bullet}}^{\ell})^\vee)\setminus\{V_{J^{\ell}}^{\infty}\}$, it follows from \ref{it: x loc an St vanishing 1} of Lemma~\ref{lem: x loc an St vanishing} (by taking $C_{w,I}=W$, $x=1$ and $I_0=I_1=J^{\ell}$ in \emph{loc.cit.}) that
\begin{equation}\label{equ: Hodge Ext transfer vanishing}
\mathrm{Ext}_{G}^k(W,V_{J^{\ell}}^{\rm{an}})=0
\end{equation}
for $k\geq 0$. By a d\'evissage with respect to $\mathrm{JH}_{G}(V/V_{J^{\ell}}^{\infty})$ we have $\mathrm{Ext}_{G}^k(V/V_{J^{\ell}}^{\infty},V_{J^{\ell}}^{\rm{an}})=0$ for $k\geq 0$, and in particular the short exact sequence $0\rightarrow V_{J^{\ell}}^{\infty}\rightarrow V\rightarrow V/V_{J^{\ell}}^{\infty}\rightarrow 0$ induces an isomorphism
\[\Hom_{G}(V,V_{J^{\ell}}^{\rm{an}})\buildrel\sim\over\longrightarrow \Hom_{G}(V_{J^{\ell}}^{\infty},V_{J^{\ell}}^{\rm{an}})\]
which together with $\mathrm{soc}_{G}(V_{J^{\ell}}^{\rm{an}})=V_{J^{\ell}}^{\infty}$ and (\ref{equ: Hodge Ext transfer 1}) gives (\ref{equ: Hodge Ext transfer dim}).
As $\mathrm{ker}(d_{\Omega^{\bullet}}^{\ell})^\vee$ has cosocle (see the description of the partial-order on $\mathrm{JH}_{G}(\mathrm{ker}(d_{\Omega^{\bullet}}^{\ell})^\vee)$ given above)
\[C_{w_{n-1,n-\ell},[1,n-2]}=C_{x_{0,\ell}},\] 
so does $V$ and thus $V$ must be the unique subrepresentation of $V_{J^{\ell}}^{\rm{an}}$ with cosocle $C_{x_{0,\ell}}$ (using $[V_{J^{\ell}}^{\rm{an}}:C_{x_{0,\ell}}]=1$ from Lemma~\ref{lem: JH mult PS}) which is $\tld{C}_{x_{0,\ell}}$ by definition. In particular, this determines the partially-ordered set $\mathrm{JH}_{G}(\tld{C}_{x_{0,\ell}})$ explicitly (which equals (\ref{equ: Hodge Ext transfer JH})).
Note that each $W\in\mathrm{JH}_{G}(\mathrm{ker}(q))$ has the form $C_{w_{j,n-\ell},[1,j-1]}$ for some $1\leq j<n-\ell$, and that we have $J^{\ell'}=[1,n-1-\ell']\not\subseteq [1,j-1]$ for each $\ell'\leq \ell$. Hence, we deduce from \ref{it: Ext with sm 1} of Lemma~\ref{lem: Ext with sm} that
\[\mathrm{Ext}_{G}^k(V_{J^{\ell'}}^{\infty},W)=0\]
for each $W\in \mathrm{JH}_{G}(\mathrm{ker}(q))$ and $k\geq 0$. This together with a d\'evissage with respect to $\mathrm{JH}_{G}(\mathrm{ker}(q))$ and the short exact sequence $0\rightarrow \mathrm{ker}(q)\rightarrow \mathrm{ker}(d_{\Omega^{\bullet}}^{\ell})^\vee\rightarrow V\rightarrow 0$ implies that the surjection $q$ induces the isomorphism (\ref{equ: Hodge Ext transfer}).
Finally, since we have shown that $\mathrm{JH}_{G}(\tld{C}_{x_{0,\ell}})=\mathrm{JH}_{G}(V)$ as described in (\ref{equ: Hodge Ext transfer JH}), we see that $x\in\Gamma^{J^{\ell'}\setminus J^{\ell}}$ (namely $x\in\Gamma$ with $\mathrm{Supp}(x)\subseteq [n-\ell,n-1-\ell']$) satisfies $C_{x}^{J^{\ell}}\in\mathrm{JH}_{G}(\tld{C}_{x_{0,\ell}})$ if and only if $x=x_{\ell'',\ell}$ for some $\ell'\leq \ell''\leq \ell$. This together with Lemma~\ref{lem: Ext subquotient} finishes the proof of the last claim of this lemma.
\end{proof}

\begin{lem}\label{lem: dR End}
For each $0\leq \ell\leq k\leq n-1$, the maps $\sigma_{\geq \ell}(\tau^{\leq k}\Omega^{\bullet})\rightarrow \tau^{\leq k}\Omega^{\bullet}\rightarrow \Omega^{\bullet}$ induce the following commutative diagram of isomorphisms
\begin{equation}\label{equ: dR End}
\xymatrix{
\mathrm{End}_{\cM(G)}(\tau^{\leq k}\Omega^{\bullet}) \ar^{\sim}[rr] \ar^{\wr}[d] & & \Hom_{\cM(G)}(\tau^{\leq k}\Omega^{\bullet},\Omega^{\bullet}) \ar^{\wr}[d]\\
\Hom_{\cM(G)}(\sigma_{\geq \ell}(\tau^{\leq k}\Omega^{\bullet}),\tau^{\leq k}\Omega^{\bullet}) \ar^{\sim}[rr] & & \Hom_{\cM(G)}(\sigma_{\geq \ell}(\tau^{\leq k}\Omega^{\bullet}),\Omega^{\bullet})
}.
\end{equation}
\end{lem}
\begin{proof}
The splitting $s$ from (\ref{equ: abstract dR splitting}) induces a splitting
\[s_{\leq k}: \tau^{\leq k}\Omega^{\bullet}\buildrel\sim\over\longrightarrow \bigoplus_{\ell'=0}^k H^{\ell'}(\Omega^{\bullet})\cong \bigoplus_{\ell'=0}^k (V_{J^{\ell'}}^{\infty})^\vee,\]
which further induces the following isomorphism
\begin{equation}\label{equ: dR End splitting}
\Hom_{\cM(G)}(\tau^{\leq k}\Omega^{\bullet},\Omega^{\bullet})\buildrel\sim\over\longrightarrow \Hom_{\cM(G)}(\bigoplus_{\ell'=0}^k H^{\ell'}(\Omega^{\bullet})[-\ell'],\bigoplus_{\ell'=0}^{n-1}H^{\ell'}(\Omega^{\bullet})[-\ell']), f \mapsto s\circ f\circ s_{\leq k}^{-1}.
\end{equation}
As $\Hom_{\cM(G)}(H^{\ell'}(\Omega^{\bullet})[-\ell'],H^{\ell}(\Omega^{\bullet})[-\ell])=0$ for each $\ell'<\ell$, we deduce by a d\'evissage that the map $\bigoplus_{\ell'=0}^k H^{\ell'}(\Omega^{\bullet})[-\ell']\rightarrow \bigoplus_{\ell'=0}^{n-1}H^{\ell'}(\Omega^{\bullet})[-\ell']$ induces an isomorphism
\begin{equation}\label{equ: dR End truncation}
\mathrm{End}_{\cM(G)}(\bigoplus_{\ell'=0}^k H^{\ell'}(\Omega^{\bullet})[-\ell'])\buildrel\sim\over\longrightarrow \Hom_{\cM(G)}(\bigoplus_{\ell'=0}^k H^{\ell'}(\Omega^{\bullet})[-\ell'], \bigoplus_{\ell'=0}^{n-1}H^{\ell'}(\Omega^{\bullet})[-\ell'])
\end{equation}
with both sides having dimension $\frac{(k+1)(k+2)}{2}$ (using Lemma~\ref{lem: Ext between sm}).
The isomorphism (\ref{equ: dR End truncation}) together with the splitting $s$ and $s_{\leq k}$ gives the top horizontal isomorphism of (\ref{equ: dR End}).
Now we consider the maps $\mathrm{ker}(d^k_{\Omega^{\bullet}})[-k]\rightarrow \sigma_{\geq \ell}(\tau^{\leq k}\Omega^{\bullet})\rightarrow \tau^{\leq k}\Omega^{\bullet}$ which induce maps
\begin{equation}\label{equ: dR End composition}
\Hom_{\cM(G)}(\tau^{\leq k}\Omega^{\bullet},\Omega^{\bullet})\rightarrow \Hom_{\cM(G)}(\sigma_{\geq \ell}(\tau^{\leq k}\Omega^{\bullet}),\Omega^{\bullet})\rightarrow \Hom_{\cM(G)}(\mathrm{ker}(d^k_{\Omega^{\bullet}})[-k],\Omega^{\bullet}).
\end{equation}
The splitting $s$ induces an isomorphism
\[\Hom_{\cM(G)}(\mathrm{ker}(d^k_{\Omega^{\bullet}})[-k],\Omega^{\bullet})\cong \bigoplus_{\ell'=0}^{n-1}\Hom_{\cM(G)}(\mathrm{ker}(d^k_{\Omega^{\bullet}})[-k],(V_{J^{\ell'}}^{\infty})^\vee[-\ell'])\cong \bigoplus_{\ell'=0}^{n-1}\mathrm{Ext}_{G}^{k-\ell'}(V_{J^{\ell'}}^{\infty},(\mathrm{ker}(d^k_{\Omega^{\bullet}}))^\vee),\]
which together with Lemma~\ref{lem: Hodge Ext transfer} implies that
\[\Dim_E \Hom_{\cM(G)}(\mathrm{ker}(d^k_{\Omega^{\bullet}})[-k],\Omega^{\bullet})=\sum_{\ell'=0}^k \Dim_E \mathrm{Ext}_{G}^{k-\ell'}(V_{J^{\ell'}}^{\infty},(\mathrm{ker}(d^k_{\Omega^{\bullet}}))^\vee)=\frac{(k+1)(k+2)}{2}.\]
Since the first term and the third term of (\ref{equ: dR End composition}) share the same dimension, to prove that the maps in (\ref{equ: dR End composition}) are isomorphisms, it suffices to show that they are injections.
Note that we have distinguished triangles $\sigma_{\geq \ell}(\tau^{\leq k}\Omega^{\bullet})\rightarrow \tau^{\leq k}\Omega^{\bullet}\rightarrow \sigma^{\leq \ell-1}\Omega^{\bullet}\rightarrow$ and $\mathrm{ker}(d^k_{\Omega^{\bullet}})[-k]\rightarrow\sigma_{\geq \ell}(\tau^{\leq k}\Omega^{\bullet})\rightarrow \sigma^{\leq k-1}(\sigma_{\geq \ell}\Omega^{\bullet})\rightarrow $. Hence, upon using a d\'evissage and using the splitting $s$, it suffices to show that
\begin{equation}\label{equ: dR Hom vanishing}
\Hom_{\cM(G)}(\Omega^{\ell}[-\ell], H^{\ell'}(\Omega^{\bullet})[-\ell'])\cong \mathrm{Ext}_{G}^{\ell-\ell'}(V_{J^{\ell'}}^{\infty},(\Omega^{\ell})^\vee)=0
\end{equation}
for each $0\leq \ell\leq n-2$ and $0\leq \ell'\leq n-1$. The vanishing (\ref{equ: dR Hom vanishing}) is obvious if $\ell'>\ell$. When $\ell'\leq \ell$, (\ref{equ: dR Hom vanishing}) follows from Lemma~\ref{lem: dR Ext vanishing}.
The proof is thus finished.
\end{proof}

For $0\leq \ell'\leq \ell\leq n-1$, note that $x_{\ell',\ell}$ is a maximal element in $\Gamma^{J^{\ell'}\setminus J^{\ell}}$ under both the partial-order $\unlhd$ and the Bruhat order, and thus we have the following maps
\begin{equation}\label{equ: Hodge position Coxeter grade}
\mathrm{Fil}_{x_{\ell',\ell}}(\mathbf{E}_{\ell',\ell})\hookrightarrow\mathbf{E}_{\ell',\ell}\twoheadrightarrow \mathrm{gr}_{x_{\ell',\ell}}(\mathbf{E}_{\ell',\ell})
\end{equation}
whose composition is again a surjection.
Under Lemma~\ref{lem: Hodge Ext transfer}, we know that the composition
\begin{equation}\label{equ: Hodge position composition}
\mathrm{ker}(d_{\Omega^{\bullet}}^{\ell})[-\ell]\rightarrow \sigma_{\geq \ell}\Omega^{\bullet}\rightarrow \Omega^{\bullet}\buildrel s_{\ell'}\over\longrightarrow H^{\ell'}(\Omega^{\bullet})[-\ell']\cong (V_{J^{\ell'}}^{\infty})^\vee[-\ell']
\end{equation}
corresponds to an element $\theta_{\ell',\ell}$ of $\mathrm{Fil}_{x_{\ell',\ell}}(\mathbf{E}_{\ell',\ell})$.
\begin{lem}\label{lem: Hodge position}
Let $0\leq \ell'\leq\ell\leq n-1$. We have the following results.
\begin{enumerate}[label=(\roman*)]
\item \label{it: Hodge position 1} The element $\theta_{\ell',\ell}\in\mathrm{Fil}_{x_{\ell',\ell}}(\mathbf{E}_{\ell',\ell})$ has non-zero image in $\mathrm{gr}_{x_{\ell',\ell}}(\mathbf{E}_{\ell',\ell})$ under (\ref{equ: Hodge position Coxeter grade}).
\item \label{it: Hodge position 2} We have
\[\mathbf{E}_{\ell',\ell}=E\theta_{\ell',\ell}\oplus \mathbf{E}_{\ell',\ell}^{<},\]
or equivalently the element $\theta_{\ell',\ell}\in\mathrm{Fil}_{x_{\ell',\ell}}(\mathbf{E}_{\ell',\ell})$ has non-zero image under the composition of
\begin{equation}\label{equ: Hodge position surjection}
\mathrm{Fil}_{x_{\ell',\ell}}(\mathbf{E}_{\ell',\ell})\rightarrow \mathbf{E}_{\ell',\ell}\twoheadrightarrow \overline{\mathbf{E}}_{\ell',\ell}=\mathbf{E}_{\ell',\ell}/\mathbf{E}_{\ell',\ell}^{<}.
\end{equation}
\end{enumerate}
\end{lem}
\begin{proof}
If $\ell=\ell'$, there is nothing to prove. We assume in the rest of the proof that $\ell'<\ell$.
Note that we have $\mathrm{ker}(d_{\Omega^{\bullet}}^{\ell})[-\ell]=\sigma_{\geq \ell}(\tau^{\leq \ell}\Omega^{\bullet})$. By taking $k=\ell$ in (\ref{equ: dR End}), we obtain isomorphisms
\begin{equation}\label{equ: differential ker Ext}
\xymatrix{
\mathrm{End}_{\cM(G)}(\tau^{\leq \ell}\Omega^{\bullet}) \ar^{\sim}[rr] \ar^{\wr}[d] & & \Hom_{\cM(G)}(\tau^{\leq \ell}\Omega^{\bullet},\Omega^{\bullet}) \ar^{\wr}[d]\\ 
\Hom_{\cM(G)}(\mathrm{ker}(d_{\Omega^{\bullet}}^{\ell})[-\ell],\tau^{\leq \ell}\Omega^{\bullet}) \ar^{\sim}[rr] & & \Hom_{\cM(G)}(\mathrm{ker}(d_{\Omega^{\bullet}}^{\ell})[-\ell],\Omega^{\bullet})
}.
\end{equation}
For each $0\leq \ell'\leq \ell''\leq \ell$, the map
\[ N^{\ell''-\ell'}: \bigoplus_{k=0}^{\ell}(V_{J^{k}}^{\infty})^\vee[-k]\rightarrow \bigoplus_{k=0}^{\ell}(V_{J^{k}}^{\infty})^\vee[-k]\]
restricts to a map $N^{\ell''-\ell'}: (V_{J^{\ell''}}^{\infty})^\vee[-\ell'']\rightarrow (V_{J^{\ell'}}^{\infty})^\vee[-\ell']$ and thus $N^{\ell''-\ell'}\circ s_{\ell''}$ is an element of $\Hom_{\cM(G)}(\tau^{\leq \ell}\Omega^{\bullet}, (V_{J^{\ell'}}^{\infty})^\vee[-\ell'])$.
The splitting $s_{\leq \ell}$ from (\ref{equ: abstract dR splitting truncation}) induces an isomorphism
\[\mathrm{End}_{\cM(G)}(\tau^{\leq \ell}\Omega^{\bullet})\buildrel\sim\over\longrightarrow \bigoplus_{\ell'=0}^{\ell}\Hom_{\cM(G)}(\tau^{\leq \ell}\Omega^{\bullet},(V_{J^{\ell'}}^{\infty})^\vee[-\ell']),\]
and (using Lemma~\ref{lem: Ext between sm}) the set
\begin{equation}\label{equ: Hodge position basis}
\{N^{\ell''-\ell'}\circ s_{\ell''}\mid \ell'\leq \ell''\leq \ell\}
\end{equation}
forms a basis of $\Hom_{\cM(G)}(\tau^{\leq \ell}\Omega^{\bullet},(V_{J^{\ell'}}^{\infty})^\vee[-\ell'])$ for each $0\leq \ell'\leq \ell$. In particular, we have
\[s_{\ell'}\in \Hom_{\cM(G)}(\tau^{\leq \ell}\Omega^{\bullet},(V_{J^{\ell'}}^{\infty})^\vee[-\ell'])\setminus N(\Hom_{\cM(G)}(\tau^{\leq \ell}\Omega^{\bullet},(V_{J^{\ell'+1}}^{\infty})^\vee[-\ell'-1])),\]
which under the isomorphisms (\ref{equ: differential ker Ext}) become
\begin{equation}\label{equ: Hodge avoid subspace}
s_{\ell'}\in \Hom_{\cM(G)}(\mathrm{ker}(d_{\Omega^{\bullet}}^{\ell})[-\ell],(V_{J^{\ell'}}^{\infty})^\vee[-\ell'])\setminus N(\Hom_{\cM(G)}(\mathrm{ker}(d_{\Omega^{\bullet}}^{\ell})[-\ell],(V_{J^{\ell'+1}}^{\infty})^\vee[-\ell'-1])).
\end{equation}
By Lemma~\ref{lem: Hodge Ext transfer} we know that the (unique up to scalar) non-zero map $\mathrm{ker}(d_{\Omega^{\bullet}}^{\ell})^\vee\rightarrow V_{J^{\ell}}^{\rm{an}}$ induces an isomorphism
\[\Hom_{\cM(G)}(\mathrm{ker}(d_{\Omega^{\bullet}}^{\ell})[-\ell],(V_{J^{\ell'}}^{\infty})^\vee[-\ell'])\buildrel\sim\over\longrightarrow \mathrm{Fil}_{x_{\ell',\ell}}(\mathbf{E}_{\ell',\ell}),\]
and
\begin{multline}\label{equ: Hodge position N image}
N(\Hom_{\cM(G)}(\mathrm{ker}(d_{\Omega^{\bullet}}^{\ell})[-\ell],(V_{J^{\ell'+1}}^{\infty})^\vee[-\ell'-1]))\\
\buildrel\sim\over\longrightarrow
\kappa_{\ell',\ell}^{\ell'+1}(\mathrm{Fil}_{1}(\mathbf{E}_{\ell',\ell'+1}\otimes_E\mathrm{Fil}_{x_{\ell'+1,\ell}}(\mathbf{E}_{\ell'+1,\ell}))=\mathrm{Fil}_{x_{\ell'+1,\ell}}(\mathbf{E}_{\ell',\ell}),
\end{multline}
with the first isomorphism from the definition of $N$ and the last equality from Lemma~\ref{lem: cup with sm}. Consequently, the condition (\ref{equ: Hodge avoid subspace}) becomes
\begin{equation}\label{equ: Hodge generic position}
\theta_{\ell',\ell}\in \mathrm{Fil}_{x_{\ell',\ell}}(\mathbf{E}_{\ell',\ell})\setminus \mathrm{Fil}_{x_{\ell'+1,\ell}}(\mathbf{E}_{\ell',\ell}).
\end{equation}
As $x'\unlhd x_{\ell',\ell}$ if and only if either $x'=x_{\ell',\ell}$ or $x'\unlhd x_{\ell'+1,\ell}$, we see that $\mathrm{Fil}_{x_{\ell'+1,\ell}}(\mathbf{E}_{\ell',\ell})$ has codimension $1$ in $\mathrm{Fil}_{x_{\ell',\ell}}(\mathbf{E}_{\ell',\ell})$ by Theorem~\ref{thm: coxeter filtration}, and thus (\ref{equ: Hodge generic position}) is equivalent to
\begin{equation}\label{equ: Hodge position span}
\mathrm{Fil}_{x_{\ell',\ell}}(\mathbf{E}_{\ell',\ell})=E\theta_{\ell',\ell}\oplus \mathrm{Fil}_{x_{\ell'+1,\ell}}(\mathbf{E}_{\ell',\ell}),
\end{equation}
and in particular the image of $\theta_{\ell',\ell}$ under (\ref{equ: Hodge position Coxeter grade}) is non-zero. As the composition of (\ref{equ: Hodge position surjection}) is surjective by Theorem~\ref{thm: cup top grade} yet the composition of
\[\mathrm{Fil}_{x_{\ell'+1,\ell}}(\mathbf{E}_{\ell',\ell})\rightarrow \mathbf{E}_{\ell',\ell}\twoheadrightarrow \overline{\mathbf{E}}_{\ell',\ell}=\mathbf{E}_{\ell',\ell}/\mathbf{E}_{\ell',\ell}^{<}\]
is zero (see the last equality in (\ref{equ: Hodge position N image})), we deduce from (\ref{equ: Hodge position span}) that the image of $\theta_{\ell',\ell}$ under the composition of (\ref{equ: Hodge position surjection}) is non-zero.
\end{proof}

We consider the following composition
\begin{multline}\label{equ: Hodge image embedding}
\Hom_{\cM(G)}((\mathrm{St}_G^{\rm{an}})^\vee[1-n],\mathrm{ker}(d_{\Omega^{\bullet}}^{\ell})[-\ell])
\rightarrow \Hom_{\cM(G)}((\mathrm{St}_G^{\rm{an}})^\vee[1-n],\Omega^{\bullet})\\
\buildrel\sim\over\longrightarrow \bigoplus_{k=0}^{n-1}\Hom_{\cM(G)}((\mathrm{St}_G^{\rm{an}})^\vee[1-n],(V_{J^k}^{\infty})^\vee[-k])\cong \bigoplus_{k=0}^{n-1}\mathbf{E}_{k,n-1}
\end{multline}
with the first map being an embedding by (\ref{equ: Hodge seq 0}) and the second map induced from (\ref{equ: abstract dR splitting}). 
For each $0\leq \ell\leq n-1$, we recall the map $\kappa^{\ell}$ from (\ref{equ: phi N total cup}).
\begin{lem}\label{lem: Hodge image}
Let $0\leq \ell\leq n-1$. We have the following results.
\begin{enumerate}[label=(\roman*)]
\item \label{it: Hodge image 1}
Then there exists an $E$-line
\[L_{\ell}\subseteq \bigoplus_{k=0}^{\ell}\mathrm{Fil}_{x_{k,\ell}}(\mathbf{E}_{k,\ell})\subseteq \bigoplus_{k=0}^{\ell}\mathbf{E}_{k,\ell}\]
such that the composition of (\ref{equ: Hodge image embedding}) has image
\begin{equation}\label{equ: Hodge image}
\kappa^{\ell}(L_{\ell}\otimes_E\mathbf{E}_{\ell,n-1})\subseteq \bigoplus_{k=0}^{\ell}\mathbf{E}_{k,n-1}\subseteq \bigoplus_{k=0}^{n-1}\mathbf{E}_{k,n-1}.
\end{equation}
\item \label{it: Hodge image 2} For each $0\leq k\leq \ell$, the image of $L_{\ell}$ under the projection $\mathrm{Fil}_{x_{k,\ell}}(\mathbf{E}_{x_{k,\ell}})\twoheadrightarrow \mathrm{gr}_{x_{k,\ell}}(\mathbf{E}_{x_{k,\ell}})$ and the image of $L_{\ell}$ under the composition of
    \[\mathrm{Fil}_{x_{k,\ell}}(\mathbf{E}_{x_{k,\ell}})\hookrightarrow\mathbf{E}_{k,\ell}\twoheadrightarrow \overline{\mathbf{E}}_{k,\ell}=\mathbf{E}_{k,\ell}/\mathbf{E}_{k,\ell}^{<}\]
    are non-zero.
\end{enumerate}
\end{lem}
\begin{proof}
Using Lemma~\ref{lem: Hodge Ext transfer}, we define $\theta_{\ell}$ as the element in
\begin{multline}\label{equ: Hodge image line}
\Hom_{\cM(G)}(\mathrm{ker}(d_{\Omega^{\bullet}}^{\ell})[-\ell],\Omega^{\bullet})\buildrel\sim\over\longrightarrow \bigoplus_{k=0}^{n-1}\Hom_{\cM(G)}(\mathrm{ker}(d_{\Omega^{\bullet}}^{\ell})[-\ell], H^k(\Omega^{\bullet})[-k])\\
\cong \bigoplus_{k=0}^{n-1}\mathrm{Ext}_{G}^{\ell-k}(V_{J^k}^{\infty},\mathrm{ker}(d_{\Omega^{\bullet}}^{\ell})^\vee)
=\bigoplus_{k=0}^{\ell}\mathrm{Ext}_{G}^{\ell-k}(V_{J^k}^{\infty},\mathrm{ker}(d_{\Omega^{\bullet}}^{\ell})^\vee)\hookrightarrow \bigoplus_{k=0}^{\ell}\mathbf{E}_{k,\ell}
\end{multline}
given the composition of $\mathrm{ker}(d_{\Omega^{\bullet}}^{\ell})[-\ell]\rightarrow \sigma_{\geq \ell}\Omega^{\bullet}\rightarrow \Omega^{\bullet}$, and then set $L_{\ell}\defeq E\theta_{\ell}$.
Hence, the image of the composition of (\ref{equ: Hodge image embedding}) equals the image of the cup product map
\begin{multline}\label{equ: Hodge image cup}
\Hom_{\cM(G)}((\mathrm{St}_G^{\rm{an}})^\vee[1-n],\mathrm{ker}(d_{\Omega^{\bullet}}^{\ell})[-\ell])\otimes_E L_{\ell}\\ \buildrel\cup\over\longrightarrow \bigoplus_{k=0}^{\ell}\Hom_{\cM(G)}((\mathrm{St}_G^{\rm{an}})^\vee[1-n],(V_{J^k})^\vee[-k])
=\bigoplus_{k=0}^{\ell}\mathbf{E}_{k,n-1}.
\end{multline}
For $0\leq k\leq \ell$, we define $\theta_{k,\ell}\in\mathrm{Fil}_{x_{k,\ell}}(\mathbf{E}_{k,\ell})\subseteq \mathbf{E}_{k,\ell}$ as the projection of $\theta_{\ell}$ to $\mathbf{E}_{k,\ell}$, and then set $L_{k,\ell}\defeq E\theta_{k,\ell}\subseteq \mathrm{Fil}_{x_{k,\ell}}(\mathbf{E}_{k,\ell})$.
It is clear that
\[\theta_{\ell}=\sum_{k=0}^{\ell}\theta_{k,\ell}\]
and $\theta_{k,\ell}$ recovers the element (\ref{equ: Hodge position composition}) with $k$ replacing $\ell'$ in \emph{loc.cit.}.
Recall from Lemma~\ref{lem: Hodge Ext transfer} that $\Dim_E\Hom_G(\mathrm{ker}(d_{\Omega^{\bullet}}^{\ell})^\vee, V_{J^{\ell}}^{\rm{an}})=1$. Note that $\mathrm{Ext}_G^k(W,\mathrm{St}_{G}^{\rm{an}})=0$ for $k\geq 0$ and each \[W\in(\mathrm{JH}_{G}(V_{J^{\ell}}^{\rm{an}})\cup\mathrm{JH}_{G}(\mathrm{ker}(d_{\Omega^{\bullet}}^{\ell})^\vee))\setminus\{V_{J^{\ell}}^{\infty}\}\] upon applying \ref{it: x loc an St vanishing 1} of Lemma~\ref{lem: x loc an St vanishing} by taking $C_{w,I}=W$ (with $w\neq 1$), $x=1$, $I_0=\emptyset$ and $I_1=\Delta$ in \emph{loc.cit.}. Hence, a d\'evissage shows that the unique up to scalar non-zero map $\mathrm{ker}(d_{\Omega^{\bullet}}^{\ell})^\vee\rightarrow V_{J^{\ell}}^{\rm{an}}$ induces an isomorphism
\[\mathbf{E}_{\ell,n-1}=\mathrm{Ext}_G^{n-1-\ell}(V_{J^{\ell}}^{\rm{an}},\mathrm{St}_G^{\rm{an}})\buildrel\sim\over\longrightarrow \mathrm{Ext}_G^{n-1-\ell}(\mathrm{ker}(d_{\Omega^{\bullet}}^{\ell})^\vee,\mathrm{St}_G^{\rm{an}}),\]
and thus the cup product map (\ref{equ: Hodge image cup}) factors through
\begin{multline*}
\Hom_{\cM(G)}((\mathrm{St}_G^{\rm{an}})^\vee[1-n],\mathrm{ker}(d_{\Omega^{\bullet}}^{\ell})[-\ell])\otimes_E L_{\ell}\\ 
\cong \Hom_{\cM(G)}((\mathrm{St}_G^{\rm{an}})^\vee[1-n],(V_{J^{\ell}}^{\rm{an}})^{\vee}[-\ell])\otimes_E\Hom_{\cM(G)}((V_{J^{\ell}}^{\rm{an}})^{\vee}[-\ell],\mathrm{ker}(d_{\Omega^{\bullet}}^{\ell})[-\ell])\otimes_E L_{\ell}\\
\buildrel\sim\over\longrightarrow \mathbf{E}_{\ell,n-1}\otimes_E L_{\ell}
\buildrel\kappa^{\ell}\over\longrightarrow \bigoplus_{k=0}^{\ell}\mathbf{E}_{k,n-1}
\end{multline*}
where we abuse $L_{\ell}$ for its image in $\bigoplus_{k=0}^{\ell}\mathbf{E}_{k,\ell}$ under (\ref{equ: Hodge image line}).
This finishes the proof of \ref{it: Hodge image 1}.
Now that $L_{k,\ell}=E\theta_{k,\ell}$ is the projection of $L_{\ell}=E\theta_{\ell}$ to $\mathbf{E}_{k,\ell}$, \ref{it: Hodge image 2} follows from Lemma~\ref{lem: Hodge position} (upon replacing $\ell'$ in \emph{loc.cit.} with $k$ here).
\end{proof}

Recall from Lemma~\ref{lem: Hodge Ext transfer} (with $\ell=n-1$ in \emph{loc.cit.}) that we have a unique (up to scalars) embedding
\begin{equation}\label{equ: St Fil max}
(\Omega^{n-1})^{\vee}\hookrightarrow \mathrm{St}_{G}^{\rm{an}}
\end{equation}
under which we may identify $(\Omega^{n-1})^{\vee}$ with the subrepresentation $\tld{C}_{x_{1,n-1}}^{\emptyset}\subseteq V_{\emptyset}^{\rm{an}}=\mathrm{St}_{G}^{\rm{an}}$.
\begin{lem}\label{lem: Fil max image}
We have the following results.
\begin{enumerate}[label=(\roman*)]
\item \label{it: Fil max image 1} For each $0\leq k\leq n-1$, the embedding (\ref{equ: St Fil max}) induces an embedding
\begin{equation}\label{equ: St Fil max sub}
\mathrm{Ext}_{G}^{n-1-k}(V_{J^{k}}^{\infty},(\Omega^{n-1})^{\vee})\hookrightarrow \mathrm{Ext}_{G}^{n-1-k}(V_{J^{k}}^{\infty},\mathrm{St}_{G}^{\rm{an}})=\mathbf{E}_{k,n-1}
\end{equation}
whose image is $\mathrm{Fil}_{x_{k,n-1}}(\mathbf{E}_{k,n-1})$.
\item \label{it: Fil max image 2} The embedding (\ref{equ: St Fil max}) induces a strict embedding between filtered $(\varphi,N)$-modules
\begin{equation}\label{equ: Fil max phi N embedding}
D^{\flat}\defeq \Hom_{\cM(G)}(\Omega^{n-1}[1-n],\Omega^{\bullet})\hookrightarrow D\defeq \Hom_{\cM(G)}((\mathrm{St}_{G}^{\rm{an}})^{\vee}[1-n],\Omega^{\bullet}).
\end{equation}
Moreover, $D^{\flat}$ is the minimal $(\varphi,N)$-module that contains $\mathrm{Fil}_{H}^{n-1}(D)$.
\end{enumerate}
\end{lem}
\begin{proof}
\ref{it: Fil max image 1} follows directly from Lemma~\ref{lem: Hodge Ext transfer} with $\ell=n-1$ and $\ell'=k$ in \emph{loc.cit.}.

We prove \ref{it: Fil max image 2}.\\
Using the abstract splitting (\ref{equ: abstract dR splitting}), we see that the embedding (\ref{equ: St Fil max}) induces an embedding of the form (\ref{equ: Fil max phi N embedding}) which corresponds to
\begin{equation}\label{equ: Fil max phi N embedding split}
\bigoplus_{k=0}^{n-1}\mathrm{Fil}_{x_{k,n-1}}(\mathbf{E}_{k,n-1})\hookrightarrow \bigoplus_{k=0}^{n-1}\mathbf{E}_{k,n-1}.
\end{equation}
under (\ref{equ: abstract dR splitting}). In particular, the underlying $(\varphi,N)$-module of $D^{\flat}$ is the minimal $(\varphi,N)$-module of $D$ that contains $\mathrm{Fil}_{H}^{n-1}(D)$ by \ref{it: Fil max sub phi N 1} of Lemma~\ref{lem: Fil max sub phi N}. Now that the embedding (\ref{equ: St Fil max}) clearly induces an isomorphism
\begin{multline*}
\mathrm{Fil}_{H}^{n-1}(D^{\flat})=\Hom_{\cM(G)}(\Omega^{n-1}[1-n],\Omega^{n-1}[1-n])\\
\buildrel\sim\over\longrightarrow \Hom_{\cM(G)}((\mathrm{St}_{G}^{\rm{an}})^{\vee}[1-n],\Omega^{n-1}[1-n])=\mathrm{Fil}_{H}^{n-1}(D),
\end{multline*}
we conclude from \ref{it: Fil max sub phi N 2} of Lemma~\ref{lem: Fil max sub phi N} that Hodge filtration of $D^{\flat}$ must coincide with the submodule filtration induced from the Hodge filtration on $D$, which satisfies $\mathrm{Fil}_{H}^{0}(D^{\flat})=D^{\flat}$, $\mathrm{Fil}_{H}^{n}(D^{\flat})=0$ and $\mathrm{Fil}_{H}^{\ell}(D^{\flat})=\mathrm{Fil}_{H}^{n-1}(D)$ for each $1\leq \ell \leq n-1$. In particular, the embedding (\ref{equ: Fil max phi N embedding}) must be strict.
\end{proof}

\begin{thm}\label{thm: main filtered phi N}
Let $n\geq 2$ with $G=\mathrm{PGL}_{n}(K)$. Then we have the following commutative diagram of maps (up to scalars) between filtered $(\varphi,N)$-modules
\begin{equation}\label{equ: main filtered phi N}
\xymatrix{
\Hom_{\cM(G)}(\Omega^{n-1}[1-n], \Omega^\bullet) \ar^{\sim}[r] \ar@{^{(}->}[d] & \mathrm{D}_{\rm{dR}}(\cD_{n-1}^{\flat}\otimes_{\cR}\cR(\theta^{1-n})) \ar@{^{(}->}[d]\\
\Hom_{\cM(G)}((\mathrm{St}_G^{\rm{an}})^\vee[1-n], \Omega^\bullet) \ar^{\sim}[r] & \mathrm{D}_{\rm{dR}}(\cD_{n-1}\otimes_{\cR}\cR(\theta^{1-n}))
}
\end{equation}
with both the horizontal maps being isomorphisms, the LHS vertical map induced from the embedding $(\Omega^{n-1})^{\vee}\hookrightarrow \mathrm{St}_{G}^{\rm{an}}$ (see (\ref{equ: St Fil max})), and the RHS vertical map induced from the embedding $\cD_{n-1}^{\flat}\hookrightarrow \cD_{n-1}$ (see \ref{it: Galois sub flat 1} of Lemma~\ref{lem: Galois sub flat}).
\end{thm}
\begin{proof}
Recall from (\ref{equ: Hodge seq 5}), Lemma~\ref{lem: Hodge image} and Definition~\ref{def: universal phi N} that $D\defeq \Hom_{\cM(G)}((\mathrm{St}_G^{\rm{an}})^\vee[1-n], \Omega^\bullet)$ is universal special with level $n-1$.
We then deduce from \ref{it: Fil max image 2} of Lemma~\ref{lem: Fil max image} that $D^{\flat}\defeq \Hom_{\cM(G)}(\Omega^{n-1}[1-n], \Omega^\bullet)$ coincides with the minimal $(\varphi,N)$-submodule of $D$ containing $\mathrm{Fil}_{H}^{n-1}(D)$ endowed with the induced filtration from $D$. These together with Proposition~\ref{prop: Galois dR Ext exact} and the isomorphism $S_{n-1,n-1}\cong D^{\flat}$ between filtered $(\varphi,N)$-modules (see (\ref{equ: Fil max basis})) give the commutative diagram (\ref{equ: main filtered phi N}).
\end{proof}

\section{Higher $\mathscr{L}$-invariants}\label{sec: total}
In this section, we introduce a definition of (higher) $\mathscr{L}$-invariants for $G=\mathrm{PGL}_{n}(K)$ and propose a candidate for its explicit relation to the Fontaine-Mazur $\mathscr{L}$-invariants (see Question~\ref{ques: BS FM match}).
\subsection{Breuil-Schraen $\mathscr{L}$-invariant}\label{subsec: BS inv}
We define what we call \emph{Breuil-Schraen $\mathscr{L}$-invariants} and study its moduli space (see Theorem~\ref{thm: moduli of inv}).

For each $I\subseteq \Delta$, we consider the following element (see (\ref{equ: explicit cup generator}))
\[x_{I,I}=\cup_{\al\in I}x_{\al}^{\infty} \in \mathbf{E}_{I}\]
and define $\mathbf{E}_{I}^{\infty}\subseteq \mathbf{E}_{I}$ as its $E$-span. By successively applying Lemma~\ref{lem: cup with sm} and using \ref{it: tau and coxeter 1} of Lemma~\ref{lem: tau and coxeter}, we deduce that
\begin{equation}\label{equ: sm line}
\mathrm{Fil}_{1}(\mathbf{E}_{I})=\mathbf{E}_{I}^{\infty}=\tau^{I}(\mathbf{E}_{I}).
\end{equation}
For each $I_0\subseteq\Delta$, we define
\begin{equation}\label{equ: E sm saturation}
\tld{\mathbf{E}}_{I_0}\defeq \kappa_{\Delta\setminus I_0,I_0}(\mathbf{E}_{\Delta\setminus I_0}^{\infty}\otimes_E \mathbf{E}_{I_0})\subseteq \mathbf{E}_{\Delta}
\end{equation}
and note that $\tld{\mathbf{E}}_{I_0}=\tau^{\Delta\setminus I_0}(\mathbf{E}_{\Delta})$ admits a basis of the form (see \ref{it: tau and coxeter 1} of Lemma~\ref{lem: tau and coxeter} and Remark~\ref{rem: tau cup basis})
\[\{x_{S,I}\mid S\in \cS_{\Delta}, \Delta\setminus I_0\subseteq I\}.\]
It is clear that 
\[\{\tld{\mathbf{E}}_{I_0}=\tau^{\Delta\setminus I_0}(\mathbf{E}_{\Delta})\}_{I_0\subseteq\Delta}\] 
forms an increasing filtration on $\mathbf{E}_{\Delta}$.
The map
\[W_{I_0} \mapsto \kappa_{\Delta\setminus I_0,I_0}(\mathbf{E}_{\Delta\setminus I_0}^{\infty}\otimes_{E}W_{I_0})\]
defines a natural bijection between the set of $E$-subspaces of $\mathbf{E}_{I_0}$ and that of $\tld{\mathbf{E}}_{I_0}$.
The formula $\kappa_{\Delta\setminus I_0,I_0}(x_{\Delta\setminus I_0,\Delta\setminus I_0}\otimes-)$ defines an isomorphism
\begin{equation}\label{equ: total cup with sm}
\iota_{I_0}:\mathbf{E}_{I_0}\buildrel\sim\over\longrightarrow\tld{\mathbf{E}}_{I_0}, x\mapsto \kappa_{\Delta\setminus I_0,I_0}(x_{\Delta\setminus I_0,\Delta\setminus I_0}\otimes x).
\end{equation}

We write $\mathbf{E}_{n-1}\defeq \mathbf{E}_{\Delta}$ and
\begin{equation}\label{equ: pure sm cup}
\mathbf{E}_{n-1}^{\infty}\defeq \tld{\mathbf{E}}_{\emptyset}=\mathbf{E}_{\Delta}^{\infty}
\end{equation}
for short.
Let $W\subseteq \mathbf{E}_{n-1}$ be a subspace.
Then for each $I_0\subseteq \Delta$, there exists a unique $E$-subspace $W_{I_0}\subseteq \mathbf{E}_{I_0}$ such that
\[W\cap\tld{\mathbf{E}}_{I_0}=\kappa_{\Delta\setminus I_0,I_0}(\mathbf{E}_{\Delta\setminus I_0}^{\infty}\otimes_E W_{I_0})\]
or equivalently $W_{I_0}=\iota_{I_0}^{-1}(W\cap\tld{\mathbf{E}}_{I_0})$.
If $W$ has codimension $1$ in $\mathbf{E}_{n-1}$ and $W\cap \mathbf{E}_{n-1}^{\infty}=0$ (see (\ref{equ: pure sm cup}) for $\mathbf{E}_{n-1}^{\infty}$), then $W\cap\tld{\mathbf{E}}_{I_0}$ has codimension $1$ in $\tld{\mathbf{E}}_{I_0}$ and thus $W_{I_0}$ has codimension $1$ in $\mathbf{E}_{I_0}$, with $W_{I_0}\cap \mathbf{E}_{I_0}^{\infty}=0$ for each $I_0\subseteq \Delta$.

Now we are ready for the following definition.
\begin{defn}\label{def: total inv}
A codimension $1$ $E$-subspace $W\subseteq \mathbf{E}_{n-1}$ is called a \emph{Breuil-Schraen $\mathscr{L}$-invariant} if the following conditions hold
\begin{enumerate}[label=(\roman*)]
\item \label{it: total 1} $W\cap \mathbf{E}_{n-1}^{\infty}=0$ and thus $W_{I_0}$ is a codimension $1$ $E$-subspace of $\mathbf{E}_{I_0}$ satisfying $W_{I_0}\cap \mathbf{E}_{I_0}^{\infty}=0$ for each $I_0\subseteq \Delta$;
\item \label{it: total 2} the composition
\begin{equation}\label{equ: total inv cup composition}
\mathbf{E}_{J}\otimes_E \mathbf{E}_{J'}\buildrel\kappa_{J,J'}\over\longrightarrow \mathbf{E}_{J\sqcup J'}\twoheadrightarrow \mathbf{E}_{J\sqcup J'}/W_{J\sqcup J'}
\end{equation}
factors through an isomorphism between $1$-dimensional $E$-vector spaces \[(\mathbf{E}_{J}/W_{J})\otimes_E(\mathbf{E}_{J'}/W_{J'})\buildrel\sim\over\longrightarrow\mathbf{E}_{J\sqcup J'}/W_{J\sqcup J'}\] for each pair of disjoint subsets $J,J'\subseteq \Delta$.
\end{enumerate}
\end{defn}

\begin{rem}\label{rem: auto L inv}
Our definition of Breuil-Schraen $\mathscr{L}$-invariants in Definition~\ref{def: total inv} is directly motivated by Gehrmann's work on \emph{automorphic $\mathscr{L}$-invariants} in \cite[\S 3.3]{Geh21}. See \S \ref{subsubsec: intro auto L inv} for more details.
\end{rem}

\begin{lem}\label{lem: total inv cond transfer}
Let $W\subseteq \mathbf{E}_{n-1}$ be a codimension $1$ $E$-subspace that satisfies \ref{it: total 1} of Definition~\ref{def: total inv}. We have the following results.
\begin{enumerate}[label=(\roman*)]
\item \label{it: total inv cond 1} The inclusion $\mathbf{E}_{I_0}^{\infty}\subseteq \mathbf{E}_{I_0}$ induces an isomorphism between $1$-dimensional $E$-vector spaces
    \[\mathbf{E}_{I_0}^{\infty}\buildrel\sim\over\longrightarrow \mathbf{E}_{I_0}/W_{I_0}\]
    for each $I_0\subseteq \Delta$.
\item \label{it: total inv cond 2} Assume that
\begin{equation}\label{equ: total inv cond inclusion}
\kappa_{J,J'}(W_{J}\otimes_E\mathbf{E}_{J'}+\mathbf{E}_{J}\otimes_E W_{J'})\subseteq W_{J\sqcup J'}
\end{equation}
for each pair of disjoint subsets $J,J'\subseteq \Delta$. Then \ref{it: total 2} of Definition~\ref{def: total inv} holds.
\end{enumerate}
\end{lem}
\begin{proof}
\ref{it: total inv cond 1} is clear from \ref{it: total 1} of Definition~\ref{def: total inv}.

We prove \ref{it: total inv cond 2}.\\
It follows from (\ref{equ: total inv cond inclusion}) that (\ref{equ: total inv cup composition}) factors through
\begin{equation}\label{equ: total inv cond factor}
(\mathbf{E}_{J}/W_{J})\otimes_E(\mathbf{E}_{J'}/W_{J'})\rightarrow\mathbf{E}_{J\sqcup J'}/W_{J\sqcup J'}.
\end{equation}
Now that the cup product map
\[\kappa_{J,J'}: \mathbf{E}_{J}\otimes_E \mathbf{E}_{J'}\rightarrow \mathbf{E}_{J\sqcup J'}\]
induces the following isomorphism between $1$-dimensional $E$-vector spaces
\[\mathbf{E}_{J}^{\infty}\otimes_E \mathbf{E}_{J'}^{\infty}\buildrel\sim\over\longrightarrow \mathbf{E}_{J\sqcup J'}^{\infty},\]
we conclude using \ref{it: total inv cond 1} that (\ref{equ: total inv cond factor}) is also an isomorphism between $1$-dimensional $E$-vector spaces (for each pair of disjoint subsets $J,J'\subseteq \Delta$), which gives \ref{it: total 2} of Definition~\ref{def: total inv}.
\end{proof}

Note that codimension $1$ $E$-subspaces $W\subseteq \mathbf{E}_{n-1}$ satisfying Definition~\ref{def: total inv} are parameterized by a Zariski locally closed subvariety $\cB\cS\subseteq \mathbb{P}(\mathbf{E}_{n-1})$.

For each $\al\in\Phi^+$, we write
\[x_{\al}^{\infty}\defeq x_{I_{\al},I_{\al}}=\cup_{\beta\in I_{\al}}x_{\beta}^{\infty}\in\mathbf{E}_{I_{\al}}^{\infty}\]
for short.
Let $W\subseteq \mathbf{E}_{n-1}$ be a Breuil-Schraen $\mathscr{L}$-invariant in the sense of Definition~\ref{def: total inv}. For each $\al\in\Phi^+$, it follows from \ref{it: total inv cond 1} of Lemma~\ref{lem: total inv cond transfer} that the inclusion $Ex_{\al}^{\infty}=\mathbf{E}_{I_{\al}}^{\infty}\subseteq \mathbf{E}_{I_{\al}}$ induces an isomorphism
\[Ex_{\al}^{\infty}\buildrel\sim\over\longrightarrow \mathbf{E}_{I_{\al}}/W_{I_{\al}}.\]
In particular, there exists a unique $\mathscr{L}_{\al}\in E$ such that the image of $x_{\al}$ in $\mathbf{E}_{I_{\al}}/W_{I_{\al}}$ equals that of $\mathscr{L}_{\al}x_{\al}^{\infty}$, or equivalently $x_{\al}-\mathscr{L}_{\al}x_{\al}^{\infty}\in W_{I_{\al}}$.
We have thus defined a map
\begin{equation}\label{equ: abstract L inv}
\cB\cS(E)\cong \fu^+(E):~W\mapsto (\mathscr{L}_{\al})_{\al\in\Phi^+},
\end{equation}
where $U^+$ is the unipotent radical of $B^+$ and $\fu^+$ is its associated $E$-Lie algebra.

Recall that $U^+$ is the unipotent radical of $B^+$ with $\fu^+$ being its associated $E$-Lie algebra.
\begin{thm}\label{thm: moduli of inv}
The map (\ref{equ: abstract L inv}) is a bijection.
\end{thm}
\begin{proof}
It suffices to construct the inverse of (\ref{equ: abstract L inv}).
Let $(\mathscr{L}_{\al})_{\al\in\Phi^+}\in \fu^+(E)$.
For each $I_0\subseteq \Delta$, we write
\begin{equation}\label{equ: total basis index}
\Sigma_{I_0}\defeq \{(S,I)\mid S\in \cS_{I_0}, I\subseteq S\cap\Delta\}
\end{equation}
for short and recall from Proposition~\ref{prop: total cup basis} that $\{x_{S,I}\}_{(S,I)\in \Sigma_{I_0}}$ forms a basis of $\mathbf{E}_{I_0}$.
The set $\Sigma_{I_0}$ admits a natural partial order that $(S,I)\leq (S',I')$ if and only if each element of $S'$ is a sum of elements in $S$ and that $I'\subseteq I$.
We write
\[\tld{\Sigma}_{I_0}\defeq \{(S,I)\mid S\in \cS_{\Delta}, I\cup I_0=\Delta\}\subseteq \Sigma_{\Delta}\]
and note that we have a natural bijection
\[\Sigma_{I_0}\buildrel\sim\over\longrightarrow \tld{\Sigma}_{I_0}: (S,I)\mapsto (S\sqcup(\Delta\setminus I_0),I\sqcup(\Delta\setminus I_0))\]
which is compatible with (\ref{equ: total cup with sm}) (by considering the induced bijection between bases).
For each $(S,I)\in \Sigma_{I_0}$, we define
\[y_{S,I}\defeq \big(\cup_{\beta\in I}x_{\beta}^{\infty}\big)\cup\big(\cup_{\beta\in S\setminus I}(x_{\beta}^{\infty}-\mathscr{L}_{\beta}x_{\beta}^{\infty})\big)\in \mathbf{E}_{I_0}\]
and $y_{S,I}=x_{S,I}+\ast$ with $\ast$ being a linear combination of $x_{S',I'}$ for various $(S',I')<(S,I)$. Consequently, $\{y_{S,I}\}_{(S,I)\in\Sigma_{I_0}}$ differs from $\{x_{S,I}\}_{(S,I)\in\Sigma_{I_0}}$ by a unipotent matrix and thus gives another basis of $\mathbf{E}_{I_0}$.
In particular,
\begin{equation}\label{equ: pair BS codim 1}
W_{I_0}\defeq \bigoplus_{(S,I)\in \Sigma_{I_0}\setminus\{(I_0,I_0)\}}Ey_{S,I}
\end{equation}
is a codimension $1$ $E$-subspace of $\mathbf{E}_{I_0}$ that satisfies $\mathbf{E}_{I_0}=W_{I_0}+Ex_{\al}^{\infty}$.
Now we check that $W\defeq W_{\Delta}$ is a Breuil-Schraen $\mathscr{L}$-invariant in the sense of Definition~\ref{def: total inv}.
Similar argument as above shows that $\{y_{S,I}\}_{(S,I)\in\tld{\Sigma}_{I_0}}$ is a basis of $\tld{\mathbf{E}}_{I_0}$ and $\{y_{S,I}\}_{(S,I)\in\tld{\Sigma}_{I_0}\setminus\{(\Delta,\Delta)\}}$ is a basis of $\iota_{I_0}(W_{I_0})$, which together with $\{y_{S,I}\}_{(S,I)\in\Sigma_{\Delta}}$ being a basis of $\mathbf{E}_{n-1}$ and $\{y_{S,I}\}_{(S,I)\in\Sigma_{\Delta}\setminus\{(\Delta,\Delta)\}}$ being a basis of $W$ implies that $\iota_{I_0}(W_{I_0})=W\cap \tld{\mathbf{E}}_{I_0}$ and thus (\ref{equ: pair BS codim 1}) is compatible with our definition of $W_{I_0}$ before Definition~\ref{def: total inv} (and \ref{it: total 1} of Definition~\ref{def: total inv} clearly holds). To check \ref{it: total 2} of Definition~\ref{def: total inv}, it suffices to prove (\ref{equ: total inv cond inclusion}) thanks to \ref{it: total inv cond 2} of Lemma~\ref{lem: total inv cond transfer}.
But then (\ref{equ: total inv cond inclusion}) follows directly from the fact that $\{y_{S,I}\}_{(S,I)\in\Sigma_{J}}$ (resp.~$\{y_{S,I}\}_{(S,I)\in\Sigma_{J}\setminus\{(J,J)\}}$) is a basis of $\mathbf{E}_{J}$ (resp.~of $W_{J}$) as well as similar facts with $J'$ and $J\sqcup J'$ replacing $J$.
The proof is thus finished.
\end{proof}

\begin{rem}\label{rem: relation to FM L inv}
Our symbol $(\mathscr{L}_{\al})_{\al\in\Phi^+}$ in (\ref{equ: abstract L inv}) might suggest that we should be able to relate Breuil-Schraen $\mathscr{L}$-invariants to the Fontaine-Mazur $\mathscr{L}$-invariants defined for semi-stable $n$-dimensional Galois representations whose associated $(\varphi,N)$-module satisfies $N^{n-1}\neq 0$. Understanding such a relation is essentially equivalent to specifying a particular choice of $x_{\al}$ for each $\al\in\Phi^+$, which we do not investigate here in Theorem~\ref{thm: moduli of inv}. We will return to this point in \S \ref{subsec: map candidate} (see the discussion before Question~\ref{ques: BS FM match}).
\end{rem}

\subsection{Isomorphism between $\mathbf{E}_{n}$ and $\cE_{n}$}\label{subsec: map candidate}
In this section, we construct an isomorphism
\begin{equation}\label{equ: toy regulator}
r_{n}: \mathbf{E}_{n}\buildrel\sim\over\longrightarrow \cE_{n}
\end{equation}
for each $n\geq 1$ under Hypothesis~\ref{hypo: log projection}.

Recall from (\ref{equ: auto total filtration}) (upon replacing $n-1$ in \emph{loc.cit.} with $n$) that $\mathbf{E}_{n}$ admits a canonical decreasing filtration $\mathrm{Fil}_{W}^k(\mathbf{E}_{n})$ that satisfies
\begin{equation}\label{equ: auto Ext filtration grade 1}
\mathrm{gr}_{W}^{n-1}(\mathbf{E}_{n})=\mathbf{E}_{n-1}\otimes_E \mathbf{E}_{1}\cong \mathbf{E}_{n-1}\otimes_E \Hom(K^\times,E)
\end{equation}
and
\begin{equation}\label{equ: auto Ext filtration grade 2}
\mathrm{gr}_{W}^{k}(\mathbf{E}_{n})=\mathbf{E}_{k}\otimes_E\overline{\mathbf{E}}_{n-k} \cong \mathbf{E}_{k-1}\otimes_E P^{2(n-k)-1}(\fg_{n-k}).
\end{equation}
for each $0\leq k<n-1$.
Upon replacing $n-1$ with $n$ in Corollary~\ref{cor: sym sharp decomposition}, we have a unique surjection
\begin{equation}\label{equ: auto pass to sym}
\mathbf{E}_{n}\twoheadrightarrow\mathrm{Sym}^{n}(\mathbf{E}_{1})
\end{equation}
whose kernel is $\mathbf{E}_{n}^{\sharp}$ and admits the natural map $\mathrm{Sym}^{n}(\mathbf{E}_{1})\hookrightarrow \mathbf{E}_{n}$ as a section.

Recall from the discussion below Proposition~\ref{prop: Galois dR Ext exact} that $\cE_{n}$ admits a canonical decreasing filtration $\mathrm{Fil}_{W}^{\bullet}(\cE_{n})$ with graded piece
\[
\mathrm{gr}_{W}^k(\cE_{n})=\cE_{k}\otimes_E \mathrm{Ext}_g^1(\cR(\theta^{k}),\cR(\theta^{n}))^\vee
\]
for each $0\leq k\leq n-1$.
Recall from (\ref{equ: universal to sym section}) we have a unique surjection
\begin{equation}\label{equ: Galois pass to sym}
\cE_{n}\twoheadrightarrow\mathrm{Sym}^{n}(\cE_{1})
\end{equation}
whose kernel is $\cE_{n}^{\sharp}$ and admits the natural map $\mathrm{Sym}^{n}(\cE_{1})\hookrightarrow \cE_{n}$ as a section.

Recall from \cite[Lem.~1.19]{Ding17} that we have canonical isomorphisms
\begin{equation}\label{equ: log regulator}
r_{1}: \mathbf{E}_{1}=\mathrm{gr}_{W}^{0}(\mathbf{E}_{1})=\Hom(K^\times,E)\buildrel\sim\over\longrightarrow \mathrm{Ext}_{g}^1(\cR(1_{K^\times}),\cR(\delta))^{\vee}=\mathrm{gr}_{W}^{0}(\cE_{1})=\cE_{1},
\end{equation}
which induce an isomorphism
\begin{equation}\label{equ: log regulator sym}
\mathrm{Sym}^{n}(r_{1}): \mathrm{Sym}^{n}(\mathbf{E}_{1})\buildrel\sim\over\longrightarrow \mathrm{Sym}^{n}(\cE_{1})
\end{equation}
for each $n\geq 1$.

We recall the $E$-subspaces $\mathbf{E}_{n}^{\sharp},\mathbf{E}_{n}^{\flat}\subseteq \mathbf{E}_{n}$ from Proposition~\ref{prop: sharp flat intersection} (upon replacing the index $n-1$ in \emph{loc.cit.} with $n$, see also (\ref{equ: sharp subspace full}) and (\ref{equ: flat subspace})).
\begin{thm}\label{thm: regulator candidate}
Assume that Hypothesis~\ref{hypo: log projection} holds for each $n\geq 2$.
Then there exists at most one collection of isomorphisms $\{r_{n}: \mathbf{E}_{n}\rightarrow \cE_{n}\}_{n\geq 1}$ such that $r_{1}$ is given by (\ref{equ: log regulator}) and the following conditions hold.
\begin{enumerate}[label=(\roman*)]
\item \label{it: regulator 1} For each $n\geq 2$, we have the following commutative diagram
\begin{equation}\label{equ: pass to sym diagram}
\xymatrix{
\mathbf{E}_{n} \ar^{r_{n}}[rr] \ar[d] & & \cE_{n} \ar[d]\\
\mathrm{Sym}^{n}(\mathbf{E}_{1}) \ar^{\mathrm{Sym}^{n}(r_{1})}[rr] & & \mathrm{Sym}^{n}(\cE_{1})
}
\end{equation}
with the vertical maps from (\ref{equ: auto pass to sym}) and (\ref{equ: Galois pass to sym}).
In particular, we have $r_{n}(\mathbf{E}_{n}^{\sharp})=\cE_{n}^{\sharp}$ for each $n\geq 1$.
\item \label{it: regulator 2} For each $m,n\geq 1$, we have the following commutative diagram
\begin{equation}\label{equ: regulator cup}
\xymatrix{
\mathbf{E}_{m}\otimes_E\mathbf{E}_{n} \ar^{r_{m}\otimes_E r_{n}}[rr] \ar[d] & & \cE_{m}\otimes_E\cE_{n} \ar[d]\\
\mathbf{E}_{m+n} \ar^{r_{m+n}}[rr] & & \cE_{m+n}
}
\end{equation}
with the vertical maps from (\ref{equ: auto stable cup}) and (\ref{equ: Galois tensor sub Ext embedding}).
\item \label{it: regulator 3} For each $n\geq 2$, we have $r_{n}(\mathbf{E}_{n}^{\flat})=\cE_{n}^{\flat}$.
\end{enumerate}
\end{thm}
\begin{proof}
Following Proposition~\ref{prop: sharp flat intersection} (upon replacing $n-1$ in \emph{loc.cit.} with $n$), we write
\[\mathbf{E}_{n}^{\plog}=E\plog_{\{1\}}\cup\cdots\cup\plog_{\{n\}}\subseteq\mathbf{E}_{n}\]
for each $n\geq 1$, using (\ref{equ: explicit cup generator}).
Assuming Hypothesis~\ref{hypo: log projection} for $n\geq 1$, we recall from \ref{it: sharp flat intersection 1} of Proposition~\ref{prop: sharp flat intersection} that
\[\mathbf{E}_{n}=\mathbf{E}_{n}^{\flat}\oplus\mathbf{E}_{n}^{\sharp}\]
and that
\begin{equation}\label{equ: auto diamond}
\mathbf{E}_{n}^{\diamond}\defeq \mathbf{E}_{n}^{\flat}\cap (\mathbf{E}_{n}^{\sharp}+\mathbf{E}_{n}^{\plog})
\end{equation}
is $1$ dimensional for each $n\geq 2$.
We write $\Phi^+_{\infty}$ for the set of pairs $(i,j)$ for each $j>i\geq 1$, identified with the union of the set of positive roots $\Phi_{n}^{+}$ (bounded by $j\leq n+1$) of $G_{n+1}$, with inclusions
\[\Phi^+_{1}\subseteq \Phi^+_{2}\subseteq \cdots\subseteq \Phi^+_{n}\subseteq \cdots\]
induced from the rational maps
\[G_{2}\dashrightarrow G_{3}\dashrightarrow \cdots\dashrightarrow G_{n}\dashrightarrow \cdots.\]
We write $\Delta_{\infty}\defeq \bigcup_{n=1}^{\infty}\Delta_{n}\subseteq \Phi^+_{\infty}$ for the subset of simple roots, namely those $(i,j)$ with $j=i+1$.
We write $\cS_{n}$ for the set of partitions of $(1,n+1)$ into elements of $\Phi^+_{\infty}$. So each $S\in\cS_{n}$ corresponds to the partition of $\Delta_{n}=[1,n]$ into subintervals $\{I_{\al}\mid \al\in S\}$. We can uniquely write $S$ as $\{\al_{1},\dots,\al_{\#S}\}$ so that $I_{\al_{s}}<I_{\al_{t}}$ for each $1\leq s<t\leq \#S$. We write $m_{s}\defeq \#I_{\al_{s}}$ for each $1\leq s\leq \#S$ so that $\sum_{s=1}^{\#S}m_s=n$.
For each $n\geq 1$, $S\in \cS_{n}$ (with associated ordered tuple $m_1,\dots,m_{\#S}$) and $I\subseteq S\cap\Delta_{n}$, we define $\mathbf{E}_{n,S,I}$ as the image of
\[\mathbf{E}_{m_1}^{\ast}\otimes_E\cdots \otimes_E\mathbf{E}_{m_{\#S}}^{\ast}\]
under the map (see (\ref{equ: multi cup map}))
\[\mathbf{E}_{m_1}\otimes_E\cdots \otimes_E\mathbf{E}_{m_{\#S}} \rightarrow \mathbf{E}_{n},\]
where $\mathbf{E}_{m_{s}}^{\ast}=\mathbf{E}_{m_{s}}^{\diamond}$ (resp.~$\mathbf{E}_{m_{s}}^{\ast}=\mathbf{E}_{1}^{\plog}$, resp.~$\mathbf{E}_{m_{s}}^{\ast}=\mathbf{E}_{1}^{\infty}$) if $m_s\geq 2$ (resp.~if $m_s=1$ and $\al_s\notin I$, resp.~if $m_s=1$ and $\al_s\in I$).
In other words, given $n\geq 1$, we have defined a $1$-dimensional $E$-subspace $\mathbf{E}_{n,S,I}\subseteq \mathbf{E}_{n}$ for each $S\in \cS_{n}$ and $I\subseteq S\cap\Delta_{n}$ which altogether satisfy
\[\mathbf{E}_{n}\cong \bigoplus_{S\subseteq \cS_{n}, I\subseteq S\cap\Delta_{n}}\mathbf{E}_{n,S,I}.\]
Similarly, we can define a $1$-dimensional $E$-subspace $\cE_{n,S,I}\subseteq \cE_{n}$ for each $S\in \cS_{n}$ and $I\subseteq S\cap\Delta_{n}$ which altogether satisfy
\[\cE_{n}\cong \bigoplus_{S\subseteq \cS_{n}, I\subseteq S\cap\Delta_{n}}\cE_{n,S,I}.\]

In the rest of the proof, we show that \ref{it: regulator 1}, \ref{it: regulator 2} and \ref{it: regulator 3} together force $r_{n}$ to have the form $\bigoplus_{S\in\cS_{n},I\subseteq S\cap\Delta_{n}}r_{n,S,I}$ for some uniquely determined isomorphisms
\begin{equation}\label{equ: regular match line}
r_{n,S,I}: \mathbf{E}_{n,S,I}\buildrel\sim\over\longrightarrow \cE_{n,S,I}
\end{equation}
between $1$ dimensional $E$-vector spaces.

Note that \ref{it: regulator 1} determines $r_{1}$ which restricts to isomorphisms $\mathbf{E}_{1}^{\infty}\buildrel\sim\over\longrightarrow \cE_{1}^{\infty}$ and $\mathbf{E}_{1}^{\plog}\buildrel\sim\over\longrightarrow \cE_{1}^{\plog}$.
This together with \ref{it: regulator 2} implies that $r_{n}$ must restrict to an (uniquely determined) isomorphism
\begin{equation}\label{equ: log regulator n}
\mathbf{E}_{n}^{\plog}\buildrel\sim\over\longrightarrow \cE_{n}^{\plog}
\end{equation}
for each $n\geq 1$.

Now that $r_{n}$ restricts to an isomorphism $\mathbf{E}_{n}^{\flat}\buildrel\sim\over\longrightarrow \cE_{n}^{\flat}$ and an isomorphism $\mathbf{E}_{n}^{\sharp}\buildrel\sim\over\longrightarrow \cE_{n}^{\sharp}$ for each $n\geq 2$ by \ref{it: regulator 3} and \ref{it: regulator 1}, we deduce from (\ref{equ: log regulator n}), (\ref{equ: auto diamond}) and Lemma~\ref{lem: Galois Ext decomposition} that $r_{n}$ restricts to an isomorphism
\begin{equation}\label{equ: diamond regulator n}
r_{n}^{\diamond}: \mathbf{E}_{n}^{\diamond}\buildrel\sim\over\longrightarrow \cE_{n}^{\diamond}
\end{equation}
between $1$ dimension $E$-vector spaces for each $n\geq 2$. The isomorphism $r_{n}^{\diamond}$ is uniquely determined by \ref{it: regulator 1} and the following commutative diagram
\[
\xymatrix{
\mathbf{E}_{n}^{\diamond} \ar^{r_{n}^{\diamond}}[rr] \ar[d] & & \cE_{n}^{\diamond} \ar[d]\\
\mathrm{Sym}^{n}(\mathbf{E}_{1}) \ar^{\mathrm{Sym}^{n}(r_{1})}[rr] & & \mathrm{Sym}^{n}(\cE_{1})
}
\]
with the LHS anf RHS vertical maps being embeddings by Proposition~\ref{prop: sharp flat intersection} and Proposition~\ref{prop: sharp flat decomposition} respectively.
For each $n\geq 1$, $S\in \cS_{n}$ and $I\subseteq S\cap\Delta_{n}$, since we have a uniquely determined map
\[r_{m_{s}}^{\ast}: \mathbf{E}_{m_{s}}^{\ast}\buildrel\sim\over\longrightarrow \cE_{m_{s}}^{\ast}\]
for each $1\leq s\leq \#S$ by the discussion above, we conclude from \ref{it: regulator 2} that $r_{n}$ must restrict to a uniquely determined map (\ref{equ: regular match line}). The proof is thus finished.
\end{proof}

Given $\{r_{n}:\mathbf{E}_{n}\buildrel\sim\over\longrightarrow\cE_{n}\}$ as described by Theorem~\ref{thm: regulator candidate}, it is clear from \ref{it: regulator 2} of \emph{loc.cit.} that $r_{n}$ restricts to an isomorphism \[\mathrm{Fil}_{W}^{k}(\mathbf{E}_{n})=\cong \bigoplus_{S\subseteq \cS_{n}, I\subseteq S\cap\Delta_{n},\#S>k}\mathbf{E}_{n,S,I}\buildrel\sim\over\longrightarrow \cong \bigoplus_{S\subseteq \cS_{n}, I\subseteq S\cap\Delta_{n},\#S>k}\cE_{n,S,I}\cong \mathrm{Fil}_{W}^{k}(\cE_{n})\] 
for each $n\geq 1$ and $0\leq k\leq n$.
Recall from \cite{Sou81} and \cite[\S~1.3]{HK11} that we have a canonical isomorphism between $1$ dimensional $E$-vector spaces
\begin{equation}\label{equ: HKS regulator}
r_{\rm{HKS}}: \mathrm{gr}_{W}^{0}(\mathbf{E}_{n})=P^{2n-1}(\fg_n)\buildrel\sim\over\longrightarrow \mathrm{Ext}_{g}^1(\cR(1_{K^\times}),\cR(\delta^{n}))^{\vee}=\mathrm{gr}_{W}^{0}(\cE_{n})
\end{equation}
for each $n\geq 2$.
We thus arrive at the following natural question.
\begin{ques}\label{ques: HKS regulator}
Given $\{r_{n}:\mathbf{E}_{n}\buildrel\sim\over\longrightarrow\cE_{n}\}$ as described by Theorem~\ref{thm: regulator candidate}. Do we have the following commutative diagram
\begin{equation}\label{equ: HKS regulator diagram}
\xymatrix{
\mathbf{E}_{n} \ar^{r_{n}}[rrr] \ar@{->>}[d] & & & \cE_{n} \ar@{->>}[d]\\
\mathrm{gr}_{W}^{0}(\mathbf{E}_{n}) \ar@{=}[d] & & &\mathrm{gr}_{W}^{0}(\cE_{n}) \ar@{=}[d]\\
P^{2n-1}(\fg_n) \ar^{r_{\rm{HKS}}}[rrr] & & & \mathrm{Ext}_{g}^1(\cR(1_{K^\times}),\cR(\delta^{n}))^{\vee}
}
\end{equation}
for each $n\geq 2$?
\end{ques}

\subsection{Relation to the Fontaine-Mazur $\mathscr{L}$-invariants}
Closely related to Theorem~\ref{thm: regulator candidate}, we are now ready to propose a candidate on how to match the Breuil-Schraen $\mathscr{L}$-invariants (see Definition~\ref{def: total inv}) with the Fontaine-Mazur $\mathscr{L}$-invariants.
For each $m\geq 1$, we write
\[x_{m}^{\infty}\defeq \val_{\{1\}}\cup\cdots\cup\val_{\{m\}}\in \mathbf{E}_{m}^{\infty}\subseteq \mathbf{E}_{m}\]
for short.
For each $m\geq 1$, the composition of the following maps (with convention $\mathbf{E}_{1}^{\sharp}=0$)
\[\mathbf{E}_{m}^{\diamond}=\mathbf{E}_{m}^{\flat}\cap(\mathbf{E}_{m}^{\sharp}+\mathbf{E}_{m}^{\plog})\subseteq \mathbf{E}_{m}\twoheadrightarrow \mathbf{E}_{m}/\mathbf{E}_{m}^{\sharp}\buildrel\sim\over\longrightarrow \mathrm{Sym}^{m}(\mathbf{E}_{1})\]
factors through 
\begin{equation}\label{equ: auto diamond log}
\mathbf{E}_{m}^{\diamond}\buildrel\sim\over\longrightarrow (\mathbf{E}_{m}^{\sharp}+\mathbf{E}_{m}^{\plog})/\mathbf{E}_{m}^{\sharp} \hookrightarrow \mathrm{Sym}^{m}(\mathbf{E}_{1}),
\end{equation}
and thus there exist a unique element $x_{m}\in \mathbf{E}_{m}^{\diamond}$ which share the same image in $\mathrm{Sym}^{m}(\mathbf{E}_{1})$ as 
\[x_{m,\plog}\defeq \plog_{1}\cup\cdots\cup\plog_{m}\in\mathbf{E}_{m}^{\plog}\subseteq \mathbf{E}_{m}.\]
We use the convention $x_{0}\defeq 1\in E=\mathbf{E}_{0}$.
For each $\ell\geq 0$, we consider the element
\[x_{\leq \ell}\defeq \sum_{k=0}^{\ell}x_{k}\in \bigoplus_{k=0}^{\ell}\mathbf{E}_{k}^{\diamond}\subseteq \mathbf{E}_{\leq \ell}\defeq \bigoplus_{k=0}^{\ell}\mathbf{E}_{k}.\]
For each $\al\in\Phi^{+}_{\infty}=\bigcup_{m\geq 1}\Phi^{+}_{m}$, we choose $x_{\al}\in\mathbf{E}_{I_{\al}}$ as the unique element corresponding to $x_{\#I_{\al}}$ under $\mathbf{E}_{I_{\al}}\cong \mathbf{E}_{\#I_{\al}}$ (see (\ref{equ: Tits interval reduction})). 
This way, for each $m\geq 0$, we could define $x_{S,I}\in\mathbf{E}_{m}$ for each $S\in\cS_{m}$ and $I\subseteq S\cap\Delta_{m}$, so that $\mathbf{E}_{m}$ admits a basis of the form
\[\{x_{S,I}\mid S\in\cS_{m}, I\subseteq S\cap\Delta_{m}\}.\]

We fix from now a choice of $n\geq 2$.
Given $0\leq \ell\leq n-1$, the following map (see (\ref{equ: auto stable cup}))
\[\mathbf{E}_{n-1-\ell}\otimes_E\mathbf{E}_{k}\hookrightarrow \mathbf{E}_{n-1-\ell+k}\]
for each $0\leq k\leq \ell$ induces the follow maps
\begin{equation}\label{equ: auto package cup}
\mathbf{E}_{n-1-\ell}\otimes_E\mathbf{E}_{\leq \ell}\hookrightarrow \bigoplus_{k=0}^{\ell}\mathbf{E}_{n-1-\ell+k}\hookrightarrow\mathbf{E}_{\leq n-1}.
\end{equation}
We use the shortened notation $\kappa_{\leq\ell}$ for the composition of (\ref{equ: auto package cup}).
Now we endow $\mathbf{E}_{\leq n-1}$ with an explicit structure of a filtered $(\varphi,N)$-module as below.
\begin{itemize}
\item For each $0\leq m\leq n-1$, $\varphi$ acts on $\mathbf{E}_{m}$ by $p^{-m}$.
\item We have $N(\mathbf{E}_{n-1})=0$. For each $0\leq m\leq n-2$, $N$ sends each $x\in\mathbf{E}_{m}$ to the image of $x\otimes_E\val\in \mathbf{E}_{m}\otimes_E\mathbf{E}_{1}$ in $\mathbf{E}_{m+1}$.
\item For each $0\leq m\leq n-1$, we have
\begin{equation}\label{equ: auto package filtration}
\mathrm{Fil}_{H}^{-m}(\mathbf{E}_{\leq n-1})\cong\sum_{\ell=m}^{n-1}\kappa_{\leq\ell}(\mathbf{E}_{n-1-\ell}\otimes_E Ex_{\leq \ell}).
\end{equation}
\end{itemize}
It is not difficult to check from Definition~\ref{def: universal phi N}) and the discussion right below it that $\mathbf{E}_{\leq n-1}\langle1-n\rangle[1-n]$ is universal special with level $n-1$ and the sum (\ref{equ: auto package filtration}) is a direct sum, which together with Proposition~\ref{prop: Galois dR Ext exact} implies that we have
\[\mathbf{E}_{\leq n-1}\cong D_{\rm{dR}}(\cD_{n-1})\]
as filtered $(\varphi,N)$-modules.

Now we consider a Breuil-Schraen $\mathscr{L}$-invariant $W\subseteq\mathbf{E}_{n-1}$ which induces $W_{I}\subseteq\mathbf{E}_{I}$ for each $I\subseteq\Delta=\Delta_{n-1}$ and in particular a Breuil-Schraen $\mathscr{L}$-invariant $W_{m}\defeq W_{[1,m]}\subseteq \mathbf{E}_{[1,m]}\cong\mathbf{E}_{m}$ for each $1\leq m\leq n-1$ (with convention $W_{0}\defeq 0$).
We write 
\[W_{\leq n-1}\defeq \bigoplus_{m=0}^{n-1}W_{m}\subseteq \bigoplus_{m=0}^{n-1}\mathbf{E}_{m}=\mathbf{E}_{\leq n-1}\]
for short and endow $\mathbf{E}_{\leq n-1}/W_{\leq n-1}$ with the filtration induced from that of $\mathbf{E}_{\leq n-1}$.
For each $\al\in\Phi^{+}=\Phi^{+}_{n-1}$, we deduce from Theorem~\ref{thm: moduli of inv} that there exists a unique $\mathscr{L}_{\al}=\mathscr{L}_{i,j}\in E$ (with $\al=(i,j)$ for some $1\leq i<j\leq n$) such that $y_{\al}\defeq x_{\al}-\mathscr{L}_{\al}x_{\al}^{\infty}\in W_{I_{\al}}$. 
Moreover, defining $y_{S,I}$ for each $S\in\cS_{m}$, $I\subseteq S\cap\Delta_{m}$ and $0\leq m\leq n$ similarly as in Theorem~\ref{thm: moduli of inv}, we know that $W_{m}$ admits a basis of the form
\[\{y_{S,I}\mid S\in\cS_{m}, \Delta_{m}\neq I\subseteq S\cap\Delta_{m}\}\]
for each $0\leq m\leq n-1$.
In particular, we see that
\[N^{m}(x_{0})=x_{m}^{\infty}\in\mathbf{E}_{m}^{\infty}\]
induces a basis of $\mathbf{E}_{m}/W_{m}$ for each $0\leq m\leq n-1$, and thus
\[\{N^{m}(x_0)\mid 0\leq m\leq n-1\}=\{x_{m}^{\infty}\mid 0\leq m\leq n-1\}\]
induces a basis of $\mathbf{E}_{\leq n-1}/W_{\leq n-1}$. Moreover, the image of
\[\kappa_{\leq \ell}(\mathbf{E}_{n-1-\ell}\otimes_E Ex_{\leq \ell})\]
under the surjection $\mathbf{E}_{\leq n-1}\twoheadrightarrow \mathbf{E}_{\leq n-1}/W_{\leq n-1}$ is spanned by the image of
\[x_{n-1-\ell}^{\infty}\otimes_E(x_{0}^{\infty}+\sum_{k=1}^{\ell}\mathscr{L}_{n-\ell,n-\ell+k}x_{k}^{\infty})\in\mathbf{E}_{n-1-\ell}\otimes_E\mathbf{E}_{\leq \ell}\]
under the composition of
\[\mathbf{E}_{n-1-\ell}\otimes_E\mathbf{E}_{\leq \ell}\buildrel \kappa_{\leq \ell}\over\longrightarrow \mathbf{E}_{\leq n-1}\twoheadrightarrow \mathbf{E}_{\leq n-1}/W_{\leq n-1},\]
which is nothing but
\[x_{n-1-\ell}^{\infty}+\sum_{k=1}^{\ell}\mathscr{L}_{n-\ell,n-\ell+k}x_{n-1-\ell+k}^{\infty}
=N^{n-1-\ell}(x_0)+\sum_{k=n-\ell}^{n-1}\mathscr{L}_{n-\ell,k}N^{k}(x_0).\]
We conclude that, for each $0\leq m\leq n-1$, the $E$-vector space
\[\mathrm{Fil}_{H}^{-m}(\mathbf{E}_{\leq n-1}/W_{\leq n-1})\defeq (\mathrm{Fil}_{H}^{-m}(\mathbf{E}_{\leq n-1})+W_{\leq n-1})/W_{\leq n-1}\]
admits a basis of the form
\[\{N^{n-1-\ell}(x_0)+\sum_{k=n-\ell}^{n-1}\mathscr{L}_{n-\ell,k}N^{k}(x_0)\mid m\leq \ell\leq n-1\}.\]
This is saying nothing but the fact that the Fontaine-Mazur $\mathscr{L}$-invariants of the filtered $(\varphi,N)$-module $\mathbf{E}_{\leq n-1}/W_{\leq n-1}$ are exactly given by $\{\mathscr{L}_{\al}\mid \al\in\Phi^{+}\}$.
Finally, we assume that $K=\Q_p$ and set
\[\rho_{W}\defeq ((\mathbf{E}_{\leq n-1}/W_{\leq n-1})\otimes_{\Q_p}\mathbf{B}_{\rm{st}})^{\varphi=1,N=0}.\]
Via the construction of an explicit filtered $(\varphi,N)$-module $\mathbf{E}_{\leq n-1}$, we obtain a map that associate with each Breuil-Schraen $\mathscr{L}$-invariant $W$ a semi-stable $p$-adic continuous Galois representation $\rho_{W}:\mathrm{Gal}_{\Q_p}\rightarrow\mathrm{GL}_{n}(E)$ which is characterized by the isomorphism $D_{\rm{st}}(\rho)\cong \mathbf{E}_{\leq n-1}/W_{\leq n-1}$ between filtered $(\varphi,N)$-modules.
We arrive at the following question.
\begin{ques}\label{ques: BS FM match}
Does the matching
\begin{equation}\label{equ: BS FM matching}
W \leftrightsquigarrow \rho_{W}
\end{equation}
satisfy certain form of $p$-adic local-global compatibility?
\end{ques}
Our construction of the explicit filtered $(\varphi,N)$-module $\mathbf{E}_{\leq n-1}$ and thus the matching (\ref{equ: BS FM matching}) is closely related to the following conjecture on the explicit description of Hyodo-Kato splitting (see \cite{CI10}, \cite{BdS16} and \cite{V24} for results in this direction).
Recall the de Rham complex $\Omega^{\bullet}=\Omega^{\bullet}(\bH)$ from (\ref{equ: Drinfeld dR}) whose differential maps are denoted by $d^{\bullet}$ for short. Recall from 
\cite{GK05} that the inclusions of $\varphi=p^{i}$ eigenspaces for $0\leq i\leq n-1$ induce an isomorphism (in the derived category of $E$-vector spaces)
\begin{equation}\label{equ: match HK splitting}
\bigoplus_{i=0}^{n-1}H^{i}_{\rm{HK}}(\bH)[-i]\buildrel\sim\over\longrightarrow R\Gamma_{\rm{HK}}(\bH)
\end{equation}
whose inverse together with the Hyodo-Kato isomorphism (which depends on the choice of $\plog$)
\begin{equation}\label{equ: match HK isom}
\Omega^{\bullet}=R\Gamma_{\rm{dR}}(\bH)\cong R\Gamma_{\rm{HK}}(\bH)
\end{equation}
gives a map
\begin{equation}\label{equ: match HK projection}
s_{\rm{HK}}^{i}: \Omega^{\bullet}=R\Gamma_{\rm{dR}}(\bH)\rightarrow H^{i}_{\rm{HK}}(\bH)[-i]
\end{equation}
in the derived category of $E$-vector spaces for each $0\leq i\leq n-1$. Note that (\ref{equ: match HK projection}) is actually a map in $\cM(G)$. For each $0\leq i\leq \ell\leq n-1$, we define $s_{\rm{HK}}^{i,\ell}$ as the composition of the following maps
\begin{equation}\label{equ: HK BdS splitting}
\mathrm{ker}(d^{\ell})[-\ell]\rightarrow \sigma_{\geq \ell}(\Omega^{\bullet})\rightarrow \Omega^{\bullet}\buildrel s_{\rm{HK}}^{i}\over\longrightarrow H^{i}_{\rm{HK}}(\bH)[-i].
\end{equation}
\begin{conj}\label{conj: HK BdS}
There exists a choice of an isomorphism between $D(G)$-modules
\begin{equation}\label{equ: HK rep isom}
H^{i}_{\rm{HK}}(\bH)\cong (V_{[1,n-1-i]}^{\infty})^{\vee}
\end{equation}
for each $0\leq i\leq n-1$ such that the following conditions hold.
\begin{enumerate}[label=(\roman*)]
\item \label{it: HK BdS 1} For each $1\leq i\leq n-1$, we have the following commutative diagram in $\cM(G)$
\begin{equation}\label{equ: HK BdS diagram}
\xymatrix{
H^{i}_{\rm{HK}}(\bH)[-i] \ar^{\sim}[r] \ar^{N}[d] & (V_{[1,n-1-i]}^{\infty})^{\vee}[-i] \ar^{\val}[d]\\
H^{i-1}_{\rm{HK}}(\bH)[1-i] \ar^{\sim}[r] & (V_{[1,n-i]}^{\infty})^{\vee}[1-i]
}
\end{equation}
with the horizontal isomorphisms from (\ref{equ: HK rep isom}), the LHS vertical map being the monodromy operator, and the RHS vertical map given by the element 
\[\val\in\mathrm{Ext}_{G}^{1}(V_{[1,n-i]}^{\infty},V_{[1,n-1-i]}^{\infty})=\Hom_{\cM(G)}((V_{[1,n-1-i]}^{\infty})^{\vee}[-i],(V_{[1,n-i]}^{\infty})^{\vee}[1-i]).\]
\item \label{it: HK BdS 2} For each $0\leq i\leq \ell\leq n-1$, we have the following commutative diagram in $\cM(G)$
\begin{equation}\label{equ: HK BdS Fil diagram}
\xymatrix{
\mathrm{ker}(d^{\ell})[-\ell] \ar^{s_{\rm{HK}}^{i,\ell}}[r] \ar[d] & H^{i}_{\rm{HK}}(\bH)[-i] \ar^{\wr}[d]\\
(V_{[1,n-1-\ell]}^{\rm{an}})^{\vee}[-\ell] \ar^{x_{i,\ell}}[r] & (V_{[1,n-1-i]}^{\infty})^{\vee}[-i]
}.
\end{equation}
with the top horizontal map being the composition of (\ref{equ: HK BdS splitting}), the right vertical isomorphism from (\ref{equ: HK rep isom}), the bottom horizontal map given by the element 
\[x_{i,\ell}\in \mathrm{Ext}_{G}^{i-\ell}(V_{[1,n-1-i]}^{\infty},V_{[1,n-1-\ell]}^{\rm{an}})=\Hom_{\cM(G)}((V_{[1,n-1-\ell]}^{\rm{an}})^{\vee}[-\ell],(V_{[1,n-1-i]}^{\infty})^{\vee}[-i]),\]
and the left vertical map being the unique map whose shift fits into the following commutative diagram
\begin{equation}\label{equ: HK BdS Fil scalar}
\xymatrix{
\mathrm{ker}(d^{\ell}) \ar[r] \ar[d] & H^{\ell}(\Omega^{\bullet}) \ar^{\wr}[d]\\
(V_{[1,n-1-\ell]}^{\rm{an}})^{\vee} \ar[r] & (V_{[1,n-1-\ell]}^{\infty})^{\vee}
}.
\end{equation}
whose RHS vertical isomorphism is given by the composition of (\ref{equ: HK rep isom}) with $H^{\ell}(\Omega^{\bullet})\cong H^{\ell}_{\rm{HK}}(\bH)$ from (\ref{equ: match HK isom}).
\end{enumerate}
\end{conj}

\newpage
\bibliographystyle{alpha}

\end{document}